\pdfoutput=1
\documentclass[12pt,oneside,openany]{book}
\usepackage[T1]{fontenc}
\usepackage{lmodern}

\usepackage{amsrefs}
\BibSpec{article}{%
+{}{\PrintAuthors}                {author}
+{,}{ \textit}                    {title}
+{.}{ }                           {part}
+{:}{ }                           {subtitle}
+{,}{ \PrintContributions}        {contribution}
+{.}{ \PrintPartials}             {partial}
+{,}{ }                           {journal}
+{}{ \textbf}                     {volume}
+{}{ \parenthesize}               {date}
+{,}{ }                           {pages}
+{,}{ }                           {status}
+{,}{ \PrintDOI}                  {doi}
+{.}{ Available at \eprint}       {eprint}
+{}{}                             {transition}
+{;}{ \PrintReviews}              {review}
+{.}{ }                           {note}
}

\usepackage{amsmath, amsfonts, amssymb, amsthm}
\usepackage{enumerate}
\usepackage{graphicx}
\usepackage{array}
\usepackage{wrapfig}
\usepackage{needspace}
\usepackage{mathrsfs}
\usepackage{xcolor}
\usepackage{longtable,array,needspace}
\usepackage[noautomatic]{imakeidx}
\makeindex[title=Subject index,intoc,columns=2,columnsep=24pt]
\newcommand{\SubjectIndexLayout}{%
\small
\let\SubjectIndexOriginalItem\item
\def\item{\needspace{3\baselineskip}\SubjectIndexOriginalItem}%
}
\indexsetup{othercode=\SubjectIndexLayout}

\makeatletter
\newcommand{\SubjectIndexAt}[2]{%
\expandafter\gdef\csname ham@index@#1\endcsname{#2}%
}
\newcount\SubjectIndexAnchorCount
\newcommand{\SubjectIndexLabel}[1]{%
\@ifundefined{ham@index@#1}{}{%
\@ifundefined{ham@index@done@#1}{%
\global\advance\SubjectIndexAnchorCount by 1\relax
\expandafter\gdef\csname ham@index@done@#1\endcsname{}%
\csname ham@index@#1\endcsname
}{}%
}%
}
\makeatother
\SubjectIndexAt{C:finite-gaussian-operators}{%
\index{Gaussian operators!finite variance estimates}%
}
\SubjectIndexAt{C:frequency-time}{%
\index{frequency estimates}%
}
\SubjectIndexAt{C:intrinsic-spatial-bracket}{%
\index{variational constant!spatial-energy bracket}%
}
\SubjectIndexAt{C:static-moment-constant-continuity}{%
\index{static limit!moment constants}%
}
\SubjectIndexAt{C:three-finite-images}{%
\index{covariance comparison!linear images and sums}%
\index{anisotropic covariances}%
}
\SubjectIndexAt{C:white-relative-width}{%
\index{spatial white noise!relative bound width}%
}
\SubjectIndexAt{Conj:intrinsic-even-moments}{%
\index{higher moments!even-order conjecture}%
}
\SubjectIndexAt{Conj:positive-control-compactness}{%
\index{optimizer!compactness conjecture}%
}
\SubjectIndexAt{Conj:white-first-correction}{%
\index{static limit!first-correction conjecture}%
}
\SubjectIndexAt{E:B-envelope}{%
\index{pairings!counting}%
}
\SubjectIndexAt{E:HAM}{%
\index{hyperbolic Anderson model}%
\index{initial data!constant displacement}%
\index{Anderson model|see{hyperbolic Anderson model}}%
}
\SubjectIndexAt{E:compact-common-trial}{%
\index{ground states!common trial function}%
}
\SubjectIndexAt{E:constant-phase-space}{%
\index{initial data!constant displacement}%
\index{wave phase!extended constant-data space}%
\index{constant initial data}%
\index{distributional phase space}%
}
\SubjectIndexAt{E:cov}{%
\index{covariance!temporal}%
\index{covariance!amplitude}%
\index{amplitude of noise|see{covariance, amplitude}}%
\index{fractional temporal noise|see{covariance, temporal}}%
}
\SubjectIndexAt{E:cp-input}{%
\index{numerical range}%
\index{Crouzeix--Palencia theorem}%
\index{functional calculus!bounded operators}%
}
\SubjectIndexAt{E:cutoff-peetre}{%
\index{Peetre's inequality}%
}
\SubjectIndexAt{E:damped-multiplier-comparison}{%
\index{Laplace transforms!damped wave norms}%
}
\SubjectIndexAt{E:damped-norm-definitions}{%
\index{Laplace transforms!damped wave norms}%
\index{Green function|see{Green measure}}%
\index{resolvent estimates}%
\index{wave kernel|see{Green measure}}%
\index{sine propagator}%
}
\SubjectIndexAt{E:degree-one-clock-loss}{%
\index{Laplace transforms!cross-pair discrepancy}%
}
\SubjectIndexAt{E:dyadic-time-majorant}{%
\index{dyadic intervals}%
\index{covariance!dyadic temporal majorant}%
}
\SubjectIndexAt{E:exact-fast-equation}{%
\index{small-noise normalization}%
}
\SubjectIndexAt{E:expansion}{%
\index{Stratonovich integrals!ordered expansion}%
\index{pairings!ordered chain kernels}%
\index{Green measure}%
}
\SubjectIndexAt{E:finite-gram-score}{%
\index{Gram matrices!nonnegative entries}%
}
\SubjectIndexAt{E:finite-hermitian-gaussian}{%
\index{Gaussian matrices!Hermitian dilation}%
}
\SubjectIndexAt{E:full-upper-bound}{%
\index{scaling exponent}%
}
\SubjectIndexAt{E:future-centered-log-wave}{%
\index{logarithmic wavefields}%
}
\SubjectIndexAt{E:future-log-wave}{%
\index{logarithmic wavefields}%
}
\SubjectIndexAt{E:future-tail-rate}{%
\index{large deviations!conditional upper-tail rate}%
}
\SubjectIndexAt{E:future-upper-tail}{%
\index{large deviations!conditional upper-tail rate}%
\index{exponential tilting!available moment orders}%
}
\SubjectIndexAt{E:gram-active-duality}{%
\index{active subsets}%
}
\SubjectIndexAt{E:gram-dual-coercivity}{%
\index{coercivity!finite dual problem}%
}
\SubjectIndexAt{E:intrinsic-positive-formula}{%
\index{Gaussian cost!variance-half normalization}%
\index{Gaussian action|see{Gaussian cost}}%
}
\SubjectIndexAt{E:kappa}{%
\index{scaling exponent}%
}
\SubjectIndexAt{E:masked-spectral-convention}{%
\index{covariance!spectral measure}%
\index{Fourier transform!normalization}%
}
\SubjectIndexAt{E:raw-block-square-identity}{%
\index{Gaussian cost!original covariance norm}%
}
\SubjectIndexAt{E:spatial-physical-pair}{%
\index{covariance!spatial measure}%
}
\SubjectIndexAt{E:temporal-equilibrium-mass}{%
\index{temporal equilibrium density!mass}%
}
\SubjectIndexAt{E:threshold}{%
\index{admissibility threshold}%
}
\SubjectIndexAt{E:tube-cameron-space}{%
\index{Cameron--Martin space}%
}
\SubjectIndexAt{E:variance-free-matrix}{%
\index{sine propagator!zero frequency}%
}
\SubjectIndexAt{E:variational-positive-class}{%
\index{positive controls!admissible class}%
\index{positive controls!nonseparable tests}%
\index{nonseparable controls|see{positive controls, nonseparable tests}}%
}
\SubjectIndexAt{Fig:active-subset-optimization}{%
\index{active subsets}%
}
\SubjectIndexAt{Fig:chapter-roadmap}{%
\index{chapter roadmap}%
}
\SubjectIndexAt{Fig:dyadic-cover}{%
\index{dyadic intervals}%
}
\SubjectIndexAt{Fig:gaussian-enlargement}{%
\index{Gaussian enlargement}%
}
\SubjectIndexAt{Fig:global-control-blocks}{%
\index{time blocks!correlated inputs}%
\index{time blocks!terminal interval}%
}
\SubjectIndexAt{Fig:hermite-cutoff}{%
\index{Hermite degree}%
}
\SubjectIndexAt{Fig:interval-end-metrics}{%
\index{phase norms!changes at interval endpoints}%
}
\SubjectIndexAt{Fig:nonseparable-controls}{%
\index{positive controls!nonseparable tests}%
}
\SubjectIndexAt{Fig:order-of-limits}{%
\index{regularization!order of limits}%
}
\SubjectIndexAt{Fig:positive-gram}{%
\index{Gram matrices!singular and zero cases}%
}
\SubjectIndexAt{Fig:walsh-guitar}{%
\index{Walsh's guitar example}%
}
\SubjectIndexAt{H:homog}{%
\index{covariance!homogeneous}%
\index{covariance!spatial measure}%
\index{covariance!spectral measure}%
\index{covariance!positive definiteness}%
\index{homogeneity!measure scaling}%
\index{homogeneous covariance|see{covariance, homogeneous}}%
}
\SubjectIndexAt{H:spatial-maximum}{%
\index{convolution bounds!translated spatial maximum}%
}
\SubjectIndexAt{H:three-endpoints}{%
\index{covariance comparison!three-dimensional hypothesis}%
}
\SubjectIndexAt{L:alpha-one-convolution}{%
\index{convolution bounds!degree-one case}%
}
\SubjectIndexAt{L:cameron-fractional-products}{%
\index{fractional Sobolev norms}%
}
\SubjectIndexAt{L:cameron-supported-tensor}{%
\index{covariance space!supported tensor product}%
\index{tensor products!covariance spaces}%
}
\SubjectIndexAt{L:compact-fractional-phase}{%
\index{fractional Sobolev norms}%
}
\SubjectIndexAt{L:compact-positive-order}{%
\index{trace inequalities}%
\index{compact operators!positive-order domination}%
}
\SubjectIndexAt{L:compact-raw-profile}{%
\index{Gaussian mean!original control and variance}%
}
\SubjectIndexAt{L:complex-phase}{%
\index{wave phase!complex norm bound}%
\index{numerical range}%
\index{resolvent estimates}%
}
\SubjectIndexAt{L:constant-phase-free}{%
\index{initial data!velocity}%
\index{constant initial data}%
\index{finite propagation}%
\index{Green measure}%
}
\SubjectIndexAt{L:constant-source-solution}{%
\index{mild solution!regularized identification}%
\index{physical wave!identification}%
}
\SubjectIndexAt{L:cutoff-factorization}{%
\index{compact operators!cutoff factorization}%
\index{Schur's test}%
\index{periodic extension!auxiliary cutoff factorization}%
}
\SubjectIndexAt{L:driver-covariance-coefficients}{%
\index{covariance coefficients}%
}
\SubjectIndexAt{L:driver-dominated-tails}{%
\index{operator tails!compact domination}%
}
\SubjectIndexAt{L:driver-holder-limit}{%
\index{Gaussian operator paths!Holder regularity@Hölder regularity}%
}
\SubjectIndexAt{L:driver-interval-limit}{%
\index{Gaussian operator paths!fixed-interval limit}%
}
\SubjectIndexAt{L:endpoint-closure}{%
\index{covariance comparison!linear images and sums}%
}
\SubjectIndexAt{L:equilibrium-gamma-beta}{%
\index{Gamma function!recurrence and duplication}%
\index{Beta function}%
}
\SubjectIndexAt{L:exact-canonical-scaling}{%
\index{homogeneity!moment scaling}%
\index{scaling!exact law of the canonical solution}%
}
\SubjectIndexAt{L:finite-word-cutoff}{%
\index{Wiener chaos!iteration order versus chaos order}%
\index{Fock space!degree cutoffs}%
\index{Hermite degree}%
}
\SubjectIndexAt{L:fixed-split}{%
\index{coefficient estimates!degree allocation}%
}
\SubjectIndexAt{L:fock-form}{%
\index{Gaussian creation and annihilation operators}%
\index{Fock space!degree cutoffs}%
}
\SubjectIndexAt{L:gaussian-scalar-powers}{%
\index{Young's inequality!scalar powers}%
}
\SubjectIndexAt{L:gaussian-trace-bounds}{%
\index{trace inequalities}%
}
\SubjectIndexAt{L:gram-dual-minimum}{%
\index{dual optimization!positive minimizer}%
}
\SubjectIndexAt{L:gram-realization}{%
\index{Gram matrices!realization}%
}
\SubjectIndexAt{L:higher-riesz}{%
\index{Riesz covariance!spatial maxima}%
}
\SubjectIndexAt{L:interval-complex}{%
\index{interval fields}%
\index{phase norms!changes at interval endpoints}%
}
\SubjectIndexAt{L:interval-contractions}{%
\index{Gaussian contractions}%
}
\SubjectIndexAt{L:inventory-product}{%
\index{pairings!counting}%
\index{coefficient estimates!degree allocation}%
}
\SubjectIndexAt{L:joint-fixed-order}{%
\index{regularization!joint source widths}%
}
\SubjectIndexAt{L:masked-chain-comparison}{%
\index{mixed orders}%
\index{masked solutions!mixed-order comparison}%
}
\SubjectIndexAt{L:masked-continuous-wave}{%
\index{sine propagator!zero frequency}%
\index{sine propagator}%
}
\SubjectIndexAt{L:masked-temporal-modulation}{%
\index{temporal modulation}%
}
\SubjectIndexAt{L:masked-time-completion}{%
\index{covariance!temporal}%
\index{covariance space!null tests and completion}%
}
\SubjectIndexAt{L:masked-trig-products}{%
\index{trigonometric covariance products}%
}
\SubjectIndexAt{L:occupation-rayleigh-selection}{%
\index{occupation profiles!near-Rayleigh selection}%
}
\SubjectIndexAt{L:pairing-product-envelope}{%
\index{pairings!finite product bound}%
}
\SubjectIndexAt{L:phase-bessel}{%
\index{Bessel inequality}%
}
\SubjectIndexAt{L:phase-contractions}{%
\index{Gaussian contractions}%
}
\SubjectIndexAt{L:physical-spectral-energy}{%
\index{spatial energy!physical and spectral forms}%
\index{measure-valued covariance!energy identification}%
}
\SubjectIndexAt{L:planar-endpoints}{%
\index{covariance!positive definiteness}%
\index{covariance comparison!planar endpoints}%
\index{anisotropic covariances}%
}
\SubjectIndexAt{L:plane-endpoints}{%
\index{covariance comparison!plane-supported endpoints}%
}
\SubjectIndexAt{L:point-fast-sources}{%
\index{regularization!order of limits}%
\index{source convergence!point moments}%
}
\SubjectIndexAt{L:positive-series-type}{%
\index{positive coefficient series!type}%
}
\SubjectIndexAt{L:positive-uniform-optimizer}{%
\index{optimizer!uniform norm bound and smoothing}%
}
\SubjectIndexAt{L:product-endpoints}{%
\index{covariance comparison!product endpoints}%
}
\SubjectIndexAt{L:raw-block-indexed-cone}{%
\index{spectral gap}%
\index{time blocks!positive cone}%
}
\SubjectIndexAt{L:raw-block-inputs}{%
\index{time blocks!correlated inputs}%
}
\SubjectIndexAt{L:raw-block-restart}{%
\index{initial data!velocity}%
\index{time blocks!composition and restart}%
}
\SubjectIndexAt{L:raw-block-uniform-error}{%
\index{time blocks!uniform error}%
}
\SubjectIndexAt{L:signed-limiting-multiplier}{%
\index{signed potentials!limiting multiplier}%
}
\SubjectIndexAt{L:signed-score-envelope}{%
\index{signed controls!upper variational expression}%
}
\SubjectIndexAt{L:signed-spectral-wave}{%
\index{signed potentials!spectral wave}%
\index{functional calculus!signed static wave}%
}
\SubjectIndexAt{L:source-finite-comparison}{%
\index{Gaussian comparison!finite linear images}%
}
\SubjectIndexAt{L:source-pointfield-regularity}{%
\index{source convergence!regularized random fields}%
}
\SubjectIndexAt{L:source-retaining-isometry}{%
\index{noise!preservation of the original process}%
\index{noise!auxiliary independent extension}%
\index{feature isometry}%
}
\SubjectIndexAt{L:source-spatial-heat}{%
\index{heat smoothing}%
}
\SubjectIndexAt{L:spatial-fourier-bound}{%
\index{spatial energy!deterministic Fourier bound}%
}
\SubjectIndexAt{L:spatial-physical-pair}{%
\index{spatial pairings!measure covariance}%
}
\SubjectIndexAt{L:spatial-wick-measures}{%
\index{Wick measures!spatial limit}%
}
\SubjectIndexAt{L:static-four}{%
\index{finite-order checks!static calibration}%
}
\SubjectIndexAt{L:symplectic-free}{%
\index{symplectic form}%
\index{free wave!norm lower bound}%
}
\SubjectIndexAt{L:temporal-equilibrium}{%
\index{temporal equilibrium density}%
\index{Riesz equilibrium density|see{temporal equilibrium density}}%
}
\SubjectIndexAt{L:terminal-wick}{%
\index{pairings!empty chains}%
\index{Wick formula!terminal chains}%
}
\SubjectIndexAt{L:translated-covariance-CS}{%
\index{Cauchy--Schwarz inequality!translated covariance}%
}
\SubjectIndexAt{L:tube-closed-geometry}{%
\index{compactness!images of Hilbert balls}%
}
\SubjectIndexAt{L:tube-half-action}{%
\index{Gaussian cost!original covariance norm}%
}
\SubjectIndexAt{L:tube-rayleigh-tail}{%
\index{Rayleigh integral!tail estimate}%
}
\SubjectIndexAt{L:uniform-mean-cutoff}{%
\index{physical masks}%
\index{spatial cutoff!uniform mean bounds}%
\index{cutoff|see{spatial cutoff}}%
}
\SubjectIndexAt{L:uniform-mean-spatial-tails}{%
\index{spatial cutoff!heat and frequency tails}%
}
\SubjectIndexAt{L:uniform-mean-tensor}{%
\index{tensor products!Hilbert-valued extension}%
}
\SubjectIndexAt{L:uniform-product-budgets}{%
\index{operator products!uniform growth bounds}%
}
\SubjectIndexAt{L:uniform-series-tail}{%
\index{series tails!uniform summability}%
}
\SubjectIndexAt{L:variational-control-dilation}{%
\index{homogeneity!measure scaling}%
\index{positive controls!dilation}%
}
\SubjectIndexAt{L:variational-heat-envelope}{%
\index{heat smoothing}%
}
\SubjectIndexAt{L:vertex-crude}{%
\index{regularization!uniform vertex bound}%
}
\SubjectIndexAt{L:white-cdf-square}{%
\index{Gagliardo--Nirenberg inequality!sharp one-dimensional constant}%
}
\SubjectIndexAt{L:white-clock-laplace}{%
\index{Laplace transforms!clock comparison}%
\index{Brownian clocks!low-order transforms}%
}
\SubjectIndexAt{L:young-upper-growth}{%
\index{Young integration!perturbation bounds}%
\index{Young integration!inverse evolution}%
}
\SubjectIndexAt{P:Jensen}{%
\index{first moment!Jensen comparison}%
\index{Jensen's inequality}%
}
\SubjectIndexAt{P:actual-operator-driver}{%
\index{Gaussian operator paths!construction}%
}
\SubjectIndexAt{P:actual-wave-form}{%
\index{wave phase!displacement and velocity}%
}
\SubjectIndexAt{P:cameron-actual-embedding}{%
\index{Cameron map!block embedding}%
}
\SubjectIndexAt{P:cameron-full-mean}{%
\index{Gaussian mean!completed space-time map}%
}
\SubjectIndexAt{P:cameron-spatial-map}{%
\index{Gaussian mean!compact spatial map}%
}
\SubjectIndexAt{P:cameron-temporal-map}{%
\index{Schur's test}%
\index{Gaussian mean!compact temporal map}%
\index{Young's inequality!convolution}%
}
\SubjectIndexAt{P:clock-copy-obstruction}{%
\index{Brownian clocks!copying obstruction}%
}
\SubjectIndexAt{P:compact-E-wave}{%
\index{spectral projections!complementary evolution}%
}
\SubjectIndexAt{P:compact-F-complement}{%
\index{spectral projections!complementary evolution}%
}
\SubjectIndexAt{P:compact-F-projection}{%
\index{spectral projections!growing mode}%
}
\SubjectIndexAt{P:compact-ground-family}{%
\index{ground states!compact potential family}%
\index{spectral gap}%
\index{ground state|see{ground states}}%
\index{ground states}%
}
\SubjectIndexAt{P:compact-original-common-trial}{%
\index{ground states!common trial function}%
}
\SubjectIndexAt{P:constant-gaussian-flow}{%
\index{Young integration!Gaussian flow moments}%
}
\SubjectIndexAt{P:constant-young-linear}{%
\index{Young integration!linear evolution}%
}
\SubjectIndexAt{P:cutoff-boundary-tests}{%
\index{compact operators!whole-space obstruction}%
}
\SubjectIndexAt{P:cutoff-singular-values}{%
\index{Schatten norms}%
\index{singular values!whole-space cutoffs}%
}
\SubjectIndexAt{P:d1-upper}{%
\index{direct upper bounds!one dimension}%
}
\SubjectIndexAt{P:envelope}{%
\index{pairings!uniform-estimate obstruction}%
\index{factorial bounds}%
}
\SubjectIndexAt{P:finite-gaussian-matrix}{%
\index{Gaussian matrices!moment bounds}%
\index{Schatten norms}%
\index{Khintchine inequality!Gaussian matrices}%
\index{noncommutative Khintchine inequality|see{Khintchine inequality, Gaussian matrices}}%
}
\SubjectIndexAt{P:frequency-row}{%
\index{frequency estimates}%
}
\SubjectIndexAt{P:gram-entrywise-comparison}{%
\index{Gram matrices!entrywise comparison}%
}
\SubjectIndexAt{P:gram-positive-optimizer}{%
\index{Gram matrices!singular and zero cases}%
\index{dual optimization!finite-subset identity}%
\index{positive cone!Gram maximizer}%
}
\SubjectIndexAt{P:higher-obstructions}{%
\index{convolution bounds!obstructions}%
}
\SubjectIndexAt{P:higher-transfer}{%
\index{direct upper bounds!conditional higher exponents}%
}
\SubjectIndexAt{P:interval-coefficients}{%
\index{interval fields}%
\index{coefficient estimates!interval products}%
\index{factorial bounds}%
}
\SubjectIndexAt{P:localized-finite-variance}{%
\index{Gaussian operators!finite variance estimates}%
}
\SubjectIndexAt{P:low-alpha-upper}{%
\index{direct upper bounds!spatial exponent at most one}%
}
\SubjectIndexAt{P:masked-auxiliary-norm}{%
\index{physical masks}%
}
\SubjectIndexAt{P:masked-canonical-transfer}{%
\index{source convergence!point moments}%
}
\SubjectIndexAt{P:masked-noise-completion}{%
\index{covariance space!null tests and completion}%
\index{noise-null quotient|see{covariance space, null tests and completion}}%
}
\SubjectIndexAt{P:masked-quotient-obstruction}{%
\index{spatial cutoff!covariance quotient obstruction}%
}
\SubjectIndexAt{P:masked-source-moments}{%
\index{masked solutions!tested moments}%
}
\SubjectIndexAt{P:masked-spectral-growth}{%
\index{spectral measure!translated mass bounds}%
}
\SubjectIndexAt{P:mild-construction}{%
\index{mild solution!construction}%
\index{Stratonovich interpretation}%
}
\SubjectIndexAt{P:mixed}{%
\index{Stratonovich integrals!self-pairings}%
\index{second moment!complete mixed expansion}%
\index{mixed orders}%
\index{pairings!self-pairings}%
\index{pairings!cross-pairings}%
\index{pairings!empty chains}%
\index{pairings!ordered chain kernels}%
\index{self-contractions|see{pairings, self-pairings}}%
\index{Wick contractions|see{pairings}}%
}
\SubjectIndexAt{P:moment-growth-logconvex}{%
\index{static noise!temporal interpolation}%
\index{log-convexity!moment growth type}%
}
\SubjectIndexAt{P:near-static-mean-gap}{%
\index{first moment!near-static comparison}%
\index{intermittency!near-static mean gap}%
}
\SubjectIndexAt{P:nogo-01}{%
\index{no-go results!the common-parameter majorant has the wrong scale}%
}
\SubjectIndexAt{P:nogo-02}{%
\index{no-go results!no uniform cubic factorial denominator per diagram}%
}
\SubjectIndexAt{P:nogo-03}{%
\index{no-go results!failure of the all-time comparison with unit prefactor}%
}
\SubjectIndexAt{P:nogo-04}{%
\index{no-go results!divergence of the first coefficient ratio}%
}
\SubjectIndexAt{P:nogo-05}{%
\index{no-go results!support obstruction to pointwise concatenation}%
}
\SubjectIndexAt{P:nogo-06}{%
\index{no-go results!the geometry-free majorant}%
}
\SubjectIndexAt{P:nogo-07}{%
\index{no-go results!large label span does not imply a cubic denominator}%
}
\SubjectIndexAt{P:nogo-08}{%
\index{no-go results!a spatial maximum need not be finite}%
}
\SubjectIndexAt{P:nogo-09}{%
\index{no-go results!radial lowering can destroy positive definiteness}%
}
\SubjectIndexAt{P:nogo-10}{%
\index{no-go results!no finite spatial-supremum pullout in the merge lens}%
}
\SubjectIndexAt{P:nogo-11}{%
\index{no-go results!a scalar wave moment is not log-convex in variance}%
}
\SubjectIndexAt{P:nogo-12}{%
\index{no-go results!failure of a heat--Mellin endpoint recipe}%
}
\SubjectIndexAt{P:nogo-13}{%
\index{no-go results!no uniform isotropic envelope for a plane covariance}%
}
\SubjectIndexAt{P:nogo-14}{%
\index{no-go results!two independent wave blocks do not multiply scalar moments}%
}
\SubjectIndexAt{P:nogo-15}{%
\index{no-go results!the first variance coefficient need not be nonnegative}%
}
\SubjectIndexAt{P:nogo-16}{%
\index{no-go results!a deterministic surrogate can lose the annealed factor}%
}
\SubjectIndexAt{P:nogo-17}{%
\index{no-go results!constant-correlation clock copying misses a positive factor}%
}
\SubjectIndexAt{P:nogo-18}{%
\index{no-go results!the shift-erased envelope is not below the wave transform}%
}
\SubjectIndexAt{P:nogo-19}{%
\index{no-go results!a positive semidefinite residual can lower the wave mean}%
}
\SubjectIndexAt{P:nogo-20}{%
\index{no-go results!a positive control need not have proportional covariance profiles}%
}
\SubjectIndexAt{P:nogo-21}{%
\index{no-go results!a same-noise Hilbert--Schmidt obstruction}%
}
\SubjectIndexAt{P:nogo-22}{%
\index{no-go results!a mask that separates noise-null representatives}%
}
\SubjectIndexAt{P:nogo-23}{%
\index{no-go results!a one-sided variance bound fails}%
}
\SubjectIndexAt{P:nogo-24}{%
\index{no-go results!translations defeat whole-space compactness}%
}
\SubjectIndexAt{P:nogo-25}{%
\index{no-go results!escape from every fixed smooth projection}%
}
\SubjectIndexAt{P:nogo-26}{%
\index{no-go results!a pure Young power and a maximum from an average}%
}
\SubjectIndexAt{P:nogo-27}{%
\index{no-go results!the temporal concentration threshold}%
}
\SubjectIndexAt{P:nogo-28}{%
\index{no-go results!an occupation ceiling need not improve a comparison value}%
}
\SubjectIndexAt{P:nogo-29}{%
\index{no-go results!empty sets, concavity and matrix order}%
}
\SubjectIndexAt{P:nogo-30}{%
\index{no-go results!escaping variational scores}%
}
\SubjectIndexAt{P:nogo-31}{%
\index{no-go results!growth does not give every-phase lower bounds or local cutoff consistency}%
}
\SubjectIndexAt{P:nogo-32}{%
\index{no-go results!the point feature is outside the closed noise image}%
}
\SubjectIndexAt{P:nogo-33}{%
\index{no-go results!three invalid proof steps}%
}
\SubjectIndexAt{P:occupation-physical-class}{%
\index{occupation profiles!admissible class}%
}
\SubjectIndexAt{P:original-point-transfer}{%
\index{averaging!squares versus square of an average}%
\index{point moments!finite-propagation transfer}%
}
\SubjectIndexAt{P:original-signed-point-upper}{%
\index{signed controls!upper variational expression}%
}
\SubjectIndexAt{P:planar-upper}{%
\index{covariance comparison!planar endpoints}%
}
\SubjectIndexAt{P:point-source-cone}{%
\index{finite propagation}%
\index{averaging!squares versus square of an average}%
\index{backward light cone|see{finite propagation}}%
}
\SubjectIndexAt{P:positive-score-cutoff}{%
\index{cutoff removal!optimized variational value}%
}
\SubjectIndexAt{P:positive-signed-comparison}{%
\index{signed controls!comparison with positive controls}%
}
\SubjectIndexAt{P:radius-actual-moments}{%
\index{Gaussian operator paths!radius dependence}%
\index{spatial cutoff!explicit radius dependence}%
}
\SubjectIndexAt{P:radius-compact-envelopes}{%
\index{spatial cutoff!explicit radius dependence}%
}
\SubjectIndexAt{P:radius-mean-anchor}{%
\index{spatial cutoff!growing-radius mean}%
}
\SubjectIndexAt{P:raw-block-gaussian-law}{%
\index{Cameron--Martin shift}%
\index{Gaussian change of measure}%
\index{Cameron--Martin formula}%
}
\SubjectIndexAt{P:raw-block-original-lower}{%
\index{Gaussian cost!variance-half normalization}%
\index{lower bound!fixed Gaussian control}%
}
\SubjectIndexAt{P:raw-occupation-energy}{%
\index{occupation profiles!original covariance energy}%
}
\SubjectIndexAt{P:riesz-small}{%
\index{finite-order checks!Riesz covariance}%
\index{Riesz covariance!finite coefficients}%
}
\SubjectIndexAt{P:signed-closed-form}{%
\index{signed potentials!closed quadratic form}%
\index{quadratic forms!closed-form representation}%
}
\SubjectIndexAt{P:signed-moment-bulk-tail}{%
\index{Gaussian tails!bulk-tail decomposition}%
}
\SubjectIndexAt{P:signed-operator-upper}{%
\index{order of limits!same endpoint and block length}%
}
\SubjectIndexAt{P:signed-physical-wave}{%
\index{physical wave!identification}%
}
\SubjectIndexAt{P:source-canonical-convergence}{%
\index{noise!preservation of the original process}%
\index{source convergence!regularized random fields}%
}
\SubjectIndexAt{P:source-gaussian-transfer}{%
\index{Gaussian comparison!source transfer}%
\index{source convergence!operator paths}%
}
\SubjectIndexAt{P:source-ordinary-response}{%
\index{source convergence!ordinary response}%
}
\SubjectIndexAt{P:spatial-energy-map}{%
\index{spatial energy!variational definition}%
}
\SubjectIndexAt{P:spatial-four-linear}{%
\index{spatial pairings!four-linear extension}%
}
\SubjectIndexAt{P:spatial-smooth-ratio}{%
\index{spatial ratio}%
}
\SubjectIndexAt{P:static-holder-upper}{%
\index{static noise!temporal interpolation}%
\index{Riesz covariance!static interpolation}%
}
\SubjectIndexAt{P:static-minorant-lower}{%
\index{static limit!temporal-minimum comparison}%
}
\SubjectIndexAt{P:static-phase-all-orders}{%
\index{static noise!all-order phase bounds}%
\index{Duhamel expansion}%
}
\SubjectIndexAt{P:static-phase-variance}{%
\index{wave phase!displacement and velocity}%
\index{static noise!first variance estimate}%
}
\SubjectIndexAt{P:symplectic-raw}{%
\index{symplectic conservation!limiting wave}%
}
\SubjectIndexAt{P:symplectic-source}{%
\index{symplectic conservation!regularized sources}%
}
\SubjectIndexAt{P:tube-rayleigh-transfer}{%
\index{Gaussian enlargement}%
\index{Rayleigh integral!Gaussian comparison}%
}
\SubjectIndexAt{P:tube-wiener-realization}{%
\index{isonormal Gaussian process}%
\index{Gaussian support!first-chaos realization}%
}
\SubjectIndexAt{P:uniform-actual-cameron}{%
\index{Cameron map!uniform radius bound}%
}
\SubjectIndexAt{P:uniform-arbitrary-residual}{%
\index{operator products!arbitrary residual paths}%
}
\SubjectIndexAt{P:uniform-mean-product}{%
\index{time blocks!terminal interval}%
\index{operator products!finite-mesh means}%
}
\SubjectIndexAt{P:uniform-mean-regulator}{%
\index{spatial cutoff!uniform mean bounds}%
}
\SubjectIndexAt{P:uniform-rate-averaging}{%
\index{averaging!completed mean fields}%
}
\SubjectIndexAt{P:variational-positive-scores}{%
\index{lower bound!nonnegative controls}%
\index{positive controls!coefficient homogeneity}%
}
\SubjectIndexAt{P:variational-remove-mask}{%
\index{Rayleigh quotient!mask removal}%
}
\SubjectIndexAt{P:variational-unmasked-rate}{%
\index{Rayleigh quotient!unmasked rate}%
}
\SubjectIndexAt{P:wave-square-columns}{%
\index{variance bounds!column and adjoint-column maps}%
}
\SubjectIndexAt{P:wave-two-variances}{%
\index{variance bounds!input and output directions}%
}
\SubjectIndexAt{P:white-four}{%
\index{pairings!self-pairings}%
\index{pairings!cross-pairings}%
\index{finite-order checks!spatial white noise}%
}
\SubjectIndexAt{P:white-spatial-energy}{%
\index{spatial energy!white-noise evaluation}%
\index{spatial white noise!energy maximizer}%
\index{hyperbolic secant profile}%
\index{Sz.-Nagy inequality|see{Gagliardo--Nirenberg inequality}}%
}
\SubjectIndexAt{P:young-upper-multiplicative}{%
\index{Young integration!perturbation bounds}%
}
\SubjectIndexAt{R:G0-verdict}{%
\index{pairings!uniform-estimate obstruction}%
}
\SubjectIndexAt{R:gaussian-one-variance-test}{%
\index{variance bounds!input and output directions}%
}
\SubjectIndexAt{R:riesz-static}{%
\index{finite-order checks!static calibration}%
}
\SubjectIndexAt{R:three-endpoint-scope}{%
\index{source normalization!first-moment comparison}%
}
\SubjectIndexAt{S:actual-operator-driver}{%
\index{isonormal Gaussian process}%
}
\SubjectIndexAt{S:book-no-go}{%
\index{no-go results}%
}
\SubjectIndexAt{S:future-real-orders}{%
\index{higher moments!real-order extension}%
}
\SubjectIndexAt{S:future-rough-noise-extension}{%
\index{rough noise!proposed spectral extensions}%
}
\SubjectIndexAt{S:future-universality}{%
\index{Rytov transformation}%
\index{KPZ fixed point!comparison and limitations}%
\index{universality!open scaling question}%
}
\SubjectIndexAt{S:gaussian-isoperimetry-input}{%
\index{Gaussian isoperimetric theorem}%
\index{Borell--Sudakov--Tsirelson theorem|see{Gaussian isoperimetric theorem}}%
}
\SubjectIndexAt{S:literature-context}{%
\index{Stratonovich interpretation!comparison with Skorohod}%
\index{Skorohod interpretation!comparison with Stratonovich}%
\index{Wiener chaos!iteration order versus chaos order}%
\index{chaos expansion|see{Wiener chaos}}%
\index{Feynman--Kac representations!comparison of stochastic products}%
\index{Stratonovich correction|see{Stratonovich interpretation}}%
\index{Stratonovich interpretation}%
}
\SubjectIndexAt{S:source-driver-convergence}{%
\index{regularized sources|see{source convergence}}%
}
\SubjectIndexAt{S:theory-outlook}{%
\index{variational methods!possible extensions}%
}
\SubjectIndexAt{T:first}{%
\index{first moment!source bounds}%
\index{source normalization!first-moment comparison}%
}
\SubjectIndexAt{T:full-upper}{%
\index{second moment!independent upper bound}%
\index{moment growth}%
\index{adjacent dyadic grids}%
\index{one-third trick|see{adjacent dyadic grids}}%
}
\SubjectIndexAt{T:intrinsic-occupation-ceiling}{%
\index{variational constant!occupation ceiling}%
}
\SubjectIndexAt{T:intrinsic-second-moment}{%
\index{hyperbolic Anderson model}%
\index{covariance!amplitude}%
\index{covariance!zero case}%
\index{second moment!intrinsic limiting constant}%
\index{variational constant!intrinsic second moment}%
\index{moment growth}%
}
\SubjectIndexAt{T:intro-second-moment}{%
\index{covariance!zero case}%
\index{second moment!intrinsic limiting constant}%
}
\SubjectIndexAt{T:intro-upper}{%
\index{second moment!independent upper bound}%
}
\SubjectIndexAt{T:white-moment-bounds}{%
\index{second moment!explicit spatial-white bounds}%
\index{spatial white noise!explicit moment bounds}%
\index{white noise|see{spatial white noise}}%
}
\SubjectIndexAt{Tab:future-moments}{%
\index{point moments!higher-order obstruction}%
\index{Sobolev embeddings!higher-moment ranges}%
}

\AddToHookWithArguments{label}[ham/subject-index]{\SubjectIndexLabel{#1}}
\AtEndDocument{%
\typeout{SUBJECT INDEX: \the\SubjectIndexAnchorCount\space curated anchors emitted.}%
}

\usepackage[colorlinks=true]{hyperref}
\AtBeginDocument{%
\renewcommand{\eprint}[1]{%
\href{https://arxiv.org/abs/#1}{\nolinkurl{#1}}%
}%
}

\usepackage{framed}
\usepackage{tikz}
\usepackage[breakable,skins]{tcolorbox}
\usetikzlibrary{arrows.meta,calc,positioning,decorations.pathreplacing}
\definecolor{hamblue}{RGB}{36,75,105}
\definecolor{hamteal}{RGB}{30,112,106}
\definecolor{hamgold}{RGB}{156,101,20}
\definecolor{hamvideo}{RGB}{178,34,34}
\newcommand{\VideoPlaylist}{\url{https://www.youtube.com/playlist?list=PLWiXjeUrQ62k}}
\newcommand{\VideoEp}[2]{\href{https://youtu.be/#2}{\textcolor{white}{\textbf{$\blacktriangleright$\,#1}}}}
\newcommand{\VideoRef}[2]{\href{https://youtu.be/#2}{\textcolor{hamvideo}{$\blacktriangleright$\,#1}}}
\tikzset{
ham/box/.style={draw=hamblue, fill=hamblue!4, rounded corners=2pt,
align=center, font=\small, inner sep=5pt, minimum height=9mm},
ham/arrow/.style={-{Latex[length=2mm]}, semithick, draw=hamblue},
ham/pair/.style={dashed, thick, draw=hamteal},
ham/dot/.style={circle, draw=hamblue, fill=white, inner sep=1.8pt},
ham/note/.style={font=\small, align=center},
ham/brace/.style={decorate, decoration={brace, amplitude=4pt}, semithick},
ham/video/.style={fill=hamvideo, rounded corners=1.5pt, inner sep=1.8pt, font=\scriptsize}
}

\newcommand{\aiCommits}{2{,}387}
\newcommand{\aiFirstDay}{October~27}
\newcommand{\aiLastDay}{November~16, 2026}
\newcommand{\aiDays}{21}
\newcommand{\aiReviewLastDay}{November~12}
\newcommand{\aiRuns}{627}
\newcommand{\aiFollowups}{43}
\newcommand{\aiRunsDatedByFolder}{114}
\newcommand{\aiReviewFolders}{731}
\newcommand{\aiFoldersWithRun}{434}
\newcommand{\aiReviewTableRows}{%
Claude (Opus 5.5) & 244 & 228 & 16 & 147 & 6 & 77 & 21.2 \\
DeepSeek V4 Pro & 152 & 147 & 5 & 56 & 0 & 35 & 11.7 \\
Kimi K3 & 107 & 91 & 16 & 81 & 0 & 4 & 8.0 \\
GPT (Codex) & 124 & 124 & 0 & 41 & 70 & 19 & 14.5 \\
\midrule
Total & 627 & 590 & 37 & 325 & 76 & 135 & \\
}
\newcommand{\aiHoursCoverage}{Claude (Opus 5.5) 244 of 244; DeepSeek V4 Pro 89 of 152; Kimi K3 22 of 107; GPT (Codex) 124 of 124}
\newcommand{\aiRunsKimiTimed}{22}

\newcommand{\aiRunsFail}{76}

\newcommand{\aiRunsKimi}{107}

\newcommand{\aiFailGPT}{70}
\newcommand{\aiCodexSessions}{88}
\newcommand{\aiCodexInB}{22.62}
\newcommand{\aiCodexCachedPct}{98.9}
\newcommand{\aiCodexOutM}{50.7}
\newcommand{\aiCodexReasonM}{29.7}
\newcommand{\aiClaudeCalls}{32{,}229}
\newcommand{\aiClaudeOutM}{39.5}
\newcommand{\aiClaudeCacheReadB}{11.28}
\newcommand{\aiDeepSeekSessions}{237}
\newcommand{\aiDeepSeekInM}{75.9}
\newcommand{\aiDeepSeekCacheM}{188.8}
\newcommand{\aiDeepSeekOutM}{2.08}
\newcommand{\aiDeepSeekCost}{72.56}
\newcommand{\aiUsageSnapshot}{2026-11-16 16:48Z}
\newcommand{\aiWorkerToDeepSeekRatio}{85}
\newcommand{\aiNogoTotal}{44}
\newcommand{\aiNogoMath}{33}
\newcommand{\aiNogoRoutes}{18}
\newcommand{\aiNogoCounter}{15}
\newcommand{\aiNogoSource}{4}
\newcommand{\aiNogoMethod}{6}
\newcommand{\aiNogoProcess}{1}

\newcommand{\aiNogoTotalFirst}{42}
\newcommand{\aiNogoInBookFirst}{8}
\newcommand{\aiNogoReviewedFirst}{14}
\newcommand{\aiNogoRecordedFirst}{11}
\newcommand{\aiNogoSnapshot}{November~16, 2026}
\newcommand{\aiLaterAttempts}{118}
\newcommand{\aiLaterInLedger}{25}
\newcommand{\aiLaterFirstDay}{November~12}
\newcommand{\aiLaterLastDay}{November~15}
\newcommand{\aiTimelineFigure}{%
\begin{tikzpicture}[x=1cm,y=1cm]
  \fill[aiClaude] (0.080,0.000) rectangle (0.491,0.317);
  \fill[aiClaude] (0.651,0.000) rectangle (1.063,0.872);
  \fill[aiDeepSeek] (0.651,0.872) rectangle (1.063,0.991);
  \fill[aiClaude] (1.223,0.000) rectangle (1.634,0.872);
  \fill[aiDeepSeek] (1.223,0.872) rectangle (1.634,1.070);
  \fill[aiClaude] (1.794,0.000) rectangle (2.206,1.585);
  \fill[aiClaude] (2.366,0.000) rectangle (2.777,1.981);
  \fill[aiClaude] (2.937,0.000) rectangle (3.349,1.268);
  \fill[aiDeepSeek] (2.937,1.268) rectangle (3.349,1.308);
  \fill[aiKimi] (2.937,1.308) rectangle (3.349,1.347);
  \fill[aiDeepSeek] (3.509,0.000) rectangle (3.920,2.100);
  \fill[aiKimi] (3.509,2.100) rectangle (3.920,2.140);
  \fill[aiDeepSeek] (4.080,0.000) rectangle (4.491,0.594);
  \fill[aiDeepSeek] (4.651,0.000) rectangle (5.063,0.317);
  \fill[aiDeepSeek] (5.223,0.000) rectangle (5.634,0.357);
  \fill[aiClaude] (5.794,0.000) rectangle (6.206,0.158);
  \fill[aiDeepSeek] (5.794,0.158) rectangle (6.206,0.277);
  \fill[aiKimi] (5.794,0.277) rectangle (6.206,0.396);
  \fill[aiClaude] (6.366,0.000) rectangle (6.777,2.258);
  \fill[aiDeepSeek] (6.366,2.258) rectangle (6.777,3.566);
  \fill[aiKimi] (6.366,3.566) rectangle (6.777,3.804);
  \fill[aiDeepSeek] (6.937,0.000) rectangle (7.349,0.357);
  \fill[aiKimi] (6.937,0.357) rectangle (7.349,3.328);
  \fill[aiClaude] (7.509,0.000) rectangle (7.920,0.357);
  \fill[aiDeepSeek] (7.509,0.357) rectangle (7.920,0.475);
  \fill[aiKimi] (7.509,0.475) rectangle (7.920,0.674);
  \fill[aiDeepSeek] (8.080,0.000) rectangle (8.491,0.040);
  \fill[aiKimi] (8.080,0.040) rectangle (8.491,0.119);
  \fill[aiGPT] (8.080,0.119) rectangle (8.491,0.634);
  \fill[aiGPT] (8.651,0.000) rectangle (9.063,1.109);
  \fill[aiDeepSeek] (9.223,0.000) rectangle (9.458,0.357);
  \fill[aiKimi] (9.223,0.357) rectangle (9.458,0.911);
  \fill[aiGPT] (9.223,0.911) rectangle (9.458,4.200);
  \draw[aiCommits,fill=aiCommits!15] (9.498,0.000) rectangle (9.634,0.277);
  \draw[aiCommits,fill=aiCommits!15] (9.794,0.000) rectangle (10.206,1.030);
  \draw[aiCommits,fill=aiCommits!15] (10.366,0.000) rectangle (10.777,3.130);
  \draw[aiCommits,fill=aiCommits!15] (10.937,0.000) rectangle (11.349,0.238);
  \draw[aiCommits,thick] (0.286,0.196) -- (0.857,0.980) -- (1.429,1.078) -- (2.000,1.722) -- (2.571,2.450) -- (3.143,1.610) -- (3.714,1.470) -- (4.286,0.644) -- (4.857,1.190) -- (5.429,0.574) -- (6.000,0.518) -- (6.571,2.240) -- (7.143,0.924) -- (7.714,0.924) -- (8.286,1.722) -- (8.857,1.792) -- (9.429,1.932) -- (10.000,2.394) -- (10.571,3.206) -- (11.143,1.652) -- (11.714,4.200);
  \fill[aiCommits] (0.286,0.196) circle (1.3pt);
  \node[font=\scriptsize,anchor=north] at (0.286,0) {Oct~27};
  \fill[aiCommits] (0.857,0.980) circle (1.3pt);
  \fill[aiCommits] (1.429,1.078) circle (1.3pt);
  \node[font=\scriptsize,anchor=north] at (1.429,0) {Oct~29};
  \fill[aiCommits] (2.000,1.722) circle (1.3pt);
  \fill[aiCommits] (2.571,2.450) circle (1.3pt);
  \node[font=\scriptsize,anchor=north] at (2.571,0) {Oct~31};
  \fill[aiCommits] (3.143,1.610) circle (1.3pt);
  \fill[aiCommits] (3.714,1.470) circle (1.3pt);
  \node[font=\scriptsize,anchor=north] at (3.714,0) {Nov~2};
  \fill[aiCommits] (4.286,0.644) circle (1.3pt);
  \fill[aiCommits] (4.857,1.190) circle (1.3pt);
  \node[font=\scriptsize,anchor=north] at (4.857,0) {Nov~4};
  \fill[aiCommits] (5.429,0.574) circle (1.3pt);
  \fill[aiCommits] (6.000,0.518) circle (1.3pt);
  \node[font=\scriptsize,anchor=north] at (6.000,0) {Nov~6};
  \fill[aiCommits] (6.571,2.240) circle (1.3pt);
  \fill[aiCommits] (7.143,0.924) circle (1.3pt);
  \node[font=\scriptsize,anchor=north] at (7.143,0) {Nov~8};
  \fill[aiCommits] (7.714,0.924) circle (1.3pt);
  \fill[aiCommits] (8.286,1.722) circle (1.3pt);
  \node[font=\scriptsize,anchor=north] at (8.286,0) {Nov~10};
  \fill[aiCommits] (8.857,1.792) circle (1.3pt);
  \fill[aiCommits] (9.429,1.932) circle (1.3pt);
  \node[font=\scriptsize,anchor=north] at (9.429,0) {Nov~12};
  \fill[aiCommits] (10.000,2.394) circle (1.3pt);
  \fill[aiCommits] (10.571,3.206) circle (1.3pt);
  \node[font=\scriptsize,anchor=north] at (10.571,0) {Nov~14};
  \fill[aiCommits] (11.143,1.652) circle (1.3pt);
  \fill[aiCommits] (11.714,4.200) circle (1.3pt);
  \node[font=\scriptsize,anchor=north] at (11.714,0) {Nov~16};
  \draw (0,0) -- (12.000,0);
  \draw (0,0.000) -- (-0.08,0.000) node[left,font=\scriptsize] {0};
  \draw (0,0.792) -- (-0.08,0.792) node[left,font=\scriptsize] {20};
  \draw (0,1.585) -- (-0.08,1.585) node[left,font=\scriptsize] {40};
  \draw (0,2.377) -- (-0.08,2.377) node[left,font=\scriptsize] {60};
  \draw (0,3.170) -- (-0.08,3.170) node[left,font=\scriptsize] {80};
  \draw (0,3.962) -- (-0.08,3.962) node[left,font=\scriptsize] {100};
  \draw (12.000,0.000) -- (12.080,0.000) node[right,font=\scriptsize,aiCommits] {0};
  \draw (12.000,0.700) -- (12.080,0.700) node[right,font=\scriptsize,aiCommits] {50};
  \draw (12.000,1.400) -- (12.080,1.400) node[right,font=\scriptsize,aiCommits] {100};
  \draw (12.000,2.100) -- (12.080,2.100) node[right,font=\scriptsize,aiCommits] {150};
  \draw (12.000,2.800) -- (12.080,2.800) node[right,font=\scriptsize,aiCommits] {200};
  \draw (12.000,3.500) -- (12.080,3.500) node[right,font=\scriptsize,aiCommits] {250};
  \draw (12.000,4.200) -- (12.080,4.200) node[right,font=\scriptsize,aiCommits] {300};
  \draw (0,0) -- (0,4.200);
  \draw (12.000,0) -- (12.000,4.200);
  \node[rotate=90,font=\scriptsize,anchor=south] at (-0.75,2.100) {review runs};
  \node[rotate=-90,font=\scriptsize,anchor=south,aiCommits] at (12.850,2.100) {commits};
  \draw[densely dotted,aiMilestone] (0.286,0) -- (0.286,4.480);
  \draw[densely dotted,aiMilestone] (0.857,0) -- (0.857,4.480);
  \draw[densely dotted,aiMilestone] (3.143,0) -- (3.143,4.480);
  \draw[densely dotted,aiMilestone] (5.429,0) -- (5.429,4.480);
  \draw[densely dotted,aiMilestone] (7.714,0) -- (7.714,4.480);
  \draw[densely dotted,aiMilestone] (8.286,0) -- (8.286,4.480);
  \draw[densely dotted,aiMilestone] (8.857,0) -- (8.857,4.480);
  \draw[densely dotted,aiMilestone] (11.143,0) -- (11.143,4.480);
  \node[circle,fill=aiMilestone,text=white,inner sep=0pt,minimum size=0.34cm,font=\tiny\bfseries] at (0.286,4.650) {1};
  \node[circle,fill=aiMilestone,text=white,inner sep=0pt,minimum size=0.34cm,font=\tiny\bfseries] at (0.286,5.070) {2};
  \node[circle,fill=aiMilestone,text=white,inner sep=0pt,minimum size=0.34cm,font=\tiny\bfseries] at (0.286,5.490) {3};
  \node[circle,fill=aiMilestone,text=white,inner sep=0pt,minimum size=0.34cm,font=\tiny\bfseries] at (0.857,4.650) {4};
  \node[circle,fill=aiMilestone,text=white,inner sep=0pt,minimum size=0.34cm,font=\tiny\bfseries] at (3.143,4.650) {5};
  \node[circle,fill=aiMilestone,text=white,inner sep=0pt,minimum size=0.34cm,font=\tiny\bfseries] at (3.143,5.070) {6};
  \node[circle,fill=aiMilestone,text=white,inner sep=0pt,minimum size=0.34cm,font=\tiny\bfseries] at (3.143,5.490) {7};
  \node[circle,fill=aiMilestone,text=white,inner sep=0pt,minimum size=0.34cm,font=\tiny\bfseries] at (5.429,4.650) {8};
  \node[circle,fill=aiMilestone,text=white,inner sep=0pt,minimum size=0.34cm,font=\tiny\bfseries] at (7.714,4.650) {9};
  \node[circle,fill=aiMilestone,text=white,inner sep=0pt,minimum size=0.34cm,font=\tiny\bfseries] at (8.286,4.650) {10};
  \node[circle,fill=aiMilestone,text=white,inner sep=0pt,minimum size=0.34cm,font=\tiny\bfseries] at (8.286,5.070) {11};
  \node[circle,fill=aiMilestone,text=white,inner sep=0pt,minimum size=0.34cm,font=\tiny\bfseries] at (8.286,5.490) {12};
  \node[circle,fill=aiMilestone,text=white,inner sep=0pt,minimum size=0.34cm,font=\tiny\bfseries] at (8.857,4.650) {13};
  \node[circle,fill=aiMilestone,text=white,inner sep=0pt,minimum size=0.34cm,font=\tiny\bfseries] at (8.857,5.070) {14};
  \node[circle,fill=aiMilestone,text=white,inner sep=0pt,minimum size=0.34cm,font=\tiny\bfseries] at (11.143,4.650) {15};
  \fill[aiClaude] (0.00,-0.95) rectangle (0.25,-0.75);
  \node[right,font=\scriptsize] at (0.25,-0.85) {Claude (Opus 5.5)};
  \fill[aiDeepSeek] (3.00,-0.95) rectangle (3.25,-0.75);
  \node[right,font=\scriptsize] at (3.25,-0.85) {DeepSeek V4 Pro};
  \fill[aiKimi] (5.90,-0.95) rectangle (6.15,-0.75);
  \node[right,font=\scriptsize] at (6.15,-0.85) {Kimi K3};
  \fill[aiGPT] (7.70,-0.95) rectangle (7.95,-0.75);
  \node[right,font=\scriptsize] at (7.95,-0.85) {GPT (Codex)};
  \draw[aiCommits,fill=aiCommits!15] (10.00,-0.95) rectangle (10.25,-0.75);
  \node[right,font=\scriptsize] at (10.25,-0.85) {recovered attempts};
\end{tikzpicture}}
\newcommand{\aiMilestoneKey}{%
\textbf{1}~Oct~26: Project set up: notes scaffold and source manifest \\
\textbf{2}~Oct~27: Gate G0: one bound for every pairing cannot reach the target scale \\
\textbf{3}~Oct~27: Gate G1: all fifteen order-four white-noise terms computed \\
\textbf{4}~Oct~28: Upper bound at the correct time scale established after review \\
\textbf{5}~Oct~31: Goal revised: derive the actual constant, not match a predicted value \\
\textbf{6}~Nov~1: Variational formula for the second-moment constant integrated \\
\textbf{7}~Nov~1: DeepSeek and Kimi become the reviewers \\
\textbf{8}~Nov~5: Reading edition of about twenty pages written \\
\textbf{9}~Nov~9: Lean checks reach the actual Gaussian objects \\
\textbf{10}~Nov~9: Fifth whole review: PASS WITH NOTES at MEDIUM, not accepted \\
\textbf{11}~Nov~9: Lean companion recorded \\
\textbf{12}~Nov~10: Technical text organized as a book \\
\textbf{13}~Nov~11: Appendix E added: the AI work in numbers \\
\textbf{14}~Nov~11: Lean correspondence gates closed through wave L4 \\
\textbf{15}~Nov~15: Final publication-facing PDF on master \\
}
\newcommand{\aiMilestones}{15}
\newcommand{\aiInfraBookStart}{October~26, 2026}
\newcommand{\aiInfraBibStart}{December 2021}
\newcommand{\aiInfraSourcesStart}{January 2026}
\newcommand{\aiInfraWikiStart}{April 2026}
\newcommand{\aiInfraSkillsWhen}{in the same month}
\newcommand{\aiRereadFiles}{33}
\newcommand{\aiRereadProofs}{249}
\newcommand{\aiRereadChanged}{208}
\newcommand{\aiRereadUnchanged}{41}
\newcommand{\aiRereadStatementsChanged}{3}
\newcommand{\aiRereadReaderEditionProofs}{4}
\newcommand{\aiRereadCalls}{59}
\newcommand{\aiRereadReports}{58}
\newcommand{\aiRereadNoReport}{1}
\newcommand{\aiRereadPass}{44}
\newcommand{\aiRereadPassNotes}{14}
\newcommand{\aiRereadPairsSent}{364}
\newcommand{\aiRereadNeverSent}{16}
\newcommand{\aiRereadCurrent}{38}
\newcommand{\aiRereadCurrentPass}{28}
\newcommand{\aiRereadCurrentPassNotes}{10}
\newcommand{\aiRereadEarlier}{21}
\newcommand{\aiRereadComplete}{249}
\newcommand{\aiRereadChecked}{226}
\newcommand{\aiRereadConditional}{9}
\newcommand{\aiRereadConditionalGroups}{8}
\newcommand{\aiRereadPartial}{14}
\newcommand{\aiRereadRestarted}{4}
\newcommand{\aiRereadRestartedHeld}{3}

\newcommand{\aiRereadRestartedCurrentProofs}{10}
\newcommand{\aiRereadTokens}{5{,}042{,}639}
\newcommand{\aiRereadTokensCached}{207{,}872}
\newcommand{\aiRereadCostCounter}{5.84}
\newcommand{\aiRereadDay}{November~16, 2026}
\newcommand{\aiRereadSource}{cdef8c83}
\newcommand{\aiRereadReaderRuns}{11}
\newcommand{\aiRereadStopPoints}{53}
\newcommand{\aiLeanReCalls}{3}
\newcommand{\aiLeanReModel}{GPT-6.1-Sol}
\newcommand{\aiLeanReDay}{November~16, 2026}
\newcommand{\aiLeanReFirstMajor}{7}
\newcommand{\aiLeanReFirstMinor}{2}
\newcommand{\aiLeanReSecondMajor}{2}
\newcommand{\aiLeanReThirdMinor}{1}
\newcommand{\aiLeanReMinutes}{25}
\newcommand{\aiLeanReInput}{11{,}775{,}176}
\newcommand{\aiLeanReCached}{11{,}261{,}952}
\newcommand{\aiLeanReOutput}{59{,}443}
\newcommand{\aiLeanReReasoning}{20{,}596}

\newcommand{\aiCloseoutTable}{%
\begin{tabular}{@{}lrrrrr@{}}
  \toprule
  Outcome & \shortstack[r]{DeepSeek\\\texttt{deepseek-{}v4-{}pro}} & \shortstack[r]{Kimi\\\texttt{k3}} & \shortstack[r]{Model not\\recorded} & \shortstack[r]{No model\\call} & All \\
  \midrule
  Accepted at the stated scope & 3 & 2 & 0 & 0 & 5 \\
  Not accepted & 9 & 0 & 1 & 0 & 10 \\
  Invalid report & 2 & 6 & 0 & 0 & 8 \\
  Recovery refused & 0 & 0 & 0 & 1 & 1 \\
  Recovered & 0 & 0 & 0 & 4 & 4 \\
  Preflight (no call) & 0 & 0 & 0 & 3 & 3 \\
  Continuation (no new vote) & 0 & 2 & 0 & 0 & 2 \\
  Timed out & 0 & 1 & 0 & 0 & 1 \\
  Running at snapshot & 0 & 0 & 1 & 0 & 1 \\
  \midrule
  Records & 14 & 11 & 2 & 8 & 35 \\
  \bottomrule
\end{tabular}}
\newcommand{\aiCloseoutAccountingKey}{%
\item[S1:] The usage object as the tool reported it; no wider total is implied.
\item[S2:] The final step of a panel run only, not the whole session.
\item[S3:] One response of the model; no session total is implied.
\item[S4:] The final step only; the usage of an earlier abandoned attempt is unknown.
}
\newcommand{\aiCloseoutHeaderNote}{ For the 5 records with scopes~S2 and~S4 the reported counts cover only the final step, not the whole session.}
\newcommand{\aiCloseoutTokenRows}{%
\texttt{deepseek-{}v4-{}pro} (DeepSeek), format A, scope~S1 & \texttt{cache.read} 0, \texttt{cache.write} 0, \texttt{input} 1{,}986{,}239, \texttt{output} 33{,}007, \texttt{reasoning} 121{,}948, \texttt{total} 2{,}141{,}194 & 5 of 14 \\
\texttt{deepseek-{}v4-{}pro} (DeepSeek), format A, scope~S3 & \texttt{cache.read} 0, \texttt{cache.write} 0, \texttt{input} 1{,}154{,}041, \texttt{output} 61{,}540, \texttt{reasoning} 193{,}545, \texttt{total} 1{,}409{,}126 & 5 of 14 \\
\texttt{deepseek-{}v4-{}pro} (DeepSeek), format A, scope~S4 (header/final-step only) & \texttt{cache.read} 15{,}488, \texttt{cache.write} 0, \texttt{input} 34, \texttt{output} 12{,}327, \texttt{reasoning} 60{,}815, \texttt{total} 88{,}664 & 1 of 14 \\
\texttt{deepseek-{}v4-{}pro} (DeepSeek), format D, scope~S2 (header/final-step only) & \texttt{input} 19{,}811, \texttt{output} 10{,}253, \texttt{reasoning} 34{,}463 & 3 of 14 \\
\texttt{k3} (Kimi), format B, scope~S1 & \texttt{inputCacheCreation} 0, \texttt{inputCacheRead} 5{,}405{,}065, \texttt{inputOther} 264{,}438, \texttt{output} 100{,}424 & 1 of 9 \\
\texttt{k3} (Kimi), format C, scope~S1 & \texttt{inputCacheCreation} 0, \texttt{inputCacheRead} 75{,}707{,}586, \texttt{inputOther} 2{,}503{,}780, \texttt{output} 493{,}885 & 8 of 9 \\
Model not recorded, format D, scope~S2 (header/final-step only) & \texttt{input} 26{,}040, \texttt{output} 2{,}515, \texttt{reasoning} 3{,}621 & 1 of 2 \\
Model not recorded, scope~S1 & no reported value & 1 of 2 \\
}
\newcommand{\aiCloseoutFormatKey}{%
\item[A:] \path{{cache:{read:number,write:number},input:number,output:number,reasoning:number,total:number}}
\item[B:] \path{list[{inputCacheCreation:number,inputCacheRead:number,inputOther:number,output:number};{inputCacheCreation:number,inputOther:number,output:number}]}
\item[C:] \path{list[{inputCacheCreation:number,inputCacheRead:number,inputOther:number,output:number}]}
\item[D:] \path{{input:number,output:number,reasoning:number}}
}
\newcommand{\aiCloseoutCost}{\texttt{deepseek-{}v4-{}pro} (DeepSeek), format A, scope~S1: 1.52 (currency not recorded; 4 records; known reported components only); \texttt{deepseek-{}v4-{}pro} (DeepSeek), format D, scope~S2: 0.12 (currency not recorded; 3 records; known reported components only); \texttt{deepseek-{}v4-{}pro} (DeepSeek), format A, scope~S3: 0.87 (currency not recorded; 4 records; known reported components only); Model not recorded, format D, scope~S2: 0.04 (currency not recorded; 1 record; known reported components only)}
\newcommand{\aiCloseoutRunningNote}{; historical running records retain their snapshot status (one record, whose provisional outcome is not accepted)}
\newcommand{\aiCloseoutHistoricalRecords}{29}
\newcommand{\aiCloseoutRecentRecords}{6}
\newcommand{\aiCloseoutRecentWords}{Six}
\newcommand{\aiCloseoutLocalWords}{Four}
\newcommand{\aiCloseoutLocalPassWords}{three}
\newcommand{\aiCloseoutLocalNotesWords}{one}
\newcommand{\aiCloseoutLocalAcceptedWords}{Two}
\newcommand{\aiCloseoutLocalPartialWords}{two}
\newcommand{\aiCloseoutRecords}{35}
\newcommand{\aiCloseoutUsageRecords}{24}
\newcommand{\aiCloseoutDenominator}{25}
\newcommand{\aiCloseoutSource}{\path{experiment/ledger-updates/2026-11-15-cold-closeout/run-ledger-v4.json}}
\newcommand{\aiCloseoutSourceCommit}{c2359a4f}
\newcommand{\aiCloseoutSourceBlob}{eab20562}
\newcommand{\aiCloseoutCorrection}{November~15, 2026}
\newcommand{\aiCloseoutMonth}{November 2026}
\newcommand{\aiCloseoutSnapshot}{November~15, 2026}
\newcommand{\aiEvidenceAttempts}{118}
\newcommand{\aiEvidenceBaseAttempts}{25}
\newcommand{\aiEvidenceExtraAttempts}{93}
\newcommand{\aiEvidenceFollowups}{4}
\newcommand{\aiEvidenceRecordedStarts}{6}
\newcommand{\aiEvidenceClientClocks}{3}
\newcommand{\aiEvidenceStartRows}{%
\texttt{C001} & \texttt{2026-{}11-{}13 20:47:11.227388} & Runner launch \\
\texttt{C006} & \texttt{2026-{}11-{}14 00:00:37.645321} & Runner launch \\
\texttt{C007} & \texttt{2026-{}11-{}14 00:22:13.811729} & Runner launch \\
\texttt{C009} & \texttt{2026-{}11-{}14 00:32:32.408329} & Runner launch \\
\texttt{C010} & \texttt{2026-{}11-{}14 00:52:24.214003} & Runner launch \\
\texttt{C011} & \texttt{2026-{}11-{}14 01:19:53.631078} & CLI receipt start \\
}
\newcommand{\aiEvidenceDecisionDate}{November~15, 2026}
\newcommand{\aiEvidenceInventoryPath}{\path{experiment/ledger-updates/2026-11-15-evidence-update/review-attempts.json}}
\newcommand{\aiEvidenceUpdatesPath}{\path{experiment/ledger-updates/2026-11-15-evidence-update/observations.json}}

\newcommand{\lcCommit}{c8d8e07}

\newcommand{\lcRows}{101}
\newcommand{\lcSummary}{With the manuscript's hypotheses, 48 labelled results are formalized in full and 16 with a narrower conclusion; \lcRef{Proposition}{P:tube-rayleigh-transfer}{10.2.3} in full with two inputs supplied; steps, parts or narrower forms of 34 further labelled results; the entries distinguish supplied hypotheses, where present, from analytic work left outside Lean; \lcRef{Theorem}{T:full-upper}{3.3.8} and \lcRef{Theorem}{T:intrinsic-second-moment}{11.2.6} are outside Lean.}
\newcommand{\lcUnlabelled}{2239}
\newcommand{\lcGroups}{five}
\newcommand{\lcTables}{%
\lcCategory{Formalized with the paper's own hypotheses; nothing supplied.}{48}
\begin{lcTable}
  \lcRef{Equation}{E:kappa}{1.2.4} & The growth exponent kappa
\end{lcTable}
\begin{lcTable}
  \lcRef{Lemma}{L:damped-homogeneous-norms}{2.3.1} & Homogeneous damped norms and temporal tails
\end{lcTable}
\begin{lcTable}
  \lcRef{Lemma}{L:phase-contractions}{3.1.1} & Two static contraction bounds
\end{lcTable}
\begin{lcTable}
  \lcRef{Lemma}{L:finite-word-cutoff}{3.2.3} & Exact cutoffs for finite Gaussian words
\end{lcTable}
\begin{lcTable}
  \lcRef{Lemma}{L:inventory-product}{3.3.3} & A coefficient bound for a product with nonnegative coefficients
\end{lcTable}
\begin{lcTable}
  \lcRef{Lemma}{L:terminal-wick}{3.3.6} & The finite terminal-chain Wick formula
\end{lcTable}
\begin{lcTable}
  \lcRef{Lemma}{L:uniform-series-tail}{3.3.7} & Uniform series tails from scaled moments
\end{lcTable}
\begin{lcTable}
  \lcRef{Proposition}{P:masked-spectral-growth}{4.1.1} & Translated spectral mass
\end{lcTable}
\begin{lcTable}
  \lcRef{Proposition}{P:masked-auxiliary-norm}{4.1.2} & Auxiliary mask norm
\end{lcTable}
\begin{lcTable}
  \lcRef{Lemma}{L:masked-time-completion}{4.1.3} & Temporal covariance completion
\end{lcTable}
\begin{lcTable}
  \lcRef{Lemma}{L:masked-temporal-modulation}{4.1.6} & Modulated amplitudes
\end{lcTable}
\begin{lcTable}
  \lcRef{Equation}{E:masked-time-schur}{4.1.7} & Schur bound for the temporal form
\end{lcTable}
\begin{lcTable}
  \lcRef{Lemma}{L:masked-trig-products}{4.1.7} & Correlated trigonometric products
\end{lcTable}
\begin{lcTable}
  \lcRef{Lemma}{L:masked-continuous-wave}{4.1.8} & Continuous sine factors
\end{lcTable}
\begin{lcTable}
  \lcRef{Proposition}{P:frequency-row}{4.2.1} & Uniform scalar row
\end{lcTable}
\begin{lcTable}
  \lcRef{Corollary}{C:frequency-time}{4.2.2} & A time exponent above one half
\end{lcTable}
\begin{lcTable}
  \lcRef{Lemma}{L:phase-bessel}{4.2.4} & Hilbert-valued finite Bessel inequality
\end{lcTable}
\begin{lcTable}
  \lcRef{Lemma}{L:gaussian-scalar-powers}{4.3.1} & Scalar power estimates
\end{lcTable}
\begin{lcTable}
  \lcRef{Lemma}{L:gaussian-trace-bounds}{4.3.2} & Two Hermitian trace bounds
\end{lcTable}
\begin{lcTable}
  \lcRef{Proposition}{P:finite-gaussian-matrix}{4.3.3} & Finite Gaussian matrix moment
\end{lcTable}
\begin{lcTable}
  \lcRef{Remark}{R:gaussian-one-variance-test}{4.3.4} & Both variance directions are necessary
\end{lcTable}
\begin{lcTable}
  \lcRef{Equation}{E:driver-parameters}{5.1.1} & driver parameters
\end{lcTable}
\begin{lcTable}
  \lcRef{Lemma}{L:source-finite-comparison}{5.2.1} & Finite Gaussian linear images
\end{lcTable}
\begin{lcTable}
  \lcRef{Lemma}{L:source-spatial-heat}{5.2.3} & Original spatial space and heat
\end{lcTable}
\begin{lcTable}
  \lcRef{Lemma}{L:variational-heat-envelope}{7.3.1} & Gaussian envelope for the original covariance
\end{lcTable}
\begin{lcTable}
  \lcRef{Lemma}{L:variational-global-mean}{7.3.2} & Global regularity of the original mean
\end{lcTable}
\begin{lcTable}
  \lcRef{Lemma}{L:equilibrium-gamma-beta}{8.1.1} & Gamma and Beta identities
\end{lcTable}
\begin{lcTable}
  \lcRef{Lemma}{L:temporal-equilibrium}{8.1.2} & Exact potential and mass
\end{lcTable}
\begin{lcTable}
  \lcRef{Lemma}{L:spatial-physical-pair}{8.2.1} & The physical two-point integral
\end{lcTable}
\begin{lcTable}
  \lcRef{Lemma}{L:spatial-fourier-bound}{8.2.3} & Deterministic Fourier bound for the spatial energy
\end{lcTable}
\begin{lcTable}
  \lcRef{Proposition}{P:spatial-four-linear}{8.2.5} & The physical four-linear extension
\end{lcTable}
\begin{lcTable}
  \lcRef{Lemma}{L:occupation-kernel-row}{8.2.7} & Exact row maximum
\end{lcTable}
\begin{lcTable}
  \lcRef{Proposition}{P:occupation-physical-class}{8.2.8} & The full physical occupation class
\end{lcTable}
\begin{lcTable}
  \lcRef{Lemma}{L:tube-half-action}{10.2.8} & Supremum with the Gaussian quadratic cost
\end{lcTable}
\begin{lcTable}
  \lcRef{Lemma}{L:gram-scalar-algebra}{11.1.1} & Scalar and Hilbert identities
\end{lcTable}
\begin{lcTable}
  \lcRef{Proposition}{P:raw-occupation-energy}{11.3.2} & Original occupation energy and dilation
\end{lcTable}
\begin{lcTable}
  \lcRef{Corollary}{C:white-relative-width}{11.4.4} & Uniform relative width of the bound formulas
\end{lcTable}
\begin{lcTable}
  \lcRef{Proposition}{P:white-four}{A.1.1} & White-space order four
\end{lcTable}
\begin{lcTable}
  \lcRef{Lemma}{L:d1-spatial}{B.1.1} & Spatial estimate
\end{lcTable}
\begin{lcTable}
  \lcRef{Lemma}{L:low-alpha-convolution}{B.1.3} & Low-alpha convolution bound
\end{lcTable}
\begin{lcTable}
  \lcRef{Lemma}{L:alpha-one-convolution}{B.1.4} & Degree-one convolution bound
\end{lcTable}
\begin{lcTable}
  \lcRef{Lemma}{L:higher-riesz}{B.2.3} & Riesz spatial maxima
\end{lcTable}
\begin{lcTable}
  \lcRef{Lemma}{L:product-endpoints}{C.2.1} & Strict product endpoints
\end{lcTable}
\begin{lcTable}
  \lcRef{Lemma}{L:plane-endpoints}{C.2.3} & Plane-supported endpoints
\end{lcTable}
\begin{lcTable}
  \lcRef{Lemma}{L:endpoint-closure}{C.2.4} & Linear images and finite positive sums
\end{lcTable}
\begin{lcTable}
  \lcRef{Lemma}{L:positive-series-type}{D.1.2} & Type of a positive coefficient series
\end{lcTable}
\begin{lcTable}
  \lcRef{Lemma}{L:white-clock-laplace}{D.2.1} & White-space low-order Laplace transforms
\end{lcTable}
\begin{lcTable}
  \lcRef{Proposition}{P:clock-copy-obstruction}{D.2.2} & Obstruction to unweighted constant-correlation clock copying
\end{lcTable}

\lcCategory{With the paper's hypotheses; Lean proves a narrower conclusion.}{16}
\lcEntry{\lcRef{Lemma}{L:fock-form}{3.2.2}}{Gaussian cutoff form bound}
\lcPart{Lean proves}{With the lemma's hypotheses and nothing supplied, \lcRef{Equation}{E:fock-sandwich}{3.2.5}, in coefficient form and as an operator-norm bound with $A_0\ge0$ and $B$ bounded self-adjoint on the cutoff Hilbert space, uniformly in $\epsilon$, $m$ and the band, for finitely many real modes, Hermite degree at most $N$ in the $m$ multiplied coordinates, on a Fourier-side model.}
\lcPart{Not yet formalized}{The Fourier-side model is not identified with the physical-space operators.}
\lcEntryEnd
\lcEntry{\lcRef{Lemma}{L:interval-contractions}{3.3.1}}{Dimension-independent Gaussian contractions}
\lcPart{Lean proves}{With the lemma's hypotheses and nothing supplied, \lcRef{Equation}{E:interval-column}{3.3.1}, its band-compressed form, the phase-insertion bound $\sqrt{R_\alpha}\,\rho^{(\alpha-2)/2}$, the sandwich bound $R_\alpha\rho^{\alpha-4}$ with self-adjointness and \lcRef{Equation}{E:interval-fock}{3.3.2}, for finitely many real modes, on a Fourier-side model, without tensoring by another Hilbert factor.}
\lcPart{Not yet formalized}{The Fourier-side model is not identified with the physical-space operators.}
\lcEntryEnd
\lcEntry{\lcRef{Proposition}{P:masked-quotient-obstruction}{4.1.5}}{A physical mask need not preserve noise-null tests}
\lcPart{Lean proves}{The proposition, with its hypotheses, for the covariance functional: $\E W(F)^2$ is read as the covariance functional of the test $F$.}
\lcPart{Not yet formalized}{The noise $W$ is not constructed, and the identification of $\E W(F)^2$ with the variance of the noise is not formalized.}
\lcEntryEnd
\lcEntry{\lcRef{Proposition}{P:variational-unmasked-rate}{7.3.3}}{The unmasked Rayleigh rate}
\lcPart{Lean proves}{With the section's hypotheses and nothing supplied, the statements about the Rayleigh supremum \lcRef{Equation}{E:variational-real-rayleigh}{7.3.4}: for $\lambda$ defined as its square root, a uniform positive floor, the Lipschitz bounds and $I_\Psi$ positive and finite.}
\lcPart{Not yet formalized}{The operator statements ($H_s$ self-adjoint, form domain $H^1$, spectral bottom $E_s$, $\lambda=\sqrt{-E_s}$) are not formalized.}
\lcEntryEnd
\lcEntry{\lcRef{Lemma}{L:variational-control-dilation}{7.3.5}}{Radon and Rayleigh dilation}
\lcPart{Lean proves}{The bijection and the four identities with the section's hypotheses (\lcRef{Hypothesis}{H:homog}{1.2.1}, $\Gamma\neq0$, $a=\alpha_0$) and nothing supplied, for the variance in its integral form and the rate defined as the square root of the real Rayleigh supremum \lcRef{Equation}{E:variational-real-rayleigh}{7.3.4}.}
\lcPart{Not yet formalized}{The identification of the variance and the rate with $\E W(\Psi)^2$ (\lcRef{Lemma}{L:compact-raw-profile}{6.3.6}) and $\sqrt{-E_s}$ (\lcRef{Proposition}{P:variational-unmasked-rate}{7.3.3}) is not formalized.}
\lcEntryEnd
\lcEntry{\lcRef{Proposition}{P:variational-positive-scores}{7.3.6}}{Intrinsic variational lower bounds}
\lcPart{Lean proves}{The proposition with the section's hypotheses and nothing supplied, including the boundedness of the value sets and the supremum statements ($S_c^{\mathrm{pos}}$ finite, positive, equal to $c^{(4-\alpha)/(3-\alpha)}S_1^{\mathrm{pos}}$), with $\lambda$ the square root of the real Rayleigh supremum and $v_\Psi$ the cost in its integral form.}
\lcPart{Not yet formalized}{$Q_2\ge S_2^{\mathrm{pos}}$ is not formalized.}
\lcEntryEnd
\lcEntry{\lcRef{Lemma}{L:spatial-static-controls}{8.2.2}}{Smooth temporal controls}
\lcPart{Lean proves}{$\lambda$ is the square root of the supremum of the real Rayleigh quotient, and $I_\Psi$, $v_\Psi$ are given by their integral formulas; \lcRef{Equation}{E:spatial-static-score}{8.2.2} is proved both in the form that every number below the bound is exceeded by some admissible score and in its supremum form, with nothing supplied.}
\lcPart{Not yet formalized}{The identifications of $\lambda$ with $\sqrt{-E_s}$ and of $v_\Psi$ with $\E W(\Psi)^2$ are not formalized.}
\lcEntryEnd
\lcEntry{\lcRef{Proposition}{P:spatial-smooth-ratio}{8.2.4}}{The finite smooth spatial ratio}
\lcPart{Lean proves}{\lcRef{Equation}{E:spatial-smooth-score}{8.2.4}, with nothing supplied, both in the form that every $y$ below the bound is exceeded by some $I_\Psi-v_\Psi/2$ and as a bound on the supremum, with $\lambda$ the square root of the real Rayleigh supremum and $v_\Psi$ given by its integral formula. $R_0<\infty$ is proved by \lcRef{Lemma}{L:spatial-fourier-bound}{8.2.3}.}
\lcPart{Not yet formalized}{The identification of $\lambda$ and $v_\Psi$ with $\sqrt{-E_s}$ and $\E W(\Psi)^2$ is not formalized.}
\lcEntryEnd
\lcEntry{\lcRef{Proposition}{P:spatial-energy-map}{8.2.6}}{The spatial energy and static lower variational expression}
\lcPart{Lean proves}{\lcRef{Equation}{E:spatial-static-energy-score}{8.2.10}, with nothing supplied, both in an existential form (every $y$ below $M_a^{-1/k}F_0(\mathcal{E}_{\mathrm{sp}})$ is exceeded by some $I_\Psi-v_\Psi/2$) and in its supremum form, with $\lambda$ the square root of the Rayleigh supremum \lcRef{Equation}{E:variational-real-rayleigh}{7.3.4} and $v_\Psi$ by its integral formula.}
\lcPart{Not yet formalized}{The identification of $\lambda$ and $v_\Psi$ with $\sqrt{-E_s}$ and $\E W(\Psi)^2$ is not formalized.}
\lcEntryEnd
\lcEntry{\lcRef{Lemma}{L:tube-rayleigh-tail}{10.2.7}}{Exact Rayleigh tail bound}
\lcPart{Lean proves}{The bound with integrability and the finite head integral.}
\lcPart{Not yet formalized}{The ``In particular'' application to the coarse-envelope bound of another proposition is not formalized.}
\lcEntryEnd
\lcEntry{\lcRef{Lemma}{L:positive-raw-embedding}{11.2.1}}{Inclusion into the covariance space and the smooth mean}
\lcPart{Lean proves}{\lcRef{Equation}{E:positive-bochner-inclusion}{11.2.1} for every real Hilbert $G$ ($H_t\otimes G$ the kernel completion of $L^2(I;G)$, M5), the inclusion of $\mathcal{C}_{\mathrm{pos}}$ in $H_I$ with $\|\Psi\|^2=v_\Psi$, the variance identity in integral form and the row and Schur bounds.}
\lcPart{Not yet formalized}{The mean identity of \lcRef{Equation}{E:positive-raw-mean}{11.2.2} and the canonical versions are not formalized.}
\lcEntryEnd
\lcEntry{\lcRef{Lemma}{L:positive-physical-gram}{11.2.2}}{Finite covariance vectors and entrywise order}
\lcPart{Lean proves}{With the section's hypotheses and nothing supplied, that $k_i^R$ and $k_i^\infty$ lie in $H_I$ (M5), that both Grams factor into the temporal pairing and $B_\Gamma$, and that both are positive semidefinite.}
\lcPart{Not yet formalized}{The label row \lcRef{Equation}{E:positive-label-row}{11.2.4} ($\mathcal{T}_R$, $v_L$ not constructed) and the entrywise order are not formalized.}
\lcEntryEnd
\lcEntry{\lcRef{Lemma}{L:occupation-raw-square}{11.3.1}}{Full-profile squares in the original space}
\lcPart{Lean proves}{With the section's hypotheses and nothing supplied, $Q:H^1\to H_{\mathrm{sp}}$ with $\langle Q(f),Q(h)\rangle=B_\Gamma(f^2,h^2)$, continuous on $H^1$ (M6), and for fields of $\mathcal{O}$ modelled by joint representatives $\rho_g\in H_I$ (M5) both identities of \lcRef{Equation}{E:occupation-raw-identities}{11.3.2} and $\langle\Psi,\rho_g\rangle\le\sqrt{v_\Psi A(g)}$.}
\lcPart{Not yet formalized}{These identities and this bound are not restated on the class $\mathcal{O}$ itself.}
\lcEntryEnd
\lcEntry{\lcRef{Lemma}{L:occupation-rayleigh-selection}{11.3.3}}{Countable Rayleigh near-occupations}
\lcPart{Lean proves}{\lcRef{Equation}{E:occupation-rayleigh-selection}{11.3.6} with the section's hypotheses and nothing supplied, with $g_\varepsilon$ in the class $\mathcal{O}$ itself (strongly $H^1$-measurable fields, read as a.e.-strongly measurable); $\lambda$ is the square root of the real Rayleigh supremum and $v_\Psi$ the integral-form cost.}
\lcPart{Not yet formalized}{The identification of $\lambda$ and $v_\Psi$ with $\sqrt{-E_s}$ and $\mathbb{E}W(\Psi)^2$ is not formalized.}
\lcEntryEnd
\lcEntry{\lcRef{Theorem}{T:intrinsic-occupation-ceiling}{11.3.4}}{The computed occupation ceiling}
\lcPart{Lean proves}{With the section's hypotheses and nothing supplied, for $\mathfrak{c}>0$: $\mathfrak{c}I_\Psi-v_\Psi/2\le U_{\mathfrak{c}}$ for nonnegative smooth controls, $S_{\mathfrak{c}}^{\mathrm{pos}}\le U_{\mathfrak{c}}$ with the exact $C_{\rm occ}$, and, for every $\zeta>0$, $U_{\mathfrak{c}}=\mathfrak{c}^m(2\zeta)^{\alpha/(2k)}F_0(\mathcal{E}_\zeta^{\mathcal{O}})$, over the class $\mathcal{O}$ itself; $\lambda$ and $v_\Psi$ in Rayleigh and integral form.}
\lcPart{Not yet formalized}{$C_2(\Gamma,a)\le U_2$ is not formalized.}
\lcEntryEnd
\lcEntry{\lcRef{Corollary}{C:intrinsic-spatial-bracket}{11.3.5}}{A common spatial-energy bracket}
\lcPart{Lean proves}{With the section's hypotheses and nothing supplied, \lcRef{Equation}{E:intrinsic-spatial-bracket}{11.3.9} and the comparison $S_2^{\mathrm{pos}}\le2^mF_0(\mathcal{E}_{1/2}^{\mathcal{O}})<2^mF_a(\mathcal{E}_{1/2}^{\mathcal{O}})$ for the variational value $S_2^{\mathrm{pos}}$ in place of $C_2(\Gamma,a)$ (identified by \lcRef{Theorem}{T:intrinsic-second-moment}{11.2.6}), with $\mathcal{E}_{1/2}^{\mathcal{O}}$ over the class $\mathcal{O}$ itself; $\lambda$ is the square root of the real Rayleigh supremum.}
\lcEntryEnd

\lcCategory{Every conclusion proved; named analytic inputs supplied as hypotheses.}{1}
\lcEntry{\lcRef{Proposition}{P:tube-rayleigh-transfer}{10.2.3}}{Comparison with the Rayleigh integral, including jumps and zero law}
\lcInputs{Inputs supplied}
\lcInput{1}{\texttt{TubeAttained}: every tube $\mathcal{T}_r$ is closed and $\{d_\delta\le r\}=\mathcal{T}_r$ (first part of \lcRef{Equation}{E:tube-attained-radius}{10.2.4} of \lcRef{Lemma}{L:tube-closed-geometry}{10.2.2}, proved there from the compactness of $\mathcal{M}_\delta$; nonzero law only)}
\lcInput{2}{\texttt{TubeTail}: $\mathbb{P}(d_\delta(X)>r)\le e^{-r^2/(2\delta)}$ for $r\ge0$ (second part of \lcRef{Equation}{E:tube-attained-radius}{10.2.4}, proved there from the imported Gaussian isoperimetric theorem of Borell and Sudakov--Tsirelson; nonzero law only)}
\lcEntryEnd

\lcCategory{Steps checked, not the result itself; analytic inputs either supplied as hypotheses or left outside Lean.}{34}
\lcEntry{\lcRef{Lemma}{L:white-test}{2.2.1}}{White-space test}
\lcPart{Lean proves}{With nothing supplied, every value of \lcRef{Equation}{E:white-test}{2.2.1} and $A_1(t)=2c_0t^{3-\alpha_0}$ for the diagram integrals of the model M4 (white noise, $d=1$, $\alpha_0<1$); the statements about $\E X_1^2$, $\E X_2$ and the contribution to \lcRef{Equation}{E:full-second}{2.1.3} follow from the supplied mixed-moment expansion of \lcRef{Proposition}{P:mixed}{2.1.1} for random variables standing for the chaos terms of the solution.}
\lcPart{Not yet formalized}{The solution is not constructed.}
\lcInputs{Input supplied}
\lcInput{1}{\texttt{MixedMoment}: the statement of \lcRef{Proposition}{P:mixed}{2.1.1} (\lcRef{Equation}{E:mixed-moment}{2.1.2} and \lcRef{Equation}{E:full-second}{2.1.3}): the diagram limits exist and are finite, $\E[X_mX_n]=\sum_DJ_{m,n,D}(t)$ with the products integrable, $u^2$ integrable and $\E u^2=\sum_kA_k(t)$; true for the solution under \lcRef{Hypothesis}{H:homog}{1.2.1}, not proved}
\lcEntryEnd
\lcEntry{\lcRef{Proposition}{P:mild-construction}{2.6.2}}{Mild construction from the fixed-order limits}
\lcPart{Lean proves}{Only the multiplication bound \lcRef{Equation}{E:occupied-multiplication}{2.6.4} of the proof, for $F$ a real polynomial of total degree at most $n$ in finitely many independent standard Gaussians and $H$ a real combination of them.}
\lcPart{Not yet formalized}{The extension of this bound to the closed space of degree at most $n$, the almost-sure identification and the mild-solution statement \lcRef{Equation}{E:mild-construction}{2.6.3} are not formalized.}
\lcInputs{Inputs supplied}
\lcInput{1}{the closure step of \lcRef{Equation}{E:occupied-multiplication}{2.6.4} (approximation in $\mathscr{P}_{\le n}$, the $L^2$ Cauchy property of the products and the almost-sure identification of the limit with $HF$); not formalized}
\lcInput{2}{the fixed-order limits $X_n$ in $\mathscr{P}_{\le n}$ (in the proof of \lcRef{Proposition}{P:mild-construction}{2.6.2}) and the proof of \lcRef{Equation}{E:mild-construction}{2.6.3}; not formalized}
\lcEntryEnd
\lcEntry{\lcRef{Lemma}{L:complex-phase}{3.2.1}}{A complex phase bound}
\lcPart{Lean proves}{With the lemma's hypotheses, the numerical-range and scalar steps on $W(L)$ with nothing supplied, and, with the inputs below, the bounds $3$, $6(1+\lambda h)$, $6$ (times $e^{\lambda h}$) for $C_h(L)$, $F_h(L)$, $H_h(L)$, the root comparisons $7/6$, $6/5$ and the bound $100e^{2\lambda h}$ for the weighted block map.}
\lcPart{Not yet formalized}{The identification of the weighted block map with $\mathcal{M}_z(h)$ is not formalized.}
\lcInputs{Inputs supplied}
\lcInput{1}{\lcRef{Equation}{E:cp-input}{3.2.1} (Crouzeix--Palencia, imported in the manuscript) for $C_h$, $F_h$, $H_h$ at $L=A_0-zB$ with its holomorphy precondition: operators of norm at most $(1+\sqrt2)M$ whenever the scalar function is bounded by $M$ on $W(L)$}
\lcInput{2}{the inverse-root integral representation in the proof of \lcRef{Lemma}{L:complex-phase}{3.2.1} for $L^*=A_0-\overline{z}B$ (bounded holomorphic calculus): $Q'-1=\pi^{-1}\int_0^\infty J(t)\ud t$ with $J$ strongly measurable and $\|J(t)\|\le(|z|d/(1-\delta_a))t^{-1/2}\kappa^{a-3}$}
\lcEntryEnd
\lcEntry{\lcRef{Proposition}{P:static-phase-all-orders}{3.2.4}}{All regularized static phase orders}
\lcPart{Lean proves}{The Cauchy coefficient estimate, the arithmetic at $\lambda=n/h$ and the series bound of the proof, and from them both displays for abstract Banach-space quantities, with the circle bound (or the inputs of \lcRef{Lemma}{L:complex-phase}{3.2.1} and \lcRef{Equation}{E:fock-sandwich}{3.2.5}), entireness, the coefficient identification, the norm comparisons, the cutoff removals and the free bound supplied.}
\lcPart{Not yet formalized}{Phases are not constructed.}
\lcInputs{Inputs supplied}
\lcInput{1}{the bound $100e^{2\lambda h}$ for the phase propagator on the coupling disc of radius $\delta_{\alpha/2}\lambda^{2-\alpha/2}/d'$, for every $d'\ge d_n$, $d'>0$ (in the proof of \lcRef{Proposition}{P:static-phase-all-orders}{3.2.4}), supplied in \texttt{StaticPhaseCauchy}}
\lcInput{2}{in \texttt{static\_phase\_all\_orders\_cp} the circle bound is derived from \lcRef{Lemma}{L:complex-phase}{3.2.1} and three named inputs: \texttt{FockSandwich} (\lcRef{Equation}{E:fock-sandwich}{3.2.5} in operator form), \texttt{CPInputs} (\lcRef{Equation}{E:cp-input}{3.2.1}, inverse-root integral) and \texttt{PropagatorBlock}, with positive $A_0$ and self-adjoint $B$ at every cutoff (supplied, with the named inputs)}
\lcInput{3}{in \texttt{FockSandwichModel.static\_phase\_all\_orders\_model} (wave L12) \texttt{FockSandwich}, $A_0\ge0$ and $B=B^*$ are proved on the Fourier-side model, not identified with physical space}
\lcInput{4}{the added restriction of \texttt{static\_phase\_all\_orders\_model}: an orthonormal real family $\psi_m$ of size $m$ in $L^2(\mu_\epsilon)$ for every Gaussian cutoff $m$ (so $\mu\ne0$; the case $\mu=0$, where every positive-order term vanishes, is not covered by that wrapper)}
\lcInput{5}{the phase entries are entire in the coupling (in the proof of \lcRef{Proposition}{P:static-phase-all-orders}{3.2.4})}
\lcInput{6}{the identification of the $n$-th coefficient, applied to $S_b(f,g)$ times the vacuum, with the cutoff $n$-th Duhamel phase via \lcRef{Lemma}{L:finite-word-cutoff}{3.2.3}, with $\|z\|_{1/h}\le n\|z\|_{n/h}$; \lcRef{Lemma}{L:finite-word-cutoff}{3.2.3} is proved with several coordinates per color (\texttt{FiniteWordCutoff}), but the identification of the phase coefficients is not formalized}
\lcInput{7}{$\|S_b(f,g)\|_{n/h}\le\|(f,g)\|_{1/h}$ (in the proof of \lcRef{Proposition}{P:static-phase-all-orders}{3.2.4})}
\lcInput{8}{the removal of the band and Gaussian cutoffs, as convergence of the norms (in the proof of \lcRef{Proposition}{P:static-phase-all-orders}{3.2.4})}
\lcInput{9}{the free bound $\|Z_0(h)\|_{1/h}\le\sqrt e\,\|(f,g)\|_{1/h}$ (in the proof of \lcRef{Proposition}{P:static-phase-all-orders}{3.2.4})}
\lcEntryEnd
\lcEntry{\lcRef{Corollary}{C:static-interaction}{3.2.5}}{Summable static interaction bounds}
\lcPart{Lean proves}{\lcRef{Equation}{E:static-block-series}{3.2.10} in full; \lcRef{Equation}{E:static-interaction}{3.2.9}, with the identification of the interaction map with $\sum_n\sigma^nY_n$ and the order bounds derived from the primitive inputs below; and \lcRef{Equation}{E:independent-static-blocks}{3.2.11} from the per-block bound on adapted phases; the Gaussian parity, the first-variance bounds and that per-block bound are supplied.}
\lcInputs{Inputs supplied}
\lcInput{1}{for the order bounds (\lcRef{Equation}{E:all-static-phase-orders}{3.2.7}, $n\ge1$; circle bound from \lcRef{Lemma}{L:complex-phase}{3.2.1}): \texttt{FockSandwich} (\lcRef{Equation}{E:fock-sandwich}{3.2.5} in operator form), \texttt{CPInputs} (\lcRef{Equation}{E:cp-input}{3.2.1} with the inverse-root integral) and \texttt{PropagatorBlock} (the weighted block map as $\mathcal{M}_z(h)$), with positive $A_0$ and self-adjoint $B$ at every cutoff (supplied, with the named inputs, in \texttt{StaticInteractionCompleted})}
\lcInput{2}{\texttt{FockSandwich}, $A_0\ge0$ and $B=B^*$ are proved in wave L12 on the Fourier-side model (\texttt{FockSandwichModel.static\_interaction\_model\_paper}), with the added restriction of an orthonormal real family $\psi_m$ of size $m$ in $L^2(\mu_\epsilon)$ for every Gaussian cutoff $m$ (so $\mu\ne0$; the case $\mu=0$, where every positive-order term vanishes, is not covered by that wrapper)}
\lcInput{3}{the remaining inputs of \lcRef{Proposition}{P:static-phase-all-orders}{3.2.4} for the order bounds: entireness in the coupling, the coefficient identification, the band and the cutoff removals (see that row)}
\lcInput{4}{the free estimate $\|\mathcal{S}_0(\mp h)\|_\lambda\le\sqrt e$, the group law $\mathcal{S}_0(-h)\mathcal{S}_0(h)=1$, $Z_0=\mathcal{S}_0(h)(f,g)$ and the norm comparisons $\|Z_n\|_\lambda\le(\lambda h)^{-1}\|Z_n\|_{1/h}$, $\|(f,g)\|_{1/h}\le\|(f,g)\|_\lambda$ (in the proof of \lcRef{Corollary}{C:static-interaction}{3.2.5}); an abstract complete real inner product space stands for $L^2(\Omega;X_\lambda)$}
\lcInput{5}{Gaussian symmetry, $\E\langle Y_r,Y_s\rangle_\lambda=0$ for $r+s$ odd (in the proof of \lcRef{Corollary}{C:static-interaction}{3.2.5})}
\lcInput{6}{the first-variance bounds of \lcRef{Proposition}{P:static-phase-variance}{3.1.2}, $\E\|Y_1\|^2\le e$ times \lcRef{Equation}{E:phase-born}{3.1.6} and $|\E\langle Y_0,Y_2\rangle|\le N\sqrt e$ times \lcRef{Equation}{E:phase-mean}{3.1.7} (in the proof of \lcRef{Corollary}{C:static-interaction}{3.2.5})}
\lcInput{7}{for \lcRef{Equation}{E:independent-static-blocks}{3.2.11}, the per-block bound $\|\mathcal{I}_jy\|^2\le(1+\mathcal{E}(h_j,\sigma_j,\lambda))\|y\|^2$ on phases adapted to the past (\lcRef{Equation}{E:static-interaction}{3.2.9} with density and tensoring with the independent past, in the proof of \lcRef{Corollary}{C:static-interaction}{3.2.5}), adaptation being an abstract predicate preserved by each block, and $\|\mathcal{S}_0(h_j)\|^2\le e^{\lambda h_j}$}
\lcEntryEnd
\lcEntry{\lcRef{Lemma}{L:interval-complex}{3.3.2}}{Complex sums and changing phase weights}
\lcPart{Lean proves}{With the lemma's hypotheses, the adjoint sandwich bound and the numerical-range and scalar steps on $W(A_0-V)$ with nothing supplied; with the inputs below, the bound $100e^{2\Lambda h}$ for the weighted block map; and \lcRef{Equation}{E:interval-metric}{3.3.3} for given finite jump data and, for the scalar conversion factors, along any monotone time partition.}
\lcPart{Not yet formalized}{The weighted block map is not identified with the propagator in $X_\Lambda$.}
\lcInputs{Inputs supplied}
\lcInput{1}{\lcRef{Equation}{E:cp-input}{3.2.1} (Crouzeix--Palencia) for $C_h$, $F_h$, $H_h$ at $L=A_0-V$ with weight $\Lambda$: operators of norm at most $(1+\sqrt2)M$ whenever the scalar function is bounded by $M$ on $W(L)$}
\lcInput{2}{the inverse-root integral representation in the proof of \lcRef{Lemma}{L:interval-complex}{3.3.2} for $L^*=A_0-V^*$: $Q'-1=\pi^{-1}\int_0^\infty J(t)\ud t$ with $J$ strongly measurable and $\|J(t)\|\le(D/(1-\delta_a))t^{-1/2}\kappa^{a-3}$}
\lcEntryEnd
\lcEntry{\lcRef{Lemma}{L:fixed-split}{3.3.4}}{A fixed degree-split estimate}
\lcPart{Lean proves}{\lcRef{Equation}{E:fixed-split-bound}{3.3.6} with $C_{\rm split}=100^8e^8\max\{1,12\delta^{-2}\}$ for any $\delta>0$, on a finite Fock model of the cutoff Gaussian space with abstract spatial operators, for coefficients of an evolution; the segment bound, the Dyson cube identity, the invariance of the cross subspaces, entire coefficients and \lcRef{Equation}{E:interval-column}{3.3.1} are supplied.}
\lcPart{Not yet formalized}{The evolution is not constructed.}
\lcInputs{Inputs supplied}
\lcInput{1}{\lcRef{Equation}{E:interval-column}{3.3.1} for the modes of each color at every $\rho>0$, as the named input \texttt{FockColumn} of \texttt{FockForm}}
\lcInput{2}{the segment bound of \lcRef{Lemma}{L:interval-complex}{3.3.2} for the background propagators, $\|P_z(t,s)\|\le100e^{2\Lambda(t-s)}$ on every interval free of background endpoints, as applied in \lcRef{Equation}{E:split-segment-product}{3.3.9} (\texttt{SegmentBound}); its operator theory is not formalized}
\lcInput{3}{\lcRef{Equation}{E:interval-cross-cube}{3.3.13} on the torus of the background radii, with the integrability of the word on each simplex (\texttt{DysonCrossCube}); variation of constants is not formalized}
\lcInput{4}{the invariance of every cross subspace $\mathscr{Z}_S$ under the background propagator (in the proof of \lcRef{Lemma}{L:fixed-split}{3.3.4}; \texttt{BackgroundPreserves}) and the absolute convergence of the ordinary coefficients on every polydisc (in the proof of \lcRef{Lemma}{L:fixed-split}{3.3.4}; \texttt{EntireCoefficients}); the cutoff evolution, its propagators and the coefficients $U_{\boldsymbol r}v$ are not constructed}
\lcEntryEnd
\lcEntry{\lcRef{Proposition}{P:interval-coefficients}{3.3.5}}{Uniform interval coefficients}
\lcPart{Lean proves}{\lcRef{Equation}{E:interval-coefficient-bound}{3.3.19} (length substitution, split count, product-coefficient identity), assembled from \lcRef{Equation}{E:interval-energy-splits}{3.3.22}, with \lcRef{Lemma}{L:inventory-product}{3.3.3} proved, \lcRef{Lemma}{L:fixed-split}{3.3.4} assembled on its supplied inputs, $Q_{r,q}$ identified with the cutoff inner product (real part used) from the supplied Duhamel word expansion, for $R_\alpha,\sigma_i>0$, and the sum to \lcRef{Equation}{E:interval-phase-bound}{3.3.20}.}
\lcInputs{Inputs supplied}
\lcInput{1}{\lcRef{Equation}{E:interval-energy-splits}{3.3.22} combined with the reflection parity (in the proof of \lcRef{Proposition}{P:interval-coefficients}{3.3.5}), the coefficient as the sum over degree vectors $\boldsymbol k$ ($\sum_ik_i=k$) and splits $r_i+q_i=2k_i$ of $Q_{\boldsymbol r,\boldsymbol q}$}
\lcInput{2}{the inputs of \lcRef{Lemma}{L:fixed-split}{3.3.4} for every split (\lcRef{Equation}{E:interval-column}{3.3.1}, the segment bound of \lcRef{Lemma}{L:interval-complex}{3.3.2}, \lcRef{Equation}{E:interval-cross-cube}{3.3.13}, the invariance of the cross subspaces and entire coefficients; see the row of \lcRef{Lemma}{L:fixed-split}{3.3.4}), with the interval endpoints $(a_i,b_i)$ in $[0,T]$, $b_i-a_i\le h_i$, and any $\delta>0$ in $C_{\rm split}=100^8e^8\max\{1,12\delta^{-2}\}$}
\lcInput{3}{for \lcRef{Equation}{E:interval-phase-bound}{3.3.20}, the power-series identity at $s=1$, the vanishing of odd coefficients and the zero-order bound $e\|v\|^2$ (all in the proof of \lcRef{Proposition}{P:interval-coefficients}{3.3.5})}
\lcInput{4}{the Duhamel word expansion of the ordinary coefficients (in the proof of \lcRef{Proposition}{P:interval-coefficients}{3.3.5}; \texttt{DuhamelWords}): $U_{\boldsymbol r}v$ and the cutoff coefficient $c_{\boldsymbol r}$ are the same finite sum of integrals of words with $r_i$ insertions of color $i$, uncut and cut; the evolutions are not constructed}
\lcInput{5}{from \texttt{DuhamelWords} Lean proves $|Q_{\boldsymbol r,\boldsymbol q}|\le|\langle c_{\boldsymbol r},c_{\boldsymbol q}\rangle|$, which the earlier step row assumed as the added hypothesis \texttt{hQ} (accepted by Le for now on 2026-11-12); $Q_{\boldsymbol r,\boldsymbol q}$ is the real part of $\langle U_{\boldsymbol r}v,U_{\boldsymbol q}v\rangle_{1/T}$}
\lcInput{6}{the reduction to $\sigma_i>0$ (zero noise and sign change, in the proof of \lcRef{Proposition}{P:interval-coefficients}{3.3.5}); not formalized}
\lcEntryEnd
\lcEntry{\lcRef{Proposition}{P:masked-noise-completion}{4.1.4}}{Completion on the same noise}
\lcPart{Lean proves}{Only the finite Gaussian comparison: for $A$ and $B-A$ positive semidefinite, every seminorm $N$ and $p\ge1$, $\int N^p\ud\mathcal{N}(0,A)\le\int N^p\ud\mathcal{N}(0,B)$ with both moments finite.}
\lcPart{Not yet formalized}{The identification of $A$ and $B$ with the masked covariance and auxiliary Gram matrices, \lcRef{Equation}{E:masked-actual-noise-norm}{4.1.8}, the contraction extension and joint Gaussianity are not formalized.}
\lcInputs{Inputs supplied}
\lcInput{1}{the covariance matrix $A$ of the masked evaluations $W(\zeta_RF_1),\ldots,W(\zeta_RF_m)$ and the auxiliary Gram matrix $B$ of $N_R$ (in the statement and the proof of \lcRef{Proposition}{P:masked-noise-completion}{4.1.4}), which need the noise $W$ and the completion of the test space}
\lcInput{2}{$A\le B$ in positive semidefinite order, from $\E W(\zeta_RF)^2\le N_R(F)^2$ (\lcRef{Proposition}{P:masked-auxiliary-norm}{4.1.2} with $H=H_a^{\mathbb{C}}$ and \lcRef{Equation}{E:masked-actual-noise-norm}{4.1.8}, in the proof of \lcRef{Proposition}{P:masked-noise-completion}{4.1.4})}
\lcInput{3}{the joint centered Gaussian law of the masked evaluations, so that their law is $\mathcal{N}(0,A)$ (in the statement and the proof of \lcRef{Proposition}{P:masked-noise-completion}{4.1.4})}
\lcEntryEnd
\lcEntry{\lcRef{Proposition}{P:actual-wave-form}{4.2.3}}{Form and phase kernel}
\lcPart{Lean proves}{Only the finite, algebraic or scalar steps of this result listed in the companion's statement ledger.}
\lcPart{Not yet formalized}{The result itself.}
\lcPart{Assumptions and limits}{The checked steps take no analytic input beyond the paper's hypotheses: the Lean statements are definitions and a matrix identity for arbitrary real parameters. The reading of the parameters as cosines and sines of a phase and as a covariance coefficient, the continuous extension at zero frequency, and the Bessel-weighted change of coordinates and the weight $1/\sqrt{\langle\xi\rangle\langle\eta\rangle}$ that lead to \lcRef{Equation}{E:variance-normalized-kernel}{4.2.17} are not hypotheses of these statements; they are outside Lean (see the statement-ledger entry \texttt{DimensionalMatrix}).}
\lcEntryEnd
\lcEntry{\lcRef{Proposition}{P:wave-two-variances}{4.2.5}}{Both variance directions}
\lcPart{Lean proves}{A bound on an abstract variance tied to $\E|\ell_R|^2$ only through the named input \lcRef{Equation}{E:variance-normalized-covariance}{4.2.18}. The spectral core (the sums of the $\nu_R$-integrals of \lcRef{Equation}{E:variance-normalized-covariance}{4.2.18}, both directions, the $h^{2\beta}$ form) in full with no input.}
\lcPart{Not yet formalized}{$\ell_R$, the spaces $E_r$ and $J$ are not constructed.}
\lcInputs{Inputs supplied}
\lcInput{1}{\lcRef{Equation}{E:variance-normalized-covariance}{4.2.18} (last claim of \lcRef{Proposition}{P:actual-wave-form}{4.2.3}, needs the construction of the noise), supplied as \texttt{hlink}, \texttt{ofReal(var R t0 h z w)} $\le$ the $\nu_R$ lower integral of $\|\texttt{coef}(Jz,Jw)\|^2$, the manuscript's display with $U=Jz$, $V=Jw$}
\lcInput{2}{$J$ maps the phases into $L^2(\R^d;\mathbb{C}^2)$ preserving inner products (\texttt{hJL2}, \texttt{hJ}; at the manuscript's definition of $J$), weaker than the manuscript's unitary $J$}
\lcInput{3}{the identity $\|z\|_{E_{s-\gamma}}^2=\|\langle\xi\rangle^{-\gamma}Jz\|_2^2$ (\texttt{hnormSG}; in the proof of \lcRef{Proposition}{P:wave-two-variances}{4.2.5})}
\lcEntryEnd
\lcEntry{\lcRef{Lemma}{L:driver-interval-limit}{5.1.8}}{Fixed-interval operator construction}
\lcPart{Lean proves}{Only the finite, algebraic or scalar steps of this result listed in the companion's statement ledger.}
\lcPart{Not yet formalized}{The result itself.}
\lcPart{Assumptions and limits}{The checked step takes no analytic input beyond the paper's hypotheses: the Lean statement assumes only $C\ge0$, $h>0$ and $L\ge0$. The variance budget \lcRef{Equation}{E:driver-short-variance-budget}{5.1.12} and the Schatten estimate it feeds are not hypotheses of this statement; they are outside Lean (see the statement-ledger entry \texttt{InteractionDriver}).}
\lcEntryEnd
\lcEntry{\lcRef{Lemma}{L:driver-holder-limit}{5.1.9}}{Hölder path limit}
\lcPart{Lean proves}{Only the finite, algebraic or scalar steps of this result listed in the companion's statement ledger.}
\lcPart{Not yet formalized}{The result itself.}
\lcPart{Assumptions and limits}{The checked steps take no analytic input beyond the paper's hypotheses: the Lean statements assume only order conditions on the real exponents. The Gaussian moment bounds for increments, the dyadic chaining argument and the H\"older interpolation are not hypotheses of these statements; they are outside Lean (see the statement-ledger entry \texttt{InteractionDriver}).}
\lcEntryEnd
\lcEntry{\lcRef{Lemma}{L:exact-canonical-scaling}{7.1.1}}{Canonical law scaling}
\lcPart{Lean proves}{Only the finite, algebraic or scalar steps of this result listed in the companion's statement ledger.}
\lcPart{Not yet formalized}{The result itself.}
\lcPart{Assumptions and limits}{The checked steps take no analytic input beyond the paper's hypotheses: the Lean statements assume only $T>0$ and $\alpha<3$. The scaling of the Gaussian noise and its covariance, the equality in law of the rescaled solution, the scaling of the Green measure and the extension to the full series and its second moment are not hypotheses of these statements; they are outside Lean (see the statement-ledger entry \texttt{DirectEvent}).}
\lcEntryEnd
\lcEntry{\lcRef{Lemma}{L:raw-block-inputs}{7.2.1}}{Coupled Gaussian block inputs}
\lcPart{Lean proves}{Only the finite, algebraic or scalar steps of this result listed in the companion's statement ledger.}
\lcPart{Not yet formalized}{The result itself.}
\lcPart{Assumptions and limits}{The checked step takes no analytic input beyond the paper's hypotheses: the Lean statement assumes only $\delta>0$. The rescaled noise and its mean shift are not hypotheses of this statement; they are outside Lean (see the statement-ledger entry \texttt{JointTilt}).}
\lcEntryEnd
\lcEntry{\lcRef{Lemma}{L:raw-block-indexed-cone}{7.2.4}}{Indexed cone and terminal displacement}
\lcPart{Lean proves}{Every conclusion, for an abstract ground family (flows, vectors, duals, rates and observable, with the properties of the compact mean family supplied as inputs).}
\lcPart{Not yet formalized}{The identification of the abstract ground family with the physical wave family of the compact mean family is not formalized.}
\lcInputs{Inputs supplied}
\lcInput{1}{\texttt{GroundFamily} data fields \texttt{P}, \texttt{e}, \texttt{ell}, \texttt{lam} (\lcRef{Equation}{E:compact-F-block}{6.3.11}, \lcRef{Equation}{E:compact-F-projection}{6.3.13}, \lcRef{Proposition}{P:compact-ground-family}{6.3.1}): flows $P_\sigma(t)$, vectors $e_\sigma$, bounded linear duals $\ell_\sigma$, rates $\lambda_\sigma$, for $\sigma\in I=[0,1]$, $t\ge0$ ($\Pi_\sigma=e_\sigma\ell_\sigma$ and $Q_\sigma=I-\Pi_\sigma$ are derived definitions \texttt{GroundFamily.Pi}, \texttt{GroundFamily.Q}, not fields)}
\lcInput{2}{\texttt{GroundFamily} data fields \texttt{e0}, \texttt{obs} (\lcRef{Section}{S:actual-raw-block-lower}{7.2}): the data $e_0\in F$ and the bounded linear observable $O_\phi:F\to\R$}
\lcInput{3}{\texttt{GroundFamily} fields \texttt{d0}, \texttt{lam\_pos}, \texttt{d0\_pos} (\lcRef{Proposition}{P:compact-ground-family}{6.3.1}, \lcRef{Equation}{E:compact-rate-gap}{6.3.4}): a gap $d_0>0$ and $\lambda_\sigma>0$ on $I$ (weakening of $\sqrt\varepsilon\le\lambda_\sigma$)}
\lcInput{4}{\texttt{GroundFamily} fields \texttt{ell\_e}, \texttt{P\_Pi}, \texttt{Pi\_P} (\lcRef{Proposition}{P:compact-F-projection}{6.3.4}, \lcRef{Equation}{E:compact-F-eigen}{6.3.14}): $\ell_\sigma(e_\sigma)=1$ and $P_\sigma(t)\Pi_\sigma=\Pi_\sigma P_\sigma(t)=e^{\lambda_\sigma t}\Pi_\sigma$ on $I$ for $t\ge0$}
\lcInput{5}{\texttt{GroundFamily} field \texttt{normalized\_complement} (\lcRef{Proposition}{P:compact-F-complement}{6.3.5}, \lcRef{Equation}{E:compact-normalized-complement}{6.3.17}): for every $0<\eta<d_0$ a $C_\eta$ with $\|e^{-\lambda_\sigma t}P_\sigma(t)Q_\sigma\|\le C_\eta e^{-\eta t}$ uniformly over $\sigma\in I$, $t\ge0$}
\lcInput{6}{\texttt{GroundFamily} fields \texttt{e\_cont}, \texttt{ell\_cont} (\lcRef{Proposition}{P:compact-F-projection}{6.3.4}): norm continuity of $\sigma\mapsto e_\sigma$ and $\sigma\mapsto\ell_\sigma$ on $I$}
\lcInput{7}{\texttt{GroundFamily} fields \texttt{margin\_ell}, \texttt{margin\_obs} (\lcRef{Equation}{E:compact-positive-margins}{6.3.15}): $\inf_I\ell_\sigma(e_0)>0$ and $\inf_IO_\phi(e_\sigma)>0$}
\lcEntryEnd
\lcEntry{\lcRef{Lemma}{L:raw-block-restart}{7.2.5}}{Composition over time blocks and regularized products}
\lcPart{Lean proves}{Only the finite, algebraic or scalar steps of this result listed in the companion's statement ledger.}
\lcPart{Not yet formalized}{The result itself.}
\lcInputs{Input supplied}
\lcInput{1}{The Lean statement takes as hypotheses: (1) the recurrence $d(n+1)\le Cd(n)+C^ne_{\mathrm{src}}$ for every $n$, for the difference $d(n)$ of the first $n$ ordered products (\texttt{hstep}); in the proof of \lcRef{Lemma}{L:raw-block-restart}{7.2.5} it comes from operator-norm submultiplicativity for the product factors, after the finite factor lists are extended by identities; (2) $d(0)\le0$ (\texttt{hzero}). The conditions $C\ge1$ and $e_{\mathrm{src}}\ge0$ are parameter conditions. The moments of $C$ and of the error are outside Lean (see the statement-ledger entry \texttt{ActualEvent}).}
\lcEntryEnd
\lcEntry{\lcRef{Proposition}{P:raw-block-original-lower}{7.2.7}}{Original point lower rate at the Gaussian quadratic cost}
\lcPart{Lean proves}{Only the finite, algebraic or scalar steps of this result listed in the companion's statement ledger.}
\lcPart{Not yet formalized}{The result itself.}
\lcPart{Assumptions and limits}{The checked step takes no analytic input beyond the paper's hypotheses: the Lean statement assumes only $\delta>0$. The Gaussian change of law, the event on which the penalty bound holds and the limit $\delta\to0$ are not hypotheses of this statement; they are outside Lean (see the statement-ledger entry \texttt{ActualEvent}).}
\lcEntryEnd
\lcEntry{\lcRef{Proposition}{P:variational-remove-mask}{7.3.4}}{Monotone rates and the original point moment}
\lcPart{Lean proves}{The mask family and every rate statement with the section's hypotheses and nothing supplied, for $\lambda_{R,s}$ defined as the square root of the real Rayleigh supremum of $m_Rh_s$; \lcRef{Equation}{E:variational-unmasked-point-lower}{7.3.5}, derived from \lcRef{Equation}{E:raw-block-original-lower}{7.2.14} as a supplied input, for a real function in place of the point moment.}
\lcPart{Not yet formalized}{Compact ground rates are not formalized.}
\lcInputs{Input supplied}
\lcInput{1}{\lcRef{Equation}{E:raw-block-original-lower}{7.2.14} of \lcRef{Proposition}{P:raw-block-original-lower}{7.2.7} for the masks $m_R$, all large $R$ and one fixed $x$, with the compact ground rate of $m_Rh_s$ read as its Rayleigh rate and $v_\Psi$ as its integral form; used only by the point-moment display}
\lcEntryEnd
\lcEntry{\lcRef{Proposition}{P:radius-mean-anchor}{9.2.5}}{Averaged block mean tends to zero}
\lcPart{Lean proves}{Every conclusion, for abstract random paths with the supplied moment bound.}
\lcPart{Not yet formalized}{The operator paths $M_R(Y_m)$ and the residual noises are not constructed, so the identification of these paths with the manuscript's block paths is not formalized.}
\lcInputs{Inputs supplied}
\lcInput{1}{the moment bound of \lcRef{Corollary}{C:radius-covariance-residuals}{9.2.4} (from \lcRef{Proposition}{P:radius-actual-moments}{9.2.3}) on the unstopped block paths: almost surely $C^\theta$, $p$-th norm power integrable, $(\E\|M_R(Y_m)\|^p)^{1/p}\le C\sqrt p\,R^{d/Q}$}
\lcInput{2}{measurability of the $C^\theta$ norm of the stopped block paths (implicit in the manuscript)}
\lcEntryEnd
\lcEntry{\lcRef{Proposition}{P:signed-closed-form}{9.3.2}}{Closed form and quarter-power phase}
\lcPart{Lean proves}{For an abstract bounded Hermitian pairing on Mathlib's $H^1$ with \lcRef{Equation}{E:signed-multiplier-form}{9.3.1} supplied, \lcRef{Equation}{E:signed-form-bounds}{9.3.3}, the form version of $A\ge1$, the step $\|u\|_{H^{1/2}}^2\le bz$ and the bounded resolvent of the form (self-adjoint, positive, injective, norm at most $1$, unique weak solution).}
\lcPart{Not yet formalized}{$H_V$, $D(A^{1/2})=H^1$, \lcRef{Equation}{E:signed-quarter-phase}{9.3.4} and $H^{-1/2}$ are not formalized.}
\lcInputs{Inputs supplied}
\lcInput{1}{\texttt{SignedMultiplierForm P k}: \lcRef{Equation}{E:signed-multiplier-form}{9.3.1}, the conclusion of \lcRef{Lemma}{L:signed-limiting-multiplier}{9.3.1} (not formalized), for an abstract bounded operator $P$ on $H^1$ standing for the pairing $v_V$; the multiplier $V$ is not constructed}
\lcInput{2}{the self-adjointness of $P$: the Hermitian pairing of the proof of \lcRef{Lemma}{L:signed-limiting-multiplier}{9.3.1}}
\lcEntryEnd
\lcEntry{\lcRef{Lemma}{L:symplectic-free}{10.1.1}}{Free wave form and compact Young baseline}
\lcPart{Lean proves}{Both compact-baseline displays on an abstract nonzero Hilbert space with a supplied free group, isometric velocity reversal and weakly vanishing translations, for a Young evolution given only by convergence of its left Young sums; and $W_t^{\mathsf{T}}J_{\rm symp}W_t=J_{\rm symp}$ per frequency and the force skewness for the $L^2$ pairing.}
\lcPart{Not yet formalized}{$J_E$ and the preservation of $\omega$ on $E$ are not formalized.}
\lcInputs{Inputs supplied}
\lcInput{1}{the free wave group on $E=H^{1/2}\times H^{-1/2}$ with its group law, the isometric velocity reversal with $S_E(-h)=\mathcal{R}S_E(h)\mathcal{R}$, and the spatial translations, isometric, commuting with $S_E(h)$ and weakly null as $|y|\to\infty$ (in the proof of \lcRef{Lemma}{L:symplectic-free}{10.1.1}); $E$ is an abstract nonzero real Hilbert space}
\lcInput{2}{\lcRef{Proposition}{P:constant-young-linear}{6.2.4} on $\mathcal{L}(E)$ (in the proof of \lcRef{Lemma}{L:symplectic-free}{10.1.1}): the left Young sums of $X$ over partitions of $[s,t]$ converge in operator norm to $U_X(t,s)-I$; supplied in the weaker form of convergence along some nontrivial filter of finite node lists in $[s,t]$ (no ordering, endpoint or mesh condition)}
\lcEntryEnd
\lcEntry{\lcRef{Proposition}{P:gram-positive-optimizer}{11.1.3}}{Positive conic maximizer}
\lcPart{Lean proves}{Only the finite, algebraic or scalar steps of this result listed in the companion's statement ledger.}
\lcPart{Not yet formalized}{The result itself.}
\lcPart{Assumptions and limits}{The checked step takes no analytic input beyond the paper's hypotheses: the Lean statement assumes only positive weights and active labels, $\langle\Phi,k_i\rangle-b_i\ge0$. Dual coercivity, the existence and uniqueness of the dual minimizer, attainment and the maximum over active sets are outside Lean (see the statement-ledger entry \texttt{PositiveGram}).}
\lcEntryEnd
\lcEntry{\lcRef{Proposition}{P:gram-entrywise-comparison}{11.1.5}}{Entrywise Gram comparison}
\lcPart{Lean proves}{Only the finite, algebraic or scalar steps of this result listed in the companion's statement ledger.}
\lcPart{Not yet formalized}{The result itself.}
\lcPart{Assumptions and limits}{The checked step takes no analytic input beyond the paper's hypotheses: the Lean statement assumes only nonnegative weights and the entrywise order $C_{ij}\le D_{ij}$, a hypothesis of \lcRef{Proposition}{P:gram-entrywise-comparison}{11.1.5}. The Gram realization (\lcRef{Lemma}{L:gram-realization}{11.1.4}) and the optimization over weights are outside Lean (see the statement-ledger entry \texttt{PositiveGram}).}
\lcEntryEnd
\lcEntry{\lcRef{Proposition}{P:white-spatial-energy}{11.4.2}}{Exact white spatial energy}
\lcPart{Lean proves}{Only the finite, algebraic or scalar steps of this result listed in the companion's statement ledger.}
\lcPart{Not yet formalized}{The result itself.}
\lcInputs{Input supplied}
\lcInput{1}{The Lean statements take as hypothesis: (1) the inequality \lcRef{Equation}{E:white-cdf-square}{11.4.1} of \lcRef{Lemma}{L:white-cdf-square}{11.4.1}, $0\le b-2zA+z^2/3$ for every real $z$, with $A=\int_\R f^4$ and $b=\int_\R(f')^2$ for a real unit $f\in H^1(\R)$ (hypothesis \texttt{h} of \texttt{WhiteCDF.cdf\_square\_bound}). The bound $3A^2\le b$ that \texttt{WhiteCDF.energy\_upper} assumes is the conclusion of \texttt{cdf\_square\_bound}. The white-model statements assume only $c\ge0$, and $0<a<1$ for \lcRef{Hypothesis}{H:homog}{1.2.1}. The variational energy, its domain and its attainment are outside Lean (see the statement-ledger entries \texttt{WhiteCDF} and \texttt{WhiteModel}).}
\lcEntryEnd
\lcEntry{\lcRef{Theorem}{T:white-moment-bounds}{11.4.3}}{Explicit white moment bounds}
\lcPart{Lean proves}{\lcRef{Equation}{E:white-moment-bounds}{11.4.2} for the variational value $S_2^{\mathrm{pos}}$ in place of $C_2(c,a)$ (identified by \lcRef{Theorem}{T:intrinsic-second-moment}{11.2.6}), for $c>0$ and $0<a<1$, with the value $\mathcal{E}_{\mathrm{sp}}=c^2/6$ of \lcRef{Proposition}{P:white-spatial-energy}{11.4.2} supplied (each bound uses one inequality of it), and that $M_a$, $K_0(a)$ and both bound formulas are positive.}
\lcPart{Not yet formalized}{The time exponent sentence and the case $c=0$ are not formalized.}
\lcInputs{Input supplied}
\lcInput{1}{\lcRef{Proposition}{P:white-spatial-energy}{11.4.2}: the spatial energy of $c\delta_0$ is $c^2/6$ for $c>0$; the lower bound uses only $c^2/6\le\mathcal{E}_{\mathrm{sp}}$ (the sech profile), the upper bound only $\mathcal{E}_{\mathrm{sp}}\le c^2/6$ (the CDF square, \lcRef{Equation}{E:white-cdf-square}{11.4.1}); not proved}
\lcEntryEnd
\lcEntry{\lcRef{Proposition}{P:riesz-small}{A.1.3}}{Riesz finite coefficients}
\lcPart{Lean proves}{With nothing supplied, \lcRef{Equation}{E:riesz-two}{A.1.7} and $A_1=J_{1,1}+2J_{2,0}$ for the diagram integrals of M4, and \lcRef{Equation}{E:riesz-four}{A.1.8} with finiteness for all splits and pairings; the order-two statements about $\E X_2$, $\E X_1^2$ follow from the supplied mixed-moment expansion of \lcRef{Proposition}{P:mixed}{2.1.1}.}
\lcPart{Not yet formalized}{The solution is not constructed. That the total coefficient uses the multiplicities of \lcRef{Equation}{E:total-four}{A.1.3} is not formalized.}
\lcInputs{Input supplied}
\lcInput{1}{\texttt{MixedMoment}: the statement of \lcRef{Proposition}{P:mixed}{2.1.1}: the diagram limits exist and are finite, $\E[X_mX_n]=\sum_DJ_{m,n,D}(t)$ with the products integrable, $u^2$ integrable and $\E u^2=\sum_kA_k(t)$; true for the solution under \lcRef{Hypothesis}{H:homog}{1.2.1}, not proved (order-two wrapper \texttt{riesz\_small\_paper'} only)}
\lcEntryEnd
\lcEntry{\lcRef{Proposition}{P:d1-upper}{B.1.2}}{One-dimensional upper bound}
\lcPart{Lean proves}{\lcRef{Equation}{E:d1-coefficient}{B.1.5} on the model M4 at every noise width, uniformly, with nothing supplied, and for the limit coefficients (Riesz case: the supplied convergence of every diagram as the noise width tends to zero (\lcRef{Lemma}{L:joint-fixed-order}{2.5.2})); the form $\sum_m\E[X_mX_{2k-m}]$ and \lcRef{Equation}{E:d1-upper}{B.1.6} hold for random variables with the supplied mixed-moment expansion of \lcRef{Proposition}{P:mixed}{2.1.1}.}
\lcPart{Not yet formalized}{The solution is not constructed.}
\lcInputs{Inputs supplied}
\lcInput{1}{\texttt{MixedMoment}: the statement of \lcRef{Proposition}{P:mixed}{2.1.1} (with its proof): every diagram limit exists and is finite, $\E[X_mX_n]=\sum_DJ_{m,n,D}(t)$, $u^2$ integrable and $\E u^2=\sum_kA_k(t)$; TRUE for the solution under \lcRef{Hypothesis}{H:homog}{1.2.1}; about random variables standing for the unconstructed solution and chaos terms}
\lcInput{2}{\texttt{DiagramLimits} (Riesz case only): every model diagram converges to \texttt{diagramLim} as $\epsilon\to0+$; part of \lcRef{Lemma}{L:joint-fixed-order}{2.5.2}, invoked in the proof of \lcRef{Proposition}{P:d1-upper}{B.1.2}; TRUE under \lcRef{Hypothesis}{H:homog}{1.2.1}; proved for white noise (\texttt{D1Upper.diagramLimits\_white}) and implied by \texttt{MixedMoment}}
\lcEntryEnd
\lcEntry{\lcRef{Proposition}{P:low-alpha-upper}{B.1.5}}{Upper bound for alpha \textless{}= 1}
\lcPart{Lean proves}{Under \lcRef{Hypothesis}{H:homog}{1.2.1} with $d\in\{1,2,3\}$, $\alpha\le1$ and $V>0$, both inequalities of \lcRef{Equation}{E:low-alpha-coefficient}{B.1.13} on the diagram model at every noise width $\epsilon>0$ with nothing supplied; the first inequality for $A_k(t)=\sum_m\E[X_mX_{2k-m}]$ and $\E u^2\le\exp(C_{\rm dir}^{\rm gen}t^\kappa)$ for random variables with the supplied mixed-moment expansion of \lcRef{Proposition}{P:mixed}{2.1.1}.}
\lcPart{Not yet formalized}{The solution is not constructed.}
\lcInputs{Input supplied}
\lcInput{1}{\texttt{MixedMoment}: the statement of \lcRef{Proposition}{P:mixed}{2.1.1}: the diagram limits exist and are finite (the manuscript proves them by \lcRef{Lemma}{L:joint-fixed-order}{2.5.2} and the $\eta$-limit paragraph in the proof of \lcRef{Proposition}{P:low-alpha-upper}{B.1.5}), $\E[X_mX_n]$ is the sum of the diagram limits, $\E u^2=\sum_kA_k(t)$; TRUE for the solution under \lcRef{Hypothesis}{H:homog}{1.2.1}}
\lcEntryEnd
\lcEntry{\lcRef{Corollary}{C:d1-first-calibration}{B.1.7}}{First-moment calibration}
\lcPart{Lean proves}{$C_{1,*}=C_{2,*}/2^{(4-\alpha)/(3-\alpha)}$ and $2C_{1,*}<C_{2,*}$ with nothing supplied, \lcRef{Equation}{E:d1-first-calibration}{B.1.14} for random variables with the supplied mixed-moment expansion of \lcRef{Proposition}{P:mixed}{2.1.1} and the supplied first-moment series of the chaos expansion, and the conditional rate clause for an abstract first moment.}
\lcPart{Not yet formalized}{The solution is not constructed. The \lcRef{Proposition}{P:Jensen}{1.4.4} comparison is commentary.}
\lcInputs{Inputs supplied}
\lcInput{1}{\texttt{MixedMoment}: the statement of \lcRef{Proposition}{P:mixed}{2.1.1}; here it gives $\E X_n=\sum_DJ_{n,0,D}(t)$ (with $X_0=1$), the odd terms zero and the diagram limits; TRUE for the solution under \lcRef{Hypothesis}{H:homog}{1.2.1}}
\lcInput{2}{\texttt{FirstMomentSeries}: $X_0=1$, $u$ integrable and $\E u=\sum_n\E X_n$ (the definition of $X_n$ at the start of \lcRef{Section}{S:G0}{2} and the proof of \lcRef{Corollary}{C:d1-first-calibration}{B.1.7}, ``the $L^2$ expansion permits taking expectations of partial sums''); TRUE for the solution and its chaos terms}
\lcEntryEnd
\lcEntry{\lcRef{Corollary}{C:three-product-riesz}{C.2.6}}{Ordinary product and Riesz covariances in dimension three}
\lcPart{Lean proves}{With the corollary's hypotheses, that the product and Riesz targets are nonzero covariances with scaling $3-\alpha$ whose endpoints satisfy \lcRef{Hypothesis}{H:three-endpoints}{C.2.2}, and, from the supplied \lcRef{Proposition}{P:three-transfer}{C.2.5}, the coefficient, moment and asymptotic bounds for an abstract model of $A_k(t)$, $\E u^2$ and $\mathcal{M}_{\rm st}$; the asymptotic bound is for every spectral measure of the low endpoint.}
\lcInputs{Input supplied}
\lcInput{1}{\texttt{ThreeTransfer}: the exact statement of \lcRef{Proposition}{P:three-transfer}{C.2.5} (proved in the manuscript) for an abstract \texttt{ThreeMomentModel} at the midpoint $\vartheta$; true for the manuscript's $A_k(t)$, $\E u^2(t,x)$ and $\mathcal{M}_{\rm st}$, which the package does not construct}
\lcEntryEnd
\lcEntry{\lcRef{Corollary}{C:three-finite-images}{C.2.7}}{Finite sums of linear images}
\lcPart{Lean proves}{With the corollary's hypotheses, that $\sum_iw_i(T_i)_\#\Gamma_i$ is a nonzero covariance with scaling $3-\alpha$ whose endpoints satisfy \lcRef{Hypothesis}{H:three-endpoints}{C.2.2}, and, from the supplied \lcRef{Proposition}{P:three-transfer}{C.2.5}, the coefficient, moment and asymptotic bounds for an abstract model of $A_k(t)$, $\E u^2$ and $\mathcal{M}_{\rm st}$; the asymptotic bound is for every spectral measure of the low endpoint.}
\lcInputs{Input supplied}
\lcInput{1}{\texttt{ThreeTransfer}: the exact statement of \lcRef{Proposition}{P:three-transfer}{C.2.5} (proved in the manuscript) for an abstract \texttt{ThreeMomentModel} at the midpoint $\vartheta$; true for the manuscript's $A_k(t)$, $\E u^2(t,x)$ and $\mathcal{M}_{\rm st}$, which the package does not construct}
\lcEntryEnd
\lcEntry{\lcRef{Proposition}{P:moment-growth-logconvex}{D.1.7}}{Log-convexity of the upper growth type}
\lcPart{Lean proves}{Every conclusion, for an abstract moment function with the supplied coefficient, interpolation, static-law, minorant and upper-bound facts.}
\lcPart{Not yet formalized}{The moments $\E u^{(a)}(T,x)^p$ of the solution are not constructed, so their identification with that function is not formalized.}
\lcInputs{Inputs supplied}
\lcInput{1}{\lcRef{Equation}{E:static-ratio-input}{D.1.1} (Chen--Hu static law): $0<L_0$ and $T^{-\kappa_{\rm st}}\log\E U(T,x)^p\to p^{\kappa_{\rm st}}L_0$ for $p=1,2$ (fields \texttt{L0\_pos}, \texttt{static\_law})}
\lcInput{2}{\lcRef{Equation}{E:full-moment-positive-series}{D.1.11} of \lcRef{Lemma}{L:full-moment-coefficient-holder}{D.1.6} with $c_{p,0}(a)=1$ (\lcRef{Equation}{E:full-moment-coefficients}{D.1.10}), nonnegative coefficients, summable for $t>0$, $0\le a<a_{\max}$ (field \texttt{series})}
\lcInput{3}{\lcRef{Equation}{E:full-moment-coefficient-holder}{D.1.12} of \lcRef{Lemma}{L:full-moment-coefficient-holder}{D.1.6} (field \texttt{holder})}
\lcInput{4}{\lcRef{Equation}{E:static-minorant-rate}{D.1.6} of \lcRef{Proposition}{P:static-minorant-lower}{D.1.3}, as: every $y<p^{\kappa_{\rm st}}L_0$ is eventually below $T^{-\kappa(a)}\log\E u^{(a)}(T,x)^p$, $0<a<a_{\max}$ (field \texttt{minorant})}
\lcInput{5}{\lcRef{Equation}{E:full-upper-bound}{3.3.28} of \lcRef{Theorem}{T:full-upper}{3.3.8} at temporal exponent $a$: $\E u^{(a)}(T,x)^2\le e(T+2)^d\exp(CT^{\kappa(a)})$ for $T\ge1$, one finite $C$ (field \texttt{full\_upper})}
\lcInput{6}{Cauchy--Schwarz $(\E u^{(a)}(T,x))^2\le\E u^{(a)}(T,x)^2$ (field \texttt{cauchy\_schwarz}; in the proof of \lcRef{Proposition}{P:moment-growth-logconvex}{D.1.7})}
\lcEntryEnd
\lcEntry{\lcRef{Corollary}{C:static-moment-constant-continuity}{D.1.8}}{Recovery of the static moment constant}
\lcPart{Lean proves}{Every conclusion, for an abstract moment function with the supplied facts of the logconvexity proposition and the static minorant and law.}
\lcPart{Not yet formalized}{The moments of the solution are not constructed, so their identification with that function is not formalized.}
\lcInputs{Inputs supplied}
\lcInput{1}{\lcRef{Equation}{E:static-minorant-rate}{D.1.6} of \lcRef{Proposition}{P:static-minorant-lower}{D.1.3} (field \texttt{minorant}) and $\overline{\Lambda}_p(0)=C_p$ from \lcRef{Equation}{E:static-ratio-input}{D.1.1} (fields \texttt{L0\_pos}, \texttt{static\_law})}
\lcInput{2}{the inputs of \lcRef{Proposition}{P:moment-growth-logconvex}{D.1.7}: \lcRef{Equation}{E:full-moment-positive-series}{D.1.11} and \lcRef{Equation}{E:full-moment-coefficient-holder}{D.1.12} of \lcRef{Lemma}{L:full-moment-coefficient-holder}{D.1.6} (fields \texttt{series}, \texttt{holder}), \lcRef{Equation}{E:full-upper-bound}{3.3.28} of \lcRef{Theorem}{T:full-upper}{3.3.8} (field \texttt{full\_upper}), Cauchy--Schwarz (field \texttt{cauchy\_schwarz})}
\lcEntryEnd

\lcCategory{Not formalized.}{2}
\begin{lcTable}
  \lcRef{Theorem}{T:full-upper}{3.3.8} & Full upper scale
\end{lcTable}
\begin{lcTable}
  \lcRef{Theorem}{T:intrinsic-second-moment}{11.2.6} & Intrinsic original second-moment constant
\end{lcTable}

}
\newcommand{\lcLeanPkgVersion}{0.29.0}
\newcommand{\lcClaimsChecked}{4502}
\newcommand{\lcComparatorTheorems}{3653}
\newcommand{\lcHaveReceipts}{1}
\newcommand{\lcVerified}{2261}
\newcommand{\lcOpen}{0}
\newcommand{\lcHaveLedger}{1}
\newcommand{\lcCompletionTable}{%
\begin{tabular}{@{}lrr@{}}
  \toprule
  Status & Results & Share \\
  \midrule
  Formalized with the paper's own hypotheses & 44 & 20\% \\
  With the paper's hypotheses; narrower conclusion & 16 & 7\% \\
  Every conclusion proved; inputs supplied & 1 & 0\% \\
  Steps checked; inputs supplied & 34 & 15\% \\
  Main theorems outside Lean & 2 & 1\% \\
  Not formalized & 123 & 56\% \\
  \midrule
  Total & 220 & 100\% \\
  \bottomrule
\end{tabular}
}
\newcommand{\lcMilestoneTable}{%
\begin{tabular}{@{}llrrr@{}}
  \toprule
  Stage & Date & Claims & Theorems & Results covered \\
  \midrule
  \texttt{L0-companion-created} & November~9 & 88 & 82 & 25 \\
  \texttt{ledger-gate} & November~10 & 88 & 82 & 28 \\
  \texttt{L1-gate} & November~10 & 117 & 104 & 28 \\
  \texttt{L2-gate} & November~10 & 223 & 188 & 29 \\
  \texttt{L3-gate} & November~10 & 557 & 470 & 33 \\
  \texttt{L4-gate} & November~11 & 724 & 611 & 40 \\
  \texttt{reanchor-975feb3} & November~11 & 724 & 611 & 40 \\
  \texttt{L5-gate} & November~11 & 972 & 810 & 46 \\
  \texttt{L6-gate} & November~11 & 1463 & 1210 & 59 \\
  \texttt{L7-gate} & November~12 & 1731 & 1430 & 67 \\
  \texttt{L8-gate} & November~12 & 2108 & 1747 & 71 \\
  \texttt{L9-gate} & November~12 & 2468 & 2050 & 80 \\
  \texttt{L10-gate} & November~13 & 2757 & 2268 & 82 \\
  \texttt{L11-gate} & November~13 & 2955 & 2423 & 82 \\
  \texttt{L12-gate} & November~13 & 3107 & 2541 & 85 \\
  \texttt{L13-gate} & November~13 & 4504 & 3654 & 101 \\
  \texttt{L14-gate} & November~13 & 4504 & 3654 & 101 \\
  \bottomrule
\end{tabular}
}
\newcommand{\lcNumberedTotal}{220}
\newcommand{\lcNotFormalized}{123}
\newcommand{\lcNotFormalizedShare}{56}
\newcommand{\lcOutsideNumbered}{two}
\newcommand{\lcSameFullNumbered}{44}
\newcommand{\lcSameFullAll}{48}
\newcommand{\lcOtherRows}{4}
\newcommand{\lcFirstStageDate}{November~9}
\newcommand{\lcRecordTotal}{220}
\newcommand{\lcAddedSinceRecord}{0}
\newcommand{\lcHaveCompletion}{1}
\newcommand{\lcModelKey}{%
\item[M1] In dimension three the wave Green measure is defined through Mathlib's sphere measure; that it equals the surface measure used in the text is not proved in Lean.
\item[M2] The Gaussian space over finitely many coordinates is represented by normalized Hermite coefficient vectors; that Bochner integrals in it agree with the coordinatewise integrals is not proved in Lean.
\item[M3] The smoothed covariance $p_\varepsilon*\Gamma$ is taken as the measure with density $\gamma_\varepsilon$; that it equals Mathlib's convolution of measures is not proved in Lean.
\item[M4] The two-chain diagram integrals are the deterministic integrals of a diagram model; their identification with the joint limit in the text at arbitrary vertexwise widths is not proved in Lean (the positions form is proved in dimension one).
\item[M5] The completion spaces $H_t\otimes G$, $H_t$, $H_I$ and $H_\Gamma$ are realized as $L^2$ completions under the kernel and covariance forms; no tensor-product identification and no isometry between $H_{\mathrm{sp}}$ and $H_\Gamma$ is proved in Lean.
\item[M6] $H^1$ is the graph of the weak gradient in $L^2\times L^2$, with the topology of Mathlib's $H^1$; that every real element of Mathlib's $H^1$ arises from it is not proved in Lean.
}

\newcommand{\vsDate}{November~16, 2026}

\newcommand{\vsRows}{253}
\newcommand{\vsLeanAny}{95}
\newcommand{\vsOutside}{two}
\newcommand{\vsHumanChecked}{0}

\newcommand{\vsTable}{%
Numbered results & 31 & 72 & 82 & 35 & 33 & 253 \\
\midrule
\multicolumn{7}{@{}l}{\emph{Lean}} \\
\quad Formalized in full & 6 & 16 & 11 & 11 & 0 & 44 \\
\quad Formalized, narrower conclusion & 2 & 4 & 10 & 0 & 0 & 16 \\
\quad Every conclusion, inputs supplied & 0 & 0 & 1 & 0 & 0 & 1 \\
\quad Steps checked; inputs supplied or left outside Lean & 8 & 11 & 7 & 8 & 0 & 34 \\
\quad Outside Lean & 1 & 0 & 1 & 0 & 0 & 2 \\
\quad No Lean claim & 14 & 41 & 52 & 16 & 33 & 156 \\
\midrule
\multicolumn{7}{@{}l}{\emph{Checked by a person}} \\
\quad Checked & 0 & 0 & 0 & 0 & 0 & 0 \\
\quad Disputed & 0 & 0 & 0 & 0 & 0 & 0 \\
\quad Not yet checked & 31 & 72 & 82 & 35 & 33 & 253 \\
}

\newcommand{\cpResults}{168}
\newcommand{\cpProved}{180}
\newcommand{\cpFirstChapter}{2}
\newcommand{\cpLastChapter}{11}
\newcommand{\cpNotReached}{12}
\newcommand{\cpSecondMoment}{153}
\newcommand{\cpLeanFull}{29}
\newcommand{\cpLeanNarrower}{15}
\newcommand{\cpLeanSupplied}{one}
\newcommand{\cpLeanPartial}{23}
\newcommand{\cpLeanNone}{100}
\newcommand{\cpLinksAdded}{28}
\newcommand{\cpLinksRemoved}{2}
\newcommand{\cpReadByHand}{52}
\newcommand{\cpTheoremRows}{%
Theorem~\ref{T:full-upper} & Full upper scale & 20 & 6 & 2 & 0 & 5 & 7 \\
Theorem~\ref{T:intrinsic-second-moment} & Intrinsic original second-moment constant & 153 & 25 & 8 & 1 & 21 & 98 \\
Theorem~\ref{T:intrinsic-occupation-ceiling} & The computed occupation ceiling & 157 & 26 & 10 & 1 & 21 & 99 \\
Theorem~\ref{T:white-moment-bounds} & Explicit white moment bounds & 167 & 29 & 15 & 1 & 22 & 100 \\
}
\newcommand{\cpChapterRows}{%
2 & 11 & 9 & 1 & 0 & 0 & 1 & 7 \\
3 & 15 & 12 & 5 & 2 & 0 & 4 & 1 \\
4 & 24 & 22 & 12 & 0 & 0 & 3 & 7 \\
5 & 16 & 16 & 2 & 0 & 0 & 2 & 12 \\
6 & 18 & 18 & 0 & 0 & 0 & 0 & 18 \\
7 & 14 & 14 & 2 & 3 & 0 & 6 & 3 \\
8 & 19 & 17 & 5 & 3 & 0 & 0 & 9 \\
9 & 14 & 14 & 0 & 0 & 0 & 2 & 12 \\
10 & 29 & 27 & 0 & 1 & 1 & 1 & 24 \\
11 & 20 & 19 & 2 & 6 & 0 & 4 & 7 \\
\midrule
All & 180 & 168 & 29 & 15 & 1 & 23 & 100 \\
}
\newcommand{\cpRows}{%
\cpRow{Proposition~\ref{P:mixed}}{Mixed terms}{\ref{L:pairing-product-envelope}, \ref{L:spatial-wick-measures}, \ref{L:joint-fixed-order}, \ref{L:crude-corrected}, \ref{L:vertex-crude}, \ref{P:mild-construction}}{40}{---}
\cpRow{Lemma~\ref{L:damped-homogeneous-norms}}{Homogeneous damped norms and temporal tails}{(\ref{L:translated-covariance-CS})}{50}{full}
\cpRow{Lemma~\ref{L:translated-covariance-CS}}{Translated covariance Cauchy--Schwarz}{---}{50}{---}
\cpRow{Lemma~\ref{L:pairing-product-envelope}}{Pairing product envelope}{\ref{L:translated-covariance-CS}}{49}{---}
\cpRow{Lemma~\ref{L:spatial-wick-measures}}{Positive spatial Wick-measure limit}{---}{49}{---}
\cpRow{Lemma~\ref{L:joint-fixed-order}}{Joint diagram limits and fixed-order construction}{\ref{L:damped-homogeneous-norms}, \ref{L:pairing-product-envelope}, \ref{L:spatial-wick-measures}}{48}{---}
\cpRow{Lemma~\ref{L:crude-corrected}}{Corrected crude estimate}{\ref{L:damped-homogeneous-norms}, \ref{L:pairing-product-envelope}, \ref{L:joint-fixed-order}}{47}{---}
\cpRow{Lemma~\ref{L:vertex-crude}}{Uniform vertexwise regularization bound}{\ref{L:damped-homogeneous-norms}, \ref{L:pairing-product-envelope}, \ref{L:joint-fixed-order}, \ref{L:crude-corrected} (\ref{P:mild-construction})}{46}{---}
\cpRow{Proposition~\ref{P:mild-construction}}{Mild construction from the fixed-order limits}{\ref{L:joint-fixed-order}, \ref{L:vertex-crude}}{41}{partial}
\cpRow{Lemma~\ref{L:phase-contractions}}{Two static contraction bounds}{---}{12}{full}
\cpRow{Lemma~\ref{L:complex-phase}}{A complex phase bound}{---}{11}{partial}
\cpRow{Lemma~\ref{L:fock-form}}{Gaussian cutoff form bound}{\ref{L:phase-contractions}}{11}{narrower}
\cpRow{Lemma~\ref{L:finite-word-cutoff}}{Exact cutoffs for finite Gaussian words}{(\ref{L:fock-form})}{10}{full}
\cpRow{Lemma~\ref{L:interval-contractions}}{Dimension-independent Gaussian contractions}{\ref{L:phase-contractions}, \ref{L:fock-form}}{10}{narrower}
\cpRow{Lemma~\ref{L:interval-complex}}{Complex sums and changing phase weights}{\ref{L:complex-phase}}{10}{partial}
\cpRow{Lemma~\ref{L:inventory-product}}{A coefficient bound for a product with nonnegative coefficients}{(\ref{P:interval-coefficients})}{9}{full}
\cpRow{Lemma~\ref{L:fixed-split}}{A fixed degree-split estimate}{\ref{L:complex-phase}, \ref{L:finite-word-cutoff}, \ref{L:interval-contractions}, \ref{L:interval-complex}}{9}{partial}
\cpRow{Proposition~\ref{P:interval-coefficients}}{Uniform interval coefficients}{\ref{L:finite-word-cutoff}, \ref{L:inventory-product}, \ref{L:fixed-split}}{8}{partial}
\cpRow{Lemma~\ref{L:terminal-wick}}{The finite terminal-chain Wick formula}{---}{9}{full}
\cpRow{Lemma~\ref{L:uniform-series-tail}}{Uniform series tails from scaled moments}{\ref{L:terminal-wick} (\ref{P:interval-coefficients})}{8}{full}
\cpRow{Theorem~\ref{T:full-upper}}{Full upper scale}{\ref{P:mixed}, \ref{L:joint-fixed-order}, \ref{L:crude-corrected}, \ref{P:mild-construction}, \ref{L:fixed-split}, \ref{P:interval-coefficients}, \ref{L:terminal-wick}, \ref{L:uniform-series-tail}}{7}{outside}
\cpRow{Proposition~\ref{P:masked-spectral-growth}}{Translated spectral mass}{---}{77}{full}
\cpRow{Proposition~\ref{P:masked-auxiliary-norm}}{Auxiliary mask norm}{\ref{P:masked-spectral-growth}}{76}{full}
\cpRow{Lemma~\ref{L:masked-time-completion}}{Temporal covariance completion}{---}{81}{full}
\cpRow{Proposition~\ref{P:masked-noise-completion}}{Completion on the same noise}{\ref{P:masked-auxiliary-norm}, \ref{L:masked-time-completion}}{75}{partial}
\cpRow{Lemma~\ref{L:masked-temporal-modulation}}{Modulated amplitudes}{\ref{L:masked-time-completion}}{76}{full}
\cpRow{Lemma~\ref{L:masked-trig-products}}{Correlated trigonometric products}{\ref{L:masked-temporal-modulation}}{75}{full}
\cpRow{Lemma~\ref{L:masked-continuous-wave}}{Continuous sine factors}{\ref{L:masked-time-completion}, \ref{L:masked-trig-products}}{74}{full}
\cpRow{Proposition~\ref{P:frequency-row}}{Uniform scalar row}{(\ref{L:masked-trig-products})}{69}{full}
\cpRow{Corollary~\ref{C:frequency-time}}{A time exponent above one half}{\ref{P:frequency-row}}{68}{full}
\cpRow{Proposition~\ref{P:actual-wave-form}}{Form and phase kernel}{\ref{P:masked-auxiliary-norm}, \ref{L:masked-time-completion}, \ref{P:masked-noise-completion}, \ref{L:masked-continuous-wave}}{73}{partial}
\cpRow{Lemma~\ref{L:phase-bessel}}{Hilbert-valued finite Bessel inequality}{(\ref{P:actual-wave-form}, \ref{P:wave-two-variances})}{68}{full}
\cpRow{Proposition~\ref{P:wave-two-variances}}{Both variance directions}{\ref{P:masked-auxiliary-norm}, \ref{L:masked-time-completion}, \ref{L:masked-continuous-wave}, \ref{P:frequency-row}, \ref{C:frequency-time}, \ref{P:actual-wave-form}, \ref{L:phase-bessel}}{67}{partial}
\cpRow{Proposition~\ref{P:wave-square-columns}}{Column and adjoint-column maps}{\ref{P:wave-two-variances}}{66}{---}
\cpRow{Lemma~\ref{L:gaussian-scalar-powers}}{Scalar power estimates}{---}{69}{full}
\cpRow{Lemma~\ref{L:gaussian-trace-bounds}}{Two Hermitian trace bounds}{\ref{L:gaussian-scalar-powers}}{68}{full}
\cpRow{Proposition~\ref{P:finite-gaussian-matrix}}{Finite Gaussian matrix moment}{\ref{L:gaussian-scalar-powers}, \ref{L:gaussian-trace-bounds}}{67}{full}
\cpRow{Lemma~\ref{L:cutoff-factorization}}{Cutoff factorization}{---}{77}{---}
\cpRow{Proposition~\ref{P:cutoff-singular-values}}{Whole-space cutoff singular values}{\ref{L:cutoff-factorization}}{72}{---}
\cpRow{Proposition~\ref{P:compact-two-envelopes}}{Both compact envelopes}{\ref{L:cutoff-factorization}, \ref{P:cutoff-singular-values}}{70}{---}
\cpRow{Lemma~\ref{L:compact-positive-order}}{Compact positive order and trace}{---}{69}{---}
\cpRow{Proposition~\ref{P:localized-finite-variance}}{Localized finite variance bounds}{\ref{P:compact-two-envelopes}, \ref{L:compact-positive-order}}{67}{---}
\cpRow{Corollary~\ref{C:finite-gaussian-operators}}{Finite Gaussian operators}{\ref{P:finite-gaussian-matrix}, \ref{L:compact-positive-order}, \ref{P:localized-finite-variance}}{66}{---}
\cpRow{Proposition~\ref{P:actual-operator-driver}}{Compact operator path on the original noise space}{\ref{L:driver-noise-basis}, \ref{L:driver-interval-limit}, \ref{L:driver-holder-limit}}{8}{---}
\cpRow{Lemma~\ref{L:driver-noise-basis}}{Basis in the prescribed quotient}{\ref{L:masked-time-completion}}{72}{---}
\cpRow{Lemma~\ref{L:driver-covariance-coefficients}}{Covariance coefficients}{\ref{P:actual-wave-form}, \ref{L:driver-noise-basis}}{70}{---}
\cpRow{Lemma~\ref{L:driver-gain-localization}}{Deterministic coefficient operators}{\ref{P:actual-wave-form}, \ref{L:driver-covariance-coefficients}}{69}{---}
\cpRow{Lemma~\ref{L:driver-actual-variances}}{Both variances and their tails}{\ref{P:wave-two-variances}, \ref{P:wave-square-columns}, \ref{P:compact-two-envelopes}, \ref{L:driver-noise-basis}, \ref{L:driver-gain-localization}}{65}{---}
\cpRow{Lemma~\ref{L:driver-dominated-tails}}{Compact domination of positive tails}{\ref{L:compact-positive-order} (\ref{L:driver-actual-variances})}{65}{---}
\cpRow{Lemma~\ref{L:driver-gaussian-tail-average}}{Retained finite Gaussian variance powers}{\ref{P:finite-gaussian-matrix}, \ref{L:compact-positive-order}, \ref{C:finite-gaussian-operators}}{65}{---}
\cpRow{Lemma~\ref{L:driver-interval-limit}}{Fixed-interval operator construction}{\ref{P:localized-finite-variance}, \ref{C:finite-gaussian-operators}, \ref{L:driver-noise-basis}, \ref{L:driver-actual-variances}, \ref{L:driver-dominated-tails}, \ref{L:driver-gaussian-tail-average}}{64}{partial}
\cpRow{Lemma~\ref{L:driver-holder-limit}}{H\"older path limit}{\ref{C:frequency-time}, \ref{P:wave-two-variances}, \ref{C:finite-gaussian-operators}, \ref{L:driver-gain-localization}, \ref{L:driver-actual-variances}, \ref{L:driver-interval-limit} (\ref{P:actual-operator-driver})}{63}{partial}
\cpRow{Lemma~\ref{L:source-finite-comparison}}{Finite Gaussian linear images}{---}{61}{full}
\cpRow{Proposition~\ref{P:source-gaussian-transfer}}{Gaussian transfer and stability}{\ref{L:driver-interval-limit}, \ref{L:driver-holder-limit}, \ref{L:source-finite-comparison}}{60}{---}
\cpRow{Lemma~\ref{L:source-spatial-heat}}{Original spatial space and heat}{\ref{P:masked-noise-completion} (\ref{P:source-gaussian-transfer})}{60}{full}
\cpRow{Lemma~\ref{L:source-retaining-isometry}}{Feature isometry preserving the given noise}{\ref{L:driver-noise-basis}}{60}{---}
\cpRow{Proposition~\ref{P:source-canonical-convergence}}{Canonical sources on the extended space containing the given noise}{\ref{P:source-gaussian-transfer}, \ref{L:source-spatial-heat}, \ref{L:source-retaining-isometry}}{59}{---}
\cpRow{Lemma~\ref{L:source-pointfield-regularity}}{Regularity for each fixed cutoff}{\ref{P:source-canonical-convergence}}{42}{---}
\cpRow{Proposition~\ref{P:source-ordinary-response}}{Ordinary source response}{\ref{P:actual-wave-form}, \ref{L:driver-gain-localization}, \ref{P:source-canonical-convergence}}{47}{---}
\cpRow{Lemma~\ref{L:masked-chain-comparison}}{Chains and mixed-order comparison}{\ref{L:joint-fixed-order}, \ref{L:vertex-crude}, \ref{L:source-pointfield-regularity} (\ref{P:source-canonical-convergence})}{41}{---}
\cpRow{Proposition~\ref{P:masked-source-moments}}{Ordinary solutions and tested moments}{\ref{L:vertex-crude}, \ref{L:masked-chain-comparison}}{40}{---}
\cpRow{Proposition~\ref{P:masked-canonical-transfer}}{Canonical point-moment limit and conditional transfer}{\ref{P:mixed}, \ref{L:joint-fixed-order}, \ref{P:masked-source-moments}}{38}{---}
\cpRow{Lemma~\ref{L:constant-phase-free}}{Constants and finite propagation}{\ref{L:driver-gain-localization}}{47}{---}
\cpRow{Proposition~\ref{P:constant-phase-driver}}{Extended operator path}{\ref{L:driver-gain-localization}, \ref{P:source-gaussian-transfer}, \ref{P:source-canonical-convergence}}{48}{---}
\cpRow{Lemma~\ref{L:constant-source-response}}{Ordinary response on all of $F$}{\ref{P:source-ordinary-response}, \ref{L:constant-phase-free}}{44}{---}
\cpRow{Proposition~\ref{P:constant-young-linear}}{Linear Young evolution}{---}{46}{---}
\cpRow{Proposition~\ref{P:constant-gaussian-flow}}{Gaussian Young and physical flows}{\ref{P:constant-phase-driver}, \ref{P:constant-young-linear}}{38}{---}
\cpRow{Lemma~\ref{L:constant-source-solution}}{Identification with the ordinary source equation}{\ref{P:masked-source-moments}, \ref{L:constant-source-response}}{38}{---}
\cpRow{Proposition~\ref{P:constant-raw-average}}{Original average and original point moment}{\ref{P:masked-source-moments}, \ref{P:masked-canonical-transfer}, \ref{P:constant-gaussian-flow}, \ref{L:constant-source-solution}}{35}{---}
\cpRow{Proposition~\ref{P:constant-compact-base}}{Compact deterministic base}{\ref{L:constant-source-response}, \ref{P:constant-raw-average}}{34}{---}
\cpRow{Proposition~\ref{P:compact-ground-family}}{Ground states and a uniform rate gap}{---}{32}{---}
\cpRow{Lemma~\ref{L:compact-fractional-phase}}{Uniform fractional phase norms}{(\ref{P:compact-ground-family})}{29}{---}
\cpRow{Proposition~\ref{P:compact-E-wave}}{The complementary wave on $E$}{\ref{P:compact-ground-family}, \ref{L:compact-fractional-phase}}{28}{---}
\cpRow{Proposition~\ref{P:compact-F-projection}}{The growing projection on the wave phase space}{\ref{P:constant-raw-average}, \ref{P:constant-compact-base}, \ref{P:compact-E-wave}}{27}{---}
\cpRow{Proposition~\ref{P:compact-F-complement}}{Complementary wave evolution and parameter dependence}{\ref{P:compact-E-wave}, \ref{P:compact-F-projection} (\ref{L:raw-block-indexed-cone})}{15}{---}
\cpRow{Lemma~\ref{L:compact-raw-profile}}{The original mean and variance}{---}{45}{---}
\cpRow{Proposition~\ref{P:compact-original-common-trial}}{Binding with the original control and cost fixed}{\ref{P:compact-ground-family}, \ref{L:compact-raw-profile}}{13}{---}
\cpRow{Lemma~\ref{L:exact-canonical-scaling}}{Canonical law scaling}{\ref{L:joint-fixed-order}, \ref{L:vertex-crude} (\ref{P:raw-block-point-comparison})}{15}{partial}
\cpRow{Lemma~\ref{L:raw-block-inputs}}{Coupled Gaussian block inputs}{\ref{L:driver-gain-localization}, \ref{P:constant-phase-driver}}{39}{partial}
\cpRow{Proposition~\ref{P:raw-block-gaussian-law}}{Change of law on the whole norm-Borel space}{\ref{L:raw-block-inputs}}{15}{---}
\cpRow{Lemma~\ref{L:raw-block-uniform-error}}{Uniform block error}{\ref{P:constant-young-linear}, \ref{P:constant-compact-base}, \ref{P:compact-F-projection}, \ref{L:raw-block-inputs}}{14}{---}
\cpRow{Lemma~\ref{L:raw-block-indexed-cone}}{Indexed cone and terminal displacement}{\ref{P:constant-raw-average}, \ref{P:compact-F-projection}, \ref{P:compact-F-complement} (\ref{L:raw-block-uniform-error})}{14}{partial}
\cpRow{Lemma~\ref{L:raw-block-restart}}{Composition over time blocks and regularized products}{\ref{P:source-gaussian-transfer}, \ref{L:source-retaining-isometry}, \ref{P:source-canonical-convergence}, \ref{L:constant-source-response}, \ref{P:constant-gaussian-flow}, \ref{P:constant-compact-base}, \ref{L:raw-block-inputs} (\ref{P:raw-block-point-comparison})}{22}{partial}
\cpRow{Proposition~\ref{P:raw-block-point-comparison}}{The square observable and canonical comparison}{\ref{P:masked-source-moments}, \ref{P:masked-canonical-transfer}, \ref{P:constant-gaussian-flow}, \ref{L:constant-source-solution}, \ref{P:constant-raw-average}, \ref{L:exact-canonical-scaling}, \ref{P:raw-block-gaussian-law}, \ref{L:raw-block-restart}}{14}{---}
\cpRow{Proposition~\ref{P:raw-block-original-lower}}{Original point lower rate at the Gaussian quadratic cost}{\ref{P:constant-young-linear}, \ref{L:raw-block-inputs}, \ref{L:raw-block-uniform-error}, \ref{L:raw-block-indexed-cone}, \ref{P:raw-block-point-comparison}}{13}{partial}
\cpRow{Lemma~\ref{L:variational-heat-envelope}}{Gaussian envelope for the original covariance}{(\ref{L:variational-global-mean})}{14}{full}
\cpRow{Lemma~\ref{L:variational-global-mean}}{Global regularity of the original mean}{\ref{L:variational-heat-envelope}}{13}{full}
\cpRow{Proposition~\ref{P:variational-unmasked-rate}}{The unmasked Rayleigh rate}{\ref{P:compact-original-common-trial}, \ref{L:variational-global-mean}}{12}{narrower}
\cpRow{Proposition~\ref{P:variational-remove-mask}}{Monotone rates and the original point moment}{\ref{P:compact-ground-family}, \ref{P:compact-original-common-trial}, \ref{P:raw-block-original-lower}, \ref{P:variational-unmasked-rate}}{10}{partial}
\cpRow{Lemma~\ref{L:variational-control-dilation}}{Radon and Rayleigh dilation}{(\ref{P:variational-remove-mask}, \ref{P:variational-positive-scores})}{7}{narrower}
\cpRow{Proposition~\ref{P:variational-positive-scores}}{Intrinsic variational lower bounds}{\ref{T:full-upper}, \ref{L:compact-raw-profile}, \ref{P:variational-remove-mask}, \ref{L:variational-control-dilation}}{6}{narrower}
\cpRow{Lemma~\ref{L:equilibrium-gamma-beta}}{Gamma and Beta identities}{(\ref{L:temporal-equilibrium})}{6}{full}
\cpRow{Lemma~\ref{L:temporal-equilibrium}}{Exact potential and mass}{\ref{L:equilibrium-gamma-beta}}{5}{full}
\cpRow{Lemma~\ref{L:spatial-physical-pair}}{The physical two-point integral}{(\ref{P:variational-positive-scores})}{39}{full}
\cpRow{Lemma~\ref{L:spatial-static-controls}}{Smooth temporal controls}{\ref{P:variational-unmasked-rate}, \ref{P:variational-remove-mask}, \ref{L:temporal-equilibrium}}{4}{narrower}
\cpRow{Lemma~\ref{L:spatial-fourier-bound}}{Deterministic Fourier bound for the spatial energy}{\ref{L:damped-homogeneous-norms}}{4}{full}
\cpRow{Proposition~\ref{P:spatial-smooth-ratio}}{The finite smooth spatial ratio}{\ref{L:spatial-static-controls}, \ref{L:spatial-fourier-bound}}{3}{narrower}
\cpRow{Proposition~\ref{P:spatial-four-linear}}{The physical four-linear extension}{\ref{L:compact-raw-profile}, \ref{L:spatial-physical-pair} (\ref{P:occupation-physical-class})}{38}{full}
\cpRow{Proposition~\ref{P:spatial-energy-map}}{The spatial energy and static lower variational expression}{\ref{P:spatial-smooth-ratio}, \ref{P:spatial-four-linear}}{2}{narrower}
\cpRow{Lemma~\ref{L:cameron-fractional-products}}{Fractional covariance products}{\ref{L:spatial-physical-pair}, \ref{P:spatial-four-linear}}{32}{---}
\cpRow{Proposition~\ref{P:cameron-spatial-map}}{Compact spatial mean map}{\ref{P:cutoff-singular-values}, \ref{L:cameron-fractional-products}}{26}{---}
\cpRow{Lemma~\ref{L:cameron-supported-tensor}}{Supported original tensor space}{(\ref{P:cameron-spatial-map})}{30}{---}
\cpRow{Proposition~\ref{P:cameron-temporal-map}}{The compact temporal mean, without trace class}{(\ref{P:cameron-full-mean})}{29}{---}
\cpRow{Proposition~\ref{P:cameron-full-mean}}{The compact time-space mean}{\ref{L:constant-source-response}, \ref{P:cameron-spatial-map}, \ref{L:cameron-supported-tensor}, \ref{P:cameron-temporal-map}}{25}{---}
\cpRow{Lemma~\ref{L:cameron-actual-residual}}{Supported residual paths}{\ref{P:source-gaussian-transfer}, \ref{P:constant-phase-driver}, \ref{L:raw-block-inputs}}{25}{---}
\cpRow{Proposition~\ref{P:cameron-actual-embedding}}{Cameron covariance and averaged bound}{\ref{P:source-gaussian-transfer}, \ref{L:raw-block-inputs}, \ref{P:cameron-spatial-map}, \ref{P:cameron-full-mean}, \ref{L:cameron-actual-residual}}{24}{---}
\cpRow{Lemma~\ref{L:young-upper-growth}}{Refined growth and uniform inverse control}{\ref{P:constant-young-linear}}{17}{---}
\cpRow{Proposition~\ref{P:young-upper-multiplicative}}{Multiplicative upper perturbation}{\ref{L:young-upper-growth}}{10}{---}
\cpRow{Lemma~\ref{L:uniform-mean-cutoff}}{Uniform cutoff and physical mean}{\ref{L:cameron-fractional-products}}{24}{---}
\cpRow{Lemma~\ref{L:uniform-mean-spatial-tails}}{Uniform heat and frequency tails}{\ref{L:uniform-mean-cutoff} (\ref{P:uniform-mean-regulator})}{23}{---}
\cpRow{Lemma~\ref{L:uniform-mean-tensor}}{Dimension-free Hilbert tensor extension}{(\ref{P:uniform-mean-regulator})}{27}{---}
\cpRow{Proposition~\ref{P:uniform-mean-regulator}}{Uniform approximation of the completed mean}{\ref{P:cameron-temporal-map}, \ref{L:uniform-mean-cutoff}, \ref{L:uniform-mean-spatial-tails}, \ref{L:uniform-mean-tensor}}{22}{---}
\cpRow{Proposition~\ref{P:uniform-actual-cameron}}{Uniform-radius Cameron bound}{\ref{P:source-gaussian-transfer}, \ref{L:cameron-actual-residual}, \ref{P:cameron-actual-embedding}, \ref{P:uniform-mean-regulator} (\ref{L:raw-block-inputs})}{12}{---}
\cpRow{Lemma~\ref{L:radius-outer-maps}}{Uniform rescaled outer maps}{\ref{L:cutoff-factorization}}{13}{---}
\cpRow{Proposition~\ref{P:radius-compact-envelopes}}{Radius-explicit compact envelopes}{\ref{P:cutoff-singular-values}, \ref{P:compact-two-envelopes}, \ref{L:radius-outer-maps}}{12}{---}
\cpRow{Proposition~\ref{P:radius-actual-moments}}{Radius-explicit path moments}{\ref{P:wave-two-variances}, \ref{C:finite-gaussian-operators}, \ref{L:driver-noise-basis}, \ref{L:driver-gain-localization}, \ref{L:driver-actual-variances}, \ref{L:driver-holder-limit}, \ref{P:radius-compact-envelopes} (\ref{P:actual-operator-driver})}{11}{---}
\cpRow{Corollary~\ref{C:radius-covariance-residuals}}{Deterministic covariance residuals}{\ref{P:source-gaussian-transfer}, \ref{L:raw-block-inputs}, \ref{P:radius-actual-moments} (\ref{P:radius-mean-anchor})}{10}{---}
\cpRow{Proposition~\ref{P:radius-mean-anchor}}{Averaged block mean tends to zero}{\ref{C:radius-covariance-residuals}}{9}{partial}
\cpRow{Lemma~\ref{L:signed-limiting-multiplier}}{The limiting multiplier}{\ref{L:cutoff-factorization}}{19}{---}
\cpRow{Proposition~\ref{P:signed-closed-form}}{Closed form and quarter-power phase}{\ref{L:signed-limiting-multiplier}, \ref{L:signed-spectral-wave}}{18}{partial}
\cpRow{Lemma~\ref{L:signed-spectral-wave}}{The bottom and spectral wave}{\ref{L:constant-phase-free}, \ref{P:signed-closed-form}}{18}{---}
\cpRow{Proposition~\ref{P:signed-physical-wave}}{Identification with the physical wave}{\ref{P:compact-F-projection}, \ref{P:cameron-actual-embedding}, \ref{P:signed-closed-form}, \ref{L:signed-spectral-wave}}{16}{---}
\cpRow{Lemma~\ref{L:symplectic-free}}{Free wave form and compact Young baseline}{\ref{P:constant-young-linear} (\ref{P:source-gaussian-transfer}, \ref{L:source-retaining-isometry}, \ref{L:source-pointfield-regularity}, \ref{P:source-ordinary-response}, \ref{P:constant-phase-driver})}{22}{partial}
\cpRow{Lemma~\ref{L:symplectic-source-mean}}{Covariance means and source primitives}{\ref{P:source-gaussian-transfer}, \ref{P:source-canonical-convergence}, \ref{P:source-ordinary-response}}{9}{---}
\cpRow{Proposition~\ref{P:symplectic-source}}{Source conservation}{\ref{L:symplectic-free}, \ref{L:symplectic-source-mean}}{8}{---}
\cpRow{Proposition~\ref{P:symplectic-block-mean}}{The exact completed block mean}{\ref{P:constant-phase-driver}, \ref{P:cameron-full-mean}, \ref{P:cameron-actual-embedding}, \ref{L:symplectic-free}}{9}{---}
\cpRow{Proposition~\ref{P:tube-wiener-realization}}{First Wiener chaos realization on the support}{---}{15}{---}
\cpRow{Lemma~\ref{L:tube-closed-geometry}}{Compact images and attained tube radii}{\ref{P:tube-wiener-realization}}{13}{---}
\cpRow{Proposition~\ref{P:tube-rayleigh-transfer}}{Comparison with the Rayleigh integral, including jumps and zero law}{\ref{L:tube-closed-geometry}}{12}{supplied}
\cpRow{Proposition~\ref{P:tube-actual-law}}{Eligibility of the E law}{\ref{L:cameron-supported-tensor}, \ref{P:cameron-full-mean}, \ref{L:cameron-actual-residual}, \ref{P:cameron-actual-embedding}, \ref{P:uniform-mean-regulator}, \ref{C:radius-covariance-residuals}, \ref{P:radius-mean-anchor}, \ref{P:symplectic-block-mean}, \ref{P:tube-wiener-realization}}{8}{---}
\cpRow{Lemma~\ref{L:tube-actual-endpoint}}{The endpoint functional}{\ref{L:raw-block-restart}, \ref{L:symplectic-free}}{12}{---}
\cpRow{Proposition~\ref{P:tube-coarse-envelope}}{Coarse tube envelope and operator moment}{\ref{P:constant-young-linear}, \ref{L:young-upper-growth}, \ref{P:uniform-actual-cameron}, \ref{P:tube-rayleigh-transfer}, \ref{L:tube-actual-endpoint}}{11}{---}
\cpRow{Lemma~\ref{L:tube-rayleigh-tail}}{Exact Rayleigh tail bound}{\ref{P:tube-coarse-envelope}}{9}{narrower}
\cpRow{Lemma~\ref{L:uniform-E-lift}}{Unique potential and fixed-radius lift}{---}{17}{---}
\cpRow{Proposition~\ref{P:uniform-E-spectral}}{Spectral rate in E}{\ref{P:constant-young-linear}, \ref{P:signed-physical-wave}, \ref{L:uniform-E-lift}}{15}{---}
\cpRow{Proposition~\ref{P:uniform-rate-averaging}}{Uniform completed-field averaging}{\ref{P:cameron-temporal-map}, \ref{L:uniform-mean-tensor}, \ref{P:uniform-E-spectral}}{13}{---}
\cpRow{Lemma~\ref{L:uniform-product-budgets}}{Explicit physical primitive and growth bounds}{\ref{L:young-upper-growth}, \ref{P:uniform-mean-regulator}, \ref{L:symplectic-free}}{11}{---}
\cpRow{Proposition~\ref{P:uniform-mean-product}}{Finite-mesh physical mean-product bound}{\ref{P:young-upper-multiplicative}, \ref{P:uniform-E-spectral}, \ref{L:uniform-product-budgets}}{9}{---}
\cpRow{Corollary~\ref{C:uniform-polynomial-balls}}{Every polynomial control radius}{\ref{P:uniform-rate-averaging}, \ref{P:uniform-mean-product}}{8}{---}
\cpRow{Proposition~\ref{P:uniform-arbitrary-residual}}{Residual comparison and the integrated rate}{\ref{P:young-upper-multiplicative}, \ref{L:uniform-product-budgets}, \ref{C:uniform-polynomial-balls}}{7}{---}
\cpRow{Lemma~\ref{L:point-actual-finiteness}}{Restriction and fixed-parameter finiteness}{\ref{P:constant-phase-driver}, \ref{L:symplectic-free}, \ref{P:tube-actual-law}, \ref{L:tube-actual-endpoint}, \ref{P:tube-coarse-envelope}, \ref{L:tube-rayleigh-tail}}{7}{---}
\cpRow{Lemma~\ref{L:point-fast-sources}}{The same regularized fast noises preserving the given noise}{\ref{P:actual-operator-driver}, \ref{P:source-gaussian-transfer}, \ref{P:source-ordinary-response}, \ref{P:masked-canonical-transfer}, \ref{P:constant-gaussian-flow}, \ref{L:constant-source-solution}, \ref{P:symplectic-source}}{7}{---}
\cpRow{Lemma~\ref{L:point-compact-source}}{Compact-data source identification}{---}{7}{---}
\cpRow{Proposition~\ref{P:point-source-cone}}{Cone identity and an average of squares}{\ref{P:masked-source-moments}, \ref{L:constant-phase-free}, \ref{L:point-fast-sources}, \ref{L:point-compact-source}}{6}{---}
\cpRow{Proposition~\ref{P:original-point-transfer}}{Exact original point transfer and polynomial cost}{\ref{P:masked-canonical-transfer}, \ref{L:exact-canonical-scaling}, \ref{P:raw-block-point-comparison}, \ref{P:raw-block-original-lower}, \ref{L:symplectic-free}, \ref{L:point-actual-finiteness}, \ref{L:point-fast-sources}, \ref{P:point-source-cone}}{5}{---}
\cpRow{Lemma~\ref{L:signed-score-envelope}}{Signed variational expression and quantitative radial envelope}{\ref{P:tube-coarse-envelope}, \ref{P:uniform-rate-averaging}, \ref{L:uniform-product-budgets} (\ref{P:tube-actual-law}, \ref{L:tube-actual-endpoint})}{7}{---}
\cpRow{Proposition~\ref{P:signed-moment-bulk-tail}}{Finite-delta bulk and exact tail}{\ref{L:tube-rayleigh-tail}, \ref{P:uniform-arbitrary-residual}, \ref{L:signed-score-envelope} (\ref{P:tube-coarse-envelope})}{6}{---}
\cpRow{Proposition~\ref{P:signed-operator-upper}}{Two limits on the same moment}{\ref{L:point-actual-finiteness}, \ref{L:signed-score-envelope}, \ref{P:signed-moment-bulk-tail}}{5}{---}
\cpRow{Proposition~\ref{P:original-signed-point-upper}}{Pointwise upper bound from the signed variational expression}{\ref{P:raw-block-original-lower}, \ref{P:original-point-transfer}, \ref{P:signed-operator-upper}}{4}{---}
\cpRow{Lemma~\ref{L:gram-scalar-algebra}}{Scalar and Hilbert identities}{---}{11}{full}
\cpRow{Lemma~\ref{L:gram-dual-minimum}}{Interior dual minimum}{---}{11}{---}
\cpRow{Proposition~\ref{P:gram-positive-optimizer}}{Positive conic maximizer}{\ref{L:gram-scalar-algebra}, \ref{L:gram-dual-minimum}}{10}{partial}
\cpRow{Lemma~\ref{L:gram-realization}}{Real Gram realization}{\ref{P:gram-positive-optimizer}}{9}{---}
\cpRow{Proposition~\ref{P:gram-entrywise-comparison}}{Entrywise Gram comparison}{\ref{P:gram-positive-optimizer}, \ref{L:gram-realization}}{8}{partial}
\cpRow{Lemma~\ref{L:positive-raw-embedding}}{Inclusion into the covariance space and the smooth mean}{\ref{L:cameron-supported-tensor}, \ref{L:uniform-mean-tensor} (\ref{L:uniform-E-lift})}{10}{narrower}
\cpRow{Lemma~\ref{L:positive-physical-gram}}{Finite covariance vectors and entrywise order}{\ref{P:uniform-mean-regulator}, \ref{L:uniform-E-lift}, \ref{P:gram-entrywise-comparison}, \ref{L:positive-raw-embedding} (\ref{P:variational-positive-scores}, \ref{P:original-signed-point-upper}, \ref{T:intrinsic-second-moment})}{7}{narrower}
\cpRow{Lemma~\ref{L:positive-uniform-optimizer}}{Uniform optimizer norm and positive smoothing}{\ref{P:uniform-rate-averaging}, \ref{P:gram-positive-optimizer}, \ref{L:positive-raw-embedding}, \ref{L:positive-physical-gram}}{6}{---}
\cpRow{Proposition~\ref{P:positive-signed-comparison}}{Signed variational values bounded by the variational value over nonnegative controls}{\ref{P:variational-positive-scores}, \ref{P:uniform-E-spectral}, \ref{P:uniform-rate-averaging}, \ref{L:positive-raw-embedding}, \ref{L:positive-physical-gram}, \ref{L:positive-uniform-optimizer}}{5}{---}
\cpRow{Proposition~\ref{P:positive-score-cutoff}}{Cutoff removal for the optimized variational expression}{\ref{P:variational-unmasked-rate}, \ref{L:positive-raw-embedding}, \ref{P:positive-signed-comparison}}{4}{---}
\cpRow{Theorem~\ref{T:intrinsic-second-moment}}{Intrinsic original second-moment constant}{\ref{P:raw-block-original-lower}, \ref{P:variational-unmasked-rate}, \ref{P:variational-positive-scores}, \ref{P:original-signed-point-upper}, \ref{P:positive-score-cutoff}}{3}{outside}
\cpRow{Lemma~\ref{L:occupation-raw-square}}{Full-profile squares in the original space}{\ref{P:spatial-four-linear}, \ref{L:positive-raw-embedding} (\ref{P:occupation-physical-class})}{4}{narrower}
\cpRow{Proposition~\ref{P:raw-occupation-energy}}{Original occupation energy and dilation}{\ref{P:spatial-four-linear}}{3}{full}
\cpRow{Lemma~\ref{L:occupation-rayleigh-selection}}{Countable Rayleigh near-occupations}{\ref{P:variational-unmasked-rate}, \ref{L:occupation-raw-square}}{3}{narrower}
\cpRow{Theorem~\ref{T:intrinsic-occupation-ceiling}}{The computed occupation ceiling}{\ref{T:intrinsic-second-moment}, \ref{P:raw-occupation-energy}, \ref{L:occupation-rayleigh-selection} (\ref{P:spatial-energy-map})}{2}{narrower}
\cpRow{Corollary~\ref{C:intrinsic-spatial-bracket}}{A common spatial-energy bracket}{\ref{L:temporal-equilibrium}, \ref{P:spatial-energy-map}, \ref{T:intrinsic-second-moment}, \ref{T:intrinsic-occupation-ceiling}}{1}{narrower}
\cpRow{Lemma~\ref{L:white-cdf-square}}{A square identity for the distribution function}{---}{2}{---}
\cpRow{Proposition~\ref{P:white-spatial-energy}}{Exact white spatial energy}{\ref{P:spatial-four-linear}, \ref{L:white-cdf-square} (\ref{P:spatial-energy-map})}{1}{partial}
\cpRow{Theorem~\ref{T:white-moment-bounds}}{Explicit white moment bounds}{\ref{P:spatial-energy-map}, \ref{T:intrinsic-second-moment}, \ref{C:intrinsic-spatial-bracket}, \ref{P:white-spatial-energy}}{0}{partial}
}
\newcommand{\cpNotReachedList}{%
Lemma~\ref{L:white-test},
Proposition~\ref{P:envelope},
Proposition~\ref{P:static-phase-variance},
Proposition~\ref{P:static-phase-all-orders},
Corollary~\ref{C:static-interaction},
Proposition~\ref{P:masked-quotient-obstruction},
Proposition~\ref{P:cutoff-boundary-tests},
Lemma~\ref{L:occupation-kernel-row},
Proposition~\ref{P:occupation-physical-class},
Proposition~\ref{P:symplectic-raw},
Lemma~\ref{L:tube-half-action},
Corollary~\ref{C:white-relative-width}%
}

\newcommand{\cpRowFormat}[5]{\noindent{\footnotesize\setlength{\tabcolsep}{3pt}%
\begin{tabular}{@{}>{\raggedright\arraybackslash}p{0.18\textwidth}>{\raggedright\arraybackslash}p{0.27\textwidth}>{\raggedright\arraybackslash}p{0.31\textwidth}>{\raggedleft\arraybackslash}p{0.07\textwidth}>{\raggedright\arraybackslash}p{0.09\textwidth}@{}}#1
  & #2 & #3 & #4 & #5\end{tabular}}\par}
\newcommand{\cpRow}[5]{\cpRowFormat{#1}{#2}{#3}{#4}{#5}}
\newcommand{\cpHead}{\cpRowFormat{\textbf{Result}}{\textbf{Title}}{\textbf{Rests directly on}}{\textbf{Used by}}{\textbf{Lean}}\vspace{2pt}\hrule\vspace{3pt}}
\makeatletter
\def\lc@equation{Equation}
\newif\iflc@lead
\newcommand{\lcRef}[3]{\def\lc@kind{#1}%
\iflc@lead\global\lc@leadfalse
\@ifundefined{LeanTag@#2}{}{\leavevmode\hypertarget{lean-entry:#2}{}}\fi
\ifx\lc@kind\lc@equation Equation~\eqref{#2}\else #1~\ref{#2}\fi}
\makeatother
\expandafter\def\csname LeanTag@C:d1-first-calibration\endcsname{Lean: partial, with inputs}
\expandafter\def\csname LeanTarget@C:d1-first-calibration\endcsname{lean-entry:C:d1-first-calibration}
\expandafter\def\csname LeanShort@C:d1-first-calibration\endcsname{partial}
\expandafter\def\csname LeanTag@C:frequency-time\endcsname{Lean: full statement}
\expandafter\def\csname LeanTarget@C:frequency-time\endcsname{lean-entry:C:frequency-time}
\expandafter\def\csname LeanShort@C:frequency-time\endcsname{full}
\expandafter\def\csname LeanTag@C:intrinsic-spatial-bracket\endcsname{Lean: narrower statement}
\expandafter\def\csname LeanTarget@C:intrinsic-spatial-bracket\endcsname{lean-entry:C:intrinsic-spatial-bracket}
\expandafter\def\csname LeanShort@C:intrinsic-spatial-bracket\endcsname{narrower}
\expandafter\def\csname LeanTag@C:static-interaction\endcsname{Lean: partial, with inputs}
\expandafter\def\csname LeanTarget@C:static-interaction\endcsname{lean-entry:C:static-interaction}
\expandafter\def\csname LeanShort@C:static-interaction\endcsname{partial}
\expandafter\def\csname LeanTag@C:static-moment-constant-continuity\endcsname{Lean: partial, with inputs}
\expandafter\def\csname LeanTarget@C:static-moment-constant-continuity\endcsname{lean-entry:C:static-moment-constant-continuity}
\expandafter\def\csname LeanShort@C:static-moment-constant-continuity\endcsname{partial}
\expandafter\def\csname LeanTag@C:three-finite-images\endcsname{Lean: partial, with inputs}
\expandafter\def\csname LeanTarget@C:three-finite-images\endcsname{lean-entry:C:three-finite-images}
\expandafter\def\csname LeanShort@C:three-finite-images\endcsname{partial}
\expandafter\def\csname LeanTag@C:three-product-riesz\endcsname{Lean: partial, with inputs}
\expandafter\def\csname LeanTarget@C:three-product-riesz\endcsname{lean-entry:C:three-product-riesz}
\expandafter\def\csname LeanShort@C:three-product-riesz\endcsname{partial}
\expandafter\def\csname LeanTag@C:white-relative-width\endcsname{Lean: full statement}
\expandafter\def\csname LeanTarget@C:white-relative-width\endcsname{lean-entry:C:white-relative-width}
\expandafter\def\csname LeanShort@C:white-relative-width\endcsname{full}
\expandafter\def\csname LeanTag@E:driver-parameters\endcsname{Lean: full statement}
\expandafter\def\csname LeanTarget@E:driver-parameters\endcsname{lean-entry:E:driver-parameters}
\expandafter\def\csname LeanShort@E:driver-parameters\endcsname{full}
\expandafter\def\csname LeanTag@E:kappa\endcsname{Lean: full statement}
\expandafter\def\csname LeanTarget@E:kappa\endcsname{lean-entry:E:kappa}
\expandafter\def\csname LeanShort@E:kappa\endcsname{full}
\expandafter\def\csname LeanTag@E:masked-time-schur\endcsname{Lean: full statement}
\expandafter\def\csname LeanTarget@E:masked-time-schur\endcsname{lean-entry:E:masked-time-schur}
\expandafter\def\csname LeanShort@E:masked-time-schur\endcsname{full}
\expandafter\def\csname LeanTag@L:alpha-one-convolution\endcsname{Lean: full statement}
\expandafter\def\csname LeanTarget@L:alpha-one-convolution\endcsname{lean-entry:L:alpha-one-convolution}
\expandafter\def\csname LeanShort@L:alpha-one-convolution\endcsname{full}
\expandafter\def\csname LeanTag@L:complex-phase\endcsname{Lean: partial, with inputs}
\expandafter\def\csname LeanTarget@L:complex-phase\endcsname{lean-entry:L:complex-phase}
\expandafter\def\csname LeanShort@L:complex-phase\endcsname{partial}
\expandafter\def\csname LeanTag@L:d1-spatial\endcsname{Lean: full statement}
\expandafter\def\csname LeanTarget@L:d1-spatial\endcsname{lean-entry:L:d1-spatial}
\expandafter\def\csname LeanShort@L:d1-spatial\endcsname{full}
\expandafter\def\csname LeanTag@L:damped-homogeneous-norms\endcsname{Lean: full statement}
\expandafter\def\csname LeanTarget@L:damped-homogeneous-norms\endcsname{lean-entry:L:damped-homogeneous-norms}
\expandafter\def\csname LeanShort@L:damped-homogeneous-norms\endcsname{full}
\expandafter\def\csname LeanTag@L:driver-holder-limit\endcsname{Lean: partial, with inputs}
\expandafter\def\csname LeanTarget@L:driver-holder-limit\endcsname{lean-entry:L:driver-holder-limit}
\expandafter\def\csname LeanShort@L:driver-holder-limit\endcsname{partial}
\expandafter\def\csname LeanTag@L:driver-interval-limit\endcsname{Lean: partial, with inputs}
\expandafter\def\csname LeanTarget@L:driver-interval-limit\endcsname{lean-entry:L:driver-interval-limit}
\expandafter\def\csname LeanShort@L:driver-interval-limit\endcsname{partial}
\expandafter\def\csname LeanTag@L:endpoint-closure\endcsname{Lean: full statement}
\expandafter\def\csname LeanTarget@L:endpoint-closure\endcsname{lean-entry:L:endpoint-closure}
\expandafter\def\csname LeanShort@L:endpoint-closure\endcsname{full}
\expandafter\def\csname LeanTag@L:equilibrium-gamma-beta\endcsname{Lean: full statement}
\expandafter\def\csname LeanTarget@L:equilibrium-gamma-beta\endcsname{lean-entry:L:equilibrium-gamma-beta}
\expandafter\def\csname LeanShort@L:equilibrium-gamma-beta\endcsname{full}
\expandafter\def\csname LeanTag@L:exact-canonical-scaling\endcsname{Lean: partial, with inputs}
\expandafter\def\csname LeanTarget@L:exact-canonical-scaling\endcsname{lean-entry:L:exact-canonical-scaling}
\expandafter\def\csname LeanShort@L:exact-canonical-scaling\endcsname{partial}
\expandafter\def\csname LeanTag@L:finite-word-cutoff\endcsname{Lean: full statement}
\expandafter\def\csname LeanTarget@L:finite-word-cutoff\endcsname{lean-entry:L:finite-word-cutoff}
\expandafter\def\csname LeanShort@L:finite-word-cutoff\endcsname{full}
\expandafter\def\csname LeanTag@L:fixed-split\endcsname{Lean: partial, with inputs}
\expandafter\def\csname LeanTarget@L:fixed-split\endcsname{lean-entry:L:fixed-split}
\expandafter\def\csname LeanShort@L:fixed-split\endcsname{partial}
\expandafter\def\csname LeanTag@L:fock-form\endcsname{Lean: narrower statement}
\expandafter\def\csname LeanTarget@L:fock-form\endcsname{lean-entry:L:fock-form}
\expandafter\def\csname LeanShort@L:fock-form\endcsname{narrower}
\expandafter\def\csname LeanTag@L:gaussian-scalar-powers\endcsname{Lean: full statement}
\expandafter\def\csname LeanTarget@L:gaussian-scalar-powers\endcsname{lean-entry:L:gaussian-scalar-powers}
\expandafter\def\csname LeanShort@L:gaussian-scalar-powers\endcsname{full}
\expandafter\def\csname LeanTag@L:gaussian-trace-bounds\endcsname{Lean: full statement}
\expandafter\def\csname LeanTarget@L:gaussian-trace-bounds\endcsname{lean-entry:L:gaussian-trace-bounds}
\expandafter\def\csname LeanShort@L:gaussian-trace-bounds\endcsname{full}
\expandafter\def\csname LeanTag@L:gram-scalar-algebra\endcsname{Lean: full statement}
\expandafter\def\csname LeanTarget@L:gram-scalar-algebra\endcsname{lean-entry:L:gram-scalar-algebra}
\expandafter\def\csname LeanShort@L:gram-scalar-algebra\endcsname{full}
\expandafter\def\csname LeanTag@L:higher-riesz\endcsname{Lean: full statement}
\expandafter\def\csname LeanTarget@L:higher-riesz\endcsname{lean-entry:L:higher-riesz}
\expandafter\def\csname LeanShort@L:higher-riesz\endcsname{full}
\expandafter\def\csname LeanTag@L:interval-complex\endcsname{Lean: partial, with inputs}
\expandafter\def\csname LeanTarget@L:interval-complex\endcsname{lean-entry:L:interval-complex}
\expandafter\def\csname LeanShort@L:interval-complex\endcsname{partial}
\expandafter\def\csname LeanTag@L:interval-contractions\endcsname{Lean: narrower statement}
\expandafter\def\csname LeanTarget@L:interval-contractions\endcsname{lean-entry:L:interval-contractions}
\expandafter\def\csname LeanShort@L:interval-contractions\endcsname{narrower}
\expandafter\def\csname LeanTag@L:inventory-product\endcsname{Lean: full statement}
\expandafter\def\csname LeanTarget@L:inventory-product\endcsname{lean-entry:L:inventory-product}
\expandafter\def\csname LeanShort@L:inventory-product\endcsname{full}
\expandafter\def\csname LeanTag@L:low-alpha-convolution\endcsname{Lean: full statement}
\expandafter\def\csname LeanTarget@L:low-alpha-convolution\endcsname{lean-entry:L:low-alpha-convolution}
\expandafter\def\csname LeanShort@L:low-alpha-convolution\endcsname{full}
\expandafter\def\csname LeanTag@L:masked-continuous-wave\endcsname{Lean: full statement}
\expandafter\def\csname LeanTarget@L:masked-continuous-wave\endcsname{lean-entry:L:masked-continuous-wave}
\expandafter\def\csname LeanShort@L:masked-continuous-wave\endcsname{full}
\expandafter\def\csname LeanTag@L:masked-temporal-modulation\endcsname{Lean: full statement}
\expandafter\def\csname LeanTarget@L:masked-temporal-modulation\endcsname{lean-entry:L:masked-temporal-modulation}
\expandafter\def\csname LeanShort@L:masked-temporal-modulation\endcsname{full}
\expandafter\def\csname LeanTag@L:masked-time-completion\endcsname{Lean: full statement}
\expandafter\def\csname LeanTarget@L:masked-time-completion\endcsname{lean-entry:L:masked-time-completion}
\expandafter\def\csname LeanShort@L:masked-time-completion\endcsname{full}
\expandafter\def\csname LeanTag@L:masked-trig-products\endcsname{Lean: full statement}
\expandafter\def\csname LeanTarget@L:masked-trig-products\endcsname{lean-entry:L:masked-trig-products}
\expandafter\def\csname LeanShort@L:masked-trig-products\endcsname{full}
\expandafter\def\csname LeanTag@L:occupation-kernel-row\endcsname{Lean: full statement}
\expandafter\def\csname LeanTarget@L:occupation-kernel-row\endcsname{lean-entry:L:occupation-kernel-row}
\expandafter\def\csname LeanShort@L:occupation-kernel-row\endcsname{full}
\expandafter\def\csname LeanTag@L:occupation-raw-square\endcsname{Lean: narrower statement}
\expandafter\def\csname LeanTarget@L:occupation-raw-square\endcsname{lean-entry:L:occupation-raw-square}
\expandafter\def\csname LeanShort@L:occupation-raw-square\endcsname{narrower}
\expandafter\def\csname LeanTag@L:occupation-rayleigh-selection\endcsname{Lean: narrower statement}
\expandafter\def\csname LeanTarget@L:occupation-rayleigh-selection\endcsname{lean-entry:L:occupation-rayleigh-selection}
\expandafter\def\csname LeanShort@L:occupation-rayleigh-selection\endcsname{narrower}
\expandafter\def\csname LeanTag@L:phase-bessel\endcsname{Lean: full statement}
\expandafter\def\csname LeanTarget@L:phase-bessel\endcsname{lean-entry:L:phase-bessel}
\expandafter\def\csname LeanShort@L:phase-bessel\endcsname{full}
\expandafter\def\csname LeanTag@L:phase-contractions\endcsname{Lean: full statement}
\expandafter\def\csname LeanTarget@L:phase-contractions\endcsname{lean-entry:L:phase-contractions}
\expandafter\def\csname LeanShort@L:phase-contractions\endcsname{full}
\expandafter\def\csname LeanTag@L:plane-endpoints\endcsname{Lean: full statement}
\expandafter\def\csname LeanTarget@L:plane-endpoints\endcsname{lean-entry:L:plane-endpoints}
\expandafter\def\csname LeanShort@L:plane-endpoints\endcsname{full}
\expandafter\def\csname LeanTag@L:positive-physical-gram\endcsname{Lean: narrower statement}
\expandafter\def\csname LeanTarget@L:positive-physical-gram\endcsname{lean-entry:L:positive-physical-gram}
\expandafter\def\csname LeanShort@L:positive-physical-gram\endcsname{narrower}
\expandafter\def\csname LeanTag@L:positive-raw-embedding\endcsname{Lean: narrower statement}
\expandafter\def\csname LeanTarget@L:positive-raw-embedding\endcsname{lean-entry:L:positive-raw-embedding}
\expandafter\def\csname LeanShort@L:positive-raw-embedding\endcsname{narrower}
\expandafter\def\csname LeanTag@L:positive-series-type\endcsname{Lean: full statement}
\expandafter\def\csname LeanTarget@L:positive-series-type\endcsname{lean-entry:L:positive-series-type}
\expandafter\def\csname LeanShort@L:positive-series-type\endcsname{full}
\expandafter\def\csname LeanTag@L:product-endpoints\endcsname{Lean: full statement}
\expandafter\def\csname LeanTarget@L:product-endpoints\endcsname{lean-entry:L:product-endpoints}
\expandafter\def\csname LeanShort@L:product-endpoints\endcsname{full}
\expandafter\def\csname LeanTag@L:raw-block-indexed-cone\endcsname{Lean: partial, with inputs}
\expandafter\def\csname LeanTarget@L:raw-block-indexed-cone\endcsname{lean-entry:L:raw-block-indexed-cone}
\expandafter\def\csname LeanShort@L:raw-block-indexed-cone\endcsname{partial}
\expandafter\def\csname LeanTag@L:raw-block-inputs\endcsname{Lean: partial, with inputs}
\expandafter\def\csname LeanTarget@L:raw-block-inputs\endcsname{lean-entry:L:raw-block-inputs}
\expandafter\def\csname LeanShort@L:raw-block-inputs\endcsname{partial}
\expandafter\def\csname LeanTag@L:raw-block-restart\endcsname{Lean: partial, with inputs}
\expandafter\def\csname LeanTarget@L:raw-block-restart\endcsname{lean-entry:L:raw-block-restart}
\expandafter\def\csname LeanShort@L:raw-block-restart\endcsname{partial}
\expandafter\def\csname LeanTag@L:source-finite-comparison\endcsname{Lean: full statement}
\expandafter\def\csname LeanTarget@L:source-finite-comparison\endcsname{lean-entry:L:source-finite-comparison}
\expandafter\def\csname LeanShort@L:source-finite-comparison\endcsname{full}
\expandafter\def\csname LeanTag@L:source-spatial-heat\endcsname{Lean: full statement}
\expandafter\def\csname LeanTarget@L:source-spatial-heat\endcsname{lean-entry:L:source-spatial-heat}
\expandafter\def\csname LeanShort@L:source-spatial-heat\endcsname{full}
\expandafter\def\csname LeanTag@L:spatial-fourier-bound\endcsname{Lean: full statement}
\expandafter\def\csname LeanTarget@L:spatial-fourier-bound\endcsname{lean-entry:L:spatial-fourier-bound}
\expandafter\def\csname LeanShort@L:spatial-fourier-bound\endcsname{full}
\expandafter\def\csname LeanTag@L:spatial-physical-pair\endcsname{Lean: full statement}
\expandafter\def\csname LeanTarget@L:spatial-physical-pair\endcsname{lean-entry:L:spatial-physical-pair}
\expandafter\def\csname LeanShort@L:spatial-physical-pair\endcsname{full}
\expandafter\def\csname LeanTag@L:spatial-static-controls\endcsname{Lean: narrower statement}
\expandafter\def\csname LeanTarget@L:spatial-static-controls\endcsname{lean-entry:L:spatial-static-controls}
\expandafter\def\csname LeanShort@L:spatial-static-controls\endcsname{narrower}
\expandafter\def\csname LeanTag@L:symplectic-free\endcsname{Lean: partial, with inputs}
\expandafter\def\csname LeanTarget@L:symplectic-free\endcsname{lean-entry:L:symplectic-free}
\expandafter\def\csname LeanShort@L:symplectic-free\endcsname{partial}
\expandafter\def\csname LeanTag@L:temporal-equilibrium\endcsname{Lean: full statement}
\expandafter\def\csname LeanTarget@L:temporal-equilibrium\endcsname{lean-entry:L:temporal-equilibrium}
\expandafter\def\csname LeanShort@L:temporal-equilibrium\endcsname{full}
\expandafter\def\csname LeanTag@L:terminal-wick\endcsname{Lean: full statement}
\expandafter\def\csname LeanTarget@L:terminal-wick\endcsname{lean-entry:L:terminal-wick}
\expandafter\def\csname LeanShort@L:terminal-wick\endcsname{full}
\expandafter\def\csname LeanTag@L:tube-half-action\endcsname{Lean: full statement}
\expandafter\def\csname LeanTarget@L:tube-half-action\endcsname{lean-entry:L:tube-half-action}
\expandafter\def\csname LeanShort@L:tube-half-action\endcsname{full}
\expandafter\def\csname LeanTag@L:tube-rayleigh-tail\endcsname{Lean: narrower statement}
\expandafter\def\csname LeanTarget@L:tube-rayleigh-tail\endcsname{lean-entry:L:tube-rayleigh-tail}
\expandafter\def\csname LeanShort@L:tube-rayleigh-tail\endcsname{narrower}
\expandafter\def\csname LeanTag@L:uniform-series-tail\endcsname{Lean: full statement}
\expandafter\def\csname LeanTarget@L:uniform-series-tail\endcsname{lean-entry:L:uniform-series-tail}
\expandafter\def\csname LeanShort@L:uniform-series-tail\endcsname{full}
\expandafter\def\csname LeanTag@L:variational-control-dilation\endcsname{Lean: narrower statement}
\expandafter\def\csname LeanTarget@L:variational-control-dilation\endcsname{lean-entry:L:variational-control-dilation}
\expandafter\def\csname LeanShort@L:variational-control-dilation\endcsname{narrower}
\expandafter\def\csname LeanTag@L:variational-global-mean\endcsname{Lean: full statement}
\expandafter\def\csname LeanTarget@L:variational-global-mean\endcsname{lean-entry:L:variational-global-mean}
\expandafter\def\csname LeanShort@L:variational-global-mean\endcsname{full}
\expandafter\def\csname LeanTag@L:variational-heat-envelope\endcsname{Lean: full statement}
\expandafter\def\csname LeanTarget@L:variational-heat-envelope\endcsname{lean-entry:L:variational-heat-envelope}
\expandafter\def\csname LeanShort@L:variational-heat-envelope\endcsname{full}
\expandafter\def\csname LeanTag@L:white-clock-laplace\endcsname{Lean: full statement}
\expandafter\def\csname LeanTarget@L:white-clock-laplace\endcsname{lean-entry:L:white-clock-laplace}
\expandafter\def\csname LeanShort@L:white-clock-laplace\endcsname{full}
\expandafter\def\csname LeanTag@L:white-test\endcsname{Lean: partial, with inputs}
\expandafter\def\csname LeanTarget@L:white-test\endcsname{lean-entry:L:white-test}
\expandafter\def\csname LeanShort@L:white-test\endcsname{partial}
\expandafter\def\csname LeanTag@P:actual-wave-form\endcsname{Lean: partial, with inputs}
\expandafter\def\csname LeanTarget@P:actual-wave-form\endcsname{lean-entry:P:actual-wave-form}
\expandafter\def\csname LeanShort@P:actual-wave-form\endcsname{partial}
\expandafter\def\csname LeanTag@P:clock-copy-obstruction\endcsname{Lean: full statement}
\expandafter\def\csname LeanTarget@P:clock-copy-obstruction\endcsname{lean-entry:P:clock-copy-obstruction}
\expandafter\def\csname LeanShort@P:clock-copy-obstruction\endcsname{full}
\expandafter\def\csname LeanTag@P:d1-upper\endcsname{Lean: partial, with inputs}
\expandafter\def\csname LeanTarget@P:d1-upper\endcsname{lean-entry:P:d1-upper}
\expandafter\def\csname LeanShort@P:d1-upper\endcsname{partial}
\expandafter\def\csname LeanTag@P:finite-gaussian-matrix\endcsname{Lean: full statement}
\expandafter\def\csname LeanTarget@P:finite-gaussian-matrix\endcsname{lean-entry:P:finite-gaussian-matrix}
\expandafter\def\csname LeanShort@P:finite-gaussian-matrix\endcsname{full}
\expandafter\def\csname LeanTag@P:frequency-row\endcsname{Lean: full statement}
\expandafter\def\csname LeanTarget@P:frequency-row\endcsname{lean-entry:P:frequency-row}
\expandafter\def\csname LeanShort@P:frequency-row\endcsname{full}
\expandafter\def\csname LeanTag@P:gram-entrywise-comparison\endcsname{Lean: partial, with inputs}
\expandafter\def\csname LeanTarget@P:gram-entrywise-comparison\endcsname{lean-entry:P:gram-entrywise-comparison}
\expandafter\def\csname LeanShort@P:gram-entrywise-comparison\endcsname{partial}
\expandafter\def\csname LeanTag@P:gram-positive-optimizer\endcsname{Lean: partial, with inputs}
\expandafter\def\csname LeanTarget@P:gram-positive-optimizer\endcsname{lean-entry:P:gram-positive-optimizer}
\expandafter\def\csname LeanShort@P:gram-positive-optimizer\endcsname{partial}
\expandafter\def\csname LeanTag@P:interval-coefficients\endcsname{Lean: partial, with inputs}
\expandafter\def\csname LeanTarget@P:interval-coefficients\endcsname{lean-entry:P:interval-coefficients}
\expandafter\def\csname LeanShort@P:interval-coefficients\endcsname{partial}
\expandafter\def\csname LeanTag@P:low-alpha-upper\endcsname{Lean: partial, with inputs}
\expandafter\def\csname LeanTarget@P:low-alpha-upper\endcsname{lean-entry:P:low-alpha-upper}
\expandafter\def\csname LeanShort@P:low-alpha-upper\endcsname{partial}
\expandafter\def\csname LeanTag@P:masked-auxiliary-norm\endcsname{Lean: full statement}
\expandafter\def\csname LeanTarget@P:masked-auxiliary-norm\endcsname{lean-entry:P:masked-auxiliary-norm}
\expandafter\def\csname LeanShort@P:masked-auxiliary-norm\endcsname{full}
\expandafter\def\csname LeanTag@P:masked-noise-completion\endcsname{Lean: partial, with inputs}
\expandafter\def\csname LeanTarget@P:masked-noise-completion\endcsname{lean-entry:P:masked-noise-completion}
\expandafter\def\csname LeanShort@P:masked-noise-completion\endcsname{partial}
\expandafter\def\csname LeanTag@P:masked-quotient-obstruction\endcsname{Lean: narrower statement}
\expandafter\def\csname LeanTarget@P:masked-quotient-obstruction\endcsname{lean-entry:P:masked-quotient-obstruction}
\expandafter\def\csname LeanShort@P:masked-quotient-obstruction\endcsname{narrower}
\expandafter\def\csname LeanTag@P:masked-spectral-growth\endcsname{Lean: full statement}
\expandafter\def\csname LeanTarget@P:masked-spectral-growth\endcsname{lean-entry:P:masked-spectral-growth}
\expandafter\def\csname LeanShort@P:masked-spectral-growth\endcsname{full}
\expandafter\def\csname LeanTag@P:mild-construction\endcsname{Lean: partial, with inputs}
\expandafter\def\csname LeanTarget@P:mild-construction\endcsname{lean-entry:P:mild-construction}
\expandafter\def\csname LeanShort@P:mild-construction\endcsname{partial}
\expandafter\def\csname LeanTag@P:moment-growth-logconvex\endcsname{Lean: partial, with inputs}
\expandafter\def\csname LeanTarget@P:moment-growth-logconvex\endcsname{lean-entry:P:moment-growth-logconvex}
\expandafter\def\csname LeanShort@P:moment-growth-logconvex\endcsname{partial}
\expandafter\def\csname LeanTag@P:occupation-physical-class\endcsname{Lean: full statement}
\expandafter\def\csname LeanTarget@P:occupation-physical-class\endcsname{lean-entry:P:occupation-physical-class}
\expandafter\def\csname LeanShort@P:occupation-physical-class\endcsname{full}
\expandafter\def\csname LeanTag@P:radius-mean-anchor\endcsname{Lean: partial, with inputs}
\expandafter\def\csname LeanTarget@P:radius-mean-anchor\endcsname{lean-entry:P:radius-mean-anchor}
\expandafter\def\csname LeanShort@P:radius-mean-anchor\endcsname{partial}
\expandafter\def\csname LeanTag@P:raw-block-original-lower\endcsname{Lean: partial, with inputs}
\expandafter\def\csname LeanTarget@P:raw-block-original-lower\endcsname{lean-entry:P:raw-block-original-lower}
\expandafter\def\csname LeanShort@P:raw-block-original-lower\endcsname{partial}
\expandafter\def\csname LeanTag@P:raw-occupation-energy\endcsname{Lean: full statement}
\expandafter\def\csname LeanTarget@P:raw-occupation-energy\endcsname{lean-entry:P:raw-occupation-energy}
\expandafter\def\csname LeanShort@P:raw-occupation-energy\endcsname{full}
\expandafter\def\csname LeanTag@P:riesz-small\endcsname{Lean: partial, with inputs}
\expandafter\def\csname LeanTarget@P:riesz-small\endcsname{lean-entry:P:riesz-small}
\expandafter\def\csname LeanShort@P:riesz-small\endcsname{partial}
\expandafter\def\csname LeanTag@P:signed-closed-form\endcsname{Lean: partial, with inputs}
\expandafter\def\csname LeanTarget@P:signed-closed-form\endcsname{lean-entry:P:signed-closed-form}
\expandafter\def\csname LeanShort@P:signed-closed-form\endcsname{partial}
\expandafter\def\csname LeanTag@P:spatial-energy-map\endcsname{Lean: narrower statement}
\expandafter\def\csname LeanTarget@P:spatial-energy-map\endcsname{lean-entry:P:spatial-energy-map}
\expandafter\def\csname LeanShort@P:spatial-energy-map\endcsname{narrower}
\expandafter\def\csname LeanTag@P:spatial-four-linear\endcsname{Lean: full statement}
\expandafter\def\csname LeanTarget@P:spatial-four-linear\endcsname{lean-entry:P:spatial-four-linear}
\expandafter\def\csname LeanShort@P:spatial-four-linear\endcsname{full}
\expandafter\def\csname LeanTag@P:spatial-smooth-ratio\endcsname{Lean: narrower statement}
\expandafter\def\csname LeanTarget@P:spatial-smooth-ratio\endcsname{lean-entry:P:spatial-smooth-ratio}
\expandafter\def\csname LeanShort@P:spatial-smooth-ratio\endcsname{narrower}
\expandafter\def\csname LeanTag@P:static-phase-all-orders\endcsname{Lean: partial, with inputs}
\expandafter\def\csname LeanTarget@P:static-phase-all-orders\endcsname{lean-entry:P:static-phase-all-orders}
\expandafter\def\csname LeanShort@P:static-phase-all-orders\endcsname{partial}
\expandafter\def\csname LeanTag@P:tube-rayleigh-transfer\endcsname{Lean: supplied inputs}
\expandafter\def\csname LeanTarget@P:tube-rayleigh-transfer\endcsname{lean-entry:P:tube-rayleigh-transfer}
\expandafter\def\csname LeanShort@P:tube-rayleigh-transfer\endcsname{supplied}
\expandafter\def\csname LeanTag@P:variational-positive-scores\endcsname{Lean: narrower statement}
\expandafter\def\csname LeanTarget@P:variational-positive-scores\endcsname{lean-entry:P:variational-positive-scores}
\expandafter\def\csname LeanShort@P:variational-positive-scores\endcsname{narrower}
\expandafter\def\csname LeanTag@P:variational-remove-mask\endcsname{Lean: partial, with inputs}
\expandafter\def\csname LeanTarget@P:variational-remove-mask\endcsname{lean-entry:P:variational-remove-mask}
\expandafter\def\csname LeanShort@P:variational-remove-mask\endcsname{partial}
\expandafter\def\csname LeanTag@P:variational-unmasked-rate\endcsname{Lean: narrower statement}
\expandafter\def\csname LeanTarget@P:variational-unmasked-rate\endcsname{lean-entry:P:variational-unmasked-rate}
\expandafter\def\csname LeanShort@P:variational-unmasked-rate\endcsname{narrower}
\expandafter\def\csname LeanTag@P:wave-two-variances\endcsname{Lean: partial, with inputs}
\expandafter\def\csname LeanTarget@P:wave-two-variances\endcsname{lean-entry:P:wave-two-variances}
\expandafter\def\csname LeanShort@P:wave-two-variances\endcsname{partial}
\expandafter\def\csname LeanTag@P:white-four\endcsname{Lean: full statement}
\expandafter\def\csname LeanTarget@P:white-four\endcsname{lean-entry:P:white-four}
\expandafter\def\csname LeanShort@P:white-four\endcsname{full}
\expandafter\def\csname LeanTag@P:white-spatial-energy\endcsname{Lean: partial, with inputs}
\expandafter\def\csname LeanTarget@P:white-spatial-energy\endcsname{lean-entry:P:white-spatial-energy}
\expandafter\def\csname LeanShort@P:white-spatial-energy\endcsname{partial}
\expandafter\def\csname LeanTag@R:gaussian-one-variance-test\endcsname{Lean: full statement}
\expandafter\def\csname LeanTarget@R:gaussian-one-variance-test\endcsname{lean-entry:R:gaussian-one-variance-test}
\expandafter\def\csname LeanShort@R:gaussian-one-variance-test\endcsname{full}
\expandafter\def\csname LeanTag@T:intrinsic-occupation-ceiling\endcsname{Lean: narrower statement}
\expandafter\def\csname LeanTarget@T:intrinsic-occupation-ceiling\endcsname{lean-entry:T:intrinsic-occupation-ceiling}
\expandafter\def\csname LeanShort@T:intrinsic-occupation-ceiling\endcsname{narrower}
\expandafter\def\csname LeanTag@T:white-moment-bounds\endcsname{Lean: partial, with inputs}
\expandafter\def\csname LeanTarget@T:white-moment-bounds\endcsname{lean-entry:T:white-moment-bounds}
\expandafter\def\csname LeanShort@T:white-moment-bounds\endcsname{partial}

\makeatletter
\let\lean@pending\@empty
\newcommand{\LeanMarker}[1]{\gdef\lean@pending{#1}}
\newcommand{\lean@margin}[1]{%
\expandafter\ifx\csname LeanTag@#1\endcsname\@empty\else
\makebox[0pt][r]{\smash{\normalfont\footnotesize
\begin{tabular}[t]{@{}r@{}}%
  \hyperlink{\csname LeanTarget@#1\endcsname}{Lean\vphantom{Lp}}\\[-2pt]
  \hyperlink{\csname LeanTarget@#1\endcsname}{\csname
  LeanShort@#1\endcsname\vphantom{Lp}}%
\end{tabular}}\hskip 0.9em}%
\fi}
\newcommand{\LeanEqMarker}[1]{%
\expandafter\ifx\csname LeanTag@#1\endcsname\@empty\else
\mathopen{\hbox to\z@{\tikz[remember picture,overlay]
\node[anchor=base east,inner sep=0pt,xshift=\dimexpr 1in+\oddsidemargin-10.8pt\relax]
at (current page.west |- 0,0) {\normalfont\footnotesize
\begin{tabular}[t]{@{}r@{}}%
  \hyperlink{\csname LeanTarget@#1\endcsname}{Lean\vphantom{Lp}}\\[-2pt]
  \hyperlink{\csname LeanTarget@#1\endcsname}{\csname
  LeanShort@#1\endcsname\vphantom{Lp}}%
\end{tabular}};\hss}}%
\fi}
\let\lean@thmhead@plain\thmhead@plain
\renewcommand{\thmhead@plain}[3]{%
\ifx\lean@pending\@empty\else
\@ifundefined{LeanTag@\lean@pending}{\PackageError{leanmarker}{No marker
data for \lean@pending}{Regenerate notes/lean-markers-data.tex.}}{}%
\lean@margin{\lean@pending}%
\global\let\lean@pending\@empty
\fi
\lean@thmhead@plain{#1}{#2}{#3}}
\let\thmhead\thmhead@plain
\makeatother
\makeatletter
\def\lc@nodot#1.\@nil{#1}
\newcommand{\lcCategory}[2]{\subsection[\lc@nodot#1\@nil]{\raggedright\lc@nodot#1\@nil~(#2)}}
\makeatother
\makeatletter
\newenvironment{lcTable}{\global\lc@leadtrue\noindent\small\begin{tabular}{@{}>{\raggedright\arraybackslash}p{0.27\textwidth}>{\raggedright\arraybackslash}p{0.70\textwidth}@{}}}{\end{tabular}\par}
\newcommand{\lcEntry}[2]{\par\Needspace{4\baselineskip}\smallskip
\global\lc@leadtrue\noindent\textbf{#1.}\enspace\emph{#2.}\par\nopagebreak}

\newcommand{\lcPart}[2]{{\small\leftskip=1.5em\noindent\textit{#1:} #2\par}}
\newcommand{\lcInputs}[1]{{\small\leftskip=1.5em\noindent\textit{#1:}\par
\nopagebreak}}
\newcommand{\lcInput}[2]{{\small\leftskip=3.3em\noindent
\llap{(#1)\enspace}#2\par}}
\newcommand{\lcEntryEnd}{\par\smallskip}
\makeatother
\usepackage{url}
\numberwithin{equation}{section}

\newtheorem{theorem}{Theorem}[section]
\newtheorem{lemma}[theorem]{Lemma}
\newtheorem{proposition}[theorem]{Proposition}
\newtheorem{corollary}[theorem]{Corollary}
\newtheorem{conjecture}[theorem]{Conjecture}
\newtheoremstyle{nogostyle}{0pt}{0pt}{\normalfont}{}
{\bfseries\color{hamblue}}{.}{0.5em}{}
\theoremstyle{nogostyle}
\newtheorem{nogo}{No-go}
\renewcommand{\thenogo}{\ifnum\value{nogo}<10 0\fi\arabic{nogo}}
\tcolorboxenvironment{nogo}{
enhanced jigsaw,
breakable,
colback=hamblue!6!white,
colframe=hamblue!32!white,
boxrule=0.35pt,
arc=2mm,
boxsep=0pt,
left=6pt,
right=6pt,
top=5pt,
bottom=5pt,
before skip=8pt,
after skip=8pt
}
\theoremstyle{definition}

\newtheorem{hypothesis}[theorem]{Hypothesis}
\newtheorem{remark}[theorem]{Remark}

\newcommand{\R}{\mathbb{R}}
\newcommand{\N}{\mathbb{N}}
\newcommand{\E}{\mathbb{E}}
\newcommand{\W}{\dot W}
\newcommand{\ud}{\,\mathrm{d}}
\newcommand{\cE}{\mathcal{E}}

\newcommand{\ProofRef}[2]{#1~\ref{#2}}
\newcommand{\GuideRef}[2]{#1~\ref{#2}}

\usepackage{microtype}
\usepackage{fancyhdr}
\usepackage{titlesec}
\usepackage{booktabs}
\fancypagestyle{plain}{%
\fancyhf{}
\fancyfoot[C]{\thepage}

}
\titleformat{\chapter}[display]
{\normalfont\huge\bfseries\raggedright\color{hamblue}}
{\normalfont\large\scshape\chaptertitlename\ \thechapter}{10pt}{}
\titleformat{name=\chapter,numberless}[display]
{\normalfont\huge\bfseries\raggedright\color{hamblue}}{}{0pt}{}
\titlespacing*{\chapter}{0pt}{-12pt}{28pt}
\titleformat{\section}
{\normalfont\Large\bfseries\raggedright\color{hamblue}}{\thesection}{0.75em}{}
\titleformat{\subsection}
{\normalfont\large\bfseries}{\thesubsection}{0.75em}{}
\hypersetup{bookmarksnumbered=true, linkcolor=hamblue,
citecolor=hamteal, urlcolor=hamblue}
\newenvironment{abstract}{%
\chapter*{Abstract}
\addcontentsline{toc}{chapter}{Abstract}
}{}
\makeatletter
\renewcommand*{\l@subsection}{\@dottedtocline{2}{3.8em}{4em}}
\renewcommand{\maketitle}{%
\begin{titlepage}
  \thispdfpagelabel{Cover}
  \raggedright
  \vspace*{0.8cm}
  {\small\sffamily\color{hamblue}\textls[120]{STOCHASTIC WAVE EQUATIONS}\par}
  \vspace{0.4cm}
  {\color{hamblue}\rule{\textwidth}{0.8pt}\par}
  \vspace{1.5cm}
  {\fontsize{32}{39}\selectfont\bfseries\color{hamblue}\@title\par}
  \vspace{0.7cm}
  {\Large Second Moment and Exact Asymptotics\par}
  \vspace{0.45cm}
  {\large\itshape An experiment in doing mathematics with AI\par}
  \vspace{1.1cm}
  {\LARGE\@author\par}
  \vspace{0.35cm}
  {\normalsize Department of Mathematics and Statistics\par}
  \vspace{0.12cm}
  {\large Auburn University\par}
  \vfill
  {\large Construction, moment bounds\par}
  {\large and the variational constant\par}
  \vspace{0.6cm}
  {\color{hamgold}\rule{0.14\textwidth}{1.2pt}\par}
  \vspace{0.6cm}
  {\small Version 1.0.0\hfill\@date\par}
\end{titlepage}
}
\renewcommand{\listoffigures}{%
\chapter*{Illustrations}
\addcontentsline{toc}{chapter}{Illustrations}
\markboth{Illustrations}{Illustrations}
\@starttoc{lof}
}
\makeatother

\begin{document}
\frontmatter
\hypersetup{pageanchor=false}

\title{Stratonovich Waves\\
Driven by Space-Time\\
Gaussian Noise}
\author{Le Chen}
\date{November~16, 2026}
\hypersetup{pdftitle={Stratonovich Waves Driven by Space-Time Gaussian Noise: Second Moment and Exact Asymptotics},
pdfauthor={Le Chen}, pdfsubject={Stratonovich moment asymptotics}}
\maketitle
\clearpage
% Imprint page: follows the cover and, like it, carries no page number, so
% the Preface still starts on page i.
\begin{titlepage}
  \thispdfpagelabel{Imprint}
  \raggedright
  \vspace*{\fill}
  {\small
  \textbf{Stratonovich Waves Driven by Space-Time Gaussian Noise:\\
  Second Moment and Exact Asymptotics}\par
  \smallskip
  Le Chen\par
  \bigskip
  Version 1.0.0, November~16, 2026\par
  \bigskip
  This work is licensed under the Creative Commons Attribution 4.0
  International licence (CC BY 4.0):\\
  \url{https://creativecommons.org/licenses/by/4.0/}\par
  \bigskip
  Source and development record:\\
  \url{https://github.com/chenle02/stratonovich-ham-second-moment}\par
  }
\end{titlepage}

\hypersetup{pageanchor=true}

\chapter*{Preface}
\addcontentsline{toc}{chapter}{Preface}
\markboth{Preface}{Preface}

This book studies the second moment of the Stratonovich hyperbolic Anderson
model driven by a Gaussian noise that is fractional in time. For this model,
X.~Chen~\cite{chen:author2026:time-dependency} constructed the solution and
showed that $\E u^2(t,x) \le e^{Ct^2}$ for large $t$. We ask two questions. What
is the correct growth scale of $\log \E u^2(t,x)$, and what determines the
constant in front of it? The scale is $t^\kappa$ with $\kappa =
(4-\alpha-\alpha_0)/(3-\alpha)$, and the limit is given by an explicit
variational formula. The first moment is better understood because its exact
representation uses the time order of a single chain. The second moment is an
integral over two chains whose points interleave, and most of this book is about
what replaces that representation. The precise statements are in
Chapter~\ref{S:book-setting}.

\section*{How this book was written}

This book was drafted together with Amy, our AI research system. The name began
as a typo. In an email to a coauthor in the summer of 2026, the author meant to
write ``my AI Army'' but missed the letter r, so Army became Amy. The coauthor
started calling the system Amy, and the name stayed. Amy is not a single language model.
She is a research environment that we built: written rules for how the work is
done (how each step is organized, how reviews are requested, what counts as
evidence, when to stop), a research library (a bibliography, a collection of
papers and their \LaTeX\ sources, and a wiki on stochastic partial differential
equations), and interfaces through which a language model can use all of these.
For this book, Amy was run with models from four families: Anthropic Claude,
OpenAI GPT (through Codex), DeepSeek and Kimi. Amy is model-agnostic by design:
any of these can be replaced. We have used Qwen through its API in earlier
projects, though not in this one. Running open-weight models on our own
machines, such as the open-weight Qwen models, is the next thing to test.

One way to picture this is as a car. The author built the car and keeps it up:
the written rules, a library of more than a thousand skills learned from earlier
work, the research library, and the interfaces. The bibliography is
\href{https://github.com/chenle02/SPDEs-Bib}{SPDEs-Bib}, and the wiki is built
on Andrej Karpathy's LLM Wiki pattern~\cite{karpathy:26:llm}: a knowledge base
that the language models write and maintain. The models are the fuel; they can
be replaced, and this book used four. The author is the driver.

Of the parts of this car, the wiki is the one we expect to last. It is the
author's second brain. One way to see it: the papers are the source code, a
language model is the compiler, and the wiki is the compiled form, written by
models for models to read. A person can read it too. Unlike a compiler, a model
can misread, which is why a result is checked against its source before it is
used. This is also how a general model becomes a specialist. No model was
trained for this book. What Amy knows about stochastic partial differential
equations is in the wiki, not in the model, and it reflects what the author has
read and how he connects it. Change the model, and the specialist is still
there. The wiki is a graph: its pages are linked, and Amy finds what bears on a
question by searching the text and by ranking pages along their links with
PageRank, the method Google first used to rank web pages. The models will be
replaced. So, in time, will the rules, skills and interfaces around them, which
engineers now call the harness. What has been read and linked stays. The wiki
grows with every project, and the more it is used, the more it knows.
Connections that nobody planned may then appear between pages written far
apart. Drawn as a graph, it looks like a starry night.

\begin{figure}[!htb]
  \centering
  \resizebox{\ifdim\width>\textwidth\textwidth\else\width\fi}{!}{%
  \input{notes/figures/amy-observatory-drawing}}

  \smallskip
  \parbox{0.9\textwidth}{\centering\footnotesize The wiki as a night sky.
  Language models compile papers and projects into it, one star for each page.
  The car steers by it, and what the car finds comes back as new
  projects.\par}
\end{figure}

The name Amy is used in two ways. In the narrow sense, Amy is the car: the
rules, skills and interfaces, with a model in the tank, that the author drives
on one project. In the broad sense, Amy is everything in the picture: the car,
the tower that turns papers into wiki pages, and the sky they fill. When we say
that Amy drafted a proof, we mean the car. When we speak of what Amy knows, we
mean the whole.

References and results are looked up in the bibliography and the wiki instead
of being recalled, which guards against invented citations and misremembered
theorems, and it is what makes checking every cited result against its source
practical. Beyond the mathematics, the choices were ours: what goes into the
book and what stays out, how it is built, the standard of presentation, what is
claimed, what is stated only as a conjecture, and what is left for future work.
Terence Tao has said that mathematics is entering an era of proof
abundance~\cite{tao:26:math}; when proofs are plentiful, choices like these
shape a book as much as any proof does.

\begin{figure}[!tp]
  \centering
  \includegraphics[width=0.6\textwidth]{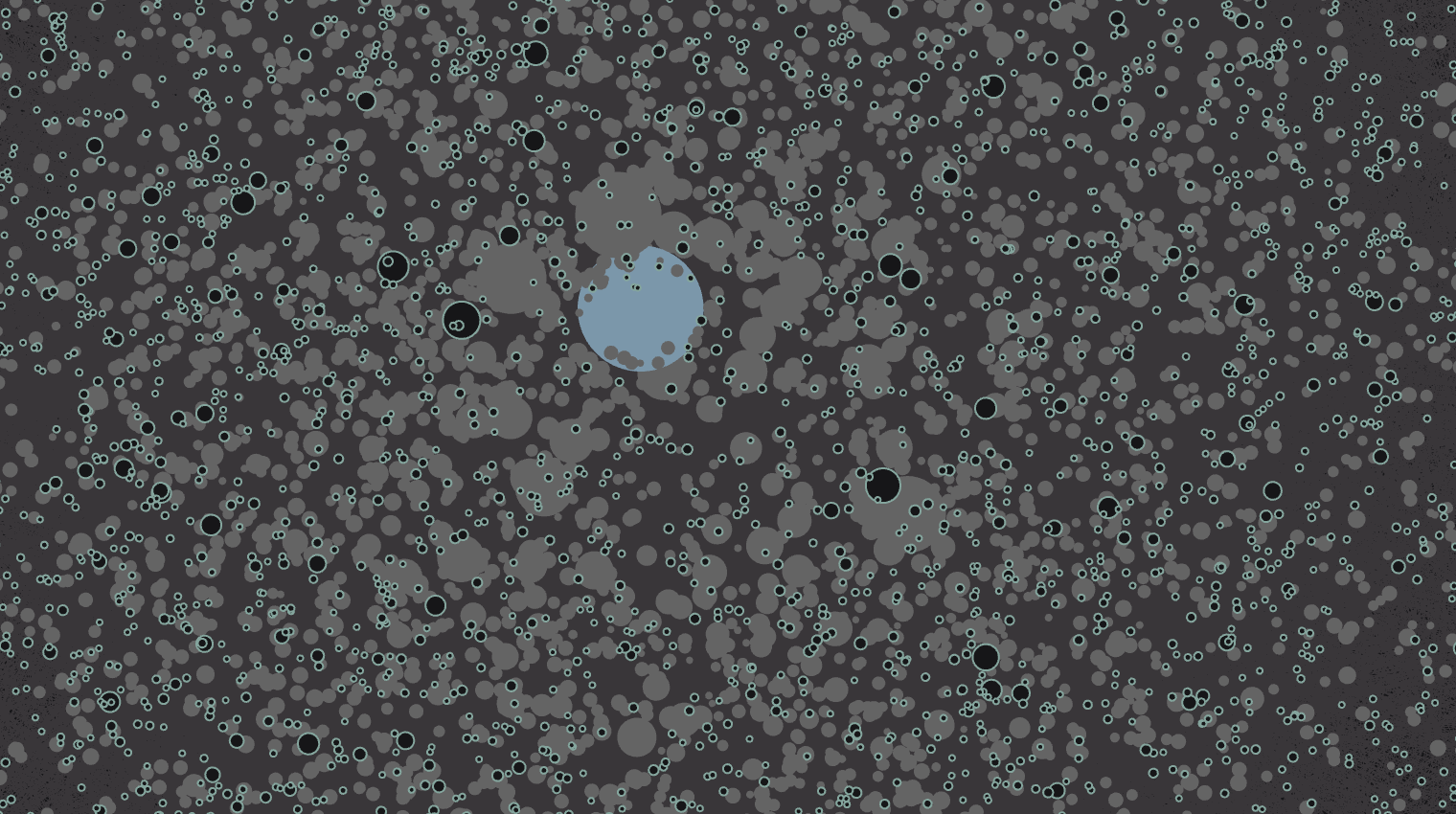}

  \smallskip
  \parbox{0.9\textwidth}{\centering\footnotesize The SPDEs-wiki as a graph: each dot is one of its more than
  4,000 pages, each a Markdown file the models wrote from a paper, a concept or a
  person they read about. With SPDEs-Bib, it is Amy's long-term memory, which
  helps her not to hallucinate. Both are parts of the infrastructure, with
  agent-agnostic skills and the command-line tools built around them.\par}
\end{figure}

The division of labour was as follows. Amy developed the mathematics: the models
working within her drafted the proofs, the \LaTeX\ source and the bookkeeping.
Our role was to drive the work: we chose the problem, set and revised the goals,
decided which approaches to pursue or abandon, chose the reviewers and decided
which results to keep. Each new proof was sent to a
reviewer from a different model family than its author, and every finding was
answered in writing. We treat an argument
produced by AI with more suspicion than one from a trusted colleague, and we ask
where each of its ideas comes from: every cited result is checked against its
source before it is used. Hairer recommends both practices in his advice on
writing mathematics with AI~\hyperlink{pp:P1}{[P1]}. The proof assistant Lean was also used: it checks some
results in full under the book's own hypotheses and the finite and algebraic
steps of others; analytic inputs that enter Lean as hypotheses are not verified
there, and two of the four main theorems, the full upper scale
(Theorem~\ref{T:full-upper}) and the second-moment constant
(Theorem~\ref{T:intrinsic-second-moment}), are outside Lean.

Mathematics is after what is true, and with AI the bottleneck has moved: it is
no longer finding an argument but understanding it. Galileo's telescope showed
Saturn's strange shape half a century before Huygens could say it was a ring;
the new instruments, too, reach results faster than anyone can understand them.
This book is such an observation, made with the checks described above and with
its development record, failed attempts included, open to view. The present
author, Le Chen, set up the instrument and ran it, and is responsible for what
it reports; later versions may have more authors, as described below. We will
look into every objection, and correct, qualify or withdraw any claim that does
not hold up, recording each change in the public history. We have not read every
line: the technical chapters have not yet been read in full by a human, and the
public table records, result by result, what was checked and how. Understanding
what we have seen is work for many eyes, and it is where we expect the next
ideas to come from. The eyepiece is open to anyone, and every check, correction
or explanation is recorded and credited by name. Amy will take part in this work
too, as she did in writing the book: under our guidance she can answer
mathematical questions, check arguments and repair them. She is not an author,
any more than Galileo's telescope was, because she cannot be held responsible
for what she produces; that responsibility rests with a person. But she did far
more than a telescope does, so we acknowledge her here by name, and the AI-use
statement records in detail which models did what.

\begin{figure}[!ht]
  \centering
  \includegraphics[width=\textwidth]{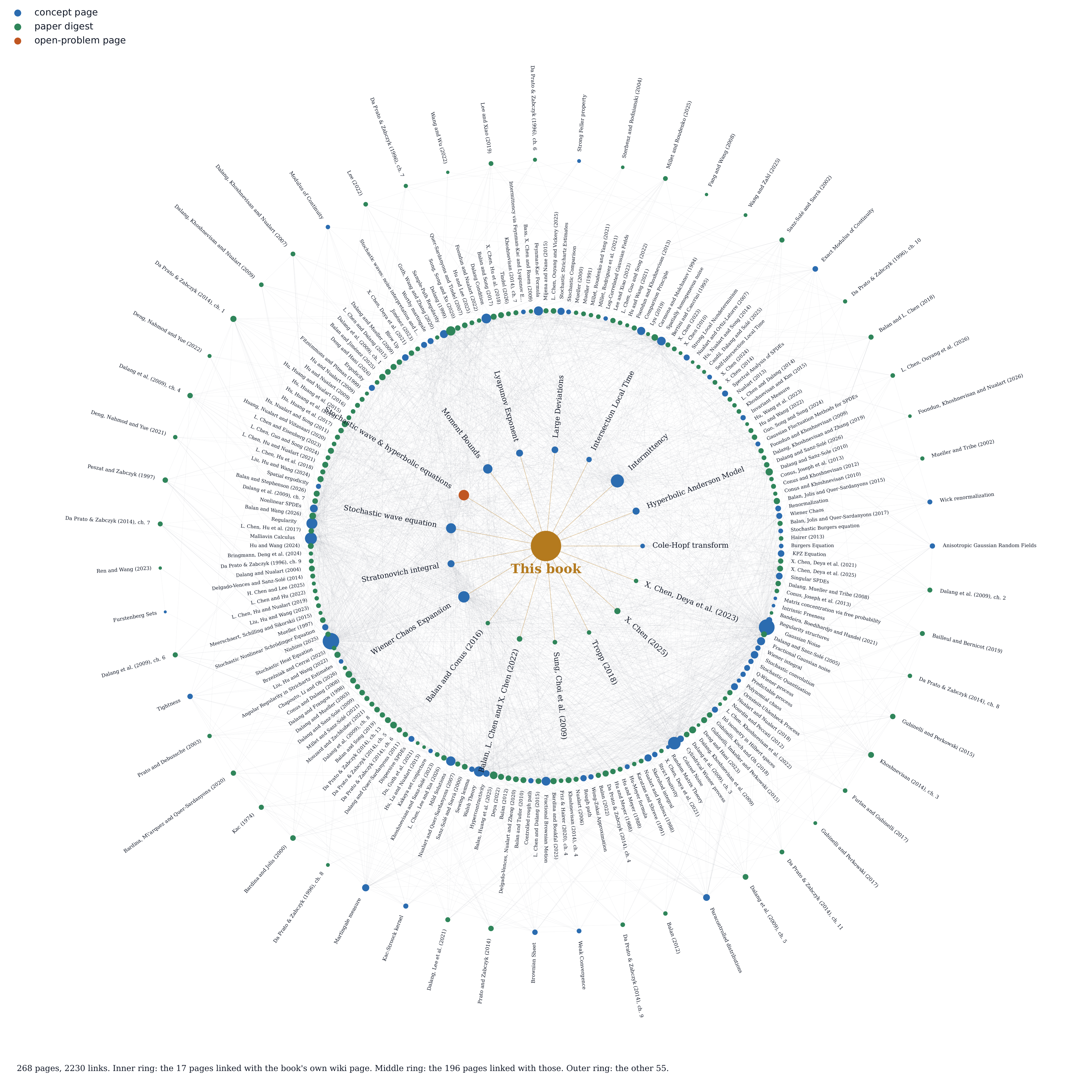}

  \smallskip
  \parbox{0.9\textwidth}{\centering\footnotesize A close-up of the graph shown
  earlier in this preface: the wiki around the stochastic wave equation, with this
  book in the centre. Each dot is one page of the wiki, a Markdown file that
  the models write and maintain from what they have read: a digest of a paper,
  or a page on a concept or on open problems. Lines are links between pages.
  Shown are 268 pages and the 2,230 links among them, on three rings: the
  pages linked with this book's project page, the pages linked with those, and
  the rest of this selection. Drawn from the wiki's link index, November~16,
  2026; see \path{media/spdes-wiki-graph/}.\par}
\end{figure}

\clearpage
\section*{An experiment in doing mathematics}

\begin{wrapfigure}{r}{0.29\textwidth}
  \centering
  \vspace{-\baselineskip}
  \includegraphics[width=0.27\textwidth]{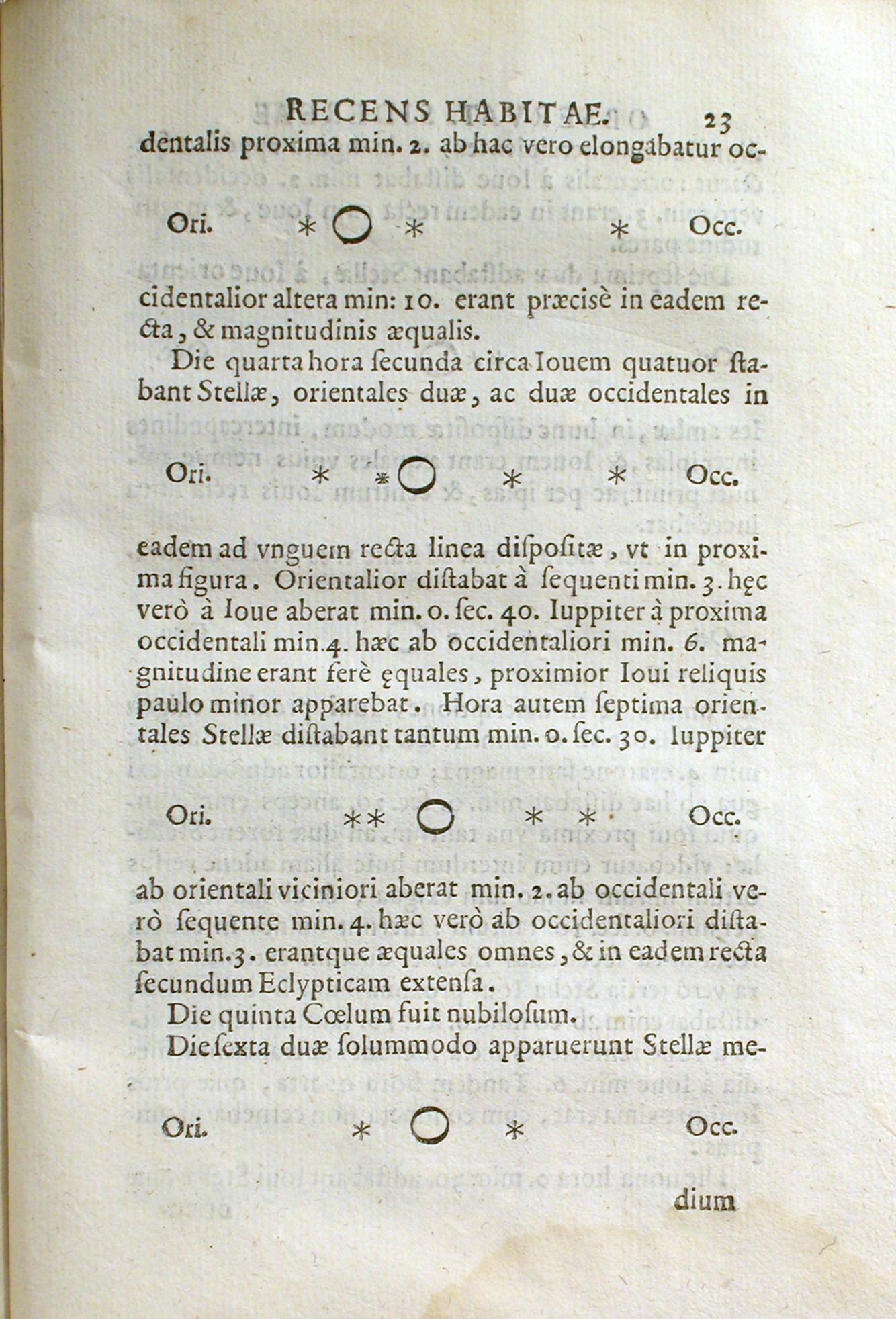}

  {\footnotesize Galileo, \emph{Sidereus Nuncius} (1610): Jupiter (the circle)
  and its moons (the stars) on several nights.\par}

  \medskip
  \includegraphics[width=0.27\textwidth]{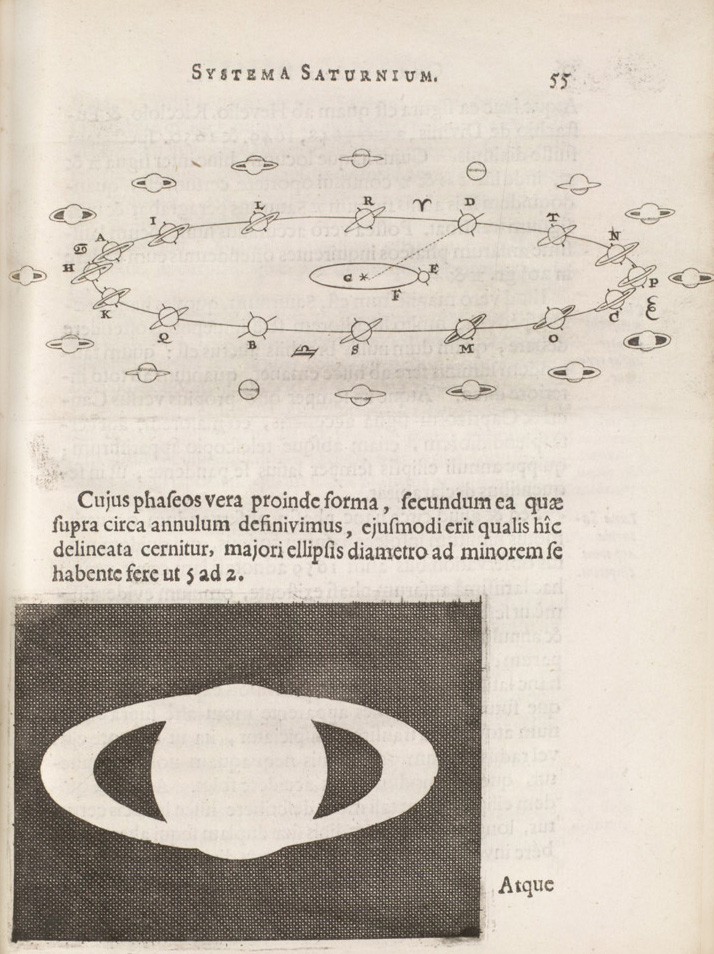}

  {\footnotesize Huygens, \emph{Systema Saturnium} (1659): Saturn's ring as
  seen from Earth around its orbit.\par}

  \smallskip {\scriptsize Scans: University of Oklahoma Libraries; Smithsonian
  Libraries. Public domain; see \path{media/historical-astronomy/}.\par}
\end{wrapfigure} In 1610 Galileo pointed a telescope that magnified about twenty
times at the sky. Within about a week he saw four moons moving around Jupiter,
which no one had seen before~\cite{galilei:89:sidereus}. A few months later he
looked at Saturn and reported a planet made of three bodies. Two years later the
two side bodies disappeared, and he could not explain why. What he saw was real,
but its meaning was harder to understand. Christiaan Huygens later explained
these changing appearances by a ring around Saturn; he published the explanation
in 1659~\cite{vanhelden:95:saturn}. An AI system is a new instrument of this
kind. It lets us see arguments that we could not see before, but people still
have to check what it shows and decide what it means. A telescope finds things;
the astronomer still has to check them, name them, put them in order and explain
them. Mourrat~\hyperlink{pp:P5}{[P5]} suggests that mathematicians will become
naturalists in this sense.

This book is also an experiment in how such mathematics can be checked,
shared and understood. It has four parts: opening the arguments to community
review, preserving the development and the failed attempts, distinguishing what
has been checked, and asking what general ideas can be extracted from the work.
The development, the writing and much of the checking were done with AI. We also
use AI in the search for an explanation, as described below. Human understanding
and public discussion remain aims of the experiment. The manuscript and its
filtered development record are on GitHub. Each result has a row in a public
table, and readers can check a proof, dispute it, repair it or explain it
through a pull request, or extend the Lean verification and the tools behind the
book. Every merged contribution is acknowledged by name, and a substantial
contribution, to the mathematics, the Lean verification or the tools, leads to
coauthorship of the next version.
Help with \mbox{de-slopping} is especially welcome: turning AI-written
arguments into ones a person can read and trust. And whether it is the
mathematics or the methodology, we would be happy to see this work inspire your
own. This is
version 1.0.0 of the book. Later versions follow one rule, written in the
repository: the last number changes for corrections of wording and form, the
middle number for changes to the mathematics that leave the main theorems as
stated, and the first number for a new version on arXiv, which comes with any
change to a main theorem, a correction that affects one, or a substantially
simpler proof. Similar proposals have
appeared on the blog \emph{Proofs and Prompts}: an open platform where
AI-generated results gradually acquire ``collective mathematical confidence''
on a ``Git Hub'' of proofs~\hyperlink{pp:P4}{[P4]}, open work groups that
explain important preprints in a Git repository~\hyperlink{pp:P11}{[P11]}, a
community-owned repository of AI proof reports~\hyperlink{pp:P12}{[P12]}, teams
of researchers who work to understand AI-generated
proofs~\hyperlink{pp:P14}{[P14]}, and credit for works at the level of a
doctoral thesis, defended in public with open reviews~\hyperlink{pp:P10}{[P10]}.

We also treat the development record as data for an experiment in doing
mathematics with AI. The public repository preserves a filtered history:
proof revisions, failed approaches with their reasons, Lean records, review
summaries, selected original reports, and marked copies whose output format
was repaired. It does not contain every reviewer input or response. Readers
cannot reconstruct every earlier state from the public repository alone.

Withheld material includes other authors' texts and images, private model
reasoning, system and account details, correspondence and planning notes,
collaborators' unpublished work, assessments of other authors' manuscripts
awaiting release review, and evidence archives. The
\href{https://github.com/chenle02/stratonovich-ham-second-moment/blob/v1.0.0/REDACTIONS.md}{redaction record}
gives the reasons. The original history and review records are retained
privately. The released material lets readers trace the construction,
criticism and repair it records, and study where the process failed.
Several posts on
\emph{Proofs and Prompts} argue in the same direction: Kretschmer calls a
detailed disclosure of AI use ``a worthy scientific
unit''~\hyperlink{pp:P7}{[P7]}, Kaipel proposes that every interaction with AI
that influenced a manuscript be recorded in full and made
available~\hyperlink{pp:P8}{[P8]}, Lewin~\hyperlink{pp:P12}{[P12]} compares AI
proofs with experimental data released with documentation for others to
analyze, and Gonen Cohen~\hyperlink{pp:P2}{[P2]} discusses releasing the full
prompting history to improve reproducibility.

Review arguments themselves require checking: No-go~\ref{P:nogo-33}
records three incorrect steps proposed in reviews and the calculations
that exposed them.

We also treat an approach that does not work as a result.
Chapter~\ref{S:book-no-go} gives numbered statements and proofs for the
mathematical obstructions, explaining why each approach was attractive,
where it fails and what it does not rule out. Readers can use these results
without repeating the failed attempts. This echoes
Morris's call for more accounts of failed attempts~\hyperlink{pp:P9}{[P9]},
Duminil-Copin's account of how failed attempts on one conjecture brought him
collaborators and ideas that he later used elsewhere~\hyperlink{pp:P3}{[P3]},
and Tao's remark that negative results are rarely reported, so that what one
hears about AI in mathematics is skewed towards success~\cite{tao:26:post}.

For the same reason, the book keeps its claims separate. Proved results, imported
statements, conditional consequences and open questions are distinguished
throughout. Lerario makes the same point about mathematics in general: it should
always be clear how far a result has been checked, and an argument checked step
by step by a proof assistant such as Lean should not have the same standing as a
heuristic analogy~\hyperlink{pp:P13}{[P13]}. A review by a language model is a
check, not a proof, and Lean checks only part of the book; two of the four main theorems are outside it. At the time of writing, the technical chapters have not yet been
checked line by line by a human reader; the public table records how far each
result has been checked.

\section*{What can we learn from the proof?}

A long argument produced with AI leaves another task: to understand which
ideas make it work and whether they can be used elsewhere. This book offers a
small experiment in that direction. We ask Amy to return to the arguments she
helped produce, including failed approaches, and to look for recurring
constructions, indispensable assumptions and possible general principles. We use
her proposals to guide further questions and checks. The aim is to involve AI in
the search for an explanation as well as in the production of an argument. A
proposed principle should explain why several steps belong together and help us
treat a second problem with less repeated work; whether it does is for people
to judge. Madore~\hyperlink{pp:P6}{[P6]} describes open problems as milestones
toward understanding and new theories, and argues that understanding is a human
activity that computers cannot perform for us.
Lerario~\hyperlink{pp:P13}{[P13]}, in ``Qualitative growth'' and ``Continuity
of research,'' discusses organizing separate results through common principles
and allows that theories and research directions could originate in artificial
systems.

Section~\ref{S:theory-outlook} records the first candidates from this exercise.
It separates the classical tools, the structure visible in this proof and the
extensions that remain to be tested. The length of the computation does not
establish that a new theory has emerged. The useful test is what the proposed
principle explains, what hypotheses it needs and what it lets us do next.

A further aim is a much shorter proof. The argument in this book is long
because it was built step by step: each obstruction was met by a new estimate,
and each estimate was kept once it had been checked. Once the ideas that make
the proof work are understood, many of these steps should merge or disappear.
Several estimates may follow from one principle, some detours may turn out to
be unnecessary, and the sharp constant may need only part of the present
machinery. We expect this simplification to come after, and through, human
understanding of the argument, with AI used to test each proposed shortcut
against the full proof. No reduction in length is claimed in this version.
Hairer makes a similar point: in digesting an argument produced by AI, one may
well realize that it can be presented in a clearer and simpler
way~\hyperlink{pp:P1}{[P1]}.

\section*{How to read this book}

The reader should have a background in graduate probability and analysis.
Familiarity with Gaussian processes, Fourier analysis, Sobolev spaces and
bounded operators will be useful, and the later chapters also use spectral
theory and Young integration. A separate reading edition gives a short account
of the model, the main statements and the proof strategy, with links to the
detailed arguments here. Readers who want only the upper bound at the correct
scale can read Chapters~\ref{S:book-setting}--\ref{S:G0} and then
Chapter~\ref{S:book-independent-upper}. The proof of the sharp constant occupies
Chapters~\ref{S:book-gaussian-foundations}--\ref{S:book-constant}. The guide to
the proof after these front pages maps the dependencies, and each chapter opens
with its purpose and the earlier results it needs. The appendices are not
prerequisites for the main theorems. They keep low-order calculations, shorter
estimates under additional assumptions, covariance comparisons and an
obstruction to a proposed clock representation, which show what a successful
argument must retain. Chapter~\ref{S:book-no-go} gathers these and other
obstructions by mathematical cause, with links back to the successful
arguments. The final three appendices give the AI work in numbers
(Appendix~\ref{S:appendix-ai-statistics}), what Lean checks
(Appendix~\ref{S:appendix-lean-checks}) and how far each result has been
checked (Appendix~\ref{S:appendix-check-status}).

\medskip
\noindent\textbf{Where to find everything.} The source and development record
are in the book's repository at release 1.0.0. The verification ledger there
records, for each numbered result, how far it has
been checked and how a reader can check it. The Lean companion, at release
0.10.0 (DOI~\href{https://doi.org/10.5281/zenodo.24295453}{10.5281/zenodo.24295453}), is described in
Appendix~\ref{S:appendix-lean-checks}.
\begin{center}\footnotesize
\href{https://github.com/chenle02/stratonovich-ham-second-moment/tree/v1.0.0}{\nolinkurl{github.com/chenle02/stratonovich-ham-second-moment/tree/v1.0.0}}\\
\href{https://github.com/chenle02/stratonovich-ham-second-moment/blob/v1.0.0/VERIFICATION.md}{\nolinkurl{github.com/chenle02/stratonovich-ham-second-moment/blob/v1.0.0/VERIFICATION.md}}\\
\href{https://github.com/chenle02/stratonovich-ham-second-moment-lean/tree/v0.10.0}{\nolinkurl{github.com/chenle02/stratonovich-ham-second-moment-lean/tree/v0.10.0}}
\end{center}

\section*{Acknowledgements}

We thank Robert Dalang, Aur\'elien Deya and Samy Tindel for their comments on
an early unpublished version of this work, and Raluca Balan, Yaozhong Hu, Jian
Song and Xiong Wang for their interest and correspondence.

The author was supported in part by the U.S. National Science Foundation under
Grant \mbox{DMS-2246850} and CAREER Grant \mbox{DMS-2443823}, and by Simons
Foundation Collaboration Grant for Mathematicians 959981. We thank OpenAI for
providing access to its AI tools through the ChatGPT for Academic Researchers
program.

\medskip
\begin{flushright}
Le Chen, with Amy\\
Auburn, Alabama\\
November 2026
\end{flushright}

% Keep the whole list together: start it on the next page if fewer than
% about 20 lines remain.
\Needspace{20\baselineskip}
\section*{Posts from \emph{Proofs and Prompts}}

In the text these posts are marked [P1]--[P14]. They are listed here in order of
date, not in the bibliography, and each title is a link to the post. They
include proposals we try to put into practice and criticisms we seek to address;
citing a post does not imply that its author endorses this experiment.
\begin{itemize}
  \item[\hypertarget{pp:P1}{[P1]}] Martin Hairer,
    \href{https://proofsandprompts.com/2026/08/07/writing-mathematics-in-the-age-of-ai/}{\emph{Writing mathematics in the age of AI}},
    August~7, 2026.
  \item[\hypertarget{pp:P2}{[P2]}] Segev Gonen Cohen,
    \href{https://proofsandprompts.com/2026/08/11/the-question-of-reasoning-traces/}{\emph{The question of reasoning traces}},
    August~11, 2026.
  \item[\hypertarget{pp:P3}{[P3]}] Hugo Duminil-Copin,
    \href{https://proofsandprompts.com/2026/08/30/care-for-a-little-more-ai/}{\emph{Care for a little more AI?}},
    August~30, 2026.
  \item[\hypertarget{pp:P4}{[P4]}] Marco Trombetti,
    \href{https://proofsandprompts.com/2026/09/22/proofforum-keeping-ai-generated-mathematics-human/}{\emph{ProofForum: Keeping AI-generated mathematics human}},
    September~22, 2026.
  \item[\hypertarget{pp:P5}{[P5]}] Jean-Christophe Mourrat,
    \href{https://proofsandprompts.com/2026/10/08/why-i-chose-to-do-math-and-how-im-grappling-with-ai-advances/}{\emph{Why I chose to do math, and how I'm grappling with AI advances}},
    October~8, 2026.
  \item[\hypertarget{pp:P6}{[P6]}] David A. Madore,
    \href{https://proofsandprompts.com/2026/10/20/letter-to-a-worried-student/}{\emph{Letter to a worried student}},
    October~20, 2026.
  \item[\hypertarget{pp:P7}{[P7]}] William Kretschmer,
    \href{https://proofsandprompts.com/2026/10/27/on-ai-disclosures/}{\emph{On AI disclosures}},
    October~27, 2026.
  \item[\hypertarget{pp:P8}{[P8]}] Maximilian Kaipel,
    \href{https://proofsandprompts.com/2026/10/28/why-and-how-to-transparently-disclose-ai-use-in-mathematics-research/}{\emph{Why and how to transparently disclose AI use in mathematics research}},
    October~28, 2026.
  \item[\hypertarget{pp:P9}{[P9]}] Patrick Morris,
    \href{https://proofsandprompts.com/2026/10/30/lets-stop-writing-papers/}{\emph{Let's stop writing papers\ldots}},
    October~30, 2026.
  \item[\hypertarget{pp:P10}{[P10]}] Andrew Ng,
    \href{https://proofsandprompts.com/2026/11/01/a-new-system-of-accreditation/}{\emph{A new system of accreditation}},
    November~1, 2026.
  \item[\hypertarget{pp:P11}{[P11]}] Alexander Van Werde,
    \href{https://proofsandprompts.com/2026/11/01/collaborative-work-groups-on-exposition-problems/}{\emph{Collaborative work groups on exposition problems}},
    November~1, 2026.
  \item[\hypertarget{pp:P12}{[P12]}] Mathieu Lewin,
    \href{https://proofsandprompts.com/2026/11/02/please-stop-writing-papers-there-is-another-way/}{\emph{Please STOP writing papers! (there is another way)}},
    November~2, 2026.
  \item[\hypertarget{pp:P13}{[P13]}] Antonio Lerario,
    \href{https://proofsandprompts.com/2026/11/09/the-shape-of-mathematical-knowledge/}{\emph{The shape of mathematical knowledge}},
    November~9, 2026.
  \item[\hypertarget{pp:P14}{[P14]}] Jonathan Noel,
    \href{https://proofsandprompts.com/2026/11/10/im-giving-my-theorems-away/}{\emph{I'm giving my theorems away}},
    November~10, 2026.
\end{itemize}

\begin{abstract}
  Motivated by Xia Chen's moment
  conjecture~\cite{chen:author2026:time-dependency}*{(1.22)}, we determine the
  exact second-moment asymptotics through an intrinsic variational formula.
  We study the Stratonovich hyperbolic Anderson model with temporal covariance
  $|s-t|^{-\alpha_0}$ and an even, nonnegative, positive-definite tempered
  Radon spatial covariance $\Gamma$ homogeneous of order $-\alpha$. We assume
  $d\in\{1,2,3\}$, $0<\alpha_0<1$, $0<\alpha<d$ or $\alpha=d=1$, and
  $\alpha+\alpha_0<2$. For initial displacement one and velocity zero, our
  main result, Theorem~\ref{T:intro-second-moment}, is the exact limit
  \[
    \lim_{T\to\infty}T^{-\kappa}\log\E u^2(T,x)
    =C_2(\Gamma,\alpha_0),
    \qquad \kappa=\frac{4-\alpha-\alpha_0}{3-\alpha},
  \]
  for every fixed $x\in\R^d$, with
  \[
    C_2(\Gamma,\alpha_0)
    =\sup_{\substack{\Psi\in C_c^\infty((0,1)\times\R^d;\R)\\
    \Psi\ge0,\ \Psi\ne0}}
    \left\{2\int_0^1\sqrt{\mathcal{R}_\Psi(s)}\ud s
    -\frac{v_\Psi}{2}\right\}.
  \]
  Here the covariance mean, variance and Rayleigh value are
  \[
    \begin{aligned}
      h_{\Psi,s}(x)
      &=\int_0^1|s-r|^{-\alpha_0}(\Gamma*\Psi_r)(x)\ud r,\\
      v_\Psi&=\int_0^1\int_{\R^d}\Psi(s,x)h_{\Psi,s}(x)\ud x\ud s,\\
      \mathcal{R}_\Psi(s)
      &=\sup_{\substack{f\in H^1(\R^d;\R)\\\|f\|_2=1}}
      \left\{\int_{\R^d}h_{\Psi,s}(x)f(x)^2\ud x
      -\|\nabla f\|_2^2\right\}.
    \end{aligned}
  \]
  The spatial norms are $L^2$ norms, and $\Psi_r=\Psi(r,\cdot)$. The
  variational class includes controls that do not separate in time and space.
  The formula retains the original covariance amplitude, kinetic coefficient
  one and Gaussian cost $v_\Psi/2$. Its value is positive and finite when
  $\Gamma\ne0$; for $\Gamma=0$, the solution is one and the constant is zero.

  An independent upper bound establishes finiteness at the scale $T^\kappa$.
  The identification of the constant combines Gaussian phase operators,
  convergence from the regularized noise, estimates uniform in the spatial
  localization radius, and finite optimization with nonnegative Gram entries.
  We also obtain bounds for the same constant through spatial energies and
  occupation measures, including explicit bounds for spatial white noise in
  one dimension.
\end{abstract}

\tableofcontents
\listoffigures

\chapter*{A guide to the proof}
\addcontentsline{toc}{chapter}{A guide to the proof}
\markboth{A guide to the proof}{A guide to the proof}

The reading edition states the main results and explains the route through
this volume. The chapters below supply the proofs in their dependency order.
The appendices retain the direct estimates, finite-order calculations and
alternative approaches, together with their original hypotheses and limitations.
They are useful comparisons; they are not all inputs to the main argument. Each
chapter opens with its motivating question, the earlier definitions or results
needed to enter it, and an outline of the argument. Readers may use these
openings to choose a starting point and follow the linked prerequisites when
needed.

How much to read depends on the time at hand.

\begin{center}
\begin{tabular}{@{}p{0.20\textwidth}p{0.71\textwidth}@{}}
  \toprule
  Time & What to read \\
  \midrule
  Ten minutes & The two main theorems in Section~\ref{S:targets}, the reason for the exponent in Section~\ref{S:why-kappa}, and the map in Figure~\ref{Fig:chapter-roadmap}. \\
  One evening & The reading edition, or Chapter~\ref{S:book-setting} followed by the opening of each later chapter. The first two videos cover Chapter~\ref{S:book-setting}. \\
  The full proof & Chapters~\ref{S:G0}--\ref{S:book-constant} in order. Each opening says which results to keep and which proofs can wait on a first reading. \\
  \bottomrule
\end{tabular}
\end{center}

Lemma~\ref{L:white-test} gives concrete variance and mean calculations.
No-go~\ref{P:nogo-05} shows geometrically why merging two chains can violate
wave support.

\begin{center}
\begin{tabular}{@{}p{0.20\textwidth}p{0.71\textwidth}@{}}
  \toprule
  Chapters & Purpose \\
  \midrule
  \ref{S:book-setting}--\ref{S:G0} & Model, literature and construction of one solution. \\
  \ref{S:book-independent-upper} & A finite upper bound at the correct time scale. \\
  \ref{S:book-gaussian-foundations}--\ref{S:book-wave-evolution} & Covariance estimates, operator paths and wave evolution. \\
  \ref{S:book-lower} & Exact scaling and the lower bound for fixed controls. \\
  \ref{S:book-means}--\ref{S:book-growing-radius} & Spatial means and estimates uniform in the localization radius. \\
  \ref{S:book-point-upper}--\ref{S:book-constant} & The point upper bound and the variational constant. \\
  \ref{S:book-no-go} & Obstructions and failed approaches: every mathematical no-go, and what it does not rule out. \\
  \midrule
  \ref{S:appendix-finite-order}--\ref{S:appendix-static-clock} & Other estimates, consistency checks and alternative approaches. \\
  \ref{S:appendix-ai-statistics} & The AI work in numbers: who did the work, the reviews, and how many approaches failed. \\
  \ref{S:appendix-lean-checks} & What the Lean companion checks, result by result. \\
  \ref{S:appendix-check-status} & How far each result has been checked: by AI review, by Lean and by people. \\
  \bottomrule
\end{tabular}
\end{center}

Figure~\ref{Fig:chapter-roadmap} gives a map of all twelve chapters and the
four mathematical appendices. The reader links to the labelled statements in this companion. The numbering
includes the chapter and section, so a statement can be located
without searching the entire volume. The figures are explanatory schematics;
their captions identify the precise results and the limits of the drawings.

\paragraph{Video introduction.}
An introduction to this book in eleven videos, each about twenty minutes and
narrated by a synthetic voice, is on YouTube at \VideoPlaylist. We made the
series with Amy: she drafted the scripts and animations from this book, every
displayed formula is copied from the statement that carries its label, and each
script was reviewed against the manuscript by a model family different from the
one that drafted it. The table
lists what each episode covers; each episode number is a link, and the red
badges in Figure~\ref{Fig:chapter-roadmap} lead to the same videos.
Episodes~1--10 follow the draft at version \texttt{38d67ec} and episode~11
the draft at version \texttt{0f553e08}, so statement numbers quoted in them
may differ from this version. The videos explain the argument; they are not
part of it or of its review.

\begin{center}
\begin{tabular}{@{}p{0.07\textwidth}p{0.36\textwidth}p{0.48\textwidth}@{}}
  \toprule
  Ep. & Title & Covers \\
  \midrule
  \VideoRef{1}{ktqRCapNP2c} & A guitar in a sandstorm & Chapter~\ref{S:book-setting}: the model, its history, the exponent $\kappa$. \\
  \VideoRef{2}{YTP_G_anCaE} & The answer: a growth rate and a variational constant & Chapter~\ref{S:book-setting}: the two main theorems and the constant. \\
  \VideoRef{3}{Glay6Q_jsLM} & Two chains and their pairings & Chapter~\ref{S:G0}. \\
  \VideoRef{4}{bUyxCpboios} & The correct time scale & Chapter~\ref{S:book-independent-upper}. \\
  \VideoRef{5}{2jv8nuqgfBY} & A Gaussian toolkit for wave fields & Chapters~\ref{S:book-gaussian-foundations}--\ref{S:book-wave-evolution}. \\
  \VideoRef{6}{VWemsluFZiA} & The lower bound: steering the noise & Chapters~\ref{S:book-lower}--\ref{S:book-means}. \\
  \VideoRef{7}{Z35qRxt05kA} & The upper bound: no control does better & Chapters~\ref{S:book-growing-radius}--\ref{S:book-point-upper}. \\
  \VideoRef{8}{RllmgM_zB48} & The constant and what it tells us & Chapter~\ref{S:book-constant} and the open questions of Chapter~\ref{S:book-setting}. \\
  \VideoRef{9}{eLz9Hd20ZT0} & Roads not taken & Appendices~\ref{S:appendix-finite-order}--\ref{S:appendix-static-clock}. \\
  \VideoRef{10}{ZwJM8BL3KSQ} & How this book was written & The preface and Appendices~\ref{S:appendix-ai-statistics}--\ref{S:appendix-check-status}. \\
  \VideoRef{11}{QzPhuFUUDT0} & Approaches that do not work & Chapter~\ref{S:book-no-go}. \\
  \bottomrule
\end{tabular}
\end{center}

\mainmatter

\chapter{Introduction}\label{S:book-setting}

A wave is driven by a random force that remembers its past. How fast can
its second moment grow? The answer depends on the covariance and on what the
stochastic product means. This chapter introduces the Stratonovich model and
states the answer in two parts: the scale of the growth, and the constant
at that scale. Both parts hold in space dimensions one, two and three.

We begin with a guitar in a sandstorm (Figure~\ref{Fig:walsh-guitar})
and a short history. Section~\ref{S:model} states the model, the
covariance and the initial data. Section~\ref{S:known} records what was
known, including the difference between the author and arXiv versions
of the source. Section~\ref{S:targets} states the two main theorems, and
Section~\ref{S:why-kappa} explains the scale. Figure~\ref{Fig:chapter-roadmap}
maps the rest of the book.

No earlier chapter is needed. On a first reading, the literature
comparison in Section~\ref{S:literature-context} and the outlook
in Sections~\ref{S:future-work} and~\ref{S:theory-outlook} can
wait. A reader entering later should start with the model and
initial data in Section~\ref{S:model}, Hypothesis~\ref{H:homog} and
the exponent in~\eqref{E:kappa}. Episodes~\VideoRef{1}{ktqRCapNP2c}
and~\VideoRef{2}{YTP_G_anCaE} of the video series follow this chapter. The
series retells the book in eleven short episodes, and each chapter guide
names its episode.

\Needspace{0.44\textheight}
\section{A guitar in a sandstorm: forty years since Walsh}\label{S:wave-history}

\input{notes/figures/walsh-guitar}

Imagine leaving a guitar outside during a sandstorm. Each grain gives a string
an unpredictable impulse. The string carries the disturbance away, while new
grains strike elsewhere. How can we describe the motion when the individual
impacts are too numerous to follow? Walsh used this picture to introduce a
stochastic wave equation in his Saint-Flour lectures, delivered in 1984 and
published in 1986~\cite{walsh:86:introduction}*{p.~282, (1.14)}. Forty years
after their publication, the picture still provides a useful way into the
subject. We take these lectures as a starting point for this account. Earlier
work on stochastic wave equations lies outside its scope.

\paragraph{1986: from random impacts to a mathematical equation.}
A physical guitar has finite strings and endpoints, but Walsh makes a deliberate
simplification in the worked example: the string is infinite and initially at
rest, so boundary conditions do not obscure the random forcing. Centering the
impacts and replacing their fluctuations by space-time white noise gives an
additively forced wave equation. The differential notation is suggestive, but
the solution need not possess the derivatives appearing in it. The integral
formulation gives the equation its meaning. We must first specify what we mean
by a solution and the space in which it is defined. The guitar is a motivation
for this model, not a derivation of every stochastic wave equation considered
below.

\paragraph{1999: how far does the integration theory reach?}
The one-dimensional string is only the beginning. In higher dimensions the wave
propagator may itself be a measure, and spatially correlated noise calls for a
wider integration framework. Dalang's extension of the martingale-measure
integral~\cite{dalang:99:extending} made it possible to treat SPDEs with more
singular Green kernels under spectral integrability conditions. The question now
includes the interaction of propagation, spatial dimension and noise covariance.
These conditions determine whether the equation admits a random-field solution.

\paragraph{The 2000s and 2010s: from existence to growth.}
Once a solution is available, we can ask how its fluctuations develop over long
times. Multiplicative noise changes the problem: the forcing depends on the
current solution, so the disturbance and its random environment interact.
Moments measure the size of the fluctuations and the contribution of rare, large
values. Dalang--Mueller--Tribe~\cite{dalang.mueller.ea:08:feynman-kac-type}
developed probabilistic moment representations using Poisson processes and
random flights. Dalang--Mueller~\cite{dalang.mueller:09:intermittency} obtained
intermittency estimates for a three-dimensional model under their spatial
covariance assumptions. These results illustrate how the wave equation has its
own probabilistic methods, while the form of the noise continues to matter.

For one-dimensional space-time white noise,
Chen--Dalang~\cite{chen.dalang:15:moment} developed moment bounds and explicit
second-moment formulas, including rough initial data. Temporally correlated
forcing brings another change: random impacts at distinct times are no longer
independent. Balan--Conus~\cite{balan.conus:16:intermittency} studied the
resulting Skorohod moment growth under their covariance assumptions. For
white-in-time noise, Balan--Song~\cite{balan.song:19:second} related hyperbolic
and parabolic second-moment Laplace transforms and obtained an actual second
Lyapunov limit in their Riesz case. Construction, bounds of the correct growth
order and exact limits are successive questions, but a theorem answering one
need not answer the others.

\paragraph{The 2020s: exact rates and the meaning of multiplication.}
The same formal multiplicative equation can describe different stochastic
models. Skorohod integration uses Wick products; ordinary Stratonovich
approximation retains the self-contractions generated by multiplication of
smooth Gaussian fields. This distinction is central when the noise is static or
correlated in time. For time-independent Skorohod noise, Balan--Le Chen--Xia
Chen~\cite{balan.chen.ea:22:exact} obtained exact moment asymptotics.
Chen--Eisenberg~\cite{chen.eisenberg:23:interpolating} treated a related
fractional differential family, with time order below two and static noise.
Chen--Guo--Song~\cite{chen.guo.ea:24:moments} obtained exact second moments for
a white-noise fractional SPDE family and matching higher-moment growth orders in
its wave specialization. The latter are growth orders, not sharp higher-moment
constants.

On the Stratonovich side, Balan~\cite{balan:22:stratonovich} constructed
static-noise solutions and compared chaos expansions, while Xia
Chen--Hu~\cite{chen.hu:25:hyperbolic} established static homogeneous moment
asymptotics. For fractional temporal correlation, Xia Chen's author
version~\cite{chen:author2026:time-dependency} gives two-sided first-moment
bounds and a coarse second-moment upper bound. The source editions and all
model, data and normalization restrictions are detailed in
Section~\ref{S:literature-context}.

\paragraph{The present question.}
These works lead to three questions. How do we construct a solution? At what
rate do its moments grow? Can we identify the leading constant? Our problem is
the second moment of a multiplicative Stratonovich wave equation with fractional
temporal correlation and homogeneous spatial covariance. We give an upper bound
at the correct scale and a variational expression for the second-moment limit.
We do not obtain a numerical closed form or a sharp formula for every moment
order. The scientific review status is recorded separately. The following
sections specify the model and explain why its self-pairings and correlated
forcing require more than a transfer of the Skorohod results.

\section{Model and standing assumptions}\label{S:model}

We study the hyperbolic Anderson model
\cite{chen:author2026:time-dependency}*{(1.1)}
\begin{equation}\label{E:HAM}
  \frac{\partial^2 u}{\partial t^2}(t,x)=\Delta u(t,x)+\W(t,x)u(t,x),
  \qquad u(0,x)=1,\quad \frac{\partial u}{\partial t}(0,x)=0,
\end{equation}
for $(t,x) \in \R_+\times\R^d$, $d \in \{1,2,3\}$, interpreted in the
Stratonovich sense (limit of classical wave equations driven by smoothed noises
$\W_{\epsilon,\delta}$). The noise is centered Gaussian with
\begin{equation}\label{E:cov}
  \operatorname{Cov}\big(\W(s,x),\W(t,y)\big)=|s-t|^{-\alpha_0}\gamma(x-y),
  \qquad 0<\alpha_0<1,
\end{equation}
where $\gamma \ge 0$ is nonnegative definite with spectral measure $\mu$.

\begin{hypothesis}[Homogeneous covariance]\label{H:homog}
  Let $d \in \{1,2,3\}$ and $0 < \alpha_0 < 1$. The spatial covariance is a
  nonnegative, even, positive-definite tempered Radon measure $\Gamma$ on
  $\R^d$, with positive even tempered Radon spectral measure $\mu$ in the
  convention
  \[
    \Gamma=\int e^{i\xi\cdot x}\mu(\ud\xi)
    \quad\text{in the sense of tempered distributions}.
  \]
  The symbol $\gamma$ in~\eqref{E:cov} denotes this covariance as a function
  when one exists, and as a distribution otherwise. For every Borel set $A$ and
  $c > 0$, $\Gamma(cA) = c^{d-\alpha}\Gamma(A)$, and
  \begin{equation}\label{E:threshold}
    \alpha_0+\alpha<2,
  \end{equation}
  with $0 < \alpha < d$, or $\alpha = d = 1$. In particular, the last case
  includes spatial white noise $\Gamma = \delta_0$.
\end{hypothesis}

\paragraph{Zero covariance.}
Hypothesis~\ref{H:homog} allows $\Gamma = 0$ by design. In this case the
centered Gaussian noise is $\W = 0$, and the solution with initial data
$(1,0)$ is $u(t,x) = 1$. Thus $\E|u(t,x)|^p = 1$ for every $p > 0$, and
all normalized logarithmic moment rates are zero. Results asserting strict
positivity of covariance-dependent quantities impose $\Gamma \ne 0$ locally.

Throughout we write
\begin{equation}\label{E:kappa}
  \LeanEqMarker{E:kappa}
  \kappa:=\frac{4-\alpha-\alpha_0}{3-\alpha}\in(1,2).
\end{equation}

\paragraph{Contribution of the present work.}
We address the second moment of the original Stratonovich solution with
fractional temporal covariance throughout Hypothesis~\ref{H:homog}. The
independent upper bound reaches the scale $t^\kappa$, rather than only the
coarse quadratic exponential scale in the source. The main identification gives
an intrinsic variational expression for the limiting constant over all
nonnegative smooth compact spacetime controls, without a separability
restriction. It retains the covariance amplitude, kinetic coefficient one and a
Gaussian cost equal to one half of the variance, and includes the zero
covariance case. The upper and lower comparisons refer to the same canonical
solution. The occupation and spatial-energy estimates then give bounds on this
constant; evaluating the white spatial energy is not a numerical evaluation of
the full moment supremum. The earlier results discussed below concern different
stochastic interpretations or temporal covariances and do not supply this
identification by substitution.

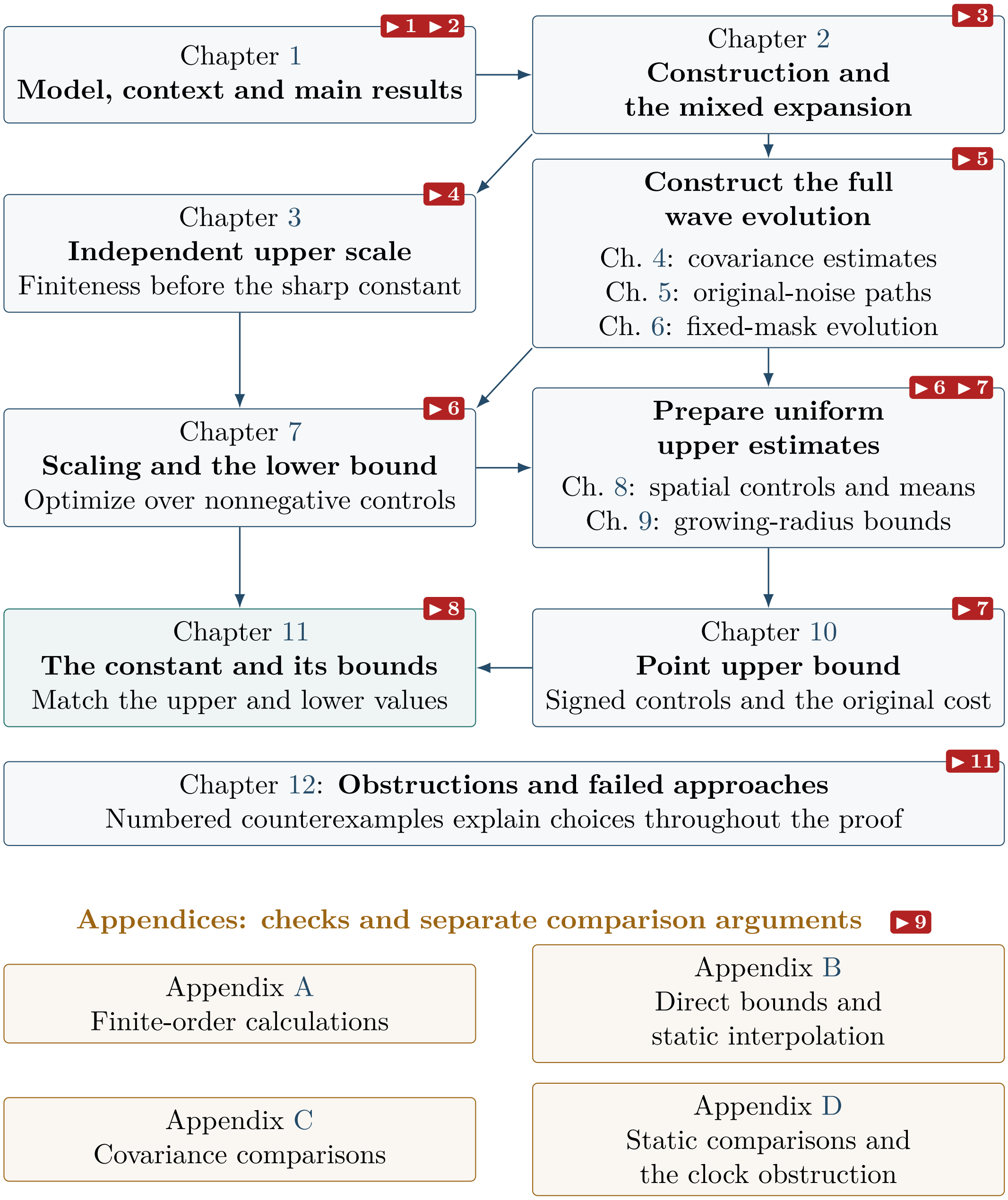
\begin{figure}
  \centering
  \begin{tikzpicture}[x=1cm,y=1cm,
    route/.style={ham/box,text width=6.25cm,minimum height=1.35cm},
    synthesis/.style={route,draw=hamteal,fill=hamteal!7},
    auxiliary/.style={route,draw=hamgold,fill=hamgold!5,
    minimum height=1.1cm}]
    \node[route] (setting) at (0,0)
    {\GuideRef{Chapter}{S:book-setting}\\
    \textbf{Model, context and main results}};
    \node[route] (construction) at (7.4,0)
    {\GuideRef{Chapter}{S:G0}\\
    \textbf{Construction and the mixed expansion}};
    \node[route] (upper) at (0,-2.5)
    {\GuideRef{Chapter}{S:book-independent-upper}\\
    \textbf{Independent upper scale}\\
    Finiteness before the sharp constant};
    \node[route,minimum height=2.1cm] (phase) at (7.4,-2.5)
    {\textbf{Construct the full wave evolution}\\[3pt]
    \GuideRef{Ch.}{S:book-gaussian-foundations}: covariance estimates\\
    \GuideRef{Ch.}{S:book-operator-paths}: original-noise paths\\
    \GuideRef{Ch.}{S:book-wave-evolution}: fixed-mask evolution};
    \node[route] (lower) at (0,-5.5)
    {\GuideRef{Chapter}{S:book-lower}\\
    \textbf{Scaling and the lower bound}\\
    Optimize over nonnegative controls};
    \node[route,minimum height=1.65cm] (means) at (7.4,-5.5)
    {\textbf{Prepare uniform upper estimates}\\[3pt]
    \GuideRef{Ch.}{S:book-means}: spatial controls and means\\
    \GuideRef{Ch.}{S:book-growing-radius}: growing-radius bounds};
    \node[synthesis] (constant) at (0,-8.3)
    {\GuideRef{Chapter}{S:book-constant}\\
    \textbf{The constant and its bounds}\\
    Match the upper and lower values};
    \node[route] (point) at (7.4,-8.3)
    {\GuideRef{Chapter}{S:book-point-upper}\\
    \textbf{Point upper bound}\\
    Signed controls and the original cost};
    \draw[ham/arrow] (setting)--(construction);
    \draw[ham/arrow] (construction.south west)--(upper.north east);
    \draw[ham/arrow] (construction)--(phase);
    \draw[ham/arrow] (upper)--(lower);
    \draw[ham/arrow] (phase.south west)--(lower.north east);
    \draw[ham/arrow] (phase)--(means);
    \draw[ham/arrow] (lower)--(means);
    \draw[ham/arrow] (lower)--(constant);
    \draw[ham/arrow] (means)--(point);
    \draw[ham/arrow] (point)--(constant);
    \node[route,text width=13.65cm,minimum height=1.1cm] (nogo) at (3.7,-10.2)
    {\GuideRef{Chapter}{S:book-no-go}: \textbf{Obstructions and failed approaches}\\
    Numbered counterexamples explain choices throughout the proof};
    % Video badges: each links to the episode that covers the box.
    \node[ham/video] at ([xshift=-0.75cm]setting.north east)
    {\VideoEp{1}{ktqRCapNP2c}\enspace\VideoEp{2}{YTP_G_anCaE}};
    \node[ham/video] at ([xshift=-0.45cm]construction.north east) {\VideoEp{3}{Glay6Q_jsLM}};
    \node[ham/video] at ([xshift=-0.45cm]upper.north east) {\VideoEp{4}{bUyxCpboios}};
    \node[ham/video] at ([xshift=-0.45cm]phase.north east) {\VideoEp{5}{2jv8nuqgfBY}};
    \node[ham/video] at ([xshift=-0.45cm]lower.north east) {\VideoEp{6}{VWemsluFZiA}};
    \node[ham/video] at ([xshift=-0.75cm]means.north east)
    {\VideoEp{6}{VWemsluFZiA}\enspace\VideoEp{7}{Z35qRxt05kA}};
    \node[ham/video] at ([xshift=-0.45cm]constant.north east) {\VideoEp{8}{RllmgM_zB48}};
    \node[ham/video] at ([xshift=-0.45cm]point.north east) {\VideoEp{7}{Z35qRxt05kA}};
    \node[ham/video] at ([xshift=-0.45cm]nogo.north east) {\VideoEp{11}{QzPhuFUUDT0}};
    \node[ham/note,text=hamgold,text width=13.8cm] at (3.7,-11.85)
    {\textbf{Appendices: checks and separate comparison arguments}\quad
    \tikz[baseline=(v.base)]\node[ham/video] (v) {\VideoEp{9}{eLz9Hd20ZT0}};};
    \node[auxiliary] (finitechecks) at (0,-13)
    {\GuideRef{Appendix}{S:appendix-finite-order}\\
    Finite-order calculations};
    \node[auxiliary] (direct) at (7.4,-13)
    {\GuideRef{Appendix}{S:appendix-direct}\\
    Direct bounds and static interpolation};
    \node[auxiliary] (covariance) at (0,-14.9)
    {\GuideRef{Appendix}{S:appendix-covariance}\\
    Covariance comparisons};
    \node[auxiliary] (clock) at (7.4,-14.9)
    {\GuideRef{Appendix}{S:appendix-static-clock}\\
    Static comparisons and the clock obstruction};
  \end{tikzpicture}
  \caption[Chapter and appendix roadmap]{A reading map for the whole book.
    Blue arrows select principal inputs, not every dependency. In particular,
    the independent upper scale supplies finiteness for the variational lower
    value, which is also used in the spatial estimates. The upper argument
    retains the earlier wave construction and the radius estimates; it does
    not follow from the lower bound alone. The green box joins the two moment
    bounds. Chapter~\ref{S:book-no-go} explains failed alternatives to these steps.
    Gold boxes identify auxiliary checks and comparisons with their
    own hypotheses. No arrow from an appendix into the main proof is asserted.
    Every chapter and appendix label is a link to its opening guide. Red badges
    link to the episodes of the video introduction that cover each box; the
    table in the guide to the proof lists them.}
  \label{Fig:chapter-roadmap}
\end{figure}

Figure~\ref{Fig:proof-map} summarizes the main dependencies and identifies
the results used at each step. Figure~\ref{Fig:chapter-roadmap} places these
steps within the chapters and distinguishes the auxiliary appendices.

\paragraph{Main result and dependencies.}
Theorem~\ref{T:intrinsic-second-moment} identifies the original point-moment
limit for every covariance in Hypothesis~\ref{H:homog}. Its full variational
domain and normalization are displayed in~\eqref{E:intrinsic-positive-formula}.
Theorem~\ref{T:intrinsic-occupation-ceiling} and
Corollary~\ref{C:intrinsic-spatial-bracket} give derived bounds;
Section~\ref{S:intrinsic-white-evaluation} evaluates the white spatial energy
and its resulting moment bounds. These bounds need not be attained at the moment
constant.

The independent earlier Theorem~\ref{T:full-upper} establishes finite upper
growth at scale $t^\kappa$ and precedes its uses proving variational finiteness.
It uses the fixed-order construction proved in Lemma~\ref{L:joint-fixed-order},
the Gaussian and complex phase estimates in
Sections~\ref{S:static-phase}--\ref{S:static-all-orders}, and the interval
coefficient and bounds for the finite family of time intervals in
Section~\ref{S:interval-upper}. It does not use the alternative covariance
interpolation arguments in
Sections~\ref{S:higher-alpha}--\ref{S:energy-fidelity} as inputs. Its uniform
geometric tail is proved here from scaled finite-model moments; the older
mixed-series and crude estimates in Chapter~\ref{S:G0} are not needed for that
moment-limit passage. The separate verification of the mild equation is proved
in Proposition~\ref{P:mild-construction}, using the corrected vertexwise bound
of Lemma~\ref{L:vertex-crude}. The alternative arguments are retained to record
their separate hypotheses and limitations. Section~\ref{S:d1-upper} also gives
explicit constants for $\alpha \le 1$; the full-class theorem asserts a finite
uniform constant, without giving a numerically optimized value. The later proof
constructs the operator and regularized noises preserving the original noise in
Sections~\ref{S:actual-operator-driver}--\ref{S:constant-phase-young}, then
proves the same-solution lower bound and positive-control variational
expressions in
Sections~\ref{S:actual-raw-block-lower}--\ref{S:unmasked-variational}. The
Gaussian upper argument culminates in
Sections~\ref{S:original-point-transfer}--\ref{S:original-signed-score-upper}.
Finite Gram optimization and its physical application then match the upper and
lower values in Section~\ref{S:intrinsic-second-moment}. This dependency order
does not use that final identification to prove its own finiteness.

\paragraph{Notation and source-dependent statements.}
The target spatial degree is $\alpha$, whereas $\chi$ denotes the degree of an
auxiliary static covariance in
Sections~\ref{S:static-holder}--\ref{S:energy-fidelity}. The coefficient power
is $\beta = 4-\alpha-\alpha_0$; the factorial power in
Section~\ref{S:interval-upper} is $\nu = 3-\alpha$. The local exponent $b =
2-\alpha$ in Section~\ref{S:d1-upper} is different from coordinate powers $b_j$
in Section~\ref{S:three-endpoints}. In Section~\ref{S:higher-alpha}, $\delta =
\alpha-1$ is an exponent; when used for temporal noise smoothing, $\delta$ is a
cutoff width. The covariance measure $\Gamma(\ud x)$ and the Euler Gamma
function $\Gamma(q)$ are distinguished by their arguments. The local $p =
2-\alpha_0$ in Section~\ref{S:interval-upper} is an interval-length decay
exponent; it is not the moment order $p$ in Conjecture~\ref{Conj:allp}. The
unconditional first-to-second moment comparison is Proposition~\ref{P:Jensen}.
Any numerical consequence involving the source's $L_*$ explicitly assumes a
valid first-moment lower bound with that constant. No correction of its printed
normalization is adopted. When a radius is given for topological support, the
closed ball is meant unless a strict open inclusion is specified. In particular,
the spatial cutoffs have support in $\overline{B}_{2R}$ and the wave Green
measures have support in $\overline{B}_t$, including the sphere measure in
dimension three. Plateau and causal-cone arguments retain their strict margins.

\paragraph{Which source statements are inputs.}
The moment estimate uses the bounded-operator Crouzeix--Palencia theorem
in~\eqref{E:cp-input}. The fixed-order construction is supplied here by
Lemmas~\ref{L:spatial-wick-measures}--\ref{L:joint-fixed-order}, using
Lemma~\ref{L:pairing-product-envelope}, which proves the finite product envelope
used here, and the direct norm bounds in Lemma~\ref{L:damped-homogeneous-norms}.
The Case~2 covariance step is proved independently in
Lemma~\ref{L:translated-covariance-CS}; the precise source-display
qualifications and regularized application are in
Remark~\ref{R:source-case2-time-shift}. Thus, neither Theorem~4.1 nor Lemma~4.4
nor Appendix~A.1 of the source is an analytic construction input. The
finite-width Gaussian field and multiple Stratonovich integral use its
definitions and regularization conventions, with Green smoothing removed before
noise smoothing.

Proposition~\ref{P:mild-construction} proves that the resulting series satisfies
its onefold Stratonovich mild equation, using the corrected vertexwise bound. It
does not import the author's Section~5(ii) argument or (5.3)--(5.4). No
uniqueness assertion is made. The uniform tail passage for the target-scale
moment estimate is proved in Lemma~\ref{L:uniform-series-tail}, independently of
the crude bound and the source's series estimate. The mixed-width and
order-limit passages for the mild equation are supplied here as well.

The first-moment formulas in Theorem~\ref{T:first} are quoted source statements,
with the printed-normalization qualification there; they are not premises of the
full upper theorem. Chen--Hu's static coefficient estimates and Chen's 2023
variational finiteness theorem are used only in the explicitly identified
auxiliary sections. In particular, our uses of the latter cite Theorem~1.1, not
Corollary~1.2. Thus, neither a disputed first-moment constant nor an auxiliary
static coefficient theorem enters the interval-field upper bound.

\section{Skorohod and Stratonovich wave equations}\label{S:literature-context}

The meaning of the noise product is part of the model. The same formal wave
operator and covariance do not specify the same equation under the Skorohod
and Stratonovich interpretations. We review both interpretations below,
separating existence results, moment-growth bounds and exact asymptotics.
All comparisons retain the hypotheses of the cited papers.

\subsection{The two products and their consequences}

The Skorohod integral is the divergence operator of Malliavin calculus. For the
linear equation with deterministic initial data, its iterated integrals form
orthogonal homogeneous Wiener chaoses. In the convergent expansion, every
positive-order term has mean zero, so the mean solution is the deterministic
free wave. Orthogonality reduces its second moment to a sum of squared chaos
norms. For adapted integrands and white-in-time noise, this interpretation
agrees with the It\^o integral in the corresponding integration framework;
see Walsh~\cite{walsh:86:introduction} and Dalang~\cite{dalang:99:extending}.
Temporal correlation removes the martingale structure, even when the
Skorohod integral remains available.

The Stratonovich construction considered here takes ordinary products of
regularized Gaussian noise and then removes the regularization. Wick's
formula retains contractions within an iterated term as well as contractions
between different terms. These self-pairings change the mean; terms of unequal
iteration orders can contribute to the second moment. This is why the
Skorohod second-moment isometry does not directly solve the present problem.
It does not mean that a Stratonovich solution has no orthogonal chaos
expansion: trace terms can be regrouped into Wiener chaoses, as in
Balan~\cite{balan:22:stratonovich}. The iteration order and the final Wiener
chaos order are different indices.

For regular noise, ordinary multiplication connects the equation to
deterministic wave evolution in a random potential. Stratonovich models are
therefore natural for specified smooth-noise approximations. For singular noise, convergence, the approximation scheme and any
required counterterms must be proved. There is no universal conversion formula
that identifies the two wave equations. In particular, an adapted
white-in-time Skorohod solution being an It\^o solution does not identify it
with an unspecified smooth-noise limit. The three-dimensional approximation
analysis of Delgado-Vences--Sanz-Sol\'e
\cite{delgado-vences.sanz-sole:14:approximation} shows why
regularization requires care.

\begin{figure}
\centering
\begin{tikzpicture}[x=1cm,y=1cm]
  \node[ham/box,minimum width=5.2cm] (hamSkoroProduct) at (0,0)
  {Skorohod / Wick product\\$u\diamond\dot W$};
  \node[ham/box,minimum width=5.2cm] (hamStratoProduct) at (6.5,0)
  {Ordinary Stratonovich construction\\regularize first, then multiply};
  \node[ham/box,minimum width=5.2cm] (hamSkoroIntegral) at (0,-1.8)
  {Divergence stochastic integral\\and Wiener-chaos expansion};
  \node[ham/box,minimum width=5.2cm] (hamStratoLimit) at (6.5,-1.8)
  {Classical random wave equations\\followed by the specified limit};
  \draw[ham/arrow] (hamSkoroProduct)--(hamSkoroIntegral);
  \draw[ham/arrow] (hamStratoProduct)--(hamStratoLimit);
  \node[ham/note] at (0,-3.6)
  {Adapted white-in-time integrands:\\
  Skorohod agrees with It\^o\\
  under the integration hypotheses.};
  \node[ham/note] at (6.5,-3.6)
  {Singular noise requires\\
  a specified approximation scheme\\
  and a convergence theorem.\\
  Correction terms may be needed.};
  \draw[hamblue!35] (-2.8,-5.0)--(9.3,-5.0);
  \node[ham/box,minimum width=10.9cm] at (3.25,-6.0)
  {The same formal wave operator and noise covariance\\
  do not identify the solutions or their moment constants.};
\end{tikzpicture}
\caption[Two meanings of the noise product]{The stochastic product is part of the model. The arrows run only
within each construction; there is no conversion arrow between the columns.
For the linear Skorohod equation with deterministic data, iteration order
agrees with orthogonal chaos order. Ordinary Stratonovich iteration retains
self-contractions and unequal-order moment terms, although the final
solution can still be regrouped into orthogonal Wiener chaoses. See
Balan~\cite{balan:22:stratonovich} for a comparison of the expansions.
An adapted It\^o identification does not identify an unspecified
smooth-noise approximation.}
\label{Fig:literature-interpretations}
\end{figure}
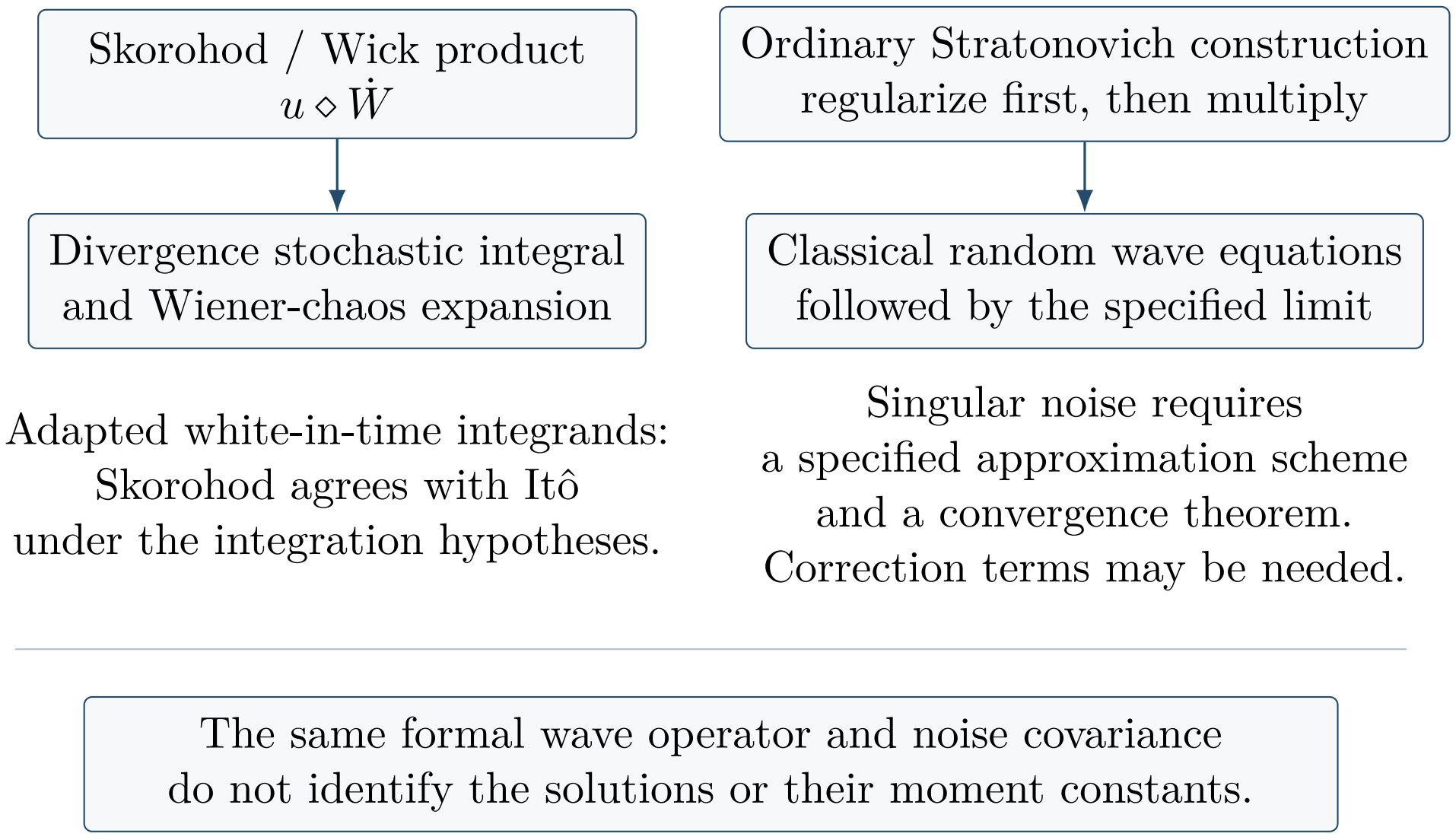

\subsection{White-in-time noise}

White-in-time wave equations provide a useful starting point because adapted
stochastic integration is available. Conus--Dalang
\cite{conus.dalang:08:non-linear} develop constructions in higher dimensions
under their noise and initial-data assumptions. Dalang--Mueller--Tribe
\cite{dalang.mueller.ea:08:feynman-kac-type} obtain Feynman--Kac-type
representations using Poisson processes and random flights. Thus, the absence
of the usual Brownian heat-equation representation does not imply the absence
of all probabilistic moment representations for wave equations.

Chen--Dalang~\cite{chen.dalang:15:moment}*{Theorem~2.3 and Corollary~1.1}
derive moment bounds and explicit second-moment formulas for the
one-dimensional space-time white-noise equation, including rough initial
data. Chen--Guo--Song~\cite{chen.guo.ea:24:moments}*{Theorem~1.1} study a
fractional SPDE family driven by space-time white noise. Their second moments
and second Lyapunov exponents are explicit. Their wave specialization gives
matching upper and lower $p^{3/2}$ moment-growth orders under the stated data
conditions~\cite{chen.guo.ea:24:moments}*{Proposition~5.10; (1.10)--(1.11),
pp.~4257--4258}. The higher-moment constants are not optimal, and these bounds
do not assert an exact higher-moment leading constant or a general
higher-moment limit.

Balan--Song~\cite{balan.song:19:second}*{Theorems~1.1 and 1.3} relate the
second-moment Laplace transforms of hyperbolic and parabolic Anderson models.
Their general result identifies an upper second Lyapunov exponent; their
Riesz specialization proves an actual limit. These conclusions should not be
interchanged. In three dimensions with a different spatial covariance class,
Dalang--Mueller~\cite{dalang.mueller:09:intermittency} establish upper and
lower moment-growth bounds of order $p^{4/3}$, with the lower bound at even
orders and with separate hypotheses on the two bounds. The power of the moment
order therefore depends on the spatial covariance, not only on the word
``wave''.

Rough spatial noise is another distinct direction. Balan--Jolis--Quer-Sardanyons
\cite{balan.jolis.ea:17:intermittency} treat one-dimensional spatial fractional
noise with Hurst index between $1/4$ and $1/2$; Liu--Hu--Wang
\cite{liu.hu.ea:22:nonlinear} develop nonlinear well-posedness in this range
under their coefficient assumptions. Such distributional spatial covariances
are outside the nonnegative spatial-measure class of the present paper. Their
results do not extend the scope of Hypothesis~\ref{H:homog} by analogy.

\subsection{Time-independent noise}

For the static Skorohod equation with initial data $(1,0)$, Balan--Le Chen--Xia
Chen~\cite{balan.chen.ea:22:exact}*{Theorem~1.6} obtain exact large-time and
high-moment asymptotics with a variational constant. Their homogeneous spatial
covariances cover dimensions at most three; their global spatial white-noise
case covers dimensions one and two. Three-dimensional white noise belongs to
a separate finite-time discussion, not this global large-time assertion.
Chen--Eisenberg~\cite{chen.eisenberg:23:interpolating}*{Theorems~1.6--1.7}
study fractional differential operators with static Skorohod noise. Their time
order lies strictly below two. Fractional time differentiation is distinct
from temporal correlation of the noise and does not silently include the
ordinary wave endpoint.

For ordinary static Stratonovich products, Balan
\cite{balan:22:stratonovich} constructs solutions in dimensions one and two
under a spectral condition stronger than Dalang's condition, and gives both
a probabilistic representation and a comparison of chaos expansions.
Xia Chen--Hu~\cite{chen.hu:25:hyperbolic} treat dimensions at most three under
Dalang's condition and establish homogeneous moment asymptotics at positive
integer orders. Their published variational expressions require the
normalization qualifications recorded in this manuscript; we do not import
a numerical white-noise constant from them.

There is a concrete difference in admissible roughness for these unrenormalized
static constructions. For the common flat data and homogeneous function
covariances, the Skorohod range in the preceding exact-asymptotics result
allows spatial homogeneity exponent below three, subject to its dimension
restriction, whereas the ordinary Stratonovich range requires exponent below
two. This comparison concerns those solution classes. Renormalized
Anderson-Hamiltonian wave equations, such as those in
Gubinelli--Ugurcan--Zachhuber~\cite{gubinelli.ugurcan.ea:20:semilinear}, are a
related but different construction; their existence does not identify an
ordinary Stratonovich product with a Skorohod product.

\subsection{Fractional temporal correlation}

On the Skorohod side, Balan--Song~\cite{balan.song:17:hyperbolic} develop the
space-time homogeneous Gaussian-noise construction. Its higher-dimensional
statements have additional spectral hypotheses, including absolute continuity
of the spectral measure in the relevant results. Xia Chen--Deya--Song--Tindel
\cite{chen.deya.ea:25:solving} sharpen the spectral solvability analysis. For
homogeneous spatial covariance and temporal singularity exponent $\alpha_0$,
the resulting noise threshold is $\alpha + \alpha_0 < 3$, under the paper's
stated data conditions. The sufficiency statement in dimensions two and three
imposes Fourier integrability on the initial displacement; the necessity
statement uses a positive lower bound on the free wave. These are different
hypotheses and directions, not an unconditional common-data equivalence with
the present construction. The static endpoint also needs separate treatment.

Balan--Conus~\cite{balan.conus:16:intermittency} give moment bounds for
fractional-in-time noise under their spatial covariance conditions, including
a second-moment lower bound in dimensions at most three at the same time-growth
order as the upper bound for positive constant initial displacement and
nonnegative constant initial velocity. Hu--Wang~\cite{hu.wang:24:matching}*{Theorem~3.7}
obtain matching higher-moment growth orders under their covariance and
Green-kernel hypotheses. Neither conclusion by itself identifies a sharp
leading constant. Space-time white noise is not an automatic specialization
of their functional spatial lower-bound hypotheses.

On the Stratonovich side, Xia Chen--Deya--Song--Tindel
\cite{chen.deya.ea:23:hyperbolic} construct solutions in dimensions one and two
using Strichartz estimates and weighted function spaces. Their sufficient
regularity ranges and initial-data assumptions differ from the full
homogeneous range considered here. Xia Chen's current author version
\cite{chen:author2026:time-dependency} treats fractional temporal correlation
in dimensions one through three, with homogeneous condition
$\alpha + \alpha_0 < 2$. It proves two-sided first-moment asymptotic bounds and
a coarse quadratic-exponential second-moment upper bound. The stronger exact
first-moment statement of arXiv:2510.01412v1 is not the statement of this
edition; the version convention is specified in Section~\ref{S:known}.

\begin{figure}
  \centering
  % Shared vector drawing for the report and reading edition.
\begin{tikzpicture}[x=1cm,y=1cm,
  lit/text/.style={anchor=north west,align=flush left,text width=5.35cm,
  font=\small,inner sep=0pt,
  execute at begin node={\hyphenpenalty=10000\relax}},
  lit/tag/.style={anchor=west,font=\small\bfseries,inner sep=0pt},
  lit/regime/.style={align=center,font=\small,inner sep=0pt}]
  % The two colored columns classify products, not solution inclusions.
  \fill[hamblue!3] (1.65,.35) rectangle (7.35,-11.75);
  \fill[hamteal!3] (7.75,.35) rectangle (13.45,-11.75);
  \fill[hamblue] (1.65,.35) rectangle (7.35,1.15);
  \fill[hamteal] (7.75,.35) rectangle (13.45,1.15);
  \node[text=white,font=\bfseries] at (4.5,.75) {Skorohod};
  \node[text=white,font=\bfseries] at (10.6,.75) {Ordinary Stratonovich};
  \node[lit/regime,text=hamblue!80!black] at (.55,.75)
  {Temporal\\covariance};

  % White time: the spike represents a delta covariance, not a noise path.
  \node[lit/regime] at (.55,-.25) {White\\in time};
  \draw[gray] (-.05,-1.4)--(1.15,-1.4);
  \draw[-{Latex[length=1.6mm]},thick,hamgold]
  (.55,-1.4)--(.55,-.85);
  \node[lit/regime] at (.55,-1.8) {$\delta_0$};
  \node[lit/tag,text=hamblue] at (1.9,-.1) {Exact second moments};
  \node[lit/text] at (1.9,-.45)
  {Chen--Dalang~\cite{chen.dalang:15:moment}\\
  Space-time white noise, $d = 1$.\\[4pt]
  Chen--Guo--Song~\cite{chen.guo.ea:24:moments}\\
  Exact $p = 2$; wave higher-$p$\\
  growth orders $p^{3/2}$.};
  \node[lit/tag,text=hamteal] at (8,-.1) {Approximation matters};
  \node[lit/text] at (8,-.45)
  {No source claim is listed here.\\[4pt]
  Specify smoothing and convergence.\\[4pt]
  No It\^o equality or openness\\
  is inferred.};

  \draw[hamblue!18] (1.9,-3.2)--(7.1,-3.2);
  \draw[hamteal!22] (8,-3.2)--(13.2,-3.2);
  % Static time: one constant covariance profile.
  \node[lit/regime] at (.55,-3.65) {Static\\noise};
  \draw[gray] (-.05,-4.8)--(1.15,-4.8);
  \draw[thick,hamgold] (-.05,-4.25)--(1.15,-4.25);
  \node[lit/regime] at (.55,-5.2) {$1$};
  \node[lit/tag,text=hamblue] at (1.9,-3.7) {Exact rates};
  \node[lit/text] at (1.9,-4.05)
  {Balan--Le Chen--Xia Chen\\
  \cite{balan.chen.ea:22:exact}: $p \ge 2$.\\
  Homogeneous: $d \le 3$;\\
  spatial white: $d \le 2$.\\[4pt]
  Chen--Eisenberg~\cite{chen.eisenberg:23:interpolating}\\
  Exact rates; static noise;\\
  fractional time order $< 2$.};
  \node[lit/tag,text=hamteal] at (8,-3.7) {Construction / exact rates};
  \node[lit/text] at (8,-4.05)
  {Balan~\cite{balan:22:stratonovich}: $d = 1,2$ construction.\\
  Stronger spectral condition.\\[4pt]
  Xia Chen--Hu~\cite{chen.hu:25:hyperbolic}\\
  Exact rates, integer $p \ge 1$;\\
  Dalang and homogeneity hypotheses.};

  \draw[hamblue!18] (1.9,-7.7)--(7.1,-7.7);
  \draw[hamteal!22] (8,-7.7)--(13.2,-7.7);
  % Fractional time: schematic singular, decreasing covariance.
  \node[lit/regime] at (.55,-8.45) {Fractional\\in time};
  \draw[gray] (-.05,-9.65)--(1.15,-9.65);
  \draw[thick,hamgold]
  (.05,-8.85) .. controls (.18,-9.35) and (.65,-9.45) .. (1.15,-9.52);
  \node[lit/regime] at (.55,-10.1) {$|t|^{-\alpha_0}$};
  \node[lit/tag,text=hamblue] at (1.9,-8.1) {Bounds / construction};
  \node[lit/text] at (1.9,-8.45)
  {Balan--Conus; Hu--Wang\\
  \cites{balan.conus:16:intermittency,hu.wang:24:matching}\\
  Growth orders; no sharp constant.\\[4pt]
  Chen--Deya--Song--Tindel\\
  \cite{chen.deya.ea:25:solving}: construction;\\
  data hypotheses depend on $d$.};
  \node[lit/tag,text=hamteal] at (8,-8.1) {First-moment bounds};
  \node[lit/text] at (8,-8.45)
  {Xia Chen, author version\\
  \cite{chen:author2026:time-dependency}\\
  $d = 1,2,3$, $\alpha + \alpha_0 < 2$.\\[4pt]
  Two-sided first-moment bounds;\\
  the exact limit is predicted.};

  % A separate manuscript contribution, not a transferred literature claim.
  \fill[hamgold!9,rounded corners=2pt]
  (1.65,-12.1) rectangle (13.45,-13.65);
  \draw[hamgold,line width=2pt] (1.65,-12.1)--(1.65,-13.65);
  \node[lit/tag,text=hamgold!85!black] at (1.95,-12.45)
  {Present manuscript \quad Fractional-time Stratonovich, $p = 2$};
  \node[anchor=north west,align=left,text width=10.9cm,
  font=\small,inner sep=0pt] at (1.95,-12.75)
  {Intrinsic variational formula for the second-moment limit.\\
  Written proof; independent review remains pending.};
\end{tikzpicture}
  \caption[Source-qualified literature comparison]{Selected results, organized
    by temporal covariance and stochastic product. The small profiles are
    schematic covariance sketches, not sample noise paths or limiting arrows.
    Construction, growth orders and exact rates are different conclusions.
    Every source retains its spatial covariance, dimension and data hypotheses;
    the Chen--Guo--Song higher-moment constants are nonoptimal, and the
    higher-dimensional Skorohod construction retains its data qualifications.
    Chen--Hu supplies no numerical constant here. Xia Chen's author
    version gives first-moment bounds, unlike the stronger arXiv-v1 wording.
    The highlighted manuscript claim has a separate review status. The columns
    assert neither inclusions nor a universal conversion between products.}
  \label{Fig:literature-results-map}
\end{figure}

\subsection{Advantages, limitations and the present contribution}

For second-moment analysis, the main advantage of the Skorohod formulation is
orthogonality: it removes unequal-chaos-order cross terms and makes isometries
and hypercontractivity directly available. In some static covariance classes
it also permits rougher noise than ordinary Stratonovich multiplication. The
price is that Wick-ordered multiplication specifies a different stochastic model;
it cannot generally be substituted for an ordinary smooth random potential.
With temporally correlated noise, it does not restore adapted martingale
methods.

The Stratonovich formulation retains the interactions generated by ordinary
multiplication and can connect to deterministic wave and Gaussian variational
methods. Its cost is control of traces, self-pairings and unequal iteration
orders, together with potentially stronger conditions for an unrenormalized
limit. Neither interpretation supplies a general wave maximum principle.
Both admit chaos methods and, in suitable settings, probabilistic
representations. The appropriate choice is determined by the model and the
specified approximation, rather than by a universal preference for one
calculus.

The present contribution concerns the temporally correlated Stratonovich
problem. Theorem~\ref{T:full-upper} gives the full homogeneous-class upper bound
at the scale $t^\kappa$. Theorem~\ref{T:intrinsic-second-moment} identifies the
second-moment limit through nonnegative compact space-time controls, retaining
the original covariance amplitude, unit kinetic coefficient and global
half-variance cost. This intrinsic expression is not a numerical closed form,
a sharp all-moment formula or an identification with the source's
first-moment variational constants. These are the claims of the written
proofs; the independent review status is recorded separately from this
literature comparison.

Related work in progress with Raluca Balan and Jinxin Wang treats three
Skorohod projects: noise colored in time and space in one dimension; colored
noise in two dimensions; and white-in-time noise in dimensions one, two and
three. We mention these scopes to locate the present work within the broader
program. No unpublished result from them is used here, and no equality of
their constants with the present Stratonovich constant is asserted.

\section{What is known}\label{S:known}

Pinpoints to Xia Chen's time-dependent Stratonovich paper refer to his
public author version~\cite{chen:author2026:time-dependency}.
Theorem~1.2 differs from the statement in arXiv:2510.01412v1
\cite{chen:25:time-dependency}; we therefore cite that theorem from the
author version.

\subsection{Stratonovich expansion and moments}

The solution is the Stratonovich expansion
\cite{chen:author2026:time-dependency}*{(2.13)--(2.16)}
\begin{equation}\label{E:expansion}
  \begin{split}
    u(t,x)&=\sum_{n=0}^\infty S_n\big(g_n(\cdot,t,x)\big),\\
    g_n(s_1,\dots,s_n,x_1,\dots,x_n,t,x)
    &=\Big(\prod_{l=1}^nG(s_{l+1}-s_l,x_{l+1}-x_l)\Big)
    \mathbf{1}_{[0,t]^n_<}(s_1,\dots,s_n),
  \end{split}
\end{equation}
with $s_{n+1} = t$, $x_{n+1} = x$, $G$ the fundamental solution of the wave
equation, and $[0,t]^n_ < = \{0<s_1<\dots<s_n<t\}$. Odd expected Stratonovich
integrals vanish and even ones are given by Wick's formula
\cite{chen:author2026:time-dependency}*{(2.19)--(2.21)}, a sum over pair
partitions $\mathcal{D} \in \Pi_n$ of $\{1,\dots,2n\}$ that includes
self-pairings within one chain. The second moment of one ordered term is the
two-chain pairing integral, stated
in~\cite{chen:author2026:time-dependency}*{Theorem~4.1, (4.3)} and proved here
under Hypothesis~\ref{H:homog} in Proposition~\ref{P:mixed}:
\begin{equation}\label{E:Sn2}
  \begin{split}
    \E\big[S_n(g_n(\cdot,t,x))^2\big]
    &=\sum_{\mathcal{D}\in\Pi_n}\int_{[0,t]^n_<\times[0,t]^n_<}\ud s
    \int_{(\R^d)^n\times(\R^d)^n}\ud x
    \Big(\prod_{\rho=1}^2\prod_{l\in I_\rho}
    G(s_l-s_{p(l)},x_l-x_{p(l)})\Big)\\
    &\qquad\times\prod_{(j,k)\in\mathcal{D}}|s_j-s_k|^{-\alpha_0}\gamma(x_j-x_k),
  \end{split}
\end{equation}
with $I_1 = \{1,\dots,n\}$ and $I_2 = \{n+1,\dots,2n\}$. Here $p(l)$ is the
preceding label in its own chain; at the first vertex of either chain use the
common root $(s_{p(l)},x_{p(l)}) = (0,0)$. Thus, the second chain does not
continue the first one. This spells out the source's conventions in Theorem~4.1.

The full mixed-term expansion is justified in Proposition~\ref{P:mixed}.

\subsection{The universal two-chain bound}

\cite{chen:author2026:time-dependency}*{Theorem~3.1, (3.6), p.~15} states an
estimate for such two-chain pairing integrals with nonnegative kernels $G_l$ by
products of the norms $\|G_l\|^{(0)}$, $\|G_l\|^{(1)}_{j,k}$,
$\|G_l\|^{(2)}_{j,k}$ of~\cite{chen:author2026:time-dependency}*{(3.3)--(3.5)},
over a partition $\{Q_0,Q_1,Q_2\}$ of $\{1,\dots,2n\}$ with $\#Q_0 = \#Q_1$ and
$\#Q_2$ even. As stated, the theorem uses finite spatial spectral measures
$\mu_{j,k}$ with nonnegative covariance functions and temporal kernels
$\Gamma(\alpha_0)^{-1}\int_{A_{j,k}}e^{-\lambda|r|}
\lambda^{\alpha_0-1}\ud\lambda$ for measurable $A_{j,k} \subset \R_+$. The
application to the full covariance therefore requires regularization; this is
carried out in Lemma~\ref{L:crude-corrected}. The theorem includes
cross-pairings between the two chains. Our finite maximum formulation is proved
in Lemma~\ref{L:pairing-product-envelope}; a maximizing assignment may depend on
the kernels. Lemma~\ref{L:translated-covariance-CS} supplies its Case~2 estimate
with the predecessor-time difference retained, as explained in
Remark~\ref{R:source-case2-time-shift}. In the source it is used only for
Stratonovich integrability (Section~4) and, via Lemma~A.1 with $G_l(t,x) =
e^{-\theta t}G(t,x)$, in Laplace variables.

\subsection{Known moment bounds}

The following statement records the source as printed. Its variational constants
require normalization clarification; see Remark~\ref{R:three-endpoint-scope}. We
retain the printed formulas without adopting a correction. The upper-bound
proofs below do not use the first-moment upper constant.

\begin{theorem}[\cite{chen:author2026:time-dependency}*{Theorem~1.2, (1.10)--(1.13)}]\label{T:first}
  \emph{Source statement as printed.} Under its covariance and homogeneity
  assumptions~\cite{chen:author2026:time-dependency}*{(1.2), (1.3), (1.8), (1.9)},
  the author version states
  \[
    L_*\le\liminf_{t\to\infty}t^{-\kappa}\log\E u(t,x)
    \le\limsup_{t\to\infty}t^{-\kappa}\log\E u(t,x)\le L^*,
  \]
  where, for $\cE^\sharp \in \{\cE_*,\cE\}$,
  \[
    L(\cE^\sharp):=(3-\alpha)\bigg(
    \frac{(4-\alpha-2\alpha_0)^{\frac{4-\alpha-2\alpha_0}{2}}}
    {2^3(4-\alpha-\alpha_0)^{4-\alpha-\alpha_0}}
    \Big(\frac{2\cE^\sharp}{2-\alpha}\Big)^{\frac{2-\alpha}{2}}
    \bigg)^{\frac1{3-\alpha}},
    \qquad L_*=L(\cE_*),\quad L^*=L(\cE),
  \]
  and $\cE_*,\cE$ are the variational quantities
  \cite{chen:author2026:time-dependency}*{(1.11)--(1.13)}. The strict source
  comparison $\cE_* < \cE$ is asserted only for $\Gamma \ne 0$; it is not
  asserted at $\Gamma = 0$, where both quantities vanish. These printed
  normalizations are quoted, not adopted as established inputs. Numerical
  calibrations using $L_*$ below explicitly assume the intended first-moment
  lower bound.
\end{theorem}

The source writes its spatial covariance as a function and explicitly includes
one-dimensional spatial white noise. We do not claim an extension of this
quoted result to every singular Radon covariance in Hypothesis~\ref{H:homog}.

The source predicts that the limit exists and equals $L^*$
\cite{chen:author2026:time-dependency}*{p.~4}; this is not proved.
The source also suggests that the constant in its lower estimate may not be
sharp, relating the difficulty to the complex-valued intersection-local-time
representation~\cite{chen:author2026:time-dependency}*{Remark~8.1, (8.8)}.

\begin{theorem}[\cite{chen:author2026:time-dependency}*{Remark~5.1, (5.8)}]\label{T:crude}
  $\E u^2(t,x) \le e^{Ct^2}$ for large $t$.
\end{theorem}

\begin{conjecture}[\cite{chen:author2026:time-dependency}*{(1.22)}]\label{Conj:allp}
  For $p = 1,2,\dots$,
  \[
    \lim_{t\to\infty}t^{-\kappa}\log\E|u(t,x)|^p
    =p^{\frac{4-\alpha}{3-\alpha}}L^*.
  \]
\end{conjecture}

Since $p^{(4-\alpha)/(3-\alpha)}$ grows faster than $p$, the conjectured
second-moment constant $2^{(4-\alpha)/(3-\alpha)}L^*$ exceeds $2L^*$. The source
explains why its method stops at $p = 1$
\cite{chen:author2026:time-dependency}*{Remark~6.3}: the representation of
Theorem~6.1 there uses the monotone order $s_1 \le \dots \le s_{2n}$ of a single
chain.

\begin{proposition}[Jensen lower bound]\label{P:Jensen}
  Under Hypothesis~\ref{H:homog}, $\E u(t,x) \ge 1$ and $\E u^2(t,x) \ge (\E
  u(t,x))^2$. Consequently, for every real $\ell$,
  \[
    \liminf_{t\to\infty}t^{-\kappa}\log\E u(t,x)\ge\ell
    \quad\Longrightarrow\quad
    \liminf_{t\to\infty}t^{-\kappa}\log\E u^2(t,x)\ge2\ell.
  \]
  The specialization $\ell = L_*$ is conditional on the intended first-moment
  lower bound with that constant. It is not an unconditional adoption of the
  printed normalization in Theorem~\ref{T:first}.
\end{proposition}

\begin{proof}
  For each fixed $t,x$, Lemma~\ref{L:vertex-crude} below gives
  $\sum_n\|X_n(t,x)\|_2 < \infty$; its proof does not use the present
  proposition. Thus, $\sum_n\E|X_n(t,x)| < \infty$, and the mean of the
  canonical series is the sum of its means. By~\eqref{E:expansion}
  and~\cite{chen:author2026:time-dependency}*{(2.20)--(2.21)}, $\E u(t,x) =
  \sum_{n\ge0}\E S_{2n}(g_{2n}(\cdot,t,x)) \ge 1$, since every term is a pairing
  integral of a nonnegative integrand and $S_0 = 1$. By the Cauchy--Schwarz
  inequality, $\E u^2(t,x) \ge (\E u(t,x))^2$, so $t^{-\kappa}\log\E u^2(t,x)
  \ge 2\,t^{-\kappa}\log\E u(t,x)$. Taking lower limits proves the implication
  for every valid $\ell$. Substituting $L_*$ requires its stated first-moment
  premise. This proves Proposition~\ref{P:Jensen}.
\end{proof}

\section{Main results}\label{S:targets}

We state the two main results for the model in Section~\ref{S:model}. Throughout
this section, Hypothesis~\ref{H:homog} holds, $u$ is the canonical Stratonovich
expansion~\eqref{E:expansion} with initial data $(1,0)$, and $\kappa$ is the
exponent in~\eqref{E:kappa}. The spectral measure $\mu$ has the normalization
specified in Hypothesis~\ref{H:homog}. Write $B_1$ for the unit ball in
$\R^d$, and put
\begin{equation}\label{E:intro-spectral-mass}
  R_\alpha=\frac{\alpha\pi\mu(B_1)}{2\sin(\pi\alpha/2)}.
\end{equation}

\begin{theorem}[Upper bound at the correct scale]\label{T:intro-upper}
  \label{G:upper}
  For every admissible dimension $d \in \{1,2,3\}$, there is a finite constant
  $C = C(\alpha,\alpha_0)$ such that, for all $T \ge 1$ and $x \in \R^d$,
  \begin{equation}\label{E:intro-upper-bound}
    \E u^2(T,x)
    \le e(T+2)^d
    \exp\left\{C R_\alpha^{1/(3-\alpha)}T^\kappa\right\}.
  \end{equation}
  In particular, for every fixed $x \in \R^d$,
  \[
    \limsup_{T\to\infty}T^{-\kappa}\log\E u^2(T,x)<\infty.
  \]
\end{theorem}

The proof is given in Theorem~\ref{T:full-upper}. It establishes finiteness
independently of the variational identification below. The more explicit bounds
for $0 < \alpha \le 1$ are retained in Section~\ref{S:d1-upper}.

To state the limiting constant, let $I = (0,1)$ and define
\begin{equation}\label{E:intro-control-class}
  \mathcal{C}_{\mathrm{pos}}
  =\left\{\Psi\in C_c^\infty(I\times\R^d;\R):
  \Psi\ge0,\ \Psi\ne0\right\}.
\end{equation}

No-go~\ref{P:nogo-20} exhibits controls whose covariance profiles do not
separate into a scalar time factor and one spatial function. This does not
decide whether a smaller family suffices for optimization.

This class includes controls that do not separate into a function of time and
a function of space. For $\Psi \in \mathcal{C}_{\mathrm{pos}}$, its covariance
mean and variance are
\begin{align}
  h_{\Psi,s}(x)
  &=\int_I|s-t|^{-\alpha_0}
  \int_{\R^d}\Psi(t,x-y)\Gamma(\ud y)\ud t,
  \label{E:intro-control-mean}\\
  v_\Psi
  &=\int_{I\times\R^d}\Psi(s,x)h_{\Psi,s}(x)\ud x\ud s.
  \label{E:intro-control-variance}
\end{align}
For $s \in I$, define the nonnegative Rayleigh value
\begin{equation}\label{E:intro-control-rate}
  \mathcal{R}_\Psi(s)
  =\sup_{\substack{f\in H^1(\R^d;\R)\\\|f\|_2=1}}
  \left\{\int_{\R^d}h_{\Psi,s}(x)f(x)^2\ud x
  -\|\nabla f\|_2^2\right\}.
\end{equation}
Here $H^1(\R^d;\R)$ is the real Sobolev space of square-integrable functions
with square-integrable first weak derivatives, and $\|\cdot\|_2$ denotes the
spatial $L^2$ norm. The spatial
covariance retains its original amplitude in these definitions.

\begin{theorem}[Second-moment limit]\label{T:intro-second-moment}
  \label{G:sharp}
  For every fixed $x \in \R^d$,
  \begin{equation}\label{E:intro-second-moment-limit}
    \lim_{T\to\infty}T^{-\kappa}\log\E u^2(T,x)
    =C_2(\Gamma,\alpha_0),
  \end{equation}
  where
  \begin{equation}\label{E:intro-second-moment-constant}
    C_2(\Gamma,\alpha_0)
    =\sup_{\Psi\in\mathcal{C}_{\mathrm{pos}}}
    \left\{2\int_I\sqrt{\mathcal{R}_\Psi(s)}\ud s
    -\frac{v_\Psi}{2}\right\}.
  \end{equation}
  If $\Gamma \ne 0$, then $0 < C_2(\Gamma,\alpha_0) < \infty$.
  If $\Gamma = 0$, then $u(T,x) = 1$ and $C_2(0,\alpha_0) = 0$.
\end{theorem}

The proof is given in Theorem~\ref{T:intrinsic-second-moment}. The variational
formula uses the full class of nonnegative smooth compact controls. Its kinetic
coefficient is one, and its Gaussian cost is one half of the original variance.
It does not identify an auxiliary variational value with the first-moment limit
or with the constants in Chen's conjecture.
Theorem~\ref{T:intrinsic-occupation-ceiling},
Corollary~\ref{C:intrinsic-spatial-bracket} and
Theorem~\ref{T:white-moment-bounds} give bounds for this same constant. In
particular, the explicit spatial white-noise energy yields moment bounds,
not a closed-form evaluation of the full supremum.

\section{Why the scale is \texorpdfstring{$t^\kappa$}{t to the kappa}}\label{S:why-kappa}

This section explains, without proof, where the exponent $\kappa$ and the form
of $C_2(\Gamma,\alpha_0)$ in Theorem~\ref{T:intro-second-moment} come from.

\begin{figure}[htbp]
  \centering
  % Exponent only; open markers denote excluded parameter endpoints.
\begin{tikzpicture}[x=4.7cm,y=2.7cm,font=\normalsize]
  \draw[ham/arrow] (0,1)--(1.06,1);
  \draw[ham/arrow] (0,1)--(0,2.04) node[above] {$\kappa$};
  \foreach \xx/\lab in {0/0,.5/{1/2},1/1} {
  \draw[hamblue] (\xx,1)--(\xx,.96);
  \node[anchor=north] at (\xx,.94) {$\lab$};
  }
  \foreach \yy/\lab in {1/1,{3/2}/{3/2},{5/3}/{5/3},2/2} {
  \draw[hamblue] (-.025,{\yy})--(0,{\yy});
  \node[anchor=east] at (-.04,{\yy}) {$\lab$};
  }
  \draw[hamblue] (-.025,{7/5})--(0,{7/5});
  \draw[hamblue,line width=.45pt]
  (-.025,{7/5})--(-.055,{7/5})--(-.055,1.30);
  \node[anchor=east,inner sep=0pt] at (-.07,1.30) {$7/5$};
  \node[align=center,text width=4.1cm] at (.5,1.93)
  {Static limit at left ends:\\ $(4-\alpha)/(3-\alpha)$};
  \draw[hamblue,line width=1pt] (0,{7/5})--(1,1);
  \draw[hamteal,line width=1pt] (0,{3/2})--(1,1);
  \draw[hamgold,line width=1pt] (0,{5/3})--(.5,{4/3});
  \draw[hamgold!45,densely dashed,line width=.8pt]
  (.5,{4/3})--(1,1);
  \foreach \yy/\colour in {{7/5}/hamblue,{3/2}/hamteal,{5/3}/hamgold} {
  \filldraw[draw=\colour,fill=white,line width=.8pt] (0,{\yy}) circle (1.5pt);
  }
  \draw[hamgold!65,dotted,line width=.45pt]
  (.5,{4/3})--(.56,1.40);
  \node[anchor=south west,inner sep=0pt,text=hamgold] at (.56,1.40) {$4/3$};
  \filldraw[draw=hamgold,fill=white,line width=.8pt] (.5,{4/3}) circle (1.5pt);
  \filldraw[draw=black!65,fill=white,line width=.8pt] (1,1) circle (1.5pt);
  \draw[hamblue,line width=1pt] (1.12,1.96)--(1.27,1.96);
  \node[anchor=west] at (1.32,1.96) {$\alpha=1/2$};
  \draw[hamteal,line width=1pt] (1.12,1.75)--(1.27,1.75);
  \node[anchor=west] at (1.32,1.75) {$\alpha=1$};
  \node[anchor=west] at (1.12,1.56) {white in space};
  \node[anchor=west] at (1.12,1.39) {$d=1$};
  \draw[hamgold,line width=1pt] (1.12,1.17)--(1.27,1.17);
  \node[anchor=west] at (1.32,1.17) {$\alpha=3/2$};
  \node[anchor=west,align=left] at (1.12,.87) {dashed: outside\\ hypothesis};
  \node at (.5,.61) {temporal exponent $\alpha_0$};
  \draw[ham/arrow] (.27,.40)--(.10,.40);
  \node[anchor=west,inner sep=0pt] at (.31,.40) {longer memory};
\end{tikzpicture}
  \caption[The exponent $\kappa$ against the temporal exponent]{The exponent
    $\kappa=1+(1-\alpha_0)/(3-\alpha)$ of~\eqref{E:kappa}, drawn against the
    temporal exponent $\alpha_0$ for three spatial exponents $\alpha$. Longer
    memory in time (smaller $\alpha_0$) and a more singular spatial covariance
    (larger $\alpha$) both raise $\kappa$. As $\alpha_0\downarrow0$ each line
    tends to the static exponent $(4-\alpha)/(3-\alpha)$. The end points are
    excluded: Hypothesis~\ref{H:homog} requires $0<\alpha_0<1$ and
    $\alpha+\alpha_0<2$, so the line for $\alpha=3/2$ stops at $\alpha_0=1/2$.
    The picture shows the exponent only; it says nothing about the constant
    $C_2(\Gamma,\alpha_0)$.}
  \label{Fig:kappa-fan}
\end{figure}
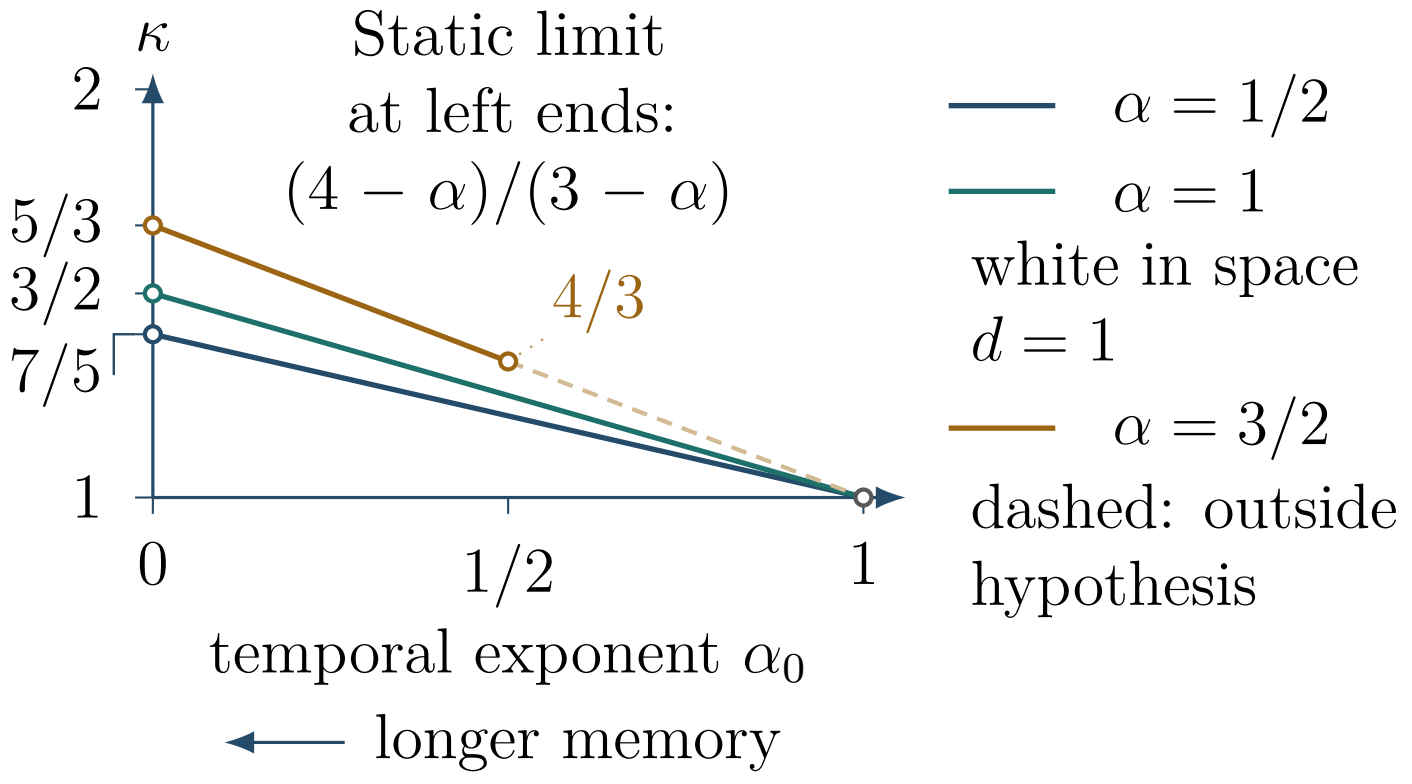

The first step is an exact scaling, proved in
Lemma~\ref{L:exact-canonical-scaling}. Stretch time by a factor $T$ and space
by $r_T=T^{(\alpha_0-1)/(3-\alpha)}$. Because the temporal covariance
$|s-t|^{-\alpha_0}$ and the spatial covariance $\Gamma$ are homogeneous
(Hypothesis~\ref{H:homog}), the stretched noise has the law of the original
noise multiplied by a constant. With $\delta=T^{-\kappa}$, the stretched
solution has the law of the solution $V_\delta$ of
\[
  \delta^2\partial_t^2V_\delta=\Delta V_\delta+\sqrt\delta\,\W V_\delta,
  \qquad V_\delta(0,\cdot)=1,\quad \partial_tV_\delta(0,\cdot)=0,
\]
and
\[
  T^{-\kappa}\log\E|u(T,0)|^2=\delta\log\E|V_\delta(1,0)|^2
\]
(see~\eqref{E:exact-fast-equation} and~\eqref{E:exact-moment-scaling}). With
this choice of $\kappa$, the time derivative carries the factor $\delta^2$ and
the noise the factor $\sqrt\delta$. The growth of the second moment over long
times thus becomes a question about small noise: does
$\delta\log\E|V_\delta(1,0)|^2$ have a limit as $\delta\to0$?

Heuristically it does, because a gain and a cost of the same order $1/\delta$
balance. Shift the noise toward a deterministic potential, chosen through a
nonnegative control $\Psi$, so that the potential seen by the wave at time $s$
is $h_{\Psi,s}$ of~\eqref{E:intro-control-mean}. In the scaled equation this
needs a shift of order $1/\sqrt\delta$, and its Gaussian cost is $v_\Psi/(2\delta)$,
one half of the variance~\eqref{E:intro-control-variance} divided by $\delta$.
Freezing the potential at time $s$ gives the equation
$\delta^2a''=\mathcal{R}_\Psi(s)\,a$ for the amplitude of the leading mode,
whose growing solution has exponential rate $\sqrt{\mathcal{R}_\Psi(s)}/\delta$.
Heuristically, the square of the wave then grows like
$\exp\bigl(2\int_I\sqrt{\mathcal{R}_\Psi(s)}\ud s/\delta\bigr)$. Gain minus
cost, multiplied by $\delta$, is
$2\int_I\sqrt{\mathcal{R}_\Psi(s)}\ud s-v_\Psi/2$, and taking the supremum over $\Psi\in\mathcal{C}_{\mathrm{pos}}$ gives
$C_2(\Gamma,\alpha_0)$. Returning to the original variables,
$\log\E u^2(T,x)$ grows like $C_2(\Gamma,\alpha_0)\,T^\kappa$.

This argument only explains the answer. A proof needs a lower bound for each
control, which computes the cost using the original covariance, retaining correlations between different time intervals
(Chapter~\ref{S:book-lower}), and an upper bound showing that no other
fluctuation of the noise does better (Chapters~\ref{S:book-means}--\ref{S:book-point-upper});
Chapter~\ref{S:book-constant} matches the two. Section~\ref{S:strategy} describes
these steps.

\section{Consistency checks}\label{S:consistency}

These comparisons place the second-moment theorem alongside known bounds and
Chen's conjecture. They are not inputs to the proof of
Theorem~\ref{T:intro-second-moment}.

\paragraph{Time-independent noise.}
As $\alpha_0\downarrow0$, the exponent $\kappa$ tends to
$(4-\alpha)/(3-\alpha)$, the exponent of the time-independent Stratonovich
model found by Chen and Hu~\cite{chen.hu:25:hyperbolic}*{Theorem~1.2}. This is
a formal boundary limit: the large-time limit and the limit
$\alpha_0\downarrow0$ are not interchanged (Remark~\ref{R:G1-scope}).

\paragraph{The first moment.}
If $\liminf_{t\to\infty}t^{-\kappa}\log\E u(t,x)\ge\ell$, then Jensen's
inequality $\E u^2(t,x)\ge(\E u(t,x))^2$ gives the corresponding second-moment
lower limit at least $2\ell$ (Proposition~\ref{P:Jensen}). Hence
$C_2(\Gamma,\alpha_0)\ge2\ell$ in Theorem~\ref{T:intro-second-moment}. Taking
$\ell=L_*$ of Theorem~\ref{T:first} is conditional on the intended first-moment
lower bound with that constant; its printed normalization is not adopted
unconditionally.

\paragraph{The Skorohod model.}
For Riesz covariances with $a=\alpha<d$ and $\alpha<2$, set
$H=1-\alpha_0/2>1/2$. The Skorohod upper bound of Balan and
Conus~\cite{balan.conus:16:intermittency}*{Theorem~2.1(b)} has the time
exponent $(2H+2-a)/(3-a)=\kappa$. The comparison concerns time exponents only:
the solution concepts and covariance normalizations differ, and no constants
are compared. The Stratonovich moment expansion also includes self-pairings,
that is, pairs of noise factors from the same ordered integral, which are
absent from the Skorohod expansion (Remark~\ref{R:G1-scope}).

\paragraph{Chen's conjecture.}
For $p=1,2,\dots$, Chen conjectures that $t^{-\kappa}\log\E|u(t,x)|^p$ tends to
$p^{(4-\alpha)/(3-\alpha)}L^*$ (Conjecture~\ref{Conj:allp}). For $p=2$ the
conjectured value is $2^{(4-\alpha)/(3-\alpha)}L^*$. When $L^*>0$, this exceeds
$2L^*$, the lower bound that Jensen's inequality would give if the first-moment
rate were $L^*$; at zero covariance both values vanish.
Theorem~\ref{T:intro-second-moment} does not decide this conjecture: its
constant $C_2(\Gamma,\alpha_0)$ is given by its own variational formula and is
not identified with the conjectured constant or with the first-moment limit.
The conjectured value is a point of comparison, not a target of the proof.

\section{Strategy of the proof}\label{S:strategy}

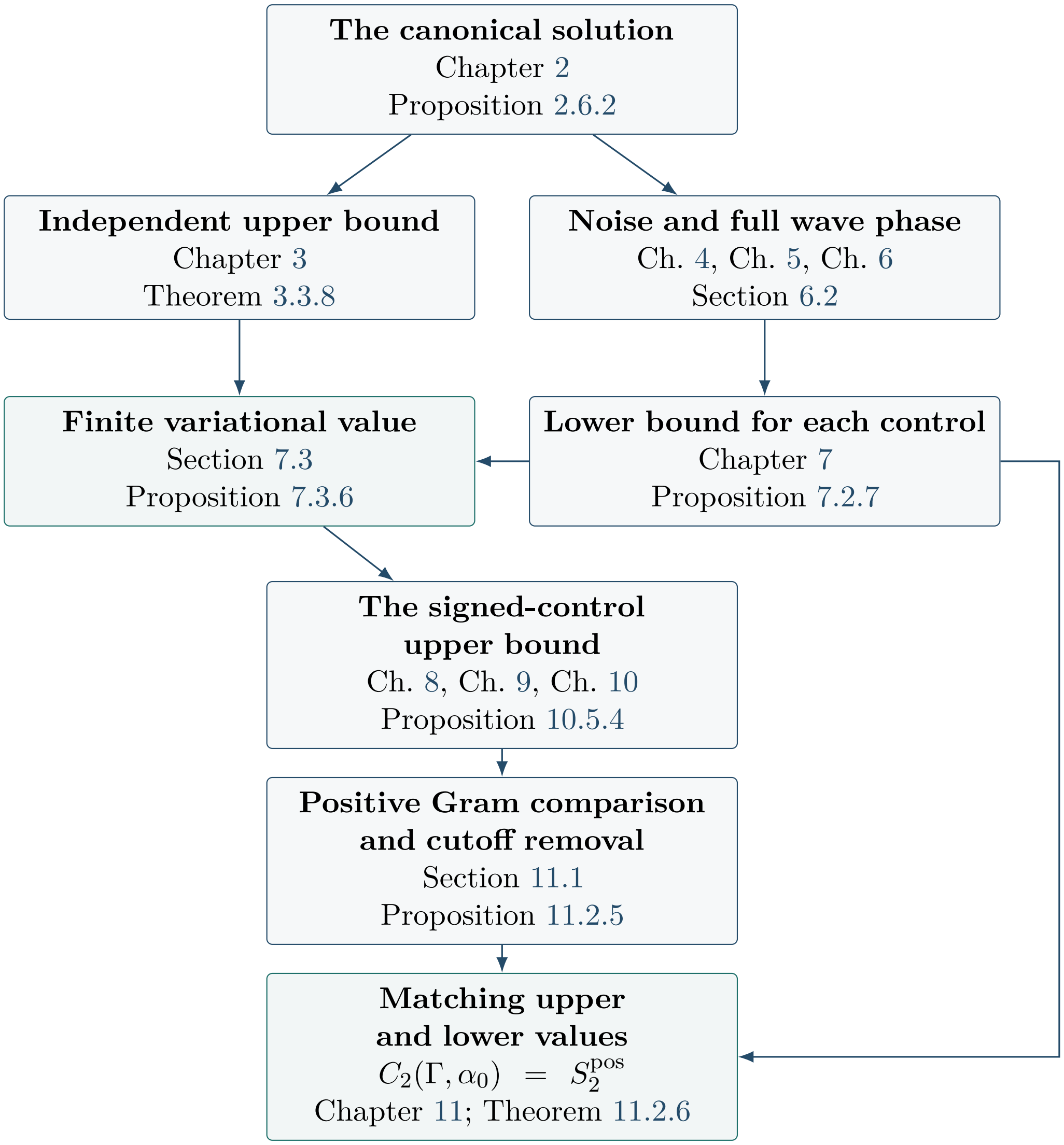
\begin{figure}
  \centering
  \begin{tikzpicture}[x=1cm,y=1cm,
    destination/.style={ham/box,text width=5.65cm,minimum height=1.55cm},
    result/.style={destination,draw=hamteal,fill=hamteal!6}]
    \node[destination] (construction) at (0,0)
    {\textbf{The canonical solution}\\
    \GuideRef{Chapter}{S:G0}\\
    \GuideRef{Proposition}{P:mild-construction}};
    \node[destination] (upper) at (-3.35,-2.4)
    {\textbf{Independent upper bound}\\
    \GuideRef{Chapter}{S:book-independent-upper}\\
    \GuideRef{Theorem}{T:full-upper}};
    \node[destination] (phase) at (3.35,-2.4)
    {\textbf{Noise and full wave phase}\\
    \GuideRef{Ch.}{S:book-gaussian-foundations},
    \GuideRef{Ch.}{S:book-operator-paths},
    \GuideRef{Ch.}{S:book-wave-evolution}\\
    \GuideRef{Section}{S:constant-phase-young}};
    \node[result] (finite) at (-3.35,-5)
    {\textbf{Finite variational value}\\
    \GuideRef{Section}{S:unmasked-variational}\\
    \GuideRef{Proposition}{P:variational-positive-scores}};
    \node[destination] (lower) at (3.35,-5)
    {\textbf{Lower bound for each control}\\
    \GuideRef{Chapter}{S:book-lower}\\
    \GuideRef{Proposition}{P:raw-block-original-lower}};
    \node[destination] (signed) at (0,-7.6)
    {\textbf{The signed-control upper bound}\\
    \GuideRef{Ch.}{S:book-means},
    \GuideRef{Ch.}{S:book-growing-radius},
    \GuideRef{Ch.}{S:book-point-upper}\\
    \GuideRef{Proposition}{P:original-signed-point-upper}};
    \node[destination] (comparison) at (0,-10.1)
    {\textbf{Positive Gram comparison}\\
    \textbf{and cutoff removal}\\
    \GuideRef{Section}{S:finite-positive-gram}\\
    \GuideRef{Proposition}{P:positive-score-cutoff}};
    \node[result] (limit) at (0,-12.6)
    {\textbf{Matching upper and lower values}\\
    $C_2(\Gamma,\alpha_0)=S_2^{\mathrm{pos}}$\\
    \GuideRef{Chapter}{S:book-constant};
    \GuideRef{Theorem}{T:intrinsic-second-moment}};
    \draw[ham/arrow] (construction)--(upper);
    \draw[ham/arrow] (construction)--(phase);
    \draw[ham/arrow] (upper)--(finite);
    \draw[ham/arrow] (phase)--(lower);
    \draw[ham/arrow] (lower)--(finite);
    \draw[ham/arrow] (finite)--(signed);
    \draw[ham/arrow] (signed)--(comparison);
    \draw[ham/arrow] (comparison)--(limit);
    \draw[ham/arrow] (lower.east)--++(.75,0) |- (limit.east);
  \end{tikzpicture}
  \caption[Proof dependencies and destinations]{The principal proof
    dependencies, with links to chapters and results in the technical companion.
    The independent
    upper bound and the fixed-control lower bound first give a finite
    variational lower value. Further operator and Gaussian estimates give a
    signed-control upper bound for the same point moment. Positive Gram
    comparison and cutoff removal identify the two values. Arrows show
    principal uses; each box also has the hypotheses and earlier inputs stated
    at its linked destination. They do not assert independence of successive
    time blocks. The zero-covariance case is included in the final theorem.}
  \label{Fig:proof-map}
\end{figure}

The second moment is a convergent sum of two-chain pairing integrals,
including empty chains and all relative time orders. Stratonovich integrals
retain pairings within either chain as well as between the chains, so products
of different orders can contribute (Proposition~\ref{P:mixed}). Estimating
every pairing by one common bound gives a finite moment but loses the target
growth scale (Proposition~\ref{P:envelope}). Chapter~\ref{S:book-no-go}
discusses this and the other approaches that do not work. The proof proceeds in
four steps, shown in Figure~\ref{Fig:proof-map}.

\paragraph{One solution and an upper bound at the right scale.}
Chapter~\ref{S:G0} derives the complete pairing expansion and constructs the
solution used throughout (Proposition~\ref{P:mild-construction}).
Chapter~\ref{S:book-independent-upper} then proves an upper bound at the scale
$t^\kappa$ for the full covariance class, before the limiting constant is known
(Theorem~\ref{T:full-upper}). The wave equation forces this argument to carry
both displacement and velocity across time intervals. The constant of this
bound is not claimed to be optimal.

\paragraph{A lower bound for each control.}
A lower bound comes from shifting the Gaussian noise toward a deterministic
potential that makes the wave grow, and comparing the gain with the cost of the
shift. Because the noise is correlated in time, the cost of a control is
charged in the original covariance space, not added over time intervals as if
they were independent. Chapters~\ref{S:book-gaussian-foundations}--\ref{S:book-wave-evolution}
prepare the tools: covariance estimates for a potential multiplied by a smooth spatial cutoff of compact support, a
random operator path driven by the noise of the original equation, and the
evolution on the displacement--velocity phase space, identified with the
physical wave. Chapter~\ref{S:book-lower} converts long time into an
equivalent small-noise problem by equality in law
(Section~\ref{S:exact-canonical-scaling}), for nonzero covariance, proves a lower bound for each fixed control in the
class defined in Section~\ref{S:targets}, with an eligible fixed spatial mask
(Proposition~\ref{P:raw-block-original-lower}), and then removes that mask and
optimizes over controls to obtain a variational lower value
(Section~\ref{S:unmasked-variational}). For zero covariance, the solution is
identically one and the rate is zero. Together with the upper bound of the
first step, this value is finite (Proposition~\ref{P:variational-positive-scores}).

\paragraph{A matching upper bound.}
The upper bound must account for all Gaussian fluctuations, not only the
favourable shifts, and it must control the moment at one point while the region
reached by the wave grows with time. Chapter~\ref{S:book-means} relates limits of smooth controls in the norm determined by the noise covariance to the deterministic means they induce. Chapter~\ref{S:book-growing-radius} separates the estimates that are
uniform in the radius of the spatial cutoff from those whose dependence on the radius must be shown, and constructs the wave for completed static mean potentials of either sign at a fixed spatial cutoff. Chapter~\ref{S:book-point-upper} reduces the probabilistic estimate to
growth under deterministic controls, returns to the original point moment through finite propagation and a spatial average of squared solutions, controls the Gaussian tail, and obtains an upper
variational expression over signed controls
(Proposition~\ref{P:original-signed-point-upper}).

\paragraph{Identifying the constant.}
The two bounds now concern the same point moment, but over different controls:
nonnegative smooth controls below, signed controls in the covariance-norm completion just described above.
Chapter~\ref{S:book-constant} shows that the two optimizations have the same
value. It first solves a finite optimization problem whose matrix of pairwise covariance inner products has nonnegative entries (Section~\ref{S:finite-positive-gram}), then applies it to
the physical means and removes the spatial cutoff
(Proposition~\ref{P:positive-score-cutoff}). Theorem~\ref{T:intrinsic-second-moment}
states the resulting limit. The remaining sections give bounds for its constant
and, for spatial white noise in one dimension, explicit moment bounds; they do
not evaluate the full supremum.

\section{Future work}\label{S:future-work}

The second-moment argument suggests questions about higher moments, optimizing
controls and scaling limits. The following subsections separate precise
conjectures, conditional consequences and exploratory questions. The proposed
extensions require further proof.
We retain the initial data $(1,0)$ throughout this section. Except in
Subsection~\ref{S:future-rough-noise-extension}, we also retain
Hypothesis~\ref{H:homog} and the exponent $\kappa$ from~\eqref{E:kappa}.

\subsection{Higher moments and pointwise estimates}\label{S:future-higher-moments}
\begin{conjecture}[Higher even-moment asymptotics]
  \label{Conj:intrinsic-even-moments}
  Under Hypothesis~\ref{H:homog}, let $u$ be the canonical Stratonovich
  solution with initial data $(1,0)$. For each fixed even integer $p \ge 4$
  and each fixed $x \in \R^d$,
  \begin{equation}\label{E:future-even-moments}
    \lim_{t\to\infty}t^{-\kappa}\log\E|u(t,x)|^p
    =\left(\frac{p}{2}\right)^{\frac{4-\alpha}{3-\alpha}}
    C_2(\Gamma,\alpha_0).
  \end{equation}
  Here $C_2$ is the intrinsic constant in
  Theorem~\ref{T:intro-second-moment}.
\end{conjecture}

The case $p = 2$ is the theorem already stated. For zero covariance, the
conjectured identity reduces to zero on both sides. The open assertion concerns
higher even moments for nonzero covariance.

The candidate comes from the variational expression in
Theorem~\ref{T:intro-second-moment}. If a Gaussian shift produces an integrated
deterministic growth rate $I_\Psi$, its contribution to the $p$th moment should
be compared with the same Gaussian cost $v_\Psi/2$. This suggests optimizing
$pI_\Psi-v_\Psi/2$ over the full positive-control class. The homogeneity of
these variational values is already proved for every positive coefficient in
Proposition~\ref{P:variational-positive-scores}. What remains is their
identification with the higher point-moment limits. Formula
\eqref{E:future-even-moments} uses the intrinsic second-moment constant; it does
not require identifying that constant with the one predicted in the source.

Even orders allow us to use $|u|^p = u^p$ and Gaussian pairing expansions
without assuming that the solution is nonnegative. A $p$-chain expansion must
retain self-pairings, cross-pairings, unequal chain lengths and empty chains.
The scaling identity and the Gaussian change of measure offer a different
starting point that may avoid estimating all these pairings separately.
The lower-bound construction and finite Gram optimization also suggest
extensions with coefficient $p$. Higher-moment convergence, uniform
integrability and the comparison with the original point observable must
still be justified.

\paragraph{Point moments and spatial regularity.}
The main obstacle in the present upper-bound method occurs in
Proposition~\ref{P:original-point-transfer}. Stationarity and finite
propagation identify the pointwise second moment with an expected spatial
$L^2$ norm. The displacement component of the wave phase controls this norm.
For a $p$th moment, the corresponding argument would require a spatial $L^p$
bound. Taking a higher power of an $L^2$ bound does not supply it.

The existing displacement space is $H^{1/2}(\R^d)$. Its embeddings into
$L^p$ on a fixed ball give the ranges in Table~\ref{Tab:future-moments};
see~\cite{di-nezza.palatucci.ea:12:hitchhikers}*{Section~6}.
The table describes the present proof method. The conjecture covers every
admissible dimension, including orders beyond these embedding ranges.

\begin{table}[htbp]
  \centering
  \caption{Higher moments through the current wave phase. The Sobolev ranges
    are known; the proposed moment extensions still require proof.}
  \label{Tab:future-moments}
  \begin{tabular}{@{}p{0.07\textwidth}p{0.15\textwidth}p{0.20\textwidth}p{0.45\textwidth}@{}}
    \toprule
    $d$ & Available $L^p$ range & First target & Remaining step \\
    \midrule
    $1$ & $2 \le p < \infty$ & $p = 4$, then fixed even $p$ &
    Higher-moment source limits and point comparison; matching Gaussian upper
    and shifted lower bounds. \\[6pt]
    $2$ & $2 \le p \le 4$ & $p = 4$ &
    The same probabilistic steps at the critical embedding exponent. Larger
    even orders need additional spatial control. \\[6pt]
    $3$ & $2 \le p \le 3$ & A point estimate for $p = 4$ &
    A new comparison or stronger spatial estimates before applying the
    higher-moment Gaussian argument. \\
    \bottomrule
  \end{tabular}
\end{table}

The failure of this embedding route at a higher order is not a nonexistence
statement for the moment. Any stronger spatial estimate must be proved under
the original covariance assumptions. Even within the displayed ranges,
the embedding alone does not establish the conjectured limit.

\paragraph{A first test and further questions.}
The fourth moment in one spatial dimension is a useful first test. One would
first establish the fourth-moment point comparison and the required source
limits. The Gaussian upper argument would then need to control the fourth
power of the evolution norm, with bounds uniform in the growing spatial
radius. Matching this bound with the shifted lower estimate would test the
prediction $C_4 = 2^{(4-\alpha)/(3-\alpha)}C_2$. The validity of the Sobolev
embedding alone does not complete these steps.

Even a proof for every fixed even $p$ would leave two separate questions:
absolute moments of other orders, and estimates uniform as $p\to\infty$.
Interpolation between adjacent even moments need not identify the exact
constant at an intermediate order. Uniform large-$p$ estimates also require
tracking the dependence on $p$ throughout the construction, source limits and
Gaussian tail bounds. These questions may clarify which parts of the present
method apply to moment growth beyond the second moment.

\subsection{Real orders and upper-tail deviations}\label{S:future-real-orders}
A stronger question is whether~\eqref{E:future-even-moments} holds for every
fixed real $p \ge 1$. The variational homogeneity allows real coefficients,
and the Gaussian change of measure applies to absolute powers. However, the
comparison between masked and original moments in
Proposition~\ref{P:masked-canonical-transfer} uses nonnegative Gaussian pairing
terms. Such an expansion does not directly represent a noninteger absolute
power. A replacement comparison and the required moment convergence are
therefore needed. Interpolation between even orders gives bounds, but does
not in general identify the exact constant at an intermediate order. At
$p = 1$, the proposed quantity is $\E|u|$, which must be distinguished from
the signed mean $\E u$ in the source. We do not assume positivity of the wave
solution. Thus, the real-order extension is plausible from the control
formula, but remains a separate open question.

For a nonzero covariance, suppose that the real-order extension has been
established for every $p \ge 1$. For this discussion, set
\[
  m=\frac{4-\alpha}{3-\alpha}>1,\qquad
  A=2^{-m}C_2(\Gamma,\alpha_0)>0,\qquad z_*=mA.
\]
Then the limiting logarithmic moment function on this range is $Ap^m$.
The corresponding sharp upper-tail consequence is
\begin{equation}\label{E:future-upper-tail}
  \lim_{t\to\infty}t^{-\kappa}
  \log\mathbb{P}\big\{|u(t,x)|\ge e^{zt^\kappa}\big\}
  =-\mathcal{I}(z),\qquad z\ge z_*,
\end{equation}
for each fixed $x$, where
\begin{equation}\label{E:future-tail-rate}
  \mathcal{I}(z)
  =(m-1)m^{-m/(m-1)}A^{-1/(m-1)}z^{4-\alpha}.
\end{equation}
The threshold comes from the maximizing moment order
$p_z=(z/(mA))^{1/(m-1)}$. For $z>z_*$ this order lies in $(1,\infty)$,
where the hypothesized moment limits are available on both sides of $p_z$.
The exponential-tilting argument uses this local information for the lower
bound and the moment inequality for the upper bound; the endpoint follows
from thresholds above $z_*$ by monotonicity. This is the moment-method
application described in~\cite{dembo.zeitouni:98:large}*{Section~2.3}.
We do not invoke a full G\"artner--Ellis theorem on an unverified domain of
moment orders. Formula~\eqref{E:future-upper-tail} is conditional on the
real-order extension, not a large-deviation result established here for the
wave equation.

For $0<z<z_*$, optimizing only over $p\ge1$ gives the upper bound
\[
  \limsup_{t\to\infty}t^{-\kappa}
  \log\mathbb{P}\big\{|u(t,x)|\ge e^{zt^\kappa}\big\}
  \le-\max\{0,z-A\}.
\]
Those moment orders do not determine the sharp rate on this remaining range.
If the same moment formula were proved for every real $p>0$, the sharp
upper-tail expression would extend to every $z>0$. Neither version supplies
a lower-tail principle for $\log|u|$, where small values and zeros require
separate analysis. For spatial white noise in one dimension,
$m=3/2$ and the predicted upper-tail rate is
$\mathcal{I}(z)=32z^3/(27C_2^2)$; with real orders only above one, its stated
range is $z\ge3C_2/(4\sqrt2)$. The time speed and $C_2$ still depend on the
temporal covariance.

\subsection{Compactness of optimizing controls}\label{S:future-optimizers}
The variational formula raises a question about the controls responsible for
large moment growth. Its admissible tests are smooth and compactly supported,
so attainment in that class should not be presumed. We first ask for a weaker
form of compactness in the one-dimensional spatial-white case.

Fix $\Gamma = c\delta_0$, with $c > 0$, and write $a = \alpha_0 \in (0,1)$.
Let $H_a$ be the completion of real smooth compactly supported tests on
$I\times\R$, where $I = (0,1)$, modulo covariance-null tests, under the norm
\[
  \|\Phi\|_{H_a}^2
  =c\int_{I^2}|s-t|^{-a}\int_\R
  \Phi(s,x)\Phi(t,x)\ud x\ud s\ud t.
\]
For $\Psi \in \mathcal{C}_{\mathrm{pos}}$
of~\eqref{E:intro-control-class}, put
\[
  J_a(\Psi)=2\int_I\sqrt{\mathcal{R}_\Psi(s)}\ud s
  -\frac12\|\Psi\|_{H_a}^2,
  \qquad (\tau_y\Psi)(s,x)=\Psi(s,x+y).
\]
The Rayleigh value is defined in~\eqref{E:intro-control-rate}.

\begin{conjecture}[Compactness of an optimizing sequence]
  \label{Conj:positive-control-compactness}
  For every $c > 0$ and $0 < a < 1$, there are
  $\Psi_j \in \mathcal{C}_{\mathrm{pos}}$, spatial translations $y_j \in \R$
  and a nonzero $\Phi_* \in H_a$ such that
  \[
    J_a(\Psi_j)\longrightarrow C_2(c\delta_0,a),
    \qquad
    \|\tau_{y_j}\Psi_j-\Phi_*\|_{H_a}\longrightarrow0.
  \]
\end{conjecture}

The conjecture asks for one strongly convergent optimizing sequence after
spatial translation. It asserts neither compactness of every optimizing
sequence nor a smooth compactly supported maximizer. Identifying the Rayleigh
functional at the limit would be a further step. Another question is whether
controls that separate into a time factor and a space factor give the full
supremum. Comparing such controls with sums of two products is a first test;
we do not assume either equality or strict improvement.

\subsection{Sharper bounds near the static case}\label{S:future-static-correction}
The lower and upper formulas in Theorem~\ref{T:white-moment-bounds} have the
same limit as $a\downarrow0$. A more precise question concerns the first
correction to this limit.

\begin{conjecture}[First correction in the spatial-white case]
  \label{Conj:white-first-correction}
  For each fixed $c > 0$,
  \[
    \lim_{a\downarrow0}\frac1a
    \left(\frac{3C_2(c\delta_0,a)}{2\sqrt c}-1\right)
    =\frac34.
  \]
\end{conjecture}

The coefficient is suggested by formally expanding a temporally constant
trial and the temporal covariance operator around the constant kernel.
Justifying the trial approximation and a matching upper estimate remains part
of the question. The common zeroth-order limit of the existing bounds does
not establish this first-order correction.

For a fixed $0 < a < 1$, another question is whether both inequalities in
Theorem~\ref{T:white-moment-bounds} are strict. Varying the temporal equilibrium
trial may improve the lower formula. On the upper side, the norm of the
temporal covariance operator may improve the current row-integral estimate.
A separate numerical objective is a certified enclosure of $C_2$ for a fixed
$a$, with explicit errors from truncation and discretization. Optimizing a
finite family of controls alone gives lower bounds and does not certify the
continuum supremum from above.

\subsection{Logarithmic wavefields and possible scaling limits}
\label{S:future-universality}

A logarithmic change of variables provides another way to view random wave
propagation. This idea has classical precedents in the Rytov transformation,
which represents a complex wavefield through its logarithmic amplitude and
phase. For example, Sung et al.~\cite{sung.choi.ea:09:optical}*{Section~2}
substitute $U=e^\psi$ into a Helmholtz equation before applying a
weak-scattering approximation. The substitution and the subsequent
approximation are distinct steps. That frequency-domain setting does not
supply a scaling limit for the time-dependent model considered here.

For the smooth equation $u_{tt}=\Delta u+Vu$, set $h=\log|u|$ on a region where
$u$ is nonzero. The ordinary chain rule gives
\begin{equation}\label{E:future-log-wave}
  h_{tt}+(h_t)^2=\Delta h+|\nabla h|^2+V.
\end{equation}
Thus, the logarithmic wave equation retains both a second time derivative and
a quadratic temporal derivative. In the stochastic model, this calculation
applies first to smooth noise approximations. Passing through zeros or to the
singular-noise limit requires a separate argument; we do not assert a global
logarithmic equation for the canonical solution.

A possible mechanism for an effective first-order time evolution can be seen
formally by writing $h=b(t)+\varphi$, where $b$ is a smooth deterministic
centering. Then
\begin{equation}\label{E:future-centered-log-wave}
  \varphi_{tt}+2b'(t)\varphi_t+(\varphi_t)^2
  =\Delta\varphi+|\nabla\varphi|^2+V-b''(t)-(b'(t))^2.
\end{equation}
If $b'$ approached a positive constant and the other temporal terms on the
left became negligible on a fluctuation scale, the term $2b'\varphi_t$ could
produce an effective parabolic evolution. Neither condition follows from
moment growth. In particular, logarithmic moments need not give the typical
centering, and the long temporal memory of the noise must still be accounted
for. This is a heuristic mechanism to investigate.

The ordinary one-dimensional KPZ comparison concerns the centered logarithm
of the positive stochastic heat solution. Its height, space and time scaling
exponents have ratio $1:2:3$; convergence to the KPZ fixed point is a statement about
random fields and their evolution~\cite{wu:26:kpz}*{Theorems~1.6 and~1.8}. The moment normalization
in the present wave problem controls rare amplification. It does not identify
the typical fluctuation size or spatial correlation length. The cubic
upper-tail prediction in the spatial-white case also gives no evidence by
itself for KPZ behavior: the usual KPZ upper-tail power is $3/2$, while a
cubic power appears in a lower-tail regime~\cite{corwin.ghosal:20:kpz}.
Different rare-event scales do not, by themselves, exclude a common limit in
another centered fluctuation regime.

The open question is whether a specified logarithmic amplitude observable
$H$ admits deterministic functions $b$, $S$ and $L$ for which
\[
  \frac{H(Tt,L(T)x)-b(Tt)}{S(T)}
\]
converges to a nondegenerate random field as $T\to\infty$. A logarithm of a
local spatial $L^2$ amplitude is one possible observable that reduces
sensitivity to isolated pointwise zeros. It can itself vanish, so its domain
must also be specified; it differs from $\log|u|$. The scales and topology must
be specified, and a claim of universality would require stability under
changes of regularization and microscopic covariance details. Temporal
correlations and sensitivity to the initial velocity would help determine
whether a scalar amplitude describes the limiting evolution or whether
additional variables remain necessary. We leave this as an exploratory
question, without prescribing exponents or a new universality class.

\subsection{Rough noise beyond the present covariance class}
\label{S:future-rough-noise-extension}

We ask whether we can construct the Stratonovich wave equation with initial
data $(1,0)$, and determine its moment growth, beyond
Hypothesis~\ref{H:homog}.
We use $d = 1$ to illustrate these questions and consider Gaussian noises with
product spectral measure $\mu_0\otimes\mu$, where
\[
  \mu_0(\ud\tau) = c_0|\tau|^{\alpha_0-1}\ud\tau,
  \qquad
  \mu(\ud\xi) = c_1|\xi|^{\alpha-1}\ud\xi,
  \qquad c_0, c_1 > 0,
\]
with $\alpha_0, \alpha > 0$ and $\alpha_0 + \alpha < 2$. The corresponding Hurst
parameters are $H_0 = 1 - \alpha_0/2$ and $H = 1 - \alpha/2$. Thus, the
proposed condition is $H_0 + H > 1$. We propose this region for investigation,
without asserting a solvability criterion. The
physical covariance in the rough regions is understood as a distribution and
need not be a nonnegative Radon measure.

After matching the covariance normalization, the current theorem covers
$0 < \alpha_0 < 1$, $0 < \alpha < 1$, and the spatial white noise segment
$\alpha = 1$, $0 < \alpha_0 < 1$.
Figure~\ref{Fig:rough-noise-region} separates this range from the two proposed
regions: rough space with $\alpha > 1$, and rough time with $\alpha_0 > 1$. The
temporal white noise segment $\alpha_0 = 1$, $0 < \alpha < 1$ is outside the
current theorem. Its integration convention and its relation to the present
variational formula require separate treatment. On either white noise line
we mean a Dirac covariance with
its specified normalization, not the literal function $|r|^{-1}$. The point
$(\alpha,\alpha_0) = (1,1)$ lies on the excluded boundary $\alpha_0 + \alpha =
2$; this exclusion is not a nonexistence statement.

\usetikzlibrary{patterns}
% Requires TikZ libraries arrows.meta and patterns.
\begin{figure}[htbp]
  \centering
  % Palette matches the book's hamblue/hamgold without changing shared styles.
\begingroup%
\definecolor{rnBlue}{RGB}{36,75,105}%
\definecolor{rnGold}{RGB}{156,101,20}%
\definecolor{rnGray}{RGB}{105,110,115}%
\begin{tikzpicture}[
  x=3.8cm,y=3.8cm,
  font=\small,
  region/.style={align=center,fill=white,fill opacity=.94,
  text opacity=1,inner sep=2pt,rounded corners=1pt},
  legend/.style={anchor=west,align=left,text width=4.4cm},
  axis/.style={-{Latex[length=2.1mm]},line width=.7pt},
  threshold/.style={draw=rnGray!45,densely dotted,line width=.8pt}
  ]
  % Three disjoint interiors within the proposed open triangle.
  \fill[rnBlue!16] (0,0) rectangle (1,1);
  \fill[rnGold!7] (1,0)--(2,0)--(1,1)--cycle;
  \fill[pattern=north east lines,pattern color=rnGold!38]
  (1,0)--(2,0)--(1,1)--cycle;
  \fill[rnGold!7] (0,1)--(1,1)--(0,2)--cycle;
  \fill[pattern=north west lines,pattern color=rnGold!38]
  (0,1)--(1,1)--(0,2)--cycle;

  % White-noise thresholds continue across the plot; exterior is subdued.
  \draw[threshold] (1,1)--(1,2.2);
  \draw[threshold] (1,1)--(2.2,1);
  \draw[rnBlue,line width=1.7pt] (1,0)--(1,1);
  \draw[rnGold,line width=1.35pt,dash pattern=on 5pt off 3pt]
  (0,1)--(1,1);
  \draw[rnGray,line width=1pt,dash pattern=on 4pt off 3pt]
  (0,2)--(2,0)
  node[pos=.25,above=6pt,sloped,text=rnGray]
  {$\alpha + \alpha_0 = 2$};

  \draw[axis] (0,0)--(2.2,0) node[right=5pt] {$\alpha$};
  \draw[axis] (0,0)--(0,2.2) node[above=5pt] {$\alpha_0$};
  \node[below left=3pt] at (0,0) {$0$};
  \foreach \v in {1,2} {
  \draw (\v,0)--++(0,-.025) node[below=4pt] {$\v$};
  \draw (0,\v)--++(-.025,0) node[left=4pt] {$\v$};
  }

  \node[region,text=rnBlue,text width=2.75cm] at (.49,.46)
  {\textbf{Current theorem}\\regular region\\
  $0 < \alpha < 1$, $0 < \alpha_0 < 1$};
  \node[region,text=rnGold,text width=2.4cm] at (1.35,.17)
  {\mbox{Rough space}\\\mbox{open question}};
  \node[region,text=rnGold,text width=2.4cm] at (.35,1.16)
  {\mbox{Rough time}\\\mbox{open question}};
  \node[text=rnGray,align=center,text width=3.0cm] at (1.67,1.79)
  {Rough in both\\time and space\\outside the proposed region};
  \node[text=rnGray,anchor=west,align=left,text width=3.0cm]
  at (1.13,1.20) {Both white $(1,1)$\\excluded};
  \draw[rnGray,line width=.6pt] (1.02,1.02)--(1.10,1.12);

  % Open markers: neither axes nor the diagonal belong to the range.
  \foreach \p in {(0,0),(0,1),(1,0),(0,2),(2,0)} {
  \filldraw[fill=white,draw=rnGray,line width=.7pt] \p circle[radius=2.2pt];
  }
  \filldraw[fill=white,draw=rnGray,line width=1.2pt]
  (1,1) circle[radius=3.4pt];

  % Legend uses fill, hatching and line style as well as color.
  \node[legend,font=\normalsize\bfseries,text=rnBlue] at (2.39,2.02)
  {Current theorem};
  \filldraw[fill=rnBlue!16,draw=rnBlue!50]
  (2.39,1.745) rectangle (2.57,1.855);
  \node[legend] at (2.66,1.80) {Regular region};
  \draw[rnBlue,line width=1.7pt] (2.39,1.52)--(2.57,1.52);
  \node[legend] at (2.66,1.52)
  {White noise in space\\$\alpha = 1$, $0 < \alpha_0 < 1$};

  \node[legend,font=\normalsize\bfseries,text=rnGold] at (2.39,1.20)
  {Extension questions};
  \filldraw[fill=rnGold!7,draw=rnGold!60]
  (2.39,.865) rectangle (2.57,.975);
  \fill[pattern=north east lines,pattern color=rnGold!55]
  (2.39,.865) rectangle (2.57,.975);
  \node[legend] at (2.66,.92) {Rough space: $\alpha > 1$\\rough time: $\alpha_0 > 1$};
  \draw[rnGold,line width=1.35pt,dash pattern=on 5pt off 3pt]
  (2.39,.58)--(2.57,.58);
  \node[legend] at (2.66,.58)
  {White noise in time\\Endpoint question};

  \draw[rnGray,line width=1pt,dash pattern=on 4pt off 3pt]
  (2.39,.29)--(2.57,.29);
  \node[legend] at (2.66,.29) {Excluded diagonal};
  \draw[threshold] (2.39,.07)--(2.57,.07);
  \node[legend] at (2.66,.07) {White noise lines\\(continued)};
\end{tikzpicture}%
\endgroup%
  \caption[Regular, white and rough noise in one dimension]{Parameter schematic for $d = 1$. The candidate range for
    investigation is $\alpha > 0$, $\alpha_0 > 0$, $\alpha + \alpha_0 < 2$. Blue
    marks the current theorem under its standing covariance assumptions: the
    open regular region, together with the segment white in space, $\alpha = 1$
    and $0 < \alpha_0 < 1$. The hatched rough regions and the dashed segment
    white in time, $\alpha_0 = 1$ and $0 < \alpha < 1$, are questions about
    extending the present results. The axes and the diagonal are excluded,
    including $(1,1)$. The region rough in both time and space does not meet the
    triangle. No solvability claim is made outside the triangle. White in space
    or time means that the corresponding covariance factor is a multiple of a
    Dirac mass. For $d \ge 2$, white noise in space lies on $\alpha = d$,
    outside the open positive candidate triangle.}
  \label{Fig:rough-noise-region}
\end{figure}
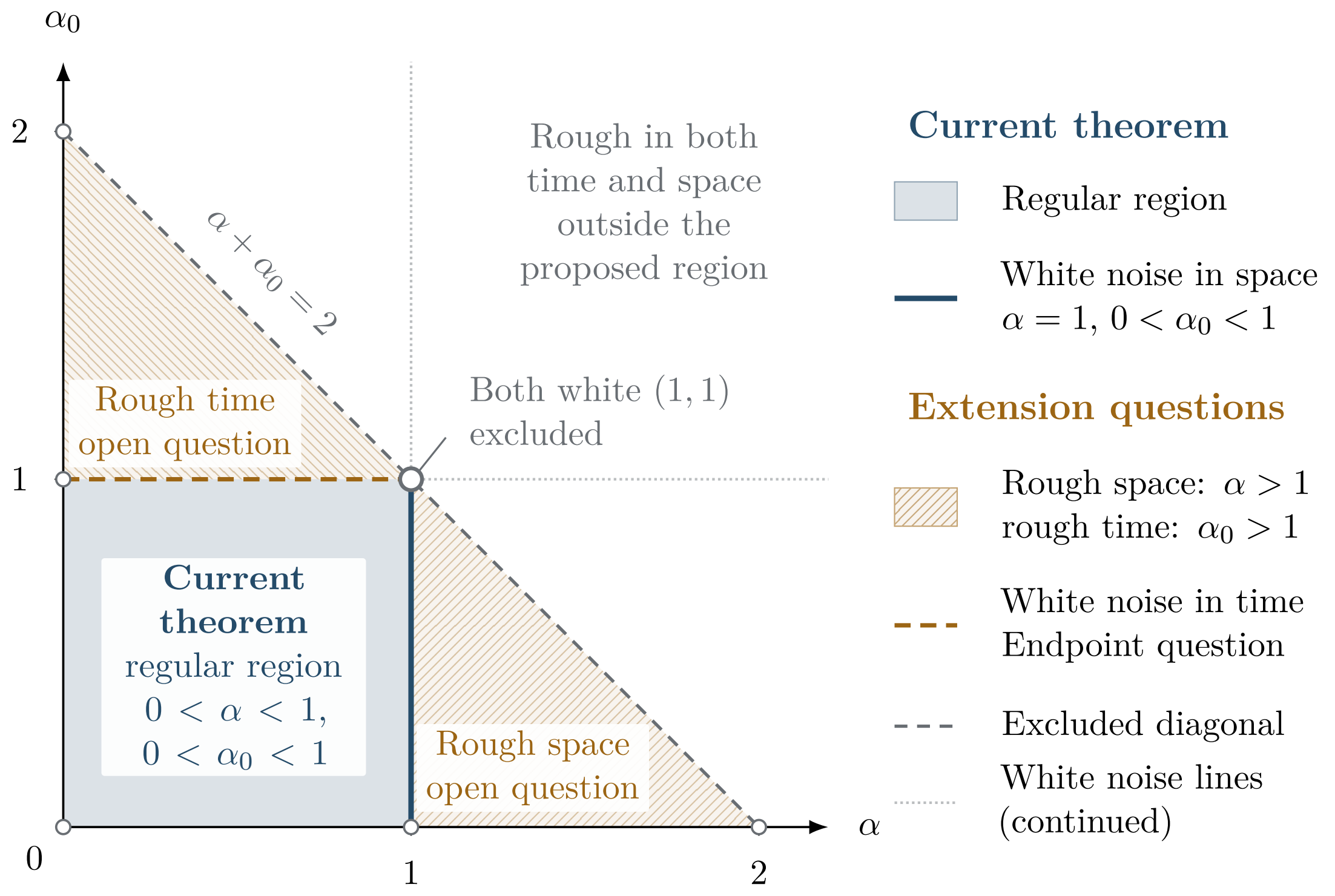

Rough spatial covariance is already used in related equations. Balan, Chen and
Ma~\cite{balan.chen.ea:22:parabolic}*{Theorem~1.1} construct the Skorohod
solution of the parabolic Anderson model for their class of initial data,
with $H_0 \in (1/2, 1)$, $0 < H < 1/2$ and
$H_0 + H > 3/4$. In the notation above, this gives $0 < \alpha_0 < 1$, $1 <
\alpha < 2$ and $\alpha_0 + \alpha < 5/2$. Their frequency estimates and bounds
for ratios of Gamma functions suggest methods to investigate. Their estimates
for the heat kernel and orthogonal chaos expansion do not directly treat the
oscillatory wave kernel or the self-pairings in an ordinary Stratonovich
expansion. Related results for wave equations with rough spatial noise are
given by Song, Song and Xu~\cite{song.song.ea:20:fractional}*{Theorem~3.2} and
Liu, Hu and Wang~\cite{liu.hu.ea:22:nonlinear}; their parameter ranges and
definitions of the stochastic integral must also be retained.

The present proof uses nonnegative physical covariance in its pairing bounds
and comparisons of masked moments. It also uses Young integration to construct
the wave evolution. These arguments require further work for the proposed
spectral model. Their limitations do not by themselves imply that the model is
ill posed. At each regularization, noise enters the velocity equation, while
displacement is its time integral. This structure may permit a different
construction. Even after a construction, the sharp moment asymptotics would
remain a separate question. The rough temporal range requires particular
care. The formal continuation of the expression in~\eqref{E:kappa} gives
\[
  \frac{4-\alpha-\alpha_0}{3-\alpha}
  = 1+\frac{1-\alpha_0}{3-\alpha}.
\]
It equals one at $\alpha_0 = 1$ and falls below one when $\alpha_0 > 1$.
For $\alpha_0 \ge 1$, the factor $\delta^{(1-\alpha_0)/2}$ used in the
current block error estimate no longer tends to zero. This factor still
vanishes in the proposed rough spatial range with $\alpha_0 < 1$. We leave
both the appropriate moment scale and its limiting constant open in these
extensions.

\subsection{Exact first-moment asymptotics}
\label{S:future-first-moment}

For a nonzero spatial covariance and a fixed positive temporal exponent
$\alpha_0$, can we determine the exact large-time growth of the signed mean
$\E u(t,x)$? This question is distinct
from determining the growth of $\E|u(t,x)|$. The current author version of
Xia Chen's paper gives two-sided bounds for the signed mean and predicts
that the upper bound is sharp;
see~\cite{chen:author2026:time-dependency}*{Theorem~1.2,
(1.10), the discussion following (1.14), and Remark~8.1}. The second-moment
result of this book does not identify either first-moment constant.

The positive-control variational problem provides a quantity to investigate.
Proposition~\ref{P:variational-positive-scores} and
Theorem~\ref{T:intrinsic-second-moment} give
\[
  S_1^{\mathrm{pos}}=2^{-m}C_2,
  \qquad m=\frac{4-\alpha}{3-\alpha}.
\]
This identity concerns control suprema. It does not prove that either
$t^{-\kappa}\log\E u(t,x)$ or $t^{-\kappa}\log\E|u(t,x)|$ converges to
$S_1^{\mathrm{pos}}$. Before proposing such an identification, we must check it
against the actual mean bounds, the static limit and the source's covariance
normalization. A variational characterization of the first mean may require
a different functional.

The main missing lower-bound argument concerns cancellation. The second-moment
proof restricts a nonnegative squared observable to a favorable event. For a
signed observable, the contribution from the complementary event cannot be
discarded. Positivity of the unconditional Wick coefficients does not supply
this comparison. Absolute values remove that particular obstruction, but a
first absolute-moment argument still needs a comparison between the localized
observable and the original solution, together with a matching upper bound.
The completed quadratic comparisons do not provide these statements.

There is a useful boundary calibration. At $\alpha_0=0$, Chen--Hu's static
first-mean asymptotic is an imported result; see~\eqref{E:static-ratio-input}.
Corollary~\ref{C:static-moment-constant-continuity} shows that the signed
first-mean lower and upper growth rates approach that constant as
$\alpha_0\downarrow0$, after the large-time liminf and limsup have been taken.
It does not establish an exact time limit at a fixed positive $\alpha_0$.
A first step is to test a proposed first-moment formula on the small-order
pairing coefficients and this static calibration, then isolate the needed
cancellation or localization estimate. We leave the existence and variational
identification of both first-moment limits open.

\section{Toward a general variational method for random evolution equations}
\label{S:theory-outlook}

The preceding section asks how to extend the results. Here we ask what the
proof itself teaches us. Several long arguments preserve the same three
objects: the displacement--velocity evolution, the covariance of the original
noise, and a variational reward for deterministic amplification. Their
connection suggests a possible method for other random evolution equations.
These questions arose by asking Amy to revisit the calculations and failed
approaches, then checking her proposed abstractions against the arguments.
The present book supplies a detailed case study; a general theorem would
require assumptions that can be verified independently in a second model.

\Needspace{5\baselineskip}
\paragraph{Compose the evolution and retain the joint Gaussian cost.}
For a wave equation, displacement at the end of a time block does not determine
the next block: velocity is also needed. Lemma~\ref{L:raw-block-restart}
composes the full phase-space propagators and identifies them through limits
of the same regularized noise. The Gaussian change of
law in Proposition~\ref{P:raw-block-gaussian-law} acts on that common noise.
It does not replace correlated blocks by independent copies.

Even two Gaussian coordinates show why this distinction matters. If their
covariance matrix is
\[
  C_\rho=\begin{pmatrix}1&\rho\\\rho&1\end{pmatrix},\qquad |\rho|<1,
\]
then a shift by $h=(q,q)$ has Cameron--Martin cost
\[
  \frac12 h^{\mathsf T}C_\rho^{-1}h=\frac{q^2}{1+\rho}.
\]
Indeed, $C_\rho h=(1+\rho)h$ and $|h|^2=2q^2$. The cost is $q^2$ for
independent coordinates and tends to $q^2/2$ as $\rho\uparrow1$.
Discarding correlation can therefore change the leading variational cost.
The Gaussian shift formula is classical~\cite{janson:97:gaussian}*{Theorem~14.1};
the work here is to connect it to the physical wave evolution and its terminal
observable, with estimates uniform in the necessary limits.

This suggests a precise task: find sufficient conditions under which a
phase-space evolution, controlled by one Gaussian field, has matching
long-time upper and lower moment bounds with the same deterministic reward
and the same covariance cost. The conditions must include convergence from
the original regularized equation, control of products on growing intervals,
and passage from phase-space estimates to the chosen observable.
Proposition~\ref{P:clock-copy-obstruction} rules out one unweighted
constant-correlation Brownian-clock substitution; it does not rule out every
probabilistic representation. The pathwise composition used here supplies
neither a conditional-expectation identity nor a general martingale calculus.

\paragraph{Separate the finite optimization from the wave estimates.}
Propositions~\ref{P:gram-positive-optimizer}
and~\ref{P:gram-entrywise-comparison} already apply in any real Hilbert space.
For the square-root reward in~\eqref{E:finite-gram-score}, nonnegative Gram
entries permit a positive conic maximizer and an entrywise comparison of
optimized values. Singular matrices and inactive labels are included. This
is the most explicit piece of the argument that can be stated without wave
notation. Its physical application still needs the covariance and spectral
estimates in Section~\ref{S:intrinsic-second-moment}.

A useful test would replace the square-root reward by the growth rate of a
damped equation. For the scalar equation $u''+2\eta u'=vu$, with $\eta>0$
and $v\ge0$, substituting $u(t)=e^{rt}$ gives
$r^2+2\eta r-v=0$. Its larger root is
\[
  r_\eta(v)=\sqrt{\eta^2+v}-\eta.
\]
Unlike $\sqrt v$, this reward has a finite derivative at zero and is not
homogeneous in $v$ when $\eta$ is fixed. The active-set analysis and dilation
argument therefore need to be checked again. A damped-wave transfer would
also require new bounds for the physical propagator and noise approximation;
the scalar calculation alone supplies none of these bounds.

\paragraph{From an optimized value to the shape of large events.}
The constant $C_2$ identifies an optimized reward, but does not show that the
noise or solution, conditioned on large amplification, is close to an
optimizing profile. Conjecture~\ref{Conj:positive-control-compactness} asks
for compactness of one optimizing sequence after spatial translation.
A further question is whether a suitable conditioned or moment-weighted law
concentrates near such profiles. This needs a specified observable and
scaling, tightness of the relevant laws, and upper bounds away from the
optimizing set. Neither the moment limit nor the compactness conjecture alone
provides these ingredients. This question concerns the shape of rare large
events; its relation to typical fluctuations remains open, as discussed in
Section~\ref{S:future-universality}.

The first test of a broader method should be a second model for which one
abstract statement removes a substantial part of the proof. A damped wave
with correlated Gaussian forcing is a candidate. The Gaussian geometry and
finite optimization should be separated from the estimates specific to its
propagator. If those estimates require rebuilding nearly every argument,
then the current work offers an instructive organization of a proof, while
the proposed general method remains to be found.

\chapter{Construction and two-chain estimates}\label{S:G0}

\begin{figure}
  \centering
  \begin{tikzpicture}[x=1cm,y=1cm]
    \foreach \hamPanelOffset/\hamPanelTitle in {0/{Cross-pairs},4.25/{Self-pairs},8.5/{Unequal orders}} {
    \node[ham/note] at (\hamPanelOffset+1.7,1.05) {\hamPanelTitle};
    }
    \begin{scope}
      \node[ham/dot] (a1) at (.6,0) {};
      \node[ham/dot] (a2) at (2.8,0) {};
      \node[ham/dot] (b1) at (.6,-1.05) {};
      \node[ham/dot] (b2) at (2.8,-1.05) {};
      \draw[ham/arrow] (a1)--(a2);
      \draw[ham/arrow] (b1)--(b2);
      \draw[ham/pair] (a1)--(b2);
      \draw[ham/pair] (a2)--(b1);
      \node[ham/note] at (1.7,-1.92) {$m=n=2$};
    \end{scope}
    \begin{scope}[xshift=4.25cm]
      \node[ham/dot] (c1) at (.6,0) {};
      \node[ham/dot] (c2) at (2.8,0) {};
      \node[ham/dot] (d1) at (.6,-1.05) {};
      \node[ham/dot] (d2) at (2.8,-1.05) {};
      \draw[ham/arrow] (c1)--(c2);
      \draw[ham/arrow] (d1)--(d2);
      \draw[ham/pair] (c1) to[bend left=32] (c2);
      \draw[ham/pair] (d1) to[bend right=32] (d2);
      \node[ham/note] at (1.7,-1.92) {$m=n=2$};
    \end{scope}
    \begin{scope}[xshift=8.5cm]
      \node[ham/dot] (e1) at (1.7,0) {};
      \node[ham/dot] (f1) at (.35,-1.05) {};
      \node[ham/dot] (f2) at (1.7,-1.05) {};
      \node[ham/dot] (f3) at (3.05,-1.05) {};
      \draw[ham/arrow] (f1)--(f2);
      \draw[ham/arrow] (f2)--(f3);
      \draw[ham/pair] (e1)--(f2);
      \draw[ham/pair] (f1) to[bend right=28] (f3);
      \node[ham/note] at (1.7,-1.92) {$m=1,\ n=3$};
    \end{scope}
    \draw[ham/arrow] (.5,-2.65)--(1.25,-2.65);
    \node[anchor=west,font=\small] at (1.4,-2.65) {order within a chain};
    \draw[ham/pair] (6.4,-2.65)--(7.15,-2.65);
    \node[anchor=west,font=\small] at (7.3,-2.65) {one covariance pairing};
  \end{tikzpicture}
  \caption[Pairings of two ordered chains]{Three permitted pairing
    contributions. Each chain has its own time order; no common interleaving is
    imposed by the picture. Stratonovich products include self-pairs and unequal
    orders, as well as cross-pairs. These examples do not list every pairing.
    See \ProofRef{Proposition}{P:mixed}.}
  \label{Fig:two-chain-pairings}
\end{figure}
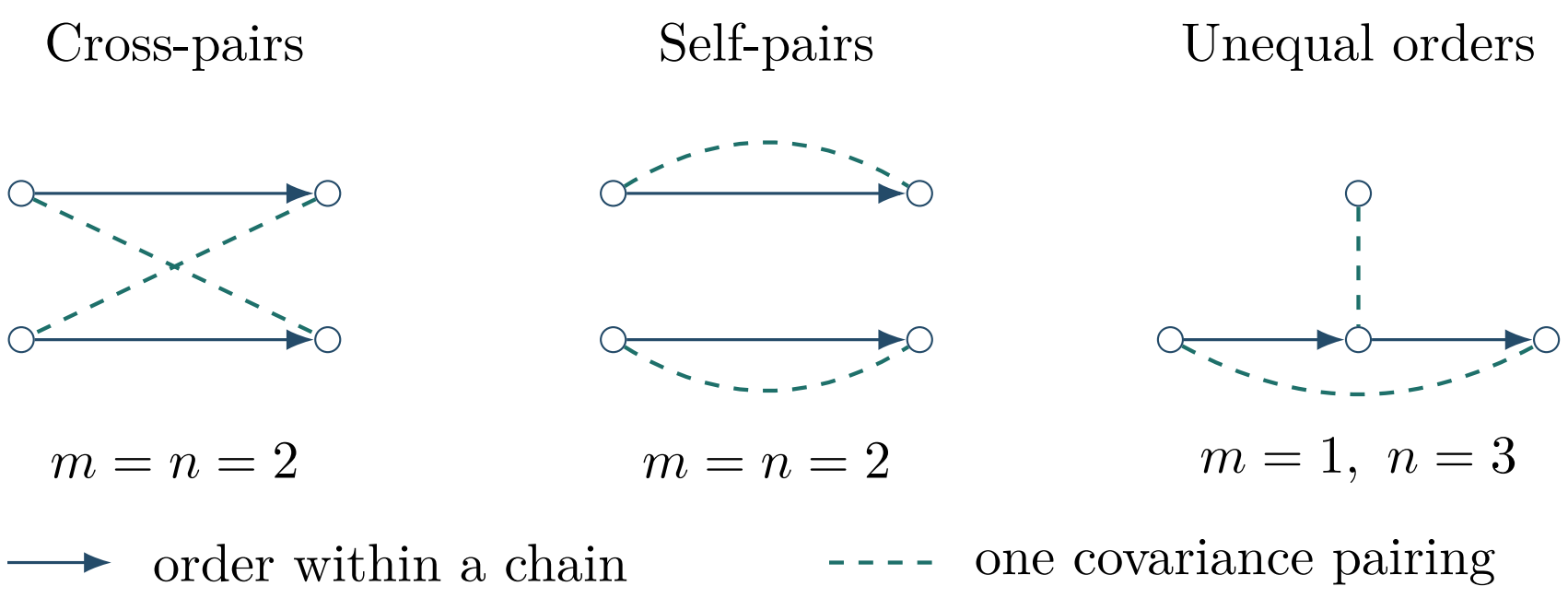

Which solution do we mean, and which terms does its second moment contain? A
Stratonovich integral keeps the pairings inside one ordered chain as well as
the pairings between two chains (Figure~\ref{Fig:two-chain-pairings}). So
products of different orders need not have zero expectation, and dropping
them would change the moment under study.

We work with the model of Section~\ref{S:model}. Proposition~\ref{P:mixed}
gives the complete formula for the mixed
moments. Lemma~\ref{L:joint-fixed-order} constructs each fixed-order term. The
summable bounds of Lemma~\ref{L:vertex-crude} then construct their series, and
Proposition~\ref{P:mild-construction} identifies it as a mild solution. These
results give the construction and mixed-moment formula used later.

The last section asks how far one common bound for every pairing can go. It
still loses the correct growth scale (Proposition~\ref{P:envelope}).
Chapters~\ref{S:book-independent-upper}--\ref{S:book-constant}
do not use that proposition, so it can wait on a first
reading. Episode~\VideoRef{3}{Glay6Q_jsLM} of the video series follows
this chapter.

Put $X_n(t,x) = S_n(g_n(\cdot,t,x))$, with $X_0 = 1$, and
\[
  \beta=4-\alpha-\alpha_0>2,\qquad
  P_k=\#\Pi_k=\frac{(2k)!}{2^k k!}.
\]
The count follows by pairing the smallest remaining label with one of $2k-1$
labels: $P_0 = 1$ and $P_k = (2k-1)P_{k-1}$, which gives the displayed formula
by induction. All constants below may depend on the fixed covariance and its
parameters, but not on the orders or on $t$.

\section{The complete two-chain expansion}

For $m+n = 2k$, label the two chains by $I_1 = \{1,\ldots,m\}$ and $I_2 =
\{m+1,\ldots,m+n\}$. A self-pairing has both endpoints in one chain; a
cross-pairing has one endpoint in each chain. Both are included in $\Pi_k$. Let
$p(l)$ denote the preceding vertex in its own chain; at each first vertex set
$(s_{p(l)},y_{p(l)}) = (0,0)$. Thus, the second chain has its own root, rather
than the last vertex of the first chain. For $D \in \Pi_k$, set
\begin{equation}\label{E:mixed-diagram}
  \begin{split}
    J_{m,n,D}(t)
    ={}&\int_{[0,t]^m_<\times[0,t]^n_<}\ud\mathbf{s}
    \int_{(\R^d)^{m+n}}\ud\mathbf{y}\,
    \prod_{\rho=1}^2\prod_{l\in I_\rho}
    G(s_l-s_{p(l)},y_l-y_{p(l)})\\
    &\hspace{20mm}\times
    \prod_{(j,h)\in D}|s_j-s_h|^{-\alpha_0}\gamma(y_j-y_h).
  \end{split}
\end{equation}
An empty chain contributes the factor one. For $d = 3$, the spatial Green
kernels are measures; for spatial white noise the covariance is a delta measure.
Integrals in~\eqref{E:mixed-diagram} use the generalized-kernel and
regularization conventions of the source, not pointwise values of these
measures. Lemmas~\ref{L:spatial-wick-measures}--\ref{L:joint-fixed-order}
construct them with the convention
of~\cite{chen:author2026:time-dependency}*{Remark~2.3, p.~11, and the Green
regularization on pp.~13--14}.

\begin{proposition}[Mixed terms]\label{P:mixed}
  Under Hypothesis~\ref{H:homog},
  \begin{equation}\label{E:mixed-moment}
    \E[X_m(t,x)X_n(t,x)]
    =\begin{cases}
      \displaystyle\sum_{D\in\Pi_k}J_{m,n,D}(t),&m+n=2k,\\
      0,&m+n\text{ odd},
    \end{cases}
  \end{equation}
  where the case $m = n = 0$ equals one. In particular,
  \begin{equation}\label{E:full-second}
    \E u^2(t,x)=\sum_{k=0}^\infty\sum_{m=0}^{2k}
    \sum_{D\in\Pi_k}J_{m,2k-m,D}(t).
  \end{equation}
  This is a convergent series of nonnegative terms.
\end{proposition}

\begin{proof}
  We first pass the finite Wick formula to the fixed-order limits. We then
  reverse the two terminal chains separately and identify their initial-chain
  diagrams. Finally, norm summability and nonnegativity justify the complete
  second-moment series.

  \medskip
  \noindent\textit{Step 1.~}We identify the fixed-order limits without using
  the infinite series.
  Lemma~\ref{L:joint-fixed-order} below constructs each $X_n$ and each
  individual two-chain diagram at arbitrary vertexwise widths, including unequal
  and empty chains. Its proof uses the positive spatial-measure limit and
  Lemma~\ref{L:pairing-product-envelope}; it uses neither the present
  proposition nor any infinite solution series.

  Here $Y_m,Y_n$ denote the positive-width terminal integrals
  $Y_m^{\boldsymbol a},Y_n^{\boldsymbol a}$ with the respective vertexwise
  widths suppressed; their joint width limit is the $L^2$ limit in
  Lemma~\ref{L:joint-fixed-order}.
  At positive noise widths Wick's formula gives zero at odd total order and the
  finite pairing sum at even total order. The Green-width limit is taken first,
  as in the construction preceding Lemma~\ref{L:spatial-wick-measures}. At fixed
  $m,n$, joint $L^2$ convergence of the two factors gives product convergence in
  $L^1$:
  \[
    \|Y_mY_n-X_mX_n\|_1
    \le\|Y_m-X_m\|_2\|Y_n\|_2
    +\|X_m\|_2\|Y_n-X_n\|_2\longrightarrow0.
  \]
  A zero-order factor stays equal to one. Each individual Wick diagram also has
  the finite nonnegative limit in Lemma~\ref{L:joint-fixed-order}.

  \medskip
  \noindent\textit{Step 2.~}We reverse the two chains and relabel the
  pairings.
  To identify its initial-chain form, for terminal variables in the first chain
  set
  \[
    s'_j=t-s_{m+1-j},\qquad y'_j=x-y_{m+1-j},\qquad1\le j\le m,
  \]
  and use the corresponding length-$n$ transformation on the second chain. Each
  transformed chain has initial root $(0,0)$. For example, $G(t-s_m,x-y_m)$
  becomes $G(s'_1,y'_1)$; every other factor becomes the corresponding
  successive increment. The differences in each covariance change sign, absorbed
  by evenness. The vertex involution is
  \[
    \rho(j)=m+1-j\ (1\le j\le m),\qquad
    \rho(m+j)=m+n+1-j\ (1\le j\le n).
  \]
  It induces the pairing bijection $D\mapsto\{\{\rho(a),\rho(b)\}:\{a,b\}\in
  D\}$. The same transformation holds before removing noise smoothing. Passing
  to the finitely many diagram limits proves \eqref{E:mixed-moment} with
  precisely~\eqref{E:mixed-diagram}.

  \medskip
  \noindent\textit{Step 3.~}We sum the mixed expectations using absolute norm
  convergence.
  Lemmas~\ref{L:crude-corrected} and~\ref{L:vertex-crude} independently prove
  fixed-$t$ norm summability of the $X_n$, using the fixed-order construction
  just cited and Lemma~\ref{L:pairing-product-envelope}.
  Proposition~\ref{P:mild-construction} identifies their series as a
  Stratonovich mild solution. None of those proofs uses the present mixed-term
  identity. Hence, $U_N = \sum_{n=0}^N X_n \to u$ in $L^2$, and
  \[
    |\E U_N^2-\E u^2|
    \le\|U_N-u\|_2(\|U_N\|_2+\|u\|_2)\longrightarrow0.
  \]
  Each finite square expands over the square of indices $0\le m,n\le N$.

  The finite expansions of $\E U_N^2$ increase, since every mixed expectation
  computed above is nonnegative. Their limit is the full double sum, independent
  of its enumeration. Regrouping by total order gives~\eqref{E:full-second}.
  This proves Proposition~\ref{P:mixed}.
\end{proof}

The two chain orders must remain distinct: chronological concatenation can
violate wave support, and a difference of block labels need not measure
chronological separation. No-gos~\ref{P:nogo-05} and~\ref{P:nogo-07}
give explicit positive-measure witnesses.

\section{A small-order check and a source uniformity issue}

\LeanMarker{L:white-test}
\begin{lemma}[White-space test]\label{L:white-test}
  For $d = \alpha = 1$ and $\gamma = \delta_0$, put $c_0 =
  ((1-\alpha_0)(2-\alpha_0)(3-\alpha_0))^{-1}$. Then
  \begin{equation}\label{E:white-test}
    \E X_1(t,x)^2=c_0t^{3-\alpha_0},\qquad
    \E X_2(t,x)=\tfrac12c_0t^{3-\alpha_0}.
  \end{equation}
  In particular, the total-order-two contribution to~\eqref{E:full-second} is
  $2c_0t^{3-\alpha_0}$.
\end{lemma}

\begin{proof}
  Here $G(s,y) = \frac12\mathbf1_{\{|y|<s\}}$. Therefore, $\int_\R
  G(s,y)G(r,y)\ud y = \min(s,r)/2$. The $(1,1)$ integral is
  \[
    \int_0^t\ud s\int_0^s r(s-r)^{-\alpha_0}\ud r
    =\frac{t^{3-\alpha_0}}
    {(1-\alpha_0)(2-\alpha_0)(3-\alpha_0)}.
  \]
  For the $(0,2)$ integral, the delta covariance sets the two spatial vertices
  equal. Write $0 < r < s < t$ for their ordered time variables. The first
  factor $G_r$ has spatial mass $r$, while $G_{s-r}(0) = 1/2$. Its time integral
  is half the preceding display. The $(2,0)$ term is identical. This proves
  Lemma~\ref{L:white-test}.
\end{proof}

\begin{remark}[Source fidelity]\label{R:source-uniformity}
  The author version's Remark~5.1 refers to (4.4), but (4.4), on p.~22, is a
  smoothing limit. The relevant bound is (5.4), on p.~35, $\E X_n(t,x)^2 \le
  C^nt^{2n}/n!$, with the preceding sentence asserting that $C$ is independent
  of $n$ and $t$. Taken for all $t > 0$, this uniformity
  contradicts~\eqref{E:white-test} at $n = 1$. The proof on pp.~34--35 sets
  $\theta = n/t$. Its first-norm estimate uses (A.1), which is finite for every
  positive $\theta$ but does not give a constant uniform as $\theta\downarrow0$.
  Its second-norm estimate uses the large-$\theta$ bound (A.3) without retaining
  that restriction. Both losses of uniformity matter when $n/t$ is small. This
  is a source-internal uniformity discrepancy; it does not refute Remark~5.1's
  final large-time conclusion. The next lemma repairs this particular deduction
  without importing the asserted uniformity of (5.4).
\end{remark}

\section{Recovering the quadratic-exponential bound}

Write $G_\theta(s,y) = e^{-\theta s}G(s,y)$ and $\widehat{G}_s(\xi) =
\sin(s|\xi|)/|\xi|$, with value $s$ at $\xi = 0$. We use the norms in the
source's (3.3)--(3.5), with the full covariance, but give their bounds directly.
Write $\Gamma_{\mathrm{E}}$ for Euler's Gamma function in the following
definitions, to distinguish it from the covariance measure:
\begin{equation}\label{E:damped-norm-definitions}
  \begin{split}
    N_0(\theta)&=\int_0^\infty e^{-\theta s}G_s(\R^d)\ud s,\\
    N_1(\theta)&=\frac1{\Gamma_{\mathrm{E}}(\alpha_0)}
    \int_0^\infty\lambda^{\alpha_0-1}\ud\lambda
    \int\frac{\mu(\ud\xi)}{(\theta+\lambda)^2+|\xi|^2},\\
    N_2(\theta)^2&=\int_0^\infty\int_0^\infty
    e^{-\theta(s+t)}|s-t|^{-\alpha_0}
    \int\widehat{G}_s(\xi)\overline{\widehat{G}_t(\xi)}
    \mu(\ud\xi)\ud s\ud t.
  \end{split}
\end{equation}
The second expression uses the elementary wave Laplace transform $\int_0^\infty
e^{-as}\widehat{G}_s(\xi)\ud s = (a^2+|\xi|^2)^{-1}$. The last expression is the
space-time covariance quadratic form; in physical variables it pairs
$G_\theta(s,x)$ with $G_\theta(t,y)$. Thus, the second Green factor in the
printed source (3.5) is read with the independent variable $y$; see
Lemma~\ref{L:translated-covariance-CS} and
Remark~\ref{R:source-case2-time-shift} for the independent justification and the
exact source qualification.

\LeanMarker{L:damped-homogeneous-norms}
\begin{lemma}[Homogeneous damped norms and temporal tails]
  \label{L:damped-homogeneous-norms}
  Under Hypothesis~\ref{H:homog}, the integrals
  in~\eqref{E:damped-norm-definitions} are finite for every $\theta > 0$, with
  $N_2(\theta)^2 \ge 0$, and
  \[
    N_0(\theta)=\theta^{-2},\qquad
    N_1(\theta)=\theta^{\alpha+\alpha_0-2}N_1(1),\qquad
    N_2(\theta)^2=\theta^{-\beta}N_2(1)^2,
    \quad\beta=4-\alpha-\alpha_0.
  \]
  Let
  \[
    K_M(r)=\frac1{\Gamma_{\mathrm{E}}(\alpha_0)}
    \int_M^\infty e^{-\lambda|r|}\lambda^{\alpha_0-1}\ud\lambda.
  \]
  Denote by $N_1^{(M)}(1)$ the first norm with $\lambda > M$, and by
  $N_2^{(M)}(1)$ the second norm with $K_M(s-t)$ in place of
  $|s-t|^{-\alpha_0}$. Then
  \begin{equation}\label{E:damped-tail-norms}
    N_1^{(M)}(1)\longrightarrow0,\qquad
    N_2^{(M)}(1)\longrightarrow0\quad(M\to\infty).
  \end{equation}
  More precisely, let $m:\R^d\to[0,1]$ be measurable and let
  $S\subset S'\subset(0,\infty)$ be measurable. Define $N_{1,m,S}(\theta)$
  by replacing $\mu$ with $m\mu$ and restricting the $\lambda$ integral in
  $N_1$ to $S$. Define $N_{2,m,S}(\theta)^2$ by replacing $\mu$ with $m\mu$
  and the temporal kernel in $N_2^2$ with
  \[
    K_S(r)=\frac1{\Gamma_{\mathrm{E}}(\alpha_0)}
    \int_S e^{-\lambda|r|}\lambda^{\alpha_0-1}\ud\lambda.
  \]
  Then, for every $\theta>0$,
  \begin{equation}\label{E:damped-multiplier-comparison}
    \begin{split}
      0\le N_{1,m,S}(\theta)&\le N_{1,1,S'}(\theta),\\
      0\le N_{2,m,S}(\theta)^2&\le N_{2,1,S'}(\theta)^2.
    \end{split}
  \end{equation}
  These are comparisons of spectral expressions, including full and tail
  domains. An arbitrary $m$ need not preserve Hypothesis~\ref{H:homog}; the
  scaling identities above are asserted for the original covariance only.
  Spatial Gaussian smoothing and Gaussian smoothing of the Green kernel
  give multipliers in this class. The separate mass identity for $N_0$ is
  unchanged by Gaussian smoothing of the Green kernel.
\end{lemma}

\begin{proof}
  We use radial spectral masses to bound the first norm and absolute wave
  estimates to control the second. Scaling gives their exact homogeneous
  powers. Tail convergence and multiplier comparisons are then proved from
  nonnegative spectral expressions after time integration.

  \medskip
  \noindent\textit{Step 1.~}We bound the first norm by the homogeneous
  resolvent integral.
  Scaling the Fourier representation of $\Gamma$ and using Fourier uniqueness
  gives $\mu(cA) = c^\alpha\mu(A)$. Radon finiteness and $\alpha > 0$ then imply
  $\mu(\{0\}) = 0$. In particular, $\mu(B_r) = r^\alpha V$, where $V = \mu(B_1)
  < \infty$. Therefore, for $a > 0$,
  \[
    \int\frac{\mu(\ud\xi)}{a^2+|\xi|^2}
    =\alpha V a^{\alpha-2}
    \int_0^\infty\frac{u^{\alpha-1}}{1+u^2}\ud u<\infty.
  \]
  The last integral is finite since $0 < \alpha < 2$. Substitution into $N_1$
  reduces its finiteness to
  $\int_0^\infty\lambda^{\alpha_0-1}(\theta+\lambda)^{\alpha-2} \ud\lambda <
  \infty$: the conditions at zero and infinity are $\alpha_0 > 0$ and
  $\alpha+\alpha_0 < 2$.

  \medskip
  \noindent\textit{Step 2.~}We prove absolute integrability and positivity of
  the second norm.
  For $N_2$ use $|\widehat{G}_s(\xi)| \le \min(s,|\xi|^{-1})$. On $|\xi| \le 1$,
  absolute integration is bounded by
  \[
    V\int_0^\infty\int_0^\infty
    e^{-\theta(s+t)}|s-t|^{-\alpha_0}st\ud s\ud t<\infty.
  \]
  On $|\xi| > 1$, use $\int_{|\xi|>1}|\xi|^{-2}\mu(\ud\xi) = \alpha
  V/(2-\alpha)$ and the same time integral without $st$. Both time integrals are
  finite by splitting into $s \le t$ and $t \le s$, then using $u = \min(s,t)$
  and $v = |s-t|$; the only singular factor is $v^{-\alpha_0}$, with $\alpha_0 <
  1$. This proves absolute integrability of the displayed spectral formula. For
  each $\xi$, its time integral is nonnegative because the temporal kernel is
  positive definite. Indeed, it is a positive Laplace mixture of the kernels
  $e^{-\lambda|s-t|}$. This proves $N_2(\theta)^2 \ge 0$.

  \medskip
  \noindent\textit{Step 3.~}We compute the homogeneous powers.
  The mass identity $G_s(\R^d) = s$ proves the formula for $N_0$. In $N_1$
  change $\lambda = \theta u$ and $\xi = \theta\eta$. In $N_2$ change $s =
  u/\theta$, $t = v/\theta$ and $\xi = \theta\eta$, using
  $\widehat{G}_{u/\theta}(\theta\eta) = \theta^{-1}\widehat{G}_u(\eta)$. These
  changes give exactly the stated powers.

  For $N_1$, the temporal measure contributes $\theta^{\alpha_0}$, the spatial
  measure contributes $\theta^\alpha$, and the resolvent contributes
  $\theta^{-2}$. For $N_2^2$, the two time Jacobians and two wave factors
  contribute $\theta^{-4}$, while the temporal kernel and spatial measure
  contribute $\theta^{\alpha_0+\alpha}$.

  \medskip
  \noindent\textit{Step 4.~}We remove temporal tails using the absolute
  majorant.
  The first tail in~\eqref{E:damped-tail-norms} is a tail of the finite
  nonnegative $N_1(1)$ integral. For the second, $0 \le K_M(r) \le
  |r|^{-\alpha_0}$ and $K_M(r) \to 0$ for $r \ne 0$. The absolute majorant just
  proved permits dominated convergence in $(s,t,\xi)$, giving $N_2^{(M)}(1)^2
  \to 0$. The tail kernel is positive definite by the same Laplace-mixture
  argument, so its quadratic form is nonnegative.

  \medskip
  \noindent\textit{Step 5.~}We compare nonnegative spectral quadratic forms.
  For the multiplier comparisons, define
  \[
    B_\theta(\lambda,\xi)=\int_0^\infty\int_0^\infty
    e^{-\theta(s+t)}e^{-\lambda|s-t|}
    \widehat G_s(\xi)\overline{\widehat G_t(\xi)}\ud s\ud t.
  \]
  For each $\lambda>0$ and $\xi$, the positive definiteness of
  $e^{-\lambda|s-t|}$ gives $B_\theta(\lambda,\xi)\ge0$.
  The absolute-integrability bound above, $0\le m\le1$ and
  $0\le K_S\le|\cdot|^{-\alpha_0}$ justify changing the integration order:
  \[
    N_{2,m,S}(\theta)^2
    =\frac1{\Gamma_{\mathrm{E}}(\alpha_0)}
    \int_S\lambda^{\alpha_0-1}
    \int m(\xi)B_\theta(\lambda,\xi)\mu(\ud\xi)\ud\lambda.
  \]
  Nonnegativity now proves the second comparison
  in~\eqref{E:damped-multiplier-comparison}. It is used after time integration,
  not pointwise on the oscillating product of wave Fourier kernels. The first
  comparison follows directly from the nonnegative integrand defining $N_1$.
  The choices $S=(0,\infty)$ or $S=(M,\infty)$ give the full or tail bounds.
  Gaussian spatial and Green smoothing contribute products of multipliers in
  $[0,1]$ depending only on $\xi$. Gaussian Green smoothing also preserves
  $G_s(\R^d)=s$, proving its separate assertion for $N_0$. This proves
  Lemma~\ref{L:damped-homogeneous-norms}.
\end{proof}

\paragraph{Covariance pairings of translated kernels.}
Fix $d \in \{1,2,3\}$. Let $\mu$ be a finite positive Borel measure on $\R^d$
such that $\gamma(z) = \int e^{i\xi\cdot z}\mu(\ud\xi)$ is real and nonnegative
for every $z$. Thus, $\gamma$ is bounded, continuous and even, with $0 \le
\gamma(z) \le \mu(\R^d)$. Let $\nu$ be a positive Borel measure on $(0,\infty)$
with $\nu((0,L]) < \infty$ for every $L > 0$, and put
\[
  \kappa_L(t)=\int_{(0,L]}e^{-\lambda|t|}\nu(\ud\lambda),\qquad
  C_L((t,x),(s,y))=\kappa_L(t-s)\gamma(x-y).
\]
For complex $L^1(\R\times\R^d)$ functions $f,h$, define
\begin{equation}\label{E:translated-covariance-form}
  Q_L(f,h)=\int_{\R^2}\int_{\R^{2d}}
  f(t,x)\overline{h(s,y)}\kappa_L(t-s)\gamma(x-y)
  \ud x\ud y\ud t\ud s.
\end{equation}
This integral is absolutely convergent, with $|Q_L(f,h)| \le
\nu((0,L])\mu(\R^d)\|f\|_1\|h\|_1$. For nonnegative $f,h$, define $Q(f,h) =
\lim_{L\to\infty}Q_L(f,h)$ as an extended nonnegative number. This positive cap
limit also defines the pairing when the full temporal kernel is singular on the
diagonal. The spatial covariance is fixed throughout this definition. For a
nonnegative Borel function $G$ on $(0,\infty)\times\R^d$, let $G^0$ be its zero
extension to nonpositive times and put
\[
  N_0(G)=\int_0^\infty\int_{\R^d}G(t,x)\ud x\ud t,\qquad
  N_C(G)^2=Q(G^0,G^0).
\]
For $(A,a) \in \R\times\R^d$, write $T_{A,a}f(t,x) = f(t-A,x-a)$. All translated
pairings are over the full time line.

\begin{lemma}[Translated covariance Cauchy--Schwarz]
  \label{L:translated-covariance-CS}
  Under the preceding assumptions on $\mu$ and $\nu$, let $G_1,G_2$ be
  nonnegative Borel functions with $N_0(G_i),N_C(G_i) < \infty$, $i = 1,2$. For
  arbitrary origins $(A,a),(B,b) \in \R\times\R^d$,
  \begin{equation}\label{E:translated-covariance-CS}
    0\le Q(T_{A,a}G_1^0,T_{B,b}G_2^0)\le N_C(G_1)N_C(G_2),
  \end{equation}
  and $Q(T_{A,a}G_i^0,T_{A,a}G_i^0) = N_C(G_i)^2$ for $i = 1,2$. The coefficient
  one is independent of the origins, dimension, cap and covariance masses. Zero
  norms and zero measures are included.
\end{lemma}

\begin{proof}
  At a finite temporal cap we represent the stationary covariance as an inner
  product. We use its Cauchy--Schwarz inequality together with translation
  invariance of each self-pairing. Positive truncation and cap limits then
  extend the result to the stated kernels.

  \medskip
  \noindent\textit{Step 1.~}We represent the capped covariance by
  deterministic Hilbert features.
  We first prove positivity and Cauchy--Schwarz for the stationary kernel $C_L$
  at fixed finite $L$. For $\lambda > 0$, define
  \[
    h_\lambda(t,r)=\sqrt{2\lambda}\,e^{-\lambda(t-r)}
    \mathbf1_{\{r\le t\}}.
  \]
  Direct integration gives
  \[
    \int_\R h_\lambda(t,r)h_\lambda(s,r)\ud r
    =2\lambda e^{-\lambda(t+s)}
    \int_{-\infty}^{\min(t,s)}e^{2\lambda r}\ud r
    =e^{-\lambda|t-s|}.
  \]
  In particular, $\int h_\lambda(t,r)^2\ud r = 1$. On the complex Hilbert space
  \[
    \mathcal{H}_L=L^2((0,L]\times\R^d\times\R,
    \nu(\ud\lambda)\mu(\ud\xi)\ud r),
  \]
  with inner product linear in its first entry, put $\Phi_L(t,x)(\lambda,\xi,r)
  = h_\lambda(t,r)e^{i\xi\cdot x}$. Then
  \[
    \|\Phi_L(t,x)\|^2=m_L:=\nu((0,L])\mu(\R^d),\qquad
    \langle\Phi_L(t,x),\Phi_L(s,y)\rangle=C_L((t,x),(s,y)).
  \]
  The scalar integrals are absolutely convergent: Cauchy--Schwarz in $r$ bounds
  $\int|h_\lambda(t,r)h_\lambda(s,r)|\ud r$ by one, and both remaining measures
  are finite on the cap.

  \medskip
  \noindent\textit{Step 2.~}We identify the covariance integral as a Hilbert
  pairing.
  The map $(t,x)\mapsto\Phi_L(t,x)$ is norm-continuous, since its squared
  increment is $2m_L-2\operatorname{Re}C_L((t,x),(s,y))$. Continuity of $C_L$
  follows by dominated convergence against the finite measures. Thus, the map is
  strongly measurable, has separable range and is bounded. For $f \in L^1$, the
  Bochner integral $F_L(f) = \int f(t,x)\Phi_L(t,x)\ud t\ud x$ exists. Bochner
  Fubini, justified by $m_L\|f\|_1\|h\|_1$, gives $Q_L(f,h) = \langle
  F_L(f),F_L(h)\rangle$. Therefore,
  \[
    |Q_L(f,h)|\le Q_L(f,f)^{1/2}Q_L(h,h)^{1/2}.
  \]
  This uses the original stationary covariance, rather than positivity of an
  independently shifted cross kernel. The feature space is deterministic and
  introduces no new Gaussian field.

  \medskip
  \noindent\textit{Step 3.~}We compare different translates through their
  unchanged self norms.
  Changing both integration variables in~\eqref{E:translated-covariance-form} by
  the same time and spatial translation gives
  \[
    Q_L(T_{A,a}f,T_{A,a}h)=Q_L(f,h).
  \]
  Both differences remain unchanged. The full time line and Lebesgue translation
  invariance justify this equality without an origin factor. Apply
  Cauchy--Schwarz to the two different translates and this equality separately
  to their self pairings. This proves the capped form
  of~\eqref{E:translated-covariance-CS}.

  \medskip
  \noindent\textit{Step 4.~}We pass the nonnegative kernel and temporal-cap
  limits.
  To make the passage for possibly unbounded nonnegative kernels explicit, set
  \[
    G_{i,M}(t,x)=\min(G_i(t,x),M)
    \mathbf1_{\{0<t\le M,\ |x|\le M\}},\qquad M\in\N.
  \]
  These functions increase to $G_i$ almost everywhere. Their zero extensions and
  translates are bounded compactly supported $L^1$ functions. Since $C_L \ge 0$,
  monotone convergence passes the capped inequality and common-translation
  identity to $G_i$ as $M \to \infty$. Equivalently, the finite-cap passage
  follows from $L^1$ convergence and the absolute bound
  after~\eqref{E:translated-covariance-form}.

  Finally, positivity of $\nu$ and nonnegativity of $\gamma$ give
  \[
    Q_L(G_i^0,G_i^0)\uparrow N_C(G_i)^2,\qquad
    Q_L(T_{A,a}G_1^0,T_{B,b}G_2^0)
    \uparrow Q(T_{A,a}G_1^0,T_{B,b}G_2^0).
  \]
  The finite self norms bound the latter limit by $N_C(G_1)N_C(G_2)$. The
  common-translation identity passes to the limit in the same way. No division
  by a norm is used, so zero norms cause no exception. No spatial Gaussian
  mollifier is removed in this argument. This proves
  Lemma~\ref{L:translated-covariance-CS}.
\end{proof}

\begin{remark}[Source Case 2 and the regularized inputs]
  \label{R:source-case2-time-shift}
  \emph{(i)} In Case~2 of the proof
  of~\cite{chen:author2026:time-dependency}*{Theorem~3.1, p.~21}, the final
  vertices belong to two nonempty chains of lengths $n_1,n_2$ with $n_1+n_2 =
  2n$. Fix their predecessor origins separately:
  \[
    (A,a)=
    \begin{cases}
      (0,0),&n_1=1,\\
      (s_{n_1-1},x_{n_1-1}),&n_1>1,
    \end{cases}
    \qquad
    (B,b)=
    \begin{cases}
      (0,0),&n_2=1,\\
      (s_{2n-1},x_{2n-1}),&n_2>1.
    \end{cases}
  \]
  The second chain's first vertex uses its own root. An empty chain has no
  interchain Case~2 pair. The final-pair integral is, using evenness to fix the
  orientation,
  \[
    \begin{split}
      I&=\int_A^\infty\int_B^\infty\int_{\R^{2d}}
      G_1(t-A,x-a)\kappa(t-s)\gamma(x-y)G_2(s-B,y-b)
      \ud x\ud y\ud s\ud t\\
      &=Q(T_{A,a}G_1^0,T_{B,b}G_2^0),
    \end{split}
  \]
  where the full temporal expression means the positive cap limit.
  Lemma~\ref{L:translated-covariance-CS} gives $0 \le I \le N_C(G_1)N_C(G_2)$
  uniformly in the predecessors, so the bound can be inserted into the remaining
  nonnegative integrals in the source's induction. If only $u = t-A$ and $v =
  s-B$ are substituted, the integrand has factors
  \begin{equation}\label{E:translated-covariance-offset}
    G_1(u,x-a)\kappa(u-v+A-B)\gamma(x-y)G_2(v,y-b).
  \end{equation}
  The time difference $A-B$ must remain. The source's first Case~2 equality on
  physical p.~21 displays $\kappa(u-v)$ after this change of time variables,
  which is not an identity for general $A \ne B$. Its spatial origins remain in
  the two tests, so $\gamma(x-y)$ is correct in those absolute spatial
  variables. Changing also $x = X+a$, $y = Y+b$ gives $\gamma(X-Y+a-b)$. We use
  the original integral and the proved inequality, not the offset-free equality
  or positivity of a shifted cross kernel.

  \emph{(ii)} The norm $N_C(G)$ pairs $G(t,x)$ with $G(s,y)$
  in~\eqref{E:translated-covariance-form}, with independent spatial variables.
  The source's printed (3.5), physical p.~14, repeats $x$ in its second Green
  factor. Our quadratic form states explicitly the convention used
  in~\eqref{E:damped-norm-definitions}; it is not a literal quotation of that
  display.

  \emph{(iii)} We use the common-prefix source convention
  of~\cite{chen:author2026:time-dependency}*{(2.11), p.~9, and (2.24), p.~13}.
  Write $Y_{\epsilon,\delta}$ for its positive-width fields. Independence of
  disjoint Laplace-frequency increments makes the cross covariance of two
  prefixes their common part. Their upper limits are $\delta_j^{-1}$ and
  $\delta_h^{-1}$, whose minimum is $(\delta_j\vee\delta_h)^{-1}$. Spatial heat
  convolution on the two fields gives $p_{\epsilon_j}*p_{\epsilon_h} =
  p_{\epsilon_j+\epsilon_h}$. Thus, with the original covariance amplitude,
  \begin{equation}\label{E:pairing-source-covariance}
    \begin{split}
      \E[Y_{\epsilon_j,\delta_j}(s,x)
      Y_{\epsilon_h,\delta_h}(t,y)]
      &=\kappa_{\delta_j\vee\delta_h}(s-t)
      (p_{\epsilon_j+\epsilon_h}*\Gamma)(x-y),\\
      \kappa_\delta(u)&=\frac1{\Gamma_{\mathrm{E}}(a)}
      \int_0^{\delta^{-1}}e^{-\lambda|u|}\lambda^{a-1}\ud\lambda,
      \qquad a=\alpha_0.
    \end{split}
  \end{equation}
  For generalized $\Gamma$, the spatial identity follows by testing against heat
  functions, or by their spectral multipliers; it requires no covariance
  density. The endpoint zero in the Laplace interval has zero Lebesgue measure,
  so $[0,L]$ and $(0,L]$ give the same integral. This verifies the covariance
  used on each vertexwise edge. It does not construct the coupled family from
  the given $W$ alone. Its common-field realization, any prescribed auxiliary
  Gaussian extension, and its joint source limits remain the separate hypotheses
  and deductions of
  Lemmas~\ref{L:spatial-wick-measures}--\ref{L:joint-fixed-order}; no new
  point-field measurability or uniqueness claim is made here. At an edge $e =
  \{j,h\}$, retain positive spatial widths $\epsilon_j,\epsilon_h$ and use
  \[
    \mu_e(\ud\xi)=e^{-(\epsilon_j+\epsilon_h)|\xi|^2/2}\mu(\ud\xi),
    \qquad \gamma_e=p_{\epsilon_j+\epsilon_h}*\Gamma\ge0.
  \]
  For $b = \epsilon_j+\epsilon_h > 0$, the radial ball masses $\mu(B_r) =
  Vr^\alpha$ from Lemma~\ref{L:damped-homogeneous-norms} give
  \[
    \mu_e(\R^d)=\alpha Vb^{-\alpha/2}
    \int_0^\infty u^{\alpha-1}e^{-u^2/2}\ud u<\infty.
  \]
  This uses ball masses rather than a radial density, and establishes finiteness
  only. It does not replace the later uniform norm bounds. Nonnegativity of the
  heat convolution follows from the nonnegative Radon measure $\Gamma$. This
  spatial covariance stays fixed during the lemma's application. For a
  measurable temporal domain $A_e \subset (0,\infty)$, use
  \[
    \nu_e(\ud\lambda)=\mathbf1_{A_e}(\lambda)
    \lambda^{\alpha_0-1}\ud\lambda/\Gamma_{\mathrm{E}}(\alpha_0).
  \]
  Since $\alpha_0 > 0$, it is finite on every $(0,L]$. This covers both bounded
  regularization domains and marked-tail domains after adding the finite cap,
  with the Gamma normalization already included. The tested kernels are
  $e^{-\theta t}G_t$, with positive Green smoothing in dimension three when a
  function is required. They have $N_0 = \theta^{-2}$, and
  Lemma~\ref{L:damped-homogeneous-norms} supplies the finite full and tail
  second-norm bounds and first-norm bounds. These are precisely the edge norm
  hypotheses used in Lemma~\ref{L:pairing-product-envelope}. The Case~2 estimate
  contributes exactly its two second norms, with coefficient one. It changes
  neither a $Q_2$ factor nor any coefficient estimate and adds no factor
  depending on $A+B$. Temporal cap removal uses positivity; Green-width removal
  at fixed positive spatial widths and the subsequent spatial-noise limit retain
  their separate weak-limit arguments in Lemmas~\ref{L:joint-fixed-order}
  and~\ref{L:spatial-wick-measures}. No pointwise monotonicity of spatial
  Gaussian smoothing or general unregularized spatial-distribution version is
  asserted.

  \emph{(iv)} We prove the product estimate needed here in
  Lemma~\ref{L:pairing-product-envelope}, rather than importing the remaining
  displayed induction. On p.~16, the empty-first-chain base caption interchanges
  the mass and first-norm kernel labels; their chronological order gives mass on
  the first kernel. On p.~18, (3.9) omits the factor two when a fused kernel
  occupies a first-norm slot. Retaining that factor and using the mass times
  first-norm convolution bound below preserves the half-bound needed in Case~3.
  The p.~19 mass-slot partition uses a union of the two excluded labels, not the
  printed intersection; its first-slot edge count excludes the fused slot. The
  p.~21 Case~2 offset is handled by the preceding lemma. These are explicit
  qualifications of the printed formulas, not corrections attributed to the
  source or a claim that its intended product bound is false. The finite maximum
  below is sufficient for all our uses; we assert no partition independent of
  the kernel family.
\end{remark}

\section{A finite product bound for two ordered chains}

Fix $d \ge 1$, $0 < a < 1$ and a matching $D$ on $2n$ labels. For each edge $e
\in D$, let $\mu_e$ be finite and positive, with $\gamma_e(x) = \int
e^{i\xi\cdot x}\mu_e(\ud\xi) \ge 0$. For a measurable $A_e \subset (0,\infty)$
put
\[
  \nu_e(\ud\lambda)=\Gamma_{\mathrm{E}}(a)^{-1}
  \mathbf1_{A_e}(\lambda)\lambda^{a-1}\ud\lambda,
  \qquad \kappa_e(t)=\int e^{-\lambda|t|}\nu_e(\ud\lambda).
\]
For nonnegative $G \in L^1((0,\infty)\times\R^d)$ define
\begin{equation}\label{E:pairing-edge-norms}
  \begin{split}
    N_0(G)&=\int G(t,x)\ud t\ud x,\\
    N_{1,e}(G)&=\int\left|\int e^{-\lambda t+i\xi\cdot x}
    G(t,x)\ud t\ud x\right|\nu_e(\ud\lambda)\mu_e(\ud\xi),\\
    N_{2,e}(G)^2&=\int G(s,x)G(t,y)
    \kappa_e(s-t)\gamma_e(x-y)\ud s\ud x\ud t\ud y.
  \end{split}
\end{equation}
Time integrations in these definitions are positive. Denote the inner
Laplace--Fourier transform by $\mathcal{L} G(\lambda,\xi)$. The last integral
means its nonnegative temporal cap limit, as in the preceding covariance
pairing. In particular, its spatial variables are independent.

Let $I_1 = \{1,\ldots,n_1\}$ and $I_2 = \{n_1+1,\ldots,2n\}$, where $n_1+n_2 =
2n$ and either chain may be empty. Each chain is ordered in its own time
variables. Write $p(l)$ for its predecessor, using $(s_{p(l)},x_{p(l)}) = (0,0)$
at the first vertex of each chain. Set
\begin{equation}\label{E:pairing-diagram-general}
  \begin{split}
    I_D(\boldsymbol G)=\int
    &\prod_{\rho=1}^2\prod_{l\in I_\rho}
    G_l(s_l-s_{p(l)},x_l-x_{p(l)})\\
    &\times\prod_{e=\{j,k\}\in D}
    \kappa_e(s_j-s_k)\gamma_e(x_j-x_k)
    \ud\boldsymbol s\ud\boldsymbol x .
  \end{split}
\end{equation}
The domain is the product of the two positive ordered time simplices and all
spatial variables.

Call $(Q,\iota)$ admissible if $Q_0,Q_1,Q_2$ partition the kernel labels, $|Q_0|
= |Q_1|$, and the map $\iota:Q_1\cup Q_2 \to D$ uses each edge once on $Q_1$ or
twice on $Q_2$, never both. There is no incidence requirement between $l$ and
$\iota(l)$. These assignments form a finite nonempty set $\mathcal{A}_D$. Put
\begin{equation}\label{E:pairing-product-general}
  \mathcal{P}_{Q,\iota}(\boldsymbol G)=
  \prod_{l\in Q_0}N_0(G_l)
  \prod_{l\in Q_1}2N_{1,\iota(l)}(G_l)
  \prod_{l\in Q_2}N_{2,\iota(l)}(G_l).
\end{equation}
Thus, $|Q_2|$ is even and $|Q_0|+|Q_2|/2 = n$.

\begin{lemma}[Pairing product envelope]\label{L:pairing-product-envelope}
  Suppose every $G_l$ is nonnegative and integrable, and
  $N_{1,e}(G_l),N_{2,e}(G_l) < \infty$ for every $l$ and $e \in D$. Then
  \begin{equation}\label{E:pairing-product-envelope}
    I_D(\boldsymbol G)\le
    \max_{(Q,\iota)\in\mathcal{A}_D}
    \mathcal{P}_{Q,\iota}(\boldsymbol G).
  \end{equation}
  A maximizing assignment may depend on the kernels and edge data. No single
  assignment valid for all kernel families is asserted. For $n = 0$, both sides
  are one. Zero kernels, measures and norms are allowed.
\end{lemma}

\begin{proof}
  We establish convolution bounds for fused kernels, then prove a regional
  half-bound by induction. The endpoint configurations either give that bound
  directly, give two second norms, or split into two half-bounds. The finite
  maximum permits different assignments on the two regions.

  \medskip
  \noindent\textit{Step 1.~}We work at finite temporal cap and control fused
  kernels.
  First cap every temporal domain at $L < \infty$; its measure is finite since
  $a > 0$. The stationary covariance pairing in
  Lemma~\ref{L:translated-covariance-CS} has norm $N_{2,e}$. Its feature
  $h_\lambda(t,r)e^{i\xi\cdot x}$ and its proof apply in any finite dimension:
  their only spatial input is the finite measure $\mu_e$. Hence, its
  coefficient-one Cauchy--Schwarz inequality and invariance under simultaneous
  time and space translation apply here. All Fourier and Bochner interchanges at
  this cap are justified by the finite covariance mass times the relevant $L^1$
  masses.

  For one-sided spacetime convolution $F = G*H$, positivity and the
  Laplace--Fourier product give
  \begin{equation}\label{E:pairing-convolution-norms}
    \begin{split}
      N_0(F)&=N_0(G)N_0(H),\\
      N_{1,e}(F)&\le N_0(G)N_{1,e}(H),\\
      N_{2,e}(F)&\le N_0(G)N_{2,e}(H).
    \end{split}
  \end{equation}
  Indeed, $|\mathcal{L} G(\lambda,\xi)| \le N_0(G)$ for $\lambda \ge 0$ proves
  the first-norm estimate. For the second, zero-extend $H$ in time. Its
  translates have the same covariance Hilbert norm. Express $G*H$ as their
  Bochner integral against $G$ and apply Minkowski. Thus, every fused kernel
  satisfies all remaining edge hypotheses. No division by a norm is used.

  \medskip
  \noindent\textit{Step 2.~}We prove the regional half-bound and its
  admissible assignments.
  We use induction on $n$, beginning with the empty product at $n = 0$. With two
  vertices on one chain the integral is at most $N_0(G_{\rm
  first})N_{1,e}(G_{\rm last})$; with one on each chain it is at most their two
  second norms by covariance Cauchy--Schwarz. These also verify both
  nonempty-chain orientations at $n = 1$. Assume the full bound at $n-1$.

  We first prove a half-bound on any region where an endpoint $v$, its
  predecessor $p$ and its partner $a$ satisfy $s_a \le s_p$. Require the region
  indicator to be independent of $(s_v,x_v)$. Write $e_0 = \{a,v\}$. Since
  $\kappa_{e_0}$ is decreasing on positive arguments, replace
  $\kappa_{e_0}(s_v-s_a)$ by $\kappa_{e_0}(s_v-s_p)$. Integrating the endpoint
  and taking the modulus of its spatial Fourier phase gives the factor
  \[
    \int |\mathcal{L} G_v(\lambda,\xi)|
    \nu_{e_0}(\ud\lambda)\mu_{e_0}(\ud\xi)=N_{1,e_0}(G_v).
  \]
  Now drop the region indicator by positivity. Remove $v$ and its edge. Its
  partner $a$ is unpaired. If $a$ is terminal in its remaining chain, integrate
  it to $N_0(G_a)$ and apply induction. The result is one half of an admissible
  product, adding $a$ to $Q_0$ and $v$ to $Q_1$ with edge $e_0$.

  Otherwise, let $b$ be the immediate successor of $a$ on its own chain.
  Integration at $a$ gives exactly $F = G_a*G_b$ at $b$, also when $a$ was the
  first vertex. Apply induction to the two reduced chains and
  $D\setminus\{e_0\}$. Mark the slot $i = 0,1,2$ occupied by $F$ in a maximizing
  reduced product. Let $R_0,R_1,R_2$ be its other slots and $B_{\rm rest}$ their
  product, retaining every first-slot factor two. The new slots can be chosen as
  follows:
  \[
    \begin{array}{c|ccc}
      i&Q_0&Q_1&Q_2\\ \hline
      0&R_0\cup\{a,b\}&R_1\cup\{v\}&R_2\\
      1&R_0\cup\{a\}&R_1\cup\{b,v\}&R_2\\
      2&R_0\cup\{a\}&R_1\cup\{v\}&R_2\cup\{b\}.
    \end{array}
  \]
  For $i = 0$, $|R_0|+1 = |R_1|$ and $|R_2|$ is even. Mass multiplicativity
  gives one half of the displayed new product. For $i = 1$, $|R_0| = |R_1|+1$
  and $|R_2|$ is even. If $e_1$ is the inherited edge, the term and its bound
  are
  \begin{align*}
    &N_{1,e_0}(G_v)\,[2N_{1,e_1}(F)]B_{\rm rest}\\
    &\qquad\le\tfrac12N_0(G_a)
    [2N_{1,e_1}(G_b)][2N_{1,e_0}(G_v)]B_{\rm rest}.
  \end{align*}
  The inherited edge stays first exactly once; $e_0$ is first once. No incidence
  condition is needed. For $i = 2$, $|R_0| = |R_1|$ and $|R_2|$ is odd; $e_1$
  already occurs once there. The last inequality
  in~\eqref{E:pairing-convolution-norms} puts its other second norm on $G_b$, a
  mass on $G_a$ and $2N_{1,e_0}(G_v)$ on $v$, again with a factor one half.
  Every old edge retains its type and multiplicity. The displayed counts give
  $|Q_0| = |Q_1|$ and even $|Q_2|$ in all three cases. This proves the regional
  half-bound, using only the full induction hypothesis, not an assumed stronger
  induction.

  \medskip
  \noindent\textit{Step 3.~}We treat same-chain pairs and the pair of terminal
  endpoints.
  If an endpoint is paired to a vertex in its own chain, its partner is at or
  before its predecessor. Hence, this half-bound applies on the entire domain.
  This also handles an empty chain. If both endpoints $u = n_1$ and $v = 2n$ are
  paired to each other, condition on their separate predecessors $(A,z),(B,w)$.
  The endpoint integral keeps the temporal argument $r-q+A-B$ in relative times.
  Here $r=s_u-A$ and $q=s_v-B$ are relative times;
  Lemma~\ref{L:translated-covariance-CS} applies with origins $(A,z),(B,w)$.
  The symbol $q$ is reused as a vertex label only in the next step.
  The original stationary covariance Cauchy--Schwarz inequality bounds it by
  $N_{2,e}(G_u)N_{2,e}(G_v)$; translate each self pairing separately. Induction
  and adding $u,v$ to $Q_2$ prove the full bound. There is no origin factor or
  discarded predecessor-time offset.

  \medskip
  \noindent\textit{Step 4.~}We split the two distinct cross-pair
  configurations.
  It remains that both endpoint pairs are cross and distinct. Let $a$ be the
  partner of $v$, $b$ the partner of $u$, and put $p = u-1$, $q = v-1$. Then $a
  < u$ and $b < v$. On $\{s_a\le s_q\}$ remove $v$ by the regional half-bound.
  Its indicator is independent of $s_v,x_v$. On the complementary region the
  chain orders give $s_p \ge s_a > s_q \ge s_b$, so remove $u$ instead; that
  indicator is independent of $s_u,x_u$. The removed partner has a successor on
  its own chain in both cases. The two bounds are $P/2$ and $P'/2$ for possibly
  different admissible products of the original labels. Therefore, their sum is
  at most $(P+P')/2 \le \max(P,P')$. This
  proves~\eqref{E:pairing-product-envelope} at the finite cap.

  \medskip
  \noindent\textit{Step 5.~}We remove the finite cap by positivity.
  Finally, increase $L$ to infinity. The diagram integrands and all norm factors
  increase to their declared values, by positivity and monotone convergence.
  There are only finitely many admissible products, so their maximum converges
  to the maximum of their limits. No kernel-independent or measurable maximizing
  choice is needed. Zero norms and measures require no exception; the positive
  cap convention avoids any undefined product on a singular time diagonal. This
  proves Lemma~\ref{L:pairing-product-envelope}.
\end{proof}

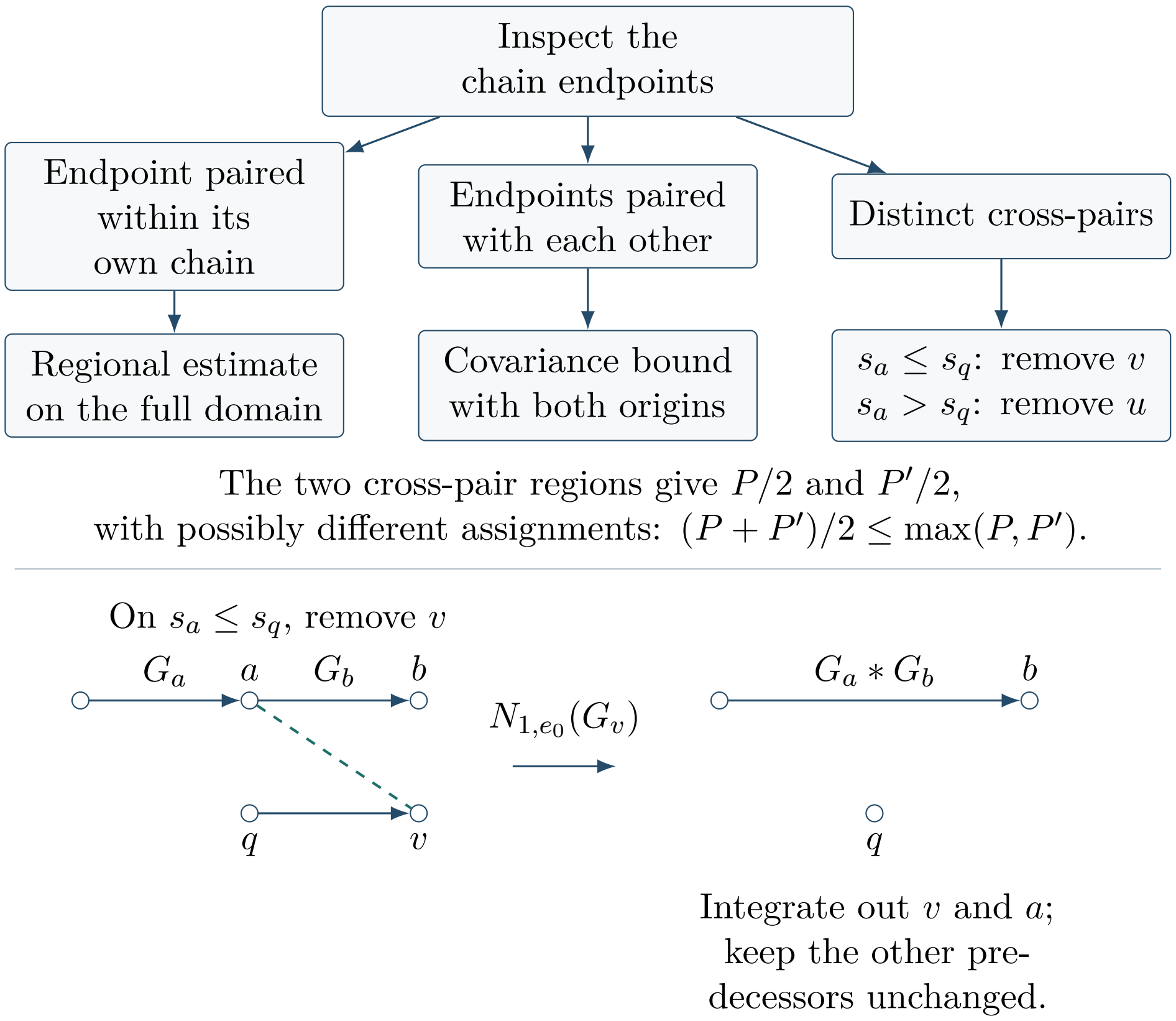
\begin{figure}
\centering
\begin{tikzpicture}[x=1cm,y=1cm]
  \node[ham/box,text width=5.3cm] (hamRoot) at (5,0)
  {Inspect the\\chain endpoints};
  \node[ham/box,text width=3.25cm] (hamSelf) at (.6,-1.65)
  {Endpoint paired\\within its own chain};
  \node[ham/box,text width=3.25cm] (hamTogether) at (5,-1.65)
  {Endpoints paired\\with each other};
  \node[ham/box,text width=3.25cm] (hamDistinct) at (9.4,-1.65)
  {Distinct cross-pairs};
  \draw[ham/arrow] (hamRoot)--(hamSelf);
  \draw[ham/arrow] (hamRoot)--(hamTogether);
  \draw[ham/arrow] (hamRoot)--(hamDistinct);
  \node[ham/box,text width=3.25cm] (hamHalf) at (.6,-3.45)
  {Regional estimate\\on the full domain};
  \node[ham/box,text width=3.25cm] (hamCS) at (5,-3.45)
  {Covariance bound\\with both origins};
  \node[ham/box,text width=3.25cm] (hamRegions) at (9.4,-3.45)
  {$s_a\le s_q$: remove $v$\\$s_a>s_q$: remove $u$};
  \draw[ham/arrow] (hamSelf)--(hamHalf);
  \draw[ham/arrow] (hamTogether)--(hamCS);
  \draw[ham/arrow] (hamDistinct)--(hamRegions);
  \node[ham/note,text width=10.5cm] at (5,-4.75)
  {The two cross-pair regions give $P/2$ and $P'/2$,\\
  with possibly different assignments: $(P+P')/2\le\max(P,P')$.};
  \draw[hamblue!35] (-1.1,-5.4)--(11.1,-5.4);
  \node[ham/note] at (1.7,-5.95) {On $s_a\le s_q$, remove $v$};
  \node[ham/dot] (hamBefore) at (-.4,-6.8) {};
  \node[ham/dot] (hamA) at (1.4,-6.8) {};
  \node[ham/dot] (hamB) at (3.2,-6.8) {};
  \node[ham/dot] (hamQ) at (1.4,-8.0) {};
  \node[ham/dot] (hamV) at (3.2,-8.0) {};
  \draw[ham/arrow] (hamBefore)--node[above,font=\small] {$G_a$} (hamA);
  \draw[ham/arrow] (hamA)--node[above,font=\small] {$G_b$} (hamB);
  \draw[ham/arrow] (hamQ)--(hamV);
  \draw[ham/pair] (hamA)--(hamV);
  \node[font=\small,above=2pt] at (hamA) {$a$};
  \node[font=\small,above=2pt] at (hamB) {$b$};
  \node[font=\small,below=2pt] at (hamQ) {$q$};
  \node[font=\small,below=2pt] at (hamV) {$v$};
  \draw[ham/arrow] (4.2,-7.5)--(5.3,-7.5);
  \node[font=\small] at (4.75,-7.0) {$N_{1,e_0}(G_v)$};
  \node[ham/dot] (hamAfterBefore) at (6.4,-6.8) {};
  \node[ham/dot] (hamAfterB) at (9.7,-6.8) {};
  \node[ham/dot] at (8.05,-8.0) {};
  \draw[ham/arrow] (hamAfterBefore)--node[above,font=\small] {$G_a*G_b$} (hamAfterB);
  \node[font=\small,above=2pt] at (hamAfterB) {$b$};
  \node[font=\small,below=2pt] at (8.05,-8.0) {$q$};
  \node[ham/note,text width=5.2cm] at (8.05,-9.5)
  {Integrate out $v$ and $a$;\\keep the other predecessors unchanged.};
\end{tikzpicture}
\caption[Endpoint elimination in two ordered chains]{The endpoint cases in the finite two-chain product bound. The empty
matching is the induction base; if only one chain is nonempty, use the
within-chain case.
Here $u,v$ are the chain endpoints, $a$ is paired with $v$,
and $q$ precedes $v$ in the distinct-cross-pair case. The lower sketch
shows one local fragment used in the regional fusion step, with
$b$ the immediate successor of $a$; other pairings are suppressed. Deleting $v$ on $s_a\le s_q$ extracts
$N_{1,e_0}(G_v)$, where $e_0=\{a,v\}$. If the unpaired $a$ is terminal,
integrate it to $N_0(G_a)$ instead; otherwise fuse its two consecutive
kernels by one-sided spacetime convolution. The mutually paired endpoint case retains
both predecessor origins, including their time difference. The argument
chooses a finite maximum of products; it does not produce one assignment
valid for all kernel families. See
\ProofRef{Lemma}{L:pairing-product-envelope}.}
\label{Fig:endpoint-elimination}
\end{figure}

Use the joint Gaussian regularizations of the source's (2.11). Write $a =
(\epsilon,\delta) > 0$ for one pair of widths and $\boldsymbol a =
(a_1,\ldots,a_n)$ for vertexwise widths. Let $Y_n^{\boldsymbol a}(t,x)$ be the
ordered terminal integral with kernel $g_n$ in~\eqref{E:expansion} and field
$W_{\epsilon_j,\delta_j}(s_j,y_j)$ at vertex $j$; put $Y_0 = 1$. At fixed
positive widths these are Gaussian-polynomial Bochner integrals. Indeed, Wick's
formula bounds their product integrands in $L^2$ by $\sqrt{(2n-1)!!}\,V^{n/2}$,
where $V$ is the largest of the finitely many field variances. The terminal
Green measure has mass $t^{2n}/(2n)!$, by induction using $G_r(\R^d) = r$. In
dimension three these integrals use the Green measures: products and linear
pushforwards of $p_\eta*G$ converge against the bounded continuous covariance
products. Mixed moments with two Green widths and the same mass bound give the
$L^2$ limit as $\eta\downarrow0$, before any noise width is removed.

\section{Spatial Wick measures and fixed-order construction}

Put $\gamma_b = p_b*\Gamma$, where $p_b$ is the Brownian heat kernel. Fix $m+n =
2k \ge 2$, the two chains $I_1,I_2$ above and $D \in \Pi_k$. For a gap vector
$\boldsymbol r \in (0,\infty)^{2k}$, let $s_l$ be the sum of the gaps up to
vertex $l$ in its own chain. Similarly, for spatial increments $z_l$, let $x_l$
be their cumulative sum in that chain. For positive spatial widths
$\boldsymbol\epsilon$ define
\[
  J_{\boldsymbol\epsilon,D}(\boldsymbol r)
  =\int\prod_{l=1}^{2k}G_{r_l}(\ud z_l)
  \prod_{\{j,h\}\in D}
  \gamma_{\epsilon_j+\epsilon_h}(x_j-x_h),\qquad
  \rho_{\boldsymbol\epsilon,D}(\ud\boldsymbol r)
  =J_{\boldsymbol\epsilon,D}(\boldsymbol r)\ud\boldsymbol r.
\]
These are nonnegative measures. They are locally finite because the smoothed
covariances are bounded and the Green masses are $r_l$. Write
$\nu_{\boldsymbol\epsilon,D} = e^{-\sum_l r_l} \rho_{\boldsymbol\epsilon,D}$,
and let $\ell_i$ be the sum of the gaps in chain $i$, with value zero for an
empty chain.

\begin{lemma}[Positive spatial Wick-measure limit]\label{L:spatial-wick-measures}
  Under Hypothesis~\ref{H:homog}, the finite measures
  $\nu_{\boldsymbol\epsilon,D}$ converge weakly, as all spatial widths tend to
  zero at arbitrary relative rates, to a finite measure $\nu_D$. The locally
  finite measure $\rho_D = e^{\sum_l r_l}\nu_D$ satisfies
  \[
    \rho_D(cA)=c^{(4-\alpha)k}\rho_D(A),\qquad c>0.
  \]
  It gives zero mass to every face $r_l = 0$ and every face $\ell_i = T > 0$.
  Consequently, for every bounded continuous $\varphi$ and $T_1,T_2 > 0$,
  \begin{equation}\label{E:spatial-wick-limit}
    \int_{\ell_1<T_1,\,\ell_2<T_2}\varphi\ud\rho_{\boldsymbol\epsilon,D}
    \longrightarrow
    \int_{\ell_1<T_1,\,\ell_2<T_2}\varphi\ud\rho_D.
  \end{equation}
\end{lemma}

\begin{proof}
  We compute positive Laplace transforms and uniformly bound their masses
  through Brownian occupation integrals. We then recover a unique finite
  measure limit, exclude boundary mass, and use scaling to treat the
  finite-horizon faces.

  \medskip
  \noindent\textit{Step 1.~}We obtain the positive Laplace transform.
  Orient each pair $e = \{j,h\}$ as $j < h$, and let $\eta_l(\boldsymbol\xi)$ be
  the coefficient of $z_l$ in $\sum_{e=\{j,h\}}\xi_e\cdot(x_j-x_h)$. Finite
  heat-weighted spectral measures and the finite damped Green masses justify
  Fourier and Laplace Fubini. The wave Laplace transform gives, for $\lambda_l
  \ge 0$,
  \begin{equation}\label{E:spatial-positive-laplace}
    \begin{split}
      L_{\boldsymbol\epsilon,D}(\boldsymbol\lambda)
      &:=\int e^{-\sum_l\lambda_l r_l}\ud\nu_{\boldsymbol\epsilon,D}\\
      &=\int\exp\left\{-\frac12\sum_{e=\{j,h\}}
      (\epsilon_j+\epsilon_h)|\xi_e|^2\right\}
      \prod_{l=1}^{2k}\frac1{(1+\lambda_l)^2+|\eta_l(\boldsymbol\xi)|^2}
      \prod_{e\in D}\mu(\ud\xi_e).
    \end{split}
  \end{equation}
  The integrand is nonnegative. Thus, the transforms increase when any spatial
  width decreases. This is transform monotonicity, not pointwise monotonicity of
  the physical covariances.

  \medskip
  \noindent\textit{Step 2.~}We bound Brownian occupation moments uniformly in
  the smoothing.
  We prove a uniform bound on their total masses. With $V = \mu(B_1)$, Fourier
  triangle and radial homogeneity give
  \begin{equation}\label{E:construction-heat-bound}
    0\le(p_t*\gamma_b)(x)\le A t^{-\alpha/2},\qquad t,b>0,
    \quad A=\alpha V\int_0^\infty u^{\alpha-1}e^{-u^2/2}\ud u<\infty.
  \end{equation}
  Indeed, its Fourier multiplier is $e^{-(t+b)|\xi|^2/2}$, whose integral is
  $A(t+b)^{-\alpha/2} \le At^{-\alpha/2}$. For a Brownian motion $B$ put $S_b(T)
  = \int_0^T\int_0^T\gamma_b(B_s-B_t)\ud s\ud t$. For every integer $j \ge 1$,
  \begin{equation}\label{E:construction-brownian-moments}
    \E_B S_b(T)^j
    \le\frac{(2j)![A\Gamma_{\mathrm{E}}(1-\alpha/2)]^j}
    {\Gamma_{\mathrm{E}}((2-\alpha/2)j+1)}T^{(2-\alpha/2)j}.
  \end{equation}
  To verify this, order the $2j$ time variables in the moment integral.
  Condition on all but the last remaining Brownian vertex. Its pair factor has
  conditional mean at most $A$ times the last remaining time gap to the power
  $-\alpha/2$, by~\eqref{E:construction-heat-bound}. Remove this vertex and its
  partner. The remaining factors do not involve either removed vertex, and their
  marginal law is again the Brownian law at the remaining times. Repeat this
  procedure. The $j$ selected right endpoints are distinct; each selected gap is
  at least the original adjacent gap at its right endpoint. On each ordered
  simplex the resulting bound is therefore $A^j$ times a product of $j$ distinct
  adjacent gaps to the power $-\alpha/2$. Integration over the simplex gives
  $\Gamma_{\mathrm{E}}(1-\alpha/2)^j
  T^{(2-\alpha/2)j}/\Gamma_{\mathrm{E}}((2-\alpha/2)j+1)$. There are $(2j)!$
  possible orders. This proves the bound uniformly in $b$; the zeroth moment is
  one.

  \medskip
  \noindent\textit{Step 3.~}We bound the two-chain transform by finite
  Brownian moments.
  Let $\epsilon_* = \min_l\epsilon_l$. The positive formula
  \eqref{E:spatial-positive-laplace} bounds $L_{\boldsymbol\epsilon,D}(0)$ by
  the transform with every vertex width equal to $\epsilon_*$. The identity of
  finite positive measures
  \[
    \int_0^\infty e^{-at}G_t(\ud x)\ud t
    =\frac12\int_0^\infty e^{-a^2u/2}p_u(x)\ud u\ud x
  \]
  follows by their equal Fourier transforms $(a^2+|\xi|^2)^{-1}$. Apply it to
  the $2k$ Green factors at unit damping. This writes the common-width transform
  as $2^{-2k}$ times a two-Brownian-chain integral with gap weights $e^{-\sum_l
  u_l/2}$ and covariance $\gamma_{2\epsilon_*}$ on every pair. If $\tau_i$ is
  the last Brownian time in chain $i$, insert
  \[
    1=\frac14\int_{\tau_1}^\infty\int_{\tau_2}^\infty
    e^{-((a-\tau_1)+(b-\tau_2))/2}\ud a\ud b.
  \]
  This also holds for an empty chain with $\tau_i = 0$. Positivity permits
  extension of the ordered Brownian time domains to their full cubes. At fixed
  paths the time integrations then factor into one occupation integral per pair.
  For independent Brownian paths, spectral Cauchy--Schwarz gives
  \[
    0\le A_{12}:=\int_0^a\int_0^b
    \gamma_{2\epsilon_*}(B_1(s)-B_2(t))\ud s\ud t
    \le\sqrt{S_{2\epsilon_*}^{(1)}(a)S_{2\epsilon_*}^{(2)}(b)}.
  \]
  If the pairing has $v_1,v_2$ self-pairs and $w$ cross-pairs, with $v_1+v_2+w =
  k$, its occupation product is at most
  \[
    S_1^{v_1+w/2}S_2^{v_2+w/2}\le(1+S_1)^k(1+S_2)^k.
  \]
  By~\eqref{E:construction-brownian-moments}, each expected factor is bounded by
  a fixed finite sum $M_k$ of nonnegative powers of its horizon, independently
  of $\epsilon_*$. Independence of the two paths now gives
  \[
    \sup_{\boldsymbol\epsilon>0}L_{\boldsymbol\epsilon,D}(0)
    \le2^{-2k-2}\left(\int_0^\infty e^{-a/2}M_k(a)\ud a\right)^2<\infty.
  \]

  \medskip
  \noindent\textit{Step 4.~}We identify the finite-measure limit and rule out
  escape of mass.
  Monotone convergence in~\eqref{E:spatial-positive-laplace} yields finite
  limits $L_D(\boldsymbol\lambda)$ for all $\boldsymbol\lambda \ge 0$. The
  limiting integrand at zero is integrable, so dominated convergence shows
  continuity at zero. Here is the finite-measure Laplace argument, including
  possible loss of mass. Push the measures to $[0,1]^{2k}$ by $y_l = e^{-r_l}$.
  Uniformly bounded masses on this compact cube have weak cluster limits. Their
  integer monomial moments are the prescribed $L_D$ values, so polynomial
  density determines a unique cluster limit. For $c > 0$ the integral of
  $\prod_l y_l^c$ is $L_D(c,\ldots,c)$. As $c\downarrow0$, continuity at zero
  and monotone convergence show that the cluster limit has no mass on the faces
  where any $y_l = 0$. Closed neighborhoods of the $y$-zero faces have limiting
  mass tending to zero; the portmanteau inequality bounds the approximating tail
  masses by these same quantities. Pulling back therefore proves weak
  convergence on $[0,\infty)^{2k}$. Bounded continuous test functions are
  handled by cutting off a neighborhood of those zero faces and then using the
  vanishing tail mass. Every sequence of widths tending to zero has this same
  limit, so the convergence is joint and independent of their relative rates.

  \medskip
  \noindent\textit{Step 5.~}We exclude zero gaps and positive horizon faces.
  Sending just $\lambda_l \to \infty$ in the limiting positive integral
  \eqref{E:spatial-positive-laplace} gives zero by domination by its integrable
  zero-parameter integrand. On the measure side this limit is $\nu_D\{r_l=0\}$.
  Thus, every coordinate-zero face has zero mass.

  Heat and wave scaling give, for Borel sets,
  \[
    \rho_{\boldsymbol\epsilon,D}(cA)
    =c^{(4-\alpha)k}\rho_{\boldsymbol\epsilon/c^2,D}(A).
  \]
  This follows from $G_{cr}(c\ud z) = cG_r(\ud z)$, $\gamma_b(cx) =
  c^{-\alpha}\gamma_{b/c^2}(x)$, and the $2k$ time Jacobians. Passing to vague
  limits proves the stated homogeneity of $\rho_D$. If a compact part $A$ of
  $\{\ell_i=T\}$, $T > 0$, had positive mass, its dilates $cA$ at distinct $c
  \in [1,2]$ would be disjoint, have a common positive lower mass bound, and lie
  in one compact set. This contradicts local finiteness. Exhausting the face by
  compacts proves it has zero mass. For an empty chain the positive-level face
  is empty. The finite-horizon region is bounded, and its boundary is contained
  in these horizon faces and the coordinate-zero faces. Weak convergence with a
  compact cutoff, and the bounded continuous multiplier $e^{\sum r_l}\varphi$
  there, proves \eqref{E:spatial-wick-limit}. This proves
  Lemma~\ref{L:spatial-wick-measures}.
\end{proof}

Lemma~\ref{L:spatial-wick-measures} removes only the spatial widths. We now
remove the temporal widths too and construct each fixed-order term.

\begin{lemma}[Joint diagram limits and fixed-order construction]
  \label{L:joint-fixed-order}
  Under Hypothesis~\ref{H:homog}, every two-chain diagram with $m+n = 2k \ge 2$
  has a finite joint limit as all its positive vertexwise spatial and temporal
  widths tend to zero, at arbitrary relative rates. This includes unequal and
  empty chains. For every $n \ge 1,t > 0,x$, $Y_n^{\boldsymbol a}(t,x)$ has a
  unique joint $L^2$ limit, denoted $X_n(t,x)$; it is the multiple Stratonovich
  integral of $g_n$. Green smoothing in dimension three is removed before noise
  smoothing.
\end{lemma}

\begin{proof}
  We first remove spatial smoothing at a fixed bounded temporal cutoff. A
  marked-edge product estimate makes the temporal error uniformly small in all
  spatial widths. This gives joint diagram limits, and the mixed second
  moments then make each fixed-order random integral Cauchy in $L^2$.

  \medskip
  \noindent\textit{Step 1.~}We take the fixed-temporal-cutoff spatial limit.
  Use the gap variables and spatial measures of
  Lemma~\ref{L:spatial-wick-measures}. With $s_j = s_j(\boldsymbol r)$ define,
  on $\ell_1 < T_1,\ell_2 < T_2$,
  \[
    F_{\boldsymbol\epsilon,\boldsymbol\delta,D}
    =\int\prod_{\{j,h\}\in D}
    \gamma^0_{\delta_j\vee\delta_h}(s_j-s_h)
    \ud\rho_{\boldsymbol\epsilon,D}.
  \]
  Here $\gamma^0_\delta$ denotes the capped kernel $\kappa_\delta$
  in~\eqref{E:pairing-source-covariance}; off the time diagonal,
  $|r|^{-\alpha_0}-\gamma^0_\eta(r)=K_{1/\eta}(r)$ is its omitted Laplace
  tail.
  At a common temporal width $\eta > 0$ the multiplier is bounded and
  continuous. Hence, the spatial limit exists at every $T_1,T_2 > 0$
  by~\eqref{E:spatial-wick-limit}.

  \medskip
  \noindent\textit{Step 2.~}We bound the temporal remainder uniformly in
  spatial widths.
  We next bound removal of that temporal cutoff uniformly in all spatial widths.
  Write $F_{\boldsymbol\epsilon,0,D}$ for the integral with full factors
  $|s_j-s_h|^{-\alpha_0}$. It is finite at positive spatial widths: bound the
  spatial covariances by their finite suprema and the Green masses by their
  horizon bounds, then integrate the disjoint-pair time factors on the full time
  cube. Each factor $|s-r|^{-\alpha_0}$ is integrable since $\alpha_0 < 1$. The
  elementary difference-of-products inequality bounds
  $F_{\boldsymbol\epsilon,0,D}-F_{\boldsymbol\epsilon,\eta,D}$ by a sum of $k$
  nonnegative integrals, each with one marked time covariance $K_{1/\eta}$ and
  all other time covariances full.

  Apply Lemma~\ref{L:pairing-product-envelope} to the damped kernels
  $e^{-s}G_s$. At first smooth the Green kernels in dimension three and cap the
  temporal Laplace domains above. Each spatial edge measure is a finite
  Gaussian-weighted $\mu$ with nonnegative physical covariance, and each
  temporal domain is a Laplace interval or its restriction. These are exactly
  the hypotheses of that lemma. The finite-horizon exponential loss is at most
  $e^{T_1+T_2}$. Lemma~\ref{L:damped-homogeneous-norms} bounds all unmarked
  norms independently of the smoothing parameters. The marked edge appears once
  in the first-norm product or twice in the second-norm product in every
  admissible assignment. Consequently,
  \begin{equation}\label{E:joint-temporal-tail}
    0\le F_{\boldsymbol\epsilon,0,D}
    -F_{\boldsymbol\epsilon,\eta,D}
    \le C_{k,T_1,T_2}
    \left(N_1^{(1/\eta)}(1)+N_2^{(1/\eta)}(1)^2\right),
  \end{equation}
  uniformly in all spatial widths. Constants may depend on the fixed covariance
  and $k,T_1,T_2$; no uniformity in these quantities is needed. Let the upper
  temporal caps tend to infinity by positivity, and then remove Green smoothing
  at these fixed spatial widths. For the latter passage, the bounded spatial
  covariance products, weak convergence of Green measures and the integrable
  time-cube majorant just used give convergence. The bound survives both
  passages. Its right side tends to zero as $\eta\downarrow0$.

  \medskip
  \noindent\textit{Step 3.~}We squeeze all vertexwise cutoff choices to the
  same diagram limit.
  Denote the fixed-$\eta$ spatial limit by $b_D(\eta)$. These numbers increase
  as $\eta$ decreases, and~\eqref{E:joint-temporal-tail} makes them Cauchy with
  a finite limit $b_D$. If every $\delta_j < \eta$, positivity gives
  \[
    F_{\boldsymbol\epsilon,\eta,D}
    \le F_{\boldsymbol\epsilon,\boldsymbol\delta,D}
    \le F_{\boldsymbol\epsilon,0,D}.
  \]
  Fix $\eta$, pass to the spatial limit in the lower term, and use the uniform
  remainder in~\eqref{E:joint-temporal-tail}; then let $\eta\downarrow0$. This
  proves the asserted joint limit $b_D$ with arbitrary spatial/temporal rates.
  If one chain is empty its horizon contributes no variables, and the same proof
  applies.

  \medskip
  \noindent\textit{Step 4.~}We obtain the joint mean-square limit from the
  mixed moments.
  At positive widths Wick's formula gives the mixed second moments of two copies
  of $Y_n$ as the finite sum of these diagrams, after reversing both terminal
  chains at their common horizon $t$. Evenness of the covariances absorbs the
  sign changes. Every such mixed moment, including either diagonal moment, tends
  to the same finite sum $Q_n = \sum_{D\in\Pi_n}b_D$. Hence,
  \[
    \|Y_n^{\boldsymbol a}-Y_n^{\boldsymbol b}\|_2^2
    \longrightarrow Q_n+Q_n-2Q_n=0.
  \]
  Completeness of $L^2$ gives a unique limit independent of all width ratios. At
  fixed positive noise widths the Green-width limit was already taken in the
  terminal construction, so the result has the required Green-before-noise
  convention. Restricting to common noise widths gives precisely the multiple
  Stratonovich definition for $g_n$. This proves
  Lemma~\ref{L:joint-fixed-order}.
\end{proof}

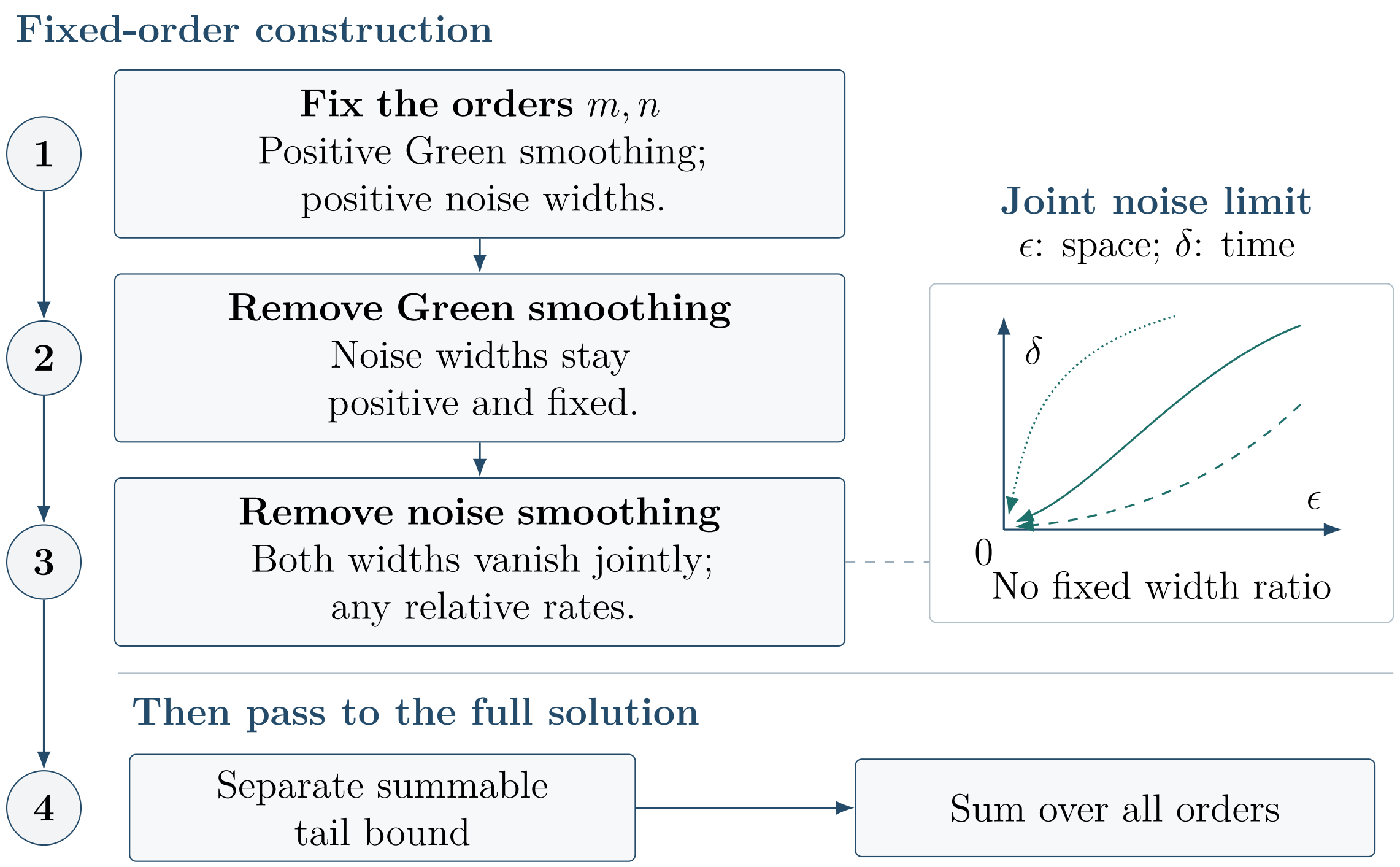
\begin{figure}
  \centering
  % Schematic of the construction order; paths show no convergence rates.
\begin{tikzpicture}[x=1cm,y=1cm,font=\normalsize,
  stage/.style={ham/box,font=\normalsize,text width=7.45cm,
  minimum height=1.65cm,inner sep=6pt},
  number/.style={circle,draw=hamblue,fill=hamblue!6,
  minimum size=8mm,inner sep=0pt,font=\normalsize\bfseries}]
  \node[anchor=west,font=\normalsize\bfseries,text=hamblue]
  at (0,0.55) {Fixed-order construction};
  \node[number] (n1) at (0.45,-0.8) {1};
  \node[stage] (start) at (5.15,-0.8)
  {\textbf{Fix the orders $m,n$}\\
  Positive Green smoothing;\\
  positive noise widths.};
  \node[number] (n2) at (0.45,-3.0) {2};
  \node[stage] (green) at (5.15,-3.0)
  {\textbf{Remove Green smoothing}\\
  Noise widths stay\\
  positive and fixed.};
  \node[number] (n3) at (0.45,-5.2) {3};
  \node[stage] (noise) at (5.15,-5.2)
  {\textbf{Remove noise smoothing}\\
  Both widths vanish jointly;\\
  any relative rates.};
  \draw[ham/arrow] (start.south)--(green.north);
  \draw[ham/arrow] (green.south)--(noise.north);
  \draw[ham/arrow] (n1.south)--(n2.north);
  \draw[ham/arrow] (n2.south)--(n3.north);

  \node[font=\normalsize\bfseries,text=hamblue] at (12.45,-1.3)
  {Joint noise limit};
  \node at (12.45,-1.8) {$\epsilon$: space; $\delta$: time};
  \draw[hamblue!35,rounded corners=2pt] (10.0,-2.2)
  rectangle (15.0,-5.85);
  % Positive spatial and temporal widths approach the same origin along
  % different schematic routes. No rate, scale or covariance is encoded.
  \coordinate (origin) at (10.8,-4.85);
  \draw[ham/arrow] (origin)--(14.45,-4.85)
  node[above left=3pt] {$\epsilon$};
  \draw[ham/arrow] (origin)--(10.8,-2.55)
  node[below right=3pt] {$\delta$};
  \node[anchor=north east,inner sep=3pt] at (origin) {$0$};
  \draw[ham/arrow,draw=hamteal] (14.0,-2.65)
  .. controls (12.8,-3.1) and (11.8,-4.4) .. (10.92,-4.77);
  \draw[ham/arrow,draw=hamteal,dashed] (14.0,-3.5)
  .. controls (12.9,-4.55) and (11.7,-4.75) .. (10.92,-4.82);
  \draw[ham/arrow,draw=hamteal,densely dotted] (12.65,-2.55)
  .. controls (11.2,-2.95) and (11.05,-3.8) .. (10.85,-4.7);
  \node[align=center] at (12.5,-5.45) {No fixed width ratio};
  \draw[hamblue!45,dashed] (noise.east)--(10,-5.2);

  \node[number] (n4) at (0.45,-7.85) {4};
  \draw[ham/arrow] (n3.south)--(n4.north);
  \draw[hamblue!35] (1.25,-6.4)--(15,-6.4);
  \node[anchor=west,font=\normalsize\bfseries,text=hamblue]
  at (1.25,-6.85) {Then pass to the full solution};
  \node[ham/box,font=\normalsize,text width=5.1cm,
  minimum height=1.05cm] (tail) at (4.1,-7.85)
  {Separate summable tail bound};
  \node[ham/box,font=\normalsize,text width=5.25cm,
  minimum height=1.05cm] (sum) at (12.0,-7.85)
  {Sum over all orders};
  \draw[ham/arrow] (tail.east)--(sum.west);
\end{tikzpicture}
  \caption[The order of construction limits]{The order of the construction
    limits. The Green smoothing is removed while the external noise widths
    remain positive. Fixed-order convergence and a separate summable bound then
    identify the full solution. See \ProofRef{Lemma}{L:joint-fixed-order} and
    \ProofRef{Lemma}{L:uniform-series-tail}. Later source limits are taken at
    fixed radius before the large-time limit.}
  \label{Fig:order-of-limits}
\end{figure}

Since $\beta > 2$ and $\alpha+\alpha_0 < 2$,
Lemma~\ref{L:damped-homogeneous-norms} gives constants $C_0 < \infty$ and
$\theta_0 = 1$ such that, for $\theta \ge \theta_0$,
\begin{equation}\label{E:crude-norms}
  N_0(\theta)=\theta^{-2},\qquad
  N_1(\theta)\le C_0,\qquad N_2(\theta)\le C_0\theta^{-1}.
\end{equation}

\begin{lemma}[Corrected crude estimate]\label{L:crude-corrected}
  Under Hypothesis~\ref{H:homog}, for suitable $C < \infty$ and $\theta_0 > 0$,
  \begin{equation}\label{E:crude-corrected}
    \E X_n(t,x)^2\le e^{2\theta_0t}\frac{C^nt^{2n}}{n!}
    \quad(n\ge1,t>0),\qquad
    \E u^2(t,x)\le2e^{2\theta_0t+2Ct^2}.
  \end{equation}
\end{lemma}

\begin{proof}
  We first damp the two chains and apply the finite pairing envelope uniformly
  in the regularizations. We then choose the damping parameter as a function
  of the order and sum the resulting norm bounds.

  \medskip
  \noindent\textit{Step 1.~}We obtain the damped pairing estimate and remove
  regularizations in order.
  On the two finite chains, the product of the exponential weights telescopes to
  $e^{-\theta(s_n+s_{2n})} \ge e^{-2\theta t}$. Extend the nonnegative integral
  to the two infinite time simplices and apply
  Lemma~\ref{L:pairing-product-envelope}. For every admissible product, put $q_i
  = \#Q_i$. Then $q_0 = q_1$ and $q_0+q_2/2 = n$. Hence,~\eqref{E:crude-norms}
  bounds each pairing by $e^{2\theta t}C^n\theta^{-2n}$. Summing the $P_n$
  pairings gives
  \begin{equation}\label{E:theta-envelope}
    \E X_n(t,x)^2\le P_n e^{2\theta t}C^n\theta^{-2n},
    \qquad\theta\ge\theta_0.
  \end{equation}
  Equations (1.5) and (2.24) here belong
  to~\cite{chen:author2026:time-dependency}*{(1.5), p.~3, and (2.24), p.~13}.
  Under Hypothesis~\ref{H:homog}, the ball growth $\mu(B_r)=Vr^\alpha$ and
  $\alpha+\alpha_0<2$ give the integrability in (1.5), while
  Remark~\ref{R:source-case2-time-shift}(iii) directly proves finiteness of
  the Gaussian-weighted measure.
  For the finite-spectral-measure hypothesis in
  Lemma~\ref{L:pairing-product-envelope}, first apply it with the spatially
  smoothed covariance and the temporal truncations of (2.24). The smoothed
  spectral measure is finite: the Gaussian multiplier is bounded by a constant
  times $(1+|\xi|^2)^{-(2-\alpha_0)/2}$, whose integral is finite by (1.5). This
  argument requires no radial density for the spectral measure. The zero norm is
  unchanged; the first and second norms are bounded by
  Lemma~\ref{L:damped-homogeneous-norms}: their Gaussian multipliers can be
  removed in the upper bounds for the nonnegative transforms and covariance
  quadratic forms. The bound is uniform in these regularizations. In $d = 3$,
  first use $G_\eta = p_\eta*G$, as defined on pp.~13--14 of the source. Its
  spatial Fourier multiplier is bounded by one and its spatial mass equals that
  of $G$, so the same bounds are uniform in $\eta$ as well. For each fixed noise
  regularization $(\epsilon,\delta)$, first let $\eta\downarrow0$ in the
  source's $L^2$ definition. The squared moments converge, so the inequality
  survives this first limit. Then let $(\epsilon,\delta) \to 0$ using
  Lemma~\ref{L:joint-fixed-order}. In $d = 1,2$ only this latter limit is
  needed. These limits transfer a uniform inequality; no pointwise delta
  covariance bound is being used.

  \medskip
  \noindent\textit{Step 2.~}We optimize a permitted damping parameter and sum
  the order bounds.
  Choose $\theta = \theta_0+n/t$ in~\eqref{E:theta-envelope}. Since $(2n)! \le
  (2n)^{2n}$, we have $P_n n^{-2n} \le 2^n/n!$. Absorb $2^ne^{2n}$ in $C^n$ to
  obtain the first bound in~\eqref{E:crude-corrected}. Its square roots are
  summable for each fixed $t$, by the ratio test. Thus, the norm series
  converges independently of Proposition~\ref{P:mixed}. By Minkowski and
  Cauchy--Schwarz,
  \[
    \begin{split}
      \|u(t,x)\|_2
      &\le e^{\theta_0t}\sum_{n\ge0}\frac{(\sqrt C t)^n}{\sqrt{n!}},\\
      \left(\sum_{n\ge0}\frac{z^n}{\sqrt{n!}}\right)^2
      &\le\left(\sum_{n\ge0}2^{-n}\right)
      \left(\sum_{n\ge0}\frac{(2z^2)^n}{n!}\right)
      =2e^{2z^2}.
    \end{split}
  \]
  This proves Lemma~\ref{L:crude-corrected}.
\end{proof}

\section{The mild equation with the corrected bound}

Lemma~\ref{L:crude-corrected} states bounds for the limits $X_n$ and their
sum. We now bound the regularized terms uniformly in all positive vertexwise
widths; Proposition~\ref{P:mild-construction} uses this bound.

\begin{lemma}[Uniform vertexwise regularization bound]\label{L:vertex-crude}
  Under Hypothesis~\ref{H:homog}, there are $C < \infty$ and $\theta_0 > 0$ such
  that, for all positive vertexwise widths and $n \ge 1$,
  \begin{equation}\label{E:vertex-crude}
    \|Y_n^{\boldsymbol a}(t,x)\|_2\le q_n(t),\qquad
    q_n(t):=e^{\theta_0t}\frac{C^{n/2}t^n}{\sqrt{n!}}.
  \end{equation}
  The bound also holds for $X_n(t,x)$, with $q_0(t) = e^{\theta_0t} \ge 1$. For
  each $T < \infty$,
  \begin{equation}\label{E:vertex-summability}
    \sum_{n\ge0}\sqrt{n+1}\,q_n(T)<\infty.
  \end{equation}
\end{lemma}

\begin{proof}
  We retain the separate width on every vertex while applying the edgewise
  product envelope. Uniform edge norms then give the factorial bound, which
  passes to the joint limit and supplies the weighted summability needed for
  multiplication.

  \medskip
  \noindent\textit{Step 1.~}We identify all edge regularizations and bound
  their norms.
  Wick's formula \cite{chen:author2026:time-dependency}*{(2.19), p.~12}, with
  the joint covariance (2.11), gives the square of $Y_n^{\boldsymbol a}$, with
  the widths copied to the second chain. Its edge covariances are
  \[
    \gamma^0_{\delta_j\vee\delta_l}(s_j-s_l)
    \gamma_{\epsilon_j+\epsilon_l}(y_j-y_l).
  \]
  On each terminal chain use $s'_j = t-s_{n+1-j}$ and $y'_j = x-y_{n+1-j}$. This
  reverses the chain and only permutes its vertex widths; covariance differences
  change sign and the kernels are even. For each pairing apply
  Lemma~\ref{L:pairing-product-envelope} to $e^{-\theta r}G_r$, first with Green
  smoothing if $d = 3$. The lemma permits a separate finite spectral measure and
  temporal truncation on each edge. These hypotheses hold here: the spatial
  spectral measure on an edge is
  $e^{-(\epsilon_j+\epsilon_l)|\xi|^2/2}\mu(\ud\xi)$, which is finite by (1.5),
  and its temporal Laplace domain is $[0,(\delta_j\vee\delta_l)^{-1}]$.

  Condition (1.5) is the source condition just identified in the preceding
  proof; under Hypothesis~\ref{H:homog}, finiteness of each Gaussian-weighted
  edge measure follows directly from
  Remark~\ref{R:source-case2-time-shift}(iii).
  All edge norms obey~\eqref{E:crude-norms}, uniformly in their widths. For the
  first norm the Laplace--Fourier transform of the damped Green kernel is
  positive, and the spatial Gaussian multiplier and temporal truncation can only
  decrease its nonnegative spectral integral. For the second norm first enlarge
  the temporal covariance to $|s-r|^{-\alpha_0}$, using positivity of the Green
  kernels and the smoothed spatial covariance. Its space-time spectral
  expression is an integral of the squared absolute transform against
  nonnegative measures, so removal of the spatial Gaussian multiplier increases
  it. Green smoothing has multiplier at most one and unchanged mass, so the same
  bounds hold before its removal. These assertions follow from
  Lemma~\ref{L:damped-homogeneous-norms}, rather than an imported normalization
  or large-parameter formula in the source's Appendix. The constants are
  independent of all the vertex widths.

  \medskip
  \noindent\textit{Step 2.~}We count the norm factors and optimize the
  damping.
  Since $\#Q_0 = \#Q_1$ and $\#Q_0+\#Q_2/2 = n$ for every admissible product,
  exactly as in~\eqref{E:theta-envelope} we obtain the following bound. Indeed,
  put $D = \max\{1,2C_0\}$. The norm constants contribute at most
  $D^{\#Q_1+\#Q_2} \le D^{2n}$, so we may take $C_1 = D^2$:
  \[
    \|Y_n^{\boldsymbol a}(t,x)\|_2^2
    \le P_n e^{2\theta t}C_1^n\theta^{-2n},\qquad\theta\ge\theta_0.
  \]
  Let the Green width tend to zero at these fixed positive noise widths, and
  choose $\theta = \theta_0+n/t$. The inequality $P_n n^{-2n} \le 2^n/n!$ then
  proves \eqref{E:vertex-crude}, after increasing $C$.
  Lemma~\ref{L:joint-fixed-order} gives the joint noise-width limit $X_n$; the
  same bound passes to it. Finally,
  \[
    \frac{\sqrt{n+2}\,q_{n+1}(T)}
    {\sqrt{n+1}\,q_n(T)}
    =\sqrt C\,T\frac{\sqrt{n+2}}{n+1}\longrightarrow0.
  \]
  The ratio test proves~\eqref{E:vertex-summability}. This proves
  Lemma~\ref{L:vertex-crude}.
\end{proof}

We now show that the series of the fixed-order terms solves the mild
equation.

\LeanMarker{P:mild-construction}
\begin{proposition}[Mild construction from the fixed-order limits]
  \label{P:mild-construction}
  Under Hypothesis~\ref{H:homog}, the series $u(t,x) = \sum_{n\ge0}X_n(t,x)$ is
  a mild solution in~\cite{chen:author2026:time-dependency}*{Definition~2.1 and
  (2.12), pp.~9--10}. More precisely, for each $t > 0$ and $x \in \R^d$,
  \begin{equation}\label{E:mild-construction}
    \lim_{\epsilon,\delta\downarrow0}
    \int_0^t\int_{\R^d}G(t-s,x-\ud y)
    u(s,y)W_{\epsilon,\delta}(s,y)\ud s
    =u(t,x)-1\quad\text{in }L^2(\Omega).
  \end{equation}
  The limit allows arbitrary relative rates of the two widths. In dimension
  three the Green-width limit precedes the noise-width limit. No uniqueness
  assertion is made.
\end{proposition}

For $d = 3$ we read the source's generalized-integrand definition using its
Green regularization, with the random factor held fixed. The proof below
establishes that limit for the product $uW_{\epsilon,\delta}$.

\begin{proof}
  We first control multiplication by one Gaussian variable on bounded
  polynomial degrees. We then identify a partially regularized terminal
  vertex, respecting the Green-first convention. Joint fixed-order convergence
  and uniform order tails finally give the full mild integral.

  \medskip
  \noindent\textit{Step 1.~}We obtain the Gaussian polynomial multiplication
  bound.
  Let $\mathscr{P}_{\le n}$ be the closed Gaussian polynomial space of total
  degree at most $n$, generated by all the coupled source fields. It includes
  all lower-degree contractions. The fixed-width terminal integrals lie in this
  space, by their Bochner construction; their fixed-order $L^2$ limits $X_n$ do
  too.

  We first record an elementary multiplication bound on these spaces. If $H$ is
  any centered linear Gaussian variable and $F \in \mathscr{P}_{\le n}$, then
  \begin{equation}\label{E:occupied-multiplication}
    \|HF\|_2\le2\sqrt{n+1}\,\|H\|_2\|F\|_2.
  \end{equation}
  If $H = 0$ the assertion is immediate. Otherwise take $H/\|H\|_2$ as the first
  coordinate of an orthonormal Gaussian basis. On normalized Hermite polynomials
  its creation and annihilation operators have norms at most $\sqrt{n+1}$ and
  $\sqrt n$ on total degrees at most $n$: the squared coefficient sums have
  respective weights $a_1+1 \le n+1$ and $a_1 \le n$. Their sum proves the bound
  for finite-coordinate polynomials. Approximation in the degree-closed space
  proves it for $F$ as well. Indeed, the products of an approximating sequence
  are $L^2$ Cauchy by this bound; subsequences converging almost surely identify
  their limit with $HF$. This also proves the closedness of multiplication
  needed below. The argument allows $F$ to depend on $H$.

  \medskip
  \noindent\textit{Step 2.~}We obtain measurable fields and spatial
  mean-square continuity.
  For $t>0$, divide every ordered time and spatial Green increment by $t$,
  using $G_{tr}(t\,\ud z)=tG_r(\ud z)$; the time Jacobians and Green masses
  contribute $t^{2n}$, leaving a fixed unit-time chain domain with $(t,x)$ in
  the Gaussian arguments. At $t=0$ all positive-order terms vanish.
  Lemma~\ref{L:vertex-crude} proves absolute $L^2$ convergence of the defining
  series for $u$, locally uniformly in time in the displayed norm bounds. All
  fields are strongly measurable as $L^2(\Omega)$ vectors. At positive widths,
  rescale the ordered times and Green increments to fixed unit-time domains. The
  Gaussian evaluations are jointly measurable, and their bounded moments and the
  finite Green masses give Bochner domination on bounded time intervals.
  Parameter integration therefore gives joint strong measurability in $(s,y)$.
  Pointwise $L^2$ limits along a fixed countable sequence of widths give this
  property for $X_n$, then for the series. For later Green limits we also need
  spatial continuity at each fixed time. Positive-width terminal terms are
  $L^2$-continuous in space, by continuity of their bounded smoothed covariance
  products and dominated integration against the finite Green measures. Joint
  spatial stationarity of the coupled regularizations shows that
  \[
    \|Y_n^{\boldsymbol b}(s,y)-X_n(s,y)\|_2
    =\|Y_n^{\boldsymbol b}(s,0)-X_n(s,0)\|_2.
  \]
  Thus, their convergence is uniform in $y$ for fixed $s$. Consequently,
  $X_n(s,\cdot)$, and by summability $u(s,\cdot)$, have continuous $L^2$
  representatives.

  \medskip
  \noindent\textit{Step 3.~}We identify the partially regularized terminal
  integral.
  Fix an external width $a = (\epsilon,\delta) > 0$ and write $W_a =
  W_{\epsilon,\delta}$, with variance $V_a =
  \gamma^0_\delta(0)\gamma_{2\epsilon}(0) < \infty$. For $n \ge 1$ define the
  Bochner integral
  \[
    I_{n,a}(t,x)=\int_0^t\int_{\R^d}G(t-s,x-\ud y)
    W_a(s,y)X_{n-1}(s,y)\ud s.
  \]
  It exists because~\eqref{E:occupied-multiplication} gives
  \begin{equation}\label{E:fixed-external-domination}
    \|W_a(s,y)X_{n-1}(s,y)\|_2
    \le2\sqrt n\sqrt{V_a}\,q_{n-1}(T),\qquad0\le s\le T,
  \end{equation}
  and the spatial Green mass is $t-s$. At fixed $a$, multiplication of the inner
  regularizations by $W_a$ preserves their $L^2$ convergence, by the same bound
  on degree at most $n-1$. Their norms have the uniform domination just used.
  Dominated Bochner integration therefore gives
  \begin{equation}\label{E:partial-vertex-limit}
    I_{n,a}(t,x)=\lim_{\boldsymbol b\downarrow0}
    Y_n^{(\boldsymbol b,a)}(t,x)\quad\text{in }L^2.
  \end{equation}
  Here $a$ is at the last, terminal vertex and $\boldsymbol b$ at the preceding
  $n-1$ vertices. At positive widths this identity before passage to the limit
  is ordinary terminal-kernel Fubini; its moment bound was established before
  Lemma~\ref{L:vertex-crude}. For $n = 1$ there are no inner widths and the
  identity is exact. In particular,~\eqref{E:vertex-crude} passes to this inner
  limit:
  \begin{equation}\label{E:mixed-vertex-bound}
    \sup_{a>0}\|I_{n,a}(t,x)\|_2\le q_n(t).
  \end{equation}

  \medskip
  \noindent\textit{Step 4.~}We retain the dimension-three Green-limit
  convention.
  We spell out the dimension-three convention in these steps. The product
  $W_a(s,y)X_{n-1}(s,y)$ is $L^2$-continuous in $y$.
  Apply~\eqref{E:occupied-multiplication} first to
  $W_a(s,y)(X_{n-1}(s,y)-X_{n-1}(s,y'))$ and then to
  $(W_a(s,y)-W_a(s,y'))X_{n-1}(s,y')$. The two factors tend to zero in the
  required norms, by the continuity just proved and the smoothed covariance of
  $W_a$. The bound~\eqref{E:fixed-external-domination} is uniform in $y$.
  Therefore, convolution against $p_\eta*G_{t-s}$ converges in $L^2$ to
  convolution against $G_{t-s}$: on a compact set approximate this continuous
  vector function by a finite continuous partition of unity times fixed vector
  values, use scalar weak convergence, and control the complement by tightness
  and the uniform bound. The mass $t-s$ and~\eqref{E:fixed-external-domination}
  permit dominated time integration. This proves the Green-width limit at fixed
  external widths. The same argument applies before the inner noise limit; all
  $Y_n^{(\boldsymbol b,a)}$ already use their Green limit at positive noise
  widths. Thus,~\eqref{E:partial-vertex-limit} does not interchange the source's
  Green and noise limits.

  \medskip
  \noindent\textit{Step 5.~}We remove the external width using joint
  convergence.
  For each fixed $n$, the joint fixed-order convergence in
  Lemma~\ref{L:joint-fixed-order} implies
  \begin{equation}\label{E:external-vertex-limit}
    \lim_{a\downarrow0}I_{n,a}(t,x)=X_n(t,x)\quad\text{in }L^2.
  \end{equation}
  To see the role of joint convergence, fix $\zeta > 0$ and choose $r > 0$ so
  that every vertex width below $r$ gives $\|Y_n^{(\boldsymbol b,a)}-X_n\|_2 <
  \zeta$. Fix any external pair $a$ below $r$ and let the inner widths tend to
  zero in~\eqref{E:partial-vertex-limit}. Then $\|I_{n,a}-X_n\|_2 \le \zeta$.
  This proves the asserted limit with arbitrary relative rates; it does not
  assume a partial-limit theorem.

  \medskip
  \noindent\textit{Step 6.~}We sum the terminal integrals and remove the final
  uniform tail.
  It remains to identify the full product integral, rather than only its
  finite-order terms. At fixed $a$, \eqref{E:occupied-multiplication}
  and~\eqref{E:vertex-summability} imply
  \[
    \sum_{n\ge1}\|W_a(s,y)X_{n-1}(s,y)\|_2
    \le2\sqrt{V_a}\sum_{n\ge1}\sqrt n\,q_{n-1}(T)<\infty,
    \qquad0\le s\le T.
  \]
  Multiplication is closed: partial sums of $u$ converge in $L^2$, and their
  products converge by this estimate, so an almost surely convergent subsequence
  identifies the product limit with $W_a u$. Dominated Green integration and
  norm summability now prove
  \[
    I_a(t,x):=\int_0^t\int_{\R^d}G(t-s,x-\ud y)
    W_a(s,y)u(s,y)\ud s
    =\sum_{n\ge1}I_{n,a}(t,x)\quad\text{in }L^2.
  \]
  The sum of the product integrands is spatially $L^2$-continuous by uniform
  summability, so its Green-width limit is justified in exactly the same way.
  Finally, for every $N$,
  \[
    \|I_a(t,x)-(u(t,x)-1)\|_2
    \le\sum_{n=1}^N\|I_{n,a}(t,x)-X_n(t,x)\|_2
    +2\sum_{n>N}q_n(t).
  \]
  First let $a\downarrow0$ using~\eqref{E:external-vertex-limit}, then let $N
  \to \infty$ using~\eqref{E:vertex-summability}. This
  proves~\eqref{E:mild-construction}, hence the mild equation in the stated
  definition. The source's (5.3)--(5.4) and its Section~5(ii) identification
  proof have not been used. This proves Proposition~\ref{P:mild-construction}.
\end{proof}

Thus, the exponent two enters through replacing the paired norm products by
$C\theta^{-2}$, followed by the pairing count and the order summation. The sum
over unequal chain lengths is not responsible for this loss: the same loss
already occurs in the diagonal estimate and Minkowski's inequality.

\section{What homogeneity recovers, and what the bound still loses}

Under homogeneity the norms scale exactly in $\theta$. We now use this
scaling in place of~\eqref{E:crude-norms}.

\begin{proposition}[The uniform pairing estimate]\label{P:envelope}
  Under Hypothesis~\ref{H:homog}, there is $K > 0$ such that
  \begin{equation}\label{E:homogeneous-envelope}
    \E[X_m(t,x)X_n(t,x)]
    \le P_k K^k e^{2\theta t}\theta^{-\beta k},
    \qquad m+n=2k\ge2,\quad\theta>0.
  \end{equation}
  Define the numerical majorant produced by applying this estimate separately at
  every total order by
  \begin{equation}\label{E:B-envelope}
    B_K(t)=1+\sum_{k\ge1}(2k+1)P_kK^k
    \inf_{\theta>0}\{e^{2\theta t}\theta^{-\beta k}\}.
  \end{equation}
  Then $\E u^2(t,x) \le B_K(t)$ and, for some $c,C > 0$ and all large $t$,
  \begin{equation}\label{E:envelope-growth}
    ct^{\rho}\le\log B_K(t)\le Ct^{\rho},\qquad
    \rho=\frac{\beta}{\beta-1}
    =\frac{4-\alpha-\alpha_0}{3-\alpha-\alpha_0}>\kappa.
  \end{equation}
  The lower bound here concerns the numerical majorant $B_K$, not the second
  moment.
\end{proposition}

\begin{proof}
  We use homogeneous norm scaling to bound each total order, then estimate the
  resulting numerical series from above and below. The lower estimate concerns
  that majorant only.

  \medskip
  \noindent\textit{Step 1.~}We obtain the homogeneous pairing majorant.
  The wave kernel scales as $G(cs,c\ud y) = cG(s,\ud y)$, and the covariance has
  degree $-\alpha$. Changing all time and spatial variables by $\theta^{-1}$ in
  the norms yields
  \begin{equation}\label{E:homogeneous-norms}
    N_0(\theta)=\theta^{-2},\quad
    N_1(\theta)=N_1(1)\theta^{\alpha+\alpha_0-2},\quad
    N_2(\theta)^2=N_2(1)^2\theta^{-\beta}.
  \end{equation}
  All constants and scalings follow directly from
  Lemma~\ref{L:damped-homogeneous-norms}. The products $2N_0N_1$ and $N_2^2$
  therefore both have degree $-\beta$. Choose $K \ge \max\{2N_1(1),N_2(1)^2\}$,
  increasing it to a positive constant if needed.
  Lemma~\ref{L:pairing-product-envelope} and $q_0+q_2/2 = k$ give
  $K^k\theta^{-\beta k}$ for every admissible product. The finite horizons
  contribute at most $e^{2\theta t}$. Apply this first to regularized moments.
  In $d = 3$, fix $(\epsilon,\delta)$ and let $\eta\downarrow0$ in each factor's
  $L^2$ limit; Cauchy--Schwarz gives $L^1$ convergence of their product. Then
  let $(\epsilon,\delta) \to 0$ in the same way, using
  Lemma~\ref{L:joint-fixed-order} for each chain. In $d = 1,2$ only the noise
  limit is required. The estimates are uniform at both stages, by the multiplier
  bounds in Lemma~\ref{L:crude-corrected}. This
  proves~\eqref{E:homogeneous-envelope}. Empty-chain cases satisfy the same
  bound (using $e^{2\theta t}$ is harmless). Proposition~\ref{P:mixed} now gives
  the majorant.

  \medskip
  \noindent\textit{Step 2.~}We determine the growth scale of the numerical
  series.
  The infimum in~\eqref{E:B-envelope} is attained at $\theta = \beta k/(2t)$ and
  equals $e^{\beta k}(2t/(\beta k))^{\beta k}$. Put $q = \beta-1 > 1$. The
  elementary estimates
  \[
    k!\le P_k\le2^k k!,\qquad
    (k/e)^k\le k!\le k^k,\qquad 2k+1\le3^k
  \]
  show that, for fixed positive $a,b$ depending only on $K,\beta$,
  \[
    \sum_{k\ge0}\frac{(at^\beta)^k}{(k!)^q}
    \le B_K(t)\le
    \sum_{k\ge0}\frac{(bt^\beta)^k}{(k!)^q}.
  \]
  Indeed, the additional factor $k!^\beta/k^{\beta k}$ lies between $e^{-\beta
  k}$ and one and can be absorbed into $a^k,b^k$. For $z > 0$, the right-hand
  series is at most $(\sum_{k\ge0}z^{k/q}/k!)^q = e^{qz^{1/q}}$, by the
  elementary $\ell^q$ versus $\ell^1$ inequality for nonnegative sequences. For
  a lower bound choose $k = \lfloor z^{1/q}/2\rfloor$ and use $k! \le k^k$; its
  single summand is at least $2^{qk}$. These estimates
  prove~\eqref{E:envelope-growth}. Finally, $\beta-1 = 3-\alpha-\alpha_0 <
  3-\alpha$, so $\rho > \kappa$. This proves Proposition~\ref{P:envelope}.
\end{proof}

\begin{remark}[Limitation of the uniform pairing estimate]\label{R:G0-verdict}
  Since $\beta > 2$, we have $1 < \rho < 2$. Thus, Proposition~\ref{P:envelope}
  improves the crude quadratic time scale under homogeneity. It does not itself
  prove Theorem~\ref{T:intro-upper}; Section~\ref{S:interval-upper} uses a different
  argument. Theorem~3.1 used through its uniform pairing estimate is \emph{not
  sufficient to prove Theorem~\ref{T:intro-upper}}: even with exact homogeneous norm scaling, the
  resulting majorant has exponent $\rho > \kappa$. The order coefficient has
  factorial denominator $(k!)^{3-\alpha-\alpha_0}$, whereas the desired scale
  requires $(k!)^{3-\alpha}$; the missing factor is $(k!)^{\alpha_0}$, up to
  exponential factors in $k$. Counting the $2k+1$ chain-length splits does not
  cause this loss, since it can be absorbed into an exponential factor. This
  conclusion does not rule out a more refined application retaining temporal
  order or grouping pairings before estimating them, and does not disprove the
  conjecture. In the formal boundary limit $\alpha_0\downarrow0$, both $\rho$
  and $\kappa$ tend to $(4-\alpha)/(3-\alpha)$. This necessary time-independent
  calibration does not distinguish them and asserts no endpoint theorem or
  uniformity of constants in $\alpha_0$. At $p = 1$ use the empty-chain case
  with $e^{\theta t}$ in place of $e^{2\theta t}$; minimizing at $\theta = \beta
  k/t$ changes only a constant to the power $k$, leaving the same factorial
  denominator $(k!)^{3-\alpha-\alpha_0}$. It therefore also loses the factorial
  power $\alpha_0$ relative to $\beta/\kappa = 3-\alpha$, the power
  corresponding to the target series scale. This is a comparison of scales, not
  a coefficient estimate asserted by the source's Theorem~1.2. The loss is not
  specific to the second moment. Both the crude and homogeneous upper bounds are
  consistent with the Jensen lower bound, being at larger time scales. No
  identification with the Skorohod second moment is used: its exclusion of
  self-pairings changes the combinatorics, already visible in the $(0,2)$ and
  $(2,0)$ terms of Lemma~\ref{L:white-test}.
\end{remark}

No-go~\ref{P:nogo-01} isolates this scale obstruction. The adjacent-diagram
calculation in No-go~\ref{P:nogo-02} also explains why one cannot repair it
by assigning a cubic-factorial estimate to every individual pairing.

\paragraph{What we now have.}
One solution, built from its expansion and shown to solve the mild equation
(Proposition~\ref{P:mild-construction}), with a formula for every mixed moment
(Proposition~\ref{P:mixed}). A single bound for all pairings misses the scale
$T^\kappa$; the next chapter reaches it.

\chapter{An independent upper bound}\label{S:book-independent-upper}

Chapter~\ref{S:G0} shows that the second moment is finite. How fast can
it grow? This chapter gives a first answer, before the exact constant is
known: no faster than $\exp(CT^\kappa)$, up to a polynomial factor. This is
Theorem~\ref{T:full-upper}, and it proves Theorem~\ref{T:intro-upper}. Its
constant is not claimed to be optimal.

\begin{figure}
  \centering
  \begin{tikzpicture}[x=1cm,y=1cm]
    \draw[ham/brace] (0,.62)--(3,.62);
    \node[font=\small,anchor=south] at (1.5,.77) {$h_j$};
    \foreach \hamGridRow/\hamPanelShift/\hamGridName in {0/0/{0},1/1/{1/3},2/2/{2/3}} {
    \draw[hamblue] (0,-\hamGridRow)--(11.4,-\hamGridRow);
    \foreach \hamGridIndex in {0,1,2,3} {
    \pgfmathsetmacro{\hamGridBoundary}{3*\hamGridIndex+\hamPanelShift}
    \draw (\hamGridBoundary,-\hamGridRow-.12)--(\hamGridBoundary,-\hamGridRow+.12);
    }
    \node[anchor=east,font=\small] at (-.25,-\hamGridRow) {$t=\hamGridName$};
    }
    \draw[hamteal,line width=3pt] (1,-1)--(4,-1);
    \draw[hamteal] (1,-1.2)--(1,-.8);
    \draw[hamteal] (4,-1.2)--(4,-.8);
    \draw[dashed] (2.8,.2)--(2.8,-2.55);
    \draw[dashed] (3.4,.2)--(3.4,-2.55);
    \node[anchor=north,font=\small] at (2.8,-2.6) {$s$};
    \node[anchor=north,font=\small] at (3.4,-2.6) {$r$};
    \node[ham/note,text width=5.4cm] at (7.4,-2.95)
    {$3|s-r|<h_j\le6|s-r|$\\both points in one interval};
  \end{tikzpicture}
  \caption[Three shifted interval grids]{The three grids at one selected scale.
    Their boundaries are separated by $h_j/3>|s-r|$, so at least one grid has no
    boundary between $s$ and $r$. The full nested construction uses endpoints
    $h_j(n+(-1)^jt)$ for $t=0,1/3,2/3$; the shift alternates with the level. The
    resulting infinite sum dominates $|s-r|^{-a}$ off the diagonal. This is not
    a claim that an arbitrary finite partial sum already dominates it. See
    \ProofRef{Theorem}{T:full-upper}.}
  \label{Fig:dyadic-cover}
\end{figure}
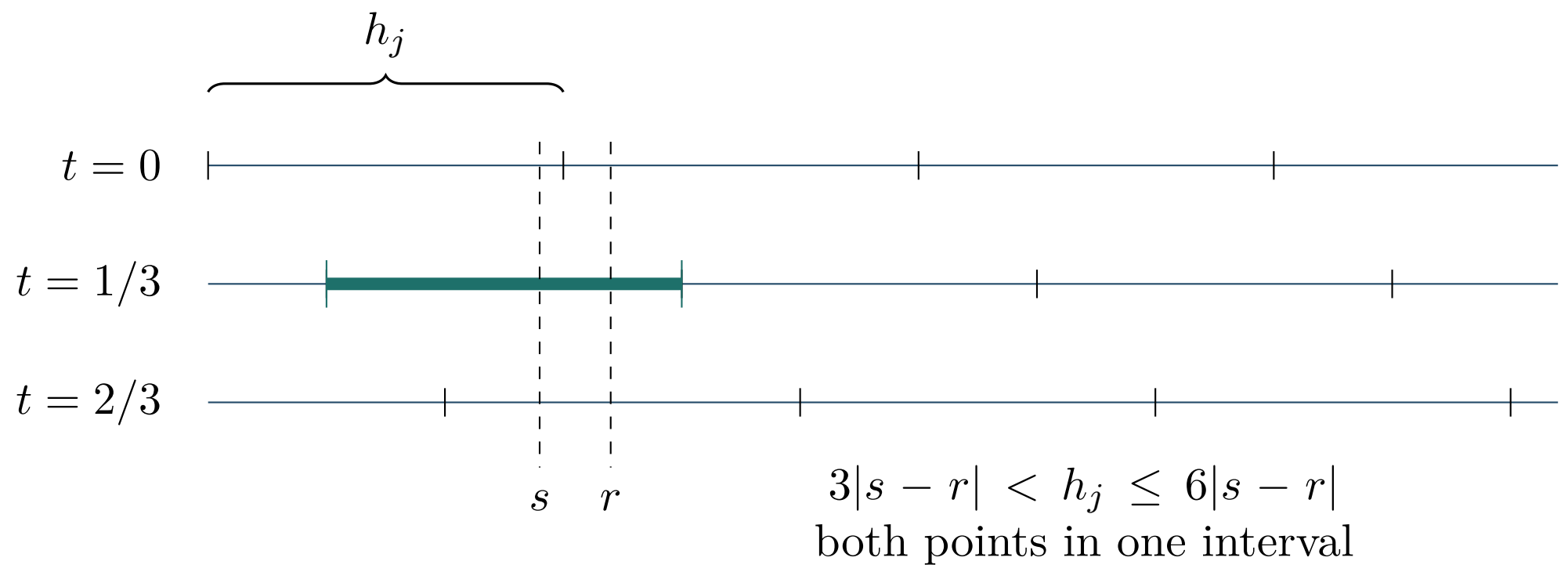

The natural plan is to cut time into blocks and multiply the bounds. A wave
does not allow this easily. Each block must hand on both displacement
and velocity (No-go~\ref{P:nogo-14}), and the noise in one block is
correlated with the noise in the others. Sections~\ref{S:static-phase}
and~\ref{S:static-all-orders} show how far independent blocks
go. Section~\ref{S:interval-upper} then replaces the blocks by Gaussian fields
on overlapping intervals; Figure~\ref{Fig:dyadic-cover} is the picture to
keep in mind.

The model prerequisites are Hypothesis~\ref{H:homog} and the
construction of Chapter~\ref{S:G0}. On a first reading, skip the proofs of
Propositions~\ref{P:static-phase-variance} and~\ref{P:static-phase-all-orders}
and of Corollary~\ref{C:static-interaction}: Theorem~\ref{T:full-upper}
uses only the lemmas of those sections and one calculation, the free-phase
energy bound at the end of the proof of
Proposition~\ref{P:static-phase-variance}.
Episode~\VideoRef{4}{bUyxCpboios} of the video series follows this chapter.

\section{A static phase estimate at first variance order}
\label{S:static-phase}

A time-block argument would have to carry both position and velocity. We
establish a uniform estimate for the first variance coefficient with arbitrary
deterministic initial phase. The covariance may be any member of the original
three-dimensional class; no endpoint representation is assumed. The estimate
below is not an all-order propagator bound.
No-go~\ref{P:nogo-14} shows why scalar solution-value moments cannot be
multiplied across blocks. No-go~\ref{P:nogo-15} shows why both the first
random term and the second-order mean enter the phase energy coefficient.

Let $0 < \alpha < 2$, and let $\Gamma$ be a nonnegative, even, positive-definite
tempered Radon covariance on $\R^3$, satisfying $\Gamma(rA) =
r^{3-\alpha}\Gamma(A)$. Its spectral convention is $\Gamma = \int e^{i\xi\cdot
x}\mu(\ud\xi)$. Put
\begin{equation}\label{E:phase-constants}
  \begin{split}
    V&=\mu(B_1),\qquad M=\Gamma(B_1),\qquad
    R_\alpha=\frac{\alpha\pi V}{2\sin(\pi\alpha/2)},\\
    c_\Gamma&=\frac{(3-\alpha)M}{4\pi},\qquad
    b_\alpha=2^\alpha+\frac{3^{2-\alpha}}{2-\alpha},\qquad
    D_\Gamma=c_\Gamma b_\alpha.
  \end{split}
\end{equation}
Here $B_1$ is the open unit ball in the relevant variable. The physical mass $M$
and the spectral mass $V$ are not identified. Let $p_\epsilon$ be the Brownian
Gaussian density with covariance $\epsilon I$, and set $\Gamma_\epsilon =
p_{2\epsilon}*\Gamma$.

\LeanMarker{L:phase-contractions}
\begin{lemma}[Two static contraction bounds]\label{L:phase-contractions}
  For $\lambda > 0$, $q \in \R^3$, $a > 0$, and $\epsilon > 0$,
  \begin{align}
    \int_{\R^3}\frac{\mu(\ud\xi)}{\lambda^2+|q+\xi|^2}
    &\le R_\alpha\lambda^{\alpha-2},
    \label{E:phase-resolvent}\\
    \left|T_\epsilon(a,q)\right|
    &\le D_\Gamma a^{1-\alpha},
    \label{E:phase-trace}\\
    T_\epsilon(a,q)
    &:=\int_{\R^3}e^{iq\cdot x}\Gamma_\epsilon(x)G_a(\ud x)
    =\int_{\R^3}\widehat{G}_a(q+\xi)e^{-\epsilon|\xi|^2}
    \mu(\ud\xi).\nonumber
  \end{align}
\end{lemma}

\begin{proof}
  We first compare translated heat transforms and integrate them to obtain the
  resolvent bound. For the trace bound, positivity removes the Fourier phase
  and radial integration reduces the calculation to a spherical convolution.

  \medskip
  \noindent\textit{Step 1.~}We bound the shifted resolvent.
  Homogeneity and Fourier uniqueness give $\mu(B_r) = Vr^\alpha$. In particular,
  there is no origin atom. For $s > 0$, Schwartz Fourier pairing yields
  \[
    \int e^{-s|q+\xi|^2}\mu(\ud\xi)
    =\int p_{2s}(x)e^{iq\cdot x}\Gamma(\ud x)
    \le\int p_{2s}(x)\Gamma(\ud x)
    =\int e^{-s|\xi|^2}\mu(\ud\xi).
  \]
  The integrals are finite. Evenness makes the middle integral real, and
  nonnegativity of $\Gamma$ gives the inequality. Multiply by $e^{-\lambda^2s}$
  and integrate in $s$. Nonnegative Tonelli on the spectral side gives the
  resolvent comparison. At $q = 0$ its value is
  \[
    \alpha V\int_0^\infty\frac{r^{\alpha-1}}{\lambda^2+r^2}\ud r
    =R_\alpha\lambda^{\alpha-2}.
  \]
  This proves~\eqref{E:phase-resolvent}.

  \medskip
  \noindent\textit{Step 2.~}We bound the regularized trace using the positive
  spherical integral.
  Since $\mu$ is tempered, its Gaussian weighting is finite. Fourier Fubini
  therefore gives the identity for $T_\epsilon$. The measure
  $\Gamma_\epsilon(x)G_a(\ud x)$ is finite and nonnegative, whence
  $|T_\epsilon(a,q)| \le T_\epsilon(a,0)$. The radial distribution function of
  $\Gamma$ is $\Gamma(B_r) = Mr^{3-\alpha}$. Thus, for every nonnegative radial
  Schwartz function $F$,
  \[
    \int F(x)\Gamma(\ud x)
    =(3-\alpha)M\int_0^\infty F(r)r^{2-\alpha}\ud r
    =c_\Gamma\int_{\R^3}F(x)|x|^{-\alpha}\ud x.
  \]
  Apply this to $F = G_a*p_{2\epsilon}$, which is radial and Schwartz.
  Nonnegative Fubini gives
  \[
    T_\epsilon(a,0)
    =c_\Gamma[p_{2\epsilon}*(G_a*|\cdot|^{-\alpha})](0)
    \le c_\Gamma\sup_z(G_a*|\cdot|^{-\alpha})(z).
  \]
  Write $r = |z|$. The spherical wave measure has mass $a$, and direct
  integration over its normal coordinate gives, for $r > 0$,
  \[
    (G_a*|\cdot|^{-\alpha})(z)
    =\frac{(a+r)^{2-\alpha}-|a-r|^{2-\alpha}}
    {2r(2-\alpha)}.
  \]
  For $r < a/2$, distances from $z$ to the sphere are at least $a/2$, giving the
  bound $2^\alpha a^{1-\alpha}$. For $r > 2a$, they are at least $a$, giving
  $a^{1-\alpha}$. On $a/2 \le r \le 2a$, the displayed formula is at most
  $3^{2-\alpha}a^{1-\alpha}/(2-\alpha)$. These bounds also cover the endpoints
  by continuity and give~\eqref{E:phase-trace}. This proves
  Lemma~\ref{L:phase-contractions}.
\end{proof}

Let $A_\lambda = \lambda^2-\Delta$, and equip deterministic phase pairs with the
norm
\begin{equation}\label{E:phase-norm}
  \|(f,g)\|_\lambda^2
  =\|f\|_2^2+\|A_\lambda^{-1/2}g\|_2^2.
\end{equation}
Let $C_t$ be the Fourier multiplier $\cos(t|\xi|)$. On one static isonormal
Gaussian space with covariance $\Gamma$, use $W_\epsilon(x) =
W(p_\epsilon(x-\cdot))$. Its covariance is exactly $\Gamma_\epsilon(x-y)$, with
spectral measure $\mu_\epsilon(\ud\xi) = e^{-\epsilon|\xi|^2}\mu(\ud\xi)$. For
real Schwartz $f,g$, define the free position $w_0(t) = C_tf+G_t*g$ and its
phase $Z_0(t) = (w_0(t),\dot{w}_0(t))$. The next two regularized static Duhamel
phases are
\begin{equation}\label{E:phase-duhamel}
  Z_{j+1}(t)
  =\int_0^t
  \begin{pmatrix}
    G_{t-s}*(W_\epsilon Z_{j,u}(s))\\
    C_{t-s}(W_\epsilon Z_{j,u}(s))
  \end{pmatrix}\ud s,\qquad j=0,1,
\end{equation}
where $Z_{j,u}$ is the position component. At fixed $\epsilon > 0$ these are
Gaussian-polynomial Bochner integrals in the indicated $L^2$ phase space.

\begin{proposition}[First variance estimate for arbitrary phase]
  \label{P:static-phase-variance}
  Let $N = \|(f,g)\|_\lambda$. Uniformly in $\epsilon > 0$,
  \begin{align}
    \E\|Z_1(h)\|_\lambda^2
    &\le R_\alpha h^2\lambda^{\alpha-2}
    (1+\lambda^2h^2)^2N^2,\label{E:phase-born}\\
    \|\E Z_2(h)\|_\lambda
    &\le\frac{D_\Gamma(1+\lambda^2h^2)h^{3-\alpha}}
    {\lambda(2-\alpha)(3-\alpha)}N.
    \label{E:phase-mean}
  \end{align}
  In particular, put $\lambda = h^{-1}$ and let $K_\epsilon$ be the coefficient
  of $\sigma^2$ in the finite polynomial $\E\|Z_0(h)+\sigma
  Z_1(h)+\sigma^2Z_2(h)\|_\lambda^2$. Then
  \begin{equation}\label{E:phase-variance}
    |K_\epsilon|
    \le\left(4R_\alpha+
    \frac{4\sqrt e\,D_\Gamma}{(2-\alpha)(3-\alpha)}\right)
    h^{4-\alpha}\|(f,g)\|_{1/h}^2.
  \end{equation}
\end{proposition}

\begin{proof}
  We estimate the first random phase term and the second-order mean
  separately. We then retain both contributions to the variance coefficient
  and insert the common phase scale.

  \medskip
  \noindent\textit{Step 1.~}We estimate the first random phase term.
  We use the unnormalized spatial Fourier transform and Plancherel measure
  $(2\pi)^{-3}\ud\eta$. At fixed positive smoothing, Gaussian second and fourth
  moments make~\eqref{E:phase-duhamel} integrable: $G_t$ has finite positive
  mass, $C_t$ is an $L^2$ contraction, and $\E|W_\epsilon(x)W_\epsilon(y)|^2 \le
  3\Gamma_\epsilon(0)^2$. These facts, Minkowski and time integration justify
  the Gaussian and Fourier Fubini steps below first for Schwartz data.

  Set $a_s(\eta) = \widehat{w}_0(s,\eta)$ and
  \[
    e_\lambda(\eta)=|\widehat{f}(\eta)|^2+
    \frac{|\widehat{g}(\eta)|^2}{\lambda^2+|\eta|^2}.
  \]
  Cauchy--Schwarz on the two initial coordinates gives
  \[
    |a_s(\eta)|^2\le(1+\lambda^2s^2)e_\lambda(\eta),
  \]
  because the squared norm of their coefficient vector is
  $\cos^2(sq)+(1+\lambda^2/q^2)\sin^2(sq) \le 1+\lambda^2s^2$, with the
  continuous value at $q = 0$. For the two-component output multiplier put
  \[
    k_\tau(r)=\left(\frac{\sin(\tau r)}r,
    \frac{\cos(\tau r)}{\sqrt{\lambda^2+r^2}}\right).
  \]
  Its Euclidean squared norm is
  \[
    |k_\tau(r)|^2
    =\frac{1+\lambda^2\sin^2(\tau r)/r^2}{\lambda^2+r^2}
    \le\frac{1+\lambda^2h^2}{\lambda^2+r^2},\qquad 0\le\tau\le h.
  \]
  Expanding the covariance of the two copies of $W_\epsilon$, followed by
  spatial Plancherel, gives exactly
  \[
    \E\|Z_1(h)\|_\lambda^2
    =\frac1{(2\pi)^3}\int\ud\eta\int\mu_\epsilon(\ud\xi)
    \left|\int_0^h k_{h-s}(|\eta+\xi|)a_s(\eta)\ud s\right|^2.
  \]
  Cauchy--Schwarz in time bounds the last squared norm by
  \[
    \frac{h^2(1+\lambda^2h^2)^2}
    {\lambda^2+|\eta+\xi|^2}e_\lambda(\eta).
  \]
  Now $\mu_\epsilon \le \mu$ and Lemma~\ref{L:phase-contractions}
  give~\eqref{E:phase-born}.

  \medskip
  \noindent\textit{Step 2.~}We estimate the contracted second-order mean.
  The contraction inside the mean of $Z_2$ is the convolution operator with
  kernel $G_{s-r}(\ud x)\Gamma_\epsilon(x)$. Indeed, this follows directly from
  $\E[W_\epsilon(x)W_\epsilon(y)] = \Gamma_\epsilon(x-y)$ in the two Duhamel
  integrals. Its Fourier multiplier is $T_\epsilon(s-r,\eta)$, up to the
  immaterial even sign. Its $L^2$ operator norm is therefore at most
  $D_\Gamma(s-r)^{1-\alpha}$. Under the isometry
  $(u,v)\mapsto(u,A_\lambda^{-1/2}v)$, $\|k_\tau(|D|)\|_{L^2\to L^2\oplus L^2}
  \le \sqrt{1+\lambda^2h^2}/\lambda$ and $\|w_0(r)\|_2 \le
  \sqrt{1+\lambda^2h^2}N$. Minkowski and
  \[
    \int_0^h\ud s\int_0^s(s-r)^{1-\alpha}\ud r
    =\frac{h^{3-\alpha}}{(2-\alpha)(3-\alpha)}
  \]
  prove~\eqref{E:phase-mean}.

  \medskip
  \noindent\textit{Step 3.~}We combine both terms in the phase variance
  coefficient.
  Finally, the free phase obeys $\|Z_0(h)\|_\lambda \le e^{\lambda h/2}N$. To
  see this, for Schwartz free data differentiate its squared norm. The $-\Delta$
  terms cancel, leaving $2\lambda^2\operatorname{Re}\langle
  A_\lambda^{-1}\dot{w}_0,w_0\rangle$, whose absolute value is at most
  $\lambda\|Z_0\|_\lambda^2$. Gronwall gives the stated bound. The required
  coefficient is exactly
  \[
    K_\epsilon=\E\|Z_1(h)\|_\lambda^2+
    2\operatorname{Re}\langle Z_0(h),\E Z_2(h)\rangle_\lambda.
  \]
  Substituting $\lambda = h^{-1}$ into the two estimates proves
  \eqref{E:phase-variance}. This proves
  Proposition~\ref{P:static-phase-variance}.
\end{proof}

The second term in $K_\epsilon$ must be retained; neither its sign nor
monotonicity of the full phase energy in covariance is assumed. The estimates
control both terms at this one variance order. They do not control higher
Gaussian contractions, products of independent time blocks, overlapping noise
levels or the temporal regularization limit. Thus, they do not yet establish
Theorem~\ref{T:intro-upper}.

\section{All static orders with arbitrary initial phase}
\label{S:static-all-orders}

We now seek all static Duhamel orders, retaining the full spatial covariance
class of Section~\ref{S:static-phase}. The temporal covariance is not included
in this section. Static estimates alone do not identify the evolution for
overlapping intervals; that step is treated in the next section.

We use one external operator theorem. From Crouzeix--Palencia
\cite{crouzeix.palencia:source2017:numerical}*{author preprint v1, (3), p.~1,
and Theorem~6, pp.~5--6}, if $T$ is a bounded operator on a complex Hilbert
space and $F$ is holomorphic on a neighborhood of $\overline{W(T)}$, then
\begin{equation}\label{E:cp-input}
  \|F(T)\|\le(1+\sqrt2)\sup_{w\in W(T)}|F(w)|,\qquad
  W(T)=\{\langle Tu,u\rangle:\|u\|=1\}.
\end{equation}
The theorem is imported, not proved here. All operators to which we apply it
below are bounded; no unbounded-operator extension is assumed.

\LeanMarker{L:complex-phase}
\begin{lemma}[A complex phase bound]\label{L:complex-phase}
  Let $H$ be a complex Hilbert space, $A_0$ a bounded nonnegative selfadjoint
  operator, and $B$ a bounded selfadjoint operator. Write $A_\rho = \rho^2+A_0$.
  Suppose that, for $0 < a < 1$ and $d > 0$,
  \begin{equation}\label{E:abstract-sandwich}
    \|A_\rho^{-1/2}BA_\rho^{-1/2}\|\le d\rho^{a-2},
    \qquad \rho>0.
  \end{equation}
  Put
  \[
    I_a=\int_0^\infty s^{-1/2}(1+s)^{(a-3)/2}\ud s,\qquad
    \delta_a=\min\{1/4,\pi/(8I_a)\}.
  \]
  For $\lambda > 0$ and $|z| \le \delta_a\lambda^{2-a}/d$, let $L = A_0-zB$. The
  phase propagator of $u''+Lu = 0$ satisfies
  \begin{equation}\label{E:complex-phase-bound}
    \|\mathcal{M}_z(h)\|_{X_\lambda\to X_\lambda}
    \le100e^{2\lambda h},\qquad
    \|(f,g)\|_{X_\lambda}^2
    =\|f\|^2+\|A_\lambda^{-1/2}g\|^2,\quad h\ge0.
  \end{equation}
  The constant is independent of $H$, $\|A_0\|$ and $\|B\|$.
\end{lemma}

\begin{proof}
  We first compare the perturbed and reference square roots through
  resolvents. We then enclose the numerical range and bound the three scalar
  functions used in the phase matrix. The bounded-operator theorem applies to
  those functions, after which the ordered root factors give the phase
  estimate.

  \medskip
  \noindent\textit{Step 1.~}We obtain inverse and square-root comparisons.
  Set $\delta = \delta_a$ and $L_\lambda = L+\lambda^2$. For $\kappa =
  (\lambda^2+t)^{1/2}$,
  \[
    L_\lambda+t
    =A_\kappa^{1/2}(1-zK_\kappa)A_\kappa^{1/2},\qquad
    K_\kappa=A_\kappa^{-1/2}BA_\kappa^{-1/2}.
  \]
  Since $a-2 < 0$, $\|zK_\kappa\| \le \delta$. Thus, these operators are
  invertible and
  \[
    \|A_\kappa^{1/2}(L_\lambda+t)^{-1}\|
    \le\frac1{(1-\delta)\kappa}.
  \]
  Also $\operatorname{Re}\langle L_\lambda u,u\rangle \ge (1-\delta)\langle
  A_\lambda u,u\rangle$.

  Principal inverse square roots are defined on the right half-plane. The scalar
  identity $w^{-1/2} = \pi^{-1}\int_0^\infty t^{-1/2}(w+t)^{-1}\ud t$, and
  bounded holomorphic calculus, give
  \[
    Q_z-1
    =\frac1\pi\int_0^\infty t^{-1/2}A_\lambda^{1/2}
    \big[(L_\lambda+t)^{-1}-(A_\lambda+t)^{-1}\big]\ud t,
    \qquad Q_z=A_\lambda^{1/2}L_\lambda^{-1/2}.
  \]
  The integrals converge in operator norm. Indeed, the resolvent identity,
  $\|A_\lambda^{1/2}A_\kappa^{-1/2}\| \le 1$, and the following ordered
  factorization give the required integrand bound:
  \[
    \begin{split}
      &A_\lambda^{1/2}
      \big[(L_\lambda+t)^{-1}-A_\kappa^{-1}\big]\\
      &\qquad=A_\lambda^{1/2}A_\kappa^{-1/2}
      (1-zK_\kappa)^{-1}zK_\kappa A_\kappa^{-1/2}.
    \end{split}
  \]
  The four factors have norms at most $1$, $(1-\delta)^{-1}$, $|z|d\kappa^{a-2}$
  and $\kappa^{-1}$, respectively. Thus, the integrand without $t^{-1/2}/\pi$ is
  bounded by
  \[
    \frac{|z|d}{1-\delta}\kappa^{a-3}.
  \]
  Consequently,
  \[
    \|Q_z-1\|
    \le\frac{|z|d\lambda^{a-2}I_a}{\pi(1-\delta)}
    \le\frac16.
  \]
  In particular, $Q_z$ and $Q_z^{-1}$ have norm at most two. Apply the same
  argument to $L^* = A_0-\overline zB$ and take adjoints. It follows that all
  four operators
  \begin{equation}\label{E:phase-root-comparisons}
    A_\lambda^{1/2}L_\lambda^{-1/2},\quad
    L_\lambda^{-1/2}A_\lambda^{1/2},\quad
    L_\lambda^{1/2}A_\lambda^{-1/2},\quad
    A_\lambda^{-1/2}L_\lambda^{1/2}
  \end{equation}
  have norm at most two. This is a quantitative comparison at the bounded
  cutoffs, not an assumed square-root domain identity.

  \medskip
  \noindent\textit{Step 2.~}We enclose the numerical range and its square
  roots.
  We next control the numerical range. For $\|u\| = 1$ put $s = \|A_0^{1/2}u\|$.
  Optimizing~\eqref{E:abstract-sandwich} at $\rho = s$ gives $|\langle
  Bu,u\rangle| \le 2ds^a$. If $s = 0$, the same conclusion follows by letting
  $\rho\downarrow0$. Separately, evaluating~\eqref{E:abstract-sandwich} at $\rho
  = \lambda$ gives
  \[
    |\langle Bu,u\rangle|
    \le d\lambda^{a-2}(s^2+\lambda^2),\qquad
    |z\langle Bu,u\rangle|\le\delta(s^2+\lambda^2).
  \]
  This is the estimate for the real part; the optimized estimate controls the
  imaginary part. Thus, for $w = x+iy = \langle Lu,u\rangle$,
  \[
    x+\lambda^2\ge(1-\delta)(s^2+\lambda^2),\qquad
    |y|\le2\delta\lambda^{2-a}s^a.
  \]
  Since $s^{2a}\lambda^{2-2a} \le s^2+\lambda^2$,
  \[
    y^2\le\frac{4\delta^2}{1-\delta}\lambda^2(x+\lambda^2)
    \le\frac13\lambda^2(x+\lambda^2),\qquad
    x\ge-\lambda^2/4.
  \]
  Squaring the positive quantity $x+2\lambda^2$ reduces its comparison with
  $|w|$ to $y^2\le4\lambda^2(x+\lambda^2)$, while $|\operatorname{Im}\sqrt
  w|^2=(|w|-x)/2$.
  Hence, $|w| \le x+2\lambda^2$, and either square root of $w$ has
  $|\operatorname{Im}\sqrt w| \le \lambda$.

  \medskip
  \noindent\textit{Step 3.~}We apply the numerical-range theorem to
  holomorphic scalar functions.
  The scalar functions
  \[
    C_h(w)=\cos(h\sqrt w),\qquad
    F_h(w)=\sqrt{w+\lambda^2}\frac{\sin(h\sqrt w)}{\sqrt w},\qquad
    H_h(w)=\frac{w}{\sqrt{w+\lambda^2}}
    \frac{\sin(h\sqrt w)}{\sqrt w}
  \]
  are holomorphic on a neighborhood of $\overline{W(L)}$. The values at $w = 0$
  are removable; the only square-root branch required is that of $w+\lambda^2$,
  whose real part is positive. On the numerical range,
  \[
    |C_h(w)|\le e^{\lambda h},\qquad
    |F_h(w)|\le2(1+\lambda h)e^{\lambda h},\qquad
    |H_h(w)|\le2e^{\lambda h}.
  \]
  Indeed, $|\operatorname{Im}\sqrt w| \le \lambda$ gives $|\cos(h\sqrt
  w)|,|\sin(h\sqrt w)| \le e^{\lambda h}$. The entire identity $\sin(h\sqrt
  w)/\sqrt w = \int_0^h\cos(r\sqrt w)\ud r$ also gives the bound $he^{\lambda
  h}$, including $w = 0$. If $|w| \le \lambda^2/2$, then
  \[
    |F_h(w)|\le\sqrt{3/2}\,\lambda h e^{\lambda h}.
  \]
  If $|w| > \lambda^2/2$, then $|w+\lambda^2|/|w| \le 3$, so $|F_h(w)| \le
  \sqrt3 e^{\lambda h}$. These bounds imply the displayed estimate for $F_h$.
  Finally,
  \[
    \frac{|w|}{|w+\lambda^2|}
    \le1+\frac{\lambda^2}{|w+\lambda^2|}
    \le1+(1-\delta)^{-1}\le7/3,
  \]
  and $H_h(w) = \sqrt w\sin(h\sqrt w)/\sqrt{w+\lambda^2}$ gives its stated
  bound. All identities have their removable values at zero. Applying the
  imported~\eqref{E:cp-input} gives bounds $3e^{\lambda h}$, $6(1+\lambda
  h)e^{\lambda h}$ and $6e^{\lambda h}$ for these three operator functions.

  \medskip
  \noindent\textit{Step 4.~}We assemble the phase blocks in their stated
  noncommuting order.
  Under the isometry $(u,v)\mapsto(u,A_\lambda^{-1/2}v)$, the four phase blocks
  are
  \[
    C_h(L),\quad
    F_h(L)L_\lambda^{-1/2}A_\lambda^{1/2},\quad
    -A_\lambda^{-1/2}L_\lambda^{1/2}H_h(L),\quad
    A_\lambda^{-1/2}L_\lambda^{1/2}C_h(L)
    L_\lambda^{-1/2}A_\lambda^{1/2}.
  \]
  This is the usual matrix
  $\left(\begin{smallmatrix}C_h(L)&G_h(L)\\-LG_h(L)&C_h(L)\end{smallmatrix}\right)$,
  with $G_h(w) = \sin(h\sqrt w)/\sqrt w$, after weighting. In the order
  displayed, the four block norms are at most
  \[
    3e^{\lambda h},\qquad
    12(1+\lambda h)e^{\lambda h},\qquad
    12e^{\lambda h},\qquad12e^{\lambda h}.
  \]
  For example, the last block has two root-comparison factors of norm at most
  two, with $C_h(L)$ between them. Only functions of $L$ commute; no commutation
  with $A_\lambda$ is used in these ordered products. A block operator on
  $H\oplus H$ has norm at most the sum of its four block norms, by the triangle
  inequality and the norm-one coordinate inclusions and projections. Therefore,
  its norm is at most
  \[
    (39+12\lambda h)e^{\lambda h}
    \le39(1+\lambda h)e^{\lambda h}\le100e^{2\lambda h}.
  \]
  The phase operator itself is the exponential of the bounded block generator
  with rows $(0,1)$ and $(-L,0)$. It is entire in $z$; the intermediate square
  roots are used to estimate it, not to impose a branch-dependent definition of
  that entire operator. This proves Lemma~\ref{L:complex-phase}.
\end{proof}

We use the standard Gaussian creation and annihilation algebra; see
Janson~\cite{janson:97:gaussian}*{Section~13.2, Theorem~13.7, pp.~200--201}.
The following calculation records the normalization used for the finite
occupation spaces in this argument.

\LeanMarker{L:fock-form}
\begin{lemma}[Gaussian cutoff form bound]\label{L:fock-form}
  Retain the covariance and $R_\alpha$ of Section~\ref{S:static-phase}. Fix
  $\epsilon > 0$. Represent the smoothed real field on independent standard
  Gaussian coordinates as $W_\epsilon(x) = \sum_{i\ge1}g_i e_i(x)$ in
  $L^2(\Omega)$. Let $P_N$ be projection onto Gaussian Hermite degrees at most
  $N$, and $S_b$ spatial Fourier projection onto $|\eta| \le b$. For finite $m$,
  on their joint range set
  \[
    B_{N,m,b}=P_NS_b\Big(\sum_{i=1}^m g_i M_{e_i}\Big)S_bP_N,\qquad
    A_\rho=\rho^2-\Delta .
  \]
  These are bounded selfadjoint operators and, uniformly in $\epsilon,m,b$,
  \begin{equation}\label{E:fock-sandwich}
    \|A_\rho^{-1/2}B_{N,m,b}A_\rho^{-1/2}\|
    \le d_N\rho^{\alpha/2-2},\qquad
    d_N=2\sqrt{N+1}\sqrt{R_\alpha}.
  \end{equation}
\end{lemma}

\begin{proof}
  We first estimate the column of sandwiched spatial multiplication operators.
  We then compose it with the Gaussian creation row, retain the annihilation
  adjoint, and truncate by total degree.

  \medskip
  \noindent\textit{Step 1.~}We bound the spatial covariance column.
  Choose an orthonormal basis of the real isonormal space, which is separable,
  and write $e_i(x) = \E[W_\epsilon(x)g_i]$. Then $\sum_i e_i(x)^2 =
  \Gamma_\epsilon(0) = :v_\epsilon$, and $|e_i(x)| \le \sqrt{v_\epsilon}$. Thus,
  each multiplication $M_{e_i}$ is bounded and selfadjoint. Put $H_i =
  A_\rho^{-1/2}M_{e_i}A_\rho^{-1/2}$ on spatial $L^2$. Gaussian isometry and
  spatial Plancherel give
  \[
    \sum_i\|H_i f\|_2^2
    =\frac1{(2\pi)^3}\int
    \frac{|\widehat{f}(\eta)|^2}{\rho^2+|\eta|^2}
    J_{\rho,\epsilon}(\eta)\ud\eta,\qquad
    J_{\rho,\epsilon}(\eta)
    =\int\frac{e^{-\epsilon|\xi|^2}\mu(\ud\xi)}
    {\rho^2+|\eta+\xi|^2}.
  \]
  This is the covariance calculation for $A_\rho^{-1/2}W_\epsilon
  A_\rho^{-1/2}f$. By Lemma~\ref{L:phase-contractions}, $J_{\rho,\epsilon} \le
  R_\alpha\rho^{\alpha-2}$. Therefore, $\sum_i H_i^2 \le
  R_\alpha\rho^{\alpha-4}I$ in quadratic forms. Finite sums have the same upper
  bound. Compression to $S_b$ preserves it: for selfadjoint $H_i$,
  $(S_bH_iS_b)^2 \le S_bH_i^2S_b$ on its range. Both row and column variances
  are the same, since these operators are selfadjoint.

  The compression difference $S_bH_i^2S_b-(S_bH_iS_b)^2=S_bH_i(I-S_b)H_iS_b$
  is positive, so summing preserves the variance bound.
  \medskip
  \noindent\textit{Step 2.~}We use the Hermite creation and annihilation
  bounds.
  For completeness, use normalized Hermite polynomials, so multiplication by
  $g_i$ is $a_i^*+a_i$, with $g_i h_k =
  \sqrt{k_i+1}\,h_{k+e_i}+\sqrt{k_i}\,h_{k-e_i}$. On degree $k$, the row map
  $(v_i)\mapsto\sum_i a_i^*v_i$ has norm $\sqrt{k+1}$: its product with its
  adjoint is $\sum_i a_i^*a_i = (k+1)I$ on degree $k+1$. Composing with the
  column of spatial operators $S_bH_iS_b$ bounds the creation part by
  $\sqrt{k+1}\sqrt{R_\alpha}\rho^{(\alpha-4)/2}$. Distinct input degrees have
  orthogonal output degrees. The annihilation part is its adjoint, with the same
  bound after truncation. Adding them on degrees at most $N$ proves
  \eqref{E:fock-sandwich}. Boundedness of $-\Delta$ on $S_b$ is explicit; no
  infinite-frequency operator theorem was used. This proves
  Lemma~\ref{L:fock-form}.
\end{proof}

\subsection{Finite Gaussian products and their degree cutoffs}

Lemma~\ref{L:fock-form} bounds the field only on Hermite degrees at most $N$.
We now check that this cutoff does not change the finite words that we need.

\LeanMarker{L:finite-word-cutoff}
\begin{lemma}[Exact cutoffs for finite Gaussian words]\label{L:finite-word-cutoff}
  Let $g_{ij}$, $1 \le i \le M$, $1 \le j \le m_i$, be independent standard real
  Gaussian coordinates, grouped by color $i$, and let $\mathscr{X}$ be a complex
  Hilbert space. On Gaussian polynomials with coefficients in $\mathscr{X}$ set
  \[
    J_i=\sum_{j=1}^{m_i}M_{g_{ij}}\otimes K_{ij},
  \]
  where each $K_{ij}$ is bounded on $\mathscr{X}$. Let $F_0,\ldots,F_n$ be
  bounded operators acting only on $\mathscr{X}$. Fix a word $i_1,\ldots,i_n$,
  and let $r_i$ be the number of its occurrences of color $i$. If $N_i \ge r_i$
  and $P$ projects each color onto total Hermite degree at most $N_i$, then, for
  every $v \in \mathscr{X}$,
  \begin{equation}\label{E:finite-word-cutoff}
    \begin{split}
      &(I\otimes F_n)(PJ_{i_n}P)(I\otimes F_{n-1})\cdots
      (PJ_{i_1}P)(I\otimes F_0)(1\otimes v)\\
      &\qquad=(I\otimes F_n)J_{i_n}(I\otimes F_{n-1})\cdots
      J_{i_1}(I\otimes F_0)(1\otimes v).
    \end{split}
  \end{equation}
  The equality holds for finite sums and for Bochner integrals of such words
  whenever the displayed vector integrals exist.
\end{lemma}

\begin{proof}
  In normalized Hermite coordinates, multiplication by $g_{ij}$ maps
  $h_{\boldsymbol a}$ to $\sqrt{a_{ij}+1}\,h_{\boldsymbol a+e_{ij}}
  +\sqrt{a_{ij}}\,h_{\boldsymbol a-e_{ij}}$; the second term is zero if $a_{ij}
  = 0$. Thus, it changes the total degree of color $i$ by one and leaves all
  other color degrees unchanged. The operators $I\otimes F_l$ leave every
  Gaussian degree unchanged. After the first $l$ insertions, the uncompressed
  word therefore has color-$i$ degree at most $r_i(l) = |\{j\le l:i_j=i\}| \le
  r_i \le N_i$. This is an induction starting from degree zero in every color.
  Consequently, $P$ is the identity both before and after every insertion along
  that word. Induction on $l$ now proves \eqref{E:finite-word-cutoff}. The
  lowering term has been retained: the assertion includes creation followed by
  annihilation and all resulting self-contractions. Linearity and equality of
  the vector integrands prove the last assertion. This is an equality of
  selected finite words, not of full propagators at a fixed cutoff. This proves
  Lemma~\ref{L:finite-word-cutoff}.
\end{proof}

\begin{figure}
  \centering
  \begin{tikzpicture}[x=1.7cm,y=.62cm]
    \fill[hamblue!4] (0,0)--(4,4)--(4,0)--cycle;
    \draw[ham/arrow] (-.15,0)--(4.5,0) node[right] {insertions};
    \draw[ham/arrow] (0,-.25)--(0,4.85) node[above] {degree};
    \draw[dashed,hamgold,thick] (0,4)--(4.25,4);
    \node[anchor=west,font=\small] at (4.3,4) {$N=4$};
    \draw[dotted,thick] (0,0)--(1,1)--(2,2)--(3,3)--(4,4);
    \node[font=\small,fill=white] at (2.6,3.25) {$r(\ell)=\ell$};
    \draw[hamblue,line width=1.2pt]
    (0,0)--(1,1)--(2,2)--(3,1)--(4,2);
    \foreach \hamPointX/\hamPointY in {0/0,1/1,2/2,3/1,4/2} {
    \filldraw[draw=hamblue,fill=white] (\hamPointX,\hamPointY) circle (2pt);
    }
    \foreach \hamPointX in {1,2,3,4} {\node[font=\small,anchor=north] at (\hamPointX,-.18) {$\hamPointX$};}
    \node[ham/note,text width=8cm] at (2,-2.25)
    {Every prefix stays below its occurrence count.\\
    Projection changes none of these selected finite words.};
  \end{tikzpicture}
  \caption[Exact cutoffs for selected Gaussian words]{One possible degree path
    for four multiplications in a single Gaussian color. Multiplication is a sum
    of raising and lowering terms; the solid line illustrates one term, not the
    entire polynomial. For each color, every prefix has degree at most its
    occurrence count, which is at most the chosen cutoff. The result concerns
    selected finite words acting on the degree-zero vector, not entire
    propagators at a fixed cutoff. See \ProofRef{Lemma}{L:finite-word-cutoff}.}
  \label{Fig:hermite-cutoff}
\end{figure}
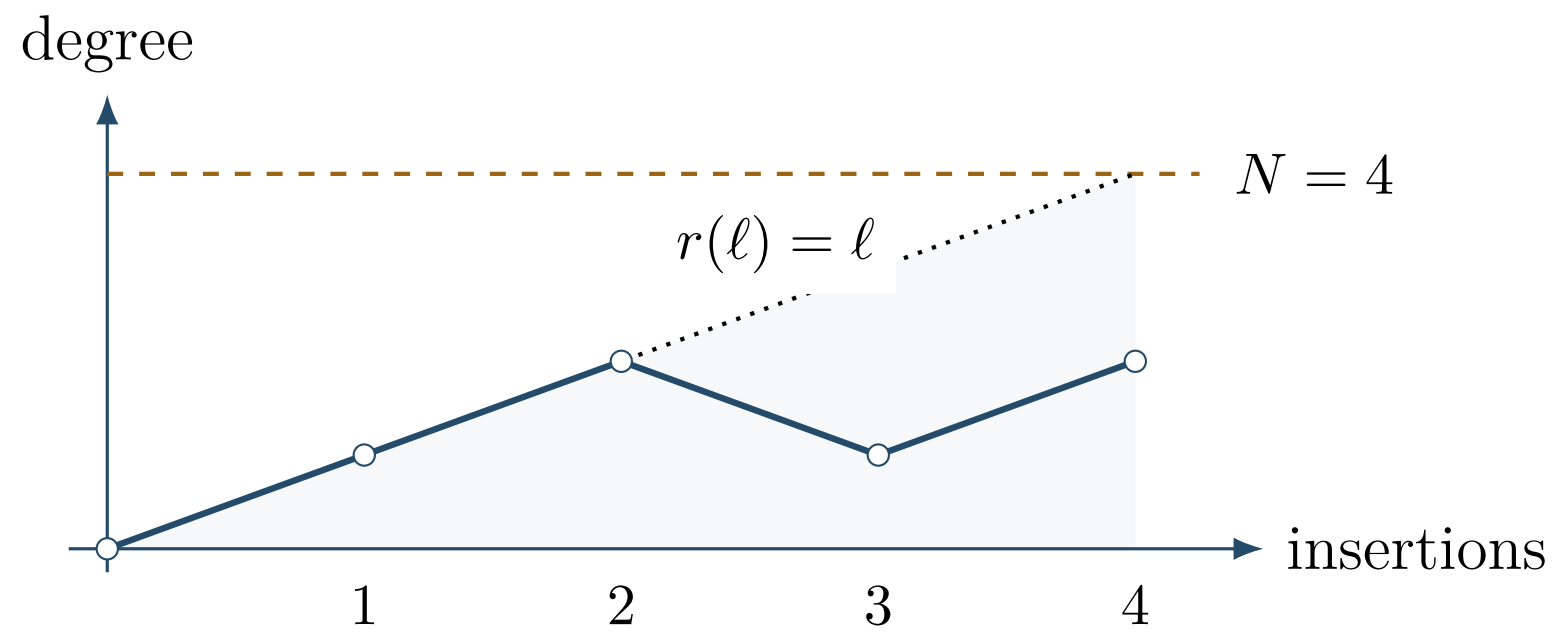

\subsection{Uniform bounds for all orders}

Proposition~\ref{P:static-phase-variance} controls only the first variance
coefficient. We now combine Lemmas~\ref{L:complex-phase}, \ref{L:fock-form}
and~\ref{L:finite-word-cutoff} to bound every static Duhamel order.

\LeanMarker{P:static-phase-all-orders}
\begin{proposition}[All regularized static phase orders]
  \label{P:static-phase-all-orders}
  Retain the full covariance class of Section~\ref{S:static-phase}. For real
  Schwartz $f,g$, let $Z_n^\epsilon(h)$ be the $n$th regularized static Duhamel
  phase, with the recursion \eqref{E:phase-duhamel} for every $n \ge 0$. This is
  an ordinary Gaussian polynomial of degree at most $n$, including
  self-contractions; it is not a pure Wick chaos. There is an explicit finite
  constant
  \[
    A_\alpha=\frac{400\sqrt2\,e^2}{\delta_{\alpha/2}}
  \]
  such that, for $n \ge 1$, uniformly in $\epsilon > 0$,
  \begin{equation}\label{E:all-static-phase-orders}
    \|Z_n^\epsilon(h)\|_{L^2(\Omega;X_{1/h})}
    \le
    \frac{(A_\alpha\sqrt{R_\alpha})^n
    h^{(4-\alpha)n/2}}{(n!)^{(3-\alpha)/2}}
    \|(f,g)\|_{1/h}.
  \end{equation}
  At every fixed positive smoothing and real coupling $\sigma$,
  $\sum_{n\ge0}\sigma^nZ_n^\epsilon(h)$ converges in $L^2$ phase. Its squared
  norm is at most
  \begin{equation}\label{E:all-static-phase-sum}
    (\sqrt e+2)^2
    \exp\left\{6A_\alpha^{2/(3-\alpha)}
    (|\sigma|^2R_\alpha)^{1/(3-\alpha)}
    h^{(4-\alpha)/(3-\alpha)}\right\}\|(f,g)\|_{1/h}^2.
  \end{equation}
\end{proposition}

\begin{proof}
  We estimate the coupling coefficients at finite cutoffs, choose the phase
  scale by order, and then remove spatial and Gaussian cutoffs at fixed
  smoothing. Finally, we sum the resulting factorial estimates.

  \medskip
  \noindent\textit{Step 1.~}We apply the complex phase bound at bounded
  cutoffs.
  If $R_\alpha = 0$, spectral homogeneity gives $\mu = 0$, and every
  positive-order term is zero. Assume $R_\alpha > 0$. First work at finite
  Gaussian cutoff $m$, degree cutoff $N$ and spatial band $b$, with initial
  phase $S_b(f,g)$ times the Gaussian vacuum $1$. All relevant operators are
  bounded. Lemmas~\ref{L:fock-form} and~\ref{L:complex-phase}, with $a =
  \alpha/2$, give a bound $100e^{2\lambda h}$ on the complex coupling circle of
  radius $\delta_{\alpha/2}\lambda^{2-\alpha/2}/d_N$.

  The four phase entries are entire in the coupling, by the power series for
  $\cos(h\sqrt L)$ and $\sin(h\sqrt L)/\sqrt L$ of a bounded $L$. Cauchy's
  formula for their operator-valued series bounds the $n$th coefficient by
  \[
    100e^{2\lambda h}
    \left(\frac{d_N}
    {\delta_{\alpha/2}\lambda^{2-\alpha/2}}\right)^n
    \|S_b(f,g)\|_\lambda .
  \]
  The coefficient is exactly the $n$th Duhamel phase from the vacuum. Apply
  Lemma~\ref{L:finite-word-cutoff} with one color, $N = n$, the spatial phase
  Hilbert space as $\mathscr{X}$, and its free wave operators as the $F_l$. It
  identifies each cutoff Duhamel word exactly with its uncompressed word,
  including self-contractions. Dominated time integration preserves this
  identity.

  \medskip
  \noindent\textit{Step 2.~}We optimize the phase scale and convert powers to
  factorials.
  Choose $\lambda = n/h$. For any phase $z$, $\|z\|_{1/h} \le n\|z\|_{n/h}$ and
  $\|(f,g)\|_{n/h} \le \|(f,g)\|_{1/h}$. Also $d_n \le
  2\sqrt{2n}\sqrt{R_\alpha}$. The preceding bound is therefore at most
  \[
    100n\left(\frac{2\sqrt2 e^2}{\delta_{\alpha/2}}\right)^n
    R_\alpha^{n/2}h^{(4-\alpha)n/2}
    n^{-(3-\alpha)n/2}\|(f,g)\|_{1/h}.
  \]
  Since $100n \le 200^n$ and $n! \le n^n$, this gives
  \eqref{E:all-static-phase-orders} at the cutoffs.

  \medskip
  \noindent\textit{Step 3.~}We remove the cutoffs for each fixed word.
  We justify removal of those cutoffs at fixed $\epsilon$ and $n$. On Gaussian
  degrees at most $N$, the unweighted field multiplication has norm at most
  $2\sqrt{N+1}\sqrt{v_\epsilon}$, uniformly in $m$, by the same
  creation/annihilation argument applied pointwise in $x$. It follows that the
  noise phase insertion $(u,v)\mapsto(0,W_\epsilon u)$, between the required
  finite degree levels, is bounded on $X_\lambda$ by this constant divided by
  $\lambda$. Free phase evolution is bounded by $e^{\lambda h/2}$. For finite
  $m$, $S_b \to I$ strongly and commutes with the free propagator. Each fixed
  Duhamel word therefore converges strongly as $b \to \infty$, by dominated time
  integration.

  Next let $m \to \infty$ on the full Gaussian space. The truncated field
  multiplications converge strongly between finite degree levels. To see this,
  first take a vector depending on finitely many Gaussian coordinates. The
  annihilation tail is then zero, while the creation tail is bounded by a
  constant times
  \[
    \left(\int\sum_{i>m}e_i(x)^2\|F(x)\|^2\ud x\right)^{1/2},
  \]
  which tends to zero by dominated convergence, since the sum is bounded by
  $v_\epsilon$. Such vectors are dense; the uniform multiplication bound extends
  strong convergence to all vectors in those degree levels. Dominated time
  integration again gives convergence of every fixed word. The limit is
  precisely $Z_n^\epsilon$, so the estimate survives both limits. We assert no
  zero-smoothing phase limit in this argument.

  \medskip
  \noindent\textit{Step 4.~}We sum all regularized static orders.
  Finally, put $\beta = (3-\alpha)/2 \in (1/2,3/2)$ and $x =
  A_\alpha|\sigma|\sqrt{R_\alpha}h^{(4-\alpha)/2}$. The elementary series
  satisfies
  \[
    \sum_{n\ge0}\frac{x^n}{(n!)^\beta}
    \le2e^{3x^{1/\beta}}.
  \]
  For $\beta \ge 1$, this follows from $\sum a_n^\beta \le (\sum a_n)^\beta$
  with $a_n = x^{n/\beta}/n!$. For $\beta < 1$, apply H\"older to
  $((2x^{1/\beta})^n/n!)^\beta (2^{-\beta n/(1-\beta)})^{1-\beta}$. The second
  geometric sum is at most two. Since $\|Z_0(h)\|_{1/h} \le \sqrt
  e\|(f,g)\|_{1/h}$, absolute $L^2$ summation and squaring give
  \eqref{E:all-static-phase-sum}. This proves
  Proposition~\ref{P:static-phase-all-orders}.
\end{proof}

\begin{figure}[!tb]
  \centering
  % Correlation pattern only: the cited product estimate uses static block fields.
\begin{tikzpicture}[x=1cm,y=1cm,font=\normalsize]
  \foreach \panel/\shift in {0/0,1/7.1} {
  \begin{scope}[xshift=\shift cm]
    \foreach \i in {0,...,23} {
    \foreach \j in {0,...,23} {
    \pgfmathtruncatemacro{\same}{int(\i/6)==int(\j/6)}
    \pgfmathsetmacro{\shade}{15+85*(1/(abs(\i-\j)+.5)^.5-.2063)/(1.4143-.2063)}
    \ifnum\panel=0
    \fill[hamblue!\shade!white] (\i*.2,\j*.2) rectangle ++(.201,.201);
    \else
    \ifnum\same=1
    \fill[hamblue!\shade!white] (\i*.2,\j*.2) rectangle ++(.201,.201);
    \fi
    \fi
    }
    }
    \draw[hamblue] (0,0) rectangle (4.8,4.8);
    \foreach \k in {0,...,3} {
    \draw[hamgold,dashed,thick] (\k*1.2,\k*1.2) rectangle ++(1.2,1.2);
    }
    \draw[ham/arrow] (0,-.12)--(5,-.12) node[right] {$s$};
    \draw[ham/arrow] (-.12,0)--(-.12,5) node[above] {$r$};
    \node[anchor=north] at (0,-.22) {$0$};
    \node[anchor=north] at (4.8,-.22) {$T$};
    \node[anchor=east] at (-.22,4.8) {$T$};
    \node at (2.4,-.85) {\ifnum\panel=0 Cross-block correlation\else Cross-block covariance: zero\fi};
  \end{scope}
  }
  \node[font=\normalsize\bfseries] at (2.4,6.2) {Original time kernel};
  \node at (2.4,5.65) {$|s-r|^{-\alpha_0}$};
  \node[font=\normalsize\bfseries] at (9.5,6.2) {Independent-block pattern};
  \node at (9.5,5.65) {Only diagonal blocks remain};
  \node at (7.7,3) {$0$};
  \node at (10.1,.6) {$0$};
  \draw[ham/arrow] (5.15,2.4)--(6.65,2.4);
  \node[align=center,text width=2cm] at (5.9,3.25) {separate\\ fields};
  \node at (5.95,-1.5)
  {$\alpha_0=1/2$ shown; shading is capped near the diagonal.};
\end{tikzpicture}
  \caption[The noise remembers across time blocks]{Left: the time covariance
    $|s-r|^{-\alpha_0}$ of the noise on $[0,T]^2$, shown for $\alpha_0 = 1/2$ with
    the shading capped near the diagonal; the dashed squares are the time blocks.
    Right: which pairs of times stay correlated when each block carries its own
    independent field, as in Corollary~\ref{C:static-interaction}: only pairs in
    the same block. The shading inside the blocks repeats the original kernel for
    comparison; in the corollary each block's field is static, so its covariance
    is constant within a block. The product bound~\eqref{E:independent-static-blocks}
    holds for such independent blocks; it is not a bound for the original noise,
    whose correlation across blocks it drops.}
  \label{Fig:noise-memory}
\end{figure}
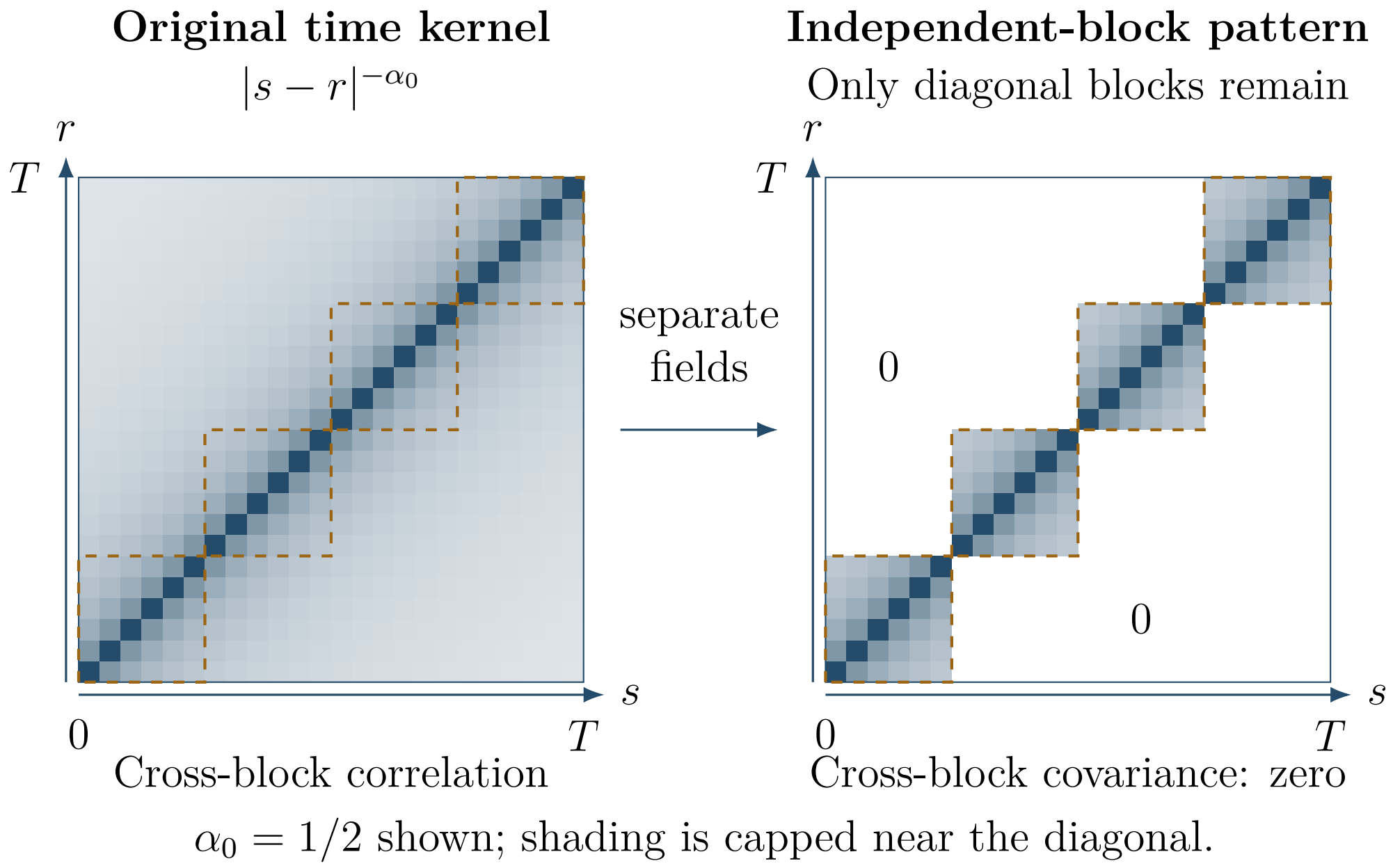

\LeanMarker{C:static-interaction}
\begin{corollary}[Summable static interaction bounds]
  \label{C:static-interaction}
  Let $\mathcal{S}_0(t)$ be free phase evolution and, for $\lambda h \le 1$,
  define the static interaction map by its convergent series
  \[
    \mathcal{I}^\epsilon_\sigma(h)(f,g)
    =\mathcal{S}_0(-h)\sum_{n\ge0}\sigma^n Z_n^\epsilon(h).
  \]
  Put $\gamma = 3-\alpha$, $B_\alpha = 3eA_\alpha^2$ and
  \[
    \begin{split}
      L_\alpha(h,\lambda)
      &=4eR_\alpha h^2\lambda^{\alpha-2}
      +\frac{4\sqrt eD_\Gamma}{(2-\alpha)(3-\alpha)}
      \frac{h^{3-\alpha}}{\lambda},\\
      x&=B_\alpha\sigma^2R_\alpha h^{4-\alpha},\\
      \mathcal{E}(h,\sigma,\lambda)
      &=\sigma^2L_\alpha(h,\lambda)
      +(\lambda h)^{-2}x^2e^{\gamma x^{1/\gamma}}.
    \end{split}
  \]
  Uniformly in $\epsilon > 0$,
  \begin{equation}\label{E:static-interaction}
    \E\|\mathcal{I}^\epsilon_\sigma(h)(f,g)\|_\lambda^2
    \le[1+\mathcal{E}(h,\sigma,\lambda)]\|(f,g)\|_\lambda^2.
  \end{equation}
  If $0 < \alpha_0 < 1$ and $\alpha+\alpha_0 < 2$, then, for $h_j =
  2^{-j}/\lambda$, $\sigma_j^2 = h_j^{-\alpha_0}$ and $T > 0$,
  \begin{equation}\label{E:static-block-series}
    \sum_{j\ge0}\frac{T}{h_j}\mathcal{E}(h_j,\sigma_j,\lambda)<\infty.
  \end{equation}
  For a finite single partition, let $\mathcal{M}_j(h_j)$ be the static phase
  maps with mutually independent copies of the smoothed field, real couplings
  $\sigma_j$, and $\lambda h_j \le 1$. They have continuous $L^2$ extensions to
  independent random initial phases. Then
  \begin{equation}\label{E:independent-static-blocks}
    \E\|\mathcal{M}_K(h_K)\cdots\mathcal{M}_1(h_1)(f,g)\|_\lambda^2
    \le
    \exp\left\{\lambda\sum_{j=1}^K h_j+
    \sum_{j=1}^K\mathcal{E}(h_j,\sigma_j,\lambda)\right\}
    \|(f,g)\|_\lambda^2.
  \end{equation}
  This concerns independent successive blocks on one partition.
  Formula~\eqref{E:static-block-series} only sums individual upper bounds across
  dyadic levels; it does not identify a product with an evolution driven by
  overlapping or repeated fields.
\end{corollary}

\begin{proof}
  We keep the complete first variance coefficient separate from the
  higher-order series. We then check dyadic summability and use independence
  only for successive blocks on one finite partition.

  \medskip
  \noindent\textit{Step 1.~}We retain both first-variance contributions.
  Write $N = \|(f,g)\|_\lambda$ and $Y_n = \mathcal{S}_0(-h)Z_n^\epsilon(h)$.
  The free estimate from Section~\ref{S:static-phase} holds for both time
  directions, since its differentiated energy was bounded in absolute value.
  Thus, $\|\mathcal{S}_0(-h)\|_\lambda \le \sqrt e$. Since $Y_0 = (f,g)$,
  Proposition~\ref{P:static-phase-variance} bounds the complete first variance
  coefficient of $\E\|\sum\sigma^nY_n\|_\lambda^2$ in absolute value by
  $L_\alpha(h,\lambda)N^2$. In detail, its first random squared norm is bounded
  by $e$ times~\eqref{E:phase-born}, and its second-order mean cross term by
  $2N\sqrt e$ times~\eqref{E:phase-mean}.

  \medskip
  \noindent\textit{Step 2.~}We bound the higher-order interaction series.
  For $n \ge 1$, the norm comparisons when $\lambda h \le 1$
  and~\eqref{E:all-static-phase-orders} give
  \[
    \|Y_n\|_{L^2(\Omega;X_\lambda)}
    \le\frac{\sqrt e}{\lambda h}
    \frac{(A_\alpha\sqrt{R_\alpha})^n h^{(4-\alpha)n/2}}
    {(n!)^{\gamma/2}}N.
  \]
  Gaussian symmetry makes $\E\langle Y_r,Y_s\rangle_\lambda = 0$ if $r+s$ is
  odd. All double series below converge absolutely by these bounds. For $r+s =
  2k$, the elementary factorial inequality $r!s! \ge (k!)^2$ shows that the sum
  of absolute inner products is at most
  \[
    \frac{e(2k+1)}{(\lambda h)^2}
    \frac{(A_\alpha^2R_\alpha)^k h^{(4-\alpha)k}}
    {(k!)^\gamma}N^2,\qquad k\ge1.
  \]
  This includes $r = 0$ or $s = 0$, since $\lambda h \le 1$ and $\|Y_0\| = N$.
  The bound $e(2k+1) \le (3e)^k$ therefore controls all orders $k \ge 2$ by
  \[
    \frac1{(\lambda h)^2}
    \sum_{k\ge2}\frac{x^k}{(k!)^\gamma}N^2
    \le\frac{x^2e^{\gamma x^{1/\gamma}}}{(\lambda h)^2}N^2.
  \]
  Here $\gamma > 1$, and $\sum_{j\ge0}x^j/(j!)^\gamma \le e^{\gamma
  x^{1/\gamma}}$. Keeping the sharper first coefficient separately proves
  \eqref{E:static-interaction}.

  \medskip
  \noindent\textit{Step 3.~}We sum the dyadic budgets and estimate the
  independent-block product.
  Finally, set $\beta = 4-\alpha-\alpha_0 > 2$. After multiplication by $T/h_j$,
  the two first-coefficient terms in~\eqref{E:static-block-series} are constant
  multiples of
  \[
    T\lambda^{\alpha-2}h_j^{1-\alpha_0},\qquad
    \frac T\lambda h_j^{2-\alpha-\alpha_0}.
  \]
  The higher-order term is
  \[
    \frac{T(B_\alpha R_\alpha)^2}{\lambda^2}
    h_j^{2\beta-3}
    \exp\{\gamma(B_\alpha R_\alpha)^{1/\gamma}h_j^{\beta/\gamma}\}.
  \]
  The exponential factor is bounded by its value at $h_0$. All three powers are
  positive, so their geometric series converge. For the finite partition, write
  $\mathcal{M}_j(h_j) = \mathcal{S}_0(h_j)\mathcal{I}_j(h_j)$.
  Bound~\eqref{E:static-interaction} and $\|\mathcal{S}_0(h_j)\|_\lambda^2 \le
  e^{\lambda h_j}$ give a bounded linear map from deterministic $X_\lambda$ to
  $L^2(\Omega_j;X_\lambda)$, initially on Schwartz phases and then by density.
  Tensoring with the past probability space extends this bound to random initial
  phases independent of the $j$th field. Iteration gives the product of
  $e^{\lambda h_j}(1+\mathcal{E}(h_j,\sigma_j,\lambda))$. The inequality $1+y
  \le e^y$ for $y \ge 0$ proves \eqref{E:independent-static-blocks}. This proves
  Corollary~\ref{C:static-interaction}.
\end{proof}

The independent-block estimate does not use log-convexity in additive
covariance amplitudes. That unrestricted interpolation principle already
fails for a scalar Gaussian wave; see No-go~\ref{P:nogo-11}.

These estimates concern individual static phase and interaction maps. Their
summability alone does not identify the evolution driven by overlapping interval
fields. The coefficient argument in Section~\ref{S:interval-upper} supplies the
separate estimate needed for that evolution.

\section{Interval fields and the full upper bound}
\label{S:interval-upper}

The intrinsic constant is identified independently in
Theorem~\ref{T:intrinsic-second-moment}.

Write
\[
  \nu=3-\alpha>1,\qquad
  \beta=4-\alpha-\alpha_0>2,\qquad p=2-\alpha_0>1.
\]
Thus, $\beta-\nu = 1-\alpha_0 > 0$. In this section the spatial dimension is any
admissible $d \in \{1,2,3\}$. The spatial covariance is the nonnegative, even,
positive-definite tempered Radon measure $\Gamma$, with scaling $\Gamma(rA) =
r^{d-\alpha}\Gamma(A)$ and spectral convention $\Gamma = \int e^{i\xi\cdot
x}\mu(\ud\xi)$. This includes $\Gamma = \delta_0$ when $d = \alpha = 1$. Set
\[
  R_\alpha=\frac{\alpha\pi\mu(B_1)}{2\sin(\pi\alpha/2)},\qquad
  A_\rho=\rho^2-\Delta,\qquad
  \|(f,g)\|_\rho^2=\|f\|_2^2+\|A_\rho^{-1/2}g\|_2^2.
\]
Complex Hilbert inner products are linear in the first slot.

\subsection{Contractions and complex phase estimates}

The resolvent estimate in Lemma~\ref{L:phase-contractions} and the Gaussian
cutoff bound in Lemma~\ref{L:fock-form} were stated for a covariance on
$\R^3$. We now derive the needed column and finite-degree bounds for every
admissible $d$.

\LeanMarker{L:interval-contractions}
\begin{lemma}[Dimension-independent Gaussian contractions]
  \label{L:interval-contractions}
  Uniformly in positive spatial smoothing and in finite Gaussian-mode and
  spatial Fourier-band cutoffs,
  \begin{equation}\label{E:interval-column}
    \sum_m\|A_\rho^{-1/2}M_{e_m}f\|_2^2
    \le R_\alpha\rho^{\alpha-2}\|f\|_2^2.
  \end{equation}
  Here $W_\epsilon(x) = \sum_m g_me_m(x)$ is an orthonormal-coordinate
  representation of $W(p_\epsilon(x-\cdot))$, whose covariance is
  $\Gamma_\epsilon = p_{2\epsilon}*\Gamma$. Consequently, on a Gaussian color
  that is initially in its vacuum, the phase insertion
  $(u,v)\mapsto(0,W_\epsilon u)$ has norm at most
  $\sqrt{R_\alpha}\rho^{(\alpha-2)/2}$ into its Gaussian $L^2$ phase. With
  Hermite degree of this color truncated at $N$, its sandwiched field has
  operator norm at most
  \begin{equation}\label{E:interval-fock}
    2\sqrt{N+1}\sqrt{R_\alpha}\rho^{\alpha/2-2}.
  \end{equation}
  Both statements remain valid tensoring with any other Hilbert factor.
\end{lemma}

\begin{proof}
  We derive the spatial covariance column estimate in every admissible
  dimension. The vacuum insertion uses that column directly, while occupied
  Hermite degrees require the creation and annihilation estimate on the
  sandwiched operators.

  \medskip
  \noindent\textit{Step 1.~}We compare translated resolvents using physical
  covariance positivity.
  Fourier uniqueness gives $\mu(B_r) = r^\alpha\mu(B_1)$. The Gaussian
  comparison in Lemma~\ref{L:phase-contractions} uses no dimension-specific
  identity:
  \[
    \int e^{-s|\eta+\xi|^2}\mu(\ud\xi)
    =\int p_{2s}^{(d)}(x)e^{i\eta\cdot x}\Gamma(\ud x)
    \le\int e^{-s|\xi|^2}\mu(\ud\xi).
  \]
  Evenness gives a real integral, and physical nonnegativity gives the
  inequality. Integrating against $e^{-\rho^2s}\ud s$ and using radial spectral
  homogeneity yields
  \[
    \sup_\eta\int\frac{\mu(\ud\xi)}{\rho^2+|\eta+\xi|^2}
    \le R_\alpha\rho^{\alpha-2}.
  \]
  No exclusion of a physical origin atom is used.

  \medskip
  \noindent\textit{Step 2.~}We obtain the column bound and the new-color
  insertion estimate.
  Gaussian isometry and unnormalized spatial Plancherel give
  \[
    \sum_m\|A_\rho^{-1/2}M_{e_m}f\|_2^2
    =\frac1{(2\pi)^d}\int|\widehat{f}(\eta)|^2
    \left(\int\frac{e^{-\epsilon|\xi|^2}\mu(\ud\xi)}
    {\rho^2+|\eta+\xi|^2}\right)\ud\eta .
  \]
  This proves~\eqref{E:interval-column}. A finite mode sum is a positive
  sub-sum. For a Fourier-band projection $S_b$, which commutes with $A_\rho$,
  use $S_bA_\rho^{-1}S_b \le A_\rho^{-1}$ to obtain the same inequality for
  $S_bM_{e_m}S_b$. The phase assertion on a new color follows because
  $\sum_mg_m(0,M_{e_m}u)$ has squared norm given by the left side
  of~\eqref{E:interval-column}, and its input is independent of this color.

  \medskip
  \noindent\textit{Step 3.~}We obtain the finite-degree bound and retain other
  Hilbert factors.
  Substituting $A_\rho^{-1/2}f$ gives
  \[
    \sum_m(A_\rho^{-1/2}M_{e_m}A_\rho^{-1/2})^2
    \le R_\alpha\rho^{\alpha-4}I.
  \]
  These spatial sandwiches are selfadjoint, so their row and column variances
  agree. The normalized Hermite creation-row norm from degree $j$ to degree
  $j+1$ is $\sqrt{j+1}$, and annihilation is its adjoint, as proved in
  Lemma~\ref{L:fock-form}. Composition with the spatial column and addition of
  the two parts prove~\eqref{E:interval-fock}. The same positive operator
  identities hold with an identity on every other Hilbert factor. The
  calculation used only the radial mass $\mu(B_1)$, homogeneity and physical
  nonnegativity. Thus, $R_\alpha$ is the only covariance quantity in these
  estimates. No angular bound or unsmoothed multiplication norm has been
  absorbed into a constant. This proves Lemma~\ref{L:interval-contractions}.
\end{proof}

\LeanMarker{L:interval-complex}
\begin{lemma}[Complex sums and changing phase weights]
  \label{L:interval-complex}
  Let $A_0$ be bounded, nonnegative and selfadjoint, and let $V$ be bounded on
  the same complex Hilbert space. Write $A_\rho = \rho^2+A_0$. The complex phase
  conclusion of Lemma~\ref{L:complex-phase} holds for $L = A_0-V$ whenever $V$
  is bounded, not necessarily selfadjoint,
  \[
    \|A_\rho^{-1/2}VA_\rho^{-1/2}\|\le D\rho^{a-2},\qquad \rho>0,
    \quad 0<a<1,\quad D\le\delta_a\Lambda^{2-a}.
  \]
  Its bound is $100e^{2\Lambda h}$ in $X_\Lambda$. If a piecewise-constant
  potential is the sum of connected-interval fields, assign a positive weight
  $\omega_i$ to each field present and put $\Lambda(t) =
  \lambda+\sum_i\omega_i\mathbf1_{I_i}(t)$. The total loss from changing between
  these phase norms, starting and ending in $X_\lambda$, is at most
  \begin{equation}\label{E:interval-metric}
    \prod_i(1+\omega_i/\lambda).
  \end{equation}
\end{lemma}

\begin{proof}
  We repeat the resolvent and numerical-range estimates for a general bounded
  complex potential, retaining adjoint comparisons. We then count the
  phase-norm losses at jumps of the active interval set.

  \medskip
  \noindent\textit{Step 1.~}We obtain the four ordered root comparisons
  without selfadjointness.
  The adjoint has the identical sandwich bound, since
  \[
    A_\rho^{-1/2}V^*A_\rho^{-1/2}
    =(A_\rho^{-1/2}VA_\rho^{-1/2})^*.
  \]
  The form estimate at $\rho = \Lambda$ gives
  $\operatorname{Re}\langle(L+\Lambda^2)u,u\rangle \ge (1-\delta_a)\langle
  A_\Lambda u,u\rangle$. The shifted operator and its adjoint thus have spectra
  in the right half-plane, where their principal inverse roots are defined. Put
  $K_\kappa = A_\kappa^{-1/2}VA_\kappa^{-1/2}$ and $\kappa =
  (\Lambda^2+t)^{1/2}$. The resolvent difference is exactly
  \[
    (L+\Lambda^2+t)^{-1}-A_\kappa^{-1}
    =A_\kappa^{-1/2}(1-K_\kappa)^{-1}K_\kappa A_\kappa^{-1/2}.
  \]
  Its product on the left with $A_\Lambda^{1/2}$ has norm at most
  $D\kappa^{a-3}/(1-\delta_a)$. The inverse-root integral in
  Lemma~\ref{L:complex-phase} therefore gives the same $1/6$ bound. More
  explicitly, put $L_\Lambda = L+\Lambda^2$ and $Q =
  A_\Lambda^{1/2}L_\Lambda^{-1/2}$. Then $\|Q-I\| \le 1/6$, so $\|Q\| \le 7/6$
  and $\|Q^{-1}\| \le 6/5$ by its Neumann series. Taking adjoints under the
  norm-convergent resolvent integral gives $(L_\Lambda^{-1/2})^* =
  (L_\Lambda^*)^{-1/2}$. Apply the same estimate to $L^* = A_0-V^*$ and take
  adjoints. This bounds $L_\Lambda^{-1/2}A_\Lambda^{1/2}$ by $7/6$ and its
  inverse $A_\Lambda^{-1/2}L_\Lambda^{1/2}$ by $6/5$. Together with $Q$ and
  $Q^{-1} = L_\Lambda^{1/2}A_\Lambda^{-1/2}$, these are all four ordered
  comparisons in~\eqref{E:phase-root-comparisons}, each at most two.

  \medskip
  \noindent\textit{Step 2.~}We retain the numerical-range and phase-block
  bounds.
  For a unit vector $u$ put $s = \|A_0^{1/2}u\|$ and $\langle Lu,u\rangle =
  x+iy$. The form estimate at $\rho = \Lambda$ and, separately, at $\rho = s$
  gives
  \[
    |\langle Vu,u\rangle|\le D\Lambda^{a-2}(s^2+\Lambda^2),
    \qquad |\langle Vu,u\rangle|\le2Ds^a.
  \]
  At $s = 0$ the second bound follows by $\rho\downarrow0$. Thus, $x+\Lambda^2
  \ge (1-\delta_a)(s^2+\Lambda^2)$ and $|y| \le 2\delta_a\Lambda^{2-a}s^a$.
  These are the two inequalities used for the numerical-range parabola in
  Lemma~\ref{L:complex-phase}. Its scalar bounds and the bounded-operator
  Crouzeix--Palencia theorem apply to $L$ without selfadjointness. The four
  weighted block formulas there use only functions of $L$ and the root
  comparisons just proved. Their norms are at most $3e^{\Lambda h}$,
  $12(1+\Lambda h)e^{\Lambda h}$, $12e^{\Lambda h}$ and $12e^{\Lambda h}$. Their
  sum is at most $100e^{2\Lambda h}$. This uses neither selfadjointness of $V$
  nor commutation of $A_\Lambda$ with $L$.

  \medskip
  \noindent\textit{Step 3.~}We charge each norm loss only when an interval
  ends.
  For a common phase vector,
  \[
    \|z\|_{\rho_{\rm new}}
    \le\max\{1,\rho_{\rm old}/\rho_{\rm new}\}\|z\|_{\rho_{\rm old}}.
  \]
  At a jump, write $E$ for ending fields, $S$ for starting fields, and $R$ for
  those remaining. The ratio is at most
  \[
    \frac{\lambda+\sum_R\omega_i+\sum_E\omega_i}
    {\lambda+\sum_R\omega_i+\sum_S\omega_i}
    \le1+\sum_E\omega_i/\lambda
    \le\prod_E(1+\omega_i/\lambda).
  \]
  An initial increase of the weight has no loss. Count a final return to
  $X_\lambda$ as an ending jump at $T$. Each connected interval ends once,
  giving~\eqref{E:interval-metric}. This proves Lemma~\ref{L:interval-complex}.
\end{proof}

\subsection{Summing the coefficient bounds}

We first record an elementary bound for power series. The proof of
Proposition~\ref{P:interval-coefficients} uses it to sum over the ways a total
order $k$ splits among the interval fields.

\LeanMarker{L:inventory-product}
\begin{lemma}[A coefficient bound for a product with nonnegative coefficients]\label{L:inventory-product}
  For $k \ge 1$, define
  \[
    v_i=i^{-\beta}+k^{-1}i^{-(\beta-1)}
    +k^{\alpha-2}i^{-p},\qquad
    b_{i,l}=\frac{i^{-\beta l}(1+li/k)^2}{(l!)^\nu},\quad l\ge2.
  \]
  There is a finite $C = C(\alpha,\alpha_0)$, independent of $M,k$, such that
  \begin{equation}\label{E:inventory-product}
    [z^k]\prod_{i=1}^M\left(1+v_i z+\sum_{l\ge2}b_{i,l}z^l\right)
    \le\frac{C^k}{(k!)^\nu}.
  \end{equation}
\end{lemma}

\begin{proof}
  We evaluate the positive generating product at one order-dependent radius.
  We separate the third singleton term, then estimate the coarse and fine
  interval indices. Positivity converts the product bound back into a
  coefficient bound.

  \medskip
  \noindent\textit{Step 1.~}We separate and bound the third singleton
  contribution.
  Evaluate the positive product at $z = \theta k^\nu$, with $\theta = 1/4$. Put
  $P_i = \theta k^\nu i^{-\beta}$ and $Q_i = \theta k i^{-p}$. The third
  singleton term is $Q_i$ because $\nu+\alpha-2 = 1$. The inequality $1+A+B \le
  (1+A)(1+B)$ separates the factors $1+Q_i$. For $p > 1$,
  \[
    \sum_{i\ge1}\log(1+A i^{-p})
    \le\int_0^\infty\log(1+A x^{-p})\ud x
    =A^{1/p}\int_0^\infty\log(1+x^{-p})\ud x.
  \]
  Near zero, $\log(1+x^{-p})$ is bounded by a constant plus $p|\log x|$; at
  infinity it is at most $x^{-p}$, which is integrable for $p>1$.
  The last integral is finite. Hence, their logarithmic contribution is at most
  $Ck^{1/p} \le Ck$.

  \medskip
  \noindent\textit{Step 2.~}We estimate the coarse factors by the factorial
  series.
  For $x\ge0$, the sum of the $\nu$th powers of the nonnegative terms
  $x^{l/\nu}/l!$ is at most the $\nu$th power of their sum, and that sum is
  $e^{x^{1/\nu}}$.
  If $i \le k$, the remaining factor is at most $S_\nu(4P_i) =
  \sum_{l\ge0}(4P_i)^l/(l!)^\nu$: its singleton is at most $2P_i$, and
  $(1+li/k)^2 \le (1+l)^2 \le 4^l$. Since $\nu > 1$, $S_\nu(x) \le e^{\nu
  x^{1/\nu}}$. Therefore,
  \[
    \sum_{i\le k}\log S_\nu(4P_i)
    \le\nu(4\theta)^{1/\nu}k\sum_{i\ge1}i^{-\beta/\nu}
    \le Ck,
  \]
  because $\beta/\nu > 1$.

  \medskip
  \noindent\textit{Step 3.~}We sum the fine-index tails and extract the
  coefficient.
  If $i > k$, then $P_i \le \theta$. The remaining singleton tail is at most
  \[
    \theta\left(\frac1{\beta-1}+\frac1{\beta-2}\right)
    k^{\nu+1-\beta}=Ck^{\alpha_0}.
  \]
  For $l\ge2$, $0\le P\le1/4$ and $y\ge0$, the bounds $(1+ly)^2\le2l^2(1+y^2)$
  and $P^l\le P^2 4^{-(l-2)}$ leave a summable geometric series with a
  quadratic factor, independent of $P,y$.
  Also, uniformly for $0 \le P \le 1/4$ and $y \ge 0$,
  $\sum_{l\ge2}P^l(1+ly)^2/(l!)^\nu \le CP^2(1+y^2)$. Summing this bound gives
  \[
    Ck^{2\nu}\left(\frac{k^{1-2\beta}}{2\beta-1}
    +\frac{k^{-2}k^{3-2\beta}}{2\beta-3}\right)
    =Ck^{2\alpha_0-1}\le Ck.
  \]
  Use $\log(1+y) \le y$ for these fine factors. The entire product is
  consequently at most $e^{Ck}$. Positivity of its coefficients bounds the $k$th
  coefficient by $\theta^{-k}k^{-\nu k}e^{Ck} \le C^k/(k!)^\nu$. This proves
  Lemma~\ref{L:inventory-product}.
\end{proof}

Let $\mathcal{F}$ be a finite family of connected intervals $I_i$ meeting
$[0,T]$, with lengths $h_i$ sorted decreasingly. Let $W_i^\epsilon$ be
independent copies of the smoothed unit field, and suppose
\begin{equation}\label{E:interval-inventory}
  h_i\le L T/i,\qquad \sigma_i^2\le c h_i^{-\alpha_0}.
\end{equation}
Write $\mathcal{U}_{\mathcal{F}}^\epsilon(T;s)$ for phase evolution with
potential $s\sum_i\sigma_i\mathbf1_{I_i}(t)W_i^\epsilon(x)$.

\begin{figure}
\centering
\begin{tikzpicture}[x=1cm,y=1cm]
  \node[ham/note,text width=11.5cm] at (4.6,1.65)
  {For each interval, $r_i+q_i=2k_i$: count its occurrences on the two copies.};
  \node[ham/box,text width=3.1cm] (hamCross) at (.4,-.2)
  {One on each copy\\$k_i=1$, $(r_i,q_i)=(1,1)$\\$i\in\mathcal C$};
  \node[ham/box,text width=3.1cm] (hamSame) at (4.6,-.2)
  {Two on one copy\\$k_i=1$, $(r_i,q_i)=(2,0)$\\or $(0,2)$};
  \node[ham/box,text width=3.1cm] (hamHigher) at (8.8,-.2)
  {Higher degree\\$k_i\ge2$\\all splits retained};
  \node[ham/note,text width=11.5cm] at (4.6,-2.35)
  {The pair counts $k_i$ are then allocated across the intervals.\\
  Example: $k=2$ and $M=3$ give exactly six weak compositions.};
  \foreach \hamAllocationX/\hamAllocationTuple in
  {.4/{(2,0,0)},4.6/{(0,2,0)},8.8/{(0,0,2)}} {
  \node[ham/box,text width=2.7cm] at (\hamAllocationX,-3.6) {$\hamAllocationTuple$};
  }
  \foreach \hamAllocationX/\hamAllocationTuple in
  {.4/{(1,1,0)},4.6/{(1,0,1)},8.8/{(0,1,1)}} {
  \node[ham/box,text width=2.7cm] at (\hamAllocationX,-4.75) {$\hamAllocationTuple$};
  }
  \node[ham/note,text width=11.8cm] at (4.6,-6.1)
  {$[z^2]\prod_{i=1}^3\left(\sum_{l\ge0}a_{i,l}z^l\right)
  =\sum_{k_1+k_2+k_3=2}a_{1,k_1}a_{2,k_2}a_{3,k_3}$.\\
  Each displayed allocation contributes once.};
\end{tikzpicture}
\caption[Degree splits and pair-count allocations]{Two different counts must be kept separate. First split each
interval's even total degree between two copies built from the same Gaussian
fields; only distinct colors are independent. For $k_i=l$, there are
$2l+1$ possible splits. Cross singletons use
one insertion on each copy; a singleton confined to one copy is estimated
with the other background terms and has one ending factor, not two.
Then sum the pair counts over interval
indices. Ordinary positive-product coefficient extraction counts weak
compositions, without an extra multinomial coefficient. The displayed
$a_{i,l}$ are generic coefficients, not additional model parameters.
See \ProofRef{Lemma}{L:fixed-split} and
\ProofRef{Lemma}{L:inventory-product}.}
\label{Fig:interval-allocations}
\end{figure}

\LeanMarker{L:fixed-split}
\begin{lemma}[A fixed degree-split estimate]\label{L:fixed-split}
  Retain the finite interval fields above, with $R_\alpha > 0$ and $\sigma_i >
  0$. Fix integers $r_i,q_i,k_i \ge 0$ with $r_i+q_i = 2k_i$, and put $k =
  \sum_i k_i \ge 1$, $\lambda = k/T$. Represent $W_i^\epsilon(x) =
  \sum_jg_{ij}e_{ij}(x)$ as in Lemma~\ref{L:interval-contractions}. Retain $m_i$
  Gaussian modes, a spatial Fourier band $S_b$, and Hermite cutoff $N_i = 2k_i$
  on color $i$. Write $P = \bigotimes_iP_{N_i}$ and set
  \[
    B_i=P S_b\left(\sum_{j=1}^{m_i}g_{ij}M_{e_{ij}}\right)S_bP,
    \qquad V_{\boldsymbol z}(t)
    =\sum_i z_i\sigma_i\mathbf1_{I_i}(t)B_i.
  \]
  Let $\mathcal{M}^{\rm cut}_{\boldsymbol z}(T)$ be the phase evolution with
  this bounded potential and $A_0 = -\Delta$ on the band. Define $U_{\boldsymbol
  r}v = [\boldsymbol z^{\boldsymbol r}] \mathcal{M}^{\rm cut}_{\boldsymbol
  z}(T)v$, with the projected deterministic initial phase in the Gaussian
  vacuum. Put
  \[
    \begin{gathered}
      \mathcal{I}=\{i:k_i>0\},\qquad
      \mathcal{C}=\{i:k_i=1,\ r_i=q_i=1\},\\
      \mathcal{B}=\mathcal{I}\setminus\mathcal{C},\qquad
      \ell_i=\mathbf1_{\{r_i>0\}}+\mathbf1_{\{q_i>0\}}.
    \end{gathered}
  \]
  There is a finite $C_{\rm split} = C_{\rm split}(\alpha)$, uniform in all
  these cutoffs and in the finite interval family, such that
  \begin{equation}\label{E:fixed-split-bound}
    \begin{split}
      |\langle U_{\boldsymbol r}v,U_{\boldsymbol q}v
      \rangle_{L^2(\Omega;X_{1/T})}|
      \le{}&k^2 C_{\rm split}^k R_\alpha^k\|v\|_{1/T}^2
      \prod_{i\in\mathcal{C}}
      (\sigma_i^2h_i^2\lambda^{\alpha-2})\\
      &\quad\times\prod_{i\in\mathcal{B}}
      \frac{(\sigma_i^2h_i^{4-\alpha})^{k_i}}{(k_i!)^\nu}
      (1+k_i/(\lambda h_i))^{\ell_i}.
    \end{split}
  \end{equation}
  In particular, a singleton used twice on one copy has $\ell_i = 1$; a higher
  color has at most two metric factors.
\end{lemma}

\begin{proof}
  We estimate each vector coefficient on bounded cutoffs before pairing the
  two copies. We bound background segments on one common partition, count
  every change of phase norm once, and insert each cross color only from its
  vacuum. A cube integral and normalized vector Cauchy formula then give the
  coefficient with no extra permutation factor.

  \medskip
  \noindent\textit{Step 1.~}We set the active colors and background amplitude
  radii.
  We give a bound for each vector coefficient before taking their inner product.
  Every operator below acts on the finite Gaussian-mode, Hermite and
  spatial-band Hilbert space in the statement. In particular, it is bounded. Its
  piecewise-constant evolution is a finite product of block exponentials, hence
  is entire in every amplitude. All coefficient extractions in this proof are
  ordinary power-series extractions.

  An amplitude with $k_i = 0$ is evaluated at zero on both copies. If $k_i > 0$
  but $I_i\cap[0,T]$ has zero length, every coefficient with a positive power of
  that amplitude vanishes. The claimed inner product is then zero. We may
  therefore assume that every active interval has a positive-length intersection
  $I_i^T = (a_i,b_i)$ with $(0,T)$; endpoint conventions do not affect the
  evolution or its time integrals. Its length is at most $h_i$.

  Write $a = \alpha/2$, and use the positive constant $\delta_a$ in
  Lemma~\ref{L:complex-phase}. For every $i \in \mathcal{B}$ put
  \begin{equation}\label{E:split-radii}
    \omega_i=\frac{k_i}{h_i},\qquad
    d_i=2\sqrt{2k_i+1}\sqrt{R_\alpha},\qquad
    \varrho_i=\frac{\delta_a\omega_i^{2-a}}{\sigma_i d_i}.
  \end{equation}
  Fix one of the two copies. Its background set is $\mathcal{A}_r =
  \{i\in\mathcal{B}:r_i>0\}$ or $\mathcal{A}_q = \{i\in\mathcal{B}:q_i>0\}$;
  denote that set by $\mathcal{A}$. Evaluate every unused amplitude at zero.
  Colors in $\mathcal{C}$ have degree one on this copy. Put $m = |\mathcal{C}|$.
  Thus, there are at most $k$ background and cross colors altogether.

  \medskip
  \noindent\textit{Step 2.~}We bound background evolution on the common
  chronological partition.
  \emph{Background segments and the common partition.} Fix $|z_i| = \varrho_i$
  for $i \in \mathcal{A}$. With the cross amplitudes zero, the potential is
  $V_z(t) = \sum_{i\in\mathcal{A}}z_i\sigma_i\mathbf1_{I_i^T}(t)B_i$. Write
  $P_z(t,s)$ for its evolution, and set
  \[
    \Lambda(t)=\lambda+\sum_{i\in\mathcal{A}}
    \omega_i\mathbf1_{I_i^T}(t).
  \]
  For each $\rho > 0$, Lemma~\ref{L:interval-contractions} gives
  \[
    \|A_\rho^{-1/2}V_z(t)A_\rho^{-1/2}\|
    \le D(t)\rho^{a-2},\qquad
    D(t)=\sum_{i\in\mathcal{A}}|z_i|\sigma_i d_i
    \mathbf1_{I_i^T}(t)
    \le\delta_a\Lambda(t)^{2-a}.
  \]
  The last step follows from~\eqref{E:split-radii} and $\sum x_i^{2-a} \le (\sum
  x_i)^{2-a}$ for nonnegative $x_i$. The spatial resolvents here include the
  identity on every Gaussian factor. This is a deterministic operator inequality
  on their joint Hilbert space; the input can depend on the background fields.

  For the moment fix a tuple $(t_i)_{i\in\mathcal{C}}$ with $t_i \in I_i^T$, all
  distinct and away from the background endpoints. The excluded tuples have
  Lebesgue measure zero. Form the common partition
  \begin{equation}\label{E:split-partition}
    0=\tau_0<\tau_1<\cdots<\tau_M=T
  \end{equation}
  from these times and all background endpoints in $(0,T)$. On
  $(\tau_{j-1},\tau_j)$ let $\Lambda_j$ be the constant value of $\Lambda$, and
  put $P_j = P_z(\tau_j,\tau_{j-1})$. There are at most $2|\mathcal{A}|+m$
  interior partition points, so $M \le 2|\mathcal{A}|+m+1 \le 3k+1$.
  Lemma~\ref{L:interval-complex}, with its adjoint bound as stated there, gives
  \begin{equation}\label{E:split-segment-product}
    \begin{split}
      \prod_{j=1}^M\|P_j\|_{\mathcal{H}_{\Lambda_j}\to
      \mathcal{H}_{\Lambda_j}}
      &\le100^M\exp\left(2\sum_{j=1}^M
      \Lambda_j(\tau_j-\tau_{j-1})\right)\\
      &\le100^{3k+1}e^{4k}=:M_k,
      \qquad \mathcal{H}_\rho=L^2(\Omega;X_\rho).
    \end{split}
  \end{equation}
  Indeed, the sum in the exponential is exactly $\lambda
  T+\sum_{i\in\mathcal{A}}\omega_i(b_i-a_i) \le k+\sum_{i\in\mathcal{A}}k_i \le
  2k$. Insertion times refine a partition; they do not repeat this integral.

  \medskip
  \noindent\textit{Step 3.~}We count the phase-metric changes.
  \emph{Every identity between phase norms.} Put $\Lambda_0 = \Lambda_{M+1} =
  \lambda$. At $\tau_j$ let $\mathscr{I}_j$ be the identity from
  $\mathcal{H}_{\Lambda_j}$ to $\mathcal{H}_{\Lambda_{j+1}}$, for $0 \le j \le
  M$. Its norm is at most
  \[
    q_j:=\max\{1,\Lambda_j/\Lambda_{j+1}\}.
  \]
  This follows directly from the two components of the phase norm and the
  spectral inequality $A_{\Lambda_{j+1}}^{-1} \le q_j^2 A_{\Lambda_j}^{-1}$. At
  an interior background endpoint, let $E_j,S_j,R_j$ denote the ending, starting
  and remaining background colors. Then
  \[
    \frac{\Lambda_j}{\Lambda_{j+1}}
    =\frac{\lambda+\sum_{R_j}\omega_i+\sum_{E_j}\omega_i}
    {\lambda+\sum_{R_j}\omega_i+\sum_{S_j}\omega_i}
    \le1+\sum_{E_j}\omega_i/\lambda
    \le\prod_{i\in E_j}(1+\omega_i/\lambda).
  \]
  The last product is at least one, so it also bounds $q_j$. At an
  insertion-only cut the heights agree and $q_j = 1$. Initially $q_0 = 1$. At
  $T$ the same inequality charges all intervals with $b_i = T$ as ending there.
  Each connected $I_i^T$ ends at precisely one such cut. Consequently,
  \begin{equation}\label{E:interval-copy-metric}
    \prod_{j=0}^M q_j
    \le D_{\mathcal{A}}:=
    \prod_{i\in\mathcal{A}}(1+k_i/(\lambda h_i)).
  \end{equation}
  This product includes the initial and final conversions. There is no new copy
  of $D_{\mathcal{A}}$ at an insertion time or an overlapping interval's start.

  \begin{figure}
\centering
\begin{tikzpicture}[x=.96cm,y=1cm]
  \draw[ham/arrow] (-.3,0)--(12.2,0) node[right] {$t$};
  \draw[hamblue,line width=2pt] (1,3.0)--(7,3.0);
  \draw[hamteal,line width=2pt,dashed] (3,2.45)--(9,2.45);
  \draw[hamgold,line width=2pt,densely dotted] (5,1.9)--(11,1.9);
  \node[anchor=west,font=\small] at (7.2,3.0) {$I_1^T$};
  \node[anchor=west,font=\small] at (9.2,2.45) {$I_2^T$};
  \node[anchor=west,font=\small] at (11.2,1.9) {$I_3^T$};
  \draw[hamblue,semithick]
  (0,.3)--(1,.3)--(1,.55)--(3,.55)--(3,.8)--(5,.8)
  --(5,1.05)--(7,1.05)--(7,.8)--(9,.8)--(9,.55)
  --(11,.55)--(11,.3)--(12,.3);
  \node[font=\small,anchor=east] at (-.15,.8) {$\Lambda(t)$};
  \foreach \hamStart in {1,3,5} {
  \draw[hamblue!40,dashed] (\hamStart,-.2)--(\hamStart,1.1);
  }
  \foreach \hamEnd in {7,9,11} {
  \draw[hamgold,dashed] (\hamEnd,-.2)--(\hamEnd,1.1);
  }
  \filldraw[draw=hamblue,fill=white] (4,0) circle (2pt);
  \node[ham/note,text width=3.5cm] at (2.9,-1.1)
  {Cuts with no ending interval:\\$q_j=1$};
  \node[ham/note,text width=2.1cm] at (7,-1.1) {$q_j\le1+\omega_1/\lambda$};
  \node[ham/note,text width=2.1cm] at (9,-1.1) {$q_j\le1+\omega_2/\lambda$};
  \node[ham/note,text width=2.1cm] at (11,-1.1) {$q_j\le1+\omega_3/\lambda$};
  \node[ham/box,text width=10.3cm] at (6,-3.05)
  {$\displaystyle\prod_jq_j\le
  \prod_{i\in\mathcal A}(1+\omega_i/\lambda)$\\
  Each connected interval is charged at its end, once.};
\end{tikzpicture}
\caption[One phase-norm charge at each interval end]{Changes of phase norm are charged at interval endpoints.
The staircase illustrates $\Lambda(t)=\lambda+
\sum_{i\in\mathcal A}\omega_i\mathbf1_{I_i^T}(t)$ for three overlapping
intervals; its heights are schematic. The conversion satisfies
$q_j=\max\{1,\Lambda_j/\Lambda_{j+1}\}$. A cut containing only starts or insertions adds no factor. At a shared ending
cut, charge all ending intervals together; the final conversion at $T$ is
also included. This prevents repeating the full product at every insertion.
See \ProofRef{Lemma}{L:fixed-split}.}
\label{Fig:interval-end-metrics}
\end{figure}
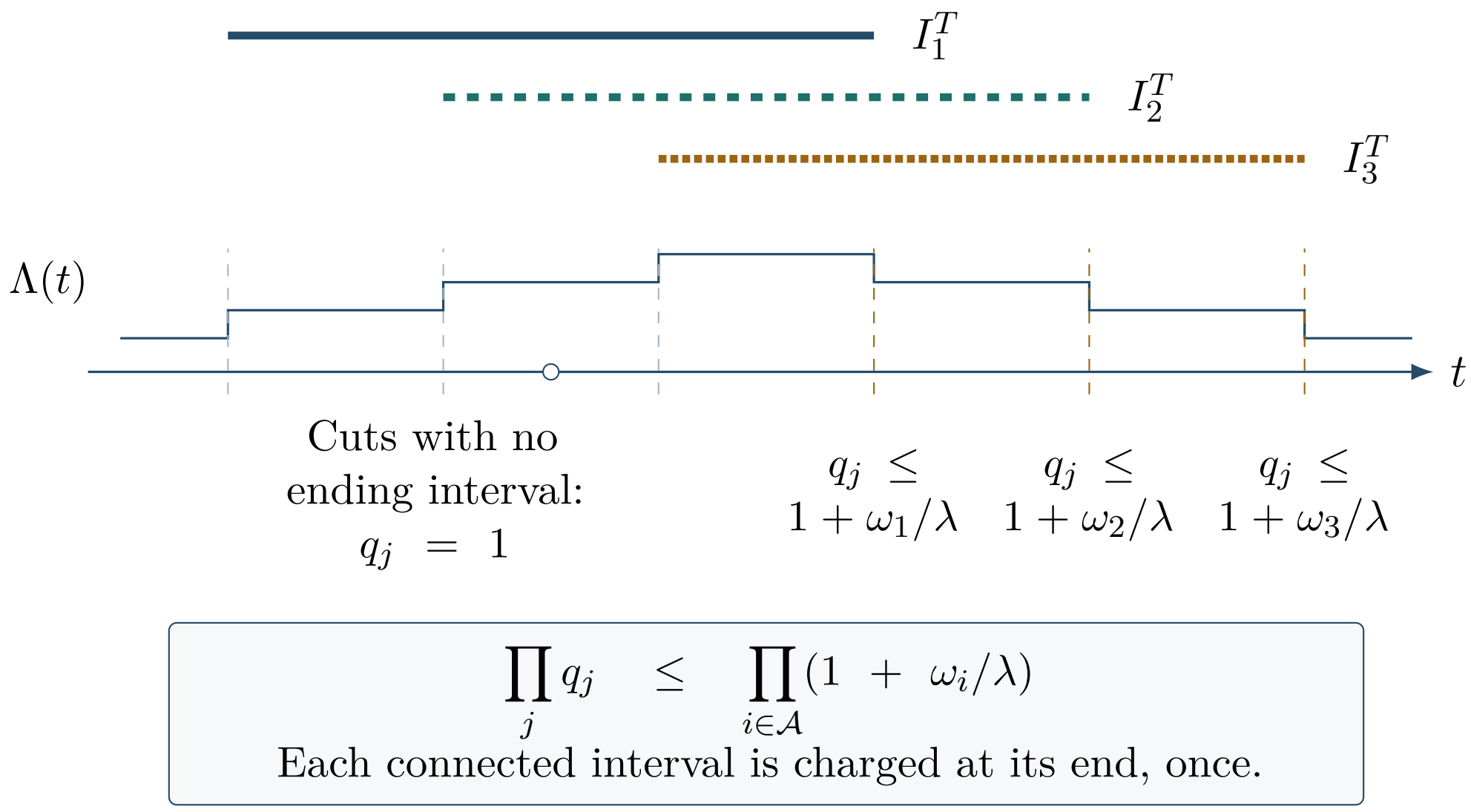

  \medskip
  \noindent\textit{Step 4.~}We estimate every cross-color insertion on its
  actual tensor subspace.
  \emph{The cross-color subspaces and the word norm.} Let $\mathscr{G}_i^{(0)}$
  be the vacuum of color $i$ and $\mathscr{G}_i^{(1)}$ its first Gaussian
  degree. For $S \subset \mathcal{C}$, use the tensor subspace
  \[
    \mathscr{Z}_S=
    \bigotimes_{i\in\mathcal{A}}\mathscr{G}_i^{\le N_i}
    \otimes\bigotimes_{i\in S}\mathscr{G}_i^{(1)}
    \otimes\bigotimes_{i\in\mathcal{C}\setminus S}
    \mathscr{G}_i^{(0)},
  \]
  with every omitted Gaussian color in its vacuum and with the common spatial
  band understood. Each $P_j$ preserves the corresponding phase space over
  $\mathscr{Z}_S$: its potential acts only on colors in $\mathcal{A}$, and its
  spatial operators act identically on all cross Gaussian factors. The
  identities $\mathscr{I}_j$ preserve the same tensor subspace.

  Define $J_i(u,w) = (0,\sigma_i B_i u)$ for $i \in \mathcal{C}$. If $i\notin
  S$, multiplication by its field takes its vacuum into first Gaussian degree.
  The Hermite cutoff $N_i = 2$ does not alter this step. Thus, $J_i$ maps the
  phase over $\mathscr{Z}_S$ into the phase over $\mathscr{Z}_{S\cup\{i\}}$.
  Collect all other Gaussian factors into $\mathscr{R}_i$. Every input has the
  form $1_i\otimes(u,w)$, with $(u,w) \in \mathscr{R}_i\otimes X_\rho$. The
  cutoffs on other colors are identities on its retained degrees, which
  multiplication by the color-$i$ field leaves unchanged. Hence,
  \[
    J_i(1_i\otimes(u,w))
    =\sigma_i\sum_{j=1}^{m_i}g_{ij}\otimes
    (0,S_bM_{e_{ij}}S_bu).
  \]
  Orthogonality $\E[g_{ij}g_{il}] = \mathbf1_{\{j=l\}}$ gives
  \[
    \begin{split}
      \|J_i(1_i\otimes(u,w))\|_{\mathcal{H}_\rho}^2
      &=\sigma_i^2\sum_{j=1}^{m_i}
      \|A_\rho^{-1/2}S_bM_{e_{ij}}S_bu
      \|_{\mathscr{R}_i\otimes L^2_x}^2\\
      &\le\sigma_i^2R_\alpha\rho^{\alpha-2}
      \|u\|_{\mathscr{R}_i\otimes L^2_x}^2\\
      &\le\sigma_i^2R_\alpha\rho^{\alpha-2}
      \|1_i\otimes(u,w)\|_{\mathcal{H}_\rho}^2.
    \end{split}
  \]
  The first inequality is~\eqref{E:interval-column} tensored with
  $\mathscr{R}_i$ and allows arbitrary dependence on earlier colors. Thus, at
  the insertion time,
  \begin{equation}\label{E:split-fresh-insertion}
    \|J_i\|\le\sigma_i\sqrt{R_\alpha}
    \Lambda(t_i)^{(\alpha-2)/2}
    \le c_i:=\sigma_i\sqrt{R_\alpha}\lambda^{(\alpha-2)/2}.
  \end{equation}
  Starting with $\mathscr{Z}_\varnothing$, these subspaces describe every
  partial word. They justify each new-color estimate even though the current
  vector depends on the background and earlier cross colors. No occupied-color
  annihilation is estimated by a vacuum bound.

  Equivalently, apply Lemma~\ref{L:finite-word-cutoff} to the cross colors, with
  $\mathscr{X}$ containing the bounded background Gaussian space and the spatial
  phase space. Background segments act only on $\mathscr{X}$, and each cross
  color occurs once. Its cutoff is an exact identity on every prefix of the
  selected word. This does not identify an untruncated background propagator
  with a fixed cutoff.

  At an interior cut $\tau_j$ let $B_j^{\rm cut} = J_i$ if $\tau_j = t_i$ and
  let $B_j^{\rm cut}$ be the identity otherwise. Apply $B_j^{\rm cut}$ after
  $\mathscr{I}_j$; at an insertion cut the two heights are equal. The
  chronological word, as a map from $\mathcal{H}_\lambda$ back to
  $\mathcal{H}_\lambda$, is exactly
  \[
    \mathcal{W}_{\mathcal{A}}(z,\boldsymbol t)
    =\mathscr{I}_M P_M B_{M-1}^{\rm cut}\mathscr{I}_{M-1}P_{M-1}
    \cdots B_1^{\rm cut}\mathscr{I}_1P_1\mathscr{I}_0.
  \]
  For $M = 1$ this means $\mathscr{I}_1P_1\mathscr{I}_0$. The rightmost factor
  acts first; no spatial operators are commuted. Multiplying their norms, using
  the subspaces just specified, gives the explicit bound
  \begin{equation}\label{E:split-word-product}
    \begin{split}
      \|\mathcal{W}_{\mathcal{A}}(z,\boldsymbol t)v\|_\lambda
      &\le\left(\prod_{j=1}^M\|P_j\|\right)
      \left(\prod_{j=0}^M q_j\right)
      \left(\prod_{i\in\mathcal{C}}c_i\right)\|v\|_\lambda\\
      &\le M_kD_{\mathcal{A}}\|v\|_\lambda
      \prod_{i\in\mathcal{C}}c_i.
    \end{split}
  \end{equation}
  The background norm bound is on the joint Gaussian Hilbert space; it is not
  obtained by conditioning on successive fields.

  \medskip
  \noindent\textit{Step 5.~}We integrate over the insertion cube and extract
  vector coefficients.
  \emph{The insertion cube and vector Cauchy formula.} Give the cross colors
  amplitudes $y_i$ and define the vector
  \[
    \Phi_{\mathcal{A}}(z)
    :=\left[\prod_{i\in\mathcal{C}}y_i\right]
    \mathcal{M}^{\rm cut}_{(z,y,0)}(T)v.
  \]
  Variation of constants about $P_z$ and iteration give a Dyson series. A
  monomial of degree one in every $y_i$ uses exactly $m$ insertions, with each
  distinct cross color used once. Hence,
  \begin{equation}\label{E:interval-cross-cube}
    \Phi_{\mathcal{A}}(z)
    =\sum_{\pi}\int_{D_\pi}
    \mathcal{W}_{\mathcal{A}}(z,\boldsymbol t)v\,\ud\boldsymbol t,
    \qquad
    D_\pi=\left\{\boldsymbol t\in\prod_{i\in\mathcal{C}}I_i^T:
    t_{\pi_1}<\cdots<t_{\pi_m}\right\}.
  \end{equation}
  Each ordered Dyson integral has coefficient one. The permutation simplices
  $D_\pi$ partition the cube up to null diagonals. In particular,
  \[
    \sum_\pi |D_\pi|=\prod_{i\in\mathcal{C}}|I_i^T|
    \le\prod_{i\in\mathcal{C}}h_i.
  \]
  Minkowski applied to~\eqref{E:interval-cross-cube} and the single-word
  bound~\eqref{E:split-word-product} therefore gives
  \begin{equation}\label{E:interval-cross-vector}
    \|\Phi_{\mathcal{A}}(z)\|_\lambda
    \le M_kD_{\mathcal{A}}\|v\|_\lambda
    \prod_{i\in\mathcal{C}}c_i h_i.
  \end{equation}
  There is no factor $m!$ or $1/m!$ in this estimate. For $m = 0$, the cube has
  unit mass and $\Phi_{\mathcal{A}}(z) = P_z(T,0)v$.

  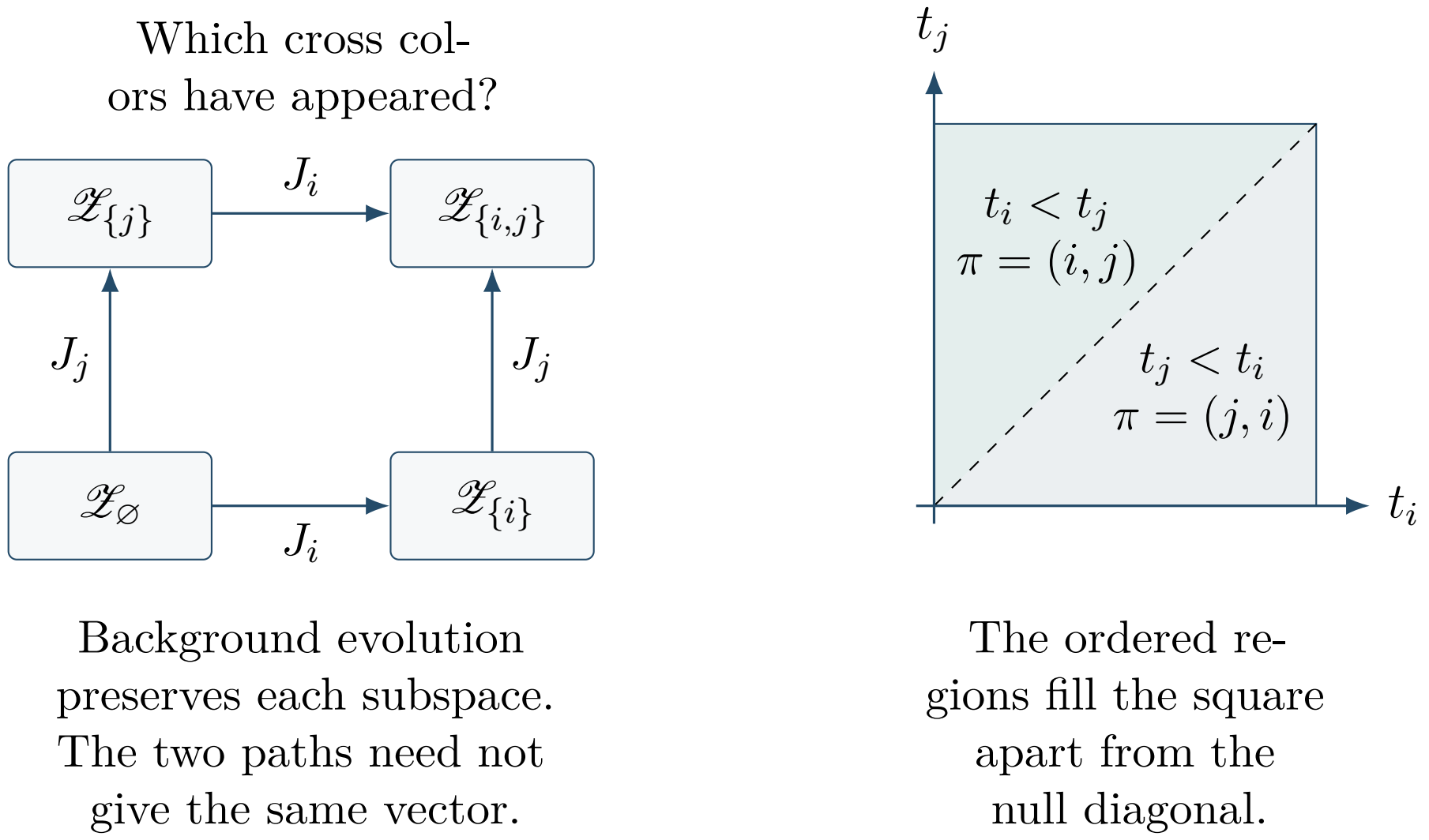
\begin{figure}
\centering
\begin{tikzpicture}[x=1cm,y=1cm]
  \node[ham/note,text width=4.8cm] at (1.6,3.65)
  {Which cross colors have appeared?};
  \node[ham/box,text width=1.35cm] (hamEmpty) at (0,0) {$\mathscr Z_\varnothing$};
  \node[ham/box,text width=1.35cm] (hamI) at (3.2,0) {$\mathscr Z_{\{i\}}$};
  \node[ham/box,text width=1.35cm] (hamJ) at (0,2.45) {$\mathscr Z_{\{j\}}$};
  \node[ham/box,text width=1.35cm] (hamIJ) at (3.2,2.45) {$\mathscr Z_{\{i,j\}}$};
  \draw[ham/arrow] (hamEmpty)--node[below,font=\small] {$J_i$} (hamI);
  \draw[ham/arrow] (hamEmpty)--node[left,font=\small] {$J_j$} (hamJ);
  \draw[ham/arrow] (hamI)--node[right,font=\small] {$J_j$} (hamIJ);
  \draw[ham/arrow] (hamJ)--node[above,font=\small] {$J_i$} (hamIJ);
  \node[ham/note,text width=4.8cm] at (1.6,-1.85)
  {Background evolution preserves each subspace.\\
  The two paths need not give the same vector.};
  \begin{scope}[xshift=6.9cm]
    \fill[hamblue!9] (0,0)--(3.2,0)--(3.2,3.2)--cycle;
    \fill[hamteal!12] (0,0)--(0,3.2)--(3.2,3.2)--cycle;
    \draw[hamblue] (0,0) rectangle (3.2,3.2);
    \draw[dashed] (0,0)--(3.2,3.2);
    \draw[ham/arrow] (-.15,0)--(3.65,0) node[right] {$t_i$};
    \draw[ham/arrow] (0,-.15)--(0,3.65) node[above] {$t_j$};
    \node[ham/note] at (.95,2.25) {$t_i<t_j$\\$\pi=(i,j)$};
    \node[ham/note] at (2.25,.95) {$t_j<t_i$\\$\pi=(j,i)$};
    \node[ham/note,text width=4.3cm] at (1.6,-1.85)
    {The ordered regions fill the square\\apart from the null diagonal.};
  \end{scope}
\end{tikzpicture}
\caption[Cross-color subsets and chronological regions]{For two distinct cross colors, each insertion starts from the
vacuum in its own unused color. The left square records subspaces only;
it is not a commutative operator diagram. The right picture uses identical
time windows for clarity. In general the ordered regions partition
$I_i^T\times I_j^T$ up to the equal-time diagonal. Each chronological
operator word has coefficient one, and the total volume is
$|I_i^T||I_j^T|$, without an extra $2!$ or $1/2!$. Spatial operators are
not reordered. See \ProofRef{Lemma}{L:fixed-split}.}
\label{Fig:cross-colour-cube}
\end{figure}

  The vector $\Phi_{\mathcal{A}}$ is entire in the background amplitudes. On the
  first copy its coefficient is given by the normalized formula
  \begin{equation}\label{E:split-vector-cauchy}
    \begin{split}
      U_{\boldsymbol r}v
      ={}&\left(\prod_{i\in\mathcal{A}_r}\varrho_i^{-r_i}\right)
      \int_{[0,2\pi]^{\mathcal{A}_r}}
      e^{-i\sum_{i\in\mathcal{A}_r}r_i\theta_i}
      \Phi_{\mathcal{A}_r}((\varrho_i e^{i\theta_i})_i)
      \prod_{i\in\mathcal{A}_r}\frac{\ud\theta_i}{2\pi}.
    \end{split}
  \end{equation}
  This is precisely the desired coefficient because each cross exponent equals
  one and every omitted exponent is zero. The normalized angular measure has
  total mass one. Thus,~\eqref{E:interval-cross-vector} bounds $\|U_{\boldsymbol
  r}v\|_\lambda$ with the additional factor
  $\prod_{i\in\mathcal{A}_r}\varrho_i^{-r_i}$, and the identical formula bounds
  the second copy using $\mathcal{A}_q,q_i$. No holomorphic formula is applied
  to a squared norm.

  \medskip
  \noindent\textit{Step 6.~}We pair the two coefficients and account for all
  constants.
  \emph{The inner product and all constants.} Since $\lambda = k/T$ and $k \ge
  1$,
  \[
    \|w\|_{1/T}\le k\|w\|_\lambda,\qquad
    \|v\|_\lambda\le\|v\|_{1/T}.
  \]
  Cauchy--Schwarz between the two vectors on the same Gaussian space now gives,
  before any constants are absorbed,
  \begin{equation}\label{E:split-before-absorption}
    \begin{split}
      |\langle U_{\boldsymbol r}v,U_{\boldsymbol q}v
      \rangle_{\mathcal{H}_{1/T}}|
      \le{}&k^2 M_k^2D_{\mathcal{A}_r}D_{\mathcal{A}_q}
      \|v\|_{1/T}^2
      \prod_{i\in\mathcal{B}}\varrho_i^{-2k_i}\\
      &\quad\times\prod_{i\in\mathcal{C}}
      (R_\alpha\sigma_i^2h_i^2\lambda^{\alpha-2}).
    \end{split}
  \end{equation}
  The radius exponent is $r_i+q_i = 2k_i$. The metric product is exactly
  \[
    D_{\mathcal{A}_r}D_{\mathcal{A}_q}
    =\prod_{i\in\mathcal{B}}(1+k_i/(\lambda h_i))^{\ell_i}.
  \]
  In particular, a singleton on only one copy has one ending charge, not two.
  For every $i \in \mathcal{B}$, $d_i^2 \le 12k_iR_\alpha$ gives
  \begin{equation}\label{E:split-radius-factorial}
    \begin{split}
      \varrho_i^{-2k_i}
      &=\left(\frac{\sigma_i^2d_i^2}{\delta_a^2}\right)^{k_i}
      \left(\frac{h_i}{k_i}\right)^{(4-\alpha)k_i}\\
      &\le(12\delta_a^{-2})^{k_i}
      (R_\alpha\sigma_i^2h_i^{4-\alpha})^{k_i}
      k_i^{-(3-\alpha)k_i}\\
      &\le(12\delta_a^{-2})^{k_i}
      \frac{(R_\alpha\sigma_i^2h_i^{4-\alpha})^{k_i}}
      {(k_i!)^{3-\alpha}}.
    \end{split}
  \end{equation}
  The last step uses $k_i! \le k_i^{k_i}$. Finally,
  \[
    M_k^2=100^{6k+2}e^{8k}\le(100^8e^8)^k,
    \qquad |\mathcal{C}|+\sum_{i\in\mathcal{B}}k_i=k.
  \]
  Substituting these inequalities in~\eqref{E:split-before-absorption}
  proves~\eqref{E:fixed-split-bound} with the explicit admissible choice
  \begin{equation}\label{E:split-explicit-constant}
    C_{\rm split}=100^8e^8
    \max\{1,12\delta_{\alpha/2}^{-2}\}.
  \end{equation}
  This constant depends only on $\alpha$. No permutation, Gaussian-mode, factor
  depending on the number of intervals or further degree factor remains. The
  projected initial phase has norm at most the displayed original $\|v\|_{1/T}$.
  This proves Lemma~\ref{L:fixed-split}.
\end{proof}

\LeanMarker{P:interval-coefficients}
\begin{proposition}[Uniform interval coefficients]\label{P:interval-coefficients}
  There is $C = C(\alpha,\alpha_0,L,c)$ such that, for real Schwartz phase $v$,
  every $k \ge 1$ satisfies
  \begin{equation}\label{E:interval-coefficient-bound}
    \left|[s^{2k}]\E\|\mathcal{U}_{\mathcal{F}}^\epsilon(T;s)v\|_{1/T}^2\right|
    \le\frac{C^k R_\alpha^k T^{\beta k}}{(k!)^\nu}\|v\|_{1/T}^2,
  \end{equation}
  uniformly in $\epsilon > 0$ and the finite family. In particular,
  \begin{equation}\label{E:interval-phase-bound}
    \E\|\mathcal{U}_{\mathcal{F}}^\epsilon(T;1)v\|_{1/T}^2
    \le e\exp\{C R_\alpha^{1/\nu}T^{\beta/\nu}\}\|v\|_{1/T}^2.
  \end{equation}
\end{proposition}

\begin{proof}
  We first construct the ordinary coefficient series at fixed smoothing and
  interval family. We retain only even total degree in each color and use
  exact finite-word cutoffs. The fixed-split estimate and sorted interval
  lengths reduce the uniform sum to the inventory product, after which
  separate cutoff limits identify and sum the full series.

  \medskip
  \noindent\textit{Step 1.~}We set up the bounded cutoffs and entire vector
  series.
  Zero noise is trivial; omit any zero $\sigma_i$. Enlarge the interval-length
  constant $L$ to $\max\{1,L\}$ if necessary. The hypotheses remain true and the
  final constant may depend on $L$, so throughout this proof we may assume $L
  \ge 1$. Replacing $(\sigma_i,W_i^\epsilon)$ by
  $(|\sigma_i|,\operatorname{sgn}(\sigma_i)W_i^\epsilon)$ leaves the potential
  unchanged and preserves the joint law of the independent real Gaussian fields.
  Thus, we may assume $\sigma_i > 0$. We first work with finitely many Gaussian
  coordinates of each field and a bounded spatial Fourier band. Denote their
  common phase Hilbert space by $\mathcal{H}_\rho = L^2(\Omega;X_\rho)$. The
  initial data are projected to the band. This contracts every $X_\rho$ norm; we
  keep the original $\|v\|_\rho$ on the right of the estimates. The projections
  converge strongly to the original data.

  \emph{Ordinary coefficients and their existence.} Fix $\lambda > 0$
  temporarily. Give interval $i$ its own amplitude $z_i$ and write
  \[
    \mathcal{U}(T;\boldsymbol z)v
    =\sum_{\boldsymbol r\in\N_0^{\mathcal{F}}}
    \boldsymbol z^{\boldsymbol r}U_{\boldsymbol r}v.
  \]
  These are ordinary power-series coefficients, without factorials incorporated
  in their definition. At fixed smoothing and finite interval family, put
  $v_{\max} = \sup_t\sum_i\sigma_i^2\mathbf1_{I_i}(t) \Gamma_\epsilon(0) <
  \infty$. On Gaussian degrees at most $n$, multiplication by the instantaneous
  field has norm at most $2\sqrt{n+1}\sqrt{v_{\max}}$. The free-phase energy
  calculation at the end of the proof of
  Proposition~\ref{P:static-phase-variance}, which uses no dimension-specific
  or covariance assumption, gives the bound $e^{\lambda h/2}$ in $X_\lambda$
  also on each spatial Fourier band. A phase insertion
  $(u,w)\mapsto(0,Wu)$ costs at most its unweighted multiplication norm divided
  by $\lambda$. Thus, on $\max_i|z_i| \le z_*$ the total order-$n$ phase is
  bounded by
  \begin{equation}\label{E:interval-entire-majorant}
    e^{\lambda T/2}
    \frac{(2z_*\sqrt{n+1}\sqrt{v_{\max}}T/\lambda)^n}{n!}
    \|v\|_\lambda .
  \end{equation}
  The series of these bounds converges for every finite $z_*$. Its locally
  uniform convergence defines an entire vector function. We also verify
  coefficientwise absolute convergence. Fix $z_* > 0$ and put
  \[
    K_{2z_*}=\sup_{\max_i|z_i|\le2z_*}
    \|\mathcal{U}(T;\boldsymbol z)v\|_{\mathcal{H}_{1/T}}<\infty,
  \]
  using the majorant with its temporary parameter $\lambda = 1/T$. Normalized
  vector Cauchy integration on this larger polydisc gives
  \[
    \|U_{\boldsymbol r}v\|_{\mathcal{H}_{1/T}}
    \le K_{2z_*}(2z_*)^{-|\boldsymbol r|},\qquad
    \sum_{\boldsymbol r}z_*^{|\boldsymbol r|}
    \|U_{\boldsymbol r}v\|_{\mathcal{H}_{1/T}}
    \le2^{|\mathcal{F}|}K_{2z_*}.
  \]
  The sum is the product of $|\mathcal{F}|$ geometric series. Cauchy--Schwarz
  therefore bounds the absolute energy double series by
  $4^{|\mathcal{F}|}K_{2z_*}^2$. This finite factor justifies the coefficient
  identities at fixed family and smoothing; it does not enter the uniform
  coefficient estimate below.

  \medskip
  \noindent\textit{Step 2.~}We identify the energy coefficient and its color
  splits.
  Set $Q_{\boldsymbol r,\boldsymbol q} = \langle U_{\boldsymbol
  r}v,U_{\boldsymbol q}v \rangle_{\mathcal{H}_{1/T}}$. Both vectors use the same
  Gaussian fields. Independence is between distinct colors, not between the two
  copies in this inner product. For real $s$, substitution $z_i = s$ gives the
  exact identity
  \begin{equation}\label{E:interval-energy-splits}
    [s^{2k}]\E\|\mathcal{U}(T;s)v\|_{1/T}^2
    =\sum_{|\boldsymbol r|+|\boldsymbol q|=2k}
    Q_{\boldsymbol r,\boldsymbol q}.
  \end{equation}
  Reflection of Gaussian color $i$ multiplies $U_{\boldsymbol r}v$ by
  $(-1)^{r_i}$. Hence, $Q_{\boldsymbol r,\boldsymbol q} = 0$ unless $r_i+q_i =
  2k_i$ for every $i$, where $\sum_i k_i = k$. If $k_i = 0$, both coefficients
  are obtained by setting $z_i = 0$; this removes that field exactly. There are
  at most $k$ active colors. For fixed $\boldsymbol k$, the number of splits is
  $\prod_i(2k_i+1) \le 3^k$.

  \medskip
  \noindent\textit{Step 3.~}We retain the selected coefficients exactly at
  finite Hermite degrees.
  For this fixed $\boldsymbol k$, truncate the Hermite degree of color $i$ at
  $N_i = 2k_i$. A word contributing to $U_{\boldsymbol r}$ has $r_i \le 2k_i$
  insertions of color $i$, and the other copy has $q_i \le 2k_i$.
  Lemma~\ref{L:finite-word-cutoff}, with the spatial phase space as
  $\mathscr{X}$ and free phase factors between insertions, therefore preserves
  both selected vector coefficients exactly. Their Duhamel time integrals and
  finite word sums preserve the equality. Creation followed by annihilation is
  retained; there is no Hermite-degree limit in this coefficient calculation.

  \medskip
  \noindent\textit{Step 4.~}We convert the split estimate into the interval
  inventory weights.
  For the uniform coefficient estimate put $\lambda = k/T$. Apply
  Lemma~\ref{L:fixed-split} to each selected coefficient, retaining the original
  sorted interval labels in the final sum. Its factor $\ell_i$ is one for a
  singleton mean and at most two for every higher color. The separate cross
  factor applies only at $(r_i,q_i) = (1,1)$. Absorb $k^2C_{\rm split}^k3^kc^k$
  into $C^k$, using $k^2 \le 4^k$. After $\sigma_i^2 \le c h_i^{-\alpha_0}$, the
  resulting local weights, with this common per-total-order constant understood,
  are
  \[
    \begin{array}{c|c}
      \text{case}&\text{weight in energy}\\ \hline
      (2,0)\text{ or }(0,2)&R_\alpha h_i^\beta
      (1+1/(\lambda h_i))\\
      (1,1)&R_\alpha\lambda^{\alpha-2}h_i^{2-\alpha_0}\\
      k_i=l\ge2&R_\alpha^l h_i^{\beta l}
      (1+l/(\lambda h_i))^2/(l!)^\nu .
    \end{array}
  \]
  The first metric factor occurs on only one copy; the last costs at most one
  factor on each copy. To justify the length substitution, expand the entire
  higher-degree weight, before its positive constant:
  \[
    h^{\beta l}(1+l/(\lambda h))^2
    =h^{\beta l}+\frac{2l}{\lambda}h^{\beta l-1}
    +\frac{l^2}{\lambda^2}h^{\beta l-2}.
  \]
  Since $\beta > 2$ and $l \ge 2$, all three powers are positive. The singleton
  mean likewise has powers $\beta,\beta-1$, and the cross singleton has power $p
  = 2-\alpha_0 > 0$. Thus, each whole local weight increases with $h$; we do not
  bound the inverse-length factor separately. At $h = LT/i$, $1+l/(\lambda h) =
  1+li/(kL) \le 1+li/k$, since $L \ge 1$. The remaining powers of $L$ are
  absorbed in the allowed $C^k$. For a singleton let $m_i,c_i$ be the first and
  second table weights. Every one of its splits is bounded by $m_i+c_i$. The
  higher-degree table weight is already uniform over its splits. Hence, the
  previously counted $\prod_i(2k_i+1) \le 3^k$ bounds all split choices once,
  using these combined singleton weights; that is the factor absorbed into $C^k$
  above. Substituting $h_i \le LT/i$ and $\lambda = k/T$ therefore extracts
  $C^kR_\alpha^kT^{\beta k}$, leaving $v_i$ at $k_i = 1$ and $b_{i,l}$ at $k_i =
  l \ge 2$ from Lemma~\ref{L:inventory-product}.

  \medskip
  \noindent\textit{Step 5.~}We sum ordinary coefficients without a multinomial
  factor.
  There is no further multinomial factor: for nonnegative weights $w_{i,0} = 1$,
  ordinary multiplication of power series gives exactly
  \[
    \sum_{\sum_i k_i=k}\prod_i w_{i,k_i}
    =[z^k]\prod_i\left(\sum_{l\ge0}w_{i,l}z^l\right).
  \]
  Apply this identity to the sum in~\eqref{E:interval-energy-splits}, using the
  two singleton alternatives together in $v_i$. Lemma~\ref{L:inventory-product}
  now proves \eqref{E:interval-coefficient-bound} at the finite cutoffs.

  \medskip
  \noindent\textit{Step 6.~}We remove cutoffs at fixed coefficients and sum
  the uniform bound.
  The closed operator used below is multiplication by the measurable smoothed
  Gaussian field on its natural Gaussian--spatial $L^2$ domain. Mode sums
  agree with this multiplication on finite-coordinate inputs in the limit;
  closedness then identifies the graph limit as those inputs are approximated
  in $L^2$.
  \emph{Identification and summation.} Keep positive smoothing, finite interval
  family and a fixed coefficient. The finite-degree multiplication bound used
  in~\eqref{E:interval-entire-majorant} is uniform in mode and band cutoffs. On
  finite-coordinate Gaussian vectors the omitted annihilation tail is zero and
  the creation tail tends to zero by its omitted pointwise variances and
  dominated convergence. Density in each finite degree extends this strong
  convergence. Closedness identifies the limit with multiplication by the
  measurable smoothed field. Fourier projections converge strongly and commute
  with free evolution. Dominated time integration thus identifies every fixed
  Duhamel word. The selected coefficients of the local Hermite truncation were
  already equal to those untruncated words. The full series is constructed
  separately using~\eqref{E:interval-entire-majorant}, choosing a sufficient
  degree cutoff anew at each order. A fixed-degree truncated propagator is not
  identified with the full evolution. The entire majorant gives convergence
  uniformly on compact amplitude sets. The coefficients estimated above are
  therefore those of the finite-interval model, with no zero-smoothing limit at
  this step.

  The nonuniform number $v_{\max}$ was used only for this identification. The
  uniform coefficient bound contains no such number. Its constants come from
  $\delta_{\alpha/2}$, the fixed phase bound, the interval-length bound lemma
  and $L,c$; covariance enters only through $R_\alpha^k$. The zero-order energy
  is at most $e\|v\|_{1/T}^2$. Put $x = CR_\alpha T^\beta$ and $Y = e^{\nu
  x^{1/\nu}} \ge 1$. Since $\nu > 1$, $\sum_{k\ge0}x^k/(k!)^\nu \le Y$. Summing
  the coefficient bounds in absolute value gives
  \[
    \E\|\mathcal{U}(T;1)v\|_{1/T}^2
    \le(e+Y-1)\|v\|_{1/T}^2\le eY\|v\|_{1/T}^2.
  \]
  Absorbing $\nu C^{1/\nu}$ into the parameter-only constant proves
  \eqref{E:interval-phase-bound}. This proves
  Proposition~\ref{P:interval-coefficients}.
\end{proof}

\subsection{Finite Gaussian models and their terminal kernels}

For $n \ge 1$, define the positive terminal-chain measure
\begin{equation}\label{E:terminal-chain-measure}
  \mathsf{g}_n^{T,x}(\ud\boldsymbol s,\ud\boldsymbol y)
  =\mathbf1_{[0,T]^n_<}(\boldsymbol s)\ud\boldsymbol s
  \prod_{j=1}^n
  G_{s_{j+1}-s_j}(y_{j+1}-\ud y_j),
  \qquad (s_{n+1},y_{n+1})=(T,x).
\end{equation}
The spatial product is read backwards from $x$. Explicitly, it is the
pushforward of $\prod_{j=1}^nG_{s_{j+1}-s_j}(\ud w_j)$ under $y_j =
x-\sum_{l=j}^n w_l$. This defines the measure also in dimension three. In
dimensions one and two its density is the ordered kernel $g_n$
in~\eqref{E:expansion}. Set $\mathsf{g}_0^{T,x}$ equal to unit mass on the empty
tuple. This definition introduces no kernel symmetrization.

\LeanMarker{L:terminal-wick}
\begin{lemma}[The finite terminal-chain Wick formula]\label{L:terminal-wick}
  Let $\mathcal{W}(s,y)$ be a centered real Gaussian field on $[0,T]\times\R^d$,
  $d = 1,2,3$, with covariance $\mathcal{K}(s,y;r,z)$. Suppose that $V =
  \sup_{s,y}\E|\mathcal{W}(s,y)|^2 < \infty$, and that the field is jointly
  continuous in mean square on each open cell of a finite partition of $[0,T]$,
  times $\R^d$. Define $Y_0 = 1$ and
  \[
    Y_n(T,x)=\int\prod_{j=1}^n\mathcal{W}(s_j,y_j)
    \,\mathsf{g}_n^{T,x}(\ud\boldsymbol s,\ud\boldsymbol y).
  \]
  These Bochner integrals exist in $L^2(\Omega)$ and satisfy
  \begin{equation}\label{E:terminal-term-majorant}
    \|Y_n(T,x)\|_2
    \le\sqrt{(2n-1)!!}\,V^{n/2}\frac{T^{2n}}{(2n)!},
    \qquad n\ge1.
  \end{equation}
  For $m+n = 2k$, concatenate $\boldsymbol\tau =
  (s_1,\ldots,s_m,r_1,\ldots,r_n)$ and $\boldsymbol\eta =
  (y_1,\ldots,y_m,z_1,\ldots,z_n)$. Then
  \begin{equation}\label{E:terminal-wick-formula}
    \E[Y_m(T,x)Y_n(T,x)]
    =\sum_{D\in\Pi_k}\mathcal{T}_{m,n,D}^{\mathcal{K}}(T,x),
  \end{equation}
  where, with both chains terminating at $(T,x)$,
  \begin{equation}\label{E:terminal-pairing-integral}
    \begin{split}
      \mathcal{T}_{m,n,D}^{\mathcal{K}}
      ={}&\int_{[0,T]^m_<\times[0,T]^n_<}
      \ud\boldsymbol s\,\ud\boldsymbol r
      \int_{(\R^d)^m\times(\R^d)^n}\\
      &\quad\times\prod_{j=1}^m
      G_{s_{j+1}-s_j}(y_{j+1}-\ud y_j)
      \prod_{j=1}^n
      G_{r_{j+1}-r_j}(z_{j+1}-\ud z_j)\\
      &\quad\times\prod_{\{a,b\}\in D}
      \mathcal{K}(\tau_a,\eta_a;\tau_b,\eta_b).
    \end{split}
  \end{equation}
  Each pairing has coefficient one. If $m+n$ is odd, the expectation is zero.
  Empty chains and the case $m = n = 0$ have the usual empty-product meaning. If
  $\mathcal{K} \ge 0$, every integral above is nonnegative. Consequently,
  entrywise comparison $0 \le \mathcal{K}_1 \le \mathcal{K}_2$ compares all
  finite partial-order second moments in the same direction. The series
  $\sum_{n\ge0}r^nY_n$ converges absolutely in $L^2$ for every real or complex
  $r$.
\end{lemma}

\begin{proof}
  We establish the finite chain mass and Gaussian product integrability first.
  We then derive the labelled Wick formula and check symmetrization and Green
  smoothing. Positivity compares finite moments, while the explicit ratio test
  sums every fixed regularized series.

  \medskip
  \noindent\textit{Step 1.~}We prove Bochner integrability of each terminal
  chain.
  The mass of $G_h$ is $h$. If $M_n(T)$ is the total mass of
  $\mathsf{g}_n^{T,x}$, integration at the last vertex gives $M_0 = 1$ and
  $M_n(T) = \int_0^T(T-s)M_{n-1}(s)\ud s$. Induction therefore yields $M_n(T) =
  T^{2n}/(2n)!$. Gaussian moments and H\"older's inequality give
  \[
    \left\|\prod_{j=1}^n\mathcal{W}(s_j,y_j)\right\|_2
    \le\prod_{j=1}^n\|\mathcal{W}(s_j,y_j)\|_{2n}
    \le\sqrt{(2n-1)!!}\,V^{n/2}.
  \]
  Mean-square continuity of Gaussian variables implies their continuity in every
  fixed finite $L^p$. Hence, the product is strongly measurable in $L^2$ on the
  stated pieces. The finite kernel mass and the displayed bound prove Bochner
  integrability and~\eqref{E:terminal-term-majorant}.

  \medskip
  \noindent\textit{Step 2.~}We derive the coefficient-one Wick formula.
  Fubini now passes expectation through the product of the two finite kernel
  measures. For the $m+n$ Gaussian variables at their labelled vertices,
  differentiate
  \[
    \E\exp\left(\sum_i t_i\mathcal{W}(\tau_i,\eta_i)\right)
    =\exp\left(\frac12\sum_{i,j}t_i t_j
    \mathcal{K}(\tau_i,\eta_i;\tau_j,\eta_j)\right)
  \]
  once in each $t_i$ at zero. Odd total degree gives zero. At degree $2k$ the
  result is the sum over unordered pair partitions, with one covariance product
  per pairing. Orientations and orderings of the $k$ pairs cancel the
  exponential's $2^k k!$ denominator. Integration
  gives~\eqref{E:terminal-wick-formula} with precisely the domains and
  coefficients in~\eqref{E:terminal-pairing-integral}.

  \medskip
  \noindent\textit{Step 3.~}We retain the ordered convention and the
  Green-width limit.
  Symmetrizing a chain measure does not change its stochastic integral: the
  product of its $n$ scalar field values is invariant under relabeling. Thus,
  averaging the $n!$ pushforwards of $\mathsf{g}_n^{T,x}$ with coefficient
  $1/n!$ gives the same $Y_n$. There is no additional $n!$ factor
  in~\eqref{E:terminal-wick-formula}, or $m!n!$ in a two-chain product. This
  identifies the ordinary ordered regularized integrals and their symmetrized
  convention in~\eqref{E:expansion}.

  In dimension three one may first replace $G$ by $p_\eta*G$. Those positive
  measures have the same masses and converge weakly to $G$. Off the finitely
  many time cut points, all finite Gaussian moments are bounded continuous
  functions of the spatial variables. Dominated convergence for their finite
  Wick sums proves the same limit for the mixed second moments with two Green
  widths. It gives $L^2$ convergence to the Bochner integral above, exactly the
  fixed-noise Green convention used in the source.

  \medskip
  \noindent\textit{Step 4.~}We compare positive diagrams and sum the
  regularized orders.
  If $\mathcal{K} \ge 0$, positivity of the Green measures makes each covariance
  product nonnegative. Comparing them proves the stated finite-moment
  monotonicity, including self-pairs and unequal chains. Finally, the ratio of
  successive right sides of~\eqref{E:terminal-term-majorant}, with an additional
  factor $|r|^n$, is $|r|\sqrt V T^2/((2n+2)\sqrt{2n+1})$, which tends to zero.
  The zero-variance case is immediate. This proves Lemma~\ref{L:terminal-wick}.
\end{proof}

\subsection{Removing the regularization}

Proposition~\ref{P:interval-coefficients} holds uniformly in the smoothing
$\epsilon > 0$, but it concerns regularized fields only. We still have to pass
to the unregularized solution.

\LeanMarker{L:uniform-series-tail}
\begin{lemma}[Uniform series tails from scaled moments]\label{L:uniform-series-tail}
  Fix $(T,x)$. For each positive regularization pair $(\epsilon,\delta)$, let
  $Y_n^{\epsilon,\delta}$ be the terms of Lemma~\ref{L:terminal-wick} for a
  field with nonnegative covariance. Write $u^{\epsilon,\delta}(r) =
  \sum_{n\ge0}r^nY_n^{\epsilon,\delta}$. Suppose that, for some $r > 1$ and
  finite $B_r$,
  \[
    \sup_{\epsilon,\delta>0}
    \E|u^{\epsilon,\delta}(r)|^2\le B_r,
    \qquad
    Y_n^{\epsilon,\delta}\longrightarrow X_n
    \text{ in }L^2\text{ for each fixed }n.
  \]
  Then $u = \sum_{n\ge0}X_n$ converges absolutely in $L^2$,
  $u^{\epsilon,\delta}(1) \to u$ in $L^2$, and
  \begin{equation}\label{E:uniform-series-tail}
    \|Y_n^{\epsilon,\delta}\|_2,\ \|X_n\|_2
    \le B_r^{1/2}r^{-n},\qquad
    \left\|\sum_{n>N}Y_n^{\epsilon,\delta}\right\|_2,
    \left\|\sum_{n>N}X_n\right\|_2
    \le\frac{B_r^{1/2}r^{-(N+1)}}{1-r^{-1}}.
  \end{equation}
  If also $\E|u^{\epsilon,\delta}(1)|^2 \le B_1$ uniformly, then $\E|u|^2 \le
  B_1$.
\end{lemma}

\begin{proof}
  Lemma~\ref{L:terminal-wick} gives absolute $L^2$ convergence of every
  regularized series and nonnegative mixed moments. Its scaled second moment
  therefore contains the nonnegative term $r^{2n}\E|Y_n^{\epsilon,\delta}|^2$,
  so $r^{2n}\E|Y_n^{\epsilon,\delta}|^2 \le B_r$. Taking the fixed-order $L^2$
  limit proves the same bound for $X_n$. The geometric series and the triangle
  inequality give all estimates in~\eqref{E:uniform-series-tail}, including
  absolute convergence.

  For each $N$, write both full series as their first $N+1$ terms and their
  tails. Their distance is at most
  \[
    \frac{2B_r^{1/2}r^{-(N+1)}}{1-r^{-1}}
    +\left\|\sum_{n=0}^N(Y_n^{\epsilon,\delta}-X_n)\right\|_2.
  \]
  First let $\epsilon,\delta\downarrow0$ at fixed $N$, and then let $N \to
  \infty$. This proves convergence of the full regularized solution; it was not
  assumed. Continuity of the $L^2$ norm passes the bound $B_1$ to $u$. This
  proves Lemma~\ref{L:uniform-series-tail}.
\end{proof}

\subsection{Temporal comparison and the full upper bound}

The use of finitely many shifted dyadic grids is classical. For the covering
principle and its antecedents, see
Lerner~\cite{lerner:source2013:simple}*{Section~2.2, p.~3162}.
We give the construction, nesting check and covariance comparison needed for
the time intervals below.

\begin{theorem}[Full upper scale]\label{T:full-upper}
  Under Hypothesis~\ref{H:homog}, let $u$ be the canonical expansion
  solution~\eqref{E:expansion}, constructed as a Stratonovich mild solution in
  Proposition~\ref{P:mild-construction}. Put
  \[
    R_\alpha=\frac{\alpha\pi\mu(B_1)}{2\sin(\pi\alpha/2)}.
  \]
  For every admissible $d = 1,2,3$, there is a finite $C = C(\alpha,\alpha_0)$
  such that
  \begin{equation}\label{E:full-upper-bound}
    \E u^2(T,x)
    \le e(T+2)^d
    \exp\{C R_\alpha^{1/(3-\alpha)}
    T^{(4-\alpha-\alpha_0)/(3-\alpha)}\},\qquad T\ge1.
  \end{equation}
  In particular, this proves Theorem~\ref{T:intro-upper}.
\end{theorem}

\begin{proof}
  We dominate the temporal covariance by finite-depth shifted interval models
  with a controlled length inventory. We transfer their phase bound to
  constant data on a fixed observation ball, compare the positive terminal
  Wick integrals, and remove noise regularization using a separate uniform
  series-tail bound. All identifications precede the large-time conclusion.

  \medskip
  \noindent\textit{Step 1.~}We construct the dyadic temporal majorant and its
  length inventory.
  Zero noise is trivial. Choose $H_0$ with $6T \le H_0 < 12T$. At lengths $h_j =
  H_0 2^{-j}$, use the three grids with endpoints $h_j(n+(-1)^j t)$, $n \in
  \mathbb{Z}$, $t = 0,1/3,2/3$. They are nested: the difference between the
  parent and child shifts in child units is an integer, since $3t \in
  \mathbb{Z}$. For $s \ne r$ in $[0,T]$, choose a scale $3|s-r| < h_j \le
  6|s-r|$. Boundaries of the three grids at this scale have mutual spacing
  $h_j/3 > |s-r|$. Hence, at least one grid interval contains both points,
  including their endpoints after choosing a boundary-free grid. It follows that
  \begin{equation}\label{E:dyadic-time-majorant}
    |s-r|^{-\alpha_0}
    \le6^{\alpha_0}\sum_{\text{three grids},\,|I|\le H_0}
    |I|^{-\alpha_0}\mathbf1_I(s)\mathbf1_I(r)=:K_\infty(s,r).
  \end{equation}
  For context, the sum is at most $C_{\alpha_0}|s-r|^{-\alpha_0}$ off the
  diagonal: only lengths at least $|s-r|$ contribute, with at most one interval
  per grid and length. Endpoint values do not affect the time integrals.

  Keep intervals meeting $[0,T]$ through depth $J$, and use independent static
  smoothed fields with $\sigma_I^2 = 6^{\alpha_0}|I|^{-\alpha_0}$. Their time
  covariance is the corresponding partial sum $K_J$. The decreasing interval
  lengths obey $h_i \le 78T/i$. Indeed, at length $h_j$ there are at most
  $3(T/h_j+2)$ intervals. Summing levels $0,\ldots,j$ bounds their cumulative
  count by
  \[
    6T/h_j+6(j+1)\le6T/h_j+6\,2^j\le78T/h_j,
  \]
  since $H_0 < 12T$ and $2^j = H_0/h_j$. Each index at level $j$ is at most that
  cumulative count, which gives the asserted decreasing-length bound. Thus, one
  may take $L = 78$ and $c = 6^{\alpha_0}$ in
  Proposition~\ref{P:interval-coefficients}.
  Proposition~\ref{P:interval-coefficients} applies uniformly in $J,\epsilon$.

  \medskip
  \noindent\textit{Step 2.~}We localize constant initial data on the full
  backward cone.
  For the kept intervals, $W_{J,\epsilon}(s,y)=\sum_{I\in\mathcal
  F_J}\sigma_I1_I(s)W_I^\epsilon(y)$, and $u_{J,\epsilon}$ is the series of
  its terminal terms from Lemma~\ref{L:terminal-wick}. For the compact data
  $(f_T,0)$ chosen below, $u_{f_T,J,\epsilon}$ denotes the position component
  of the phase evolution in Proposition~\ref{P:interval-coefficients}.
  To pass from phase data to constant initial data, choose $f_T \in C_c^\infty$
  with $0 \le f_T \le 1$, equal to one on $B_{T+1}$ and supported in $B_{T+2}$,
  and set the velocity to zero. For $|y|+s < T+1$, $G_s*f_T(y) = s$, so
  $C_sf_T(y) = 1$. The support of $G_t$ is contained in $B_t$. Along a chain
  ending at $x \in B_1$ at time $T$, each earlier point $(s,y)$ satisfies $|y|+s
  \le |x|+T < T+1$. Thus, its free initial-data factor is one. Induction on
  Duhamel order identifies the position generated from $(f_T,0)$ at time $T$ on
  $B_1$ with the constant-data model $u_{J,\epsilon}$. This identity can be
  taken in Gaussian $L^2$ term by term. At fixed $J,\epsilon$, each finite-order
  terminal integral is continuous in $x$ as an $L^2(\Omega)$ vector: write its
  Green variables as displacements from $x$, and use continuity of
  $\Gamma_\epsilon$ in its mixed Wick moments. Bounded covariance products and
  the finite terminal Green mass permit dominated convergence. These measurable
  representatives identify the $L^2(\Omega;L^2_x)$ position on $B_1$ with the
  same pointwise constant-data terms. The summable bound below passes this
  identity to their series. At finite $J,\epsilon$, pointwise convergence for
  constant data follows from the bound
  $\sqrt{(2n-1)!!}\,v_{\max}^{n/2}T^{2n}/(2n)!$ on the $n$th position term. The
  model is spatially stationary. Consequently,
  \[
    |B_1|\E|u_{J,\epsilon}(T,0)|^2
    \le\E\|u_{f_T,J,\epsilon}(T)\|_2^2
    \le e\,|B_{T+2}|\exp\{C R_\alpha^{1/\nu}T^{\beta/\nu}\}.
  \]

  The final inequality above is~\eqref{E:interval-phase-bound} with
  $v=(f_T,0)$: the position $L^2$ norm is at most its phase norm, and
  $\|(f_T,0)\|_{1/T}^2=\|f_T\|_2^2\le|B_{T+2}|$ because the velocity is zero
  and $0\le f_T\le1$.
  \medskip
  \noindent\textit{Step 3.~}We identify the finite interval-model Wick sums
  before comparison.
  We now identify the interval-model Wick sums before comparing them. Its
  centered Gaussian noise has covariance
  \[
    \E[W_{J,\epsilon}(s,y)W_{J,\epsilon}(r,z)]
    =\sum_{I\in\mathcal{F}_J}\sigma_I^2
    \mathbf1_I(s)\mathbf1_I(r)\Gamma_\epsilon(y-z)
    =K_J(s,r)\Gamma_\epsilon(y-z).
  \]
  Lemma~\ref{L:terminal-wick} applies: at fixed $J,\epsilon$ the field has
  bounded variance, is mean-square continuous in space and is piecewise
  continuous in time. Its exact terminal kernels are
  \eqref{E:terminal-chain-measure}, and every term
  of~\eqref{E:terminal-pairing-integral} has covariance factor
  $K_J(s_j,s_l)\Gamma_\epsilon(y_j-y_l)$, with coefficient one. This comparison
  is made before any time reversal; $K_J$ need not be stationary in time. The
  position kernels and constant initial data are nonnegative. Thus, every such
  integral, including self-pairings and unequal chain lengths, is nonnegative.
  The pointwise $L^2$ convergence proved above identifies the full moment with
  their nonnegative sum. In particular, the moment of every finite partial-order
  sum is at most the full moment, even though samplewise partial sums need not
  increase.

  \medskip
  \noindent\textit{Step 4.~}We compare the source model through positive
  terminal diagrams.
  From the current author version
  \cite{chen:author2026:time-dependency}*{(2.23)--(2.24), p.~13}, the common
  noise regularization has covariance
  $\gamma_\delta^0(s-r)\Gamma_\epsilon(y-z)$, where
  \[
    0\le\gamma_\delta^0(t)
    =\Gamma(\alpha_0)^{-1}\int_0^{\delta^{-1}}
    e^{-q|t|}q^{\alpha_0-1}\ud q
    \le |t|^{-\alpha_0}.
  \]
  Here $\Gamma(\alpha_0)$ is the Euler Gamma function. At fixed positive
  $\epsilon,\delta$, first take the source's Green-width limit. The bounded
  continuous covariances, weak convergence of the positive Green measures and
  the dominated finite Wick calculation in Lemma~\ref{L:terminal-wick} identify
  its classical measure-kernel Duhamel terms. For the source's continuous field
  modification, these Bochner terms equal the pathwise smoothed integrals:
  writing $\Phi_n$ for the product of field values, the Gaussian bound in that
  lemma gives $\int\|\Phi_n\|_2\ud\mathsf{g}_n^{T,x} < \infty$. Tonelli gives
  pathwise integrability and Minkowski puts the integral in $L^2$. Fubini
  against every $L^2$ test variable and Bochner duality identify the two
  versions. The same argument applies before the Green-width limit in dimension
  three.

  Write $U_N^{\epsilon,\delta}$ for the sum of source terms of orders
  $0,\ldots,N$ after that Green limit. Write $V_{N,J}^\epsilon$ for the same
  partial-order sum of the finite interval model. Both moments are finite sums
  over exactly the same terminal-chain kernels and pairings; only the covariance
  factors differ. By~\eqref{E:dyadic-time-majorant} each source covariance
  factor is bounded by its $K_\infty$ counterpart. Since $K_J\uparrow K_\infty$,
  nonnegative Tonelli at fixed $N$ yields
  \begin{equation}\label{E:interval-source-comparison}
    \E|U_N^{\epsilon,\delta}(T,x)|^2
    \le\lim_{J\to\infty}\E|V_{N,J}^\epsilon(T,x)|^2
    \le\sup_J\E|u_{J,\epsilon}(T,x)|^2\le B(T),
  \end{equation}
  where $B(T)$ is the right side of~\eqref{E:full-upper-bound}. The middle
  inequality uses the positive pairing sum just identified. No Gaussian field of
  covariance $K_\infty\Gamma_\epsilon$ is constructed; $K_\infty$ is used only
  inside these monotone limits of positive integrals. This is entrywise
  comparison of Wick integrals, not covariance order for arbitrary initial
  phases or a coupling of random solutions.

  \medskip
  \noindent\textit{Step 5.~}We obtain a scaled source bound that controls all
  order tails.
  The same finite-model argument works after multiplying the source noise by any
  real amplitude $a > 0$. This replaces $\Gamma$ by $a^2\Gamma$ and $R_\alpha$
  by $a^2R_\alpha$, without changing its homogeneity degree or any
  parameter-only constant. Thus,
  \begin{equation}\label{E:regularized-scaled-moment}
    \E|u^{\epsilon,\delta}(a;T,x)|^2
    \le B(T;a):=e(T+2)^d
    \exp\{C(a^2R_\alpha)^{1/\nu}T^\kappa\}.
  \end{equation}
  Here the full regularized series exists by Lemma~\ref{L:terminal-wick}; it is
  the ordinary Duhamel series with its order-$n$ term multiplied by $a^n$.
  Taking the finite partial-order cutoff to infinity
  in~\eqref{E:interval-source-comparison} is justified at positive
  regularization by that lemma, not by a theorem about the original unsmoothed
  solution.

  \medskip
  \noindent\textit{Step 6.~}We identify fixed orders and remove the full
  source regularization.
  Lemma~\ref{L:joint-fixed-order} supplies the fixed-order construction: each
  $Y_n^{\epsilon,\delta}$ converges in $L^2$ to $X_n = S_n(g_n)$ as
  $\epsilon,\delta\downarrow0$, with arbitrary relative rates. The bounded-noise
  Green limit and ordered-kernel convention needed to match these terms are
  already verified in Lemma~\ref{L:terminal-wick}.

  Apply Lemma~\ref{L:uniform-series-tail} with $r = 2$ and $B_r = B(T;2)$. It
  proves absolute $L^2$ convergence of the series defining $u$
  in~\eqref{E:expansion} and convergence of the full regularized solution to
  that series. Taking the limit in~\eqref{E:regularized-scaled-moment} at $a =
  1$ proves the required bound. This series-limit argument is independent of
  Proposition~\ref{P:mixed}, Lemma~\ref{L:crude-corrected}, and the source's
  Section~5 convergence proof. No first-moment asymptotic or printed
  normalization is used. The fixed-order construction just used is proved here
  from the universal pairing bound and the homogeneous norm estimates; neither
  the source's Theorem~4.1 nor its Lemma~4.4 is imported.
  Proposition~\ref{P:mild-construction} separately proves that this same series
  satisfies the source's Definition~2.1: its onefold Stratonovich integral
  equals $u-1$. That construction uses the corrected mixed-vertex crude bound
  for its product and order limits. This is a separate dependency from the
  sharp-scale moment estimate just proved; the source's Section~5(ii) proof is
  not imported. The crude constants used for this identification may depend on
  the fixed covariance; they do not enter~\eqref{E:full-upper-bound}. No
  uniqueness assertion is made.

  \medskip
  \noindent\textit{Step 7.~}We record the distinct cutoff limits and the
  parameter-only constant.
  \emph{Cutoff dependence and the order of limits.} The preceding argument uses
  the following distinct statements.
  \begin{enumerate}[(i)]
    \item At fixed total order $2k$, the coefficient estimates are uniform
      in the Gaussian-mode and spatial-band cutoffs, in the finite interval
      family, and in positive spatial smoothing. The local Hermite cutoff $N_i =
      2k_i$ preserves each selected coefficient exactly; it is an algebraic
      choice, not a degree limit. Lemma~\ref{L:fixed-split} and the
      interval-length bound contain no smoothing-dependent variance.
    \item At fixed finite interval family and $\epsilon > 0$, mode and band
      limits identify these coefficients with the Duhamel words. The
      entire-series majorant in Proposition~\ref{P:interval-coefficients} may
      depend on this family and $\epsilon$, through $v_{\max}$. It is used for
      identification, not as a uniform moment estimate. The uniform coefficient
      bound passes to these limits. Localization and positive covariance
      comparison are applied only after removing these cutoffs; positivity of a
      finite-mode covariance is not assumed.
    \item For the resulting finite-depth interval model, the moment bound
      is uniform in $J,\epsilon$. Fix positive source widths $\epsilon,\delta$
      and finite partial order $N$. The Green-width limit is taken first. The
      subsequent $J \to \infty$ comparison
      in~\eqref{E:interval-source-comparison} is monotone convergence of finite
      positive Wick integrals. At these fixed noise widths, $N \to \infty$ is
      justified by Lemma~\ref{L:terminal-wick}, yielding
      \eqref{E:regularized-scaled-moment}. There is no infinite-depth Gaussian
      field or simultaneous mode-depth limit in this step.
    \item At each fixed $T,x,\Gamma$ and order $n$,
      Lemma~\ref{L:joint-fixed-order} gives qualitative $L^2$ convergence as
      $\epsilon,\delta\downarrow0$, with arbitrary relative rates. The bound
      $B(T;2)$ gives tails uniform in these two widths.
      Lemma~\ref{L:uniform-series-tail} first uses their limit at fixed partial
      order, then sends that order to infinity. The final moment bound comes
      from $B(T;1)$, not the larger tail bound $B(T;2)$.
    \item All identifications above are for fixed $T,x,\Gamma$; no
      convergence rate uniform in $T$ or over covariance shapes is claimed. They
      prove~\eqref{E:full-upper-bound} for every fixed $T \ge 1$ with the same
      parameter-only constant. Only after this estimate is proved is $T \to
      \infty$ used for the growth-rate conclusion.
  \end{enumerate}
  The constants come from $L = 78$, $c = 6^{\alpha_0}$ and the parameter-only
  constants of Proposition~\ref{P:interval-coefficients}; covariance enters only
  through $R_\alpha$. The volume ratio is $|B_{T+2}|/|B_1| = (T+2)^d$. This
  proves Theorem~\ref{T:full-upper}.
\end{proof}

\begin{remark}[Claim and calibration boundary]\label{R:full-upper-boundary}
  Theorem~\ref{T:full-upper} asserts a nonsharp upper growth rate. It does not
  identify the sharp constant or require a correction of the printed
  first-moment normalization as an input. Its exponent tends to
  $(4-\alpha)/(3-\alpha)$ as $\alpha_0\downarrow0$. Cauchy--Schwarz and the
  unconditional Jensen inequality remain compatible with the nonsharp upper
  constant; any numeric $L_*$ calibration retains the source-normalization
  qualification stated earlier. The Skorohod exponent is compared, not used to
  identify constants: both Gaussian creation and annihilation were retained
  here. The intrinsic identification in Theorem~\ref{T:intro-second-moment} is proved in
  Theorem~\ref{T:intrinsic-second-moment}, using a separate argument.
\end{remark}

\paragraph{What we now have.}
For every covariance in Hypothesis~\ref{H:homog}, the second moment is at most
$\exp(CT^\kappa)$ times a polynomial in $T$ (Theorem~\ref{T:full-upper}). The
constant $C$ is not claimed sharp; the next eight chapters identify the exact
asymptotic constant.

\chapter{Covariance estimates and Gaussian approximation}\label{S:book-gaussian-foundations}

A rough Gaussian potential is not a bounded multiplication operator. So what
does it do to a wave? We need estimates that hold three things together:
the covariance, wave propagation and spatial localization. There is one
trap. An operator and its adjoint need two variance estimates, and one
direction alone does not give the later Gaussian bounds.

\begin{figure}[htbp]
  \centering
  \input{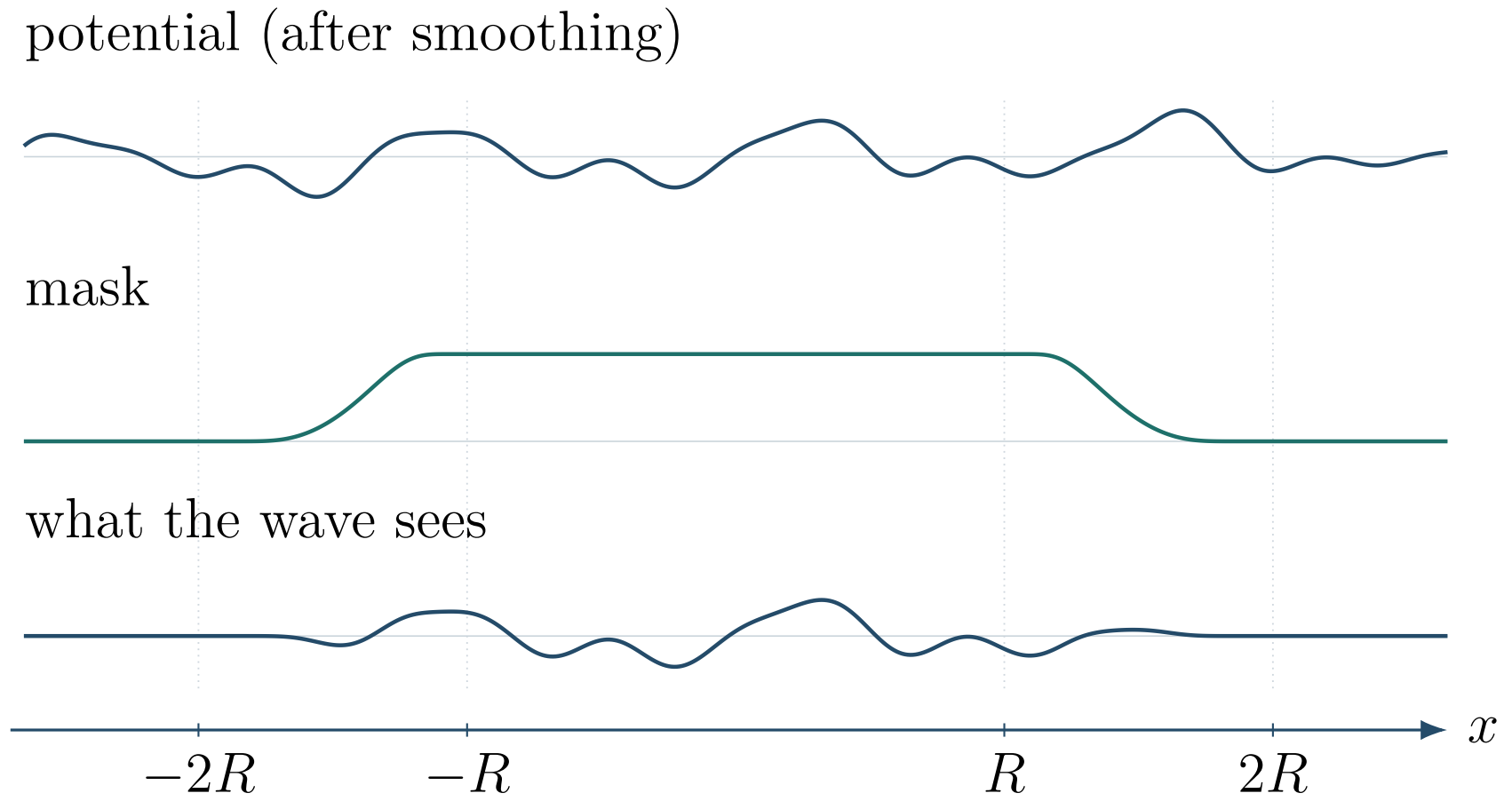}
  \caption[A physical mask as a window on the potential]{A physical mask,
    equal to one on a ball and zero outside a larger ball, and its product
    with the potential. The mask multiplies the potential in space. The
    covariance of the product has a mask factor at each spatial argument. The
    original noise and the covariance norm used for control costs remain
    fixed. A schematic in one space dimension: the unsmoothed potential is a
    distribution, not a curve.}
  \label{Fig:mask-window}
\end{figure}

Figure~\ref{Fig:mask-window} shows a window on the potential: a physical
mask is a smooth function of compact support that multiplies it in
space. Section~\ref{S:masked-covariance} bounds the covariance behind
it. Section~\ref{S:wave-variance} carries the bounds to the full wave
phase. Section~\ref{S:compact-gaussian} adds finite Gaussian matrix estimates
and the compactness from spatial cutoffs.

The assumptions are those of Hypothesis~\ref{H:homog}. On a
first reading, keep Proposition~\ref{P:wave-two-variances}
and Corollary~\ref{C:finite-gaussian-operators} and let the
proofs wait. The Gaussian operator-series estimates here concern
finite noise sums; fixed-vector Gaussian column limits are also
constructed. Radius-explicit operator-path moment bounds come in
Chapter~\ref{S:book-growing-radius}. Episode~\VideoRef{5}{2jv8nuqgfBY}
of the video series covers this chapter and the next two.

\section{Covariance bounds for physically masked wave fields}
\label{S:masked-covariance}

We establish the covariance estimates needed for the operator construction.
Throughout this section, Hypothesis~\ref{H:homog} holds, $a = \alpha_0$, and
$\Gamma \ne 0$. For zero spatial covariance, the spatial spectral measures and
noise variances vanish; the purely temporal estimates remain valid independently
of $\Gamma$. We use the unitary Fourier transform
\[
  \widehat{f}(\xi) = c_d\int_{\R^d}e^{-i\xi\cdot x}f(x)\ud x,
  \qquad c_d = (2\pi)^{-d/2}.
\]
The spectral convention in Hypothesis~\ref{H:homog} means that, for every
complex Schwartz test $\varphi$,
\begin{equation}\label{E:masked-spectral-convention}
  \int_{\R^d}\varphi(x)\Gamma(\ud x)
  = \int_{\R^d}\left(\int_{\R^d}e^{i\xi\cdot x}
  \varphi(x)\ud x\right)\mu(\ud\xi).
\end{equation}
This is an identity of distribution pairings. We will not interchange two
infinite measures in the associated oscillatory expression.

A positive tempered Radon measure has polynomial mass growth. Indeed, choose a
nonnegative compact smooth test that is at least one on the unit ball, and apply
a defining Schwartz seminorm bound for the distribution to its dilations. Its
mass on a ball of radius $R \ge 1$ is bounded by a power of $R$. Summation over
dyadic shells then gives integrability of $(1+|x|)^{-N}$ for some $N$, and hence
absolute integrability of every Schwartz test. These facts apply to both
$\Gamma$ and $\mu$. Smooth cutoffs converge in the Schwartz topology, so these
integrals agree with the tempered pairings.

\subsection{Translated spectral growth and physical masks}

\LeanMarker{P:masked-spectral-growth}
\begin{proposition}[Translated spectral mass]\label{P:masked-spectral-growth}
  Put
  \[
    M_\Gamma = \int_{\R^d}e^{-|x|^2/4}\Gamma(\ud x),
    \qquad C_\Gamma = \frac{eM_\Gamma}{(2\sqrt{\pi})^d}.
  \]
  Then $0 < M_\Gamma < \infty$, and for $y \in \R^d$ and $r > 0$,
  \begin{align}
    \int_{\R^d}e^{-|\xi-y|^2/r^2}\mu(\ud\xi)
    &= \left(\frac{r}{2\sqrt{\pi}}\right)^d
    \int_{\R^d}e^{-r^2|x|^2/4}e^{-iy\cdot x}\Gamma(\ud x),
    \label{E:masked-gaussian-transform}\\
    0 \le \int_{\R^d}e^{-|\xi-y|^2/r^2}\mu(\ud\xi)
    &\le \frac{r^\alpha M_\Gamma}{(2\sqrt{\pi})^d},
    \qquad \mu(B(y,r)) \le C_\Gamma r^\alpha.
    \label{E:masked-translated-ball}
  \end{align}
  The ball may be open or closed.
\end{proposition}

\begin{proof}
  We first compute the Gaussian transform in the stated Fourier convention. We
  then use positivity and homogeneity to bound it, and compare the Gaussian
  with the indicator of a ball.

  \medskip
  \noindent\textit{Step 1.~}We compute the transform and its constant.
  Differentiating the Gaussian Fourier integral below is justified by its
  integrable first moment, and Gaussian decay makes the integration-by-parts
  boundary terms vanish.
  Temperedness gives finiteness of $M_\Gamma$. Its positivity follows from
  $\Gamma \ne 0$. To check the Fourier constant, set $J(b) = \int_\R
  e^{-Ax^2}e^{ibx}\ud x$ for $A > 0$. Squaring the real Gaussian integral and
  using polar coordinates gives $J(0) = \sqrt{\pi/A}$. Differentiation and
  integration by parts give $J'(b) = -bJ(b)/(2A)$. Therefore, $J(b) =
  \sqrt{\pi/A}e^{-b^2/(4A)}$. Taking products over the coordinates gives
  \[
    \int_{\R^d}e^{-r^2|x|^2/4}e^{i(\xi-y)\cdot x}\ud x
    = \left(\frac{2\sqrt{\pi}}r\right)^d e^{-|\xi-y|^2/r^2}.
  \]
  Apply~\eqref{E:masked-spectral-convention} to the single Schwartz test
  $e^{-r^2|x|^2/4}e^{-iy\cdot x}$ to obtain \eqref{E:masked-gaussian-transform}.

  \medskip
  \noindent\textit{Step 2.~}We use positivity and homogeneity.
  Its left side is nonnegative. Taking the absolute value on the right removes
  the complex phase because $\Gamma$ is nonnegative. The pushforward of $\Gamma$
  by $x \mapsto rx$ assigns mass $r^{\alpha-d}\Gamma(A)$ to a Borel set $A$.
  Thus,
  \[
    \int e^{-r^2|x|^2/4}\Gamma(\ud x) = r^{\alpha-d}M_\Gamma.
  \]
  Multiplying the last equality by the prefactor $r^d/(2\sqrt{\pi})^d$ gives
  $r^\alpha M_\Gamma/(2\sqrt{\pi})^d$.
  This proves the Gaussian bound. On either ball $B(y,r)$ the integrand
  $e^{-|\xi-y|^2/r^2}$ is at least $e^{-1}$, proving the mass bound. This proves
  Proposition~\ref{P:masked-spectral-growth}.
\end{proof}

Choose and fix a real nonzero Schwartz function $\zeta$. For $R > 0$ put
$\zeta_R(x) = \zeta(x/R)$, $k_R = \widehat{\zeta_R}$ and $L =
\|\widehat{\zeta}\|_1$. Then $k_R(u) = R^d\widehat{\zeta}(Ru)$ and $\|k_R\|_1 =
L > 0$. Define
\begin{equation}\label{E:masked-auxiliary-measure}
  \nu_R(\ud\eta) = \rho_R(\eta)\ud\eta,
  \qquad \rho_R(\eta) = c_d^2L\int |k_R(\xi-\eta)|\mu(\ud\xi).
\end{equation}

\LeanMarker{P:masked-auxiliary-norm}
\begin{proposition}[Auxiliary mask norm]\label{P:masked-auxiliary-norm}
  The measure $\nu_R$ is even, tempered and Radon, with
  \begin{equation}\label{E:masked-auxiliary-growth}
    \nu_R(B(y,r)) \le c_d^2L^2C_\Gamma r^\alpha
  \end{equation}
  for all $y \in \R^d$ and $r > 0$, uniformly in $R > 0$. If $H$ is a separable
  complex Hilbert space and $F$ is an $H$-valued Schwartz function, then
  \begin{equation}\label{E:masked-hilbert-fourier}
    \int \|\widehat{\zeta_RF}(\xi)\|_H^2\mu(\ud\xi)
    \le \int \|\widehat{F}(\eta)\|_H^2\nu_R(\ud\eta).
  \end{equation}
  For a real $H$, the same assertion holds after complexification.
\end{proposition}

\begin{proof}
  We first establish the properties and translated growth bound of the
  auxiliary measure. We then prove the Fourier norm inequality by a weighted
  Cauchy--Schwarz estimate.

  \medskip
  \noindent\textit{Step 1.~}We examine the density and its ball masses.
  For fixed $R$, translating the argument of a Schwartz function over a
  compact set preserves a common rapid-decay bound. Choose the decay power
  large enough to be integrable against the polynomially growing measure
  $\mu$.
  Schwartz decay of $k_R$ and polynomial growth of $\mu$ show that
  $\rho_R(\eta)$ is finite at every $\eta$. On compact sets of $\eta$, the
  integrands have a common integrable bound. Dominated convergence therefore
  makes $\rho_R$ continuous, and $\nu_R$ is Radon. Since $\zeta$ is real,
  $k_R(-u) = \overline{k_R(u)}$. Evenness of $\mu$ gives evenness of $\rho_R$.
  To bound a ball mass, write the density as its nonnegative integral and
  exchange the two integrations. For fixed $\xi$, set $u = \xi-\eta$. The
  condition $\eta \in B(y,r)$ becomes $\xi \in B(y+u,r)$ when the integrations
  are exchanged again.
  Tonelli's theorem gives
  \[
    \nu_R(B(y,r)) = c_d^2L\int |k_R(u)|\mu(B(y+u,r))\ud u.
  \]
  Proposition~\ref{P:masked-spectral-growth} and $\|k_R\|_1 = L$ give
  \eqref{E:masked-auxiliary-growth}. Polynomial growth on centered balls proves
  temperedness by the same dyadic argument as above.

  \medskip
  \noindent\textit{Step 2.~}We justify the Hilbert-valued convolution
  identity.
  The finite-dimensional orthogonal projections used below have norm at most
  one, so their Fourier integrals have the same Schwartz majorants.
  Bochner Fourier integration gives
  \[
    \widehat{\zeta_RF}(\xi)
    = c_d\int k_R(\xi-\eta)\widehat{F}(\eta)\ud\eta.
  \]
  The factor $c_d$ is the product formula for the unitary Fourier transform.
  We justify the same identity for vector-valued functions by approximation,
  so that the convolution is a Bochner integral.
  For scalar Schwartz functions this follows by Fourier inversion and absolute
  integrability. It follows next for finite-rank Hilbert-valued functions. For
  general $F$, take projections onto the successive spans of a countable
  orthonormal basis of its separable range closure. They converge pointwise in
  norm. The Schwartz bounds on $F$ dominate the Fourier integrals, proving the
  identity by dominated convergence. In a finite-dimensional space the full
  identity projection suffices; the zero space is immediate. The triangle
  inequality and weighted Cauchy--Schwarz yield
  \[
    \|\widehat{\zeta_RF}(\xi)\|_H^2
    \le c_d^2L\int |k_R(\xi-\eta)|\|\widehat{F}(\eta)\|_H^2\ud\eta.
  \]

  \medskip
  \noindent\textit{Step 3.~}We estimate the convolution in norm and integrate.
  Integrating against $\mu$ and applying Tonelli's theorem proves
  \eqref{E:masked-hilbert-fourier}. Its right side is finite by Schwartz decay
  and temperedness of $\nu_R$. This proves
  Proposition~\ref{P:masked-auxiliary-norm}.
\end{proof}

\subsection{Temporal completion and the original noise}

The noise in~\eqref{E:cov} is also correlated in time, with kernel
$|t-s|^{-a}$. We now complete the temporal tests in that form.

\LeanMarker{L:masked-time-completion}
\begin{lemma}[Temporal covariance completion]\label{L:masked-time-completion}
  Let $T > 0$. The form
  \[
    \langle f,g\rangle_a
    = \int_0^T\int_0^T f(t)g(s)|t-s|^{-a}\ud t\ud s
  \]
  is positive semidefinite on real compact smooth tests. Its null-quotient
  completion $H_a$ is separable. The inclusion of tests extends continuously
  from $L^2(0,T)$, with
  \begin{equation}\label{E:masked-time-schur}
    \LeanEqMarker{E:masked-time-schur}
    \|f\|_a^2 \le \frac{2T^{1-a}}{1-a}\|f\|_2^2.
  \end{equation}
  In its complexification the norm is given by the same double integral with
  $f(t)\overline{f(s)}$. In particular, bounded measurable amplitudes on a
  subinterval, extended by zero, belong to this completion.
\end{lemma}

\begin{proof}
  We prove positivity by expressing the temporal kernel as a positive mixture
  of square forms. A uniform row-integral bound then gives the continuous
  inclusion and the completion statements.

  \medskip
  \noindent\textit{Step 1.~}We prove positivity of the covariance form.
  For fixed $t$, the kernel row integral is
  $[t^{1-a}+(T-t)^{1-a}]/(1-a)\le2T^{1-a}/(1-a)$; symmetry gives the same
  bound for columns.
  Write $\Gamma_E$ for Euler's Gamma function. Away from the null diagonal,
  \[
    |t-s|^{-a} = \frac1{\Gamma_E(a)}\int_0^\infty
    \lambda^{a-1}e^{-\lambda|t-s|}\ud\lambda.
  \]
  This follows by substitution in the defining Gamma integral. Moreover,
  \[
    e^{-\lambda|t-s|}
    = \int_{-\infty}^{\min(t,s)}2\lambda
    e^{-\lambda(t-u)}e^{-\lambda(s-u)}\ud u.
  \]
  For each $\lambda > 0$ its quadratic form is an integral of a square. All
  exchanges for compact smooth tests are absolutely integrable because the
  kernel row integral is at most $2T^{1-a}/(1-a)$. The same bound and $2|xy| \le
  |x|^2+|y|^2$ prove~\eqref{E:masked-time-schur}. A countable $L^2$-dense family
  of compact smooth tests is consequently dense in $H_a$, proving separability.

  \medskip
  \noindent\textit{Step 2.~}We explain the row bound and the completion.
  The complex quadratic form is the sum of those of the real and imaginary
  parts, since symmetry cancels the imaginary cross term. Thus, the same bound
  holds after complexification. Smooth density in $L^2$, together with the Schur
  bound on the absolute bilinear form, extends the inclusion and the integral
  identity to $L^2$. This includes the stated interval amplitudes. This proves
  Lemma~\ref{L:masked-time-completion}.
\end{proof}

Let $\mathcal{D}$ consist of finite sums $F(t,x) = \sum_j\tau_j(t)f_j(x)$ with
real $\tau_j \in C_c^\infty(0,T)$ and real Schwartz $f_j$. These give
$H_a$-valued Schwartz functions of $x$. Let $W$ denote the original real
isonormal noise in~\eqref{E:cov}.

\LeanMarker{P:masked-noise-completion}
\begin{proposition}[Completion on the same noise]\label{P:masked-noise-completion}
  The space $\mathcal{D}$ belongs to the original test-noise completion and
  \begin{equation}\label{E:masked-actual-noise-norm}
    \E W(F)^2 = (2\pi)^d\int
    \|\widehat{F}(\xi)\|_{H_a^{\mathbb{C}}}^2\mu(\ud\xi).
  \end{equation}
  Define the auxiliary seminorm
  \[
    N_R(F)^2 = (2\pi)^d\int
    \|\widehat{F}(\eta)\|_{H_a^{\mathbb{C}}}^2\nu_R(\ud\eta).
  \]
  Then $F \mapsto W(\zeta_RF)$ extends as a contraction from the null-quotient
  completion for $N_R$ into the original $L^2(\Omega)$. The extension is
  centered and jointly Gaussian. For finitely many real $F_1,\ldots,F_m \in
  \mathcal{D}$, their masked covariance matrix is bounded in positive
  semidefinite order by the matrix of auxiliary inner products. This gives
  covariance comparison for every finite linear Gaussian image and every convex
  norm power $p \ge 1$.
\end{proposition}

\begin{proof}
  We first identify the spectral norm of the original noise. We then extend
  masked evaluation by the auxiliary bound, and finally compare finite
  Gaussian vectors by adding independent Gaussian noise.

  \medskip
  \noindent\textit{Step 1.~}We identify the original covariance norm.
  For real spatial Schwartz $f,g$, put $h(z) = \int f(y+z)g(y)\ud y$. Fourier
  inversion and Parseval give
  \[
    h(z) = \int e^{i\xi\cdot z}\widehat{f}(\xi)
    \overline{\widehat{g}(\xi)}\ud\xi.
  \]
  The product of the two Fourier transforms is Schwartz, so the inverse
  Fourier identity can be applied again. Evaluating its Fourier transform with
  the positive sign reflects the frequency and contributes the factor
  $(2\pi)^d$.
  Therefore,
  \[
    \int e^{i\xi\cdot z}h(z)\ud z
    = (2\pi)^d\widehat{f}(-\xi)\overline{\widehat{g}(-\xi)}.
  \]
  Apply~\eqref{E:masked-spectral-convention} to the single test $h$ and reflect
  $\xi$ using evenness of $\mu$. The spatial covariance is thus
  $(2\pi)^d\int\widehat{f}(\xi)\overline{\widehat{g}(\xi)}\mu(\ud\xi)$. Summing
  over the finitely many temporal factors gives
  \eqref{E:masked-actual-noise-norm}. Schwartz decay, temperedness and temporal
  integrability justify these pairings. If $W$ is initially given only on
  compact spatial tests, smooth expanding cutoffs converge in the Schwartz
  topology. The spectral form is continuous in that topology, admitting
  $\mathcal{D}$ in the same completion.

  \medskip
  \noindent\textit{Step 2.~}We extend masked evaluation on the same
  probability space.
  Proposition~\ref{P:masked-auxiliary-norm}, applied with $H =
  H_a^{\mathbb{C}}$, gives $\E W(\zeta_RF)^2 \le N_R(F)^2$. It annihilates
  auxiliary-null tests. An auxiliary Cauchy sequence thus gives a Cauchy
  sequence of noise evaluations in $L^2(\Omega)$, proving the contraction
  extension and its uniqueness. Every finite linear combination of the limits is
  an $L^2$ limit of centered Gaussian variables. Their variances converge, so
  their characteristic functions prove that the limits remain jointly centered
  Gaussian.

  \medskip
  \noindent\textit{Step 3.~}We compare finite Gaussian images.
  The comparison uses conditional Jensen: adding an independent centered
  Gaussian vector leaves the conditional mean of the evaluation vector
  unchanged, and the $p$th power of a norm is convex for $p\ge1$.
  For the finite family let $A$ be the covariance matrix and $B$ the auxiliary
  Gram matrix. Applying the norm inequality to $\sum_i b_iF_i$ gives
  $b^{\mathsf{T}}Ab \le b^{\mathsf{T}}Bb$ for every real vector $b$.
  Diagonalizing $B-A$ constructs a centered Gaussian vector $Z$ of that
  covariance, independent of the evaluation vector $X$, on a product space. Then
  $X+Z$ has covariance $B$. Conditional Jensen gives $\E\|X\|^p \le \E\|X+Z\|^p$
  for every norm and $p \ge 1$. Finite Gaussian moments ensure finiteness of
  both sides. This proves Proposition~\ref{P:masked-noise-completion}.
\end{proof}

The auxiliary completion depends on $R$ and $\zeta$. Its construction does not
define multiplication by $\zeta_R$ on the original noise quotient. Nor does the
finite linear Gaussian comparison assert a covariance comparison for nonlinear
wave solutions. The auxiliary and original test norms need not dominate one
another. No spatial representative for every auxiliary class, nonnegative
physical inverse covariance for $\nu_R$, or uniform density bound in $R$ is
asserted.

\LeanMarker{P:masked-quotient-obstruction}
\begin{proposition}[A physical mask need not preserve noise-null tests]
  \label{P:masked-quotient-obstruction}
  In the admissible case $d = 2$, $\alpha = 1$, $0 < a < 1$ and $\Gamma =
  \delta_0(\ud x_1)\ud x_2$, there is a real compact smooth mask $\zeta$ such
  that no finite $C$ satisfies $\E W(\zeta F)^2 \le C\E W(F)^2$ for every
  compact smooth test $F$.
\end{proposition}

\begin{proof}
  We exhibit a test with zero unmasked variance and positive masked variance.
  This rules out any finite comparison constant.

  \medskip
  \noindent\textit{Step 1.~}We identify the covariance in this example.
  Fourier inversion on Schwartz tests gives $\mu =
  (2\pi)^{-1}\ud\xi_1\delta_0(\ud\xi_2)$ in the stated spectral convention. The
  measure $\Gamma$ is nonnegative, even, positive definite and tempered, with
  $\Gamma(cA) = c\Gamma(A)$. It therefore satisfies Hypothesis~\ref{H:homog}
  with $\alpha = 1$. This is not spatial white noise in dimension two. For real
  compact smooth $g$, its spatial covariance is, by ordinary Fubini,
  \[
    V(g) = \int_\R\left(\int_\R g(x_1,x_2)\ud x_2\right)^2\ud x_1.
  \]

  \medskip
  \noindent\textit{Step 2.~}We construct the mask and a noise-null test.
  Choose real even $\chi \in C_c^\infty(\R)$ equal to one on $[-1,1]$ and zero
  outside $[-2,2]$. For example, on $1 < r = |x| < 2$ use
  \[
    \chi(x) = \frac{e^{-1/(2-r)}}{e^{-1/(2-r)}+e^{-1/(r-1)}}.
  \]
  The endpoint exponentials vanish to every order, so the plateau pastes
  smoothly. Let $I_\chi = \int(\chi')^2 > 0$ and $J_\chi = \int\chi^4 > 0$.
  Positivity of $I_\chi$ follows because $\chi$ is not constant, and $J_\chi \ge
  2$. Put $g(x) = \chi(x_1)\chi''(x_2)$ and $\zeta(x) = \chi(x_1)\chi(x_2)$.
  Their unmasked marginal vanishes since $\int\chi'' = 0$, whereas integration
  by parts gives the masked marginal $-I_\chi\chi(x_1)^2$. Thus, $V(g) = 0$ and
  $V(\zeta g) = J_\chi I_\chi^2 > 0$. For a nonzero nonnegative $\tau \in
  C_c^\infty(0,T)$, set $F(t,x) = \tau(t)g(x)$. The temporal factor
  $\|\tau\|_a^2$ is finite and strictly positive because its nonnegative
  integrand is positive on a set of positive two-dimensional measure.
  Consequently, the full unmasked variance is zero and the masked variance is
  positive. This proves Proposition~\ref{P:masked-quotient-obstruction}.
\end{proof}

No-go~\ref{P:nogo-22} records the quotient interpretation of this witness:
the mask separates two representatives of the same noise-null class.
This is why we mask tests before constructing the auxiliary completion.

\subsection{Temporal modulation and continuous wave factors}

Write $\langle x\rangle = \sqrt{1+x^2}$. Let $I = [s,s+h] \subset [0,T]$, with
$0 < h \le 1$. Norms below are those of the complexified temporal completion in
Lemma~\ref{L:masked-time-completion}.

\LeanMarker{L:masked-temporal-modulation}
\begin{lemma}[Modulated amplitudes]\label{L:masked-temporal-modulation}
  Let $A \in C^1(I;\mathbb{C})$, $M = \sup_I|A|$ and $D = \sup_I|A'|$. For every
  real $\omega$,
  \begin{equation}\label{E:masked-temporal-modulation}
    \|1_IA(t)e^{i\omega t}\|_a^2
    \le C_0(a)h^{2-a}(M+hD)^2\langle h\omega\rangle^{a-1},
    \qquad C_0(a) = \frac{16}{1-a}.
  \end{equation}
  The bound is uniform in $T,s,h,A$ and $\omega$.
\end{lemma}

\begin{proof}
  We express the covariance through the time difference and rescale the
  interval to unit length. At low frequency we use absolute values; at high
  frequency we separate the singular endpoint before integrating by parts.

  \medskip
  \noindent\textit{Step 1.~}We express and rescale the covariance.
  Put $c(v) = \int_s^{s+h-v}A(u+v)\overline{A(u)}\ud u$ for $0 \le v \le h$.
  Differentiation at the moving endpoint gives
  \[
    c'(v) = -A(s+h)\overline{A(s+h-v)}
    +\int_s^{s+h-v}A'(u+v)\overline{A(u)}\ud u.
  \]
  Thus, $c(h) = 0$, $|c(v)| \le hM^2$ and $|c'(v)| \le M^2+hMD$. Splitting the
  covariance square into its two conjugate triangles gives
  \[
    \|1_IA(t)e^{i\omega t}\|_a^2
    = 2\operatorname{Re}\int_0^h v^{-a}e^{i\omega v}c(v)\ud v
    = 2h^{2-a}\operatorname{Re}\int_0^1 y^{-a}e^{i\omega hy}f(y)\ud y,
  \]
  where $f(y) = c(hy)/h$. Set $K = M+hD$. Then $f(1) = 0$ and $|f|,|f'| \le
  K^2$. If $K = 0$ the assertion is immediate. If $|\omega h| \le 1$, the last
  integral is at most $K^2/(1-a)$ in absolute value, and $\langle\omega
  h\rangle^{a-1} \ge 2^{(a-1)/2}$. The claimed constant is sufficient in this
  case.

  The low-frequency constant check is $2/(1-a)\le C_0(a)2^{(a-1)/2}$, after
  including the factor $2h^{2-a}$.
  \medskip
  \noindent\textit{Step 2.~}We separate the singular endpoint at high
  frequency.
  Now suppose $k = |\omega h| > 1$. On $[0,1/k]$ the integral is at most
  $K^2k^{a-1}/(1-a)$. On $[1/k,1]$, integrate the exponential by parts. The
  endpoint at one vanishes, while the other endpoint contributes at most
  $K^2k^{a-1}$. The derivative integral, divided by $k$, is bounded by
  \[
    \frac{K^2}k\left(a\int_{1/k}^1 y^{-a-1}\ud y
    +\int_{1/k}^1 y^{-a}\ud y\right)
    \le K^2 k^{a-1}\left(1+\frac1{1-a}\right).
  \]
  Here $k^{-1} \le k^{a-1}$ since $a > 0$. The total is at most
  $4K^2k^{a-1}/(1-a)$. Finally, $k^{a-1} \le \sqrt{2}\langle\omega
  h\rangle^{a-1}$ and $8\sqrt{2} \le 16$. Including the factor $2h^{2-a}$ proves
  the bound for either sign of $\omega$. This proves
  Lemma~\ref{L:masked-temporal-modulation}.
\end{proof}

Lemma~\ref{L:masked-temporal-modulation} treats one exponential mode. A
product of two sines or cosines is a sum of four such modes, which we bound
next.

\LeanMarker{L:masked-trig-products}
\begin{lemma}[Correlated trigonometric products]\label{L:masked-trig-products}
  Let $r,\rho \ge 0$. Let $P_r$ be either $\cos(rt)$ or $\sin(rt)$, and let
  $Q_\rho$ be either $\cos(\rho t)$ or $\sin(\rho t)$. With the amplitude and
  constant above,
  \begin{align}
    \|1_IAP_rQ_\rho\|_a^2
    &\le C_0(a)h^{2-a}(M+hD)^2\mathcal{Q}_h(r,\rho),
    \label{E:masked-trig-products}\\
    \mathcal{Q}_h(r,\rho)
    &:= \langle h(r-\rho)\rangle^{a-1}
    +\langle h(r+\rho)\rangle^{a-1}.
    \label{E:masked-frequency-pair}
  \end{align}
\end{lemma}

\begin{proof}
  Expanding the two trigonometric factors gives four complex exponential modes
  with frequencies $\sigma r+\tau\rho$, $\sigma,\tau \in \{-1,1\}$, and
  coefficients of absolute value $1/4$. In a Hilbert space,
  \[
    \left\|\sum_j c_jv_j\right\|^2
    \le \left(\sum_j|c_j|\right)\sum_j|c_j|\|v_j\|^2.
  \]
  This Hilbert-space inequality follows from the triangle inequality followed
  by scalar Cauchy--Schwarz applied to $|c_j|^{1/2}$ and $|c_j|^{1/2}\|v_j\|$.

  The first coefficient sum is one. Apply
  Lemma~\ref{L:masked-temporal-modulation} to each mode. Two frequencies have
  absolute value $|r-\rho|$ and two have absolute value $r+\rho$. This gives one
  half of the asserted right side, and hence the bound. Repeated or zero
  frequencies are permitted; no orthogonality or independence of the modes is
  used. This proves Lemma~\ref{L:masked-trig-products}.
\end{proof}

\LeanMarker{L:masked-continuous-wave}
\begin{lemma}[Continuous sine factors]\label{L:masked-continuous-wave}
  Define $c_r(t) = \cos(rt)$ and
  \[
    b_r(t) = \langle r\rangle\frac{\sin(rt)}r\quad(r > 0),
    \qquad b_0(t) = t.
  \]
  For $F_r \in \{c_r,b_r\}$ and $G_\rho \in \{c_\rho,b_\rho\}$,
  \begin{equation}\label{E:masked-continuous-wave}
    \|1_IF_rG_\rho\|_a^2
    \le C_1(a)(1+T)^4h^{2-a}\mathcal{Q}_h(r,\rho),
    \qquad C_1(a) = 4\cdot3^{1-a}C_0(a).
  \end{equation}
\end{lemma}

\begin{proof}
  We replace low-frequency sine quotients by averages of cosines. The
  trigonometric estimate applies to each integrand; a comparison of shifted
  frequency weights then gives the stated constant.

  \medskip
  \noindent\textit{Step 1.~}We represent each factor without a singular
  quotient.
  For $r \ge 1$, $b_r$ is a sine times a coefficient at most $\sqrt{2}$. For $0
  \le r \le 1$, the fundamental theorem of calculus gives
  \[
    b_r(t) = t\langle r\rangle\int_0^1\cos(vrt)\ud v.
  \]
  For $r > 0$, integrate the derivative of $v \mapsto \sin(vrt)$ from zero to
  one and divide by $r$. At $r = 0$, the right side is $t$, in agreement with
  the definition.
  Each low-frequency sine factor contributes one factor $t$ and one averaging
  parameter; each high-frequency sine contributes a coefficient of absolute
  value at most $\sqrt2$, while a cosine contributes neither.
  This includes $r = 0$. Hence, a product $F_rG_\rho$ is a scalar coefficient of
  absolute value at most two times a normalized integral of
  $t^jP_{r'}(t)Q_{\rho'}(t)$, where $j \in \{0,1,2\}$ and $|r'-r|,|\rho'-\rho|
  \le 1$. The integral is over zero, one or two copies of $[0,1]$; with no
  averages, use a one-point probability space. Its parameter dependence is
  continuous in $L^2(0,T)$ by dominated convergence.
  Lemma~\ref{L:masked-time-completion} makes it continuous in $H_a$ as well.
  Thus, the averages are Bochner integrals, and their $H_a$ images agree with
  the pointwise averages by continuity of the inclusion. Hilbert Jensen bounds
  the squared norm of an average by the average of the squared norms.

  \medskip
  \noindent\textit{Step 2.~}We bound each integrand and compare frequency
  weights.
  Apply Lemma~\ref{L:masked-trig-products} to each integrand with $A(t) = t^j$.
  For $j = 0,1,2$, respectively, $M+hD$ is at most $1$, $T+1$ and $T^2+2T$. In
  every case $(M+hD)^2 \le (1+T)^4$. For real $x,e$, the triangle inequality and
  its reverse imply
  \[
    \frac{\langle x\rangle}{1+|e|}
    \le \langle x+e\rangle \le (1+|e|)\langle x\rangle.
  \]
  For the upper bracket comparison, view $\langle x\rangle$ as the Euclidean
  norm of $(1,x)$ and use the triangle inequality. Applying that same upper
  bound with $x+e$ and $-e$ gives the lower comparison.
  The changes in $h(r\pm\rho)$ have magnitude at most two, since $h \le 1$.
  Because $a-1 < 0$, each altered bracket to the power $a-1$ is at most
  $3^{1-a}$ times the original one. Apply these comparisons, Hilbert Jensen and
  the squared scalar coefficient bound four to obtain
  \eqref{E:masked-continuous-wave}. This proves
  Lemma~\ref{L:masked-continuous-wave}.
\end{proof}

These estimates concern covariance forms and temporal factors. Their application
to displacement and velocity is in Section~\ref{S:wave-variance}; the subsequent
operator, source and moment results require their separate proofs. They do not
follow from covariance domination alone.

\section{Frequency estimates and both wave variance directions}
\label{S:wave-variance}

We now apply the covariance and temporal estimates of
Section~\ref{S:masked-covariance} to the full wave phase. We retain
Hypothesis~\ref{H:homog}, write $a = \alpha_0$, and choose and fix
\begin{equation}\label{E:variance-gamma-range}
  0 < \gamma < \min\left(\frac12,\frac{2-\alpha}{2}\right).
\end{equation}
The interval is nonempty because $\alpha+a < 2$. We use $\langle\xi\rangle =
\sqrt{1+|\xi|^2}$ and set
\begin{equation}\label{E:variance-scale-profile}
  B_\alpha(h) =
  \begin{cases}
    1, & \alpha < 1,\\
    1+\log(1/h), & \alpha = 1,\\
    h^{1-\alpha}, & \alpha > 1,
  \end{cases}
  \qquad 0 < h \le 1.
\end{equation}

\subsection{A scalar row estimate for a measure}

Let $\nu$ be a positive Radon measure satisfying
\begin{equation}\label{E:frequency-translated-growth}
  \nu(B(y,r)) \le C r^\alpha,
  \qquad y \in \R^d,\quad r > 0,
\end{equation}
where $C \ge 0$. If the estimate is initially known for open balls, apply it at
radius $r+\epsilon$ and let $\epsilon \downarrow 0$ to obtain the closed-ball
estimate. Neither evenness nor a density is assumed here.
Using~\eqref{E:masked-frequency-pair}, define
\begin{align}
  \Phi_h(\xi,\eta)
  &:= \frac{\mathcal{Q}_h(|\xi|,|\eta|)}
  {\langle\xi\rangle\langle\eta\rangle},
  \label{E:frequency-row-kernel}\\
  J_h(\xi)
  &:= \langle\xi\rangle^{2\gamma}
  \int_{\R^d}\Phi_h(\xi,\xi-\lambda)\nu(\ud\lambda).
  \label{E:frequency-row-integral}
\end{align}
These are initially nonnegative extended integrals. In particular, the notation
denotes a pushforward of $\nu$, not integration against an assumed spectral
density.

\LeanMarker{P:frequency-row}
\begin{proposition}[Uniform scalar row]\label{P:frequency-row}
  There is a finite $K = K(\alpha,a,\gamma)$, independent of $d,\nu,C,h$ and
  $\xi$, such that
  \begin{equation}\label{E:frequency-row-bound}
    \sup_{\xi\in\R^d}J_h(\xi) \le K C B_\alpha(h),
    \qquad 0 < h \le 1.
  \end{equation}
\end{proposition}

\begin{proof}
  We first bound the integral of the inverse Bessel weight on a ball. We then
  split the row integral into frequencies near and far from the origin
  relative to the fixed output frequency. The near bound uses the choice of
  $\gamma$; the far bound uses temporal decay and $\alpha+a < 2$.

  \medskip
  \noindent\textit{Step 1.~}We bound the weighted mass on a ball.
  If $C = 0$, a countable covering by finite balls gives $\nu = 0$, and the
  result is immediate. Suppose $C > 0$. For each $\xi$, let $\nu_\xi$ be the
  pushforward under $\lambda \mapsto \eta = \xi-\lambda$. Then $\nu_\xi(B(0,R))
  \le CR^\alpha$. For $q \ge 1$, put
  \[
    A_\xi(q) = \int_{|\eta|\le5q}\langle\eta\rangle^{-1}
    \nu_\xi(\ud\eta).
  \]
  Choose $N = \lceil\log_2(5q)\rceil$, so $5q \le 2^N < 10q$. Splitting into the
  unit ball and the annuli $2^j < |\eta| \le 2^{j+1}$ gives
  \[
    A_\xi(q) \le C\left(1+2^\alpha
    \sum_{j=0}^{N-1}2^{j(\alpha-1)}\right).
  \]
  The geometric sum is bounded when $\alpha < 1$, is $N$ when $\alpha = 1$, and
  is at most $2^{N(\alpha-1)}/(2^{\alpha-1}-1)$ when $\alpha > 1$. Hence,
  \begin{equation}\label{E:frequency-inverse-bessel}
    A_\xi(q) \le K_\alpha C
    \begin{cases}
      1, & \alpha < 1,\\
      1+\log q, & \alpha = 1,\\
      q^{\alpha-1}, & \alpha > 1,
    \end{cases}
  \end{equation}
  with a constant depending only on $\alpha$.

  \medskip
  \noindent\textit{Step 2.~}We estimate the near region.
  Set $q = \max(1,|\xi|)$ and $\epsilon = 1-2\gamma > 0$. On $|\lambda| \le 4q$
  we have $|\eta| \le 5q$. Both bracket powers in $\mathcal{Q}_h$ are at most
  one, since $a-1 < 0$. The contribution of this region to $J_h(\xi)$ is at most
  \[
    2\langle\xi\rangle^{-\epsilon}A_\xi(q)
    \le 2q^{-\epsilon}A_\xi(q).
  \]
  If $\alpha < 1$, this is uniformly bounded
  by~\eqref{E:frequency-inverse-bessel}. If $\alpha = 1$, the bound follows from
  \[
    e^{-\epsilon u}(1+u) \le 1+\frac1{e\epsilon},
    \qquad u = \log q \ge 0,
  \]
  since $u e^{-\epsilon u}$ has maximum $1/(e\epsilon)$. If $\alpha > 1$, the
  remaining power is $q^{\alpha+2\gamma-2}$, whose exponent is negative
  by~\eqref{E:variance-gamma-range}. Thus, the near contribution is at most a
  constant times $C$ in all three cases.

  \medskip
  \noindent\textit{Step 3.~}We use temporal decay in the far region.
  On the far region write $r = |\lambda| > 4q$. Then
  \[
    |\eta| \ge \frac{3r}4,\qquad
    \big||\eta|-|\xi|\big| \ge \frac r2,\qquad
    |\eta|+|\xi| \ge \frac r2.
  \]
  Indeed, $|\xi| \le q < r/4$. The reverse triangle inequality gives $|\eta| =
  |\xi-\lambda| \ge r-|\xi| > 3r/4$, and then $|\eta|-|\xi| > r/2$. The sum of
  the two magnitudes is at least their difference. These bounds apply to both
  frequency arguments in $\mathcal{Q}_h$.
  Each oscillatory bracket is at least $\langle hr\rangle/2$, whereas
  $\langle\eta\rangle^{-1} \le 4/(3r)$ and $\langle\xi\rangle^{-\epsilon} \le
  1$. Consequently, the far contribution is at most
  \[
    \frac43\,2^{2-a}\int_{|\lambda|>4q}
    |\lambda|^{-1}\langle h|\lambda|\rangle^{a-1}\nu(\ud\lambda).
  \]
  Raising the lower bound on each bracket to $a-1 < 0$ gives an upper bound
  with coefficient $2^{1-a}$. There are two brackets, so their sum contributes
  $2^{2-a}$. The factor $4/3$ comes from the inverse Bessel weight of $\eta$.
  Enlarge the domain to $|\lambda| > 1$ and split it into dyadic annuli.
  Applying~\eqref{E:frequency-translated-growth} to their masses bounds this
  expression by a constant depending only on $\alpha,a$ times $C S(h)$, where
  \[
    S(h) = \sum_{j=0}^\infty
    2^{j(\alpha-1)}\langle h2^j\rangle^{a-1}.
  \]
  On $2^j < |\lambda| \le 2^{j+1}$, the inverse radius is at most $2^{-j}$ and
  the decreasing bracket power is at most $\langle h2^j\rangle^{a-1}$. The
  mass is at most $C2^{(j+1)\alpha}$. Multiplying these bounds gives the $j$th
  term of $S(h)$ up to a constant depending only on $\alpha$.

  \medskip
  \noindent\textit{Step 4.~}We sum below and above the temporal scale.
  Choose $L = \lfloor\log_2(1/h)\rfloor \ge 0$. For $j \le L$ the bracket power
  is at most one. For $j \ge L+1$ it is at most $(h2^j)^{a-1}$. Since
  $\alpha+a-2 < 0$, the tail is
  \[
    \frac{h^{a-1}2^{(L+1)(\alpha+a-2)}}{1-2^{\alpha+a-2}}
    \le \frac{h^{1-\alpha}}{1-2^{\alpha+a-2}}.
  \]
  The inequality uses $2^{L+1} > 1/h$ and the negative exponent. The finite sum
  is bounded by $1/(1-2^{\alpha-1})$ for $\alpha < 1$, by $1+\log_2(1/h)$ for
  $\alpha = 1$, and by
  \[
    \frac{2^{\alpha-1}h^{1-\alpha}}{2^{\alpha-1}-1}
  \]
  for $\alpha > 1$. Since $h^{1-\alpha} \le 1$ in the first case, we obtain
  $S(h) \le K_{\alpha,a}B_\alpha(h)$. Combining the near and far bounds, and
  using $B_\alpha(h) \ge 1$, proves the assertion and finiteness of the original
  integral. All constants arose from numerical shell estimates and are
  independent of $d$. This proves Proposition~\ref{P:frequency-row}.
\end{proof}

For $\alpha = 1$ the profile $B_\alpha(h)$ in
Proposition~\ref{P:frequency-row} carries a logarithm. After multiplication
by $h^{2-a}$, we absorb it into a power of $h$.

\LeanMarker{C:frequency-time}
\begin{corollary}[A time exponent above one half]\label{C:frequency-time}
  For every
  \begin{equation}\label{E:variance-beta-range}
    \frac12 < \beta < \min\left(1-\frac a2,
    \frac{3-a-\alpha}{2}\right),
  \end{equation}
  there is a finite $K_\beta = K_\beta(\alpha,a,\gamma,\beta)$ such that
  \begin{equation}\label{E:frequency-time-bound}
    h^{2-a}\sup_\xi J_h(\xi) \le K_\beta C h^{2\beta},
    \qquad 0 < h \le 1.
  \end{equation}
\end{corollary}

\begin{proof}
  Both upper endpoints in~\eqref{E:variance-beta-range} exceed one half, by $a <
  1$ and $\alpha+a < 2$. If $\alpha < 1$, use $2-a > 2\beta$ in
  Proposition~\ref{P:frequency-row}. If $\alpha > 1$, use $3-a-\alpha > 2\beta$.
  For $\alpha = 1$, set $\delta = 2-a-2\beta > 0$ and $u = \log(1/h)$. Then
  \[
    h^\delta(1+\log(1/h))
    = e^{-\delta u}(1+u) \le 1+\frac1{e\delta}.
  \]
  For the displayed bound, $u \ge 0$, the exponential term is at most one, and
  $u e^{-\delta u}$ has maximum $1/(e\delta)$. Multiplying the display by
  $h^{2\beta}$ gives the required estimate at $\alpha = 1$.
  This absorbs the critical logarithm. No endpoint value of $\beta$ is used.
  This proves Corollary~\ref{C:frequency-time}.
\end{proof}

\subsection{Forms on the full phase}

Choose and fix $T > 0$. Set $\Lambda = (1-\Delta)^{1/2}$, $s = 1/2$ and $E_r =
H^r(\R^d)\times H^{r-1}(\R^d)$. On the complexification of $E = E_s$, use inner
products linear in the first slot. The map
\[
  Jz(\xi) = \left(\langle\xi\rangle^{1/2}\widehat{z}_1(\xi),
  \langle\xi\rangle^{-1/2}\widehat{z}_2(\xi)\right)
\]
is unitary onto $L^2(\R^d;\mathbb{C}^2)$. The full massless wave group has
Fourier matrix
\begin{equation}\label{E:variance-free-matrix}
  S(t,\xi) =
  \begin{pmatrix}
    \cos(t|\xi|) & \sin(t|\xi|)/|\xi|\\
    -|\xi|\sin(t|\xi|) & \cos(t|\xi|)
  \end{pmatrix},
\end{equation}
where the sine quotient is $t$ at zero. Let $B(b)(u,v) = (0,bu)$. Choose and fix
a real even nonnegative $\zeta \in C_c^\infty(\R^d)$ equal to one on $B_1$ and
zero off $B_2$, and put $\zeta_R(x) = \zeta(x/R)$ for $R > 0$. For $I =
[t_0,t_0+h] \subset [0,T]$, $0 < h \le 1$, and Schwartz phases $z,w$, put
\begin{align}
  v_z(t) &= [S(t)z]_1,
  \label{E:variance-displacement}\\
  G_w(t) &= -\frac{\sin(t|D|)}{|D|}\Lambda w_1
  +\cos(t|D|)\Lambda^{-1}w_2,
  \label{E:variance-dual-wave}\\
  F_{z,w}(t,x) &= 1_I(t)v_z(t,x)\overline{G_w(t,x)},
  \qquad \ell_R(I;z,w) = W_{\mathbb{C}}(\zeta_RF_{z,w}),
  \label{E:variance-actual-form}
\end{align}
where $W_{\mathbb{C}}(F) = W(\operatorname{Re}F)+iW(\operatorname{Im}F)$. There
is no observation cutoff in this definition.

We first assume $\Gamma \ne 0$. If $\Gamma = 0$, all forms and the column maps
below are zero, so their bounds hold directly. For the nonzero case, let $\nu_R$
be~\eqref{E:masked-auxiliary-measure} and set
\begin{equation}\label{E:variance-actual-growth-constant}
  L_\zeta = \|\widehat{\zeta}\|_1,\qquad
  C_\nu = c_d^2L_\zeta^2C_\Gamma.
\end{equation}
Proposition~\ref{P:masked-auxiliary-norm} gives a positive even tempered Radon
measure with $\nu_R(B(y,r)) \le C_\nu r^\alpha$, uniformly in $R$.

\LeanMarker{P:actual-wave-form}
\begin{proposition}[Form and phase kernel]\label{P:actual-wave-form}
  The tests $\zeta_RF_{z,w}$ in~\eqref{E:variance-actual-form} belong to the
  original noise completion. The form $\ell_R$ is linear in $z$ and conjugate
  linear in $w$, and uses the same $W$ as the model. For a smooth deterministic
  coefficient $b$ with common compact spatial support and bounded spatial
  derivatives on $[0,T]$,
  \begin{equation}\label{E:variance-physical-adjoint}
    \langle S(-t)B(b(t))S(t)z,w\rangle_E
    = \int b(t,x)v_z(t,x)\overline{G_w(t,x)}\ud x.
  \end{equation}
  For a real compact smooth $\chi$, multiplication $D_\chi$ on both phase
  entries has adjoint
  \begin{equation}\label{E:variance-cutoff-adjoint}
    D_\chi^{*,E}w
    = \left(\Lambda^{-1}\chi\Lambda w_1,
    \Lambda\chi\Lambda^{-1}w_2\right).
  \end{equation}
  With $U = Jz$, $V = Jw$, $r = |\xi|$, $\rho = |\eta|$ and the continuous
  functions $c_r,b_r$ of Lemma~\ref{L:masked-continuous-wave}, set
  \begin{equation}\label{E:variance-normalized-kernel}
    \mathcal{K}(t;\xi,\eta)
    = \frac1{\sqrt{\langle\xi\rangle\langle\eta\rangle}}
    \begin{pmatrix}
      -b_rc_\rho & -b_rb_\rho\\
      c_rc_\rho & c_rb_\rho
    \end{pmatrix}(t).
  \end{equation}
  Then
  \begin{equation}\label{E:variance-normalized-covariance}
    \E|\ell_R(I;z,w)|^2
    \le \int\left\|\int\sum_{i,j=1}^2
    1_I\mathcal{K}_{ij}(\cdot;\eta+\lambda,\eta)
    U_j(\eta)\overline{V_i(\eta+\lambda)}\ud\eta
    \right\|_{H_a^{\mathbb{C}}}^2\nu_R(\ud\lambda).
  \end{equation}
\end{proposition}

\begin{proof}
  We first admit the masked products in the original noise completion and
  identify the same limits in the auxiliary norm. We then extend the
  comparison to complex tests, compute the phase adjoints, and write the
  Fourier kernel with its normalization.

  \medskip
  \noindent\textit{Step 1.~}We admit the masked products in the original
  completion.
  The wave multipliers are smooth at zero through the continuous sine quotient,
  with polynomially bounded frequency derivatives on finite time intervals.
  Bessel multipliers preserve Schwartz space. Thus, $v_z,G_w$ and their product
  are continuous in all Schwartz seminorms on $[0,T]$. For tests supported in a
  fixed compact $K$, the absolute spatial pairing is bounded by
  \[
    |K|\Gamma(K-K)\|F(t)\|_\infty\|H(u)\|_\infty.
  \]
  Indeed, only $y \in K-K$ contributes to $\int\Gamma(\ud y)\int
  F(t,x)\overline{H(u,x-y)}\ud x$. Combine this with the temporal Schur bound
  and approximate the interval indicator and the continuous time dependence by
  smooth interior temporal tests, keeping the spatial support of $\zeta_R$. This
  admits the masked product in the original completion in every allowed
  dimension.

  \medskip
  \noindent\textit{Step 2.~}We compare the limits in the two completions.
  We also identify its auxiliary and spectral limits. Approximate a
  Schwartz-valued function on $I$ by time partitions with frozen spatial values,
  and then replace the finitely many indicators by compact smooth temporal
  functions. A diagonal sequence $F_n$ of finite temporal tensors converges in
  $L^2$ time in every Schwartz seminorm. Fourier integration by parts and the
  spatial weighted supremum bounds give, for sufficiently large fixed $N$,
  \[
    \|\widehat{F_n-F}(\xi)\|_{L^2(0,T)}
    \le \epsilon_n\langle\xi\rangle^{-N},\qquad \epsilon_n \to 0.
  \]
  For a fixed seminorm, continuity on the compact time interval makes a fine
  time partition uniformly accurate. Each frozen value is Schwartz; replacing
  finitely many indicators in $L^2$ time controls the seminorm error by a
  finite sum. A diagonal choice over the countable seminorms gives the
  sequence, and Fourier integration by parts turns those errors into the
  displayed frequency decay.
  Choose $N$ so that the square weight is integrable against both tempered
  measures $\mu$ and $\nu_R$. Lemma~\ref{L:masked-time-completion} gives
  convergence in their spectral $H_a$ norms. Multiplication by the fixed
  $\zeta_R$ preserves these $L^2$ Schwartz-seminorm limits, so the masked
  differences converge in the $\mu$ spectral norm too. The preceding
  compact-support estimate also gives convergence in the original completion.
  Hence, the auxiliary extension and original-noise limit agree on these
  products. Constants for this admission argument may depend on $R$; they are
  not the uniform variance constants below.

  \medskip
  \noindent\textit{Step 3.~}We pass from real tests to complex tests.
  For real tensor tests $f,g$, the Hermitian symmetry of their spatial
  transforms gives
  \[
    A(-\xi) = \overline{A(\xi)},\qquad
    A(\xi) = \langle\widehat{f}(\xi),\widehat{g}(\xi)\rangle_{H_a^{\mathbb{C}}}.
  \]
  Integration against either even $\mu$ or even $\nu_R$ cancels the imaginary
  mixed term. Expanding the squared norm of $\widehat{f}+i\widehat{g}$, then
  using the preceding approximation, extends
  Proposition~\ref{P:masked-noise-completion} to the present complex tests:
  \begin{equation}\label{E:variance-complex-noise}
    \E|W_{\mathbb{C}}(\zeta_RF)|^2
    = (2\pi)^d\int\|\widehat{\zeta_RF}(\xi)\|_{H_a^{\mathbb{C}}}^2\mu(\ud\xi)
    \le (2\pi)^d\int\|\widehat{F}(\lambda)\|_{H_a^{\mathbb{C}}}^2\nu_R(\ud\lambda).
  \end{equation}
  This is a bound for the masked noise, not a declaration that its law is
  stationary. Complex linearity of $W_{\mathbb{C}}$ and linearity of $v_z,G_w$
  prove the stated sesquilinearity.

  \medskip
  \noindent\textit{Step 4.~}We compute the physical adjoints.
  The displacement and velocity pairings carry $\Lambda$ and $\Lambda^{-1}$,
  respectively; moving these weights onto $w$ in $S(-t)(0,b(t)v_z(t))$ gives
  the negative sine and cosine terms defining $G_w$.
  The first-response integral below means $\int_I S(-t)B(b(t))S(t)z\,\ud t$ in
  $E$; its existence allows the pointwise identity to be integrated, but is
  not needed for the pointwise adjoint calculation itself.
  Smooth multiplication is bounded on each Sobolev space: Fourier convolution,
  the inequality $\langle\xi\rangle^r\langle\eta\rangle^{-r} \le
  2^{|r|/2}\langle\xi-\eta\rangle^{|r|}$ and weighted Schwartz decay reduce this
  to the $L^1*L^2$ bound. The free matrix is strongly continuous and bounded on
  $E$ on compact time intervals, directly by its weighted Fourier entries. Thus,
  the ordinary first-response integral exists on states. Since $S(-t)(0,f) =
  (-\sin(t|D|)f/|D|,\cos(t|D|)f)$, moving the real self-adjoint multipliers and
  the two Sobolev weights onto $w$ proves \eqref{E:variance-physical-adjoint}.
  The observed form is $\ell_R(I;z,D_\chi^{*,E}w)$: applying $D_\chi$ to an
  output moves its phase adjoint onto the test vector $w$.
  The same weighted duality gives \eqref{E:variance-cutoff-adjoint}. Inserting
  this adjoint in $G_w$ is the correctly observed form; $D_\chi$ need not be
  self-adjoint on $E$. This deterministic identity does not identify a
  positive-width source with the original noise or assert a source limit.

  \medskip
  \noindent\textit{Step 5.~}We compute the normalized covariance kernel.
  The definitions give exactly
  \begin{align*}
    \widehat{v}_z(t,\eta)
    &= \langle\eta\rangle^{-1/2}
    [c_\rho(t)U_1(\eta)+b_\rho(t)U_2(\eta)],\\
    \widehat{G}_w(t,\xi)
    &= \langle\xi\rangle^{-1/2}
    [-b_r(t)V_1(\xi)+c_r(t)V_2(\xi)].
  \end{align*}
  In the first line, $\eta$ is the input frequency and $U_1,U_2$ are its two
  weighted coordinates. In the second, $\xi$ is the frequency paired with
  $V_1,V_2$. Multiplying the first line by the complex conjugate of the second
  gives four terms, indexed by input component $j$ and output component $i$.
  The wave factors are real, so conjugation affects the coordinates but not
  these factors.
  Their product gives~\eqref{E:variance-normalized-kernel}. At zero frequencies
  one uses the continuous sine quotient, not division by zero; at $\xi = \eta =
  0$ the matrix is
  $\left(\begin{smallmatrix}-t&-t^2\\1&t\end{smallmatrix}\right)$. Fourier
  multiplication with the conjugate yields
  \[
    \widehat{F}_{z,w}(t,-\lambda)
    = c_d1_I(t)\int\sum_{i,j=1}^2\mathcal{K}_{ij}(t;\eta+\lambda,\eta)
    U_j(\eta)\overline{V_i(\eta+\lambda)}\ud\eta.
  \]
  The product with the conjugate has frequency $\eta-\xi$. Evaluating it at
  $-\lambda$ therefore sets $\xi = \eta+\lambda$, explaining the shift in
  every kernel entry and in $V_i$.
  The spatial integral is absolutely convergent uniformly on the finite time
  interval by Schwartz decay. The Schur bound also makes it a Bochner
  $H_a^{\mathbb{C}}$ integral. Reflect $\lambda$ using evenness of $\nu_R$
  in~\eqref{E:variance-complex-noise}. Since $(2\pi)^dc_d^2 = 1$, this
  proves~\eqref{E:variance-normalized-covariance} with no remaining Fourier
  factor. This proves Proposition~\ref{P:actual-wave-form}.
\end{proof}

\subsection{Bessel's inequality and the two variances}

Proposition~\ref{P:actual-wave-form} bounds the variance of one form
$\ell_R(I;z,w)$. To sum it over an orthonormal family of phases, we need a
Bessel inequality with values in $H_a^{\mathbb{C}}$.
Proposition~\ref{P:wave-two-variances} uses it in both directions.

\LeanMarker{L:phase-bessel}
\begin{lemma}[Hilbert-valued finite Bessel inequality]\label{L:phase-bessel}
  Let $H$ be a separable complex Hilbert space, let $\{u_k\}$ be a finite
  orthonormal family in $L^2(\R^d;\mathbb{C}^2)$, and let $g = (g_1,g_2) \in
  L^2(\R^d;H^2)$. Then
  \begin{equation}\label{E:phase-bessel}
    \sum_k\left\|\int\sum_{j=1}^2u_{k,j}(\eta)g_j(\eta)\ud\eta\right\|_H^2
    \le \int\sum_{j=1}^2\|g_j(\eta)\|_H^2\ud\eta.
  \end{equation}
  The same bound holds with $\overline{u_{k,j}}$ in place of $u_{k,j}$.
\end{lemma}

\begin{proof}
  The nonconjugated integrand requires the conjugate orthonormal family in
  scalar Bessel, because the scalar inner product conjugates its second
  argument; complex conjugation preserves orthonormality.
  Expand $g$ in a countable orthonormal basis of $H$. On each scalar coefficient
  apply scalar Bessel to the orthonormal family $\{\overline{u_k}\}$, or to
  $\{u_k\}$ for the conjugated version. Scalar finite Bessel follows by
  subtracting an orthogonal finite projection and expanding its nonnegative
  squared norm. Sum over the $H$ basis and use Parseval and nonnegative Tonelli.
  The Bochner integrals exist by $L^2$ Cauchy--Schwarz. A finite-dimensional $H$
  uses its finite basis, and the zero space is immediate. No independence is
  involved. This proves Lemma~\ref{L:phase-bessel}.
\end{proof}

\LeanMarker{P:wave-two-variances}
\begin{proposition}[Both variance directions]\label{P:wave-two-variances}
  There is a finite $C_T$, depending on the fixed model, mask,
  $a,\alpha,\gamma,T$, but not on $R,I$ or the orthonormal family, such that for
  every finite complex $E$-orthonormal Schwartz family $\{e_k\}$ and every
  Schwartz phase $z$,
  \begin{align}
    \sum_k\E|\ell_R(I;e_k,z)|^2
    &\le C_T h^{2-a}B_\alpha(h)\|z\|_{E_{s-\gamma}}^2,
    \label{E:wave-input-variance}\\
    \sum_k\E|\ell_R(I;z,e_k)|^2
    &\le C_T h^{2-a}B_\alpha(h)\|z\|_{E_{s-\gamma}}^2.
    \label{E:wave-output-variance}
  \end{align}
  For every $\beta$ in~\eqref{E:variance-beta-range}, either right side may be
  replaced by $C_{T,\beta}h^{2\beta}\|z\|_{E_{s-\gamma}}^2$.
\end{proposition}

\begin{proof}
  We combine the temporal bound for each matrix entry with the scalar row
  estimate. We apply Hilbert-valued Bessel first to the input family and then
  to the output family, keeping both variance directions and the same
  constant.

  \medskip
  \noindent\textit{Step 1.~}We bound the temporal kernels and their rows.
  Each entry in~\eqref{E:variance-normalized-kernel} is a product of two
  continuous wave factors times $(\langle\xi\rangle\langle\eta\rangle)^{-1/2}$,
  with a possible minus sign. Lemma~\ref{L:masked-continuous-wave} therefore
  gives
  \begin{equation}\label{E:variance-kernel-time-bound}
    \|1_I\mathcal{K}_{ij}(\cdot;\xi,\eta)\|_{H_a^{\mathbb{C}}}^2
    \le A_T h^{2-a}\Phi_h(\xi,\eta),
    \qquad A_T = C_1(a)(1+T)^4.
  \end{equation}
  The positive even Radon measure $\nu_R$ satisfies the growth bound
  \eqref{E:variance-actual-growth-constant}. Thus, all hypotheses of
  Proposition~\ref{P:frequency-row} hold, and
  \begin{equation}\label{E:variance-actual-row-bound}
    \int\Phi_h(\xi,\xi-\lambda)\nu_R(\ud\lambda)
    \le K C_\nu B_\alpha(h)\langle\xi\rangle^{-2\gamma}.
  \end{equation}
  No nonnegative physical inverse covariance of $\nu_R$ is needed.

  \medskip
  \noindent\textit{Step 2.~}We check the hypotheses of Hilbert-valued Bessel.
  The kernels used in Lemma~\ref{L:phase-bessel} belong to its stated $L^2$
  space for each fixed $\lambda$. Indeed, $|c_r(t)| \le 1$ and $|b_r(t)| \le
  \sqrt{2}\max(1,T)$ on $[0,T]$, by the high- and low-frequency representations.
  The kernels are continuous into $L^2(0,T)$ in their frequency variables, hence
  into $H_a^{\mathbb{C}}$ by Schur. Separability gives strong measurability.
  Multiplication by shifted Schwartz components of $U$ or $V$ then gives the
  required $L^2$ integrability. Later $\lambda$ integrals are of nonnegative
  functions; their finiteness follows from the estimates below.

  \medskip
  \noindent\textit{Step 3.~}We sum over the input family.
  First fix $V = Jz$ and let $U_k = Je_k$ be the input family.
  In~\eqref{E:variance-normalized-covariance} set
  \[
    g_j^\lambda(\eta) = \sum_{i=1}^2
    1_I\mathcal{K}_{ij}(\cdot;\eta+\lambda,\eta)
    \overline{V_i(\eta+\lambda)}.
  \]
  Lemma~\ref{L:phase-bessel} and Cauchy--Schwarz over the two output components
  give
  \begin{align*}
    \sum_j\|g_j^\lambda(\eta)\|_{H_a^{\mathbb{C}}}^2
    &\le 2\sum_{i,j}|V_i(\eta+\lambda)|^2
    \|1_I\mathcal{K}_{ij}(\cdot;\eta+\lambda,\eta)\|_{H_a^{\mathbb{C}}}^2\\
    &\le 4A_T h^{2-a}\Phi_h(\eta+\lambda,\eta)
    \sum_i|V_i(\eta+\lambda)|^2.
  \end{align*}
  Under $\xi=\eta+\lambda$, the factor $\sum_i|V_i(\xi)|^2$ is independent of
  $\lambda$, leaving the row integral as its coefficient.
  Change variables $\xi = \eta+\lambda$ in the Lebesgue integral for each
  $\lambda$, apply Tonelli and~\eqref{E:variance-actual-row-bound}, and use
  $\|z\|_{E_{s-\gamma}} = \|\langle\xi\rangle^{-\gamma}Jz\|_2$. This
  proves~\eqref{E:wave-input-variance} with $C_T = 4A_TKC_\nu$.

  \medskip
  \noindent\textit{Step 4.~}We sum over the output family.
  Next fix $U = Jz$ and let $V_k = Je_k$ be the output family. After $\eta =
  \xi-\lambda$ in~\eqref{E:variance-normalized-covariance}, put
  \[
    f_i^\lambda(\xi) = \sum_{j=1}^2
    1_I\mathcal{K}_{ij}(\cdot;\xi,\xi-\lambda)U_j(\xi-\lambda).
  \]
  Apply the conjugated version of Lemma~\ref{L:phase-bessel}. Cauchy--Schwarz
  over the two input components bounds its right side by
  \[
    4A_T h^{2-a}\int\Phi_h(\xi,\xi-\lambda)
    \sum_j|U_j(\xi-\lambda)|^2\ud\xi.
  \]
  After the change of variables and Tonelli, the row at fixed $\eta$ is
  $\int\Phi_h(\eta+\lambda,\eta)\nu_R(\ud\lambda)$. Reflect $\lambda$ using the
  evenness of $\nu_R$ and then use symmetry of $\Phi_h$. This row becomes
  $\int\Phi_h(\eta,\eta-\lambda)\nu_R(\ud\lambda)$, which is bounded
  by~\eqref{E:variance-actual-row-bound} with $\xi = \eta$. This proves
  \eqref{E:wave-output-variance} with the same $C_T$.

  \medskip
  \noindent\textit{Step 5.~}We obtain the time-power bounds and their
  uniformity.
  Use Corollary~\ref{C:frequency-time} for $h^{2-a}$ times the row integral at
  the two applications above to obtain the asserted $h^{2\beta}$ bounds,
  including strict absorption of the critical logarithm. The constants are
  uniform in $R,I$ and family size. The $R$-dependent admission estimate was not
  used as this uniform variance bound. This proves
  Proposition~\ref{P:wave-two-variances}.
\end{proof}

\subsection{Mean-square columns on the original probability space}

Write $E_r^{\mathrm{real}}$ for the real phase space in dimension $d$.

\begin{proposition}[Column and adjoint-column maps]\label{P:wave-square-columns}
  On the same probability space there are deterministic bounded real linear maps
  \[
    m_{R,I},n_{R,I}:E_{s-\gamma}^{\mathrm{real}}\longrightarrow
    L^2(\Omega;E_s^{\mathrm{real}}),
  \]
  with squared bounds $C_T h^{2-a}B_\alpha(h)$, such that, for every fixed real
  Schwartz pair $z,w$,
  \begin{equation}\label{E:variance-square-columns}
    \langle m_{R,I}z,w\rangle_E
    = \ell_R(I;z,w) = \langle z,n_{R,I}w\rangle_E
    \quad\text{in }L^2(\Omega).
  \end{equation}
  Fixed-vector values and finite families of values are centered Gaussian.
\end{proposition}

\begin{proof}
  We choose a real Schwartz orthonormal basis and construct the two columns by
  summing their scalar noise coefficients. The two variance bounds give
  convergence and continuous extension. We then identify the pairings and the
  Gaussian laws of the limits.

  \medskip
  \noindent\textit{Step 1.~}We choose the basis and dense input class.
  Choose a countable real Schwartz orthonormal basis $\{e_k\}$ of $E_s$. Such a
  basis follows by Gram--Schmidt from a countable dense Schwartz family,
  discarding zero residuals. Schwartz density follows by smooth Fourier cutoff
  approximation in the weighted $L^2$ norms: approximate by compact smooth
  Fourier coordinates, replace $g(\xi)$ by $(g(\xi)+\overline{g(-\xi)})/2$, and
  apply the inverse Bessel/Fourier maps. The resulting phases are real Schwartz,
  and Gram--Schmidt preserves that class. The same density holds in
  $E_{s-\gamma}$.

  \medskip
  \noindent\textit{Step 2.~}We construct the two columns in mean square.
  For fixed real Schwartz $z$, define
  \[
    m_N(z) = \sum_{k\le N}\ell_R(I;z,e_k)e_k,
    \qquad n_N(z) = \sum_{k\le N}\ell_R(I;e_k,z)e_k.
  \]
  The output-variance estimate controls $m_N(z)$, and the input-variance
  estimate controls $n_N(z)$.
  Both are finite Gaussian sums on the original $W$. By
  Proposition~\ref{P:wave-two-variances}, their coefficient expected squares
  have finite sum. Orthogonality gives the squared tail norms as the
  corresponding sums of expected squares. Thus, both sequences converge in
  $L^2(\Omega;E_s)$ to $m(z),n(z)$ with the asserted bounds, without
  independence of the coefficients. Finite-sum measurability and separable $E_s$
  give strong measurability of the limits. The bounded real linear maps extend
  uniquely from Schwartz phases to $E_{s-\gamma}$.

  \medskip
  \noindent\textit{Step 3.~}We identify the pairings and remove the basis
  choice.
  If $w_N$ is the $E_s$ projection of a Schwartz $w$ onto the first $N$ basis
  vectors, then $\langle m_N(z),w\rangle_E = \ell_R(I;z,w_N)$. The singleton
  form of~\eqref{E:wave-output-variance}, after normalizing $w-w_N$ in $E_s$,
  gives
  \[
    \E|\ell_R(I;z,w-w_N)|^2
    \le C_T h^{2-a}B_\alpha(h)\|z\|_{E_{s-\gamma}}^2\|w-w_N\|_{E_s}^2.
  \]
  The residual is Schwartz, since the projection is a finite linear
  combination of Schwartz basis vectors. If it is nonzero, divide it by its
  $E_s$ norm to make a singleton orthonormal family and apply the variance
  bound. Multiplication by the squared norm gives the displayed estimate.
  The zero residual is immediate. The right side tends to zero, giving the first
  identity in~\eqref{E:variance-square-columns}. For the second, project $z$
  instead and use the singleton input bound \eqref{E:wave-input-variance}. A
  fixed countable Schwartz orthonormal basis and Parseval show basis
  independence of the constructed limits.

  \medskip
  \noindent\textit{Step 4.~}We identify the Gaussian laws after taking limits.
  Continuous linear functionals of the limits, and finite linear combinations of
  such functionals, are $L^2$ limits of centered Gaussian variables with
  convergent variances. Their characteristic functions therefore converge to
  those of centered Gaussians. The same reasoning after the bounded input
  extension proves joint Gaussianity for arbitrary fixed inputs and finite
  families of the two maps. This proves Proposition~\ref{P:wave-square-columns}.
\end{proof}

These are maps into $L^2(\Omega;E_s)$ for fixed intervals. No common uncountable
almost-sure version, random bounded operator, Hilbert--Schmidt realization, time
path, source convergence or new solution-moment bound is inferred by this
section alone. Both variance directions are retained in
Section~\ref{S:actual-operator-driver}.

\section{Finite Gaussian moments and whole-space Sobolev compactness}
\label{S:compact-gaussian}

We prove the finite Gaussian and compact-cutoff estimates needed to construct
the wave operator path. All traces and Schatten norms in this section are
unnormalized. The spatial states remain on $\R^d$, $d \in \{1,2,3\}$. Here $p$
is an auxiliary Schatten and Gaussian-moment index. The auxiliary periodic space
below is used only to factor a compactly supported multiplication operator. No
infinite Gaussian series or time path is constructed in this section.

\subsection{Finite Gaussian matrices with both variance powers}

\LeanMarker{L:gaussian-scalar-powers}
\begin{lemma}[Scalar power estimates]\label{L:gaussian-scalar-powers}
  For $0 < \theta < 1$ and $u,v \ge 0$,
  \begin{equation}\label{E:gaussian-weighted-power}
    u^\theta v^{1-\theta} \le \theta u+(1-\theta)v.
  \end{equation}
  Consequently, if $L > 0$ and $0 \le q \le L$ are integers, then for real
  $x,y$,
  \begin{equation}\label{E:gaussian-integer-power}
    |x|^q|y|^{L-q}
    \le \frac qL|x|^L+\frac{L-q}L|y|^L.
  \end{equation}
  For nonnegative finite sequences, the usual H\"older inequality holds with
  exponents $m$ and $m/(m-1)$, $m > 1$. If $Z \ge 0$ has finite expectation $F$
  and $0 < r < 1$, then $\E Z^r \le F^r$.
\end{lemma}

\begin{proof}
  We first prove the weighted scalar inequality and its integer-power
  consequence. We then use normalization to obtain the two remaining bounds.

  \medskip
  \noindent\textit{Step 1.~}We prove the two scalar inequalities.
  If $v > 0$, divide~\eqref{E:gaussian-weighted-power} by $v$ and put $t = u/v$.
  After division, the desired inequality is $t^\theta \le \theta t+1-\theta$.
  Thus, it suffices to show that the difference of the right and left sides is
  nonnegative.
  The function $\theta t+1-\theta-t^\theta$ has derivative
  $\theta(1-t^{\theta-1})$. Since $\theta-1 < 0$, this derivative is negative
  for $0 < t < 1$ and positive for $t > 1$. The function therefore has minimum
  zero at $t = 1$. Continuity covers $t = 0$; if $v = 0$ the inequality is
  immediate.

  For $0 < q < L$, apply this with $u = |x|^L$, $v = |y|^L$ and $\theta =
  q/L$. The left side becomes $|x|^q|y|^{L-q}$, and the coefficients on the
  right are $q/L$ and $(L-q)/L$. The endpoint cases $q = 0,L$ are identities,
  with the ordinary integer convention $a^0 = 1$.

  \medskip
  \noindent\textit{Step 2.~}We prove the finite sequence inequality.
  To prove finite H\"older, divide the two sequences by their respective
  $\ell^m$ and $\ell^{m/(m-1)}$ norms. If either norm vanishes, the result is
  immediate. Apply~\eqref{E:gaussian-weighted-power} with $\theta = 1/m$ to the
  corresponding normalized powers and sum; the sum of products is at most one.

  \medskip
  \noindent\textit{Step 3.~}We prove the expectation bound.
  Finally, $F = 0$ gives $Z = 0$ almost surely. If $F >
  0$, apply \eqref{E:gaussian-weighted-power} to $u = Z/F$, $v = 1$ and $\theta
  = r$, and integrate. The right side has expectation one.
  This proves
  Lemma~\ref{L:gaussian-scalar-powers}.
\end{proof}

\LeanMarker{L:gaussian-trace-bounds}
\begin{lemma}[Two Hermitian trace bounds]\label{L:gaussian-trace-bounds}
  Let $B,Y$ be finite Hermitian matrices, $m \ge 2$ an integer, and $L = 2m-2$.
  For integers $0 \le q \le L$,
  \begin{equation}\label{E:gaussian-trace-word}
    \operatorname{Re}\operatorname{Tr}(BY^qBY^{L-q})
    \le \operatorname{Tr}(B^2|Y|^L).
  \end{equation}
  If $D \ge 0$ is Hermitian, then
  \begin{equation}\label{E:gaussian-trace-holder}
    \operatorname{Tr}(D|Y|^{2m-2})
    \le (\operatorname{Tr}D^m)^{1/m}
    (\operatorname{Tr}|Y|^{2m})^{(m-1)/m}.
  \end{equation}
\end{lemma}

\begin{proof}
  We work in an eigenbasis of $Y$ throughout. The first bound follows by
  comparing scalar powers; the second follows by finite H\"older.

  \medskip
  \noindent\textit{Step 1.~}We bound the trace word.
  Diagonalize $Y$, with real eigenvalues $\lambda_i$, and write $b_{ij}$ for the
  entries of $B$ in that basis. Hermitian symmetry gives $b_{ji} =
  \overline{b_{ij}}$. The trace word is the real sum
  \[
    \sum_{i,j}|b_{ij}|^2\lambda_j^q\lambda_i^{L-q}.
  \]
  To obtain this sum, the $i$th diagonal entry of the product is summed over
  the intermediate index $j$. The diagonal powers of $Y$ contribute
  $\lambda_j^q$ and $\lambda_i^{L-q}$, and $b_{ij}b_{ji} = |b_{ij}|^2$. All
  these terms are real.
  Bound each product by its absolute value and use
  \eqref{E:gaussian-integer-power}. After interchanging $i,j$, the two resulting
  sums agree because $|b_{ij}| = |b_{ji}|$. Each equals
  $\sum_{i,j}|b_{ij}|^2|\lambda_i|^L = \operatorname{Tr}(B^2|Y|^L)$. The two
  weights sum to one, proving \eqref{E:gaussian-trace-word}. No commutation of
  $B$ and $Y$ is used.

  \medskip
  \noindent\textit{Step 2.~}We bound the trace involving $D$.
  In the same basis put $d_i = D_{ii} \ge 0$. If $\theta_j$ are the nonnegative
  eigenvalues of $D$, a unitary change of basis gives $d_i =
  \sum_j|U_{ij}|^2\theta_j$. Both row and column sums of $|U_{ij}|^2$ equal one.
  Finite convexity therefore yields
  \[
    \sum_i d_i^m \le \sum_{i,j}|U_{ij}|^2\theta_j^m
    = \operatorname{Tr}D^m.
  \]
  The row sums provide the convex weights for each $d_i$. The column sums
  then give the equality after summing over $i$.
  For completeness, the convexity step follows from $t^m \ge s^m+ms^{m-1}(t-s)$
  for $s,t \ge 0$, obtained by differentiating the difference. Average this
  inequality with $s$ equal to the weighted mean.

  Now apply finite H\"older from
  Lemma~\ref{L:gaussian-scalar-powers} to $\sum_i d_i|\lambda_i|^{2m-2}$. Since
  $(2m-2)m/(m-1) = 2m$, we obtain the factors in
  \eqref{E:gaussian-trace-holder}.
  This proves the claimed inequality,
  including zero sums. This proves Lemma~\ref{L:gaussian-trace-bounds}.
\end{proof}

The estimate below is a finite-dimensional Gaussian noncommutative Khintchine
inequality. For its origins in the work of Lust-Piquard and a Gaussian
integration-by-parts proof, see Tropp~\cite{tropp:source2018:second-order}*{arXiv:1504.05919v2,
Proposition~2.3 and Section~7}. We include the proof and Hermitian dilation
to fix the normalization and retain both variance operators.

\LeanMarker{P:finite-gaussian-matrix}
\begin{proposition}[Finite Gaussian matrix moment]\label{P:finite-gaussian-matrix}
  Let $A_1,\ldots,A_N$ be complex $o$-by-$n$ matrices, with $o,n,N \ge 1$, and
  let $g_k$ be independent real standard Gaussians. Put
  \[
    X=\sum_k g_kA_k,\qquad V_{\rm out}=\sum_k A_kA_k^*,
    \qquad V_{\rm in}=\sum_k A_k^*A_k.
  \]
  For even $p = 2m \ge 2$, with $\|A\|_{S_p}^p = \operatorname{Tr}((A^*A)^m)$,
  \begin{align}
    (\E\|X\|_{S_p}^p)^{1/p}
    &\le \sqrt{p-1}\left(
    \frac{\operatorname{Tr}V_{\rm out}^m+
    \operatorname{Tr}V_{\rm in}^m}{2}\right)^{1/p}
    \label{E:finite-gaussian-average}\\
    &\le \sqrt{p-1}\max\left(
    \|V_{\rm out}^{1/2}\|_{S_p},\|V_{\rm in}^{1/2}\|_{S_p}\right).
    \label{E:finite-gaussian-maximum}
  \end{align}
  The constants are independent of $o,n,N$. The first inequality is an equality
  at $p = 2$. Empty or zero families give the zero identity.
\end{proposition}

\begin{proof}
  We first prove the moment bound for Hermitian matrices by Gaussian
  integration by parts. We then use Hermitian dilation to treat rectangular
  matrices and keep track of the factor two.

  \medskip
  \noindent\textit{Step 1.~}We justify Gaussian integration by parts.
  First suppose $B_k$ are Hermitian and put $Y = \sum_k g_kB_k$ and $D = \sum_k
  B_k^2$. All polynomial moments are finite because Gaussian decay dominates
  every fixed power. Componentwise integration by parts against the product
  Gaussian density gives $\E[g_kP(g)] = \E[\partial_kP(g)]$, with zero boundary
  terms. It applies to complex polynomials through their real and imaginary
  parts; absolute polynomial integrability justifies Fubini.

  For a fixed value of the other Gaussian coordinates, the derivative of the
  density in the $k$th coordinate is minus that coordinate times the density.
  Moving this derivative onto $P$ gives the positive sign in the identity. At
  both ends of the real line, the polynomial times the Gaussian density tends
  to zero, so no boundary term remains.

  \medskip
  \noindent\textit{Step 2.~}We estimate the Hermitian moment.
  For $m \ge 2$, put $F = \E\operatorname{Tr}Y^{2m}$ and $L = 2m-2$. Write
  $\operatorname{Tr}Y^{2m} = \sum_k g_k \operatorname{Tr}(B_kY^{2m-1})$ and
  differentiate the $2m-1$ factors in their order. This gives exactly
  \[
    F=\sum_k\sum_{q=0}^{L}
    \E\operatorname{Tr}(B_kY^qB_kY^{L-q}).
  \]
  Here $\partial_kY = B_k$. In the product rule, $q$ counts the factors $Y$
  before the differentiated factor and $L-q$ counts those after it. The
  leading $B_k$ comes from the expression for the trace before
  differentiation. This explains both the order of the matrices and the range
  $0 \le q \le L$.
  The eigenbasis calculation above shows these traces are real. The
  differentiated polynomial is real by conjugation and trace cyclicity. No
  derivative of a spectral absolute value or measurable choice of an eigenbasis
  is needed. The two trace bounds give
  \begin{align*}
    F&\le (2m-1)\E\operatorname{Tr}(D|Y|^{2m-2})\\
    &\le (2m-1)(\operatorname{Tr}D^m)^{1/m}
    \E(\operatorname{Tr}|Y|^{2m})^{(m-1)/m}\\
    &\le (2m-1)(\operatorname{Tr}D^m)^{1/m}F^{(m-1)/m}.
  \end{align*}
  The first inequality bounds each trace word and sums the $2m-1$
  positions of the derivative. The second applies the trace H\"older bound.
  Here $\operatorname{Tr}|Y|^{2m} = \operatorname{Tr}Y^{2m}$, and the last
  inequality is the expectation bound in Lemma~\ref{L:gaussian-scalar-powers}.
  The matrix $D$ is positive because each $B_k^2$ is positive. Also,
  $\operatorname{Tr}|Y|^{2m}$ is nonnegative and has finite expectation $F$.
  Thus, the expectation bound applies with $r = (m-1)/m$, which lies strictly
  between zero and one.
  If $F = 0$ the conclusion is immediate. Otherwise divide by its positive
  fractional power and raise to $m$ to obtain
  \begin{equation}\label{E:finite-hermitian-gaussian}
    \E\operatorname{Tr}Y^{2m} \le (2m-1)^m\operatorname{Tr}D^m.
  \end{equation}
  For $m = 1$, independence gives directly $\E\operatorname{Tr}Y^2 =
  \sum_k\operatorname{Tr}B_k^2 = \operatorname{Tr}D$, proving the same formula
  as equality.

  \medskip
  \noindent\textit{Step 3.~}We pass to rectangular matrices.
  For rectangular $A_k$, use the Hermitian dilation
  \[
    B_k=\begin{pmatrix}0&A_k\\A_k^*&0\end{pmatrix},\qquad
    Y=\begin{pmatrix}0&X\\X^*&0\end{pmatrix},\qquad
    D=\operatorname{diag}(V_{\rm out},V_{\rm in}).
  \]
  Each $B_k$ is Hermitian. Squaring it gives diagonal blocks $A_kA_k^*$ and
  $A_k^*A_k$, so its sum of squares is exactly the displayed $D$. Likewise,
  $Y^2$ has diagonal blocks $XX^*$ and $X^*X$. Its $m$th power therefore has
  blocks $(XX^*)^m$ and $(X^*X)^m$.
  Finite entrywise summation proves $\operatorname{Tr}(UV) =
  \operatorname{Tr}(VU)$ for compatible rectangular matrices. Since $(XX^*)^m =
  X(X^*X)^{m-1}X^*$, it follows that
  \[
    \operatorname{Tr}Y^{2m}=2\operatorname{Tr}((X^*X)^m),\qquad
    \operatorname{Tr}D^m=\operatorname{Tr}V_{\rm out}^m+
    \operatorname{Tr}V_{\rm in}^m.
  \]
  Substitute in~\eqref{E:finite-hermitian-gaussian}, divide by two and take the
  $p$th root. This proves~\eqref{E:finite-gaussian-average}. For positive $V$,
  spectral diagonalization gives $\|V^{1/2}\|_{S_p}^p = \operatorname{Tr}V^m$,
  and the average of two nonnegative numbers is bounded by their maximum. At $p
  = 2$, $\E\operatorname{Tr}X^*X = \operatorname{Tr}V_{\rm in} =
  \operatorname{Tr}V_{\rm out}$, proving the claimed equality. This proves
  Proposition~\ref{P:finite-gaussian-matrix}.
\end{proof}

Both variance directions are needed: the row-matrix example in
No-go~\ref{P:nogo-23} rules out a dimension-independent estimate that
keeps only one of them.

\LeanMarker{R:gaussian-one-variance-test}
\begin{remark}[Both variance directions are necessary]
  \label{R:gaussian-one-variance-test}
  Let $n = 16$ and $A_k = e_k^{\mathsf{T}}$ be the coordinate row matrices. The
  only nonzero singular value of $X$ is $(\sum_k g_k^2)^{1/2}$. Gaussian
  integration by parts gives $\E g^2 = 1$, $\E g^4 = 3$, and hence
  \[
    \E\|X\|_{S_4}^4=3n+n(n-1)=n(n+2)=288.
  \]
  Here $V_{\rm out} = n$ and $V_{\rm in} = I_n$. Keeping only the input variance
  in~\eqref{E:finite-gaussian-maximum} would give fourth power
  $9\operatorname{Tr}I_n^2 = 144$, which is false. Transposing the family
  refutes the corresponding output-only bound. The average and maximum bounds
  have fourth powers $1224$ and $2304$, respectively. This finite example
  asserts nothing about a specialized wave estimate.
\end{remark}

\subsection{Whole-space cutoffs through an auxiliary periodic space}

Use the unitary Fourier transform with $c_d = (2\pi)^{-d/2}$ and $\Lambda =
(1-\Delta)^{1/2}$. For real $r$, the norm on $H^r(\R^d)$ is
$\|\langle\xi\rangle^r\widehat{u}\|_2$. The deterministic cutoff arguments allow
any smoothing order $\gamma > 0$. Fix real $\sigma$ and a real compact smooth
$\chi$. Singular values $s_j$, $j \ge 1$, are counted with multiplicity, with
zeros appended when necessary. For compact $K$, put $\|K\|_{S_p}^p = \sum_j
s_j(K)^p$.

\begin{lemma}[Cutoff factorization]\label{L:cutoff-factorization}
  Choose a cube $Q_L = (-L/2,L/2)^d$ containing $\operatorname{supp}\chi$ in its
  interior, and a real $\psi \in C_c^\infty(Q_L)$ equal to one near
  $\operatorname{supp}\chi$. On the auxiliary torus $\mathbb{T}_L$ use Lebesgue
  measure and the orthonormal basis
  \[
    e_k(x)=L^{-d/2}e^{i\kappa_k\cdot x},\qquad
    \kappa_k=2\pi k/L,\qquad k \in \mathbb{Z}^d,
  \]
  with $\|v\|_{H^r(\mathbb{T}_L)}^2 = \sum_k\langle\kappa_k\rangle^{2r}|v_k|^2$.
  Periodization of $\psi u$ and multiplication of a periodic function by $\chi$
  extend to bounded maps
  \begin{equation}\label{E:cutoff-periodic-maps}
    \mathcal{P}_\psi:H^r(\R^d)\to H^r(\mathbb{T}_L),\qquad
    \mathcal{E}_\chi:H^r(\mathbb{T}_L)\to H^r(\R^d)
  \end{equation}
  for every real $r$. If $\mathcal{I}_\gamma$ denotes the identity inclusion
  from $H^{\sigma+\gamma}(\mathbb{T}_L)$ to $H^\sigma(\mathbb{T}_L)$, then
  \begin{equation}\label{E:cutoff-whole-space-recovery}
    \mathcal{E}_\chi\mathcal{I}_\gamma\mathcal{P}_\psi u=\chi u,
    \qquad u \in H^{\sigma+\gamma}(\R^d).
  \end{equation}
\end{lemma}

\begin{proof}
  We have to show two things. First, the two maps
  in~\eqref{E:cutoff-periodic-maps} are bounded for every real $r$. Second, the
  identity~\eqref{E:cutoff-whole-space-recovery} holds on
  $H^{\sigma+\gamma}(\R^d)$. In both cases we argue first on a dense class,
  namely Schwartz functions on $\R^d$ and finite torus polynomials on
  $\mathbb{T}_L$, and then pass to the whole space by density. The proof has
  five steps. Step~1 collects three elementary facts. Step~2 recalls the mixed
  Schur bound. Step~3 writes the two maps as kernel operators in the Fourier
  variables and bounds them. Step~4 extends the maps by density. Step~5 proves
  the identity~\eqref{E:cutoff-whole-space-recovery}.

  \medskip
  \noindent\textit{Step 1.~}We collect three elementary facts. The first fact
  compares the weights at
  two points. For $x,y \in \R^d$ and real $r$, the elementary inequality
  $\langle x\rangle \le \sqrt{2}\langle y\rangle\langle x-y\rangle$ gives
  \begin{equation}\label{E:cutoff-peetre}
    \left(\frac{\langle x\rangle}{\langle y\rangle}\right)^r
    \le 2^{|r|/2}\langle x-y\rangle^{|r|}.
  \end{equation}
  For $r \ge 0$ we raise the elementary inequality to the power $r$. For
  negative $r$, we interchange $x,y$ before raising to the power $|r|$.

  The second fact is that integration by parts gives rapid decay of the Fourier
  transform of every compact smooth cutoff.

  The third fact is a lattice sum that is bounded uniformly in its shift. For
  any lattice spacing $a > 0$ and $N > d$,
  \begin{equation}\label{E:cutoff-lattice-sum}
    \sup_y\sum_{k\in\mathbb{Z}^d}\langle y-ak\rangle^{-N}<\infty.
  \end{equation}
  Indeed, the sum does not change when $y$ is shifted by a lattice vector. We
  translate $y$ into a centered lattice cell, where $|y| \le a\sqrt d/2$. For
  $|k| \ge \sqrt d$ the triangle inequality gives $|y-ak| \ge a|k|-|y|$, and
  we have $|y-ak| \ge a|k|/2$. We treat the terms with $|k| < \sqrt d$ and the
  tail $|k| \ge \sqrt d$ separately. The remaining terms, those with $|k| <
  \sqrt d$, are finite: each is at most one, and their number depends only on
  $d$. The tail converges by dyadic lattice counting, since there are at most
  $C_d2^{jd}$ lattice points in the $j$th shell $2^j \le |k| < 2^{j+1}$ in
  $k$. More precisely, for $|k| \ge \sqrt d$, the preceding inequality gives
  \[
    \langle y-ak\rangle^{-N}
    \le \max(1,2/a)^N\langle k\rangle^{-N}.
  \]
  Here we use $\langle y-ak\rangle \ge \langle ak/2\rangle \ge
  \min(1,a/2)\langle k\rangle$. The comparison-series terms $\langle
  k\rangle^{-N}$ are bounded by $2^{-jN}$ on each dyadic $k$-shell. Hence, the
  comparison series is bounded by the geometric sum $\sum_{j \ge 0}
  C_d2^{j(d-N)}$, which converges because $N > d$. Thus, the resulting bound
  carries constants depending on $a,d,N$, as
  does~\eqref{E:cutoff-lattice-sum}.

  \medskip
  \noindent\textit{Step 2.~}We recall the mixed Schur bound. We use the mixed
  Schur bound for a
  kernel whose absolute row and column integrals are bounded by $M_1,M_2$. It
  states that the squared operator norm of the kernel operator is at most
  $M_1M_2$. To prove it, we apply the weighted Cauchy--Schwarz inequality at
  each output point and then Tonelli's theorem. The same proof works with
  counting measure on one side and Lebesgue measure on the other. For counting
  measure inputs, use finitely supported sequences; for Lebesgue measure
  inputs, use smooth compactly supported functions. The estimate extends to
  the corresponding $L^2$ spaces by density. For the kernels below, the same
  Cauchy--Schwarz estimate on truncated inputs identifies the extension with
  the stated integral or sum.

  \medskip
  \noindent\textit{Step 3.~}We prove the bounds
  in~\eqref{E:cutoff-periodic-maps}. We write each map as a kernel operator
  between the weighted Fourier variables and apply Step~2. We begin with
  $\mathcal{P}_\psi$. For Schwartz $u$, the $k$th coefficient of its cutoff
  periodization is
  \[
    (\mathcal{P}_\psi u)_k
    = L^{-d/2}\int\psi(x)u(x)e^{-i\kappa_k\cdot x}\ud x
    = L^{-d/2}\int\widehat{\psi}(\kappa_k-\xi)\widehat{u}(\xi)\ud\xi.
  \]
  In the second equality, the $c_d$ from Fourier inversion cancels the
  $c_d^{-1}$ in the ordinary Fourier integral of $\psi$. The norms of
  $H^r(\R^d)$ and $H^r(\mathbb{T}_L)$ carry the weights $\langle\xi\rangle^r$
  and $\langle\kappa_k\rangle^r$. After weighting, that is, as a map from
  $\langle\xi\rangle^r\widehat{u}(\xi)$ to
  $\langle\kappa_k\rangle^r(\mathcal{P}_\psi u)_k$, its exact kernel is
  \[
    L^{-d/2}\left(\frac{\langle\kappa_k\rangle}{\langle\xi\rangle}\right)^r
    \widehat{\psi}(\kappa_k-\xi).
  \]
  We turn to $\mathcal{E}_\chi$. For a finite torus polynomial $v$, the
  normalized Fourier kernel of $\mathcal{E}_\chi v$ is
  \[
    L^{-d/2}\left(\frac{\langle\xi\rangle}{\langle\kappa_k\rangle}\right)^r
    \widehat{\chi}(\xi-\kappa_k).
  \]
  There is no additional $c_d$ in this kernel.

  We now bound both kernels. By~\eqref{E:cutoff-peetre}, each is bounded in
  absolute value by a rapidly decreasing function of its frequency difference,
  times a fixed constant. For the first kernel this bound is
  \[
    L^{-d/2}2^{|r|/2}\langle\kappa_k-\xi\rangle^{|r|}
    \left|\widehat{\psi}(\kappa_k-\xi)\right|,
  \]
  and for the second kernel it is the same expression with
  $\widehat{\chi}(\xi-\kappa_k)$ in place of $\widehat{\psi}(\kappa_k-\xi)$.
  These functions of the frequency difference decrease rapidly by the second
  fact of Step~1. One Schur quantity is its $L^1$ integral and the other a
  uniformly translated lattice sum, bounded by~\eqref{E:cutoff-lattice-sum}.
  Indeed, the integral over $\xi$ at fixed $k$ is the $L^1$ integral of the
  bounding function and does not depend on $k$. The sum over $k$ at fixed
  $\xi$ is a lattice sum with spacing $a = 2\pi/L$ and shift $y = \xi$. We
  bound the rapidly decreasing function by a constant times
  $\langle\cdot\rangle^{-N}$ for some $N > d$,
  and~\eqref{E:cutoff-lattice-sum} bounds the sum uniformly in $\xi$. This
  proves~\eqref{E:cutoff-periodic-maps}, including negative $r$.

  \medskip
  \noindent\textit{Step 4.~}We extend the two maps by density. We check that
  the two classes used in
  Step~3 are dense. Smooth compact approximation in weighted Fourier $L^2$
  proves Schwartz density on $\R^d$. Weighted sequence truncation proves
  polynomial density on the torus. Hence, these maps extend uniquely. Their
  distributional meanings follow by testing against compact smooth functions
  and density: each periodic coefficient pairs the distribution with a compact
  smooth test. No point sampling of a rough Fourier function is used. The same
  full-Lebesgue Schur argument, with Lebesgue measure on both sides, proves
  smooth compact multiplication bounded on every whole-space $H^r$.

  \medskip
  \noindent\textit{Step 5.~}We prove the
  identity~\eqref{E:cutoff-whole-space-recovery}. We first
  prove the identity for Schwartz $u$. For Schwartz $u$, the periodic Fourier
  series of $\psi u$ converges smoothly on the cube. Its coefficients are the
  numbers $(\mathcal{P}_\psi u)_k$ of Step~3, and its partial sums are finite
  torus polynomials, on which $\mathcal{E}_\chi$ is multiplication by $\chi$.
  Rapid decay of the Fourier coefficients gives convergence of these partial
  sums in $H^{\sigma+\gamma}(\mathbb{T}_L)$, hence also in
  $H^\sigma(\mathbb{T}_L)$. Boundedness of $\mathcal{E}_\chi$ gives convergence
  of their images to $\mathcal{E}_\chi\mathcal{I}_\gamma\mathcal{P}_\psi u$
  in $H^\sigma(\R^d)$. Their smooth convergence on the support of $\chi$ also
  gives convergence of the products with $\chi$ to $\chi\psi u$ as
  distributions. Since Sobolev convergence implies distributional
  convergence for every real $\sigma$, these two limits agree.
  Since $\chi\psi = \chi$ and $\chi$ is supported strictly inside the cube,
  multiplying this series by $\chi$ recovers $\chi u$ with no alias term.
  Since $\mathcal{I}_\gamma$ is the identity inclusion, this
  is~\eqref{E:cutoff-whole-space-recovery} for Schwartz $u$.

  It remains to pass to all of $H^{\sigma+\gamma}(\R^d)$. Both sides
  of~\eqref{E:cutoff-whole-space-recovery} are bounded from
  $H^{\sigma+\gamma}(\R^d)$ to $H^\sigma(\R^d)$, using the just-proved maps,
  the inclusion and smooth multiplication. In detail, the left side uses
  $\mathcal{P}_\psi$ with $r = \sigma+\gamma$, the inclusion
  $\mathcal{I}_\gamma$, which is bounded because $\gamma > 0$, and
  $\mathcal{E}_\chi$ with $r = \sigma$. The right side uses the multiplication
  by $\chi$ of Step~4 with $r = \sigma$, together with the inclusion of
  $H^{\sigma+\gamma}(\R^d)$ in $H^\sigma(\R^d)$. Since Schwartz functions are
  dense in $H^{\sigma+\gamma}(\R^d)$ by Step~4, density
  proves~\eqref{E:cutoff-whole-space-recovery} on that full domain. This
  proves Lemma~\ref{L:cutoff-factorization}.
\end{proof}

\begin{proposition}[Whole-space cutoff singular values]\label{P:cutoff-singular-values}
  For real $\sigma$ and $\gamma > 0$, the initially Schwartz-defined operator
  $K_{\chi,\sigma,\gamma} = \Lambda^\sigma\chi\Lambda^{-\sigma-\gamma}$ extends
  uniquely to a compact operator on $L^2(\R^d)$, and
  \begin{equation}\label{E:cutoff-singular-values}
    s_j(K_{\chi,\sigma,\gamma})
    \le C_{\chi,\sigma,\gamma,d}j^{-\gamma/d},\qquad j \ge 1.
  \end{equation}
  Consequently, $K_{\chi,\sigma,\gamma} \in S_p$ whenever $p \ge 1$ and $p\gamma
  > d$.
\end{proposition}

\begin{proof}
  We use the cutoff factorization to approximate the operator by finite-rank
  operators. Their ranks and approximation errors will give the singular-value
  bound, which we then sum to obtain Schatten membership.

  \medskip
  \noindent\textit{Step 1.~}We construct finite-rank approximations.
  The zero cutoff is immediate. Fix one cube and cutoff pair $L,\psi$ as in
  Lemma~\ref{L:cutoff-factorization}. Its recovery identity gives
  \[
    K_{\chi,\sigma,\gamma}
    =\Lambda^\sigma\mathcal{E}_\chi\mathcal{I}_\gamma
    \mathcal{P}_\psi\Lambda^{-\sigma-\gamma}.
  \]
  The cutoff maps are bounded by the preceding lemma. Thus, the only
  factor we need to approximate is the torus inclusion.
  The rightmost Bessel map takes $L^2(\R^d)$ unitarily to
  $H^{\sigma+\gamma}(\R^d)$, and the leftmost takes $H^\sigma(\R^d)$ unitarily
  to $L^2(\R^d)$. The bounded factorization agrees with the original operator
  on Schwartz functions. Their density gives the unique bounded extension.
  In the weighted torus coordinates, $\mathcal{I}_\gamma$ is diagonal with
  entries $\langle\kappa_k\rangle^{-\gamma}$. Let $Q_n$ retain $|k|_\infty \le
  n$. Its rank is $(2n+1)^d$, and
  \[
    \|(1-Q_n)\mathcal{I}_\gamma\|
    \le [1+(2\pi/L)^2(n+1)^2]^{-\gamma/2}\longrightarrow0.
  \]
  The weighted input coefficient has weight
  $\langle\kappa_k\rangle^{\sigma+\gamma}$ and the output coefficient has
  weight $\langle\kappa_k\rangle^\sigma$. The output weight divided by
  the input weight is
  $\langle\kappa_k\rangle^{-\gamma}$. Outside the retained cube, at least one
  integer coordinate has absolute value at least $n+1$, so $|\kappa_k| \ge
  (2\pi/L)(n+1)$. Taking the supremum of the remaining diagonal entries proves
  the norm bound.
  Replacing the inclusion by $Q_n\mathcal{I}_\gamma$ gives finite-rank $K_n$
  converging in norm to $K_{\chi,\sigma,\gamma}$. This proves compactness:
  finite-dimensional images have finite nets of each positive radius, and norm
  approximation transfers those nets to the limiting image of the unit ball.

  \medskip
  \noindent\textit{Step 2.~}We compare approximation rank with
  singular-value index.
  Put
  \[
    M=\|\mathcal{E}_\chi:H^\sigma(\mathbb{T}_L)\to H^\sigma(\R^d)\|
    \|\mathcal{P}_\psi:H^{\sigma+\gamma}(\R^d)
    \to H^{\sigma+\gamma}(\mathbb{T}_L)\|.
  \]
  The outer Bessel maps are unitary. For $j \ge 3^d$, choose $n = \lfloor
  j^{1/d}/4\rfloor$. Then
  \[
    (2n+1)^d\le(j^{1/d}/2+1)^d<j,\qquad
    \|K-K_n\|\le M(2L/\pi)^\gamma j^{-\gamma/d}.
  \]
  The strict rank inequality lets us use $K_n$ to bound the $j$th singular
  value. The finitely many smaller indices will be treated by the operator
  norm.
  The second estimate follows from $n+1 > j^{1/d}/4$. For a compact operator,
  every rank less than $j$ approximation $F$ satisfies $s_j(K) \le \|K-F\|$: the
  span of the first $j$ right singular vectors contains a unit vector in $\ker
  F$, on which $\|K\cdot\| \ge s_j(K)$. If $s_j(K) = 0$ the inequality is
  immediate. For $j < 3^d$ use $s_j(K) \le \|K\| \le M$. Thus, one may take
  \[
    C_{\chi,\sigma,\gamma,d}
    =M\max(3^\gamma,(2L/\pi)^\gamma).
  \]

  \medskip
  \noindent\textit{Step 3.~}We sum the singular-value bound.
  Finally, $\sum_j j^{-p\gamma/d}$ converges for $p\gamma > d$.
  Raise the singular-value bound to $p$ and sum over $j$. The resulting series
  bounds $\sum_j s_j(K)^p$, which is the defining Schatten norm power. This
  proves
  Proposition~\ref{P:cutoff-singular-values}.
\end{proof}

No uniform support-radius constant or endpoint necessity is asserted in
Proposition~\ref{P:cutoff-singular-values}. The periodic inclusion was only a
factor of the recovered whole-space operator.

\subsection{Balanced adjoints and localized variance envelopes}

For $E = H^{1/2}(\R^d)\times H^{-1/2}(\R^d)$ use the physical unitary
coordinates
\[
  J_0(u,v)=(\Lambda^{1/2}u,\Lambda^{-1/2}v),\qquad
  H=L^2(\R^d;\mathbb{C}^2).
\]
The frequency map $J$ in Section~\ref{S:wave-variance} is
$(\mathcal{F}\oplus\mathcal{F})J_0$. Let $D_\chi$ multiply both phase entries
and define
\begin{align}
  P_\chi&=J_0D_\chi J_0^{-1},\qquad
  D_\gamma=\operatorname{diag}(\Lambda^{-\gamma},\Lambda^{-\gamma}),
  \label{E:compact-balanced-cutoff}\\
  K_{\rm out}&=P_\chi D_\gamma,\qquad
  K_{\rm in}=P_\chi^*D_\gamma.
  \label{E:compact-two-factors}
\end{align}

\begin{proposition}[Both compact envelopes]\label{P:compact-two-envelopes}
  The operator $P_\chi$ is bounded. For the unitary swap $U(v,w) = (w,v)$,
  \[
    P_\chi^*=UP_\chi U,\qquad K_{\rm in}=UK_{\rm out}U.
  \]
  Both factors are compact and have equal $S_p$ norms for $p \ge 1$, $p\gamma >
  d$. For $p \ge 2$, $p\gamma > d$, the positive compact operators
  \[
    B_{\rm out}=P_\chi D_\gamma^2P_\chi^*=K_{\rm out}K_{\rm out}^*,
    \qquad B_{\rm in}=P_\chi^*D_\gamma^2P_\chi=K_{\rm in}K_{\rm in}^*
  \]
  satisfy
  \begin{equation}\label{E:compact-envelope-traces}
    \operatorname{Tr}B_{\rm out}^{p/2}
    =\operatorname{Tr}B_{\rm in}^{p/2}=\|K_{\rm out}\|_{S_p}^p.
  \end{equation}
\end{proposition}

\begin{proof}
  We first identify the adjoint of the balanced cutoff. The swap identity
  then reduces both factors to the same singular values, and these values
  give the trace identities.

  \medskip
  \noindent\textit{Step 1.~}We identify the adjoint.
  Put $A_\sigma = \Lambda^\sigma\chi\Lambda^{-\sigma}$. The full-Lebesgue Schur
  argument in Lemma~\ref{L:cutoff-factorization} shows boundedness for every
  real $\sigma$. On Schwartz functions, self-adjointness of real multiplication
  and of the Bessel powers gives $A_\sigma^* = A_{-\sigma}$. Density and
  boundedness extend this to the Hilbert adjoints. No evenness of $\chi$ is
  needed. Thus,
  \[
    P_\chi=\operatorname{diag}(A_{1/2},A_{-1/2}),\qquad
    P_\chi^*=\operatorname{diag}(A_{-1/2},A_{1/2})=UP_\chi U.
  \]

  \medskip
  \noindent\textit{Step 2.~}We obtain compactness and Schatten membership.
  The swap commutes with $D_\gamma$, giving the claimed factor identity. The
  blocks of $K_{\rm out}$ are $K_{\chi,1/2,\gamma}$ and $K_{\chi,-1/2,\gamma}$.
  Proposition~\ref{P:cutoff-singular-values} gives compactness and $S_p$
  membership. The singular values of their direct sum form the multiset union,
  so
  \[
    \|K_{\rm out}\|_{S_p}^p
    =\|K_{\chi,1/2,\gamma}\|_{S_p}^p+
    \|K_{\chi,-1/2,\gamma}\|_{S_p}^p.
  \]
  For a block diagonal operator, its adjoint times itself is block diagonal as
  well. The positive eigenvalues are therefore those of the two blocks,
  counted with multiplicity. Taking their square roots gives the asserted
  multiset union of singular values.
  Unitary conjugation preserves them.

  \medskip
  \noindent\textit{Step 3.~}We compute the two traces.
  The products $K_{\rm out}K_{\rm out}^*$ and $K_{\rm in}K_{\rm in}^*$ are
  positive and compact. Their positive eigenvalues are the squared singular
  values of the corresponding factors.
  Finally, the positive eigenvalues of
  $KK^*$ and $K^*K$ agree with multiplicity: on an eigenspace with eigenvalue
  $\lambda > 0$, the maps $K/\sqrt\lambda$ and $K^*/\sqrt\lambda$ are inverse
  isometries. Indeed, $K$ takes the $\lambda$-eigenspace of $K^*K$ into that
  of $KK^*$, because applying $KK^*$ after $K$ equals applying $K$ after
  $K^*K$. On this eigenspace, $K$ multiplies squared norms by $\lambda$. The
  corresponding statement for $K^*$ gives the inverse map after division by
  $\sqrt\lambda$.

  Zero eigenvalues contribute nothing to the trace. Raising each squared
  singular value to $p/2$ gives its $p$th power. The sum is finite by the
  Schatten membership already proved, and the two factors have the same
  singular values. Thus, summing gives~\eqref{E:compact-envelope-traces}. This
  proves
  Proposition~\ref{P:compact-two-envelopes}.
\end{proof}

We use the usual positive square root of a bounded positive operator and the
spectral decomposition of compact positive operators. The needed order argument
is recorded explicitly.

\begin{lemma}[Compact positive order and trace]\label{L:compact-positive-order}
  If $A,B$ are bounded self-adjoint operators with $0 \le A \le B$ and $B$
  compact positive, then $A$ is compact and its ordered eigenvalues satisfy
  $\lambda_j(A) \le \lambda_j(B)$, including zero padding. Consequently, for $t
  \ge 1$ with finite $\operatorname{Tr}B^t$,
  \begin{equation}\label{E:compact-positive-trace}
    \operatorname{Tr}A^t\le\operatorname{Tr}B^t.
  \end{equation}
\end{lemma}

\begin{proof}
  We first prove compactness using the square root of $A$. Once compactness is
  known, we compare Rayleigh quotients to order the eigenvalues, and then sum
  their powers.

  \medskip
  \noindent\textit{Step 1.~}We prove compactness.
  Let $Q_n$ be finite spectral projections of $B$ whose tails tend to zero in
  norm. Then
  \[
    \|(1-Q_n)A^{1/2}\|^2
    =\|(1-Q_n)A(1-Q_n)\|
    \le\|(1-Q_n)B(1-Q_n)\|\longrightarrow0.
  \]
  The equality is the Hilbert-space identity that the squared norm of an
  operator equals the norm of its product with its adjoint. The inequality
  follows by applying $A \le B$ between the two projections. These products
  are positive, and order between positive operators implies order of their
  norms.
  The finite-rank approximants here are $Q_nA^{1/2}$; the displayed norm
  is precisely their error.
  Thus, $A^{1/2}$ is a norm limit of finite-rank operators and is compact, as is
  $A$, since $A = A^{1/2}A^{1/2}$ is a product of a compact and a
  bounded operator.

  \medskip
  \noindent\textit{Step 2.~}We compare the eigenvalues and traces.
  For compact positive operators, the max-min formula for $\lambda_j$
  follows in an eigenbasis: the span of the first $j$ eigenvectors gives the
  lower bound for the supremum over $j$-dimensional subspaces of their minimum
  Rayleigh quotient, while every such subspace meets the orthogonal complement
  of the first $j-1$ eigenvectors, giving the upper bound. Compare the Rayleigh
  quotients of $A$ and $B$ to obtain the eigenvalue inequality.
  On each fixed $j$-dimensional subspace, every unit vector has its Rayleigh
  quotient for $A$ bounded by that for $B$. Taking the infimum within that
  subspace, and then the supremum over such subspaces, preserves the
  inequality. The eigenbasis argument also covers zero eigenvalues: if fewer
  than $j$ positive eigenvalues exist, every such subspace meets their
  orthogonal complement and the infimum there is zero.

  Since taking a positive scalar power preserves order, the eigenvalue
  comparison gives $\lambda_j(A)^t \le \lambda_j(B)^t$. Sum first
  over finitely many indices and then pass to the increasing limit. The finite
  trace of $B^t$ bounds that limit. This proves
  \eqref{E:compact-positive-trace}. No inequality $A^t
  \le B^t$ is asserted. This proves Lemma~\ref{L:compact-positive-order}.
\end{proof}

\begin{proposition}[Localized finite variance bounds]\label{P:localized-finite-variance}
  Suppose finitely many bounded operators $Z_\ell$ on $H$ satisfy
  \begin{equation}\label{E:localized-finite-hypotheses}
    P_\chi Z_\ell=Z_\ell=Z_\ell P_\chi,\qquad
    V_{\rm out}:=\sum_\ell Z_\ell Z_\ell^*\le CD_\gamma^2,\qquad
    V_{\rm in}:=\sum_\ell Z_\ell^*Z_\ell\le CD_\gamma^2,
  \end{equation}
  with $C \ge 0$. If $p \ge 2$ and $p\gamma > d$, every $Z_\ell$ is compact in
  $S_p$, the two variances are compact in $S_{p/2}$, and
  \begin{align}
    V_{\rm out}&\le CB_{\rm out},\qquad V_{\rm in}\le CB_{\rm in},
    \label{E:localized-two-envelopes}\\
    \operatorname{Tr}V_{\rm out}^{p/2},\quad
    \operatorname{Tr}V_{\rm in}^{p/2}
    &\le C^{p/2}\|K_{\rm out}\|_{S_p}^p.
    \label{E:localized-trace-bounds}
  \end{align}
\end{proposition}

\begin{proof}
  We first use the localization identities to replace the two variance
  bounds by compact envelopes. We then apply the same positive-order argument
  to each individual operator.

  \medskip
  \noindent\textit{Step 1.~}We bound the two variances.
  Left localization gives $V_{\rm out} = P_\chi V_{\rm out}P_\chi^*$. Apply
  positive congruence to the output bound
  in~\eqref{E:localized-finite-hypotheses} to get $V_{\rm out} \le CP_\chi
  D_\gamma^2P_\chi^* = CB_{\rm out}$. Right localization similarly gives $V_{\rm
  in} = P_\chi^*V_{\rm in}P_\chi \le CB_{\rm in}$.
  Positive congruence preserves order: one tests the original inequality
  on the adjoint of the multiplying operator applied to each vector.
  The exponent $t=p/2$ is at least one, and $p\gamma>d$ makes the envelope
  trace powers finite; scaling an envelope by $C$ multiplies its trace power
  by $C^{p/2}$.
  Proposition~\ref{P:compact-two-envelopes} and
  Lemma~\ref{L:compact-positive-order} prove compactness and the two trace
  bounds with $t = p/2$. If $C = 0$, both variances and every $Z_\ell$ vanish,
  so all assertions are immediate.

  \medskip
  \noindent\textit{Step 2.~}We obtain compactness of each summand.
  Every summand of the output variance is positive. Subtracting one summand
  leaves a positive sum, which gives the order comparison needed below.
  For each $\ell$, $Z_\ell Z_\ell^* \le V_{\rm out}$ is compact by the same
  lemma. Its finite spectral projections $Q_n$ give $\|(1-Q_n)Z_\ell\|^2 =
  \|(1-Q_n)Z_\ell Z_\ell^*(1-Q_n)\| \to 0$. Hence, $Z_\ell$ is compact. Its
  singular values squared are the nonzero eigenvalues of $Z_\ell Z_\ell^*$.
  The approximating operators are $Q_nZ_\ell$, whose ranges are finite
  dimensional. The displayed error tends to zero, proving their norm
  convergence. The positive-order lemma also bounds each eigenvalue of
  $Z_\ell Z_\ell^*$ by the corresponding one of $V_{\rm out}$. Taking
  $p/2$ powers converts the former to $p$th powers of singular values.
  Comparing with $V_{\rm out}$ and summing proves $Z_\ell \in S_p$ with
  $\|Z_\ell\|_{S_p}^p \le \operatorname{Tr}V_{\rm out}^{p/2}$. This proves
  Proposition~\ref{P:localized-finite-variance}.
\end{proof}

The two localization identities are hypotheses here; a noise coefficient family
has not yet been shown to satisfy them. In particular, $P_\chi$ is not assumed
self-adjoint or a projection.

\begin{corollary}[Finite Gaussian operators]\label{C:finite-gaussian-operators}
  Under~\eqref{E:localized-finite-hypotheses}, let $g_\ell$ be independent real
  standard Gaussians. For even $p \ge 2$ with $p\gamma > d$, the finite sum $X =
  \sum_\ell g_\ell Z_\ell$ satisfies
  \begin{equation}\label{E:finite-gaussian-operator-bound}
    (\E\|X\|_{S_p}^p)^{1/p}
    \le\sqrt{p-1}\sqrt C\|K_{\rm out}\|_{S_p}.
  \end{equation}
  Its extended norm is measurable, with finite expectation as displayed.
\end{corollary}

\begin{proof}
  We compress the operators to finite matrices, where the Gaussian moment
  bound is already available. We then pass to the original operators using
  convergence of singular values and two applications of Fatou's lemma.

  \medskip
  \noindent\textit{Step 1.~}We bound the variances of the compressions.
  Choose finite-rank orthogonal projections $\Pi_n$ from a countable orthonormal
  basis of $H$, with $\Pi_n \to I$ strongly. Put $Z_{\ell,n} = \Pi_nZ_\ell\Pi_n$
  and $X_n = \sum_\ell g_\ell Z_{\ell,n}$, viewed as matrices on
  $\operatorname{ran}\Pi_n$. Their output variance satisfies
  \[
    \sum_\ell\Pi_nZ_\ell\Pi_nZ_\ell^*\Pi_n
    \le\Pi_nV_{\rm out}\Pi_n,
  \]
  since the difference is $\sum_\ell\Pi_nZ_\ell(1-\Pi_n)Z_\ell^*\Pi_n \ge 0$.
  Each term in this difference is positive: its quadratic form is the squared
  norm of $(1-\Pi_n)Z_\ell^*\Pi_n$ applied to the test vector. Thus, inserting
  the middle projection can only decrease the variance.
  Lemma~\ref{L:compact-positive-order}, with $t=p/2$, first bounds each
  compressed variance trace by the corresponding trace of $\Pi_nV\Pi_n$; the
  compression comparison below and~\eqref{E:localized-trace-bounds} then bound
  it by $C^{p/2}\|K_{\rm out}\|_{S_p}^p$.
  The input inequality is identical with adjoints interchanged. Compact max-min
  also gives $\lambda_j(\Pi_nV\Pi_n) \le \lambda_j(V)$ for positive compact $V$,
  since Rayleigh subspaces inside $\operatorname{ran}\Pi_n$ form a subset of all
  subspaces. Thus, both variance trace powers of these matrices are bounded by
  $C^{p/2}\|K_{\rm out}\|_{S_p}^p$.

  \medskip
  \noindent\textit{Step 2.~}We apply the finite matrix estimate.
  Proposition~\ref{P:finite-gaussian-matrix} applies to these finite complex
  matrices, the same real independent Gaussians and the even $p$. Taking the
  average of the two already bounded trace powers gives
  \begin{equation}\label{E:finite-gaussian-compression-bound}
    \E\|X_n\|_{S_p}^p
    \le(p-1)^{p/2}C^{p/2}\|K_{\rm out}\|_{S_p}^p.
  \end{equation}
  No factor two is added by this average.

  \medskip
  \noindent\textit{Step 3.~}We pass to the limit and check measurability.
  This passage is pointwise in the Gaussian coordinates until the final
  application of Fatou in probability.
  For each Gaussian coordinate vector, $X$ is a finite sum of compact operators.
  Strong convergence of $\Pi_n$ implies $\Pi_nX\Pi_n \to X$ in operator norm:
  approximate $X$ by a finite-rank $F$, use uniform strong convergence on the
  finite-dimensional ranges of $F$ and $F^*$, and control the remainder using
  $\|\Pi_n\| \le 1$.

  Singular values are norm-Lipschitz. Indeed, their formula
  \[
    s_j(A)=\inf_{\dim L<j}
    \sup_{v\in L^\perp,\ \|v\|=1}\|Av\|
  \]
  follows by choosing $L$ to be the span of the first $j-1$ right singular
  vectors for the upper bound. For the lower bound, if $s_j(A) > 0$, the span
  of the first $j$ right singular vectors contains a unit vector perpendicular
  to any $L$ of dimension less than $j$. On that vector, the norm of $A$ is at
  least $s_j(A)$. If $s_j(A) = 0$, the lower bound is immediate.
  Apply the triangle inequality inside the supremum in this formula, then
  take the infimum over $L$. This gives $s_j(A) \le s_j(B)+\|A-B\|$,
  and the reverse bound by interchanging $A,B$. Hence, $s_j(X_n) \to s_j(X)$ for
  each $j$. Fatou for the nonnegative singular-value sum gives $\|X\|_{S_p}^p
  \le \liminf_n\|X_n\|_{S_p}^p$, allowing extended infinity. These variables are
  measurable as follows. Each $Z_\ell$ is bounded, and there are finitely many
  Gaussian coordinates, so changing those coordinates changes $X$ continuously
  in operator norm. The Lipschitz bound makes each singular value continuous
  in these coordinates. Partial sums of their nonnegative $p$th powers are
  measurable, and their increasing limit is measurable, including when it is
  infinite. Fatou in probability
  and~\eqref{E:finite-gaussian-compression-bound} prove the bound and
  almost-sure finiteness.
  No $S_p$-Bochner convergence is
  assumed. This proves
  Corollary~\ref{C:finite-gaussian-operators}.
\end{proof}

\subsection{A nonzero witness and the need for a compact cutoff}

Proposition~\ref{P:localized-finite-variance} and
Corollary~\ref{C:finite-gaussian-operators} assume
\eqref{E:localized-finite-hypotheses}. We check that these hypotheses admit a
nonzero operator, and then test what is lost without the compact cutoff.

\begin{proposition}[Whole-space cutoff tests]\label{P:cutoff-boundary-tests}
  If $\chi$ equals one on a nonempty open set, there is a nonzero rank-one
  operator satisfying~\eqref{E:localized-finite-hypotheses} with a finite common
  $C$. On the other hand, $\Lambda^{-\gamma}$ on $L^2(\R^d)$ is not compact for
  any $\gamma > 0$.
\end{proposition}

\begin{proof}
  We treat the two assertions separately. A rank-one operator built from test
  functions supported on the plateau provides the witness. Translates of a
  fixed function show why smoothing alone does not give compactness on the
  whole space.

  \medskip
  \noindent\textit{Step 1.~}We construct the localized witness.
  Choose nonzero real $v,f \in C_c^\infty$ supported inside the plateau, and put
  $\widetilde{v} = (\Lambda^{1/2}v,0)$ and $\widetilde{w} =
  (\Lambda^{-1/2}f,0)$. Both are Schwartz. The block formulas give
  $P_\chi\widetilde{v} = \widetilde{v}$ and $P_\chi^*\widetilde{w} =
  \widetilde{w}$. Therefore, $Zx = \langle x,\widetilde{w}\rangle\widetilde{v}$
  satisfies $P_\chi Z = Z = ZP_\chi$. Its two variances are
  $\|\widetilde{w}\|^2\widetilde{v}\otimes\widetilde{v}^*$ and
  $\|\widetilde{v}\|^2\widetilde{w}\otimes\widetilde{w}^*$. For any Schwartz $a$
  and any $x \in H$,
  \[
    |\langle x,a\rangle|
    =|\langle D_\gamma x,D_\gamma^{-1}a\rangle|
    \le\|D_\gamma x\|\|D_\gamma^{-1}a\|.
  \]
  The last inequality is Cauchy--Schwarz. The preceding equality uses
  self-adjointness of the multiplier and its inverse on Schwartz functions.
  The inverse multiplier is applied only to Schwartz $a$, so it belongs to $H$.
  Apply the inequality with $a = \widetilde{v},\widetilde{w}$ to obtain the two
  variance bounds with
  \[
    C=\max\left(
    \|\widetilde{w}\|^2\|D_\gamma^{-1}\widetilde{v}\|^2,
    \|\widetilde{v}\|^2\|D_\gamma^{-1}\widetilde{w}\|^2\right)<\infty.
  \]
  This witness uses a plateau; no nonzero family is inferred for every compact
  cutoff or identified with a source family.

  \medskip
  \noindent\textit{Step 2.~}We test compactness by translation.
  For noncompactness, take $f$ with nonzero Gaussian Fourier transform. Then $g
  = \Lambda^{-\gamma}f$ is a nonzero Schwartz function. The images of translates
  of $f$ are the same translates of $g$. With $T_hg(x) = g(x-h)$, Plancherel
  gives
  \[
    \langle T_hg,g\rangle
    =\int e^{-i\xi\cdot h}|\widehat{g}(\xi)|^2\ud\xi
    =c_d^{-1}\widehat{|\widehat{g}|^2}(h).
  \]
  The Fourier multiplier commutes with translations, so the input
  sequence and its image have the same translation parameters.
  This tends to zero because $|\widehat{g}|^2$ is Schwartz. Choose a large
  spacing so that any two distinct translates in the resulting sequence have
  image distance at least $\|g\|$. The input norms are identical, but the images
  have no Cauchy subsequence. This disproves compactness. This proves
  Proposition~\ref{P:cutoff-boundary-tests}.
\end{proof}

The translation obstruction is restated in No-go~\ref{P:nogo-24}:
frequency decay alone does not replace compact spatial localization.

Thus, the choice $\chi = 1$ cannot be substituted into the compact-cutoff
argument. No necessity at the Schatten endpoint $p\gamma = d$ is asserted. The
results concern finite matrices and conditional finite operator sums.
Coefficient localization, compact variance tails, an infinite Gaussian series
and time regularity are proved in Section~\ref{S:actual-operator-driver}; source
and Young identification are proved in
Sections~\ref{S:source-driver-convergence}--\ref{S:constant-phase-young}. They
are not consequences of the finite comparison alone.

\paragraph{What we now have.}
Variance estimates for an operator and for its adjoint
(Proposition~\ref{P:wave-two-variances}), and bounds for finite Gaussian
operator sums (Corollary~\ref{C:finite-gaussian-operators}). The next chapter
builds the random operator path itself.

\chapter{The original noise and operator paths}\label{S:book-operator-paths}

\begin{figure}
  \centering
  % One common extension: the auxiliary channel vanishes on the original range.
\begin{tikzpicture}[x=1cm,y=1cm,font=\normalsize,
  input/.style={ham/box,font=\normalsize,text width=5.9cm,
  minimum height=1.45cm,inner sep=7pt},
  case/.style={ham/box,font=\normalsize,text width=5.9cm,
  minimum height=2.4cm,inner sep=7pt}]
  \node[input] (w) at (3.5,0)
  {\strut\textbf{Given noise $W$}\\\strut Retained literally};
  \node[input,draw=hamteal,fill=hamteal!4] (j) at (11.3,0)
  {\strut\textbf{$J$ independent of $W$}\\\strut One auxiliary for all cutoffs};
  \node[ham/box,font=\normalsize,text width=13.7cm,
  minimum height=2cm,inner sep=7pt] (iso) at (7.4,-2.2)
  {\strut\textbf{One common Gaussian field}\\[3pt]
  \strut $\mathcal I(g)=W(\Phi^{-1}Pg)+J((I-P)g)$\\[3pt]
  \strut $P$ projects onto the original range $\Phi H_W$};
  \draw[ham/arrow] (w.south)--(w.south |- iso.north);
  \draw[ham/arrow,draw=hamteal] (j.south)--(j.south |- iso.north);

  \node[case] (cutoff) at (3.5,-4.9)
  {\strut\textbf{Positive cutoff widths}\\[5pt]
  \strut $Y_\rho(F)=\mathcal I(Q_\rho\Phi F)$\\[5pt]
  \strut The $J$ term may be nonzero.};
  \node[case,fill=hamblue!7] (limit) at (11.3,-4.9)
  {\strut\textbf{Limit on the original range}\\[5pt]
  \strut $(I-P)\Phi F=0$\\[5pt]
  \strut $\mathcal I(\Phi F)=W(F)$};
  \draw[ham/arrow] (cutoff.north |- iso.south)--(cutoff.north);
  \draw[ham/arrow] (cutoff.east)--(limit.west)
  node[midway,above=5pt] {$\rho\downarrow0$};
  \node[align=center] at (7.4,-6.65)
  {The limit retains the original random variable $W(F)$.};
\end{tikzpicture}
  \caption[The original and auxiliary Gaussian noises]{One coupling for every
    cutoff. The sources may depend on the auxiliary noise $J$, but the limiting
    noise and operator path retain the given $W$. A new independent copy is not
    chosen at each cutoff. Here $P$ projects onto $\Phi H_W$; it is not a time
    cutoff. See \ProofRef{Lemma}{L:source-retaining-isometry} and
    \ProofRef{Proposition}{P:source-canonical-convergence}.}
  \label{Fig:common-noise}
\end{figure}
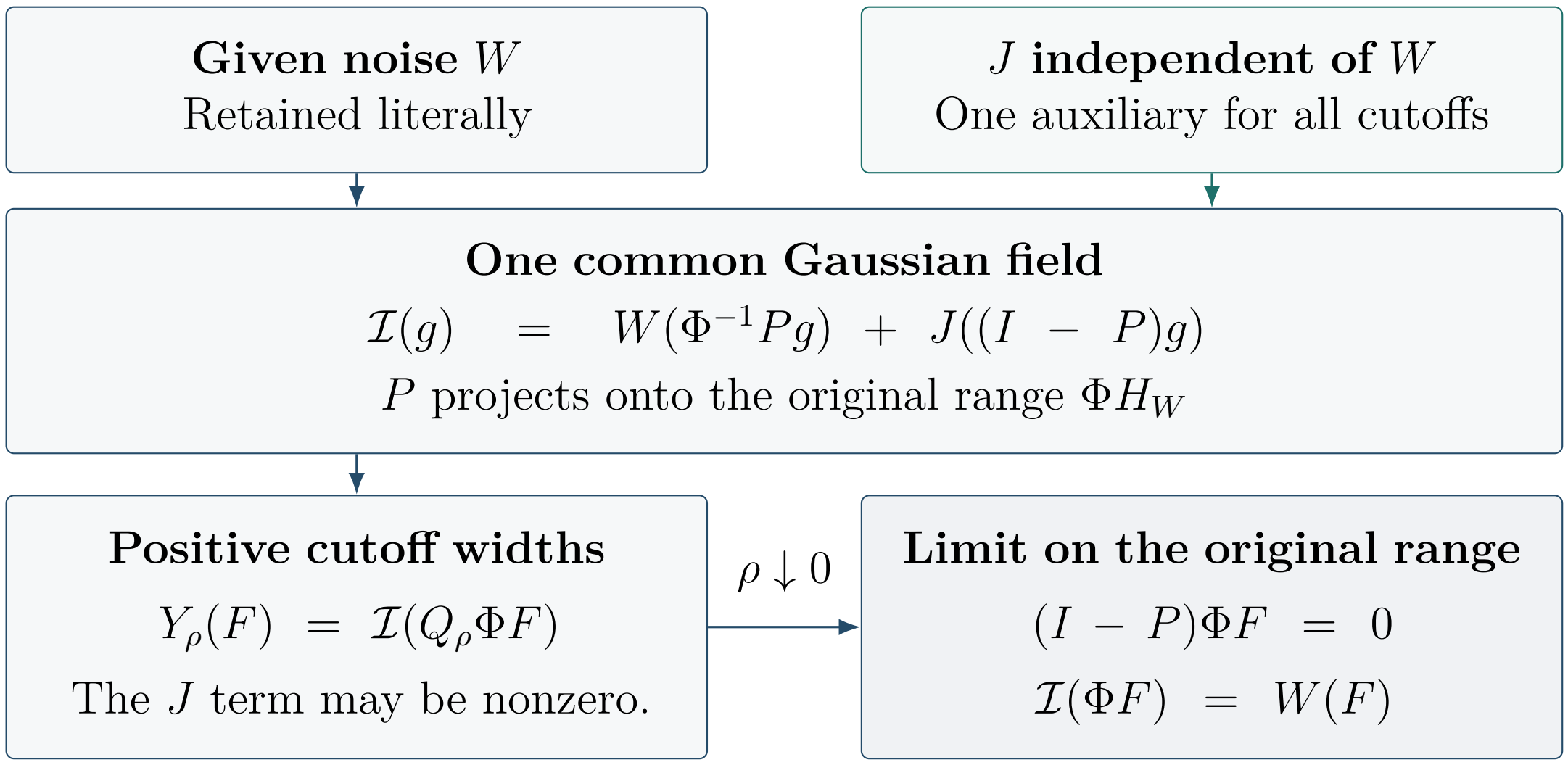

The covariance estimates must now produce something random: an operator path
driven by the noise of the original equation. Another Gaussian field with a
convenient covariance would not do, because the later comparisons are made
on the original noise (Figure~\ref{Fig:common-noise}).

The wave phase is a displacement in $H^{1/2}$ and a velocity
in $H^{-1/2}$. For a fixed spatial mask and a finite time horizon,
Section~\ref{S:actual-operator-driver} builds the integrated interaction path
on this space. Proposition~\ref{P:actual-operator-driver} states its regularity
and how it is approximated. Section~\ref{S:source-driver-convergence} matches
the regularized sources to the original noise and proves convergence in the
path norm. The construction keeps the joint noise information that Gaussian
changes of law will need.

The chapter uses the two variance directions and the compact
Gaussian estimates of Chapter~\ref{S:book-gaussian-foundations}. On
a first reading, keep Propositions~\ref{P:actual-operator-driver}
and~\ref{P:source-gaussian-transfer} and let the proofs wait. The
dependence of the operator-path bounds on the mask radius is tracked in
Chapter~\ref{S:book-growing-radius}. Episode~\VideoRef{5}{2jv8nuqgfBY}
of the video series covers this chapter.

\section{The Gaussian operator path}
\label{S:actual-operator-driver}

We now construct the integrated operator path in the interaction representation
on the same probability space as the original noise. Retain
Hypothesis~\ref{H:homog}, set $a = \alpha_0$, $s = 1/2$ and $E_r =
H^r(\R^d)\times H^{r-1}(\R^d)$, and initially use real phases $E = E_s$. Let
$H_W$ be the real completion of $C_c^\infty((0,\infty)\times\R^d)$ modulo its
covariance-null tests. The given isonormal process on it is $W$. If $\Gamma =
0$, the operator path and all covariance forms are zero, so the conclusions
below hold directly. Assume $\Gamma \ne 0$ in the construction.

Choose and fix $R,T > 0$ and a real even nonnegative $\zeta \in
C_c^\infty(\R^d)$ equal to one on $B_1$ and zero off $B_2$. The physical mask is
$\zeta_R(x) = \zeta(x/R)$. The free group $S(t)$ and forms $\ell_R$ are those of
Section~\ref{S:wave-variance}; both phase entries and the continuous
zero-frequency sine quotient are retained. Fix parameters
\begin{equation}\label{E:driver-parameters}
  \LeanEqMarker{E:driver-parameters}
  0<\gamma<\min\left(\frac12,\frac{2-\alpha}{2}\right),\qquad
  \frac12<\theta<\beta<\min\left(1-\frac a2,
  \frac{3-a-\alpha}{2}\right).
\end{equation}
The strict intervals are nonempty. The $C^\theta$ norm is the supremum norm plus
the $\theta$-H\"older seminorm.

\begin{proposition}[Compact operator path on the original noise space]\label{P:actual-operator-driver}
  There is an orthonormal basis $\{f_j\}$ of $H_W$ with real smooth compact test
  representatives, independent standard real Gaussians $g_j = W(f_j)$, and
  deterministic operators $A_j(0,t)$ on $E$ such that
  \begin{equation}\label{E:driver-path-convergence}
    X_N(t):=\sum_{j\le N}g_jA_j(0,t)\longrightarrow M_R(t)
    \quad\text{in }L^p(\Omega;C^\theta([0,T];\mathcal{B}(E)))
  \end{equation}
  for every finite $p \ge 1$. The limit is a centered Gaussian path on the
  original probability space, starts at zero, and has compact values at every
  time on one probability-one path set. For each fixed real Schwartz pair $z,w$
  and $t \in [0,T]$,
  \begin{equation}\label{E:driver-raw-row}
    \langle M_R(t)z,w\rangle_E=\ell_R((0,t);z,w)
    \quad\text{in }L^2(\Omega).
  \end{equation}
  For long intervals, use the completed covariance test, equivalently the
  additive unit-interval construction justified below. Different compact smooth
  test bases give indistinguishable paths on the same $W$. For constants $C,c >
  0$ depending on the fixed data but not the basis or $p$,
  \begin{equation}\label{E:driver-gaussian-norm}
    (\E\|M_R\|_{C^\theta}^p)^{1/p}\le C\sqrt p,\qquad
    \E\exp(c\|M_R\|_{C^\theta}^2)<\infty.
  \end{equation}
  Its complexification has the corresponding sesquilinear covariance identity.
  Fixed-interval Schatten moments hold at the admissible even indices
  in~\eqref{E:driver-fixed-schatten}. Arbitrary $p$ above refers to
  $\mathcal{B}(E)$, not arbitrary Schatten indices or a Schatten-valued path.
\end{proposition}

\subsection{A test basis and explicit covariance coefficients}

\begin{lemma}[Basis in the prescribed quotient]\label{L:driver-noise-basis}
  The nonzero space $H_W$ is separable and infinite-dimensional. It has a real
  orthonormal basis with compact smooth representatives, whose evaluations by
  $W$ are independent standard real Gaussians.
\end{lemma}

\begin{proof}
  We first obtain a countable dense collection of tests in the prescribed
  covariance norm. We then show that the quotient is infinite-dimensional and
  apply Gram--Schmidt within that quotient.

  \medskip
  \noindent\textit{Step 1.~}We prove separability of the test completion.
  For tests supported in $[0,L]\times K$, with $K$ compact, the physical
  covariance and the temporal Schur bound give the absolute estimate
  \begin{equation}\label{E:driver-test-continuity}
    |\langle F,G\rangle_{H_W}|
    \le |K|\Gamma(K-K)\frac{2L^{1-a}}{1-a}
    \|\|F(t)\|_\infty\|_{L^2(0,L)}
    \|\|G(t)\|_\infty\|_{L^2(0,L)}.
  \end{equation}
  The spatial covariance is bounded by the product of the two spatial supremum
  norms, the volume $|K|$ and the mass $\Gamma(K-K)$. The temporal Schur
  estimate then bounds the integral of this product by the two $L^2$ time
  norms. Thus, uniform approximation on a fixed box implies approximation in
  the covariance seminorm.
  Only spatial differences in $K-K$ contribute. On each rational compact box
  inside positive time, rational step approximations followed by compact smooth
  mollifiers at rational scales give a countable uniformly dense test
  collection, supported in a fixed slightly larger box. The estimate proves
  noise-norm approximation. The union over increasing boxes is dense in the test
  completion. This proves separability without assuming the covariance form is
  injective.

  \medskip
  \noindent\textit{Step 2.~}We show that infinitely many temporal tests
  survive the quotient.
  Choose $L \ge 1$ with $0 < \Gamma(B_L) < \infty$ and nonnegative compact
  smooth $q \ge 1$ on $B_{2L+1}$. Then
  \[
    Q_{\rm sp}:=\int\Gamma(\ud y)\int q(x)q(x-y)\ud x
    \ge |B_{L+1}|\Gamma(B_L)>0,
  \]
  and the compact-support estimate makes it finite. For every real compact
  smooth positive-time $\tau$, the tensor identity is
  \[
    \|\tau\otimes q\|_{H_W}^2=Q_{\rm sp}\|\tau\|_a^2.
  \]
  For $y \in B_L$ and $x \in B_{L+1}$, both $x$ and $x-y$ lie in $B_{2L+1}$.
  Each factor in the integrand is therefore at least one on this region.
  Here the temporal norm is computed on any finite interval containing the
  support of $\tau$. It is strictly positive for signed nonzero $\tau$. Indeed,
  zero-extend this deterministic test and use the positive Laplace factorization
  from Lemma~\ref{L:masked-time-completion} to write
  \begin{align*}
    \|\tau\|_a^2
    &=\frac1{\Gamma_E(a)}\int_0^\infty\lambda^{a-1}
    E_\lambda(\tau)\ud\lambda,\\
    E_\lambda(\tau)
    &=\int_\R 2\lambda\left|\int_r^\infty
    e^{-\lambda(t-r)}\tau(t)\ud t\right|^2\ud r.
  \end{align*}
  Signed exchanges are absolutely integrable on the test support, by the
  temporal Schur estimate. If the norm vanished, some $\lambda > 0$ would have
  $E_\lambda = 0$. Its inner integral $v(r)$ is continuous and satisfies $v' =
  \lambda v-\tau$. Vanishing almost everywhere makes it identically zero and
  forces the derivative $v'$ to vanish as well. The identity $v' = \lambda
  v-\tau$ then
  forces $\tau = 0$. Thus, infinitely many linearly independent compact time
  bumps inject into $H_W$ after multiplication by $q$.

  \medskip
  \noindent\textit{Step 3.~}We construct the basis and its Gaussian
  coordinates.
  Apply Gram--Schmidt to a countable dense real test collection, discarding zero
  residuals in the given null quotient. Infinitude gives an infinite basis; each
  representative is a finite linear combination of compact smooth tests. The
  isonormal process depends only on the quotient class, so a zero residual
  really is a zero noise variable. The surviving normalized vectors have
  covariance equal to their Hilbert inner product. Thus, the
  isonormal identity gives $\E g_jg_k = \delta_{jk}$. Joint Gaussian
  characteristic functions then factor, giving independence and unit variances.
  No physical null direction is removed beyond the prescribed quotient, and no
  $H_W$-valued white Gaussian vector is used. This proves
  Lemma~\ref{L:driver-noise-basis}.
\end{proof}

The operator-norm topology is essential here: in dimensions two and three,
No-go~\ref{P:nogo-21} rules out a compatible square-integrable
Hilbert--Schmidt realization on a same-index Sobolev phase pair,
for suitable coefficient masks and a sufficiently large observation cutoff.

\begin{lemma}[Covariance coefficients]\label{L:driver-covariance-coefficients}
  For each basis representative put
  \begin{equation}\label{E:driver-covariance-coefficient}
    H_{j,R}(t,x)=\zeta_R(x)\int_{u>0}|t-u|^{-a}
    \int\Gamma(\ud y)f_j(u,x-y)\ud u.
  \end{equation}
  These functions are real smooth on finite time-space compacts, supported
  spatially in $B_{2R}$, with bounded derivatives there. For smooth compact real
  positive-time $F$, and also for the admitted wave tests
  of~\eqref{E:variance-actual-form},
  \begin{equation}\label{E:driver-literal-covariance}
    \E[W_{\mathbb{C}}(\zeta_RF)g_j]
    =\int H_{j,R}(t,x)F(t,x)\ud t\ud x.
  \end{equation}
  The coefficient on $t \ge 0$ depends only on the $H_W$ class of $f_j$.
\end{lemma}

\begin{proof}
  We move derivatives onto the smooth basis test to avoid differentiating the
  singular covariance kernel. We then prove the covariance identity by
  approximation and show that the coefficient depends only on the noise class.

  \medskip
  \noindent\textit{Step 1.~}We prove smoothness of the deterministic
  coefficient.
  Extend $f_j$ smoothly by zero in time and substitute $v = t-u$. The unmasked
  expression is $\int_\R |v|^{-a}\int\Gamma(\ud y)f_j(t-v,x-y)\ud v$. On each
  compact $(t,x)$ set, the integrand involves only bounded $v$ and a common
  compact set of $y$. Local Radon mass and local integrability of $|v|^{-a}$
  provide an integrable majorant for each derivative of the test on these
  fixed sets. The singular kernel is left unchanged under differentiation.
  Thus, every mixed derivative can move onto the smooth test. Dominated
  integration proves smoothness and the derivative bounds. Multiplication by
  $\zeta_R$ gives the stated support. Zero extension is only of a deterministic
  test; no negative-time noise is introduced.

  \medskip
  \noindent\textit{Step 2.~}We identify the coefficient for tests and their
  limits.
  The physical covariance formula and absolute Fubini give
  \eqref{E:driver-literal-covariance} first for real compact tests. The compact
  difference set in $\Gamma$ and integrable temporal singularity justify the
  exchange. For wave tests use the common-support original approximations of
  Proposition~\ref{P:actual-wave-form}, converging in $L^2$ time in the spatial
  supremum norm after masking. The left side converges by noise Cauchy--Schwarz
  and the right side by the coefficient's boundedness on that support. Real
  and
  imaginary parts are treated separately because the original noise is real.
  Combining the two identities gives the complex evaluation defined in the
  preceding section. If two representatives give the same
  class, their Gaussian difference is zero. Testing the covariance identity
  against every compact smooth $F$ makes the difference of their masked
  coefficients zero as a distribution for $t > 0$. Smoothness, then continuity
  at zero, proves the representative assertion on $t \ge 0$. This proves
  Lemma~\ref{L:driver-covariance-coefficients}.
\end{proof}

\subsection{Gain, adjoints and localization}

\begin{lemma}[Deterministic coefficient operators]\label{L:driver-gain-localization}
  For $I \subset [0,T]$, the statewise integral
  \begin{equation}\label{E:driver-coefficient-operator}
    A_j(I)z=\int_I S(-t)B(H_{j,R}(t))S(t)z\ud t
  \end{equation}
  satisfies, for every real $r$,
  \begin{equation}\label{E:driver-gain-one}
    \|A_j(I)\|_{\mathcal{B}(E_r,E_{r+1})}
    \le C_{j,r,R,T}|I|.
  \end{equation}
  Its $E_s$ adjoint has the gain $E_{2s-r-1} \to E_{2s-r}$, so both operators
  act from $E_{s-\gamma}$ to $E$. Choose real compact smooth $\chi = 1$ on
  $B_{2R+T+1}$. On every $E_r$,
  \begin{equation}\label{E:driver-actual-localization}
    D_\chi A_j(I)=A_j(I)=A_j(I)D_\chi.
  \end{equation}
  For each fixed Schwartz pair,
  \begin{equation}\label{E:driver-covariance-row}
    \E[g_j\ell_R(I;z,w)]=\langle A_j(I)z,w\rangle_E.
  \end{equation}
\end{lemma}

\begin{proof}
  We establish the gain in Sobolev order and its adjoint version first. Finite
  propagation then gives both localization identities. Finally, the physical
  pairing identifies the Gaussian coefficient of each row.

  \medskip
  \noindent\textit{Step 1.~}We prove the gain and the adjoint gain.
  The full massless wave is strongly continuous and bounded on every $E_r$ on
  finite time intervals. Its normalized off-diagonal Fourier entries have sizes
  at most $C_T$ and one, respectively; continuity and domination prove the
  strong assertion. Multiplication by $H_{j,R}(t)$ is bounded on every $H^r$,
  continuously in $t$, by weighted Fourier convolution and its compact smooth
  derivative bounds. Hence, $B(H_{j,R}(t)):E_r \to E_{r+1}$ is bounded: its
  second entry is a multiplier of the displacement in $H^r$. The integrand
  in~\eqref{E:driver-coefficient-operator} is strongly continuous with values in
  $E_{r+1}$ for each input in $E_r$, and uniformly bounded in that operator
  norm. Statewise Bochner integration proves~\eqref{E:driver-gain-one},
  including the Lipschitz operator-norm primitive $A_j(0,t)$. Norm continuity of
  $S(t)$ is not needed.

  Under the $E_s$ pairing, the dual of $E_r$ is isometrically $E_{2s-r}$,
  directly by the two Bessel weights. Dualizing the bound with a gain of one
  derivative therefore gives the asserted adjoint gain. It agrees on Schwartz
  phases with the bounded $E_s$ adjoint. Since $\gamma < 1$, gain one at every
  index implies the two stated maps from $E_{s-\gamma}$ to $E$.

  \medskip
  \noindent\textit{Step 2.~}We prove both localization identities by finite
  propagation.
  For smooth data the free-wave energy on a shrinking ball has nonpositive
  boundary flux: $v\partial_nu-(v^2+|\nabla u|^2)/2 \le 0$. Vanishing initial
  displacement and velocity give zero velocity and gradient in the cone;
  integration of velocity from the zero displacement gives zero displacement
  there as well. Time reversal gives the backward assertion. Thus, both entries
  of $S(\pm t)$ enlarge supports by at most $|t|$. Density and testing against
  compact smooth functions extend this support assertion to Sobolev phases. In
  particular, for input $(1-D_\chi)z$, approximate $z$ by Schwartz phases and
  retain the factor $1-D_\chi$; the approximants vanish on $B_{2R+T+1}$. Their
  free displacement vanishes on $B_{2R}$ throughout $[0,T]$, and the same is
  true in the distributional limit. They contribute nothing to the integral. Its
  output is supported in $B_{2R+T}$, by the support of $H_{j,R}$ and backward
  finite propagation. Since $\chi = 1$ on a neighborhood of this ball, both
  identities in~\eqref{E:driver-actual-localization} follow. This uses both
  entries, not a scalar restart or a positivity theorem.

  \medskip
  \noindent\textit{Step 3.~}We identify the covariance row.
  Finally, apply the physical Sobolev pairing
  \eqref{E:variance-physical-adjoint} and~\eqref{E:driver-literal-covariance} to
  the state integral. They give \eqref{E:driver-covariance-row}.
  In~\eqref{E:driver-covariance-coefficient}, the mask is after convolution;
  replacing it by a mask inside $\Gamma*f_j$ would change this identity. This
  proves Lemma~\ref{L:driver-gain-localization}.
\end{proof}

\subsection{Finite variances and first Wiener chaos tails}

Put $J_0(u,v) = (\Lambda^{1/2}u,\Lambda^{-1/2}v)$, $H_0 = L^2(\R^d;\R^2)$, $a_j
= J_0A_j(I)J_0^{-1}$ and $P_\chi = J_0D_\chi J_0^{-1}$. The operator $D_\gamma =
\operatorname{diag}(\Lambda^{-\gamma},\Lambda^{-\gamma})$ acts on these physical
coordinates. Its two compact envelopes $B_{\rm out},B_{\rm in}$ are those of
Proposition~\ref{P:compact-two-envelopes}, after complexification when Schatten
norms are used. For $I \subset [0,T]$ with $0 < h = |I| \le 1$, set
\begin{equation}\label{E:driver-short-variance-budget}
  b_I=C_T h^{2-a}B_\alpha(h),
\end{equation}
where $C_T$ is the uniform-in-$R$ constant in
Proposition~\ref{P:wave-two-variances}.

\begin{lemma}[Both variances and their tails]\label{L:driver-actual-variances}
  For every $N$,
  \begin{equation}\label{E:driver-two-variances}
    V_{{\rm out},N}:=\sum_{j\le N}a_ja_j^*\le b_ID_\gamma^2,
    \qquad V_{{\rm in},N}:=\sum_{j\le N}a_j^*a_j\le b_ID_\gamma^2.
  \end{equation}
  The finite families satisfy both localizations by $P_\chi$. There are bounded
  positive limits $V_{{\rm out},\infty},V_{{\rm in},\infty}$ with positive tails
  $T_{{\rm out},N},T_{{\rm in},N} \to 0$ strongly, and
  \begin{equation}\label{E:driver-compact-tail-envelopes}
    T_{{\rm out},N}\le b_IB_{\rm out},\qquad
    T_{{\rm in},N}\le b_IB_{\rm in}.
  \end{equation}
  For each fixed $z \in E_{s-\gamma}$, the real column maps satisfy
  \begin{equation}\label{E:driver-column-expansions}
    m_{R,I}(z)=\sum_jg_jA_j(I)z,\qquad
    n_{R,I}(z)=\sum_jg_jA_j(I)^{*,E}z
    \quad\text{in }L^2(\Omega;E).
  \end{equation}
\end{lemma}

\begin{proof}
  We first compare finite noise projections with the original covariance
  forms. Summing over a phase basis gives the two variances. We then identify
  the column series and use their positive tails to obtain compact envelopes.

  \medskip
  \noindent\textit{Step 1.~}We project scalar rows onto the noise basis.
  By~\eqref{E:driver-covariance-row}, for real Schwartz $z,w$,
  \[
    \langle X_N(I)z,w\rangle_E
    =\sum_{j\le N}g_j\E[g_j\ell_R(I;z,w)]
    =\mathcal{P}_N\ell_R(I;z,w),
  \]
  where $\mathcal{P}_N$ is the orthogonal projection in $L^2(\Omega)$ onto the
  span of $g_1,\ldots,g_N$. In particular,
  \[
    \sum_{j\le N}|\langle A_j(I)z,w\rangle_E|^2
    =\E|\langle X_N(I)z,w\rangle_E|^2
    \le\E|\ell_R(I;z,w)|^2.
  \]
  The variables $g_j$ are orthonormal in $L^2(\Omega)$. Thus, the squared norm
  of this finite orthogonal projection is the sum of the squared coefficients
  and cannot exceed the squared norm of the original row. This is the
  contraction used in both directions below.
  These are projections onto the noise-test basis, not source regularizations.
  The two variance bounds of Proposition~\ref{P:wave-two-variances}, with the
  mean-square columns of Proposition~\ref{P:wave-square-columns}, apply to these
  same covariance forms and the bound~\eqref{E:driver-short-variance-budget}.

  \medskip
  \noindent\textit{Step 2.~}We sum the row bounds in both phase directions.
  Choose a real Schwartz orthonormal basis $\{e_k\}$ of $E$. Fix real Schwartz
  $w$. For every finite phase truncation $K$,
  \[
    \begin{split}
      \sum_{k\le K}\sum_{j\le N}
      |\langle A_j(I)e_k,w\rangle_E|^2
      &=\sum_{k\le K}\E|\langle X_N(I)e_k,w\rangle_E|^2\\
      &\le\sum_{k\le K}\E|\ell_R(I;e_k,w)|^2
      \le b_I\|w\|_{E_{s-\gamma}}^2,
    \end{split}
  \]
  by~\eqref{E:wave-input-variance}. Let $K \to \infty$. Parseval in the phase
  space and nonnegative monotone convergence, with $N$ still finite, give
  \[
    \sum_{j\le N}\|A_j(I)^{*,E}w\|_E^2
    =\E\|X_N(I)^{*,E}w\|_E^2
    \le b_I\|w\|_{E_{s-\gamma}}^2.
  \]
  Conversely, for fixed real Schwartz $z$, use \eqref{E:wave-output-variance}
  with the output basis varying:
  \[
    \begin{split}
      \sum_{k\le K}\sum_{j\le N}
      |\langle A_j(I)z,e_k\rangle_E|^2
      &=\sum_{k\le K}\E|\langle X_N(I)z,e_k\rangle_E|^2\\
      &\le\sum_{k\le K}\E|\ell_R(I;z,e_k)|^2
      \le b_I\|z\|_{E_{s-\gamma}}^2.
    \end{split}
  \]
  The same phase Parseval limit gives
  \[
    \sum_{j\le N}\|A_j(I)z\|_E^2
    =\E\|X_N(I)z\|_E^2
    \le b_I\|z\|_{E_{s-\gamma}}^2.
  \]
  Thus, both directions use their respective Hilbert-Bessel bounds; neither
  follows from retaining only one variance.

  \medskip
  \noindent\textit{Step 3.~}We extend the finite bounds to the weaker input
  space.
  These inequalities extend at each finite $N$ to the full $E_{s-\gamma}$ input
  space. Indeed, Lemma~\ref{L:driver-gain-localization} gives $A_j:E_{s-\gamma}
  \to E_{s+1-\gamma} \subset E$. Its adjoint gain, at $r = s+\gamma-1$, gives
  the same map for $A_j^{*,E}$. The inclusion holds since $\gamma < 1$.
  Approximate the input by Schwartz phases and pass each of the finitely many
  operator terms to the limit. The resulting constant is still $b_I$,
  independent of $N$; individual $j$-dependent gain constants were not summed.
  Restricting to $E$ and using $\|z\|_{E_{s-\gamma}} = \|D_\gamma J_0z\|_{H_0}$
  gives \eqref{E:driver-two-variances} as two operator quadratic inequalities.
  In physical coordinates, the sum of squared adjoint norms is the quadratic
  form of $\sum_{j\le N}a_ja_j^*$, and the sum of squared direct norms is that
  of $\sum_{j\le N}a_j^*a_j$. The weaker input norm is the quadratic form of
  $D_\gamma^2$. This converts the two scalar inequalities into the displayed
  operator bounds.
  Both localizations follow from~\eqref{E:driver-actual-localization}. No
  self-adjointness or projection property of $P_\chi$ is assumed.

  \medskip
  \noindent\textit{Step 4.~}We identify the first-chaos expansions of the
  columns.
  The coordinates of $m_{R,I}(z)$, for Schwartz $z$, are the original covariance
  forms $\ell_R(I;z,e_k)$ and belong to the closed first Wiener chaos $W(H_W)$.
  Their coefficients against $g_j$ are $\langle A_jz,e_k\rangle_E$. Parseval in
  that chaos, followed by Parseval in $E$ and positive Tonelli, gives
  \[
    \sum_j\|A_jz\|_E^2=\E\|m_{R,I}(z)\|_E^2,
    \qquad m_{R,I}(z)=\sum_jg_jA_jz\quad\text{in }L^2(\Omega;E).
  \]
  For each phase coordinate, completeness of the noise basis gives convergence
  in the scalar first chaos. Summing the squared coefficient tails over the
  phase basis gives the vector-valued mean-square error. Positive Tonelli
  permits this summation, and the finite total column variance forces the
  tails to zero.
  The same argument for the adjoint column gives the other expansion. Membership
  in this first chaos comes from linear noise tests, not merely Gaussianity. To
  extend to a fixed $E_{s-\gamma}$ input, approximate by Schwartz phases. The
  finite inequalities extend to that space by the gain bounds and have a
  constant uniform in $N$. Apply them to the input difference, then use the
  Schwartz expansion and continuity of $m,n$. This
  proves~\eqref{E:driver-column-expansions} on the full claimed input domain.

  \medskip
  \noindent\textit{Step 5.~}We pass to positive variance limits and their
  compact envelopes.
  The convergent positive quadratic sums, restricted to $E$ and transported by
  $J_0$, define the bounded positive operators $V_{\cdot,\infty}$ by
  polarization and Riesz. They retain the bounds $b_ID_\gamma^2$. Their tails
  have quadratic forms tending to zero. Since they are uniformly bounded and
  positive, $\|T_Nv\|^2 \le \|T_N\|\langle T_Nv,v\rangle \to 0$. Passing to the
  limit in the finite localization identities gives
  \[
    V_{{\rm out},\infty}=P_\chi V_{{\rm out},\infty}P_\chi^*,\qquad
    V_{{\rm in},\infty}=P_\chi^*V_{{\rm in},\infty}P_\chi,
  \]
  and the same identities for their tails. Positive congruence applied to $T_N
  \le b_ID_\gamma^2$ proves~\eqref{E:driver-compact-tail-envelopes}.

  Complexification preserves these statements and constants. The Hilbert norms
  split into real and imaginary squared norms; real operators extend
  isometrically. A real orthonormal basis complexifies to a complex orthonormal
  basis, preserving singular values and multiplicities. The complex column maps
  have the same bounds and sesquilinear row identities. This proves
  Lemma~\ref{L:driver-actual-variances}.
\end{proof}

Lemma~\ref{L:driver-actual-variances} gives tails
$T_{{\rm out},N},T_{{\rm in},N}$ that tend to zero only strongly. They lie
below the compact envelopes in~\eqref{E:driver-compact-tail-envelopes}. We
now use this compact domination to improve the mode of convergence.

\begin{lemma}[Compact domination of positive tails]\label{L:driver-dominated-tails}
  If $0 \le T_N \le B$, with positive compact $B \in S_m$, $m \ge 1$, and $T_N
  \to 0$ weakly, then
  \begin{equation}\label{E:driver-trace-tails}
    \|T_N\|\longrightarrow0,\qquad \operatorname{Tr}T_N^m\longrightarrow0.
  \end{equation}
\end{lemma}

\begin{proof}
  Choose a finite spectral projection $Q$ of $B$ with $\|(1-Q)B(1-Q)\| <
  \epsilon$. Positive quadratic Cauchy--Schwarz gives, for unit $x$,
  \[
    \langle T_Nx,x\rangle
    \le2\langle T_NQx,Qx\rangle
    +2\langle T_N(1-Q)x,(1-Q)x\rangle
    \le2\|QT_NQ\|+2\epsilon.
  \]
  Weak convergence makes the finite-dimensional matrix $QT_NQ$ converge in norm
  to zero. Letting $\epsilon\downarrow0$ proves norm convergence.
  Lemma~\ref{L:compact-positive-order} gives $\lambda_j(T_N) \le \lambda_j(B)$,
  while each $\lambda_j(T_N) \to 0$ by norm convergence. Dominated summation of
  $\lambda_j(B)^m$ proves the trace limit. No order comparison of the operator
  powers is used. This proves Lemma~\ref{L:driver-dominated-tails}.
\end{proof}

Trace powers of variances must still be turned into moments of Gaussian sums.
Proposition~\ref{P:finite-gaussian-matrix} keeps the average of the two trace
powers, but only for matrices. Corollary~\ref{C:finite-gaussian-operators}
bounds compact sums only through the envelope norm $\|K_{\rm out}\|_{S_p}$.

\begin{lemma}[Retained finite Gaussian variance powers]\label{L:driver-gaussian-tail-average}
  Let $Z_\ell$ be finitely many compact operators whose two variances $V_{\rm
  out},V_{\rm in}$ belong to $S_{q/2}$, for an even $q \ge 2$. With independent
  standard real $g_\ell$,
  \begin{equation}\label{E:driver-retained-gaussian-average}
    \E\left\|\sum_\ell g_\ell Z_\ell\right\|_{S_q}^q
    \le(q-1)^{q/2}\frac{\operatorname{Tr}V_{\rm out}^{q/2}
    +\operatorname{Tr}V_{\rm in}^{q/2}}2.
  \end{equation}
\end{lemma}

\begin{proof}
  We retain the two variance trace powers while compressing to finite
  matrices. Compactness then allows passage to the original finite Gaussian
  sum, using Fatou twice.

  \medskip
  \noindent\textit{Step 1.~}We apply the matrix estimate to the compressed
  sum.
  Compress by finite-rank orthogonal $Q_k \to I$ strongly. The compressed output
  variance is at most $Q_kV_{\rm out}Q_k$, since their difference is $\sum_\ell
  Q_kZ_\ell(1-Q_k)Z_\ell^*Q_k \ge 0$; the input inequality is the same with
  adjoints interchanged. Positive-order and compression max-min compare their
  trace powers with the uncompressed ones.
  Apply
  Proposition~\ref{P:finite-gaussian-matrix} on the finite range, retaining the
  average. For each finite Gaussian coordinate vector, compactness gives
  operator-norm convergence of the compressed sum to the original sum.

  \medskip
  \noindent\textit{Step 2.~}We pass to the uncompressed sum.
  The singular-value norm-Lipschitz argument in the proof of
  Corollary~\ref{C:finite-gaussian-operators}, then Fatou over singular values
  and probability, proves the displayed bound. The extended norm is measurable
  by the same finite-coordinate continuity argument. No infinite Gaussian series
  or $S_q$-Bochner convergence is assumed. This proves
  Lemma~\ref{L:driver-gaussian-tail-average}.
\end{proof}

\subsection{Interval limits on the given probability space}

The two previous lemmas control tails of the series $\sum_jg_jA_j(I)$ through
its two variances. We now fix one interval and use them to pass to the limit in
$N$, without leaving the original probability space.

\LeanMarker{L:driver-interval-limit}
\begin{lemma}[Fixed-interval operator construction]\label{L:driver-interval-limit}
  For each fixed short interval $I$, $X_N(I)$ converges in every finite
  $L^p(\Omega;\mathcal{B}(E))$ to a compact-valued centered Gaussian operator
  $X_R(I)$ on the original probability space. For every admissible even $q \ge
  2$, $q\gamma > d$,
  \begin{equation}\label{E:driver-fixed-schatten}
    (\E\|X_R(I)\|_{S_q}^q)^{1/q}
    \le\sqrt{q-1}\sqrt{C_T}\,h^{1-a/2}B_\alpha(h)^{1/2}
    \|K_\chi\|_{S_q},\qquad K_\chi=P_\chi D_\gamma.
  \end{equation}
  Its fixed Schwartz rows equal the covariance forms in $L^2$.
\end{lemma}

\begin{proof}
  We first show that the finite sums are Cauchy using the compact variance
  tails. We construct their limit on the original probability space, identify
  its rows, and finally extend additively to long intervals.

  \medskip
  \noindent\textit{Step 1.~}We obtain uniform finite moments and vanishing
  tails.
  The finite families $a_j$, $j \le N$, satisfy every hypothesis of
  Proposition~\ref{P:localized-finite-variance} and
  Corollary~\ref{C:finite-gaussian-operators}: both localizations, both weighted
  variances and the independent real $g_j$ have been proved in
  Lemma~\ref{L:driver-actual-variances}. Thus, the coefficients are compact, and
  their finite Gaussian sums satisfy~\eqref{E:driver-fixed-schatten}, uniformly
  in $N$. The fixed cutoff norm may depend on $R,T$.

  Fix an admissible even $q \ge 2$, $q\gamma > d$. Both positive compact
  envelopes in~\eqref{E:driver-compact-tail-envelopes} belong to $S_{q/2}$.
  Lemma~\ref{L:driver-dominated-tails} therefore gives
  \[
    \operatorname{Tr}T_{{\rm out},N}^{q/2}\longrightarrow0,
    \qquad
    \operatorname{Tr}T_{{\rm in},N}^{q/2}\longrightarrow0.
  \]
  For $M > N$, put $Y_{N,M} = \sum_{j=N+1}^Mg_jA_j(I)$. After the unitary $J_0$
  transport, its two finite variances are
  \[
    V_{{\rm out},M}-V_{{\rm out},N}\le T_{{\rm out},N},\qquad
    V_{{\rm in},M}-V_{{\rm in},N}\le T_{{\rm in},N}.
  \]
  Lemma~\ref{L:driver-gaussian-tail-average} and positive eigenvalue max-min
  comparison now give the uniform tail estimate
  \[
    \sup_{M>N}\E\|Y_{N,M}\|_{S_q}^q
    \le(q-1)^{q/2}
    \frac{\operatorname{Tr}T_{{\rm out},N}^{q/2}
    +\operatorname{Tr}T_{{\rm in},N}^{q/2}}2
    \longrightarrow0.
  \]
  For each $M > N$, the two variances of the finite tail lie below the
  corresponding infinite tails. The positive-order comparison bounds their
  trace powers by the right side, which no longer depends on $M$. This is why
  the estimate proves a Cauchy property rather than only boundedness of the
  partial sums.
  The average of the two trace powers is retained. The comparison uses
  eigenvalues, not operator monotonicity of the $q/2$ power. Since $\|Y_{N,M}\|
  \le \|Y_{N,M}\|_{S_q}$, this proves $L^q(\Omega;\mathcal{B}(E))$ Cauchy
  convergence; a uniform moment bound alone would not suffice.

  \medskip
  \noindent\textit{Step 2.~}We construct the operator limit on the prescribed
  probability space.
  All finite sums lie in the separable operator-norm closed real span
  $\mathcal{F}_I$ of the countably many deterministic $A_j(I)$. Choose $N_k$
  increasing so that $\|X_{N_{k+1}}(I)-X_{N_k}(I)\|_{L^q\mathcal{B}(E)} \le
  2^{-k}$. Their $L^1$ norms are summable. Tonelli therefore shows that the
  pointwise operator-norm differences sum absolutely almost surely. Completeness
  gives a limit in $\mathcal{F}_I$ on the original probability space; assign it
  zero on the exceptional null set. It is strongly measurable as a limit of
  separably valued finite Gaussian sums. Fatou and the Cauchy property give
  convergence of the whole sequence in $L^q$ to this limit $X_R(I)$. No new
  noise or probability space has been introduced. Limits constructed at
  different admissible even indices agree by uniqueness of convergence in
  probability. For each requested finite $p \ge 1$, choose an admissible even $q
  \ge p$ and use Lyapunov to obtain $L^p$ operator-norm convergence to that same
  limit. One fixed Schatten index is not used to infer all moment orders.

  \medskip
  \noindent\textit{Step 3.~}We identify compactness, moments and Gaussianity
  of the limit.
  The compact operators are norm closed, so an almost surely norm-convergent
  subsequence makes $X_R(I)$ compact valued. Singular-value Fatou applied to the
  finite uniform bound gives \eqref{E:driver-fixed-schatten} at each admissible
  even $q$. Every continuous real linear functional of $X_R(I)$ is an $L^2$
  limit of centered finite Gaussian sums; convergence of means, variances and
  characteristic functions makes it centered Gaussian. The same applies to
  finite linear combinations of these functionals, proving joint Gaussianity on
  the given space.

  \medskip
  \noindent\textit{Step 4.~}We identify the rows as actual noise variables.
  For each fixed real Schwartz $z,w$, the scalar rows $\langle
  X_N(I)z,w\rangle_E$ are the orthogonal projections of the completed covariance
  form $\ell_R(I;z,w)$ onto the complete noise basis. That form lies in the
  original closed first Wiener chaos, so these projections converge to it in
  $L^2(\Omega)$. Operator-norm $L^2$ convergence also gives convergence to
  $\langle X_R(I)z,w\rangle_E$. Uniqueness of the $L^2$ limit identifies these
  two variables on the same space, not merely their laws.

  \medskip
  \noindent\textit{Step 5.~}We extend the construction by finite additivity.
  For a fixed long interval, choose a finite partition into intervals of length
  at most one. Their indicator tests add to the long-interval test except at
  finitely many endpoints. The masked wave products have common compact spatial
  support, and the temporal estimate~\eqref{E:driver-test-continuity} makes
  those endpoint changes have zero original covariance norm. Covariance
  completion and linearity give exact additivity of the forms. The deterministic
  statewise coefficient integrals, hence every $X_N$, add on the same partition.
  Passing their finite-sum $L^p$ limits gives the corresponding operator
  additivity and covariance identity for matrix elements, independently of the
  partition by common refinement. These are statements at each fixed interval
  and partition; no uncountable intersection of independently selected null sets
  is used. No long-interval version of the short-$h$ bound is asserted. The
  Schatten conclusion is a fixed-interval moment at admissible even indices, not
  a Schatten-valued path or arbitrary-index claim. This proves
  Lemma~\ref{L:driver-interval-limit}.
\end{proof}

\subsection{Time-path convergence and Gaussian norm bounds}

Lemma~\ref{L:driver-interval-limit} works one interval at a time.
Proposition~\ref{P:actual-operator-driver} asks for more: convergence in the
$C^\theta$ path norm and compact values at every time on one path event.

\LeanMarker{L:driver-holder-limit}
\begin{lemma}[H\"older path limit]\label{L:driver-holder-limit}
  The finite paths $X_N(t) = X_N((0,t))$ converge in
  $L^p(\Omega;C^\theta([0,T];\mathcal{B}(E)))$ for every finite $p \ge 1$, in a
  separable closed path span on the same probability space. The limit starts at
  zero and has compact values for all times on one path event. For sufficiently
  large even $q_0$, its $C^\theta$ norm has $L^q$ bound $C\sqrt q$, uniformly
  over even $q \ge q_0$ and the basis.
\end{lemma}

\begin{proof}
  We first obtain a uniform moment bound on the dyadic increments in a
  stronger Hölder norm. Finite-grid convergence then gives uniform
  convergence, and interpolation gives the desired Hölder convergence. We
  construct the limit in a separable path span.

  \medskip
  \noindent\textit{Step 1.~}We bound increments and the dyadic sum.
  Corollary~\ref{C:frequency-time} bounds the budget $b_I$ by
  $C|t-u|^{2\beta}$ for $I=(u,t)$; substitute this
  in~\eqref{E:driver-two-variances} and apply
  Corollary~\ref{C:finite-gaussian-operators} to the finite increment sum.
  Choose $\theta < \eta < \beta$ and an even $q_0 \ge 2$ with $q_0\gamma > d$
  and $\beta-\eta-1/q_0 > 0$. For $q \ge q_0$, the singular-value sequence gives
  $\|K_\chi\|_{S_q} \le \|K_\chi\|_{S_{q_0}}$: if the latter is nonzero,
  normalize by it, so every singular value is at most one, and sum its $q$th
  powers. The zero case is immediate. Apply the finite Gaussian operator
  estimate with the time-refined variances from
  Proposition~\ref{P:wave-two-variances}. This gives, uniformly in $N$ and even
  $q \ge q_0$,
  \begin{equation}\label{E:driver-increment-moments}
    \|X_N(t)-X_N(u)\|_{L^q(\mathcal{B}(E))}
    \le C\sqrt q\,|t-u|^\beta,\qquad |t-u|\le1.
  \end{equation}
  The constant includes the fixed $S_{q_0}$ cutoff norm and does not depend on
  $q$. Longer increments are finite sums of unit increments, with constants
  depending on $T$.

  Let $D_n(F)$ be the maximum of adjacent increments on the grid $kT2^{-n}$, $0
  \le k \le 2^n$. At fine levels $T2^{-n} \le 1$, the elementary bound for a
  finite maximum yields
  \[
    \|D_n(X_N)\|_{L^q}
    \le 2^{n/q}C\sqrt q\,(T2^{-n})^\beta.
  \]
  The $q$th power of a maximum is at most the sum of the $q$th powers. There
  are $2^n$ adjacent increments, so taking the $q$th root contributes
  $2^{n/q}$. Each increment has length $T2^{-n}$, which yields the displayed
  estimate.
  No independence of these increments is used. The finitely many coarse levels
  are controlled by unit partitions. Minkowski and a geometric sum therefore
  give
  \begin{equation}\label{E:driver-dyadic-moment-sum}
    \sup_N\left\|\sum_{n\ge0}2^{n\eta}D_n(X_N)\right\|_{L^q}
    \le C_\eta\sqrt q,
  \end{equation}
  since the fine-level ratio is at most $2^{-(\beta-\eta-1/q_0)} < 1$.

  \medskip
  \noindent\textit{Step 2.~}We control the stronger Hölder norm by dyadic
  increments.
  For a continuous Banach-valued $F$ starting at zero, this weighted dyadic sum
  controls its $C^\eta$ norm. Approximate two times by dyadic points at a scale
  comparable to their separation. The two grid points differ by at most two
  adjacent increments, and each endpoint's finer binary corrections are bounded
  by the tail sum of $D_n$. Multiplying by the separation to the power $-\eta$
  bounds that total by the weighted sum, up to a constant depending on $T,\eta$.
  Continuity identifies the endpoint limits and also controls the supremum norm
  from $F(0) = 0$. The finite paths are Lipschitz by~\eqref{E:driver-gain-one}.
  Thus, \eqref{E:driver-dyadic-moment-sum} gives uniform $L^qC^\eta$ bounds, and
  the triangle inequality gives the same bound for their differences.

  \medskip
  \noindent\textit{Step 3.~}We prove convergence in the required Hölder norm.
  On a fixed finite time grid, Lemma~\ref{L:driver-interval-limit} gives $L^q$
  convergence of every point value. For continuous $F$,
  \[
    \|F\|_\infty\le\max_k\|F(t_k)\|
    +C(\operatorname{mesh})^\eta\|F\|_{C^\eta}.
  \]
  Apply this to $X_N-X_M$, send $N,M \to \infty$ at a fixed grid, and then send
  the mesh to zero. The uniform $C^\eta$ bound proves sup-norm $L^q$ Cauchy
  convergence. The deterministic interpolation inequality
  \begin{equation}\label{E:driver-holder-interpolation}
    \|F\|_{C^\theta}
    \le C\|F\|_\infty^{1-\theta/\eta}
    \|F\|_{C^\eta}^{\theta/\eta}
  \end{equation}
  follows by bounding each increment by both $2\|F\|_\infty$ and its $C^\eta$
  bound, and splitting at their balancing scale. Its supremum term is controlled
  too. Scalar H\"older in probability proves the $L^qC^\theta$ Cauchy property.
  For any requested finite $p$, choose an even $q \ge \max(p,q_0)$ and use
  Lyapunov.

  \medskip
  \noindent\textit{Step 4.~}We choose the measurable limit and one path event.
  The full H\"older/operator space need not be separable. Nevertheless, every
  finite path lies in the span of the countably many deterministic Lipschitz
  paths $A_j(0,\cdot)$. Their $C^\theta$-closed real span is a separable Banach
  space. In it the same elementary Cauchy/subsequence construction as in
  Lemma~\ref{L:driver-interval-limit} gives a strongly measurable limit $M_R$ on
  the original space. Different $p$ limits agree by uniqueness of convergence in
  probability in the path norm. An almost surely $C^\theta$-convergent
  subsequence converges uniformly in operator norm. Finite paths have compact
  values on short intervals, and on long ones by finite additivity. Therefore,
  the limit has compact values at every time on that single path set. It starts
  at zero. Passing the uniform bounds from~\eqref{E:driver-dyadic-moment-sum}
  through the path convergence proves the asserted $C\sqrt q$ bound. This proves
  Lemma~\ref{L:driver-holder-limit}.
\end{proof}

\begin{proof}[Proof of Proposition~\ref{P:actual-operator-driver}]
  We assemble the interval and path constructions, then verify basis
  independence. The uniform square-root moment bound gives the exponential
  norm estimate by summing its power series.

  \medskip
  \noindent\textit{Step 1.~}We identify the path and its covariance rows.
  Lemmas~\ref{L:driver-noise-basis}--\ref{L:driver-holder-limit} construct the
  stated finite sums and their path limit, with no new probability space.
  Evaluation at a fixed time agrees with the interval limit by uniqueness in
  probability. Its covariance identity for matrix elements is
  \eqref{E:driver-raw-row}, by Lemma~\ref{L:driver-interval-limit}, including
  the long-interval interpretation. The same additivity identifies each fixed
  increment with $X_R(I)$, proving its stated Schatten moment bound.

  The finite paths are Gaussian linear combinations of deterministic paths.
  Every continuous linear functional on their closed path span converges in
  $L^2$ to a centered Gaussian; finite linear combinations do as well. Thus,
  $M_R$ is a centered Gaussian path. The four real pair identities give its
  complexified sesquilinear row identity.

  \medskip
  \noindent\textit{Step 2.~}We compare paths built from different noise bases.
  The construction used only compact smooth representatives, orthonormality and
  completeness, so it applies to any such noise basis. Compare two constructions
  using a fixed countable dense Schwartz phase collection and rational times.
  Both rows equal the same covariance forms. A countable null-set union,
  boundedness of the operators, phase density and path continuity give equality
  at all times on one probability-one set. This proves indistinguishability
  without an assumed common event for uncountably many original evaluations.

  \medskip
  \noindent\textit{Step 3.~}We deduce all finite moments and exponential
  integrability.
  For $p \ge 1$, choose even $q \ge \max(p,q_0)$ with $q \le p+q_0+2 \le
  (q_0+3)p$. Lemma~\ref{L:driver-holder-limit} and Lyapunov give the first bound
  in~\eqref{E:driver-gaussian-norm}. Write its resulting constant as $C_1$;
  Lyapunov and $q \le (q_0+3)p$
  absorb the passage from even $q$ into the constant $C_1$. Now
  let $Y = \|M_R\|_{C^\theta}$. If $C_1 = 0$, then $Y = 0$ almost surely and the
  exponential bound is immediate. Otherwise, $\E Y^{2m} \le
  (C_1\sqrt{2m})^{2m}$. Since $m! \ge (m/e)^m$, as follows by integrating $\log
  x$ from one to $m$, the $m$th term in the exponential series is at most
  $(2ecC_1^2)^m$. For $0 < c < 1/(2eC_1^2)$, positive Tonelli and the geometric
  series prove the second bound. This proves
  Proposition~\ref{P:actual-operator-driver}.
\end{proof}

All constants here may depend on $R,T$ and the fixed admissible parameters. No
uniform-in-$R$ path-norm estimate or $C^\theta S_q$-valued path is asserted.
This constructs the local masked integrated operator on the original noise.
Source-cutoff convergence, Young evolution, canonical initial data and
solution-moment asymptotics remain separate parts of the argument.

\section{Convergence of the regularized Gaussian operator paths}
\label{S:source-driver-convergence}

We match the heat normalization, Laplace prefix and joint source law in the
author version~\cite{chen:author2026:time-dependency}*{(1.20), (1.21), (2.9),
(2.11), pp.~5--6, 8--9}. Preservation of the given noise retention, the Radon
implementation and the operator-path conclusions are proved here. No
first-moment variational theorem is an input.

Retain Hypothesis~\ref{H:homog}, the fixed $R,T,\zeta_R,E$ and parameters
of~\eqref{E:driver-parameters}, and the original noise space $H_W$. Write
$\mathscr{B}_\theta = C^\theta([0,T];\mathcal{B}(E))$. Fix one original test
basis $f_j$ and its coefficient paths $A_j(t)$ from
Proposition~\ref{P:actual-operator-driver}. Thus, $M_N = \sum_{j\le N}W(f_j)A_j
\to M_R$ in every finite $L^p(\mathscr{B}_\theta)$. Zero covariance gives zero
noises and operator paths throughout; assume $\Gamma \ne 0$ in the construction.

\subsection{Gaussian transfer in the full path norm}

\LeanMarker{L:source-finite-comparison}
\begin{lemma}[Finite Gaussian linear images]\label{L:source-finite-comparison}
  Let $v_1,\ldots,v_N$ lie in a real Banach space and let $\xi$ be a centered
  Gaussian vector with covariance bounded by $CI$, $C \ge 0$. For independent
  standard real $g_j$ and $p \ge 1$,
  \begin{equation}\label{E:source-finite-comparison}
    \left(\E\left\|\sum_j\xi_jv_j\right\|^p\right)^{1/p}
    \le\sqrt C\left(\E\left\|\sum_jg_jv_j\right\|^p\right)^{1/p}.
  \end{equation}
\end{lemma}

\begin{proof}
  The function $\varphi(x) = \|\sum_jx_jv_j\|^p$ is continuous and convex, even
  when the linear image gives only a seminorm on $\R^N$. Its Gaussian
  expectations are finite by the finite-sum bound. On an auxiliary space take an
  independent centered Gaussian $Z$ of covariance $CI-\operatorname{Cov}\xi$.
  The covariance difference is positive semidefinite by hypothesis, so its
  diagonalization constructs such a Gaussian vector even if it is degenerate.
  Independence and centering give conditional mean $\xi$ for $\xi+Z$ given
  $\xi$. Convexity therefore applies to this conditional expectation.
  Conditional Jensen gives $\E\varphi(\xi) \le \E\varphi(\xi+Z)$. The latter
  vector has the law of $\sqrt C\,g$, proving the claim. The case $C = 0$ is
  immediate. This auxiliary comparison changes neither the coupling of $\xi$ nor
  any nonlinear wave solution. This proves
  Lemma~\ref{L:source-finite-comparison}.
\end{proof}

Lemma~\ref{L:source-finite-comparison} compares two finite Gaussian sums with
the same vectors. We apply it to the coefficient paths $A_j$, keeping the
original basis $f_j$ and replacing $W(f_j)$ by the evaluations of another
Gaussian noise.

\begin{proposition}[Gaussian transfer and stability]
  \label{P:source-gaussian-transfer}
  Let $Y$ be a centered Gaussian linear noise on smooth compact positive-time
  tests, on a common coupled space, with
  \begin{equation}\label{E:source-raw-variance-bound}
    \operatorname{Var}Y(F)\le C\|F\|_{H_W}^2.
  \end{equation}
  It extends to $H_W$. The finite paths $M_{Y,N} = \sum_{j\le N}Y(f_j)A_j$
  converge in every finite $L^p(\mathscr{B}_\theta)$ to a centered Gaussian path
  $M_Y$ with compact values at all times on one path event, and
  \begin{align}
    \|M_Y-M_{Y,N}\|_{L^p\mathscr{B}_\theta}
    &\le\sqrt C\|M_R-M_N\|_{L^p\mathscr{B}_\theta},\label{E:source-transfer-tail}\\
    \|M_Y\|_{L^p\mathscr{B}_\theta}
    &\le\sqrt C\|M_R\|_{L^p\mathscr{B}_\theta}.\label{E:source-transfer-norm}
  \end{align}
  For a fixed real Schwartz pair and time,
  \begin{equation}\label{E:source-transferred-row}
    \langle M_Y(t)z,w\rangle_E=Y(h_{R,t;z,w})\quad\text{in }L^2,
  \end{equation}
  where $h_{R,t;z,w}$ is the explicit original test class
  $\zeta_R1_{(0,t)}v_zG_w$ with the adjoint from~\eqref{E:variance-dual-wave}.
  The path is basis-independent up to indistinguishability, is linear in jointly
  Gaussian noise families, and complexifies sesquilinearly.

  If $(W,Y_n)$ are jointly centered Gaussian, $Y_n$ and $D_n = W-Y_n$ have
  uniform bounds of the form~\eqref{E:source-raw-variance-bound}, and
  $\operatorname{Var}D_n(F) \to 0$ on every smooth compact test, then $M_{Y_n}
  \to M_R$ in every finite $L^p(\mathscr{B}_\theta)$.
\end{proposition}

\begin{proof}
  We first transfer finite path sums and their tails using the Gaussian
  comparison. We then identify the transferred rows on the given coupling. For
  stability, we make the common tail small before sending the finitely many
  remaining noise errors to zero.

  \medskip
  \noindent\textit{Step 1.~}We extend the noise and compare its finite
  coordinate vectors.
  The variance bound annihilates covariance-null tests and extends $Y$ as a
  bounded linear map into its $L^2$ space. Completion preserves Gaussian
  finite-dimensional laws. For each finite original orthonormal family, the
  covariance of its $Y$ evaluations is bounded by $CI$, by applying the variance
  bound to every real linear combination. Independence of $Y(f_j)$ is not
  assumed.

  \medskip
  \noindent\textit{Step 2.~}We construct the transferred path and its bounds.
  Apply Lemma~\ref{L:source-finite-comparison} with $v_j = A_j(\cdot)$ in
  $\mathscr{B}_\theta$, also to finite tails. The standard comparison vector has
  the law of the original $W(f_j)$. The already-proved original path convergence
  makes $M_{Y,N}$ Cauchy in every $L^p$. It lies in the same separable closed
  span of the countably many deterministic paths. The elementary
  Cauchy/subsequence construction of Lemma~\ref{L:driver-holder-limit} gives a
  strongly measurable limit on the space of $Y$. Different $p$ limits agree by
  convergence in probability. Passing to limits in the finite inequalities gives
  both displayed bounds. An almost surely uniformly convergent subsequence gives
  compact values for all times. Functional $L^2$ limits of Gaussian linear
  combinations prove that the path is Gaussian.

  \medskip
  \noindent\textit{Step 3.~}We identify the rows and remove dependence on the
  basis.
  For real Schwartz $z,w$, the coefficient identity is
  \[
    \langle A_j(t)z,w\rangle_E
    =\E[W(h_{R,t;z,w})W(f_j)]
    =\langle h_{R,t;z,w},f_j\rangle_{H_W}.
  \]
  Therefore, the finite transferred row is $Y(P_Nh)$, with $P_N$ the original
  $H_W$ orthogonal projection. Since $P_Nh \to h$, the noise bound gives
  $Y(P_Nh) \to Y(h)$ in $L^2$. Operator-path $L^2$ convergence
  proves~\eqref{E:source-transferred-row}. Long intervals use the same
  covariance completion and endpoint nullity as in
  Lemma~\ref{L:driver-interval-limit}. No filter is commuted through a physical
  mask or wave group. Four real pair identities give the complex assertion.
  Countably many dense Schwartz phases and rational times, followed by phase
  density and path continuity, prove basis indistinguishability.

  \medskip
  \noindent\textit{Step 4.~}We prove linearity and stability of the transfer.
  For jointly Gaussian $Y,Z$, their real linear combinations have an original
  variance bound by the $L^2$ triangle inequality. The finite sums are linear in
  the noise, hence so are the path limits. In the stability assertion $D_n$ is
  centered Gaussian by the joint hypothesis. Its uniform bound gives
  \[
    \|M_{D_n}-M_{D_n,N}\|_{L^p\mathscr{B}_\theta}
    \le\sqrt{C_D}\|M_R-M_N\|_{L^p\mathscr{B}_\theta}.
  \]
  The right side of this bound is independent of $n$ and tends to zero as $N$
  increases. We can therefore select the tail cutoff before considering
  convergence of the finite coordinate vector.
  For fixed $N$, the covariance matrix of $(D_n(f_j))_{j\le N}$ is positive and
  bounded by $s_N(n)I$, where $s_N(n) = \sum_{j\le N}\operatorname{Var}D_n(f_j)
  \to 0$. Every eigenvalue of a positive finite matrix is at most its trace,
  so the scalar $s_N(n)$ bounds that covariance matrix in positive order.
  Thus,
  the finite comparison gives $\|M_{D_n,N}\|_{L^p\mathscr{B}_\theta} \le
  \sqrt{s_N(n)}\|M_N\|_{L^p\mathscr{B}_\theta}$. Choose $N$ large for the
  uniform tail, then $n$ large for the finite head. Linearity identifies
  $M_{D_n} = M_R-M_{Y_n}$. The same argument proves a multiparameter limit
  whenever each fixed test error tends to zero, with no relation between
  parameter rates. This proves Proposition~\ref{P:source-gaussian-transfer}.
\end{proof}

\subsection{The spatial covariance space and an extension preserving the given noise}

Proposition~\ref{P:source-gaussian-transfer} reduces source-path convergence
to uniform variance domination and vanishing residual variance on each test,
for noises jointly Gaussian with $W$. We now construct such regularized
sources, beginning with the spatial covariance space and heat regularization.

\LeanMarker{L:source-spatial-heat}
\begin{lemma}[Original spatial space and heat]\label{L:source-spatial-heat}
  Let $H_\Gamma$ be the real test-null completion with
  \[
    B_\Gamma(f,g)=\int\Gamma(\ud y)\int f(x)g(x-y)\ud x
    =(2\pi)^d\int\widehat{f}(\xi)\overline{\widehat{g}(\xi)}\mu(\ud\xi).
  \]
  It is separable and includes Schwartz tests. With
  \begin{equation}\label{E:source-heat-definition}
    p_\varepsilon(x)=(2\pi\varepsilon)^{-d/2}
    e^{-|x|^2/(2\varepsilon)},\qquad C_\varepsilon f=p_\varepsilon*f,
  \end{equation}
  heat is a positive self-adjoint contraction on $H_\Gamma$, $I-C_\varepsilon$
  is a contraction, and $C_\varepsilon \to I$ strongly.
\end{lemma}

\begin{proof}
  The spectral normalization is the spatial computation in
  Proposition~\ref{P:masked-noise-completion}; its integral is real on real
  tests by evenness of $\mu$. The local bound $|B_\Gamma(f,g)| \le
  |K|\Gamma(K-K)\|f\|_\infty\|g\|_\infty$ on a common compact support gives
  separability by countable smooth approximation. Expanding smooth cutoffs of a
  Schwartz test converge in Schwartz seminorms, and its covariance correlation
  pairing converges against tempered $\Gamma$. Thus, Schwartz tests enter the
  same completion. The Fourier multiplier of $C_\varepsilon$ is
  $m_\varepsilon(\xi) = e^{-\varepsilon|\xi|^2/2}$. On dense Schwartz tests, the
  displayed norm represents heat by multiplication by the real function
  $m_\varepsilon$, which satisfies $0 \le m_\varepsilon \le 1$. The multiplier
  of $I-C_\varepsilon$ is $1-m_\varepsilon$, also between zero and one. Thus,
  the
  displayed norm gives positivity, self-adjointness and both contraction bounds.
  Dominated convergence gives strong convergence there; uniform boundedness
  extends it to the completion. This uses no physical mask operator on a noise
  quotient. This proves Lemma~\ref{L:source-spatial-heat}.
\end{proof}

Set $\nu_a(\ud\lambda) = \lambda^{a-1}\ud\lambda/\Gamma_E(a)$ and
\[
  \mathcal{G}=L^2((0,\infty)_\lambda\times\R_r,
  \nu_a(\ud\lambda)\ud r;H_\Gamma).
\]
This is a real separable Hilbert space: countable simple functions on rational
rectangles with a countable dense spatial collection suffice.

\begin{lemma}[Feature isometry preserving the given noise]\label{L:source-retaining-isometry}
  On finite smooth temporal-spatial tensors define
  \begin{equation}\label{E:source-feature-isometry}
    \Phi F(\lambda,r)=\sqrt{2\lambda}
    \int_{\substack{t>0\\t\ge r}}e^{-\lambda(t-r)}F(t,\cdot)\ud t.
  \end{equation}
  It extends to an isometry $\Phi:H_W \to \mathcal{G}$ with closed range, and
  the formula represents smooth compact tests almost everywhere. On a product
  extension retaining the original $W$, there is an isonormal $\mathcal{I}$ on
  $\mathcal{G}$ with
  \begin{equation}\label{E:source-given-W-retention}
    \mathcal{I}(\Phi h)=W(h)\quad\text{in }L^2
    \quad\text{for each }h\in H_W.
  \end{equation}
\end{lemma}

\begin{proof}
  We first prove the feature map is an isometry and identify its formula on
  completed tests. We then extend the given noise only on the orthogonal
  complement of the feature range, preserving every original evaluation.

  \medskip
  \noindent\textit{Step 1.~}We identify the feature norm with the original
  covariance norm.
  The identity
  \[
    \int_{-\infty}^{\min(t,u)}
    2\lambda e^{-\lambda(t+u-2r)}\ud r
    =e^{-\lambda|t-u|}
  \]
  and the original Laplace formula give $\|\Phi F\|_{\mathcal{G}}^2 =
  \|F\|_{H_W}^2$.
  Indeed, integrating the two exponential time features over $r$ gives the
  displayed kernel. Integrating that kernel against $\nu_a$ gives
  $|t-u|^{-a}$. The spatial inner product is the original $B_\Gamma$, so
  these two integrations reproduce the original covariance.
  All signed exchanges for finite tensors are absolutely
  integrable: compact temporal factors give integrals of
  $|\tau(t)\sigma(u)||t-u|^{-a}$, finite for $a < 1$, times fixed finite spatial
  inner products. No oscillatory infinite-measure Fubini is used. Smooth time
  partitions on one compact support approximate a general smooth spacetime test
  uniformly by finite tensors. The local covariance
  bound~\eqref{E:driver-test-continuity} proves their density in $H_W$. Hence,
  the isometry extends with closed range.

  \medskip
  \noindent\textit{Step 2.~}We identify the integral formula for general
  compact tests.
  For a general smooth compact $F$, those approximants also converge uniformly
  in spatial $H_\Gamma$ norm in time, by local Radon mass. For fixed
  $\lambda,r$, their integrals in~\eqref{E:source-feature-isometry} converge in
  $H_\Gamma$. Meanwhile the isometry gives convergence in $\mathcal{G}$, with an
  almost everywhere convergent subsequence. Its two limits agree, proving the
  stated formula. The coordinate $r$ is auxiliary; no value of original
  negative-time noise has been assigned.

  \medskip
  \noindent\textit{Step 3.~}We extend the isonormal process while retaining
  the given noise.
  Let $P$ project orthogonally onto $\Phi H_W$. On a product extension take an
  independent isonormal $J$ on $(\Phi H_W)^\perp$, built from independent real
  standard Gaussians on a finite or countable basis. Define
  \[
    \mathcal{I}(g)=W(\Phi^{-1}Pg)+J((I-P)g).
  \]
  The two summands use orthogonal feature components and independent Gaussian
  processes. Their variances are $\|Pg\|_{\mathcal G}^2$ and
  $\|(I-P)g\|_{\mathcal G}^2$, whose sum is $\|g\|_{\mathcal G}^2$. If $g =
  \Phi h$, the second component vanishes and the first is the original $W(h)$.
  Orthogonality and independence make it isonormal on $\mathcal{G}$ and
  prove~\eqref{E:source-given-W-retention} with the original random variable,
  not only its law. For the countable fixed original basis, these equalities
  hold on one common probability-one set. Thus, the original finite paths and
  their $L^p$ limits are precisely pullbacks of the already constructed original
  operator path. Cutoff fields may use $J$; they are not asserted to be
  functions of $W$ alone. This proves Lemma~\ref{L:source-retaining-isometry}.
\end{proof}

The auxiliary component can be nonzero: No-go~\ref{P:nogo-32}
constructs an orthogonal witness for every positive finite cutoff point
feature at positive time. Proposition~\ref{P:constant-raw-average} later
proves that the tested fixed-mask limit is $W$-measurable, using the
operator-path and Young-flow convergence.

\subsection{Exact joint source law and convergence}

Let $P_\eta$ multiply $\mathcal{G}$ by $1_{\{\lambda\le\eta^{-1}\}}$. Let
$C_\varepsilon$ act fiberwise on $H_\Gamma$. These commute. For $\rho =
(\varepsilon,\eta) > 0$, define
\begin{equation}\label{E:source-cutoff-noises}
  Q_\rho=C_\varepsilon P_\eta,\qquad
  Y_\rho(F)=\mathcal{I}(Q_\rho\Phi F),\qquad
  D_\rho(F)=\mathcal{I}((I-Q_\rho)\Phi F)=W(F)-Y_\rho(F).
\end{equation}
The last equality uses the coupling preserving the given noise. Define
\begin{equation}\label{E:source-covariance-functions}
  \Gamma_e(x)=(p_e*\Gamma)(x),\qquad
  K_\eta(v)=\frac1{\Gamma_E(a)}
  \int_0^{\eta^{-1}}e^{-\lambda|v|}\lambda^{a-1}\ud\lambda.
\end{equation}
The temporal operation is a Laplace prefix, not temporal convolution.

\begin{proposition}[Canonical sources on the extended space containing the given noise]
  \label{P:source-canonical-convergence}
  The noises in~\eqref{E:source-cutoff-noises} have random field realizations
  with joint covariance
  \begin{equation}\label{E:source-joint-covariance}
    \operatorname{Cov}(Y_{\varepsilon,\eta}(t,x),
    Y_{\varepsilon',\eta'}(s,y))
    =K_{\max(\eta,\eta')}(t-s)\Gamma_{\varepsilon+\varepsilon'}(x-y).
  \end{equation}
  Mean-square random field integration agrees with $Y_\rho(F)$ on tests. Both
  source and full residual obey
  \begin{equation}\label{E:source-two-dominations}
    \operatorname{Var}Y_\rho(F),\quad\operatorname{Var}D_\rho(F)
    \le\|F\|_{H_W}^2,
  \end{equation}
  and $\operatorname{Var}D_\rho(F) \to 0$ for every $F \in H_W$ as
  $\varepsilon,\eta\downarrow0$ at arbitrary relative rates. Consequently, for
  every finite $p \ge 1$,
  \begin{equation}\label{E:source-driver-limit}
    M_R^\rho:=M_{Y_\rho}\longrightarrow M_R
    \quad\text{in }L^p(\mathscr{B}_\theta),\qquad
    M_{D_\rho}=M_R-M_R^\rho.
  \end{equation}
  The source and residual have uniform local Gaussian operator-path moments and
  square-exponential norms, with constants for fixed $R,T$ and exponents.
\end{proposition}

\begin{proof}
  We compute the joint law from the point features, including the two heat
  widths and the common Laplace range. We then identify field integrals, bound
  both filters, and apply Gaussian transfer to the vanishing residual.

  \medskip
  \noindent\textit{Step 1.~}We compute the spatial covariance and
  point-feature covariance.
  Apply the spatial distribution identity to a single shifted heat test:
  \[
    \Gamma_e(x)=\int e^{i\xi\cdot x}e^{-e|\xi|^2/2}\mu(\ud\xi).
  \]
  Gaussian spectral moments of every order are finite for tempered $\mu$. Thus,
  this is smooth with bounded derivatives. Since
  $\widehat{p_\varepsilon(x-\cdot)} = c_de^{-i\xi\cdot
  x}e^{-\varepsilon|\xi|^2/2}$, the two heat factors in the spatial inner
  product cancel $(2\pi)^d$ exactly, giving
  \begin{equation}\label{E:source-heat-inner-product}
    B_\Gamma(p_\varepsilon(x-\cdot),p_{\varepsilon'}(y-\cdot))
    =\Gamma_{\varepsilon+\varepsilon'}(x-y).
  \end{equation}
  The product of the two heat multipliers is
  $e^{-(\varepsilon+\varepsilon')|\xi|^2/2}$. The product of their Fourier
  normalization factors is $c_d^2$, and $(2\pi)^dc_d^2 = 1$. This gives the
  sum of widths without any remaining Fourier factor.
  Evenness of $\mu$ removes the reflected argument.

  For each positive $\rho$ and point $(t,x)$ put
  \begin{equation}\label{E:source-point-feature}
    v_{\rho;t,x}(\lambda,r)
    =1_{\{\lambda\le\eta^{-1}\}}\sqrt{2\lambda}
    e^{-\lambda(t-r)}1_{\{r\le t\}}p_\varepsilon(x-\cdot).
  \end{equation}
  Its squared $\mathcal{G}$ norm is $K_\eta(0)\Gamma_{2\varepsilon}(0) <
  \infty$. Define $Y_\rho(t,x) = \mathcal{I}(v_{\rho;t,x})$. The common
  $\lambda$ range for two such vectors is $\lambda \le
  \min(\eta^{-1},\eta'^{-1}) = \max(\eta,\eta')^{-1}$. Integrating $r$ gives
  $e^{-\lambda|t-s|}$, and~\eqref{E:source-heat-inner-product} gives the spatial
  factor. This proves~\eqref{E:source-joint-covariance} with the prescribed
  widths and original amplitudes, also for singular Radon and white covariances.

  \medskip
  \noindent\textit{Step 2.~}We identify the mean-square field integral with
  the filtered noise.
  The representing vector is continuous in $\mathcal{G}$ norm, directly from
  this continuous covariance. Its norm is bounded on compact parameter sets.
  Hence, for a smooth compact spacetime $F$, its Bochner integral against
  $v_{\rho;t,x}$ exists. Integrate the heat translates in $x$, then the time
  feature, to obtain exactly $Q_\rho\Phi F$. Applying the bounded map
  $\mathcal{I}:\mathcal{G} \to L^2$ proves the random field-test identity on the
  same space. Common-support indicator tests follow by covariance completion and
  $L^1$ approximation; bounded point-feature norms and the domination below give
  the same $L^2$ limit.

  \medskip
  \noindent\textit{Step 3.~}We bound the source and the full residual.
  Orthogonality of the two $\lambda$ ranges gives, for $g \in \mathcal{G}$,
  \[
    \|(I-Q_\rho)g\|_{\mathcal{G}}^2
    =\|(I-P_\eta)g\|_{\mathcal{G}}^2
    +\|P_\eta(I-C_\varepsilon)g\|_{\mathcal{G}}^2
    \le\|g\|_{\mathcal{G}}^2.
  \]
  On the retained range, $I-C_\varepsilon$ is a contraction by the heat lemma.
  Hence, its squared norm is bounded by $\|P_\eta g\|_{\mathcal G}^2$. Adding
  the omitted-range norm gives at most the full squared norm, because these
  two ranges are orthogonal.
  Also $\|Q_\rho g\|_{\mathcal{G}} \le \|g\|_{\mathcal{G}}$. Together with the
  isometry these prove both variance bounds. Disjoint $\lambda$ ranges are
  independent; the full heat residual need not be independent of the source.

  \medskip
  \noindent\textit{Step 4.~}We send both regularization parameters to zero.
  The $L^2$ tail integral gives $P_\eta \to I$ strongly. Fiberwise heat
  converges strongly on simple $\mathcal{G}$ functions by
  Lemma~\ref{L:source-spatial-heat}, then on all of $\mathcal{G}$ by the
  contraction bound. The estimate
  \[
    \|(I-Q_\rho)g\|_{\mathcal{G}}
    \le\|(I-C_\varepsilon)g\|_{\mathcal{G}}
    +\|(I-P_\eta)g\|_{\mathcal{G}}\longrightarrow0
  \]
  holds at arbitrary relative rates.
  Indeed, the two terms on the right depend on separate parameters, and each
  tends to zero by strong convergence. Therefore, there is no need to compare
  the rates at which the heat width and the Laplace cutoff are removed.
  Thus, the residual variance tends to zero
  for every original $F$. Proposition~\ref{P:source-gaussian-transfer}
  gives~\eqref{E:source-driver-limit} without commuting a filter through
  $\zeta_R$ or a wave group.

  \medskip
  \noindent\textit{Step 5.~}We retain the uniform Gaussian path bounds.
  The same comparison gives $\|M_R^\rho\|_{L^p\mathscr{B}_\theta},
  \|M_{D_\rho}\|_{L^p\mathscr{B}_\theta} \le \|M_R\|_{L^p\mathscr{B}_\theta}$.
  The original $\sqrt p$ estimate and its even-moment exponential-series
  argument give uniform Gaussian operator path bounds. They are statements about
  $L^p$ equivalence classes on the extension preserving the given noise, not a
  simultaneous path event for all uncountably many cutoffs. This proves
  Proposition~\ref{P:source-canonical-convergence}.
\end{proof}

\subsection{Regularized random fields and the ordinary first response}

\begin{lemma}[Regularity for each fixed cutoff]\label{L:source-pointfield-regularity}
  For each fixed $\rho > 0$, the random field has a version continuous in time
  and smooth in space, with every spatial derivative jointly continuous on
  compact time-space rectangles.
\end{lemma}

\begin{proof}
  We first differentiate the point features and bound their increments. A
  dyadic construction gives jointly continuous versions for all spatial
  derivative fields. The feature-space fundamental theorem then identifies
  those fields as actual derivatives of one version.

  \medskip
  \noindent\textit{Step 1.~}We construct derivative fields and bound their
  increments.
  For each spatial multiindex $b$, differentiate \eqref{E:source-point-feature}
  in $x$ as a $\mathcal{G}$ vector. All derivatives exist because their spatial
  squared norms are the finite heat moments
  \[
    \int \xi^{2b}e^{-\varepsilon|\xi|^2}\mu(\ud\xi),
  \]
  and the next moments dominate difference quotients. Set $Z_b(t,x) =
  \mathcal{I}(\partial_x^bv_{\rho;t,x})$. At fixed $x$, use $1-e^{-\lambda h}
  \le \lambda h$ in the finite $\lambda$ integral; at fixed $t$, use
  $|e^{i\xi\cdot x}-e^{i\xi\cdot y}| \le |\xi||x-y|$ and the next heat moment.
  The squared triangle inequality gives
  \begin{equation}\label{E:source-derivative-increments}
    \E|Z_b(t,x)-Z_b(s,y)|^2
    \le C_{\rho,b}(|t-s|+|x-y|^2).
  \end{equation}
  For fixed cutoff, the temporal integral is over a bounded range of
  $\lambda$, so its extra factor $\lambda$ remains integrable. The heat factor
  controls the extra powers of $\xi$ needed for the spatial increment.
  Splitting a general increment into a time increment and a spatial increment
  gives the displayed sum of bounds.
  This computation uses~\eqref{E:source-joint-covariance}, including its spatial
  width, and applies to singular spectral measures without a density assumption.

  \medskip
  \noindent\textit{Step 2.~}We construct the jointly continuous versions.
  The maximum of finitely many increments has its $L^q$ norm bounded by the
  $q$th root of the sum of their $q$th moments, without independence.
  Here is a direct continuity construction. In $m = d+1$ parameter dimensions on
  a compact rectangle, even Gaussian $q$th increment moments have root bounded
  by $C\sqrt q$ times the parameter distance to the power $1/2$ at small
  distances. The dyadic grid has at most $C_m2^{nm}$ adjacent edges at level
  $n$. Their maximum has $L^q$ norm at most $C\sqrt q\,2^{-n(1/2-m/q)}$, with
  constants for this rectangle and fixed $\rho,b$. Choose $q$ large and $0 < h_0
  < 1/2-m/q$. The sum of these maxima weighted by $2^{nh_0}$ has finite $L^q$
  norm. Binary refinement, using a fixed number of adjacent edges per coordinate
  at each scale, gives a continuous extension from dyadic points. At each fixed
  parameter the extension agrees with the given Gaussian variable in $L^2$,
  by~\eqref{E:source-derivative-increments}; it is a modification. Use countably
  many nested compact rectangles and all multiindices. The versions agree on
  overlaps on countable dense sets, hence everywhere there by continuity, so
  they patch.

  \medskip
  \noindent\textit{Step 3.~}We identify the continuous fields as spatial
  derivatives.
  They are derivatives. The $\mathcal{G}$-valued fundamental theorem gives,
  along a coordinate segment,
  \[
    \partial_x^bv(t,x+he_i)-\partial_x^bv(t,x)
    =\int_0^h\partial_x^{b+e_i}v(t,x+ue_i)\ud u.
  \]
  Apply bounded linear $\mathcal{I}$ and $L^2$ integration. The continuous
  versions give the same ordinary integral by their Riemann sums, which converge
  both pathwise and in $L^2$. Enforce these identities first on countably many
  rational $t,x,h$ and multiindices, then use joint continuity to extend them to
  every segment in the compact rectangles. The ordinary fundamental theorem
  yields $\partial_iZ_b = Z_{b+e_i}$. Thus, all spatial derivatives exist and
  are jointly continuous with time. The versions are for this fixed $\rho$, not
  a common path event over uncountably many source parameters. This proves
  Lemma~\ref{L:source-pointfield-regularity}.
\end{proof}

\begin{proposition}[Ordinary source response]\label{P:source-ordinary-response}
  For each fixed $\rho > 0$, on its regularity event, put $b_\rho(t,x) =
  \zeta_R(x)Y_\rho(t,x)$. The statewise integral
  \begin{equation}\label{E:source-ordinary-response}
    L_\rho(t)z=\int_0^tS(-u)B(b_\rho(u))S(u)z\ud u
  \end{equation}
  defines an operator-norm Lipschitz path on $E$, and $L_\rho = M_R^\rho$ for
  all $t \in [0,T]$ on one probability-one set. In particular, it has the
  strongly measurable path version supplied by the Gaussian construction.
\end{proposition}

\begin{proof}
  We first construct the ordinary integral path using the samplewise smooth
  coefficient. We then identify its scalar rows with the transferred noise
  evaluations and use countability and continuity to obtain equality of paths.

  \medskip
  \noindent\textit{Step 1.~}We construct the ordinary response on its
  regularity event.
  The coefficient is compact smooth in space and continuous in time in every
  spatial derivative norm. Multiplication estimates and strong continuity of the
  full free group show that the statewise integral exists on every $E_r$. Its
  operator increments $E_r \to E_{r+1}$ are bounded by a finite sample constant
  times $|t-s|$, exactly by the statewise gain argument of
  Lemma~\ref{L:driver-gain-localization}. This proves Lipschitz continuity in
  $\mathcal{B}(E)$, without requiring norm continuity or an operator-norm
  derivative of $S(t)$.

  \medskip
  \noindent\textit{Step 2.~}We identify its rows by physical duality and test
  approximation.
  The row calculation is~\eqref{E:variance-physical-adjoint} of
  Proposition~\ref{P:actual-wave-form}, with $G_w$
  from~\eqref{E:variance-dual-wave}, applied at the fixed time $u$ to the
  spatial coefficient $b_\rho(u,\cdot)$.
  At each fixed $u$, the Sobolev duality calculation uses only spatial
  smoothness of $b_\rho(u,\cdot)$. From $S(-u)(0,f)$ and the two Bessel weights,
  for real Schwartz phases,
  \[
    \langle S(-u)B(b_\rho(u))S(u)z,w\rangle_E
    =\int b_\rho(u,x)[S(u)z]_1(x)G_w(u,x)\ud x.
  \]
  Thus, no prior hypothesis of smoothness in time is invoked for this
  instantaneous identity. For complex pairs conjugate the second factor.
  For the approximation below, Cauchy--Schwarz on the fixed finite time
  interval converts $L^2$ time error to $L^1$ error, and the common compact
  spatial support contributes only its fixed volume.
  Integrating gives the ordinary source-field integral against the specified
  wave test for the original covariance $h_{R,t;z,w}$. For smooth tests, the
  point-feature integral has already been identified with $Y_\rho(F)$. The
  continuous random field version has the same integral: compact Riemann sums
  converge pathwise and in $L^2$ to the ordinary and feature integrals,
  respectively. Approximate the indicator wave test in the covariance completion
  with common compact spatial support and $L^2$-time spatial supremum error
  tending to zero. Its $L^1$ error also tends to zero. The point-feature norms
  are bounded on that support, so their integrals converge in $\mathcal{G}$; the
  original variance domination gives the same $L^2$ limit on the test side.
  Consequently,
  \[
    \langle L_\rho(t)z,w\rangle_E
    =Y_\rho(h_{R,t;z,w})
    =\langle M_R^\rho(t)z,w\rangle_E
    \quad\text{in }L^2.
  \]

  \medskip
  \noindent\textit{Step 3.~}We obtain equality at every time on one event.
  Compare a countable dense Schwartz phase collection and rational times.
  Boundedness and the path continuity of both sides extend this equality to all
  states and all times on one probability-one set. The equality is on the
  coupled space, not only in law. The set is for fixed $\rho$; the $L^p$ limit
  in~\eqref{E:source-driver-limit} needs no simultaneous version for every
  source parameter. This proves Proposition~\ref{P:source-ordinary-response}.
\end{proof}

The regularized noise laws above match the pinned author definitions. The
author's page~9 redefinition keeps the law; the extension preserving the given
noise was proved in Lemma~\ref{L:source-retaining-isometry}. The derivative
estimate here was derived from the correct joint covariance, not from the
differing spatial index in the author's page~8 increment display. No source
correction is adopted. All limits are at fixed finite $R,T$. The nonlinear Young
flow, canonical initial datum $(1,0)$ and its solution-moment transfer remain
separate parts of the argument.

\paragraph{What we now have.}
An operator path driven by the original noise, for a fixed mask and a finite
horizon (Proposition~\ref{P:actual-operator-driver}); the regularized sources
converge to it in the path norm (Section~\ref{S:source-driver-convergence}).
The next chapter lets this path drive a wave.

\chapter{Fixed masks and wave evolution}\label{S:book-wave-evolution}

An operator path describes the forcing. Which wave does it drive? The initial
data hide a catch: the constant displacement does not lie in the Sobolev phase
of the operator estimates. So we carry the constant coordinates explicitly,
with displacement and velocity.

This chapter has no figure. Section~\ref{S:masked-moments}
compares the regularized solutions and their tested
moments. Section~\ref{S:constant-phase-young} builds the Young evolution on
the enlarged phase space (Proposition~\ref{P:constant-young-linear}) and
identifies it with the physical wave. Section~\ref{S:compact-mean-family}
studies compact families of nonnegative potentials with a common
negative trial value. Their positive ground states and uniform gap
(Proposition~\ref{P:compact-ground-family}) control growth in the lower bound.

The inputs are the paths of Chapter~\ref{S:book-operator-paths}
and the canonical expansion of Chapter~\ref{S:G0}. On a
first reading, keep Propositions~\ref{P:constant-young-linear}
and~\ref{P:compact-original-common-trial} and let the proofs wait. The
masks are fixed before source and time limits; none grows with observation
time here. Episode~\VideoRef{5}{2jv8nuqgfBY} of the video series covers
this chapter.

\section{Fixed masks and displacement moments}
\label{S:masked-moments}

We retain Hypothesis~\ref{H:homog}, data $(1,0)$, and all its original
dimensions and covariance amplitudes. Fix $\rho = (\epsilon,\eta) > 0$ and the
source $Y_\rho$ of Proposition~\ref{P:source-canonical-convergence}. For each
fixed $\rho$ use the continuous-time, smooth-space version in
Lemma~\ref{L:source-pointfield-regularity}. Source fields can use the auxiliary
extension preserving the given noise; their law and the given noise remain those
already constructed. Let $m \in C_c^\infty(\R^d)$ be real with $0 \le m \le 1$.
It multiplies the physical source after regularization.

Write $X_{\rho,n}^m(t,x)$ for the ordinary ordered displacement integral with
coefficient $mY_\rho$, and $X_{\rho,n}$ for its unmasked version; both
order-zero terms are one. We first construct these integrals. The positive Green
measures $G_d(r,\ud z)$ have mass $r$ and support in $B_r$, including the sphere
measure in $d = 3$. These are the source
kernels~\cite{chen:author2026:time-dependency}*{(2.1), (2.3)--(2.5), p.~7}; the
ordered convention and Gaussian pairing formula are
\cite{chen:author2026:time-dependency}*{(2.13)--(2.16), p.~10; (2.19), p.~12}.
We use all pair partitions, including self-pairings.

\begin{lemma}[Chains and mixed-order comparison]
  \label{L:masked-chain-comparison}
  The finite ordered integrals are Bochner $L^2$ integrals. For all $n,k \ge 0$,
  fixed $t,x,y$ and positive $\rho$,
  \begin{equation}\label{E:masked-order-comparison}
    0\le \E[X_{\rho,n}^m(t,x)X_{\rho,k}^m(t,y)]
    \le \E[X_{\rho,n}(t,x)X_{\rho,k}(t,y)].
  \end{equation}
  Both expectations vanish when $n+k$ is odd. With $q_n$ from
  Lemma~\ref{L:vertex-crude},
  \begin{equation}\label{E:masked-order-bound}
    \|X_{\rho,n}^m(t,x)\|_2\le q_n(t),\qquad
    \|X_{\rho,n}(t,x)\|_2\le q_n(t).
  \end{equation}
\end{lemma}

\begin{proof}
  We construct the ordered integrals using finite positive chain measures and
  Gaussian moment bounds. We then identify the Green-regularization limit and
  compare every Wick diagram before summing over pairings.

  \medskip
  \noindent\textit{Step 1.~}We construct the finite chain integrals.
  Let $\nu_n^{t,x}$ be the positive time-space chain measure on $0 < s_1 <
  \cdots < s_n < t$, with successive spatial factors $G_d(s_{j+1}-s_j,\ud z_j)$
  and $s_{n+1} = t$, recovering positions backwards from $x$. It uses measures
  rather than sphere densities. The order-zero measure has mass one. Integrating
  the latest vertex gives the recursion $f_n(t) = \int_0^t(t-r)f_{n-1}(r)\ud r$;
  hence,
  \begin{equation}\label{E:masked-chain-mass}
    \nu_n^{t,x}(1)=\frac{t^{2n}}{(2n)!}.
  \end{equation}
  In the covariance and mask factors below, $z_i$ denotes the vertex position
  recovered backwards from the terminal point, not the Green increment used to
  construct the chain measure.
  Starting from mass one at order zero, the recursion integrates the mass of
  the last Green factor, $t-r$, against the previous mass. Substituting
  $f_{n-1}(r) = r^{2n-2}/(2n-2)!$ gives the displayed mass by elementary
  integration. In particular, the measure is finite at each fixed order.
  Measurability follows by representing $G_d(r)$ as $r$ times the dilation by
  $r$ of its unit-time probability measure.

  The covariance is
  \[
    C_\rho(s,z;r,w)=K_\eta(s-r)\Gamma_{2\epsilon}(z-w),
    \qquad 0\le C_\rho\le c_\rho
    :=K_\eta(0)\Gamma_{2\epsilon}(0)<\infty.
  \]
  Positivity uses the nonnegative Radon measure $\Gamma$ and positive heat and
  Laplace kernels. The upper bound also follows by Gaussian Cauchy--Schwarz.
  H\"older's inequality and scalar Gaussian moments bound the absolute
  expectation of any fixed number of source factors, uniformly over their
  vertices. Applied to $2n$ factors this bounds the $L^2$ norm of each $n$-fold
  product. The field is measurable and continuous, so integration against the
  finite chain measure defines a Bochner $L^2$ integral and permits stochastic
  Fubini between two chains.

  The bound on the $L^2$ norm of each product is uniform over the integration
  variables at fixed $\rho,n$. Its integral against the finite chain measure
  is therefore finite. This is the Bochner integrability condition, and
  applying the same argument to two products gives the absolute integrability
  needed to exchange expectation and the two chain integrals.

  \medskip
  \noindent\textit{Step 2.~}We identify the integrals obtained by Green
  smoothing.
  These are the canonical regularized terminal integrals at common noise widths.
  The unmasked $X_{\rho,n}(t,x)$ is $Y_n^{\boldsymbol a}(t,x)$ from
  Lemma~\ref{L:vertex-crude} with every vertexwise noise width in $\boldsymbol
  a$ set equal to $\rho$, so its bound is $q_n(t)$.
  Here is the Green-width identification needed in dimension three. At fixed
  positive $\rho$, replace every Green measure by its convolution with a
  positive heat approximate identity. At fixed times the measures converge
  weakly, preserving masses. Every Wick covariance product is bounded by a power
  of $c_\rho$ and is spatially continuous by
  Lemma~\ref{L:source-pointfield-regularity}. Thus, its spatial integrals
  converge. The chain mass integrands from~\eqref{E:masked-chain-mass} give an
  integrable time majorant. Square and cross moments between smoothed and
  unsmoothed integrals converge to the same diagram sum; therefore, the random
  integrals converge in $L^2$. This is precisely the Green-before-noise limit in
  Lemma~\ref{L:joint-fixed-order}. It also applies in dimensions one and two if
  Green smoothing is used. No samplewise source bound on heat-kernel tails is
  needed. Consequently, Lemma~\ref{L:vertex-crude} applies to the unmasked
  integrals.

  \medskip
  \noindent\textit{Step 3.~}We compare each masked Wick diagram with its
  unmasked diagram.
  If $n+k$ is odd, Gaussian centering gives zero. If $n+k = 2l$, each pair
  partition $D$ contributes the nonnegative finite integral
  \[
    \int \prod_{\{i,j\}\in D}
    C_\rho(s_i,z_i;s_j,z_j)
    \ud(\nu_n^{t,x}\otimes\nu_k^{t,y}).
  \]
  The masked integrand has the additional factor $\prod_{i=1}^{2l}m(z_i) \in
  [0,1]$. Every vertex appears exactly once, irrespective of its chain or
  partner. Summing gives \eqref{E:masked-order-comparison}, including unequal
  and empty chains and every self-pairing. At $n = k = 0$ both sides are one.
  Taking $k = n$, $y = x$ and then applying the unmasked bound gives
  \eqref{E:masked-order-bound}. This is a statement about expected displacement
  diagrams; no samplewise solution order or Gaussian covariance order is used.
  This proves Lemma~\ref{L:masked-chain-comparison}.
\end{proof}

\begin{proposition}[Ordinary solutions and tested moments]
  \label{P:masked-source-moments}
  The sums $u_\rho^m = \sum_{n\ge0}X_{\rho,n}^m$ and $u_\rho =
  \sum_{n\ge0}X_{\rho,n}$ converge absolutely in $L^2$ and have locally
  continuous versions that are the unique locally continuous ordinary mild
  solutions with coefficients $mY_\rho$ and $Y_\rho$, respectively, and data
  $(1,0)$. For every nonnegative $\phi \in C_c^\infty(\R^d)$,
  \begin{align}
    0&\le\E[u_\rho^m(t,x)u_\rho^m(t,y)]
    \le\E[u_\rho(t,x)u_\rho(t,y)],\label{E:masked-two-point}\\
    \E\left|\int\phi(x)u_\rho^m(t,x)\ud x\right|^2
    &\le \|\phi\|_1^2\,\E|u_\rho(t,0)|^2.
    \label{E:masked-average-source}
  \end{align}
  The spatial integral is a Bochner $L^2$ integral.
\end{proposition}

\begin{proof}
  We first identify the mean-square series with the continuous Picard solution
  on each compact observation region. We then pass the mixed-order comparison
  through the series and integrate it against the nonnegative spatial test.

  \medskip
  \noindent\textit{Step 1.~}We construct and identify the ordinary mild
  solutions.
  Put $Q_T = \sum_{n\ge0}q_n(T) < \infty$, using~\eqref{E:vertex-summability}.
  Lemma~\ref{L:masked-chain-comparison} proves absolute $L^2$ convergence
  uniformly in space and widths. To identify the sums pathwise, fix $\rho$, a
  compact observation set and a finite horizon $T$. On the source regularity
  event the coefficient is bounded by a finite sample-dependent $M$ on the
  compact backward-cone enlargement. Finite propagation
  and~\eqref{E:masked-chain-mass} bound each Picard term there by
  $M^nT^{2n}/(2n)!$. This series converges uniformly. Each term is continuous by
  the unit-time Green representation and dominated integration. Summing the
  Volterra recursion gives the ordinary mild equation with initial displacement
  one and velocity zero. Iterating the difference equation between two
  continuous solutions bounds its value on the observation set by $B
  M^nT^{2n}/(2n)!$, where $B$ bounds their difference on the larger compact
  cone. Letting $n \to \infty$ proves uniqueness.

  \medskip
  \noindent\textit{Step 2.~}We identify the continuous version with the
  Bochner series.
  Ordinary Fubini and the absolute product integrability in the preceding lemma
  identify each Bochner term with the pathwise chain integral at a fixed
  terminal point, simultaneously over countably many orders. An almost surely
  convergent subsequence of the $L^2$ partial sums then identifies its limit
  with the continuous Picard version at that point. These versions define the
  solution for all terminal points at this fixed $\rho$; no common event over
  uncountably many widths is needed.

  \medskip
  \noindent\textit{Step 3.~}We sum the expected two-point comparison.
  For either masked or unmasked terms, Cauchy--Schwarz gives
  \[
    \sum_{n,k\ge0}\E|X_{\rho,n}^{(m)}(t,x)
    X_{\rho,k}^{(m)}(t,y)|\le Q_t^2.
  \]
  For the unmasked stationary field, Cauchy--Schwarz bounds its two-point
  expectation by the geometric mean of two equal point second moments.
  Products of partial sums converge in $L^1$. Thus, the preceding lemma passes
  termwise to~\eqref{E:masked-two-point}. The unmasked source law is spatially
  stationary, and its finite Picard integrals commute with translations. Their
  $L^2$ limits inherit stationarity, so $\E[u_\rho(t,x)u_\rho(t,y)] \le
  \E|u_\rho(t,0)|^2$. Stationarity is not required of the masked solution.

  \medskip
  \noindent\textit{Step 4.~}We form the tested spatial average.
  Translating every recovered vertex position by $x$ identifies $\nu_n^{t,x}$
  with the translate of $\nu_n^{t,0}$; on this fixed increment measure,
  varying $x$ changes only the arguments of the source factors and, in the
  masked case, of $m$.
  Each finite Bochner term is $L^2$-continuous in $x$ by source increment
  continuity on compact cones and dominated Gaussian product moments. Uniform
  $q_n$ tails pass this property to the full sum. Its separable $L^2$ image and
  bound $Q_t$ give the Bochner average. Absolute Fubini is justified by
  $\|\phi\|_1^2Q_t^2$. Integrating the two-point inequality against
  $\phi(x)\phi(y)$ proves~\eqref{E:masked-average-source}. This proves
  Proposition~\ref{P:masked-source-moments}.
\end{proof}

\begin{proposition}[Canonical point-moment limit and conditional transfer]
  \label{P:masked-canonical-transfer}
  At every fixed $t$, at arbitrary relative rates of the two widths,
  \begin{equation}\label{E:masked-unmasked-moment-limit}
    \lim_{\epsilon,\eta\downarrow0}\E|u_\rho(t,0)|^2
    =\E|u(t,0)|^2,
  \end{equation}
  where $u$ is the original canonical Stratonovich solution. If the masked
  averages in~\eqref{E:masked-average-source} converge in $L^2$ to
  $Z_m(t,\phi)$, either at arbitrary rates or along a specified sequence with
  both widths tending to zero, then
  \begin{equation}\label{E:masked-conditional-transfer}
    \E|Z_m(t,\phi)|^2\le\|\phi\|_1^2\,\E|u(t,0)|^2.
  \end{equation}
\end{proposition}

\begin{proof}
  We first pass to the canonical limit at each pair of chain orders. A
  summable double majorant then gives the full point moment. The conditional
  masked transfer follows by passing to the limit in the tested source
  inequality.

  \medskip
  \noindent\textit{Step 1.~}We identify the limit for fixed orders.
  The relabelling $r$ below reverses the vertex order separately within each
  chain.
  At every fixed pair of orders the common-width source law is the vertexwise
  law in Lemma~\ref{L:joint-fixed-order}. At terminal $(t,0)$ reverse each chain
  by $s'_j = t-s_{n+1-j}$, $z'_j = -z_{n+1-j}$, and likewise for the other
  chain. Even covariances preserve the reversed differences. Relabeling vertices
  sends each pairing $\{i,j\}$ to $\{r(i),r(j)\}$, giving the exact
  initial-chain convention in Proposition~\ref{P:mixed}, also for an empty
  chain. The preceding Green-width identification removes Green smoothing first.
  Lemma~\ref{L:joint-fixed-order} then gives the canonical limit of each diagram
  at arbitrary relative noise-width rates. Odd mixed moments stay zero. Hence,
  at fixed $n,k$,
  \[
    \E[X_{\rho,n}(t,0)X_{\rho,k}(t,0)]
    \longrightarrow \E[X_n(t,0)X_k(t,0)].
  \]

  \medskip
  \noindent\textit{Step 2.~}We pass from finitely many orders to the full
  moment.
  All absolute mixed moments, before and after the limit, are bounded by
  $q_n(t)q_k(t)$. Given $\delta > 0$, choose a finite square of orders with the
  double majorant outside it less than $\delta$. On the square the finitely many
  limits hold jointly; off it both tails are bounded by $\delta$.
  Proposition~\ref{P:mixed} identifies the canonical sum, proving
  \eqref{E:masked-unmasked-moment-limit}. This argument uses the single-width
  diagram laws, without assigning a new pathwise realization of the entire
  canonical solution family.

  \medskip
  \noindent\textit{Step 3.~}We transfer the bound under the stated
  masked-limit premise.
  Finally, $L^2$ convergence of the masked averages implies convergence of their
  second moments by
  \[
    |\E|B|^2-\E|A|^2|
    \le\|B-A\|_2(\|B\|_2+\|A\|_2).
  \]
  The inequality follows by factoring the difference of squared norms and
  applying Cauchy--Schwarz. Convergence in $L^2$ also bounds both norms on the
  right, so this difference tends to zero. Thus, both sides of the source
  inequality have the limits needed for the comparison.
  Take the same source limit in~\eqref{E:masked-average-source} and
  use~\eqref{E:masked-unmasked-moment-limit}. Existence of the masked limit
  remains an explicit premise here; the next section constructs it on a phase
  containing the constant data. This proves
  Proposition~\ref{P:masked-canonical-transfer}.
\end{proof}

\section{Young evolution with constant initial data}
\label{S:constant-phase-young}

We use the fixed $R,T,E,\gamma,\theta,\beta$ of~\eqref{E:driver-parameters}, and
choose the allowed even mask $m = \zeta_R$ with $0 \le m \le 1$. All constants
in this section can depend on these fixed data, the model and the auxiliary
physical cutoffs. There is no uniform estimate as $R$ or $T$ grows.

\subsection{The distributional phase and its free evolution}

Define the real Hilbert space
\begin{equation}\label{E:constant-phase-space}
  F=E\oplus\R^2,\qquad
  \|(e,c_1,c_2)\|_F^2=\|e\|_E^2+c_1^2+c_2^2,\qquad
  e_0=(0,1,0).
\end{equation}
Its distributional interpretation is $(e,c_1,c_2)\mapsto(e_1+c_1,e_2+c_2)$. Let
$i:E \to F$ denote $e\mapsto(e,0,0)$ and let $\pi_C:F \to \R^2$ be the
constant-coordinate projection.

\begin{lemma}[Constants and finite propagation]\label{L:constant-phase-free}
  The distributional interpretation of $F$ is continuous and injective. The
  space is complete and separable, and $iE$ is closed and proper. For every real
  $\psi \in C_c^\infty(\R^d)$,
  \begin{equation}\label{E:constant-phase-cutoff}
    L_\psi(e,c_1,c_2)=\psi e+(c_1\psi,c_2\psi)
  \end{equation}
  is bounded from $F$ to $E$. The free wave group on $F$ is
  \begin{equation}\label{E:constant-phase-free}
    S_F(t)(e,c_1,c_2)=(S(t)e,c_1+tc_2,c_2).
  \end{equation}
  It is strongly continuous and bounded on finite horizons and has finite
  propagation for the total distributional data, including local cancellation
  between $e$ and the constants.
\end{lemma}

\begin{proof}
  We first separate constant distributions from the Sobolev phase and prove
  the cutoff bound. We then identify the free group as convolution on the
  total distributional data, which gives finite propagation even when the two
  parts cancel locally.

  \medskip
  \noindent\textit{Step 1.~}We prove injectivity and the basic Hilbert-space
  properties.
  A nonzero constant has Fourier transform a nonzero multiple of $\delta_0$.
  Every $H^\sigma$ Fourier transform is locally $L^2$ near zero, since its
  weight is bounded above and below there. A locally $L^2$ function $f$ cannot
  represent $\delta_0$: testing against $\varphi(\xi/\varepsilon)$, with
  $\varphi(0) = 1$ and $\varphi$ compactly supported in $B_1$, bounds its
  integral by
  \[
    \|f\|_{L^2(B_\varepsilon)}
    \varepsilon^{d/2}\|\varphi\|_2\longrightarrow0,
  \]
  whereas $\delta_0$ gives one. Thus, constants lie in neither $H^{1/2}$ nor
  $H^{-1/2}$; the interpretation is injective. Sobolev duality and the two test
  integrals prove continuity into distributions, uniformly on bounded families
  of compact tests. Completeness, separability and the closed proper inclusion
  follow from the Hilbert direct sum.

  \medskip
  \noindent\textit{Step 2.~}We bound localization from the extended phase to
  the Sobolev phase.
  Multiplication by $\psi$ is bounded on both Sobolev entries. Indeed, the
  Peetre inequalities give $\langle\xi\rangle^\sigma\langle\eta\rangle^{-\sigma}
  \le 2^{|\sigma|/2}\langle\xi-\eta\rangle^{|\sigma|}$. In the unitary Fourier
  convention its weighted kernel is therefore bounded by convolution with the
  $L^1$ function $2^{|\sigma|/2}c_d|\widehat{\psi}(q)|\langle
  q\rangle^{|\sigma|}$, where $c_d = (2\pi)^{-d/2}$. The $L^1*L^2$ bound applies
  at $\sigma = \pm1/2$; the two constant columns are also in $E$. This proves
  boundedness of~\eqref{E:constant-phase-cutoff}.

  \medskip
  \noindent\textit{Step 3.~}We identify the free group and finite propagation.
  The finite-dimensional block in~\eqref{E:constant-phase-free} and the strong,
  finite-horizon bound for $S(t)$ from Lemma~\ref{L:driver-gain-localization}
  prove the group assertions. To verify its physical meaning, for $t > 0$
  consider the matrix of compactly supported distributions
  \[
    \begin{pmatrix}
      \partial_tG_d(t)&G_d(t)\\
      \Delta G_d(t)&\partial_tG_d(t)
    \end{pmatrix}.
  \]
  The source Green convention gives $\widehat{G}_d(t)(\xi) =
  c_d\sin(t|\xi|)/|\xi|$; hence its convolution multiplier is the sine quotient,
  and the differentiated multipliers are $\cos(t|\xi|)$ and
  $-|\xi|\sin(t|\xi|)$. The matrix agrees with $S(t)$ on Schwartz phases.
  Convolution by a compact distribution is continuous on tempered distributions:
  its transpose convolves a Schwartz test with the reflected compact
  distribution, continuously in all Schwartz seminorms. Schwartz approximation
  in $E$ extends the identity to $E$. The masses of the four entries are
  $1,t,0,1$, so constants give exactly~\eqref{E:constant-phase-free}. Time
  reversal handles $t < 0$ with the same support bound. The support property
  concerns the total distribution, not its chosen decomposition. Each entry is
  supported in $B_{|t|}$, also for the
  three-dimensional sphere kernel. If the total initial distribution vanishes on
  an enlarged region, each translated kernel tests zero there, so the evolved
  total distribution vanishes on the smaller cone. This argument permits
  cancellation between its Sobolev and constant parts. This proves
  Lemma~\ref{L:constant-phase-free}.
\end{proof}

\subsection{The two columns for constant initial data}

Choose real $\chi \in C_c^\infty$, equal to one on $B_{2R+T+1}$ and supported in
$B_{2R+T+2}$, and choose real compact smooth $\psi = 1$ on a neighborhood of
$\operatorname{supp}\chi$. Write $M_E = M_R$ and $M_E^\rho = M_{Y_\rho}$ for the
paths already constructed, with the same original deterministic coefficients.

\begin{proposition}[Extended operator path]\label{P:constant-phase-driver}
  As identities of continuous operator paths,
  \begin{equation}\label{E:constant-phase-localization}
    D_\chi M_E=M_E=M_ED_\chi,
  \end{equation}
  and the same holds for each fixed $\rho$. Define
  \begin{equation}\label{E:constant-phase-driver}
    A_E=M_EL_\psi,\qquad A=iA_E,\qquad
    A_E^\rho=M_E^\rho L_\psi,\qquad A^\rho=iA_E^\rho.
  \end{equation}
  These are independent of $\psi$, have compact values, and satisfy $Ai = iM_E$,
  $\pi_CA = 0$ and their source counterparts. In every finite $L^pC^\theta$,
  $A^\rho \to A$ at arbitrary relative source-width rates. Their norms have
  uniform $\sqrt p$ moment bounds and a common finite square-exponential moment.
  The original $A$ is measurable with respect to the given $W$.
\end{proposition}

\begin{proof}
  We pass both localizations from the deterministic coefficients to the path
  limits. Right localization then removes the choice of the auxiliary cutoff,
  while bounded composition transfers convergence, compactness and moment
  estimates to the larger phase.

  \medskip
  \noindent\textit{Step 1.~}We pass the localization identities to the path
  limits.
  Lemma~\ref{L:driver-gain-localization} gives both localizations for every
  deterministic $A_j(0,t)$. The original finite sums retain them.
  Proposition~\ref{P:source-gaussian-transfer} constructs the source finite sums
  with exactly these coefficients and correlated $Y_\rho(f_j)$. Their $C^\theta$
  limits preserve both identities. A uniformly convergent subsequence, or
  rational times and continuous paths,
  gives~\eqref{E:constant-phase-localization} for all $t$ on one event. Use a
  separate event for each fixed $\rho$.

  \medskip
  \noindent\textit{Step 2.~}We define the extended path independently of the
  auxiliary cutoff.
  For two choices $\psi,\psi'$, $(L_\psi-L_{\psi'})z \in E$ and
  $D_\chi(L_\psi-L_{\psi'})z = 0$. The right localization kills this difference.
  For $z = ie$, $D_\chi D_\psi e = D_\chi e$ and hence $Ai = iM_E$. The constant
  columns are precisely $iM_E(\psi,0)$ and $iM_E(0,\psi)$; they are retained.
  The range of $i$ has zero
  constant coordinates, which also gives $\pi_C A = 0$.

  \medskip
  \noindent\textit{Step 3.~}We transfer the path bounds and measurability.
  Bounded composition preserves compactness and transfers
  \eqref{E:source-driver-limit} and the Gaussian norm bounds
  to~\eqref{E:constant-phase-driver}. If $\|X\|_{L^p} \le C\sqrt p$, then $\E
  X^{2n} \le C^{2n}(2n)^n$. Since $n! \ge (n/e)^n$, the exponential series gives
  $\E e^{cX^2} < \infty$ uniformly when $2ecC^2 < 1$. Finally, $A$ is a
  deterministic bounded transform of the original $W$-measurable operator path.
  This proves Proposition~\ref{P:constant-phase-driver}.
\end{proof}

\begin{lemma}[Ordinary response on all of $F$]\label{L:constant-source-response}
  For fixed $\rho$ put $b(r) = mY_\rho(r)$ and define
  \[
    B_F(b)(e,c_1,c_2)=(0,b(e_1+c_1))\in E.
  \]
  Then
  \begin{equation}\label{E:constant-source-response}
    A_E^\rho(t)z=\int_0^tS(-r)B_F(b(r))S_F(r)z\ud r
  \end{equation}
  as a statewise strong integral. Its primitive is Lipschitz in
  $\mathcal{B}(F,E)$ and the identity holds for all states and times on one
  event for this fixed $\rho$.
\end{lemma}

\begin{proof}
  We first bound the ordinary coefficient operator on the full phase,
  including the constant displacement. We then use finite propagation to
  replace each initial state by its localized Sobolev representative near the
  coefficient support.

  \medskip
  \noindent\textit{Step 1.~}We construct the statewise integral and its
  Lipschitz bound.
  The source regularity gives one compact spatial support for $b$, continuity in
  $L^\infty$ and $L^2$, and bounded spatial derivatives on the horizon. Since
  $e_1 \in L^2$ and $\|\cdot\|_{H^{-1/2}} \le \|\cdot\|_2$,
  \[
    \|B_F(b)\|_{\mathcal{B}(F,E)}
    \le\|b\|_\infty+\|b\|_2.
  \]
  Strong continuity and finite-horizon bounds of both free groups make the
  integrand strongly continuous on each state and uniformly bounded as an
  operator. Integrating its bound proves operator-norm Lipschitz increments
  without assuming a Banach-valued operator derivative.

  \medskip
  \noindent\textit{Step 2.~}We identify the ordinary response for every
  extended state.
  The total initial distributions of $z$ and $iL_\psi z$ agree on the entire
  enlarged ball where $\psi = 1$. Lemma~\ref{L:constant-phase-free} gives
  equality of the free displacements near $\operatorname{supp}m$ for $0 \le r
  \le T$. Multiplication by $b(r)$ kills the difference. Apply
  Proposition~\ref{P:source-ordinary-response} to the $E$ state $L_\psi z$ to
  obtain~\eqref{E:constant-source-response}. Separability and rational times,
  followed by boundedness and continuity, give its path identity. This proves
  the constant columns directly; no density of $E$ in $F$ is used. This proves
  Lemma~\ref{L:constant-source-response}.
\end{proof}

\subsection{The exact linear Young equation}

We use Young's integration theory~\cite{young:36:inequality} in its
Banach-space form. The estimate is stated by
Friz and Hairer~\cite{friz.hairer:20:course}*{Sections~1.4--1.5, 4.1, (4.3), pp.~13, 61--62}.
For $D,V$ with exponent $\theta > 1/2$ and bounded bilinear multiplication, left
composition gives
\begin{equation}\label{E:constant-young-estimate}
  \left\|\int_s^t\ud D(r)V(r)-D(s,t)V(s)\right\|
  \le C_\theta[D]_\theta[V]_\theta|t-s|^{2\theta}.
\end{equation}
The Banach spaces here are $\mathcal{B}(F)$, and multiplication is bounded with
constant one. We derive the linear solution and stability estimates below rather
than invoke a nonlinear bounded-vector-field theorem.

\begin{proposition}[Linear Young evolution]\label{P:constant-young-linear}
  For every $D \in C^\theta([0,T];\mathcal{B}(F))$ with $D(0) = 0$ there is a
  unique $C^\theta$ solution, in the terminal variable, of
  \begin{equation}\label{E:constant-young-equation}
    U_D(t,s)=I+\int_s^t\ud D(r)U_D(r,s).
  \end{equation}
  It has the composition law, invertible values and joint operator-norm
  continuity on $0 \le s \le t \le T$. The map $D\mapsto U_D$ is locally
  Lipschitz into the uniform time-triangle norm, and
  \begin{equation}\label{E:constant-young-growth}
    \sup_{s\le t\le T}\|U_D(t,s)\|
    \le C_T\exp\big(C_T(1+\|D\|_{C^\theta})^{1/\theta}\big).
  \end{equation}
\end{proposition}

\begin{proof}
  We solve the equation by contraction on short intervals in a scaled Hölder
  norm. We compare two such solutions on the same partition, then patch the
  intervals to obtain the global bound, composition law and invertibility.

  \medskip
  \noindent\textit{Step 1.~}We construct the solution on a short interval.
  On an interval $[a,a+\ell]$ use the complete scaled norm $N_\ell(V) =
  \sup\|V\|+\ell^\theta[V]_\theta$ and put $I_D(V)(t) = \int_a^t\ud D\,V$.
  Applying~\eqref{E:constant-young-estimate} to every subinterval gives
  \begin{align*}
    \sup\|I_D(V)\|
    &\le[D]_\theta\ell^\theta\sup\|V\|
    +C_\theta[D]_\theta\ell^{2\theta}[V]_\theta,\\
    [I_D(V)]_\theta
    &\le[D]_\theta\sup\|V\|
    +C_\theta[D]_\theta\ell^\theta[V]_\theta.
  \end{align*}
  For the supremum estimate, use the Young remainder bound between $a$ and
  $t$. For the seminorm estimate, apply it to each subinterval and divide by
  its length to the power $\theta$. Multiplying the second estimate by
  $\ell^\theta$ and adding it to the first gives the scaled-norm estimate.
  Thus, with $c_\theta = 2\max(1,C_\theta)$,
  \begin{equation}\label{E:constant-young-local}
    N_\ell(I_D(V))
    \le c_\theta[D]_\theta\ell^\theta N_\ell(V).
  \end{equation}
  Take $K \ge 1$, $[D]_\theta \le K$, and choose $\ell$ with $c_\theta
  K\ell^\theta \le 1/4$. The map $V\mapsto V_a+I_D(V)$ is a contraction on the
  complete affine space of paths starting at $V_a$. Its solution satisfies
  $N_\ell(V) \le (4/3)\|V_a\| \le 2\|V_a\|$.

  Subtracting two images of the fixed-point map removes the initial value, so
  its Lipschitz constant is at most $1/4$. The image starts at $V_a$ because
  the integral vanishes at the left endpoint. Finally, the fixed-point
  equation gives $N_\ell(V) \le \|V_a\|+N_\ell(V)/4$, which yields the factor
  $4/3$.

  \medskip
  \noindent\textit{Step 2.~}We compare solutions driven by two paths.
  To compare $D,D'$, require both seminorms at most $K$ and write $\delta =
  \|D-D'\|_{C^\theta}$, $W = V-V'$. Subtracting gives $W =
  W_a+I_D(W)+I_{D-D'}(V')$. Absorbing the factor $1/4$ yields
  \begin{equation}\label{E:constant-young-difference}
    N_\ell(W)\le2\|W_a\|
    +2c_\theta\ell^\theta\delta N_\ell(V').
  \end{equation}
  In the difference equation, the first integral contains the unknown
  difference $W$ and has coefficient at most $1/4$. Move it to the left. The
  remaining integral is bounded by the same local estimate with driver $D-D'$,
  whose seminorm is at most $\delta$. Replacing $4/3$ by two gives the
  displayed bound.
  Cover $[s,T]$ by at most $n_{K,T}$ such pieces. If the second initial norm is
  $B$, its norm on piece $j$, indexed from zero, is at most $2^{j+1}B$. With
  $E_j$ the difference at its left endpoint, induction gives
  \[
    E_{j+1}\le2E_j+2^{j+2}c_\theta\ell^\theta\delta B,\qquad
    E_j\le2^jE_0+j2^{j+1}c_\theta\ell^\theta\delta B.
  \]
  Equation~\eqref{E:constant-young-difference} gives the same bound inside each
  piece; a shorter final piece only reduces the coefficient. For evolutions,
  $E_0 = 0$ and $B = 1$. The number of pieces is uniform in $s,t$, proving local
  Lipschitz dependence in the full triangle norm.

  \medskip
  \noindent\textit{Step 3.~}We patch the intervals and prove the flow
  properties.
  Now choose $K = 1+\|D\|_{C^\theta}$. At most
  $C_T(1+\|D\|_{C^\theta}^{1/\theta})$ pieces suffice. Each multiplies the norm
  by at most two, giving \eqref{E:constant-young-growth}. The product of these
  factors is bounded by an exponential in the number of pieces.   Uniqueness
  patches the
  solutions; right multiplication by the preceding evolution gives the
  composition law. On a small piece~\eqref{E:constant-young-local} gives
  $\|U_D-I\|_\infty \le 1/3$, so a Neumann series proves local invertibility,
  then products give global invertibility. Finally, $U_D(t,s) =
  U_D(t,0)U_D(s,0)^{-1}$ proves joint norm continuity, since $U_D(\cdot,0)$ is
  continuous and inversion is continuous on invertible operators. This proves
  Proposition~\ref{P:constant-young-linear}.
\end{proof}

Proposition~\ref{P:constant-young-linear} concerns one deterministic operator
path. We now apply it to the random paths $A$ and $A^\rho$
of~\eqref{E:constant-phase-driver} and pass to the physical flow on $F$.

\begin{proposition}[Gaussian Young and physical flows]
  \label{P:constant-gaussian-flow}
  Apply Proposition~\ref{P:constant-young-linear} to $A,A^\rho$. The resulting
  flows are strongly measurable in uniform triangle norm, have every finite norm
  moment and satisfy
  \begin{equation}\label{E:constant-flow-limit}
    U_{A^\rho}\longrightarrow U_A
    \quad\hbox{in every finite }L^p
    \left(\sup_{0\le s\le t\le T}\|\cdot\|_{\mathcal{B}(F)}\right).
  \end{equation}
  They preserve the interaction constant coordinates:
  \begin{equation}\label{E:constant-flow-range}
    \pi_CU_A(t,s)=\pi_C,\qquad
    (U_A(t,s)-I)F\subset iE,
  \end{equation}
  and likewise for the source. The physical flow is
  \begin{equation}\label{E:constant-physical-flow}
    P_A(t,s)=S_F(t)U_A(t,s)S_F(-s).
  \end{equation}
  It has composition and strong continuity on states, physical constant block
  $(c_1,c_2)\mapsto(c_1+(t-s)c_2,c_2)$, all finite moments and the analogous
  uniform $L^p$ source convergence. Its difference from the free flow has range
  in $iE$.
\end{proposition}

\begin{proof}
  We first obtain measurability and the constant-coordinate identities. Local
  stability gives convergence in probability after restricting to bounded
  driver norms. Gaussian tails give the moments needed to strengthen this to
  $L^p$ convergence, and bounded free-group composition transfers it to the
  physical flow.

  \medskip
  \noindent\textit{Step 1.~}We obtain measurable flows and their range
  identities.
  The original paths lie in separable closed path spans by their Gaussian
  construction. Bounded composition in~\eqref{E:constant-phase-driver} and the
  continuous Young map preserve strong measurability into
  $C(\Delta_T;\mathcal{B}(F))$, where $\Delta_T = \{(t,s):0\le s\le t\le T\}$
  has the uniform norm. Indeed, the continuous image of the separable space of
  operator paths has separable closure in this triangle space. Applying $\pi_C$
  to Young Riemann sums gives the first identity
  in~\eqref{E:constant-flow-range}. Each increment applied to any state lies in
  $iE$; closedness of this subspace passes the property through sums and limits,
  giving the second identity.

  \medskip
  \noindent\textit{Step 2.~}We prove convergence in probability of the Young
  flows.
  Write $X_\rho = \|A^\rho\|_{C^\theta}$ and $X_0 = \|A\|_{C^\theta}$.
  Proposition~\ref{P:constant-phase-driver} supplies $c_0 > 0$ and $H_0 <
  \infty$ such that
  \[
    \sup_\rho\E e^{c_0X_\rho^2}\le H_0,\qquad
    \E e^{c_0X_0^2}\le H_0.
  \]
  For $K \ge 1$, let $L_K$ be the Lipschitz constant of
  Proposition~\ref{P:constant-young-linear} on the ball of operator paths of
  radius $K$, in the whole time-triangle norm. On $\{X_\rho\le K,\ X_0\le K\}$
  it gives
  \[
    \|U_{A^\rho}-U_A\|_\infty
    \le L_K\|A^\rho-A\|_{C^\theta}.
  \]
  Consequently, for $\varepsilon > 0$ and any finite $q \ge 1$, Markov and the
  union bound give
  \[
    \mathbb{P}(\|U_{A^\rho}-U_A\|_\infty>\varepsilon)
    \le2H_0e^{-c_0K^2}
    +(L_K/\varepsilon)^q\E\|A^\rho-A\|_{C^\theta}^q.
  \]
  The exceptional event is contained in the union where either driver norm
  exceeds $K$. Markov applied to each square exponential bounds this union by
  the first term. On its complement, the local Lipschitz estimate and the
  ordinary moment form of Markov give the second term.
  At fixed $K$ the second term tends to zero at arbitrary relative source-width
  rates. Only then let $K \to \infty$ to obtain convergence in probability. No
  independence of the original, source and residual operator paths is required.

  \medskip
  \noindent\textit{Step 3.~}We obtain flow moments and strengthen convergence
  to $L^p$.
  Put $\nu = 1/\theta < 2$. For each finite $r \ge 1$ and $c > 0$, there is
  $C_{r,c} < \infty$ such that
  \[
    rC_T(1+x)^\nu\le cx^2+C_{r,c},\qquad x\ge0.
  \]
  This follows because $(1+x)^\nu/x^2 \to 0$ as $x \to \infty$, while the
  remainder is bounded on compact intervals. Use $c \le c_0$
  and~\eqref{E:constant-young-growth} to obtain
  \[
    \sup_\rho\E\|U_{A^\rho}\|_\infty^r<\infty,\qquad
    \E\|U_A\|_\infty^r<\infty.
  \]
  In particular, for fixed finite $p \ge 1$ and $V_\rho =
  \|U_{A^\rho}-U_A\|_\infty$, the $2p$ moments of $V_\rho$ are uniformly
  bounded. For every $\varepsilon > 0$,
  \[
    \E V_\rho^p
    \le\varepsilon^p+
    (\E V_\rho^{2p})^{1/2}
    \mathbb{P}(V_\rho>\varepsilon)^{1/2}.
  \]
  Split the expectation into $V_\rho \le \varepsilon$ and its complement. On
  the first event the integrand is at most $\varepsilon^p$. On the second,
  Cauchy--Schwarz pairs $V_\rho^p$ with the indicator, giving the displayed
  probability factor. The uniform $2p$ moment bound therefore supplies uniform
  integrability.
  The second term tends to zero by convergence in probability; then let
  $\varepsilon\downarrow0$. This proves \eqref{E:constant-flow-limit}. The
  square-exponential bound used here is an operator path bound; no
  square-exponential flow bound is needed.

  \medskip
  \noindent\textit{Step 4.~}We transfer the conclusions to the physical flow.
  The free group is strongly continuous and bounded on the horizon. Define
  $\mathcal{B}_\infty(\Delta_T;\mathcal{B}(F))$ to be the Banach space of all
  bounded operator families on $\Delta_T$. The deterministic map
  \[
    \mathcal{L}(V)(t,s)=S_F(t)V(t,s)S_F(-s)
  \]
  is bounded linear from $C(\Delta_T;\mathcal{B}(F))$ into this bounded-family
  space, with norm at most $C_T^2$, where $C_T \ge \sup_{|r|\le T}\|S_F(r)\|$.
  Hence, it preserves strong measurability on the separable ranges and gives, in
  particular,
  \[
    \|P_{A^\rho}-P_A\|_\infty
    \le C_T^2\|U_{A^\rho}-U_A\|_\infty.
  \]
  For the range identity below, the rightmost free-group factor changes
  the input, while the leftmost factor preserves $iE$. Thus, composing
  $U_A(t,s)-I$ with these factors preserves its range inclusion in $iE$.
  All finite moments and uniform $L^p$ convergence therefore transfer to the
  physical flows. This does not place those flows in the norm-continuous
  triangle space or assert operator-norm continuity of the free group. For each
  fixed state, joint continuity follows by splitting a physical-flow difference
  into the norm-continuous $U$ difference and two strong free-group differences,
  using the uniform group bound. Group cancellation gives composition. The
  finite constant block $J_t(c_1,c_2) = (c_1+tc_2,c_2)$ satisfies $J_tJ_{-s} =
  J_{t-s}$, which gives the asserted physical constant coordinates from
  $\pi_CU_A = \pi_C$. Finally,
  \[
    P_A(t,s)-S_F(t-s)
    =S_F(t)(U_A(t,s)-I)S_F(-s)
  \]
  has range in $iE$, since $S_F(t)iE = iE$. This proves
  Proposition~\ref{P:constant-gaussian-flow}.
\end{proof}

\subsection{Ordinary source solutions and the original tested observable}

So far $P_{A^\rho}$ is defined through Young integrals, not through the source
equation of Proposition~\ref{P:masked-source-moments}. We now connect the two
at a fixed source width $\rho$.

\begin{lemma}[Identification with the ordinary source equation]
  \label{L:constant-source-solution}
  For fixed $\rho$, the displacement distribution of $P_{A^\rho}(t,0)e_0$ is the
  locally continuous ordinary source solution $u_\rho^m(t,\cdot)$ from
  Proposition~\ref{P:masked-source-moments}.
\end{lemma}

\begin{proof}
  We first identify the Young integral with the ordinary interaction integral
  for fixed regularization. We then compare the phase Picard terms with the
  scalar Green-chain terms by distributional testing and pass to their limits.

  \medskip
  \noindent\textit{Step 1.~}We identify the statewise Young equation with the
  ordinary equation.
  In~\eqref{E:constant-source-response} let $C(r) = S(-r)B_F(b(r))S_F(r)$. The
  primitive is $A_E^\rho$. For left endpoints $r_j = t_j$, the Young sums
  applied to a state $v$ are
  \[
    \sum_j\int_{t_j}^{t_{j+1}}iC(r)U_{A^\rho}(r_j,s)v\ud r.
  \]
  Uniform continuity of $U$ and the uniform bound for $C$ make these converge to
  the strong integral with $U(r,s)v$. Thus, the Young equation is precisely the
  ordinary interaction integral on states. Multiplying by free groups gives the
  Duhamel equation for displacement and velocity. Its unique $C([0,T];F)$
  solution follows by the bounded-coefficient Picard argument with the strongly
  continuous, bounded free group. No operator-norm derivative of $S_F$ is used.

  \medskip
  \noindent\textit{Step 2.~}We construct the phase Picard series.
  More explicitly, put $W_0(t) = e_0$ and, for $n \ge 1$,
  \[
    W_n(t)=\int_0^t S_F(t-r)iB_F(b(r))W_{n-1}(r)\ud r.
  \]
  The ordered $n$-fold time simplex has volume $T^n/n!$.
  If $B_T = \sup_r\|B_F(b(r))\|_{\mathcal{B}(F,E)}$, induction on the time
  simplex gives $\sup_t\|W_n(t)\|_F \le C_T^{n+1}B_T^nT^n/n!$. We take $C_T \ge
  \max(1,\sup_{|t|\le T}\|S_F(t)\|)$, which also supplies the bound at order
  zero. The series converges in $C([0,T];F)$. The Duhamel map there has norm at
  most $C_TB_TT$, so its sum solves the same integral equation and equals the
  ordinary phase by uniqueness.

  \medskip
  \noindent\textit{Step 3.~}We identify the displacement of each term and then
  of the sum.
  Let $X_0(t,x) = 1$ and $X_n(t) = \int_0^tG_d(t-r)*(b(r)X_{n-1}(r))\ud r$. At
  fixed $n$ these are continuous in time and smooth in space. Derivatives fall
  on the smooth compact factor $bX_{n-1}$; Green mass $t-r$ bounds the
  integrals, while the unit-time measure representation gives time continuity.
  This is induction using only finitely many bounded derivatives for each fixed
  derivative order. At $n = 0$ the phase displacement is $X_0$. If the equality
  holds at $n-1$, multiplication by smooth $b$ identifies the forcing with
  $(0,f_r)$, where $f_r = b(r)X_{n-1}(r)$ is smooth and compactly supported. For
  a compact smooth test $\varphi$, continuous distributional testing of the
  phase integral yields
  \[
    \langle[W_n(t)]_1,\varphi\rangle
    =\int_0^t\int f_r(y)
    \int G_d(t-r,\ud z)\varphi(y+z)\ud y\ud r
    =\langle X_n(t),\varphi\rangle.
  \]
  Absolute Fubini follows from $\|\varphi\|_\infty\|f_r\|_1(t-r)$, integrable on
  the horizon because $f_r$ has one compact support and continuous $L^1$ norm.
  This also covers sphere measures in $d = 3$. Iterating the scalar recursion
  gives exactly $X_{\rho,n}^m$, with chronological order and all physical mask
  factors. The scalar compact-cone factorial bound gives locally uniform
  convergence to $u_\rho^m$, whereas the phase series converges in $F$ and hence
  in distributions. Their displacement limits agree. This proves
  Lemma~\ref{L:constant-source-solution}.
\end{proof}

Lemma~\ref{L:constant-source-solution} holds at each fixed $\rho$. To reach
the original noise we test against a fixed $\phi$ and use the source
limit~\eqref{E:constant-flow-limit}.

\begin{proposition}[Original average and original point moment]
  \label{P:constant-raw-average}
  Let real $\phi \in C_c^\infty(\R^d)$ be nonnegative with integral one. The
  bounded functional and original random variable
  \begin{equation}\label{E:constant-raw-average}
    O_\phi(e,c_1,c_2)=\langle e_1,\phi\rangle+c_1,\qquad
    Z_m(t,\phi)=O_\phi P_A(t,0)e_0
  \end{equation}
  are well defined, and $Z_m$ is $W$-measurable. At arbitrary relative
  source-width rates, for every finite $p$,
  \begin{equation}\label{E:constant-tested-limit}
    \int\phi(x)u_\rho^m(t,x)\ud x
    \longrightarrow Z_m(t,\phi)\quad\text{in }L^p.
  \end{equation}
  In particular,
  \begin{equation}\label{E:constant-original-point-transfer}
    \E|Z_m(t,\phi)|^2\le\E|u(t,0)|^2.
  \end{equation}
\end{proposition}

\begin{proof}
  We first identify the ordinary tested source observable with the phase
  functional. We then show that it is the same random variable as the Bochner
  average. This allows the flow convergence to supply the masked-limit premise
  for the original moment comparison.

  \medskip
  \noindent\textit{Step 1.~}We identify the phase functional and the ordinary
  source average.
  Sobolev duality gives $\|O_\phi\| \le (\|\phi\|_{H^{-1/2}}^2+1)^{1/2}$. The
  original operator path is $W$-measurable, and the Young and deterministic maps
  above are continuous or bounded measurable on their separable ranges. Thus,
  $Z_m$ is $W$-measurable and does not depend on the auxiliary $J$ used in
  source realizations. At each fixed $\rho$,
  Lemma~\ref{L:constant-source-solution} identifies the displacement as a
  distribution with the ordinary source solution. Testing by this same $\phi$,
  with $\int\phi = 1$, therefore gives the almost sure identity
  \[
    O_\phi P_{A^\rho}(t,0)e_0
    =C_\phi:=\int\phi(x)u_\rho^m(t,x)\ud x.
  \]
  This uses the displacement pairing and its constant displacement coordinate
  in~\eqref{E:constant-raw-average}; it takes no point trace of an $H^{1/2}$
  state. All identities here are at one fixed $\rho$, so no uncountable
  intersection of source-version events is used.

  \medskip
  \noindent\textit{Step 2.~}We identify the ordinary and Bochner averages in
  $L^2$.
  We verify that $C_\phi$ equals the Bochner $L^2$ average $B_\phi$ from
  Proposition~\ref{P:masked-source-moments}. That proposition gives the
  spatially uniform bound $\|u_\rho^m(t,x)\|_2 \le Q_t$, with $Q_t$ also
  independent of source widths. Scalar Fubini gives
  \[
    \E\int|\phi(x)u_\rho^m(t,x)|\ud x
    \le Q_t\|\phi\|_1<\infty.
  \]
  Thus, the ordinary spatial integral exists almost surely. Weighted
  Cauchy--Schwarz and scalar Fubini also give
  \[
    \E|C_\phi|^2
    \le\|\phi\|_1\int|\phi(x)|\E|u_\rho^m(t,x)|^2\ud x
    \le Q_t^2\|\phi\|_1^2.
  \]
  For every bounded real random variable $\Xi$ on the fixed source space,
  $V\mapsto\E[\Xi V]$ is a continuous functional on $L^2$. It therefore commutes
  with the Bochner integral, while scalar Fubini is justified by $\|\Xi\|_\infty
  Q_t\|\phi\|_1$. Consequently,
  \[
    \E[\Xi B_\phi]
    =\int\phi(x)\E[\Xi u_\rho^m(t,x)]\ud x
    =\E[\Xi C_\phi].
  \]
  Both averages are real $L^2$ variables. Taking $\Xi =
  \operatorname{sign}(B_\phi-C_\phi)$ yields $\E|B_\phi-C_\phi| = 0$, so the two
  averages agree almost surely.

  \medskip
  \noindent\textit{Step 3.~}We pass to the original tested observable and
  moment bound.
  Proposition~\ref{P:constant-gaussian-flow} now gives, for every finite $p$,
  \[
    \|B_\phi-Z_m(t,\phi)\|_{L^p}
    \le\|O_\phi\|\|e_0\|
    \big\|\|P_{A^\rho}-P_A\|_\infty\big\|_{L^p}
    \longrightarrow0.
  \]
  This is~\eqref{E:constant-tested-limit}, including its $L^2$ case, for the
  masked Bochner average, with arbitrary relative source-width rates. It
  supplies the conditional convergence premise of
  Proposition~\ref{P:masked-canonical-transfer} at the same fixed mask, horizon,
  original covariance and observable. Since $\phi \ge 0$ and $\|\phi\|_1 = 1$,
  that proposition gives \eqref{E:constant-original-point-transfer}. Zero
  covariance gives the free constant solution and equality. All source limits
  are at fixed mask and horizon, before time asymptotics. This proves
  Proposition~\ref{P:constant-raw-average}.
\end{proof}

So far the operator path comes from the noise alone.
Section~\ref{S:compact-mean-family} also needs a deterministic potential.

\begin{proposition}[Compact deterministic base]\label{P:constant-compact-base}
  Let $h \in C([0,T];C_c^\infty(K))$ be real for one fixed compact $K$, with all
  spatial derivatives continuous on that support. The statewise primitive
  \[
    A_h(t)=\int_0^t iS(-r)B_F(h(r))S_F(r)\ud r
  \]
  is Lipschitz in $\mathcal{B}(F)$, takes values with range in $iE$, and has a
  bound $C_T\sup_r(\|h(r)\|_\infty+\|h(r)\|_2)$. The flow conclusions above,
  including moments, source convergence and ordinary-source identification, hold
  for $A_h+gA$ and $A_h+gA^\rho$, uniformly on bounded displayed deterministic
  bounds and bounded real $g$. No Cameron-shift identity is asserted here.
\end{proposition}

\begin{proof}
  We bound the deterministic primitive and then add it to the Gaussian driver.
  Its bounded shift preserves the required exponential control, so the
  preceding stability and ordinary-source arguments apply.

  \medskip
  \noindent\textit{Step 1.~}We bound the shifted driver and its moments.
  The estimate for $B_F$ and the strong free groups prove the primitive and its
  Lipschitz bound exactly as in Lemma~\ref{L:constant-source-response}. If
  $\|A_h\|_{C^\theta} \le H$, $|g| \le G$ and $X = \|A\|_{C^\theta}$, then
  $\|A_h+gA\|_{C^\theta} \le H+GX$. If $\E e^{cX^2}$ is bounded, choose $c' \le
  c/(2(1+G^2))$ to obtain
  \[
    \E e^{c'\|A_h+gA\|_{C^\theta}^2}
    \le e^{2c'H^2}\E e^{2c'G^2X^2}<\infty.
  \]
  The elementary inequality $(H+GX)^2 \le 2H^2+2G^2X^2$ gives the displayed
  exponential comparison. The choice of $c'$ ensures that the coefficient of
  $X^2$ stays below the available Gaussian exponential constant. The
  deterministic factor is uniform when $H$ and $G$ are bounded.
  Alternatively $pC(1+H+Gx)^{1/\theta} \le cx^2+C_{p,H,G,c}$ transfers the Young
  growth estimate directly to all finite moments.

  \medskip
  \noindent\textit{Step 2.~}We transfer convergence and the ordinary-source
  identification.
  Source operator-path convergence, local Young stability and uniform
  integrability give the same
  flow convergence. The state-integral and Picard-to-chain arguments apply to
  $h+gmY_\rho$, whose spatial support lies in one fixed compact union. These are
  individual estimates and versions, without an uncountable common-parameter
  event. Only $h = 0$, $g = 1$ enters
  \eqref{E:constant-original-point-transfer}. This proves
  Proposition~\ref{P:constant-compact-base}.
\end{proof}

\section{Compact mean families and the spectral gap for wave evolution}
\label{S:compact-mean-family}

Let $I = [0,1]$, let $K_x \subset \R^d$ be compact, and let $\sigma\mapsto
V_\sigma$ be continuous from $I$ into $C_c^\infty(K_x)$, real and nonnegative.
The index $\sigma$ is a family parameter, distinct from physical propagation
time $t$. Assume that one real $w \in H^1(\R^d)$ with $\|w\|_2 = 1$ satisfies
\begin{equation}\label{E:compact-common-trial}
  \|\nabla w\|_2^2-\int V_\sigma w^2\le-\varepsilon
  \quad\hbox{for all }\sigma\in I,\qquad \varepsilon>0.
\end{equation}
Put $M = \sup_I\|V_\sigma\|_\infty$ and $H_\sigma = -\Delta-V_\sigma$. The
phases $E,F$, inclusion $i$, data $e_0$ and observable $O_\phi$ are exactly
those of Section~\ref{S:constant-phase-young}.
Proposition~\ref{P:constant-compact-base}, with deterministic coefficient $h =
V_\sigma$ and noise amplitude zero, gives the ordinary physical flow
$P_\sigma(t)$ on $F$. Its $E$ restriction is $P_\sigma^E(t)$. Uniqueness of the
Duhamel solution makes finite-horizon constructions agree on overlaps; integer
horizons and autonomy give the same semigroup for every $t \ge 0$.

We import the following statements from the \emph{second edition}
of~\cite{teschl:14:mathematical}: Theorem~10.2, p.~258, gives the self-adjoint
$H^2$ operator, lower bound, essential spectrum $[0,\infty)$ and compact smooth
core for real compact smooth potentials in the dimensions considered here.
Theorem~10.13, p.~273, gives a simple positive ground state \emph{if the
spectral bottom is an eigenvalue}. Theorem~4.12 and Corollary~4.13,
pp.~140--141, give max--min and form monotonicity, with the essential threshold
used after discrete eigenvalues end and with the quadratic-form domain allowed.
Uniformity, fractional norms and all applications to $F$ are proved below.

\subsection{The ground family and the global phase}

\begin{proposition}[Ground states and a uniform rate gap]
  \label{P:compact-ground-family}
  The common form domain of $H_\sigma$ is $H^1$ with form $q_\sigma(u) =
  \|\nabla u\|_2^2-\int V_\sigma|u|^2$. Its bottom $E_1(\sigma)$ is a simple
  eigenvalue with a unique real positive-a.e. normalized ground state
  $p_\sigma$. Let $E_2(\sigma)$ be the second max--min level, set to zero when
  there is no second negative eigenvalue. Then
  \begin{align}
    -M&\le E_1(\sigma)\le-\varepsilon,\qquad
    E_1(\sigma)<E_2(\sigma)\le0,\label{E:compact-energy-levels}\\
    |E_j(\sigma)-E_j(\tau)|
    &\le\|V_\sigma-V_\tau\|_\infty,\qquad j=1,2.
    \label{E:compact-level-continuity}
  \end{align}
  The map $\sigma\mapsto p_\sigma$ is continuous into $H^2$ and bounded there.
  With $\lambda_\sigma = \sqrt{-E_1(\sigma)}$ and $\mu_\sigma =
  \sqrt{-E_2(\sigma)}$,
  \begin{equation}\label{E:compact-rate-gap}
    \sqrt\varepsilon\le\lambda_\sigma\le\sqrt M,\qquad
    d_0:=\inf_I(\lambda_\sigma-\mu_\sigma)>0.
  \end{equation}
\end{proposition}

\begin{proof}
  We first identify the common form domain and place the spectral bottom below
  the essential threshold. We then use simplicity and continuity to obtain a
  uniform gap, and compare ground states to prove their continuity in $H^2$.

  \medskip
  \noindent\textit{Step 1.~}We identify the common form domain.
  The imported Theorem~10.2 applies to $-V_\sigma$, which is real, compactly
  supported and in $L^2$ for $d \le 3$. It gives the stated self-adjoint
  operator on $H^2$ and essential threshold zero. Write $c = 1+M$ and $A_\sigma
  = H_\sigma+c \ge 1$. On $H^2$, integration by parts and the spectral theorem
  give
  \[
    \|A_\sigma^{1/2}u\|_2^2
    =q_\sigma(u)+c\|u\|_2^2,\qquad
    \|u\|_{H^1}^2\le\|A_\sigma^{1/2}u\|_2^2
    \le c\|u\|_{H^1}^2.
  \]
  Since $0 \le V_\sigma \le M$ and $c = 1+M$, the shifted potential term
  $c-V_\sigma$ lies between one and $c$. This gives the two form bounds on
  $H^2$. The issue is then to identify their completions, rather than assume
  the form domain from the formal expression.
  Spectral cutoffs $1_{[1,n]}(A_\sigma)u$ show that $H^2$ is a core for the
  half-power form norm. Fourier cutoffs show it is dense in $H^1$. The two
  equivalent completions are therefore $D(A_\sigma^{1/2})$ and $H^1$. Both norms
  dominate $L^2$ and give the same $L^2$ limit, so these are equal as subspaces,
  with the displayed form. This verifies the form domain used below.

  \medskip
  \noindent\textit{Step 2.~}We obtain the isolated simple ground state.
  The form lower bound and~\eqref{E:compact-common-trial} give $-M \le E_1 \le
  -\varepsilon$ by the imported form max--min principle. The bottom belongs to
  the spectrum and lies below the essential threshold, so it is an isolated
  eigenvalue of finite multiplicity. Only now apply the imported Theorem~10.13:
  it is simple with a positive ground state. Conjugation preserves the real
  operator; positivity and unit $L^2$ norm fix its real representative.

  \medskip
  \noindent\textit{Step 3.~}We make the spectral and rate gaps uniform.
  Simplicity, isolation and threshold zero give $E_1 < E_2 \le 0$. For every
  common unit form test, $|q_\sigma(u)-q_\tau(u)| \le
  \|V_\sigma-V_\tau\|_\infty$. Taking the same infima over orthogonal
  complements and suprema in the source max--min formula proves
  \eqref{E:compact-level-continuity}, including passage of $E_2$ to the
  threshold. Compactness gives $\Delta_*: = \min_I(E_2-E_1) > 0$. Thus,
  \[
    \lambda_\sigma-\mu_\sigma
    =\frac{E_2(\sigma)-E_1(\sigma)}
    {\lambda_\sigma+\mu_\sigma}
    \ge\frac{\Delta_*}{2\sqrt M}>0.
  \]
  Here $M \ge \varepsilon > 0$, and no division by $\mu_\sigma$ occurs.

  \medskip
  \noindent\textit{Step 4.~}We prove continuity of the positive normalized
  ground state.
  Put $d_V = \|V_\sigma-V_\tau\|_\infty$ and decompose $p_\tau = a p_\sigma+r$,
  with $r\perp p_\sigma$ and $0 < a \le 1$. The variational bound gives
  $q_\sigma(p_\tau) \le E_1(\sigma)+2d_V$, whereas the complementary spectral
  bound gives $q_\sigma(p_\tau) \ge E_1(\sigma)+\Delta_*\|r\|_2^2$. Hence,
  \begin{equation}\label{E:compact-ground-continuity}
    \|p_\tau-p_\sigma\|_2^2
    =2(1-a)\le2(1-a^2)\le4d_V/\Delta_*.
  \end{equation}
  Positivity of both ground states gives $a = \langle p_\tau,p_\sigma\rangle >
  0$. Orthogonality gives $a^2+\|r\|_2^2 = 1$. The upper form bound follows by
  changing the potential and then the lowest eigenvalue, each at cost at most
  $d_V$. The complementary lower bound uses the uniform gap on $r$. Combining
  them gives $\|r\|_2^2 \le 2d_V/\Delta_*$ and hence the displayed distance
  estimate.
  The eigenfunction equation is
  \[
    (1-\Delta)p_\sigma
    =(1-\lambda_\sigma^2+V_\sigma)p_\sigma.
  \]
  Its right side is continuous in $L^2$ and bounded there by $1+2M$. The inverse
  of $1-\Delta$ is an isometry from $L^2$ to the Bessel $H^2$ space, proving the
  asserted $H^2$ continuity and bound. This proves
  Proposition~\ref{P:compact-ground-family}.
\end{proof}

The rate gap~\eqref{E:compact-rate-gap} is a statement about the spectrum of
$H_\sigma$. The wave on $E$ is measured in the fixed norm of $E$, so we first
compare the two.

\begin{lemma}[Uniform fractional phase norms]\label{L:compact-fractional-phase}
  Let $B = 1-\Delta$, $c = 1+M$ and $A_\sigma = H_\sigma+c$. The quarter-power
  domain is $H^{1/2}$ and
  \begin{equation}\label{E:compact-quarter-power}
    \|B^{1/4}u\|_2\le\|A_\sigma^{1/4}u\|_2
    \le c^{1/4}\|B^{1/4}u\|_2.
  \end{equation}
  The norm
  \begin{equation}\label{E:compact-spectral-phase-norm}
    \|(u,z)\|_{E_\sigma}^2
    =\|A_\sigma^{1/4}u\|_2^2
    +\|A_\sigma^{-1/4}z\|_2^2
  \end{equation}
  is uniformly equivalent to the fixed norm on $E = H^{1/2}\times H^{-1/2}$,
  interpreting the negative entry by completion.
\end{lemma}

\begin{proof}
  We prove the needed fractional order comparison from resolvents. The
  resulting quarter-power norm comparison identifies the positive phase entry,
  and duality identifies the negative entry with uniform constants.

  \medskip
  \noindent\textit{Step 1.~}We compare the quarter-power norms and their
  domains.
  In the resolvent supremum below, increasing the quadratic form decreases
  every term being maximized over the common form domain, reversing the order.
  On the common form domain $H^1$, $B \le A_\sigma \le cB$ and all three
  operators are at least one. For positive self-adjoint $T,S$ with common form
  domain, $T \le S$ implies $(T+t)^{-1} \ge (S+t)^{-1}$ for $t > 0$. Indeed,
  \[
    \langle f,(T+t)^{-1}f\rangle
    =\sup_u\{2\operatorname{Re}\langle f,u\rangle
    -q_T(u)-t\|u\|_2^2\};
  \]
  completing the square in the closed form norm proves the identity, and
  comparison of the suprema proves the order. Consequently, $I-t(T+t)^{-1} \le
  I-t(S+t)^{-1}$. The scalar identity
  \[
    \sqrt x=\frac1\pi\int_0^\infty
    t^{-1/2}\frac{x}{x+t}\ud t
  \]
  follows by $t = xr^2$ and $\int_0^\infty(1+r^2)^{-1}\ud r = \pi/2$. Spectral
  integration and positive Tonelli, first over finite intervals and then by
  monotone convergence, give $T^{1/2} \le S^{1/2}$ as forms, including domain
  inclusions. Apply this to $B \le A_\sigma \le cB$ to obtain
  \eqref{E:compact-quarter-power}: the form of $A_\sigma^{1/2}$ is the squared
  norm of $A_\sigma^{1/4}$. The lower and upper comparisons give both domain
  inclusions, so they also give
  equality of the quarter-power domains with
  $H^{1/2}$.

  \medskip
  \noindent\textit{Step 2.~}We identify the negative phase norm by duality.
  Under the $L^2$ pairing the dual norm of $\|A_\sigma^{1/4}u\|_2$ at $z \in
  L^2$ is $\|A_\sigma^{-1/4}z\|_2$. Cauchy--Schwarz gives the upper bound, and a
  scalar multiple of $A_\sigma^{-1/2}z$ attains it. Comparison
  of these dual norms preserves the comparison with reversed inequalities. The
  constants remain uniform because the positive norm constants depend only on
  $c$. Completion
  and~\eqref{E:compact-quarter-power} identify this negative space with
  $H^{-1/2}$, uniformly in $\sigma$. Adding the squared norms proves the claim.
  This proves Lemma~\ref{L:compact-fractional-phase}.
\end{proof}

\begin{proposition}[The complementary wave on $E$]\label{P:compact-E-wave}
  Put
  \[
    e_\sigma^E=(p_\sigma,\lambda_\sigma p_\sigma),\qquad
    \ell_\sigma^E(u,z)
    =\frac{\langle p_\sigma,u\rangle+
    \langle p_\sigma,z\rangle/\lambda_\sigma}{2},
    \qquad \Pi_\sigma^E=e_\sigma^E\ell_\sigma^E.
  \]
  These are bounded real objects, continuous in their respective norms, with
  $\ell_\sigma^E(e_\sigma^E) = 1$, and
  \begin{align}
    P_\sigma^E(t)e_\sigma^E
    &=e^{\lambda_\sigma t}e_\sigma^E,\qquad
    \ell_\sigma^EP_\sigma^E(t)
    =e^{\lambda_\sigma t}\ell_\sigma^E,\label{E:compact-E-eigen}\\
    \|P_\sigma^E(t)(I-\Pi_\sigma^E)\|_{\mathcal{B}(E)}
    &\le C(1+t)e^{\mu_\sigma t},\qquad t\ge0.
    \label{E:compact-E-complement}
  \end{align}
  The constant is uniform over $\sigma \in I$.
\end{proposition}

\begin{proof}
  We estimate the spectral wave matrix in the uniformly equivalent phase norm,
  separating the complementary spatial spectrum from the two ground-state
  branches. We then identify this spectral matrix with the ordinary physical
  flow.

  \medskip
  \noindent\textit{Step 1.~}We write the spectral wave and its complementary
  spatial bound.
  Define real spectral functions $C_t(h),S_t(h)$ by $\cos(t\sqrt h),\sin(t\sqrt
  h)/\sqrt h$ for $h > 0$, by $1,t$ at zero, and by
  $\cosh(t\sqrt{-h}),\sinh(t\sqrt{-h})/\sqrt{-h}$ for $h < 0$. The spectral wave
  matrix is
  \begin{equation}\label{E:compact-spectral-wave}
    \begin{pmatrix}
      C_t(H_\sigma)&S_t(H_\sigma)\\
      -H_\sigma S_t(H_\sigma)&C_t(H_\sigma)
    \end{pmatrix}.
  \end{equation}
  Real potentials, their resolvents and spectral resolutions commute with
  conjugation. All these real functions and positive $A_\sigma$ powers preserve
  the real phases.

  The spatial projection $Q_\sigma^{\mathrm{sp}} = I-|p_\sigma\rangle\langle
  p_\sigma|$ commutes with $H_\sigma,A_\sigma$ and is contractive on both
  entries of~\eqref{E:compact-spectral-phase-norm}. On its support $h \ge E_2 =
  -\mu_\sigma^2$, the phase isometry gives scalar matrix entries
  \[
    C_t(h),\quad \sqrt{h+c}\,S_t(h),\quad
    -hS_t(h)/\sqrt{h+c},\quad C_t(h).
  \]
  The isometry multiplies the displacement by $(h+c)^{1/4}$ and the velocity
  by $(h+c)^{-1/4}$. Conjugating the two off-diagonal entries by these weights
  produces the displayed factors $\sqrt{h+c}$ and its reciprocal. Thus, an
  ordinary scalar matrix bound gives an operator bound in the spectral phase
  norm.
  For $h \ge 0$, use $|S_t(h)| \le \min(t,h^{-1/2})$ and split $h \le 1$ from $h
  \ge 1$; the lower off-diagonal entry is at most one. For $-\mu_\sigma^2 \le h
  < 0$, use $|C_t(h)| \le e^{\mu_\sigma t}$ and $|S_t(h)| \le t e^{\mu_\sigma
  t}$, with $|h| \le M$ and $1 \le h+c \le c$. Thus, the complementary spatial
  matrix has norm at most $C(1+t)e^{\mu_\sigma t}$, including $\mu_\sigma = 0$.

  \medskip
  \noindent\textit{Step 2.~}We retain the decaying ground branch in the phase
  complement.
  The two spatial ground branches are $(p_\sigma,\pm\lambda_\sigma p_\sigma)$.
  Besides $\Pi_\sigma^E$, their projections include the decaying branch
  \[
    (p_\sigma,-\lambda_\sigma p_\sigma)
    \frac{\langle p_\sigma,u\rangle
    -\langle p_\sigma,z\rangle/\lambda_\sigma}{2}.
  \]
  Their sum is the spatial ground projection on both phase entries. Both
  projection norms are uniformly bounded in the spectral phase norm because
  $\sqrt\varepsilon \le \lambda_\sigma \le \sqrt M$ and $1 \le
  c-\lambda_\sigma^2 \le c$. The decaying branch has factor $e^{-\lambda_\sigma
  t}$ and remains in $I-\Pi_\sigma^E$. The fractional norm equivalence now
  proves the stated complement estimate for the spectral matrix.

  \medskip
  \noindent\textit{Step 3.~}We identify the spectral matrix with the ordinary
  flow.
  We identify this matrix with the physical flow. The same scalar bounds at
  graph levels $(A_\sigma u,A_\sigma^{1/2}z)$ give continuous solutions in
  $H^2\times H^1$. These graph norms are uniformly equivalent to the Sobolev
  norms: the $H^1$ case was proved above, and $A_\sigma-B = M-V_\sigma$ is
  bounded by $M$ while $A_\sigma,B \ge 1$. Spectral cutoffs give $u_t = z$ in
  $H^1$ and $z_t = -H_\sigma u$ in $L^2$, with continuous derivatives there.
  The
  free Fourier group is used at these graph levels so that the
  interaction-picture derivative is justified in $H^1\times L^2$. This
  derivative is not an operator-norm derivative on the whole phase. The
  free Fourier group on $H^1\times L^2$ has, after its Bessel isometry, the same
  bounded weighted entries and is strongly continuous by dominated Fourier
  convergence. Its generator is $(z,\Delta u)$ on $H^2\times H^1$, as the
  Fourier multipliers show. For compact smooth initial data, differentiating the
  free interaction picture at these graph levels gives Duhamel with force
  $(0,V_\sigma u)$. The identity also holds in $E$. Multiplication by $V_\sigma$
  maps $H^{1/2}$ into $L^2$ and then $H^{-1/2}$ with bound $M$, so Picard
  uniqueness identifies the spectral solution with the ordinary $E$ flow on
  these data. Density within $E$ and finite-horizon boundedness extend equality
  to all $E$. No density of $E$ in $F$ is used. Finally, the ground branches
  give~\eqref{E:compact-E-eigen}; their pairings are continuous because
  $p_\sigma \in H^2$. Continuity and uniform bounds of the displayed vectors and
  duals follow from Proposition~\ref{P:compact-ground-family}. This proves
  Proposition~\ref{P:compact-E-wave}.
\end{proof}

\subsection{The growing projection and its complementary wave evolution}

For the constant coordinates write $C(t)(c_1,c_2) = (c_1+tc_2,c_2)$ and
$D_\sigma(c_1,c_2) = (0,V_\sigma c_1) \in E$.
The data $e_0 = (0,1,0)$ lie entirely in these constant coordinates, outside
$iE$. So we follow the flow on all of $F$, where $D_\sigma$ couples the
constants to $E$.

\begin{proposition}[The growing projection on the wave phase space]
  \label{P:compact-F-projection}
  The physical flow satisfies
  \begin{equation}\label{E:compact-F-block}
    P_\sigma(t)(e,c)
    =\left(P_\sigma^E(t)e+
    \int_0^tP_\sigma^E(t-r)D_\sigma C(r)c\ud r,\ C(t)c\right).
  \end{equation}
  Set
  \begin{align}
    f_\sigma&=\frac{\int V_\sigma p_\sigma}{2\lambda_\sigma}>0,
    &\ell_\sigma^C(c)&=
    f_\sigma\left(\frac{c_1}{\lambda_\sigma}
    +\frac{c_2}{\lambda_\sigma^2}\right),
    \label{E:compact-constant-dual}\\
    e_\sigma&=ie_\sigma^E,
    &\ell_\sigma(e,c)&=\ell_\sigma^E(e)+\ell_\sigma^C(c),\qquad
    \Pi_\sigma=e_\sigma\ell_\sigma,\quad Q_\sigma=I-\Pi_\sigma.
    \label{E:compact-F-projection}
  \end{align}
  These vectors, duals and rank-one projections are norm-continuous and
  uniformly bounded, and
  \begin{equation}\label{E:compact-F-eigen}
    \ell_\sigma(e_\sigma)=1,\qquad
    P_\sigma(t)\Pi_\sigma=\Pi_\sigma P_\sigma(t)
    =e^{\lambda_\sigma t}\Pi_\sigma.
  \end{equation}
  For every fixed nonnegative smooth compact $\phi$ of integral one,
  \begin{equation}\label{E:compact-positive-margins}
    \inf_I\ell_\sigma(e_0)>0,\qquad
    \inf_I O_\phi(e_\sigma)>0.
  \end{equation}
\end{proposition}

\begin{proof}
  We derive the block formula by variation of constants in $E$. We then
  compute the constant-coordinate contribution to the growing dual, combine
  its finite-time contribution with an exponential tail, and obtain uniform
  positive margins by compactness.

  \medskip
  \noindent\textit{Step 1.~}We derive the block formula for the full phase.
  The constant coordinates obey $C(t)$ because the potential force has range in
  $iE$. To justify variation of constants on the $E$ component, let $g$ be a
  continuous $E$ force and substitute the ordinary free-Duhamel equation for
  $P_\sigma^E(t-r)g(r)$ in its time integral. With $B_\sigma^E(u,z) =
  (0,V_\sigma u)$, strong $E$-valued Fubini gives
  \begin{align*}
    &\int_0^t\int_r^t S(t-s)B_\sigma^E
    P_\sigma^E(s-r)g(r)\ud s\ud r\\
    &\qquad=\int_0^tS(t-s)B_\sigma^E
    \left[\int_0^sP_\sigma^E(s-r)g(r)\ud r\right]\ud s.
  \end{align*}
  The integration region is the triangle $0 \le r \le s \le t$. Reversing the
  order changes the outer variable from $r$ to $s$, with $0 \le r \le s$ in
  the inner integral. This is exactly the forcing term in the free-Duhamel
  equation for the proposed $E$ component.
  All norms are bounded on the finite triangle, so this exchange is absolute in
  $E$. Thus, $Y(t) = P_\sigma^E(t)e+\int_0^tP_\sigma^E(t-r)g(r)\ud r$ solves
  free Duhamel with force $B_\sigma^EY+g$. Take $g(r) = D_\sigma C(r)c$ and
  adjoin the constant block. This is the same $F$ Duhamel equation as the
  ordinary flow; uniqueness proves~\eqref{E:compact-F-block}. Its finite-horizon
  bounds are uniform since the family has bounded $L^\infty$ and $L^2$ norms on
  one compact support.

  \medskip
  \noindent\textit{Step 2.~}We compute the growing dual on the constant
  coordinates.
  The integral defining $f_\sigma$ is finite by $L^2$ Cauchy--Schwarz and
  positive because $V_\sigma$ is nonzero by the common trial and $p_\sigma > 0$
  a.e. Its coefficients are continuous and uniformly bounded by the $H^2$
  continuity and positive lower bound for the rate; no $p_\sigma \in L^1$
  assumption is used. Directly,
  \begin{align*}
    \ell_\sigma^ED_\sigma C(r)c&=f_\sigma(c_1+rc_2),\\
    \ell_\sigma^C(c)
    &=\int_0^\infty e^{-\lambda_\sigma r}
    f_\sigma(c_1+rc_2)\ud r,\\
    \ell_\sigma^C(C(t)c)
    &=\int_t^\infty e^{\lambda_\sigma(t-r)}
    f_\sigma(c_1+rc_2)\ud r.
  \end{align*}
  The first integral uses $\int_0^\infty e^{-\lambda_\sigma r}\ud r =
  1/\lambda_\sigma$ and $\int_0^\infty r e^{-\lambda_\sigma r}\ud r =
  1/\lambda_\sigma^2$. The last identity follows by translating the
  integration variable by $t$ and using $C(t)c = (c_1+tc_2,c_2)$.
  These integrals converge absolutely. Applying $\ell_\sigma$
  to~\eqref{E:compact-F-block}, the finite integral and this tail combine to
  give $\ell_\sigma P_\sigma(t) = e^{\lambda_\sigma t}\ell_\sigma$. The ground
  mode also gives $P_\sigma(t)e_\sigma = e^{\lambda_\sigma t}e_\sigma$ and
  $\ell_\sigma(e_\sigma) = 1$, proving~\eqref{E:compact-F-eigen}. Bounded
  bilinear multiplication and the preceding continuous coefficients prove norm
  continuity of all objects. Finally, $\ell_\sigma(e_0) =
  f_\sigma/\lambda_\sigma > 0$ and $O_\phi(e_\sigma) = \int\phi p_\sigma > 0$.
  Nonnegativity and positive integral of $\phi$ give a set of positive measure
  on which $\phi > 0$. Since $p_\sigma > 0$ almost everywhere, its pairing
  with $\phi$ is strictly positive.
  Both are continuous on compact $I$, proving
  \eqref{E:compact-positive-margins}. Positivity is in integrals, without a
  point trace or a one-dimensional ODE. This proves
  Proposition~\ref{P:compact-F-projection}.
\end{proof}

Proposition~\ref{P:compact-F-projection} isolates the growing part
$\Pi_\sigma$. It remains to bound the flow on the complement $Q_\sigma$. The
proof of Lemma~\ref{L:raw-block-indexed-cone} uses the normalized form of this
bound.

\begin{proposition}[Complementary wave evolution and parameter dependence]
  \label{P:compact-F-complement}
  Uniformly over $\sigma$ and $t \ge 0$,
  \begin{equation}\label{E:compact-F-complement}
    \|P_\sigma(t)Q_\sigma\|_{\mathcal{B}(F)}
    \le C(1+t)^3e^{\mu_\sigma t}.
  \end{equation}
  For every fixed $0 < \eta < d_0$,
  \begin{equation}\label{E:compact-normalized-complement}
    \|e^{-\lambda_\sigma t}P_\sigma(t)Q_\sigma\|
    \le C_\eta e^{-\eta t}.
  \end{equation}
  For every finite $T_0$,
  \begin{equation}\label{E:compact-parameter-flow}
    \sup_{0\le t\le T_0}\|P_\sigma(t)-P_\tau(t)\|_{\mathcal{B}(F)}
    \le C_{T_0}
    \big(\|V_\sigma-V_\tau\|_\infty+\|V_\sigma-V_\tau\|_2\big).
  \end{equation}
\end{proposition}

\begin{proof}
  We use the block formula to keep track of the growing response to the
  constant coordinates when subtracting the projection. We estimate the
  remaining integral and absorb its polynomial factor into the rate gap.
  Finally, we compare two potentials through their Duhamel equations.

  \medskip
  \noindent\textit{Step 1.~}We bound the complementary flow on the full phase.
  The growing forcing component equals
  \[
    e_\sigma^E\int_0^t e^{\lambda_\sigma(t-r)}
    f_\sigma(c_1+rc_2)\ud r
    =e_\sigma^E\bigl(e^{\lambda_\sigma t}\ell_\sigma^C(c)
    -\ell_\sigma^C(C(t)c)\bigr).
  \]
  Its first term cancels the constant-coordinate contribution of
  $e^{\lambda_\sigma t}\Pi_\sigma(e,c)$, leaving the negative term in the
  complement formula below.
  Write $Q_\sigma^E = I-\Pi_\sigma^E$. Split $D_\sigma$ into its growing and
  complementary parts in~\eqref{E:compact-F-block} and use the two tail
  integrals in the preceding proof. The exact $E$ component of
  $P_\sigma(t)Q_\sigma(e,c)$ is
  \begin{equation}\label{E:compact-complement-forcing}
    P_\sigma^E(t)Q_\sigma^Ee
    +\int_0^tP_\sigma^E(t-r)Q_\sigma^ED_\sigma C(r)c\ud r
    -e_\sigma^E\ell_\sigma^C(C(t)c),
  \end{equation}
  and the constant component is $C(t)c$. The final subtraction retains the
  growing response to constants in the projection. Since $\|Q_\sigma^ED_\sigma
  C(r)\| \le C(1+r)$, the integral in~\eqref{E:compact-complement-forcing} has
  norm at most
  \[
    Ce^{\mu_\sigma t}|c|
    \int_0^t(1+t-r)(1+r)\ud r
    =Ce^{\mu_\sigma t}(t+t^2+t^3/6)|c|.
  \]
  The integrand has one factor $1+t-r$ from the complementary $E$ flow and one
  factor $1+r$ from the free constant block. Its exponential is bounded by
  $e^{\mu_\sigma t}$ because $\mu_\sigma \ge 0$. Integrating the product of
  the two linear factors gives the displayed cubic polynomial.
  The first term is bounded by $C(1+t)e^{\mu_\sigma t}\|e\|_E$; the final term
  and constant component by $C(1+t)|c|$. Since $\mu_\sigma \ge 0$ this proves
  \eqref{E:compact-F-complement}. The boundedness of $(1+t)^3e^{-(d_0-\eta)t}$
  proves \eqref{E:compact-normalized-complement}. In particular, the
  constant-velocity polynomial is included in this estimate.

  \medskip
  \noindent\textit{Step 2.~}We prove finite-horizon dependence on the
  potential.
  Put $\mathcal{B}_\sigma = iB_F(V_\sigma)$ and $d_{\sigma,\tau} =
  \|V_\sigma-V_\tau\|_\infty+ \|V_\sigma-V_\tau\|_2$. The bound for $B_F$ gives
  $\|\mathcal{B}_\sigma-\mathcal{B}_\tau\| \le d_{\sigma,\tau}$. Subtract the
  two free-Duhamel equations on states. Their difference is the sum of the
  integral with $S_F(t-r)\mathcal{B}_\sigma(P_\sigma(r)-P_\tau(r))$ and that
  with $S_F(t-r)(\mathcal{B}_\sigma-\mathcal{B}_\tau)P_\tau(r)$. The uniform
  finite-horizon bounds control the second integral by a constant times
  $d_{\sigma,\tau}$. The first integral is bounded by a constant times the
  time integral of the unknown flow difference. Thus, the
  finite-horizon bounds and scalar Gronwall prove
  \eqref{E:compact-parameter-flow}. The operator norm is measurable:
  separability of $F$ expresses it as the supremum over a countable dense set of
  unit states, each with continuous image. Hence, Gronwall applies without a
  $\mathcal{B}(F)$-valued derivative or operator-norm continuity of the free
  group. This proves Proposition~\ref{P:compact-F-complement}.
\end{proof}

\subsection{The original nonnegative control and a common trial}

So far $V_\sigma$ was any continuous family
satisfying~\eqref{E:compact-common-trial}. We now produce such a family from a
nonnegative test $\Psi$ and the original covariance $\Gamma$, in two steps.

\begin{lemma}[The original mean and variance]\label{L:compact-raw-profile}
  Let $\Gamma \ne 0$ satisfy Hypothesis~\ref{H:homog}, let $a = \alpha_0$, and
  let real nonnegative nonzero $\Psi \in C_c^\infty((0,1)\times\R^d)$ be
  arbitrary, without a separability assumption. Define
  \begin{equation}\label{E:compact-raw-mean}
    h_\sigma(x)=\int_0^1|\sigma-u|^{-a}
    (\Gamma*\Psi(u,\cdot))(x)\ud u,\qquad
    \overline{q}(x)=\int_0^1\Psi(u,x)\ud u,\qquad Z=W(\Psi).
  \end{equation}
  Then $h$ is smooth in time and space on compact sets, and
  \begin{equation}\label{E:compact-raw-cost}
    0<v_\Psi:=\E Z^2
    =\int_0^1\int\Psi(\sigma,x)h_\sigma(x)\ud x\ud\sigma<\infty.
  \end{equation}
  Against every smooth compact positive-time test $F_0$,
  \[
    \E[W(\Psi)W(F_0)]=\int F_0(s,x)h_s(x)\ud x\ud s.
  \]
  Moreover, on $I$,
  \begin{equation}\label{E:compact-positive-mean-floor}
    h_\sigma\ge\Gamma*\overline{q},\qquad
    v_\Psi\ge Q(\overline{q})
    :=\int\Gamma(\ud z)\int\overline{q}(x)\overline{q}(x-z)\ud x>0.
  \end{equation}
  For every fixed mask, $V_\sigma = m_Rh_\sigma$ is a smooth compact family,
  Lipschitz in $\sigma$ in $L^\infty$ and $L^2$.
\end{lemma}

\begin{proof}
  We identify the original covariance mean and show its cost is finite.
  Positivity and homogeneity give a nonnegative, nonzero time-independent lower
  profile and a strictly positive covariance cost. We then differentiate the
  test, rather than the singular time kernel, to obtain regularity.

  \medskip
  \noindent\textit{Step 1.~}We identify the mean and bound its covariance
  cost.
  Choose $A \ge 1$ containing the spatial support of $\Psi$ in $B_A$. The
  original covariance identity and absolute Fubini give the stated mean pairing:
  the two compact tests restrict spatial differences to a compact set of finite
  $\Gamma$ mass, and the time singularity is locally integrable since $a < 1$.
  This identifies the linear covariance mean without evaluating noise at a
  point. The same calculation gives~\eqref{E:compact-raw-cost}. Its spatial
  integrand is bounded in integral by
  $\|\Psi\|_\infty^2|B_A|\Gamma(\overline{B}_{2A})$, and
  $\int_0^1\int_0^1|\sigma-u|^{-a}\ud u\ud\sigma = 2/((1-a)(2-a))$, proving
  finiteness.

  \medskip
  \noindent\textit{Step 2.~}We obtain strict positivity and the uniform mean
  floor.
  Nonnegative nonzero smooth $\Psi$ has positive integral, so $\overline{q} \ne
  0$. Nonzero $\Gamma$ charges some ball, and homogeneity gives $\Gamma(B_1) =
  r^{\alpha-d}\Gamma(B_r) > 0$. Thus, every origin-centered ball $B_r$, $r > 0$,
  has positive finite mass. The continuous compact autocorrelation
  $\int\overline{q}(x)\overline{q}(x-z)\ud x$ is positive at zero and bounded
  below on some ball. This proves $0 < Q(\overline{q}) < \infty$. Since
  $|\sigma-u|^{-a} \ge 1$ on the unit time square, positive Tonelli gives both
  inequalities in~\eqref{E:compact-positive-mean-floor}.

  \medskip
  \noindent\textit{Step 3.~}We prove smoothness and the fixed-mask Lipschitz
  bounds.
  Extend $\Psi$ by zero in time and write
  \[
    h_s(x)=\int_\R |r|^{-a}
    \int\Gamma(\ud z)\Psi(s-r,x-z)\ud r.
  \]
  On compact $s,x$ sets the relevant $r,z$ range over fixed compact sets;
  $|r|^{-a}$ is integrable there. Every derivative falls on the smooth test, not
  on the singular time kernel. Dominated differentiation gives smoothness and
  continuity of all derivatives on compacts. Multiplication by $m_R$ now gives
  the claimed compact family and its fixed-$R$ Lipschitz bounds. This proves
  Lemma~\ref{L:compact-raw-profile}.
\end{proof}

Lemma~\ref{L:compact-raw-profile} gives a smooth compact family for every
fixed mask, but not yet the common negative
trial~\eqref{E:compact-common-trial}.

\begin{proposition}[Binding with the original control and cost fixed]
  \label{P:compact-original-common-trial}
  For every $\Psi$ in Lemma~\ref{L:compact-raw-profile}, a sufficiently large
  finite mask $m_R$ gives $V_\sigma = m_Rh_\sigma$ a common negative
  trial~\eqref{E:compact-common-trial}. Thus, all preceding ground, projection,
  positive-margin and complementary-wave conclusions apply. The rate is
  Lipschitz on $I$ and has $\lambda_0 > 0$. The original test $\Psi$, original
  variable $Z$ and variance $v_\Psi$ are unchanged by this choice. The mask can
  additionally cover any fixed compact variational test supports.
\end{proposition}

\begin{proof}
  We dilate a fixed trial function and choose the mask to contain its support.
  Homogeneity gives a potential contribution of order $L^{-\alpha}$, which
  dominates the kinetic contribution of order $L^{-2}$. We then apply the
  ground-family estimates at this fixed mask.

  \medskip
  \noindent\textit{Step 1.~}We bound the potential contribution of the dilated
  trial.
  Choose a real unit $L^2$ function $\psi \in C_c^\infty(B_2)$ with $\psi \ge
  c_0 > 0$ on $B_1$, and put $\psi_L(x) = L^{-d/2}\psi(x/L)$, $D_\psi =
  \|\nabla\psi\|_2^2$, where $L \ge \max(1,2A)$. Choose the even mask with $R
  \ge 2L$, so $V_\sigma = h_\sigma$ on $\operatorname{supp}\psi_L$. This is the
  allowed real mask $0 \le m_R \le 1$, equal to one on $B_R$ and with support in
  $\overline{B}_{2R}$. Positive Tonelli
  and~\eqref{E:compact-positive-mean-floor} give
  \begin{align}
    \int V_\sigma\psi_L^2
    &\ge c_0^2L^{-d}\int_{B_L}(\Gamma*\overline{q})(x)\ud x\nonumber\\
    &=c_0^2L^{-d}\int\overline{q}(y)\Gamma(B_L-y)\ud y\nonumber\\
    &\ge c_0^2L^{-d}\Gamma(B_{L/2})\int\overline{q}
    =C_\Gamma L^{-\alpha},\label{E:compact-trial-potential}\\
    C_\Gamma&=c_0^22^{\alpha-d}\Gamma(B_1)\int\overline{q}>0.\nonumber
  \end{align}
  On $B_L$, the dilated trial has squared value at least $c_0^2L^{-d}$.

  Homogeneity gives $\Gamma(B_{L/2})=(L/2)^{d-\alpha}\Gamma(B_1)$, yielding
  the displayed constant.
  Here $B_{L/2} \subset B_L-y$ for $|y| \le A$ and $L \ge 2A$. All covariance
  masses in the argument are compact and finite. The kinetic energy is $D_\psi
  L^{-2}$. Since $\alpha < 2$, choose $L$ large enough that $\varepsilon_R: =
  C_\Gamma L^{-\alpha}-D_\psi L^{-2} > 0$. Then
  \[
    \|\nabla\psi_L\|_2^2-\int V_\sigma\psi_L^2
    \le-\varepsilon_R\quad\hbox{for every }\sigma\in I.
  \]

  \medskip
  \noindent\textit{Step 2.~}We fix a negative trial and obtain the rate
  regularity.
  Fix this finite $R$ before invoking Proposition~\ref{P:compact-ground-family}.
  Its level bound and the mean Lipschitz estimate give
  \[
    |\lambda_\sigma-\lambda_\tau|
    \le\frac{\|V_\sigma-V_\tau\|_\infty}{2\sqrt{\varepsilon_R}}
    \le C_R|\sigma-\tau|.
  \]
  The lowest eigenvalue remains at most $-\varepsilon_R$ for every parameter.
  Factoring the difference of its negative square roots gives a denominator
  $\lambda_\sigma+\lambda_\tau \ge 2\sqrt{\varepsilon_R}$. The level
  comparison then gives the displayed Lipschitz estimate.
  The positive mean integrates the entire time test, so $\lambda_0 > 0$ even
  though $\Psi(0,\cdot) = 0$. Only the trial and mask have grown; the original
  test, its variable and~\eqref{E:compact-raw-cost} were fixed first. Enlarging
  $R$ to contain additional fixed compact test supports preserves the trial and
  gives equality of masked and unmasked potential energies on those supports.
  This includes singular invariant-direction covariances and the one-dimensional
  white case. It is not an all-strength binding assertion for one fixed compact
  potential. Masks are fixed before source or time asymptotics; a later supremum
  over controls must follow their individual bounds. This proves
  Proposition~\ref{P:compact-original-common-trial}.
\end{proof}

\paragraph{What we now have.}
The Young evolution and its physical-wave identification
(Proposition~\ref{P:constant-gaussian-flow} and
Lemma~\ref{L:constant-source-solution}), together with positive ground states
and a uniform gap for compact potential families satisfying the common
negative trial condition (Proposition~\ref{P:compact-ground-family}). The next
chapter uses both for the lower bound.

\chapter{Scaling and the point lower bound}\label{S:book-lower}

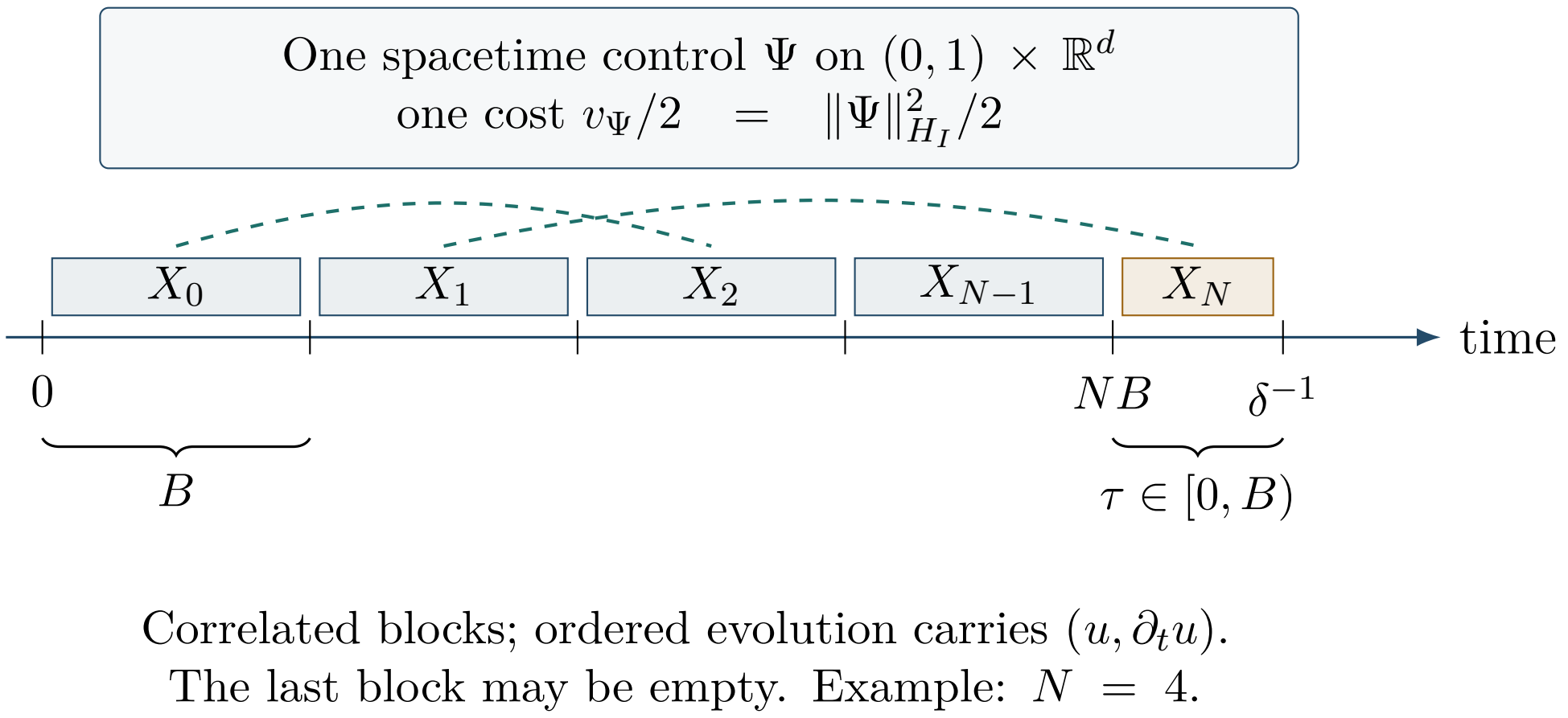
\begin{figure}
  \centering
  \begin{tikzpicture}[x=1cm,y=1cm]
    \node[ham/box,text width=9.5cm] at (5.4,2.05)
    {One spacetime control $\Psi$ on $(0,1)\times\R^d$\\
    one cost $v_\Psi/2=\|\Psi\|_{H_I}^2/2$};
    \draw[ham/arrow] (-.3,0)--(11.5,0) node[right] {time};
    \foreach \hamPointX in {0,2.2,4.4,6.6,8.8} {\draw (\hamPointX,-.14)--(\hamPointX,.14);}
    \draw (10.2,-.14)--(10.2,.14);
    \foreach \hamPointX/\hamBlockIndex in {1.1/0,3.3/1,5.5/2,7.7/{N-1}} {
    \draw[fill=hamblue!9,draw=hamblue] (\hamPointX-1.02,.18) rectangle (\hamPointX+1.02,.65);
    \node[font=\small] at (\hamPointX,.42) {$X_{\hamBlockIndex}$};
    }
    \draw[fill=hamgold!12,draw=hamgold] (8.88,.18) rectangle (10.12,.65);
    \node[font=\small] at (9.5,.42) {$X_N$};
    \node[font=\small,anchor=north] at (0,-.18) {$0$};
    \node[font=\small,anchor=north] at (8.8,-.18) {$NB$};
    \node[font=\small,anchor=north] at (10.2,-.18) {$\delta^{-1}$};
    \draw[ham/brace] (2.2,-.83)--(0,-.83);
    \node[font=\small,anchor=north] at (1.1,-.98) {$B$};
    \draw[ham/brace] (10.2,-.83)--(8.8,-.83);
    \node[font=\small,anchor=north] at (9.5,-.98) {$\tau\in[0,B)$};
    \draw[ham/pair] (1.1,.75) to[bend left=16] (5.5,.75);
    \draw[ham/pair] (3.3,.75) to[bend left=12] (9.5,.75);
    \node[ham/note,text width=10cm] at (5.25,-2.65)
    {Correlated blocks; ordered evolution carries $(u,\partial_tu)$.\\
    The last block may be empty. Example: $N=4$.};
  \end{tikzpicture}
  \caption[Correlated blocks and one Gaussian cost]{The true partition in fast
    time. Here $N=\lfloor1/(\delta B)\rfloor$ and $\tau=\delta^{-1}-NB$. All
    blocks come from the same noise and a single control on the whole time
    interval. The change-of-law density has cost $v_\Psi/(2\delta)$; the
    displayed $v_\Psi/2$ is its normalized cost. The proof does not multiply
    independent block probabilities or add separate Gaussian costs. See
    \ProofRef{Lemma}{L:raw-block-inputs} and
    \ProofRef{Proposition}{P:tube-actual-law}.}
  \label{Fig:global-control-blocks}
\end{figure}

To bound a moment from below, push the Gaussian noise toward a deterministic
potential that makes the wave grow, and pay for the push. The noise is
correlated in time, so the cost of one control is charged once, in the
original covariance space. It is not a sum over independent time blocks
(Figure~\ref{Fig:global-control-blocks}).

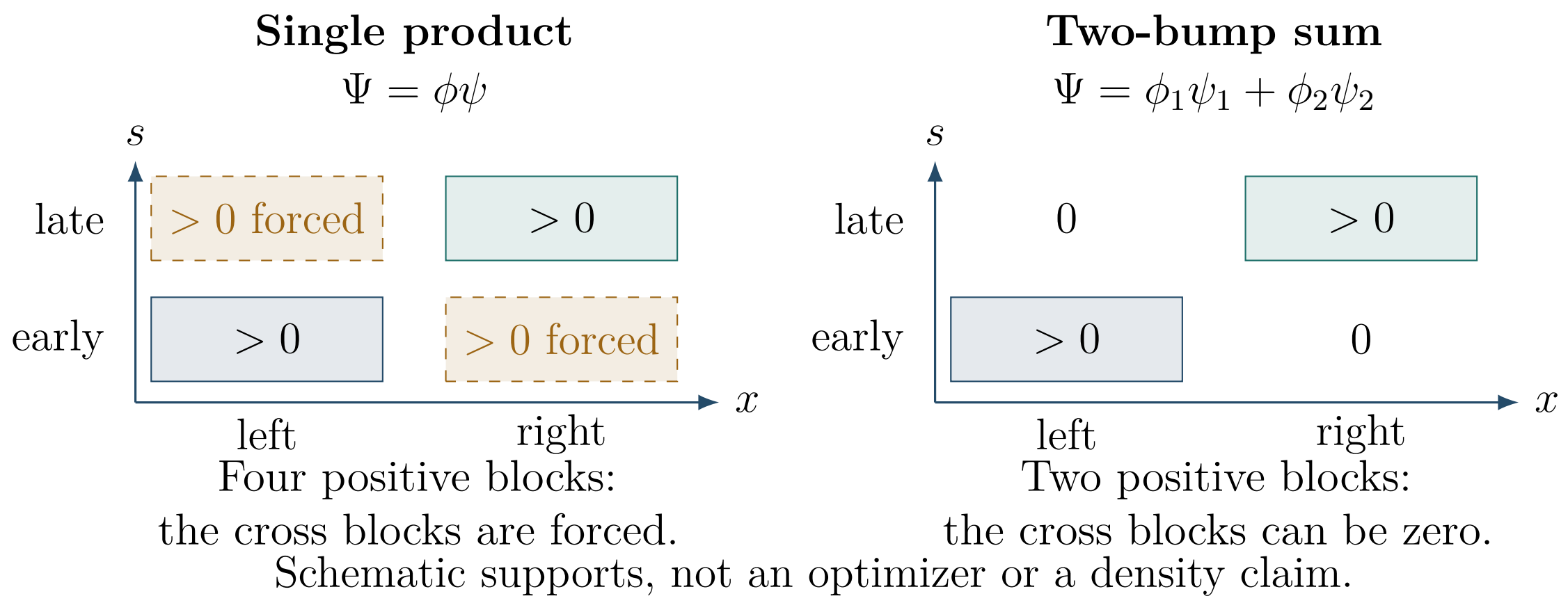
\begin{figure}
  \centering
  \begin{tikzpicture}[x=1cm,y=1cm,font=\normalsize]
    % Identical panel widths, baselines, box sizes and type sizes.
    \foreach \shift in {0,7.6} {
    \begin{scope}[xshift=\shift cm]
      \draw[ham/arrow] (.85,0)--(6.4,0) node[right] {$x$};
      \draw[ham/arrow] (.85,0)--(.85,2.3) node[above] {$s$};
      \node[anchor=east] at (.7,.6) {early};
      \node[anchor=east] at (.7,1.75) {late};
      \node at (2.1,-.3) {left};
      \node at (4.9,-.3) {right};
      \filldraw[draw=hamblue,fill=hamblue!12]
      (1,.2) rectangle (3.2,1);
      \filldraw[draw=hamteal,fill=hamteal!12]
      (3.8,1.35) rectangle (6,2.15);
      \node at (2.1,.6) {$>0$};
      \node at (4.9,1.75) {$>0$};
    \end{scope}
    }
    \node[font=\normalsize\bfseries] at (3.5,3.5) {Single product};
    \node[font=\normalsize\bfseries] at (11.1,3.5) {Two-bump sum};
    \node at (3.5,2.95) {$\Psi=\phi\psi$};
    \node at (11.1,2.95) {$\Psi=\phi_1\psi_1+\phi_2\psi_2$};
    \foreach \x/\y in {4.9/.6,2.1/1.75} {
    \filldraw[draw=hamgold,dashed,fill=hamgold!12]
    (\x-1.1,\y-.4) rectangle (\x+1.1,\y+.4);
    \node[text=hamgold] at (\x,\y) {$>0$ forced};
    }
    \node at (12.5,.6) {$0$};
    \node at (9.7,1.75) {$0$};
    \node[align=center,text width=6.6cm] at (3.5,-.95)
    {Four positive blocks:\\ the cross blocks are forced.};
    \node[align=center,text width=6.6cm] at (11.1,-.95)
    {Two positive blocks:\\ the cross blocks can be zero.};
    \node at (7.3,-1.65)
    {Schematic supports, not an optimizer or a density claim.};
  \end{tikzpicture}
  \caption[Controls that need not factor]{The full control class allows spatial
    profiles to change with time, as in this sum of two smooth nonnegative product bumps that no single
    time--space product can represent. The full class is infinite-dimensional;
    see \ProofRef{Theorem}{T:intrinsic-second-moment}.}
  \label{Fig:nonseparable-controls}
\end{figure}

Section~\ref{S:exact-canonical-scaling} turns long time into small
noise by an equality in law. Section~\ref{S:actual-raw-block-lower}
proves, block by block in time, a pointwise lower bound for each
nonnegative smooth control (Proposition~\ref{P:raw-block-original-lower}).
Section~\ref{S:unmasked-variational} removes the spatial mask and optimizes
over controls, including those that do not factor into time and space
(Figure~\ref{Fig:nonseparable-controls}).

We start from the canonical solution of Chapter~\ref{S:G0} and the
phase and ground states of Chapter~\ref{S:book-wave-evolution}. On
a first reading, keep Lemma~\ref{L:exact-canonical-scaling} and
Proposition~\ref{P:variational-positive-scores}, and let the proofs wait. The
result is a variational lower value; only Chapter~\ref{S:book-constant}
proves that it equals the upper value. Episode~\VideoRef{6}{VWemsluFZiA}
of the video series covers this chapter and the next.

\section{Exact scaling of the canonical solution}
\label{S:exact-canonical-scaling}

We first rewrite the second moment at a large time $T$ as a second moment at
time one, for an equation with the small parameter $\delta = T^{-\kappa}$. The
proof of Proposition~\ref{P:raw-block-point-comparison} uses this for the final
equality in~\eqref{E:raw-block-original-comparison}.

\LeanMarker{L:exact-canonical-scaling}
\begin{lemma}[Canonical law scaling]\label{L:exact-canonical-scaling}
  Under Hypothesis~\ref{H:homog}, put $a = \alpha_0$, $\kappa =
  (4-\alpha-a)/(3-\alpha)$, and, for $T > 0$,
  \[
    \delta=T^{-\kappa},\qquad
    r_T=T\delta=T^{(a-1)/(3-\alpha)}.
  \]
  The law of $u(Tt,r_Tx)$ is the canonical expansion solution $V_\delta(t,x)$
  with data $(1,0)$ of
  \begin{equation}\label{E:exact-fast-equation}
    \delta^2\partial_t^2V_\delta
    =\Delta V_\delta+\sqrt\delta\,\W V_\delta,
  \end{equation}
  where $\W$ has the original covariance. In particular,
  \begin{equation}\label{E:exact-moment-scaling}
    T^{-\kappa}\log\E|u(T,0)|^2
    =\delta\log\E|V_\delta(1,0)|^2.
  \end{equation}
  The assertion is equality in law; it assumes no pathwise uniqueness for a full
  rough solution.
\end{lemma}

\begin{proof}
  We first scale the regularized noise and the Green measure. We check that
  the coefficient at each ordered vertex is the one in the fast equation. We
  then pass the finite-order law identities to the canonical series at fixed
  dilation.

  \medskip
  \noindent\textit{Step 1.~}We scale the source covariance and the Green
  measure.
  For the scalar powers below, substituting $r_T=T\delta$ and
  $(3-\alpha)\kappa=4-\alpha-a$ with $\delta=T^{-\kappa}$ gives
  $T^{1-a/2}r_T^{1-\alpha/2}=\delta^{-1/2}$; multiplying by $\delta$ gives the
  companion identity used for the noise coefficient $\sqrt\delta$.
  At each vertex the time substitution contributes the Jacobian $T$, the
  Green-measure scaling contributes $r_T$, and the single Gaussian source
  factor contributes the square root $T^{-a/2}r_T^{-\alpha/2}$ of its
  covariance amplitude.
  At positive widths, time and space dilation of the Gaussian source gives
  covariance amplitude $T^{-a}r_T^{-\alpha}$. This follows distributionally from
  the homogeneous Radon covariance, including spatial white noise. Its source
  cutoffs transform by $K_\eta(Ts) = T^{-a}K_{\eta/T}(s)$ and
  $\epsilon\mapsto\epsilon/r_T^2$; the joint maximum-time and summed-space width
  laws are preserved at every vertex. The Green measure has the scaling
  \[
    G_d(T\tau,r_T\ud z)
    =r_TG_d(\tau/\delta,\ud z).
  \]
  It follows from its degree $1-d$ distributional scaling, also for the sphere
  measure in dimension three. Each terminal Duhamel vertex thus contributes
  $Tr_T\,T^{-a/2}r_T^{-\alpha/2}$. Substitution of the displayed positive real
  powers gives
  \[
    T^{-a/2}r_T^{2-\alpha/2}=\sqrt\delta,\qquad
    T^{1-a/2}r_T^{1-\alpha/2}=\delta^{-1/2}.
  \]

  \medskip
  \noindent\textit{Step 2.~}We match the scaled vertex with the fast-equation
  vertex.
  The Green kernel for the forcing in~\eqref{E:exact-fast-equation} is
  $\delta^{-1}G_d((t-s)/\delta)$; after multiplication by $\sqrt\delta$ it gives
  exactly that vertex factor. Therefore, every finite list of fixed-order
  coefficients has the same joint law under the dilation, including all
  self-pairings, unequal orders and empty chains. At fixed $T,r_T$, remove Green
  smoothing first in dimension three, with its width divided by $r_T^2$ on the
  scaled side. Lemma~\ref{L:joint-fixed-order} and the absolute $L^2$
  summability in Lemma~\ref{L:vertex-crude} pass the identities to the canonical
  series and its second moment. No large-time and source limit are interchanged.
  Finally, evaluate the law
  identity at $t=1,x=0$ and use $T^{-\kappa}=\delta$ to obtain the logarithmic
  moment identity.
  This proves the law assertion and~\eqref{E:exact-moment-scaling}. This proves
  Lemma~\ref{L:exact-canonical-scaling}.
\end{proof}

\section{Time blocks for the original noise and the pointwise lower bound}
\label{S:actual-raw-block-lower}

Fix a nonzero covariance from Hypothesis~\ref{H:homog} and any real nonnegative
nonzero $\Psi \in C_c^\infty((0,1)\times\R^d)$. Retain the \emph{original} $Z =
W(\Psi)$, $v = v_\Psi > 0$ and mean $h_s$ of Lemma~\ref{L:compact-raw-profile}.
Choose one finite even mask $m_R \in [0,1]$, equal to one on $B_R$ and supported
in $\overline{B}_{2R}$, as in Proposition~\ref{P:compact-original-common-trial}.
Then $V_s = m_Rh_s$ has the ground family, rate $\lambda_s$ and full $F$
projection of Section~\ref{S:compact-mean-family}. Set
\[
  I_R=\int_0^1\lambda_s\ud s,\qquad
  Q_2=\liminf_{T\to\infty}T^{-\kappa}\log\E|u(T,0)|^2.
\]
We will prove $Q_2 \ge 2I_R-v/2$, for this same control and covariance. Fix a
real nonnegative smooth compact $\phi$ of integral one, the bounded observable
$O_\phi$, and the data $e_0 = (0,1,0) \in F$. Both displacement and velocity are
retained. This is the same physical profile $m_R = \zeta_R$ used in the operator
path and compact mean; no mask is changed in a covariance identity. The
semiscaled $W$ below is the given-law noise. The comparison to the unscaled
solution uses equality in law, not samplewise identities on one probability
space for every $T$.

\subsection{The common original noise and its exact Gaussian shift}

Choose an integer $B \ge 1$, later fixed by the cone estimate, and take $0 <
\delta < \min(1,1/B)$. Define
\begin{equation}\label{E:raw-block-clocks}
  N=\left\lfloor\frac1{\delta B}\right\rfloor\ge1,\quad
  s_m=m\delta B\ (0\le m\le N),\quad
  \tau=\delta^{-1}-NB\in[0,B),\quad
  g_\delta=\delta^{(1-a)/2}.
\end{equation}
For a smooth compact positive-fast-time test $f$, extended by zero, put
\begin{equation}\label{E:raw-block-pullback}
  T_{\delta,m}f(s,x)
  =\delta^{a/2-1}f((s-s_m)/\delta,x),\qquad
  W_{\delta,m}(f)=W(T_{\delta,m}f).
\end{equation}
Write $W_{\mathrm{fast}} = W_{\delta,0}$. Use the same admissible
$\gamma,\theta,\beta$ as~\eqref{E:driver-parameters}. Write $A_R(W_{\delta,m})$
for the $F$ integrated operator obtained by the bounded extension
\eqref{E:constant-phase-driver}, on the local horizon $[0,B]$, and set
\begin{equation}\label{E:raw-block-drivers}
  X_m=g_\delta A_R(W_{\delta,m}),\qquad
  H_m(u)=\int_0^uS_F(-r)iB_F(V_{s_m+\delta r})S_F(r)\ud r.
\end{equation}
The convolution formula for $h$ extends smoothly to $s \in [0,2]$, covering
these means. Ground templates will only use $s_m \in [0,1]$.

\LeanMarker{L:raw-block-inputs}
\begin{lemma}[Coupled Gaussian block inputs]\label{L:raw-block-inputs}
  Each $T_{\delta,m}$ is an isometry of the original test completion, respecting
  covariance-null tests. All $W_{\delta,m}$ and $Z$ are on the same given noise
  space, without an independence assertion. Let $\mathcal{K}(F)$ have the
  operator norm and let $\mathscr{B}_0 = h_0^\theta([0,B];\mathcal{K}(F))$ be
  the $C^\theta$ closure of zero-initial polygonal compact-operator paths. It is
  a separable Banach space. In $\mathscr{B}_\delta = \mathscr{B}_0^{N+1}$ with
  maximum norm, $(X,Z)$, $X = (X_m)_{m=0}^N$, is jointly centered Gaussian and
  \begin{equation}\label{E:raw-block-covariance}
    b:=\E[XZ]=\sqrt\delta\,H,\qquad H=(H_m)_{m=0}^N.
  \end{equation}
  Moreover,
  \begin{equation}\label{E:raw-block-moments}
    \|X_m\|_{L^p\mathscr{B}_0}
    \le C_{R,B}\sqrt p\,g_\delta,\qquad p\ge1.
  \end{equation}
\end{lemma}

\begin{proof}
  We first check the time pullbacks in the original covariance norm. We then
  place the operator paths in a separable space and compute their joint
  covariance with the single original Gaussian variable.

  \medskip
  \noindent\textit{Step 1.~}We check the common-space pullback identities.
  Substituting $s = s_m+\delta u$, $t = s_m+\delta r$ in the original covariance
  gives the factor $\delta^{a-2}\delta^2\delta^{-a} = 1$. Thus, the maps extend
  isometrically through the null quotient. Direct test substitution gives
  $W_{\delta,m}(f) = W_{\mathrm{fast}}(f(\cdot-mB,\cdot))$. These are explicit
  common-space identities.

  \medskip
  \noindent\textit{Step 2.~}We identify a separable path space for the
  Gaussian variables.
  The pullback isometry gives $W_{\delta,m}$ the same Gaussian law as $W$, so
  Proposition~\ref{P:constant-phase-driver} supplies the path convergence and
  Gaussian norm bound on $[0,B]$, with the same $C_{R,B}$ for every
  $m,\delta$.
  Finite-rank rational matrices in a countable real Hilbert basis are dense in
  $\mathcal{K}(F)$: for compact $T$, $P_nTP_n \to T$ in norm by compactness of
  $T,T^*$, and finite entries admit rational approximation. Polygonal paths with
  rational breakpoints and such matrices give separability of $\mathscr{B}_0$.
  To approximate breakpoints, use uniform convergence, uniform Lipschitz bounds
  and $[f]_\theta \le (2\|f\|_\infty)^{1-\theta}\operatorname{Lip}(f)^\theta$.
  The same inequality puts every zero-initial Lipschitz compact-valued path in
  $\mathscr{B}_0$ by polygonal approximation. The finite original coefficient
  paths $iA_jL_\psi$ have these properties. Their Gaussian sums converge in
  $L^pC^\theta$ and hence in the closed space $\mathscr{B}_0$. Their laws and
  the parent Gaussian norm bound give~\eqref{E:raw-block-moments}. Every finite
  joint coefficient list with $Z$ is Gaussian on the same $W$. Scalar continuous
  functionals pass through the $L^2$ limits, preserving Gaussianity by
  convergence of variances and characteristic functions. Thus, $(X,Z)$ is the
  asserted Banach Gaussian variable. Cauchy--Schwarz gives $\E\|XZ\| \le
  \|X\|_2\sqrt v < \infty$, so its Bochner covariance exists.

  \medskip
  \noindent\textit{Step 3.~}We compute the covariance direction of the
  operator paths.
  The time substitution in the mean pairing contributes a Jacobian $\delta$,
  changing the pullback factor $\delta^{a/2-1}$ to $\delta^{a/2}$.
  The original covariance and the test pullback give
  \begin{equation}\label{E:raw-block-mean-jacobian}
    \E[W_{\delta,m}(f)Z]
    =\delta^{a/2}\int f(u,x)h_{s_m+\delta u}(x)\ud x\ud u.
  \end{equation}
  Absolute Fubini uses compact supports, local Radon mass and $a < 1$. This
  identity extends to the explicit matrix-element tests for the operator path on
  $E$: their time indicator is approximated by smooth bounded tests on one
  compact support, converging in the original covariance norm by dominated
  convergence against $|u-r|^{-a}\Gamma$; their mean pairings converge as well.
  Insert the matrix-element test $\zeta_R1_{(0,u)}v_zG_w$
  in~\eqref{E:raw-block-mean-jacobian}: the mask multiplying the mean is
  $V_{s_m+\delta r}$, and~\eqref{E:driver-covariance-row} identifies the
  resulting integral with the row of $H_m$.
  For constant inputs, the free displacement from $(c_1\psi,c_2\psi)$ equals
  $c_1+rc_2$ on $B_{2R}$ for $r\le B$, since $\psi=1$ on the enlarged ball;
  multiplying by $V_{s_m+\delta r}$ therefore gives the same constant-column
  forcing in $H_m$.
  The Sobolev adjoint is the one in~\eqref{E:driver-covariance-row}. Hence, the
  path covariance is $g_\delta\delta^{a/2}H_m = \sqrt\delta\,H_m$ on real
  Schwartz $E$ pairs. For constant inputs the two columns are the same $L_\psi$
  inputs $(\psi,0)$ and $(0,\psi)$, and total finite propagation identifies
  their deterministic response. The span of $E$ and both constant coordinates is
  dense in $F$; $E$ alone is not. Boundedness, rational times and path
  continuity prove~\eqref{E:raw-block-covariance} on all $F$. In particular, $H
  = b/\sqrt\delta \in \mathscr{B}_\delta$. This proves
  Lemma~\ref{L:raw-block-inputs}.
\end{proof}

The next identity is the covariance-direction Cameron--Martin formula;
see Janson~\cite{janson:97:gaussian}*{Theorem~14.1(v), pp.~216--217}.
We prove it on the norm-Borel space used here to verify the measurable
realization and its deterministic translate.

\begin{proposition}[Change of law on the whole norm-Borel space]
  \label{P:raw-block-gaussian-law}
  For every nonnegative norm-Borel $j$ on all of $\mathscr{B}_\delta$,
  \begin{equation}\label{E:raw-block-gaussian-law}
    \E j(X)=
    \E\left[e^{-Z/\sqrt\delta-v/(2\delta)}j(X+H)\right].
  \end{equation}
  The equality is valid also as an extended nonnegative integral.
\end{proposition}

\begin{proof}
  We prove the change of measure first for finitely many continuous
  coordinates, allowing singular covariance matrices. We then show that these
  coordinates generate the norm-Borel sets and extend the identity to
  nonnegative functions.

  \medskip
  \noindent\textit{Step 1.~}We shift the finite Gaussian coordinate vectors.
  For real $s$, the density $e^{-sZ-s^2v/2}$ defines an equivalent probability
  measure $Q_s$ by completing the square in the $N(0,v)$ density. For any finite
  list $\ell$ of continuous real coordinates on $\mathscr{B}_\delta$, the
  Gaussian vector $(\ell(X),Z)$ is equal in law to $(TG,z\cdot G)$ for a finite
  standard Gaussian vector $G$. A positive covariance square root gives this
  representation even when singular; no inverse is used. Here $Tz = \ell(b)$ and
  $|z|^2 = v$. Substitution $G' = G+sz$ in its density gives
  $\E_{Q_s}f(\ell(X+sb)) = \E f(\ell(X))$ for every bounded Borel $f$.

  \medskip
  \noindent\textit{Step 2.~}We pass from coordinates to the full norm-Borel
  space.
  The matrix entries of each block in a real Hilbert basis at rational times,
  together with $0,B$, generate the norm-Borel sigma-field. To see the
  operator-norm point, for a countable dense set $\mathcal{D}$ of unit vectors
  obtained from normalized rational finite basis combinations, $\|T\| =
  \sup_{u,w\in\mathcal{D}}|\langle Tu,w\rangle|$. These pairings are finite
  linear combinations of matrix entries. Supremum and H\"older norms of
  continuous paths are suprema over rational times and pairs. Thus, distance to
  every fixed element of $\mathscr{B}_\delta$ is coordinate-measurable. Open
  balls are measurable by rational closed-ball approximation, and separability
  makes every open set a countable union of balls from dense centers. Conversely
  the coordinates are continuous. The cylinder sets form a generating pi-system;
  uniqueness of probability measures extends the finite-coordinate equality to
  all norm-Borel sets. Simple functions and monotone convergence give the
  identity for nonnegative $j$. Finally, take $s = \delta^{-1/2}$ and use $sb =
  H$ from~\eqref{E:raw-block-covariance}. The map and its deterministic
  translate are defined on the entire input space; there is no translated
  exceptional-domain assumption. This proves
  Proposition~\ref{P:raw-block-gaussian-law}.
\end{proof}

The identity retains the random shifted input $X+H$. Replacing its
expected contribution by evaluation at the deterministic mean can lose
the exponential constant, even in a scalar Gaussian calculation; see
No-go~\ref{P:nogo-16}. Positive semidefinite residual covariance does not
justify discarding that residual: No-go~\ref{P:nogo-19} gives an admissible
parent model whose rank-one projection has a larger small-coupling wave mean.

\subsection{Small noisy blocks and a cone in the displacement-velocity space}

For any $x \in \mathscr{B}_0$ define $D_x(u) = S_F(u)U_x(u,0)$ by the
everywhere-defined linear Young map of
Proposition~\ref{P:constant-young-linear}. Write $D_m^0 = D_{X_m}$ and $D_m^+ =
D_{X_m+H_m}$. The latter retains both the mean and the residual noise.

The shifted blocks $D_m^+$ are random. We now show that their largest distance
from the deterministic $P_{s_m}$ of~\eqref{E:compact-F-block}, over all blocks
$m \le N$, tends to zero in every finite $L^q$.

\begin{lemma}[Uniform block error]\label{L:raw-block-uniform-error}
  Let $M_\delta = \max_{m\le N}\|X_m\|_{\mathscr{B}_0}$, including its factor
  $g_\delta$. Then $M_\delta \to 0$ in every finite $L^q$, and, for some $c >
  0$,
  \begin{equation}\label{E:raw-block-maximum-exponential}
    \sup_\delta\E e^{cM_\delta^2}<\infty.
  \end{equation}
  With the deterministic $P_s$ in~\eqref{E:compact-F-block},
  \begin{equation}\label{E:raw-block-operator-error}
    E_{\mathrm{op},\delta}:=
    \max_{m\le N}\sup_{0\le u\le B}
    \|D_m^+(u)-P_{s_m}(u)\|_{\mathcal{B}(F)}
    \longrightarrow0
    \quad\hbox{in every finite }L^q.
  \end{equation}
\end{lemma}

\begin{proof}
  We control the largest scaled Gaussian block by a moment bound and its
  exponential series. We then compare shifted drivers on a fixed ball and use
  uniform higher moments to pass from convergence in probability to
  convergence in every finite $L^q$.

  \medskip
  \noindent\textit{Step 1.~}We bound the scaled maximum and its exponential
  moment.
  The maximum and~\eqref{E:raw-block-moments}, with $N+1 \le 2/\delta$, give
  \[
    \|M_\delta\|_p
    \le (N+1)^{1/p}C\sqrt p\,g_\delta
    \le C'\sqrt p\,\delta^{(1-a)/2-1/p}.
  \]
  For a finite family, the $p$th power of the maximum is at most the sum of
  the $p$th powers. This gives the factor $(N+1)^{1/p}$ without independence.
  For fixed $q$, choose $p > \max(q,2/(1-a))$ to prove the first limit. No
  independence is needed. Choose an integer $n_0$ with $n_0(1-a) \ge 1$. For $n
  \ge n_0$,
  \[
    \E M_\delta^{2n}
    \le2C^{2n}(2n)^n\delta^{n(1-a)-1}
    \le2C^{2n}(2n)^n.
  \]
  For $1 \le n < n_0$, Lyapunov bounds the moment by the appropriate power of
  the $2n_0$ moment. These finitely many terms are bounded. For the remaining
  exponential series, $n! \ge (n/e)^n$ gives a geometric tail $2\sum_{n\ge
  n_0}(2ecC^2)^n$ for small $c$. This
  proves~\eqref{E:raw-block-maximum-exponential}; no such assertion is made for
  the unscaled maximum.

  \medskip
  \noindent\textit{Step 2.~}We compare each noisy block with its frozen
  deterministic flow.
  Let $H_m^{\mathrm{stat}}$ be the primitive for constant potential $V_{s_m}$.
  The primitive estimate is Proposition~\ref{P:constant-compact-base} on
  horizon $B$, applied to $V_{s_m+\delta r}$, $V_{s_m}$ and their difference,
  in the $L^\infty$ and $L^2$ coefficient norms.
  The frozen flow is $P_{s_m}(u)=S_F(u)U_{H_m^{\mathrm{stat}}}(u,0)$;
  Proposition~\ref{P:constant-young-linear} compares it with the driver
  $X_m+H_m$, whose distance from $H_m^{\mathrm{stat}}$ is at most
  $M_\delta+C\delta$ on the fixed driver ball.
  Compact mean smoothness on $[0,2]$ and the statewise primitive bound give
  \[
    \sup_m\max(\|H_m\|_{C^\theta},
    \|H_m^{\mathrm{stat}}\|_{C^\theta})\le H_0,\qquad
    \sup_m\|H_m-H_m^{\mathrm{stat}}\|_{C^\theta}\le C\delta.
  \]
  Indeed, the potential difference is at most $C\delta B$ in $L^\infty$ and
  $L^2$; integration of the bounded free groups gives its Lipschitz primitive
  bound on this fixed horizon. On $\{M_\delta\le1\}$, local Young stability on
  the fixed ball of radius $H_0+1$ gives $E_{\mathrm{op},\delta} \le
  C(M_\delta+\delta)$. This proves convergence in probability. The global Young
  bound also gives
  \[
    \max_m\sup_{u\le B}\|D_m^+(u)\|
    \le C_B\exp\big(C_B(1+H_0+M_\delta)^{1/\theta}\big).
  \]
  Since $1/\theta < 2$, every finite power of this bound is dominated by a
  constant times the exponential in~\eqref{E:raw-block-maximum-exponential}.
  Uniform higher moments yield uniform integrability, proving
  \eqref{E:raw-block-operator-error}. Suprema are measurable: on a countable
  dense set of unit states each difference is continuous in $u$, so rational
  times and those states suffice. Operator-norm continuity of the physical free
  group is not used. Only $s_m \le 1$ supplies the ground templates. This proves
  Lemma~\ref{L:raw-block-uniform-error}.
\end{proof}

Lemma~\ref{L:raw-block-uniform-error} bounds the distance of each shifted block
from $P_{s_m}$. It says nothing about the product of the blocks. The next lemma
is deterministic and holds for every $N \ge 1$. If each factor $D_m$ is within
$\epsilon_{\mathrm{op}}$ of $P_{\sigma_m}$, and adjacent indices $\sigma_m$ are
within $\omega$, then the product satisfies the lower bound
\eqref{E:raw-block-positive-product}.

\LeanMarker{L:raw-block-indexed-cone}
\begin{lemma}[Indexed cone and terminal displacement]
  \label{L:raw-block-indexed-cone}
  Fix $0 < \varrho < 1/2$. There are an integer $B \ge 1$ and positive
  $\omega,\epsilon_{\mathrm{op}},c_*$ such that the following holds for every $N
  \ge 1$. If $\sigma_0,\ldots,\sigma_N \in I$ have adjacent distances at most
  $\omega$ and real bounded operator families $D_m(u)$, $0 \le u \le B$, satisfy
  \[
    \max_{m\le N}\sup_{0\le u\le B}
    \|D_m(u)-P_{\sigma_m}(u)\|\le\epsilon_{\mathrm{op}},
  \]
  then, for $0 \le u \le B$,
  \begin{equation}\label{E:raw-block-positive-product}
    O_\phi D_N(u)D_{N-1}(B)\cdots D_0(B)e_0
    \ge c_*(1-\varrho)^{N-1}
    \exp\left(B\sum_{m<N}\lambda_{\sigma_m}
    +u\lambda_{\sigma_N}\right)>0.
  \end{equation}
  No positivity, composition or independence assumption on the $D_m$ is
  required.
\end{lemma}

\begin{proof}
  We choose the cone aperture before the block length and perturbation
  tolerances. We prove entrance, rebase at each changing index, and finally
  bound the terminal observable even when the last interval is short.

  \medskip
  \noindent\textit{Step 1.~}We choose the margins in an order independent of
  the number of blocks.
  Proposition~\ref{P:compact-F-projection}
  and~\eqref{E:compact-positive-margins} give the positive margins below and
  norm continuity of $e_s,\ell_s,\Pi_s$; compactness of $I$ makes their bounds
  finite and their continuity uniform.
  From~\eqref{E:compact-normalized-complement} fix $C_{\mathrm{gap}} \ge 1$,
  $\eta_{\mathrm{gap}} > 0$ so the normalized complementary norm is at most
  $C_{\mathrm{gap}}e^{-\eta_{\mathrm{gap}}t}$. Define positive margins and
  finite bounds
  \[
    a_*=\inf_I\ell_s(e_0),\quad j_*=\inf_I O_\phi(e_s),\quad
    M_e=\max(1,\sup_I\|e_s\|),\quad
    M_\ell=\max(1,\sup_I\|\ell_s\|).
  \]
  Also put $M_P = \max(1,\sup_I\|\Pi_s\|)$, $M_Q = 1+M_P$, $M_z = M_e+1$, $E_0 =
  \|e_0\| = 1$ and $H_{\mathrm{obs}} = \max(1,\|O_\phi\|)$. Choose the aperture
  \emph{first}: $\vartheta = \min(1,j_*/(16H_{\mathrm{obs}}C_{\mathrm{gap}})) >
  0$. Next choose integer $B \ge 1$ with
  \[
    q_B:=C_{\mathrm{gap}}e^{-\eta_{\mathrm{gap}}B}
    \le\min\left(\frac18,\frac{\vartheta a_*}{8M_QE_0}\right).
  \]
  Choose $\nu > 0$ such that $\nu M_z \le \min(\varrho/2,\vartheta,1/4)$.
  Uniform continuity of $\ell_s,\Pi_s$ supplies $\omega > 0$ for which adjacent
  distances at most $\omega$ give both norm differences at most $\nu$. Finally,
  choose $\epsilon_{\mathrm{op}} > 0$ satisfying all five smallness conditions
  \begin{align*}
    M_\ell\epsilon_{\mathrm{op}}M_z&\le\varrho/2,&
    M_Q\epsilon_{\mathrm{op}}M_z&\le\vartheta/4,\\
    M_\ell\epsilon_{\mathrm{op}}E_0&\le a_*/2,&
    M_Q\epsilon_{\mathrm{op}}E_0&\le\vartheta a_*/8,&
    H_{\mathrm{obs}}\epsilon_{\mathrm{op}}M_z&\le j_*/8.
  \end{align*}
  All denominators are positive. This order of choices uses the normalized
  complement rather than a growing terminal norm.

  \medskip
  \noindent\textit{Step 2.~}We prove entrance into the first cone.
  The identities in~\eqref{E:compact-F-eigen} give
  $\ell_sP_s(B)=e^{\lambda_sB}\ell_s$ and $Q_sP_s(B)=P_s(B)Q_s$, so the
  unperturbed normalized entrance amplitude is $\ell_s(e_0)$ and its residual
  is at most $q_BM_QE_0$; perturbation errors use $e^{-\lambda_sB}\le1$.
  The cone at $s$ consists of $z = Ae_s+r$, with $A = \ell_s(z) > 0$, $\ell_s(r)
  = 0$ and $\|r\|_F \le \vartheta A$. In particular, $r = Q_sz$ and $\|z\| \le
  M_zA$. For entrance, $z_1 = D_0(B)e_0$, with $b = \sigma_0$, has normalized
  amplitude at least $a_*-M_\ell\epsilon_{\mathrm{op}}E_0 \ge a_*/2$. Its
  normalized residual is at most $q_BM_QE_0+M_Q\epsilon_{\mathrm{op}}E_0 \le
  \vartheta a_*/4$. Thus, it enters the cone with residual ratio at most
  $\vartheta/2$ and amplitude at least $e^{\lambda_bB}a_*/2$.

  \medskip
  \noindent\textit{Step 3.~}We rebase the cone and propagate it through each
  full block.
  For an old-cone state $z = Ae_b+r$ and an adjacent index $c$, rebase before
  using the complementary estimate:
  \[
    \ell_c(z)\ge(1-\nu M_z)A,\qquad
    \|Q_cz\|\le\|Q_bz\|+\|(\Pi_b-\Pi_c)z\|
    \le(\vartheta+\nu M_z)A.
  \]
  Since $Q_c-Q_b=\Pi_b-\Pi_c$, the change in the residual is controlled by the
  projection difference times $\|z\|$. Likewise, $\ell_c(z)-\ell_b(z)$ is
  controlled by the dual difference times $\|z\|$. Using $\|z\|\le M_zA$ gives
  the two rebasing errors in the display.
  For the next full block $D$ close to $P_c(B)$, put $x =
  e^{-\lambda_cB}\ell_c(Dz)$ and $y = e^{-\lambda_cB}\|Q_cDz\|$. Since
  $\lambda_c > 0$, the perturbation and smallness conditions give
  \[
    x\ge(1-\nu M_z-M_\ell\epsilon_{\mathrm{op}}M_z)A
    \ge(1-\varrho)A>A/2,
  \]
  and
  \[
    y\le
    [q_B(\vartheta+\nu M_z)+M_Q\epsilon_{\mathrm{op}}M_z]A
    \le\vartheta A/2.
  \]
  The last inequality uses $q_B \le 1/8$, $0 \le \nu M_z \le \vartheta$ and
  $M_Q\epsilon_{\mathrm{op}}M_z \le \vartheta/4$. Together with $\varrho < 1/2$,
  these give $y \le \vartheta x$. These are precisely the finite residual and
  cone-conversion margins; the old residual was not assumed to lie in the new
  kernel. Each full block therefore preserves the new cone and multiplies its
  amplitude by at least $e^{\lambda_cB}(1-\varrho)$.

  \medskip
  \noindent\textit{Step 4.~}We bound the observable on the final partial
  block.
  For the terminal $0 \le u \le B$, rebase once more to $c = \sigma_N$. The same
  bounds give $\ell_c(z) \ge 3A/4$ and $\|Q_cz\| \le 2\vartheta A$.
  Consequently,
  \[
    e^{-\lambda_cu}O_\phi(P_c(u)z)/A
    \ge3j_*/4-2H_{\mathrm{obs}}C_{\mathrm{gap}}\vartheta.
  \]
  The normalized growing component contributes at least $3j_*/4$. The rebased
  residual has norm at most $2\vartheta A$, and the normalized complementary
  flow is at most $C_{\mathrm{gap}}$ for every nonnegative terminal time.
  Applying the observable norm gives the subtracted residual term, uniformly
  even at $u=0$.
  The perturbation changes this by at most
  $H_{\mathrm{obs}}\epsilon_{\mathrm{op}}M_z \le j_*/8$. Since
  $16H_{\mathrm{obs}}C_{\mathrm{gap}}\vartheta \le j_*$, the result is at least
  $j_*/2$. Combine entrance, all $N-1$ intermediate full blocks and this
  terminal estimate to obtain~\eqref{E:raw-block-positive-product}, with $c_* =
  a_*j_*/4$. The induction is empty for $N = 1$ and the estimate includes $u =
  0$. This proves Lemma~\ref{L:raw-block-indexed-cone}.
\end{proof}

\subsection{Source restart and comparison with the original point}

So far the blocks $D_m^0$ and $D_m^+$ are defined one at a time, from the local
operator paths $X_m$ and $X_m+H_m$. We now show that each family forms the
successive pieces of its corresponding global flow, and that its finite
products are source limits in $L^2$ at fixed $\delta,N$. The proof of
Proposition~\ref{P:raw-block-point-comparison} uses this to verify the
source-limit premise of the masked comparison.

\LeanMarker{L:raw-block-restart}
\begin{lemma}[Composition over time blocks and regularized products]
  \label{L:raw-block-restart}
  Fix $\delta,B,N$ and let $L = \delta^{-1}$, $L_{\mathrm{ext}} = (N+1)B \ge L$.
  On this global fast horizon define $A_{\mathrm{glob}} = g_\delta
  A_R(W_{\mathrm{fast}})$. Its physical flow obeys
  \begin{equation}\label{E:raw-block-restart}
    P_{\mathrm{glob}}(mB+u,mB)=D_m^0(u),\qquad 0\le u\le B.
  \end{equation}
  The analogous identity holds with the deterministic mean $V_{\delta t}$ added
  globally and local operator paths $X_m+H_m$. Their finite block products,
  including the last remainder $\tau$, are source limits in $L^2$, at this fixed
  $\delta,N$.
\end{lemma}

\begin{proof}
  We use a single regularized slow noise to realize every fast and local
  source together. An ordinary change of variables gives the restart identity.
  We pass it through the fixed-horizon source limit and control the finite
  products by a telescoping estimate.

  \medskip
  \noindent\textit{Step 1.~}We construct coupled global and local
  regularizations.
  Take one regularized slow noise preserving the given noise $Y_{\epsilon,\eta}$
  on the extension of the given $W$ constructed in
  Section~\ref{S:source-driver-convergence}. Define the fast and local noises by
  \[
    Y_{\mathrm{fast}}(f)=Y_{\epsilon,\eta}(T_{\delta,0}f),\qquad
    Y_m(f)=Y_{\epsilon,\eta}(T_{\delta,m}f).
  \]
  The pullback isometry and~\eqref{E:source-two-dominations} give source and
  residual variance bounds by $\|f\|_{H_W}^2$. The residual tends to zero on
  each transformed smooth compact test at fixed $\delta,m$. Joint Gaussianity is
  retained. Proposition~\ref{P:source-gaussian-transfer} therefore gives the
  $L^pC^\theta$ operator-path limits on $[0,L_{\mathrm{ext}}]$ and on $[0,B]$,
  including both constant columns.

  The field-test integral identity of
  Proposition~\ref{P:source-canonical-convergence} passes through these test
  pullbacks by change of variables; the proof of
  Lemma~\ref{L:constant-source-response} then identifies the ordinary
  primitives using spatial smoothness and total finite propagation.
  Their random fields are $Y_{\mathrm{fast}}(u,x) =
  \delta^{a/2}Y_{\epsilon,\eta}(\delta u,x)$ and $Y_m(u,x) =
  Y_{\mathrm{fast}}(mB+u,x)$. Indeed, these follow by change of variables in
  their test integrals. The temporal law is $K_{\eta/\delta}$ because
  $\delta^aK_\eta(\delta r) = K_{\eta/\delta}(r)$; the spatial heat law and all
  joint width rules are unchanged. They are continuous in time and smooth in
  space. Their ordinary primitives equal their Gaussian operator paths: their
  matrix entries are $Y_m(h_R)$, by the indicator approximation used
  in~\eqref{E:raw-block-mean-jacobian} and the point-test integral identity.
  Total finite propagation extends this identity to both constant columns. This
  does not assume a separately chosen extension preserving the given noise
  automatically commutes with time dilation.

  \medskip
  \noindent\textit{Step 2.~}We identify the restart relation before and after
  removing the source widths.
  Changing the ordinary source integration variable and composing the free
  groups gives the exact identity on $F$
  \[
    S_F(mB)\big[A_{\mathrm{glob}}^{\mathrm{src}}(mB+u)
    -A_{\mathrm{glob}}^{\mathrm{src}}(mB)\big]
    S_F(-mB)=X_m^{\mathrm{src}}(u).
  \]
  No physical mask is commuted through a wave group. Restriction, difference and
  conjugation are bounded maps on these fixed-horizon $C^\theta$ spaces. Passing
  the source limits gives the same identity for $A_{\mathrm{glob}},X_m$. For the
  finitely many $m$, a convergent subsequence or rational times and path
  continuity gives simultaneous identities in $u$. Conjugation norms need only
  be finite at this fixed $\delta$.

  \medskip
  \noindent\textit{Step 3.~}We identify the local Young equations and compose
  the blocks.
  Conjugate the global interaction Young equation on $[mB,mB+u]$ by $S_F(mB)$.
  Bounded linearity of its Riemann sums and that original identity give the
  local equation driven by $X_m$. Uniqueness identifies it with $U_{X_m}$,
  proving \eqref{E:raw-block-restart}. The ordinary deterministic primitive for
  $V_{\delta t}$ has the identical cocycle with local $H_m$. Composition now
  identifies the full product, with terminal $\tau$.

  \medskip
  \noindent\textit{Step 4.~}We pass finite products through the source limit.
  For the unshifted factors use Proposition~\ref{P:constant-gaussian-flow} on
  the fixed horizon $[0,L_{\mathrm{ext}}]$; for the shifted factors use
  Proposition~\ref{P:constant-compact-base} with $g=g_\delta$ and the fixed
  mean $V_{\delta t}$, followed by the restart identity.
  For completeness, source convergence of these finite products requires no
  growing-product stability assertion. At fixed $\delta,N$ let $C \ge 1$ be the
  maximum norm of their source and original factors, and let $e_{\mathrm{src}}$
  be the maximum factor difference. The parent Young convergence gives every
  finite moment of $C$ uniformly in source widths and $e_{\mathrm{src}} \to 0$
  in every finite $L^p$. If $d(n)$ is the difference of the first $n$ ordered
  products, then
  \[
    d(n+1)\le C d(n)+C^n e_{\mathrm{src}},\qquad d(0)=0.
  \]
  Subtract the next pair of products by adding and subtracting the source
  factor times the old partial product. One term is at most $C$ times the
  previous difference; the other is at most $e_{\mathrm{src}}$ times the norm
  of an old product of $n$ factors, bounded by $C^n$. This gives the displayed
  recurrence without commuting factors.
  Extend both finite sequences by identities so the recurrence holds for all
  natural $n$. The finite error-accumulation estimate gives $d(N+1) \le
  (N+1)C^{N+1}e_{\mathrm{src}}$. H\"older, using $C$ in $L^{4(N+1)}$ and the
  error in $L^4$, proves $L^2$ convergence. This is only a fixed-$\delta$ source
  limit. Shifted source products may use the fixed original mean plus the
  regularized noise; no finite-source Cameron mean is substituted. This proves
  Lemma~\ref{L:raw-block-restart}.
\end{proof}

\begin{proposition}[The square observable and canonical comparison]
  \label{P:raw-block-point-comparison}
  Define the signed averages without and with the deterministic mean
  \[
    U_{\delta,B}
    =O_\phi D_N^0(\tau)D_{N-1}^0(B)\cdots D_0^0(B)e_0,
    \qquad
    Y_{\delta,B}
    =O_\phi D_N^+(\tau)D_{N-1}^+(B)\cdots D_0^+(B)e_0.
  \]
  They have all finite unweighted moments, and
  \begin{align}
    \E U_{\delta,B}^2
    &=\E\left[e^{-Z/\sqrt\delta-v/(2\delta)}
    Y_{\delta,B}^2\right],
    \label{E:raw-block-square-identity}\\
    \E U_{\delta,B}^2
    &\le\E|V_\delta(1,0)|^2
    =\E|u(T,0)|^2,\qquad \delta=T^{-\kappa}.
    \label{E:raw-block-original-comparison}
  \end{align}
  The variance $v$ in the density is the original slow original variance.
\end{proposition}

\begin{proof}
  We first apply the exact Gaussian change of law to the squared observable on
  the whole path space. We then verify the source-limit premise for the
  unshifted masked product and compare it with the canonical point moment at
  the scaled amplitude.

  \medskip
  \noindent\textit{Step 1.~}We apply the Gaussian identity to an
  everywhere-defined squared observable.
  On all of $\mathscr{B}_\delta$ put
  \[
    j_\delta(x)=
    |O_\phi D_{x_N}(\tau)D_{x_{N-1}}(B)\cdots D_{x_0}(B)e_0|^2.
  \]
  The everywhere-defined continuous Young map, fixed-time bounded free groups,
  finite products and bounded $O_\phi$ make $j_\delta$ continuous and
  nonnegative on the whole input space. Apply
  Proposition~\ref{P:raw-block-gaussian-law} directly to these original paths.
  Since $X+H$ is its exact translate, this proves
  \eqref{E:raw-block-square-identity} with the single original $Z = W(\Psi)$ and
  its variance $v$. Finite products of the local Young flows have every finite
  unweighted moment by their Gaussian path bounds and H\"older, at fixed
  $\delta,N$. The Gaussian identity itself did not require prior integrability
  of $j_\delta$. No source convergence under its exponential density is used.

  \medskip
  \noindent\textit{Step 2.~}We identify the amplitude-rescaled original model.
  To apply the masked comparison at the amplitude, put $\Gamma_\delta =
  g_\delta^2\Gamma$ and $W' = g_\delta W_{\mathrm{fast}}$. This is a nonzero
  member of the same admissible covariance class, with unchanged degree $\alpha$
  and temporal parameter $a$. Its integrated operator is $g_\delta$ times the
  original one: the original matrix forms are linear in noise, and
  basis-independent Gaussian transfer and the bounded $F$ extension retain both
  constant columns. The $n$th source chain is multiplied by $g_\delta^n$. Thus,
  this is exactly the model without the deterministic mean at that amplitude.

  \medskip
  \noindent\textit{Step 3.~}We pass the unshifted source product to its masked
  limit.
  We verify the source-limit premise of the masked comparison. Keep $\delta,B,N$
  fixed and use the single regularized noise preserving the given noise of
  Lemma~\ref{L:raw-block-restart}. For $m < N$ let $F_m = D_m^0(B)$ and put $F_N
  = D_N^0(\tau)$; write $F_{m,\rho}$ for the corresponding ordinary source
  factors, where both source widths in $\rho$ tend to zero.
  Lemma~\ref{L:raw-block-restart} and Proposition~\ref{P:constant-gaussian-flow}
  give factor convergence in every finite $L^q$ and uniform bounds for every
  finite moment of their norms. These statements are on the fixed finite
  horizons $[0,(N+1)B]$ and $[0,B]$, and include both constant-response columns.
  Set $n = N+1$ and
  \[
    C_\rho=\max\left(1,\max_{m\le N}\|F_{m,\rho}\|,
    \max_{m\le N}\|F_m\|\right),\qquad
    e_\rho=\max_{m\le N}\|F_{m,\rho}-F_m\|.
  \]
  Then $\sup_\rho\|C_\rho\|_q < \infty$ for every finite $q$, and $\|e_\rho\|_q
  \to 0$. Telescoping the ordered products gives
  \[
    \|F_{N,\rho}\cdots F_{0,\rho}-F_N\cdots F_0\|
    \le n C_\rho^n e_\rho.
  \]
  Each term in the telescoping difference contains one factor error and at
  most $n-1$ factors bounded by $C_\rho$. There are $n$ such terms. Since
  $C_\rho\ge1$, replacing $C_\rho^{n-1}$ by $C_\rho^n$ gives the displayed
  bound.
  Therefore, H\"older with exponents $4n$ and $4$ gives
  \[
    \big\|F_{N,\rho}\cdots F_{0,\rho}-F_N\cdots F_0\big\|_{L^2}
    \le n\|C_\rho\|_{L^{4n}}^n\|e_\rho\|_{L^4}
    \longrightarrow0.
  \]
  No uniform estimate in the growing number $N$ is needed: this source limit is
  completed at each fixed $\delta,N$. The physical source identification in
  Lemma~\ref{L:constant-source-solution}, with the same $m_R$ and initial state
  $e_0$, identifies the source product with the tested ordinary masked solution
  at fast time $L = 1/\delta$. Boundedness of $O_\phi$ and $\|e_0\| = 1$ thus
  identify $U_{\delta,B}$ as its $L^2$ limit. The sources $g_\delta
  Y_{\mathrm{fast}}$ have the correct law, contraction and original-test limit
  by Lemma~\ref{L:raw-block-restart}.

  \medskip
  \noindent\textit{Step 4.~}We apply the masked moment comparison at the fixed
  amplitude.
  At positive source widths, $0 \le m_R \le 1$ multiplies each nonnegative
  displacement Wick integrand by the product of its vertex masks. Positive Green
  measures and covariance retain this comparison for unequal, self and empty
  chains as well. Proposition~\ref{P:masked-source-moments} gives the tested
  comparison, and the just-proved $L^2$ convergence discharges the premise of
  Proposition~\ref{P:masked-canonical-transfer} at $\Gamma_\delta$. Its
  canonical unmasked moment limit uses the fixed-order joint source limits and
  the absolutely summable order bounds. Since $\phi \ge 0$ and $\int\phi = 1$,
  these facts give
  \[
    \E U_{\delta,B}^2
    \le\E|u_{\Gamma_\delta}(L,0)|^2.
  \]
  Only the source products without the deterministic mean were needed for this
  passage. The shifted original product is already the everywhere-defined
  $j_\delta(X+H)$ in the Gaussian identity; no tilted source limit or change of
  the original slow $\Psi,Z,v$ has been made.

  \medskip
  \noindent\textit{Step 5.~}We identify the canonical fast-time moment with
  the large-time moment.
  The ordinary change from fast $u$ to slow $t = \delta u$, with phase
  $(V,\delta V_t)$, transforms the unmasked source equation into the regularized
  version of~\eqref{E:exact-fast-equation}, with slow terminal time one. Indeed,
  $g_\delta Y_{\mathrm{fast}}(u) = \sqrt\delta\,Y_{\epsilon,\eta}(\delta u)$,
  and its fast cutoff is $\eta/\delta$. The canonical fixed-order and absolutely
  convergent series construction passes this law identity at fixed $\delta$ to
  the original moment, preserving every mixed, unequal, self and empty chain.
  Lemma~\ref{L:exact-canonical-scaling} gives the final equality
  in~\eqref{E:raw-block-original-comparison}. The physical scale $r_T = T\delta$
  is distinct from the fixed mask radius $R$. The mask need not contain the
  growing backward cone: the comparison is a moment inequality, not equality of
  masked and full fields. All source limits preceded this large-time comparison.
  This proves Proposition~\ref{P:raw-block-point-comparison}.
\end{proof}

Proposition~\ref{P:raw-block-point-comparison} bounds $\E|u(T,0)|^2$ from below
by the weighted second moment of $Y_{\delta,B}$
in~\eqref{E:raw-block-square-identity}. We still need a lower bound for
$Y_{\delta,B}$ itself. The proof below obtains one on an event whose probability
tends to one. It uses the estimates of Lemma~\ref{L:raw-block-uniform-error} and
the cone bound of Lemma~\ref{L:raw-block-indexed-cone}.

\LeanMarker{P:raw-block-original-lower}
\begin{proposition}[Original point lower rate at the Gaussian quadratic cost]
  \label{P:raw-block-original-lower}
  For every original positive $\Psi$ and eligible fixed mask $R$ chosen above,
  and every fixed $x \in \R^d$,
  \begin{equation}\label{E:raw-block-original-lower}
    \liminf_{T\to\infty}T^{-\kappa}\log\E|u(T,x)|^2
    \ge 2\int_0^1\lambda_s\ud s-\frac{v_\Psi}{2}.
  \end{equation}
\end{proposition}

\begin{proof}
  We fix the cone constants before taking the small-parameter limit. A single
  high-probability event controls every block and the original Gaussian
  variable. On that event the cone bound and a Riemann sum give the growth
  rate, while the linear Gaussian penalty vanishes after normalization.

  \medskip
  \noindent\textit{Step 1.~}We fix the rate tolerance and bound the coupled
  block maximum.
  Fix $e_{\mathrm{rate}} > 0$. Choose $0 < \varrho < 1/2$ with $\log(1-\varrho)
  > -e_{\mathrm{rate}}/3$, then choose the cone constants, including an integer
  $B$ fixed independently of $\delta$, from
  Lemma~\ref{L:raw-block-indexed-cone}. The control $\Psi$ and eligible mask $R$
  remain fixed throughout.

  Choose a finite $p > 2/(1-a)$ and put $\gamma_p = p(1-a)/2-1 > 0$.
  From~\eqref{E:raw-block-moments} and $N+1 \le 2/\delta$, Markov and the union
  bound give, for $\eta > 0$,
  \[
    \begin{split}
      \mathbb{P}(M_\delta>\eta)
      &\le\eta^{-p}\sum_{m=0}^N
      \E\|X_m\|_{C^\theta}^p\\
      &\le A_p\eta^{-p}\delta^{\gamma_p},\qquad
      A_p=2(C_{R,B}\sqrt p)^p<\infty.
    \end{split}
  \]
  This estimate concerns the coupled blocks; no independence is needed. At these
  fixed data $A_p$ is independent of $\delta$.

  \medskip
  \noindent\textit{Step 2.~}We choose a deterministic block tolerance and one
  common good event.
  We also specify the deterministic tolerance used with this bound. The
  primitive estimates in Lemma~\ref{L:raw-block-uniform-error} give, for the
  constant-potential primitives $H_m^{\mathrm{stat}}$,
  \[
    \max_m\max(\|H_m\|_{C^\theta},
    \|H_m^{\mathrm{stat}}\|_{C^\theta})\le H_0,\qquad
    \max_m\|H_m-H_m^{\mathrm{stat}}\|_{C^\theta}\le C_B\delta.
  \]
  Apply Proposition~\ref{P:constant-young-linear} on the fixed ball of operator
  paths of radius $H_0+1$ and use the bounded free groups on $[0,B]$. There is
  $L_B \ge 1$, depending on these fixed $R,B,\Psi$ and the admissible exponents,
  such that on $\{M_\delta\le1\}$,
  \[
    E_{\mathrm{op},\delta}\le L_B(M_\delta+\delta).
  \]
  Set $\eta = \min(1,\epsilon_{\mathrm{op}}/(2L_B)) > 0$. If $\delta \le
  \epsilon_{\mathrm{op}}/(2L_B)$, then $M_\delta \le \eta$ implies
  $E_{\mathrm{op},\delta} \le \epsilon_{\mathrm{op}}$. Thus, for the measurable
  event on the original common space
  \[
    E_\delta=\{E_{\mathrm{op},\delta}\le\epsilon_{\mathrm{op}}\}
    \cap\{|Z|\le\delta^{-1/4}\},\qquad
    p_\delta=\mathbb{P}(E_\delta),
  \]
  Chebyshev for the same single $Z = W(\Psi)$, with $\E Z^2 = v$, gives
  \[
    p_\delta\ge1-A_p\eta^{-p}\delta^{\gamma_p}
    -v\delta^{1/2}\longrightarrow1.
  \]
  This is one union bound, with no independence between blocks or from $Z$ and
  no product of block probabilities.

  \medskip
  \noindent\textit{Step 3.~}We make the block and Riemann-sum errors smaller
  than the fixed tolerance.
  Let $C_\lambda$ be a Lipschitz bound for $s\mapsto\lambda_s$. Take $\delta$
  small enough, in addition to $0 < \delta < \min(1,1/B)$, that
  \[
    \begin{gathered}
      \delta B\le\omega,\qquad
      \delta\le\epsilon_{\mathrm{op}}/(2L_B),\qquad
      C_\lambda\delta B\le e_{\mathrm{rate}}/3,\\
      \delta|\log c_*|\le e_{\mathrm{rate}}/3,\qquad
      A_p\eta^{-p}\delta^{\gamma_p}+v\delta^{1/2}\le1/2.
    \end{gathered}
  \]
  Every constant has been fixed before this choice. The condition involving
  $C_\lambda$ is vacuous if that constant is zero. All source limits in
  Proposition~\ref{P:raw-block-point-comparison} have already been completed at
  each such fixed $\delta$.

  The indices $s_m = m\delta B$ lie in $I$ with adjacent gaps at most $\omega$.
  On $E_\delta$ the shifted full-$F$ blocks therefore satisfy
  Lemma~\ref{L:raw-block-indexed-cone}, retaining both phase entries and the
  final remainder $\tau$. In particular, $Y_{\delta,B} > 0$, and
  \[
    \delta\log Y_{\delta,B}
    \ge\delta\log c_*+\delta(N-1)\log(1-\varrho)
    +\delta B\sum_{m<N}\lambda_{s_m}
    +\delta\tau\lambda_{s_N}.
  \]
  Since $s_N+\delta\tau = 1$, the last two terms are the left Riemann sum for
  $I_R$ on the entire interval, including its last interval of length
  $\delta\tau$, with mesh at most $\delta B$. Its error is at most
  $C_\lambda\delta B$. Also $0 \le \delta(N-1) \le 1/B \le 1$, so the
  logarithmic loss is at least $\log(1-\varrho) > -e_{\mathrm{rate}}/3$. The two
  other losses are each at most $e_{\mathrm{rate}}/3$ by the choices above.
  Hence,
  \[
    Y_{\delta,B}\ge
    \exp((I_R-e_{\mathrm{rate}})/\delta)\quad\hbox{on }E_\delta,
  \]
  including the $N = 1$ and $\tau = 0$ cases of the cone estimate.

  \medskip
  \noindent\textit{Step 4.~}We use the same event in the shifted second-moment
  identity.
  The integrand in~\eqref{E:raw-block-square-identity} is nonnegative. On the
  same $E_\delta$ its linear Gaussian term satisfies $-Z/\sqrt\delta \ge
  -\delta^{-3/4}$. Restriction to this event
  and~\eqref{E:raw-block-original-comparison} therefore give the explicit bound
  for the original point moment
  \[
    \E|u(T,0)|^2\ge\E U_{\delta,B}^2
    \ge p_\delta\exp\left(
    \frac{2I_R-2e_{\mathrm{rate}}-v/2}{\delta}
    -\delta^{-3/4}\right),\qquad \delta=T^{-\kappa}.
  \]
  Here $p_\delta \ge 1/2$ for small $\delta$, so $-\delta\log2 \le \delta\log
  p_\delta \le 0$. Multiplying the logarithm of the preceding inequality by
  $\delta$ yields
  \[
    \delta\log\E|u(T,0)|^2
    \ge2I_R-2e_{\mathrm{rate}}-v/2
    -\delta^{1/4}+\delta\log p_\delta.
  \]
  Thus, the normalized linear penalty is $\delta\delta^{-3/4} = \delta^{1/4} \to
  0$, and $\delta\log p_\delta \to 0$. The quadratic Gaussian cost is precisely
  half the original variance $v$. Since $\kappa > 0$, $\delta = T^{-\kappa}$ is
  a bijection between all sufficiently large $T$ and all sufficiently small
  $\delta$. This gives the full continuous-time liminf $Q_2 \ge
  2I_R-2e_{\mathrm{rate}}-v/2$. For each $e_{\mathrm{rate}}$ the argument first
  fixed its own $B = B(e_{\mathrm{rate}})$ and then took $\delta\downarrow0$.
  The original point rate $Q_2$ is independent of $B$; only now let
  $e_{\mathrm{rate}}\downarrow0$. No choice $B = B(\delta)$ is made.

  \medskip
  \noindent\textit{Step 5.~}We extend the point bound by translation.
  Finally, spatial translation of all positive-width vertices preserves the
  original noise law, Green kernels and constant initial data. The fixed-order
  $L^2$ and absolute-series limits preserve equality of point second moments.
  Hence, zero can be replaced by every fixed $x$
  in~\eqref{E:raw-block-original-lower}. Every control and mask were fixed
  before asymptotics; a supremum over them may follow their individual
  inequalities. This proves Proposition~\ref{P:raw-block-original-lower}.
\end{proof}

\section{The unmasked Rayleigh rate and nonnegative controls}
\label{S:unmasked-variational}

Fix a nonzero covariance satisfying Hypothesis~\ref{H:homog}, including its
original amplitude, and set $a = \alpha_0$. Define the full class
\begin{equation}\label{E:variational-positive-class}
  \mathcal{C}_{\mathrm{pos}}
  =\{\Psi\in C_c^\infty((0,1)\times\R^d;\R):
  \Psi\ge0,\ \Psi\ne0\}.
\end{equation}

It includes nonseparable tests. It is nonempty: if $b(z) = e^{-1/(1-z^2)}$ on
$|z| < 1$ and zero otherwise, then $b(4s-2)\prod_{j=1}^db(x_j)$ belongs to this
class. All boundary derivatives vanish because exponential decay dominates every
power of $(1-z^2)^{-1}$. For each $\Psi$ retain the original mean $h_{\Psi,s}$
and original variance $v_\Psi = \E W(\Psi)^2$ from~\eqref{E:compact-raw-mean}--%
\eqref{E:compact-raw-cost}. Neither the test nor its covariance is changed by a
physical mask. We write $h_s$ when $\Psi$ is fixed.

\subsection{Global bounds and the unmasked spectral bottom}

We first record a bound for the covariance smoothed by the heat kernel. The
proof of Lemma~\ref{L:variational-global-mean} uses it to bound
$\Gamma*|D\Psi(t,\cdot)|$ uniformly in $t,x$.

\LeanMarker{L:variational-heat-envelope}
\begin{lemma}[Gaussian envelope for the original covariance]
  \label{L:variational-heat-envelope}
  For the normalized heat kernel $p_\epsilon(x) =
  (2\pi\epsilon)^{-d/2}e^{-|x|^2/(2\epsilon)}$,
  \begin{equation}\label{E:variational-heat-envelope}
    (\Gamma*p_\epsilon)(x)
    =\int e^{i\xi\cdot x}e^{-\epsilon|\xi|^2/2}\mu(\ud\xi),
    \qquad
    0\le(\Gamma*p_\epsilon)(x)
    \le(\Gamma*p_\epsilon)(0)<\infty.
  \end{equation}
  The value at zero is strictly positive.
\end{lemma}

\begin{proof}
  We first use temperedness to make all Gaussian pairings finite. We identify
  the convolution and spectral expressions as continuous distributions, then
  use positivity to obtain the bound by the value at zero.

  \medskip
  \noindent\textit{Step 1.~}We obtain integrable majorants for the Gaussian
  weights.
  In the dyadic summation below, mass growth $2^{jN}$ times weight decay
  $2^{-jL}$ is summable for $L>N$; Gaussian derivatives decrease faster than
  each such polynomial weight.
  A positive tempered Radon measure has polynomial ball growth. Indeed, evaluate
  its tempered bound on a nonnegative smooth cutoff equal to one on $B_1$,
  dilated by $R \ge 1$. Its finitely many Schwartz seminorms are bounded by
  $CR^N$, so positivity gives a ball-mass bound $CR^N$. Summing dyadic annuli
  shows that $(1+|x|)^{-L}$ is integrable for $L > N$. Apply this separately to
  $\Gamma$ and $\mu$. Gaussian weights are therefore integrable against both
  measures. The classical convolution on the left is finite and smooth locally
  in $x$: on compact $x$ sets its translated Gaussian derivatives are bounded by
  any fixed sufficiently fast polynomial weight in the integration variable.

  \medskip
  \noindent\textit{Step 2.~}We identify the two continuous expressions and
  bound their common value.
  Pair the convolution with a Schwartz test $\varphi$. Absolute Fubini is
  justified by these weights. The resulting pairing is that of $\Gamma$ with
  $p_\epsilon*\varphi$. Applying the spectral convention in
  Hypothesis~\ref{H:homog} and the identity $\int p_\epsilon(y)e^{-i\xi\cdot
  y}\ud y = e^{-\epsilon|\xi|^2/2}$ gives exactly the right side
  of~\eqref{E:variational-heat-envelope} as a distribution. The heat kernel is
  even, so either Fourier sign gives that same Gaussian multiplier. Its spectral
  integral is absolutely finite and continuous in $x$ by dominated convergence.
  Thus, equality of the two continuous distributions gives equality at every
  $x$, without an extra normalization factor. The physical convolution is
  nonnegative. Evenness of $\mu$ makes the spectral integral real, and its
  absolute value is bounded by its total Gaussian mass, the value at zero.
  Finally, $\Gamma \ne 0$ and $p_\epsilon > 0$ give strict positivity there.
  This proves Lemma~\ref{L:variational-heat-envelope}.
\end{proof}

\LeanMarker{L:variational-global-mean}
\begin{lemma}[Global regularity of the original mean]
  \label{L:variational-global-mean}
  Every mixed time and space derivative of $h_s$ is globally bounded in space,
  uniformly on compact time intervals. In particular, for $I = [0,1]$ there are
  finite constants $M,C_\Psi$ such that
  \begin{equation}\label{E:variational-global-mean}
    0\le h_s(x)\le M,\qquad
    \|h_s-h_t\|_\infty\le C_\Psi|s-t|,
    \qquad s,t\in I.
  \end{equation}
\end{lemma}

\begin{proof}
  Extend $\Psi$ by zero in time; its compact support inside $(0,1)$ makes this a
  smooth extension. For any mixed derivative $D\Psi$, a fixed $\epsilon > 0$ and
  compactness of its common spatial support give $|D\Psi(t,x)| \le
  C_Dp_\epsilon(x)$ for all $t,x$. Here the heat kernel has a positive minimum
  on that support. Lemma~\ref{L:variational-heat-envelope} consequently bounds
  $\Gamma*|D\Psi(t,\cdot)|$ uniformly in $t,x$. Ordinary local differentiation
  against the Radon measure gives the corresponding derivatives of the
  convolution. After the change $r = s-t$,
  \[
    h_s(x)=\int_\R |r|^{-a}
    (\Gamma*\Psi(s-r,\cdot))(x)\ud r.
  \]

  For $s$ in a compact interval the relevant $r$ lie in a fixed bounded
  interval. The preceding global bound times $|r|^{-a}$ is integrable since $a <
  1$. Every derivative thus falls on the smooth test and passes through the
  integral by dominated convergence, with the claimed global spatial bounds.
  Nonnegativity follows from positive kernels. The bounded time derivative
  proves the Lipschitz estimate. No directional decay of the mean has been
  inferred. This proves Lemma~\ref{L:variational-global-mean}.
\end{proof}

Lemma~\ref{L:variational-global-mean} shows that $h_s$ is bounded by $M$ for $s
\in I$ and is Lipschitz in $s$ in the supremum norm. We now use these two bounds
to define the rate $\lambda_{\Psi,s}$ of the potential $h_s$ without a mask and
to show that it is Lipschitz on $I$.

\LeanMarker{P:variational-unmasked-rate}
\begin{proposition}[The unmasked Rayleigh rate]
  \label{P:variational-unmasked-rate}
  On complex $L^2(\R^d)$, the operator $H_s = -\Delta-h_s$ is self-adjoint on
  $H^2$, with form domain $H^1$ and form $q_s(f) = \|\nabla f\|_2^2-\int
  h_s|f|^2$. Its spectral bottom $E_s$ is uniformly negative on $I$. The
  positive rate $\lambda_{\Psi,s} = \sqrt{-E_s}$ satisfies
  \begin{equation}\label{E:variational-real-rayleigh}
    \lambda_{\Psi,s}^2
    =\sup_{\substack{f\in H^1(\R^d;\R)\\\|f\|_2=1}}
    \left(\int h_sf^2-\|\nabla f\|_2^2\right).
  \end{equation}
  It is Lipschitz on $I$, and $I_\Psi: = \int_0^1\lambda_{\Psi,s}\ud s$ is
  positive and finite. No eigenstate at the unmasked bottom is asserted.
\end{proposition}

\begin{proof}
  We construct the unmasked self-adjoint operator using bounded perturbation
  and identify its form domain. The spectral theorem gives the real Rayleigh
  formula; the common compact trial makes its bottom uniformly negative and
  its square-root rate Lipschitz.

  \medskip
  \noindent\textit{Step 1.~}We prove self-adjointness on the stated operator
  domain.
  The Fourier Laplacian $A = -\Delta$ is self-adjoint on $H^2$, being
  multiplication by $|\xi|^2$ with its maximal domain. Multiplication $B_s =
  -h_s$ is bounded self-adjoint with norm at most $M$. Thus, $H_s = A+B_s$ is
  densely defined, symmetric and closed on $H^2$; the two graph norms are
  equivalent under a bounded perturbation. For real $b > M$,
  \[
    H_s-ib=(I+B_s(A-ib)^{-1})(A-ib).
  \]
  The first factor has a Neumann inverse because $\|B_s(A-ib)^{-1}\| \le M/b <
  1$, and the same holds with $+ib$. To check self-adjointness directly, for $u
  \in D(H_s^*)$ solve $(H_s-ib)v = (H_s^*-ib)u$ with $v \in D(H_s)$. Then $u-v
  \in \ker(H_s^*-ib) = \operatorname{Ran}(H_s+ib)^\perp = \{0\}$. Hence,
  $D(H_s^*) = D(H_s)$.

  \medskip
  \noindent\textit{Step 2.~}We identify the form domain by equivalent
  completions.
  Set $c_{\mathrm{shift}} = 1+M$ and $J_s = H_s+c_{\mathrm{shift}} \ge 1$.
  Integration by parts on $H^2$ gives
  \[
    \|J_s^{1/2}f\|_2^2=q_s(f)+c_{\mathrm{shift}}\|f\|_2^2,\qquad
    \|f\|_{H^1}^2\le\|J_s^{1/2}f\|_2^2
    \le c_{\mathrm{shift}}\|f\|_{H^1}^2.
  \]
  Spectral cutoffs $1_{[1,n]}(J_s)$ show that $H^2 = D(J_s)$ is a core for the
  half-power norm; Fourier cutoffs give its density in $H^1$. Completion in the
  equivalent norms identifies $D(J_s^{1/2})$ and $H^1$ as $L^2$ subspaces, since
  both norms dominate $L^2$. This proves the stated form domain.

  \medskip
  \noindent\textit{Step 3.~}We obtain the real Rayleigh formula without an
  eigenstate assumption.
  The spectral lower bound gives one direction of the Rayleigh identity. For the
  reverse, nonzero spectral projections in $[E_s,E_s+\epsilon]$ give unit $H^2$
  vectors with Rayleigh values at most $E_s+\epsilon$. The bottom is finite
  because the operator is bounded below and has finite Rayleigh tests. Let
  $\epsilon\downarrow0$. Real tests give the same infimum: $q_s(u+iv) =
  q_s(u)+q_s(v)$, so the quotient is the norm-weighted average of the nonzero
  real quotients. This proves the Rayleigh identity with $-E_s$ on the left. The
  following negative upper bound justifies writing that value as
  $\lambda_{\Psi,s}^2$ in~\eqref{E:variational-real-rayleigh}.

  For a complex vector $f = u+iv$ with real-valued $u,v$, its quadratic form
  and squared norm are sums of the corresponding values for $u$ and $v$. At
  least one of the nonzero components has quotient no greater than their
  norm-weighted average.

  \medskip
  \noindent\textit{Step 4.~}We obtain strict negativity and uniform rate
  continuity.
  Proposition~\ref{P:compact-original-common-trial} supplies, for this original
  $\Psi$, one unit compactly supported $w$ and a finite plateau mask with a
  common strictly negative upper bound $\epsilon_0 > 0$. Its trial has support
  in $B_{2L}$; choose its plateau radius $R_0 \ge 2L$. The unmasked form on this
  same $w$ is no larger. Consequently, $-M \le E_s \le -\epsilon_0 < 0$, and
  $\sqrt{\epsilon_0} \le \lambda_{\Psi,s} \le \sqrt M$. Comparing identical unit
  form tests gives $|E_s-E_t| \le \|h_s-h_t\|_\infty$ and therefore
  \[
    |\lambda_{\Psi,s}-\lambda_{\Psi,t}|
    \le \frac{C_\Psi}{2\sqrt{\epsilon_0}}|s-t|.
  \]
  This also proves positivity and finiteness of $I_\Psi$. No compact-potential
  spectral theorem was applied to $H_s$. This proves
  Proposition~\ref{P:variational-unmasked-rate}.
\end{proof}

\subsection{Removing masks from the numerical lower rate}

Proposition~\ref{P:raw-block-original-lower} gives the lower bound with the rate
of the masked potential $m_Rh_s$, for one fixed mask. We now let the mask radius
$R$ tend to infinity and replace this rate by the rate $\lambda_{\Psi,s}$ of
Proposition~\ref{P:variational-unmasked-rate}.

\LeanMarker{P:variational-remove-mask}
\begin{proposition}[Monotone rates and the original point moment]
  \label{P:variational-remove-mask}
  There is an increasing family of admissible radial masks $m_R$ tending to one,
  for which the compact ground rates satisfy
  $\lambda_{R,s}\uparrow\lambda_{\Psi,s}$ uniformly on $I$ as $R \to \infty$,
  with the compact ground statements used for $R \ge R_0(\Psi)$. For the
  original canonical solution and any fixed $x$,
  \begin{equation}\label{E:variational-unmasked-point-lower}
    \liminf_{T\to\infty}T^{-\kappa}\log\E|u(T,x)|^2
    \ge2I_\Psi-\frac{v_\Psi}{2}.
  \end{equation}
\end{proposition}

\begin{proof}
  We choose nested physical masks, compare their Rayleigh suprema with the
  unmasked supremum, and prove uniform convergence in the time parameter. We
  then take a numerical supremum of already established fixed-mask bounds.

  \medskip
  \noindent\textit{Step 1.~}We construct the nested masks.
  Put $F(t) = e^{-1/t}$ for $t > 0$ and zero otherwise, and define
  \[
    \chi(r)=\frac{F(2-r)}{F(2-r)+F(r-1)},\qquad
    m_R(x)=\chi(|x|/R).
  \]
  The denominator is positive. The numerator decreases and the second
  denominator term increases; thus $\chi$ is nonincreasing, equals one on $r \le
  1$, and zero on $r \ge 2$. Flatness of $F$ at zero gives smoothness. Since
  $\chi$ is constant near zero, $m_R$ is smooth also at the spatial origin.
  These even $[0,1]$ masks increase pointwise to one.

  \medskip
  \noindent\textit{Step 2.~}We pass the Rayleigh values and rates to their
  unmasked limits.
  For each fixed unit form test $f$, the mask-removal argument uses the
  integrable majorant $Mf^2$ for $m_Rh_sf^2$.
  The compact-family input is Proposition~\ref{P:compact-ground-family}, with
  the common trial supplied by
  Proposition~\ref{P:compact-original-common-trial}: its rate obeys
  $\lambda_{R,s}^2=-E_1$, the negative form bottom of $-\Delta-m_Rh_s$ on
  $H^1$.
  Choose $R_0$ so the preceding common trial lies in the plateau. Every $R \ge
  R_0$ has the same negative margin. For each fixed $R$, the compact-family
  theorem applies to $m_Rh_s$. Its real Rayleigh formula follows from the same
  real/complex quotient argument as~\eqref{E:variational-real-rayleigh}. The
  rates are continuous, positive, at most $\sqrt M$ and monotone in $R$. For
  every fixed real unit $f \in H^1$, dominated convergence gives $\int m_Rh_sf^2
  \to \int h_sf^2$. Hence,
  \[
    \sup_{R\ge R_0}\lambda_{R,s}^2
    =\sup_{\|f\|_2=1,\ f\in H^1(\R^d;\R)}
    \sup_{R\ge R_0}
    \left(\int m_Rh_sf^2-\|\nabla f\|_2^2\right)
    =\lambda_{\Psi,s}^2.
  \]
  Only two numerical suprema were interchanged. Positive square roots give the
  pointwise rate convergence. For $R_n = R_0+n$ and $\epsilon > 0$, the
  increasing open sets $\{s:\lambda_{\Psi,s}-\lambda_{R_n,s}<\epsilon\}$ cover
  compact $I$. A finite subcover lies in the largest one, proving uniform
  convergence. Monotonicity extends it to all real $R \to \infty$. Therefore,
  $\int_I\lambda_{R,s}\ud s\uparrow I_\Psi$.

  \medskip
  \noindent\textit{Step 3.~}We take the supremum of fixed-mask moment lower
  bounds.
  At each fixed $R$ use this very same mask in the operator path and mean when
  applying Proposition~\ref{P:raw-block-original-lower}. It gives the lower
  bound with $2\int_I\lambda_{R,s}\ud s-v_\Psi/2$ for the same original point
  liminf. Taking the numerical supremum over these individual bounds proves
  \eqref{E:variational-unmasked-point-lower}. There is no $R(T)$, unmasked noisy
  wave, uniform-in-mask Gaussian estimate or exchange of source and time limits.
  This proves Proposition~\ref{P:variational-remove-mask}.
\end{proof}

\subsection{Full-class dilation and finite variational values for nonnegative controls}

The bound~\eqref{E:variational-unmasked-point-lower} holds for each $\Psi \in
\mathcal{C}_{\mathrm{pos}}$ separately. Before we take the supremum over $\Psi$,
we record how $h_{\Psi,s}$, $v_\Psi$, $\lambda_{\Psi,s}$ and $I_\Psi$ change
under a spatial dilation of the control. The proof of
Proposition~\ref{P:variational-positive-scores} uses this to relate the
variational values for different coefficients $c$.

\LeanMarker{L:variational-control-dilation}
\begin{lemma}[Radon and Rayleigh dilation]\label{L:variational-control-dilation}
  For $r > 0$ the map
  \begin{equation}\label{E:variational-control-dilation}
    (D_r\Psi)(s,x)=r^{d+2-\alpha}\Psi(s,rx)
  \end{equation}
  is a bijection of $\mathcal{C}_{\mathrm{pos}}$, with inverse $D_{1/r}$. With
  the same original covariance,
  \begin{equation}\label{E:variational-dilation-values}
    h_{D_r\Psi,s}(x)=r^2h_{\Psi,s}(rx),\qquad
    v_{D_r\Psi}=r^{4-\alpha}v_\Psi,\qquad
    \lambda_{D_r\Psi,s}=r\lambda_{\Psi,s},\quad
    I_{D_r\Psi}=rI_\Psi.
  \end{equation}
\end{lemma}

\begin{proof}
  We apply the homogeneous Radon pushforward to the covariance mean and the
  Lebesgue Jacobian to its cost. We then dilate the unit Rayleigh tests to
  identify the rate and its time integral.

  \medskip
  \noindent\textit{Step 1.~}We scale the original mean and covariance cost.
  Positive dilation preserves the full smooth positive nonzero class;
  multiplication of its amplitude powers gives the stated inverse. Homogeneity
  of the physical Radon measure gives, for every nonnegative measurable $G$,
  \[
    \int G(ry)\Gamma(\ud y)
    =r^{\alpha-d}\int G(z)\Gamma(\ud z).
  \]
  Indeed, its left pushforward measure on a Borel set $A$ is $\Gamma(A/r) =
  r^{\alpha-d}\Gamma(A)$. Positive Tonelli applies this identity inside the
  temporal convolution, giving $h_{D_r\Psi,s}(x) = r^2h_{\Psi,s}(rx)$. Now use
  the \emph{Lebesgue} substitution $z = rx$ in the original pairing:
  \[
    v_{D_r\Psi}=r^{d+2-\alpha}r^2r^{-d}v_\Psi
    =r^{4-\alpha}v_\Psi.
  \]
  This distinguishes the Radon pushforward from the spatial volume factor; it
  includes singular covariances and white noise.

  \medskip
  \noindent\textit{Step 2.~}We scale the full set of Rayleigh tests.
  The map $U_rf(x) = r^{d/2}f(rx)$ is a unitary bijection on real $L^2$,
  preserving $H^1,H^2$ as sets under its inverse. Weak derivatives give
  $\|\nabla U_rf\|_2^2 = r^2\|\nabla f\|_2^2$ and $\int h_{D_r\Psi,s}(U_rf)^2 =
  r^2\int h_{\Psi,s}f^2$. The real Rayleigh bijection proves the squared rate
  identity; the rates are positive, so their square roots and integrals give
  \eqref{E:variational-dilation-values}. On $H^2$ this is also $H_{D_r\Psi,s}U_r
  = r^2U_rH_{\Psi,s}$. Time has not been rescaled. For compact masks the
  corresponding identity uses $m_R(rx) = m_{R/r}(x)$; a fixed mask radius is not
  invariant. This proves Lemma~\ref{L:variational-control-dilation}.
\end{proof}

\LeanMarker{P:variational-positive-scores}
\begin{proposition}[Intrinsic variational lower bounds]
  \label{P:variational-positive-scores}
  For $c > 0$, the set of variational values
  $\{cI_\Psi-v_\Psi/2:\Psi\in\mathcal{C}_{\mathrm{pos}}\}$ is nonempty and
  bounded above. Its finite real supremum $S_c^{\mathrm{pos}}$ is positive and
  satisfies
  \begin{equation}\label{E:variational-score-homogeneity}
    S_c^{\mathrm{pos}}=c^{(4-\alpha)/(3-\alpha)}S_1^{\mathrm{pos}},
    \qquad
    Q_2\ge S_2^{\mathrm{pos}},
  \end{equation}
  where $Q_2$ is the original point liminf. The coefficient-one value is an
  auxiliary variational value; no first-moment or source-constant identification
  is asserted.
\end{proposition}

\begin{proof}
  We first construct a positive coefficient-one value and scale the entire
  class of values. We use the independently proved point upper bound to show
  these sets are bounded above before taking their finite real suprema.

  \medskip
  \noindent\textit{Step 1.~}We exhibit a positive value.
  Write $k = 3-\alpha > 0$, $q = 4-\alpha = k+1$ and $m = q/k > 1$. Fix one
  $\Psi_0$ from the nonempty class and put $I_0 = I_{\Psi_0} > 0$, $v_0 =
  v_{\Psi_0} > 0$. Positivity of the latter is
  Lemma~\ref{L:compact-raw-profile}. For $r_0 = (I_0/v_0)^{1/k} > 0$, which need
  not be less than one, the exact real-power identity gives
  \[
    I_{D_{r_0}\Psi_0}-v_{D_{r_0}\Psi_0}/2
    =r_0I_0-r_0^qv_0/2=r_0I_0/2>0.
  \]
  The definition of $r_0$ gives $r_0^k=I_0/v_0$. Since $q=k+1$, it follows
  that $r_0^qv_0=r_0I_0$. Thus, subtracting half the cost leaves exactly half
  of the positive linear term.
  This positive value is obtained from a smooth nonnegative control.

  \medskip
  \noindent\textit{Step 2.~}We scale the whole coefficient-one set
  bijectively.
  For $c > 0$ put $r_c = c^{1/k}$ and $K_c = c^m > 0$.
  Lemma~\ref{L:variational-control-dilation} and the exact scalar powers give,
  for every $\Psi$,
  \begin{equation}\label{E:variational-score-set-scaling}
    cI_{D_{r_c}\Psi}-v_{D_{r_c}\Psi}/2
    =K_c(I_\Psi-v_\Psi/2).
  \end{equation}
  The two scalar factors agree: $cr_c=c^{1+1/k}=c^{q/k}=K_c$ and
  $r_c^q=c^{q/k}=K_c$. Consequently, dilation multiplies the linear reward and
  the quadratic cost by the same factor for this choice of $r_c$.
  Since $D_{r_c}$ is a bijection of the entire class, the whole set of
  variational values with coefficient $c$ is $K_c$ times the whole
  coefficient-one set of variational values.

  \medskip
  \noindent\textit{Step 3.~}We bound the sets using the original point upper
  estimate.
  We prove boundedness before taking finite-real suprema. The independent
  Theorem~\ref{T:full-upper}, for the \emph{same canonical solution}, gives
  \[
    \E|u(T,0)|^2\le e(T+2)^d e^{C_{\mathrm{up}}T^\kappa},
    \qquad T\ge1,\quad C_{\mathrm{up}}<\infty.
  \]
  Its prefactor is negligible: choose $\epsilon = \kappa/2 > 0$, use $\log T \le
  (T^\epsilon-1)/\epsilon$ and $\log(T+2) \le \log T+\log3$ for $T \ge 1$, and
  divide by $T^\kappa$. Thus, the original point limsup is at most
  $C_{\mathrm{up}}$. One fixed test/mask lower bound gives a finite lower bound
  for $Q_2$, so $Q_2$ is finite. By~\eqref{E:variational-unmasked-point-lower}
  every coefficient-two variational expression is at most this same $Q_2 \le
  C_{\mathrm{up}}$, uniformly in $\Psi$.
  Equation~\eqref{E:variational-score-set-scaling} at $c = 2$ consequently
  bounds all coefficient-one variational values by $C_{\mathrm{up}}/K_2$. The
  nonempty set has a finite real supremum $S_1^{\mathrm{pos}} > 0$ by the
  positive test above.

  \medskip
  \noindent\textit{Step 4.~}We take finite real suprema and retain their
  interpretation.
  For general $c$, set equality and $K_c > 0$ show that $K_cS_1^{\mathrm{pos}}$
  is the supremum: it is an upper bound, and every smaller number is exceeded by
  $K_c$ times a coefficient-one variational expression within the corresponding
  positive distance of its supremum. This proves finiteness, positivity and the
  first identity in~\eqref{E:variational-score-homogeneity}. Taking the
  numerical supremum of the coefficient-two lower bounds proves the second.
  Although $2^m > 2$, the resulting $S_2^{\mathrm{pos}} > 2S_1^{\mathrm{pos}}$
  is only a comparison of auxiliary variational values. It is not a Jensen
  improvement without a first-moment identification. This proves
  Proposition~\ref{P:variational-positive-scores}.
\end{proof}

For $\Gamma = 0$, the mean, original variance and Rayleigh rate are zero. Set
all these variational values to zero; the free canonical solution $u = 1$ has
the corresponding zero moment rate. Strict positivity above is restricted to
nonzero covariance. Matching an upper variational rate and identifying an exact
pointwise limit remain later steps.

\paragraph{What we now have.}
The lower bound $Q_2\ge S_2^{\mathrm{pos}}$
(Proposition~\ref{P:variational-positive-scores}). The upper bound takes the
next three chapters.

\chapter{Spatial controls and Gaussian means}\label{S:book-means}

The lower bound used smooth controls. The upper bound must also handle controls
completed in the Gaussian covariance norm, where pointwise formulas may no
longer exist. So we need to know how such a control, the deterministic mean
it induces, and the wave operators generated by that mean are related.

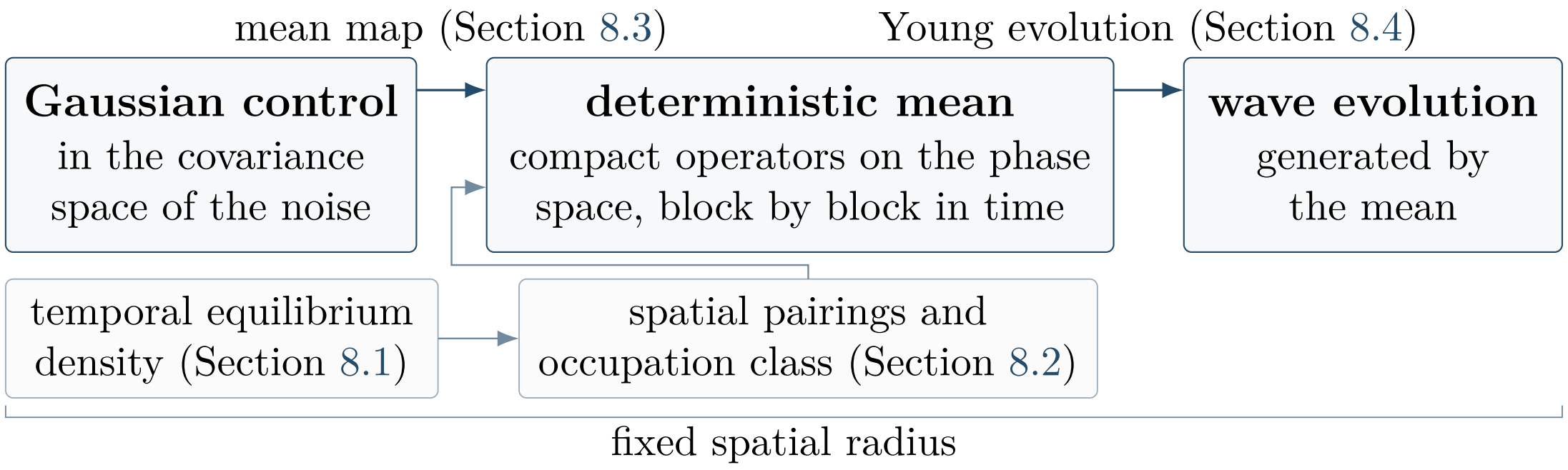
\begin{figure}[htbp]
  \centering
  \providecommand{\FlowSecEq}{\ref{S:temporal-equilibrium}}
  \providecommand{\FlowSecSp}{\ref{S:spatial-occupation-foundation}}
  \providecommand{\FlowSecMean}{\ref{S:actual-cameron-mean}}
  \providecommand{\FlowSecYoung}{\ref{S:young-upper-foundation}}
  % Construction order only; the spatial radius remains fixed.
\providecommand{\FlowSecEq}{8.1}%
\providecommand{\FlowSecSp}{8.2}%
\providecommand{\FlowSecMean}{8.3}%
\providecommand{\FlowSecYoung}{8.4}%
\begin{tikzpicture}[x=1cm,y=1cm,font=\normalsize,every node/.style={inner sep=0pt}]
  \foreach \centre/\width in {1.9/3.8,7.35/5.8,12.65/3.5} {
  \node[draw=hamblue,fill=hamblue!4,rounded corners=2pt,
  minimum width=\width cm,minimum height=1.8cm,inner sep=0pt]
  at (\centre,0) {};
  }
  \node[font=\normalsize\bfseries,anchor=base] at (1.9,.36) {Gaussian control};
  \node[font=\normalsize\bfseries,anchor=base] at (7.35,.36) {deterministic mean};
  \node[font=\normalsize\bfseries,anchor=base] at (12.65,.36) {wave evolution};
  \node[font=\small,anchor=base] at (1.9,-.13) {in the covariance};
  \node[font=\small,anchor=base] at (1.9,-.60) {space of the noise};
  \node[font=\small,anchor=base] at (7.35,-.13) {compact operators on the phase};
  \node[font=\small,anchor=base] at (7.35,-.60) {space, block by block in time};
  \node[font=\small,anchor=base] at (12.65,-.13) {generated by};
  \node[font=\small,anchor=base] at (12.65,-.60) {the mean};
  \draw[ham/arrow] (3.8,.60)--(4.45,.60);
  \draw[ham/arrow] (10.25,.60)--(10.9,.60);
  \node[font=\small,anchor=base] at (4.125,1.06) {mean map (Section~\FlowSecMean)};
  \node[font=\small,anchor=base] at (10.575,1.06) {Young evolution (Section~\FlowSecYoung)};
  % Secondary inputs: equilibrium feeds spatial pairings, not the main row directly.
  \foreach \centre/\width in {2/4,7.425/5.35} {
  \node[draw=hamblue!45,fill=hamblue!2,rounded corners=2pt,
  line width=.35pt,minimum width=\width cm,minimum height=1.1cm]
  at (\centre,-1.70) {};
  }
  \node[font=\small,anchor=base] at (2,-1.57) {temporal equilibrium};
  \node[font=\small,anchor=base] at (2,-2.04) {density (Section~\FlowSecEq)};
  \node[font=\small,anchor=base] at (7.425,-1.57) {spatial pairings and};
  \node[font=\small,anchor=base] at (7.425,-2.04)
  {occupation class (Section~\FlowSecSp)};
  \draw[ham/arrow,draw=hamblue!65,line width=.45pt]
  (4,-1.70)--(4.75,-1.70);
  \draw[ham/arrow,draw=hamblue!65,line width=.45pt]
  (7.425,-1.15)--(7.425,-1.02)--(4.125,-1.02)--(4.125,-.30)--(4.45,-.30);
  % The bracket spans the whole construction, keeping the input arrow clear.
  \draw[hamblue!65,line width=.35pt]
  (0,-2.32)--(0,-2.43)--(14.4,-2.43)--(14.4,-2.32);
  \node[font=\small,anchor=base] at (7.2,-2.78) {fixed spatial radius};
\end{tikzpicture}
  \caption[From a Gaussian control to the wave it drives]{How this chapter is
    organized. The arrows show the order of construction only, and the
    spatial radius is fixed throughout.}
  \label{Fig:control-mean-wave}
\end{figure}

Figure~\ref{Fig:control-mean-wave} is the picture to keep in
mind. Section~\ref{S:temporal-equilibrium} computes a temporal equilibrium
density. Section~\ref{S:spatial-occupation-foundation} extends the physical
spatial pairings and defines the class of time-dependent spatial profiles
used in the later energy bounds. Section~\ref{S:actual-cameron-mean}
builds the compact map from a Gaussian control to its mean on time
blocks. Section~\ref{S:young-upper-foundation} estimates perturbations of
the Young evolution while keeping the deterministic endpoint.

The chapter uses the variational value of Section~\ref{S:unmasked-variational}
and the wave phase of Section~\ref{S:constant-phase-young}. On a first reading,
keep Propositions~\ref{P:spatial-four-linear}, \ref{P:cameron-full-mean}
and~\ref{P:cameron-actual-embedding} and let the proofs wait. The spatial
radius is still fixed; the next chapter decides which bounds survive as
it grows. Episode~\VideoRef{6}{VWemsluFZiA} of the video series covers
this chapter.

\section{The temporal equilibrium density}
\label{S:temporal-equilibrium}

The density in Lemma~\ref{L:temporal-equilibrium} is the classical interval
Riesz equilibrium density. Its probability-normalized form on $[-1,1]$ is
recalled by Calef and Hardin~\cite{calef.hardin:source2008:riesz}*{p.~1}.
We prove the potential identity with the normalization used here.

We use Euler's function $\Gamma_E(x) = \int_0^\infty u^{x-1}e^{-u}\ud u$, $x >
0$, and $\operatorname{Beta}(x,y) = \int_0^1t^{x-1}(1-t)^{y-1}\ud t$, $x,y > 0$.
Both integrals are positive and finite.

\LeanMarker{L:equilibrium-gamma-beta}
\begin{lemma}[Gamma and Beta identities]\label{L:equilibrium-gamma-beta}
  For positive arguments,
  \begin{align*}
    \Gamma_E(x+1)&=x\Gamma_E(x),\\
    \operatorname{Beta}(x,y)
    &=\frac{\Gamma_E(x)\Gamma_E(y)}{\Gamma_E(x+y)},\\
    \Gamma_E(2x)
    &=\frac{2^{2x-1}\Gamma_E(x)\Gamma_E(x+1/2)}{\sqrt\pi}.
  \end{align*}
\end{lemma}

\begin{proof}
  Integration by parts gives the recurrence, with vanishing boundary terms. In
  the positive product Gamma integral, use $(u,v) = (Rt,R(1-t))$, with Jacobian
  $R$; Tonelli gives the Beta--Gamma identity. For the last identity, symmetry,
  then $t = (1+v)/2$ and $v^2 = w$, give $\operatorname{Beta}(x,x) =
  2^{1-2x}\operatorname{Beta}(1/2,x)$. Also $\Gamma_E(1/2) = \int_\R e^{-z^2}\ud
  z = \sqrt\pi$: the substitution $u = z^2$ gives the first equality, and
  squaring the positive Gaussian integral and using polar coordinates gives
  $2\pi\int_0^\infty re^{-r^2}\ud r = \pi$. Substituting these into the
  Beta--Gamma identity proves duplication. This proves
  Lemma~\ref{L:equilibrium-gamma-beta}.
\end{proof}

\LeanMarker{L:temporal-equilibrium}
\begin{lemma}[Exact potential and mass]\label{L:temporal-equilibrium}
  For $0 < a < 1$ put $r = (1-a)/2$ and, for $0 < t < 1$,
  \begin{equation}\label{E:temporal-equilibrium-density}
    f_a(t)=\frac{[t(1-t)]^{-r}}{\Gamma_E(r)\Gamma_E(1-r)}.
  \end{equation}
  This positive function belongs to $L^1(0,1)\cap L^2(0,1)$, and
  \begin{align}
    \int_0^1|s-t|^{-a}f_a(t)\ud t&=1
    &&\text{for every }s\in[0,1],\label{E:temporal-equilibrium-potential}\\
    M_a:=\int_0^1f_a(t)\ud t
    &=\frac{\Gamma_E((1+a)/2)}
    {\Gamma_E(1+a)\Gamma_E((1-a)/2)}>0.
    \label{E:temporal-equilibrium-mass}
  \end{align}
\end{lemma}

\begin{proof}
  We first justify the positive series used in the integral calculation. We
  compute the two halves of the potential through an incomplete Beta integral,
  treat the endpoints directly, and then compute the mass and integrability.

  \medskip
  \noindent\textit{Step 1.~}We justify the positive binomial expansion.
  The positive binomial series $(1-z)^{-p} = \sum_{n\ge0}(p)_nz^n/n!$, $p > 0$,
  $0 \le z < 1$, follows directly as follows. Its coefficient ratio gives
  uniform convergence on compact subintervals, including for its derivative.
  Termwise differentiation gives $(1-z)y' = py$, $y(0) = 1$; hence $(1-z)^py$ is
  constant one. Here $(p)_0 = 1$, $(p)_n = p(p+1)\cdots(p+n-1)$. All terms used
  below are nonnegative, so Tonelli applies.

  \medskip
  \noindent\textit{Step 2.~}We compute the potential at interior points and
  endpoints.
  For $0 < s < 1$ split the integral of $|s-t|^{2r-1}[t(1-t)]^{-r}$ into left
  and right parts. In the left part, $t = su$ gives
  \begin{align*}
    L(s)
    &=s^r\int_0^1(1-u)^{2r-1}u^{-r}(1-su)^{-r}\ud u\\
    &=\operatorname{Beta}(1-r,2r)s^r
    \sum_{n\ge0}
    \frac{(r)_n(1-r)_ns^n}{(1+r)_n n!}.
  \end{align*}
  The Beta ratio follows from the Gamma recurrence. Since $(r)_n/(1+r)_n =
  r/(r+n)$, expanding the incomplete Beta integral gives
  \[
    L(s)=r\operatorname{Beta}(1-r,2r)B_s(r,r),\qquad
    B_s(r,r)=\int_0^s u^{r-1}(1-u)^{r-1}\ud u.
  \]
  Expanding $(1-u)^{r-1}$ in $B_s(r,r)$ gives coefficient $(1-r)_n/n!$ and
  integral $s^{r+n}/(r+n)$. Multiplication by $r\operatorname{Beta}(1-r,2r)$
  therefore gives exactly the coefficients of $L(s)$.
  Reflection makes the right part $L(1-s)$, while symmetry gives
  $B_s(r,r)+B_{1-s}(r,r) = \operatorname{Beta}(r,r)$. The full integral is
  therefore $r\operatorname{Beta}(1-r,2r)\operatorname{Beta}(r,r) =
  \Gamma_E(r)\Gamma_E(1-r)$. All integrals are finite because $-r > -1$ and
  $2r-1 > -1$. At $s = 0$ the integral is directly $\operatorname{Beta}(r,1-r) =
  \Gamma_E(r)\Gamma_E(1-r)$; $s = 1$ follows by reflection. This proves
  \eqref{E:temporal-equilibrium-potential}, including the endpoints without a
  limiting argument there.

  \medskip
  \noindent\textit{Step 3.~}We verify both integrability classes and the
  stated mass.
  Integrability and square integrability follow from $-r > -1$ and $-2r = a-1 >
  -1$. Finally,
  \[
    M_a=
    \frac{\operatorname{Beta}(1-r,1-r)}
    {\Gamma_E(r)\Gamma_E(1-r)}
    =\frac{\Gamma_E(1-r)}{\Gamma_E(2-2r)\Gamma_E(r)},
  \]
  which is~\eqref{E:temporal-equilibrium-mass}. Its time mass has not been
  renormalized. This proves Lemma~\ref{L:temporal-equilibrium}.
\end{proof}

\section{Physical spatial pairings and the full occupation class}
\label{S:spatial-occupation-foundation}

Retain a nonzero original covariance from Hypothesis~\ref{H:homog} with its
amplitude. Set $a = \alpha_0$, $k = 3-\alpha > 1$, $q = 4-\alpha = k+1$ and $D =
q-2a > 0$. The scalar $D$ is fixed. All profiles are real. The original positive
temporal exponent and threshold imply $\alpha < 2$; this is the strict
inequality used in the scalar powers and concavity. The finite variational value
over nonnegative controls $S_1^{\mathrm{pos}}$ is that of
Proposition~\ref{P:variational-positive-scores}.

\LeanMarker{L:spatial-physical-pair}
\begin{lemma}[The physical two-point integral]\label{L:spatial-physical-pair}
  The positive measure $\nu(\ud x,\ud y) = \ud x\,\Gamma(\ud y)$ is
  sigma-finite. For nonnegative functions define
  \begin{equation}\label{E:spatial-physical-pair}
    B_\Gamma(f,h)=\int f(x)h(x-y)\nu(\ud x,\ud y)
  \end{equation}
  initially as an extended nonnegative integral, using finite Borel
  representatives of Lebesgue classes. Its value is representative independent.
  Lebesgue-a.e. convergence transfers to $\nu$-a.e. convergence at both
  positions $x$ and $x-y$.
\end{lemma}

\begin{proof}
  Local finiteness of $\Gamma$ gives sigma-finiteness on exhausting balls, and
  Lebesgue measure is sigma-finite. If $N$ is Borel and Lebesgue null, Tonelli
  applied first in $x$ gives $\nu\{x\in N\} = 0$ and $\nu\{x-y\in N\} = 0$,
  because every fixed translate of $N$ is null. An arbitrary Lebesgue null set
  has a Borel null superset. Thus, arbitrary versions are measurable in the
  $\nu$ completion with the same value, and the same observation applies to the
  exceptional sets of a.e. convergence. Signed pairings will be used only after
  absolute convergence. This proves Lemma~\ref{L:spatial-physical-pair}.
\end{proof}

Put $A_{\mathrm{sp}}(g) = B_\Gamma(g^2,g^2)$ and $b(g) = \|\nabla g\|_2^2$ on
$H^1$; initially $A_{\mathrm{sp}}$ is used only for smooth compact profiles
before the extension below. For smooth compact $z$, write $U_z = \Gamma*z$ in
its usual pointwise sense.

\LeanMarker{L:spatial-static-controls}
\begin{lemma}[Smooth temporal controls]\label{L:spatial-static-controls}
  For every nonnegative nonzero $z \in C_c^\infty(\R^d)$, let
  \[
    \lambda_z(c)^2=
    \sup_{\substack{g\in H^1(\R^d;\R)\\\|g\|_2=1}}
    \left(c\int U_zg^2-b(g)\right),\qquad c>0.
  \]
  This is positive and finite, and $\lambda_z$ is Lipschitz on every compact
  positive $c$ interval. With the exact $M_a$
  in~\eqref{E:temporal-equilibrium-mass},
  \begin{equation}\label{E:spatial-static-score}
    S_1^{\mathrm{pos}}\ge
    \lambda_z(1)-\frac{M_a}{2}B_\Gamma(z,z).
  \end{equation}
\end{lemma}

\begin{proof}
  We first obtain the static Rayleigh rate from an admissible smooth control.
  We then approximate the exact temporal equilibrium by increasing compact
  smooth controls and pass their rewards and original covariance costs to the
  limit separately.

  \medskip
  \noindent\textit{Step 1.~}We obtain the static rate and its local Lipschitz
  bound.
  Choose nonnegative nonzero $f \in C_c^\infty(0,1)$ and $s_0 = 1/2$. Let
  $k_f(s_0) = \int_0^1|s_0-t|^{-a}f(t)\ud t > 0$. The smooth test $\Psi(t,x) =
  cf(t)z(x)/k_f(s_0)$ has mean $cU_z$ at $s_0$.
  Proposition~\ref{P:variational-unmasked-rate} gives boundedness of this
  potential, the displayed real Rayleigh formula, and a positive rate for every
  $c > 0$. Comparing identical unit tests makes $\lambda_z(c)^2$ Lipschitz with
  constant $\|U_z\|_\infty$. On a positive compact interval its rate has a
  positive floor by monotonicity, giving the square-root Lipschitz bound.

  \medskip
  \noindent\textit{Step 2.~}We approximate the equilibrium density by
  admissible controls.
  Use the nonincreasing smooth cutoff $\chi$ constructed in the proof of
  Proposition~\ref{P:variational-remove-mask}. For integers $n \ge 4$, set on
  $(0,1)$
  \[
    f_n(t)=f_a(t)\chi(1/(nt))\chi(1/(n(1-t))),\qquad
    k_n(s)=\int_0^1|s-t|^{-a}f_n(t)\ud t.
  \]
  Extend $f_n$ by zero. It is smooth compact, vanishes near both endpoints, is
  nonzero, and increases to $f_a$. Thus, $\Psi_n = f_n\otimes z$ is a positive
  smooth control. Monotone convergence
  and~\eqref{E:temporal-equilibrium-potential} give $k_n(s)\uparrow1$ for every
  $s \in [0,1]$. Since the kernel is at least one on the unit square off its
  null diagonal, $0 < c_0: = \int_0^1f_4 \le k_n(s) \le 1$. The mean is
  $k_n(s)U_z$; hence its rate is $\lambda_z(k_n(s))$, converging to
  $\lambda_z(1)$ and bounded by that finite number. Dominated convergence gives
  $I_{\Psi_n} \to \lambda_z(1)$.

  \medskip
  \noindent\textit{Step 3.~}We pass the covariance cost and the score to the
  limit.
  The Gaussian cost is exactly the nonnegative Tonelli pairing
  \[
    v_{\Psi_n}=B_\Gamma(z,z)\int_0^1f_n(s)k_n(s)\ud s.
  \]
  Its time integrand increases to $f_a$, so the cost tends to
  $M_aB_\Gamma(z,z)$. The spatial factor is finite and positive: its smooth
  compact autocorrelation is positive near zero, and homogeneity scales any
  origin-centered ball of positive $\Gamma$ mass to every other positive radius.
  Such a ball exists by exhaustion. Take the limit of the individual
  inequalities $S_1^{\mathrm{pos}} \ge I_{\Psi_n}-v_{\Psi_n}/2$. The variational
  expression differences need not increase. No control depends on original large
  time and no rough stochastic evolution has been used. This proves
  Lemma~\ref{L:spatial-static-controls}.
\end{proof}

The next estimate also allows zero covariance and does not use the finiteness
of $S_1^{\mathrm{pos}}$.

\LeanMarker{L:spatial-fourier-bound}
\begin{lemma}[Deterministic Fourier bound for the spatial energy]
  \label{L:spatial-fourier-bound}
  Under Hypothesis~\ref{H:homog}, allowing $\Gamma=0$, let
  $g\in C_c^\infty(\R^d;\R)$ satisfy $\|g\|_2=1$ and put
  $b=\|\nabla g\|_2^2$. Then $b>0$ and
  \begin{equation}\label{E:spatial-fourier-bound}
    0\le A_{\mathrm{sp}}(g)
    \le 2^{1+\alpha}R_\alpha b^{\alpha/2},
  \end{equation}
  where $R_\alpha$ is defined in~\eqref{E:intro-spectral-mass}.
\end{lemma}

\begin{proof}
  We identify the spatial energy spectrally without assuming a measure
  density. A mass bound and a derivative bound on the Fourier transform give a
  resolvent majorant, whose homogeneous integral supplies the desired
  constant.

  \medskip
  \noindent\textit{Step 1.~}We identify the energy as a finite spectral
  integral.
  If $b=0$, smoothness makes $\nabla g=0$ everywhere. Thus, $g$ is constant,
  and compact support forces $g=0$, contrary to $\|g\|_2=1$. Put $Q=g^2$ and
  use the unnormalized Fourier transform
  \[
    \mathcal{F}_0Q(\xi)=\int_{\R^d}e^{-i\xi\cdot x}Q(x)\ud x.
  \]
  In the unitary convention of Section~\ref{S:masked-covariance}, this is
  $c_d^{-1}\widehat Q$, with $c_d=(2\pi)^{-d/2}$.
  The autocorrelation
  \[
    c_Q(y)=\int_{\R^d}Q(x)Q(x-y)\ud x
  \]
  is smooth and compactly supported. Apply the distributional spectral
  identity~\eqref{E:masked-spectral-convention} to this Schwartz function.
  Its inner Fourier integral is $|\mathcal{F}_0Q(\xi)|^2$, since $Q$ is real.
  Hence
  \[
    A_{\mathrm{sp}}(g)=\int c_Q(y)\Gamma(\ud y)
    =\int|\mathcal{F}_0Q(\xi)|^2\mu(\ud\xi).
  \]
  Both integrals are finite: $c_Q$ is compactly supported and $\Gamma$ is
  Radon, while $|\mathcal{F}_0Q|^2$ is Schwartz and $\mu$ is tempered.
  No density of either measure is required.

  \medskip
  \noindent\textit{Step 2.~}We bound the Fourier transform by a resolvent
  weight.
  Since $\|Q\|_1=1$, we have $|\mathcal{F}_0Q|\le1$. Integration by parts
  and the Cauchy--Schwarz inequality give
  \[
    |\xi|\,|\mathcal{F}_0Q(\xi)|
    \le\|\nabla Q\|_1
    \le2\|g\|_2\|\nabla g\|_2=2\sqrt b.
  \]
  Splitting into $|\xi|^2\le4b$ and $|\xi|^2>4b$ therefore yields
  \[
    |\mathcal{F}_0Q(\xi)|^2\le\frac{8b}{4b+|\xi|^2}.
  \]
  In the region $|\xi|^2\le4b$, the right side of the preceding display is at
  least one, so the mass bound applies. In the other region, the derivative
  bound gives $4b/|\xi|^2$, which is at most $8b/(4b+|\xi|^2)$ because
  $4b\le|\xi|^2$. Thus, the two estimates give one majorant over the whole
  spectrum.

  \medskip
  \noindent\textit{Step 3.~}We integrate the majorant and keep the original
  spectral normalization.
  The scalar integral in that resolvent calculation is
  \[
    \int_0^\infty\frac{u^{\alpha-1}}{1+u^2}\ud u
    =\frac12\operatorname{Beta}(\alpha/2,1-\alpha/2)
    =\frac{\pi}{2\sin(\pi\alpha/2)},
  \]
  by the Beta integral and Euler's reflection identity.
  In the following resolvent formula, $q$ is a free positive argument, separate
  from the fixed exponent $q=4-\alpha$ used elsewhere in this section.
  The deterministic resolvent calculation in the proof of
  Lemma~\ref{L:damped-homogeneous-norms}, together
  with~\eqref{E:intro-spectral-mass}, gives, for $q>0$,
  \[
    \mathcal{R}_\Gamma(q):=\int\frac{\mu(\ud\xi)}{q^2+|\xi|^2}
    =q^{\alpha-2}R_\alpha<\infty.
  \]
  Consequently,
  \[
    A_{\mathrm{sp}}(g)\le8b\mathcal{R}_\Gamma(2\sqrt b)
    =2^{1+\alpha}R_\alpha b^{\alpha/2}.
  \]
  Nonnegativity follows from the spectral identity. If $\Gamma=0$, then
  $\mu=0$ and both sides of~\eqref{E:spatial-fourier-bound} vanish. This proves
  Lemma~\ref{L:spatial-fourier-bound}.
\end{proof}

Lemma~\ref{L:spatial-static-controls} gives a lower bound for each nonzero,
nonnegative smooth spatial control of compact support. We now scale and
dilate a single profile and optimize. The resulting bound depends on the
profile only through the ratio of $A_{\mathrm{sp}}$ to a power of $b$.

\LeanMarker{P:spatial-smooth-ratio}
\begin{proposition}[The finite smooth spatial ratio]
  \label{P:spatial-smooth-ratio}
  For real unit $g \in C_c^\infty(\R^d)$ put $R(g) =
  A_{\mathrm{sp}}(g)/b(g)^{\alpha/2}$. Its supremum $R_0$ over these profiles is
  a positive finite real number, and
  \begin{equation}\label{E:spatial-smooth-score}
    S_1^{\mathrm{pos}}\ge
    k\left[
    \frac{R(g)\alpha^{\alpha/2}}{8M_aq^{q/2}}
    \right]^{1/k}.
  \end{equation}
\end{proposition}

\begin{proof}
  We first bound the smooth spatial ratio. We then vary the amplitude of a
  squared dilated profile so that its Rayleigh test has a prescribed positive
  value, and optimize the resulting one-variable score.

  \medskip
  \noindent\textit{Step 1.~}We prove positivity and finiteness of the ratio
  supremum.
  The gradient energy $b(g)$ is positive: a nonzero $L^2$ function with zero
  weak gradient would have Fourier transform supported only at zero, and hence
  would vanish. Its square is nonnegative nonzero smooth compact, so the
  autocorrelation argument in the proof of
  Lemma~\ref{L:spatial-static-controls} gives
  $0 < A:=A_{\mathrm{sp}}(g)<\infty$. Lemma~\ref{L:spatial-fourier-bound} gives
  $R(g)\le2^{1+\alpha}R_\alpha$ for every such profile. A normalized compact
  bump supplies a positive ratio. Thus, $0<R_0<\infty$, independently of the
  finiteness of $S_1^{\mathrm{pos}}$. It remains to prove the score inequality.
  With $U_\ell g(x)=\ell^{d/2}g(\ell x)$, the physical changes of variables give
  \begin{equation}\label{E:spatial-dilation}
    A_{\mathrm{sp}}(U_\ell g)=\ell^\alpha A,\qquad
    b(U_\ell g)=\ell^2b(g).
  \end{equation}
  The factors in the first identity are $\ell^{2d}$ from the two squares,
  $\ell^{-d}$ from $\ud x$, and $\ell^{\alpha-d}$ from the Radon pushforward.

  \medskip
  \noindent\textit{Step 2.~}We choose the amplitude and spatial scale of a
  single control.
  For $\theta,\ell > 0$ choose the smooth control $z = \theta(U_\ell g)^2$. The
  Rayleigh test $U_\ell g$ gives $\lambda_z(1)^2 \ge \theta\ell^\alpha
  A-\ell^2b$, while $B_\Gamma(z,z) = \theta^2\ell^\alpha A$. For any $L > 0$,
  put $u = \alpha/q > 0$ and choose
  \[
    \ell^2=uL^2/b,\qquad
    \theta=(1+u)L^2/(\ell^\alpha A).
  \]
  In the cost,
  substituting these same choices gives $\theta^2\ell^\alpha
  A=(1+u)^2L^4/(\ell^\alpha A)$; replacing $\ell^\alpha A$ by
  $u^{\alpha/2}L^\alpha R(g)$ yields the coefficient $C$ below.
  The radicand is $L^2$. Lemma~\ref{L:spatial-static-controls} consequently
  gives
  \[
    S_1^{\mathrm{pos}}\ge L-CL^q,\qquad
    C=\frac{M_a(1+u)^2}{2R(g)u^{\alpha/2}}>0.
  \]
  Its maximum occurs at $L^k = 1/(qC)$ and is $(k/q)(qC)^{-1/k}$. Since $1+u =
  4/q$, $1/(q^qC) = R(g)\alpha^{\alpha/2}/(8M_aq^{q/2})$,
  proving~\eqref{E:spatial-smooth-score}. This argument retains the finite
  real value $S_1^{\mathrm{pos}}$ fixed at the start of the section; only the
  proof of $R_0<\infty$ is independent of that input. This proves
  Proposition~\ref{P:spatial-smooth-ratio}.
\end{proof}

So far $A_{\mathrm{sp}}$ has been used only for smooth compact profiles.
Proposition~\ref{P:occupation-physical-class} needs it and the square pairings
for real $H^1$ profiles, so we extend them with the same constant $R_0$.

\LeanMarker{P:spatial-four-linear}
\begin{proposition}[The physical four-linear extension]
  \label{P:spatial-four-linear}
  For real $f \in H^1(\R^d)$ define $N(f) = \|f\|_2^{1-\alpha/4}\|\nabla
  f\|_2^{\alpha/4}$. Every physical square pairing is finite and
  \begin{align}
    A_{\mathrm{sp}}(f)&\le R_0N(f)^4,\label{E:spatial-homogeneous-bound}\\
    0\le B_\Gamma(f^2,h^2)
    &\le\sqrt{A_{\mathrm{sp}}(f)A_{\mathrm{sp}}(h)}.
    \label{E:spatial-square-CS}
  \end{align}
  The signed integral $C_4(f_1,f_2,f_3,f_4) = B_\Gamma(f_1f_2,f_3f_4)$ is
  absolutely convergent and jointly continuous on $(H^1)^4$, with
  \begin{equation}\label{E:spatial-four-linear}
    \int|f_1(x)f_2(x)f_3(x-y)f_4(x-y)|\ud\nu
    \le R_0\prod_{j=1}^4N(f_j).
  \end{equation}
  The ratio supremum $R_{\mathrm{sp}}$ over all real unit $H^1$ profiles equals
  $R_0$. Every nonzero $H^1$ profile has $A_{\mathrm{sp}}(f) > 0$.
\end{proposition}

\begin{proof}
  We first obtain the square-pair bound for smooth profiles and extend it by
  Fatou. We then use ordinary Cauchy–Schwarz under the physical measure to
  prove an absolutely convergent four-linear bound. That bound, rather than
  Fatou alone, gives continuity and the full-class ratio identity.

  \medskip
  \noindent\textit{Step 1.~}We establish the smooth square-pair estimates.
  For smooth compact $f$, normalization in $R_0$ proves
  \eqref{E:spatial-homogeneous-bound}, including $f = 0$. All exponents in $N$
  are positive and $N(f) \le \|f\|_{H^1}$. For smooth compact signed $Q,P$, the
  absolute correlation is compactly supported and bounded by $\|Q\|_2\|P\|_2$.
  Local finiteness of $\Gamma$ therefore makes its physical integral finite.
  Smooth signed covariance is positive: with the ordinary transform
  $\widehat{Q}_{\mathrm{ord}}(\xi) = \int e^{i\xi\cdot x}Q(x)\ud x$, the model
  convention gives
  \[
    B_\Gamma(Q,Q)=\int|\widehat{Q}_{\mathrm{ord}}(\xi)|^2\mu(\ud\xi)
    \ge0,\qquad Q\in C_c^\infty(\R^d;\R).
  \]
  Pair $\Gamma$ with the smooth compact correlation to obtain this identity;
  rapid Fourier decay and tempered $\mu$ give finiteness. Evenness of $\Gamma$
  and $(x,y)\mapsto(x-y,-y)$ give symmetry. The quadratic polynomial for $Q+tP$
  gives covariance Cauchy--Schwarz. In particular, for smooth real $f,h$,
  \[
    0\le B_\Gamma(f^2,h^2)
    \le R_0N(f)^2N(h)^2.
  \]
  The squares of signed smooth profiles are already smooth and nonnegative.

  \medskip
  \noindent\textit{Step 2.~}We approximate Sobolev profiles and obtain finite
  physical pairings.
  Smooth compact functions are dense in $H^1$. Cutoffs equal to one on the
  radius-$R$ ball, vanishing beyond $2R$ and with gradient at most $C/R$,
  converge by the function and gradient tails and the bound $C\|f\|_2/R$ for the
  extra gradient term. Mollification then converges in $H^1$ by $L^2$
  translation continuity of the function and its weak gradients. The space is
  separable through its closed isometric image $f\mapsto(f,\nabla f)$ in
  $L^2(\R^d)^{d+1}$.

  For arbitrary $f,h \in H^1$ take smooth approximants converging in $H^1$ and
  along a common subsequence Lebesgue-a.e. Lemma~\ref{L:spatial-physical-pair}
  transfers this convergence at both positions to $\nu$-a.e. convergence. Fatou
  and continuity of $N$ give $B_\Gamma(f^2,h^2) \le R_0N(f)^2N(h)^2 < \infty$,
  including the diagonal bound. At this stage Fatou supplies an inequality, not
  continuity.

  \medskip
  \noindent\textit{Step 3.~}We prove the four-linear estimate and continuity.
  Apply ordinary Cauchy--Schwarz under the positive measure $\nu$ to
  $f_1(x)f_3(x-y)$ and $f_2(x)f_4(x-y)$. The two squared integrals are finite by
  the last bound; their square root is at most $R_0\prod_jN(f_j)$. This
  proves~\eqref{E:spatial-four-linear}, absolute convergence and multilinearity.
  Telescoping a difference one argument at a time and using $N(f) \le
  \|f\|_{H^1}$ proves joint continuity. Thus, $A_{\mathrm{sp}}$ and the square
  pairings are continuous. PSD for smooth compact profiles passes to finite
  linear combinations of $H^1$ squares by this continuity, and its two-profile
  Cauchy--Schwarz gives \eqref{E:spatial-square-CS}. This completes the physical
  pairing argument beyond Fatou.

  \medskip
  \noindent\textit{Step 4.~}We identify the full ratio supremum and strict
  positivity.
  As in Lemma~\ref{L:compact-raw-profile}, nonzero $\Gamma$ gives positive
  mass to some origin-centred ball, and homogeneity rescales that ball to
  every positive radius.
  Every real unit $H^1$ profile has positive gradient energy by the Fourier
  argument above. The homogeneous bound gives $R_{\mathrm{sp}} \le R_0$, and
  smooth inclusion gives the reverse. Normalized smooth approximations preserve
  the ratio by continuity. Finally, for nonzero $f$, choose a bounded set $E$ of
  finite positive measure where $|f| \ge c > 0$. The correlation of $1_E$ is
  continuous and positive near zero by $L^2$ translation continuity. An
  origin-centered ball there has positive $\Gamma$ mass, so $A_{\mathrm{sp}}(f)
  \ge c^4B_\Gamma(1_E,1_E) > 0$. This uses a nonnegative integral, not an $H^1$
  claim for $1_E$. Finiteness was already proved. This proves
  Proposition~\ref{P:spatial-four-linear}.
\end{proof}

By Proposition~\ref{P:spatial-four-linear}, the ratio supremum is the same
over smooth and over $H^1$ profiles. We now express it through an energy and
restate the score bound in those terms.

\LeanMarker{P:spatial-energy-map}
\begin{proposition}[The spatial energy and static lower variational expression]
  \label{P:spatial-energy-map}
  Over all real unit $H^1$ profiles define the energy
  \[
    \mathcal{E}_{\mathrm{sp}}
    =\sup_{\|g\|_2=1}
    \left(A_{\mathrm{sp}}(g)-\frac{b(g)}2\right).
  \]
  Then
  \begin{equation}\label{E:spatial-energy-map}
    \mathcal{E}_{\mathrm{sp}}
    =\frac{2-\alpha}{2}\alpha^{\alpha/(2-\alpha)}
    R_{\mathrm{sp}}^{2/(2-\alpha)}
    \in(0,\infty).
  \end{equation}
  In particular,
  \begin{equation}\label{E:spatial-static-energy-score}
    S_1^{\mathrm{pos}}\ge M_a^{-1/k}F_0(\mathcal{E}_{\mathrm{sp}}),
    \qquad
    F_0(z)=k\left[
    \frac{(2z/(2-\alpha))^{(2-\alpha)/2}}{8q^{q/2}}
    \right]^{1/k},\quad z>0.
  \end{equation}
  The subscript in $F_0$ labels this scalar map; the covariance retains its
  original temporal exponent $a$.
\end{proposition}

\begin{proof}
  We optimize spatial dilation for a fixed profile, then take the supremum
  over all profiles. The resulting identity expresses the ratio constant in
  terms of the energy and lets us rewrite the static score bound.

  Physical change of variables, or smooth approximation and
  Proposition~\ref{P:spatial-four-linear}, extends \eqref{E:spatial-dilation} to
  real $H^1$ profiles. For fixed unit $g$, with $A,b > 0$, the maximum of
  $\ell^\alpha A-\ell^2b/2$ occurs at $\ell^{2-\alpha} = \alpha A/b$ and equals
  \[
    \frac{2-\alpha}{2}\alpha^{\alpha/(2-\alpha)}
    R(g)^{2/(2-\alpha)}.
  \]
  Differentiating the scalar expression gives $\alpha A\ell^{\alpha-1}-b\ell$.
  Its positive zero satisfies the stated equation. At that point,
  $b\ell^2=\alpha A\ell^\alpha$, so the value is $(1-\alpha/2)A\ell^\alpha$.
  Substitution of the critical scale gives the displayed power of $R(g)$.
  Dilation preserves the whole unit class and includes $\ell = 1$. Taking both
  numerical suprema, then profiles with ratio tending to $R_{\mathrm{sp}}$,
  proves~\eqref{E:spatial-energy-map}. No optimizer is assumed. Smooth ratios
  also approach $R_{\mathrm{sp}}$ by Proposition~\ref{P:spatial-four-linear}.
  Pass them through~\eqref{E:spatial-smooth-score} and substitute
  \[
    R_{\mathrm{sp}}\alpha^{\alpha/2}
    =\left(\frac{2\mathcal{E}_{\mathrm{sp}}}{2-\alpha}
    \right)^{(2-\alpha)/2}.
  \]
  This gives~\eqref{E:spatial-static-energy-score}. This proves
  Proposition~\ref{P:spatial-energy-map}.
\end{proof}

\subsection{The kernel row and full occupation energy}

Define, for $u \ne 0$,
\begin{equation}\label{E:occupation-damped-kernel}
  K_\star(u)=\frac1{\Gamma_E(a)}
  \int_0^\infty
  e^{-(\lambda+\lambda^2/(2D))|u|}
  \lambda^{a-1}\ud\lambda,
\end{equation}
and set $K_\star(0) = 0$ for time integration. The spatial estimates above
concern individual $H^1$ profiles. To apply them to time-dependent families,
we first compute the exact row bound for this time kernel.

\LeanMarker{L:occupation-kernel-row}
\begin{lemma}[Exact row maximum]\label{L:occupation-kernel-row}
  The kernel is even, positive off zero, locally integrable, decreasing on
  $(0,\infty)$, and $K_\star(u) \le |u|^{-a}$ off zero. Its exact row maximum on
  $[0,1]$ is
  \begin{equation}\label{E:occupation-kernel-row}
    \max_{s\in[0,1]}\int_0^1K_\star(s-t)\ud t
    =R_\star:=2\int_0^{1/2}K_\star(u)\ud u\in(0,\infty).
  \end{equation}
\end{lemma}

\begin{proof}
  Dropping the positive quadratic damping bounds the mixture by the Gamma
  integral $|u|^{-a}$. Positivity and monotonicity follow directly from the
  integrand, and $a < 1$ gives local integrability. Put $\Phi(t) =
  \int_0^tK_\star(u)\ud u$. It is concave because its integrand is decreasing on
  positive arguments. The row is $\Phi(s)+\Phi(1-s) \le 2\Phi(1/2)$, with
  equality at $s = 1/2$. The diagonal value changes no time integral. This
  proves Lemma~\ref{L:occupation-kernel-row}.
\end{proof}

Let $\mathcal{O}$ be the full class of strongly $H^1(\R^d;\R)$-measurable fields
$s\mapsto g_s$ on $(0,1)$ with $\|g_s\|_2 = 1$ a.e. and $T(g) =
\int_0^1b(g_s)\ud s < \infty$. No time smoothness, common compact support or
optimizer condition is imposed.

\LeanMarker{P:occupation-physical-class}
\begin{proposition}[The full physical occupation class]
  \label{P:occupation-physical-class}
  The pairings $(s,t)\mapsto B_\Gamma(g_s^2,g_t^2)$ are measurable,
  representative independent and equal to their physical integrals. With
  \[
    J(g)=\int_0^1\int_0^1
    K_\star(s-t)B_\Gamma(g_s^2,g_t^2)\ud s\ud t,
  \]
  every $g \in \mathcal{O}$ satisfies
  \begin{equation}\label{E:occupation-physical-bound}
    0\le J(g)\le R_\star R_{\mathrm{sp}}T(g)^{\alpha/2}<\infty.
  \end{equation}
  The energy over exactly this full class satisfies
  \begin{equation}\label{E:occupation-own-energy}
    0<\mathcal{E}_\star^{\mathcal{O}}
    :=\sup_{g\in\mathcal{O}}\big(J(g)-T(g)/2\big)
    \le R_\star^{2/(2-\alpha)}\mathcal{E}_{\mathrm{sp}}<\infty.
  \end{equation}
\end{proposition}

\begin{proof}
  We first identify measurable physical pairings for the full occupation
  class. We bound the time integral by the kernel row estimate and then
  compare its energy with the spatial energy, retaining the full class
  throughout.

  \medskip
  \noindent\textit{Step 1.~}We identify measurable representatives and their
  pairings.
  The map $s\mapsto b(g_s)$ is measurable by continuity of the squared gradient
  norm. Also $\int_0^1\|g_s\|_{H^1}^2\ud s = 1+T(g)$. Use a separably valued
  strong representative and replace exceptional null times by one fixed real
  unit compact bump. Joint measurable spatial representatives exist: approximate
  this Bochner $L^2$ field by simple $H^1$-valued functions with summable
  $L^2(0,1;H^1)$ errors, and choose measurable spatial versions for their
  countably many values. A subsequence has a measurable spacetime $L^2$ limit,
  and for a.e. $s$ its spatial $L^2$ class is $g_s$.
  Lemma~\ref{L:spatial-physical-pair} makes every pairing independent of the
  spatial versions. The continuous pairing map on the separable space $H^1\times
  H^1$ from Proposition~\ref{P:spatial-four-linear} consequently gives a
  measurable function of $(s,t)$ agreeing a.e. with the nonnegative physical
  integral.

  \medskip
  \noindent\textit{Step 2.~}We apply the spatial bound and the temporal row
  estimate.
  Set $F(s) = b(g_s)^{\alpha/4}$. Equations~\eqref{E:spatial-homogeneous-bound}
  and~\eqref{E:spatial-square-CS} give $B_\Gamma(g_s^2,g_t^2) \le
  R_{\mathrm{sp}}F(s)F(t)$. Since $0 < \alpha/2 < 1$, concavity on the unit time
  interval yields $\int_0^1F^2 = \int_0^1b(g_s)^{\alpha/2}\ud s \le
  T(g)^{\alpha/2}$. Now use $2F(s)F(t) \le F(s)^2+F(t)^2$ and the exact row
  bound \eqref{E:occupation-kernel-row}, whose column bound is the same because
  $K_\star$ is even, to obtain \eqref{E:occupation-physical-bound}. Positive
  Tonelli is valid before finiteness is known. Null-time modifications affect
  only a product-null set, so they cannot change $J$.

  \medskip
  \noindent\textit{Step 3.~}We bound the energy and show that it is positive.
  The scalar maximum of $R_\star R_{\mathrm{sp}}T^{\alpha/2}-T/2$ over $T \ge 0$
  is $R_\star^{2/(2-\alpha)}\mathcal{E}_{\mathrm{sp}}$, by the same
  differentiation as~\eqref{E:spatial-energy-map}. This bounds the occupation
  supremum from above. For positivity take a constant-in-time unit compact
  profile and its spatial dilation $U_\ell$. Its $J$ is $J_0\ell^\alpha
  A_{\mathrm{sp}}(g)$, where $J_0 = \int_0^1\int_0^1K_\star(s-t)\ud s\ud t \in
  (0,\infty)$, while $T = \ell^2b(g)$. For sufficiently small $\ell > 0$,
  $\alpha < 2$ makes $J-T/2 > 0$. Thus, the nonempty supremum is a positive
  finite real number. This is a comparison bound; no attaining profile or
  equality in it has been asserted. No spatial point trace or rough stochastic
  evolution was used. This proves Proposition~\ref{P:occupation-physical-class}.
\end{proof}

For zero covariance all spatial pairings vanish and set $R_{\mathrm{sp}} = 0$.
The spatial and occupation energies are zero suprema of $-b/2$ and $-T/2$: small
spatial dilations of one fixed unit compact profile make these energies tend to
zero, with the profile constant in time for the latter. No zero-gradient unit
state is used. The temporal identities remain valid, and the moment model is the
free case already treated.

\section{The compact mean and Cameron map on time blocks}
\label{S:actual-cameron-mean}

Retain the nonzero original covariance in Hypothesis~\ref{H:homog}, its
amplitude, $a = \alpha_0$, and the real phases $E = H^{1/2}\times H^{-1/2}$ and
$F = E\oplus\R^2$. Write $\Lambda = (1-\Delta)^{1/2}$. The form $B_\Gamma$ and
constant $R_{\mathrm{sp}}$ are those of
Section~\ref{S:spatial-occupation-foundation}. Let $H_{\mathrm{sp}}$ be the real
covariance-null completion of $C_c^\infty(\R^d)$ with inner product $B_\Gamma$,
and $H_t$ the analogous completion of $C_c^\infty(I)$, $I = (0,1)$, with time
kernel $|s-t|^{-a}$. Let $H = H_I$ be the completion of smooth tests supported
inside $I\times\R^d$. It embeds isometrically as a closed subspace of the given
$H_W$. We retain the given $W$ restricted to it.

Fix $R > 0$ and real even $\eta_R \in C_c^\infty(\R^d)$ with $0 \le \eta_R \le
1$, equal to one on the radius-$R$ ball and supported in the closed radius-$2R$
ball. The physical mask is $m_R = \eta_R^2$. Every operator path below takes its
mask argument equal to $m_R$: at this fixed radius the base profile $\zeta(y) =
m_R(Ry)$ meets the earlier hypotheses. The root $\eta_R$ factors state tests; it
does not multiply an original noise class. All compactness conclusions here fix
$R$.

\subsection{Fractional products and the compact spatial map}

The four-linear estimate~\eqref{E:spatial-four-linear} needs $H^1$ factors.
The first component of a state in $E$ lies only in $H^{1/2}$, so we first
extend that estimate to fractional Sobolev factors.

\begin{lemma}[Fractional covariance products]
  \label{L:cameron-fractional-products}
  For $\alpha/4 < s < 1/2$, products of real $H^s$ functions define elements of
  $H_{\mathrm{sp}}$, with
  \begin{equation}\label{E:cameron-fractional-product}
    \|fg\|_{H_{\mathrm{sp}}}
    \le C_s^2\sqrt{R_{\mathrm{sp}}}\|f\|_{H^s}\|g\|_{H^s}.
  \end{equation}
  Their Hilbert pairings equal their physical integrals. In particular,
  \begin{equation}\label{E:cameron-fractional-absolute}
    \int|f_1(x)f_2(x)f_3(x-y)f_4(x-y)|\ud x\,\Gamma(\ud y)
    \le C_s^4R_{\mathrm{sp}}\prod_{j=1}^4\|f_j\|_{H^s}.
  \end{equation}
\end{lemma}

\begin{proof}
  We first define covariance products for $H^1$ factors. We then decompose
  each fractional Sobolev factor into real Fourier bands and sum the absolute
  four-linear bounds, identifying the completed products with their physical
  pairings.

  \medskip
  \noindent\textit{Step 1.~}We obtain the product map on $H^1$ factors.
  First approximate real $H^1$ factors by compact smooth ones. The absolute
  four-linear bound~\eqref{E:spatial-four-linear} and covariance PSD make their
  products Cauchy in $H_{\mathrm{sp}}$, with norm at most
  $\sqrt{R_{\mathrm{sp}}}N(f)N(g)$. Expanding product differences proves
  independence of approximation; physical continuity identifies all pairings.
  PSD for products also follows from $fg = ((f+g)^2-(f-g)^2)/4$.

  \medskip
  \noindent\textit{Step 2.~}We sum the band bounds for fractional Sobolev
  factors.
  Use real radial Fourier bands $D_0 = 1_{\{|\xi|<1\}}$ and $D_j =
  1_{\{2^{j-1}\le|\xi|<2^j\}}$, $j \ge 1$. For $f \in H^s$ put $f_j = D_jf$.
  Each block is real and in $H^1$. Frequency support, Parseval and weighted
  Cauchy--Schwarz give
  \begin{equation}\label{E:cameron-dyadic-product}
    N(f_j)\le2^{j\alpha/4}\|f_j\|_2,\qquad
    \sum_{j\ge0}N(f_j)\le C_s\|f\|_{H^s}.
  \end{equation}
  On band $j$, the gradient norm is at most $2^j\|f_j\|_2$. Substituting into
  $N$ gives the first inequality. Write $2^{j\alpha/4}\|f_j\|_2$ as
  $2^{-j(s-\alpha/4)}(2^{js}\|f_j\|_2)$ and apply Cauchy--Schwarz to the two
  sequences.
  Indeed, $\sum_j2^{-2j(s-\alpha/4)} < \infty$, while $\sum_j2^{2js}\|f_j\|_2^2
  \le 2^{2s}\|f\|_{H^s}^2$. The low block has factor one. No $L^p$ multiplier
  result is used. Also $\sum_j\|f_j\|_2 < \infty$ since $s > 0$. Thus,
  $\sum_j|f_j|$ lies in $L^2$ by Minkowski and monotone convergence, and the
  series is absolutely convergent a.e. with sum $f$.

  \medskip
  \noindent\textit{Step 3.~}We identify the physical and Hilbert-space product
  series.
  Apply the absolute $H^1$ four-linear estimate to each band quadruple. The sum
  of its absolute integrals is at most
  $R_{\mathrm{sp}}\prod_{\ell=1}^4\sum_jN((f_\ell)_j) < \infty$. Tonelli and the
  null-set transfer at both spatial positions in
  Lemma~\ref{L:spatial-physical-pair} identify the physical product series,
  proving~\eqref{E:cameron-fractional-absolute}. The double series
  $\sum_{j,k}f_jg_k$ is also absolutely convergent in $H_{\mathrm{sp}}$ by the
  $H^1$ product bound. Its partial sums are products of $H^1$ functions. The
  absolute quadruple calculation identifies their limiting physical pairings and
  gives~\eqref{E:cameron-fractional-product}. Compact smooth functions are dense
  in $H^s$: approximate in Fourier space by real Schwartz functions, then cut
  each off compactly in $H^1$ and hence in $H^s$. This also proves uniqueness of
  the bilinear extension. This proves Lemma~\ref{L:cameron-fractional-products}.
\end{proof}

With these products in $H_{\mathrm{sp}}$, we can extend the spatial mean map
from smooth tests to all of $H_{\mathrm{sp}}$.

\begin{proposition}[Compact spatial mean map]\label{P:cameron-spatial-map}
  The smooth-test map
  \[
    T_{\mathrm{sp}}q=iB_F(m_R(\Gamma*q))
  \]
  extends to a compact linear map $T_{\mathrm{sp}}:H_{\mathrm{sp}} \to
  \mathcal{K}(F)$, where $\mathcal{K}(F)$ has the operator norm. For $z =
  (z_1,z_2,c_1,c_2)$ and $w = (w_1,w_2,d_1,d_2)$ set $w' = \Lambda^{-1}w_2$. Its
  exact row identity is
  \begin{equation}\label{E:cameron-spatial-row}
    \langle T_{\mathrm{sp}}qz,w\rangle_F
    =\langle q,k(z,w)\rangle_{H_{\mathrm{sp}}},
    \qquad
    k(z,w)=(\eta_Rz_1+c_1\eta_R)(\eta_Rw').
  \end{equation}
\end{proposition}

\begin{proof}
  We identify the spatial row and use compactness of localized fractional
  factors to bound it uniformly. We first extend the map with compact operator
  values, then approximate the map itself by finite-dimensional spatial
  projections.

  \medskip
  \noindent\textit{Step 1.~}We identify the physical and covariance row.
  Fix $s$ with $\alpha/4 < s < 1/2$; this interval is nonempty because the
  original threshold implies $\alpha < 2$. For smooth compact $q$, $\Gamma*q$ is
  locally smooth by differentiation under the compactly restricted Radon
  integral. The masked potential is smooth compact. The velocity pairing in $E$
  has weight $\Lambda^{-1}$; hence physical Fubini
  gives~\eqref{E:cameron-spatial-row} on smooth states. Evenness of $\Gamma$
  changes $B_\Gamma(k,q)$ to $B_\Gamma(q,k)$. For general states the localized
  products are $L^1$, and $\Gamma*q$ is bounded on their compact support. $L^2$
  state approximation passes the physical integral to the limit, while
  Lemma~\ref{L:cameron-fractional-products} passes the $H_{\mathrm{sp}}$
  pairing. This proves the row.

  \medskip
  \noindent\textit{Step 2.~}We obtain compact families of row tests.
  The velocity test satisfies $w_2\in H^{-1/2}$ and $w'=\Lambda^{-1}w_2\in
  H^{1/2}$, so it has the same localization regularity as the displacement
  test.
  Multiplication by $\eta_R$ maps bounded $H^{1/2}$ sets to precompact $H^s$
  sets by Proposition~\ref{P:cutoff-singular-values}, with gain $1/2-s > 0$.
  This applies to $z_1$ and $w'$. The $c_1\eta_R$ family is finite dimensional
  and compact. The closures of both factor families for unit $z,w$ are compact,
  so their products have compact closure $\mathcal{L}$ in $H_{\mathrm{sp}}$ by
  the fractional product bound. The row identity gives $\|T_{\mathrm{sp}}q\| \le
  \sup_{k\in\mathcal{L}}\|k\|\|q\|$. Null covariance classes therefore map to
  zero.

  \medskip
  \noindent\textit{Step 3.~}We extend the map with compact operator values.
  For smooth $q$, multiplication by its potential is compact $H^{1/2} \to
  H^{-1/2}$ by the same whole-space cutoff result with gain one. Its extra $c_1$
  column is finite rank, and the instantaneous $c_2$ column is zero. Thus, the
  smooth images lie in $\mathcal{K}(F)$, and so does the bounded extension by
  operator-norm closure.

  \medskip
  \noindent\textit{Step 4.~}We prove compactness of the map itself.
  On each fixed compact support, local Radon finiteness bounds the covariance
  norm by a constant times the supremum norm. Countable uniformly dense
  families of smooth tests on exhausting compact supports therefore give a
  countable dense set in $H_{\mathrm{sp}}$.
  Choose a complete smooth-test orthonormal basis of $H_{\mathrm{sp}}$ by
  countable compact-support approximation and Gram--Schmidt, omitting null
  classes. Let $Q_n$ project onto its first $n$ terms. The row identity yields
  \[
    \|T_{\mathrm{sp}}(I-Q_n)\|
    \le\sup_{k\in\mathcal{L}}\|(I-Q_n)k\|_{H_{\mathrm{sp}}}
    \longrightarrow0.
  \]
  Strong projection convergence is uniform on compact $\mathcal{L}$ by a
  finite-net argument and $\|I-Q_n\| \le 1$. Each $T_{\mathrm{sp}}Q_n$ has
  finite-dimensional range; their operator-norm limit is a compact map. This is
  compactness of the map, in addition to compact values. This proves
  Proposition~\ref{P:cameron-spatial-map}.
\end{proof}

\subsection{The supported tensor and temporal mean}

Proposition~\ref{P:cameron-spatial-map} treats the spatial variable alone. A
control in $H$ depends on time and space, so we first split $H$ into a time
factor and a spatial factor.

\begin{lemma}[Supported original tensor space]\label{L:cameron-supported-tensor}
  Canonically $H = H_t\otimes H_{\mathrm{sp}}$ as real Hilbert spaces. It is
  separable and admits a complete orthonormal basis $e_j$ with smooth test
  representatives supported inside $I\times\R^d$. Projection $P_J$ onto its
  first $J$ elements is also an orthogonal projection on $H_W$ by the same
  finite formula.
\end{lemma}

\begin{proof}
  On separated smooth tests, absolute Fubini makes the original inner product
  exactly the tensor product after both null quotients. Finite sums are dense:
  approximate a smooth compact spacetime test uniformly by smooth time
  partitions times its spatial slices, on one common compact support inside
  $I\times\R^d$. A uniform error $\epsilon$ on spatial compact $A$ has squared
  original covariance norm at most
  \[
    \epsilon^2|A|\Gamma(A-A)
    \int_I\int_I|s-t|^{-a}\ud s\ud t<\infty.
  \]
  This also gives countable smooth-test density on exhausting compact sets, both
  here and in $H_{\mathrm{sp}}$. Gram--Schmidt skipping zeros preserves smooth
  compact support. The isometric embedding into $H_W$ has closed image by
  completeness; its finite orthonormal projection formula remains orthogonal
  there. It is not multiplication by a time indicator and gives no independence
  of disjoint time intervals. This proves
  Lemma~\ref{L:cameron-supported-tensor}.
\end{proof}

It remains to treat the time factor $H_t$, whose kernel is $|s-t|^{-a}$.
Proposition~\ref{P:cameron-full-mean} then combines the two maps.

\begin{proposition}[The compact temporal mean, without trace class]
  \label{P:cameron-temporal-map}
  Let $K$ on $L^2(I)$ have kernel $|s-t|^{-a}$. It is compact, self-adjoint,
  strictly positive and injective. The map $\iota:f\mapsto K^{1/2}f$ extends
  isometrically from $H_t$ onto $L^2(I)$, and the temporal mean
  \begin{equation}\label{E:cameron-temporal-map}
    C_t=K^{1/2}\iota:H_t\longrightarrow L^r(I)
  \end{equation}
  is bounded and compact for every $2 < r < 2/a$. The displayed formula first
  has its $L^2(I)$ meaning; the same completed map is then shown to take values
  compactly in $L^r(I)$. The operator $K$ is not trace class and $C_t:H_t \to
  L^2(I)$ is not Hilbert--Schmidt.
\end{proposition}

\begin{proof}
  We prove compactness and strict positivity of the temporal operator before
  identifying the covariance completion. We then obtain the stronger $L^r$
  bound by duality, use interpolation for compactness, and test cell
  indicators to rule out trace class.

  \medskip
  \noindent\textit{Step 1.~}We prove positivity, injectivity and compactness
  of the kernel operator.
  Schur bounds $K$ by its finite row and column integrals. The Gamma
  representation of $|s-t|^{-a}$ and zero extension of $f \in L^2(I)$ give
  strict positivity directly: the kernel $e^{-\lambda|v|}$ has quadratic form
  \[
    \int_\R\frac{2\lambda}{\lambda^2+\xi^2}
    |\widehat{f}(\xi)|^2\ud\xi>0
    \quad(f\ne0,\ \lambda>0)
  \]
  in the unitary Fourier convention. Absolute Fubini is justified by Schur
  applied to $|f(s)f(t)|$. Integrating with positive weight
  $\lambda^{a-1}/\Gamma_E(a)$ proves positivity and injectivity. Symmetry gives
  self-adjointness. The deterministic truncated kernel $K_\epsilon =
  K1_{\{|s-t|\ge\epsilon\}}$ is square integrable, so finite separated-kernel
  approximations make it compact. The Schur norm of $K-K_\epsilon$ is at most
  $2\epsilon^{1-a}/(1-a)$, proving compactness of $K$. These truncations are not
  used as Gaussian covariances. Its positive square root is compact by the
  compact spectral theorem.

  \medskip
  \noindent\textit{Step 2.~}We identify the covariance completion with the
  square-root range.
  On tests $\|K^{1/2}f\|_2 = \|f\|_{H_t}$. Smooth-test $L^2$ density and $\ker K
  = 0$ imply that its range is dense, so the isometry extends onto $L^2(I)$. On
  tests $C_tf = Kf$, and~\eqref{E:cameron-temporal-map} first gives compactness
  into $L^2$.

  \medskip
  \noindent\textit{Step 3.~}We bound the temporal mean in $L^r$ and obtain
  compactness there.
  The restriction $r<2/a$ gives $ar/2<1$, ensuring integrability of the kernel
  power at zero.
  Extend $K_a(v) = |v|^{-a}1_{\{|v|\le1\}}$ to the line. It lies in $L^{r/2}$.
  Scalar Young convolution gives
  \[
    \|g\|_{H_t}^2\le\|K_a\|_{r/2}\|g\|_{r'}^2,\qquad
    r'=r/(r-1).
  \]
  To verify the needed Young inequality, take $u = r/2$, $v = r'$, so $u,v < r$
  and $1/u+1/v = 1+1/r$. Split the absolute convolution integrand into
  $(|f|^u|h|^v)^{1/r}$, $|f|^{1-u/r}$ and $|h|^{1-v/r}$. H\"older with exponents
  $r$, $(1/u-1/r)^{-1}$ and $(1/v-1/r)^{-1}$, followed by integration and
  Tonelli, gives $\|f*h\|_r \le \|f\|_u\|h\|_v$. Covariance Cauchy--Schwarz and
  $L^r$ duality now imply $\|Kf\|_r \le C_{a,r}\|f\|_{H_t}$ on tests. Their $Kf$
  are bounded by the integrable-kernel row, making duality legitimate before
  completion. Density extends the estimate, with $L^2$ value equal
  to~\eqref{E:cameron-temporal-map}. Choose $r < r_2 < 2/a$. Compactness into
  $L^2$ and boundedness into $L^{r_2}$ give compactness into $L^r$ by $\|u\|_r
  \le \|u\|_2^\tau\|u\|_{r_2}^{1-\tau}$ with $1/r = \tau/2+(1-\tau)/r_2$.

  \medskip
  \noindent\textit{Step 4.~}We test the trace by orthonormal cell indicators.
  Finally, the orthonormal cell indicators $N^{1/2}1_{(j/N,(j+1)/N)}$, $0 \le j
  < N$, have $K$ quadratic forms summing to $2N^a/((1-a)(2-a)) \to \infty$.
  On a cell of length $1/N$, the double kernel integral is
  $2N^{a-2}/((1-a)(2-a))$. The squared normalization of its indicator
  contributes $N$, and summing over the $N$ cells gives the stated $N^a$
  divergence. For a positive trace-class operator every such finite
  orthonormal sum is bounded by its trace, so this divergence excludes trace
  class.
  Thus, positive $K$ has infinite trace, and $K^{1/2}\iota$ is not
  Hilbert--Schmidt into $L^2$. No instantaneous $L^2$-valued Gaussian has been
  introduced. This proves Proposition~\ref{P:cameron-temporal-map}.
\end{proof}

No-go~\ref{P:nogo-27} explains the upper endpoint of this compactness
range: temporal concentration makes the mean map unbounded above
$2/\alpha_0$ and noncompact at $2/\alpha_0$. The endpoint witness does not
refute boundedness there.

\begin{proposition}[The compact time-space mean]\label{P:cameron-full-mean}
  Fix one $2 < r < 2/a$ throughout this proposition. For finite tensors with
  orthonormal spatial labels define
  \[
    \mathcal{U}\left(\sum_j f_j\otimes q_j\right)(s)
    =\sum_j(Kf_j)(s)q_j.
  \]
  It extends to a bounded map $H \to L^r(I;H_{\mathrm{sp}})$ with dimension-free
  bound. The map
  \begin{equation}\label{E:cameron-full-mean}
    T\Phi(s)=T_{\mathrm{sp}}(\mathcal{U}\Phi(s))
    :H\longrightarrow L^r(I;\mathcal{K}(F))
  \end{equation}
  is compact. For smooth supported $\Phi$, it is the multiplier
  $iB_F(m_Rh_\Phi(s))$, where $h_\Phi = |\cdot|^{-a}*(\Gamma*\Phi)$. For the
  supported projections of Lemma~\ref{L:cameron-supported-tensor},
  \begin{equation}\label{E:cameron-compact-tail}
    \epsilon_J:=\|T(I-P_J)\|_{H\to L^r(\mathcal{K}(F))}
    \longrightarrow0.
  \end{equation}
\end{proposition}

\begin{proof}
  We first prove a dimension-free temporal estimate on finite tensors and
  identify its physical mean. We then combine spatial finite-rank
  approximation with temporal compactness, and use compactness of the combined
  map to control supported projection tails.

  \medskip
  \noindent\textit{Step 1.~}We bound and extend the time-only tensor map.
  For a finite-dimensional spatial field $g$, positivity of the time kernel and
  scalar Young give
  \[
    \|g\|_{H_t\otimes H_{\mathrm{sp}}}^2
    \le\int_I\int_I|s-t|^{-a}\|g_s\|\|g_t\|\ud s\ud t
    \le C_{a,r}^2\|g\|_{L^{r'}(H_{\mathrm{sp}})}^2.
  \]
  Hilbert Cauchy--Schwarz and finite-dimensional $L^r$ duality give
  $\|\mathcal{U}\Phi\|_{L^r(H_{\mathrm{sp}})} \le C_{a,r}\|\Phi\|_H$. The
  duality can be checked by the pointwise direction
  $\|\mathcal{U}\Phi\|^{r-2}\mathcal{U}\Phi$ and smooth coefficient
  approximation in $L^{r'}$, controlled by the displayed estimate. No dimension
  factor occurs. Tensor density and completeness extend this to a strongly
  measurable Bochner field. This time-only map still contains spatial vectors in
  the covariance space. Spatial covariance is applied once, by
  $T_{\mathrm{sp}}$.

  \medskip
  \noindent\textit{Step 2.~}We identify the completed map on smooth controls.
  The multiplier estimate used here is $\|B_F(b)\|\le\|b\|_\infty+\|b\|_2$
  from Lemma~\ref{L:constant-source-response}; on the fixed support of $m_R$,
  the $L^2$ error is at most the square root of its volume times the supremum
  error.
  The physical mean identity on separated tests follows
  from~\eqref{E:cameron-spatial-row}. For a smooth spacetime test use the
  separated compact-support approximations in the tensor proof. Their uniform
  errors give uniformly small masked mean errors, using the finite time row and
  locally finite spatial covariance on the relevant compact differences. The
  multiplier bound by $L^\infty+L^2$ then gives uniform operator convergence,
  identifying the $L^r$ mean in~\eqref{E:cameron-full-mean}. General $H$
  controls denote completed operator fields, without a pointwise-potential
  assertion.

  \medskip
  \noindent\textit{Step 3.~}We combine temporal and spatial compactness.
  Let $Q_n$ be the spatial projections in the compact-map proof. The maps
  $T_n\Phi = T_{\mathrm{sp}}Q_n\mathcal{U}\Phi$ approximate $T$ in operator norm
  by the dimension-free bound and $\|T_{\mathrm{sp}}(I-Q_n)\| \to 0$. For fixed
  $n$, each spatial coordinate extraction $H \to H_t$ has norm at most one and
  its image under $C_t$ is precompact in $L^r$. Finite combination with the
  compact operators $T_{\mathrm{sp}}q_j$ makes $T_n$ compact. Hence, $T$ is
  compact. Scalar time compactness alone tensored with an infinite identity
  would not prove this.

  \medskip
  \noindent\textit{Step 4.~}We obtain the supported projection-tail bound.
  To prove~\eqref{E:cameron-compact-tail}, if it failed choose unit $v_J$ on a
  subsequence with $\|T(I-P_J)v_J\|$ bounded below. The vectors $(I-P_J)v_J$
  tend weakly to zero, since their pairings with fixed $h$ are bounded by
  $\|(I-P_J)h\|$. Compactness gives a norm-convergent image subsequence with
  weak limit zero, a contradiction. This proves
  Proposition~\ref{P:cameron-full-mean}.
\end{proof}

The radius remains fixed in this approximation. No-go~\ref{P:nogo-25}
shows that translated white-noise controls prevent these fixed smooth
projections from giving tails uniform over growing radii.

\subsection{The noise residual and its Cameron operator}

Choose integer $B \ge 1$ and $0 < \delta < \min(1,1/B)$. Put $N =
\lfloor1/(\delta B)\rfloor$, $s_m = m\delta B$, $\ell_m = B$ for $m < N$ and
$\ell_N = \delta^{-1}-NB \in [0,B)$. The intervals $I_m =
[s_m,s_m+\delta\ell_m]$ partition $[0,1]$ up to endpoints. With $Q_J = I-P_J$ on
$H_W$ and the explicit pullback~\eqref{E:raw-block-pullback}, define
\[
  Y_{m,J}(f)=W(Q_JT_{\delta,m}f),\qquad
  g_\delta=\delta^{(1-a)/2}.
\]
Fix parameters
\begin{equation}\label{E:cameron-block-parameters}
  \frac12<\theta<\min\left(1-\frac a2,\frac{3-a-\alpha}{2}\right),
  \qquad \frac1{1-\theta}<r<\frac2a,\qquad
  \frac1\theta\le p<2<r .
\end{equation}
These intervals are nonempty.

\begin{lemma}[Supported residual paths]\label{L:cameron-actual-residual}
  Gaussian transfer with mask $m_R$ and its $F$ extension constructs the
  centered Gaussian paths
  \[
    X_{m,J}(u)=g_\delta A_{m_R}(Y_{m,J})(\min(u,\ell_m)),
    \qquad 0\le u\le B.
  \]
  Their finite tuple lies in the separable space $\mathscr{B}_0^{N+1}$, where
  $\mathscr{B}_0 = h_0^\theta([0,B];\mathcal{K}(F))$. Every continuous linear
  statistic $\ell(X)$ is $W(k_\ell)$ for some $k_\ell \in H_I$. An empty last
  interval gives zero.
\end{lemma}

\begin{proof}
  We construct the stopped residual paths by Gaussian transfer. We then
  identify the supported slow-time tests in every observed row, so that each
  continuous path statistic lies in the original first chaos over $H_I$.

  \medskip
  \noindent\textit{Step 1.~}We construct the stopped paths in the separable
  little Hölder space.
  For $Y=Y_{m,J}$, the notation $A_{m_R}(Y)$ means $iM_YL_\psi$
  from~\eqref{E:constant-phase-driver}, where $M_Y$ is the transferred $E$
  path. The two localizations are~\eqref{E:constant-phase-localization},
  passed from the common deterministic coefficients to this path by bounded
  multiplication.
  Isometry of the pullback and orthogonality of $Q_J$ give
  $\operatorname{Var}Y_{m,J}(f) \le \|f\|_{H_W}^2$.
  Proposition~\ref{P:source-gaussian-transfer} therefore applies. Its finite
  sums retain both deterministic localizations, and their limits retain the true
  $F$ extension and constant columns. Choose $\theta < \beta < \beta_2$ below
  the upper threshold in~\eqref{E:cameron-block-parameters}; read the
  operator-path construction at exponent $\beta$ with auxiliary $\beta_2$.
  Constant stopping preserves compact values and $C^\beta$ bounds, hence gives
  paths in the little $C^\theta$ closure by polygonal approximation. Its
  separability is Lemma~\ref{L:raw-block-inputs}.

  \medskip
  \noindent\textit{Step 2.~}We identify the supported first-chaos tests for
  the path statistics.
  The projection $P_J$ has range in $H_I$, so subtracting it from a vector
  already in $H_I$ keeps that vector in $H_I$. Consequently, once an observed
  slow-time row is shown to belong to $H_I$, its residual row has the same
  support-space membership. This is a Hilbert-subspace statement, not
  multiplication of the original noise by a time indicator.
  As in Lemma~\ref{L:raw-block-inputs}, every observed row's slow pullback
  belongs to $H_I$: its time support is in $[0,1]$, and approximating its
  indicator inside $I$ converges in original covariance norm by dominated
  convergence against the integrable time kernel and locally finite $\Gamma$ on
  the common compact support. This includes endpoints. Schwartz $E$ states and
  both true constant columns supply the identities, then boundedness extends
  them to $F$. Let $P_H$ be the orthogonal projection onto closed $H_I$.
  Transfer for the noise $f\mapsto W(P_HQ_JT_{\delta,m}f)$ leaves each observed
  row unchanged, since its residual slow test already lies in $H_I$. Countable
  states and rational times identify the same operator paths. This
  projected-transfer series has coefficient tests in $H_I$. Applying $\ell$
  gives an $L^2$ limit in the closed isonormal image $W(H_I)$, proving the
  assertion. Original global-basis tests were not assumed individually
  supported. No time-support independence or Hilbert-valued noise is used. This
  proves Lemma~\ref{L:cameron-actual-residual}.
\end{proof}

\begin{proposition}[Cameron covariance and averaged bound]
  \label{P:cameron-actual-embedding}
  The covariance map $L_{\delta,J}\Phi = \E[XW(\Phi)]$ is
  \begin{equation}\label{E:cameron-actual-covariance}
    (L_{\delta,J}\Phi)_m(u)=\sqrt\delta
    \int_0^{\min(u,\ell_m)}
    S_F(-v)T(Q_J\Phi)(s_m+\delta v)S_F(v)\ud v ,
    \qquad \Phi\in H.
  \end{equation}
  This is a compact-operator Bochner integral. Define the exact norm and
  constants
  \begin{equation}\label{E:cameron-averaged-norm}
    \|x\|_{\delta,p}
    =\left(\delta\sum_{m=0}^N\|x_m\|_{C^\theta}^p\right)^{1/p},
    \quad G_B=\sup_{|v|\le B}\|S_F(v)\|,
    \quad C_B=G_B^2[B^{1-1/r}+B^{1-1/r-\theta}].
  \end{equation}
  Then the Gaussian Cameron embedding satisfies
  \begin{equation}\label{E:cameron-averaged-bound}
    \|L_{\delta,J}\Phi\|_{\delta,p}
    \le\sqrt\delta\,C_B(2/B)^{1/p-1/r}
    \epsilon_J\|\Phi\|_H .
  \end{equation}
  Convergence $\epsilon_J \to 0$ fixes $R$ and is independent of $\delta$.
\end{proposition}

\begin{proof}
  We compute the covariance direction first on smooth supported controls and
  justify its compact-operator integral. We then extend to all controls,
  identify the Cameron preimage norm, and sum the local primitive estimates
  with the exact time weights.

  \medskip
  \noindent\textit{Step 1.~}We compute the covariance direction and its factor
  $\sqrt\delta$.
  Insert the wave-row tests of~\eqref{E:source-transferred-row} in the scalar
  covariance below; Proposition~\ref{P:cameron-full-mean} identifies the
  resulting masked mean with $T(Q_J\Phi)$, and the row calculation of
  Lemma~\ref{L:raw-block-inputs} gives its conjugated primitive on $F$.
  For smooth supported $\Phi$, $Q_J\Phi$ is another such test. The explicit
  pullback covariance~\eqref{E:raw-block-mean-jacobian} gives
  \[
    \E[Y_{m,J}(f)W(\Phi)]
    =\delta^{a/2}\int f(v,x)
    h_{Q_J\Phi}(s_m+\delta v,x)\ud x\ud v .
  \]
  Absolute Fubini also holds for signed compact tests. Multiplying by $g_\delta$
  gives $\sqrt\delta$. Explicit rows and~\eqref{E:cameron-spatial-row} give
  \eqref{E:cameron-actual-covariance}. The constant-velocity input acquires
  displacement $v$ under $S_F(v)$; the constant displacement remains one. Both
  columns are included.

  \medskip
  \noindent\textit{Step 2.~}We justify the compact-operator Bochner integral.
  Both $S_F$ and its adjoint are strongly continuous on Hilbert $F$. In
  Bessel-unitary $E$ coordinates the Fourier matrix is
  \[
    \begin{pmatrix}
      \cos(v|\xi|)&\langle\xi\rangle\sin(v|\xi|)/|\xi|\\
      -|\xi|\sin(v|\xi|)/\langle\xi\rangle&\cos(v|\xi|)
    \end{pmatrix}.
  \]
  Its continuous zero quotient and common bound $C(1+B)$ apply also to its
  adjoint; dominated convergence proves the assertion. The constant block is
  finite dimensional. For compact $A$, conjugation $S_F(-v)AS_F(v)$ is norm
  continuous: strong convergence is uniform on compact images, and the right
  factor follows by adjoints. Joint norm continuity in $(v,A)$ and uniform
  bounds give Bochner measurability for the conjugated strongly measurable
  compact-operator field. No norm continuity of the free group itself is
  assumed.

  \medskip
  \noindent\textit{Step 3.~}We bound each stopped primitive and extend the
  covariance identity.
  Write $M = TQ_J\Phi$. H\"older and the primitive increment estimate give
  \begin{equation}\label{E:cameron-local-primitive-bound}
    \|(L_{\delta,J}\Phi)_m\|_{C^\theta}
    \le\sqrt\delta\,C_B\delta^{-1/r}\|M\|_{L^r(I_m)}.
  \end{equation}
  Hölder in the fast variable bounds an increment by $G_B^2$ times its length
  to the power $1-1/r$, multiplied by the pulled-back $L^r$ norm. The latter
  is $\delta^{-1/r}\|M\|_{L^r(I_m)}$. For the supremum use length at most $B$;
  for the $\theta$ seminorm divide by the increment length to the power
  $\theta$. The two remaining powers of $B$ are exactly those in $C_B$.
  Its primitive exponent is $1-1/r > \theta$. Simple compact-operator integrand
  approximation in $L^r$ puts the primitive in $\mathscr{B}_0$; constant
  stopping does not increase the H\"older norm. Affine pullback preserves null
  sets. For fixed $\delta,J$ the covariance map is bounded by Bochner
  Cauchy--Schwarz and finite tuple $L^2$ moments. The integral is bounded by the
  displayed estimate, so equality on smooth controls extends by density to all
  $H$. Dense states and rational times make it a full operator-path identity.

  \medskip
  \noindent\textit{Step 4.~}We identify the Cameron embedding norm.
  Let $H_{\mathrm{seen}}$ be the closed span of the $k_\ell$ from
  Lemma~\ref{L:cameron-actual-residual}. Then $\ell(L_{\delta,J}\Phi) = \langle
  k_\ell,\Phi\rangle_H$ and $\ker L_{\delta,J} = H_{\mathrm{seen}}^\perp$. The
  first Wiener chaos definition of the Cameron space is therefore
  $L_{\delta,J}(H_{\mathrm{seen}})$ with the minimum $H$ preimage norm. Extra
  directions are removed by orthogonal projection, without an inverse
  covariance. A bound in the original $H$ norm is thus the embedding bound.

  \medskip
  \noindent\textit{Step 5.~}We combine the block estimates using the true
  partition.
  Put $t_m = \|M\|_{L^r(I_m)} \ge 0$ and $n = N+1 \ge 2$. Finite power mean at
  exponent $r/p > 1$, applied to $t_m^p$, gives
  \[
    \left(\sum_m t_m^p\right)^{1/p}
    \le n^{1/p-1/r}\left(\sum_m t_m^r\right)^{1/r}.
  \]
  The true partition gives $\sum_m t_m^r = \|M\|_{L^r(I)}^r$, including an empty
  last interval, and $n \le 2/(\delta B)$. Hence, the exact $\delta$-weighted
  norm gives
  \[
    \|L_{\delta,J}\Phi\|_{\delta,p}
    \le\sqrt\delta\,C_B(\delta n)^{1/p-1/r}\|M\|_{L^r}
    \le\sqrt\delta\,C_B(2/B)^{1/p-1/r}
    \epsilon_J\|\Phi\|_H.
  \]
  The $B$ factor remains explicit; covariance and mask occur once. This proves
  Proposition~\ref{P:cameron-actual-embedding}.
\end{proof}

These conclusions fix the radius and make no uniform-radius projection tail
assertion. No concentration estimate or original moment upper bound follows from
this section alone.

\section{Young perturbations retaining the deterministic endpoint}
\label{S:young-upper-foundation}

Let $\mathcal{X}$ be any nonzero real Banach space, with its induced operator
norm on $\mathcal{B}(\mathcal{X})$. Fix $B \ge 1$, $1/2 < \theta < 1$ and $\nu =
1/\theta < 2$. All operator paths start at zero and have norm $\|D\|_\theta =
\sup_t\|D(t)\|+[D]_\theta$ on $[0,B]$. Let $U_D = I+\int\ud D\,U_D$ be the left
Young evolution. The construction in Proposition~\ref{P:constant-young-linear}
applies: its proof needs completeness and operator multiplication of norm one,
with no Hilbert-space assumption. The primary estimate is
\cite{friz.hairer:20:course}*{Section~1.5, p.~13; Section~4.1, (4.1)--(4.3),
pp.~61--62}.

\begin{lemma}[Refined growth and uniform inverse control]
  \label{L:young-upper-growth}
  With dimension-free constants,
  \begin{equation}\label{E:young-upper-refined-growth}
    \sup_{t\le B}\|U_D(t)\|
    \le\exp\{c_{\theta,B}(d+d^\nu)\},\qquad d=\|D\|_\theta .
  \end{equation}
  For $\|H\|_\theta \le K$, $K \ge 1$,
  \begin{equation}\label{E:young-upper-inverse-control}
    \|U_H\|_{C^\theta}+\|U_H^{-1}\|_{C^\theta}
    \le L_{\theta,B,K}<\infty .
  \end{equation}
\end{lemma}

\begin{proof}
  We separate small drivers, where one contraction gives a linear logarithmic
  bound, from larger drivers, where a partition gives the power bound. We then
  use a fixed partition on a bounded driver ball to control both the evolution
  and its inverse in the Hölder norm.

  \medskip
  \noindent\textit{Step 1.~}We prove the refined growth bound in the two
  driver-size regimes.
  The derivative of $x\mapsto-\log(1-x)$ is at most two on $[0,1/4]$, giving
  the linear logarithmic bound in the small-driver case below.
  Let $c_\theta = 2\max(1,C_\theta)$ for the primary Young constant. The
  complete scaled norm $N_\ell(V) = \sup\|V\|+\ell^\theta[V]_\theta$ satisfies
  $N_\ell(\int_a^\cdot\ud D\,V) \le c_\theta[D]_\theta\ell^\theta N_\ell(V)$, as
  in~\eqref{E:constant-young-local}. If $c_\theta dB^\theta \le 1/4$,
  whole-interval contraction gives $\sup\|U_D\| \le (1-c_\theta dB^\theta)^{-1}$
  and hence $\log\sup\|U_D\| \le 2c_\theta B^\theta d$. Otherwise take $n_D =
  \lceil B(4c_\theta d)^\nu\rceil$ equal pieces. Their contraction ratios are at
  most $1/4$ and local norms at most $4/3$. The expression inside the ceiling
  exceeds one, so $n_D \le 2B(4c_\theta d)^\nu$. Composition gives the $d^\nu$
  bound. At $d = 0$, $U_D = I$, proving~\eqref{E:young-upper-refined-growth}
  without a prefactor.

  \medskip
  \noindent\textit{Step 2.~}We bound the global evolution and inverse on a
  fixed driver ball.
  On the ball $\|H\|_\theta \le K$ choose one uniform equal partition with ratio
  at most $1/4$. Each local evolution $Q$ has $N_\ell(Q) \le 4/3$ and
  $\sup\|Q-I\| \le 1/3$. Its Neumann inverse is bounded by $3/2$. Products give
  uniform global sup bounds. The same scaled estimate bounds local H\"older
  seminorms. Pairs crossing a boundary at distance below $\ell$ split into at
  most two pieces; use $h_1^\theta+h_2^\theta \le 2^{1-\theta}(h_1+h_2)^\theta$.
  For larger distances use $2\sup\|U_H\|/\ell^\theta$. This gives its global
  $C^\theta$ bound. Finally,
  \[
    U_H^{-1}(t)-U_H^{-1}(s)
    =-U_H^{-1}(t)(U_H(t)-U_H(s))U_H^{-1}(s)
  \]
  bounds the inverse seminorm. This proves Lemma~\ref{L:young-upper-growth}.
\end{proof}

The linear term in~\eqref{E:young-upper-refined-growth} cannot simply be
dropped: No-go~\ref{P:nogo-26} gives a scalar commuting counterexample.
Its second observation also distinguishes an averaged block budget from
a uniform bound on every block.

\begin{proposition}[Multiplicative upper perturbation]
  \label{P:young-upper-multiplicative}
  If $K \ge 1$, $\|H\|_\theta \le K$, $x = \|X\|_\theta$ and $A \in
  \mathcal{B}(\mathcal{X})$, then
  \begin{equation}\label{E:young-upper-multiplicative}
    \|AU_{H+X}(t)\|
    \le\|AU_H(t)\|\exp\{C_{\theta,B,K}(x+x^\nu)\},
    \qquad 0\le t\le B.
  \end{equation}
  This includes noninvertible $A$ and $A = 0$ with no extra factor when $x = 0$.
\end{proposition}

\begin{proof}
  We conjugate the perturbation by the deterministic evolution to form a new
  Young driver. We verify the resulting product identity directly from
  increments, preserving the left multiplication order, and then apply the
  refined growth bound only to the perturbation driver.

  \medskip
  \noindent\textit{Step 1.~}We form the conjugated perturbation driver.
  Put $U = U_H$, $R = U^{-1}$ and $G(s)Z = R(s)ZU(s)$.
  Lemma~\ref{L:young-upper-growth} makes $G$ a $C^\theta$ path of bounded linear
  maps on $\mathcal{B}(\mathcal{X})$, uniformly bounded on the fixed ball. Thus,
  the Banach Young integral
  \[
    Z(t)=\int_0^tR(s)\ud X(s)U(s)
  \]
  exists, with $\|Z\|_\theta \le L'_{\theta,B,K}x$ and $Z_{s,t} =
  R(s)X_{s,t}U(s)+O(|t-s|^{2\theta})$. Put $V = U_Z$. Its and $U$'s Young
  remainders are $V_{s,t} = Z_{s,t}V(s)+O(|t-s|^{2\theta})$ and $U_{s,t} =
  H_{s,t}U(s)+O(|t-s|^{2\theta})$. The exact increment of $P = UV$ includes a
  cross term of order $2\theta$. Substitution in the stated left order gives
  \[
    P_{s,t}
    =H_{s,t}P(s)+U(s)R(s)X_{s,t}U(s)V(s)
    +O(|t-s|^{2\theta})
    =(H+X)_{s,t}P(s)+O(|t-s|^{2\theta}).
  \]

  \medskip
  \noindent\textit{Step 2.~}We identify the product as the full Young
  evolution.
  Each interval length to the power $2\theta$ is at most that length times the
  mesh to the power $2\theta-1$, so summing bounds the remainder by the total
  length times that mesh power.
  On any partition the remainder sum is bounded by a constant times
  $B\,\operatorname{mesh}^{2\theta-1} \to 0$. Left sums converge to the Young
  integral, so uniqueness yields $U_{H+X} = U_HU_Z$. This proves the product
  identity without assuming smooth density in the full H\"older norm.

  \medskip
  \noindent\textit{Step 3.~}We retain the deterministic endpoint factor in the
  norm bound.
  At $x=0$, the driver $Z$ is zero and $U_Z=I$.
  Apply~\eqref{E:young-upper-refined-growth} to $Z$, then left multiply by $A$
  and use submultiplicativity. No inverse or lower bound for $A$ is required.
  This proves Proposition~\ref{P:young-upper-multiplicative}.
\end{proof}

\paragraph{What we now have.}
The compact map from a Gaussian control to its mean, at a fixed spatial radius
(Propositions~\ref{P:cameron-full-mean} and~\ref{P:cameron-actual-embedding}).
The next chapter lets the radius grow.

\chapter{Bounds as the spatial radius grows}\label{S:book-growing-radius}

Finite propagation makes spatial localization useful, but the region we need
grows with the observation time (Figure~\ref{Fig:light-cone}, in the next
chapter). A bound with an uncontrolled constant in the radius is useless for
long time. This chapter sorts the estimates into two kinds: those uniform
in the radius, and stochastic ones that must display it.

Section~\ref{S:uniform-mean-regulation} gives uniform mean and Cameron bounds
on the displacement--velocity Sobolev phase; it claims nothing of the kind for
the additional constant coordinates. Section~\ref{S:radius-stochastic-anchor}
tracks the radius in the Gaussian path moments and estimates averages over
time blocks. Section~\ref{S:signed-static-wave} builds the physical wave
for completed static mean potentials of either sign, so the upper argument
keeps all its Gaussian controls.

We use the covariance and phase constructions of
Chapters~\ref{S:book-gaussian-foundations}--\ref{S:book-wave-evolution}
and the mean maps of Chapter~\ref{S:book-means}. On a first reading, keep
Propositions~\ref{P:uniform-mean-regulator} and~\ref{P:signed-physical-wave},
and let the proofs wait. The restriction to nonnegative controls enters later,
through a separate variational comparison. Episode~\VideoRef{7}{Z35qRxt05kA}
of the video series covers this chapter and the next.

\begin{center}
\small
\begin{tabular}{@{}p{0.40\textwidth}p{0.19\textwidth}p{0.33\textwidth}@{}}
  \toprule
  Estimate & Result & Dependence on the radius $R$ \\
  \midrule
  Completed mean: norm and tail & Proposition~\ref{P:uniform-mean-regulator} & none, uniform for $R\ge1$ \\
  Cameron embedding of the block paths & Proposition~\ref{P:uniform-actual-cameron} & none, uniform for $R\ge1$ \\
  \midrule
  Compact envelopes of the cutoff & Proposition~\ref{P:radius-compact-envelopes} & displayed: $\rho^{d/q}$ \\
  Moments of the operator path & Proposition~\ref{P:radius-actual-moments} & displayed: $R^{d/Q}$ \\
  \raggedright Expected weighted norm of block paths & Proposition~\ref{P:radius-mean-anchor} & $g_\delta R^{d/Q}\to0$ at $R_\delta=\delta^{-1}+2$ if $Q>2d/(1-\alpha_0)$ \\
  \bottomrule
\end{tabular}

\smallskip
\begin{minipage}{0.92\textwidth}\footnotesize
  All rows concern the Sobolev phase $E=H^{1/2}\times H^{-1/2}$, not the
  additional constant columns on $F$. In the last two rows, fix the exponents
  and even $Q$ as in~\eqref{E:radius-anchor-parameters}--\eqref{E:radius-anchor-Q}.
  The last limit fixes $B,p,Q$ before $\delta\downarrow0$.
\end{minipage}
\end{center}

\section{Mean and Cameron bounds uniform in the spatial radius}
\label{S:uniform-mean-regulation}

Retain the original nonzero covariance, amplitude and real spaces of
Section~\ref{S:actual-cameron-mean}. Fix $\alpha/4 < s < 1/2$ and put $\beta =
1/2-s > 0$. Write $A_s = C_s^2\sqrt{R_{\mathrm{sp}}}$ for the constant
in~\eqref{E:cameron-fractional-product}. Choose a real even smooth compact
$\eta$, equal to one on the unit ball, with $0 \le \eta \le 1$ and support in
the closed radius-two ball. For $R \ge 1$ set $\eta_R(x) = \eta(x/R)$ and $m_R =
\eta_R^2$. In this section the mean acts on $E = H^{1/2}\times H^{-1/2}$. The
constant data columns belong to $F$, outside $E$. All uniform-radius claims
below concern $E$. The spatial map $T_R:H_{\mathrm{sp}} \to \mathcal{K}(E)$ is
the $E$ restriction of Proposition~\ref{P:cameron-spatial-map}. Its explicit row
is
\begin{equation}\label{E:uniform-mean-row}
  \langle T_Rqz,w\rangle_E
  =\langle q,k_R(z,w)\rangle_{H_{\mathrm{sp}}},
  \qquad
  k_R(z,w)=(\eta_Rz_1)(\eta_R\Lambda^{-1}w_2).
\end{equation}
Thus, the physical mask occurs once, after covariance convolution.

\subsection{Uniform multiplier and spatial tails}

We first bound the row $k_R$ in~\eqref{E:uniform-mean-row} uniformly in
$R \ge 1$. Both factors of this row carry the cutoff $\eta_R$, so we start with
its multiplier norm on $H^\sigma$.

\begin{lemma}[Uniform cutoff and physical mean]\label{L:uniform-mean-cutoff}
  For every real $\sigma$, multiplication by $\eta_R$ on $H^\sigma$ has norm at
  most
  \[
    M_\sigma=2^{|\sigma|/2}(2\pi)^{-d/2}
    \int_{\R^d}|\widehat{\eta}(p)|\langle p\rangle^{|\sigma|}\ud p,
  \]
  independently of $R \ge 1$. With $M = M_{1/2}$ and $M_0 = A_sM^2$,
  \begin{equation}\label{E:uniform-spatial-mean}
    \|k_R(z,w)\|_{H_{\mathrm{sp}}}
    \le M_0\|z\|_E\|w\|_E,\qquad
    \|T_R\|_{H_{\mathrm{sp}}\to\mathcal{L}(E)}\le M_0 .
  \end{equation}
\end{lemma}

\begin{proof}
  Peetre's inequality bounds the weighted Fourier kernel for multiplication by
  $\eta_R$ by convolution with
  \[
    2^{|\sigma|/2}(2\pi)^{-d/2}
    |\widehat{\eta_R}(q)|\langle q\rangle^{|\sigma|}.
  \]
  Here $\widehat{\eta_R}(q) = R^d\widehat{\eta}(Rq)$. Substitution $p = Rq$ and
  $\langle p/R\rangle \le \langle p\rangle$ bound its $L^1$ norm by $M_\sigma$.
  The $L^1*L^2$ inequality proves the assertion on Schwartz functions and
  density extends it to $H^\sigma$. The map $w_2\mapsto\Lambda^{-1}w_2$ is an
  isometry from $H^{-1/2}$ to $H^{1/2}$, and $H^{1/2}\hookrightarrow H^s$ has
  norm one. Apply the cutoff bound at $\sigma = 1/2$ and
  Lemma~\ref{L:cameron-fractional-products} to the two factors
  in~\eqref{E:uniform-mean-row}. Hilbert Cauchy--Schwarz then gives the operator
  bound. This is the existing completed physical map: the row agrees on smooth
  controls and hence on their completion. Its compact values at each fixed $R$
  are retained. This proves Lemma~\ref{L:uniform-mean-cutoff}.
\end{proof}

Lemma~\ref{L:uniform-mean-cutoff} bounds $T_R$ on all of $H_{\mathrm{sp}}$.
The approximation in Proposition~\ref{P:uniform-mean-regulator} also needs
$T_R$ to be small on high frequencies, uniformly in $R \ge 1$. We obtain this
from a bound on $T_R(I-C_\varepsilon)$ for the heat contraction
$C_\varepsilon$.

\begin{lemma}[Uniform heat and frequency tails]
  \label{L:uniform-mean-spatial-tails}
  Let $C_\varepsilon$ be the covariance-Hilbert heat contraction with Fourier
  multiplier $e^{-\varepsilon|\xi|^2/2}$. There is a real orthogonal projection
  $P_K$ on $H_{\mathrm{sp}}$ corresponding to the closed frequency ball
  $\{|\xi|\le K\}$. It commutes with heat. For $K \ge 1$, $\varepsilon > 0$ and
  $R \ge 1$,
  \begin{align}
    \|T_R(I-C_\varepsilon)\|
    &\le M_h\varepsilon^{\beta/2},
    &M_h&=4A_sM^2\E|G_d|^\beta,\label{E:uniform-mean-heat}\\
    \|T_R(I-P_K)\|
    &\le M_fK^{-\beta},
    &M_f&=\frac{M_h}{1-e^{-1/2}},\label{E:uniform-mean-frequency}
  \end{align}
  where $G_d$ is a standard real $d$-dimensional normal vector. Operator norms
  are from $H_{\mathrm{sp}}$ to $\mathcal{L}(E)$. The projection $P_K$ is
  allowed to have infinite rank.
\end{lemma}

\begin{proof}
  We estimate translations of each covariance row and average them against
  heat. We then construct the frequency projection on the actual Fourier image
  and use the heat error to bound its high-frequency tail.

  \medskip
  \noindent\textit{Step 1.~}We bound translations of the localized row.
  Translations $\tau_hf(x) = f(x-h)$ extend to strongly continuous isometries of
  the $H_{\mathrm{sp}}$. Indeed, in its covariance Fourier isometry they
  multiply by $e^{-i\xi\cdot h}$; dominated convergence gives continuity on
  tests, then density and isometry give it on the completion. The real subspace
  is preserved. The inequalities $|e^{-it}-1| \le 2\min(1,|t|^\beta)$ and
  $|\xi|^\beta\langle\xi\rangle^s \le \langle\xi\rangle^{1/2}$ show that
  \[
    \|\tau_hf-f\|_{H^s}
    \le2\min(1,|h|^\beta)\|f\|_{H^{1/2}} .
  \]
  Use this for $f = \eta_Rz_1$ and $g = \eta_R\Lambda^{-1}w_2$. The product
  identity
  \[
    \tau_h(fg)-fg=(\tau_hf-f)\tau_hg+f(\tau_hg-g)
  \]
  holds in $H_{\mathrm{sp}}$, first for smooth factors and then by the
  continuous fractional product extension. Lemma~\ref{L:uniform-mean-cutoff}
  gives
  \begin{equation}\label{E:uniform-mean-row-translation}
    \|\tau_hk_R(z,w)-k_R(z,w)\|_{H_{\mathrm{sp}}}
    \le4A_sM^2\min(1,|h|^\beta)\|z\|_E\|w\|_E .
  \end{equation}

  \medskip
  \noindent\textit{Step 2.~}We average the translation bound to obtain the
  heat error.
  The Bochner integral $\int p_\varepsilon(h)\tau_hk\ud h$ exists, since the
  translation orbit is strongly continuous and bounded in the separable
  covariance Hilbert space. Its Fourier image is
  $e^{-\varepsilon|\xi|^2/2}\widehat{k}$: ordinary Gaussian integration proves
  this on tests and contraction extends it to all $k$. It is therefore the
  $C_\varepsilon$ of the source construction. Heat has variance $\varepsilon$
  here; there is no additional Fourier normalization factor. Integrating
  \eqref{E:uniform-mean-row-translation} with $h = \sqrt\varepsilon G_d$ bounds
  $\|(I-C_\varepsilon)k_R(z,w)\|$ by $M_h\varepsilon^{\beta/2}\|z\|_E\|w\|_E$.
  Here $C_\varepsilon$ is selfadjoint on $H_{\mathrm{sp}}$, since its Fourier
  multiplier is real; the row identity therefore gives
  \[
    \langle T_R(I-C_\varepsilon)qz,w\rangle_E
    =\langle q,(I-C_\varepsilon)k_R(z,w)\rangle_{H_{\mathrm{sp}}}.
  \]
  Selfadjointness in~\eqref{E:uniform-mean-row} proves
  \eqref{E:uniform-mean-heat}.

  The two terms in the product difference each use one translation error and
  one bounded factor. The fractional product bound therefore contributes
  $2A_sM^2$ twice. Averaging $\min(1,|h|^\beta)\le |h|^\beta$ with
  $h=\sqrt\varepsilon G_d$ gives $\varepsilon^{\beta/2}\E|G_d|^\beta$.

  \medskip
  \noindent\textit{Step 3.~}We construct the projection without assuming
  Fourier surjectivity.
  To construct $P_K$, let $V$ be the closed image of the covariance Fourier
  isometry into $L^2(\mu)$, including its fixed $(2\pi)^{d/2}$ factor. No
  surjectivity is assumed. Its formula extends to Schwartz functions: compact
  cutoff approximation converges in Schwartz seminorms, and temperedness of
  $\mu$ controls the Fourier error in $L^2(\mu)$. Choose real even radial smooth
  $\chi_n$, equal to one on the closed ball of radius $K$, zero outside radius
  $K+1/n$, and between zero and one. For real Schwartz $f$, the function
  $\chi_n\widehat{f}$ is the transform of another real Schwartz function and so
  lies in $V$. Dominated convergence against $|\widehat{f}|^2\mu$ gives
  $1_{\{|\xi|\le K\}}\widehat{f} \in V$. The pointwise limit includes the
  boundary sphere, even if $\mu$ charges it. Test density and contraction extend
  this assertion to all of $V$. Indicator multiplication is selfadjoint and
  idempotent on $V$; pulling it back defines the asserted real orthogonal $P_K$.
  All Fourier multipliers here commute.

  \medskip
  \noindent\textit{Step 4.~}We compare the frequency tail with a heat
  residual.
  On $|\xi| > K$ one has $1-e^{-|\xi|^2/(2K^2)} \ge 1-e^{-1/2}$, whence
  \[
    \|(I-P_K)k\|_{H_{\mathrm{sp}}}
    \le(1-e^{-1/2})^{-1}\|(I-C_{K^{-2}})k\|_{H_{\mathrm{sp}}}.
  \]
  Apply the row heat bound and selfadjointness of $P_K$
  in~\eqref{E:uniform-mean-row} to obtain \eqref{E:uniform-mean-frequency}. No
  density of $\mu$, periodic state space or finite rank was used. This proves
  Lemma~\ref{L:uniform-mean-spatial-tails}.
\end{proof}

\subsection{Tensor products and mean approximation}

We first record a tensor extension bound for maps into $L^r(I)$. The proof
of Proposition~\ref{P:uniform-mean-regulator} uses it to pass from a scalar
time map to the Hilbert-valued time map without a dimension factor.

\begin{lemma}[Dimension-free Hilbert tensor extension]
  \label{L:uniform-mean-tensor}
  If $H_0,G$ are real Hilbert spaces, $2 \le r < \infty$ and $A:H_0 \to L^r(I)$
  is bounded, then
  \[
    \|A\otimes I_G\|_{H_0\otimes G\to L^r(I;G)}\le\|A\|.
  \]
  The tensor extension need not be compact when $A$ is compact.
\end{lemma}

\begin{proof}
  We use finite Gaussian sums to estimate the Hilbert-valued tensor norm with
  the scalar operator constant. We then extend by density and check why the
  same argument does not imply compactness.

  \medskip
  \noindent\textit{Step 1.~}We compare the Gaussian norm moment with its total
  variance.
  For a finite tensor write $v = \sum_{j=1}^n f_j\otimes e_j$ with $(e_j)$
  orthonormal. Its squared norm is $\sum_j\|f_j\|^2$. Let $\gamma_j$ be
  independent standard real normals and $m_r = \E|\gamma_1|^r$. Finite
  covariance diagonalization gives
  \[
    \E\left\|\sum_j\gamma_j f_j\right\|_{H_0}^r
    \le m_r\left(\sum_j\|f_j\|_{H_0}^2\right)^{r/2}.
  \]
  In fact, if $\lambda_l$ are the covariance eigenvalues and $L =
  \sum_l\lambda_l$, the squared norm has law $\sum_l\lambda_l\zeta_l^2$ for
  independent standard $\zeta_l$. For $L > 0$, convexity at exponent $r/2 \ge 1$
  bounds $(\sum_l(\lambda_l/L)\zeta_l^2)^{r/2}$ by
  $\sum_l(\lambda_l/L)|\zeta_l|^r$. Taking expectations proves the inequality;
  $L = 0$ is immediate. Also $L = \sum_j\|f_j\|^2$.

  \medskip
  \noindent\textit{Step 2.~}We apply the scalar operator estimate and cancel
  the Gaussian moment.
  Choose measurable representatives of the finitely many $Af_j$. Scalar Gaussian
  moments and Tonelli give
  \[
    m_r\int_I\left(\sum_j|Af_j(t)|^2\right)^{r/2}\ud t
    =\E\left\|A\sum_j\gamma_jf_j\right\|_{L^r}^r
    \le\|A\|^rm_r\left(\sum_j\|f_j\|^2\right)^{r/2}.
  \]
  For each fixed $t$, the scalar sum of the $Af_j(t)$ with independent
  Gaussian coefficients has variance $\sum_j|Af_j(t)|^2$. Its $r$th absolute
  moment is that variance to the power $r/2$ times $m_r$. This explains the
  same factor $m_r$ on both sides, so no dimension-dependent factor remains
  after cancellation.
  Cancel the finite positive $m_r$ and take the $r$th root. Finite tensor
  density and Bochner completeness give the extension. The zero space is
  trivial. If $Af \ne 0$ and $G$ has an infinite orthonormal sequence $(e_j)$,
  the images of the bounded sequence $f\otimes e_j$ have pairwise distance
  $\sqrt2\|Af\|_{L^r}$. They have no convergent subsequence, even when the
  scalar map $A$ is compact. Only finite auxiliary Gaussians were used in the
  estimate. This proves Lemma~\ref{L:uniform-mean-tensor}.
\end{proof}

Recall the supported identity $H_I = H_t\otimes H_{\mathrm{sp}}$, the compact
scalar covariance map $C_t:H_t \to L^r(I)$, $2 < r < 2/a$, and its completed
tensor $U = C_t\otimes I_{H_{\mathrm{sp}}}$ from
Section~\ref{S:actual-cameron-mean}. On finite smooth tensors $U$ agrees with
the earlier $\mathcal{U}$; density gives equality of their completed maps. For
block applications fix the same $r$ as in~\eqref{E:cameron-block-parameters}.
All $c_t,\varepsilon_J$ and mean norms below use that one choice. Choose a
smooth-test orthonormal basis of $H_t$ and let $\pi_J$ be its first-$J$
orthogonal projection, with $\pi_0 = 0$. Put
\[
  c_t=\|C_t\|,\qquad
  \varepsilon_J=\|C_t(I-\pi_J)\|_{H_t\to L^r(I)},\qquad
  \mathcal{T}_R=T_RU .
\]

\begin{proposition}[Uniform approximation of the completed mean]
  \label{P:uniform-mean-regulator}
  On $H_I$ the operators
  \[
    P_{\mathrm{ret}}=\pi_J\otimes P_K,\qquad
    Q_{J,K}=I-P_{\mathrm{ret}}
  \]
  are real orthogonal projections, and $\varepsilon_J \to 0$. The mean is a
  strongly measurable $L^r(I;\mathcal{K}(E))$ field and satisfies, uniformly for
  $R \ge 1$,
  \begin{equation}\label{E:uniform-completed-mean}
    \|\mathcal{T}_R\|\le M_0c_t,\qquad
    \|\mathcal{T}_RQ_{J,K}\|
    \le M_0\varepsilon_J+M_fc_tK^{-\beta}.
  \end{equation}
  These norms have domain $H_I$ and target $L^r(I;\mathcal{L}(E))$.
\end{proposition}

\begin{proof}
  We transfer the compact scalar temporal tail to Hilbert-valued functions. We
  then decompose the residual projection into a temporal tail and a
  spatial-frequency tail, estimating these terms separately.

  \medskip
  \noindent\textit{Step 1.~}We control the temporal projection tail.
  The tensor projections are selfadjoint, idempotent and contractive on finite
  tensors and hence on their Hilbert completion. For the scalar tail, failure of
  convergence would supply unit $v_J$ with $\|C_t(I-\pi_J)v_J\|$ bounded below
  on a subsequence. But $(I-\pi_J)v_J$ tends weakly to zero, since its pairing
  with $h$ is bounded by $\|(I-\pi_J)h\|$. Compactness supplies a
  norm-convergent image subsequence; bounded linear scalar observables show its
  limit is zero, a contradiction. Lemma~\ref{L:uniform-mean-tensor} transfers
  this scalar tail norm to the Hilbert-valued time map without a dimension
  factor.

  \medskip
  \noindent\textit{Step 2.~}We decompose and bound the completed mean
  residual.
  For $\Phi\in H_I$, Lemma~\ref{L:uniform-mean-tensor} applied to $A=C_t$
  gives $\|U\Phi\|_{L^r(I;H_{\mathrm{sp}})}\le c_t\|\Phi\|_{H_I}$.
  Applying~\eqref{E:uniform-spatial-mean} pointwise in time then bounds
  $\|T_RU\Phi\|_{L^r(I;\mathcal{L}(E))}$ by $M_0c_t\|\Phi\|_{H_I}$.
  For fixed $R$, $T_R$ is bounded into the separable compact-operator space on
  $E$. Composition with the Bochner field $U\Phi$ gives the stated mean and its
  first bound. The exact decomposition is
  \[
    Q_{J,K}=(I-\pi_J)\otimes I+\pi_J\otimes(I-P_K).
  \]
  The first term is bounded by $M_0\varepsilon_J$. For the second term use the
  intertwining identity
  \[
    U(\pi_J\otimes(I-P_K))
    =(I_{L^r}\otimes(I-P_K))(C_t\pi_J\otimes I).
  \]
  It holds first on finite tensors, then by continuity.
  Apply~\eqref{E:uniform-mean-frequency} pointwise, followed by
  Lemma~\ref{L:uniform-mean-tensor} for $C_t\pi_J$, whose norm is at most $c_t$.
  This proves~\eqref{E:uniform-completed-mean}. Infinite spatial rank is allowed
  throughout; neither tensor compactness nor uniform-radius compactness has been
  asserted. This proves Proposition~\ref{P:uniform-mean-regulator}.
\end{proof}

\subsection{The block Cameron embedding for the given noise}

Use $B,\delta,N,s_m,\ell_m,I_m$ and $\theta,r,p$ exactly as
in~\eqref{E:cameron-block-parameters} and its preceding block partition. Let
$i_I:H_I \to H_W$ be the isometric closed inclusion. On the full original noise
space set
\[
  P_{\mathrm{full}}=i_IP_{\mathrm{ret}}i_I^*,\qquad
  Q_{\mathrm{full}}=I-P_{\mathrm{full}} .
\]
These are orthogonal projections on $H_W$. Define the residual fast noise and
stopped $E$ operator path by
\begin{equation}\label{E:uniform-actual-residual}
  Y_m(f)=W(Q_{\mathrm{full}}T_{\delta,m}f),\qquad
  X_m(u)=g_\delta A_{m_R}(Y_m)(\min(u,\ell_m)),\qquad
  g_\delta=\delta^{(1-a)/2}.
\end{equation}
Here $T_{\delta,m}$ is the explicit original isometry
\eqref{E:raw-block-pullback}, and $A_{m_R}$ takes its $E$ restriction. Its
isometry is proved in Lemma~\ref{L:raw-block-inputs}. The symbol
$H_{\mathrm{fast}}$ below denotes this same original covariance norm in fast
coordinates, carried isometrically into $H_W$. The $E$ little H\"older space
means the closure of Lipschitz compact-operator paths starting at zero in
$C^\theta$ norm. Compact operators on separable $E$ and polygonal approximation
make this closure separable. The last interval can be empty, in which case $X_N
= 0$.

\begin{proposition}[Uniform-radius Cameron bound]
  \label{P:uniform-actual-cameron}
  The tuple $X = (X_m)_{m=0}^N$ is a centered Gaussian variable in the separable
  little H\"older compact-operator path space on $E$. Its Bochner covariance $Lh
  = \E[XW(h)]$, $h \in H_W$, is
  \begin{equation}\label{E:uniform-actual-covariance}
    (Lh)_m(u)=\sqrt\delta\int_0^{\min(u,\ell_m)}
    S_E(-v)\mathcal{T}_R(Q_{J,K}i_I^*h)(s_m+\delta v)
    S_E(v)\ud v .
  \end{equation}
  In the exact norm $\|x\|_{\delta,p} =
  (\delta\sum_m\|x_m\|_{C^\theta}^p)^{1/p}$, its Cameron embedding has the
  uniform-radius bound
  \begin{equation}\label{E:uniform-actual-cameron-bound}
    \|Lh\|_{\delta,p}
    \le\sqrt\delta\,C_{B,E}(2/B)^{1/p-1/r}
    [M_0\varepsilon_J+M_fc_tK^{-\beta}]\|h\|_{H_W},
    \qquad C_{B,E}\le18B^{3-1/r}.
  \end{equation}
  All spatial constants are independent of $R \ge 1$.
\end{proposition}

\begin{proof}
  We construct the residual paths on the given noise, identify their supported
  covariance rows, and estimate the resulting primitive. The true block
  partition then gives the averaged Cameron bound with its exact
  radius-independent constants.

  \medskip
  \noindent\textit{Step 1.~}We construct the residual paths and their
  first-chaos coordinates.
  Orthogonality and original pullback isometry give $\operatorname{Var}Y_m(f)
  \le \|f\|_{H_{\mathrm{fast}}}^2$. Gaussian transfer,
  Proposition~\ref{P:source-gaussian-transfer}, therefore constructs the
  operator paths on the given $W$, with finite path moments at each fixed
  $R,B,\delta,J,K$. Choose an admissible H\"older exponent strictly above
  $\theta$ as in Lemma~\ref{L:cameron-actual-residual}; constant stopping puts
  the paths in the separable little $C^\theta$ closure. Finite-coefficient
  constructions and their $L^2$ limits show that every continuous linear
  coordinate is in the first chaos of the given $W$.

  \medskip
  \noindent\textit{Step 2.~}We identify the covariance direction on the
  supported subspace.
  Every explicit observed row, pulled back to slow time, belongs to $H_I$. This
  is the supported-row approximation from Lemma~\ref{L:cameron-actual-residual},
  also valid before any projection; approximate endpoints inside $I$ against the
  integrable time kernel and locally finite spatial covariance. For such a row
  $k$, exact orthogonality gives
  \[
    \langle Q_{\mathrm{full}}T_{\delta,m}k,h\rangle_{H_W}
    =\langle T_{\delta,m}k,i_IQ_{J,K}i_I^*h\rangle_{H_W}.
  \]
  In the earlier construction $P_0=0$, so
  $Y_{m,0}(f)=W(T_{\delta,m}f)$.
  The preceding pairing therefore reduces the present covariance to that
  unprojected identity with control $\Phi=Q_{J,K}i_I^*h\in H_I$.
  Directions perpendicular to $H_I$ disappear from this pairing. Apply the
  completed $J = 0$ identity \eqref{E:cameron-actual-covariance} to
  $Q_{J,K}i_I^*h$, and restrict its input and output to $E$. The free $F$ group
  preserves $E$ and restricts to $S_E$; the mean restricts to $\mathcal{T}_R$.
  This proves~\eqref{E:uniform-actual-covariance}. Its amplitude is exactly
  $g_\delta\delta^{a/2} = \sqrt\delta$, on the noise $W$. Arbitrary projected
  controls are covered by the completed mean identity; no smoothness or finite
  rank of $P_K$ is required. Compact operator values and strong continuity of
  $S_E$ and its adjoint make conjugation norm continuous on compact operators,
  as in Proposition~\ref{P:cameron-actual-embedding}. Thus, the integral is
  Bochner measurable. Fixed-parameter $L^2$ covariance continuity, dense states
  and rational times identify the entire operator path.

  \medskip
  \noindent\textit{Step 3.~}We estimate the free group and each local
  primitive.
  For explicit free constants, the Bessel-unitary Fourier matrix of $S_E(v)$ is
  the matrix displayed in that proposition. For $|v| \le B$, $B \ge 1$, the
  upper off-diagonal entry is at most $\sqrt2 B$: use $|\sin(vk)| \le Bk$ for $k
  \le 1$ and $\langle k\rangle/k \le \sqrt2$ for $k \ge 1$. The lower entry is
  at most one, as are both diagonal entries; at zero the sine quotient has value
  $v$. Its Frobenius bound is $\sqrt{3+2B^2} \le 3B$. Hence, $G_{B,E} =
  \sup_{|v|\le B}\|S_E(v)\| \le 3B$. With
  \[
    \mathcal{M}_h(t)=\mathcal{T}_R(Q_{J,K}i_I^*h)(t),\qquad
    t_m=\|\mathcal{M}_h\|_{L^r(I_m;\mathcal{L}(E))},
  \]
  the primitive increment bound and the affine time change give
  \[
    \|(Lh)_m\|_{C^\theta}
    \le\sqrt\delta\,C_{B,E}\delta^{-1/r}t_m,\qquad
    C_{B,E}=G_{B,E}^2[B^{1-1/r}+B^{1-1/r-\theta}]
    \le18B^{3-1/r}.
  \]
  Here $1-1/r > \theta$, constant stopping preserves the bound, and the empty
  block contributes zero.

  \medskip
  \noindent\textit{Step 4.~}We combine the local estimates with the block
  weights.
  Put $n = N+1 \ge 2$. The finite power-mean inequality at exponent $r/p > 1$
  with nonnegative inputs $t_m^p$ gives
  \[
    \left(\sum_m t_m^p\right)^{r/p}
    \le n^{r/p-1}\sum_m t_m^r,\qquad
    \left(\sum_m t_m^p\right)^{1/p}
    \le n^{1/p-1/r}\left(\sum_m t_m^r\right)^{1/r}.
  \]
  By~\eqref{E:uniform-completed-mean} and $\|i_I^*\|\le1$, the field
  $\mathcal{M}_h=\mathcal{T}_RQ_{J,K}i_I^*h$ satisfies
  $\|\mathcal{M}_h\|_{L^r}\le
  [M_0\varepsilon_J+M_fc_tK^{-\beta}]\|h\|_{H_W}$.
  Raising nonnegative terms uses only positive powers, including when $t_m = 0$.
  The true partition gives $\sum_m t_m^r = \|\mathcal{M}_h\|_{L^r(I)}^r$ and
  $\delta n \le 2/B$. The exact $\delta$ weight therefore yields
  \[
    \|Lh\|_{\delta,p}
    \le\sqrt\delta\,C_{B,E}(\delta n)^{1/p-1/r}\|\mathcal{M}_h\|_{L^r}
    \le\sqrt\delta\,C_{B,E}(2/B)^{1/p-1/r}
    [M_0\varepsilon_J+M_fc_tK^{-\beta}]\|h\|_{H_W}.
  \]

  The local factor $\delta^{-1/r}$ is multiplied by the averaged norm's factor
  $\delta^{1/p}$. Finite power mean contributes $n^{1/p-1/r}$, so their
  product is $(\delta n)^{1/p-1/r}$. Since $p<r$, this power is increasing,
  and $\delta n\le2/B$ gives the claimed bound. The two powers in the
  primitive constant are each at most $B^{1-1/r}$ because $B\ge1$, giving
  $18B^{3-1/r}$ from $G_{B,E}\le3B$.

  \medskip
  \noindent\textit{Step 5.~}We apply the bound to the minimum-norm Cameron
  representative.
  Finally, let $H_{\mathrm{seen}}$ be the closed span in $H_W$ of all
  coefficient tests $k_\ell$ with $\ell(X) = W(k_\ell)$. The first Wiener chaos
  Cameron identification in Proposition~\ref{P:cameron-actual-embedding} applies
  on $H_W$: $L$ kills $H_{\mathrm{seen}}^\perp$ and its Cameron space is
  $L(H_{\mathrm{seen}})$ with minimum norm of a preimage in the original
  covariance space. Thus, $\|Lh\|_{\mathrm{CM}} = \|P_{\mathrm{seen}}h\|_{H_W}$.
  Applying the bound to this minimum-norm representative gives precisely the
  Cameron embedding bound, without an inverse covariance or finite-dimensional
  density. This proves Proposition~\ref{P:uniform-actual-cameron}.
\end{proof}

The orthogonal $P_{\mathrm{full}}$ and $Q_{\mathrm{full}}$ Gaussian fields are
independent; residual time blocks can remain correlated. The retained spatial
head may have infinite rank. Heat filtering alone would not give this orthogonal
decomposition. These linear bounds do not bound the expected residual path norm
uniformly on growing masks, nor perform conditioning or integration over the
infinite sum. The remaining stochastic, product and pointwise upper estimates
require their own proofs.

\section{Bounds for the operator path with explicit radius dependence}
\label{S:radius-stochastic-anchor}

Retain every original covariance and amplitude in Hypothesis~\ref{H:homog}, with
$a = \alpha_0$, and the real phase $E = H^{1/2}\times H^{-1/2}$. Use the fixed
cutoff profile $\eta$ of Section~\ref{S:uniform-mean-regulation} and physical
mask $m_R = \eta(\cdot/R)^2$, $R \ge 1$. In the variance construction the base
mask is $\zeta = \eta^2$. The balanced unitary is $J_0(u,v) =
(\Lambda^{1/2}u,\Lambda^{-1/2}v)$. Fix a finite horizon $B \ge 1$ and parameters
\begin{equation}\label{E:radius-anchor-parameters}
  0<\gamma<\min\left(\frac12,\frac{2-\alpha}{2}\right),\qquad
  \frac12<\theta<\eta_t<\beta_t<
  \min\left(1-\frac a2,\frac{3-a-\alpha}{2}\right).
\end{equation}
The strict subcritical assumptions make these intervals nonempty. Fix one even
integer $Q \ge 2$ such that
\begin{equation}\label{E:radius-anchor-Q}
  Q\gamma>d,\qquad
  \beta_t-\eta_t-\frac1Q>0,\qquad
  Q>\frac{2d}{1-a}.
\end{equation}
Such a $Q$ exists because all three lower thresholds are finite. It stays fixed
as $\delta$ tends to zero. Constants below may depend on the model, $B,Q$, these
exponents and the fixed auxiliary cutoff profiles. They are independent of the
spatial radius, noise-basis truncation, moment index and deterministic residual
projection unless explicitly stated otherwise.

\subsection{Uniform outer factors and the exact whole-space envelope}

Fix real $\chi_0 \in C_c^\infty(\R^d)$, equal to one on $B_1$ and supported in
the closed ball of radius two. Choose real $\psi_0 \in C_c^\infty(\R^d)$, equal
to one on a neighbourhood of $\operatorname{supp}\chi_0$ and supported inside
$B_3$. For $\rho \ge 1$ put $\chi_\rho(x) = \chi_0(x/\rho)$, $\psi_\rho(x) =
\psi_0(x/\rho)$ and $L = 2\pi\rho$. Since $3 < \pi$, both cutoffs are supported
strictly inside the cube of side $L$. The auxiliary torus has unnormalized
$L^2(Q_L,\ud x)$, orthonormal modes $L^{-d/2}e^{ik\cdot x/\rho}$ and Sobolev
weights $\langle k/\rho\rangle^r$.

\begin{lemma}[Uniform rescaled outer maps]\label{L:radius-outer-maps}
  The cutoff periodization and recovery maps from
  Lemma~\ref{L:cutoff-factorization} satisfy, for every real $r$,
  \[
    \|\mathcal{P}_{\psi_\rho}:H^r(\R^d)\to H^r(\mathbb{T}_L)\|
    \le C_{\psi_0,r},\qquad
    \|\mathcal{E}_{\chi_\rho}:H^r(\mathbb{T}_L)\to H^r(\R^d)\|
    \le C_{\chi_0,r},
  \]
  independently of $\rho \ge 1$.
\end{lemma}

\begin{proof}
  We rescale the exact Fourier kernels and bound their integral and
  lattice-sum norms separately. The radius powers cancel in the mixed Schur
  estimate, after which density gives the full Sobolev maps.

  The exact weighted Fourier kernels are
  \begin{align*}
    &L^{-d/2}
    \left(\frac{\langle k/\rho\rangle}{\langle\xi\rangle}\right)^r
    \widehat{\psi_\rho}(k/\rho-\xi),\\
    &L^{-d/2}
    \left(\frac{\langle\xi\rangle}{\langle k/\rho\rangle}\right)^r
    \widehat{\chi_\rho}(\xi-k/\rho).
  \end{align*}
  There is no additional $(2\pi)^{-d/2}$: it cancels between Fourier inversion
  and the ordinary integral in the first kernel, while the torus normalization
  gives the second. For either base profile $\phi_0$, $\widehat{\phi_\rho}(t) =
  \rho^d\widehat{\phi_0}(\rho t)$. Peetre's inequality~\eqref{E:cutoff-peetre},
  also for negative $r$, and $\langle(\rho\xi-k)/\rho\rangle \le
  \langle\rho\xi-k\rangle$ bound the relevant kernel by
  \[
    C_r\rho^{d/2}F_r(k-\rho\xi),\qquad
    F_r(t)=|\widehat{\phi_0}(\pm t)|\langle t\rangle^{|r|}.
  \]
  The bracket comparison uses $\rho\ge1$; consequently, the constants
  multiplying the explicit radius powers in the two Schur bounds depend only
  on the fixed $\phi_0,r,d$.
  The nonnegative $F_r$ is integrable and bounded by arbitrarily high inverse
  powers of $\langle t\rangle$. At fixed $k$, integration in $\xi$ gives a bound
  $C\rho^{-d/2}$. At fixed $\xi$, summation in $k$ gives $C\rho^{d/2}$:
  translate $\rho\xi$ into a unit lattice cell and sum a power tail of order
  greater than $d$, as in~\eqref{E:cutoff-lattice-sum}. This lattice spacing is
  one and has no remaining $\rho$ dependence. The mixed Schur estimate bounds
  the operator norm by the square root of the product of those two quantities.
  Their radius powers cancel. For recovery the row and column roles are
  reversed.

  In the integral bound, the substitution $\rho\xi-k$ contributes $\rho^{-d}$.
  Multiplying by the kernel's factor $\rho^{d/2}$ gives $\rho^{-d/2}$.
  Bounded kernel extension from Schwartz functions and torus polynomials gives
  the full Sobolev maps by density. Their distributional meanings are the same
  cutoff maps as in Lemma~\ref{L:cutoff-factorization}; no point sampling of a
  rough Fourier function is used. This includes negative indices. This proves
  Lemma~\ref{L:radius-outer-maps}.
\end{proof}

Lemma~\ref{L:radius-outer-maps} bounds the two outer maps independently of
$\rho \ge 1$. The proof below places a diagonal inclusion on the torus between
them. We now sum the $q$th powers of the singular values of this inclusion and
take the $q$th root, which produces the factor $\rho^{d/q}$.

\begin{proposition}[Radius-explicit compact envelopes]
  \label{P:radius-compact-envelopes}
  For real $\sigma$, $\gamma > 0$ and $q \ge 1$ with $q\gamma > d$,
  \begin{equation}\label{E:radius-cutoff-Schatten}
    \|\Lambda^\sigma\chi_\rho\Lambda^{-\sigma-\gamma}\|_{S_q}
    \le C_{\sigma,\gamma,q,d,\chi_0,\psi_0}\rho^{d/q}.
  \end{equation}
  In the balanced coordinates of~\eqref{E:compact-two-factors},
  \begin{equation}\label{E:radius-two-envelopes}
    \|K_{\rm in}\|_{S_q}=\|K_{\rm out}\|_{S_q}
    \le C_q\rho^{d/q}.
  \end{equation}
\end{proposition}

\begin{proof}
  We factor the cutoff through the auxiliary torus. We compare the diagonal
  inclusion sum with a Euclidean integral and then use finite-rank
  approximations to transfer that sum to the whole-space operator.

  \medskip
  \noindent\textit{Step 1.~}We use the full-domain recovery factorization.
  The exact recovery identity, with $\chi_\rho\psi_\rho = \chi_\rho$, gives
  \[
    \Lambda^\sigma\chi_\rho\Lambda^{-\sigma-\gamma}
    =\Lambda^\sigma\mathcal{E}_{\chi_\rho}\mathcal{I}_\gamma
    \mathcal{P}_{\psi_\rho}\Lambda^{-\sigma-\gamma}.
  \]
  Here $\mathcal{I}_\gamma:H^{\sigma+\gamma}(\mathbb{T}_L) \to
  H^\sigma(\mathbb{T}_L)$ has diagonal singular values $\langle
  k/\rho\rangle^{-\gamma}$. The outer Bessel maps are unitary on their indicated
  Sobolev domains. The preceding lemma bounds both other outer maps uniformly in
  $\rho$. This is an equality on the full whole-space domain by the recovery
  proof and density, including negative $\sigma$.

  \medskip
  \noindent\textit{Step 2.~}We bound the diagonal sum by a volume integral.
  Put $u = q\gamma > d$. If $x \in k+[-1/2,1/2]^d$, then
  \[
    \langle x/\rho\rangle
    \le\sqrt2\langle k/\rho\rangle\langle(x-k)/\rho\rangle
    \le c_d'\langle k/\rho\rangle,\qquad
    c_d'=\sqrt{2(1+d/4)}.
  \]
  Integrate over the unit cubes to obtain
  \begin{equation}\label{E:radius-lattice-volume}
    \sum_{k\in\mathbb{Z}^d}\langle k/\rho\rangle^{-u}
    \le(c_d')^u\rho^d\int_{\R^d}\langle y\rangle^{-u}\ud y
    \le C_{d,u}\rho^d .
  \end{equation}
  The integral is finite: dyadic annuli outside the unit ball give a geometric
  series with terms bounded by $C_d2^{j(d-u)}$.

  The unit cubes partition Euclidean space and each has volume one. The
  preceding weight comparison bounds the constant cube value $\langle
  k/\rho\rangle^{-u}$ by $(c_d')^u\langle x/\rho\rangle^{-u}$. Summing the
  cube integrals and substituting $y=x/\rho$ gives the factor $\rho^d$. Taking
  a $q$th root later gives $\rho^{d/q}$.

  \medskip
  \noindent\textit{Step 3.~}We transfer the singular-value bound to the two
  whole-space envelopes.
  Arrange the inclusion's diagonal values in decreasing order $(\lambda_j)$.
  They tend to zero. Retaining its first $j-1$ coordinates has rank less than
  $j$ and error $\lambda_j$. Compose this approximation with the bounded outer
  factors. The singular-value approximation argument in
  Proposition~\ref{P:cutoff-singular-values} gives $s_j \le
  C_{\chi_0,\sigma}C_{\psi_0,\sigma+\gamma}\lambda_j$. Finite-rank approximation
  also proves compactness. Summing $q$th powers and
  using~\eqref{E:radius-lattice-volume} proves~\eqref{E:radius-cutoff-Schatten}.
  No separate Schatten ideal theorem is needed.

  The two blocks of $K_{\rm out} = P_{\chi_\rho}D_\gamma$ are the whole-space
  operators with $\sigma = 1/2,-1/2$. Their $q$th power norms add, proving its
  bound. Proposition~\ref{P:compact-two-envelopes} gives $K_{\rm in} = V K_{\rm
  out}V$ for the unitary coordinate swap $V(z_1,z_2) = (z_2,z_1)$, so their
  norms agree. In~\eqref{E:radius-anchor-parameters}, $\gamma < 1/2$ and
  $q\gamma > d$ force $q > 2d \ge 2$; there is no Hilbert--Schmidt conclusion.
  Neither $P_{\chi_\rho}$ nor its adjoint is assumed selfadjoint or idempotent.
  This proves Proposition~\ref{P:radius-compact-envelopes}.
\end{proof}

\subsection{The operator path and covariance contractions}

Write $M_R(W) = A_{m_R}(W)$ for the $E$ interaction operator path starting at
zero on $[0,B]$ from Proposition~\ref{P:actual-operator-driver}. The phase
remains the whole-space $E$ throughout.

\begin{proposition}[Radius-explicit path moments]
  \label{P:radius-actual-moments}
  Under~\eqref{E:radius-anchor-parameters}--\eqref{E:radius-anchor-Q}, for every
  finite $v \ge 1$,
  \begin{equation}\label{E:radius-actual-path-moments}
    \|M_R(W)\|_{L^v(C^\theta([0,B];\mathcal{L}(E)))}
    \le C_{Q,B}\sqrt v\,R^{d/Q},\qquad R\ge1.
  \end{equation}
  The constant does not depend on $v$ or $R$. This is an $E$ operator-path
  estimate; the two constant columns on $F$ require their own estimates. Any
  later projection used in a Cameron comparison must also satisfy its
  orthogonality and covariance hypotheses.
\end{proposition}

\begin{proof}
  We keep the physical noise mask fixed while choosing a larger localization
  envelope. We use only one fixed Schatten index to extract the radius factor,
  then apply the dyadic path estimate uniformly over all higher moment
  indices.

  \medskip
  \noindent\textit{Step 1.~}We choose a localization envelope covering the
  propagated support.
  Use the given smooth-test basis from Lemma~\ref{L:driver-noise-basis}, $g_j =
  W(f_j)$, and its finite paths $M_N = \sum_{j\le N}g_jA_j$. These coefficient
  paths use the same given noise as in Lemma~\ref{L:driver-holder-limit}, not
  source approximations. Choose
  \[
    \rho=2R+B+2\le(B+4)R,\qquad \chi=\chi_0(\cdot/\rho).
  \]
  The plateau covers $B_{2R+B+1}$. Thus, both explicit localization identities
  \eqref{E:driver-actual-localization} hold for this envelope. The physical
  noise mask remains $m_R$. Both envelopes $K_{\rm in},K_{\rm out}$ here use
  this same $\chi_\rho$.

  \medskip
  \noindent\textit{Step 2.~}We estimate increments using both variance bounds.
  For an interval $I \subset [0,B]$ of length $0 < h \le 1$,
  Proposition~\ref{P:wave-two-variances} and
  Lemma~\ref{L:driver-actual-variances} give both finite weighted variances with
  bound $V_{B,\mathrm{var}}h^{2-a}B_\alpha(h)$, uniformly in $R,N,I$. In
  particular, the base mask is the same fixed $\eta^2$ required by that variance
  theorem. Corollary~\ref{C:finite-gaussian-operators}, applied in the balanced
  coordinates, gives for every even $q \ge Q$
  \[
    \|M_N(I)\|_{L^q(S_q)}
    \le\sqrt{q-1}\sqrt{V_{B,\mathrm{var}}}\,
    h^{1-a/2}B_\alpha(h)^{1/2}\|K_{\rm out}\|_{S_q},
  \]
  where $B_\alpha(h)$ is $1$, $1+\log(1/h)$ or $h^{1-\alpha}$ according as
  $\alpha < 1$, $\alpha = 1$ or $\alpha > 1$. Both variances and both
  localizations have been checked; individual basis-dependent coefficient bounds
  are not summed.

  \medskip
  \noindent\textit{Step 3.~}We extract the radius factor at the fixed index.
  For $q \ge Q$, the unnormalized Schatten norms obey $\|K_{\rm out}\|_{S_q} \le
  \|K_{\rm out}\|_{S_Q}$. Indeed, normalize a nonzero singular-value sequence by
  its $\ell^Q$ norm; each entry is at most one, so its $q$th power sum is at
  most its $Q$th power sum. The zero case is immediate.
  Apply~\eqref{E:radius-two-envelopes} only at the fixed index $Q$. Use $\|A\|
  \le \|A\|_{S_q}$ and $\sqrt{q-1} \le \sqrt q$. The balanced Hilbert
  complexification has the same operator norm as its real restriction, and
  $\rho^{d/Q} \le (B+4)^{d/Q}R^{d/Q}$. Since $\beta_t$ is strictly below the
  corresponding short-time power, it follows that
  \begin{equation}\label{E:radius-finite-increment}
    \|M_N(t)-M_N(u)\|_{L^q(\mathcal{L}(E))}
    \le C_{Q,B}\sqrt q\,R^{d/Q}|t-u|^{\beta_t},
    \qquad |t-u|\le1.
  \end{equation}
  At $\alpha = 1$, the factor $\sqrt{1+\log(1/h)}$ is absorbed by the strict gap
  $1-a/2-\beta_t > 0$. For $\alpha > 1$, the power is $(3-a-\alpha)/2$. These
  power/log constants are independent of $R,N,q$. Longer increments are finite
  sums of increments of length at most one, with constants depending only on the
  fixed $B$.

  \medskip
  \noindent\textit{Step 4.~}We sum dyadic increments and pass to the path
  limit.
  Let $D_n(F)$ be the largest adjacent increment on the grid $kB2^{-n}$ of
  $[0,B]$. For fine levels $B2^{-n} \le 1$, the finite-maximum bound gives
  \[
    \|D_n(M_N)\|_{L^q}
    \le2^{n/q}C_{Q,B}\sqrt q\,R^{d/Q}(B2^{-n})^{\beta_t}.
  \]
  No independence of increments is used. The sum after multiplication by
  $2^{n\eta_t}$ converges uniformly for even $q \ge Q$, because its geometric
  ratio is at most $2^{-(\beta_t-\eta_t-1/Q)} < 1$. The finitely many coarse
  levels are controlled by unit partitions. The deterministic dyadic chaining
  argument in Lemma~\ref{L:driver-holder-limit} therefore bounds the
  $L^qC^{\eta_t}$ norms by $C_{Q,B}\sqrt q\,R^{d/Q}$. Each finite path is
  continuous and starts at zero, as required there. The fixed-horizon embedding
  $C^{\eta_t}\hookrightarrow C^\theta$ gives the same estimate in $C^\theta$.

  For each fixed $R$, Lemma~\ref{L:driver-holder-limit} identifies $M_R(W)$ as
  the all-finite-moment $C^\theta$ limit of these very $M_N$. Pass the uniform
  finite-$N$ estimate to this limit. Uniform convergence in $R$ is unnecessary:
  the limiting bound has its explicit radius dependence at every $R$. For
  arbitrary $v \ge 1$, choose an even $q \ge \max(v,Q)$ with $q \le v+Q+2 \le
  (Q+3)v$. Lyapunov's inequality gives~\eqref{E:radius-actual-path-moments},
  absorbing $\sqrt{Q+3}$ into the fixed constant. This proves
  Proposition~\ref{P:radius-actual-moments}.
\end{proof}

Proposition~\ref{P:radius-actual-moments} concerns the path $M_R(W)$ of the
noise $W$ itself. Proposition~\ref{P:radius-mean-anchor} works instead with
residual noises $Y_m$, obtained from $W$ through a deterministic orthogonal
projection. We now transfer the moment bound to such noises.

\begin{corollary}[Deterministic covariance residuals]
  \label{C:radius-covariance-residuals}
  Let $Y$ be a centered Gaussian linear noise on the coupled space, with
  $\operatorname{Var}Y(f) \le C\|f\|_{H_W}^2$, $C \ge 0$. Its transferred
  operator path satisfies
  \[
    \|M_R(Y)\|_{L^v(C^\theta)}
    \le\sqrt C\,C_{Q,B}\sqrt v\,R^{d/Q}.
  \]
  In particular, this holds with $C = 1$ for $Y_m(f) = W((I-\Pi)T_{\delta,m}f)$,
  for any deterministic orthogonal projection $\Pi$ on the slow $H_W$. Its rank
  may be infinite and it may depend on $\delta,R$.
\end{corollary}

\begin{proof}
  Proposition~\ref{P:source-gaussian-transfer} is a linear Gaussian covariance
  contraction theorem for this same physical mask and these same deterministic
  coefficient maps. It bounds the path norm by $\sqrt C$ times the reference
  isonormal operator-path norm. The reference has the canonical operator-path
  law used in Proposition~\ref{P:radius-actual-moments}; the $Y$ is not replaced
  in this comparison. The zero case gives zero directly. This proves the first
  assertion. For the second, the original isometry \eqref{E:raw-block-pullback},
  proved in Lemma~\ref{L:raw-block-inputs}, and $\|I-\Pi\| \le 1$ give
  $\operatorname{Var}Y_m(f) \le \|f\|_{H_W}^2$, where the right side is the
  original norm in fast coordinates. The joint residuals retain the given slow
  $W$ and may be correlated. They need not lie in the span of an individual fast
  noise, which is not a hypothesis of Gaussian transfer. No random
  noise-dependent projection is used. This proves
  Corollary~\ref{C:radius-covariance-residuals}.
\end{proof}

\subsection{The growing-radius expected-norm anchor}

Corollary~\ref{C:radius-covariance-residuals} bounds each residual path with
the growing factor $R^{d/Q}$. We now average over the time blocks with the
amplitude $g_\delta$. The last condition in~\eqref{E:radius-anchor-Q} makes
$g_\delta R^{d/Q}$ tend to zero at the radius $R_\delta = \delta^{-1}+2$.

\LeanMarker{P:radius-mean-anchor}
\begin{proposition}[Averaged block mean tends to zero]
  \label{P:radius-mean-anchor}
  Fix integer $B \ge 1$ and use the true block partition $B,\delta,N,\ell_m$ of
  Section~\ref{S:actual-cameron-mean}, with $0 < \delta < \min(1,1/B)$ and $N =
  \lfloor1/(\delta B)\rfloor$. For $Y_m$ as in
  Corollary~\ref{C:radius-covariance-residuals}, put
  \[
    X_m(u)=g_\delta M_R(Y_m)(\min(u,\ell_m)),\qquad
    Z_\delta=\left(\delta\sum_{m=0}^N
    \|X_m\|_{C^\theta}^p\right)^{1/p},\qquad
    g_\delta=\delta^{(1-a)/2},
  \]
  for any fixed $p \ge 1$. Then
  \begin{equation}\label{E:radius-mean-anchor}
    \E Z_\delta\le(2/B)^{1/p}C_{Q,B}\sqrt p\,
    g_\delta R^{d/Q}.
  \end{equation}
  In particular, for $R_\delta = \delta^{-1}+2$,
  \begin{equation}\label{E:radius-mean-anchor-decay}
    2\E Z_\delta
    \le2(2/B)^{1/p}C_{Q,B}\sqrt p\,3^{d/Q}
    \delta^{e_{\rm anchor}}\longrightarrow0,\qquad
    e_{\rm anchor}=\frac{1-a}{2}-\frac dQ>0 .
  \end{equation}
  Both bounds are uniform in the deterministic orthogonal projection.
\end{proposition}

\begin{proof}
  Stopping an interaction primitive starting at zero and extending it constantly
  to $B$ does not increase its supremum or H\"older seminorm. An empty last
  interval gives the zero path. Probability normalization and Lyapunov give
  \begin{align*}
    \E Z_\delta
    &\le(\E Z_\delta^p)^{1/p}
    =\left(\delta\sum_{m=0}^N
    \E\|X_m\|_{C^\theta}^p\right)^{1/p}\\
    &\le[\delta(N+1)]^{1/p}g_\delta
    C_{Q,B}\sqrt p\,R^{d/Q}.
  \end{align*}
  This is linearity for a finite sum; no time-block independence or maximum over
  blocks is used. Since $\delta B < 1$, we have $N \ge 1$ and $N+1 \le 2/(\delta
  B)$, proving~\eqref{E:radius-mean-anchor}.

  Now $R_\delta \le 3/\delta$ because $\delta < 1$. The last condition
  in~\eqref{E:radius-anchor-Q} is equivalent to $Q(1-a) > 2d$. Dividing by $2Q >
  0$ gives $e_{\rm anchor} > 0$. Substitute this radius
  into~\eqref{E:radius-mean-anchor} to
  obtain~\eqref{E:radius-mean-anchor-decay}. The estimate holds at each fixed
  $B,p,Q$ and model, uniformly in every deterministic residual projection. This
  proves Proposition~\ref{P:radius-mean-anchor}.
\end{proof}

These bounds concern the expected operator-path norm on $E$. They do not give an
$F$ estimate with constant initial-data columns. The radius substitution is
numerical; cone coverage and source-to-point transfer require their separate
arguments. No concentration, integration over an infinite sequence of Gaussian
coordinates or nonlinear moment upper bound follows from this section alone. The
constants need not be uniform as $Q,B$ tend to infinity or as model parameters
approach an endpoint. For zero covariance, the operator paths and stochastic
bounds vanish directly.

\section{Signed rough static potentials and the physical wave}
\label{S:signed-static-wave}

Fix the original nonzero covariance and a finite physical mask $m_R$ as in
Section~\ref{S:actual-cameron-mean}. For real $q \in C_c^\infty(\R^d)$ put $V_q
= m_R(\Gamma*q)$ and $L_q(u,v,c_1,c_2) = (0,V_q(u+c_1),0,0)$ on $F =
E\oplus\R^2$. Let $\mathcal{V}_R$ be the operator-norm closure of these
operators in $\mathcal{K}(F)$. It is a closed real linear space. No positivity
is imposed on $q$.

\begin{lemma}[The limiting multiplier]\label{L:signed-limiting-multiplier}
  Every $L \in \mathcal{V}_R$ determines a unique real $V \in H^{-1/2}$
  supported in $\operatorname{supp}m_R$. It has the multiplier action
  $L(u,v,c_1,c_2) = (0,Vu+c_1V,0,0)$, where multiplication extends from tests to
  $H^{1/2} \to H^{-1/2}$. For $k = \|L\|$, its symmetric pairing obeys
  \begin{equation}\label{E:signed-multiplier-form}
    |v_V(u,w)|\le k\|u\|_{H^{1/2}}\|w\|_{H^{1/2}} .
  \end{equation}
\end{lemma}

\begin{proof}
  The form $v_V(u,w)$ is the duality pairing
  $\langle Vu,w\rangle_{H^{-1/2},H^{1/2}}$.
  Once the multiplier bound $\|Vu\|_{H^{-1/2}}\le k\|u\|_{H^{1/2}}$
  is obtained below, duality gives~\eqref{E:signed-multiplier-form}.
  Choose $L_{q_n} \to L$. The constant-displacement vector has norm one, so
  $\|V_{q_n}-V_{q_j}\|_{H^{-1/2}} \le \|L_{q_n}-L_{q_j}\|$. Thus, $V_{q_n} \to
  V$ in $H^{-1/2}$. Testing outside the fixed compact mask support preserves
  that support. Every zero column of $L_{q_n}$ remains zero. Its displacement
  multipliers converge in $\mathcal{L}(H^{1/2},H^{-1/2})$ to a map of norm at
  most $k$. For compact smooth $u$, smooth multiplication is bounded on
  $H^{-1/2}$ by Lemma~\ref{L:cutoff-factorization}. Hence, $V_{q_n}u \to Vu$
  there, identifying the limiting map on tests and then uniquely on $H^{1/2}$ by
  density. The constant column recovers $V$ and proves uniqueness. On compact
  smooth $u,w$, the pairing is $\langle V,uw\rangle$, denoted $v_V(u,w)$, so it
  is symmetric. Its bound and symmetry extend by density. Complexification
  preserves the operator norm by splitting real and imaginary parts in the
  Hilbert norms. The pairing extends sesquilinearly to a Hermitian form; its
  diagonal is the sum of the real and imaginary diagonal values. This proves
  Lemma~\ref{L:signed-limiting-multiplier}.
\end{proof}

The bound~\eqref{E:signed-multiplier-form} controls the pairing $v_V$ on
$H^{1/2}$. We now use it to define a self-adjoint operator $H_V$ from the
quadratic form $q_V$ on $H^1$. We also compare the quarter power of its shift
$A$ with the $H^{1/2}$ norm.

\LeanMarker{P:signed-closed-form}
\begin{proposition}[Closed form and quarter-power phase]
  \label{P:signed-closed-form}
  Let $K \ge \max(1,\|L\|)$ and $c = 2+K^2$. On complexified $L^2$ the form
  \[
    q_V(u)=\|\nabla u\|_2^2-v_V(u,u),\qquad Q(q_V)=H^1,
  \]
  represents a unique self-adjoint $H_V$. With $A = H_V+c$ and $L_K = 3(1+K^2)$,
  \begin{align}
    A&\ge1,\qquad D(A^{1/2})=H^1,\label{E:signed-form-domain}\\
    \tfrac12\|u\|_{H^1}^2
    &\le q_V(u)+c\|u\|_2^2\le L_K\|u\|_{H^1}^2,\label{E:signed-form-bounds}\\
    D(A^{1/4})&=H^{1/2},\qquad
    2^{-1/4}\|u\|_{H^{1/2}}
    \le\|A^{1/4}u\|_2\le L_K^{1/4}\|u\|_{H^{1/2}}.
    \label{E:signed-quarter-phase}
  \end{align}
  The negative dual completion is $H^{-1/2}$. No $H^2$ operator domain is
  asserted.
\end{proposition}

\begin{proof}
  We first bound the signed potential form by the Sobolev norm and identify
  its closed representation. We then derive the quarter-power comparison from
  resolvents, and use duality to identify the negative phase.

  \medskip
  \noindent\textit{Step 1.~}We make the shifted form coercive and closed.
  Put $b = \|u\|_2$, $z = \|u\|_{H^1}$. Fourier Cauchy--Schwarz
  and~\eqref{E:signed-multiplier-form} give
  \[
    \|u\|_{H^{1/2}}^2\le bz,\qquad
    |v_V(u,u)|\le Kbz\le\tfrac12z^2+\tfrac12K^2b^2 .
  \]
  The last scalar inequality is the weighted arithmetic--geometric mean
  inequality at weights $1/2,1/2$ and nonnegative bases $z^2,(Kb)^2$, including
  zero norms. Since $q_V = z^2-b^2-v_V(u,u)$, it follows that
  \[
    q_V+cb^2\ge\tfrac12z^2+(1+K^2/2)b^2,\qquad
    q_V+cb^2\le\tfrac32z^2+(1+3K^2/2)b^2\le L_Kz^2.
  \]
  The lower bound dominates both $z^2/2$ and $b^2$. Thus, this densely defined
  Hermitian form is semibounded and its shifted norm is equivalent to the
  complete $H^1$ norm. It is closed in the definition
  of~\cite{teschl:14:mathematical}*{p.~80}.
  From~\cite{teschl:14:mathematical}*{Theorem~2.14, p.~81, (2.63)}, it
  represents a unique self-adjoint $H_V$. This representation is independent of
  the admissible shift. From that book's Theorem~3.6 and~(3.51)--(3.53), p.~109,
  applied after shifting, we obtain~\eqref{E:signed-form-domain}. Here $A \ge 1$
  is initially its form lower bound. A bounded spectral slice below one would
  produce a unit $D(A)$ test contradicting that form bound; hence it is also the
  spectral lower bound. Lemma~\ref{L:signed-spectral-wave} records the full
  bottom argument.

  \medskip
  \noindent\textit{Step 2.~}We compare fractional powers through resolvent
  order.
  Set $B_0 = 1-\Delta$. On the common form domain, $\tfrac12B_0 \le A \le
  L_KB_0$. We prove the needed square-root comparison rather than assuming an
  order rule for arbitrary operator powers. For $T_0,S_0 \ge 1$ with common form
  domain and $q_{T_0} \le q_{S_0}$, completing the square in the shifted form
  norm gives
  \[
    \langle f,(T_0+t)^{-1}f\rangle
    =\sup_u\{2\operatorname{Re}\langle f,u\rangle
    -q_{T_0}(u)-t\|u\|^2\}
    \ge\langle f,(S_0+t)^{-1}f\rangle,\quad t>0 .
  \]
  The scalar identity $\sqrt x = \pi^{-1}\int_0^\infty t^{-1/2}x/(x+t)\ud t$
  follows by $t = xv^2$. Positive spectral integration and Tonelli now imply
  $T_0^{1/2} \le S_0^{1/2}$ as forms, with their domain inclusions:
  $I-t(T_0+t)^{-1} \le I-t(S_0+t)^{-1}$. Apply this to $(B_0,2A)$ and
  $(A,L_KB_0)$. Both operators in each pair are at least one. The resulting
  inequalities give~\eqref{E:signed-quarter-phase}, including equality of its
  domains. For $v \in L^2$ the $L^2$-duality norm against $\|A^{1/4}u\|_2$ is
  $\|A^{-1/4}v\|_2$: Cauchy--Schwarz gives the upper bound, and $u = A^{-1/2}v$
  gives equality after normalization. Duality and completion therefore identify
  the negative phase. Uniqueness of the form, resolvent and spectral calculus
  makes these constructions commute with real conjugation. This proves
  Proposition~\ref{P:signed-closed-form}.
\end{proof}

Proposition~\ref{P:signed-closed-form} gives the operator $H_V$ and compares
the norms built from $A^{1/4}$ with those of $E$. We now locate the bottom of
the spectrum of $H_V$ and bound the wave that its spectral calculus defines on
$E$. The bound allows exponential growth at the rate $\lambda(V)$.

\begin{lemma}[The bottom and spectral wave]\label{L:signed-spectral-wave}
  Write $P_V^E$ for the evolution in~\eqref{E:signed-spectral-matrix}, and $S_E$
  for the free restriction from Lemma~\ref{L:constant-phase-free}. Put
  \[
    \beta(V)=\inf_{\substack{u\in H^1\\\|u\|_2=1}}q_V(u),
    \qquad \lambda(V)=\sqrt{-\beta(V)} .
  \]
  Then $\beta(V) = \inf\operatorname{spec}H_V$ and $-K^2 \le \beta(V) \le 0$.
  The spectral wave on $E$ has the bound
  \begin{equation}\label{E:signed-E-wave}
    \|P_V^E(t)\|\le8(1+K)^2(1+t)e^{\lambda(V)t},\qquad t\ge0.
  \end{equation}
\end{lemma}

\begin{proof}
  We identify the form infimum with the spectral bottom and use the compact
  support to show that it is nonpositive. We then estimate the spectral wave
  matrix in its adapted phase norm and transfer the estimate to the fixed
  Sobolev phase.

  \medskip
  \noindent\textit{Step 1.~}We locate the spectral bottom.
  For unit $u$, $z = \|u\|_{H^1} \ge 1$ gives $q_V(u) \ge z^2/2-1-K^2/2 \ge
  -(1+K^2)/2 \ge -K^2$. Conversely dilate a unit compact smooth bump and
  translate its whole support outside $\operatorname{supp}m_R$. Its potential
  pairing is zero and its kinetic energy tends to zero. Thus, $\beta \le 0$.
  Real and complex form infima agree because a complex diagonal splits into the
  two real diagonals. The spectral form formula
  in~\cite{teschl:14:mathematical}*{(3.52)--(3.53)} identifies this form infimum
  with the spectral bottom. For clarity, a bounded spectral slice below the form
  lower bound would give a unit form test violating that bound. Conversely the
  form integral is bounded below by the spectral infimum. The spectrum is
  nonempty: a unit vector's scalar spectral measure has mass one and a support
  point, and normalized projections onto arbitrarily small bounded neighborhoods
  give approximate eigenvectors, excluding a bounded inverse there. Conversely,
  off the closed PVM support $\Sigma$ the Borel function $(h-z)^{-1}$, set to
  zero off $\Sigma$, is bounded for real $z\notin\Sigma$, as is $h/(h-z)$ on
  $\Sigma$. Its calculus maps into $D(H_V)$ and is a two-sided resolvent inverse
  on the appropriate domains. Thus, the spectrum is precisely $\Sigma$. Hence,
  $\lambda \le K$ and every spectral value $h$ has $h+c \ge 2$.

  \medskip
  \noindent\textit{Step 2.~}We estimate the spectral matrix in an equivalent
  phase norm.
  Use the phase norm $\|(u,v)\|_{E_A}^2 = \|A^{1/4}u\|_2^2+\|A^{-1/4}v\|_2^2$.
  With $M_K = L_K^{1/4}$, the preceding proposition and its dual give
  $M_K^{-1}\|z\|_E \le \|z\|_{E_A} \le M_K\|z\|_E$. The condition-number cost is
  $M_K^2 = \sqrt{3(1+K^2)} \le 2(1+K)$. Define $C_t(h),S_t(h)$ by $\cos(t\sqrt
  h),\sin(t\sqrt h)/\sqrt h$ for $h > 0$, by $1,t$ at zero, and by
  $\cosh(t\sqrt{-h}),\sinh(t\sqrt{-h})/\sqrt{-h}$ for $h < 0$. In the $E_A$
  isometry the wave matrix is
  \begin{equation}\label{E:signed-spectral-matrix}
    \begin{pmatrix}
      C_t(h)&\sqrt{h+c}\,S_t(h)\\
      -hS_t(h)/\sqrt{h+c}&C_t(h)
    \end{pmatrix}.
  \end{equation}
  The isometry is $(u,v)\mapsto(A^{1/4}u,A^{-1/4}v)$ onto
  $L^2\times L^2$, with the negative-power map extended to the completed
  velocity space.
  The displayed matrix is the conjugate, initially on bounded spectral
  cutoffs, of the spectral solution
  \[
    P_V^E(t)\binom{u}{v}
    =\binom{C_t(H_V)u+S_t(H_V)v}
    {-H_VS_t(H_V)u+C_t(H_V)v}.
  \]
  For $h \ge 0$, the entries are bounded using $|C_t| \le 1$, $|S_t| \le
  \min(t,h^{-1/2})$, $\sqrt{h+c}|S_t| \le 1+\sqrt c\,t$ and $|hS_t|/\sqrt{h+c}
  \le 1$; at zero the off-diagonals are $\sqrt c\,t,0$. For $\beta \le h < 0$,
  they are bounded using $|C_t| \le e^{\lambda t}$, $|S_t| \le te^{\lambda t}$,
  $\sqrt{h+c}|S_t| \le \sqrt c\,te^{\lambda t}$ and $|hS_t|/\sqrt{h+c} \le
  Ke^{\lambda t}$. There is no division by $\lambda$. Since $\sqrt c \le 2K$,
  each entry is at most $2(1+K)(1+t)e^{\lambda t}$. Bounded Borel calculus
  from~\cite{teschl:14:mathematical}*{Theorem~3.1, pp.~102--103,
  (3.18)--(3.21)}, the two-by-two matrix norm bound and the condition-number
  cost give~\eqref{E:signed-E-wave}. Functions can be set to zero off the
  spectral support. Spectral dominated convergence gives strong continuity, and
  scalar addition identities give the group property. The same matrix and bounds
  apply on the weighted graph levels $D(A^{1/2})\times L^2$ and $D(A)\times
  D(A^{1/2})$. Spectral cutoffs justify their extensions and continuity there.
  For negative spectral values, the lower off-diagonal estimate can be read
  without dividing by $\lambda$: its absolute value is
  $\sqrt{-h}\sinh(t\sqrt{-h})/\sqrt{h+c}$. The bounds $\sqrt{-h}\le K$ and
  $h+c\ge2$ make this at most $Ke^{\lambda t}$. Thus, the estimate also
  remains meaningful when the bottom rate is zero.
  This proves Lemma~\ref{L:signed-spectral-wave}.
\end{proof}

Lemma~\ref{L:signed-spectral-wave} bounds a wave defined by spectral
calculus, and only on $E$. We still have to identify it with the physical wave
and to include the two constant columns of $F$.

\begin{proposition}[Identification with the physical wave]
  \label{P:signed-physical-wave}
  The spectral wave in Lemma~\ref{L:signed-spectral-wave} is the
  bounded-perturbation, equivalently Young, physical wave on $E$. Its extension
  on $F$ satisfies
  \begin{equation}\label{E:signed-F-wave}
    \|P_V(t)\|_{\mathcal{L}(F)}
    \le16(1+K)^4(1+t)^3e^{\lambda(V)t},\qquad t\ge0.
  \end{equation}
  Both constant response columns are included.
\end{proposition}

\begin{proof}
  We first embed the form operator domain into the free generator domain. The
  spectral solution then satisfies the physical Duhamel identity on dense
  graph data. We extend that identity and include the constant columns by
  variation of constants.

  \medskip
  \noindent\textit{Step 1.~}We obtain the graph-continuous embedding.
  The imported form representation \cite{teschl:14:mathematical}*{Theorem~2.14,
  p.~81, (2.63)}, tested against compact smooth functions, gives $H_Vu = -\Delta
  u-Vu$ in distributions for $u \in D(H_V)$. Such $u$ lies in $H^1$ and $Vu \in
  H^{-1/2}$. Consequently, $(1-\Delta)u = H_Vu+Vu+u \in H^{-1/2}$, so
  $D(H_V)\hookrightarrow H^{3/2}$. The embedding is graph-continuous:
  by~\eqref{E:signed-form-bounds}, $\|u\|_{H^1}^2 \le 2\langle Au,u\rangle \le
  \|Au\|_2^2+\|u\|_2^2$, and the distributional equation bounds the $H^{3/2}$
  norm by this graph norm.

  \medskip
  \noindent\textit{Step 2.~}We identify the free generator and apply the
  interaction identity.
  The generator of the free $E$ group is $A_0(u,v) = (v,\Delta u)$ on
  $H^{3/2}\times H^{1/2}$. Indeed, differentiate its Fourier matrix after a
  frequency cutoff, dominate by that graph norm and remove the cutoff.
  Conversely any strong generator derivative has those distributional
  components; they require $v \in H^{1/2}$ and $\Delta u \in H^{-1/2}$, hence $u
  \in H^{3/2}$. No operator-norm differentiability of the group is used.

  Start spectral data in $D(A)\times D(A^{1/2}) = D(H_V)\times H^1$. The
  preceding graph-level estimates give continuous graph trajectories and $C^1$
  trajectories in $H^1\times L^2$, with $u_t = v$, $v_t = -H_Vu$. This follows
  first on bounded spectral cutoffs and then uniformly on finite horizons in the
  graph norms by spectral dominated convergence
  and~\cite{teschl:14:mathematical}*{Theorem~3.2, p.~104, (3.29)--(3.33)}. The
  projections $1_{[1,n]}(A)$ approximate every graph datum and the same weighted
  matrix bounds control all cutoff tails on those horizons. Graph embedding
  gives continuity in $D(A_0)$ and $z_t = A_0z+L_Ez$ in $E$, where $L_E = L|_E$.
  The strong-group chain rule on this domain yields
  \[
    S_E(-t)z(t)=z(0)+\int_0^t S_E(-s)L_Ez(s)\ud s .
  \]
  Graph data are dense in $E$: spectral cutoffs make $D(A)$ dense in $D(A^{1/2})
  = H^1$, and $H^1$ is dense in both $H^{1/2}$ and $H^{-1/2}$. Boundedness on
  finite horizons extends this identity to all $E$. Volterra uniqueness by
  Picard iteration identifies the physical wave. No $D(A^{3/4}) = H^{3/2}$
  assertion is needed.

  \medskip
  \noindent\textit{Step 3.~}We identify the ordinary and Young flows on the
  extended phase.
  On $F$, bounded $L$ and strongly continuous $S_F$ similarly define a unique
  Duhamel flow. On $[0,T]$, the $n$th Volterra term has bound
  $C_T^{n+1}\|L\|^nT^n/n!$. Uniqueness makes horizons consistent and gives
  composition. Its interaction primitive is $\int_0^t S_F(-s)LS_F(s)\ud s$.
  Compactness of $L$, and strong continuity of the group and its adjoint, give
  norm-continuity of the conjugated compact operator, as in
  Proposition~\ref{P:cameron-actual-embedding}. Thus, this primitive is
  Lipschitz and its ordinary equation is exactly the Young equation.

  \medskip
  \noindent\textit{Step 4.~}We bound both constant response columns.
  For initial $(e,c_1,c_2)$ the constant coordinates are $(c_1+tc_2,c_2)$ and
  the $E$ component is
  \[
    P_V^E(t)e+\int_0^t P_V^E(t-s)(0,V(c_1+sc_2))\ud s .
  \]
  The bounded $E$-valued Fubini argument proving \eqref{E:compact-F-block}
  applies: substitute the $E$ Duhamel formula under this integral and
  interchange the absolutely bounded finite time triangle. It solves the full
  system, and uniqueness identifies it. No ground positivity or gap is used. The
  forcing norm is at most $K(1+s)|(c_1,c_2)|$. Since $\lambda \ge 0$ and
  $\int_0^t(1+t-s)(1+s)\ud s \le (1+t)^3$, the full norm is at most
  \[
    (1+t)^3e^{\lambda t}
    \{8(1+K)^2\|e\|+[8K(1+K)^2+1]|(c_1,c_2)|\}.
  \]
  The coefficient of $\|e\|$ is bounded by $8(1+K)^3$. The coefficient of the
  constant-coordinate norm is at most $9(1+K)^3$, since $K\le1+K$ and
  $1\le(1+K)^3$. Cauchy--Schwarz on these two scalar components therefore
  gives the factor $\sqrt{8^2+9^2}$ in the following estimate.
  Scalar Cauchy--Schwarz and $\sqrt{8^2+9^2} < 13$ bound this by
  $13(1+K)^3(1+t)^3e^{\lambda t}\|(e,c_1,c_2)\|_F$. This implies the stated
  sufficient bound. This proves Proposition~\ref{P:signed-physical-wave}.
\end{proof}

\paragraph{What we now have.}
Mean bounds uniform in the radius (Proposition~\ref{P:uniform-mean-regulator})
and the physical wave for signed static means
(Proposition~\ref{P:signed-physical-wave}).

\chapter{Gaussian estimates for the point upper bound}\label{S:book-point-upper}

\begin{figure}
  \centering
  \begin{tikzpicture}[x=1cm,y=1cm,font=\normalsize]
  \fill[hamblue!4] (-5.4,0) rectangle (5.4,3.4);
  \fill[white] (-2.7,0) rectangle (2.7,3.4);
  \fill[hamteal!10] (-2.5,0)--(0,2.5)--(2.5,0)--cycle;
  \draw[ham/arrow] (-5.7,0)--(5.7,0) node[right] {$y$};
  \draw[ham/arrow] (-6.1,-.25)--(-6.1,3.55) node[above] {$s$};
  \draw[thick,hamteal] (-2.5,0)--(0,2.5)--(2.5,0);
  \draw[dashed,hamblue] (-2.7,0)--(-2.7,3.25);
  \draw[dashed,hamblue] (2.7,0)--(2.7,3.25);
  \draw[hamblue] (-5.4,0)--(-5.4,3.25);
  \draw[hamblue] (5.4,0)--(5.4,3.25);
  \fill (0,2.5) circle (1.7pt);
  \node[anchor=west,fill=white,inner sep=2pt] at (.15,2.55) {$(T,x)$};
  \node[ham/note,font=\normalsize,text width=4.2cm] at (0,.65)
  {$|y-x|\le T-s$\\data: $(u,\partial_tu)$};
  \node[font=\normalsize,anchor=north] at (-2.7,-.1) {$-R$};
  \node[font=\normalsize,anchor=north] at (2.7,-.1) {$R$};
  \node[ham/note,font=\normalsize,fill=white] at (0,3.05) {mask equals one on $B_R$};
  \node[font=\normalsize,anchor=north] at (-5.4,-.1) {$-2R$};
  \node[font=\normalsize,anchor=north] at (5.4,-.1) {$2R$};
\end{tikzpicture}
  \caption[Finite propagation and spatial cutoffs]{A spatial cross-section of
    finite propagation. The entire backward cone is chosen inside the region
    where the mask equals one; the mask vanishes outside $\overline{B}_{2R}$.
    The drawing is schematic and does not prescribe a value of $R/T$. Shading
    indicates the domain of dependence, not a Green density; in dimension three
    the Green measure is supported on a sphere. In the point-moment comparison,
    $L=\delta^{-1}$, $R=L+2$ and $x\in B_1$. Both displacement and velocity must
    be passed between time intervals. See \ProofRef{Lemma}{L:raw-block-restart}
    and \ProofRef{Proposition}{P:original-point-transfer}.}
  \label{Fig:light-cone}
\end{figure}
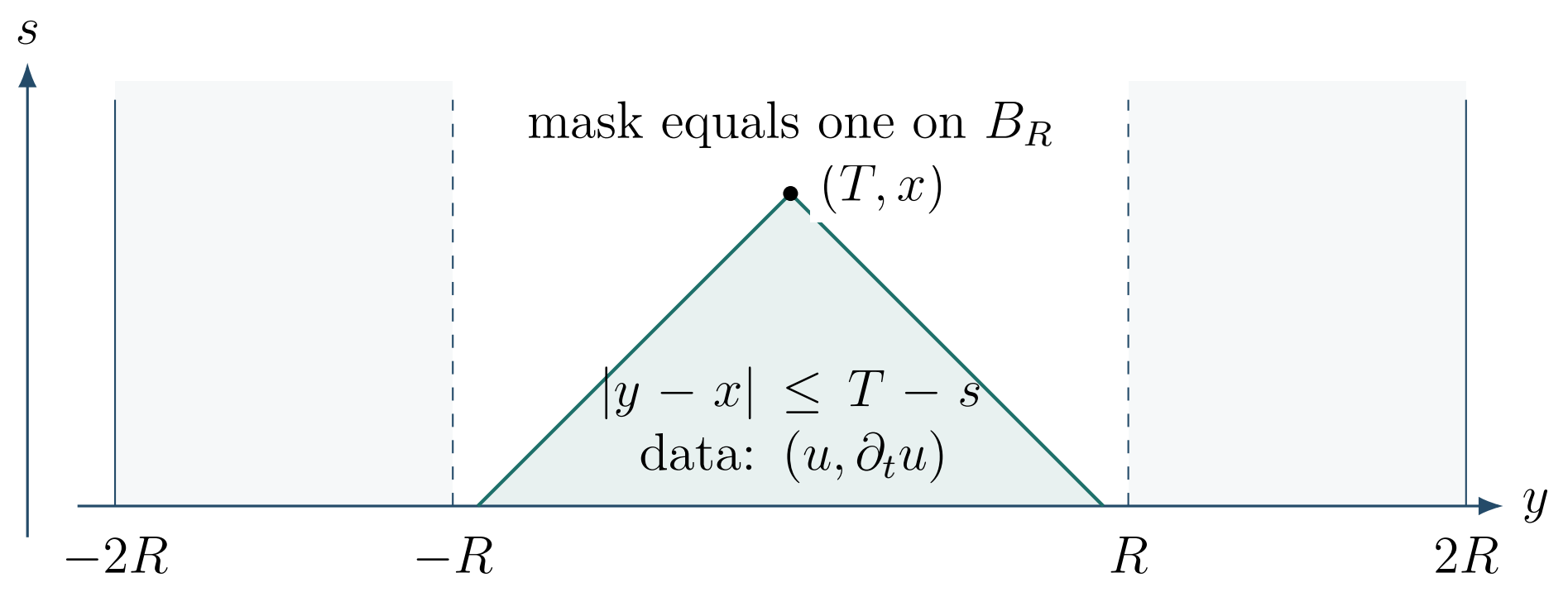

The lower bound needed only the favorable shifts of the noise. The upper
bound must account for every Gaussian fluctuation. It must also control
the moment at one point of space while the region of propagation grows
(Figure~\ref{Fig:light-cone}). The Gaussian cost stays in the original
covariance norm.

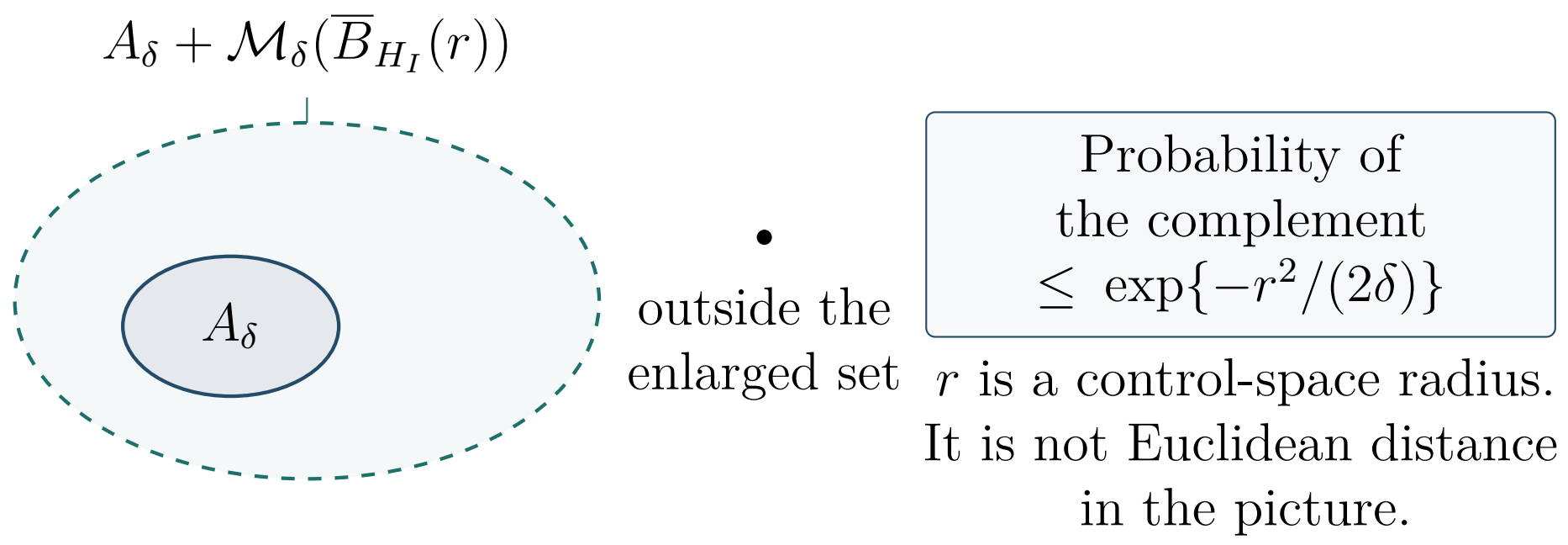
\begin{figure}
  \centering
  % Schematic geometry only; the radius belongs to the control space.
\begin{tikzpicture}[x=1cm,y=1cm,font=\normalsize]
  \filldraw[draw=hamteal,fill=hamteal!5,thick,dashed]
  (2.7,2) ellipse (2.3 and 1.4);
  \filldraw[draw=hamblue,fill=hamblue!12,thick]
  (2.1,1.8) ellipse (.85 and .55);
  \node at (2.1,1.8) {$A_\delta$};
  \node[anchor=south] at (2.7,3.65)
  {$A_\delta+\mathcal{M}_\delta(\overline{B}_{H_I}(r))$};
  \draw[hamteal] (2.7,3.6)--(2.7,3.4);
  \fill (6.3,2.5) circle (1.7pt);
  \node[align=center,anchor=north] at (6.3,2.25)
  {outside the\\enlarged set};
  \node[ham/box,font=\normalsize,text width=4.6cm] at (10.05,2.6)
  {Probability of the complement\\
  $\le \exp\{-r^2/(2\delta)\}$};
  \node[align=center,text width=5cm] at (10.05,.85)
  {$r$ is a control-space radius.\\
  It is not Euclidean distance\\in the picture.};
\end{tikzpicture}
  \caption[Gaussian enlargement by an image of a Hilbert ball]{A schematic
    finite-dimensional projection of enlargement for the nonzero Gaussian law in
    the cited lemma. The map $\mathcal{M}_\delta$ is compact, so the image of a
    closed control ball is compact; the set $A_\delta$ itself need not be
    compact. No ambient Lebesgue density or Euclidean distance formula is being
    assumed. The cost uses one global control norm. See
    \ProofRef{Lemma}{L:tube-closed-geometry} and
    \ProofRef{Proposition}{P:tube-rayleigh-transfer}.}
  \label{Fig:gaussian-enlargement}
\end{figure}

Section~\ref{S:actual-symplectic-phase} sets up the symplectic and norm
properties of the wave evolution. Section~\ref{S:actual-gaussian-tubes}
is the heart of the chapter: Gaussian enlargement turns
a probabilistic estimate into growth under deterministic
mean controls (Figure~\ref{Fig:gaussian-enlargement}).
Section~\ref{S:uniform-physical-products} makes this compatible with the
growing radius. Section~\ref{S:original-point-transfer} returns to the
original point moment. Section~\ref{S:original-signed-score-upper} controls
the Gaussian tail.

The prerequisites are the mean estimates of Chapter~\ref{S:book-means},
the radius bounds of Chapter~\ref{S:book-growing-radius} and the exact
scaling in Section~\ref{S:exact-canonical-scaling}. The chapter delivers
Proposition~\ref{P:original-signed-point-upper}: an upper variational
expression over signed controls. The next chapter matches it to the lower
value. Episode~\VideoRef{7}{Z35qRxt05kA} of the video series covers this
chapter.

\section{The symplectic form and mean evolution on time blocks}
\label{S:actual-symplectic-phase}

Retain the original model, covariance amplitude, physical mask $m_R$ and real $E
= H^{1/2}\times H^{-1/2}$. Fix finite $R,T > 0$ and admissible exponents as
in~\eqref{E:driver-parameters}. Let $M = M_R(W)$ be the from the given noise
first interaction path on $E$, and $M^\rho$ its regularized paths from
Proposition~\ref{P:source-gaussian-transfer} and the source construction. Here
$\rho$ is the positive heat/Laplace smoothing pair. The construction preserving
the given noise and ordinary rows are those of
Lemma~\ref{L:source-retaining-isometry},
Lemma~\ref{L:source-pointfield-regularity} and
Proposition~\ref{P:source-ordinary-response}. The full $F$ operator path and
Young evolution are fixed by Proposition~\ref{P:constant-phase-driver} and
Proposition~\ref{P:constant-young-linear}; their $E$ restrictions are used here.
Complexifications commute with conjugation, so these real restrictions are the
physical phases. For fixed $\Phi$ in the full positive-time original $H_W$ and
real $c$ set
\begin{equation}\label{E:symplectic-combined-driver}
  \begin{aligned}
    Z&=W(\Phi),\qquad D_\Phi=\E[ZM],\qquad
    D_\Phi^\rho=\E[ZM^\rho],\\
    D&=cM+D_\Phi,\qquad D^\rho=cM^\rho+D_\Phi^\rho .
  \end{aligned}
\end{equation}
All these are operator paths on $E$. The covariance mean is not a pointwise
original potential. Write $U_D$ for left Young evolution and $P_D(t,s) =
S_E(t)U_D(t,s)S_E(-s)$.

\LeanMarker{L:symplectic-free}
\begin{lemma}[Free wave form and compact Young baseline]
  \label{L:symplectic-free}
  There is a bounded isometric skew-adjoint $J_E$ on $E$, with $J_E^2 = -I$,
  such that
  \[
    \omega(z,w)=\langle z,J_Ew\rangle_E
    =\langle u,g\rangle_{H^{1/2},H^{-1/2}}
    -\langle f,v\rangle_{H^{1/2},H^{-1/2}},
    \quad z=(u,v),\quad w=(f,g).
  \]
  The free group preserves $\omega$. For bounded real smooth $b$, the force
  $B(b)(u,v) = (0,bu)$ obeys $\omega(B(b)z,w)+\omega(z,B(b)w) = 0$.

  Let $1/2 < \theta < 1$ and let $X \in C^\theta([0,T];\mathcal{K}(E))$ start at
  zero, with $U_X$ its left Young evolution. Then
  \begin{equation}\label{E:symplectic-compact-young-baseline}
    U_X(t,s)-I\in\mathcal{K}(E),\qquad
    \|S_E(t)U_X(t,s)S_E(-s)\|
    \ge\|S_E(t-s)\|\ge1,\qquad 0\le s\le t\le T .
  \end{equation}
  More generally, for finitely many such paths $X_j$ starting at zero on $[0,B]$
  and lengths $0 \le \ell_j \le B$, put $P_j = S_E(\ell_j)U_{X_j}(\ell_j,0)$.
  Then
  \begin{equation}\label{E:symplectic-compact-product-baseline}
    P_N\cdots P_0-S_E\Big(\sum_{j=0}^N\ell_j\Big)\in\mathcal{K}(E),
    \qquad
    \|P_N\cdots P_0\|
    \ge\left\|S_E\Big(\sum_{j=0}^N\ell_j\Big)\right\|\ge1 .
  \end{equation}
  These norm bounds require neither Gaussian nor Hamiltonian operator paths.
\end{lemma}

\begin{proof}
  We identify the symplectic form in balanced coordinates and check free-wave
  preservation. We then show that compact Young perturbations remain compact
  at the endpoint. Translation to spatial infinity gives the required norm
  lower bound, including finite products.

  \medskip
  \noindent\textit{Step 1.~}We identify the form and its preservation by the
  free group.
  The balanced map $C_{\rm bal}(u,v) = (\Lambda^{1/2}u,\Lambda^{-1/2}v)$ is an
  onto real isometry to $L^2\times L^2$. Define $J_{\rm symp}(q,p) = (p,-q)$ and
  $J_E = C_{\rm bal}^{-1}J_{\rm symp}C_{\rm bal}$. This gives all stated
  operator properties and $J_E(f,g) = (\Lambda^{-1}g,-\Lambda f)$. The two
  displayed dualities follow from the balanced $L^2$ pairing. No $L^2$ pairing
  of a general distributional velocity is taken.

  In balanced Fourier coordinates, the free matrix is
  \[
    W_t(\xi)=
    \begin{pmatrix}
      \cos(tk)&\langle\xi\rangle\sin(tk)/k\\
      -k\sin(tk)/\langle\xi\rangle&\cos(tk)
    \end{pmatrix},\qquad k=|\xi|.
  \]
  At zero it is $\left(\begin{smallmatrix}1&t\\0&1\end{smallmatrix}\right)$, by
  the continuous sine quotient. Direct matrix multiplication gives
  $W_t^{\mathsf{T}}J_{\rm symp}W_t = J_{\rm symp}$, since its determinant is
  $\cos^2(tk)+\sin^2(tk) = 1$, also at zero. The entries are real, even and
  bounded on finite time horizons. Integrating the identity against real-phase
  Fourier vectors with Hermitian symmetry is legitimate by $L^2$
  Cauchy--Schwarz. It gives $\omega(S_E(t)z,S_E(t)w) = \omega(z,w)$ on all $E$.
  Equivalently one may first use Schwartz phases and then density. No rough
  generator derivative or frequency-zero trace is used.

  \medskip
  \noindent\textit{Step 2.~}We verify force skewness and compactness of Young
  endpoint differences.
  Finally, $\|bu\|_{H^{-1/2}} \le \|bu\|_2 \le \|b\|_\infty\|u\|_{H^{1/2}}$.
  Duality for these $L^2$ functions gives the absolutely integrable cancellation
  $-\int fbu+\int ubf = 0$. This is skewness for $\omega$, not for the $E$
  Hilbert inner product. The Young construction of
  Proposition~\ref{P:constant-young-linear} applies on $\mathcal{L}(E)$: its
  proof uses completeness and operator multiplication of norm one, exactly as on
  $\mathcal{L}(F)$. For fixed $s,t$, each left Young sum
  \[
    \sum_j\big(X(t_{j+1})-X(t_j)\big)U_X(t_j,s)
  \]
  is compact. The sums converge in operator norm to $U_X(t,s)-I$, so closedness
  of the compact ideal gives the first assertion
  in~\eqref{E:symplectic-compact-young-baseline}. Bounded composition then shows
  that the physical endpoint differs from $S_E(t-s)$ by a compact operator.

  \medskip
  \noindent\textit{Step 3.~}We bound the norm of a compact perturbation by the
  free norm.
  We prove the required compact-perturbation bound for this whole-space phase.
  Translations $\tau_y$ on both entries are isometries on $E$ and commute with
  $S_E(h)$. For fixed $z,w \in E$, their translated pairing is the Fourier
  transform of an $L^1$ function: the two weighted Fourier products are
  integrable by Cauchy--Schwarz in the Bessel norms. Approximation of that
  function in $L^1$ by compact smooth functions, followed by integration by
  parts for the approximants, proves $\langle\tau_yz,w\rangle_E \to 0$ as $|y|
  \to \infty$. Thus, $\tau_yz$ tends weakly to zero. For compact $K$,
  $\|K\tau_yz\| \to 0$: otherwise compactness supplies a norm-convergent
  subsequence with nonzero norm, while bounded linear testing makes its limit
  zero. For unit $z$ it follows that
  \[
    \|S_E(h)+K\|
    \ge\liminf_{|y|\to\infty}\|(S_E(h)+K)\tau_yz\|
    \ge\|S_E(h)z\|.
  \]
  Since translations commute with the free group and preserve the norm,
  $\|S_E(h)\tau_yz\|=\|S_E(h)z\|$. The compact term tends to zero on this
  weakly vanishing family. The reverse triangle inequality therefore gives the
  second inequality in the display for each fixed unit vector, before the
  supremum over vectors is taken.
  Taking the supremum over $z$ gives $\|S_E(h)+K\| \ge \|S_E(h)\|$. The
  isometric velocity reversal $\mathcal{R}(u,v) = (u,-v)$ satisfies $S_E(-h) =
  \mathcal{R} S_E(h)\mathcal{R}$ by the full free Fourier matrix. The group law
  yields $1 \le \|S_E(-h)\|\|S_E(h)\| = \|S_E(h)\|^2$. This
  proves~\eqref{E:symplectic-compact-young-baseline}.

  \medskip
  \noindent\textit{Step 4.~}We retain the compact difference through finite
  ordered products.
  Finally, $P_j-S_E(\ell_j)$ is compact. Telescoping the finite ordered product
  against its free product gives a finite sum of bounded--compact--bounded
  products, hence a compact difference. Only the free factors are combined by
  the group law; no noisy factors are commuted. The same translation argument
  proves \eqref{E:symplectic-compact-product-baseline}. A zero length gives $P_j
  = I$. All assertions are at fixed finite horizons and products; no uniform
  compact approximation in a radius, horizon or number of factors is asserted.
  This proves Lemma~\ref{L:symplectic-free}.
\end{proof}

\begin{lemma}[Covariance means and source primitives]
  \label{L:symplectic-source-mean}
  The expectations in~\eqref{E:symplectic-combined-driver} are Bochner defined
  and
  \begin{equation}\label{E:symplectic-mean-convergence}
    \|D_\Phi^\rho-D_\Phi\|_{C^\theta}
    \le\|\Phi\|_{H_W}\|M^\rho-M\|_{L^2C^\theta}\longrightarrow0 .
  \end{equation}
  For each positive $\rho$, there is a deterministic real time-continuous,
  spatially smooth $h_\Phi^\rho$ such that
  \begin{equation}\label{E:symplectic-source-primitive}
    D_\Phi^\rho(t)
    =\int_0^t S_E(-r)B(m_Rh_\Phi^\rho(r))S_E(r)\ud r .
  \end{equation}
  This is the statewise primitive and defines a Lipschitz operator path. The
  combined source coefficient is $b^\rho = m_R(cY^\rho+h_\Phi^\rho)$.
\end{lemma}

\begin{proof}
  We bound the Bochner covariance means by path moments. At positive source
  widths, the feature representation gives a smooth deterministic coefficient,
  whose ordinary rows identify the covariance primitive.

  \medskip
  \noindent\textit{Step 1.~}We bound the covariance mean and its source error.
  Use~\eqref{E:source-driver-limit} for $M^\rho\to M$ in $L^2C^\theta$,
  and~\eqref{E:source-transfer-norm} with $C=1$
  from~\eqref{E:source-two-dominations} for
  $\|M^\rho\|_{L^2C^\theta}\le\|M\|_{L^2C^\theta}$.
  Separable Gaussian path moments give $\E\|ZM\|_{C^\theta} \le
  \|\Phi\|_{H_W}\|M\|_{L^2C^\theta}$, and the same holds for sources. Bochner
  Cauchy--Schwarz gives~\eqref{E:symplectic-mean-convergence}. Covariance
  contraction also gives the uniform source bound $\|D_\Phi^\rho\|_{C^\theta}
  \le \|\Phi\|_{H_W}\|M\|_{L^2C^\theta}$.

  \medskip
  \noindent\textit{Step 2.~}We identify the regularized deterministic
  coefficient.
  The construction preserving the given noise has an isometry $j_{\rm src}:H_W
  \to \mathcal{G}$ and an isonormal process $\mathcal{I}$ with
  $\mathcal{I}(j_{\rm src}h) = W(h)$. Its source random field is $Y^\rho(t,x) =
  \mathcal{I}(v_{\rho;t,x})$, with continuous $\mathcal{G}$-valued spatial
  derivatives at positive widths. Consequently,
  \[
    h_\Phi^\rho(t,x)=\E[W(\Phi)Y^\rho(t,x)]
    =\langle j_{\rm src}\Phi,v_{\rho;t,x}\rangle_{\mathcal{G}}
  \]
  has the claimed regularity. Bounded linear testing commutes with those
  derivatives. After multiplication by $m_R$ the coefficient has one compact
  support and bounded derivatives on $[0,T]$.

  \medskip
  \noindent\textit{Step 3.~}We identify its ordinary primitive with the source
  covariance mean.
  Proposition~\ref{P:source-ordinary-response} identifies $M^\rho$ with the
  ordinary source primitive; the approximation here extends smooth time tests
  to the indicator wave tests in~\eqref{E:source-transferred-row}, with both
  limits taken in the original covariance completion.
  On Schwartz states, the row of the right side
  of~\eqref{E:symplectic-source-primitive} is the integral of $h_\Phi^\rho$
  against the explicit masked source wave product. Ordinary Fubini is justified
  by $\E|ZY^\rho(t,x)| \le \|\Phi\|_{H_W} \sqrt{\E|Y^\rho(t,x)|^2}$ and bounded
  point variance on this fixed compact set at positive widths. The source test
  and indicator approximation extends this to the completed rows by Gaussian
  Cauchy--Schwarz and covariance contraction. Hence, this row is the row of
  $\E[ZM^\rho(t)]$. Countable dense states and rational times, then continuity,
  identify the operator path. The uniformly bounded integrand gives its
  Lipschitz bound. This identifies the combined source coefficient exactly. The
  combined source is generally noncentered and the mask occurs once. This proves
  Lemma~\ref{L:symplectic-source-mean}.
\end{proof}

\begin{proposition}[Source conservation]\label{P:symplectic-source}
  At each fixed positive $\rho$, on one probability-one set for the fixed
  $c,\Phi,R,T$, the physical source flow satisfies
  \[
    P_{D^\rho}(t,s)^*J_EP_{D^\rho}(t,s)=J_E,
    \qquad 0\le s\le t\le T .
  \]
\end{proposition}

\begin{proof}
  Intersect the source coefficient regularity event with the event identifying
  its first primitive with the ordinary primitive for all times and states. The
  source construction supplies both at this fixed $\rho$. On their intersection
  $b^\rho$ has a common compact support, is continuous in supremum norm and is
  uniformly bounded on $[0,T]$. Thus, $C_b(r) = S_E(-r)B(b^\rho(r))S_E(r)$ is
  uniformly bounded and strongly continuous on states.
  Lemma~\ref{L:symplectic-free} gives $\omega(C_b(r)z,w)+\omega(z,C_b(r)w) = 0$
  for all $r,z,w$.

  Free preservation transfers the force identity to its interaction conjugate:
  applying $S_E(r)$ to both arguments turns the two terms with $C_b(r)$ into
  the two skew-force terms with $B(b^\rho(r))$. Their sum is zero. This uses
  the symplectic form, not skew-adjointness in the phase Hilbert norm.
  The source and mean primitive identities identify the Young state equation
  with $U(t,s)e = e+\int_s^t C_b(r)U(r,s)e\ud r$. Indeed, replace each operator
  path increment in a Young sum by its statewise integral. The error is bounded
  by the uniform coefficient bound times the uniform continuity error of $U$ on
  the partition. Each state trajectory is therefore $C^1$, with derivative
  $C_b(t)U(t,s)e$, by the strong-integral difference quotient. Differentiating
  the bounded bilinear form $\omega(U(t,s)e,U(t,s)f)$ gives zero. This holds on
  the same path set for all states and time pairs, since the bounded ordinary
  evolution is defined there on all $E$. Free conjugation gives the displayed
  physical identity. No uncountable intersection of state-dependent null sets or
  derivative of the free generator on a rough state was used. This proves
  Proposition~\ref{P:symplectic-source}.
\end{proof}

Proposition~\ref{P:symplectic-source} concerns the regularized flow
$P_{D^\rho}$ at a fixed positive $\rho$. We now prove the same conservation
identity for the flow $P_D$ itself, together with a lower bound on its norm.

\begin{proposition}[Conservation, inverse and a lower norm bound]
  \label{P:symplectic-raw}
  For the fixed parameters in~\eqref{E:symplectic-combined-driver}, on one
  probability-one set,
  \begin{equation}\label{E:symplectic-raw-conservation}
    P_D(t,s)^*J_EP_D(t,s)=J_E,\qquad
    \|P_D(t,s)\|\ge1,\qquad 0\le s\le t\le T .
  \end{equation}
  The flow is invertible. For $c = 0$ these assertions are deterministic for
  every fixed $\Phi$.
\end{proposition}

\begin{proof}
  We first show that every driver increment satisfies the same closed linear
  symplectic constraint. We pass that constraint to the Gaussian path and its
  covariance mean, then prove conservation by summing Young remainder
  estimates. Invertibility comes separately from the Young construction.

  \medskip
  \noindent\textit{Step 1.~}We identify the closed constraint for the
  coefficient paths.
  We first identify a closed real linear constraint on the operator path. Write
  \[
    \mathfrak{g}_J=\{A\in\mathcal{L}(E):A^*J_E+J_EA=0\}.
  \]
  The coefficient $H_{j,R}$ in~\eqref{E:driver-covariance-coefficient} is real
  and smooth. Free conservation gives $S_E(-r)^*J_E = J_ES_E(r)$ and $J_ES_E(-r)
  = S_E(r)^*J_E$. Lemma~\ref{L:symplectic-free} therefore makes
  $S_E(-r)B(H_{j,R}(r))S_E(r)$ belong to $\mathfrak{g}_J$. Testing the statewise
  integral in~\eqref{E:driver-coefficient-operator} against two states shows
  $A_j(0,t) \in \mathfrak{g}_J$ for every $j,t$. No operator-norm derivative of
  the free group is taken. The real finite sums
  in~\eqref{E:driver-path-convergence} obey the same constraint. The map
  $A\mapsto A^*J_E+J_EA$ is bounded on the path space, so their $L^2C^\theta$
  limit $M$ obeys it at all times on one probability-one path set.

  \medskip
  \noindent\textit{Step 2.~}We pass the constraint to the combined driver.
  The mean is a Bochner expectation since $\E\|W(\Phi)M\|_{C^\theta} \le
  \|\Phi\|_{H_W}\|M\|_{L^2C^\theta} < \infty$. Bounded real linear maps commute
  with this expectation, hence $D_\Phi(t) \in \mathfrak{g}_J$ for every $t$.
  Thus, $D = cM+D_\Phi$ starts at zero and is a $C^\theta$ path in
  $\mathfrak{g}_J$ on that path set. Complex Gaussian coefficients are not used.
  All these variables remain on the given $W$.

  \medskip
  \noindent\textit{Step 3.~}We prove conservation from the Young increment
  expansion.
  The proof of Proposition~\ref{P:constant-young-linear} applies on
  $\mathcal{L}(E)$, with multiplication of norm one. Fix any $s$ and put $V(r) =
  U_D(r,s)$ for $r \ge s$. The Young remainder on $[u,v] \subset [s,T]$ gives
  \[
    V(v)-V(u)=D_{u,v}V(u)+R_{u,v},\qquad
    \|R_{u,v}\|\le C|v-u|^{2\theta}.
  \]
  The exact increment of $V^*J_EV$ has first-order term
  $V(u)^*(D_{u,v}^*J_E+J_ED_{u,v})V(u) = 0$. Its remaining terms have norm at
  most $C'|v-u|^{2\theta}$, using the bounded $C^\theta$ norm of $V$. Their sum
  on a partition is at most $C'(t-s)\operatorname{mesh}^{2\theta-1}$, which
  tends to zero. Therefore, $U_D(t,s)^*J_EU_D(t,s) = J_E$. This deterministic
  argument applies to every time pair on the same operator path; it requires no
  state-dependent null sets. Free conservation in Lemma~\ref{L:symplectic-free}
  gives \eqref{E:symplectic-raw-conservation}, including $1 = \|J_E\| \le
  \|P_D(t,s)\|^2$.

  \medskip
  \noindent\textit{Step 4.~}We retain the independent justification of
  invertibility.
  Invertibility comes from the local Young Neumann pieces and their finite
  products, followed by the free inverse; it is not inferred from conservation
  in infinite dimension. For $c = 0$ the path and the argument are deterministic
  for each $\Phi$. Zero covariance gives the free flow directly. At each fixed
  $\delta$ one may take the given reference $W\circ T_{\delta,0}$, $R =
  \delta^{-1}+2$, $T = \delta^{-1}$, $c = g_\delta$ and $\Phi = 0$. No
  uniform-radius or uniform bound for controls assertion follows. This proves
  Proposition~\ref{P:symplectic-raw}.
\end{proof}

\begin{proposition}[The exact completed block mean]\label{P:symplectic-block-mean}
  Use the true partition $B,\delta,N,s_m,\ell_m$ of
  Section~\ref{S:actual-cameron-mean}, with integer $B \ge 1$. Let $i_I:H_I \to
  H_W$ be the supported closed inclusion and $T_m = T_{\delta,m}:H_{\rm fast}
  \to H_W$ the original isometry; $H_{\rm fast}$ has the same original
  covariance in fast coordinates. For $\Phi \in H_I$ set
  \[
    W_m(f)=W(T_mf),\qquad
    \phi_m=\delta^{-a/2}T_m^*i_I\Phi,\qquad
    \mathcal{T}_E\Phi=i^*T\Phi\,i ,
  \]
  where $T$ is the full $F$ mean of Proposition~\ref{P:cameron-full-mean}.
  Define
  \begin{equation}\label{E:symplectic-block-primitive}
    H_{m,E}(\Phi)(u)=\int_0^{\min(u,\ell_m)}
    S_E(-v)\mathcal{T}_E\Phi(s_m+\delta v)S_E(v)\ud v .
  \end{equation}
  With $M_m = M_R(W_m)$ and $D_{\phi_m}(u) = \E[M_m(u)W_m(\phi_m)]$, one has
  \begin{equation}\label{E:symplectic-block-direction}
    D_{\phi_m}(\min(u,\ell_m))=H_{m,E}(\Phi)(u),\qquad 0\le u\le B .
  \end{equation}
  Consequently, the pure-mean endpoint obeys
  \[
    \|S_E(\ell_m)U_{H_{m,E}(\Phi)}(\ell_m,0)\|\ge1 .
  \]
\end{proposition}

\begin{proof}
  We define the fast covariance direction using the adjoint of the original
  isometry. We identify the stopped completed mean by bounded compression,
  then remove the orthogonal complement of the fast first chaos from its
  covariance pairing.

  \medskip
  \noindent\textit{Step 1.~}We define the fast control and its completed mean
  path.
  The original pullback~\eqref{E:raw-block-pullback} has covariance factor
  $\delta^{a-2}\delta^2\delta^{-a} = 1$, hence respects nulls and extends
  isometrically from fast to slow original space. Its Hilbert adjoint is a
  contraction. Thus, $\|\phi_m\| \le \delta^{-a/2}\|\Phi\|$ is finite at fixed
  $\delta$, without a uniform control-radius assertion. The fast noise remains
  on the same given $W$.

  The compact-conjugation and primitive estimates are those in
  Proposition~\ref{P:cameron-actual-embedding},
  including~\eqref{E:cameron-local-primitive-bound}, restricted from $F$ to
  $E$: compactness gives norm-continuous conjugation and Hölder integration
  gives exponent $1-1/r$.
  The mean $T\Phi$ has range in closed $iE$. By
  Proposition~\ref{P:cameron-full-mean}, bounded compression passes its strong
  measurability, $L^r$ bound and compact operator values to $\mathcal{T}_E\Phi
  \in L^r(I;\mathcal{K}(E))$. Choose the full parameter
  choices~\eqref{E:cameron-block-parameters}. Its intervals are nonempty because
  $\theta < 1-a/2$ and $\theta > 1/2$; $p$ only instantiates that package here.
  The compact-conjugation argument makes \eqref{E:symplectic-block-primitive} a
  Bochner integral, and the $L^r$ primitive estimate gives $C^\theta$, since
  $1-1/r > \theta$. Affine pullback preserves null sets. Its integrand is
  compact valued almost everywhere; closedness of $\mathcal{K}(E)$ under Bochner
  integration makes $H_{m,E}(\Phi)$ a compact-valued path starting at zero, also
  after constant stopping.

  \medskip
  \noindent\textit{Step 2.~}We compute the unscaled covariance by compressing
  the full phase.
  Here $A_{F,m}$ is the full $F$ interaction operator path from
  Proposition~\ref{P:constant-phase-driver}, and $H_{m,F}$ is its stopped
  completed mean primitive in~\eqref{E:cameron-actual-covariance}. The completed
  identity~\eqref{E:cameron-actual-covariance} at original projection rank zero
  gives, for $0 \le u \le \ell_m$,
  \[
    \E[g_\delta A_{F,m}(u)W(i_I\Phi)]
    =\sqrt\delta\,H_{m,F}(\Phi)(u).
  \]
  Compress through $i$ using the operator path intertwining, $S_Fi = iS_E$, and
  the mean range in $iE$ preserved by $S_F$. Bochner expectation commutes with
  this bounded compression. The constant columns are retained in the full $F$
  operator path and mean until this compression; they are inactive only on the
  $iE$ inputs. Since $g_\delta\delta^{a/2} = \sqrt\delta$ and $g_\delta > 0$,
  division gives the unscaled covariance
  \begin{equation}\label{E:symplectic-unscaled-covariance}
    \E[M_m(u)W(i_I\Phi)]
    =\delta^{a/2}H_{m,E}(\Phi)(u),\qquad u\le\ell_m .
  \end{equation}
  The factor here is $\delta^{a/2}$, not $\sqrt\delta$. Completion in the
  earlier covariance identity permits every $\Phi \in H_I$, without positivity
  or smoothness.

  \medskip
  \noindent\textit{Step 3.~}We identify the covariance direction in the fast
  noise space.
  The operator $T_mT_m^*$ is the orthogonal projection onto the original
  isometry's closed range. Thus, $W((I-T_mT_m^*)i_I\Phi)$ is orthogonal to the
  entire fast first chaos. Every continuous linear coordinate of $M_m$ is an
  $L^2$ limit of finite $W_m$-basis sums. Its covariance with this difference
  vanishes. Bochner covariance exists by $L^2$ path moments, and bounded linear
  functionals separate its values. Hence,
  \[
    \E[M_m(u)W(i_I\Phi)]
    =\E[M_m(u)W_m(T_m^*i_I\Phi)] .
  \]
  Multiply by $\delta^{-a/2}$ and use \eqref{E:symplectic-unscaled-covariance}
  to obtain \eqref{E:symplectic-block-direction} on $[0,\ell_m]$, then stop both
  sides. No unstopped equality beyond $\ell_m$ is asserted.

  \medskip
  \noindent\textit{Step 4.~}We apply the compact endpoint lower bound.
  Apply Lemma~\ref{L:symplectic-free} directly to the compact-valued path
  $H_{m,E}(\Phi)$ starting at zero on $[0,B]$. Its physical endpoint has norm at
  least $\|S_E(\ell_m)\| \ge 1$. If $\ell_m = 0$, the endpoint is the identity
  and its norm is one. The radius may be chosen $\delta^{-1}+2$ after fixing
  $\delta$; all identities are fixed-parameter instances. This proves
  Proposition~\ref{P:symplectic-block-mean}.
\end{proof}

For zero covariance the operator paths and covariance means vanish, and the
claims reduce directly to the free identity in Lemma~\ref{L:symplectic-free}.
These results give a lower bound on the operator norm and its exact completed
mean application. They provide no uniform bound on the control radius,
stochastic moment upper bound or nonlinear product ratio by themselves.

\section{Gaussian enlargement and the Rayleigh comparison}
\label{S:actual-gaussian-tubes}

We first state the Gaussian geometry with its exact covariance scale. Let $H$ be
real separable Hilbert, $W$ its given isonormal process, and $X$ a strongly
measurable centered Gaussian variable in a real separable Banach space
$\mathscr{B}$, with $\E\|X\|^2 < \infty$. Assume every continuous linear
coordinate has the form $\ell(X) = W(k_\ell)$, $k_\ell \in H$. Fix $0 < \delta <
1$ and write
\begin{equation}\label{E:tube-covariance-scale}
  Ch=\E[XW(h)]=\sqrt\delta\,\mathcal{M}_\delta h,\qquad
  H_{\rm seen}=\overline{\operatorname{span}\{k_\ell\}},\qquad
  E_0=\overline{C(H)} .
\end{equation}
Assume $\mathcal{M}_\delta:H \to \mathscr{B}$ is compact. All Hilbert balls
below are closed.

\subsection{The Wiener-space realization}

The Gaussian isoperimetric theorem that we import in
Section~\ref{S:gaussian-isoperimetry-input} is stated for an abstract Wiener
space. We first show that the law of $X$ on $E_0$ is of this kind, and we
identify its Cameron space.

\begin{proposition}[First Wiener chaos realization on the support]
  \label{P:tube-wiener-realization}
  The variable $X$ lies in $E_0$ almost surely, with a centered Radon Gaussian
  law of full support there. Its Cameron space and norm are
  \begin{equation}\label{E:tube-cameron-space}
    H_{\rm CM}=C(H_{\rm seen}),\qquad
    \|Ch\|_{\rm CM}=\|P_{\rm seen}h\|_H .
  \end{equation}
  This support triple is an abstract Wiener realization, including the
  measurable-norm condition.
\end{proposition}

\begin{proof}
  We reconstruct the Gaussian variable from its finite original-noise
  coordinates. We identify its Cameron image and support, and then verify the
  measurable-norm condition required by the imported abstract Wiener theorem.

  \medskip
  \noindent\textit{Step 1.~}We recover the variable from finite Gaussian
  regressions.
  The zero law has $C = 0$ and $E_0 = \{0\}$ and is immediate. A countable
  norming family exists in $\mathscr{B}^*$: choose a dense sequence on the unit
  sphere and apply Hahn--Banach at each point. Approximation gives $\|x\| =
  \sup_j|\ell_j(x)|$. Norms of differences from a countable dense set show that
  these coordinates generate the Borel sets. Consequently, $X$ is measurable in
  the completion of $\sigma(W(H))$. Fix an orthonormal basis $(e_j)$ of $H$ and
  set $\mathcal{G}_j = \sigma(W(e_1),\ldots,W(e_j))$. In finite dimension these
  sums and projections stabilize at the dimension. Bochner conditional
  expectations need no Radon--Nikodym property here: define them first on
  finite-valued simple variables using scalar conditional expectations, then
  extend by the $L^2$ norm contraction. Simple approximation supplies their
  scalar testing identities. Gaussian regression and the countable norming
  family give
  \begin{equation}\label{E:tube-finite-regression}
    X_j=\E[X\mid\mathcal{G}_j]
    =\sum_{l\le j}W(e_l)Ce_l,\qquad X_j\longrightarrow X
    \quad\hbox{in }L^2(\mathscr{B}).
  \end{equation}
  To justify convergence, finite-coordinate scalar events approximate all events
  in the generated sigma field: the class of events whose indicators are
  approximable in $L^2$ is a monotone class containing the finite-coordinate
  algebra. Approximate a Bochner variable by simple variables and use the
  conditional $L^2$ contraction.

  \medskip
  \noindent\textit{Step 2.~}We identify the Cameron image and the closed
  support space.
  Bounded testing gives $\ell(Ch) = \langle k_\ell,h\rangle$. Thus, $\ker C =
  H_{\rm seen}^\perp$, and $C$ is injective on $H_{\rm seen}$. The closed first
  Wiener chaos generated by $X$ is exactly $W(H_{\rm seen})$. Its covariance
  image therefore has the first Wiener chaos norm
  in~\eqref{E:tube-cameron-space}; unused directions are removed by the
  orthogonal projection, without an inverse covariance. This is the Cameron
  image under the first Wiener chaos definition. The finite sums
  in~\eqref{E:tube-finite-regression} lie in $E_0$, so $L^2$ convergence and an
  almost-sure subsequence put $X$ in that closed space almost surely.
  Functionals on $E_0$ extend to $\mathscr{B}$ by Hahn--Banach, so restriction
  preserves the scalar laws and first chaos.

  \medskip
  \noindent\textit{Step 3.~}We prove the Radon and full-support assertions.
  Its law is tight. For each radius $2^{-j}$ choose a finite union of balls of
  mass greater than $1-\varepsilon2^{-j-1}$, using a countable cover of the
  separable space. The intersection of the corresponding finite closed unions is
  closed and totally bounded and has mass at least $1-\varepsilon$. Completeness
  makes it compact. Metric Borel regularity on these compact sets, followed by
  tightness, gives a Radon law. To prove full support, the tail $X-X_j$ is
  independent of the finite Gaussian sum: every tail coordinate has label
  $(I-P_j)k_\ell$, orthogonal to $e_1,\ldots,e_j$. Joint Gaussianity and
  countable Borel coordinates give Banach tail independence. Given $h$ and
  $\varepsilon > 0$, choose $j$ with $\|CP_jh-Ch\| < \varepsilon/3$ and
  $\mathbb{P}(\|X-X_j\| < \varepsilon/3) > 0$. The finite sum has positive
  probability in every ball about $CP_jh$, even if its covariance is singular:
  the preimage under its finite linear map contains an open ball about the
  coefficient vector of $P_jh$. Independence proves positivity of every ball
  about $Ch$. Closure proves that the support is all of $E_0$.

  \medskip
  \noindent\textit{Step 4.~}We verify the finite-dimensional tail criterion
  for the measurable norm.
  We check the measurable norm explicitly. In $H_{\rm seen}$ put $V_j =
  \operatorname{span}\{P_{\rm seen}e_l:l\le j\}$. Any finite orthonormal family
  $(h_k)$ in $V_j^\perp$ is also orthogonal to $e_1,\ldots,e_j$. The same
  regression and Bochner Jensen argument gives
  \[
    \E[X-X_j\mid W(h_k),\hbox{ all }k]
    =\sum_k W(h_k)Ch_k,\qquad
    \E\left\|\sum_kW(h_k)Ch_k\right\|^2
    \le\E\|X-X_j\|^2 .
  \]
  For any $\varepsilon > 0$ choose $j$ so that the right side is less than
  $\varepsilon^3$. Markov bounds the probability that this finite sum has norm
  exceeding $\varepsilon$ by $\varepsilon$, uniformly over every such finite
  orthonormal family. Since $C:H_{\rm seen} \to H_{\rm CM}$ is isometric, the
  excluded finite-dimensional Cameron subspace is $C(V_j)$. This is the
  measurable-norm tail condition for the dense Hilbert image, whose Banach
  completion is $E_0$. The characteristic functional of the law is $\E
  e^{i\ell(X)} = \exp(-\|k_\ell\|_H^2/2)$. These facts give the stated abstract
  Wiener realization. The zero space has the direct zero-law interpretation.
  This proves Proposition~\ref{P:tube-wiener-realization}.
\end{proof}

\subsection{Closed Gaussian enlargements and the quadratic cost}
\label{S:gaussian-isoperimetry-input}

The Gaussian isoperimetric theorem is due to Borell and, independently,
Sudakov--Tsirelson~\cites{borell:75:brunn-minkowski, sudakov.cirel-son:74:extremal}.
We use its set-enlargement form from
Friz and Hairer~\cite{friz.hairer:20:course}*{Theorem~11.6, p.~190}, with footnote~3 on
p.~190; p.~191 supplies the following context. For an abstract Wiener space, a
Borel $A$ of positive mass and its Cameron unit ball $\mathcal{K}$, it states
\[
  \mu((A+t\mathcal{K})^c)
  \le\overline{\Phi}_{\rm N}(\widehat{a}+t),\qquad
  \widehat{a}=\Phi_{\rm N}^{-1}(\mu(A)),\quad t\ge0 .
\]
Here $\Phi_{\rm N}$ is the standard normal distribution function; the same page
gives $\overline{\Phi}_{\rm N}(y) \le e^{-y^2/2}$ for $y \ge 0$. This theorem is
imported, not proved here. It concerns Gaussian input sets, not Gaussianity of a
nonlinear wave. We verify Borel measurability of every enlargement we use.

\begin{lemma}[Compact images and attained tube radii]
  \label{L:tube-closed-geometry}
  Suppose the law of $X$ is nonzero and put
  \[
    b_\delta=2\E\|X\|,\qquad
    A_\delta=\{x\in E_0:\|x\|\le b_\delta\},\qquad
    \mathcal{T}_r=A_\delta+\mathcal{M}_\delta(\overline{B}_H(r)).
  \]
  Every $\mathcal{T}_r$ is closed. The extended distance $d_\delta(x) =
  \inf\{r\ge0:x\in\mathcal{T}_r\}$ satisfies
  \begin{equation}\label{E:tube-attained-radius}
    \{d_\delta\le r\}=\mathcal{T}_r,\qquad
    \mathbb{P}(d_\delta(X)>r)\le e^{-r^2/(2\delta)},\quad r\ge0 .
  \end{equation}
  In particular, $d_\delta$ is lower semicontinuous and finite almost surely.
\end{lemma}

\begin{proof}
  We first prove that compact control images make the tubes closed. Weak
  compactness also attains every finite infimal radius. We then apply Gaussian
  enlargement at the covariance scale fixed in the statement.

  \medskip
  \noindent\textit{Step 1.~}We prove closedness and attainment of finite
  radii.
  A compact map from a separable Hilbert space sends a closed bounded ball to a
  compact set, including closedness of the image. Indeed, diagonal orthonormal
  coordinates and Fatou give a weakly convergent subsequence of any bounded
  sequence, with limit still in the ball. Compactness gives norm-convergent
  image subsequences, and bounded scalar testing identifies every limit with the
  image of that weak limit. This proves sequential compactness. Closed
  $A_\delta$ plus a compact image is closed: along a converging sequence of
  sums, extract the compact component and pass the residual into the closed base
  set.

  Membership in $\mathcal{T}_r$ implies $d_\delta \le r$. Conversely, if
  $d_\delta(x) \le r$, choose $h_j$ with $\|h_j\| \le r+1/j$ and
  $x-\mathcal{M}_\delta h_j \in A_\delta$. A weak subsequence has limit $h$ of
  norm at most $r$; compactness gives $\mathcal{M}_\delta h_j \to
  \mathcal{M}_\delta h$. Closedness of the residual proves membership in
  $\mathcal{T}_r$. This also covers $r = 0$ and proves least finite-radius
  attainment. Closed sublevels give lower semicontinuity and Borel
  measurability.

  \medskip
  \noindent\textit{Step 2.~}We use the exact Cameron scale in the enlargement
  inequality.
  The nonzero law has $b_\delta > 0$, and Markov gives $\mu(A_\delta) \ge 1/2$.
  Projection onto $H_{\rm seen}$ removes all unused directions in the original
  covariance space, so~\eqref{E:tube-cameron-space} gives the exact set
  identities
  \[
    \mathcal{K}=C(\overline{B}_H(1))
    =\sqrt\delta\,\mathcal{M}_\delta(\overline{B}_H(1)),\qquad
    \mathcal{T}_r=A_\delta+(r/\sqrt\delta)\mathcal{K} .
  \]
  These are Borel enlargements by the preceding closedness. The imported theorem
  applies to the realization in Proposition~\ref{P:tube-wiener-realization}.
  Since $\widehat{a} \ge 0$, it yields~\eqref{E:tube-attained-radius}. The value
  $\widehat{a} = +\infty$ is allowed when $\mu(A_\delta) = 1$. The tail tends to
  zero, proving almost-sure finiteness. The Cameron embedding does not change
  the factor one half in the Gaussian cost. This proves
  Lemma~\ref{L:tube-closed-geometry}.
\end{proof}

Lemma~\ref{L:tube-closed-geometry} bounds the tail of the tube distance
$d_\delta(X)$. We now convert this tail bound into a bound for
$\E e^{F(X)/\delta}$, for a functional $F$ that satisfies $F \le G(r)$ on each
tube $\mathcal{T}_r$.

\LeanMarker{P:tube-rayleigh-transfer}
\begin{proposition}[Comparison with the Rayleigh integral, including jumps and zero law]
  \label{P:tube-rayleigh-transfer}
  Let $F:E_0 \to \R$ be Borel and let $G:[0,\infty) \to \R$ be Borel and
  nondecreasing. If $F(x) \le G(r)$ for every $x \in \mathcal{T}_r$ and every $r
  \ge 0$, then
  \begin{equation}\label{E:tube-rayleigh-transfer}
    \E e^{F(X)/\delta}
    \le\int_0^\infty\frac r\delta
    \exp\left(\frac{G(r)-r^2/2}{\delta}\right)\ud r .
  \end{equation}
  For zero law use $E_0 = \{0\}$, $C = \mathcal{M}_\delta = 0$; the same
  conclusion holds when $F(0) \le G(0)$.
\end{proposition}

\begin{proof}
  We compare upper-level probabilities with those of the scalar Rayleigh law,
  treating both open and closed endpoints. Layer cake then gives the
  exponential integral bound, with attained radius removing the need for a
  measurable control choice.

  The auxiliary Rayleigh variable $\mathcal{R}_\delta$ with density
  $(r/\delta)e^{-r^2/(2\delta)}$, $r > 0$, has mass one and survival function
  $e^{-r^2/(2\delta)}$ by direct integration. It is a scalar stochastic
  dominator, not a norm of a Gaussian variable in the covariance Hilbert space.
  For every level $t \ge 0$, the upper level set of $e^{G(r)/\delta}$ is an
  upper interval, possibly empty or the whole halfline. An open endpoint is
  controlled by~\eqref{E:tube-attained-radius}. At a positive closed endpoint
  $c$, let $s\uparrow c$:
  \[
    \mathbb{P}(d_\delta(X)\ge c)
    =\lim_{s\uparrow c}\mathbb{P}(d_\delta(X)>s)
    \le e^{-c^2/(2\delta)}
    =\mathbb{P}(\mathcal{R}_\delta\ge c).
  \]
  Rayleigh has no atoms, and at the closed endpoint zero both probability bounds
  are one. On the null event $\{d_\delta(X)=\infty\}$ define
  $e^{G(d_\delta(X))/\delta}$ to be one. No value of $G$ at infinity is needed
  and no expectation changes. Layer cake and Tonelli therefore order the
  expectations of $e^{G(d_\delta(X))/\delta}$ and
  $e^{G(\mathcal{R}_\delta)/\delta}$, even when $G$ jumps or expectations are
  infinite. Least-radius attainment gives $F(X) \le G(d_\delta(X))$ almost
  surely. No measurable selection of a control or residual is needed. This
  proves~\eqref{E:tube-rayleigh-transfer}. For zero law, $X = 0$ and $F(0) \le
  G(0) \le G(r)$; normalization of the Rayleigh density proves the claim
  directly. In particular, $F(0)$ need not be zero. This proves
  Proposition~\ref{P:tube-rayleigh-transfer}.
\end{proof}

\subsection{The wave tuple and its globally defined endpoint}

Retain all original covariances and amplitudes. Use the exponents and fixed even
$Q$ in~\eqref{E:radius-anchor-parameters}--\eqref{E:radius-anchor-Q}; put $\nu =
1/\theta$ and choose $(1-\theta)^{-1} < r_t < 2/a$, $\nu \le p < 2 < r_t$. All
these intervals are nonempty by the strict subcritical assumptions. Fix integer
$B \ge 1$ and the true partition $\delta,N,s_m,\ell_m$ from
Section~\ref{S:actual-cameron-mean}. Set $L_\delta = \delta^{-1}$ and eventually
$R_\delta = L_\delta+2$. On the same given slow $W$ define the full, unprojected
input
\begin{equation}\label{E:tube-actual-input}
  W_m(f)=W(T_{\delta,m}f),\qquad
  X_m(u)=g_\delta M_R(W_m)(\min(u,\ell_m)),\qquad
  g_\delta=\delta^{(1-a)/2}.
\end{equation}
Write $P_{\delta,E}^R$ for the physical $E$ flow driven by $g_\delta M_R(W\circ
T_{\delta,0})$ on the global fast horizon $[0,L_\delta]$. Let
$\mathscr{B}_\delta$ be the finite product of little
$C^\theta([0,B];\mathcal{K}(E))$ spaces of paths starting at zero, with exact
norm
\[
  \|x\|_{\delta,p}=\left(\delta\sum_{m=0}^N
  \|x_m\|_{C^\theta}^p\right)^{1/p}.
\]

\begin{proposition}[Eligibility of the E law]\label{P:tube-actual-law}
  The tuple $X_\delta = (X_m)$ satisfies the preceding Gaussian hypotheses with
  $H = H_I$ and
  \begin{align}
    (\mathcal{M}_\delta\Phi)_m(u)
    &=\int_0^{\min(u,\ell_m)}
    S_E(-v)\mathcal{T}_R\Phi(s_m+\delta v)S_E(v)\ud v,
    \label{E:tube-actual-mean}\\
    \E[X_\delta W(\Phi)]&=\sqrt\delta\,\mathcal{M}_\delta\Phi .
    \label{E:tube-actual-covariance}
  \end{align}
  Its support is $E_{0,\delta} = \overline{C_\delta(H_I)} \subset
  \mathscr{B}_\delta$, where $C_\delta\Phi = \E[X_\delta W(\Phi)] =
  \sqrt\delta\,\mathcal{M}_\delta\Phi$. The map $\mathcal{M}_\delta$ is compact
  at each fixed radius and parameters. At $R = R_\delta$,
  \begin{equation}\label{E:tube-typical-radius}
    b_\delta=2\E\|X_\delta\|_{\delta,p}
    \le C_{Q,B,p}\delta^{(1-a)/2-d/Q}\longrightarrow0 .
  \end{equation}
\end{proposition}

\begin{proof}
  We check the separable Gaussian tuple and compact completed mean at fixed
  parameters. We then identify its covariance on the supported original
  Hilbert space and apply the explicit radius moment bound to its typical
  norm.

  \medskip
  \noindent\textit{Step 1.~}We check the measurable Gaussian state space.
  Compact operators on separable Hilbert $E$ form a separable operator-norm
  space: finite-rank rational matrices in an orthonormal basis are dense.
  Polygonal paths give separability of the little H\"older closure. The
  construction is a $C^\theta$ limit of Lipschitz finite coefficient paths
  starting at zero; stopping at $\ell_m$ and extending constantly preserves this
  space and gives zero for an empty final interval. Gaussian transfer retains
  the given $W$ for all finitely many blocks. Every finite collection of linear
  statistics is an $L^2$ limit of joint first Wiener chaos Gaussians, hence
  jointly centered Gaussian. The finite tuple has finite $L^2$ norm by the
  fixed-parameter path moments in Corollary~\ref{C:radius-covariance-residuals}.
  The original Hilbert space $H_I$ is separable by
  Lemma~\ref{L:cameron-supported-tensor}. No separability of all
  $\mathcal{L}(E)$ or full $C^\theta$ is claimed.

  \medskip
  \noindent\textit{Step 2.~}We prove compactness of the fixed-parameter mean
  map.
  The $E$ mean field $\mathcal{T}_R$ is the compression of
  Proposition~\ref{P:cameron-full-mean}; it is compact into
  $L^{r_t}(I;\mathcal{K}(E))$ at fixed $R$.
  Proposition~\ref{P:uniform-mean-regulator} gives $\|\mathcal{T}_R\| \le C_{\rm
  mean} = M_0\|C_t\|$ independently of $R \ge 1$. Compact conjugation, affine
  time pullback and $1-1/r_t > \theta$ make the primitive
  in~\eqref{E:tube-actual-mean} bounded from that field space to
  $\mathscr{B}_\delta$, with constant stopping. Composition proves compactness
  of $\mathcal{M}_\delta$. There is no uniform-radius compactness assertion.

  \medskip
  \noindent\textit{Step 3.~}We identify the supported covariance and typical
  radius.
  The rank-zero identity~\eqref{E:cameron-actual-covariance}, compressed to $E$
  as in Proposition~\ref{P:symplectic-block-mean},
  gives~\eqref{E:tube-actual-covariance}. Its scale is exactly
  $g_\delta\delta^{a/2} = \sqrt\delta$. Every stopped row's slow pullback lies
  in $H_I$ by the inside-$I$ endpoint approximation of
  Lemma~\ref{L:cameron-actual-residual}. Countable matrix coordinates and
  rational times identify the tuple in the closed first Wiener chaos $W(H_I)$,
  so every linear statistic has its label there. Equivalently, apply that lemma
  at $J = 0$ on $F$ and compose any $E$-tuple statistic with the bounded
  $E$-compression. This is same-$W$ covariance, without a block-independence
  assumption or a time-indicator multiplier on the covariance space.
  Proposition~\ref{P:tube-wiener-realization} therefore applies to this exact
  law on its support. Finally, Proposition~\ref{P:radius-mean-anchor}, with zero
  residual projection, gives~\eqref{E:tube-typical-radius}. The third fixed-$Q$
  condition $Q > 2d/(1-a)$ makes this exponent positive; the first two supply
  the path moments and chaining. Fast directions used for the symplectic
  baseline do not replace this global $H_I$ action or get summed as separate
  costs. Zero covariance gives zero operator paths and means directly. This
  proves Proposition~\ref{P:tube-actual-law}.
\end{proof}

\begin{lemma}[The endpoint functional]\label{L:tube-actual-endpoint}
  For every $x \in \mathscr{B}_\delta$ define
  \[
    \begin{aligned}
      P_m(x_m)&=S_E(\ell_m)U_{x_m}(\ell_m,0),\qquad
      \mathcal{P}_\delta(x)=P_N(x_N)\cdots P_0(x_0),\\
      F_\delta(x)&=2\delta\log\max(1,\|\mathcal{P}_\delta(x)\|).
    \end{aligned}
  \]
  This is a finite continuous functional on all tuples. For the input it
  satisfies
  \begin{equation}\label{E:tube-global-endpoint}
    \mathcal{P}_\delta(X_\delta)
    =P_{\delta,E}^R(L_\delta,0),\qquad
    e^{F_\delta(X_\delta)/\delta}
    =\|P_{\delta,E}^R(L_\delta,0)\|^2 .
  \end{equation}
  These are almost-sure identities for each fixed admissible $\delta,R,B$,
  without a simultaneous event over all parameter choices.
\end{lemma}

\begin{proof}
  We define the endpoint functional on the entire tuple space, including
  inputs outside any Gaussian full-measure set. We then use the common-source
  restart identity and the compact-product lower bound to identify its value
  on the actual input.

  \medskip
  \noindent\textit{Step 1.~}We obtain continuity and the deterministic
  endpoint lower bound.
  Generic Young evolution is globally defined and locally Lipschitz in every
  $C^\theta$ operator path. Finite products, fixed-time free factors, norm and
  $\log\max(1,\cdot)$ preserve continuity. Thus, $F_\delta$ is defined even off
  every Gaussian full-measure set. Every $x_m$ in this ambient space is a
  compact-valued $C^\theta$ path starting at zero. The finite-product assertion
  of Lemma~\ref{L:symplectic-free} gives, for every such tuple,
  \[
    \|\mathcal{P}_\delta(x)\|
    \ge\left\|S_E\Big(\sum_{m=0}^N\ell_m\Big)\right\|\ge1 .
  \]

  \medskip
  \noindent\textit{Step 2.~}We identify the global endpoint from the coupled
  local factors.
  Lemma~\ref{L:raw-block-restart}, in its case without the deterministic mean,
  applies to this same full input with base physical profile $\eta^2$. To check
  the exact objects, use its single regularized slow noise preserving the given
  noise and pull it back by $T_{\delta,0}$ and $T_{\delta,m}$. At fixed $\delta$
  the global fast source random field is $\delta^{a/2}Y^\rho(\delta u,x)$ and
  its local restriction is its time shift by $mB$. Their ordinary Duhamel flows
  have the same massless free groups, amplitude $g_\delta$ and mask $m_R$.
  Source limits are taken on the fixed finite extended horizon $(N+1)B$ and on
  each local $B$ before composing the finitely many factors. The original
  cocycle from that lemma therefore gives the ordered physical $F$ product with
  the true remainder. The full operator paths and free groups intertwine through
  the closed isometry $i:E \to F$, including the source limits. Restricting the
  product and using injectivity of $i$ gives the first identity
  in~\eqref{E:tube-global-endpoint}. Both fast velocity and displacement remain
  in each factor. No mask was commuted through a free group, and the time blocks
  remain correlated on the same $W$. For zero covariance every input is zero and
  the ordered factors are $S_E(\ell_m)$, whose product is $S_E(L_\delta)$.

  The lower bound already established for every tuple removes the maximum at
  $X_\delta$. Together with the first identity, this proves the second identity
  in~\eqref{E:tube-global-endpoint}. No fixed-radius Gaussian-operator rate or
  source-to-point result is an input to this restart argument. This proves
  Lemma~\ref{L:tube-actual-endpoint}.
\end{proof}

\subsection{Uniform coarse envelope and the exact radial integral}

To apply Proposition~\ref{P:tube-rayleigh-transfer} to the endpoint functional
$F_\delta$ of Lemma~\ref{L:tube-actual-endpoint}, we need a nondecreasing Borel
function that dominates $F_\delta$ on every tube. We take the supremum of
$F_\delta$ over each tube and bound its growth in the radius, uniformly in
$\delta$ and $R \ge 1$.

\begin{proposition}[Coarse tube envelope and operator moment]
  \label{P:tube-coarse-envelope}
  On the support $E_{0,\delta}$ let $A_\delta$ be its closed ball of radius
  $b_\delta$ and put
  \[
    G_\delta(r)=
    \sup_{\substack{a\in A_\delta\\\|\Phi\|_{H_I}\le r}}
    F_\delta(a+\mathcal{M}_\delta\Phi).
  \]
  This is finite, nondecreasing and Borel. There is $C_{\rm env} < \infty$,
  independent of $\delta$ and $R \ge 1$, such that
  \begin{equation}\label{E:tube-coarse-envelope}
    G_\delta(r)\le C_{\rm env}(1+b_\delta^\nu+r^\nu),\qquad
    \E\|P_{\delta,E}^R(L_\delta,0)\|^2
    \le\int_0^\infty\frac r\delta
    e^{(G_\delta(r)-r^2/2)/\delta}\ud r .
  \end{equation}
\end{proposition}

\begin{proof}
  We first bound the sum of deterministic mean primitive norms. We combine
  this with the generic Young growth estimate and the exact weighted power
  means to obtain a finite monotone tube envelope, then apply the Rayleigh
  comparison.

  \medskip
  \noindent\textit{Step 1.~}We bound the deterministic mean primitives on the
  true partition.
  Put $G_B = \sup_{|v|\le B}\|S_E(v)\|$ and $C_{\rm pr} =
  G_B^2[B^{1-1/r_t}+B^{1-1/r_t-\theta}]$. The lower bound $G_B \ge 1$ follows
  directly from $S_E(0) = I$. Proposition~\ref{P:uniform-actual-cameron} gives
  $1 \le G_B \le 3B$, $C_{\rm pr} \le 18B^{3-1/r_t}$. Its primitive estimate,
  now without the covariance factor $\sqrt\delta$, gives
  \begin{equation}\label{E:tube-primitive-sum}
    \delta\sum_m h_m^{r_t}
    \le(C_{\rm pr}C_{\rm mean})^{r_t}\|\Phi\|_{H_I}^{r_t},
    \qquad h_m=\|(\mathcal{M}_\delta\Phi)_m\|_{C^\theta}.
  \end{equation}
  Indeed, each primitive norm is at most $C_{\rm pr}\delta^{-1/r_t}
  \|\mathcal{T}_R\Phi\|_{L^{r_t}(I_m)}$. Raising to $r_t$ cancels the $\delta$
  weight, and the true partition sums to the full field norm, including an empty
  cell.

  \medskip
  \noindent\textit{Step 2.~}We control the endpoint on every tube.
  The generic Banach Young bound supplies $c_Y = c_Y(B,\theta) > 0$ with
  $\log\|U_d(u,0)\| \le c_Y(1+\|d\|_{C^\theta}^{\nu})$ for $u \le B$. This
  follows from~\eqref{E:young-upper-refined-growth}
  and~\eqref{E:constant-young-growth} by $d \le 1+d^\nu$ and $(1+d)^\nu \le
  2^{\nu-1}(1+d^\nu)$. It applies on $\mathcal{L}(E)$ with multiplication norm
  one, without Gaussian, positivity or Hamiltonian assumptions on $d$. For $a
  \in A_\delta$ put $v_m = \|a_m\|_{C^\theta}$ and $n = N+1$. Convexity and
  submultiplicativity give
  \[
    F_\delta(a+\mathcal{M}_\delta\Phi)
    \le2\delta n(\log G_B+c_Y)
    +2^\nu c_Y\delta\sum_m(v_m^\nu+h_m^\nu).
  \]
  The upper expression is nonnegative, so it also bounds the logarithm with the
  maximum at one. The exact finite weighted power-mean inequalities are
  \begin{align*}
    \delta\sum_m v_m^\nu
    &\le(\delta n)^{1-\nu/p}\|a\|_{\delta,p}^\nu,\\
    \delta\sum_m h_m^\nu
    &\le(\delta n)^{1-\nu/r_t}
    \left(\delta\sum_m h_m^{r_t}\right)^{\nu/r_t}.
  \end{align*}
  They use exponents $p/\nu \ge 1$, $r_t/\nu > 1$ and nonnegative inputs
  $v_m^\nu,h_m^\nu$, followed by positive powers. Zero norms and empty cells
  cause no negative powers; $\delta > 0$ and $n \ge 2$, including the case $p =
  \nu$. Since $\delta n \le 2/B$, these bounds and~\eqref{E:tube-primitive-sum}
  give $F_\delta(a+\mathcal{M}_\delta\Phi) \le D_0+D_1b_\delta^\nu+D_2r^\nu$,
  where
  \begin{align*}
    D_0&=2(2/B)(\log G_B+c_Y),\\
    D_1&=2^\nu c_Y(2/B)^{1-\nu/p},\\
    D_2&=2^\nu c_Y(2/B)^{1-\nu/r_t}
    (C_{\rm pr}C_{\rm mean})^\nu .
  \end{align*}
  Hence, $C_{\rm env} = 1+D_0+D_1+D_2$ is independent of $\delta,R$ and proves
  the coarse bound for every admissible $\delta$, not just when $b_\delta$ is
  small. The supremum sets are nonempty and this bound makes their suprema
  finite. Monotonicity gives Borel measurability without asserting attainment of
  these suprema. The definition bounds $F_\delta$ on every closed tube. Apply
  Proposition~\ref{P:tube-rayleigh-transfer} and~\eqref{E:tube-global-endpoint}
  to obtain the operator-moment bound. For zero law its direct case applies,
  with no assumption that the free endpoint norm is one. This proves
  Proposition~\ref{P:tube-coarse-envelope}.
\end{proof}

\LeanMarker{L:tube-rayleigh-tail}
\begin{lemma}[Exact Rayleigh tail bound]\label{L:tube-rayleigh-tail}
  Suppose $G$ is Borel and $G(r) \le A_0+Cr^\nu$, $C \ge 0$, $1 < \nu < 2$. If
  $M > 0$ and $4C \le M^{2-\nu}$, then
  \begin{equation}\label{E:tube-rayleigh-tail}
    \int_M^\infty\frac r\delta
    e^{(G(r)-r^2/2)/\delta}\ud r
    \le4e^{(A_0-M^2/4)/\delta}.
  \end{equation}
  In particular, the operator-moment bound in~\eqref{E:tube-coarse-envelope} is
  finite at fixed parameters.
\end{lemma}

\begin{proof}
  The finite scalar cutoff inequality gives, for $r \ge M$, $Cr^\nu-r^2/2 \le
  -M^2/8-r^2/8$. Its exact conditions are $M > 0$, $r \ge M$, $C \ge 0$, $1 <
  \nu < 2$ and $4C \le M^{2-\nu}$. Indeed, $r^{2-\nu} \ge M^{2-\nu} \ge 4C$
  gives $Cr^\nu \le r^2/4$, while $r^2 \ge M^2$ splits the remaining $-r^2/4$
  into the two displayed negative terms. Analytically, $-\frac{\ud}{\ud
  r}[4e^{-r^2/(8\delta)}] = (r/\delta)e^{-r^2/(8\delta)}$. Integrating gives the
  front factor $4e^{-M^2/(8\delta)}$; multiplication by $e^{(A_0-M^2/8)/\delta}$
  proves \eqref{E:tube-rayleigh-tail}. Such a cutoff exists since $2-\nu > 0$.
  On $[0,M]$ the growth bound gives $G(r) \le A_0+CM^\nu$, so the integral is
  finite there as well. Apply this to $A_0 = C_{\rm env}(1+b_\delta^\nu)$ and $C
  = C_{\rm env}$ in the preceding proposition. No infinite-dimensional Gaussian
  density is used. This proves Lemma~\ref{L:tube-rayleigh-tail}.
\end{proof}

\LeanMarker{L:tube-half-action}
\begin{lemma}[Supremum with the Gaussian quadratic cost]\label{L:tube-half-action}
  For finite $J:H \to \R$ bounded above on every original ball, let $g_J(r) =
  \sup_{\|h\|\le r}J(h)$. In the extended reals,
  \[
    \sup_{r\ge0}\{g_J(r)-r^2/2\}
    =\sup_{h\in H}\{J(h)-\|h\|_H^2/2\}.
  \]
\end{lemma}

\begin{proof}
  For the lower inequality take $r = \|h\|$. For the upper inequality, every
  $\|h\| \le r$ gives $J(h)-r^2/2 \le J(h)-\|h\|^2/2$. Take first the ball
  supremum and then the radius supremum. Each ball is nonempty with finite
  supremum, so no infinity-minus-infinity expression occurs. No compactness,
  rank or density assumption is needed. This proves
  Lemma~\ref{L:tube-half-action}.
\end{proof}

The Gaussian comparison retains one half of the squared norm of the single
control. The coarse envelope is not the sharp spectral/product envelope. Its
application to the original canonical point moment, and the sharp envelope and
variational comparison, require their separate proofs.

\section{Uniform spectral rates and phase products}
\label{S:uniform-physical-products}

Retain every original covariance and amplitude. Set $E = H^{1/2}\times
H^{-1/2}$, $F = E\oplus\R^2$, and $m_R = \eta(\cdot/R)^2$, $R \ge 1$, with the
same fixed cutoff profile. Choose
\[
  \begin{aligned}
    &\frac12<\theta<\min(1-a/2,(3-a-\alpha)/2),\qquad \nu=1/\theta,\\
    &(1-\theta)^{-1}<r<2/a,\qquad \nu\le p<2<r .
  \end{aligned}
\]
Choose auxiliary exponents $0 < \gamma < \min(1/2,(2-\alpha)/2)$ and $\theta <
\beta < \min(1-a/2,(3-a-\alpha)/2)$ for the phase construction. The original
strict threshold makes these intervals nonempty. Write $i_I:H_I \to H_W$ for the
supported closed inclusion. For nonzero covariance the field is
\[
  \mathcal{T}_R=T_R(C_t\otimes I_{H_{\rm sp}}),\qquad
  \|\mathcal{T}_R\|_{H_I\to L^r(I;\mathcal{K}(E))}
  \le C_{\rm mean}=M_0\|C_t\| ,
\]
by Proposition~\ref{P:uniform-mean-regulator}. For zero covariance set $H_{\rm
sp} = H_I = \{0\}$ and $M_0 = C_{\rm mean} = 0$. Means and their rates are then
zero, and mean flows are free; no nonzero-noise estimate is needed for that
case. The residual comparisons below still act on the nonzero phase $E$.

\subsection{The operator class on E and its spectral constants}

For real smooth compact $q$ let $V_q = m_R(\Gamma*q)$, $M_q:H^{1/2} \to
H^{-1/2}$ be its multiplier and $L_q^E(u,v) = (0,M_qu)$. Let $\mathcal{V}_{E,R}$
be their operator-norm closure in $\mathcal{K}(E)$.

\begin{lemma}[Unique potential and fixed-radius lift]\label{L:uniform-E-lift}
  Each $L \in \mathcal{V}_{E,R}$ has the form $L(u,v) = (0,Mu)$, with $\|L\| =
  \|M\| = :k$, and determines a unique real $V \in H^{-1/2}$ supported in
  $\operatorname{supp}m_R$. The symmetric pairing satisfies
  \begin{equation}\label{E:uniform-E-pairing}
    |v_V(u,w)|\le k\|u\|_{H^{1/2}}\|w\|_{H^{1/2}} .
  \end{equation}
  It has a lift $L_F$ in the exact class of Section~\ref{S:signed-static-wave},
  with $L_Fi = iL$. Lift constants may depend on $R$ and are used only for
  identification. Here $v_V(u,w) = \langle Vu,w\rangle_{H^{-1/2},H^{1/2}}$ on
  real states, with its Hermitian extension on complex states.
\end{lemma}

\begin{proof}
  For zero covariance, $L = M = V = L_F = 0$ gives the claims directly. If
  $L_n^E \to L$, the zero columns persist and $M_n \to M$ in
  $\mathcal{L}(H^{1/2},H^{-1/2})$, with the stated norm equality. At fixed $R$
  choose real compact smooth $\chi_R = 1$ near $\operatorname{supp}m_R$. Since
  $M_n\chi_R = V_n$, $\|V_n-V_j\|_{H^{-1/2}} \le
  \|\chi_R\|_{H^{1/2}}\|M_n-M_j\|$. Thus, $V_n \to V$, with the same support by
  testing outside it. Smooth multiplication on $H^{-1/2}$ identifies $Mu = Vu$
  on compact tests, then on $H^{1/2}$ by density. The identity $M\chi_R = V$
  gives uniqueness; symmetry and~\eqref{E:uniform-E-pairing} pass to the limit.

  The original $F$ operators are $(u,v,c_1,c_2)\mapsto(0,M_nu+c_1V_n,0,0)$.
  Their differences are bounded by $(1+\|\chi_R\|_{H^{1/2}})\|M_n-M_j\|$. They
  therefore converge in the old $F$ class and give the lift. This bound is
  allowed to grow with $R$; it will enter no uniform estimate. Dense smooth
  original controls and bounded completion of $T_R$ also show that all completed
  mean fields lie in $\mathcal{V}_{E,R}$ almost everywhere. Bochner averages
  preserve this closed linear class. Complexification preserves $k$ by the
  Hilbert real/imaginary splitting. The real symmetric pairing becomes Hermitian
  and its diagonal splits into the two real diagonals. This proves
  Lemma~\ref{L:uniform-E-lift}.
\end{proof}

\begin{proposition}[Spectral rate in E]\label{P:uniform-E-spectral}
  For $L \in \mathcal{V}_{E,R}$ put
  \[
    \beta(V)=\inf_{\substack{u\in H^1\\\|u\|_2=1}}
    \{\|\nabla u\|_2^2-v_V(u,u)\},\qquad \lambda(L)=\sqrt{-\beta(V)} .
  \]
  If $K \ge \max(1,k)$, the physical wave satisfies, for $t \ge 0$,
  \begin{equation}\label{E:uniform-E-static-bound}
    -K^2\le\beta(V)=\inf\operatorname{spec}H_V\le0,\qquad
    \|P_V^E(t)\|\le8(1+K)^2(1+t)e^{\lambda(L)t}.
  \end{equation}
  If two such operators have norms at most $K$, then
  \begin{equation}\label{E:uniform-E-rate-continuity}
    |\beta(V)-\beta(W)|\le(K+2)\|L_V-L_W\|,\qquad
    |\lambda(L_V)-\lambda(L_W)|
    \le\sqrt{(K+2)\|L_V-L_W\|}.
  \end{equation}
  Every quantitative constant here uses the $E$ norm alone.
\end{proposition}

\begin{proof}
  We derive the form and spectral bounds using only the $E$ multiplier norm.
  We use the fixed-radius lift only to identify the physical flow. Finally,
  approximate minimizers give continuity of the spectral rate without assuming
  an eigenfunction at the bottom.

  \medskip
  \noindent\textit{Step 1.~}We obtain the closed form with constants measured
  on $E$.
  We derive the bounds without transferring a quantitative $F$ bound. With $b =
  \|u\|_2$, $z = \|u\|_{H^1}$, Fourier Cauchy--Schwarz gives $\|u\|_{H^{1/2}}^2
  \le bz$. Equation~\eqref{E:uniform-E-pairing} and the half-weight
  arithmetic--geometric mean inequality give $|v_V(u,u)| \le Kbz \le
  z^2/2+K^2b^2/2$. The weights are $1/2,1/2$, bases $z^2,(Kb)^2 \ge 0$; positive
  powers include zero norms. For $q_V = \|\nabla u\|_2^2-v_V(u,u)$ and $c =
  2+K^2$,
  \begin{equation}\label{E:uniform-E-form-budget}
    \tfrac12z^2+(1+K^2/2)b^2
    \le q_V(u)+cb^2
    \le\tfrac32z^2+(1+3K^2/2)b^2
    \le3(1+K^2)z^2 .
  \end{equation}
  This is a densely defined closed semibounded Hermitian form on full $H^1$: its
  shifted norm is equivalent to the complete $H^1$ norm. The lower bound
  dominates $2b^2$, since $z \ge b$ and $K \ge 1$.
  From~\cite{teschl:14:mathematical}*{Theorem~2.14, p.~81, (2.63)} and
  Theorem~3.6, p.~109, (3.51)--(3.53), it represents a unique self-adjoint $H_V$
  with $A = H_V+c \ge 2$ and $D(A^{1/2}) = H^1$. Closedness uses the definition
  on that book's p.~80. The spectral meaning of $A \ge 2$ follows from bounded
  spectral slices: a nonzero slice below two would give a unit vector in $D(A)$
  contradicting the form lower bound. The bottom argument below uses the same
  representation.

  \medskip
  \noindent\textit{Step 2.~}We compare the spectral phase with the fixed
  Sobolev phase.
  Put $B_0 = 1-\Delta$ and $L_K = 3(1+K^2)$. The common forms obey $\tfrac12B_0
  \le A \le L_KB_0$. For positive self-adjoint $T,S \ge 1$ with ordered forms,
  completing the square gives
  \[
    \langle f,(T+t)^{-1}f\rangle
    =\sup_u\{2\operatorname{Re}\langle f,u\rangle-q_T(u)-t\|u\|^2\}
    \ge\langle f,(S+t)^{-1}f\rangle .
  \]
  The identity $\sqrt x = \pi^{-1}\int_0^\infty t^{-1/2}x/(x+t)\ud t$ follows by
  $t = xv^2$. Positive spectral integration and Tonelli give square-root form
  order, with domain inclusions, from $I-t(T+t)^{-1} \le I-t(S+t)^{-1}$. Apply
  this to $(B_0,2A)$ and $(A,L_KB_0)$, both pairs at least one:
  \[
    D(A^{1/4})=H^{1/2},\qquad
    2^{-1/4}\|u\|_{H^{1/2}}\le\|A^{1/4}u\|_2
    \le L_K^{1/4}\|u\|_{H^{1/2}} .
  \]
  Under the $L^2$ pairing the dual norm is $\|A^{-1/4}v\|_2$, by Cauchy--Schwarz
  and the equality test $u = A^{-1/2}v$. Completion gives $H^{-1/2}$. Thus, with
  $M_K = L_K^{1/4}$, the phase norm $\|(u,v)\|_{E_A}^2 =
  \|A^{1/4}u\|_2^2+\|A^{-1/4}v\|_2^2$ lies between $M_K^{-2}\|(u,v)\|_E^2$ and
  $M_K^2\|(u,v)\|_E^2$. The condition cost is $M_K^2 = \sqrt{3(1+K^2)} \le
  2(1+K)$.

  \medskip
  \noindent\textit{Step 3.~}We identify and bound the spectral bottom.
  For unit $L^2$ data, $b = 1$ and $z \ge 1$. The form bound gives $q_V \ge
  z^2/2-1-K^2/2 \ge -K^2$. A dilated compact unit bump translated wholly off the
  potential support has zero potential pairing and kinetic energy tending to
  zero, so $\beta \le 0$. Complex form infima agree with real ones: dividing
  each complex diagonal by its squared norm is the convex combination of the two
  normalized real diagonals, omitting zero parts. The spectral form formulas
  identify $\beta$ with the spectral bottom by bounded spectral slices. There
  are support points because scalar spectral measures of unit vectors have mass
  one; their arbitrarily small bounded slices give approximate eigenvectors and
  the reverse form bound. If $\Sigma$ denotes the closed PVM support and real
  $z\notin\Sigma$, the Borel function $f_z(h) = (h-z)^{-1}$ on $\Sigma$, zero
  outside, is bounded, as is $hf_z(h)$. Its calculus has range in $D(H_V)$ and
  supplies both inverse identities for $H_V-z$ on their respective domains.
  Thus, no such $z$ is spectral; the approximate eigenvectors give the other
  inclusion. Hence, $\lambda \le K$ and $h+c \ge 2$ on the spectrum.

  \medskip
  \noindent\textit{Step 4.~}We estimate the spectral wave matrix.
  Define $\operatorname{Cos}_t(h),\operatorname{Sin}_t(h)$ by $\cos(t\sqrt
  h),\sin(t\sqrt h)/\sqrt h$ for $h > 0$, by $1,t$ at zero and by
  $\cosh(t\sqrt{-h}),\sinh(t\sqrt{-h})/\sqrt{-h}$ for $h < 0$. In the $E_A$
  isometry the matrix is
  \[
    \begin{pmatrix}
      \operatorname{Cos}_t(h)&\sqrt{h+c}\operatorname{Sin}_t(h)\\
      -h\operatorname{Sin}_t(h)/\sqrt{h+c}&\operatorname{Cos}_t(h)
    \end{pmatrix}.
  \]
  For $h \ge 0$, its off-diagonals have bounds $1+\sqrt c\,t$ and one; at zero
  their values are $\sqrt c\,t,0$. For $\beta \le h < 0$, the entries are
  bounded using $|\operatorname{Cos}_t| \le e^{\lambda t}$,
  $|\operatorname{Sin}_t| \le te^{\lambda t}$, upper off-diagonal at most $\sqrt
  c\,te^{\lambda t}$, and lower at most $Ke^{\lambda t}$. Since $\sqrt c \le
  2K$, each entry is at most $2(1+K)(1+t)e^{\lambda t}$. Set these entries to
  zero off the closed spectral support before applying bounded Borel calculus.
  From~\cite{teschl:14:mathematical}*{Theorems~3.1--3.2, pp.~102--104,
  (3.18)--(3.33)}, their matrix norm is at most twice that entry bound, and the
  condition cost gives \eqref{E:uniform-E-static-bound}. Spectral dominated
  convergence and scalar addition give the group.

  \medskip
  \noindent\textit{Step 5.~}We identify the physical wave without importing
  the lift norm.
  Denote the spectral matrix evolution by $P_V^E(t)$. It is the physical $E$
  wave. The fixed-$R$ lift in Lemma~\ref{L:uniform-E-lift} belongs to the
  precise class in Proposition~\ref{P:signed-physical-wave}. Its form is this
  same $q_V$ and has the same self-adjoint representative by uniqueness,
  regardless of the lift's larger shift. On dense bounded spectral cutoffs the
  two waves have the same $\operatorname{Cos},\operatorname{Sin}$ functions of
  $H_V$. Boundedness extends equality to all $E$. The old graph-core Duhamel
  identification therefore applies; no lift norm or $\chi_R$ constant entered
  the bounds above. Conjugation commutes with the real structures, so real
  restrictions are the phases. For compact $L$, finite-rank approximation and
  strong continuity of the free group and its adjoint make $v\mapsto
  S_E(-v)LS_E(v)$ norm continuous on finite intervals. Its primitive is
  Lipschitz. Its Young equation is therefore the ordinary bounded-interaction
  equation by Proposition~\ref{P:constant-young-linear}. This identifies the
  averaged static baseline without transferring any lift norm into the uniform
  bounds.

  The static baseline is the identity
  \[
    S_E(t)U_H(t,0)=P_V^E(t),\qquad
    H(u)=\int_0^u S_E(-v)LS_E(v)\ud v ,
  \]
  for the constant multiplier $L$ just considered; a cell average uses this
  identity with its frozen multiplier.

  \medskip
  \noindent\textit{Step 6.~}We compare rates by approximate form minimizers.
  Finally, take a unit $\varepsilon$-minimizer for $W$, $0 < \varepsilon \le 1$.
  Since $\beta(W) \le 0$, it has $q_W \le 1$. Thus, $z^2-1-Kz \le 1$, which
  forces $z \le K+2$. The multiplier difference changes its form by at most
  $\|L_V-L_W\|\|u\|_{H^{1/2}}^2 \le (K+2)\|L_V-L_W\|$. Take infima, let
  $\varepsilon\downarrow0$ and reverse $V,W$. The square-root inequality on
  nonnegative $-\beta$ gives \eqref{E:uniform-E-rate-continuity}. For zero
  covariance the class is $\{0\}$, $H_V = -\Delta$ and $\beta = \lambda = 0$;
  the free wave gives the claims directly. This proves
  Proposition~\ref{P:uniform-E-spectral}.
\end{proof}

\subsection{Scalar-time averaging and mean bounds}

Proposition~\ref{P:uniform-E-spectral} bounds the wave of one fixed operator
$L \in \mathcal{V}_{E,R}$ through the rate $\lambda(L)$. The mean field
$\mathcal{T}_R\Phi$ depends on the time variable. We now average it over the
cells of the partition and bound the resulting change in the integrated rate,
uniformly in $R \ge 1$.

\begin{proposition}[Uniform completed-field averaging]\label{P:uniform-rate-averaging}
  Use the true slow partition $I_m$ with mesh at most $\delta B$. Let
  $\mathcal{A}_\delta$ average over each positive-length cell; an empty final
  cell is ignored. Put $\omega_\delta = \|(\mathcal{A}_\delta-I)C_t\|_{H_t\to
  L^r(I)}$. At fixed $B$, $\omega_\delta \to 0$. For $M \ge 0$, if
  $\|\Phi\|_{H_I} \le M$, $W_R = \mathcal{T}_R\Phi$, $\overline{W}_R =
  \mathcal{A}_\delta W_R$ and $Z = C_{\rm mean}M$, then
  \begin{align}
    \|W_R\|_{L^r},\ \|\overline{W}_R\|_{L^r}&\le Z,\qquad
    \|\overline{W}_R-W_R\|_{L^r}\le M_0M\omega_\delta=:\Omega_\delta(M),
    \label{E:uniform-field-averaging}\\
    \operatorname{Rate}_R(\Phi)&=\int_I\lambda(\mathcal{T}_R\Phi(s))\ud s,
    \label{E:uniform-mean-rate}\\
    \int_I|\lambda(\overline{W}_R)-\lambda(W_R)|\ud s
    &\le\sqrt{3+2Z}\sqrt{\Omega_\delta(M)}=:\Lambda_\delta(M)\to0 .
    \label{E:uniform-rate-averaging}
  \end{align}
  These errors are uniform in $R \ge 1$. The rate is finite, nonnegative and
  continuous on $H_I$.
\end{proposition}

\begin{proof}
  We approximate the scalar time map by cell averages uniformly on its compact
  image. The dimension-free tensor bound transfers this to the completed
  operator field. A pointwise rate comparison followed by Cauchy--Schwarz
  gives the integrated rate error.

  \medskip
  \noindent\textit{Step 1.~}We obtain uniform scalar-time averaging on the
  compact image.
  Bochner averaging is a contraction on $L^r(I;\mathscr{Z})$ for any Banach
  $\mathscr{Z}$, by Jensen. For continuous fields on $[0,1]$, its error is
  bounded by the modulus of continuity at the mesh. Simple-field approximation,
  followed by continuous approximation of scalar indicators in $L^r$, proves
  density of those fields. The contraction gives strong convergence for every
  Bochner field. In the scalar case it is uniform on the compact closure of
  $C_t$'s unit-ball image, by a finite net and $\|\mathcal{A}_\delta-I\| \le 2$.
  This proves $\omega_\delta \to 0$.

  \medskip
  \noindent\textit{Step 2.~}We transfer the estimate to the full mean field.
  On finite tensors, integration commutes with the fixed bounded $T_R$:
  \[
    \mathcal{A}_\delta^{\mathcal{K}(E)}
    T_R(C_t\otimes I)\Phi
    =T_R((\mathcal{A}_\delta C_t)\otimes I)\Phi .
  \]
  The dimension-free Hilbert tensor bound in Lemma~\ref{L:uniform-mean-tensor},
  and $\|T_R\| \le M_0$, extend this identity to the completed $H_I$. Applying
  it to the difference gives~\eqref{E:uniform-field-averaging}. Averages remain
  in the same closed $E$ operator class. Only scalar-time compactness was used;
  the tensor map is not claimed compact. Nor does this assert noncompactness of
  a composite spatial family: each fixed-radius completed mean is compact.

  \medskip
  \noindent\textit{Step 3.~}We estimate the integrated rate difference.
  For any two such fields, apply~\eqref{E:uniform-E-rate-continuity} with the
  pointwise bound
  \[
    K(s)=\max(1,\|W_1(s)\|,\|W_2(s)\|).
  \]
  The bound $K(s)+2\le3+\|W_1(s)\|+\|W_2(s)\|$ controls the first factor in
  the pointwise square-root estimate.
  Cauchy--Schwarz gives
  \[
    \int_I|\lambda(W_1)-\lambda(W_2)|
    \le\sqrt{3+\|W_1\|_{L^1}+\|W_2\|_{L^1}}\,
    \sqrt{\|W_1-W_2\|_{L^1}} .
  \]
  Continuity of $\lambda$ on the class makes these quantities measurable. Since
  $|I| = 1$, this proves~\eqref{E:uniform-rate-averaging}. It also proves rate
  continuity, while $\lambda(L) \le \max(1,\|L\|)$ gives
  $\operatorname{Rate}_R(\Phi) \le 1+C_{\rm mean}\|\Phi\|$. For zero covariance
  the fields and rate errors are zero directly. This proves
  Proposition~\ref{P:uniform-rate-averaging}.
\end{proof}

For the remaining proofs use the full parameter choices and the
$B,\delta,N,\ell_m,I_m$ partition of Section~\ref{S:actual-cameron-mean}. For
$W_R = \mathcal{T}_R\Phi$ define
\[
  H_m(u)=\int_0^{\min(u,\ell_m)}
  S_E(-v)W_R(s_m+\delta v)S_E(v)\ud v,\qquad
  P_m^{\rm mean}=S_E(\ell_m)U_{H_m}(\ell_m,0).
\]

\begin{lemma}[Explicit physical primitive and growth bounds]
  \label{L:uniform-product-budgets}
  Write $G_B = \sup_{|v|\le B}\|S_E(v)\|$ and $C_B =
  G_B^2[B^{1-1/r}+B^{1-1/r-\theta}]$. Then $1 \le G_B \le 3B$, $C_B \le
  18B^{3-1/r} \le 32B^{3-1/r}$ and $C_BB^{1/r} \ge B$. The stopped primitives
  lie in the little H\"older space of paths starting at zero and
  \begin{equation}\label{E:uniform-product-primitive}
    \|H_m\|_\theta\le a_m:=C_B\delta^{-1/r}\|W_R\|_{L^r(I_m)},
    \qquad\delta\sum_m a_m^r\le C_B^r\|W_R\|_{L^r}^r .
  \end{equation}
  There is a fixed $C_{\rm gr} = C_{\rm gr}(\theta) > 0$ such that
  \begin{equation}\label{E:uniform-product-growth}
    \log\|S_E(\ell)U_D(\ell,0)\|
    \le C_{\rm gr}B(1+\|D\|_\theta^\nu),\qquad 0\le\ell\le B .
  \end{equation}
  Here $D$ is any $C^\theta([0,B];\mathcal{L}(E))$ path starting at zero; no
  Gaussian or Hamiltonian hypothesis is imposed on this growth bound. The mean
  endpoints have $\|P_m^{\rm mean}\| \ge 1$.
\end{lemma}

\begin{proof}
  We compute the free-group and primitive constants explicitly. We then count
  short contraction intervals for the generic Young growth bound and use
  compactness only for the lower bound on mean endpoint norms.

  \medskip
  \noindent\textit{Step 1.~}We estimate the free group and the physical
  primitive.
  The balanced free matrix has upper entry $\langle k\rangle\sin(vk)/k$ and
  lower $-k\sin(vk)/\langle k\rangle$. Splitting $k \le 1$ and $k \ge 1$ bounds
  the upper by $\sqrt2 B$, while the lower and diagonals have bound one. At zero
  the off-diagonals are $v,0$. The Frobenius bound gives $G_B \le \sqrt{3+2B^2}
  \le 3B$; $S_E(0) = I$ gives $G_B \ge 1$. The stated $C_B$ inequalities follow,
  including $C_BB^{1/r} = G_B^2(B+B^{1-\theta}) \ge B$. H\"older in time with
  $1-1/r > \theta$, the affine time change, compact conjugation and constant
  stopping prove \eqref{E:uniform-product-primitive}; the true partition sums
  the $r$th powers, with zero on an empty cell.

  \medskip
  \noindent\textit{Step 2.~}We make the generic Young growth constant
  independent of radius.
  For the explicit Young bound put $d = \|D\|_\theta$ and $c_\theta =
  2\max(1,C_{\rm Young})$, as in the local estimate proving
  Lemma~\ref{L:young-upper-growth}. If $c_\theta dB^\theta \le 1/4$,
  whole-interval contraction gives $\log\sup\|U_D\| \le \log(4/3) < 1/2$.
  Otherwise at most $2B(4c_\theta d)^\nu$ contraction pieces have norm at most
  $4/3$. Combine their logarithms with $\log G_B \le \log(4B)$. For example
  $C_{\rm gr} = 4+\log4+2(4c_\theta)^\nu$ gives
  \eqref{E:uniform-product-growth}, since $B \ge 1$. This is the generic algebra
  argument on $\mathcal{L}(E)$ with multiplication norm one, not an unspecified
  radius-dependent bound.

  \medskip
  \noindent\textit{Step 3.~}We retain the compact endpoint lower bound.
  Finally, Proposition~\ref{P:uniform-mean-regulator} gives $W_R \in
  L^r(I;\mathcal{K}(E))$ for each fixed $R,\Phi$. The compact conjugation and
  stopped Bochner primitive used above therefore make each $H_m$ a
  compact-valued $C^\theta$ path starting at zero. Lemma~\ref{L:symplectic-free}
  gives $\|P_m^{\rm mean}\| \ge \|S_E(\ell_m)\| \ge 1$, deterministically for
  each $\Phi$. This norm bound needs no fast covariance direction or action
  estimate. An empty endpoint is $I$; zero covariance gives the free endpoint.
  This proves Lemma~\ref{L:uniform-product-budgets}.
\end{proof}

\subsection{Full blocks, the true last cell and polynomial control balls}

\begin{proposition}[Finite-mesh physical mean-product bound]
  \label{P:uniform-mean-product}
  Fix $B$ as above, $M \ge 0$ and $K \ge 1$, and put $Z = C_{\rm mean}M$.
  Uniformly for $R \ge 1$, $\|\Phi\| \le M$,
  \begin{align}
    \delta\sum_{m=0}^N\log\|P_m^{\rm mean}\|
    -\operatorname{Rate}_R(\Phi)
    \le{}&\frac{L_{B,K}}B+
    C_{\rm gr}32^r Z^rB^{3r}(K^{-r}+K^{\nu-r})\notag\\
    &+\Lambda_\delta(M)+C(\theta,B,K)E_\delta+
    \operatorname{Last}_\delta,\label{E:uniform-mean-product}
  \end{align}
  where $L_{B,K} = \log16+4\log(1+K)+3\log(1+B)$ and
  \begin{align*}
    E_\delta&=(2/B)^{1-1/r}C_B\Omega_\delta(M)
    +(2/B)^{1-\nu/r}C_B^\nu\Omega_\delta(M)^\nu,\\
    \operatorname{Last}_\delta&=
    C_{\rm gr}B[\delta+C_B^\nu Z^\nu\delta^{1-\nu/r}].
  \end{align*}
  The last three terms vanish as $\delta\downarrow0$ at fixed $B,K,M$.
\end{proposition}

\begin{proof}
  We freeze each positive-length cell at its own mean. Good full blocks use
  the spectral baseline and the multiplicative perturbation bound; bad full
  blocks use their $L^r$ budget. We estimate the true final cell separately,
  avoiding division by its possibly small length.

  \medskip
  \noindent\textit{Step 1.~}We control the mean-freezing errors and primitive
  budgets.
  On each positive-length cell use its own Bochner average $\overline{M}_m$; set
  it to zero on an empty cell and form the constant-field primitive
  $\overline{H}_m$. Put $e_m = C_B\delta^{-1/r}
  \|W_R-\overline{W}_R\|_{L^r(I_m)}$. Jensen
  and~\eqref{E:uniform-product-primitive} give
  \[
    \|\overline{H}_m\|_\theta\le a_m,\quad
    \|H_m-\overline{H}_m\|_\theta\le e_m,\quad
    \delta\sum e_m^r\le C_B^r\Omega_\delta(M)^r,\quad
    \delta\sum a_m^r\le D_M:=C_B^rZ^r .
  \]
  Finite power mean at exponents $r$ and $r/\nu > 1$, with nonnegative inputs
  $e_m,e_m^\nu$, gives $\delta\sum_m(e_m+e_m^\nu) \le E_\delta$. Here $n = N+1
  \ge 2$, $\delta > 0$ and $\delta n \le 2/B$; all subsequent powers are
  positive and zero errors are allowed.

  \medskip
  \noindent\textit{Step 2.~}We use the spectral baseline on good full blocks.
  On full cells only, $\|\overline{M}_m\| \le (\delta
  B)^{-1/r}\|W_R\|_{L^r(I_m)} = a_m/(C_BB^{1/r}) \le a_m/B$. Split full indices
  $m < N$ into good $a_m \le K$ and bad $a_m > K$. On a good block,
  Proposition~\ref{P:young-upper-multiplicative} applies with base
  $\overline{H}_m$, perturbation $H_m-\overline{H}_m$ and the left multiplier
  $S_E(B)$. The baseline is the static wave of $\overline{M}_m$ by
  Proposition~\ref{P:uniform-E-spectral}. Here $\|\overline{M}_m\| \le K/B \le
  K$, whereas $\|\overline{H}_m\|_\theta \le a_m \le K$. The spectral estimate
  gives
  \[
    \|S_E(B)U_{\overline{H}_m}(B,0)\|
    \le8(1+K)^2(1+B)e^{B\lambda(\overline{M}_m)}.
  \]
  Consequently,
  \[
    \log\|P_m^{\rm mean}\|
    \le L_{B,K}+B\lambda(\overline{M}_m)
    +C(\theta,B,K)(e_m+e_m^\nu).
  \]
  This sufficient prefactor follows from the bound in $E$, $8(1+K)^2(1+B)$, by
  enlarging it. Since $\delta N \le 1/B$, the prefactors sum to at most
  $L_{B,K}/B$. Nonnegativity of $\lambda$ gives $\delta\sum_{\rm
  good}B\lambda(\overline{M}_m) \le \int_I\lambda(\overline{W}_R) \le
  \operatorname{Rate}_R(\Phi)+\Lambda_\delta(M)$. Bad and last rate cells can be
  dropped in this upper bound.

  \medskip
  \noindent\textit{Step 3.~}We control bad full blocks by the higher-power
  budget.
  On bad full cells~\eqref{E:uniform-product-growth} gives $C_{\rm
  gr}B(1+a_m^\nu)$. Because $a_m > K \ge 1$ and $r > \nu$, $1+a_m^\nu \le
  (K^{-r}+K^{\nu-r})a_m^r$. Their weighted sum is bounded by $C_{\rm
  gr}BC_B^rZ^r(K^{-r}+K^{\nu-r})$. The inequality $C_B \le 32B^{3-1/r}$ yields
  the displayed bad-block term.

  \medskip
  \noindent\textit{Step 4.~}We bound the final remainder directly.
  If $\ell_N = 0$, the last endpoint is $I$ and $\log\|I\| = 0$. Otherwise use
  the stopped path, not its frozen full-cell formula. Since $a_N \le
  C_B\delta^{-1/r}Z$, its weighted logarithm is at most
  $\operatorname{Last}_\delta$ by~\eqref{E:uniform-product-growth}. This tends
  to zero because $\nu < r$, even for arbitrarily short last cells; no division
  by that cell length occurs. Combining these estimates
  proves~\eqref{E:uniform-mean-product}. All constants are independent of $R$.
  This proves Proposition~\ref{P:uniform-mean-product}.
\end{proof}

\begin{corollary}[Every polynomial control radius]\label{C:uniform-polynomial-balls}
  For $c_{\rm ball} > 0$, $\kappa_{\rm ball} \ge 0$ put $M_B = c_{\rm
  ball}B^{\kappa_{\rm ball}}$. Then
  \[
    \limsup_{\substack{B\to\infty\\B\in\mathbb{N}}}
    \limsup_{\delta\downarrow0}
    \sup_{\substack{R\ge1\\\|\Phi\|\le M_B}}
    \left\{\delta\sum_m\log\|P_m^{\rm mean}\|
    -\operatorname{Rate}_R(\Phi)\right\}\le0 .
  \]
\end{corollary}

\begin{proof}
  Choose $s > r(3+\kappa_{\rm ball})/(r-\nu)$ and $K_B = B^s$. At each fixed
  $B$, the inner $\delta$ limit has a fixed finite radius and threshold, so the
  last three errors in~\eqref{E:uniform-mean-product} vanish first. Its bad
  terms become
  \[
    C_{\rm gr}32^r C_{\rm mean}^r c_{\rm ball}^r
    \{B^{r(3+\kappa_{\rm ball})-sr}
    +B^{r(3+\kappa_{\rm ball})+s(\nu-r)}\}.
  \]
  Both exponents are negative: the second by choice and the first because $s >
  3+\kappa_{\rm ball}$. Also $L_{B,K_B}/B \le [\log16+7\log2+(4s+3)\log B]/B \to
  0$. Thus, the outer integer-$B$ limit gives the claim. No diagonal choice
  $B_\delta,M_\delta$ is asserted. At fixed parameters all primitive maps are
  bounded linear into separable little H\"older products. Young local dependence
  and positive endpoint norms make the displayed logarithms continuous; the rate
  is continuous by Proposition~\ref{P:uniform-rate-averaging}. Growth bounds
  make each control-ball supremum finite. A countable dense subset of each
  closed Hilbert ball has the same supremum. Supremum over $R$ is an ordinary
  deterministic real supremum, not a random measurability or compact-family
  claim. This proves Corollary~\ref{C:uniform-polynomial-balls}.
\end{proof}

\subsection{Arbitrary small residual paths and the ordered product}

Corollary~\ref{C:uniform-polynomial-balls} compares the sum of the logarithms
$\log\|P_m^{\rm mean}\|$ with the rate, for the mean primitives $H_m$ alone. We
now add compact-operator paths $x_m$ to the $H_m$ and bound the norm of the
ordered product of the perturbed endpoints when these paths are small.

\begin{proposition}[Residual comparison and the integrated rate]
  \label{P:uniform-arbitrary-residual}
  Let $x_m$ be arbitrary little $C^\theta$ compact-operator paths starting at
  zero on $[0,B]$, put $v_m = \|x_m\|_\theta$ and $Y = (\delta\sum_m
  v_m^p)^{1/p}$. Define
  \[
    P_m^{\rm pert}=S_E(\ell_m)U_{H_m+x_m}(\ell_m,0),\qquad
    z_{\rm rat}=\prod_m\frac{\|P_m^{\rm pert}\|}{\|P_m^{\rm mean}\|}.
  \]
  For fixed integer $B \ge 1$, $M \ge 0$ and every $K_{\rm res} \ge 1$,
  \begin{align}
    \delta\log\max(1,z_{\rm rat})\le{}&
    C_{K_{\rm res}}[(2/B)^{1-1/p}Y+(2/B)^{1-\nu/p}Y^\nu]\notag\\
    &+c_BD_M(K_{\rm res}^{-r}+K_{\rm res}^{\nu-r}),
    \label{E:uniform-residual-ratio}
  \end{align}
  uniformly in $R \ge 1$, $\|\Phi\| \le M$. The finite constants may be written
  $C_{K_{\rm res}} = C(\theta,B,K_{\rm res})$ and $c_B = 2^{\nu-1}C_{\rm gr}B$.
  Here $D_M = C_B^r(C_{\rm mean}M)^r$. If $b_\delta \ge 0$ is any deterministic
  function of $\delta$ tending to zero at each fixed $B$, $M_B$ is as in
  Corollary~\ref{C:uniform-polynomial-balls}, and $\mathcal{F} =
  2\delta\log\max(1,\|P_N^{\rm pert}\cdots P_0^{\rm pert}\|)$, then
  \begin{equation}\label{E:uniform-residual-integrated-rate}
    \limsup_{\substack{B\to\infty\\B\in\mathbb{N}}}
    \limsup_{\delta\downarrow0}
    \sup_{\substack{R\ge1,\ \|\Phi\|\le M_B\\Y\le b_\delta}}
    \{\mathcal{F}-2\operatorname{Rate}_R(\Phi)\}\le0 .
  \end{equation}
  No Gaussian, Cameron or Hamiltonian condition on the residuals is imposed.
\end{proposition}

\begin{proof}
  We compare perturbed and mean endpoints separately on good and bad primitive
  blocks. We then use the mean norm lower bound to compare the ordered
  product, taking the residual threshold limit only after the small-parameter
  limit.

  \medskip
  \noindent\textit{Step 1.~}We estimate endpoint norm ratios and sum the
  residual budgets.
  All endpoint norms are positive by the Young and free inverses. On every cell
  with $a_m \le K_{\rm res}$ apply
  Proposition~\ref{P:young-upper-multiplicative} with base $H_m$, perturbation
  $x_m$ and left multiplier $S_E(\ell_m)$. It gives logarithmic ratio at most
  $C(\theta,B,K_{\rm res})(v_m+v_m^\nu)$. On cells with $a_m > K_{\rm res}$ use
  \eqref{E:uniform-product-growth}, $(a_m+v_m)^\nu \le
  2^{\nu-1}(a_m^\nu+v_m^\nu)$ and the mean denominator $\|P_m^{\rm mean}\| \ge
  1$. The logarithmic ratio is at most $c_B(1+a_m^\nu+v_m^\nu)$, with $c_B$
  independent of $K_{\rm res},R,\delta,\Phi,M$. Increase $C_{K_{\rm res}}$ to
  dominate both good and bad residual constants. Summation gives a nonnegative
  bound
  \[
    \log\max(1,z_{\rm rat})
    \le C_{K_{\rm res}}\sum_m(v_m+v_m^\nu)
    +c_B\sum_{a_m>K_{\rm res}}(1+a_m^\nu).
  \]
  The bad sum times $\delta$ is at most $D_M(K_{\rm res}^{-r}+K_{\rm
  res}^{\nu-r})$. Finite power mean at exponents $p \ge 1$ and $p/\nu \ge 1$,
  with nonnegative inputs $v_m,v_m^\nu$, gives $\delta\sum v_m \le
  (2/B)^{1-1/p}Y$ and $\delta\sum v_m^\nu \le (2/B)^{1-\nu/p}Y^\nu$. This
  includes $p = \nu$, zero entries and the true empty block, with no singular
  powers. It proves~\eqref{E:uniform-residual-ratio}.

  \medskip
  \noindent\textit{Step 2.~}We compare the ordered product and take the limits
  in order.
  The product $A_{\rm prod} = \prod_m\|P_m^{\rm mean}\|$ is at least one.
  Since $A_{\rm prod}\ge1$, the inequality $\max(1,A_{\rm prod}z_{\rm rat})\le
  A_{\rm prod}\max(1,z_{\rm rat})$ preserves the product bound when taking the
  maximum with one.
  Submultiplicativity, without commutation, gives
  \[
    \begin{aligned}
      \max(1,\|P_N^{\rm pert}\cdots P_0^{\rm pert}\|)
      &\le A_{\rm prod}\max(1,z_{\rm rat}),\\
      \mathcal{F}-2\delta\sum_m\log\|P_m^{\rm mean}\|
      &\le2\delta\log\max(1,z_{\rm rat}).
    \end{aligned}
  \]
  At fixed $B,M,K_{\rm res}$ the residual terms vanish on $Y \le b_\delta$. Let
  $\delta\downarrow0$ first and then $K_{\rm res} \to \infty$. Since $r > \nu$,
  the residual excess has limsup at most zero, uniformly in $R,\Phi$. Now take
  $M = M_B$ at each fixed $B$ and combine with
  Corollary~\ref{C:uniform-polynomial-balls}; a supremum of sums is no larger
  than the sum of the suprema. The outer integer-$B$ limit
  proves~\eqref{E:uniform-residual-integrated-rate}. The endpoint functional is
  continuous on all tuples by generic Young local dependence, including
  nonphysical residuals. Nothing was multiplied as scalar displacement moments.
  This proves Proposition~\ref{P:uniform-arbitrary-residual}.
\end{proof}

These are uniform deterministic $E$-phase product estimates. Their Gaussian
Rayleigh application, original point transfer and signed-variational value over
nonnegative controls comparison are proved separately in
Sections~\ref{S:actual-gaussian-tubes}--\ref{S:intrinsic-second-moment}.

\section{The original point moment and the growing-cone phase}
\label{S:original-point-transfer}

Retain all original covariances and amplitudes. Write $a = \alpha_0$, $\kappa =
(4-\alpha-a)/(3-\alpha) > 1$. Fix the auxiliary parameters of the preceding
Gaussian and product sections once, in their strict admissible gaps. For $0 <
\delta < 1$ put
\[
  L=\delta^{-1},\quad R=L+2,\quad T_\delta=\delta^{-1/\kappa},
  \quad g=\delta^{(1-a)/2},\quad \Gamma_\delta=g^2\Gamma .
\]
Use the same real even root $\eta$, set $\psi = \eta(\cdot/R)$, $m = \psi^2$, $z
= (\psi,0) \in E$, $D = B_1(0)$ and $v_d = |D| > 0$. Let $W_{\rm fast} = W\circ
T_{\delta,0}$ on the given slow probability space, and let $P_\delta$ be the
physical $E$ endpoint at time $L$ driven by $gM_R(W_{\rm fast})$, with physical
noise mask $m$ once. The datum $\psi$ is not the auxiliary cutoff used to extend
the operator path to constant columns on $F$.

\begin{lemma}[Restriction and fixed-parameter finiteness]
  \label{L:point-actual-finiteness}
  For nonzero covariance, the $E$ flow and its source flows are the closed
  isometric restrictions of their $F$ extensions. Their operator norms and
  differences are measurable, and
  \begin{equation}\label{E:point-fixed-finiteness}
    1\le\E\|P_\delta\|_{\mathcal{L}(E)}^2<\infty .
  \end{equation}
\end{lemma}

\begin{proof}
  We first identify the $E$ flow as a closed isometric restriction of the
  extended phase flow. We then apply the fixed-parameter Gaussian tube
  estimate and compact endpoint lower bound to this exact original endpoint.

  \medskip
  \noindent\textit{Step 1.~}We identify restrictions and measurable endpoint
  norms.
  Let $A_F = gA_{F,R}(W_{\rm fast})$ be the scaled full operator path. Linearity
  and Proposition~\ref{P:constant-phase-driver} give $A_Fi = iD_E$, with $D_E =
  gM_R(W_{\rm fast})$, and $S_F(t)i = iS_E(t)$. The extension uses a cutoff
  $\widetilde{\psi} = 1$ near its larger localization, not the datum $\psi$.
  Applying $i$ to the state Young sums and using uniqueness gives $U_F(t,s)i =
  iU_E(t,s)$, then $P_F(t,s)i = iP_E(t,s)$. The same argument applies to
  sources. Closedness and isometry imply $\|P_E\| \le \|P_F\|$ and the analogous
  difference bound; density of $E$ in $F$ is neither true nor used. The
  continuous Young map on its separable path span, followed by bounded physical
  conjugation and restriction, gives measurable endpoint norms and uniform
  time-triangle difference norms. No operator-norm continuity of the free group
  is assumed. Young Neumann pieces and the free inverse also give invertibility.

  \medskip
  \noindent\textit{Step 2.~}We prove finiteness and the lower endpoint norm
  bound.
  The clock condition is $0<\delta<\min(1,1/B)$ from
  Section~\ref{S:actual-cameron-mean}, which reduces to $0<\delta<1$ at $B=1$.
  We check a direct finiteness input, without a fixed-mask asymptotic rate. Take
  $B = 1$ in Propositions~\ref{P:tube-actual-law}
  and~\ref{P:tube-coarse-envelope} and Lemma~\ref{L:tube-rayleigh-tail}. Every
  $0 < \delta < 1$ satisfies their clock condition. The global endpoint
  in~\eqref{E:tube-global-endpoint} is exactly this $P_\delta$: the same given
  $W_{\rm fast}$, amplitude $g$, root $\eta$, one mask $m$, massless $E$ phase,
  and horizon $\delta^{-1}$. We use the same admissible $\theta$ and operator
  path version. Explicit row identities, rational times and Young uniqueness
  identify versions if a smaller common admissible exponent is used. Thus, the
  finite radial bound applies to this endpoint, not merely an equal-law copy. At
  these fixed parameters the operator path $gM_R(W_{\rm fast})$ is a
  compact-valued $C^\theta([0,L];\mathcal{K}(E))$ path starting at zero on its
  path event. Lemma~\ref{L:symplectic-free} gives $\|P_\delta\| \ge \|S_E(L)\|
  \ge 1$ almost surely, proving~\eqref{E:point-fixed-finiteness}. This use
  requires no $F$-norm rate along growing radii. Source convergence and source
  uniform integrability are checked next, separately from this original endpoint
  finiteness. This proves Lemma~\ref{L:point-actual-finiteness}.
\end{proof}

\begin{lemma}[The same regularized fast noises preserving the given noise]\label{L:point-fast-sources}
  At each fixed $\delta$, take one regularized slow noise preserving the given
  noise $Y_{\epsilon,\tau}$ and define
  \[
    Y_{\rm fast}(t,x)=\delta^{a/2}Y_{\epsilon,\tau}(\delta t,x).
  \]
  Its physical $E$ flows $P_E^{\epsilon,\tau}$ with coefficient $m gY_{\rm
  fast}$ satisfy
  \begin{equation}\label{E:point-source-flow-limit}
    \E\sup_{0\le s\le t\le L}
    \|P_E^{\epsilon,\tau}(t,s)-P_E(t,s)\|^2\longrightarrow0
    \quad\hbox{as }(\epsilon,\tau)\to(0,0).
  \end{equation}
  The unmasked field $gY_{\rm fast}$ has the canonical source law for covariance
  $\Gamma_\delta$, spatial width $\epsilon$ and temporal width $\tau/\delta$.
\end{lemma}

\begin{proof}
  We pull one regularized slow noise into fast coordinates and transfer its
  operator paths on the original coupling. We identify the ordinary response
  and compute the exact temporal width and covariance amplitude.

  \medskip
  \noindent\textit{Step 1.~}We prove convergence of the source flow at fixed
  dilation.
  Testwise this field is $Y_{\epsilon,\tau}(T_{\delta,0}f)$. The original
  isometry and source domination give source and residual variances bounded by
  $\|f\|_{H_{\rm fast}}^2$, with residual variance tending to zero at fixed
  $\delta$ on every smooth compact test. The Gaussian transfer
  Proposition~\ref{P:source-gaussian-transfer} therefore gives operator-path
  convergence in every finite $L^pC^\theta$, retaining $W_{\rm fast}$ literally.
  Multiplication by fixed $g$ and the bounded $F$ extension preserve it. The
  transfer bound and Proposition~\ref{P:actual-operator-driver} also give
  uniform source $\sqrt p$ path moments at these fixed $\delta,R,L$. The
  moment-series argument yields square-exponential path norms. Young local
  dependence gives convergence in probability on balls of operator paths; the
  growth exponent $1/\theta < 2$ and those tails give uniform flow moments of
  every finite order and uniform integrability. The proof of
  Proposition~\ref{P:constant-gaussian-flow} applies to this very pullback
  family. Restriction from $F$ then proves~\eqref{E:point-source-flow-limit}. No
  uniform-in-$\delta$ source bound or exchange of limits is asserted.

  \medskip
  \noindent\textit{Step 2.~}We identify the source primitive on the same
  states.
  This source is also the ordinary response for all $E$ data. On Schwartz states
  its first primitive has row $gY_{\rm fast}(h_R)$. The source indicator
  approximation and bounded point variances at positive widths identify it with
  the statewise integral $\int_0^t S_E(-r)B(mgY_{\rm fast}(r))S_E(r)\ud r$.
  Fubini uses $\E|ZY| \le \|Z\|_2\|Y\|_2$ when Gaussian testing is needed.
  Countable dense states, rational times and continuity give the full
  ordinary-primitive identity, as in
  Proposition~\ref{P:source-ordinary-response}.
  Lemma~\ref{L:constant-source-solution} and
  Proposition~\ref{P:symplectic-source} then identify the Young equation with
  the ordinary interaction equation on $E$.

  \medskip
  \noindent\textit{Step 3.~}We compute the transformed source law.
  In the temporal kernel $K_\tau(v) = \Gamma_{\rm E}(a)^{-1}
  \int_0^{\tau^{-1}}e^{-\lambda|v|}\lambda^{a-1}\ud\lambda$, substitution
  $\lambda' = \delta\lambda$ gives
  \[
    \delta^aK_\tau(\delta v)=K_{\tau/\delta}(v).
  \]
  The change of variable contributes $\delta^{-1}$ from $d\lambda$ and
  $\delta^{1-a}$ from $\lambda^{a-1}$. These cancel the prefactor $\delta^a$,
  while the upper cutoff becomes $\delta/\tau$. Thus, the new temporal width
  is $\tau/\delta$, with the original temporal covariance amplitude unchanged
  before multiplication by $g^2$.
  Hence, $gY_{\rm fast}$ has covariance
  $K_{\tau/\delta}(t-s)(p_{2\epsilon}*\Gamma_\delta)(x-y)$. The shared-width
  maximum and spatial heat-sum conventions persist for two sources. Positive
  multiplication by $g^2$ preserves every original covariance class, including
  spatial white noise. At fixed $\delta$, $(\epsilon,\tau/\delta) \to (0,0)$.
  The full source point-moment limit in
  Proposition~\ref{P:masked-canonical-transfer} thus applies to this
  single-source law and its continuous version; no newly reconstructed auxiliary
  family for $\Gamma_\delta$ is required. This proves
  Lemma~\ref{L:point-fast-sources}.
\end{proof}

\begin{lemma}[Compact-data source identification]\label{L:point-compact-source}
  For fixed $\delta,\epsilon,\tau$, the displacement $v_{\epsilon,\tau}(t,x) =
  (\operatorname{pr}_1P_E^{\epsilon,\tau}(t,0)z)(x)$ has a jointly measurable
  continuous version and satisfies the ordinary scalar mild equation with
  coefficient $b = mgY_{\rm fast}$ and initial data $(\psi,0)$.
\end{lemma}

\begin{proof}
  At positive source widths we solve in a sufficiently regular Sobolev phase,
  identify that solution with the $E$ flow by uniqueness, and then use Fourier
  embedding and Green testing to recover the continuous scalar mild solution.

  \medskip
  \noindent\textit{Step 1.~}We construct the high-regularity phase solution.
  On the fixed source regularity event, $b$ is smooth in space, has common
  compact support in $B_{2R}$, and is continuous in time in every bounded
  spatial derivative norm on $[0,L]$. Choose integer $k > d/2+1$ and $E_k =
  H^k\times H^{k-1}$. Its Bessel-unitary free matrix is
  \[
    \begin{pmatrix}
      \cos(t|\xi|)&\langle\xi\rangle\sin(t|\xi|)/|\xi|\\
      -|\xi|\sin(t|\xi|)/\langle\xi\rangle&\cos(t|\xi|)
    \end{pmatrix}.
  \]
  The zero sine quotient is $t$; finite-horizon bounds and dominated Fourier
  convergence give a strongly continuous group on $E_k$. Leibniz and the
  derivative-norm equivalence, obtained by expanding $(1+|\xi|^2)^{k-1}$ and
  Parseval, give
  \[
    \|b(t)f\|_{H^{k-1}}
    \le C_k\max_{|j|\le k-1}\|\partial^jb(t)\|_\infty
    \|f\|_{H^{k-1}} .
  \]
  Thus, $B(b(t))$ is a norm-continuous bounded force on $E_k$, with a finite
  uniform sample bound $B_L$. For $Z_0 = S(t)z$ and $Z_{n+1}(t) =
  \int_0^tS(t-r)B(b(r))Z_n(r)\ud r$,
  \[
    \sup_{t\le L}\|Z_n(t)\|_{E_k}
    \le C_L^{n+1}B_L^nL^n\|z\|_{E_k}/n! .
  \]
  This series converges in $C_tE_k$. Iterating a homogeneous Volterra difference
  proves uniqueness. Its continuous embedding into $E$ and ordinary $E$
  uniqueness identify the sum with $P_E^{\epsilon,\tau}(t,0)z$.

  \medskip
  \noindent\textit{Step 2.~}We obtain continuous measurable displacement and
  velocity.
  Fourier Cauchy--Schwarz gives $\|\widehat{f}\|_1 \le \|f\|_{H^k}
  (\int\langle\xi\rangle^{-2k}\ud\xi)^{1/2} < \infty$. Fourier inversion and
  dominated convergence imply continuous bounded displacement, and the same
  holds for velocity since $k-1 > d/2$. The source derivatives on the compact
  support form a measurable $C_tC_x^{k-1}$ random element, by countable
  evaluations and continuity. Picard iteration and local norm estimates make its
  phase sum measurable in $C_tE_k$; the Fourier embedding gives joint
  measurability. These high norms are used only at positive source widths.

  \medskip
  \noindent\textit{Step 3.~}We identify the scalar Green equation.
  The first component of Duhamel is, distributionally,
  \[
    v(t)=\cos(t\sqrt{-\Delta})\psi+
    \int_0^t G_d(t-r)*[b(r)v(r)]\ud r .
  \]
  Here $G_d(h)$ is the positive Green measure of mass $h$ and support in $B_h$,
  including the sphere measure for $d = 3$. The product $bv$ is continuous
  compactly supported, with finite uniform $L^1$ and supremum bounds. Against a
  compact smooth $\phi$, absolute Fubini is bounded by
  $\int_0^t\|\phi\|_\infty\|b(r)v(r)\|_1(t-r)\ud r < \infty$. Writing $G_d(h)$
  as $h$ times the dilation of its unit-time probability measure proves
  continuity of the convolution term. The free term is continuous by the Fourier
  embedding. Thus, the distributional identity agrees everywhere with the
  pointwise continuous mild versions. No Green density in $d = 3$ or original
  Sobolev point trace is used. This proves Lemma~\ref{L:point-compact-source}.
\end{proof}

The displacement $v_{\epsilon,\tau}$ of Lemma~\ref{L:point-compact-source} has
the masked coefficient $mgY_{\rm fast}$ and the data $(\psi,0)$. The unmasked
source solution has the coefficient $gY_{\rm fast}$ and the data $(1,0)$. We now
show that the two agree on a cone that contains $D$ at time $L$.

\begin{proposition}[Cone identity and an average of squares]
  \label{P:point-source-cone}
  Let $u_{\epsilon,\tau}^{\rm full}$ be the unmasked source solution with
  coefficient $gY_{\rm fast}$ and data $(1,0)$. At the fixed parameters, almost
  surely,
  \[
    v_{\epsilon,\tau}=u_{\epsilon,\tau}^{\rm full}
    \quad\hbox{on }\mathcal{C}=
    \{(t,x):0\le t\le L,\ |x|\le1+L-t\}.
  \]
  Consequently,
  \begin{equation}\label{E:point-source-square-average}
    \E\|\operatorname{pr}_1P_E^{\epsilon,\tau}(L,0)z\|_{L^2(D)}^2
    =v_d\,\E|u_{\epsilon,\tau}^{\rm full}(L,0)|^2<\infty .
  \end{equation}
\end{proposition}

\begin{proof}
  We compare the compact-data and constant-data source equations on one
  causally closed cone. Their free terms agree there by total finite
  propagation, and iteration of the homogeneous difference equation gives
  equality before taking the terminal spatial square average.

  \medskip
  \noindent\textit{Step 1.~}We identify the common source law and the causally
  closed cone.
  Proposition~\ref{P:masked-source-moments}, applied through the exact source
  law of Lemma~\ref{L:point-fast-sources}, gives the locally continuous
  measurable full solution and its $L^2$ point Picard sum. The pathwise cone
  argument uses the realized pullback field; the point second moments depend
  only on its Gaussian covariance through those Bochner chain sums. We place no
  global unmasked field in $E$ or $F$. The cone is causally closed: a
  predecessor at distance at most $t-s$ has norm at most $1+L-s$. Its spatial
  projection lies a positive distance inside the plateau of both $\psi$ and $m$.

  \medskip
  \noindent\textit{Step 2.~}We compare both initial phase components on the
  cone.
  The total initial difference $(\psi-1,0)$ is represented in $E$ plus constants
  by $(\psi,0)+(-1,0)$. Every entry of the full free phase kernel is a compactly
  supported distribution in $B_t$, and the top-left is $\partial_tG_d(t)$. Since
  $|x|+t \le 1+L = R-1$ on the cone, its entire translated kernel support lies
  in the region where this total datum vanishes.
  Lemma~\ref{L:constant-phase-free} therefore gives
  $\cos(t\sqrt{-\Delta})\psi(x) = 1$ there. Testing the $E$ and constant
  summands separately would not give this total finite-propagation statement.

  \medskip
  \noindent\textit{Step 3.~}We prove uniqueness of the two restrictions on the
  cone.
  The two mild equations have the same coefficient on the cone. For their
  difference $w$, causal closure gives the homogeneous Green integral equation
  there. Put $H = \sup_{\mathcal{C}}|w|$ and $M = \sup_{\mathcal{C}}|gY_{\rm
  fast}|$, both finite sample constants. Iteration with the positive Green
  measures gives $|w(t,x)| \le HM^nt^{2n}/(2n)!$. Indeed, the mass recursion is
  $f_0 = 1$, $f_{n+1}(t) = \int_0^t(t-r)f_n(r)\ud r = t^{2n+2}/(2n+2)!$. The
  factorial ratio forces this bound to zero on the fixed cone. No positivity of
  the sample coefficient or expectation of a global source supremum is required.

  \medskip
  \noindent\textit{Step 4.~}We take the terminal average of squares.
  At time $L$ the cone contains $D$. Tonelli and spatial stationarity of the
  unmasked single-source law give~\eqref{E:point-source-square-average}.
  Finiteness is the full source point-moment bound. This is an integral of
  squares, not a squared average or a covariance comparison with a masked model.
  Events are for each fixed $\delta,\epsilon,\tau$. This proves
  Proposition~\ref{P:point-source-cone}.
\end{proof}

Agreement on the terminal observation set alone would not suffice:
both initial phase components must agree on its backward cone.
No-go~\ref{P:nogo-31} gives a smooth free-wave counterexample when
the velocity guard is omitted, and a separate phase example explaining
why operator growth alone does not imply growth for every initial state.

Proposition~\ref{P:point-source-cone} holds at fixed positive
$\epsilon,\tau$ and for the scaled covariance $\Gamma_\delta$. We now remove the
regularization and undo the scaling to reach the original point moment
$\E|u_\Gamma(T_\delta,0)|^2$.

\begin{proposition}[Exact original point transfer and polynomial cost]
  \label{P:original-point-transfer}
  For every original model, including zero covariance,
  \begin{align}
    \E|u_\Gamma(T_\delta,0)|^2
    &=v_d^{-1}\E\|\operatorname{pr}_1P_\delta z\|_{L^2(D)}^2,
    \label{E:original-point-transfer}\\
    &\le v_d^{-1}\|\psi\|_{H^{1/2}}^2\E\|P_\delta\|^2,\notag\\
    \limsup_{T\to\infty}T^{-\kappa}\log\E|u_\Gamma(T,x)|^2
    &\le\limsup_{\delta\downarrow0}\delta\log\E\|P_\delta\|^2
    \quad\hbox{for every fixed }x .
    \label{E:original-point-rate-reduction}
  \end{align}
\end{proposition}

\begin{proof}
  We remove source regularization at fixed dilation, then use the canonical
  law scaling to identify the original point moment. The compact initial datum
  costs only a polynomial factor, which disappears after logarithmic
  normalization.

  \medskip
  \noindent\textit{Step 1.~}We pass the source square identity to the original
  operator endpoint.
  Assume nonzero covariance until the last paragraph, and fix $\delta$.
  Equation~\eqref{E:point-source-flow-limit} gives convergence of the source
  displacement in $L^2(\Omega;L^2(D))$, since $\|\operatorname{pr}_1e\|_{L^2(D)}
  \le \|e\|_E$. Expected squared norms converge by Cauchy--Schwarz applied to
  their difference times their sum. Separately, the full source point-moment
  limit \eqref{E:masked-unmasked-moment-limit} at covariance $\Gamma_\delta$
  gives convergence of the right side of~\eqref{E:point-source-square-average}
  to $v_d\E|u_{\Gamma_\delta}(L,0)|^2$. Both cutoffs $(\epsilon,\tau/\delta)$
  tend to zero at fixed $\delta$. Thus, the two limits in the equal prelimit
  moments give
  \[
    \E|u_{\Gamma_\delta}(L,0)|^2
    =v_d^{-1}\E\|\operatorname{pr}_1P_\delta z\|_{L^2(D)}^2 .
  \]

  \medskip
  \noindent\textit{Step 2.~}We identify the original moment by the canonical
  time change.
  The law normalization is the consequence of
  Lemma~\ref{L:exact-canonical-scaling} and the canonical-series time change
  proved in Proposition~\ref{P:raw-block-point-comparison}. Explicitly, for
  $\delta = T^{-\kappa}$ and $r_T = T\delta$, $V_\delta$ has the law of
  $u_\Gamma(Tt,r_Tx)$ and equation \eqref{E:exact-fast-equation}. At positive
  source widths, $U(t,x) = V_\delta(\delta t,x)$ has coefficient
  $\sqrt\delta\,Y_{\epsilon,\tau}(\delta t,x) = gY_{\rm fast}(t,x)$. The
  fixed-order and absolutely convergent canonical-series argument passes this
  time change to the point moment at fixed $\delta$. Hence,
  $\E|u_\Gamma(T_\delta,0)|^2 = \E|u_{\Gamma_\delta}(L,0)|^2$. This is a
  law/moment identity, not a pathwise identity on the same noise. It preserves
  all mixed, unequal, self and empty Stratonovich chains and the
  Green-before-noise convention in $d = 3$. The shrinking $r_T$ is not the mask
  radius $R$; no spatial Jacobian occurs at point zero. The phase/operator
  inequality now proves~\eqref{E:original-point-transfer}.

  \medskip
  \noindent\textit{Step 3.~}We control the initial-data norm and its
  logarithmic cost.
  Parseval and the inhomogeneous Bessel weights give, for $R \ge 1$,
  \[
    \|\psi\|_{H^{1/2}}^2\le\|\psi\|_{H^1}^2
    =R^d\|\eta\|_2^2+R^{d-2}\|\nabla\eta\|_2^2
    \le C_\eta R^d,\quad C_\eta=\|\eta\|_2^2+\|\nabla\eta\|_2^2>0.
  \]
  This includes the negative $d-2$ power in $d = 1$. For
  $\operatorname{DataPref}_\delta = v_d^{-1}C_\eta(\delta^{-1}+2)^d$,
  \[
    \delta\log\operatorname{DataPref}_\delta
    =\delta\log(C_\eta/v_d)+d\delta[\log(1+2\delta)-\log\delta]\to0 .
  \]
  The nontrivial term follows by $t = -\log\delta$ and $te^{-t} \to 0$, since
  $e^t \ge t^2/2$. The original point moment is at least one: its nonnegative
  full pairing expansion contains the unit empty pair of chains. It is finite
  by~\eqref{E:original-point-transfer} and~\eqref{E:point-fixed-finiteness}.
  Thus, logs are defined without subtracting infinite quantities. Since $\alpha
  < 2 < 3$, $a < 1$, we have $\kappa > 1$ and $\delta\mapsto T_\delta$ is a
  decreasing bijection $(0,1) \to (1,\infty)$ with inverse $\delta =
  T^{-\kappa}$. Taking logs and this change of parameter proves
  \eqref{E:original-point-rate-reduction} at zero. Original canonical spatial
  stationarity, proved in Proposition~\ref{P:raw-block-original-lower}, extends
  it to every fixed $x$; no stationarity of the masked operator is assumed.

  \medskip
  \noindent\textit{Step 4.~}We verify the free zero-covariance case directly.
  For zero covariance, $u = 1$, all operator paths vanish and $P_\delta =
  S_E(L)$. The same total free cone argument gives the spatial-square identity.
  Lemma~\ref{L:symplectic-free} gives $\|S_E(L)\| \ge 1$, and its Fourier matrix
  gives $\|S_E(L)\| \le C(1+L)$. Thus, $0 \le \delta\log\|S_E(\delta^{-1})\|^2
  \to 0$ and both claims hold directly without a nonzero-noise import. This
  proves Proposition~\ref{P:original-point-transfer}.
\end{proof}

\section{An upper bound from a signed variational expression}
\label{S:original-signed-score-upper}

Keep all auxiliary choices fixed as in the Gaussian tube construction, including
$\alpha/4 < s_{\rm sp} < 1/2$, the admissible $\theta,\nu,r,p$, and the three
fixed even-$Q$ conditions. The mean bound $C_{\rm mean}$ is independent of
$R,B,\delta$; it may depend on those fixed auxiliary parameters and the model.
The space $H_I$ has the original covariance $\Gamma$. Its norm is not reweighted
by $g$ and is not a sum of fast block costs. Let $R_\delta = \delta^{-1}+2$ and
$M_\delta = \E\|P_\delta\|_{\mathcal{L}(E)}^2$. For every fixed integer $B$ the
Gaussian block tuple has typical radius $b_{\delta,B} =
2\E\|X_{\delta,B}\|_{\delta,p} \to 0$. Its endpoint is the same $P_\delta$ for
every $B$ by Lemma~\ref{L:tube-actual-endpoint}, not a $B$-dependent surrogate.
Write $\mathcal{M}_{\delta,B,R}$ for the normalized mean tuple
in~\eqref{E:tube-actual-mean}, with the indicated block length and radius. Write
$\mathcal{F}_{\delta,B}$ for the everywhere-defined endpoint functional
$F_\delta$ of Lemma~\ref{L:tube-actual-endpoint}.

\begin{lemma}[Signed variational expression and quantitative radial envelope]
  \label{L:signed-score-envelope}
  For $R \ge 1$ define the variational expression
  \begin{equation}\label{E:signed-raw-score}
    \mathcal{S}_R=\sup_{\Phi\in H_I}
    \{2\operatorname{Rate}_R(\Phi)-\|\Phi\|_{H_I}^2/2\},
    \qquad V_{\rm up}=\limsup_{\delta\downarrow0}\mathcal{S}_{R_\delta}.
  \end{equation}
  Then $0 \le \mathcal{S}_R,V_{\rm up} \le 2+2C_{\rm mean}^2 < \infty$. The tube
  envelope $G_{\delta,B}$ of Proposition~\ref{P:tube-coarse-envelope} satisfies,
  uniformly in $R \ge 1$ (in particular at the radius $R_\delta$),
  \begin{equation}\label{E:signed-polynomial-envelope}
    G_{\delta,B}(\rho)\le A_0+A_1B^{3\nu}\rho^\nu
    +\eta_B b_{\delta,B}^\nu ,
  \end{equation}
  for $0 < \delta < \min(1,1/B)$ and every $\rho \ge 0$, where
  \[
    \begin{aligned}
      A_0&=4C_{\rm gr},\\
      A_1&=2^{\nu+1-\nu/r}C_{\rm gr}32^\nu C_{\rm mean}^\nu,\\
      \eta_B&=2^{\nu+1-\nu/p}C_{\rm gr}B^{\nu/p}.
    \end{aligned}
  \]
  The variational expression is not identified here with a constant in the cited
  source or a variational value over nonnegative controls.
\end{lemma}

\begin{proof}
  We bound the signed score by completing the square in the control norm. We
  then use the explicit physical Young growth bound and weighted block
  estimates to compute the radial envelope constants, keeping their dependence
  on $B$.

  \medskip
  \noindent\textit{Step 1.~}We obtain a finite uniform ceiling for the signed
  score.
  The mean is linear, so $\mathcal{T}_R0 = 0$. Its free form has nonnegative
  gradient energy and spectral bottom zero; thus $\operatorname{Rate}_R(0) = 0$.
  Proposition~\ref{P:uniform-rate-averaging} gives $\operatorname{Rate}_R(\Phi)
  \le 1+C_{\rm mean}\|\Phi\|$. With $q = \|\Phi\| \ge 0$,
  \[
    2\operatorname{Rate}_R(\Phi)-q^2/2
    \le2+2C_{\rm mean}q-q^2/2
    =2+2C_{\rm mean}^2-(q-2C_{\rm mean})^2/2 .
  \]
  The coarser inequality also follows from half-weight arithmetic--geometric
  mean at bases $q^2,(2C_{\rm mean})^2$, with nonnegative weights and bases,
  including zeros. The zero control gives the lower bound. No optimizer, finite
  rank, positive control or convergence of $\mathcal{S}_{R_\delta}$ is assumed.
  The finite ceiling is auxiliary-parameter dependent and need not be sharp.

  \medskip
  \noindent\textit{Step 2.~}We compute the block-length dependence of the tube
  envelope.
  For $\|a\|_{\delta,p} \le b_{\delta,B}$, $\|\Phi\| \le \rho$ set $v_m =
  \|a_m\|_\theta$ and $h_m = \|(\mathcal{M}_{\delta,B,R}\Phi)_m\|_\theta$. The
  quantitative physical bound~\eqref{E:uniform-product-growth}, already
  including the free factor, and convexity give
  \[
    \mathcal{F}_{\delta,B}(a+\mathcal{M}_{\delta,B,R}\Phi)
    \le2C_{\rm gr}B\delta n
    +2^\nu C_{\rm gr}B\delta\sum_m(v_m^\nu+h_m^\nu).
  \]
  Its upper expression is nonnegative, so it controls the maximum at one. The
  true primitive estimate and weighted power mean give
  \begin{align*}
    \delta\sum v_m^\nu&\le(2/B)^{1-\nu/p}b_{\delta,B}^\nu,\\
    \delta\sum h_m^\nu&\le(2/B)^{1-\nu/r}
    (32B^{3-1/r}C_{\rm mean}\rho)^\nu .
  \end{align*}
  The exponents are $p/\nu \ge 1$, $r/\nu > 1$, inputs nonnegative, $n = N+1 \ge
  2$, $\delta > 0$, and $\delta n \le 2/B$. Zeros and $p = \nu$ are legal. The
  first constant is $4C_{\rm gr}$. In the mean term the factors $B^{\nu/r}$ and
  $B^{-\nu/r}$ cancel exactly, leaving $B^{3\nu}$. These computations give the
  displayed coefficients and~\eqref{E:signed-polynomial-envelope}, uniformly for
  all radii. They use no unspecified $c_Y(B)$ or quantitative $F$ estimate. This
  proves Lemma~\ref{L:signed-score-envelope}.
\end{proof}

Lemma~\ref{L:signed-score-envelope} bounds the envelope $G_{\delta,B}$ at every
radius $\rho$. We now split the radial integral
in~\eqref{E:tube-coarse-envelope} at a radius $M$. Below $M$ we compare the
exponent with $\mathcal{S}_{R_\delta}$, and above $M$ we use
Lemma~\ref{L:tube-rayleigh-tail}.

\begin{proposition}[Finite-delta bulk and exact tail]
  \label{P:signed-moment-bulk-tail}
  For fixed integer $B \ge 1$, $0 < \delta < \min(1,1/B)$ and $M > 0$, define
  the deterministic number
  \[
    \mathcal{E}_{\delta,B,M}=
    \max\left(0,\sup_{\substack{R\ge1,\ \|\Phi\|\le M\\
    \|x\|_{\delta,p}\le b_{\delta,B}}}
    [\mathcal{F}_{\delta,B}(x+\mathcal{M}_{\delta,B,R}\Phi)
    -2\operatorname{Rate}_R(\Phi)]\right).
  \]
  It is finite. For every polynomial $M_B = c_{\rm ball}B^{\kappa_{\rm ball}}$,
  with $c_{\rm ball} > 0$, $\kappa_{\rm ball} \ge 0$, one has
  \begin{equation}\label{E:signed-error-limit}
    e_B:=\limsup_{\delta\downarrow0}\mathcal{E}_{\delta,B,M_B}
    \ \text{finite and nonnegative},\qquad e_B\longrightarrow0 .
  \end{equation}
  If $4A_1B^{3\nu} \le M^{2-\nu}$, then
  \begin{equation}\label{E:signed-finite-moment}
    M_\delta\le
    \frac{M^2}{2\delta}
    e^{(\mathcal{S}_{R_\delta}+\mathcal{E}_{\delta,B,M})/\delta}
    +4e^{(A_0+\eta_Bb_{\delta,B}^\nu-M^2/4)/\delta}.
  \end{equation}
\end{proposition}

\begin{proof}
  We first identify the vanishing deterministic error on polynomial control
  balls. We then split the exact radial integral at $M$: the bulk uses the
  signed variational score, and the tail uses the explicit subquadratic
  envelope.

  \medskip
  \noindent\textit{Step 1.~}We control the deterministic product error.
  The polynomial envelope and nonnegative rate make this supremum finite. It is
  an ordinary deterministic real number, not a random uncountable-radius
  supremum. At fixed $R$ the inner functions are continuous in controls and
  paths; no measurable maximizing selector or joint $\delta,R$ measurability is
  needed. The zero control and zero path have rate zero and
  $\mathcal{F}_{\delta,B}(0) \ge 0$, so the untruncated supremum is already
  nonnegative. The anchor $b_{\delta,B} \to 0$ at fixed $B$ satisfies
  Proposition~\ref{P:uniform-arbitrary-residual} exactly. Its polynomial-radius
  conclusion gives $\limsup_B e_B \le 0$. At fixed $B$, the polynomial envelope
  bounds the inner limsup by a finite number, and $e_B \ge 0$ then
  proves~\eqref{E:signed-error-limit}.

  \medskip
  \noindent\textit{Step 2.~}We split the radial integral and keep both
  prefactors.
  The exact Gaussian-enlargement bound is
  \[
    M_\delta\le\int_0^\infty\frac{\rho}{\delta}
    e^{(G_{\delta,B}(\rho)-\rho^2/2)/\delta}\ud\rho .
  \]
  For every $\rho \le M$, every base-ball $a$ and $\|\Phi\| \le \rho$,
  \[
    \mathcal{F}_{\delta,B}(a+\mathcal{M}_{\delta,B,R_\delta}\Phi)-\rho^2/2
    \le2\operatorname{Rate}_{R_\delta}(\Phi)-\|\Phi\|^2/2
    +\mathcal{E}_{\delta,B,M}
    \le\mathcal{S}_{R_\delta}+\mathcal{E}_{\delta,B,M}.
  \]

  Taking suprema and integrating $\rho/\delta$ over $[0,M]$ gives the first term
  of~\eqref{E:signed-finite-moment}, even if $M^2/(2\delta) < 1$. This used
  $\rho^2 \ge \|\Phi\|^2$, not equality or an optimizer. For the tail apply
  Lemma~\ref{L:tube-rayleigh-tail} with $C = A_1B^{3\nu}$ and intercept
  $A_0+\eta_Bb_{\delta,B}^\nu$. All cutoff conditions hold: $M > 0$, $\rho \ge
  M$, $C \ge 0$, $1 < \nu < 2$, $4C \le M^{2-\nu}$. Its analytic integration
  retains both the front factor four and the exponent $-M^2/4$. Adding bulk and
  tail proves the finite-delta bound. No infinite-dimensional density or Laplace
  principle is used. This proves Proposition~\ref{P:signed-moment-bulk-tail}.
\end{proof}

\begin{proposition}[Two limits on the same moment]
  \label{P:signed-operator-upper}
  The global operator moment obeys
  \[
    0\le\liminf_{\delta\downarrow0}\delta\log M_\delta
    \le\limsup_{\delta\downarrow0}\delta\log M_\delta
    \le V_{\rm up}\le2+2C_{\rm mean}^2 .
  \]
\end{proposition}

\begin{proof}
  We choose a polynomial cutoff large enough for the quadratic tail to
  dominate. We first take the small-parameter limit at fixed block length, and
  only then send the block length to infinity, keeping the same endpoint
  moment on the left.

  \medskip
  \noindent\textit{Step 1.~}We choose a polynomial tail cutoff and bound the
  logarithm of the sum.
  Choose any $c_{\rm ball} > 0$ and $\kappa_{\rm ball} > 3\nu/(2-\nu)$, and set
  $M_B = c_{\rm ball}B^{\kappa_{\rm ball}}$. The cutoff in
  Proposition~\ref{P:signed-moment-bulk-tail} holds for all sufficiently large
  integers $B$. For $A_1 = 0$ this is immediate; otherwise its side ratio is
  \[
    \frac{M_B^{2-\nu}}{4A_1B^{3\nu}}
    =\frac{c_{\rm ball}^{2-\nu}}{4A_1}
    B^{\kappa_{\rm ball}(2-\nu)-3\nu}\longrightarrow\infty .
  \]
  Both constants and the strictly positive exponent are retained. At such a
  fixed $B$, a positive sum is at most twice its larger term. Taking
  $\delta\log$ of~\eqref{E:signed-finite-moment} therefore gives
  \begin{align*}
    \delta\log M_\delta\le\delta\log2+
    \max\{&
    \mathcal{S}_{R_\delta}+\mathcal{E}_{\delta,B,M_B}
    +\delta\log(M_B^2/(2\delta)),\\
    &A_0+\eta_Bb_{\delta,B}^\nu-M_B^2/4+\delta\log4\}.
  \end{align*}

  \medskip
  \noindent\textit{Step 2.~}We take the small-parameter limit before the
  block-length limit.
  At fixed $B$, $b_{\delta,B} \to 0$. All logarithmic normalizers vanish,
  including $\delta\log(M_B^2/(2\delta)) = \delta\log(M_B^2/2)-\delta\log\delta
  \to 0$. The last limit follows from $te^{-t} \le 2/t$ for $t = -\log\delta >
  0$. Finite variational expression and error bounds permit limsup
  subadditivity, giving
  \[
    \limsup_{\delta\downarrow0}\delta\log M_\delta
    \le\max\{V_{\rm up}+e_B,\ A_0-M_B^2/4\}.
  \]

  The lower bound $M_\delta\ge1$ used here is~\eqref{E:point-fixed-finiteness}
  of Lemma~\ref{L:point-actual-finiteness}.
  The left side has no $B$: restart identifies every block product with this
  same $P_\delta$. Now let integer $B \to \infty$. We have $e_B \to 0$, $M_B \to
  \infty$ and finite $V_{\rm up} \ge 0$. This gives the asserted upper bound,
  while $M_\delta \ge 1$ gives the nonnegative lower bound. The order was
  $\delta$ first at fixed block length, ball and thresholds, then $B$; no
  diagonal or limit of $\mathcal{S}_{R_\delta}$ was assumed. This proves
  Proposition~\ref{P:signed-operator-upper}.
\end{proof}

\begin{proposition}[Pointwise upper bound from the signed variational expression]
  \label{P:original-signed-point-upper}
  For every original admissible model and every fixed $x$,
  \[
    \limsup_{T\to\infty}T^{-\kappa}
    \log\E|u_\Gamma(T,x)|^2
    \le\limsup_{\delta\downarrow0}\mathcal{S}_{\delta^{-1}+2}
    \le2+2C_{\rm mean}^2 .
  \]
  For zero covariance the two normalized moments and the signed variational
  expression have rate zero.
\end{proposition}

\begin{proof}
  Proposition~\ref{P:original-point-transfer} bounds the original point moment
  by $\operatorname{DataPref}_\delta M_\delta$, with
  $\delta\log\operatorname{DataPref}_\delta \to 0$. Its proof, which removes
  the regularization first, keeps the full chain and works with the canonical
  law, is the input here, not an original $E$ point trace. The decreasing
  bijection $T_\delta =
  \delta^{-1/\kappa}$ has inverse $\delta = T^{-\kappa}$. Take logs, normalized
  limsup and Proposition~\ref{P:signed-operator-upper}. Original canonical
  spatial stationarity from Proposition~\ref{P:raw-block-original-lower} gives
  every fixed $x$. The moment is at least one by its full expansion and finite
  by the bound, so no logarithm or difference of infinite quantities is hidden.

  For $\Gamma = 0$ set $C_{\rm mean} = 0$ directly. Then $H_I = \{0\}$,
  $\operatorname{Rate}_R = \mathcal{S}_R = V_{\rm up} = 0$, $u = 1$ and
  $P_\delta = S_E(\delta^{-1})$. The free lower and upper bounds give $0 \le
  \delta\log M_\delta \le 2\delta\log[C(1+\delta^{-1})] \to 0$. This direct
  proof needs no nonzero Gaussian or nondegenerate Cameron assumption; the free
  functional at zero may exceed zero. The formal $a\downarrow0$ exponent is
  $(4-\alpha)/(3-\alpha)$, with no new endpoint-in-$a$ theorem asserted. The
  Jensen and lower bounds remain valid. This proves
  Proposition~\ref{P:original-signed-point-upper}.
\end{proof}

The varying variational expression remains a signed variational value for the
original covariance. Its positive-control comparison, exact limit and
identification of the derived constant require the subsequent variational
arguments.

\chapter{The constant and its consequences}\label{S:book-constant}

The lower and the upper bound now concern the same point moment, but
their variational problems differ. The lower bound uses nonnegative
smooth controls (Section~\ref{S:unmasked-variational}). The
upper bound allows signed controls in a completed covariance space
(Section~\ref{S:original-signed-score-upper}). Do the two optimizations have
the same value?

Section~\ref{S:finite-positive-gram} solves a finite
optimization problem whose Gram matrix has nonnegative
entries (Figure~\ref{Fig:active-subset-optimization}).
Section~\ref{S:intrinsic-second-moment} applies this to
the physical means, removes the spatial cutoff and reaches
Theorem~\ref{T:intrinsic-second-moment}. For the limit alone, stop there.

The remaining sections study the
constant. Section~\ref{S:intrinsic-occupation-ceiling} gives occupation
and spatial-energy bounds. For spatial white noise in one dimension,
Section~\ref{S:intrinsic-white-evaluation} evaluates the spatial energy and
obtains explicit moment bounds. These are bounds: the full supremum is not
evaluated, and the auxiliary value with coefficient one is not identified
with the actual first-moment constant. Section~\ref{S:route} lists what is
established and what is open. Episode~\VideoRef{8}{RllmgM_zB48} of the video
series covers this chapter.

\section{Finite optimization with nonnegative Gram entries}
\label{S:finite-positive-gram}

Let $H$ be any real Hilbert space, $\mathcal{I}$ a finite index set, $k_i \in
H$, $\ell_i > 0$ and $b_i \ge 0$. Write $C_{ij} = \langle k_i,k_j\rangle$ and
assume $C_{ij} \ge 0$ for every pair. The Gram matrix is symmetric and positive
semidefinite, since $\sum_{i,j}a_iC_{ij}a_j = \|\sum_i a_i k_i\|^2$. Put $t_+ =
\max(t,0)$ and define
\begin{equation}\label{E:finite-gram-score}
  \mathfrak{F}_C(\Phi)
  =2\sum_{i\in\mathcal{I}}\ell_i
  \sqrt{(\langle\Phi,k_i\rangle-b_i)_+}-\|\Phi\|^2/2,
  \qquad \mathscr{S}(C)=\sup_{\Phi\in H}\mathfrak{F}_C(\Phi).
\end{equation}
Let $\mathcal{J} = \{i:k_i\ne0\} = \{i:C_{ii}>0\}$. Zero vectors contribute zero
because $b_i \ge 0$. Empty index sets, singular Grams and $H = \{0\}$ are
allowed. For nonempty $A \subset \mathcal{J}$ and $w_i > 0$, set
\begin{equation}\label{E:finite-gram-dual}
  D_A^C(w)=\sum_{i\in A}\frac{\ell_i^2}{w_i}
  -\sum_{i\in A}b_iw_i
  +\frac12\sum_{i,j\in A}w_iC_{ij}w_j,
  \qquad d_A^C=\inf_{w_i>0}D_A^C(w).
\end{equation}
Set $d_\varnothing^C = 0$.

\LeanMarker{L:gram-scalar-algebra}
\begin{lemma}[Scalar and Hilbert identities]\label{L:gram-scalar-algebra}
  For real $\ell$, $w > 0$ and $t \ge 0$,
  \begin{equation}\label{E:gram-sqrt-dual}
    2\ell\sqrt t\le\ell^2/w+wt.
  \end{equation}
  For $\Phi,z \in H$,
  \begin{equation}\label{E:gram-hilbert-square}
    \langle\Phi,z\rangle-\|\Phi\|^2/2
    =\|z\|^2/2-\|\Phi-z\|^2/2\le\|z\|^2/2.
  \end{equation}
  For finite real $w_i$, $\|\sum_iw_ik_i\|^2 = \sum_{i,j}w_iC_{ij}w_j$.
\end{lemma}

\begin{proof}
  Since $\sqrt w > 0$ and $\sqrt{wt} = \sqrt w\sqrt t$,
  \[
    0\le(\ell/\sqrt w-\sqrt{wt})^2
    =\ell^2/w+wt-2\ell\sqrt t.
  \]
  This includes $t = 0$. If $\ell,t > 0$, equality holds for $w = \ell/\sqrt t$.
  Expanding $\|\Phi-z\|^2$ proves~\eqref{E:gram-hilbert-square}. Bilinearity of
  the real inner product proves the finite Gram identity. These identities use
  neither compactness of a Hilbert ball nor an inverse covariance matrix. This
  proves Lemma~\ref{L:gram-scalar-algebra}.
\end{proof}

\begin{lemma}[Interior dual minimum]\label{L:gram-dual-minimum}
  Each nonempty $A \subset \mathcal{J}$ has a unique positive minimizer $w^A$ of
  $D_A^C$. At this minimizer,
  \begin{equation}\label{E:gram-dual-stationarity}
    \sum_{j\in A}C_{ij}w_j^A-b_i
    =\frac{\ell_i^2}{(w_i^A)^2}>0,\qquad i\in A.
  \end{equation}
\end{lemma}

\begin{proof}
  We bound each coordinate of a minimizing sequence away from both zero and
  infinity. This gives an attained interior minimum. Strict positivity of the
  Hessian then gives uniqueness and the stationarity equations.

  \medskip
  \noindent\textit{Step 1.~}We confine a nonempty sublevel to a compact
  positive box.
  Put $c_i = C_{ii} > 0$ and $K_A = \sum_{i\in A}b_i^2/c_i$. Nonnegative entries
  and positive weights give
  \[
    D_A^C(w)\ge
    \sum_i\left(\ell_i^2/w_i-b_iw_i+c_iw_i^2/2\right).
  \]
  Completing the square in $c_iw_i^2/4-b_iw_i$ yields
  \begin{equation}\label{E:gram-dual-coercivity}
    D_A^C(w)\ge
    \sum_i\ell_i^2/w_i+\frac14\sum_i c_iw_i^2-K_A.
  \end{equation}
  Thus, no coordinate can approach zero or infinity in a bounded sublevel while
  a different negative term hides its growth. Take $w^0 = (1,\ldots,1)$, $d_0 =
  D_A^C(w^0)$ and $L_A = d_0+K_A$. Equation~\eqref{E:gram-dual-coercivity} at
  $w^0$ gives $L_A \ge \sum_i\ell_i^2+\sum_i c_i/4 > 0$. Every $w$ with
  $D_A^C(w) \le d_0$ satisfies
  \[
    \ell_i^2/L_A\le w_i\le\sqrt{4L_A/c_i}.
  \]
  The sublevel is nonempty and closed inside this compact box, which lies in the
  open positive orthant. A minimizing sequence chosen in the sublevel has a
  convergent subsequence there. Continuity attains its infimum. Outside this
  sublevel the values exceed $d_0$, so the minimum is global on the entire
  positive orthant.

  \medskip
  \noindent\textit{Step 2.~}We prove uniqueness and compute the interior
  critical equations.
  The function is smooth in that orthant. Its Hessian on a nonzero vector $\xi$
  has value
  \[
    \sum_i\frac{2\ell_i^2}{w_i^3}\xi_i^2
    +\sum_{i,j}\xi_iC_{ij}\xi_j>0.
  \]
  The first term is strictly positive and the second is a squared Hilbert norm.
  Restriction to any nonconstant line segment is strictly convex, proving
  uniqueness. Each coordinate permits both signs of a small perturbation at the
  interior minimum. Differentiation gives $-\ell_i^2/w_i^2-b_i+\sum_jC_{ij}w_j =
  0$. This is~\eqref{E:gram-dual-stationarity}. This proves
  Lemma~\ref{L:gram-dual-minimum}.
\end{proof}

\LeanMarker{P:gram-positive-optimizer}
\begin{proposition}[Positive conic maximizer]\label{P:gram-positive-optimizer}
  The finite real value in~\eqref{E:finite-gram-score} satisfies
  \begin{equation}\label{E:gram-active-duality}
    \mathscr{S}(C)=\max_{A\subset\mathcal{J}}d_A^C .
  \end{equation}
  If a maximizing subset is empty, $\Phi = 0$ attains the supremum. If a
  maximizing subset $A$ is nonempty, a maximizer is $\Phi_A = \sum_{i\in A}w_i^A
  k_i$ with strictly positive coefficients. No uniqueness of the primal
  maximizer or maximizing subset is asserted.
\end{proposition}

\begin{proof}
  We bound each primal value using its own active set. We then evaluate the
  restricted objective at the corresponding dual minimizer, where both scalar
  and Hilbert inequalities are equalities, and use the nonnegative omitted
  rewards.

  \medskip
  \noindent\textit{Step 1.~}We bound the primal objective by a finite
  active-set dual value.
  For arbitrary $\Phi$, let $A(\Phi) = \{i:\langle\Phi,k_i\rangle>b_i\}$. It is
  contained in $\mathcal{J}$. All other terms in the objective vanish, including
  equality at the threshold. If $A(\Phi)$ is empty, $\mathfrak{F}_C(\Phi) =
  -\|\Phi\|^2/2 \le 0$. Otherwise, for any positive $w$ on $A = A(\Phi)$, set
  $t_i = \langle\Phi,k_i\rangle-b_i > 0$ and $z = \sum_{i\in A}w_ik_i$.
  Lemma~\ref{L:gram-scalar-algebra} gives
  \[
    \mathfrak{F}_C(\Phi)
    \le\sum_i\ell_i^2/w_i-\sum_i b_iw_i
    +\langle\Phi,z\rangle-\|\Phi\|^2/2
    \le D_A^C(w).
  \]
  This holds for every positive $w$, so $\mathfrak{F}_C(\Phi) \le d_A^C$. There
  are finitely many subsets and their values are finite by
  Lemma~\ref{L:gram-dual-minimum}. Taking the supremum bounds $\mathscr{S}(C)$
  above by their maximum.

  \medskip
  \noindent\textit{Step 2.~}We realize the maximum of the dual values by a
  nonnegative conic control.
  Conversely, $\mathfrak{F}_C(0) = 0 = d_\varnothing^C$. For nonempty $A$, take
  its dual minimizer. Equation~\eqref{E:gram-dual-stationarity} makes every
  vector in $A$ strictly active at $\Phi_A$, with
  $\sqrt{\langle\Phi_A,k_i\rangle-b_i} = \ell_i/w_i^A$. Both inequalities in the
  preceding display are then equalities for the sum restricted to $A$, since $z
  = \Phi_A$. This restricted objective equals $d_A^C$. All omitted vectors
  contribute nonnegative terms to the full objective, hence
  $\mathfrak{F}_C(\Phi_A) \ge d_A^C$, even if this value is negative. Maximizing
  over $A$, including the empty set, proves \eqref{E:gram-active-duality} and
  the asserted attainment. For a maximizing nonempty subset the upper and lower
  bounds force $\mathfrak{F}_C(\Phi_A) = d_A^C$; no outside inactivity was
  assumed. All-zero vectors give $\mathscr{S}(C) = 0$ directly. This proves
  Proposition~\ref{P:gram-positive-optimizer}.
\end{proof}

No-go~\ref{P:nogo-29} gives the three finite counterexamples behind this
formulation: the empty active set cannot be omitted, the primal objective
need not be concave, and entrywise matrix order is not positive
semidefinite order.

\begin{figure}
\centering
\begin{tikzpicture}[x=1cm,y=1cm]
  \node[ham/box,text width=7.8cm] (hamRoot) at (5.1,1.7)
  {Discard zero vectors: $\mathcal J=\{i:k_i\ne0\}$.\\
  Illustrate all subsets when $\mathcal J=\{1,2\}$.};
  \foreach \hamSubsetId/\hamSubsetX/\hamSubsetName/\hamSubsetValue in
  {empty/0/{\varnothing}/{0},first/3.4/{\{1\}}/{d_{\{1\}}^C},second/6.8/{\{2\}}/{d_{\{2\}}^C},both/10.2/{\{1,2\}}/{d_{\{1,2\}}^C}} {
  \node[ham/box,text width=2.7cm] (hamSubset\hamSubsetId) at (\hamSubsetX,-.15)
  {$A=\hamSubsetName$\\value $\hamSubsetValue$};
  \draw[ham/arrow] (hamRoot)--(hamSubset\hamSubsetId);
  }
  \node[ham/note,text width=10.6cm] at (5.1,-1.8)
  {For every nonempty $A$, minimize $D_A^C(w)$ over $w_i>0$.\\
  The weight minimizer $w^A$ is unique; $d_A^C$ need not be positive.};
  \node[ham/box,text width=8.4cm] (hamMax) at (5.1,-3.4)
  {$\mathscr S(C)=\max\{0,d_{\{1\}}^C,d_{\{2\}}^C,d_{\{1,2\}}^C\}$\\
  Take a maximizing subset; if nonempty, set $\Phi_A=\sum_{i\in A}w_i^Ak_i$.};
  \node[ham/note,text width=10.7cm] at (5.1,-5.1)
  {The omitted terms are nonnegative at $\Phi_A$.\\
  Their inactivity is not assumed in the construction.};
\end{tikzpicture}
\caption[Minimize within each subset, then maximize over subsets]{The finite-subset identity includes the empty set with value zero.
For each nonempty subset the inner problem is a strictly convex minimization
in positive weights, followed by a maximum over subsets. This does not assert
global concavity of the original objective, uniqueness of its maximizer,
or an effective algorithm. No inverse Gram matrix is used, so coincident
vectors and singular Gram matrices remain allowed. See
\ProofRef{Lemma}{L:gram-dual-minimum} and
\ProofRef{Proposition}{P:gram-positive-optimizer}.}
\label{Fig:active-subset-optimization}
\end{figure}
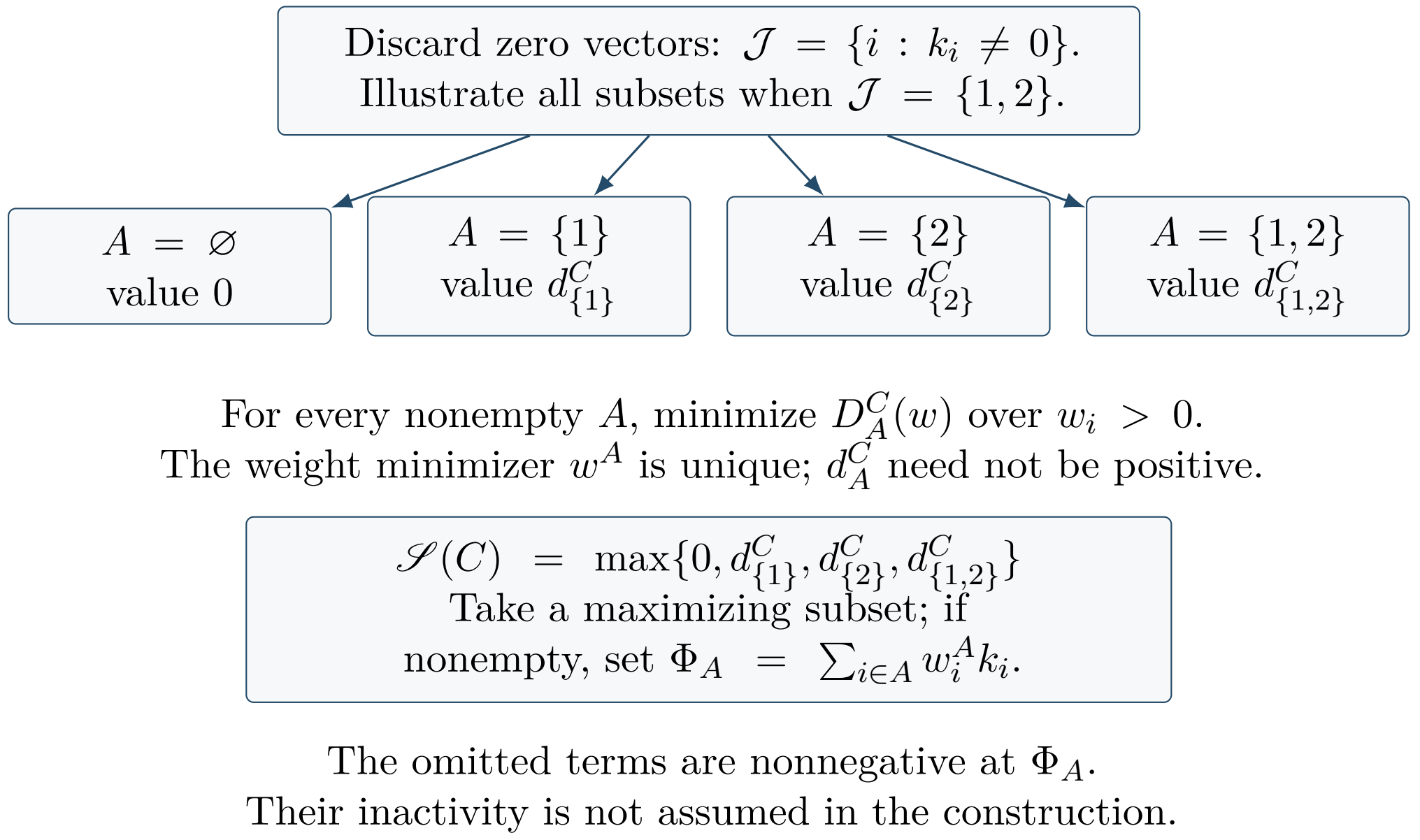

\begin{figure}
  \centering
  \begin{tikzpicture}[x=1cm,y=1cm]
    \begin{scope}
      \fill[hamteal!9] (0,0)--(4,.89655)--(4,3)--(2.1875,3)--cycle;
      \draw[ham/arrow] (0,0)--(2.9,.65) node[right] {$k_1$};
      \draw[ham/arrow] (0,0)--(1.75,2.4) node[above] {$k_2$};
      \draw[ham/arrow,hamteal,very thick] (0,0)--(3.0,2.15)
      node[anchor=west] {$\Phi_A=\sum_{i\in A}w_i^Ak_i$};
      \node[ham/note] at (1.8,-.6) {$\langle k_i,k_j\rangle\ge0,\quad w_i^A>0$};
    \end{scope}
    \begin{scope}[xshift=8.3cm]
      \draw[ham/arrow] (0,0)--(2.2,1.1) node[right] {$k_2=2k_1$};
      \draw[ham/arrow,thick] (0,0)--(1.1,.55) node[below] {$k_1$};
      \filldraw[fill=white,draw=hamgold,thick] (0,0) circle (2pt);
      \node[ham/note,text width=3.2cm] at (1.0,2.6)
      {Singular Gram matrices\\are allowed};
      \node[ham/note,text width=3.2cm] at (1,-1.45)
      {A zero vector contributes zero.\\The active set may be empty.};
    \end{scope}
  \end{tikzpicture}
  \caption[Positive combinations and singular Gram matrices]{Finite positive
    combinations, not an inverse covariance formula. For a nonempty maximizing
    active set, a maximizer has positive coefficients on its selected vectors.
    The right panel shows why linear independence must not be imposed. The
    drawing is in a finite-dimensional span and does not restrict the ambient
    Hilbert space. Neither the maximizing subset nor the maximizer is asserted
    to be unique. See \ProofRef{Proposition}{P:gram-positive-optimizer}.}
  \label{Fig:positive-gram}
\end{figure}
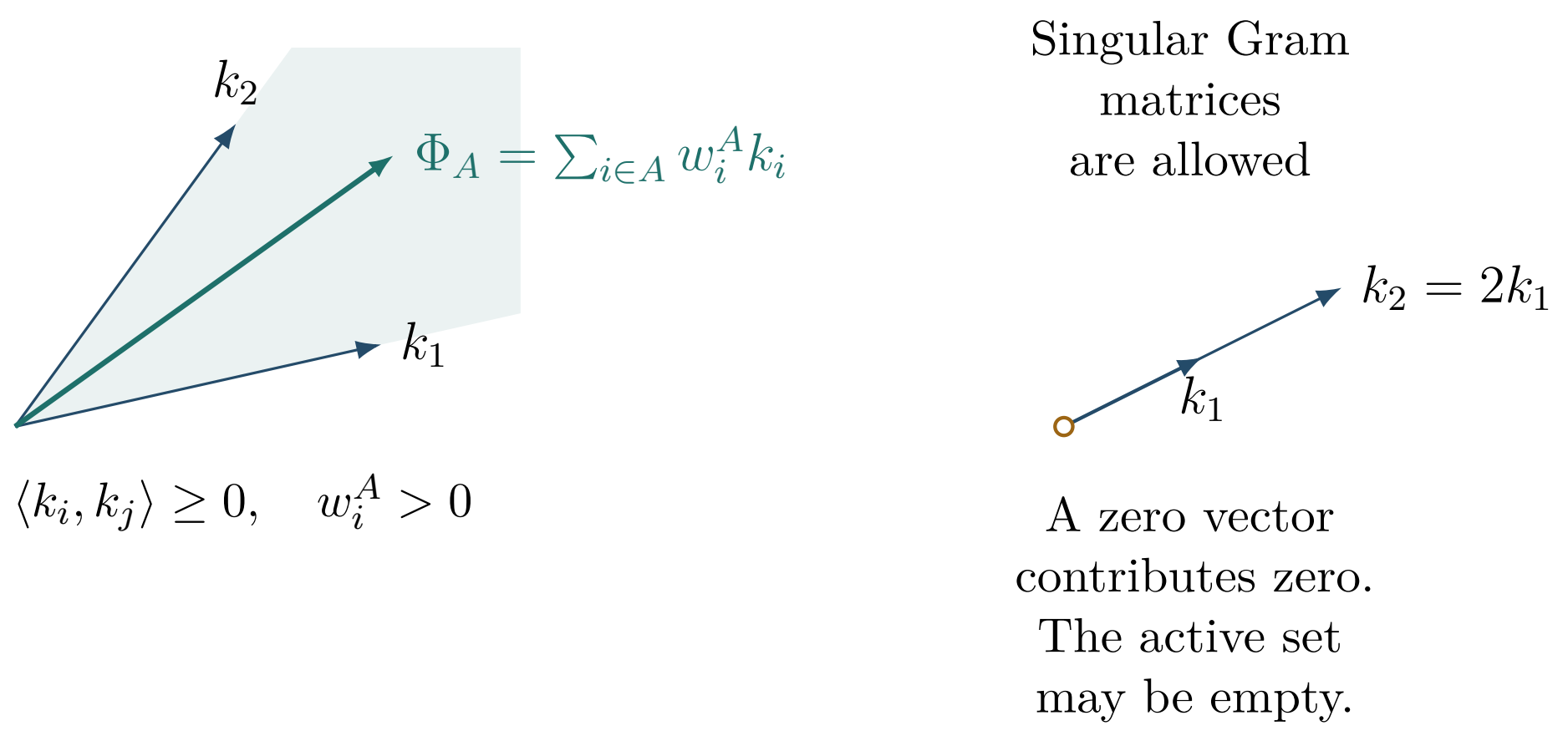

\begin{lemma}[Real Gram realization]\label{L:gram-realization}
  Every finite real symmetric positive semidefinite matrix is the Gram matrix of
  vectors in a finite-dimensional real Hilbert space. When its entries are
  nonnegative, the value $\mathscr{S}(C)$ is independent of its realization.
\end{lemma}

\begin{proof}
  Induct on the size, with the empty matrix realized in $\{0\}$. Diagonal
  entries are nonnegative. If all vanish, testing on $e_i+t e_j$ gives $2tC_{ij}
  \ge 0$ for every real $t$, hence $C = 0$, realized by zero vectors. Otherwise
  permute a positive diagonal $c$ into the first position and write $C =
  \left[\begin{smallmatrix}c&v^{\mathsf{T}}\\v&A\end{smallmatrix}\right]$. The
  Schur complement $B = A-vv^{\mathsf{T}}/c$ is positive semidefinite: test the
  original quadratic form on $(-v^{\mathsf{T}}y/c,y)$. Realize $B$ by vectors
  $u_i$ using the induction hypothesis. Then $k_1 = (\sqrt c,0)$ and $k_i =
  (v_i/\sqrt c,u_i)$ for the other indices realize $C$. Undo the permutation.
  The inner products are $\langle k_1,k_i\rangle=v_i$ and $\langle
  k_i,k_j\rangle=v_iv_j/c+B_{ij}=A_{ij}$, while $\|k_1\|^2=c$.
  Entries of $B$ need not be nonnegative; this construction requires neither
  nonnegative coordinates nor complete positivity.

  In any realization, $C_{ii} = 0$ means $k_i = 0$, so its entire row vanishes.
  Formula~\eqref{E:gram-active-duality} depends only on $C,\ell,b$ and these
  eligible indices. It therefore proves realization independence whenever the
  entries are nonnegative, which is the scope in which $\mathscr{S}(C)$ is used.
  This proves Lemma~\ref{L:gram-realization}.
\end{proof}

By Lemma~\ref{L:gram-realization}, $\mathscr{S}(C)$ depends only on the matrix
$C$ and the data $\ell,b$. We can therefore compare two Gram matrices directly.

\LeanMarker{P:gram-entrywise-comparison}
\begin{proposition}[Entrywise Gram comparison]\label{P:gram-entrywise-comparison}
  Let $C,D$ be real symmetric positive semidefinite matrices on the same finite
  index set, with $0 \le C_{ij} \le D_{ij}$. With the same $\ell_i > 0,b_i \ge
  0$,
  \[
    \mathscr{S}(C)\le\mathscr{S}(D).
  \]
  No positive-semidefinite order on $D-C$ is required.
\end{proposition}

\begin{proof}
  Lemma~\ref{L:gram-realization} defines both values. Their eligible sets
  satisfy $\mathcal{J}_C \subset \mathcal{J}_D$. For nonempty $A \subset
  \mathcal{J}_C$ and positive weights,
  \[
    D_A^D(w)-D_A^C(w)
    =\frac12\sum_{i,j\in A}w_i(D_{ij}-C_{ij})w_j\ge0.
  \]
  Taking infima gives $d_A^C \le d_A^D$. The empty values are both zero. Taking
  maxima and using the inclusion of eligible families proves the claim by
  Proposition~\ref{P:gram-positive-optimizer}. This includes the appearance of
  previously zero vectors. This proves
  Proposition~\ref{P:gram-entrywise-comparison}.
\end{proof}

\begin{remark}[Exact checks and scope]\label{R:gram-controls}
  The empty value is essential: for one vector with $C = 1$, $\ell = 1/10,b =
  4$, the nonempty dual has $D(4) = -3199/400 < 0$, whereas the primal objective
  at zero is zero. Global concavity is unavailable: for $k = \ell = b = 1$ the
  objective $2\sqrt{(\phi-1)_+}-\phi^2/2$ decreases on $[0,1]$ and has infinite
  positive right slope at one. For $k = \ell = 1,b = 0$, differentiation on
  $\phi > 0$ gives the maximum $3/2$ at $\phi = 1$; on $\phi \le 0$ the
  objective is nonpositive. Doubling the cost changes the positive critical
  point to $2^{-2/3}$ and the maximum to $3\,2^{-4/3}$. Finally, $C = I_2$, $D =
  \left[\begin{smallmatrix}1&1/2\\1/2&1\end{smallmatrix}\right]$ are nonnegative
  positive semidefinite Grams but $D-C$ has eigenvalues $1/2,-1/2$. These are
  finite analytic checks, not a covariance-order theorem for physical waves.

  The scalar, finite Gram, conditional active upper and half-quadratic entrywise
  deductions have existing local Lean certificates at their exact hypotheses.
  Dual coercivity, calculus, attainment, realization and the infimum/maximum
  argument above remain analytic proofs. No new kernel fact or formalization of
  the physical application follows merely from these certificates.
\end{remark}

\section{A variational formula for the second-moment limit}
\label{S:intrinsic-second-moment}

Retain every original covariance and amplitude. Write $a = \alpha_0$, $I =
(0,1)$, $K(s,t) = |s-t|^{-a}$ and $H_I = H_t\otimes H_{\mathrm{sp}}$ for the
original covariance. Keep the uniform $E$ mean $\mathcal{T}_R = T_R(C_t\otimes
I)$ and $C_{\rm mean} = M_0\|C_t\|$ of
Section~\ref{S:uniform-physical-products}, with exactly the fixed $r$ chosen
there. In particular, $2 < r < 2/a$. To distinguish it from the map below, write
$J_R = T_R$. For a bounded multiplier $V:H^{1/2} \to H^{-1/2}$ write
$B_E(V)(u,v) = (0,Vu)$, the $E$ operator of Lemma~\ref{L:uniform-E-lift}. The
physical mask is $m_R = \eta(\cdot/R)^2$, used once. Let
$\mathcal{C}_{\mathrm{pos}}$ be the full nonseparable class
in~\eqref{E:variational-positive-class}. For $\Psi$ in that class, retain its
unmasked mean $h_{\Psi,s}$, original variance $v_\Psi$, integrated full-gradient
rate $I_\Psi$ and variational values for nonnegative controls
$S_c^{\mathrm{pos}}$ from Section~\ref{S:unmasked-variational}. The value with
coefficient one is not identified with a first-moment limit or a constant in the
cited source. We assume $\Gamma \ne 0$ through the following comparison lemmas.
The theorem treats zero covariance directly.

\LeanMarker{L:positive-raw-embedding}
\begin{lemma}[Inclusion into the covariance space and the smooth mean]
  \label{L:positive-raw-embedding}
  For any real Hilbert space $G$ there is a bounded map $J_G:L^2(I;G) \to
  H_t\otimes G$ with
  \begin{equation}\label{E:positive-bochner-inclusion}
    \|J_GF\|^2
    =\int_{I^2}K(s,t)\langle F(s),F(t)\rangle_G\ud s\ud t
    \le K_0\|F\|_{L^2(I;G)}^2,\qquad K_0=\frac{2^a}{1-a}.
  \end{equation}
  The natural inclusion puts every $\Psi \in \mathcal{C}_{\mathrm{pos}}$ in
  $H_I$, with
  \begin{equation}\label{E:positive-raw-mean}
    \|\Psi\|_{H_I}^2=v_\Psi,\qquad
    \mathcal{T}_R\Psi(s)=B_E(m_Rh_{\Psi,s})\quad\hbox{a.e. }s.
  \end{equation}
  For each fixed smooth $\Psi$ the right side gives a canonical version for
  every $R \ge 1$ and every $s \in I$.
\end{lemma}

\begin{proof}
  We construct the covariance inclusion from the temporal Schur bound and
  identify its completed temporal mean. We then use common compact spatial
  supports to identify a smooth positive test with its original covariance
  class and physical masked mean.

  \medskip
  \noindent\textit{Step 1.~}We construct the Hilbert-valued covariance
  inclusion.
  Direct integration and concavity of $x^{1-a}$ give
  \[
    \int_I K(s,t)\ud t
    =\frac{s^{1-a}+(1-s)^{1-a}}{1-a}\le K_0 .
  \]
  For scalar $f \ge 0$, Cauchy--Schwarz against $K(s,t)\ud t$ and then Tonelli
  give $\|Kf\|_2 \le K_0\|f\|_2$. The same estimate for absolute values bounds
  the bilinear integral by $K_0\|f\|_2\|g\|_2$. Approximate scalar $L^2(I)$
  functions by real $C_c^\infty(I)$ functions. This bound makes their original
  test-null classes Cauchy in $H_t$, defines the bounded inclusion, and
  identifies its pairing with the integral. Finite $G$-valued tensors then give
  \eqref{E:positive-bochner-inclusion}. Simple fields are dense in Bochner
  $L^2(I;G)$; the same bound on differences extends $J_G$. Absolute bilinear
  continuity extends the integral identity as well.

  \medskip
  \noindent\textit{Step 2.~}We identify the temporal integral with the
  completed tensor map.
  For $G = H_{\mathrm{sp}}$, the Bochner integral $\int_IK(s,t)F(t)\ud t$ exists
  almost everywhere: its norm is bounded by $K\|F\|(s)$, which belongs to $L^2$.
  Simple-field approximation and the Schur estimate give a strongly measurable
  $L^2$ version and a bounded map of norm at most $K_0$. On smooth finite
  tensors it agrees with $(C_t\otimes I)J_G$ by the covariance definition.
  Approximation by such tensors converges in $L^2(I;G)$. The integral converges
  there by Schur; the completed tensor map converges in $L^r(I;G)$ by
  Lemma~\ref{L:uniform-mean-tensor} and~\eqref{E:positive-bochner-inclusion},
  hence also in $L^2$. Uniqueness of the $L^2$ limit identifies these versions.
  In particular, the integral belongs to $L^r$.

  \medskip
  \noindent\textit{Step 3.~}We identify the original norm of a smooth
  spacetime test.
  Choose a closed ball $B$ containing all spatial supports of $\Psi$. For real
  smooth $q$ supported there, the physical covariance gives
  \begin{equation}\label{E:positive-local-spatial-bound}
    \|q\|_{H_{\mathrm{sp}}}^2
    \le |B|\,\Gamma(B-B)\|q\|_\infty^2.
  \end{equation}
  Indeed, in $\int\ud x\int\Gamma(\ud y)q(x)q(x-y)$ only $x \in B,y \in B-B$
  contribute. Absolute values and local finiteness prove the bound, including
  white and singular Radon covariances. Thus, $t\mapsto\Psi_t$ is continuous,
  bounded and compactly supported as an $H_{\mathrm{sp}}$-valued field, and
  belongs to $L^2$. Its image under $J_{H_{\mathrm{sp}}}$ is the original test
  class of $\Psi$. Indeed, the separated compact smooth approximants in
  Lemma~\ref{L:cameron-supported-tensor} converge uniformly on one common
  compact support. The local bound above and the temporal Schur estimate give
  convergence both in $L^2(I;H_{\mathrm{sp}})$ and in the original covariance
  norm. Their finite tensor identifications therefore have this same limit.
  Formula~\eqref{E:positive-bochner-inclusion} and absolute Fubini give
  \[
    \|\Psi\|_{H_I}^2
    =\int_{I^2}K(s,t)\langle\Psi_s,\Psi_t\rangle_{H_{\mathrm{sp}}}
    \ud s\ud t
    =\int_{I\times\R^d}\Psi(s,x)h_{\Psi,s}(x)\ud x\ud s=v_\Psi.
  \]
  The compact spatial supports and~\eqref{E:positive-local-spatial-bound} make
  all these absolute integrals finite.

  \medskip
  \noindent\textit{Step 4.~}We identify its masked operator mean on all phase
  states.
  Test the spatial row~\eqref{E:uniform-mean-row} first with displacement $f$
  and dual second component $\Lambda g$, with $f,g$ compact smooth. The row
  product is $m_Rfg$, so the preceding temporal identity gives
  \[
    \left\langle\int_IK(s,t)\Psi_t\ud t,m_Rfg\right\rangle_{H_{\mathrm{sp}}}
    =\int m_Rh_{\Psi,s}fg .
  \]
  This Fubini step has two fixed compact supports and integrable time kernel.
  Both resulting operators are bounded on $E$: the left by the mean theorem, the
  right by the bounded smooth $h_{\Psi,s}$ and compact smooth mask. Dense states
  extend the identity to $E$. For this fixed $\Psi$ the temporal integral exists
  at every $s$, with norm at most $K_0\sup_t\|\Psi_t\|_{H_{\mathrm{sp}}}$.
  Applying $J_R$ gives the explicit version in~\eqref{E:positive-raw-mean} for
  every $R,s$. For each $R$ it agrees with the completed $L^r$ mean class. This
  makes no simultaneous-version assertion over all rough controls. This proves
  Lemma~\ref{L:positive-raw-embedding}.
\end{proof}

Proposition~\ref{P:original-signed-point-upper} bounds the original point
limsup by a limsup of the signed values $\mathcal{S}_R$.
Proposition~\ref{P:variational-positive-scores} bounds the original point
liminf below by $S_2^{\mathrm{pos}}$. To close this gap in
Theorem~\ref{T:intrinsic-second-moment} we compare $\mathcal{S}_R$ with
$S_2^{\mathrm{pos}}$. We start with finite Gram matrices, with and without the
mask.

\LeanMarker{L:positive-physical-gram}
\begin{lemma}[Finite covariance vectors and entrywise order]
  \label{L:positive-physical-gram}
  Take a finite partition of $I$ with nonempty cells $I_i$, lengths $\ell_i > 0$
  and $\sum_i\ell_i = 1$, and real $f_i \in C_c^\infty(\R^d)$ with $\|f_i\|_2 =
  1$. Put
  \[
    h_i=\mathbf1_{I_i}/\ell_i,\quad b_i=\|\nabla f_i\|_2^2,\quad
    k_i^R=h_i\otimes(m_Rf_i^2),\quad
    k_i^\infty=h_i\otimes f_i^2 .
  \]
  These are elements of $H_I$. If $v_L$ is the multiplier form
  of~\eqref{E:uniform-E-pairing}, then for every $\Phi \in H_I$,
  \begin{equation}\label{E:positive-label-row}
    \langle\Phi,k_i^R\rangle_{H_I}
    =\ell_i^{-1}\int_{I_i}v_{\mathcal{T}_R\Phi(s)}(f_i,f_i)\ud s .
  \end{equation}
  The Grams $C^R_{ij} = \langle k_i^R,k_j^R\rangle$ and $D_{ij} = \langle
  k_i^\infty,k_j^\infty\rangle$ satisfy $0 \le C^R_{ij} \le D_{ij}$ entrywise.
\end{lemma}

\begin{proof}
  We realize the cell labels in the original tensor Hilbert space and identify
  their rows by completion. We then compare their Gram entries through the
  nonnegative temporal and physical spatial integrands.

  \medskip
  \noindent\textit{Step 1.~}We identify the finite label rows.
  Lemma~\ref{L:positive-raw-embedding} puts $h_i$ in $H_t$. The compact smooth
  spatial factors lie in $H_{\mathrm{sp}}$ by its test definition
  and~\eqref{E:positive-local-spatial-bound}. On a smooth finite tensor,
  use~\eqref{E:uniform-mean-row} with $z = (f_i,0)$ and $w = (0,\Lambda f_i)$.
  Its row is $(\eta_Rf_i)(\eta_R\Lambda^{-1}\Lambda f_i) = m_Rf_i^2$. Temporal
  integration gives the pairing with $h_i$. This
  proves~\eqref{E:positive-label-row} on dense tensors. Its right side is a
  bounded original functional, with bound $\ell_i^{-1/r}C_{\rm
  mean}\|\Phi\|\|f_i\|_{H^{1/2}}^2$ by Hölder
  and~\eqref{E:uniform-completed-mean}. The left side is bounded by Hilbert
  Cauchy--Schwarz. Thus, the identity extends to the completed $H_I$. The form
  is the displacement form by Lemma~\ref{L:uniform-E-lift}.

  \medskip
  \noindent\textit{Step 2.~}We compare the physical Gram integrands entrywise.
  Each Gram is symmetric positive semidefinite. Its entries factor into the
  temporal pairing of nonnegative $h_i,h_j$ and the physical spatial covariance
  pairing. The kernels $K,\Gamma$ and $f_i^2,f_j^2$ are nonnegative; $0 \le m_R
  \le 1$. Tonelli therefore gives the entrywise order. The physical pairing
  already contains $\Gamma$, and the row cancellation introduces no additional
  Fourier factor. Zero labels are permitted. All hypotheses of
  Proposition~\ref{P:gram-entrywise-comparison} are met, without assuming
  $D-C^R$ is positive semidefinite. This proves
  Lemma~\ref{L:positive-physical-gram}.
\end{proof}

By Proposition~\ref{P:gram-positive-optimizer}, the unmasked finite problem has
a maximizer built from the $k_i^\infty$. Their time factors $h_i$ are scaled
indicators, so this maximizer is in general not a test in
$\mathcal{C}_{\mathrm{pos}}$. We also need a bound on the controls that holds
uniformly in $R$.

\begin{lemma}[Uniform optimizer norm and positive smoothing]
  \label{L:positive-uniform-optimizer}
  Put $M = 4C_{\rm mean}+4$. The supremum defining $\mathcal{S}_R$ may be
  restricted to $\|\Phi\| \le M$, uniformly in $R$. For every finite trial
  family in Lemma~\ref{L:positive-physical-gram}, the unmasked finite
  variational expression has a maximizer $\Phi_D$ with $\|\Phi_D\| \le M$. It is
  zero or is a positive conic combination of the $k_i^\infty$. In the latter
  case it is approximable in $H_I$ by tests in $\mathcal{C}_{\mathrm{pos}}$.
\end{lemma}

\begin{proof}
  We first restrict the signed score to a uniform control ball. We prove that
  the finite unmasked optimizer lies in that same ball before using any
  averaging estimate. Finally, its positive conic structure permits
  nonnegative smooth approximation with convergent cost and rate.

  \medskip
  \noindent\textit{Step 1.~}We restrict the signed optimization to a fixed
  ball.
  Put $C = C_{\rm mean}$ and $q = \|\Phi\|$.
  Proposition~\ref{P:uniform-rate-averaging} gives
  \[
    2\operatorname{Rate}_R(\Phi)-q^2/2
    \le2+2Cq-q^2/2<2-2q<0\quad(q>M).
  \]
  The zero control gives zero, proving the restriction.

  \medskip
  \noindent\textit{Step 2.~}We bound the finite unmasked optimizer uniformly.
  Denote the finite unmasked objective by $\mathfrak{F}_D$.
  Proposition~\ref{P:gram-positive-optimizer} supplies an attained maximum at
  zero or a positive conic sum $\Phi_D$. Choose a finite $R' \ge 1$ whose
  plateau contains all supports of $f_i$. Then $k_i^{R'} = k_i^\infty$. Write
  $W_i = \ell_i^{-1}\int_{I_i}\mathcal{T}_{R'}\Phi(s)\ud s$. The row identity
  and real Rayleigh formula imply, for every $\Phi$,
  \begin{align*}
    \mathfrak{F}_D(\Phi)
    &\le2\sum_i\ell_i\lambda(W_i)-q^2/2\\
    &\le2+2\sum_i\ell_i\|W_i\|-q^2/2
    \le2+2Cq-q^2/2 .
  \end{align*}
  For fixed $f_i$, the map $L\mapsto v_L(f_i,f_i)$ is bounded linear; hence the
  form of the Bochner average is the average form. Here $\lambda(L) \le
  \max(1,\|L\|) \le 1+\|L\|$, the norm of an average is bounded by the average
  norm, and the $L^1$ mean norm is at most its $L^r$ norm. No time-averaging
  error or concavity of the positive-part square root enters this estimate.
  Since the maximum is nonnegative, its optimizer lies in the same ball of
  radius $M$, independently of the partition, profiles and $R'$.

  \medskip
  \noindent\textit{Step 3.~}We smooth the positive temporal factors without
  changing the limit class.
  On a nonzero branch write $\Phi_D = \sum_{i\in A}w_i h_i\otimes f_i^2$, with
  $w_i > 0$. Extend $h_i$ by zero, convolve with a nonnegative compact smooth
  approximate identity, and multiply by smooth cutoffs $0 \le \chi_j \le 1$
  compactly supported in $I$, tending pointwise to one. The resulting
  nonnegative $h_{i,j} \in C_c^\infty(I)$ converge to $h_i$ in $L^2$: the error
  is bounded by the convolution error plus $\|(\chi_j-1)h_i\|_2$, which tends to
  zero by domination. Set $\Psi_j = \sum_{i\in A}w_i h_{i,j}f_i^2$. Finite
  tensor estimates and~\eqref{E:positive-bochner-inclusion} give $\Psi_j \to
  \Phi_D$ in $H_I$. The positive conic sum has positive squared norm, since at
  least one eligible diagonal is positive and all Gram entries are nonnegative.
  Thus, $\Psi_j \ne 0$ eventually, so $\Psi_j \in \mathcal{C}_{\mathrm{pos}}$.
  Costs converge. The integrated rate continuity proved in
  Proposition~\ref{P:uniform-rate-averaging} gives
  \begin{equation}\label{E:positive-rate-continuity}
    |\operatorname{Rate}_{R'}(\Psi_j)-\operatorname{Rate}_{R'}(\Phi_D)|
    \le\sqrt{3+C\|\Psi_j\|+C\|\Phi_D\|}
    \sqrt{C\|\Psi_j-\Phi_D\|}\longrightarrow0 .
  \end{equation}
  This estimate is uniform in $R'$. It asserts no uniform derivative bound on
  the smooth approximants. The zero branch will use $S_2^{\mathrm{pos}} \ge 0$
  directly. This proves Lemma~\ref{L:positive-uniform-optimizer}.
\end{proof}

Lemma~\ref{L:positive-physical-gram} orders the masked and unmasked Grams
entrywise. Lemma~\ref{L:positive-uniform-optimizer} keeps the unmasked optimizer
in a fixed ball and, when it is nonzero, approximates it by nonnegative tests.
The time averaging of Proposition~\ref{P:uniform-rate-averaging} links these
finite problems to the rate.

\begin{proposition}[Signed variational values bounded by the variational value over nonnegative controls]
  \label{P:positive-signed-comparison}
  For every $R \ge 1$,
  \begin{equation}\label{E:positive-signed-comparison}
    \mathcal{S}_R\le S_2^{\mathrm{pos}} .
  \end{equation}
\end{proposition}

\begin{proof}
  We approximate each integrated signed rate by a finite collection of smooth
  Rayleigh trials. Entrywise Gram comparison replaces the optimized masked
  value by an unmasked positive conic optimizer. Its uniform norm bound allows
  the same averaging error to be used again, after which positive smoothing
  gives the nonnegative-control bound.

  \medskip
  \noindent\textit{Step 1.~}We approximate the integrated rate uniformly on
  the fixed control ball.
  Use uniform partitions with $n \ge 2$ cells of length $1/n$. They are the $B =
  1,\delta = 1/n$ partitions in Proposition~\ref{P:uniform-rate-averaging}, with
  the empty last cell omitted. That proposition gives $e_n \to 0$, independent
  of $R$, such that
  \[
    \left|\operatorname{Rate}_R(\Phi)
    -\sum_i\ell_i\lambda\left(\ell_i^{-1}
    \int_{I_i}\mathcal{T}_R\Phi(s)\ud s\right)\right|\le e_n
    \quad(\|\Phi\|\le M).
  \]

  \medskip
  \noindent\textit{Step 2.~}We choose finitely many near-maximizing spatial
  trials.
  Fix $R,\Phi,n$, with $\|\Phi\| \le M$, and denote these averages by $W_i$.
  Real compact smooth unit trials are dense in the real $H^1$ unit sphere:
  truncate by cutoffs tending to one, whose gradient errors tend to zero, then
  mollify in $H^1$ and normalize in $L^2$. Continuity of the multiplier form in
  $H^{1/2}$ and hence $H^1$ proves
  \[
    \lambda(W_i)=
    \sup_{\substack{f\in C_c^\infty(\R^d;\R)\\\|f\|_2=1}}
    \sqrt{(v_{W_i}(f,f)-\|\nabla f\|_2^2)_+}.
  \]
  The real Rayleigh supremum is nonnegative by
  Proposition~\ref{P:uniform-E-spectral}, so the positive part does not change
  its supremum. For $\varepsilon > 0$ choose one near-maximizing $f_i$ per cell
  with error at most $\varepsilon$ in the square root. This is a finite choice,
  not a measurable continuum optimizer. Since $\sum_i\ell_i = 1$, the weighted
  trial error contributes $2\varepsilon$ to the next inequality.
  Lemma~\ref{L:positive-physical-gram} gives
  \[
    2\operatorname{Rate}_R(\Phi)-\|\Phi\|^2/2
    \le\mathfrak{F}_{C^R}(\Phi)+2e_n+2\varepsilon .
  \]

  \medskip
  \noindent\textit{Step 3.~}We compare optimized Gram values and return to an
  integrated rate.
  Optimize this finite objective and use entrywise Gram comparison:
  \[
    \mathfrak{F}_{C^R}(\Phi)
    \le\mathscr{S}(C^R)\le\mathscr{S}(D)=\mathfrak{F}_D(\Phi_D).
  \]
  This compares optimized values, not the same individual control.
  Lemma~\ref{L:positive-uniform-optimizer} proves $\|\Phi_D\| \le M$ before time
  averaging is used again. Choose $R'$ covering the finite supports as in that
  lemma. The label row, Rayleigh bound and the same uniform averaging estimate
  give
  \begin{equation}\label{E:positive-finite-optimizer-bound}
    \mathfrak{F}_D(\Phi_D)
    \le2\operatorname{Rate}_{R'}(\Phi_D)-\|\Phi_D\|^2/2+2e_n .
  \end{equation}

  \medskip
  \noindent\textit{Step 4.~}We pass through nonnegative smooth controls and
  remove the two finite errors.
  If $\Phi_D = 0$, its finite variational expression is zero. Otherwise use its
  positive smooth approximants $\Psi_j$. Lemma~\ref{L:positive-raw-embedding}
  identifies their mean with $B_E(m_{R'}h_{\Psi_j,s})$ and their Gaussian cost
  with $v_{\Psi_j}$. Since $h_{\Psi_j,s} \ge 0$ and $0 \le m_{R'} \le 1$, real
  $H^1$ form order gives $\operatorname{Rate}_{R'}(\Psi_j) \le I_{\Psi_j}$.
  Thus, $2\operatorname{Rate}_{R'}(\Psi_j)-\|\Psi_j\|^2/2 \le
  S_2^{\mathrm{pos}}$. Rate and cost convergence pass this inequality to
  $\Phi_D$. Together with~\eqref{E:positive-finite-optimizer-bound} and the
  separate zero branch this yields
  \[
    2\operatorname{Rate}_R(\Phi)-\|\Phi\|^2/2
    \le S_2^{\mathrm{pos}}+4e_n+2\varepsilon .
  \]
  All finite choices were made after fixing $n,R,\Phi$. The right side is
  independent of those choices. Let $\varepsilon\downarrow0$, take the supremum
  over the original ball, use Lemma~\ref{L:positive-uniform-optimizer}, and then
  let $n \to \infty$. This proves~\eqref{E:positive-signed-comparison}.
  Finiteness of the variational value over nonnegative controls was already
  proved from the independent canonical upper bound in
  Proposition~\ref{P:variational-positive-scores}, not from this conclusion.
  This proves Proposition~\ref{P:positive-signed-comparison}.
\end{proof}

Proposition~\ref{P:positive-signed-comparison} bounds $\mathcal{S}_R$ above for
every $R$. For the limit we still need a lower bound as the mask is removed.

\begin{proposition}[Cutoff removal for the optimized variational expression]
  \label{P:positive-score-cutoff}
  The possibly nonmonotone mask family satisfies
  \[
    \lim_{R\to\infty}\mathcal{S}_R=S_2^{\mathrm{pos}}.
  \]
\end{proposition}

\begin{proof}
  We fix one smooth positive control and use compactly supported Rayleigh
  tests to remove the mask from its rate. We then approximate the already
  finite supremum by one such fixed control; the uniform signed comparison
  supplies the opposite bound.

  Fix $\Psi \in \mathcal{C}_{\mathrm{pos}}$ before taking $R$ or moment limits.
  Use the explicit mean version of Lemma~\ref{L:positive-raw-embedding}. Its
  rate $\lambda_R(s)$ obeys $0 \le \lambda_R(s) \le \lambda_{\Psi,s}$ by real
  form order. For each real compact smooth unit $f$, $m_R = 1$ on its support
  for every sufficiently large $R$. The Rayleigh formula therefore gives
  \[
    \liminf_{R\to\infty}\lambda_R(s)^2
    \ge\int h_{\Psi,s}f^2-\|\nabla f\|_2^2.
  \]
  Compact trial density and~\eqref{E:variational-real-rayleigh} identify the
  supremum of the right side as $\lambda_{\Psi,s}^2$. Hence, $\lambda_R(s) \to
  \lambda_{\Psi,s}$ for these versions. The global smooth-source bound $H_\Psi =
  \sup_{s,x}h_{\Psi,s}(x) < \infty$ gives $\lambda_R(s) \le \sqrt{H_\Psi}$.
  Dominated convergence yields $\operatorname{Rate}_R(\Psi) \to I_\Psi$. Neither
  monotonicity of $m_R$ nor an unmasked ground eigenfunction was used.

  The original representative and cost $\|\Psi\|^2 = v_\Psi$ are independent of
  $R$. Thus, $\liminf_R\mathcal{S}_R \ge 2I_\Psi-v_\Psi/2$ for each fixed test.
  Given $\varepsilon > 0$, choose one fixed test with variational expression
  greater than $S_2^{\mathrm{pos}}-\varepsilon$; the finite supremum admits such
  a test. Take the liminf and then $\varepsilon\downarrow0$.
  Proposition~\ref{P:positive-signed-comparison} supplies the reverse limsup.
  This proves Proposition~\ref{P:positive-score-cutoff}.
\end{proof}

This passage from cutoff scores to their suprema needs the comparison
just proved. No-go~\ref{P:nogo-30} shows why pointwise convergence and
bounded control costs alone would not justify it.

\begin{theorem}[Intrinsic original second-moment constant]
  \label{T:intrinsic-second-moment}
  For every original admissible model with data $(1,0)$ and fixed $x$,
  \begin{equation}\label{E:intrinsic-second-moment-limit}
    \lim_{T\to\infty}T^{-(4-\alpha-a)/(3-\alpha)}
    \log\E|u_\Gamma(T,x)|^2
    =C_2(\Gamma,a)=S_2^{\mathrm{pos}}
    =2^{(4-\alpha)/(3-\alpha)}S_1^{\mathrm{pos}}.
  \end{equation}
  Its complete variational formula over nonnegative controls is
  \begin{align}
    C_2(\Gamma,a)
    &=\sup_{\Psi\in\mathcal{C}_{\mathrm{pos}}}
    \left\{2\int_I
    \sqrt{\sup_{\substack{f\in H^1(\R^d;\R)\\\|f\|_2=1}}
    \left(\int h_{\Psi,s}f^2-\|\nabla f\|_2^2\right)}
    \ud s-\frac{v_\Psi}{2}\right\},
    \label{E:intrinsic-positive-formula}\\
    h_{\Psi,s}(x)
    &=\int_I|s-t|^{-a}(\Gamma*\Psi_t)(x)\ud t,\qquad
    v_\Psi=\int_{I\times\R^d}\Psi(s,x)h_{\Psi,s}(x)\ud x\ud s.
    \notag
  \end{align}
  Here $\mathcal{C}_{\mathrm{pos}}$ consists of all nonzero nonnegative real
  $C_c^\infty(I\times\R^d)$ tests, including nonseparable ones. The covariance
  retains its amplitude, the kinetic term has coefficient one, and the Gaussian
  action is $v_\Psi/2$. For nonzero covariance $0 < C_2 < \infty$. For $\Gamma =
  0$ the constant and moment rate are zero.
\end{theorem}

\begin{proof}
  We compare the upper and lower rates for the same canonical solution, then
  identify the variational expression and its homogeneity. The zero-covariance
  case is handled directly.

  \medskip
  \noindent\textit{Step 1.~}We match the two finite bounds for nonzero
  covariance.
  Assume first $\Gamma \ne 0$. Proposition~\ref{P:original-signed-point-upper},
  followed by Proposition~\ref{P:positive-score-cutoff} with $R =
  \delta^{-1}+2$, bounds the original point limsup by $S_2^{\mathrm{pos}}$.
  The radius $\delta^{-1}+2$ tends to infinity, so cutoff-score convergence
  applies along this radius without a separate exchange of limits.
  Proposition~\ref{P:variational-positive-scores} bounds the original point
  liminf below by this same value. Both are for the same canonical Stratonovich
  solution, exponent, original covariance and variance-half action. Their
  identifications, made with the regularization removed first, retain all mixed,
  unequal, self and empty chains and both phase components with initial data
  $(1,0)$. The
  original spatial stationarity proved in
  Proposition~\ref{P:raw-block-original-lower} extends both bounds to every
  fixed $x$. Thus, the limit exists and equals $S_2^{\mathrm{pos}}$. Its
  homogeneity, positivity and finiteness are
  Proposition~\ref{P:variational-positive-scores}.
  Formula~\eqref{E:variational-real-rayleigh} and the definition of the original
  variance give~\eqref{E:intrinsic-positive-formula}.

  The factor relating $S_2^{\mathrm{pos}}$ to $S_1^{\mathrm{pos}}$ is the
  control-dilation factor, not an identification with a first-moment constant.

  \medskip
  \noindent\textit{Step 2.~}We verify zero covariance directly.
  If $\Gamma = 0$, $H_I = \{0\}$, all means, variances and rates vanish,
  $\mathcal{S}_R = S_c^{\mathrm{pos}} = 0$, and $u = 1$. Both formulas are
  direct zero identities. This proves Theorem~\ref{T:intrinsic-second-moment}.
\end{proof}

The variational value with coefficient one is auxiliary; the theorem does not
identify it with a first-moment limit or with Chen's printed constants. It
neither repairs a source energy normalization nor asserts a Jensen improvement.
The exponent formally tends to $(4-\alpha)/(3-\alpha)$ as $a\downarrow0$; no
interchange with the static moment limit is claimed. The explicit bounds and
full calibration are separate from this identification.

\section{Occupation bounds for the intrinsic constant}
\label{S:intrinsic-occupation-ceiling}

Retain all original models and amplitudes, put $a = \alpha_0$, $q = 4-\alpha >
2$, $k = 3-\alpha = q-1 > 1$ and $m = q/k$. The temporal kernel is $K(s,t) =
|s-t|^{-a}$ on $I = (0,1)$; its diagonal value may be set to zero. Let
$\mathcal{O}$ be the full strongly $H^1$-measurable unit-profile class defined
before Proposition~\ref{P:occupation-physical-class}. Here we use the
\emph{original} kernel, defining
\begin{equation}\label{E:raw-occupation-definitions}
  \begin{aligned}
    A(g)&=\int_{I^2}K(s,t)B_\Gamma(g_s^2,g_t^2)\ud s\ud t,\\
    C_{\rm occ}&=\sup_{g\in\mathcal{O}}
    \frac{A(g)}{T(g)^{\alpha/2}},&
    \mathcal{E}_\zeta^{\mathcal{O}}
    &=\sup_{g\in\mathcal{O}}\{A(g)-\zeta T(g)\},\qquad \zeta>0.
  \end{aligned}
\end{equation}
The scalar penalty $\zeta$ is not a spatial mask. These objects differ from the
earlier damped-kernel energy. Assume $\Gamma \ne 0$ until the direct zero case
below.

\LeanMarker{L:occupation-raw-square}
\begin{lemma}[Full-profile squares in the original space]
  \label{L:occupation-raw-square}
  There is a continuous map $Q:H^1(\R^d;\R) \to H_{\mathrm{sp}}$ with
  \[
    \langle Q(f),Q(h)\rangle_{H_{\mathrm{sp}}}
    =B_\Gamma(f^2,h^2),\qquad
    \|Q(f)\|^2=A_{\mathrm{sp}}(f).
  \]
  For $g \in \mathcal{O}$, the element $\rho_g = J_{H_{\mathrm{sp}}}(s\mapsto
  Q(g_s))$ belongs to the original $H_I$ and satisfies
  \begin{equation}\label{E:occupation-raw-identities}
    \|\rho_g\|_{H_I}^2=A(g),\qquad
    \langle\Psi,\rho_g\rangle_{H_I}
    =\int_I\int_{\R^d}h_{\Psi,s}(x)g_s(x)^2\ud x\ud s
  \end{equation}
  for every $\Psi \in \mathcal{C}_{\mathrm{pos}}$. In particular,
  \[
    \langle\Psi,\rho_g\rangle_{H_I}\le\sqrt{v_\Psi A(g)}.
  \]
\end{lemma}

\begin{proof}
  We complete the square map in the spatial covariance Hilbert space. Its norm
  bound makes a measurable occupation an $L^2$ covariance-valued field, so the
  temporal inclusion applies. We then identify the physical pairing with a
  smooth control and use Hilbert Cauchy--Schwarz.

  \medskip
  \noindent\textit{Step 1.~}We define the continuous square map by Sobolev
  approximation.
  Approximate $f$ in $H^1$ by real compact smooth $f_n$.
  Proposition~\ref{P:spatial-four-linear} gives
  \[
    \|f_n^2-f_l^2\|_{H_{\mathrm{sp}}}^2
    =C_4(f_n-f_l,f_n+f_l,f_n-f_l,f_n+f_l)
    \le R_{\mathrm{sp}}N(f_n-f_l)^2N(f_n+f_l)^2\to0.
  \]
  Here $N$ is continuous on $H^1$ and $N(f) \le \|f\|_{H^1}$. Completeness
  defines $Q(f)$; interlacing two approximating sequences proves independence.
  Passing the physical four-linear continuity and Hilbert pairings to limits
  proves the first identities. The same difference bound, first for smooth
  approximants and then by passage to the limit, proves continuity of $Q$.

  \medskip
  \noindent\textit{Step 2.~}We construct the occupation vector in the original
  covariance space.
  Thus, $F_g(s) = Q(g_s)$ is strongly measurable. Writing $b_s = \|\nabla
  g_s\|_2^2$, the unit norm and $0 < \alpha/2 < 1$ give
  \[
    \int_I\|F_g(s)\|^2\ud s
    \le R_{\mathrm{sp}}\int_I b_s^{\alpha/2}\ud s
    \le R_{\mathrm{sp}}T(g)^{\alpha/2}<\infty .
  \]
  Lemma~\ref{L:positive-raw-embedding} now constructs $\rho_g$ and identifies
  its squared norm with the physical $A(g)$. Joint measurability and
  representative independence follow from the continuous physical pairing on
  $(H^1)^2$. Absolute time integrability follows from the same Schur bound.

  \medskip
  \noindent\textit{Step 3.~}We identify its pairing with the original smooth
  control.
  Choose one compact smooth $\chi$ equal to one on all spatial supports of
  $\Psi_t$. For real $f \in H^1$, smooth approximation gives
  \[
    \langle\Psi_t,Q(f)\rangle_{H_{\mathrm{sp}}}
    =\lim_n B_\Gamma(\Psi_t,f_n^2)
    =C_4(\chi,\Psi_t,f,f)=B_\Gamma(\Psi_t,f^2).
  \]
  Proposition~\ref{P:spatial-four-linear} justifies these absolute integrals and
  limits. Polarizing~\eqref{E:positive-bochner-inclusion} gives
  \[
    \langle\Psi,\rho_g\rangle
    =\int_{I^2}K(s,t)\langle\Psi_t,Q(g_s)\rangle_{H_{\mathrm{sp}}}
    \ud t\ud s .
  \]
  The absolute integrand is bounded by the product of two $H_{\mathrm{sp}}$
  norms against $K$, integrable by Schur and $L^2$ Cauchy--Schwarz. Evenness of
  $\Gamma$ and Tonelli for the nonnegative physical integrands
  yield~\eqref{E:occupation-raw-identities}. Hilbert Cauchy--Schwarz with
  $\|\Psi\|^2 = v_\Psi$ proves the bound. No rough occupation is inserted into a
  stochastic flow. This proves Lemma~\ref{L:occupation-raw-square}.
\end{proof}

\LeanMarker{P:raw-occupation-energy}
\begin{proposition}[Original occupation energy and dilation]
  \label{P:raw-occupation-energy}
  With $K_0 = 2^a/(1-a)$,
  \[
    0<C_{\rm occ}\le K_0R_{\mathrm{sp}}<\infty .
  \]
  Moreover,
  \begin{align}
    \mathcal{E}_\zeta^{\mathcal{O}}
    &=D_\zeta C_{\rm occ}^{2/(2-\alpha)},&
    D_\zeta&=\frac{2-\alpha}{2}
    \left(\frac{\alpha}{2\zeta}\right)^{\alpha/(2-\alpha)},
    \label{E:raw-occupation-energy}\\
    C_{\rm occ}
    &=\left(\frac{2\mathcal{E}_\zeta^{\mathcal{O}}}{2-\alpha}\right)^{(2-\alpha)/2}
    \left(\frac{2\zeta}{\alpha}\right)^{\alpha/2}.
    \label{E:raw-occupation-ratio}
  \end{align}
  In particular, $\mathcal{E}_\zeta^{\mathcal{O}} > 0$ is finite and equals
  $\zeta^{-\alpha/(2-\alpha)}\mathcal{E}_1^{\mathcal{O}}$.
\end{proposition}

\begin{proof}
  We bound the occupation ratio with the original temporal row bound and show
  it is positive. We then dilate the whole occupation class and optimize the
  resulting scalar energy to recover the exact penalty dependence.

  \medskip
  \noindent\textit{Step 1.~}We prove a positive finite occupation ratio.
  A unit $H^1$ profile has positive gradient energy: a zero weak gradient makes
  it constant, and an $L^2$ constant on $\R^d$ is zero. Hence, $0 < T(g) <
  \infty$. Spatial covariance Cauchy--Schwarz gives $B_\Gamma(g_s^2,g_t^2) \le
  R_{\mathrm{sp}}b_s^{\alpha/4}b_t^{\alpha/4}$. Schur with row bound $K_0$, then
  concavity, proves
  \begin{equation}\label{E:raw-occupation-bound}
    0\le A(g)\le K_0R_{\mathrm{sp}}T(g)^{\alpha/2}.
  \end{equation}
  A fixed compact smooth unit profile, constant in time, belongs to
  $\mathcal{O}$ and has positive $A$ by Proposition~\ref{P:spatial-four-linear}.
  This proves positivity and finiteness of $C_{\rm occ}$.

  \medskip
  \noindent\textit{Step 2.~}We optimize spatial dilation over the same
  occupation class.
  The map $(U_\ell g)_s(x) = \ell^{d/2}g_s(\ell x)$ is an $H^1$ isomorphism
  preserving unit mass and strong measurability. It has inverse $U_{1/\ell}$ and
  is a bijection on $\mathcal{O}$. The physical Radon change of variables gives
  \[
    A(U_\ell g)=\ell^\alpha A(g),\qquad
    T(U_\ell g)=\ell^2T(g).
  \]
  Indeed, the two densities contribute $\ell^{2d}$, $\ud x$ contributes
  $\ell^{-d}$, and the pushforward of $\Gamma$ contributes $\ell^{\alpha-d}$.
  Approximation extends the identity to full $H^1$. Taking the union over
  $g,\ell$ leaves the class unchanged. For fixed $A,T > 0$, the derivative of
  $A\ell^\alpha-\zeta T\ell^2$ has the sign of $\alpha A-2\zeta
  T\ell^{2-\alpha}$. Its global maximum occurs at $\ell = (\alpha A/(2\zeta
  T))^{1/(2-\alpha)}$ and equals $D_\zeta[A/T^{\alpha/2}]^{2/(2-\alpha)}$. For
  $A = 0$ the supremum is zero as $\ell\downarrow0$. Bounding every ratio by
  $C_{\rm occ}$ gives one inequality in~\eqref{E:raw-occupation-energy}. Fixed
  ratios approaching its positive supremum give the other by continuity of the
  positive power. Algebra gives~\eqref{E:raw-occupation-ratio} and the penalty
  scaling. No optimizer or minimax interchange is assumed. This proves
  Proposition~\ref{P:raw-occupation-energy}.
\end{proof}

\LeanMarker{L:occupation-rayleigh-selection}
\begin{lemma}[Countable Rayleigh near-occupations]
  \label{L:occupation-rayleigh-selection}
  For fixed $\Psi \in \mathcal{C}_{\mathrm{pos}}$ and every $\varepsilon > 0$
  there is $g_\varepsilon \in \mathcal{O}$ such that
  \begin{equation}\label{E:occupation-rayleigh-selection}
    I_\Psi^2+T(g_\varepsilon)
    \le\sqrt{v_\Psi A(g_\varepsilon)}+\varepsilon .
  \end{equation}
\end{lemma}

\begin{proof}
  We choose near-maximizing Rayleigh tests from one countable dense family. A
  first-index rule makes the chosen occupation strongly measurable; the
  bounded potential controls its kinetic energy and allows the physical
  covariance pairing estimate.

  \medskip
  \noindent\textit{Step 1.~}We choose a measurable near-maximizing occupation.
  The real $H^1$ space is separable, for example by its isometric embedding
  $f\mapsto(f,\nabla f)$ into a finite product of $L^2$ spaces. Compact smooth
  density therefore supplies a countable real family. Discard zero and normalize
  each member in $L^2$ to obtain $f_j$ dense in the $L^2$-unit $H^1$ class.
  Normalization preserves convergence because the approximating $L^2$ norms tend
  to one. Put $b_j = \|\nabla f_j\|_2^2$ and $R_j(s) = \int
  h_{\Psi,s}f_j^2-b_j$. Bounded-potential form continuity and
  Proposition~\ref{P:variational-unmasked-rate} give $\lambda_{\Psi,s}^2 =
  \sup_jR_j(s)$. The canonical functions are measurable, indeed continuous. Let
  $U_j = \{s:R_j(s)>\lambda_{\Psi,s}^2-\varepsilon\}$, $A_1 = U_1$ and $A_j =
  U_j\setminus\bigcup_{i<j}U_i$. These first-index cells are Borel and cover
  $I$. Put $g_\varepsilon(s) = f_j$ on $A_j$. Finite-valued truncations using
  $f_1$ on the remaining cells converge pointwise in $H^1$, proving strong
  measurability. A null-set representative may be filled by $f_1$.

  \medskip
  \noindent\textit{Step 2.~}We bound its kinetic energy and integrate the
  Rayleigh inequality.
  With $H_\Psi = \sup_{s,x}h_{\Psi,s}(x) < \infty$, the selection implies
  \[
    b(g_\varepsilon(s))
    <\int h_{\Psi,s}g_\varepsilon(s)^2-\lambda_{\Psi,s}^2+\varepsilon
    \le H_\Psi+\varepsilon .
  \]
  Thus, it belongs to the $\mathcal{O}$. Integrating the near-Rayleigh
  inequality, using Lemma~\ref{L:occupation-raw-square}, gives
  \[
    \int_I\lambda_{\Psi,s}^2\ud s+T(g_\varepsilon)
    \le\langle\Psi,\rho_{g_\varepsilon}\rangle+\varepsilon
    \le\sqrt{v_\Psi A(g_\varepsilon)}+\varepsilon .
  \]
  Time Cauchy--Schwarz on the unit interval gives $I_\Psi^2 \le
  \int_I\lambda_{\Psi,s}^2\ud s$. This
  proves~\eqref{E:occupation-rayleigh-selection}. No eigenvector, continuum
  selector or limiting occupation is required. This proves
  Lemma~\ref{L:occupation-rayleigh-selection}.
\end{proof}

For $0 \le b < q/2$ define $F_b(0) = 0$ and
\begin{equation}\label{E:occupation-scalar-map}
  F_b(z)=k\left[
  \frac{(q-2b)^{(q-2b)/2}}{8(q-b)^{q-b}}
  \left(\frac{2z}{2-\alpha}\right)^{(2-\alpha)/2}
  \right]^{1/k},\qquad z>0.
\end{equation}
At $b = 0$ this is exactly the map in~\eqref{E:spatial-static-energy-score}. The
display parameter $b$ does not change the temporal exponent $a$.

\LeanMarker{T:intrinsic-occupation-ceiling}
\begin{theorem}[The computed occupation ceiling]\label{T:intrinsic-occupation-ceiling}
  For every variational value over nonnegative controls coefficient
  $\mathfrak{c}$,
  \begin{align}
    S_{\mathfrak{c}}^{\mathrm{pos}}
    &\le U_{\mathfrak{c}}
    :=\mathfrak{c}^m k
    \left[\frac{\alpha^{\alpha/2}C_{\rm occ}}{8q^{q/2}}\right]^{1/k}
    \label{E:occupation-score-ceiling}\\
    &=\mathfrak{c}^m(2\zeta)^{\alpha/(2k)}
    F_0(\mathcal{E}_\zeta^{\mathcal{O}}),\qquad \zeta>0.\notag
  \end{align}
  In particular, the original moment constant satisfies $C_2(\Gamma,a) \le U_2$
  for every original model. This ceiling is not asserted to be attained or equal
  to $C_2$.
\end{theorem}

\begin{proof}
  We first optimize the kinetic budget in the near-occupation inequality to
  bound the integrated rate by a power of the original variance. We then
  optimize that variance in the reward-minus-cost expression and rewrite the
  result using the occupation energy.

  \medskip
  \noindent\textit{Step 1.~}We optimize the kinetic budget before removing the
  selection error.
  Fix $\Psi$ and write $v = v_\Psi$, $C = C_{\rm occ}$ and $\beta = \alpha/4$.
  Take $T=T(g_\varepsilon)>0$ from the selected occupation; by the definition
  of $C_{\rm occ}$, $A(g_\varepsilon)\le C_{\rm occ}T^{\alpha/2}$, which
  converts its covariance term into $\sqrt{vC}\,T^\beta$.
  At the positive critical point $T_*$ of the scalar maximization below,
  $T_*=\beta\sqrt{vC}\,T_*^\beta$, so the value is
  $(1-\beta)\sqrt{vC}\,T_*^\beta$.
  Lemma~\ref{L:occupation-rayleigh-selection} gives $I_\Psi^2 \le
  \sqrt{vC}T^\beta-T+\varepsilon$ for some $T > 0$. Bound it by the supremum
  over all $T \ge 0$ before sending $\varepsilon\downarrow0$. If $v = 0$, that
  supremum is zero. If $v > 0$, the derivative $\beta\sqrt{vC}T^{\beta-1}-1$
  strictly decreases from infinity to $-1$. Its unique positive zero is $T_* =
  (\beta\sqrt{vC})^{1/(1-\beta)}$. The global maximum equals
  $(q/4)(\alpha/4)^{\alpha/q}(vC)^{2/q}$. In both cases,
  \[
    I_\Psi\le Dv^{1/q},\qquad
    D=\sqrt{q/4}(\alpha/4)^{\alpha/(2q)}C^{1/q}.
  \]

  \medskip
  \noindent\textit{Step 2.~}We optimize the variance and keep the exact
  Gaussian half-cost.
  For the constant computation below, $4^{-\alpha/2}=2^{-\alpha}$ and
  absorbing the outside factor $1/2$ into the $1/k$ power contributes
  $2^{-k}$, whose product is $2^{-3}$ because $k+\alpha=3$.
  The scalar function $\mathfrak{c} Dv^{1/q}-v/2$ vanishes at zero, tends to
  minus infinity, and has strictly decreasing derivative. Its unique positive
  maximum is at $v_* = (2\mathfrak{c} D/q)^{q/k}$, with value $kv_*/2$. Direct
  positive-power algebra gives
  \[
    (2D/q)^q=q^{-q/2}(\alpha/4)^{\alpha/2}C .
  \]
  Since $k+\alpha = 3$, the maximum is precisely $U_{\mathfrak{c}}$
  in~\eqref{E:occupation-score-ceiling}, with denominator eight. Substitution
  of~\eqref{E:raw-occupation-ratio} gives the second form, retaining the factor
  $(2\zeta)^{\alpha/(2k)}$. Take the supremum over the smooth positive-control
  class. Theorem~\ref{T:intrinsic-second-moment} identifies $C_2 =
  S_2^{\mathrm{pos}}$, proving the moment ceiling. The penalty factor cancels
  its energy scaling, so the ceiling is independent of $\zeta$. This proves
  Theorem~\ref{T:intrinsic-occupation-ceiling}.
\end{proof}

\LeanMarker{C:intrinsic-spatial-bracket}
\begin{corollary}[A common spatial-energy bracket]
  \label{C:intrinsic-spatial-bracket}
  With the same original covariance and $M_a$
  of~\eqref{E:temporal-equilibrium-mass},
  \begin{equation}\label{E:intrinsic-spatial-bracket}
    2^m M_a^{-1/k}F_0(\mathcal{E}_{\mathrm{sp}})
    \le C_2(\Gamma,a)
    \le2^m K_0^{1/k}F_0(\mathcal{E}_{\mathrm{sp}}).
  \end{equation}
  For $\Gamma \ne 0$ there is also the comparison of the energies defined here
  \[
    C_2\le2^mF_0(\mathcal{E}_{1/2}^{\mathcal{O}})
    <2^mF_a(\mathcal{E}_{1/2}^{\mathcal{O}}).
  \]
\end{corollary}

\begin{proof}
  We combine the spatial ratio estimate with the occupation ceiling for the
  common-energy bracket. We then compare the scalar maps by strict convexity,
  keeping the kinetic penalty normalization explicit.

  Insert $C_{\rm occ} \le K_0R_{\mathrm{sp}}$ into
  \eqref{E:occupation-score-ceiling} and use \eqref{E:spatial-energy-map}. This
  gives the upper bracket; the lower is \eqref{E:spatial-static-energy-score}
  multiplied by $2^m$, using Theorem~\ref{T:intrinsic-second-moment}. Both
  concern one common spatial energy, not an assumed first moment.

  For completeness put
  \[
    \Delta=\tfrac12[q\log q+(q-2a)\log(q-2a)]
    -(q-a)\log(q-a)>0.
  \]
  Here $a > 0$ and $q-2a > 2-a > 1$. The scalar inequality for the difference of
  logarithms applies at exactly these conditions, without $q \ge 3$.
  Equivalently the strict midpoint convexity of $x\log x$ on positive arguments,
  whose second derivative is $1/x$, gives the sign. Subtracting logarithms
  in~\eqref{E:occupation-scalar-map} gives
  \[
    \log(F_a(z)/F_0(z))=\Delta/k,\qquad
    \frac{\mathfrak{c}^mF_a(\mathcal{E}_\zeta^{\mathcal{O}})}{U_{\mathfrak{c}}}
    =\exp\!\left(\frac{\Delta-(\alpha/2)\log(2\zeta)}{k}\right).
  \]
  At $\zeta = 1/2$ this and the ceiling prove the strict comparison. At $\zeta =
  1$ the comparison is not universally favorable. For $\alpha = 5/4,a = 1/2,d =
  2$,
  \[
    8[\Delta-(\alpha/2)\log2]
    =\log\frac{11^{11}7^7}{2^5 9^{18}}<0,
  \]
  since $11^{11}7^7 = 234966429149994773$ and $2^5 9^{18} =
  4803028329503971872$. This contrasts computed bounds; it refutes no moment
  inequality. This proves Corollary~\ref{C:intrinsic-spatial-bracket}.
\end{proof}

The kinetic normalization matters in this comparison.
No-go~\ref{P:nogo-28} shows that changing it can reverse an ordering
between the computed ceiling and a comparison value. That example makes
no claim about the ordering of the actual moment constant.

For $\Gamma = 0$, all original means and variational values vanish and $C_2 =
0$. Here $Q = 0$, $A = 0$ and $C_{\rm occ} = 0$. The supremum of $-\zeta T(g)$
is zero: dilate a fixed unit compact profile by $\ell^{d/2}f(\ell x)$ and let
$\ell\downarrow0$. The same argument gives $\mathcal{E}_{\mathrm{sp}} = 0$ and
the earlier damped energy zero, without using a zero vector as a unit profile.
With $F_b(0) = 0$, every non-strict bound above is the zero identity. These
energies are not identified with Chen's printed energies.

\section{Spatial white noise: energy and explicit moment bounds}
\label{S:intrinsic-white-evaluation}

Take $d = \alpha = 1$ and $\Gamma = c\delta_0$, $c \ge 0$, with temporal kernel
$|s-t|^{-a}$, $0 < a < 1$, and data $(1,0)$. The parameter $c$ is covariance
mass; the coefficient of unit noise would be $\sqrt c$. No additional temporal
normalization is inserted. The wave kinetic coefficient is one. For real unit $f
\in H^1(\R)$ put $A = \int_\R f^4$ and $b = \int_\R(f')^2$. The physical
covariance is $A_{\mathrm{sp}}(f) = cA$.

The next calculation recovers the sharp one-dimensional interpolation
inequality recorded by Sz.-Nagy~\cite{sz.nagy:source1941:integralungleichungen}*{(2),
p.~64, and (4), p.~65}, now commonly called a Gagliardo--Nirenberg inequality,
and its hyperbolic-secant optimizer. We give the proof to fix the normalization
of the spatial energy.

\begin{lemma}[A square identity for the distribution function]\label{L:white-cdf-square}
  Every such $f$ has $A < \infty$. With $\mathsf{F}(x) =
  \int_{-\infty}^xf(y)^2\ud y$, for every real $z$,
  \begin{equation}\label{E:white-cdf-square}
    0\le\int_\R[f'-zf(1-2\mathsf{F})]^2
    =b-2zA+z^2/3 .
  \end{equation}
\end{lemma}

\begin{proof}
  We first obtain a continuous representative with a bounded square and zero
  limits at infinity. We then integrate the derivative of its product with the
  distribution function and expand the square, with all terms absolutely
  integrable.

  \medskip
  \noindent\textit{Step 1.~}We identify the representative and its endpoint
  decay.
  On a bounded interval the weak derivative $f'$ is $L^1$ by Cauchy--Schwarz.
  Its integral primitive has weak derivative $f'$ by Fubini. Subtracting it from
  $f$ gives a function of weak derivative zero, hence constant almost
  everywhere: test against derivatives of compact smooth primitives of all
  zero-integral test functions. The resulting locally absolutely continuous
  representatives patch on overlaps, giving one continuous representative. Thus,
  $\rho = f^2$ is locally absolutely continuous with $\rho' = 2ff' \in L^1(\R)$,
  while $\rho \ge 0$ and $\int\rho = 1$. Integrability of $\rho'$ makes $\rho$
  Cauchy on each infinite tail, so it has finite limits at both ends. Its $L^1$
  integrability forces both limits to be zero. It is therefore bounded, and $A =
  \int\rho^2 \le \|\rho\|_\infty\int\rho < \infty$.

  \medskip
  \noindent\textit{Step 2.~}We compute the two integrals in the square
  expansion.
  Now $0 \le \mathsf{F} \le 1$, $\mathsf{F}' = \rho$ continuously, and its
  limits are zero and one. The derivative of $\rho(1-2\mathsf{F})$ is
  $\rho'(1-2\mathsf{F})-2\rho^2 \in L^1$. Its endpoint limits vanish.
  Since $\rho'=2ff'$, dividing the integrated derivative identity by two gives
  the cross term with a positive sign.
  Integrating on $[-R,R]$ and sending $R \to \infty$ gives
  \[
    \int_\R f'f(1-2\mathsf{F})=A .
  \]
  The polynomial $P(y) = y-2y^2+(4/3)y^3$ has derivative $(1-2y)^2$. The chain
  rule and the endpoints of $\mathsf{F}$ give
  \[
    \int_\R f^2(1-2\mathsf{F})^2=P(1)-P(0)=1/3 .
  \]
  Absolute integrability holds since this integrand is bounded by $f^2$. For
  each real $z$, both terms in the square in~\eqref{E:white-cdf-square} are in
  $L^2$, and their cross product is absolutely integrable. Expansion and the two
  preceding identities prove the formula. Neither positivity of $f$ nor strict
  increase of its CDF is required. This proves Lemma~\ref{L:white-cdf-square}.
\end{proof}

One choice of $z$ in Lemma~\ref{L:white-cdf-square} gives an inequality
between $A$ and $b$. We use it to evaluate the spatial ratio $R_{\mathrm{sp}}$
and the energy $\mathcal{E}_{\mathrm{sp}}$ of
Proposition~\ref{P:spatial-energy-map} for white noise.

\LeanMarker{P:white-spatial-energy}
\begin{proposition}[Exact white spatial energy]\label{P:white-spatial-energy}
  For $c > 0$,
  \[
    R_{\mathrm{sp}}=c/\sqrt3,\qquad
    \mathcal{E}_{\mathrm{sp}}=c^2/6.
  \]
  Both values are zero at $c = 0$.
\end{proposition}

\begin{proof}
  We minimize the CDF square in its scalar parameter to bound the ratio. We
  then evaluate the hyperbolic-secant profiles, which attain both the ratio
  and the completed-square energy bound at the stated scale.

  \medskip
  \noindent\textit{Step 1.~}We obtain the sharp ratio upper bound.
  Substitution $z = 3A$ into~\eqref{E:white-cdf-square} gives $3A^2 \le b$. This
  is the exact premise and conclusion of the existing conditional WhiteCDF
  square deduction. Here $A > 0$, since $A = 0$ would imply $f = 0$ a.e. Also $b
  > 0$: zero derivative makes its continuous representative constant and hence
  zero in $L^2$, contradicting unit mass. Thus, $A/\sqrt b \le 1/\sqrt3$.

  \medskip
  \noindent\textit{Step 2.~}We evaluate an attaining family and the original
  covariance energy.
  For $\ell > 0$ use the real smooth $H^1$ profile $f_\ell(x) =
  \sqrt{\ell/2}\operatorname{sech}(\ell x)$, with $f_\ell' = -\ell\tanh(\ell
  x)f_\ell$. Set $t = \tanh(\ell x)$, so $\ud t = \ell\operatorname{sech}^2(\ell
  x)\ud x$. Direct integration gives
  \[
    \int f_\ell^2=\frac12\int_{-1}^1\ud t=1,\quad
    \int f_\ell^4=\frac\ell4\int_{-1}^1(1-t^2)\ud t=\ell/3,\quad
    \int(f_\ell')^2=\frac{\ell^2}2\int_{-1}^1t^2\ud t=\ell^2/3.
  \]
  The ratio is attained, proving $R_{\mathrm{sp}} = c/\sqrt3$. The physical
  Dirac measure introduces no Fourier constant or nonintegrable kernel
  $|x|^{-1}$. The existing conditional energy deduction applies to $3A^2 \le b$:
  \[
    cA-b/2\le cA-\tfrac32A^2
    =c^2/6-\tfrac32(A-c/3)^2\le c^2/6 .
  \]
  The coefficient is covariance mass $c$ in the
  physical pairing, not a noise coefficient of size $c$.
  At $\ell = c > 0$ the profile attains $c^2/6$. These are full $H^1$ energy
  profiles; compact approximation equality is
  Proposition~\ref{P:spatial-four-linear}, not a claim that the noncompact
  profile itself is a compact stochastic control. At $c = 0$ the ratio is zero
  and the energy is $\sup(-b/2) = 0$, approached by $f_\ell$ as
  $\ell\downarrow0$. No zero-scale profile is used as a unit vector. The CDF,
  integral and attainment arguments remain analytic; the two finite conditional
  deductions do not formalize them. This proves
  Proposition~\ref{P:white-spatial-energy}.
\end{proof}

Corollary~\ref{C:intrinsic-spatial-bracket} bounds $C_2$ above and below
through the spatial energy $\mathcal{E}_{\mathrm{sp}}$. With its white value,
both bounds become explicit.

\LeanMarker{T:white-moment-bounds}
\begin{theorem}[Explicit white moment bounds]\label{T:white-moment-bounds}
  For the same original point limit,
  \begin{equation}\label{E:white-moment-bounds}
    \frac23\sqrt{\frac{c}{M_a}}\le C_2(c,a)
    \le\frac23\sqrt{cK_0(a)},\qquad
    K_0(a)=\frac{2^a}{1-a},\quad
    M_a=\frac{\Gamma_E((1+a)/2)}
    {\Gamma_E(1+a)\Gamma_E((1-a)/2)} .
  \end{equation}
  Its time exponent is $(3-a)/2$. These are bounds, not a numerical evaluation
  of the intrinsic supremum.
\end{theorem}

\begin{proof}
  Here $q = 3,k = 2,m = 3/2$. Formula~\eqref{E:spatial-static-energy-score}
  gives
  \[
    F_0(z)=2\left[\frac{(2z)^{1/2}}{8\,3^{3/2}}\right]^{1/2}.
  \]
  For $c > 0$, $\sqrt{2\mathcal{E}_{\mathrm{sp}}} = c/\sqrt3$, hence
  $F_0(\mathcal{E}_{\mathrm{sp}}) = \sqrt c/(3\sqrt2)$ and
  $2^mF_0(\mathcal{E}_{\mathrm{sp}}) = (2/3)\sqrt c$. Substitution into
  Corollary~\ref{C:intrinsic-spatial-bracket}
  proves~\eqref{E:white-moment-bounds} in the original units. All Gamma factors
  are positive finite for $0 < a < 1$. At $c = 0$, $u = 1$ and $C_2 = 0$ by
  Theorem~\ref{T:intrinsic-second-moment}; both displayed bounds vanish. This
  proves Theorem~\ref{T:white-moment-bounds}.
\end{proof}

For $c>0$, both formulas in~\eqref{E:white-moment-bounds} are multiples of
$\sqrt c$. Their ratio depends only on $a$.

\LeanMarker{C:white-relative-width}
\begin{corollary}[Uniform relative width of the bound formulas]
  \label{C:white-relative-width}
  For $c > 0$, the ratio of the upper and lower formulas is
  \[
    \sqrt{K_0(a)M_a}
    \le\left[\frac{\sqrt\pi}{2\Gamma_E(5/4)^2}\right]^{1/2},
  \]
  with equality at $a = 1/2$.
\end{corollary}

\begin{proof}
  Duplication and recurrence from Lemma~\ref{L:equilibrium-gamma-beta} give
  \[
    \Gamma_E(1+a)
    =\frac{2^a\Gamma_E((1+a)/2)\Gamma_E(1+a/2)}{\sqrt\pi},
    \qquad
    (1-a)\Gamma_E((1-a)/2)=2\Gamma_E((3-a)/2).
  \]
  Therefore,
  \[
    K_0(a)M_a
    =\frac{\sqrt\pi}{2\Gamma_E(1+a/2)\Gamma_E((3-a)/2)} .
  \]
  Put $x = 1+a/2$, $y = (3-a)/2$, so $x,y > 0$ and $x+y = 5/2$. In Euler's Gamma
  integral the $L^2(0,\infty)$ functions $t^{(x-1)/2}e^{-t/2}$ and
  $t^{(y-1)/2}e^{-t/2}$ have squared norms $\Gamma_E(x),\Gamma_E(y)$ and inner
  product $\Gamma_E(5/4)$. Cauchy--Schwarz, positivity, inversion and the
  increasing square root prove the bound. At $a = 1/2$, $x = y = 5/4$, giving
  equality. Replacing $a$ by $1-a$ exchanges these two arguments, so only this
  ratio of bound formulas is symmetric, not the moment constant or exponent.
  This proves Corollary~\ref{C:white-relative-width}.
\end{proof}

This example retains the all-model theorem rather than replacing its scope. The
formal static exponent is $3/2$ as $a\downarrow0$; no interchange of long-time
and temporal-endpoint limits is asserted. Jensen holds for the same canonical
solution, while $S_1^{\mathrm{pos}}$ remains auxiliary. The read-only Skorohod
comparison at strict spatial Riesz exponent below one does not cover this
spatial-white endpoint. A separate white-time problem cannot be substituted for
$0 < a < 1$. No new Skorohod result, source correction or certified decimal
enclosure is claimed.

\section{Established steps and open questions}\label{S:route}

\begin{enumerate}[1.]
  \item \textbf{The mixed expansion and the first obstruction.}
    Chapter~\ref{S:G0} proves the complete nonnegative mixed expansion and an
    independent convergent solution series. It corrects the time-uniformity
    issue in the source's intermediate estimate. The common-per-pairing majorant
    has exponent $\rho > \kappa$ and cannot give the desired scale. This is a
    restriction on that majorant, not on the model or all uses of the source's
    Theorem~3.1.
  \item \textbf{Calibrations.}
    Section~\ref{S:G1} gives the finite white and Riesz coefficients. As
    $\alpha_0\downarrow0$, the target exponent becomes $(4-\alpha)/(3-\alpha)$,
    the static exponent. Proposition~\ref{P:Jensen} is unconditional; its
    numerical $2L_*$ specialization requires the stated first-moment premise.
    Comparisons with the Skorohod model identify time exponents only.
  \item \textbf{The full upper scale.}
    Section~\ref{S:d1-upper} gives explicit bounds when $\alpha \le 1$.
    Section~\ref{S:interval-upper} proves the full homogeneous-covariance upper
    scale by independent interval fields, local Gaussian-degree estimates and
    dyadic covariance domination. The latter proof does not assume an endpoint
    representation of the spatial covariance.
  \item \textbf{The intrinsic constant.}
    Sections~\ref{S:actual-operator-driver}--\ref{S:intrinsic-second-moment}
    construct the phase/source framework and prove matching original
    point-moment bounds. The upper variational expression uses the original
    covariance norm and its Gaussian quadratic cost. Finite Gram optimization
    and positive smoothing identify it with the full positive-control lower
    value. The auxiliary coefficient-one variational quantity is not identified
    with the first-moment limit of the solution or Chen's printed constants.
  \item \textbf{Computed bounds and the white example.}
    Section~\ref{S:intrinsic-occupation-ceiling} bounds the same constant using
    the full occupation class and a common spatial energy.
    Section~\ref{S:intrinsic-white-evaluation} evaluates that energy for $\Gamma
    = c\delta_0$ and gives explicit moment bounds and the relative width of
    their formulas. It does not evaluate the intrinsic supremum numerically or
    assert equality at either bound.
  \item \textbf{Source comparisons and scope.}
    The source first-moment normalization remains a qualified comparison
    question. Chapter~\ref{S:book-no-go} collects the failed approaches and
    explains their precise scope and the constructions that replace them.
    They are not new premises of the intrinsic-limit proof.
\end{enumerate}

\chapter{Obstructions and failed approaches}\label{S:book-no-go}
\label{S:appendix-no-go} % Retained for existing links to the former appendix.

% Stable NG-01--NG-33 records, reorganized by mathematical cause.
% Historical review grades are preserved in the provenance record.
\section{How to read the negative results}
\label{S:nogo-guide}

The proof of the main results did not follow one route from beginning to end.
Several plausible steps lost the required time scale, changed the covariance
class, discarded part of the wave state, or asserted uniformity that was not
available. This chapter explains those failures as mathematical results.
Each entry starts with the reason the approach was attractive, states a
precise obstruction, and supplies its argument. It then separates what fails
from what survives and specifies what a replacement would need to establish.

\paragraph{Where to begin.}
For pairing arguments, start with the two-chain expansion in
Proposition~\ref{P:mixed}; for covariance and phase arguments, recall
Hypothesis~\ref{H:homog} and the displacement--velocity state. The topology
section states its own spaces, and the finite examples need only Gaussian
moments and elementary optimization. Each entry gives more specific
prerequisites and points to the part of the positive proof it explains.
The five groups below can be read independently.

The stable records NG-01 to NG-33 are printed below as No-go~01 to No-go~33.
Their numbers do not change when related examples are placed together. Some
restate results proved earlier in the book; others give the retained
counterexample in full. They are records of attempted approaches, not a count
of independent discoveries. No-go~\ref{P:nogo-33} groups three elementary
errors found in review arguments.

\paragraph{The scope of a counterexample.}
Unless an entry specifies a different model, the covariance and initial data
are those of Hypothesis~\ref{H:homog} and the canonical Stratonovich equation
with data $(1,0)$. We write $a=\alpha_0$. Scalar ordinary differential
equations, finite matrices and abstract Hilbert-space examples are identified
as such. They test the general principles stated in those entries; an example
outside Hypothesis~\ref{H:homog} does not refute a theorem restricted to that
hypothesis. A numerical majorant growing too quickly is an obstruction to its
use in a proof, not a lower bound on the actual moment.

\paragraph{The status of the evidence.}
The historical record distinguishes arguments already in the book, arguments
that received a review, and author-side diagnostics. These descriptions are
not interchangeable. A review may have occurred without meeting its required
acceptance threshold, as happened for the earlier NG-09 record. The proofs
below are the present mathematical account; the original reports and their
grades remain unchanged. Where an earlier computation was only formal, the
necessary analytic justification is given explicitly rather than attributed
to that earlier computation.

\begin{center}
\begin{tabular}{@{}>{\raggedright\arraybackslash}p{0.26\textwidth}>{\raggedright\arraybackslash}p{0.41\textwidth}>{\raggedright\arraybackslash}p{0.24\textwidth}@{}}
  \toprule
  Obstruction & What the book keeps explicit & Entries \\
  \midrule
  \hyperref[S:nogo-pairings]{Pairing and coefficient estimates} & Temporal order, spatial support and
  the quantifiers on constants & 01--07, 10, 17, 18 \\
  \hyperref[S:nogo-covariance]{Covariance comparisons} & The actual covariance class and the correct
  order of comparison & 08, 09, 11--13, 19 \\
  \hyperref[S:nogo-phase]{Phase and Gaussian controls} & Velocity, the second-order mean and the
  full tilted law & 14--16, 20, 31 \\
  \hyperref[S:nogo-topology]{Operator spaces and limits} & The specified topology, radius dependence
  and control of variational suprema & 21, 22, 24--27, 30, 32 \\
  \hyperref[S:nogo-optimization]{Finite optimization and review} & Both variances, active sets and exact
  algebraic factors & 23, 28, 29, 33 \\
  \bottomrule
\end{tabular}
\end{center}

The main text cites these entries where they explain a choice in the proof.
The entries also point back to the construction or estimate that survives.
This is intended to let a reader begin with an attempted approach, understand
its obstruction, and then follow the successful replacement.

Tao has observed that negative results are underreported, which can make the
public account of AI-assisted mathematics appear more successful than the
underlying research process~\cite{tao:26:post}*{first paragraph}. The purpose here is to make
these particular failures usable. Separate source-paper qualifications and
failures of the checking process retain their own records; they are not
silently added to the 33 mathematical identifiers below.

\begin{figure}[htbp]
  \centering
  % Vector metaphor only; the cards are the book's stable No-go identifiers.
\begin{tikzpicture}[x=1cm,y=1cm,font=\normalsize]
  \fill[hamblue!5] (-6.8,-4) rectangle (6.8,0);
  \filldraw[draw=hamblue!65,fill=hamblue!15,thick]
  (-1,0)--(-5.7,-.7)--(-6.1,-1.5)--(-6.2,-2.7)
  --(-4.8,-3.6)--(-1.2,-4)--(3.8,-3.75)
  --(5.8,-2.9)--(6.15,-1.5)--(5.7,-.7)--(1,0)--cycle;
  \filldraw[draw=hamblue,fill=white,thick]
  (-1,0)--(-.6,.85)--(-.05,1.85)--(.5,1.05)--(1,0)--cycle;
  \draw[hamblue,thick] (-4.4,0)--(6.8,0);
  \node[anchor=west,fill=white,inner sep=3pt] at (-6.6,.18) {waterline};
  \node[anchor=west,text=hamteal] at (2.15,1.35)
  {\ProofRef{Theorem}{T:intro-upper}};
  \node[anchor=west,text=hamteal] at (2.15,.65)
  {\ProofRef{Theorem}{T:intro-second-moment}};
  \draw[hamteal] (.55,1.05)--(1.5,1.35)--(2,1.35);
  \draw[hamteal] (.8,.45)--(1.5,.65)--(2,.65);
  \foreach \id/\xx in {01/-4.2,02/-1.4,33/4.2} {
  \node[draw=hamblue!70,fill=white,rounded corners=1.5pt,
  minimum width=1.9cm,minimum height=.55cm,inner sep=4pt]
  at (\xx,-1.95) {\hyperref[P:nogo-\id]{No-go~\id}};
  }
  \node at (1.4,-1.95) {$\cdots$};
  \node[font=\normalsize\bfseries,text=hamblue] at (0,-3.3)
  {Chapter 12: below the waterline};
\end{tikzpicture}
  \caption[What lies below the waterline]{The results stated in Chapter~1,
    Theorems~\ref{T:intro-upper} and~\ref{T:intro-second-moment}, are what shows
    above the waterline. This chapter is what lies below it: the 33 numbered
    obstructions and failed approaches, each kept with the reason it was tried and
    the argument that stops it.}
  \label{Fig:nogo-iceberg}
\end{figure}
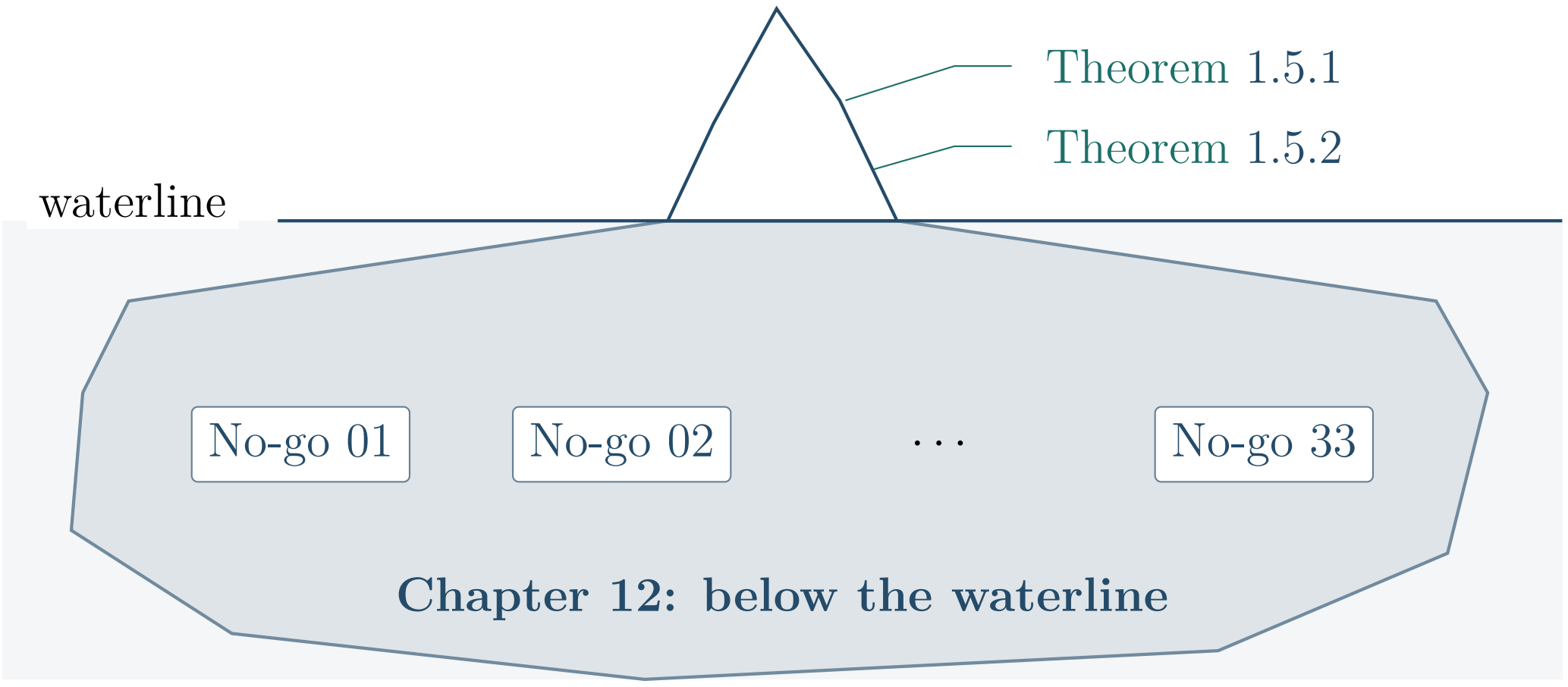

\section{Pairings, ordering, and comparison kernels}\label{S:nogo-pairings}

The second moment contains two ordered chains and all their Gaussian
pairings, including pairs within a chain. This structure offers several
natural reductions: estimate each pairing separately, compare with a single
chain, or copy the first-moment clock. The obstructions below distinguish
loss of scale from failure of an inequality. A majorant can be too large even
when every step used to obtain it is valid.

Throughout this section write $a=\alpha_0$ and
$\beta=4-\alpha-a$. Let $X_n$ denote the ordered Stratonovich term of degree
$n$, and put
\[
  a_k=\E X_{2k}(1,x),\qquad
  A_k=\sum_{m+n=2k}\E[X_m(1,x)X_n(1,x)],\qquad
  b_1=\E X_1(1,x)^2.
\]
Thus $A_1=2a_1+b_1$. Homogeneity supplies the factor $t^{\beta k}$.
The nonnegative convergent moment expansions used here are those of
Proposition~\ref{P:mixed}; they include unequal chain lengths. In dimension
one we use $G_s(z)=\frac12\mathbf1_{\{|z|<s\}}$.

\subsection{NG-01: one common estimate loses a factorial power}
\label{S:nogo-01}

Chen's author version~\cite{chen:author2026:time-dependency}*{Theorem~3.1, (3.6)}
bounds regularized two-chain integrals for arbitrary pairings under its
finite-measure and norm hypotheses. Proposition~\ref{P:envelope} gives the
application under Hypothesis~\ref{H:homog}. Applying one common norm bound
to every pairing is an attractive first attempt at the upper exponent.
The issue is what happens after optimization and summation.

\setcounter{nogo}{0}
\begin{nogo}[The common-parameter majorant has the wrong scale]
  \label{P:nogo-01}
  Under Hypothesis~\ref{H:homog}, let $P_k=(2k)!/(2^k k!)$ and $K>0$.
  The numerical majorant
  \[
    B_K(t)=1+\sum_{k\ge1}(2k+1)P_kK^k
    \inf_{\theta>0}e^{2\theta t}\theta^{-\beta k}
  \]
  satisfies $\log B_K(t)\asymp t^{\beta/(\beta-1)}$ as $t\to\infty$.
  Since $\beta/(\beta-1)>\kappa=\beta/(3-\alpha)$, this majorant cannot
  establish a finite upper constant on the scale $t^\kappa$.
\end{nogo}
\begin{proof}
  The infimum equals $e^{\beta k}(2t/(\beta k))^{\beta k}$.
  The inequalities $k!\le P_k\le2^k k!$ and
  $e^{-k}\le k!/k^k\le1$ show that $B_K$ lies between two series
  \[
    \sum_{k\ge0}\frac{(c t^\beta)^k}{(k!)^q},\qquad q=\beta-1>1,
  \]
  with different positive constants $c$. For $z>0$, this series is at most
  $\exp(qz^{1/q})$, by the $\ell^q$--$\ell^1$ inequality. Its term with
  $k=\lfloor z^{1/q}/2\rfloor$ is at least $2^{qk}$ for large $z$.
  This proves the two-sided growth. Finally, $\beta-1=3-\alpha-a<3-\alpha$.
  This proves No-go~\ref{P:nogo-01}.
\end{proof}

This restates the obstruction in Proposition~\ref{P:envelope} and
Remark~\ref{R:G0-verdict}. The lower estimate is for $B_K$, not for the
moment. The valid upper estimate improves the crude quadratic exponent,
but its factorial denominator misses $(k!)^a$. The split count $2k+1$
costs only an exponential factor in $k$; removing it does not repair the
loss. The same calculation for one chain, with $e^{\theta t}$, has the
same defect. Reopening requires an estimate retaining ordering or grouping
pairings before taking this common bound. Section~\ref{S:interval-upper}
provides such a different argument.

\subsection{NG-02: a single adjacent diagram is too large}
\label{S:nogo-02}

There are roughly an exponential factor times $k!$ pairings of $2k$
vertices. In white space, a uniform denominator $(k!)^3$ for each diagram
would therefore give the desired denominator $(k!)^2$ after summation.
The following example explains why that convenient sufficient estimate is
stronger than the problem permits.

\setcounter{nogo}{1}
\begin{nogo}[No uniform cubic factorial denominator per diagram]
  \label{P:nogo-02}
  Let $d=1$, $\Gamma=\delta_0$ and $0<a<1$. For the split $(0,2k)$
  and adjacent matching $\{(1,2),(3,4),\ldots,(2k-1,2k)\}$, its contribution is
  \[
    J_{\rm adj}(t)=
    \left[\frac{\Gamma_{\mathrm E}(1-a)}2\right]^k
    \frac{t^{(3-a)k}}{\Gamma_{\mathrm E}(1+(3-a)k)}.
  \]
  Consequently, there is no finite $C=C(a)$ bounding every diagram by
  $C^k t^{(3-a)k}/(k!)^3$ for all $k$ and $t>0$.
\end{nogo}
\begin{proof}
  We compute the adjacent diagram after its white contractions and justify
  that fixed-order limit. Its Dirichlet integral then shows why a uniform
  cubic factorial denominator is impossible.

  \medskip
  \noindent\textit{Step 1.~}We identify the adjacent-pair integral and its
  smoothing limit.
  Use root-forward increments $h_1,\ldots,h_{2k}>0$ with
  $\sum h_i<t$. Each adjacent contraction evaluates the even Green factor
  at zero, giving $1/2$, and integration of the preceding free increment
  gives mass $h_{2j-1}$. The temporal factor is $h_{2j}^{-a}$. Therefore,
  \[
    J_{\rm adj}(t)=2^{-k}
    \int_{\sum h_i<t}\prod_{j=1}^k h_{2j-1}h_{2j}^{-a}\,\ud h.
  \]
  Here delta contractions mean limits of Gaussian approximate identities.
  For completeness, fix $k$ and a strictly positive increment vector.
  After convolution the paired positions approach equality. The even Green
  factors are continuous at zero; the remaining indicator boundaries are
  finitely many hyperplanes in the $k$ free positions and have measure zero.
  The limit is thus the stated diagonal integral. All positions in the
  original chain lie in $[-t,t]$. Uniformly in the smoothing widths, replacing
  its Green indicators by one bounds the spatial integral by
  $2^{-2k}(2t)^k$: each disjoint smoothed pair has integral at most $2t$
  over $[-t,t]^2$. The temporal factors on the full time cube are a product
  of $k$ disjoint factors $|s-r|^{-a}$, each integrable because $a<1$.
  Dominated convergence justifies the fixed-order spatial limit through the
  time integration, and also temporal truncation limits.

  \medskip
  \noindent\textit{Step 2.~}We compare its growth with the proposed factorial
  bound.
  The Dirichlet simplex formula, with parameters $2$ and $1-a$ for each
  pair, now gives the asserted expression. Stirling's formula yields
  \[
    \log\big((k!)^3J_{\rm adj}(1)\big)=a k\log k+O(k).
  \]
  This exceeds $k\log C$ eventually for every fixed $C$.
  This proves No-go~\ref{P:nogo-02}.
\end{proof}

The adjacent matching is only one diagram. Its denominator is, up to
exponential factors, $(k!)^{3-a}$, which is stronger than the desired
summed denominator $(k!)^2$. Hence, this witness does not obstruct the
summed bound. A viable successor must retain the number of large diagrams,
or sum before estimating; it cannot assign the same cubic denominator to
all of them.

\subsection[NG-03: a squared-mean comparison fails]{NG-03: doubled variance does not give the proposed squared mean}
\label{S:nogo-03}

A positive mixture of exponentials of linear Gaussian functionals would
satisfy a useful comparison between the second moment and the mean at
noise amplitude $\sqrt2$. Such a representation must be established for
the wave rather than inferred from a parabolic analogy.

\setcounter{nogo}{2}
\begin{nogo}[Failure of the all-time comparison with unit prefactor]
  \label{P:nogo-03}
  For the admissible Riesz case $d=1$, $\Gamma(x)=|x|^{-1/2}$ and $a=1/2$,
  \[
    \E u(t,x)^2>\big(\E u_{\sqrt2}(t,x)\big)^2
  \]
  for all sufficiently small $t>0$, where $u_\lambda$ uses noise amplitude
  $\lambda$ and the same initial data.
\end{nogo}
\begin{proof}
  We compare the first nonconstant coefficients in the proposed inequality. We
  then bound both series remainders to turn that coefficient discrepancy into
  a strict small-time failure.

  \medskip
  \noindent\textit{Step 1.~}We compare the two first-order contributions.
  More generally, for $\Gamma(x)=c|x|^{-\alpha}$, $0<\alpha<1$,
  direct integration of the length-two chain gives
  \[
    a_1=\frac{c}{(1-\alpha)(2-\alpha-a)(3-\alpha-a)(4-\alpha-a)}.
  \]
  Indeed the spatial integral is
  $c h_1h_2^{1-\alpha}/(1-\alpha)$ and its temporal factor is $h_2^{-a}$;
  integrate over $h_1+h_2<1$. At $\alpha=a=1/2$, $c=1$, this gives $a_1=1/3$.
  For independent uniforms $U$ on $[-s,s]$ and $V$ on $[-r,r]$, $0<r\le s$,
  \[
    \E|U-V|=\frac{s}{2}+\frac{r^2}{6s}\le\frac{2s}{3}.
  \]
  Jensen's inequality for $z\mapsto z^{-\alpha}$ implies
  \[
    \int G_s(x)G_r(y)|x-y|^{-\alpha}\,\ud x\ud y
    \ge (3/2)^\alpha s^{1-\alpha}r.
  \]
  Integrating on the two temporal triangles gives
  \[
    b_1\ge
    \frac{2(3/2)^\alpha}{(1-a)(2-a)(4-\alpha-a)}
    =\frac89\sqrt{\frac32}>\frac23=2a_1.
  \]

  \medskip
  \noindent\textit{Step 2.~}We justify the strict inequality for sufficiently
  small positive time.
  The comparison is a finite-time failure, not just a formal coefficient test.
  Put $z=t^\beta$. The two nonnegative series under comparison have expansions
  \[
    \E u(t,x)^2=1+(2a_1+b_1)z+O(z^2),\qquad
    (\E u_{\sqrt2}(t,x))^2=1+4a_1z+O(z^2).
  \]
  To justify the remainders, choose any fixed positive $z_0$ at which these
  series converge, using the moment expansion also at amplitude $\sqrt2$.
  For each nonnegative coefficient sequence $(c_j)$,
  $\sum_{j\ge2}c_jz^j\le(z/z_0)^2\sum_{j\ge2}c_jz_0^j$ when
  $0\le z\le z_0$. Apply this to both series, including the nonnegative
  Cauchy product for the squared mean. Their difference therefore has the
  strictly positive leading coefficient $b_1-2a_1$.
  This proves No-go~\ref{P:nogo-03}.
\end{proof}

For a deterministic nonnegative mixture of Gaussian exponentials, the
inequality follows by Tonelli and
$\|h+k\|^2/2\le\|h\|^2+\|k\|^2$, provided the same mixture represents
all amplitudes in the comparison. That proposed mechanism is excluded by
the example. White noise, an eventual large-time comparison, or a prefactor
larger than one is not excluded. Reopening requires a restriction or a new
comparison that removes the contradicted first-coefficient requirement;
it cannot merely assert this positive-mixture representation.

\subsection[NG-04: constants near the time boundary]{NG-04: a coefficient constant cannot be uniform at the time boundary}
\label{S:nogo-04}

Comparing $A_k$ with $a_k$ retains the self-pairings of the Stratonovich
problem and would transfer a first-moment bound to the second moment.
However, a prescribed asymptotic constant must not be imposed on every
low-order coefficient without checking parameter uniformity.

\setcounter{nogo}{3}
\begin{nogo}[Divergence of the first coefficient ratio]
  \label{P:nogo-04}
  Fix $0<\alpha<1$ and $c>0$, with $\Gamma(x)=c|x|^{-\alpha}$ in $d=1$.
  Then
  \[
    \lim_{a\uparrow1}(1-a)\frac{A_1}{a_1}
    =2^{2-\alpha}(1-\alpha)>0.
  \]
  In particular no constant depending on $\alpha$ alone, including
  $2^{4-\alpha}$, can satisfy $A_k\le C^k a_k$ for every $k\ge1$
  and every $0<a<1$ at this fixed spatial covariance.
\end{nogo}
\begin{proof}
  We isolate the contribution near equal times and compute its boundary
  divergence. Dividing by the finite limiting single-chain coefficient shows
  why no constant uniform in the temporal exponent can work.

  \medskip
  \noindent\textit{Step 1.~}We identify the diagonal time contribution.
  Let
  $K(s,r)=\int G_s(x)G_r(y)c|x-y|^{-\alpha}\,\ud x\ud y$.
  This is continuous on $[0,1]^2$ by local integrability of the Riesz kernel
  and dominated convergence on a fixed square. Direct integration gives
  \[
    K(s,s)=\frac{c\,2^{1-\alpha}s^{2-\alpha}}
    {(1-\alpha)(2-\alpha)}.
  \]
  Writing $v=s-r$, define the continuous function
  $H(v)=2\int_v^1 K(s,s-v)\,\ud s$. Then
  \[
    b_1=\int_0^1v^{-a}H(v)\,\ud v,\qquad
    H(0)=\frac{c\,2^{2-\alpha}}
    {(1-\alpha)(2-\alpha)(3-\alpha)}.
  \]

  \medskip
  \noindent\textit{Step 2.~}We take the boundary limit of the coefficient
  ratio.
  For any continuous $H$ on $[0,1]$,
  $\varepsilon\int_0^1v^{-1+\varepsilon}H(v)\,\ud v\to H(0)$.
  Indeed the measures have mass one, and their mass on $[\delta,1]$ is
  $1-\delta^\varepsilon\to0$; continuity controls $[0,\delta]$.
  Thus $(1-a)b_1\to H(0)$. The formula for $a_1$ in the preceding proof gives
  \[
    a_1\longrightarrow
    \frac{c}{(1-\alpha)^2(2-\alpha)(3-\alpha)}.
  \]
  Combine these limits with $A_1=2a_1+b_1$.
  This proves No-go~\ref{P:nogo-04}.
\end{proof}

For each fixed $a<1$, a finite constant $C=C(\alpha,a)$ remains possible;
the argument only forces it to diverge as $a\uparrow1$. It also says
nothing against comparisons beginning at large order or large time with
suitable prefactors. Reopening requires explicitly allowing that parameter
dependence or removing the low orders. Historical static-endpoint tests
likewise constrain an endpoint or a uniform family, not a separately
chosen constant at an arbitrary positive $a$.

\subsection{NG-05: chronological merging breaks wave support}
\label{S:nogo-05}

A single chain has the monotone ordering needed by the first-moment
representation. Sorting two chains by time looks like a direct way to
recover it. Finite propagation prevents the resulting pointwise comparison.

\setcounter{nogo}{4}
\begin{nogo}[Support obstruction to pointwise concatenation]
  \label{P:nogo-05}
  In $d=1$, no finite constant can bound every product of two root-forward
  chain kernels by the single-chain kernel obtained by chronological merging.
  The failure occurs in the $(2,2)$ split on a set of positive measure for
  Riesz covariance and, with domestic pairs, for white covariance.
\end{nogo}
\begin{proof}
  We choose two admissible chains whose chronological merge violates a strict
  propagation constraint. We then keep this defect on sets of positive measure
  for each of the two covariance cases.

  \medskip
  \noindent\textit{Step 1.~}We exhibit the positive original product and zero
  merged product.
  Work at horizon one with root $(0,0)$. The first chain has times
  $(1/4,3/4)$ and positions $(x,x)$, $x=1/5$; the second has times
  $(1/2,7/8)$ and positions $(y,y)$, $y=-2/5$.
  The four original displacements and elapsed times are
  \[
    (1/5,1/4),\quad(0,1/2),\quad(2/5,1/2),\quad(0,3/8).
  \]
  Each displacement is strictly smaller than its elapsed time, so their
  Green product is $2^{-4}$. The merged chain instead includes the jump
  from $(1/4,1/5)$ to $(1/2,-2/5)$. Its displacement is $3/5>1/4$,
  so that Green factor vanishes. All these inequalities are strict.

  \medskip
  \noindent\textit{Step 2.~}We retain the support defect under the stated
  covariance integrations.
  For white noise choose domestic pairs within each original chain. The
  delta contractions leave precisely the two free spatial variables $x,y$.
  Open neighborhoods of the displayed four times and these two variables
  preserve both conclusions and have positive contracted measure. Gaussian
  regularization justifies this contraction as in No-go~\ref{P:nogo-02}.
  For Riesz covariance perturb all four spatial positions independently
  in small neighborhoods. Avoiding the domestic diagonals gives a set of
  positive ordinary measure on which the covariance factors are finite and
  strictly positive, while the same support defect persists.
  No multiplier can turn the zero merged product into a positive upper bound.
  This proves No-go~\ref{P:nogo-05}.
\end{proof}

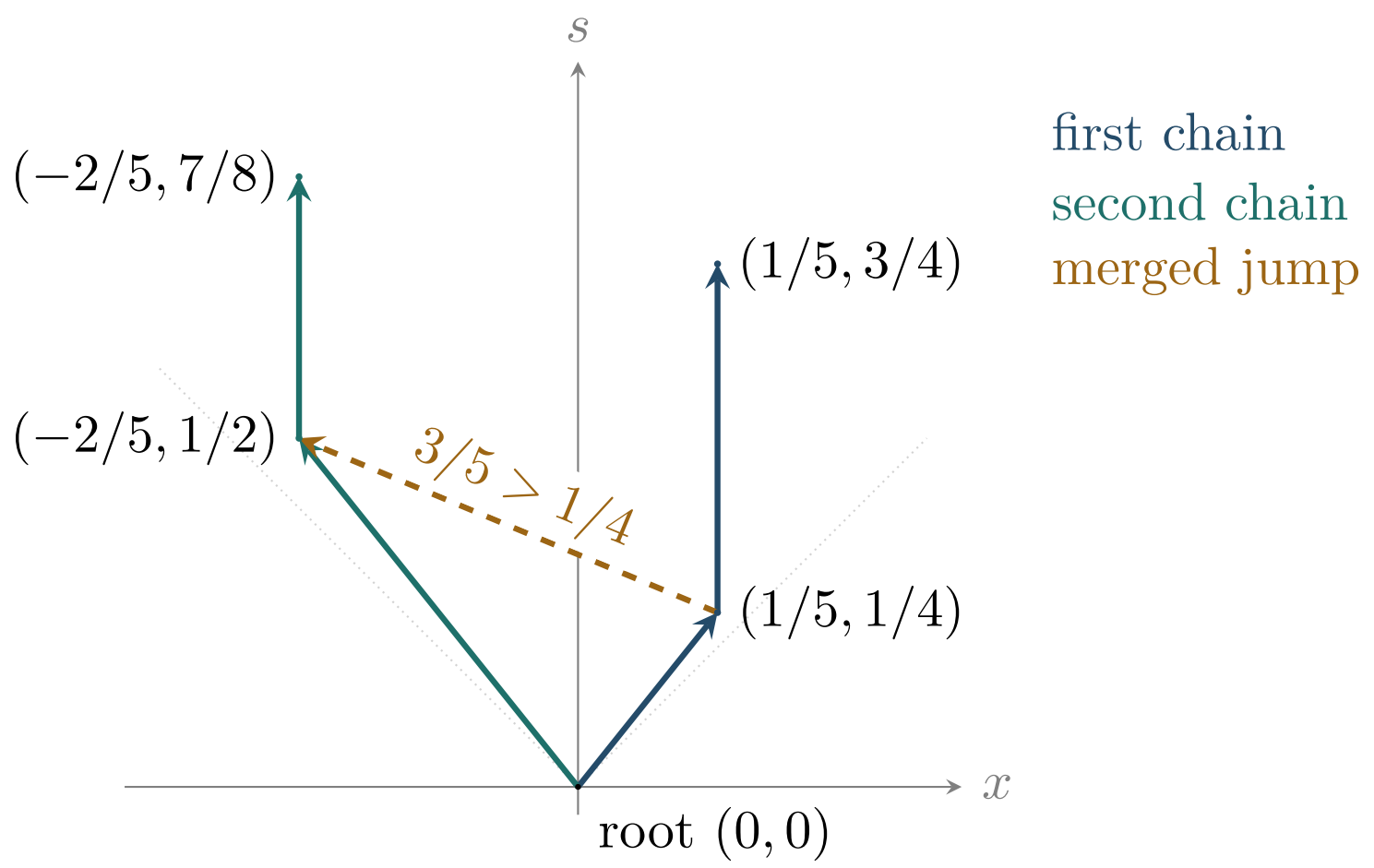
\begin{figure}[htbp]
  \centering
  \begin{tikzpicture}[x=5cm,y=5cm,>=stealth,font=\small]
    \draw[->,gray] (-.65,0) -- (.55,0) node[right] {$x$};
    \draw[->,gray] (0,-.04) -- (0,1.04) node[above] {$s$};
    \draw[gray!35,densely dotted] (-.6,.6) -- (0,0) -- (.5,.5);
    \coordinate (O) at (0,0);
    \coordinate (A) at (.2,.25);
    \coordinate (B) at (.2,.75);
    \coordinate (C) at (-.4,.5);
    \coordinate (D) at (-.4,.875);
    \draw[hamblue,very thick,->] (O) -- (A);
    \draw[hamblue,very thick,->] (A) -- (B);
    \draw[hamteal,very thick,->] (O) -- (C);
    \draw[hamteal,very thick,->] (C) -- (D);
    \draw[hamgold,very thick,dashed,->] (A) -- (C)
    node[midway,above,sloped,fill=white,inner sep=2pt]
    {$3/5>1/4$};
    \fill (O) circle (.6pt) node[below right] {root $(0,0)$};
    \fill[hamblue] (A) circle (.7pt)
    node[right,black] {$(1/5,1/4)$};
    \fill[hamblue] (B) circle (.7pt)
    node[right,black] {$(1/5,3/4)$};
    \fill[hamteal] (C) circle (.7pt)
    node[left,black] {$(-2/5,1/2)$};
    \fill[hamteal] (D) circle (.7pt)
    node[left,black] {$(-2/5,7/8)$};
    \node[hamblue,anchor=west] at (.65,.94) {first chain};
    \node[hamteal,anchor=west] at (.65,.84) {second chain};
    \node[hamgold,anchor=west] at (.65,.74) {merged jump};
  \end{tikzpicture}
  \caption[Merging two chains can violate wave support]{Root-forward
    coordinates $(x,s)$ for No-go~\ref{P:nogo-05}.
    Each solid edge lies strictly inside its propagation cone. Chronological
    merging inserts the dashed jump, whose distance $3/5$ exceeds its elapsed
    time $1/4$. The picture refutes pointwise concatenation, not an integrated
    comparison. For domestic white pairs the two free positions are \mbox{$x$ and $y$}.}
  \label{F:nogo-chain-geometry}
\end{figure}

The time coordinates here are root-forward chain coordinates; applying a
terminal-to-point convention again would change the object being tested.
The result concerns pointwise concatenation. Integrated comparisons and
comparisons of grouped Wick sums remain possible. Reopening requires
averaging, extra displacement variables, or another mechanism that preserves
the two original propagation constraints instead of imposing the false
merged one.

\subsection{NG-06: erasing geometry before summation loses the scale}
\label{S:nogo-06}

After white contractions, the temporal pairing sum is symmetric. One can
exploit that symmetry immediately by forgetting the spatial indicators.
The resulting bound is easy to sum but too weak for the desired exponent.

\setcounter{nogo}{5}
\begin{nogo}[The geometry-free majorant]
  \label{P:nogo-06}
  For $d=1$, $\Gamma=\delta_0$, $0<a<1$, replacing every wave indicator by
  one while retaining only the spatial box $[-t,t]$ gives the coefficient
  majorant
  \[
    A_k t^{(3-a)k}\le
    \frac{1}{k!}\left[\frac{2t^{3-a}}{(1-a)(2-a)}\right]^k.
  \]
  Its sum is $\exp(2t^{3-a}/((1-a)(2-a)))$, whose time exponent
  $3-a$ exceeds the target $(3-a)/2$.
\end{nogo}
\begin{proof}
  The spatial replacement gives $(2t)^k/2^{2k}=t^k/2^k$ for each pairing.
  After this replacement the sum over all pairings is symmetric in the
  $2k$ times. Its integral over a split $(m,n)$ simplex equals the cube
  integral divided by $m!n!$. The identities
  \[
    \sum_{m+n=2k}\frac1{m!n!}=\frac{2^{2k}}{(2k)!},\qquad
    P_k=\frac{(2k)!}{2^k k!},\qquad
    \int_0^t\int_0^t|s-r|^{-a}\,\ud s\ud r
    =\frac{2t^{2-a}}{(1-a)(2-a)}
  \]
  give the stated product and its exponential sum.
  This proves No-go~\ref{P:nogo-06}.
\end{proof}

This is a valid upper bound and a failure of its scale, not a moment lower
bound. Symmetry itself remains useful. The reopening condition is a spatial
estimate with enough order dependence retained before temporal
symmetrization; Proposition~\ref{P:d1-upper} retains this order dependence
and sums the temporal pairings afterwards (see
Remark~\ref{R:d1-upper-checks}).

\subsection{NG-07: block labels are not chronological ranks}
\label{S:nogo-07}

Temporal singularities favor close paired times. A proposed penalty for
long pairs must therefore measure time separation or actual chronological
rank. Labeling all vertices of one chain before those of the other creates
an artificial distance unrelated to either.

\setcounter{nogo}{6}
\begin{nogo}[Large label span does not imply a cubic denominator]
  \label{P:nogo-07}
  In white space with $0<a<1$, consider two length-$k$ chains and the matching
  $j\leftrightarrow k+j$. Although every pair has label distance $k$, its
  unit-time diagram satisfies
  \[
    J_k(1)\ge \left(\frac{3\,4^a}{128}\right)^k k^{-(3-a)k}.
  \]
  Thus no estimate assigning $C^k/(k!)^3$ to these diagrams solely on account
  of their block-label span can hold uniformly in $k$.
\end{nogo}
\begin{proof}
  Let the chain times be $s_j,r_j$, with increments $h_j,g_j$ from the
  common root. White contractions identify their paired positions, and
  integration of the successive spatial increments gives
  $2^{-k}\prod_{j=1}^k\min(h_j,g_j)$. Restrict both $s_j,r_j$ to
  $((j-1/4)/k,j/k)$. The domain has volume $(4k)^{-2k}$,
  each increment is at least $3/(4k)$, and each paired difference is at most
  $1/(4k)$. Multiplication gives the lower bound. The contraction can again
  be obtained at fixed $k$ by Gaussian smoothing; the spatial box bound and
  integrable disjoint temporal factors justify the limit.
  Stirling's formula implies
  $\log((k!)^3J_k(1))\ge a k\log k+O(k)$, contradicting any fixed $C$.
  This proves No-go~\ref{P:nogo-07}.
\end{proof}

Within each window the paired vertices are adjacent in the actual
chronological ordering. This example only refutes the label substitution;
it supplies no new shuffle estimate or Skorohod theorem. Reopening would
require a statement conditioned on actual order and an independent proof
of its summed cost. The shared pure-cross direction is recorded for later work and is not
developed here.

\subsection[NG-10: crossing covariance survives merging]{NG-10: an unweighted merge does not remove a crossing covariance}
\label{S:nogo-10}

A two-wave convolution estimate can remove an interior chain vertex when
no other factor depends on it. In a crossing pairing that vertex still
interacts with a remote partner. The following order-four calculation
identifies the missing weight before any all-order iteration is attempted.

\setcounter{nogo}{9}
\begin{nogo}[No finite spatial-supremum pullout in the merge lens]
  \label{P:nogo-10}
  In $d=2$, take $\Gamma(z)=c|z|^{-\alpha}$, $c>0$,
  $0<\alpha_0<1/2$ and $3/2\le\alpha<2-\alpha_0$.
  For the root-forward $(0,4)$ diagram
  with pairs $(1,3),(2,4)$, merging vertex two requires the weighted integral
  \[
    W=c\int_a^b(T-v)^{-\alpha_0}
    \int_{\mathbb R^2}G_{v-a}(y-x)G_{b-v}(z-y)|w-y|^{-\alpha}
    \,\ud y\ud v,
    \qquad 0<a<b<T.
  \]
  Here the letters $a,b,T$ are vertex times, so the temporal exponent is
  written $\alpha_0$. The spatial weight has no finite uniform essential
  supremum on admissible merge lenses. Hence, pulling it out by a uniform
  supremum cannot reduce $W$ to an unweighted two-wave convolution.
\end{nogo}
\begin{proof}
  Set $t=5$, $a=1$, $b=3$, $T=4$, $x=(0,0)$,
  $z=(1/2,0)$ and $w=(1/4,0)$. For $v$ near $2$, the point $y=w$
  is strictly inside both disks $|y-x|<v-a$ and $|z-y|<b-v$.
  It also satisfies the outer constraint $|w-z|<T-b$, and $|x|<a$.
  There is therefore an interior ball about $w$ in the merge lens.
  For every $M>0$ a smaller ball has positive Lebesgue measure and
  $|w-y|^{-\alpha}>M$ there. This proves infinite essential supremum;
  the strict inequalities persist under small changes of the fixed data.
  The other pair contributes $(b-a)^{-\alpha_0}c|z-x|^{-\alpha}$ outside
  the merge, whereas the two factors for pair $(2,4)$ depend on $v,y$
  and must remain in $W$. The temporal weight alone admits the bound
  $(T-v)^{-\alpha_0}\le(T-b)^{-\alpha_0}$.
  This proves No-go~\ref{P:nogo-10}.
\end{proof}

The unweighted powered-kernel convolution estimate in
\cite{bolanos-guerrero.nualart.ea:21:averaging}*{Lemma~4.3, p.~20}
assumes $1/2<q<1$ and $1\le\delta\le1/q$. It contains no covariance factor
linking the integration variable to a remote vertex and does not directly
estimate $W$ above. Infinite essential supremum does not imply divergence of $W$ or
of the diagram: the Riesz singularity is locally integrable for $\alpha<2$.
A reopening condition is a covariance-weighted convolution inequality,
or an explicit envelope integrable against the remaining diagram factors,
with the remote-gap cost retained. No such inequality is asserted here.

\subsection{NG-17: direct copying of the first-moment clock undercounts}
\label{S:nogo-17}

The source's first-moment clock representation suggests attaching one
Brownian clock to each ordered chain. A degree-two cross pair already tests
whether this copying preserves the interaction. It does not, even when
the two clocks are shared.

\setcounter{nogo}{16}
\begin{nogo}[Constant-correlation clock copying misses a positive factor]
  \label{P:nogo-17}
  Under the nonzero covariance hypotheses of
  Proposition~\ref{P:clock-copy-obstruction}, let $K_\Gamma(\theta)$ be the
  true cross-pair transform and $J_{\Gamma,\rho}(\theta)$ the direct
  unweighted candidate~\eqref{E:naive-clock-candidate}, $-1\le\rho\le1$.
  Then for every $\theta>0$,
  \[
    0<J_{\Gamma,\rho}(\theta)\le J_{\Gamma,1}(\theta)<K_\Gamma(\theta).
  \]
  If $\alpha=1$, then $J_{\Gamma,1}=(2-2^a)K_\Gamma$.
\end{nogo}
\begin{proof}
  Recall $\mathcal R_\Gamma(q)=\int(q^2+|\xi|^2)^{-1}\mu(\ud\xi)$.
  The exact transform identities proved in
  Proposition~\ref{P:clock-copy-obstruction} give
  \[
    w_\theta(\lambda)=
    \frac{\lambda^{a-2}}{2\theta\Gamma_{\mathrm E}(a)}
    [\mathcal R_\Gamma(\theta)-\mathcal R_\Gamma(\theta+\lambda)],
    \qquad
    K_\Gamma=\int_0^\infty w_\theta(\lambda)\,\ud\lambda,
  \]
  \[
    J_{\Gamma,1}=\int_0^\infty
    \frac{2\theta}{2\theta+\lambda}w_\theta(\lambda)\,\ud\lambda.
  \]
  The density is strictly positive because $\mu\ne0$ and the resolvent
  strictly decreases. It is integrable: homogeneity bounds it by multiples
  of $\lambda^{a-1}$ near zero and $\lambda^{a-2}$ at infinity.
  The multiplier is strictly between zero and one. For smaller correlations,
  $u+v-2\rho\min(u,v)\ge|u-v|$, and the positive Laplace integral
  $\Psi_a(r,v)$ in~\eqref{E:averaged-clock-kernel} decreases with $v$.
  This proves the chain of inequalities. When $\alpha=1$, insert
  $\mathcal R_\Gamma(q)=\mathcal R_\Gamma(1)/q$ and use
  $\int_0^\infty x^{a-1}/(x+c)\,\ud x
  =c^{a-1}\Gamma_{\mathrm E}(a)\Gamma_{\mathrm E}(1-a)$ to obtain the ratio.
  This proves No-go~\ref{P:nogo-17}.
\end{proof}

This is the already proved book obstruction, not a new clock impossibility
theorem. At the static endpoint the discrepancy disappears, so homogeneity
and static checks alone cannot detect it. At degree two the missing
shared-clock multiplier is exactly $1+\lambda/(2\theta)$; no higher-order
correction rule follows. Reopening requires relative-order weights or a
different coupling or representation that reproduces the actual cross-pair
transform. Remark~\ref{R:clock-copy-scope} gives the full boundary of this
statement, including the distinction between gap transforms and physical
unequal-time moments.

\subsection[NG-18: an oversized middle-cross envelope]{NG-18: erasing the middle-cross shift makes an envelope too large}
\label{S:nogo-18}

A later white-space comparison retains four vertices and a middle cross
pair. Removing the clock-overlap penalty gives a positive upper bound.
One further simplification, erasing the shift of its root kernel, also
increases the bound, but increases it beyond the comparator at an explicit
parameter point. This is an obstruction to closing that envelope, not a
counterexample to the original comparison.

\setcounter{nogo}{17}
\begin{nogo}[The shift-erased envelope is not below the wave transform]
  \label{P:nogo-18}
  Consider the white $(1,3)$ diagram with middle-cross partition $02\cdot13$,
  common gap damping one, cross rate $\lambda\ge0$ and self rate $\nu\ge0$.
  Put $b=1+\nu$, $c=1+\lambda$, $A=\lambda+\lambda^2/2$. Its actual wave
  transform $I$ and the shift-erased clock upper envelope $M_b(0;A)$ satisfy
  at $b=1,c=11/10$
  \[
    I=\frac{545}{11616},\qquad M_b(0;A)=\frac{1025}{21483},\qquad
    I-M_b(0;A)=-\frac{6005}{7562016}<0.
  \]
  Consequently the proposed pointwise closure $M_b(0;A)\le I$ is false.
\end{nogo}
\begin{proof}
  We identify the actual wave transform and the particular shift-erased
  envelope separately. We then evaluate both expressions at the stated
  rational parameters to disprove this envelope closure, without replacing the
  original clock candidate by the envelope.

  \medskip
  \noindent\textit{Step 1.~}We compute the actual middle-cross wave transform.
  We give the integral identifications so that the rational test specifies
  which bound fails. With singleton gap $u$ and long-chain gaps $v,w,z$,
  the spatial factor is
  $\int G_u(X)G_v(Y)G_w(X-Y)G_z(Y-X)\,\ud X\ud Y$.
  Set $R=u+v$, $q=u-v$. Integration in $R$, including its Jacobian,
  gives the exact wave expression
  \[
    I=\frac1{16}\int_{w,z>0}\int_{\mathbb R}
    e^{-b(w+z)-\lambda|q-w|}
    \int_0^{\min(w,z)}e^{-\max(|q|,t)}\,\ud t\ud q\ud w\ud z.
  \]
  To check the spatial reduction, the derivative of the area
  $\{|X|<u,|Y|<v,|X-Y|<C\}$ with respect to $C$ is
  $2(R-\max(C,|q|))_+$; the area vanishes at $C=0$.
  Integrate this against $e^{-R}/32$ for $R>|q|$.
  Next integrate $z\ge t$ and put $w=t+x$. Splitting the remaining $q$
  integral into $q<0$, $0\le q\le t$, and $q\ge t$ gives respectively
  \[
    \frac{1/c+1/(2b+\lambda+1)}{16b(b+\lambda)(2b+2\lambda+1)},\quad
    \frac1{16b(b+\lambda)(2b+1)(2b+1+\lambda)},
  \]
  \[
    \frac{1/(b+\lambda)+1/c}{16b(2b+1)(b+1)}.
  \]
  Adding gives
  \[
    I=\frac{2b^2+4bc+2c^2-c+1}
    {8bc(b+1)(2b+1)(b+c-1)(2b+2c-1)}.
  \]

  \medskip
  \noindent\textit{Step 2.~}We identify the root kernel and the shift factor
  that is discarded.
  For the clock envelope write $h(r)=(2\pi r)^{-1/2}$ and
  $m(s)=1/(\pi\sqrt2\sqrt{s-1/2})$, $s>1/2$.
  The two spatial pair differences have covariance matrix
  $\left(\begin{smallmatrix}u+v+w&-w\\-w&w+z\end{smallmatrix}\right)$.
  Writing $r=w+z$, its determinant is $rR+wz$, so their joint white
  density factors as $h(r)h(R+wz/r)$. After deleting the nonnegative
  clock-overlap penalty, the root kernel is
  \[
    H_\gamma(q)=\int_{|q|}^\infty e^{-R/2}h(R+\gamma)\,\ud R
    =\int_{1/2}^\infty\frac{m(s)}s
    e^{-(s-1/2)\gamma}e^{-s|q|}\,\ud s,
    \qquad \gamma=\frac{wz}{r}.
  \]
  The equality follows from the positive Gamma integral for $h$.
  The convolution of the mixture component $q\mapsto e^{-s|q|}$ with
  $q\mapsto e^{-A|q|}$, evaluated at $w$, equals $2S_{A,s}(w)/(s+A)$, where
  \[
    S_{A,s}(w)=e^{-Aw}+\int_0^w A e^{-At}e^{-s(w-t)}\,\ud t\le1.
  \]
  The inequality follows by bounding $e^{-s(w-t)}$ by one; the remaining
  integral is $1-e^{-Aw}$.
  Replacing $S$ by one is precisely the shift erasure. More generally, for
  $\eta\ge0$, replacing it by $e^{-\eta w}$ defines
  \[
    M_b(\eta;A)=\frac1{16}\int_{1/2}^\infty
    \frac{m(s)}{s(s+A)}\int_0^\infty h(r)e^{-b^2r/2}
    \int_0^r e^{-(s-1/2)w(r-w)/r-\eta w}\,\ud w\ud r\ud s.
  \]

  \medskip
  \noindent\textit{Step 3.~}We evaluate the shift-erased envelope and compare
  the rational values.
  All integrals here are positive. Set $w=rt$ and $\tau=s-1/2$.
  For $\eta,\tau\ge0$, put $p=\sqrt{b^2+2\eta}$. The elementary bridge
  integral is
  \[
    \int_0^1[b^2+2\eta t+2\tau t(1-t)]^{-3/2}\,\ud t
    =\frac{2(p+b)}{bp((p+b)^2+2\tau)}.
  \]
  It follows by integrating a quadratic power; at $\eta=\tau=0$ both sides
  are $b^{-3}$. Substitution $z=\sqrt{2s-1}$ and partial fractions give
  \begin{align*}
    M_b(\eta;A)&=\frac{p+b}{2\pi bp}
    \int_0^\infty\frac{\ud z}
    {(z^2+1)(z^2+c^2)(z^2+(p+b)^2)}\\
    &=\frac{1+c+p+b}{4bp c(c+1)(c+p+b)(1+p+b)}.
  \end{align*}
  The usual identity $\int_0^\infty(z^2+q^2)^{-1}\ud z=\pi/(2q)$
  proves the partial fraction evaluation; coincident rates follow by
  continuity. At $\eta=0$, $p=b$. Substituting the asserted rational
  parameters in these two expressions gives the displayed strict gap.
  This proves No-go~\ref{P:nogo-18}.
\end{proof}

The boundary value $\nu=0$ is a legitimate fixed-parameter transform;
continuity of both rational expressions also gives an overshoot for some
$\nu>0$ nearby. The no-go does not refute $J\le I$, where $J$ is the
original shared-clock candidate, nor does it decide an integrated
comparison using this envelope. The surviving middle-cross comparison
retains $S_{A,1/2}(w)$, instead of replacing it by one, and obtains a smaller
envelope. Reopening this particular closure requires keeping a shift cost
or proving a different integrated estimate. The historical standalone
record was not separately adjudicated; the subsequent
middle-cross analysis and this reconstructed proof do not change that
historical review status.

\section{Covariance transformations and their limits}\label{S:nogo-covariance}

Covariance comparisons offer a way to use a simpler model without estimating
all the pairings again. Three different requirements enter such an argument:
positive definiteness, local finiteness, and uniformity of the comparison.
None follows from homogeneity alone. The examples below separate these
requirements. The final example also distinguishes positive semidefinite
order from pointwise order, which matters for the wave mean.

\subsection{NG-08: the translated spatial maximum}\label{S:nogo-08}

The conditional higher-dimensional argument replaces a translated wave
convolution by a uniform spatial maximum. This would remove the dependence
on the other vertices of a pairing. The price is an additional hypothesis,
not a consequence of the covariance assumptions.

\setcounter{nogo}{7}
\begin{nogo}[A spatial maximum need not be finite]\label{P:nogo-08}
  Hypothesis~\ref{H:spatial-maximum} fails for two-dimensional Riesz
  covariance of degree $-\alpha$ when $3/2\le\alpha<2$. For each
  $1<\alpha<2$, it also fails for some nonnegative, positive definite,
  locally integrable homogeneous product covariance in dimensions two
  and three. These examples can satisfy the model condition
  $\alpha+\alpha_0<2$.
\end{nogo}
\begin{proof}
  We exhibit the boundary divergence for planar Riesz covariance and then the
  anisotropic product example in both dimensions. The smoothing limit shows
  why neither example admits the proposed uniform spatial maximum.

  \medskip
  \noindent\textit{Step 1.~}We test the translated planar Riesz convolution.
  We restate the witnesses from Proposition~\ref{P:higher-obstructions}.
  For Riesz covariance translate the singularity to $e_1$, on the boundary
  of the unit wave disk. Write $y=e_1+\rho(\cos\phi,\sin\phi)$ and restrict
  to $0<\rho<1/2$, $2\pi/3\le\phi\le4\pi/3$. Then
  $0<1-|y|^2\le2\rho$, so the wave density and polar Jacobian give
  \[
    \int G_1^{(2)}(\ud y)c|y-e_1|^{-\alpha}
    \ge \frac{c}{3\sqrt2}\int_0^{1/2}\rho^{1/2-\alpha}\ud\rho
    =\infty.
  \]
  At $\alpha=3/2$ the divergence is logarithmic.

  \medskip
  \noindent\textit{Step 2.~}We test the product covariance and its smoothed
  limits.
  For the product example choose
  \[
    a_1=\frac{1+\max\{2-\alpha,\alpha-1\}}2,
    \qquad a_2=\alpha-a_1.
  \]
  Then $0<a_1,a_2<1$ and $a_1+a_2/2>1$. The function
  $\gamma(x)=c|x_1|^{-a_1}|x_2|^{-a_2}$, constant in $x_3$ when $d=3$,
  is locally integrable and homogeneous. Each factor is a positive
  Gaussian mixture, hence the product is positive definite. With
  $u=1-y_1$ and $y_2=\sqrt{2u-u^2}\,w$, direct integration yields
  \[
    \int G_1^{(2)}(\ud y)\gamma(y-e_1)
    =\frac{c B((1-a_2)/2,1/2)}{2\pi}
    \int_0^2u^{-a_1}(2u-u^2)^{-a_2/2}\ud u=\infty.
  \]
  The projection of the three-dimensional wave measure onto its first
  two coordinates is $G_1^{(2)}$, which proves the three-dimensional claim.
  Gaussian smoothing converges off the singular sets, which have zero
  translated wave measure. Fatou's lemma therefore gives divergence in
  lower limit as the smoothing widths decrease to zero. Thus, smoothing
  cannot restore the required uniform constant. Taking
  $0<\alpha_0<2-\alpha$ ensures admissibility.
  This proves No-go~\ref{P:nogo-08}.
\end{proof}

This result excludes the spatial-supremum input, not moment finiteness or
the required time scale. The conditional theorem remains useful where its
hypothesis holds. Reopening this route for the full class would require
an integrated estimate retaining the geometry lost by the supremum;
Section~\ref{S:interval-upper} supplies a different full-class argument.

\subsection{NG-09: radial lowering in three dimensions}\label{S:nogo-09}

Multiplication by $|x|^{\alpha-1}$ gives the desired endpoint degree $-1$.
The following example shows why this degree calculation does not construct
a covariance. The original endpoint review returned MEDIUM confidence and
was not accepted at the HIGH threshold. The argument below is a complete
reconstruction; it does not change that historical disposition.

\setcounter{nogo}{8}
\begin{nogo}[Radial lowering can destroy positive definiteness]\label{P:nogo-09}
  In dimension three let
  \[
    \Gamma(x)=|x_1|^{-9/10}|x_2|^{-3/5},\qquad
    F(x)=|x|^{1/2}\Gamma(x).
  \]
  The covariance $\Gamma$ satisfies the homogeneous model assumptions with
  $\alpha=3/2$ and $\alpha_0=1/4$. The locally integrable tempered function
  $F$ is not positive definite as a distribution.
\end{nogo}
\begin{proof}
  We verify the original covariance, then show that the radially modified
  function has a heat convolution growing along the third coordinate. This
  contradicts the boundedness that positive definiteness would require.

  \medskip
  \noindent\textit{Step 1.~}We check the original covariance and the modified
  function.
  The exponents of both coordinate factors lie in $(0,1)$. Local
  integrability, temperedness, and positive definiteness of $\Gamma$
  follow from the one-dimensional Gaussian mixture formula, applied in
  each coordinate. The third factor is the constant covariance one.
  Multiplication by $|x|^{1/2}$ preserves local integrability and temperedness.

  \medskip
  \noindent\textit{Step 2.~}We obtain the contradiction after Gaussian
  convolution.
  Fix a centered Gaussian heat density $p_\sigma$. Its convolution with
  $F$ is finite at zero: the coordinate singularities are integrable and
  the remaining growth is polynomial. On the other hand, restrict the
  convolution at $Le_3$ to
  \[
    1\le y_1\le2,\qquad 1\le y_2\le2,\qquad |y_3-L|\le1.
  \]
  For $L\ge2$, the Gaussian factor $p_\sigma(Le_3-y)$ is bounded below
  there by a fixed positive constant, the coordinate powers are bounded
  below, and $|y|^{1/2}\ge(L-1)^{1/2}$. Hence,
  \[
    (p_\sigma*F)(Le_3)\ge c_\sigma\sqrt L.
  \]
  If $F$ were positive definite, its quadratic form on Schwartz tests
  would satisfy Cauchy--Schwarz. Take a Gaussian test whose autocorrelation
  is $p_\sigma$, and its translate by $x$. The two diagonal quadratic
  forms are both $(p_\sigma*F)(0)$, while the cross form is
  $(p_\sigma*F)(x)$. Consequently,
  \[
    |(p_\sigma*F)(x)|\le(p_\sigma*F)(0),\qquad x\in\mathbb R^3.
  \]
  The lower bound along $Le_3$ contradicts this inequality.
  This proves No-go~\ref{P:nogo-09}.
\end{proof}

The obstruction concerns this Euclidean radial multiplication. It leaves
coordinatewise constructions and the planar endpoint argument intact;
see Remark~\ref{R:three-endpoint-scope}. A successor must verify positive
definiteness of its endpoint, rather than infer it from degree and sign.

\subsection{NG-11: additive variance interpolation}\label{S:nogo-11}

A covariance mixture would be easy to handle if its second moment were
log-convex in the additive variance. The following scalar calculation
refutes that unrestricted principle. It is outside the strict homogeneous
model and is not a counterexample to a model-specific comparison.

\setcounter{nogo}{10}
\begin{nogo}[A scalar wave moment is not log-convex in variance]\label{P:nogo-11}
  Let $y''=Zy$, $y(0)=1$, $y'(0)=0$, where $Z$ is centered Gaussian of
  variance $A$. For $F(A)=\E[y(1)^2]$,
  \[
    F(100)^2>F(99)F(101).
  \]
  Thus, $\log F$ is not convex on $(0,\infty)$.
\end{nogo}
\begin{proof}
  We express the scalar moment as a positive series and bound its tail. The
  resulting finite rational enclosures disprove the stated log-convexity
  inequality.

  \medskip
  \noindent\textit{Step 1.~}We obtain the series and a controlled tail.
  Interpreting $\cosh\sqrt z$ by its entire power series gives
  \[
    F(A)=1+\sum_{k\ge1}t_k(A),\qquad
    t_k(A)=\frac{2^{3k-1}(2k)!}{k!(4k)!}A^k.
  \]
  Absolute termwise Gaussian integration follows from
  $\E\exp(2\sqrt{|Z|})<\infty$. All displayed terms are positive.
  Their ratio satisfies
  \[
    \frac{t_{k+1}(A)}{t_k(A)}
    =\frac{16A(2k+1)}{(4k+1)(4k+2)(4k+3)(4k+4)}
    \le\frac{16A}{(4k+1)^3}.
  \]
  Write $S_{20}=1+\sum_{k=1}^{20}t_k(A)$ and
  $q_A=16A/81^3<1$ for the three values in question. Then
  \[
    S_{20}\le F(A)\le S_{20}+\frac{t_{21}(A)}{1-q_A}.
  \]

  \medskip
  \noindent\textit{Step 2.~}We compare the rigorous rational enclosures.
  Substitution in these finite rational expressions gives the following
  outward enclosures, with common denominator $10^6$:
  \[
    \begin{array}{c|rr}
      A&10^6\text{ times lower bound}&10^6\text{ times upper bound}\\ \hline
      99&207442589&207442590\\
      100&212332413&212332414\\
      101&217304018&217304019
    \end{array}
  \]
  In particular,
  \[
    F(100)^2-F(99)F(101)
    \ge\frac{6945091633359}{10^{12}}>0.
  \]
  This proves No-go~\ref{P:nogo-11}.
\end{proof}

Positive power-series coefficients do imply convexity of
$t\mapsto\log F(e^t)$, by H\"older's inequality. Thus, geometric
interpolation of amplitudes survives. An additive-covariance argument
would need additional hypotheses and a separate proof; the scalar test
prevents using log-convexity as a general Gaussian-wave fact.

\subsection{NG-12: a heat--Mellin square need not be locally finite}
\label{S:nogo-12}

Squaring a heat-smoothed covariance preserves positive definiteness, and
a Mellin weight can prescribe homogeneity. This suggests an endpoint
recipe with no angular decomposition. Its missing requirement is local
finiteness of the infinite positive mixture.

\setcounter{nogo}{11}
\begin{nogo}[Failure of a heat--Mellin endpoint recipe]\label{P:nogo-12}
  Let $d=3$, $\alpha=11/10$, $\alpha_0=1/4$, and
  \[
    \Gamma(x)=|x_1|^{-9/10}|x_2|^{-1/5},\qquad F_s=p_s*\Gamma.
  \]
  Although $\Gamma$ is an admissible covariance, the putative endpoint
  density of degree $-\beta$, with $\beta=9/5$,
  \[
    H(x)=\int_0^\infty
    s^{(2\alpha-\beta)/2-1}F_s(x)^2\ud s
    =\int_0^\infty s^{-4/5}F_s(x)^2\ud s
  \]
  does not define a locally finite measure.
\end{nogo}
\begin{proof}
  The coordinate exponents show admissibility as in No-go~\ref{P:nogo-09}.
  Heat convolution factorizes in the first two coordinates. For
  $0<s<1/16$ and $|x_1|\le\sqrt s$, integrating the first Gaussian
  convolution over $\sqrt s\le y_1\le2\sqrt s$ gives
  \[
    (p_s^{(1)}*|\,\cdot\,|^{-9/10})(x_1)\ge c s^{-9/20}.
  \]
  For $1/2\le x_2\le1$, integrate the second convolution over
  $x_2\le y_2\le x_2+\sqrt s$; its value is bounded below by a positive
  constant. The third coordinate contributes one. On the fixed box
  $B=[-1,1]\times[1/2,1]\times[-1,1]$, the strip
  $|x_1|\le\sqrt s$ therefore gives
  \[
    \int_B F_s(x)^2\ud x\ge c s^{1/2-9/10}=c s^{-2/5}.
  \]
  Tonelli's theorem implies
  \[
    \int_B H(x)\ud x\ge c\int_0^{1/16}s^{-6/5}\ud s=\infty.
  \]
  This proves No-go~\ref{P:nogo-12}.
\end{proof}

The positive-definiteness calculation for finite mixtures is not the
problem. The infinite mixture has infinite mass on a compact set. This
excludes the displayed recipe, not endpoint interpolation itself. A
replacement must prove local finiteness, for example by controlling the
coordinate singularities before integrating the smoothing scale.

\subsection{NG-13: angular smoothing and isotropic envelopes}\label{S:nogo-13}

Smooth angular densities admit finite angular suprema. It is tempting to
approximate a singular covariance by such densities and apply the Riesz
estimate with that supremum as amplitude. The limit below shows why the
amplitudes cannot remain bounded.

\setcounter{nogo}{12}
\begin{nogo}[No uniform isotropic envelope for a plane covariance]
  \label{P:nogo-13}
  In $\mathbb R^3=\mathbb R^2\times\mathbb R$ set
  \[
    \Gamma(\ud y,\ud z)=\frac1{2\pi}|y|^{-1/2}\ud y\,\delta_0(\ud z).
  \]
  This is an admissible covariance with $\alpha=3/2$ and
  $\alpha_0=1/4$. If nonnegative measures $\Gamma_n$ converge vaguely to
  $\Gamma$ and satisfy
  \[
    \Gamma_n(\ud x)\le C_n|x|^{-3/2}\ud x,
  \]
  then $C_n\to\infty$.
\end{nogo}
\begin{proof}
  We verify that the plane-supported covariance is admissible. A bounded
  subsequence of isotropic envelope constants would force its vague limit to
  have a Lebesgue density, contradicting its planar mass.

  \medskip
  \noindent\textit{Step 1.~}We verify the singular covariance example.
  The planar Riesz factor is locally integrable and positive definite by
  its Gaussian mixture representation. The normal Dirac covariance is
  positive definite: in its quadratic form the two normal coordinates
  coincide, leaving an integral of nonnegative tangential quadratic
  forms. Thus, the product is positive definite. For every $r>0$ and
  Borel set $A$, the substitution $y=ry'$ in its defining integral gives
  $\Gamma(rA)=r^{2-1/2}\Gamma(A)=r^{3/2}\Gamma(A)$. Its finite mass on
  a radius $R$ ball is proportional to $R^{3/2}$; polynomial mass growth
  also gives temperedness.

  \medskip
  \noindent\textit{Step 2.~}We contradict bounded envelope constants by
  singular support.
  Suppose a subsequence of $C_n$ were bounded by $C<\infty$. For every
  nonnegative compactly supported continuous $f$, vague convergence gives
  \[
    \int f\ud\Gamma=\lim_n\int f\ud\Gamma_n
    \le C\int f(x)|x|^{-3/2}\ud x.
  \]
  This is measure domination, so $\Gamma$ would be absolutely continuous
  with respect to three-dimensional Lebesgue measure. But the compact
  planar annulus $\{(y,0):1\le|y|\le2\}$ has positive $\Gamma$ mass and
  zero Lebesgue measure. This contradiction excludes every bounded
  subsequence, proving the asserted limit.
  This proves No-go~\ref{P:nogo-13}.
\end{proof}

This argument applies in particular to angular smoothing whenever the
smoothed covariances converge to the plane measure. It requires no
ellipsoid representation. It rules out a uniform scalar isotropic density
envelope, not other covariance majorants or density of any proposed cone.
The reopening condition is an estimate that remains uniform without taking
an angular supremum, or a comparison adapted to the singular support.

\subsection{NG-19: discarding an independent Gaussian residual}
\label{S:nogo-19}

An orthogonal Gaussian projection produces an independent residual with
positive semidefinite covariance. This makes projection a useful change-of-law
tool. It does not give pointwise nonnegative residual covariance, and the
wave mean need not increase when that residual is restored.

\setcounter{nogo}{18}
\begin{nogo}[A positive semidefinite residual can lower the wave mean]
  \label{P:nogo-19}
  There is an admissible one-dimensional parent covariance $C$ and an
  actual rank-one Gaussian projection with covariance $P$ such that
  $C-P$ is positive semidefinite but
  \[
    \E u_C(h;1,0)<\E u_P(h;1,0)
  \]
  for all sufficiently small positive noise couplings $h$. The projection
  $P$ is not asserted to be stationary or homogeneous.
\end{nogo}
\begin{proof}
  We construct an actual Gaussian projection of the parent field and compare
  the two degree-two mean coefficients. Bounds on the higher-order tails then
  turn that coefficient gap into a strict small-coupling failure of residual
  monotonicity.

  \medskip
  \noindent\textit{Step 1.~}We construct the parent field and its rank-one
  projection.
  Put $a=99/100$ and
  \[
    C((s,x),(t,y))=(1-a)|s-t|^{-a}|x-y|^{-1/2}.
  \]
  Its spatial degree is $\alpha=1/2$, and $\alpha+a<2$. Gaussian mixtures
  of the two power kernels verify positive definiteness and local
  integrability. Let $\psi(s,x)=\mathbf1_{[0,1]}(s)\mathbf1_{[-1,1]}(x)$,
  and $V=\E W(\psi)^2$. The indicator is a legal Hilbert-space test by
  bounded smooth approximation and dominated convergence. Local integrability
  and positivity of $C$ give $0<V<\infty$. Put $Z=W(\psi)/\sqrt V$ and
  define the covariance profile by the locally integrable kernel formula
  \[
    m(s,x)=\frac1{\sqrt V}\int_0^1\int_{-1}^1
    C((s,x),(r,y))\ud y\ud r.
  \]
  For every smooth compactly supported spacetime test $f$, Fubini's
  theorem identifies $\langle m,f\rangle=\E[W(f)Z]$.
  Thus, $m$ is the covariance of the generalized field with $Z$; no
  pointwise random variable $W(s,x)$ is used. Gaussian projection gives
  $W=Zm+W^\perp$ with $W^\perp$ independent of $Z$ and
  covariance $C-P$, where $P=m\otimes m$. Cauchy--Schwarz gives its
  positive semidefiniteness.

  \medskip
  \noindent\textit{Step 2.~}We bound the projected profile and compare the
  mean coefficients.
  Direct integration gives, on $0\le s\le1$,
  \[
    E_\phi=\frac{16\sqrt2}{3},\quad V=\frac{2E_\phi}{2-a},\quad
    m(s,x)=\frac{B_a(s)F(x)}{\sqrt V},
  \]
  where $B_a(s)=s^{1-a}+(1-s)^{1-a}\ge1$ and
  $F(x)=\int_{-1}^1|x-y|^{-1/2}\ud y$. We have $F\le4$ everywhere and
  $F(x)\ge2\sqrt2$ for $|x|\le1$. Thus, $m$ is bounded and nonnegative,
  and $m(s,x)m(t,y)\ge8/V$ on the backward cone of $(1,0)$.

  For a covariance $D$ the degree-two ordinary mean coefficient is
  \[
    A_D=\int_{0<s<t<1}\int_{\mathbb R^2}
    G_{1-t}(-y)G_{t-s}(y-z)D((s,z),(t,y))
    \ud z\ud y\ud s\ud t.
  \]
  This is the single self-pair term of Proposition~\ref{P:mixed}; there
  is no additional factorial. With $G_r(x)=\mathbf1_{\{|x|<r\}}/2$,
  integration of the parent covariance yields
  \[
    A_C=\frac{1-a}{1-\alpha}
    \frac1{(2-\alpha-a)(3-\alpha-a)(4-\alpha-a)}
    =\frac{20000}{1932951}<\frac1{96}.
  \]
  For covariance one the same integral is
  $\int_{0<s<t<1}(1-t)(t-s)\ud s\ud t=1/24$. Therefore,
  \[
    A_P\ge\frac{8}{24V}
    =\frac{101}{3200\sqrt2}>\frac1{48}.
  \]

  \medskip
  \noindent\textit{Step 3.~}We pass the strict coefficient gap to the full
  means.
  We justify passage from coefficients to full means. Write the parent
  solution as $u_C(h)=\sum_{n\ge0}h^nX_n$. For any fixed $r>1$,
  finiteness of the second moment at amplified covariance $r^2C$ and
  nonnegativity of the mixed moments in Proposition~\ref{P:mixed} give
  \[
    \|X_n\|_2\le \sqrt{B_r}\,r^{-n},\qquad
    B_r=\E|u_C(r;1,0)|^2<\infty.
  \]
  Hence, the $L^2$ norm of the tail starting at $n=4$ is at most
  $\sqrt{B_r}(|h|/r)^4/(1-|h|/r)$. Odd means vanish by Gaussian sign
  symmetry, giving $\E u_C(h)=1+h^2A_C+O(h^4)$.

  For $P$, the potential is the bounded profile $Zm$. If $M=\sup|m|$,
  the deterministic Picard coefficients at time one are bounded by
  $M^n/(2n)!$. The absolute mean series for $|h|\le1$ is bounded by
  $\E\cosh\sqrt{M|Z|}<\infty$. Factoring $h^4$ from its tail and
  using symmetry therefore gives $\E u_P(h)=1+h^2A_P+O(h^4)$.
  The strict coefficient gap proves the claimed inequality.
  This proves No-go~\ref{P:nogo-19}.
\end{proof}

The same distinction persists under a variance tilt of the projected
coordinate. For $1<q<2$, its density relative to the parent law is
\[
  L_q(Z)=q^{-1/2}\exp\bigl((1-q^{-1})Z^2/2\bigr),\qquad
  \E L_q^2=[q(2-q)]^{-1/2}.
\]
The equivalent Gaussian law has covariance $Q=C+(q-1)P$. The same original
measurable solution has an $O(h^4)$ mean remainder under this law by
Cauchy--Schwarz and the preceding $L^2$ tail estimate. Its degree-two
coefficient is $A_Q=A_C+(q-1)A_P<qA_P-1/96$. This is an equivalent-law
statement, not a new stationary-model existence theorem for $Q$.
To identify its degree-two coefficient, one can use the same regularization
as the parent construction. If $K_\delta$ denotes the truncated temporal
covariance and $p_\epsilon$ the spatial heat kernel, its cross-covariance
with $Z$ is
\[
  m_{\epsilon,\delta}(s,x)
  =\frac{(1-a)\int_0^1K_\delta(s-r)\ud r}{\sqrt V}
  (p_\epsilon*F)(x).
\]
The bounds $0\le K_\delta\le|\,\cdot\,|^{-a}$ and $F\le4$ give
$0\le m_{\epsilon,\delta}\le8/\sqrt V$. These profiles converge to
$m$ at interior times. Dominated convergence on the finite-mass terminal
wave kernel identifies the added $P$ contribution. The parent degree-two
term converges in $L^2$, and $L_q\in L^2$ transfers its expectation to the
equivalent law. Thus, the fixed-order Gaussian covariance identity indeed
gives the stated $A_Q$, without invoking a homogeneous representation for
$Q$. Consequently, discarding the residual also fails in this tilted comparison.
A large-time bridge for specially chosen projections could still hold,
but it must be proved at the required scale. Positive semidefinite order
alone cannot supply it.

\section{Phase propagation and Gaussian controls}\label{S:nogo-phase}

A wave equation transports displacement and velocity together. Its scalar
solution value is not a state variable for composition. Likewise, a
Gaussian change of law transports an entire random solution, not merely
the solution evaluated at the mean of the new noise. These distinctions
explain the failures in this section. The scalar examples are diagnostics
for generic proof principles; their scope is stated separately from the
homogeneous stochastic wave model.

\subsection{NG-14: scalar moments do not compose across blocks}
\label{S:nogo-14}

Independent static time blocks look favorable for multiplying second
moments. Independence does remove correlations between their potentials,
but it does not reset the velocity to zero at a block boundary. The
original record computed formal coefficients only. We include the
integrable-series argument needed to turn that diagnostic into a
small-variance counterexample.

\setcounter{nogo}{13}
\begin{nogo}[Two independent wave blocks do not multiply scalar moments]
  \label{P:nogo-14}
  Let $Z_1,Z_2$ be independent centered Gaussians of variance $A$. Solve
  $y''=Z_1y$ on $[0,1]$ and $y''=Z_2y$ on $[1,2]$, with continuous phase
  at time one and initial phase $(y(0),y'(0))=(1,0)$. Then, as
  $A\downarrow0$,
  \[
    \E y(1)^2=1+\frac A3+O(A^2),\qquad
    \E y(2)^2=1+3A+O(A^2).
  \]
  In particular, $\E y(2)^2>(\E y(1)^2)^2$ for all sufficiently small
  positive $A$.
\end{nogo}
\begin{proof}
  We compute the full two-block phase product rather than resetting the
  velocity. We then control the Gaussian series remainders to show that the
  mismatch of scalar moments persists at small positive variance.

  \medskip
  \noindent\textit{Step 1.~}We compute the first variance coefficients of the
  one-block and two-block moments.
  Define the entire functions
  \[
    c(z)=\sum_{n\ge0}\frac{z^n}{(2n)!},\qquad
    s(z)=\sum_{n\ge0}\frac{z^n}{(2n+1)!}.
  \]
  The unit-time phase propagator and the terminal displacement are
  \[
    M(z)=\begin{pmatrix}c(z)&s(z)\\zs(z)&c(z)\end{pmatrix},\qquad
    y(2)=c(Z_2)c(Z_1)+s(Z_2)Z_1s(Z_1).
  \]
  Expanding through total degree two in $(z,w)$ gives
  \[
    c(w)c(z)+s(w)zs(z)
    =1+\frac32z+\frac12w
    +\frac5{24}z^2+\frac5{12}zw+\frac1{24}w^2
    +\text{terms of degree at least three}.
  \]
  After squaring, the coefficients of $z^2$ and $w^2$ are respectively
  $8/3$ and $1/3$. Independence and centering eliminate the $zw$ term.
  Similarly, the coefficient of $z^2$ in $c(z)^2$ is $1/3$.

  \medskip
  \noindent\textit{Step 2.~}We justify the moment expansions and their strict
  difference.
  For completeness, these expansions admit the stated remainders after
  Gaussian expectation. The series of $c$ and $s$ have nonnegative
  coefficients, and
  \[
    |c(z)|+|s(z)|\le2e^{\sqrt{|z|}}.
  \]
  Consequently, the absolute series for the square of the displayed
  two-block displacement is bounded by
  $C(1+|z|)^2e^{2\sqrt{|z|}+2\sqrt{|w|}}$.
  This is integrable at any fixed Gaussian variance $A_1>0$. Absolute
  integration allows termwise expectation. Odd powers in either variable
  vanish; every remaining term has even total degree and a nonnegative
  coefficient. Thus, the expectation is a convergent series
  $\sum_{j\ge0}b_jA^j$ with $b_j\ge0$ at $A_1$. For $0\le A\le A_1$,
  \[
    \sum_{j\ge2}b_jA^j
    \le (A/A_1)^2\sum_{j\ge2}b_jA_1^j=O(A^2).
  \]
  The same argument applies to the one-block moment. Its square equals
  $1+2A/3+O(A^2)$, whereas the two-block moment equals
  $1+3A+O(A^2)$. Their difference is $7A/3+O(A^2)$.
  This proves No-go~\ref{P:nogo-14}.
\end{proof}

\begin{figure}
  \centering
  \begin{tikzpicture}[x=1cm,y=1cm]
    \node[ham/box,text width=1.8cm] (initial) at (0,0)
    {Initial phase\\$(1,0)$};
    \node[ham/box,text width=3cm] (middle) at (4.1,0)
    {After block 1\\$(c(Z_1),Z_1s(Z_1))$};
    \node[ham/box,text width=2.6cm] (final) at (8.7,0)
    {After block 2\\$M(Z_2)M(Z_1)e_1$};
    \draw[ham/arrow] (initial)--node[above,font=\small] {$M(Z_1)$} (middle);
    \draw[ham/arrow] (middle)--node[above,font=\small] {$M(Z_2)$} (final);
    \node[ham/note,text width=10cm] at (4.4,-1.7)
    {Retained displacement: $c(Z_2)c(Z_1)+s(Z_2)Z_1s(Z_1)$.\\
    Resetting the velocity would retain only $c(Z_2)c(Z_1)$.};
    \node[ham/note,text width=10cm] at (4.4,-3)
    {True second moment: $1+3A+O(A^2)$.\\
    Product of reset-block moments: $1+\frac23A+O(A^2)$.};
  \end{tikzpicture}
  \caption[Velocity survives a block boundary]{Two independent scalar static
    Gaussian blocks, each of variance $A$. The second block acts on both
    components of the first block's phase. This is a scalar diagnostic
    outside the target homogeneous-noise model; see
    No-go~\ref{P:nogo-14}.}
  \label{Fig:nogo-phase-blocks}
\end{figure}

This scalar model has spatially constant noise and independent static
blocks, outside Hypothesis~\ref{H:homog}. It excludes the unrestricted
scalar multiplication rule, not a phase-space operator estimate. Corollary~\ref{C:static-interaction} gives a valid replacement for independent
successive blocks: it controls both phase coordinates and includes explicit
block losses whose dyadic sum is finite. Its stated scope does not identify
products across overlapping levels with the original evolution.

\subsection{NG-15: the missing second-order mean term}\label{S:nogo-15}

Keeping the first random Duhamel term gives a nonnegative squared norm.
It is not the complete coefficient of the variance in the expected phase
energy. The mean of the second Duhamel term contributes at the same order.

\setcounter{nogo}{14}
\begin{nogo}[The first variance coefficient need not be nonnegative]
  \label{P:nogo-15}
  For the scalar phase propagator $M$ in No-go~\ref{P:nogo-14} and a
  standard Gaussian $Z$, write
  \[
    \E\|M(\sigma Z)x\|^2=\|M(0)x\|^2+\sigma^2x^{\mathsf T}Kx
    +O(\sigma^4).
  \]
  Then
  \[
    K=\begin{pmatrix}4/3&4/5\\4/5&17/45\end{pmatrix},\qquad
    (1,-2)K(1,-2)^{\mathsf T}=-\frac{16}{45}.
  \]
  The matrix $K$ is the full first coefficient in the variance parameter,
  not a covariance matrix. Omitting the second-order mean term also
  underestimates this coefficient for $x=(1,0)^{\mathsf T}$.
\end{nogo}
\begin{proof}
  We compute the complete first variance matrix, including the second-order
  mean terms. We then use the integrable series remainder to obtain an actual
  phase whose expected energy decreases at small coupling.

  \medskip
  \noindent\textit{Step 1.~}We retain all terms of the first variance matrix.
  The coefficients of $M(z)=M_0+zM_1+z^2M_2+\cdots$ are
  \[
    M_0=\begin{pmatrix}1&1\\0&1\end{pmatrix},\quad
    M_1=\begin{pmatrix}1/2&1/6\\1&1/2\end{pmatrix},\quad
    M_2=\begin{pmatrix}1/24&1/120\\1/6&1/24\end{pmatrix}.
  \]
  Gaussian expectation of the squared phase norm gives
  \[
    K=M_1^{\mathsf T}M_1+M_0^{\mathsf T}M_2+M_2^{\mathsf T}M_0.
  \]
  Matrix multiplication gives the stated entries and negative value.
  For $e_1=(1,0)^{\mathsf T}$, the first random contribution is
  $\|M_1e_1\|^2=5/4$, whereas $e_1^{\mathsf T}Ke_1=4/3$.
  Thus, it is not valid even as an upper bound for the full coefficient.

  \medskip
  \noindent\textit{Step 2.~}We justify the negative small-coupling energy
  correction.
  The entire series are absolutely Gaussian integrable, including their
  squares, by the bound in No-go~\ref{P:nogo-14} with the additional
  polynomial factor from $zs(z)$. Absolute summation of the even-degree
  terms at a fixed nonzero $\sigma_1$ bounds the tail of degree at least
  four by $C\sigma^4$ for $|\sigma|\le|\sigma_1|$. Therefore, for
  $x=(1,-2)^{\mathsf T}$ the actual expected energy satisfies
  \[
    \E\|M(\sigma Z)x\|^2=5-\frac{16}{45}\sigma^2+O(\sigma^4).
  \]
  This proves No-go~\ref{P:nogo-15}.
\end{proof}

The example refutes arbitrary-phase covariance monotonicity for this
scalar flow. It says nothing against monotonicity of the original
solution-value moment with initial data $(1,0)$ in the homogeneous model.
Section~\ref{S:static-phase} retains both contributions to the coefficient.
The all-order estimate of Proposition~\ref{P:static-phase-all-orders} and
the interaction bound of Corollary~\ref{C:static-interaction} retain these
terms. Positivity of the first random variance does not control the cross
term with the mean.

\subsection{NG-16: replacing a tilted mean by a deterministic value}
\label{S:nogo-16}

The Gaussian density cost is exact. A loss may instead enter when the
random solution under the tilted law is replaced by a deterministic
surrogate. The following two calculations identify the loss without
asserting an obstruction to Gaussian shifts themselves.

\setcounter{nogo}{15}
\begin{nogo}[A deterministic surrogate can lose the annealed factor]
  \label{P:nogo-16}
  Let $A>0$, $c>0$ and $0<r<2$, and define
  \[
    J_c=\sup_{H\ge0}\{cAH^r-H^2/2\}.
  \]
  Then $J_c=c^{2/(2-r)}J_1$, whereas
  \[
    \sup_{H\ge0}\{2AH^r-H^2\}=2J_1<J_2.
  \]
  Moreover, for $Y=e^{\sigma Z}$, $Z$ standard Gaussian and $\sigma>0$,
  the full tilted-mean Cauchy--Schwarz bound recovers $\E Y^2$ exactly,
  while substituting $Y(h)=e^{\sigma h}$ loses a factor in its logarithm.
\end{nogo}
\begin{proof}
  Substitute $H=c^{1/(2-r)}K$ in the definition of $J_c$. Both terms
  acquire the factor $c^{2/(2-r)}$. The supremum $J_1$ is finite, since
  $r<2$, and positive, since for small positive $H$ the positive term
  dominates. Factoring two from the surrogate score proves the strict
  comparison because $2/(2-r)>1$.

  Let $Q_h$ be the Gaussian shift with density
  $L_h=e^{hZ-h^2/2}$ relative to the standard Gaussian law. Then
  $\E L_h^2=e^{h^2}$, so Cauchy--Schwarz gives
  \[
    \E Y^2\ge(\E_{Q_h}Y)^2e^{-h^2}.
  \]
  The exact tilted mean is $\E_{Q_h}Y=e^{\sigma^2/2+\sigma h}$.
  Maximizing its resulting exponent
  $\sigma^2+2\sigma h-h^2$ at $h=\sigma$ gives
  $e^{2\sigma^2}=\E Y^2$. Substitution of the deterministic value
  $Y(h)$ instead gives the maximum $e^{\sigma^2}$.
  This proves No-go~\ref{P:nogo-16}.
\end{proof}

These are generic scalar score calculations, not an SPDE no-go. They
identify which term is lost and show that retaining the full tilted mean
can avoid the loss. A deterministic surrogate remains useful if one
proves that its error is negligible on the intended exponential scale;
without such an estimate it cannot identify the annealed constant.

\subsection{NG-20: covariance profiles need not separate}
\label{S:nogo-20}

A fixed spatial profile times a scalar temporal factor is easier to
propagate than a general spacetime control. Such a product family is
legitimate as a restriction, but it does not represent every positive
control even in the white spatial model.

\setcounter{nogo}{19}
\begin{nogo}[A positive control need not have proportional covariance profiles]
  \label{P:nogo-20}
  Fix $0<a<1$ and one-dimensional spatial covariance
  $\Gamma=c_\Gamma\delta_0$, $c_\Gamma>0$. There is a nonnegative smooth
  compactly supported spacetime test $\psi$ for which
  \[
    h_\sigma(x)=c_\Gamma\int_0^1|\sigma-u|^{-a}\psi(u,x)\ud u,
    \qquad 0\le\sigma\le1,
  \]
  cannot be represented as $b(\sigma)U(x)$ with one fixed spatial $U$.
\end{nogo}
\begin{proof}
  Choose a nonzero nonnegative smooth $f_1$ supported in $(1/8,1/4)$,
  let $f_2(u)=f_1(1-u)$, and choose nonzero nonnegative smooth bumps
  $q_1,q_2$ supported respectively in $(-2,-1)$ and $(1,2)$. Set
  $\psi(u,x)=f_1(u)q_1(x)+f_2(u)q_2(x)$. Put
  \[
    A=\int_0^1u^{-a}f_1(u)\ud u,\qquad
    B=\int_0^1(1-u)^{-a}f_1(u)\ud u.
  \]
  Since $u<1-u$ on the support of $f_1$, we have $A>B>0$.
  Reflection in time gives
  \[
    h_0=c_\Gamma(Aq_1+Bq_2),\qquad
    h_1=c_\Gamma(Bq_1+Aq_2).
  \]
  Choose $x_i$ with $q_i(x_i)>0$. Disjoint supports imply
  \[
    h_0(x_1)h_1(x_2)-h_1(x_1)h_0(x_2)
    =c_\Gamma^2(A^2-B^2)q_1(x_1)q_2(x_2)>0.
  \]
  A representation by one $U$ would make this determinant zero.
  The white covariance has $\alpha=1$, so $1+a<2$ verifies the model
  parameter condition. This proves No-go~\ref{P:nogo-20}.
\end{proof}

This is an exact representation obstruction. It does not prove that
separable controls fail to attain a variational supremum; a theorem
reducing optimizers to a smaller family would be a different statement.
A useful weaker fact survives: because $|\sigma-u|^{-a}\ge1$ on the
unit time square, every nonnegative test has
$h_\sigma\ge\Gamma*\int_0^1\psi(u,\cdot)\ud u$.
This common lower profile does not by itself prove a phase-propagation
comparison. Reopening an exact product reduction requires extra structure
on the control, while an optimizer reduction requires a variational proof.

\subsection{NG-31: initial phase and the backward cone}\label{S:nogo-31}

Two seemingly small generalizations remove essential hypotheses: replacing
the specified initial phase by an arbitrary one, and checking agreement
of spatial cutoffs only on the terminal observation set. Simple
deterministic waves already distinguish these statements.

\setcounter{nogo}{30}
\begin{nogo}[Growth does not give every-phase lower bounds or local cutoff consistency]
  \label{P:nogo-31}
  The growth rate of a wave phase propagator need not be a lower growth
  rate for every nonzero initial phase. Also, two smooth initial-data
  masks that agree on an observation interval can give different free-wave
  displacements throughout that interval at a later time, even when their
  initial displacements agree everywhere.
\end{nogo}
\begin{proof}
  For the first assertion, let
  $B=\left(\begin{smallmatrix}0&1\\1&0\end{smallmatrix}\right)$.
  Its eigenvectors $(1,1)^{\mathsf T}$ and $(1,-1)^{\mathsf T}$ have
  eigenvalues $1$ and $-1$. Hence, $\|e^{tB}\|=e^t$ for $t\ge0$,
  but $e^{tB}(1,-1)^{\mathsf T}=e^{-t}(1,-1)^{\mathsf T}$.

  For the second assertion take the one-dimensional free wave with zero
  initial displacement and nonnegative smooth velocity $v$ supported in
  $(5,6)$, normalized by $\int v=1$. Choose smooth masks $\chi_s,\chi_l$
  between zero and one, with $\chi_s=1$ on $[-2,2]$ and support in
  $[-4,4]$, and $\chi_l=1$ on $[-14,14]$. Both masks agree on the
  observation interval $[-2,2]$. The masked velocities are nevertheless
  $\chi_sv=0$ and $\chi_lv=v$. D'Alembert's formula at time nine gives
  \[
    u_s(9,x)=0,\qquad
    u_l(9,x)=\frac12\int_{x-9}^{x+9}v(y)\ud y=\frac12,
    \qquad -2\le x\le2.
  \]
  Indeed, every interval $[x-9,x+9]$ in question contains $(5,6)$.
  This proves No-go~\ref{P:nogo-31}.
\end{proof}

The first example does not refute a terminal lower bound for the original
phase $(1,0)$. The second requires agreement of both displacement and
velocity on the full backward cone, not only on the observer. Finite
propagation with that guard survives. Thus, a successor must retain the
specified phase, or prove a suitable nonvanishing unstable component, and
must impose cutoff agreement on the relevant cone for both coordinates.

\section{Topology, localization and limits of variational arguments}
\label{S:nogo-topology}

A useful approximation must be small in the topology in which the next
argument uses it. The obstructions below distinguish a local operator from a
Hilbert--Schmidt operator, a fixed spatial cutoff from a growing one, and
pointwise convergence from convergence of variational suprema. Throughout,
$a=\alpha_0\in(0,1)$ and $\Lambda=(1-\Delta)^{1/2}$, except that the positive
number $L$ denotes the temporal feature cutoff in No-go~\ref{P:nogo-32}.

\subsection{The first driver cannot always be Hilbert--Schmidt}
\label{S:nogo-21}

The one-dimensional operator construction suggests trying the same
Hilbert--Schmidt space in higher dimensions. A compact spatial observation
removes infinite volume, but it does not remove the following high-frequency
obstruction. The relevant compatibility is with the given noise, through its
bilinear first-driver tests.

\setcounter{nogo}{20}
\begin{nogo}[A same-noise Hilbert--Schmidt obstruction]
  \label{P:nogo-21}
  Suppose that $d\in\{2,3\}$, $s\in\R$, and $0<T\le1$. Put
  $E_s=H^s(\R^d)\times H^{s-1}(\R^d)$ and let $S(t)$ be the massless wave
  phase group. For every nonzero covariance satisfying Hypothesis~\ref{H:homog},
  there are a nonnegative compact smooth physical coefficient mask $\zeta$ and
  a compact smooth observation cutoff $\chi$, equal to one on a sufficiently
  large ball, for which the localized Gaussian first-driver form has no
  compatible realization in $L^2(\Omega;\operatorname{HS}(E_s,E_s))$.
  Here compatibility means equality with the given $W$ for each fixed smooth
  phase pair, as specified in the proof.
\end{nogo}

\begin{proof}
  We first construct a deterministic coefficient response that is not
  Hilbert--Schmidt, even after the stated observation cutoff. We then realize
  that coefficient as an actual covariance mean. Compatibility with a
  mean-square Hilbert--Schmidt driver would make this response
  Hilbert--Schmidt, giving a contradiction.

  \medskip
  \noindent\textit{Step 1.~}We identify the deterministic response and its
  weighted Fourier kernel.
  First consider a deterministic coefficient $h(u,x)=k(u)H(x)$, where
  $H\in C_c^\infty(\R^d)$ is nonzero, even and nonnegative, and
  $k\in W^{1,\infty}[0,T]$ satisfies $k\ge k_0>0$. Write
  $B(b)(v,w)=(0,bv)$ and
  \[
    A_h=\int_0^T S(-u)B(h(u))S(u)\ud u,\qquad
    S(u,\xi)=\begin{pmatrix}
      \cos(u|\xi|)&\sin(u|\xi|)/|\xi|\\
      -|\xi|\sin(u|\xi|)&\cos(u|\xi|)
    \end{pmatrix}.
  \]
  The sine quotient has its continuous value at zero. The integral is a strong
  operator integral on $E_s$: $S(u)$ is strongly continuous and uniformly
  bounded on compact time intervals, and multiplication by $H$ is bounded on
  $H^s$. For the latter assertion, the weighted Fourier kernel is bounded in
  absolute value by a constant times
  $\langle\xi-\eta\rangle^{|s|}|\widehat H(\xi-\eta)|$, an integrable
  convolution kernel. This proves boundedness for every fixed real $s$.

  Normalize the Fourier coordinates by
  $\operatorname{diag}(\langle\xi\rangle^s,\langle\xi\rangle^{s-1})$.
  Using the unitary Fourier transform, the lower-left block of $A_h$ then has
  kernel
  \[
    K_s(\xi,\eta)=(2\pi)^{-d/2}
    \langle\xi\rangle^{s-1}\langle\eta\rangle^{-s}
    \widehat H(\xi-\eta)
    \int_0^T k(u)\cos(u|\xi|)\cos(u|\eta|)\ud u.
  \]
  This formula follows from both phase coordinates: before applying the weights,
  the full product, with $r=|\xi|$, $\rho=|\eta|$, is
  \[
    \widehat H(\xi-\eta)
    \begin{pmatrix}
      -\dfrac{\sin(ru)}r\cos(\rho u)&
      -\dfrac{\sin(ru)}r\dfrac{\sin(\rho u)}\rho\\
      \cos(ru)\cos(\rho u)&\cos(ru)\dfrac{\sin(\rho u)}\rho
    \end{pmatrix}.
  \]
  All Fourier and time exchanges are justified on Schwartz pairs by the rapid
  decay of $\widehat H$ and then by density.

  \medskip
  \noindent\textit{Step 2.~}We prove the high-frequency Hilbert--Schmidt
  divergence.
  Let $J_0=\int_0^T k>0$ and
  $B_0=\|k\|_\infty+T\|k'\|_\infty$. Choose $\epsilon>0$ so that
  $\epsilon\le1$, $\epsilon T\le1$ and
  $\widehat H(z)\ge\widehat H(0)/2>0$ when $|z|\le\epsilon$.
  If $|\xi-\eta|\le\epsilon$, then
  \[
    \int_0^T k(u)\cos((r-\rho)u)\ud u\ge J_0/2,
    \qquad
    \left|\int_0^T k(u)\cos((r+\rho)u)\ud u\right|
    \le B_0/(r+\rho).
  \]
  The second bound is integration by parts in the absolutely continuous
  function $k$. Product-to-sum therefore bounds the time integral in $K_s$
  from below by $J_0/8$ when
  $r\ge R_0:=\max(1,2\epsilon,4B_0/J_0)$.
  On the same region,
  \[
    \langle\xi\rangle^{s-1}\langle\eta\rangle^{-s}
    =\langle\xi\rangle^{-1}
    (\langle\xi\rangle/\langle\eta\rangle)^s
    \ge 2^{-|s|}/(\sqrt2 r).
  \]
  Consequently, $|K_s(\xi,\eta)|\ge c_s/|\xi|$ there, for $c_s>0$.
  If $P_N$ cuts off Fourier frequencies to $B_N$, its truncated kernel is
  square integrable and Parseval followed by Tonelli gives
  \[
    \|P_NK_sP_N\|_{\operatorname{HS}}^2
    \ge c_s^2|B_\epsilon||\mathbb S^{d-1}|
    \int_{R_0}^{N-\epsilon}r^{d-3}\ud r\longrightarrow\infty.
  \]
  Indeed Parseval in the second kernel variable identifies its square integral
  with the sum of squared images of an orthonormal basis. If $A_h$ were
  Hilbert--Schmidt, every block compression would have bounded
  Hilbert--Schmidt norm. The displayed divergence proves otherwise. Real
  Hilbert--Schmidt operators and their complexifications have the same norm,
  so Fourier complexification has not changed the conclusion.

  \medskip
  \noindent\textit{Step 3.~}We retain the obstruction after compact
  observation.
  Suppose that $\operatorname{supp}H\subset B_A$, and let $\chi=1$ on
  $B_{A+T+1}$. With $D_\chi$ denoting multiplication of both phase coordinates,
  finite propagation gives
  \[
    D_\chi A_h=A_h=A_hD_\chi.
  \]
  For completeness, the local wave energy on a shrinking ball has boundary
  flux $v_t\partial_n v-(v_t^2+|\nabla v|^2)/2\le0$. Zero initial data
  therefore remain zero in the interior cone; time reversal gives the same
  support assertion for $S(-u)$. Apply this first to smooth data and then use
  Sobolev continuity and distributional testing to pass to $E_s$. The
  coefficient acts only in $B_A$, and the resulting output lies in
  $B_{A+T}$. Thus the failure persists after this compact observation.

  \medskip
  \noindent\textit{Step 4.~}We realize the coefficient as a covariance mean of
  the given noise.
  We next produce $h$ as an actual covariance mean. Choose $L_0\ge1$ with
  $0<\Gamma(B_{L_0})<\infty$, and an even nonnegative $q\in C_c^\infty$
  equal to at least one on $B_{2L_0+1}$. Set $U=\Gamma*q$. Local finiteness
  of $\Gamma$ permits differentiation on every compact set, so $U$ is smooth,
  even and nonnegative, with $U(0)>0$. Moreover,
  \[
    Q=\int q(x)U(x)\ud x
    \ge |B_{L_0+1}|\Gamma(B_{L_0})>0.
  \]
  Choose an even nonnegative compact smooth $\zeta$, equal to one near zero,
  and put $H=\zeta U$. Choose $f\in C_c^\infty(0,1)$ nonnegative and
  nonzero, and set
  \[
    Z=W(f\otimes q),\qquad
    k(u)=\int_0^1|u-v|^{-a}f(v)\ud v.
  \]
  For $u\in[0,1]$, $k(u)\ge\int f>0$. Moving the derivative to the smooth
  zero extension of $f$ shows $k\in C^1[0,1]$, since $|\cdot|^{-a}$ is
  locally integrable. Also $\E Z^2=Q\int_0^1 f(u)k(u)\ud u\in(0,\infty)$.
  The mask is applied after the covariance convolution: $H=\zeta(\Gamma*q)$.

  \medskip
  \noindent\textit{Step 5.~}We state the compatible rows and derive the
  Hilbert--Schmidt contradiction.
  It remains to specify the compatible first-driver form with the correct
  Sobolev adjoint. For real Schwartz phases $z,w$ put
  \[
    v_u=[S(u)z]_1,\qquad
    G_w(u)=-\frac{\sin(u|D|)}{|D|}(\chi\Lambda^{2s}w_1)
    +\cos(u|D|)(\chi\Lambda^{2s-2}w_2),
  \]
  and define
  $F_{z,w}(u,x)=\mathbf1_{(0,T)}(u)\zeta(x)v_u(x)G_w(u,x)$.
  This test lies in the raw noise completion. In fact, on a fixed compact
  spatial support $K$, its covariance norm is bounded by a constant times
  $\|u\mapsto\|F(u)\|_\infty\|_{L^2(0,T)}$, by local Radon mass
  $\Gamma(K-K)$ and the time Schur bound for $|u-v|^{-a}$. Smooth interior
  time approximation gives membership. Define
  $\ell_\chi(z,w)=W(F_{z,w})$ on the given noise space.
  Direct transposition, with the Sobolev weights retained, gives
  \[
    G_w(u)=\Lambda^{2s-2}
    \big((D_\chi S(-u))^{*_{E_s}}w\big)_2.
  \]
  In particular, multiplication by $\chi$ has not been treated as self-adjoint
  in $E_s$. Absolute Fubini in the covariance pairing now gives
  \[
    \E[Z\ell_\chi(z,w)]
    =\int_0^T\!\int k(u)H(x)v_u(x)G_w(u,x)\ud x\ud u
    =\langle D_\chi A_hz,w\rangle_{E_s}
    =\langle A_hz,w\rangle_{E_s}.
  \]
  All spatial factors in the raw tests have compact support and the temporal
  singularity is integrable, which supplies the absolute Fubini bound.

  If $X\in L^2(\Omega;\operatorname{HS}(E_s,E_s))$ satisfied
  $\langle Xz,w\rangle_{E_s}=\ell_\chi(z,w)$ almost surely for every fixed
  Schwartz pair, then Cauchy--Schwarz would make
  $Y=\E[ZX]$ a Hilbert--Schmidt Bochner integral. Continuous matrix evaluation
  would give $\langle Yz,w\rangle=\langle A_hz,w\rangle$ for all such pairs.
  Density implies $Y=A_h$, a contradiction. This proves No-go~\ref{P:nogo-21}.
\end{proof}

This result concerns the same-index phase spaces $E_s$ and a compatible
$L^2$ Hilbert--Schmidt realization. It does not concern every cutoff: the zero
cutoff, for example, has no such obstruction. It is not a nonexistence result
for the Stratonovich solution or its moments. A route using loss of regularity
between input and output, operator norms, or a different Schatten class must
supply its own estimates; changing the fixed index $s$ alone cannot reopen
this particular route.

\subsection{A physical mask may not act on the noise quotient}
\label{S:nogo-22}

A smooth physical mask is harmless on ordinary test functions. It need not
preserve the zero classes of a covariance seminorm. The following is the
obstruction already proved in Proposition~\ref{P:masked-quotient-obstruction};
we include its witness to explain why the auxiliary completion is needed.

\setcounter{nogo}{21}
\begin{nogo}[A mask that separates noise-null representatives]
  \label{P:nogo-22}
  For $d=2$, $\alpha=1$, $0<a<1$ and
  $\Gamma=\delta_0(\ud x_1)\ud x_2$, there exist real compact smooth $F$ and
  a nonnegative smooth square mask $\zeta$ for which
  $\E W(F)^2=0$ but $\E W(\zeta F)^2>0$.
  Thus multiplication by this mask does not descend to the original noise
  quotient, and no finite uniform covariance-norm multiplier bound holds.
\end{nogo}
\begin{proof}
  Choose a nonnegative nonconstant compact smooth plateau function $\psi$ and
  put $\chi=\psi^2$. The spatial covariance seminorm is
  \[
    V(g)=\int_\R\left(\int_\R g(x_1,x_2)\ud x_2\right)^2\ud x_1.
  \]
  This also shows positive definiteness of the nonnegative even tempered
  measure $\Gamma$; its scaling is $\Gamma(cA)=c\Gamma(A)$, as required
  for $d=2$, $\alpha=1$. Put
  $g(x)=\chi(x_1)\chi''(x_2)$ and
  $\zeta(x)=\chi(x_1)\chi(x_2)=(\psi(x_1)\psi(x_2))^2$.
  Then $\int\chi''=0$, whereas integration by parts gives
  \[
    \int\zeta(x)g(x)\ud x_2=-\chi(x_1)^2\int(\chi')^2.
  \]
  Consequently, $V(g)=0$ and
  $V(\zeta g)=\int\chi^4\,\big(\int(\chi')^2\big)^2>0$.
  Tensor with a nonzero nonnegative $\tau\in C_c^\infty(0,T)$. Its temporal
  variance is finite and strictly positive, so $F=\tau\otimes g$ has the
  asserted full-noise variances. This proves No-go~\ref{P:nogo-22}.
\end{proof}

This covariance is not two-dimensional spatial white noise. One can still
mask individual tests before taking an auxiliary completion. To reopen a
multiplier argument on the original quotient, one must first prove that the
chosen mask preserves its null space and then prove the required norm bound.

\subsection{Frequency decay does not give whole-space compactness}
\label{S:nogo-24}

The compact-cutoff estimates use spatial confinement as well as frequency
smoothing. Removing the former loses compactness even for a deterministic
Bessel multiplier; see Proposition~\ref{P:cutoff-boundary-tests}.

\setcounter{nogo}{23}
\begin{nogo}[Translations defeat whole-space compactness]
  \label{P:nogo-24}
  For every $d\ge1$ and $\gamma>0$, the operator $\Lambda^{-\gamma}$ on
  $L^2(\R^d)$ is not compact.
\end{nogo}
\begin{proof}
  Take a nonzero Schwartz function $f$ with Gaussian Fourier transform and put
  $g=\Lambda^{-\gamma}f$. Then $g$ is nonzero and Schwartz. For translations
  $T_hg(x)=g(x-h)$, Plancherel gives
  \[
    \langle T_hg,g\rangle
    =\int e^{-i\xi\cdot h}|\widehat g(\xi)|^2\ud\xi\longrightarrow0
    \quad (|h|\longrightarrow\infty).
  \]
  The limit follows, for example, by integration by parts in the Fourier
  transform of the Schwartz function $|\widehat g|^2$.
  Choose a sufficiently large spacing $b$ and $h_n=nb e_1$ so that distinct
  images $T_{h_n}g$ have squared distance at least $\|g\|_2^2$.
  The inputs $T_{h_n}f$ have constant norm and
  $\Lambda^{-\gamma}T_{h_n}f=T_{h_n}g$. Their images have no Cauchy
  subsequence. This proves No-go~\ref{P:nogo-24}.
\end{proof}

Compact spatial cutoffs and the conditional Schatten estimates in the main
text survive. No necessity at a Schatten endpoint is asserted here. A
whole-space successor needs a tightness condition that prevents translation
of mass, or an explicitly different topology.

\subsection{Fixed projections do not control growing spatial masks}
\label{S:nogo-25}

At each fixed radius the Cameron mean map can be compact. Its finite-rank
approximations need not be uniform in that radius. The following witness uses
actual unit controls for one-dimensional spatial white noise. The fixed-radius
map and its compact approximation are those of
Proposition~\ref{P:cameron-full-mean}.

\setcounter{nogo}{24}
\begin{nogo}[Escape from every fixed smooth projection]
  \label{P:nogo-25}
  Let $\Gamma=\delta_0$ in $d=1$, and let $H_I$ be the closed raw-noise
  space generated by tests supported in $I\times\R$, $I=(0,1)$. Let $P_J$
  be any fixed finite-rank orthogonal projection in $H_I$ whose range is
  spanned by smooth compact tests supported in $I\times\R$. For the physical masked mean multiplier $T_R$ on
  $F=(H^{1/2}\times H^{-1/2})\oplus\R^2$, with the last coordinates
  representing constant displacement and velocity, one has, for every
  $2<r<2/a$,
  \[
    \sup_R\|T_R(I-P_J)\|_{H_I\to L^r((0,1);\mathcal K(F))}\ge c_a>0,
  \]
  where $c_a$ can be chosen independently of $J$. Here the smooth masks $m_R$
  equal one on $[-R,R]$. For a smooth control $\Psi$ with covariance mean
  $h_\Psi$, the operator sends $(v,w,c,d)$ to
  $(0,m_Rh_\Psi(s)(v+c),0,0)$.
\end{nogo}
\begin{proof}
  We translate unit controls beyond the supports of a fixed projection while
  enlarging the physical mask to contain them. A unit constant-displacement
  input then gives a nonvanishing operator lower bound, showing why
  fixed-radius compactness is not uniform in radius.

  \medskip
  \noindent\textit{Step 1.~}We construct unit controls killed by the fixed
  projection.
  Use the temporal equilibrium density $f_a$ of
  Lemma~\ref{L:temporal-equilibrium}, so that
  $\int_0^1|s-t|^{-a}f_a(t)\ud t=1$ and
  $\|f_a\|_{H_a}^2=M_a>0$.
  Choose nonnegative $q\in C_c^\infty(\R)$ with $\|q\|_2=1$, and put
  $q_n(x)=q(x-x_n)$, where $x_n\to\infty$. The controls
  \[
    \Psi_n=\frac{f_a\otimes q_n}{\sqrt{M_a}}
  \]
  have unit raw norm and time-independent covariance mean
  $q_n/\sqrt{M_a}$. Membership follows by $L^2$ approximation of $f_a$ and
  the Schur bound for the temporal kernel.
  The finitely many tests spanning $P_J$ have common compact spatial support.
  For all sufficiently large $n$, spatial whiteness gives
  $P_J\Psi_n=0$ exactly. Choose $R_n$ to contain $\operatorname{supp}q_n$
  in the plateau of $m_{R_n}$.

  \medskip
  \noindent\textit{Step 2.~}We retain a positive response inside the growing
  mask.
  The constant-displacement input $e=(0,0,1,0)$ has norm one. The phase dual
  \[
    w_n=(0,\Lambda q_n/\|q\|_{H^{1/2}},0,0)
  \]
  also has norm one. Indeed its velocity norm is $H^{-1/2}$, and translation
  preserves $\|q\|_{H^{1/2}}$. The mean multiplier sends $e$ to the velocity
  $q_n/\sqrt{M_a}$, so the $H^{-1/2}$ inner product gives, for every $s$,
  \[
    \langle T_{R_n}\Psi_n(s)e,w_n\rangle_F
    =\frac{\|q_n\|_2^2}{\sqrt{M_a}\|q\|_{H^{1/2}}}
    =:c_a>0.
  \]
  Integrating over the unit time interval proves the bound for all $r$ stated.
  This proves No-go~\ref{P:nogo-25}.
\end{proof}

\begin{figure}
  \centering
  \begin{tikzpicture}[x=1cm,y=1cm]
    \draw[ham/arrow] (-.5,0)--(10.5,0) node[right] {$x$};
    \fill[hamblue!12] (0,.1) rectangle (2.4,1.1);
    \node[ham/note] at (1.2,1.95)
    {Fixed projection $P_J$\\common spatial support};
    \draw[hamblue,thick] (0,.1) rectangle (2.4,1.1);
    \draw[hamteal,thick] (6.5,.1) .. controls (7.1,.1) and (7.1,1.1) .. (7.5,1.1)
    .. controls (7.9,1.1) and (7.9,.1) .. (8.5,.1);
    \node[ham/note] at (7.5,1.65) {Translated bump $q_n$};
    \draw[ham/arrow] (4,.6)--(5.9,.6);
    \node[ham/note] at (4.95,1.15) {$x_n\to\infty$};
    \draw[hamgold,thick] (-.3,-.3)--(-.3,-.6)--(9.5,-.6)--(9.5,-.3);
    \node[ham/note] at (4.6,-1.1) {Growing plateau $m_{R_n}=1$ retains $q_n$};
    \node[ham/note] at (4.6,-1.8)
    {$P_J\Psi_n=0$,\qquad $\|T_{R_n}\Psi_n\|\ge c_a>0$};
  \end{tikzpicture}
  \caption[Escape from a fixed spatial projection]{Schematic for
    No-go~\ref{P:nogo-25}. For one-dimensional spatial white noise, a translated
    unit control is exactly orthogonal to the finitely many compact tests of
    $P_J$, while a growing mask retains its nonzero phase functional.
    At any fixed radius sufficiently distant bumps miss the mask.}
  \label{Fig:nogo-projection-escape}
\end{figure}
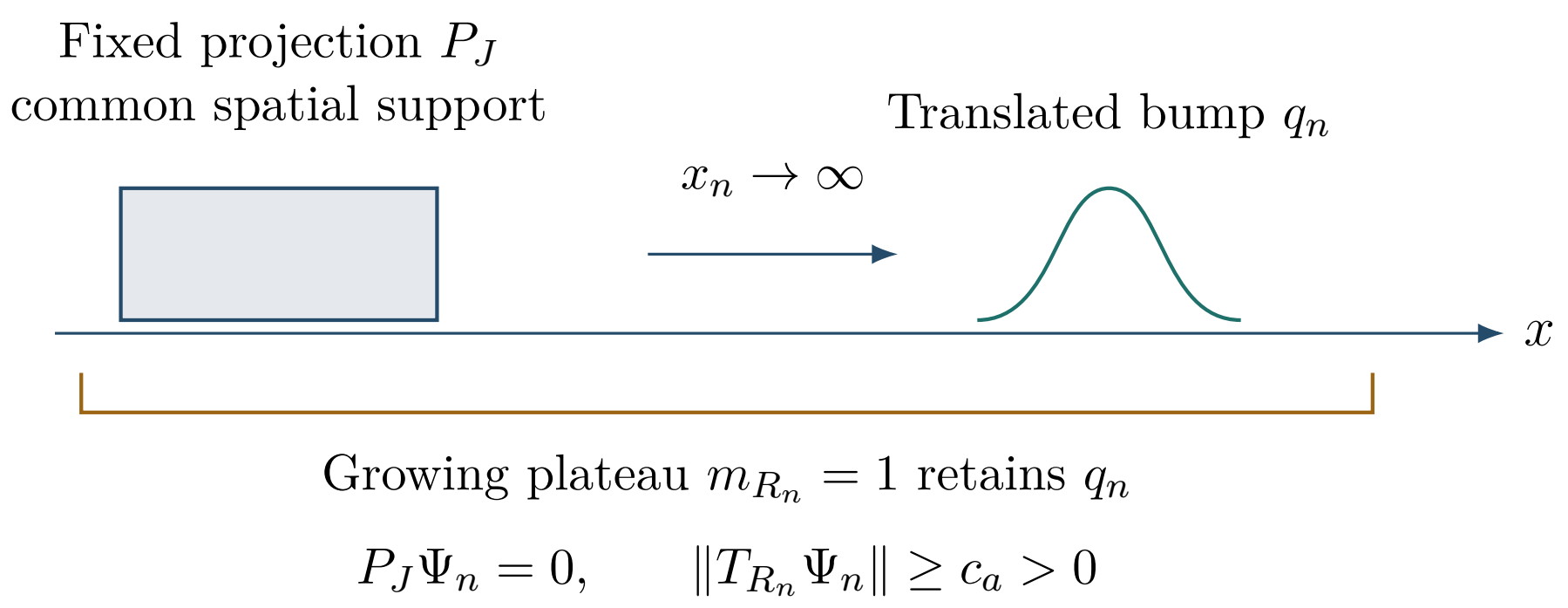

For a fixed $R$, sufficiently distant translates miss its mask, so there is
no contradiction with fixed-radius compactness. A replacement must use
projections depending on the radius with quantitative control, or a different
localization argument. The witness does not exclude every such strategy.

\subsection{Two distinct limitations of block costs}
\label{S:nogo-26}

A nonlinear block argument must handle both small perturbations and rare large
blocks. A power larger than one misses the first effect; an averaged budget
alone does not bound the largest block.

\setcounter{nogo}{25}
\begin{nogo}[A pure Young power and a maximum from an average]
  \label{P:nogo-26}
  On $[0,1]$, with $\theta=2/3$ and $\nu=1/\theta=3/2$, no finite constant
  $C>0$ gives
  $\log U(1)\le C\|D\|_{C^\theta}^{\nu}$ for all scalar paths starting at zero
  and solving $\ud U=U\ud D$, $U(0)=1$, when
  $\|D\|_{C^\theta}=\|D\|_\infty+[D]_{C^\theta}$.
  Moreover, a bound on $N^{-1}\sum_{j=1}^N b_j^r$, $r>0$, does not give a
  bound independent of $N$ on $\max_j b_j$.
\end{nogo}
\begin{proof}
  For $D(t)=xt$, $x>0$, the solution is $U(t)=e^{xt}$ and
  $\|D\|_{C^{2/3}}=2x$. Taking $x=1/(16C^2)$ gives
  $C(2x)^{3/2}=x/\sqrt2<x=\log U(1)$.
  For the second assertion, take $b_1=N^{1/r}$ and $b_j=0$ for $j>1$.
  Their averaged $r$th power is one while their maximum tends to infinity.
  This proves No-go~\ref{P:nogo-26}.
\end{proof}

A linear-plus-power cost is not refuted. For any $\tau>0$ and $\nu>1$,
$x\le\tau+\tau^{1-\nu}x^\nu$ gives an explicit alternative with slack.
Also, if $0<\nu<r$ and $N^{-1}\sum b_j^r\le M$, then H\"older gives
$N^{-1}\sum b_j^\nu\le M^{\nu/r}$, and
\[
  \frac1N\sum_{b_j>K}b_j^\nu\le K^{\nu-r}M.
\]
Thus an averaged tail argument survives. Reopening the block route requires
an actual evolution estimate with its linear term or slack, followed by an
averaged bad-block estimate. These scalar tests alone do not construct a
Cameron potential or prove a wave evolution bound.

\subsection{Concentrating temporal controls restrict integrability}
\label{S:nogo-27}

Compactness of a temporal mean map has a natural upper endpoint. A
concentrating control detects both excessive integrability and failure of
compactness at that endpoint.

\setcounter{nogo}{26}
\begin{nogo}[The temporal concentration threshold]
  \label{P:nogo-27}
  Let $I=(0,1)$, $0<a<1$, and let $H_a$ be the completion of smooth tests for
  $\|f\|_{H_a}^2=\int_I\int_I f(t)f(s)|t-s|^{-a}\ud t\ud s$.
  On smooth tests set $C_af(t)=\int_I|t-s|^{-a}f(s)\ud s$.
  The map $C_a$ is not bounded from $H_a$ to $L^r(I)$ when $r>2/a$,
  including $r=\infty$, and it is not compact into $L^{2/a}(I)$.
\end{nogo}
\begin{proof}
  We normalize a concentrating temporal test in the covariance norm. Its local
  lower bound disproves boundedness above the threshold; pointwise decay away
  from its center separately disproves compactness at the endpoint.

  \medskip
  \noindent\textit{Step 1.~}We compute the normalized concentration and its
  local norm lower bound.
  Choose $t_0\in I$ and a nonzero nonnegative
  $\phi\in C_c^\infty(-1,1)$. Put
  \[
    v_\phi=\int\!\int\phi(u)\phi(v)|u-v|^{-a}\ud u\ud v\in(0,\infty),
    \qquad
    f_\epsilon(t)=\frac{\epsilon^{a/2-1}}{\sqrt{v_\phi}}
    \phi((t-t_0)/\epsilon).
  \]
  For sufficiently small $\epsilon$, its support lies in $I$. Changing both
  variables shows $\|f_\epsilon\|_{H_a}=1$, and
  \[
    (C_af_\epsilon)(t_0+\epsilon y)
    =\frac{\epsilon^{-a/2}}{\sqrt{v_\phi}}
    \int|y-v|^{-a}\phi(v)\ud v.
  \]
  For $|y|\le1$, the last integral is at least $2^{-a}\int\phi>0$.
  Thus for finite $r$,
  $\|C_af_\epsilon\|_{L^r}\ge c_r\epsilon^{1/r-a/2}$, and the analogous
  supremum lower bound is $c\epsilon^{-a/2}$. This proves unboundedness above
  $2/a$.

  \medskip
  \noindent\textit{Step 2.~}We rule out a strongly convergent endpoint
  subsequence.
  For every fixed $t\ne t_0$, on the other hand,
  $C_af_\epsilon(t)=O_t(\epsilon^{a/2})\to0$.
  At $r=2/a$ the norm lower bound stays strictly positive. If a subsequence
  converged strongly in this $L^r$ space, a further subsequence would converge
  almost everywhere to its strong limit. The limit would have to be zero,
  contradicting the norm lower bound. This proves No-go~\ref{P:nogo-27}.
\end{proof}

The endpoint argument refutes compactness, not boundedness. Tensoring the
controls with one fixed nonzero spatial profile transfers this input
obstruction to a corresponding separated mean map. The corrected compact
range $2<r<2/a$ is consistent with the witness. A larger range requires
additional temporal regularity or restrictions that exclude concentration.

\subsection{Pointwise convergence does not pass through a supremum}
\label{S:nogo-30}

The Gaussian action cost must remain in a variational comparison. Even with
that cost, bounded near-maximizers can escape in infinitely many directions.

\setcounter{nogo}{29}
\begin{nogo}[Escaping variational scores]
  \label{P:nogo-30}
  There are continuous functions $J_n$ on the Hilbert space $\ell^2$ tending
  to zero at every point such that
  \[
    \sup_h\{J_n(h)-\tfrac12\|h\|^2\}\ge\tfrac12,
    \qquad
    \sup_h\{0-\tfrac12\|h\|^2\}=0.
  \]
  The witnesses for the first inequality have norm one.
\end{nogo}
\begin{proof}
  Let $(e_n)$ be the standard orthonormal sequence and set
  $J_n(h)=(1-2\|h-e_n\|)_+$.
  For each fixed $h\in\ell^2$, its $n$th coordinate tends to zero, and hence
  $\|h-e_n\|^2=\|h\|^2+1-2\langle h,e_n\rangle\to\|h\|^2+1\ge1$.
  Thus $J_n(h)=0$ eventually. But $J_n(e_n)=1$, so its score at $e_n$ is
  $1/2$. The limiting score is nonpositive and equals zero at $h=0$.
  This proves No-go~\ref{P:nogo-30}.
\end{proof}

A bounded action budget alone supplies no strong compactness in $\ell^2$.
To pass to an exact constant one needs a variational upper-limit argument,
for example compactness of near-maximizers in a topology with the needed
upper semicontinuity. Pointwise convergence, or compactness proved only with
an additional fixed spatial cutoff, does not supply that argument by itself.

\subsection{A cutoff point feature can lie outside the raw-noise image}
\label{S:nogo-32}

Lemma~\ref{L:source-retaining-isometry} embeds the original noise into one
fixed larger Gaussian space. Its closed image retains all original raw tests.
It does not contain every positive-cutoff point feature. The following
explicit orthogonal vector explains the role of the auxiliary noise.

\setcounter{nogo}{31}
\begin{nogo}[The point feature is outside the closed noise image]
  \label{P:nogo-32}
  Suppose that $\Gamma\ne0$ satisfies Hypothesis~\ref{H:homog}, $t_0>0$,
  $L\in(0,\infty)$ and $\epsilon>0$. In the feature space
  \[
    \mathcal G=L^2((0,\infty)_\lambda\times\R_r,
    \nu_a(\ud\lambda)\ud r;H_\Gamma),\qquad
    \nu_a(\ud\lambda)=\lambda^{a-1}\ud\lambda/\Gamma_E(a),
  \]
  let $K=\Phi H_W$ be the closed image of
  Lemma~\ref{L:source-retaining-isometry}. Set $g=p_\epsilon(x-\cdot)$ and
  \[
    v(\lambda,r)=\mathbf1_{\{\lambda\le L\}}\sqrt{2\lambda}
    e^{-\lambda(t_0-r)}\mathbf1_{\{r\le t_0\}}g.
  \]
  Then $v\notin K$. On the fixed extension
  $\mathcal I(v)=W(\Phi^{-1}Pv)+J((I-P)v)$, where $P$ projects onto $K$
  and $J$ is independent of $W$, one has
  $\operatorname{Var}(\mathcal I(v)\mid\sigma(W))>0$.
\end{nogo}
\begin{proof}
  We construct a feature-space vector orthogonal to every raw-noise test but
  not to the cutoff point feature. This separates the point feature from the
  closed noise image and gives a positive conditional variance on the fixed
  extension.

  \medskip
  \noindent\textit{Step 1.~}We construct the finite-norm witness.
  The heat profile has finite strictly positive norm:
  \[
    \|g\|_{H_\Gamma}^2=\int p_{2\epsilon}(y)\Gamma(\ud y)>0.
  \]
  Finiteness follows from temperedness and Gaussian decay, and positivity from
  nonzero nonnegative $\Gamma$. Integrating in $r$ also gives
  $\|v\|_{\mathcal G}^2=\|g\|^2\nu_a((0,L])<\infty$.
  Choose nonnegative nonzero smooth bumps $u_A,u_B$ supported in compact
  intervals $A\subset(0,L)$ and $B\subset(L,\infty)$, respectively, and put
  \[
    p(\lambda)=\frac{u_A(\lambda)}{\int u_A\ud\nu_a}
    -\frac{u_B(\lambda)}{\int u_B\ud\nu_a}.
  \]
  Then $\int p\ud\nu_a=0$ and $\int_{(0,L]}p\ud\nu_a=1$.
  Choose $\varphi\in C_c^\infty(0,\infty)$ with $\varphi(t_0)=1$, and define
  \[
    q(\lambda,r)=\frac{p(\lambda)}{\sqrt{2\lambda}}
    (\varphi'(r)+\lambda\varphi(r))g.
  \]
  Its support in $\lambda$ stays away from zero and infinity. In particular,
  \[
    \|q\|_{\mathcal G}^2=\frac{\|g\|^2}{2}
    \left(\|\varphi'\|_2^2\int\frac{p^2}{\lambda}\ud\nu_a
    +\|\varphi\|_2^2\int\lambda p^2\ud\nu_a\right)<\infty,
  \]
  where the cross term vanishes because $\int\varphi'\varphi=0$.

  \medskip
  \noindent\textit{Step 2.~}We prove orthogonality to the raw-noise image.
  For a smooth compact raw tensor $F$ the isometry is
  \[
    \Phi F(\lambda,r)=\sqrt{2\lambda}
    \int_{t>0,\ t\ge r}e^{-\lambda(t-r)}F_t\ud t.
  \]
  The identity
  \[
    \int_{-\infty}^t e^{-\lambda(t-r)}
    (\varphi'(r)+\lambda\varphi(r))\ud r=\varphi(t)
  \]
  follows by differentiating $e^{-\lambda(t-r)}\varphi(r)$ in $r$.
  All integrals paired with $q$ are absolutely integrable because its two
  supports are compact and its $\lambda$ support stays away from zero.
  Consequently,
  \[
    \langle\Phi F,q\rangle_{\mathcal G}
    =\left(\int_{t>0}\langle F_t,g\rangle_{H_\Gamma}
    \varphi(t)\ud t\right)\int p\ud\nu_a=0.
  \]

  \medskip
  \noindent\textit{Step 3.~}We separate the point feature and obtain the
  conditional variance obstruction.
  Density and continuity give $q\perp K$. The same calculation for the point
  feature instead gives
  \[
    \langle v,q\rangle_{\mathcal G}
    =\varphi(t_0)\|g\|^2\int_{(0,L]}p\ud\nu_a=\|g\|^2>0.
  \]
  Thus $(I-P)v\ne0$. Independence of $J$ and $W$ yields
  \[
    \operatorname{Var}(\mathcal I(v)\mid\sigma(W))
    =\|(I-P)v\|^2\ge\frac{|\langle v,q\rangle|^2}{\|q\|^2}>0.
  \]
  This proves No-go~\ref{P:nogo-32}.
\end{proof}

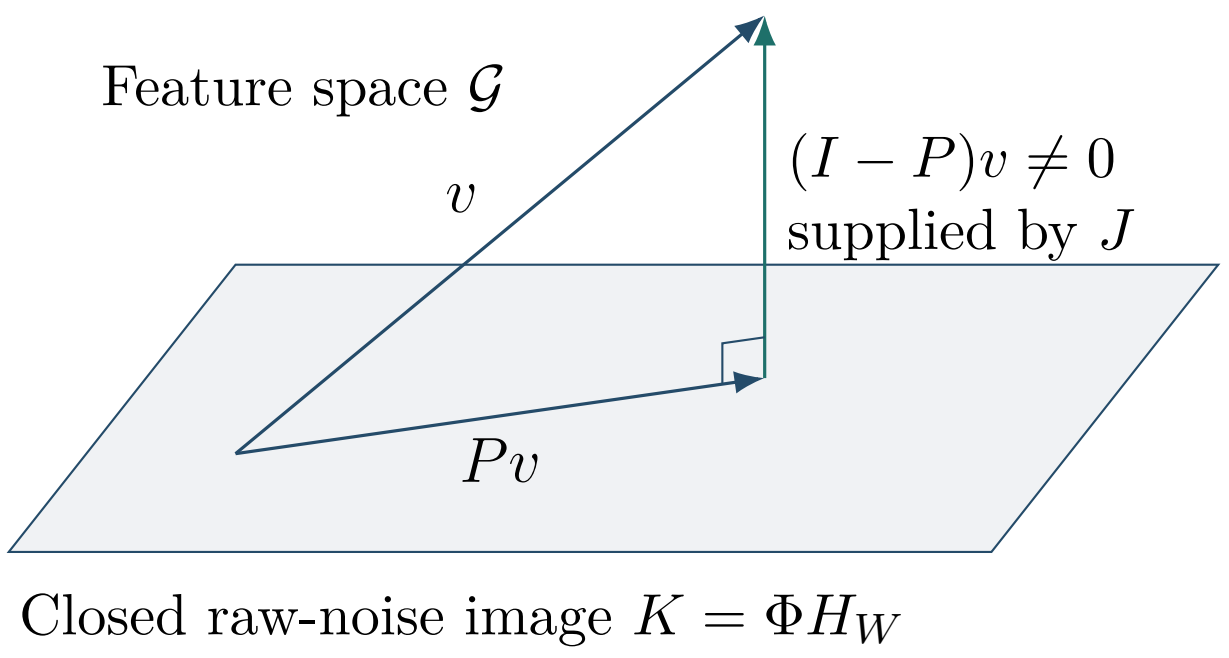
\begin{figure}
  \centering
  \begin{tikzpicture}[x=1cm,y=1cm]
    \filldraw[fill=hamblue!7,draw=hamblue] (-.5,-.5)--(6,-.5)--(7.5,1.4)--(1,1.4)--cycle;
    \coordinate (o) at (1,.15);
    \coordinate (p) at (4.5,.65);
    \coordinate (v) at (4.5,3.05);
    \draw[ham/arrow] (o)--(p) node[midway,below] {$Pv$};
    \draw[ham/arrow] (o)--(v) node[midway,above left] {$v$};
    \draw[ham/arrow,draw=hamteal] (p)--(v)
    node[ham/note,midway,right,align=left]
    {$(I-P)v\ne0$\\supplied by $J$};
    \draw[hamblue] (4.22,.61)--(4.22,.88)--(4.5,.92);
    \node[ham/note] at (2.5,-.95) {Closed raw-noise image $K=\Phi H_W$};
    \node[ham/note] at (1.45,2.55) {Feature space $\mathcal G$};
  \end{tikzpicture}
  \caption[The point feature and the closed noise image]{Hilbert-space
    schematic for No-go~\ref{P:nogo-32}, drawn in three coordinates only.
    The explicit vector $q\perp K$ has $\langle v,q\rangle>0$, so a
    positive finite cutoff point feature at $t_0>0$ has a nonzero component
    supplied by $J$. This says nothing by itself about measurability of a
    nonlinear limit with respect to $W$.}
  \label{Fig:nogo-feature-range}
\end{figure}

The conclusion concerns positive finite cutoffs and $t_0>0$ in this fixed
extension. It neither treats time zero nor proves that a nonlinear canonical
limit retains an auxiliary component. Raw-test convergence can remove such
components in a limit. The fixed-extension source construction remains valid.
For a concrete surviving limit, Proposition~\ref{P:constant-raw-average}
proves that the tested fixed-mask limit $Z_m(t,\phi)$ is $W$-measurable.
That conclusion uses the operator-path and Young-flow convergence; it is not
supplied or contradicted by this finite-feature obstruction.

\section{Finite optimization and checks on review arguments}
\label{S:nogo-optimization}

The last group concerns finite calculations. These examples are small enough
that every number can be checked directly, but the distinctions they expose
matter in the infinite-dimensional proof: the two variance directions, the
normalization of an energy, the active set of an optimization, and the exact
inequality used in a review.

\subsection{NG-23: both Gaussian variance directions}
\label{S:nogo-23}

A rectangular Gaussian matrix has an input variance and an output variance.
Keeping only one would simplify the moment estimate, but the dimensions of
the two spaces cannot be ignored. This is the counterexample behind
Remark~\ref{R:gaussian-one-variance-test}.

\setcounter{nogo}{22}
\begin{nogo}[A one-sided variance bound fails]\label{P:nogo-23}
  Let $g_1,\ldots,g_n$ be independent standard real Gaussians,
  $A_i=e_i^{\mathsf T}$ the coordinate row matrices, and
  $X=\sum_i g_iA_i$. Put $V_{\rm in}=\sum_iA_i^*A_i$ and
  $V_{\rm out}=\sum_iA_iA_i^*$. For $n=16$,
  \[
    \E\|X\|_{S_4}^4=288>144=9\operatorname{Tr}V_{\rm in}^2.
  \]
  The transposed family gives the analogous failure of an output-only bound.
  More generally, no constant independent of $n$ bounds the left side by
  $\operatorname{Tr}V_{\rm in}^2$ alone for all these row matrices.
\end{nogo}

\begin{proof}
  The only nonzero singular value of $X$ is
  $(\sum_i g_i^2)^{1/2}$. Since $\E g_i^4=3$ and
  $\E g_i^2g_j^2=1$ for $i\ne j$,
  \[
    \E\|X\|_{S_4}^4
    =\E\left(\sum_i g_i^2\right)^2
    =3n+n(n-1)=n(n+2).
  \]
  Here $V_{\rm in}=I_n$, whereas $V_{\rm out}=n$ acts on a
  one-dimensional space. The proposed bound gives $9n$, which is smaller
  than $n(n+2)$ at $n=16$. The ratio to $\operatorname{Tr}V_{\rm in}^2=n$
  is $n+2$, so no fixed constant works for all $n$. Transposition exchanges
  the variance directions and preserves the singular values. This proves
  No-go~\ref{P:nogo-23}.
\end{proof}

\paragraph{Scope and what survives.}
This is a finite-matrix counterexample, not a counterexample to a wave
estimate. Proposition~\ref{P:finite-gaussian-matrix} retains both variances.
At $n=16$, its average and maximum bounds have fourth powers $1224$ and
$2304$, respectively, and neither is contradicted by $288$.

\paragraph{What would permit a different estimate?}
A one-sided estimate would need additional structure relating the two
variance operators, with its dimension dependence stated. It cannot follow
from the Gaussian matrix form alone.

\subsection{NG-28: the normalization matters when bounds are compared}
\label{S:nogo-28}

The occupation argument gives the upper ceiling $U_2$ for $C_2$.
It is tempting to compare it with an energy expression by changing the
coefficient of the gradient penalty. The comparison must be calculated:
a value with penalty one is not interchangeable with a value with penalty
one-half. In particular, not every comparison value below is itself a proved
upper bound for the moment constant.

Retain $q=4-\alpha$, $k=3-\alpha$ and $m=q/k$, and the occupation energies
$\mathcal{E}_\zeta^{\mathcal O}$ of~\eqref{E:raw-occupation-definitions}.
For $z>0$, the scalar maps are
\[
  F_b(z)=k\left[
  \frac{(q-2b)^{(q-2b)/2}}{8(q-b)^{q-b}}
  \left(\frac{2z}{2-\alpha}\right)^{(2-\alpha)/2}
  \right]^{1/k},\qquad 0\le b<q/2.
\]
These are the maps in~\eqref{E:occupation-scalar-map}; here $b$ is only
a display parameter and $a=\alpha_0$ is the temporal exponent.

\setcounter{nogo}{27}
\begin{nogo}[An occupation ceiling need not improve a comparison value]
  \label{P:nogo-28}
  For $d=2$, $\alpha=5/4$, $a=1/2$ and every nonzero admissible covariance,
  \[
    2^mF_a(\mathcal{E}_1^{\mathcal O})<U_2.
  \]
  Thus, the reverse comparison is not valid throughout the admissible class.
\end{nogo}

\begin{proof}
  We compute the ratio of the two comparison formulas with the kinetic penalty
  retained. At the stated parameters, its logarithm is negative, disproving
  the proposed ordering of these formulas.

  \medskip
  \noindent\textit{Step 1.~}We express the comparison through its exact
  logarithmic ratio.
  The energies are finite and positive by
  Proposition~\ref{P:raw-occupation-energy}. Put
  \[
    \Delta=\frac12\{q\log q+(q-2a)\log(q-2a)\}
    -(q-a)\log(q-a).
  \]
  Subtracting the logarithms of $F_a$ and $F_0$ gives
  $\log(F_a(z)/F_0(z))=\Delta/k$. The identity in
  Theorem~\ref{T:intrinsic-occupation-ceiling} gives
  $U_2=2^m(2\zeta)^{\alpha/(2k)}F_0(\mathcal{E}_\zeta^{\mathcal O})$.
  Therefore,
  \[
    \frac{2^mF_a(\mathcal{E}_\zeta^{\mathcal O})}{U_2}
    =\exp\left\{\frac{\Delta-(\alpha/2)\log(2\zeta)}{k}\right\}.
  \]

  \medskip
  \noindent\textit{Step 2.~}We evaluate the sign at the admissible parameter
  point.
  At the stated parameters, $q=11/4$, $k=7/4$, and
  \[
    8\{\Delta-(\alpha/2)\log2\}
    =\log\frac{11^{11}7^7}{2^5 9^{18}}<0,
  \]
  because
  \[
    234966429149994773<4803028329503971872.
  \]
  This proves No-go~\ref{P:nogo-28}.
\end{proof}

\paragraph{Scope and what survives.}
This refutes an ordering of two computed quantities, not an inequality for
the actual moment. It does not decide whether $C_2$ is below the smaller
comparison value. At $\zeta=1/2$, the logarithmic difference is instead
$\Delta/k>0$, by strict convexity of $x\log x$; that is the valid strict
comparison in Corollary~\ref{C:intrinsic-spatial-bracket}.

\paragraph{What would permit a different comparison?}
Match the covariance amplitude, admissible profiles and kinetic penalty,
then prove the sign of the resulting ratio in the parameter range claimed.
Comparing the names of the energies is insufficient.

\subsection{NG-29: three safeguards in the Gram optimization}
\label{S:nogo-29}

The finite optimization of Section~\ref{S:finite-positive-gram} uses
\[
  F(\Phi)=2\sum_i\ell_i
  \sqrt{(\langle\Phi,k_i\rangle-b_i)_+}-\frac12\|\Phi\|^2,
  \qquad \ell_i>0,\quad b_i\ge0.
\]
The Gram entries $C_{ij}=\langle k_i,k_j\rangle$ are nonnegative. Three
possible simplifications fail: omitting the empty active set, invoking
global concavity, and treating entrywise order as positive semidefinite order.

\setcounter{nogo}{28}
\begin{nogo}[Empty sets, concavity and matrix order]\label{P:nogo-29}
  Each of the following failures occurs within the finite setting above.
  \begin{enumerate}
    \item With one unit vector, $\ell=1/10$ and $b=4$, the nonempty-set
      dual infimum is negative, while the primal supremum is nonnegative.
    \item With $H=\R$ and $k=\ell=b=1$, the primal objective is not concave.
    \item There are positive semidefinite Grams $C,D$ with
      $0\le C_{ij}\le D_{ij}$ for every entry but $D-C$ not positive semidefinite.
  \end{enumerate}
\end{nogo}

\begin{proof}
  In the first example the nonempty dual is
  \[
    D(w)=\frac{1}{100w}-4w+\frac{w^2}{2},\qquad w>0.
  \]
  Its infimum is at most $D(4)=-3199/400<0$, whereas $F(0)=0$.
  Dropping the empty value from the maximum of dual infima therefore cannot
  give the primal supremum.

  In the second example $F(x)=2\sqrt{(x-1)_+}-x^2/2$. Directly,
  $F(0)=F(2)=0$ but $F(1)=-1/2$, violating midpoint concavity.

  For the third example take
  \[
    C=\begin{pmatrix}1&0\\0&1\end{pmatrix},\qquad
    D=\begin{pmatrix}1&1/2\\1/2&1\end{pmatrix}.
  \]
  The eigenvalues of $D$ are $3/2$ and $1/2$, while those of $D-C$ are
  $1/2$ and $-1/2$. All entries are ordered as claimed. This proves
  No-go~\ref{P:nogo-29}.
\end{proof}

\paragraph{Scope and what survives.}
These are the checks in Remark~\ref{R:gram-controls}. They explain the
empty-set convention and active-set proof of
Proposition~\ref{P:gram-positive-optimizer}. Entrywise comparison still
works in Proposition~\ref{P:gram-entrywise-comparison}: for nonnegative
weights $w_i$, the form $\sum_{i,j}w_i(D_{ij}-C_{ij})w_j$ is nonnegative.
No order on all vectors, and no covariance monotonicity for nonlinear waves,
is asserted by that calculation.

\paragraph{What would permit a shortcut?}
A restricted concavity argument must first justify its active region.
A matrix comparison must specify the vectors on which its quadratic form
is nonnegative. Neither step follows from a convenient change of terminology.

\subsection{NG-33: check the mathematics returned by a reviewer}
\label{S:nogo-33}

A review report is also a mathematical argument and can contain an incorrect
step. This historical record groups three such checks. The manuscript had
the correct expressions; the examples below explain why the proposed review
steps could not replace them. These are not three new counterexamples to the
stochastic model.

\setcounter{nogo}{32}
\begin{nogo}[Three invalid proof steps]\label{P:nogo-33}
  The following assertions fail.
  \begin{enumerate}
    \item The free phase matrix
      $S_0(t)=\left(\begin{smallmatrix}1&t\\0&1\end{smallmatrix}\right)$
      is a contraction in the Euclidean norm for every $t\ge0$.
    \item For $x>0$, an integer $n\ge1$ and $0<\vartheta<1$, the series term
      $x^n/(n!)^\vartheta$ equals
      $A_nB_n$, where
      \[
        A_n=\frac{(2x^{1/\vartheta})^n}{n!},\qquad
        B_n=2^{-\vartheta n/(1-\vartheta)}.
      \]
    \item The denominator $k!\,k^{(2-\alpha)k}$ may always be replaced by
      the larger denominator $(k!)^3$ in the coefficient bound of
      Proposition~\ref{P:low-alpha-upper}.
  \end{enumerate}
\end{nogo}

\begin{proof}
  We test the three discarded review steps separately: free-phase
  contractivity, the unweighted series-factor identity, and the stronger
  factorial denominator. Each calculation concerns the stated inference rather
  than a new moment counterexample.

  \medskip
  \noindent\textit{Step 1.~}We disprove free-phase contractivity in the stated
  norm.
  The vector $(0,1)$ has norm one, but its image under $S_0(t)$ has norm
  $\sqrt{1+t^2}>1$ for $t>0$. To transfer this to the whole-space phase,
  use exactly the inhomogeneous Sobolev norm of
  Section~\ref{S:wave-variance}: $E=H^{1/2}\times H^{-1/2}$ and
  $J_0=\operatorname{diag}(\Lambda^{1/2},\Lambda^{-1/2})$, where
  $\Lambda=(1-\Delta)^{1/2}$. With $r=|\xi|$ and
  $\langle\xi\rangle=\sqrt{1+r^2}$, the Fourier symbol of
  $J_0S(t)J_0^{-1}$ is
  \[
    \begin{pmatrix}
      \cos(tr)&\langle\xi\rangle\sin(tr)/r\\
      -r\sin(tr)/\langle\xi\rangle&\cos(tr)
    \end{pmatrix}
    \longrightarrow S_0(t)\qquad(r\downarrow0).
  \]
  Choose normalized Fourier bumps supported in balls shrinking to zero
  and place them in the second coordinate of $J_0z$. Continuity of this
  matrix at zero shows that their image norms tend to $\sqrt{1+t^2}$.
  Thus $S(t)$ is not a contraction on this $E$. The distinction between
  $\langle\xi\rangle$ and $|\xi|$ matters: this is not a claim about the
  homogeneous energy norm. The valid estimate keeps a free-growth factor.

  \medskip
  \noindent\textit{Step 2.~}We test the proposed unweighted series product.
  The correct identity in the second assertion is
  $x^n/(n!)^\vartheta=A_n^\vartheta B_n^{1-\vartheta}$.
  At $\vartheta=3/4$, $n=1$ and $x=1$, its left side is one, whereas
  $A_1B_1=2\cdot(1/8)=1/4$. H\"older applies to the weighted product, as
  in the proof of Proposition~\ref{P:static-phase-all-orders}.

  \medskip
  \noindent\textit{Step 3.~}We test the proposed replacement of the factorial
  denominator.
  For the third assertion choose the permitted spatial-white value
  $\alpha=1$ and $k=3$. Then $k!k^{(2-\alpha)k}=162<216=(k!)^3$.
  For $M>0$, an upper bound $M/162$ does not imply the smaller upper bound
  $M/216$.
  The valid factorial weakening is instead
  \[
    k!k^{(2-\alpha)k}\ge(k!)^{3-\alpha},
  \]
  since $k!\le k^k$ and $2-\alpha>0$; at this example its right side is
  $36$. This proves No-go~\ref{P:nogo-33}.
\end{proof}

\paragraph{Scope and what survives.}
The last arithmetic check refutes an inference between numerical bounds;
it does not prove that the actual coefficient exceeds $M/216$.
All three corrected arguments remain in the main text. Their lesson is to
check an asserted derivation against its domains, factors and quantifiers,
even when the report containing it has a favorable overall verdict.

\paragraph{What would permit the discarded steps?}
Contractivity needs a suitable invariant norm or a restricted state class.
The H\"older factors need their conjugate powers. A stronger factorial
bound needs an additional estimate; it cannot be obtained by changing the
exponent in the denominator.

\appendix

\chapter{Finite-order calculations}\label{S:appendix-finite-order}

\begin{figure}
  \centering
  \begin{tikzpicture}[x=1cm,y=1cm]
    \foreach \hamPanelShift/\hamStripWidth/\hamPanelTitle in {0/.35/{Narrow strip},4.25/1.65/{Corners removed},8.5/3.15/{Full rectangle}} {
    \begin{scope}[xshift=\hamPanelShift cm]
      \node[ham/note] at (0,2.15) {\hamPanelTitle};
      \begin{scope}
        \clip (-1.55,-1.05) rectangle (1.55,1.05);
        \fill[hamteal!18] (-4,-4-\hamStripWidth)--(4,4-\hamStripWidth)--(4,4+\hamStripWidth)--(-4,-4+\hamStripWidth)--cycle;
        \draw[hamteal,densely dashed] (-4,-4-\hamStripWidth)--(4,4-\hamStripWidth);
        \draw[hamteal,densely dashed] (-4,-4+\hamStripWidth)--(4,4+\hamStripWidth);
      \end{scope}
      \draw[hamblue] (-1.55,-1.05) rectangle (1.55,1.05);
      \draw[ham/arrow] (-1.8,0)--(1.8,0) node[right] {$x$};
      \draw[ham/arrow] (0,-1.3)--(0,1.4) node[above] {$y$};
      \node[font=\small,anchor=north] at (-1.55,-1.1) {$-a$};
      \node[font=\small,anchor=north] at (1.55,-1.1) {$a$};
      \node[font=\small,anchor=east] at (-1.6,1.05) {$b$};
      \node[font=\small,anchor=east] at (-1.6,-1.05) {$-b$};
    \end{scope}
    }
    \node[font=\small] at (0,-1.95) {$0<c<a-b$};
    \node[font=\small] at (4.25,-1.95) {$a-b\le c<a+b$};
    \node[font=\small] at (8.5,-1.95) {$c\ge a+b$};
  \end{tikzpicture}
  \caption[The rectangle and strip in the fourth-order calculation]{The area
    $F(a,b,c)$ in the spatial-white order-four calculation: the shaded part is
    $[-a,a]\times[-b,b]\cap\{|x-y|<c\}$, and the dashed lines are the edges
    $|x-y|=c$ of the strip. The drawing uses $a>b>0$; here $a,b,c$
    are time-increment parameters, not covariance exponents. The three geometric
    cases explain the positive-part formula, without replacing the separate
    pairing and time-order calculations. See
    \ProofRef{Proposition}{P:white-four}.}
  \label{Fig:white-strip}
\end{figure}
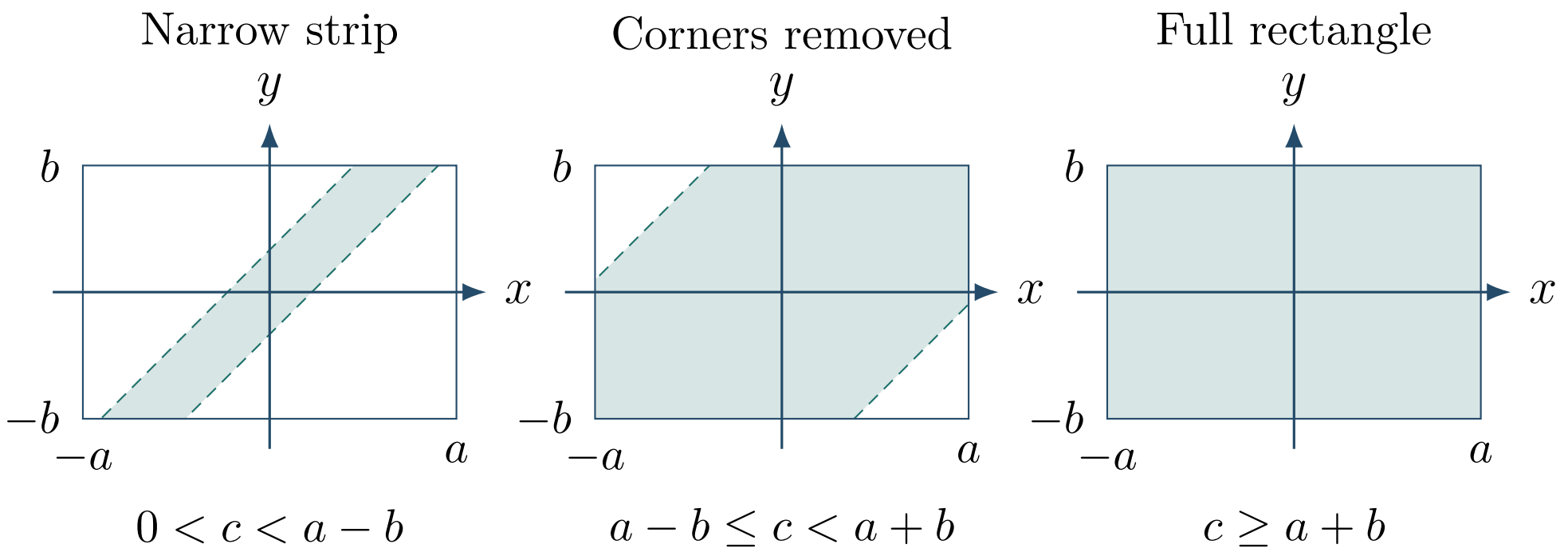

Before trusting an asymptotic argument, compute something by hand. What do the
first terms of the second moment look like? This appendix works them out in
one space dimension with initial data $(1,0)$. Figure~\ref{Fig:white-strip}
is the picture to keep in mind.

Already at small order the Stratonovich expansion shows its three
features: pairings inside one chain, pairings between the two chains,
and products of different orders. Starting from Proposition~\ref{P:mixed},
Proposition~\ref{P:white-four} computes every term of total order four for
spatial white noise, and Proposition~\ref{P:riesz-small} gives the order-two
coefficient and an integral formula for every order-four Riesz term. Here
the total order is the sum of the two Stratonovich orders; it is not a
Wiener-chaos order.

These formulas test the scaling, the pairing counts and the static
limit. Finitely many coefficients cannot determine the exponent or the
constant, and the main argument does not rely on this appendix. A reader
who wants examples first can start here after the chain notation of
Chapter~\ref{S:G0}. Episode~\VideoRef{9}{eLz9Hd20ZT0} of the video series
covers this appendix.

\section{Finite-order checks in one space dimension}\label{S:G1}

We now compute the total-order-two and total-order-four coefficients in
Proposition~\ref{P:mixed}, first for spatial white noise and then for $\gamma(x)
= c_\alpha|x|^{-\alpha}$ with $0 < \alpha < 1$ and $c_\alpha > 0$. Here $d = 1$
throughout. A total order is the sum of the two Stratonovich orders; it is not
an orthogonal Wiener-chaos order. Write
\[
  A=\{(1,2),(3,4)\},\quad B=\{(1,3),(2,4)\},\quad
  C=\{(1,4),(2,3)\}.
\]
These are the three pairings of four vertices. All finite integrals below retain
the normalization $G(s,y) = \frac12\mathbf1_{\{|y|<s\}}$ and covariance
$|s-r|^{-\alpha_0}\gamma(y-z)$.

\subsection{Spatial white noise: all order-four diagrams}

For $a,b,c \ge 0$, define
\begin{equation}\label{E:interval-area}
  F(a,b,c)=4ab-(a+b-c)_+^2+(|a-b|-c)_+^2.
\end{equation}
This is the area of $\{(x,y):|x|<a,\ |y|<b,\ |x-y|<c\}$. In particular, $F$ is
nonnegative and symmetric in its three arguments. To verify the formula,
subtract from the rectangle the two regions where $|x-y| \ge c$; when $c \ge
|a-b|$ these are two right triangles, and when $c < |a-b|$ the retained strip
has area $4\min(a,b)c$. The three formulas, including $c \ge a+b$,
give~\eqref{E:interval-area}. Symmetry also follows from swapping $x,y$ and the
determinant-one change $(x,y)\mapsto(x,x-y)$.

In the following table, $V$ is a spatial factor and $q_1,q_2$ are nonnegative
time gaps. The contribution of a row is the integral of $V(q_1q_2)^{-\alpha_0}$
over its domain, with ordinary Lebesgue measure in its four listed variables.
The domains are
\[
  \begin{array}{ll}
    (0,4):&a,b,c,d>0,\quad a+b+c+d<t,\\
    (1,3):&0<r<t,\quad a,b,c>0,\quad a+b+c<t,\\
    (2,2):&a,b,c,d>0,\quad a+b<t,\quad c+d<t.
  \end{array}
\]
The letters $a,b,c,d$ denote increments within the nonempty chains; for $(1,3)$,
$r$ is the time of the single vertex. Values on zero-gap sets do not affect
these nonnegative integrals.

\LeanMarker{P:white-four}
\begin{proposition}[White-space order four]\label{P:white-four}
  Suppose that $d = \alpha = 1$ and $\gamma = \delta_0$. Each $J_{m,n,D}(t)$
  with $m+n = 4$ is given by the following table or by exchanging the two
  chains:
  \begin{equation}\label{E:white-four-table}
    \begin{array}{c|c|c|c|c}
      (m,n)&D&V&q_1&q_2\\ \hline
      (0,4)&A&ac/4&b&d\\
      (0,4)&B&a\min(b,c,d)/4&b+c&c+d\\
      (0,4)&C&a\min(b,d)/4&b+c+d&c\\ \hline
      (1,3)&A&\min(r,a)b/4&|r-a|&c\\
      (1,3)&B&F(r,a,\min(b,c))/16&|r-a-b|&b+c\\
      (1,3)&C&F(r,a,c)/16&|r-a-b-c|&b\\ \hline
      (2,2)&A&ac/4&b&d\\
      (2,2)&B&\min(a,c)\min(b,d)/4&|a-c|&|a+b-c-d|\\
      (2,2)&C&F(a,c,\min(b,d))/16&|a-c-d|&|a+b-c|
    \end{array}
  \end{equation}
  The total-order-four contribution is
  \begin{equation}\label{E:total-four}
    2\sum_{D\in\{A,B,C\}}J_{0,4,D}(t)
    +2\sum_{D\in\{A,B,C\}}J_{1,3,D}(t)
    +\sum_{D\in\{A,B,C\}}J_{2,2,D}(t).
  \end{equation}
  Each entry scales as $t^{6-2\alpha_0}$ times its unit-time value.
\end{proposition}

\begin{proof}
  We identify the two free spatial positions for each pairing and compute the
  allowed area. We then exchange the chains and collect the time and space
  scaling factors.

  \medskip
  \noindent\textit{Step 1.~}We compute the spatial factors for each chain
  split.
  In each pairing identify its first pair's spatial positions with $x$ and the
  other pair's positions with $y$. The two delta covariances have unit mass, so
  spatial integration leaves $1/16$ times the area allowed by the four Green
  indicators. For a single chain the position lists for $A,B,C$ are,
  respectively, $(x,x,y,y)$, $(x,y,x,y)$ and $(x,y,y,x)$. Their areas are $4ac$,
  $4a\min(b,c,d)$ and $4a\min(b,d)$.

  For $(1,3)$ take the single position to be $x$. The long chain has positions
  $(x,y,y)$, $(y,x,y)$ or $(y,y,x)$ for $A,B,C$. The first area is
  $4\min(r,a)b$; the other two areas are $F(r,a,\min(b,c))$ and $F(r,a,c)$. For
  $(2,2)$ the position lists of the two chains are $(x,x),(y,y)$ for $A$,
  $(x,y),(x,y)$ for $B$, and $(x,y),(y,x)$ for $C$. The areas are $4ac$,
  $4\min(a,c)\min(b,d)$ and $F(a,c,\min(b,d))$. In every case reading the paired
  times gives the two gaps in~\eqref{E:white-four-table}. This proves each row,
  with no unlisted pairing or multiplicity. The distributional delta
  integrations can equivalently be obtained by convolution with the source's
  spatial approximate identity: for positive increments the remaining bounded
  interval indicators have only boundary discontinuities, and those boundaries
  have two-dimensional measure zero. The resulting areas converge to the
  displayed areas. They are uniformly bounded at unit time: each smoothed delta
  has mass one, so its integral over its paired coordinates in $[-1,1]^2$ is at
  most two. The four Green factors give a spatial bound $2^2/16 = 1/4$.
  Integrability of the time weights follows from the disjoint-pair bound used in
  Lemma~\ref{L:static-four} below.

  \medskip
  \noindent\textit{Step 2.~}We collect the multiplicities and the scaling.
  Exchanging chains gives the $(4,0)$ and $(3,1)$ sums without changing the
  covariance factors. This proves~\eqref{E:total-four}. Finally, scaling all
  four time variables by $t$ gives $t^4$ from their measure, $t^2$ from the area
  and $t^{-2\alpha_0}$ from the covariance factors. This proves
  Proposition~\ref{P:white-four}.
\end{proof}

\begin{lemma}[Static finite-order calibration]\label{L:static-four}
  As $\alpha_0\downarrow0$, each unit-time coefficient in
  Proposition~\ref{P:white-four} converges to the entry below:
  \begin{equation}\label{E:static-four-table}
    \begin{array}{c|ccc|c}
      (m,n)&A&B&C&\text{sum}\\ \hline
      (0,4)&1/2880&1/8640&1/5760&11/17280\\
      (1,3)&1/576&29/34560&1/640&143/34560\\
      (2,2)&1/144&1/480&31/17280&187/17280
    \end{array}
  \end{equation}
  Thus, the static total-order-four contribution is $11t^6/540$. The static
  total-order-two contribution is $t^3/3$. These are coefficients, not
  large-time asymptotic constants.
\end{lemma}

\begin{proof}
  We first justify the static limit by an integrable time-cube bound. We then
  integrate time slacks at fixed spatial positions and evaluate the remaining
  polynomial sector integrals.

  \medskip
  \noindent\textit{Step 1.~}We pass to the static coefficient by dominated
  convergence.
  On the unit-time domain all original time coordinates belong to $[0,1]$. The
  spatial factor is bounded by $1/4$, since all positions lie in $[-1,1]$ and
  the two free positions have area at most four. Fix $0 < \bar{a} < 1$. For $0
  \le \alpha_0 \le \bar{a}$, the temporal factors are bounded by the product
  with exponent $\bar{a}$. On the full time cube its integral is
  \[
    \left(\int_0^1\int_0^1|s-r|^{-\bar{a}}\ud s\ud r\right)^2
    =\left(\frac{2}{(1-\bar{a})(2-\bar{a})}\right)^2<\infty,
  \]
  because the two pairs use disjoint labels. Dominated convergence permits
  setting $\alpha_0 = 0$ in every row.

  \medskip
  \noindent\textit{Step 2.~}We reduce the coefficients to spatial sectors.
  We evaluate the resulting integrals by an independent spatial calculation. For
  fixed free positions $x,y$, let $L_1,L_2$ be the sums of absolute spatial
  increments along the two chains, starting at zero. Integrating their time
  increments first gives
  \begin{equation}\label{E:static-sector}
    J_{m,n,D}(1)\big|_{\alpha_0=0}
    =\frac{1}{16m!n!}\int_{\R^2}
    (1-L_1(x,y))_+^m(1-L_2(x,y))_+^n\ud x\ud y.
  \end{equation}
  An empty-chain factor in this formula is defined to be one. Indeed, writing
  each time increment as its required absolute spatial increment plus a
  nonnegative slack leaves a simplex of volume $(1-L_\rho)_+^{n_\rho}/n_\rho!$.
  Put $w = |x-y|$. The length pairs are
  \[
    \begin{array}{c|ccc}
      (m,n)&A&B&C\\ \hline
      (0,4)&(0,|x|+w)&(0,|x|+3w)&(0,|x|+2w)\\
      (1,3)&(|x|,|x|+w)&(|x|,|y|+2w)&(|x|,|y|+w)\\
      (2,2)&(|x|,|y|)&(|x|+w,|x|+w)&(|x|+w,|y|+w)
    \end{array}
  \]

  \medskip
  \noindent\textit{Step 3.~}We evaluate and sum the sector contributions.
  Divide the plane into $0 \le x \le y$, $0 \le y \le x$, $x \ge 0 \ge y$ and
  their reflected sectors. Use $(x,y) = (u,u+v),(u+v,u),(u,-v)$, respectively,
  where $u,v \ge 0$. Each map has determinant of absolute value one. In these
  three sectors the integrals of $(1-L_1)_+^m(1-L_2)_+^n$, before multiplication
  by $2/(16m!n!)$, are
  \begin{equation}\label{E:sector-values}
    \begin{array}{c|c|ccc}
      (m,n)&D&0\le x\le y&0\le y\le x&x\ge0\ge y\\ \hline
      (0,4)&A&1/30&1/60&1/60\\
      (0,4)&B&1/90&1/120&1/360\\
      (0,4)&C&1/60&1/90&1/180\\
      (1,3)&A&1/24&3/160&11/480\\
      (1,3)&B&1/72&3/160&11/1440\\
      (1,3)&C&1/48&1/30&1/48\\
      (2,2)&A&1/18&1/18&1/9\\
      (2,2)&B&1/30&1/60&1/60\\
      (2,2)&C&31/1440&31/1440&31/2160
    \end{array}
  \end{equation}
  All entries follow by integrating polynomials over the two linear constraints
  $L_1 \le 1$, $L_2 \le 1$. For example, $\int_{u,v\ge0}(1-pu-qv)_+^4\ud u\ud v
  = 1/(30pq)$. The last row's first entry and last entry are, explicitly,
  \[
    \begin{split}
      &\int_0^{1/2}\ud v\int_0^{1-2v}
      (1-u-v)^2(1-u-2v)^2\ud u=\frac{31}{1440},\\
      &\int_0^{1/3}\ud u\int_0^{(1-u)/2}
      (1-2u-v)^2(1-u-2v)^2\ud v\\
      &\quad+\int_{1/3}^{1/2}\ud u\int_0^{1-2u}
      (1-2u-v)^2(1-u-2v)^2\ud v=\frac{31}{2160}.
    \end{split}
  \]
  The stated lengths and linear bounds similarly specify every entry
  of~\eqref{E:sector-values}. Summing and multiplying by $2/(16m!n!)$
  proves~\eqref{E:static-four-table}. Applying the multiplicities
  in~\eqref{E:total-four} gives $11/540$. Lemma~\ref{L:white-test} at $\alpha_0
  = 0$ gives the order-two value. This proves Lemma~\ref{L:static-four}.
\end{proof}

\subsection{Riesz covariance: order two and an order-four integral}

\LeanMarker{P:riesz-small}
\begin{proposition}[Riesz finite coefficients]\label{P:riesz-small}
  Suppose $d = 1$, $0 < \alpha < 1$, $0 < \alpha_0 < 1$, and $\gamma(y) =
  c_\alpha|y|^{-\alpha}$ with $c_\alpha > 0$. With $\beta = 4-\alpha-\alpha_0$,
  the two order-two terms are
  \begin{equation}\label{E:riesz-two}
    \begin{split}
      \E X_2(t,x)
      &=\frac{c_\alpha t^\beta}
      {(1-\alpha)(2-\alpha-\alpha_0)(3-\alpha-\alpha_0)
      (4-\alpha-\alpha_0)},\\
      \E X_1(t,x)^2
      &=\frac{c_\alpha}{(1-\alpha)(2-\alpha)}
      \int_0^t v^{-\alpha_0}
      \left[\frac{(2t-v)^{3-\alpha}-v^{3-\alpha}}{2(3-\alpha)}
      -(t-v)v^{2-\alpha}\right]\ud v.
    \end{split}
  \end{equation}
  Their total is $\E X_1^2+2\E X_2$.

  The following explicit unit-domain integral gives every order-four term. For
  $m+n = 4$, let
  \[
    H_m=\{h\in(0,\infty)^4:\ \sum_{l=1}^m h_l<1,
    \ \sum_{l=m+1}^4 h_l<1\},
  \]
  with an empty sum equal to zero. For each index $l$, let $b(l)$ be the first
  index in its own nonempty chain, and define
  \[
    S_l(h)=\sum_{i=b(l)}^l h_i,\qquad
    Y_l(h,U)=\sum_{i=b(l)}^l h_iU_i,
  \]
  where $U_1,\ldots,U_4$ are independent uniform variables on $[-1,1]$. Then,
  for $D = A,B,C$,
  \begin{equation}\label{E:riesz-four}
    \begin{split}
      J_{m,n,D}(t)
      ={}&c_\alpha^2t^{2\beta}
      \int_{H_m}\left(\prod_{l=1}^4h_l\right)
      \prod_{(j,q)\in D}|S_j(h)-S_q(h)|^{-\alpha_0}\\
      &\hspace{8mm}\times
      \E_U\left[\prod_{(j,q)\in D}
      |Y_j(h,U)-Y_q(h,U)|^{-\alpha}\right]\ud h.
    \end{split}
  \end{equation}
  This is a finite nonnegative integral. It specifies all five splits and all
  three pairings without an unspecified kernel. The total coefficient uses the
  multiplicities in~\eqref{E:total-four}.
\end{proposition}

\begin{proof}
  We compute the two order-two coefficients directly. For order four, we
  express Green increments through uniform variables and check finiteness
  before using the resulting unit-domain formula.

  \medskip
  \noindent\textit{Step 1.~}We evaluate the order-two terms.
  Direct integration on the two intervals gives
  \[
    \int_{\R^2}G(s,y)G(r,z)c_\alpha|y-z|^{-\alpha}\ud y\ud z
    =\frac{c_\alpha\big((s+r)^{2-\alpha}-|s-r|^{2-\alpha}\big)}
    {2(1-\alpha)(2-\alpha)}.
  \]
  Integrate on $s > r$ and double, set $v = s-r$, and integrate $r$ over
  $(0,t-v)$. This proves the second formula in~\eqref{E:riesz-two}. In the
  single-chain term with increments $r,v$, the spatial integral is $c_\alpha
  rv^{1-\alpha}/(1-\alpha)$. Thus, the remaining integral is
  $c_\alpha(1-\alpha)^{-1}\int_{r+v<t}rv^{1-\alpha-\alpha_0}\ud r\ud v$, which
  gives its first formula by two elementary integrations.

  \medskip
  \noindent\textit{Step 2.~}We identify the four-increment formula and its
  scaling.
  For order four set each spatial increment equal to $h_lU_l$. The Green factor
  $\frac12\mathbf1_{[-h_l,h_l]}$ is $h_l$ times the law of $h_lU_l$. The map
  from increments to vertex positions is triangular within each chain with
  determinant one. Consequently, the four spatial integrals give exactly the
  product $\prod_lh_l$ and the expectation in~\eqref{E:riesz-four}; the factors
  $1/16$ are already incorporated in the four uniform probability laws. The time
  increments give the domain $H_m$ and the differences $S_j-S_q$. Scaling time
  and space by $t$ contributes $t^8$ from the eight integrations and
  $t^{-2\alpha-2\alpha_0}$ from the two covariance factors, giving $t^{2\beta}$.

  \medskip
  \noindent\textit{Step 3.~}We check integrability on the enlarged time and
  space cubes.
  Finiteness can be checked without a numerical argument. On the unit-time
  domain all vertex positions lie in $[-1,1]$. Enlarge the time domain to
  $[0,1]^4$ and the spatial domain to $[-1,1]^4$, and bound the Green indicators
  by one. The two time covariance factors involve disjoint pairs and are
  integrable since $\alpha_0 < 1$. The same is true of the two spatial factors
  since $\alpha < 1$. Tonelli's theorem now gives a finite bound. This proves
  Proposition~\ref{P:riesz-small}.
\end{proof}

\begin{remark}[Normalization and static Riesz check]\label{R:riesz-static}
  The same integrable bound with any fixed exponent $\bar{a} \in (\alpha_0,1)$
  permits $\alpha_0\downarrow0$ at these fixed orders. From~\eqref{E:riesz-two},
  the static total-order-two coefficient is
  \[
    \frac{c_\alpha2^{3-\alpha}}
    {(1-\alpha)(2-\alpha)(3-\alpha)(4-\alpha)}.
  \]
  For example, the double integral of $(s+r)^{2-\alpha}-|s-r|^{2-\alpha}$ on
  $[0,1]^2$ is $(2^{4-\alpha}-4)/((3-\alpha)(4-\alpha))$, which verifies this
  value after adding twice the single-chain coefficient. If one chooses the
  normalization $c_\alpha = (1-\alpha)/2$, this coefficient tends to $1/3$ as
  $\alpha\uparrow1$, in agreement with the white-space value. More generally, at
  each fixed $\alpha_0 \in (0,1)$ the two normalized coefficients
  in~\eqref{E:riesz-two} tend to $c_0/2$ and $c_0$, respectively. For the second
  coefficient, the bracket at $t = 1$ tends to $(1-v)^2$, the prefactor tends to
  $1/2$, and dominated convergence gives $\frac12\int_0^1v^{-\alpha_0}(1-v)^2\ud
  v = c_0$. The brackets are uniformly bounded for $\alpha$ near one, so
  $Cv^{-\alpha_0}$ is an integrable dominator. This is a finite-coefficient
  normalization check; unnormalized $|y|^{-\alpha}$ must not be identified with
  $\delta_0$.
\end{remark}

\subsection{Calibration table and limits of the finite-order test}

Let $M_2$ and $M_4$ denote only the total-order-two and total-order-four
contributions. The results above give
\[
  \begin{array}{c|c|c}
    \text{covariance and time regime}&M_2(t)&M_4(t)\\ \hline
    \delta_0,\quad0<\alpha_0<1
    &2c_0t^{3-\alpha_0}&t^{6-2\alpha_0}\,M_4(1)\\
    \delta_0,\quad\alpha_0=0
    &t^3/3&11t^6/540\\
    c_\alpha|y|^{-\alpha},\quad0<\alpha_0<1
    &t^\beta M_2(1)&t^{2\beta}M_4(1)\\
    c_\alpha|y|^{-\alpha},\quad\alpha_0=0
    &\displaystyle\frac{c_\alpha2^{3-\alpha}t^{4-\alpha}}
    {(1-\alpha)(2-\alpha)(3-\alpha)(4-\alpha)}
    &t^{8-2\alpha}M_4(1)
  \end{array}
\]
Here $c_0$ is defined in Lemma~\ref{L:white-test}; the fractional-time white
coefficient $M_4(1)$ is the explicit four-dimensional time integral
in~\eqref{E:white-four-table}, and the Riesz coefficients are specified
by~\eqref{E:riesz-two}--\eqref{E:riesz-four}. No closed form for those remaining
integrals is asserted.

\begin{remark}[Purpose of these calculations]\label{R:G1-scope}
  The static exponent from the published Chen--Hu theorem is
  $(4-\alpha)/(3-\alpha)$; see~\cite{chen.hu:25:hyperbolic}*{Theorem~1.2, (1.9),
  p.~2147}. Both $\kappa$ and the G0 exponent $\rho$ have this boundary limit.
  The finite coefficient limits above are justified separately by dominated
  convergence; they do not exchange $t \to \infty$ and $\alpha_0\downarrow0$ and
  do not verify an asymptotic constant. In particular, the disputed printed
  white-space constant in the separate Chen--Hu source audit is not used here.
  The single-chain rows agree with the first-moment expansion, and all mixed
  rows are nonnegative, consistent with the Jensen comparison. Relative to the
  Skorohod expansion, the rows with self-pairings are present here and cannot be
  dropped: already $M_2 = \E X_1^2+2\E X_2$. This is a structural comparison
  only, not a Skorohod result. Finally, any fixed truncation is a polynomial in
  $t^\beta$ with nonnegative coefficients, so its logarithm divided by
  $t^\kappa$ tends to zero. Orders two and four can falsify a normalization or
  an incorrect pairing rule, but cannot determine an all-order factorial power
  or prove either principal goal.
\end{remark}

No-gos~\ref{P:nogo-03} and~\ref{P:nogo-04} use these finite coefficients
to disprove a squared-mean comparison with unit prefactor for every time
and a coefficient constant uniform near $\alpha_0=1$, respectively.
Neither test determines the large-time constant.

\chapter{Direct upper bounds and static interpolation}\label{S:appendix-direct}

Theorem~\ref{T:full-upper} covers all of Hypothesis~\ref{H:homog}, and its
proof is long. Is there a shorter way in simpler cases? Yes: estimate the
pairing integrals of Proposition~\ref{P:mixed} directly. The price is extra
assumptions; the reward, explicit constants. This appendix has no figure.

Section~\ref{S:d1-upper} bounds the spatial factors first
and then sums the temporal pairings. It reaches the correct
scale when $\alpha\le 1$ (Proposition~\ref{P:low-alpha-upper}),
with separate constants for one-dimensional Riesz and white noise
(Proposition~\ref{P:d1-upper}). Section~\ref{S:higher-alpha} goes beyond
$\alpha=1$, but only under an extra bound on translated convolutions
(Proposition~\ref{P:higher-transfer}); Proposition~\ref{P:higher-obstructions}
shows where that bound fails. Section~\ref{S:static-holder} gives a third
route for isotropic Riesz covariance, by interpolation with static noise.

Each restriction is part of its statement. These bounds neither replace
the proof of Chapter~\ref{S:book-independent-upper} nor give the sharp
constant, and the main argument does not rely on them. On a first reading,
Section~\ref{S:d1-upper} is enough. Episode~\VideoRef{9}{eLz9Hd20ZT0} of
the video series covers this appendix.

\section{Upper bounds for spatial exponents at most one}\label{S:d1-upper}

The upper bound at scale $t^\kappa$ is proved for $0 < \alpha \le 1$ in
dimensions one, two and three under the standing assumptions. The separate
one-dimensional Riesz estimates retain their sharper constants and white
endpoint. Degree-one covariance in dimensions two and three does not mean full
spatial white noise. Section~\ref{S:interval-upper} treats the higher exponents
by a separate argument. The intrinsic constant is identified  in
Theorem~\ref{T:intrinsic-second-moment}.

We retain the two-chain expansion of Proposition~\ref{P:mixed}, but estimate its
spatial and temporal factors in a different order. The spatial estimate keeps a
power of the number of pairs. After taking a common spatial maximum, we sum all
temporal pairings before using the ordered time domains.

Throughout this section $0 < \alpha_0 < 1$. We first work in dimension one with
either of the spatial covariances
\begin{equation}\label{E:d1-covariances}
  \begin{cases}
    \gamma(y)=c_\alpha|y|^{-\alpha},&0<\alpha<1,\quad c_\alpha>0,\\
    \gamma=c_\delta\delta_0,&\alpha=1,\quad c_\delta>0.
  \end{cases}
\end{equation}
Put $b = 2-\alpha$, $q = 3-\alpha$, $\beta = 4-\alpha-\alpha_0$, and
\begin{equation}\label{E:d1-constants}
  Z=\frac{2}{(1-\alpha_0)(2-\alpha_0)},\qquad
  K_{\rm sp}=
  \begin{cases}
    \displaystyle\frac{c_\alpha 2^{2-\alpha}(1-\alpha)^{-\alpha}}
    {(2-\alpha)^{2-\alpha}},&0<\alpha<1,\\[2mm]
    c_\delta,&\alpha=1.
  \end{cases}
\end{equation}
These parameters satisfy Hypothesis~\ref{H:homog}, so the fixed-order
integrability and mixed expansion already established apply.

\subsection{A spatial bound retaining the ordered increments}

For $m+n = 2k$, label the two chains in block order as in~\eqref{E:mixed-diagram}. For each vertex $l$ write
$h_l = s_l-s_{p(l)} > 0$ and $y_l = x_l-x_{p(l)}$, with a separate root zero
in each chain. Thus, $x_l$ is the sum of increments in its own chain
up to label $l$. The spatial change of variables has determinant one.
In this section $x_l$ denotes a spatial position and $y_l$ its increment; this local
notation differs from the position variables $y_l$
in~\eqref{E:mixed-diagram}. The terminal evaluation point of $u(t,x)$
remains denoted by $x$.
Let $p_\epsilon$ be the centered Brownian Gaussian density of variance
$\epsilon$ and mass one, and set
$\gamma_\epsilon = p_\epsilon*\gamma$. Allow a separate positive
$\epsilon_e$ for every pair $e \in D$ and define
\[
  S_{D,\boldsymbol\epsilon}(h)
  =\int_{\R^{2k}}\prod_{l=1}^{2k}G(h_l,y_l)
  \prod_{(i,j)\in D}\gamma_{\epsilon_{i,j}}(x_j-x_i)\ud y,
  \qquad G(h,y)=\tfrac12\mathbf1_{\{|y|<h\}}.
\]
This is an ordinary nonnegative integral; no pointwise delta value is used.

\LeanMarker{L:d1-spatial}
\begin{lemma}[Spatial estimate]\label{L:d1-spatial}
  Suppose~\eqref{E:d1-covariances} holds, $k \ge 1$, and $D$ is any pairing of
  the two chains. Orient every pair $(i,j)$ by $i < j$ and let $O_D$ and $C_D$
  be its opening and closing label sets. Then, uniformly in
  $\boldsymbol\epsilon$,
  \begin{equation}\label{E:d1-spatial-product}
    S_{D,\boldsymbol\epsilon}(h)\le
    \begin{cases}
      \displaystyle\left(\frac{c_\alpha}{1-\alpha}\right)^k
      \prod_{i\in O_D}h_i\prod_{j\in C_D}h_j^{1-\alpha},
      &0<\alpha<1,\\[2mm]
      \displaystyle\left(\frac{c_\delta}{2}\right)^k
      \prod_{i\in O_D}h_i,&\alpha=1.
    \end{cases}
  \end{equation}
  In particular, if $\sum_{l=1}^{2k}h_l \le 2t$, then
  \begin{equation}\label{E:d1-spatial-max}
    S_{D,\boldsymbol\epsilon}(h)
    \le K_{\rm sp}^{\,k}\frac{t^{bk}}{k^{bk}}.
  \end{equation}
\end{lemma}

\begin{proof}
  We first obtain uniform convolution bounds. Descending label integration
  then removes each closing covariance and opening mass in order, and an
  increment maximum gives the final spatial bound.

  \medskip
  \noindent\textit{Step 1.~}We bound the Riesz and white convolutions.
  We first record the convolution bounds. For $0 < \alpha < 1$, the
  integral of the even decreasing function $|y|^{-\alpha}$ over an
  interval of fixed length is maximal when the interval is centered
  at zero. Therefore,
  \[
    \sup_z\int_\R G(h,y)c_\alpha|y+z|^{-\alpha}\ud y
    =\frac{c_\alpha h^{1-\alpha}}{1-\alpha}.
  \]
  By Tonelli's theorem, the integral with $\gamma_\epsilon$ equals the
  $p_\epsilon(w)\ud w$ average of the unsmoothed integrals with shifts $z-w$.
  Its nonnegative mass-one average is therefore bounded by the same supremum. In
  the white case, $\gamma_\epsilon = c_\delta p_\epsilon$, so
  \[
    \sup_z\int_\R G(h,y)\gamma_\epsilon(y+z)\ud y
    \le\frac{c_\delta}{2}.
  \]
  Both estimates are uniform in $\epsilon > 0$.

  \medskip
  \noindent\textit{Step 2.~}We integrate all increments in descending label
  order.
  Now integrate \emph{all} spatial increments in descending label order
  $2k,2k-1,\ldots,1$. Before processing label $l$, the remaining covariance
  factors are exactly those pairs whose closing labels are at most $l$; the
  bounds from the removed pairs depend only on $h$. If $l \in C_D$, its pair is
  $(i,l)$ with $i < l$ and $x_l-x_i = y_l+z$, where $z$ is independent of $y_l$.
  Every other remaining pair has both endpoints smaller than $l$, so its two
  positions are independent of $y_l$: positions involve only earlier increments
  in their own block chain. Thus, the only covariance involving $y_l$ is its own
  pair, and the convolution bound removes it. If $l \in O_D$, its partner has a
  larger label and its covariance factor has already been removed. All remaining
  covariances have endpoints smaller than $l$, so integration of $y_l$ gives
  only $\int G(h_l,y_l)\ud y_l = h_l$. This also preserves the stated invariant
  for the next label. Tonelli's theorem and the nonnegativity of all factors
  justify each step. There are $k$ labels of each type,
  proving~\eqref{E:d1-spatial-product} for every split and pairing.

  \medskip
  \noindent\textit{Step 3.~}We maximize the remaining increment product.
  For Riesz covariance, weighted arithmetic--geometric mean gives
  \[
    \prod_{i\in O_D}h_i\prod_{j\in C_D}h_j^{1-\alpha}
    \le\left(\frac{2t}{bk}\right)^{bk}
    (1-\alpha)^{(1-\alpha)k}.
  \]
  Indeed, the weights are one at the $k$ openings and $1-\alpha$ at the $k$
  closings, with total $bk$, and the sum of the increments is at most $2t$.
  Multiplication by $(c_\alpha/(1-\alpha))^k$ gives~\eqref{E:d1-spatial-max}
  with the stated $K_{\rm sp}$. In the white case, $\prod_{i\in O_D}h_i \le
  (2t/k)^k$; multiplying by $(c_\delta/2)^k$ gives the same conclusion with $b =
  1$ and $K_{\rm sp} = c_\delta$. This proves Lemma~\ref{L:d1-spatial}.
\end{proof}

\subsection{Summing temporal pairings before using time ordering}

Define the total-order coefficients
\[
  A_k(t)=\sum_{m=0}^{2k}\E[X_m(t,x)X_{2k-m}(t,x)],\qquad A_0(t)=1.
\]
They do not depend on $x$, by spatial stationarity. By
Proposition~\ref{P:mixed}, $\E u^2(t,x) = \sum_{k\ge0}A_k(t)$.

\LeanMarker{P:d1-upper}
\begin{proposition}[One-dimensional upper bound]\label{P:d1-upper}
  Under~\eqref{E:d1-covariances}, for every $k \ge 1$ and $t > 0$,
  \begin{equation}\label{E:d1-coefficient}
    A_k(t)\le
    \frac{(2ZK_{\rm sp})^k t^{\beta k}}{k!\,k^{bk}}
    \le\frac{(2ZK_{\rm sp})^k t^{\beta k}}{(k!)^q}.
  \end{equation}
  Consequently,
  \begin{equation}\label{E:d1-upper}
    \E u^2(t,x)\le\exp(C_{\rm dir}t^\kappa),\qquad
    C_{\rm dir}=q(2ZK_{\rm sp})^{1/q},\quad \kappa=\frac{\beta}{q}.
  \end{equation}
\end{proposition}

\begin{proof}
  We extract one spatial bound common to all pairings, then use symmetry of
  their temporal sum. We count all chain splits, pass the fixed-order limits
  and sum the factorial bound.

  \medskip
  \noindent\textit{Step 1.~}We separate the common spatial bound from the
  temporal pairing sum.
  We first work with the Gaussian spatial smoothing and temporal truncation
  in~\cite{chen:author2026:time-dependency}*{(2.9), (2.11), (2.23)--(2.24)}. Let
  $X_n^{\epsilon,\delta}$ be the smoothed Stratonovich term, using common
  parameters in every factor. Wick's formula expresses each mixed expectation of
  total order $2k$ as a sum over $D \in \Pi_k$. After the chain reversal in
  Proposition~\ref{P:mixed}, its spatial factor is the one in
  Lemma~\ref{L:d1-spatial}; its temporal factors are nonnegative and bounded by
  $|s_i-s_j|^{-\alpha_0}$. For times in $[0,t]^m_ < \times[0,t]^{2k-m}_ <$, the
  increments sum to the two terminal times and therefore to at most $2t$.
  Extract~\eqref{E:d1-spatial-max}, which is independent of the time vector, the
  split, the pairing, and the smoothing parameters. The remaining temporal sum
  is
  \[
    F_k(s_1,\ldots,s_{2k})
    =\sum_{D\in\Pi_k}\prod_{(i,j)\in D}|s_i-s_j|^{-\alpha_0}.
  \]
  Every permutation of the $2k$ labels bijects $\Pi_k$, so $F_k$ is fully
  symmetric. Moreover,
  \[
    \int_{[0,t]^2}|s-r|^{-\alpha_0}\ud s\ud r
    =Zt^{2-\alpha_0},\qquad
    \int_{[0,t]^{2k}}F_k(s)\ud s
    =P_k(Zt^{2-\alpha_0})^k<\infty.
  \]
  The second equality uses the disjoint endpoints of every pairing. Partitioning
  the time cube according to the order within each block now gives
  \begin{equation}\label{E:d1-time-symmetry}
    \int_{[0,t]^m_<\times[0,t]^{2k-m}_<}F_k(s)\ud s
    =\frac{P_k(Zt^{2-\alpha_0})^k}{m!(2k-m)!}.
  \end{equation}
  Ties have measure zero. Crucially, this symmetry is applied to the \emph{sum}
  of pairings after taking the common spatial bound, not to an individual
  temporally weighted diagram.

  \medskip
  \noindent\textit{Step 2.~}We count pairings and splits.
  The number of pairings obeys $P_0 = 1$ and $P_k = (2k-1)P_{k-1}$, by choosing
  the partner of the smallest label. Hence, $P_k = (2k)!/(2^kk!)$. The binomial
  identity gives
  \[
    \sum_{m=0}^{2k}\frac1{m!(2k-m)!}=\frac{2^{2k}}{(2k)!}.
  \]
  Multiplying~\eqref{E:d1-time-symmetry} by the spatial bound and summing over
  $m$ yields
  \[
    \sum_{m=0}^{2k}\E[X_m^{\epsilon,\delta}X_{2k-m}^{\epsilon,\delta}]
    \le K_{\rm sp}^{\,k}\frac{t^{bk}}{k^{bk}}
    \frac{2^{2k}}{(2k)!}
    \frac{(2k)!}{2^kk!}(Zt^{2-\alpha_0})^k
    =\frac{(2ZK_{\rm sp})^kt^{\beta k}}{k!\,k^{bk}}.
  \]
  The same proof covers empty chains, with $X_0^{\epsilon,\delta} = 1$.

  \medskip
  \noindent\textit{Step 3.~}We remove regularization and sum the full series.
  By Lemma~\ref{L:joint-fixed-order}, each smoothed factor converges in $L^2$.
  Products converge in $L^1$ by Cauchy--Schwarz; for fixed $k$ the sum over
  splits is finite. The bound is uniform in $\epsilon,\delta$, so it passes to
  $A_k(t)$. Since $k! \le k^k$ and $b > 0$, we have $k!\,k^{bk} \ge (k!)^{b+1} =
  (k!)^q$, proving both inequalities in~\eqref{E:d1-coefficient}.

  Finally, for $z \ge 0$ and $q \ge 1$,
  \[
    \sum_{k\ge0}\frac{z^k}{(k!)^q}
    =\sum_{k\ge0}\left(\frac{z^{k/q}}{k!}\right)^q
    \le\left(\sum_{k\ge0}\frac{z^{k/q}}{k!}\right)^q
    =\exp(qz^{1/q}).
  \]
  Apply this inequality with $z = 2ZK_{\rm sp}t^\beta$, adding the separately
  defined term $A_0 = 1$. Proposition~\ref{P:mixed} then
  gives~\eqref{E:d1-upper}. This proves Proposition~\ref{P:d1-upper}.
\end{proof}

\subsection{Spectral bounds for $0 < \alpha \le 1$ in all three dimensions}

We now retain Hypothesis~\ref{H:homog} with $d \in \{1,2,3\}$ and $0 < \alpha
\le 1$. The covariance is nonnegative and nonnegative definite as in
Section~\ref{S:model}; when it is a distribution, both these properties and its
homogeneity are understood distributionally. Let $\mu$ be its positive Radon
spectral measure in the convention $\gamma(x) = \int e^{i\xi\cdot
x}\mu(\ud\xi)$, and put $V = \mu(\{\xi:|\xi|\le1\})$. No isotropy or spectral
density is assumed. Write $G_h$ for the positive wave-kernel measure of mass
$h$, and $G_{h,\eta} = p_\eta*G_h$ for $\eta > 0$; set $G_{h,0} = G_h$. The
measure notation includes the sphere kernel in dimension three. The mixed-moment
definition of $A_k(t)$ is now applied to the present $d$-dimensional solution.

\LeanMarker{L:low-alpha-convolution}
\begin{lemma}[Low-$\alpha$ convolution bound]\label{L:low-alpha-convolution}
  Under the preceding assumptions with $0 < \alpha < 1$, for every $R,h > 0$,
  \begin{equation}\label{E:low-alpha-radial}
    \mu(\{\xi:|\xi|\le R\})=VR^\alpha,\qquad
    \int_{\R^d}\min\{h,|\xi|^{-1}\}\mu(\ud\xi)
    =\frac{V}{1-\alpha}h^{1-\alpha}.
  \end{equation}
  For $\epsilon > 0$ and $\eta \ge 0$,
  \begin{equation}\label{E:low-alpha-convolution}
    0\le\int_{\R^d}G_{h,\eta}(\ud y)\gamma_\epsilon(y+z)
    \le\frac{V}{1-\alpha}h^{1-\alpha},
    \qquad z\in\R^d.
  \end{equation}
  The estimate is uniform in both smoothing parameters and in $z$.
\end{lemma}

\begin{proof}
  We compute the absolute Fourier envelope from spectral ball masses. We then
  use the finite Gaussian-weighted measure to justify Fourier integration and
  bound the translated convolution.

  \medskip
  \noindent\textit{Step 1.~}We compute the homogeneous spectral envelope.
  The spectral measure of $\gamma(c\,\cdot)$ is the pushforward of $\mu$ under
  $\xi\mapsto c\xi$. Homogeneity and uniqueness of the Fourier transform on
  tempered distributions therefore give $\mu(A/c) = c^{-\alpha}\mu(A)$ for every
  Borel set $A$. Equivalently, $\mu(cA) = c^\alpha\mu(A)$. This proves the
  ball-mass formula. Since $\mu$ is Radon, $V < \infty$, and continuity from
  above gives $\mu(\{0\}) = \lim_{R\downarrow0}VR^\alpha = 0$. By the layer-cake
  identity and Tonelli's theorem,
  \[
    \int_{\R^d}\min\{h,|\xi|^{-1}\}\mu(\ud\xi)
    =\int_0^h\mu(\{\xi:|\xi|\le u^{-1}\})\ud u
    =V\int_0^h u^{-\alpha}\ud u.
  \]
  The endpoint of each indicator has no effect on its $u$ integral. Since
  $\alpha < 1$, the last expression is the finite quantity
  in~\eqref{E:low-alpha-radial}. This computation uses ball masses, not a radial
  density for $\mu$.

  \medskip
  \noindent\textit{Step 2.~}We obtain the smoothed convolution bound.
  Gaussian weighting gives a finite measure, by the polynomial ball growth of
  $\mu$. Thus, Fourier Fubini is justified for $\gamma_\epsilon$ and the finite
  measure $G_{h,\eta}$. From the source's wave-kernel
  formulas~\cite{chen:author2026:time-dependency}*{(2.3)--(2.5), p.~7},
  \[
    \widehat{G}_h(\xi)=\frac{\sin(h|\xi|)}{|\xi|},\qquad
    G_h(\R^d)=h,
  \]
  with value $h$ in the Fourier formula at zero. The Fourier multipliers of
  $p_\epsilon$ and $p_\eta$ have absolute value at most one, with the latter
  multiplier equal to one when $\eta = 0$. Consequently,
  \[
    \left|\int G_{h,\eta}(\ud y)\gamma_\epsilon(y+z)\right|
    \le\int_{\R^d}\frac{|\sin(h|\xi|)|}{|\xi|}\mu(\ud\xi)
    \le\int_{\R^d}\min\{h,|\xi|^{-1}\}\mu(\ud\xi).
  \]
  Nonnegativity of $\gamma$ as a distribution implies $\gamma_\epsilon \ge 0$,
  since Gaussian translates are nonnegative test functions; $G_{h,\eta}$ is also
  nonnegative. This proves the lower bound in~\eqref{E:low-alpha-convolution}.
  Combining the preceding displays proves its upper bound. No positivity of the
  oscillatory Fourier integrand is asserted. This proves
  Lemma~\ref{L:low-alpha-convolution}.
\end{proof}

\LeanMarker{L:alpha-one-convolution}
\begin{lemma}[Degree-one convolution bound]\label{L:alpha-one-convolution}
  Under the preceding covariance assumptions with $\alpha = 1$, for every $h >
  0$, $\epsilon > 0$, $\eta \ge 0$, and $z \in \R^d$,
  \begin{equation}\label{E:alpha-one-convolution}
    0\le\int_{\R^d}G_{h,\eta}(\ud y)\gamma_\epsilon(y+z)
    \le\frac{\pi V}{2}.
  \end{equation}
  No isotropy or angular density is required.
\end{lemma}

\begin{proof}
  We resolve the degree-one spectral measure into radial and angular parts. We
  retain the radial oscillation and evaluate it as a one-dimensional heat
  convolution, whose bound is uniform in direction.

  \medskip
  \noindent\textit{Step 1.~}We obtain the polar spectral decomposition.
  Fourier uniqueness, applied to the distributional homogeneity as above, gives
  $\mu(cA) = c\mu(A)$ and $\mu(\{0\}) = 0$. Define a finite measure on the unit
  sphere by
  \[
    \sigma(E)=\mu\{\xi:0<|\xi|\le1,\ \xi/|\xi|\in E\}.
  \]
  Its total mass is $V$. Homogeneity implies that the corresponding sector of
  radius $R$ has mass $R\sigma(E)$. Subtracting two such sectors shows that a
  polar rectangle $(r_1,r_2]\times E$ has mass $(r_2-r_1)\sigma(E)$. Uniqueness
  of measures on these generating rectangles therefore gives the polar
  decomposition $\mu(\ud\xi) = \ud r\,\sigma(\ud\theta)$, $r > 0$, $\theta \in
  \mathbb{S}^{d-1}$. The real stationary covariance is even; Fourier uniqueness
  makes $\mu$, and hence $\sigma$, even.

  \medskip
  \noindent\textit{Step 2.~}We express the convolution through its radial
  oscillatory integral.
  Use the Brownian Gaussian convention, so the combined covariance and Green
  smoothing multiplier is $e^{-ar^2}$ with $a = (\epsilon+\eta)/2 > 0$. Gaussian
  weighting makes the spectral measure finite. Fourier Fubini and evenness give
  \begin{equation}\label{E:alpha-one-polar}
    \begin{split}
      \int G_{h,\eta}(\ud y)\gamma_\epsilon(y+z)
      &=\int_{\mathbb{S}^{d-1}} I_a(h,\theta\cdot z)\sigma(\ud\theta),\\
      I_a(h,v)&=\int_0^\infty e^{-ar^2}
      \frac{\sin(hr)}r\cos(rv)\ud r.
    \end{split}
  \end{equation}
  The radial integral is absolutely convergent: its integrand is bounded near
  zero and Gaussian-decaying at infinity. Passing from the complex exponential
  to cosine introduces no factor two, since both antipodal directions are
  already included in $\sigma$.

  \medskip
  \noindent\textit{Step 3.~}We identify that integral with a bounded
  one-dimensional convolution.
  Let $g_h(v) = \tfrac12\mathbf1_{\{|v|<h\}}$ and let $p_{2a}^{(1)}(v) = (4\pi
  a)^{-1/2}e^{-v^2/(4a)}$ be the one-dimensional Gaussian density. Their Fourier
  transforms are $\sin(hr)/r$ and $e^{-ar^2}$. The one-dimensional inversion
  factor is $1/(2\pi)$ in the source convention, so
  \[
    (g_h*p_{2a}^{(1)})(v)
    =\frac1{2\pi}\int_\R e^{irv}e^{-ar^2}
    \frac{\sin(hr)}r\ud r
    =\frac1\pi I_a(h,v).
  \]
  The sign in the inversion exponential is immaterial because the integrand
  pairs evenly under $r\mapsto-r$. Since $p_{2a}^{(1)}$ has mass one and $0 \le
  g_h \le 1/2$, we obtain $0 \le I_a(h,v) \le \pi/2$. Integrating against
  $\sigma$, whose mass is $V$, proves~\eqref{E:alpha-one-convolution}. This
  retains the oscillation until the one-dimensional integral is evaluated; the
  divergent absolute Fourier envelope at $\alpha = 1$ is not used. This proves
  Lemma~\ref{L:alpha-one-convolution}.
\end{proof}

\LeanMarker{P:low-alpha-upper}
\begin{proposition}[Upper bound for $\alpha \le 1$]\label{P:low-alpha-upper}
  Suppose that Hypothesis~\ref{H:homog} and the standing covariance positivity
  assumptions hold, $d \in \{1,2,3\}$, $0 < \alpha \le 1$, and $V > 0$. Define
  \begin{equation}\label{E:low-alpha-constants}
    \begin{split}
      K_{\rm sp}^{\rm gen}
      &=\begin{cases}
        \displaystyle\frac{V2^{2-\alpha}(1-\alpha)^{-\alpha}}
        {(2-\alpha)^{2-\alpha}},&0<\alpha<1,\\[2mm]
        \pi V,&\alpha=1,
      \end{cases}\\
      C_{\rm dir}^{\rm gen}&=q(2ZK_{\rm sp}^{\rm gen})^{1/q}.
    \end{split}
  \end{equation}
  With the same $A_k$, $b$, $q$, $\beta$, and $Z$ as above, for
  $k \ge 1$ and $t > 0$,
  \begin{equation}\label{E:low-alpha-coefficient}
    A_k(t)\le\frac{(2ZK_{\rm sp}^{\rm gen})^k t^{\beta k}}
    {k!\,k^{bk}}
    \le\frac{(2ZK_{\rm sp}^{\rm gen})^k t^{\beta k}}{(k!)^q},
    \qquad
    \E u^2(t,x)\le\exp(C_{\rm dir}^{\rm gen}t^\kappa).
  \end{equation}
\end{proposition}

\begin{proof}
  We apply the descending-label spatial argument in every dimension using the
  appropriate closing bound. We then sum temporal pairings and pass the Green
  and noise limits in their prescribed order.

  \medskip
  \noindent\textit{Step 1.~}We obtain a common spatial bound in each
  dimension.
  First smooth the noise as in~\cite{chen:author2026:time-dependency}*{(2.9), (2.11), (2.23)--(2.24)}.
  In dimension three also use $G_{h,\eta}$ with $\eta > 0$, as in the
  source's pp.~13--14. The increment change of variables is now a
  triangular map on vectors in $\R^d$, still of determinant one.
  The descending-label proof of Lemma~\ref{L:d1-spatial} depends only
  on which increments enter each position, not on their dimension.
  Its invariant therefore remains valid: at a closing label there
  is exactly one covariance depending on that vector increment, and
  at an opening label its pair has already been removed.
  If $0 < \alpha < 1$, Lemma~\ref{L:low-alpha-convolution} gives the
  closing bound $Mh_l^{1-\alpha}$ with $M = V/(1-\alpha)$.
  At $\alpha = 1$, Lemma~\ref{L:alpha-one-convolution} gives the
  constant closing bound $M = \pi V/2$. Every opening gives mass
  $h_l$, uniformly in $\eta$ and all noise-smoothing parameters.
  For $\alpha < 1$, the same weighted arithmetic--geometric mean
  calculation gives
  \[
    M^k\prod_{i\in O_D}h_i\prod_{j\in C_D}h_j^{1-\alpha}
    \le (K_{\rm sp}^{\rm gen})^k\frac{t^{bk}}{k^{bk}}.
  \]
  At $\alpha = 1$, ordinary arithmetic--geometric mean on the $k$ opening
  increments gives $M^k\prod_{i\in O_D}h_i \le (2M)^kt^k/k^k$, which is the same
  bound with $b = 1$ and $K_{\rm sp}^{\rm gen} = \pi V$. Both cases use
  $\sum_lh_l \le 2t$. For $d = 1,2$ one may use the Green densities directly;
  their opening mass and closing convolution obey these same bounds.

  \medskip
  \noindent\textit{Step 2.~}We sum the temporal factors.
  Replace the nonnegative truncated temporal factors by $|s_i-s_j|^{-\alpha_0}$
  before taking their pairing sum. The temporal symmetry, cube integral and
  split count in~\eqref{E:d1-time-symmetry} do not involve the spatial
  dimension. Thus, exactly the multiplication of counts in the proof of
  Proposition~\ref{P:d1-upper} gives the first bound
  in~\eqref{E:low-alpha-coefficient} for the regularized moments.

  \medskip
  \noindent\textit{Step 3.~}We remove regularizations and sum the coefficient
  bounds.
  We spell out the limit in dimension three. Fix every noise-smoothing parameter
  and first let $\eta\downarrow0$. At each time vector, the Gaussian moment
  formula is a finite sum of bounded continuous products of the smoothed
  covariance functions, integrated against products of the finite Green
  measures. These products of measures converge weakly as $\eta\downarrow0$,
  with the same masses $\prod_l h_l$. At fixed noise smoothing the Gaussian
  covariance factors are uniformly bounded, so a constant times $\prod_lh_l$
  dominates on each bounded time simplex. Dominated convergence applies also to
  the covariance between two approximations with different $\eta$ values. The
  resulting common joint covariance limit proves the $L^2$ Cauchy property of
  each fixed-order random integral. This is the inner limit used in the source's
  definition on pp.~13--14. Products of two factors converge in $L^1$ by
  Cauchy--Schwarz, so the moment inequality survives this limit. Next send the
  noise-smoothing parameters to zero using the fixed-order $L^2$ convergence of
  Lemma~\ref{L:joint-fixed-order}, and again pass products in $L^1$. In
  dimensions one and two only the latter limit is required. The bounds at both
  stages are uniform.

  We have proved the first coefficient bound for $A_k(t)$. The factorial
  inequality and nonnegative-series argument in Proposition~\ref{P:d1-upper}
  apply unchanged, with $K_{\rm sp}^{\rm gen}$ in place of $K_{\rm sp}$ and $A_0
  = 1$. This proves Proposition~\ref{P:low-alpha-upper}.
\end{proof}

\begin{remark}[Scope of the spectral estimates]\label{R:low-alpha-scope}
  If $V = 0$, spectral homogeneity gives $\mu(B_R) = 0$ for every $R$, hence
  $\mu = 0$ and $u = 1$; the noise-free conclusion is immediate. Otherwise, all
  angular and dimensional dependence is allowed in $V$. For $0 < \alpha < 1$,
  the example $\gamma(x) = c|x_1|^{-\alpha}$ in dimensions two or three has its
  spectral measure supported on an axis; a full-dimensional density is not
  required. The absolute Fourier envelope diverges at $\alpha = 1$, where the
  separate oscillatory calculation in Lemma~\ref{L:alpha-one-convolution} is
  essential. The endpoint is not obtained by taking a limit in $K_{\rm sp}^{\rm
  gen}$ for $\alpha < 1$. For $d = 1$ and $\gamma = c_\delta\delta_0$, the
  source convention gives $\mu(\ud\xi) = c_\delta\ud\xi/(2\pi)$ and $V =
  c_\delta/\pi$. Hence, the endpoint constant $K_{\rm sp}^{\rm gen} = \pi V =
  c_\delta$ agrees exactly with Proposition~\ref{P:d1-upper}, with no extra
  antipodal factor. In dimensions two and three the endpoint covers degree-one
  homogeneous covariance, not full spatial white noise. No result for $\alpha >
  1$ follows from these estimates.
\end{remark}

\subsection{Checks on the scale, constants and scope}

\LeanMarker{C:d1-first-calibration}
\begin{corollary}[First-moment calibration]\label{C:d1-first-calibration}
  Let $K_*$ be $K_{\rm sp}$ under~\eqref{E:d1-covariances}, or $K_{\rm sp}^{\rm
  gen}$ under the hypotheses of Proposition~\ref{P:low-alpha-upper}. Put
  $C_{2,*} = q(2ZK_*)^{1/q}$. Then
  \begin{equation}\label{E:d1-first-calibration}
    \E u(t,x)\le\exp(C_{1,*}t^\kappa),\qquad
    C_{1,*}=q\left(\frac{ZK_*}{2^{3-\alpha}}\right)^{1/q}
    =\frac{C_{2,*}}{2^{(4-\alpha)/(3-\alpha)}}.
  \end{equation}
  In either nonzero case $2C_{1,*} < C_{2,*}$. Conditionally on a first-moment
  lower rate at least $L_*$, $2L_* \le 2C_{1,*} < C_{2,*}$, consistently with
  Proposition~\ref{P:Jensen}.
\end{corollary}

\begin{proof}
  For one chain of length $2k$, the increments sum to at most $t$. The spatial
  estimate in Lemma~\ref{L:d1-spatial}, or in the proof of
  Proposition~\ref{P:low-alpha-upper}, therefore applies with $t/2$ in place of
  $t$. The temporal pairing sum on one ordered simplex is
  $P_k(Zt^{2-\alpha_0})^k/(2k)!$. Using the same regularization and limit order
  as in the respective upper-bound proof gives
  \[
    \E X_{2k}(t,x)
    \le\frac{(ZK_*/2^{b+1})^k t^{\beta k}}{k!\,k^{bk}}
    \le\frac{(ZK_*/2^{3-\alpha})^k t^{\beta k}}{(k!)^q}.
  \]
  Every odd term has zero expectation. The $L^2$ expansion permits taking
  expectations of partial sums, and the resulting nonnegative even series is
  summed by the scalar inequality in the proof of Proposition~\ref{P:d1-upper}.
  This gives~\eqref{E:d1-first-calibration}. If the first-moment lower rate is
  at least $L_*$, \eqref{E:d1-first-calibration} implies $L_* \le C_{1,*}$. In
  the stated nonzero cases $K_* > 0$, hence $C_{1,*} > 0$; and
  $(4-\alpha)/(3-\alpha) > 1$ gives $2C_{1,*} < C_{2,*}$. This proves
  Corollary~\ref{C:d1-first-calibration}.
\end{proof}

\begin{remark}[Checks and comparison with the uniform pairing bound]\label{R:d1-upper-checks}
  The argument preserves the factor $k^{-(2-\alpha)k}$ from the spatial
  increment maximum. Replacing the spatial volume merely by a constant times
  $t^{(2-\alpha)k}$ would lose the factorial power needed
  in~\eqref{E:d1-coefficient}. This proof does not use the
  common-Laplace-parameter norm majorant of Proposition~\ref{P:envelope}, and it
  does not assert a stronger factorial estimate for every individual diagram.
  Instead, the temporal ordering is exploited after summing pairings
  in~\eqref{E:d1-time-symmetry}. Thus, the scoped G0 obstruction remains valid.
  For unit white covariance, the bound on the total-order-two coefficient is
  $2Z$, whereas the exact value from Lemma~\ref{L:white-test} is $2c_0$; their
  ratio is $2(3-\alpha_0) > 4$. At the static endpoint the exact order-four
  coefficient is $11/540$, while the first bound in~\eqref{E:d1-coefficient}
  gives $1/2$. These are checks of normalization, not evidence of sharp
  constants.
\end{remark}

No-go~\ref{P:nogo-06} quantifies the loss avoided here: discarding the
spatial constraints before summation gives the wrong time exponent.

\begin{remark}[Boundary and literature comparisons]\label{R:d1-comparisons}
  As $\alpha_0\downarrow0$, the time exponent tends to $(4-\alpha)/(3-\alpha)$,
  the Chen--Hu static exponent cited in Remark~\ref{R:G1-scope}. This is an
  exponent calibration, not an interchange of the static and large-time limits.
  With the Riesz normalization $c_\alpha = (1-\alpha)/2$, the displayed $K_{\rm
  sp}$ tends to one as $\alpha\uparrow1$, matching unit white covariance. For
  comparison only, the Skorohod upper bound of
  Balan--Conus~\cite{balan.conus:16:intermittency}*{Theorem~2.1(b), (8), (10),
  p.~1492} has time exponent $(2H+2-a)/(3-a)$. Setting $H = 1-\alpha_0/2$ and $a
  = \alpha$ yields the same $\kappa$. Their covariance normalization contains
  $H(2H-1)$ and they use a Skorohod mild solution. No constants or solutions are
  identified, and their theorem is not an input to the present proof. The
  current coverage of this section is $0 < \alpha \le 1$ in $d = 1,2,3$,
  including the one-dimensional white case. It does not address $1 < \alpha <
  2-\alpha_0$ in dimensions two and three.
\end{remark}

\section{Higher spatial exponents under a convolution bound}
\label{S:higher-alpha}

The convolution assumption below is additional to homogeneity, and the concrete
applications cover only the stated Riesz ranges. Section~\ref{S:interval-upper}
establishes the full general upper scale by a separate argument. This
alternative route is not an input to Theorem~\ref{T:intrinsic-second-moment}.

We retain the standing covariance assumptions and take $d \in \{2,3\}$ and $1 <
\alpha < 2-\alpha_0$. Put $\delta = \alpha-1$, $q = 3-\alpha$, and $\beta =
4-\alpha-\alpha_0$. Here $\delta$ is an exponent, not a smoothing parameter.
Write $p_\epsilon$ for the centered Gaussian density with covariance $\epsilon
I$, $\gamma_\epsilon = p_\epsilon*\gamma$, and $G_{h,\eta} = p_\eta*G_h$, with
$G_{h,0} = G_h$.

\begin{hypothesis}[Spatial maximum]\label{H:spatial-maximum}
  There is $M > 0$ such that, for every $h > 0$, $z \in \R^d$, $\epsilon > 0$,
  and $\eta \ge 0$,
  \begin{equation}\label{E:higher-convolution}
    0\le\int_{\R^d}G_{h,\eta}(\ud y)\gamma_\epsilon(y+z)
    \le M h^{-\delta}.
  \end{equation}
  The same $M$ works for all translations and both Gaussian widths. In
  particular, each pairing edge may have its own positive covariance width and
  each Green factor its own nonnegative smoothing width.
\end{hypothesis}

This uniform maximum is an additional assumption, not a consequence of
homogeneity. No-go~\ref{P:nogo-08} gives admissible covariances for which
the translated wave convolution diverges.

\begin{proposition}[Conditional upper bound]\label{P:higher-transfer}
  Suppose that Hypotheses~\ref{H:homog} and~\ref{H:spatial-maximum} hold in the
  range specified above. Choose and fix
  \[
    \frac1{1-\alpha_0}<r<\frac1\delta,\qquad
    r'=\frac r{r-1},\qquad \theta=\alpha_0r',\qquad
    a=1-\frac{\delta r}{2}.
  \]
  Define the explicit, nonsharp constants
  \begin{equation}\label{E:higher-constants}
    Z_\theta=\frac2{(1-\theta)(2-\theta)},\qquad
    D_a=e^2(e/a)^{2a},\qquad
    K=12M[\Gamma(1-\delta r)D_a]^{1/r}Z_\theta^{1/r'}.
  \end{equation}
  For the total-order coefficients $A_k(t)$ of Section~\ref{S:d1-upper},
  \begin{equation}\label{E:higher-bound}
    A_k(t)\le\frac{K^kt^{\beta k}}{(k!)^q}\quad(k\ge1),
    \qquad \E u^2(t,x)\le\exp\{qK^{1/q}t^\kappa\}.
  \end{equation}
  These estimates hold for every $t > 0$ and $x \in \R^d$.
\end{proposition}

\begin{proof}
  We retain the singular closing increments instead of replacing them by a
  spatial maximum. Hölder on the common pairing-time measure separates those
  increments from the temporal factors. Dirichlet and symmetry estimates then
  recover the required factorial power before the regularization limits.

  \medskip
  \noindent\textit{Step 1.~}We verify parameters and perform the spatial
  elimination.
  The interval for $r$ is nonempty because $\delta+\alpha_0 < 1$. Thus, $\delta
  r < 1$, $\theta < 1$, and $1/2 < a < 1$. We use the source's noise
  regularization from Chen~\cite{chen:author2026:time-dependency}*{(2.9),
  (2.11), (2.23)--(2.24)}, as in Section~\ref{S:d1-upper}. Its nonnegative
  temporal covariance is bounded by $|s_i-s_j|^{-\alpha_0}$, and its spatial
  covariance on an edge $e$ is $\gamma_{\epsilon_e+\epsilon'_e}$, where
  $\epsilon_e,\epsilon'_e > 0$ are the two endpoint noise widths. Indeed,
  $p_{\epsilon_e}*p_{\epsilon'_e} = p_{\epsilon_e+\epsilon'_e}$; equal widths
  give the source's $\gamma_{2\epsilon}$. Every such positive covariance width
  is covered by Hypothesis~\ref{H:spatial-maximum}. In dimension three first
  smooth the Green measures as prescribed on pp.~13--14 of that source.

  Fix $m+n = 2k$, $k \ge 1$, and put $\Omega_{m,n} = [0,t]^m_ < \times[0,t]^n_
  <$. Use the two separate roots and increments $h_l = s_l-s_{p(l)}$ of
  Section~\ref{S:d1-upper}. For $D \in \Pi_k$, orient its pairs by increasing
  block label and denote their opening and closing sets by $O_D,C_D$. Descending
  integration of the vector increments gives
  \[
    S_D(h)\le M^k\prod_{i\in O_D}h_i\prod_{j\in C_D}h_j^{-\delta}
    \le M^k(2t/k)^k H_D(h),\qquad
    H_D(h)=\prod_{j\in C_D}h_j^{-\delta}.
  \]
  Indeed, before label $l$ is integrated, the remaining covariances have closing
  labels at most $l$. At a closing label exactly its own pair depends on that
  increment, as $\gamma_\epsilon(y_l+z)$; all other endpoints have smaller
  labels. Hypothesis~\ref{H:spatial-maximum} removes this factor. At an opening
  its pair is already removed and the integral is the Green mass $h_l$. This
  proves the invariant in every dimension and for empty chains. Tonelli applies
  by positivity. Finally, the $k$ opening increments sum to at most $2t$, so
  ordinary arithmetic--geometric mean gives the second inequality. We have not
  applied the low-$\alpha$ conclusion of Lemma~\ref{L:d1-spatial}.

  \medskip
  \noindent\textit{Step 2.~}We separate the spatial gap powers and temporal
  pairing factors.
  Set $T_D(s) = \prod_{(i,j)\in D}|s_i-s_j|^{-\alpha_0}$. H\"older's inequality
  on counting measure on $\Pi_k$ times Lebesgue measure on $\Omega_{m,n}$ gives
  \begin{equation}\label{E:higher-holder}
    \sum_D\int_{\Omega_{m,n}}H_DT_D\ud s
    \le\left(\sum_D\int_{\Omega_{m,n}}H_D^r\ud s\right)^{1/r}
    \left(\sum_D\int_{\Omega_{m,n}}T_D^{r'}\ud s\right)^{1/r'}.
  \end{equation}
  Both integrals are finite by the following computations. If the two chains
  contain $k_1,k_2$ closing labels, then $k_1+k_2 = k$ and the simplex change to
  increments has determinant one. The Dirichlet integral, including the unused
  terminal gap with exponent zero, gives
  \[
    \int_{\Omega_{m,n}}H_D^r\ud s
    =\frac{\Gamma(1-\delta r)^k t^{(2-\delta r)k}}
    {\Gamma(1+m-\delta r k_1)\Gamma(1+n-\delta r k_2)}.
  \]
  To see the normalization, in a chain of length $m$ the $k_1$ closing
  increments have Dirichlet parameter $1-\delta r$, the other increments have
  parameter one, and the terminal gap has parameter one. Scaling gives
  $t^{m-\delta r k_1}$. Iterating the beta integral gives the displayed Gamma
  quotient. Empty chains contribute $\Gamma(1) = 1$. The two denominator
  arguments are positive and have average $1+ak$. Log convexity of $\Gamma$
  therefore bounds their product below by $\Gamma(1+ak)^2$. An elementary lower
  bound makes the constant explicit:
  \[
    \Gamma(1+ak)\ge\int_{ak}^{ak+1}x^{ak}e^{-x}\ud x
    \ge e^{-1}(a/e)^{ak}k^{ak}
    \ge e^{-1}(a/e)^{ak}(k!)^a.
  \]
  Since $k \ge 1$, its square is at least $D_a^{-k}(k!)^{2-\delta r}$. Writing
  $P_k = |\Pi_k|$, we obtain
  \[
    \sum_D\int_{\Omega_{m,n}}H_D^r\ud s
    \le\frac{P_k[\Gamma(1-\delta r)D_a]^k
    t^{(2-\delta r)k}}{(k!)^{2-\delta r}}.
  \]
  The sum $\sum_DT_D^{r'}$ is symmetric in all $2k$ time labels, because a
  permutation bijects the pairings. On the full cube each pairing has integral
  $(Z_\theta t^{2-\theta})^k$. Therefore,
  \[
    \sum_D\int_{\Omega_{m,n}}T_D^{r'}\ud s
    =\frac{P_k(Z_\theta t^{2-\theta})^k}{m!n!}
    \le\frac{P_k(Z_\theta t^{2-\theta})^k}{(k!)^2}.
  \]
  Here $m!n! \ge (k!)^2$ since $m+n = 2k$. Symmetry is used only for the pairing
  sum, not for an individual diagram. Also, $P_k = \prod_{j=1}^k(2j-1) \le
  2^kk!$ and $k^k \ge k!$. Combining these bounds with~\eqref{E:higher-holder}
  and the opening estimate proves, uniformly in all regularization parameters,
  \[
    \sum_D\int_{\Omega_{m,n}}S_D T_D\ud s
    \le\frac{\{4M[\Gamma(1-\delta r)D_a]^{1/r}
    Z_\theta^{1/r'}\}^k t^{\beta k}}{(k!)^q}.
  \]
  Indeed, the time power is $1+(2-\delta r)/r+(2-\theta)/r' = 3-\delta-\alpha_0
  = \beta$; the factorial power is $1+(2-\delta r)/r+2/r'-1 = 2-\delta = q$.
  There are $2k+1 \le 3^k$ splits, giving the asserted $K$.

  \medskip
  \noindent\textit{Step 3.~}We pass the fixed-order limits.
  We now pass to the solution using the independent fixed-order construction in
  Lemma~\ref{L:joint-fixed-order}, which applies under Hypothesis~\ref{H:homog}.
  In dimension three, fix the noise smoothing first and send the Green smoothing
  to zero. Smoothed spatial covariances are bounded and continuous. Products of
  the finite Green measures converge weakly, with masses $\prod_lh_l$. At fixed
  noise smoothing all covariance factors are bounded, so a constant times
  $\prod_lh_l$ dominates the time integrals on bounded simplexes. The same
  argument for mixed Green approximations gives a common joint covariance limit,
  and hence the fixed-order $L^2$ Cauchy property. This is the inner limit of
  the source's pp.~13--14, detailed in Proposition~\ref{P:low-alpha-upper}. Next
  remove the noise smoothing using Lemma~\ref{L:joint-fixed-order}'s fixed-order
  $L^2$ convergence. In dimension two only this latter limit is needed. At each
  stage products converge in $L^1$, since
  \[
    \|YZ-Y'Z'\|_1
    \le\|Y-Y'\|_2\|Z\|_2+\|Y'\|_2\|Z-Z'\|_2.
  \]
  Each split sum is finite and the bound is uniform. Thus, it passes to
  $A_k(t)$; zero chains remain the constant one. These convergence inputs use
  the standing assumptions, not $\alpha \le 1$.

  \medskip
  \noindent\textit{Step 4.~}We sum the full coefficient series.
  Finally, Proposition~\ref{P:mixed} supplies the nonnegative full series under
  these same hypotheses. Since $q > 1$, for $z \ge 0$,
  \[
    \sum_{k\ge0}\frac{z^k}{(k!)^q}
    =\sum_{k\ge0}\left(\frac{z^{k/q}}{k!}\right)^q
    \le\left(\sum_{k\ge0}\frac{z^{k/q}}{k!}\right)^q
    =e^{qz^{1/q}}.
  \]
  Take $z = Kt^\beta$ and use $A_0 = 1$ and $\beta/q = \kappa$. This proves
  Proposition~\ref{P:higher-transfer}.
\end{proof}

% This lemma discharges H:spatial-maximum only for its explicitly stated
% isotropic Riesz ranges.
\LeanMarker{L:higher-riesz}
\begin{lemma}[Riesz spatial maxima]\label{L:higher-riesz}
  Let $\gamma(x) = c|x|^{-\alpha}$, $c > 0$. In dimension three, for $1 < \alpha
  < 2$, Hypothesis~\ref{H:spatial-maximum} holds with the optimal uniform
  constant $M_3 = c2^{1-\alpha}/(2-\alpha)$. In dimension two, for $1 < \alpha <
  3/2$, it holds with
  \[
    M_2=c\left\{4^\alpha+
    \frac{2\alpha}{3/2-\alpha}4^{\alpha-3/2}\right\}.
  \]
\end{lemma}

\begin{proof}
  We compute the spherical maximum in dimension three and bound small
  translated ball masses in dimension two. Gaussian averaging preserves those
  bounds, and an approximate-identity limit proves the stated optimality in
  dimension three.

  \medskip
  \noindent\textit{Step 1.~}We evaluate the three-dimensional spherical
  convolution.
  From Chen~\cite{chen:author2026:time-dependency}*{(2.5), p.~7}, the wave
  measures are
  \[
    G_h^{(2)}(\ud y)=\frac{\mathbf1_{\{|y|<h\}}\ud y}
    {2\pi\sqrt{h^2-|y|^2}},\qquad
    G_h^{(3)}(\ud y)=\frac{\sigma_h(\ud y)}{4\pi h},
  \]
  where $\sigma_h$ is surface measure on the sphere of radius $h$. Both have
  mass $h$ and satisfy $G_h(hA) = hG_1(A)$. Thus, the unsmoothed convolution is
  $ch^{1-\alpha}F_d(z/h)$, where $F_d(w) = \int G_1^{(d)}(\ud
  y)|y+w|^{-\alpha}$. In dimension three spherical coordinates give, for $\rho =
  |w| > 0$,
  \[
    F_3(w)=\frac{(1+\rho)^v-|1-\rho|^v}{2v\rho},
    \qquad v=2-\alpha\in(0,1).
  \]
  For $0 < \rho < 1$ this is the average of the convex function $f(s) =
  (1+s)^{v-1}$ on $[-\rho,\rho]$. Its derivative has numerator
  $\rho[f(\rho)+f(-\rho)]-\int_{-\rho}^\rho f(s)\ud s \ge 0$, by the endpoint
  chord bound for a convex function. For $\rho > 1$, every distance
  $\sqrt{1+\rho^2+2\rho u}$, $-1 \le u \le 1$, increases with $\rho$, so its
  negative power decreases. The formula is continuous at $\rho = 1$ and at zero.
  Its maximum is therefore $F_3(e_1) = 2^{1-\alpha}/(2-\alpha)$.

  \medskip
  \noindent\textit{Step 2.~}We bound the planar convolution by local ball
  masses.
  In dimension two we claim $G_1(B(w,R)) \le 2R^{3/2}$ for $0 < R \le 1/4$.
  Rotate $w = (A,0)$, $A \ge 0$. If $A < 1/2$, the ball lies in $|y| < 3/4$,
  where the density is at most $2/(\pi\sqrt7)$; its mass is at most
  $2R^2/\sqrt7$. If $A \ge 1/2$, write $y = (x,v)$. On the ball $|v| \le R$ and
  $x \ge 1/4$. Put $L = \sqrt{1-v^2}$. Inside the disk, $1-|y|^2 = (L-x)(L+x)$
  with $L+x > 1$. Each horizontal slice has length at most $2R$, so integrating
  $(L-x)^{-1/2}$ gives at most $2\sqrt{2R}$. The $v$ range has length at most
  $2R$, giving mass at most $(2\sqrt2/\pi)R^{3/2} \le 2R^{3/2}$, as claimed.
  Layer cake now gives
  \[
    F_2(w)\le4^\alpha+
    \alpha\int_0^{1/4}R^{-\alpha-1}G_1(B(-w,R))\ud R
    \le M_2/c.
  \]
  This includes centers outside the disk and is uniform in $w$.

  \medskip
  \noindent\textit{Step 3.~}We retain smoothing uniformity and prove the
  optimality claim.
  Finally, Tonelli and the Gaussian convolution rule give
  \[
    \int G_{h,\eta}(\ud y)\gamma_\epsilon(y+z)
    =\int p_{\epsilon+\eta}(w)(G_h*\gamma)(z-w)\ud w.
  \]
  Mass-one averaging preserves both supremum bounds, for all widths in
  Hypothesis~\ref{H:spatial-maximum}. We justify the approximate identity used
  for optimality. For $f(x) = |x|^{-\lambda}$ in $\R^j$, $0 < \lambda < j$, fix
  $x \ne 0$. In the convolution $\int p_\epsilon(x-y)f(y)\ud y$, a small
  neighborhood of $x$ gives convergence by continuity. On $|y| < |x|/2$ the
  integral is bounded by a constant times
  $\epsilon^{-j/2}e^{-|x|^2/(8\epsilon)}$, since $f$ is integrable on that ball.
  On the remaining region outside the neighborhood of $x$, $f$ is bounded and
  the Gaussian tail tends to zero. Thus, $p_\epsilon*f(x) \to f(x)$ for every $x
  \ne 0$. In dimension three Fatou's lemma at a maximizing shift as
  $\epsilon\downarrow0$, $\eta = 0$, now recovers the unsmoothed value; the
  single singular point has zero sphere measure. Hence, no smaller uniform
  constant is possible. This proves Lemma~\ref{L:higher-riesz}.
\end{proof}

\begin{corollary}[Riesz moment bounds]\label{C:higher-riesz}
  For isotropic Riesz covariance $c|x|^{-\alpha}$, the bounds
  in~\eqref{E:higher-bound} hold in dimension three for $1 < \alpha <
  2-\alpha_0$, and in dimension two for $1 < \alpha < \min\{3/2,2-\alpha_0\}$,
  with $M = M_d$ in~\eqref{E:higher-constants}.
\end{corollary}

\begin{proof}
  These Riesz functions are nonnegative, homogeneous of degree $-\alpha$,
  locally integrable since $\alpha < d$, and nonnegative definite. For the last
  property, the identity
  \[
    |x|^{-\alpha}=\frac1{\Gamma(\alpha/2)}
    \int_0^\infty s^{\alpha/2-1}e^{-s|x|^2}\ud s
  \]
  represents the covariance as a positive mixture of positive definite
  Gaussians. To justify this distributionally, take a complex Schwartz function
  $\varphi$, put $\widetilde{\varphi}(x) = \overline{\varphi(-x)}$, and set
  $\psi = \varphi*\widetilde{\varphi}$. This is Schwartz and $\int
  |x|^{-\alpha}|\psi(x)|\ud x < \infty$: local integrability controls the origin
  and rapid decay controls infinity. Tonelli applied with $|\psi|$ therefore
  proves absolute integrability of the mixture against $\psi$, so Fubini
  applies. Each Gaussian has nonnegative integral against $\psi$, as follows
  from its positive Fourier transform and the factor $|\widehat{\varphi}|^2$.
  The mixture consequently has nonnegative integral against $\psi$. This proves
  nonnegative definiteness as a tempered distribution. The stated ranges satisfy
  $\alpha+\alpha_0 < 2$, and Lemma~\ref{L:higher-riesz} supplies the additional
  spatial hypothesis with all Gaussian uniformities.
  Proposition~\ref{P:higher-transfer} applies. This proves
  Corollary~\ref{C:higher-riesz}.
\end{proof}

Another possible reduction would merge two neighboring wave kernels.
No-go~\ref{P:nogo-10} explains why an unweighted convolution estimate
cannot simply discard the covariance linking a merged vertex to its remote
partner in a crossing diagram.

\begin{proposition}[Obstructions to the spatial hypothesis]
  \label{P:higher-obstructions}
  In dimension two, Riesz covariance with $3/2 \le \alpha < 2$ fails
  Hypothesis~\ref{H:spatial-maximum}. For every $1 < \alpha < 2$ there are
  nonnegative, nonnegative definite homogeneous locally integrable covariance
  functions in both dimensions two and three that fail it.
\end{proposition}

\begin{proof}
  We exhibit a divergent planar Riesz convolution at the boundary, then
  construct an anisotropic product covariance with the same obstruction in
  both dimensions. Approximate-identity limits transfer these divergences to
  the smoothed uniform-bound question.

  \medskip
  \noindent\textit{Step 1.~}We identify the divergent boundary configurations.
  For the first claim set $h = 1$, $z = -e_1$, and $y =
  e_1+\rho(\cos\phi,\sin\phi)$. On the wedge $0 < \rho < 1/2$, $2\pi/3 \le \phi
  \le 4\pi/3$, we have $0 < 1-|y|^2 \le 2\rho$. Therefore,
  \[
    \int G_1^{(2)}(\ud y)c|y-e_1|^{-\alpha}
    \ge\frac{c}{3\sqrt2}\int_0^{1/2}\rho^{1/2-\alpha}\ud\rho
    =\infty.
  \]
  This includes logarithmic divergence at $\alpha = 3/2$. For the second claim
  choose $a_1 = (1+\max\{2-\alpha,\alpha-1\})/2$ and $a_2 = \alpha-a_1$. Then $0
  < a_1,a_2 < 1$ and $a_1+a_2/2 > 1$. The function $\gamma(x) =
  c|x_1|^{-a_1}|x_2|^{-a_2}$, independent of $x_3$ in dimension three, is
  locally integrable, tempered, and homogeneous of degree $-\alpha$. The
  Gaussian mixture identity in the preceding proof applied separately to each
  coordinate shows nonnegative definiteness; products of these Gaussians and the
  constant function in the third coordinate are positive definite. More
  explicitly, for the Schwartz autocorrelation $\psi$ used in the preceding
  proof, choose $N > 1$ and bound $|\psi(x)| \le
  C_N\prod_{j=1}^d(1+|x_j|)^{-N}$. Then $\int\gamma(x)|\psi(x)|\ud x < \infty$,
  by the product of one-dimensional integrals and $a_1,a_2 < 1$. Tonelli with
  $|\psi|$ and then Fubini justify the two Gaussian-mixture integrals against
  $\psi$. Their integrands are nonnegative quadratic forms, proving the
  assertion also for the covariance independent of $x_3$. With $u = 1-y_1$ and
  then $y_2 = \sqrt{2u-u^2}\,w$, Tonelli gives
  \[
    \int G_1^{(2)}(\ud y)\gamma(y-e_1)
    =\frac{c B((1-a_2)/2,1/2)}{2\pi}
    \int_0^2u^{-a_1}(2u-u^2)^{-a_2/2}\ud u=\infty.
  \]

  \medskip
  \noindent\textit{Step 2.~}We transfer the product example to dimension three
  and remove smoothing.
  The beta factor is finite and the integral diverges at zero. The projection of
  $G_1^{(3)}$ onto its first two coordinates is $G_1^{(2)}$: the two sphere
  graphs each have surface Jacobian $(1-|y|^2)^{-1/2}$, divided by $4\pi$.
  Hence, the same divergence holds in dimension three. For radial Riesz
  covariance, the approximate-identity argument in Lemma~\ref{L:higher-riesz}
  applies. For the product covariance, Gaussian factorization and Tonelli give
  \[
    (p_\epsilon*\gamma)(x)
    =c\prod_{j=1}^2
    (p_\epsilon^{(1)}*|\,\cdot\,|^{-a_j})(x_j).
  \]
  In dimension three the third coordinate integrates to one. Each factor
  converges when $x_j \ne 0$ by the same argument with dimension one and
  exponent $a_j < 1$. This does not require the product kernel to be bounded at
  infinity; it is unbounded along its singular axes. The singular sets have zero
  shifted Green measure in both dimensions. Fatou with $\eta = 0$ implies that
  the smoothed convolution tends to infinity in lower limit as
  $\epsilon\downarrow0$. Even if both widths must be positive, take $\eta =
  \epsilon$ and use their sum in the Gaussian identity above. Thus, no uniform
  $M$ exists. This proves Proposition~\ref{P:higher-obstructions}.
\end{proof}

\begin{remark}[Calibrations and scope]\label{R:higher-calibrations}
  The printed first-moment statement and unconditional comparison are
  Theorem~\ref{T:first} and Proposition~\ref{P:Jensen}. Under the additional
  spatial hypothesis, Cauchy--Schwarz and Proposition~\ref{P:higher-transfer}
  give $\E u(t,x) \le \exp\{qK^{1/q}t^\kappa/2\}$, and hence $2L_* \le qK^{1/q}$
  if the intended first-moment lower bound holds. No sharp first- or
  second-moment constant is identified. As $\alpha_0\downarrow0$, the exponent
  tends to $(4-\alpha)/(3-\alpha)$, consistently with the Chen--Hu static
  calibration in Remark~\ref{R:d1-comparisons}; no interchange of limits is
  claimed. For the Riesz cases here, the Skorohod comparison of
  Balan--Conus~\cite{balan.conus:16:intermittency}*{Theorem~2.1(b), p.~1492;
  p.~1490} applies with $a = \alpha < d$, $\alpha < 2$, and $H = 1-\alpha_0/2 >
  1/2$. Its exponent $(2H+2-a)/(3-a)$ is $\kappa$. This comparison is restricted
  to those Riesz cases, with a different solution concept and covariance
  normalization; it is not an input here, and no equality with the
  Balan--Chen--Wang second-moment asymptotics is inferred. The obstructions
  above concern only the spatial supremum, not existence, moment finiteness, or
  failure of the target time scale. In particular, arbitrary homogeneous
  covariances remain outside this conditional argument, even when
  $\alpha+\alpha_0 < 2$.
\end{remark}

\section{Riesz bounds by interpolation with static noise}
\label{S:static-holder}

This alternative argument concerns isotropic Riesz covariance in dimensions two
and three. It does not use the preceding conditional proposition, settle general
anisotropic covariance, or identify the sharp constant. The full upper and
intrinsic formula use separate arguments.

We use the static Stratonovich model of Chen--Hu as one input and the
space-degree-one estimate of Section~\ref{S:d1-upper} as the other. The static
source allows nonnegative generalized covariances, not only Riesz functions: see
its standing assumptions (1.2)--(1.3), p.~2145, and its explicit
generalized-covariance convention (2.18), p.~2152,
in~\cite{chen.hu:25:hyperbolic}. In the latter convention every covariance
factor is heat regularized before its defining limit. A nonnegative, even,
positive-definite tempered Radon measure is such a generalized covariance. In
each application below we separately verify its homogeneity range and static
Dalang condition before using the source's quantitative coefficient estimates
(5.1) and (5.15). Positive definiteness of a weighted endpoint is proved in
Lemma~\ref{L:planar-endpoints}, or is an explicit assumption in
Hypothesis~\ref{H:three-endpoints}; weighting a measure alone is not being used
to assert that property. Chen's variational finiteness theorem supplies only
finiteness of the static constant, not a coefficient-rate theorem.

For $d \in \{2,3\}$, $c > 0$, and $0 < \lambda < \min\{2,d\}$, write
\[
  \gamma_\lambda(x)=c|x|^{-\lambda},\qquad
  \gamma_{\lambda,\epsilon}=p_\epsilon*\gamma_\lambda.
\]
All these covariances use the same amplitude $c$ and Gaussian convention. They
are nonnegative, locally integrable, and homogeneous. Their nonnegative
definiteness follows from
\[
  \gamma_\lambda(x)=\frac{c}{\Gamma(\lambda/2)}
  \int_0^\infty s^{\lambda/2-1}e^{-s|x|^2}\ud s.
\]
Indeed, each Gaussian is positive definite. For a Schwartz autocorrelation
$\psi$, the integral $\int |x|^{-\lambda}|\psi(x)|\ud x$ is finite by $\lambda <
d$ and rapid decay. Tonelli with $|\psi|$, followed by Fubini, therefore
justifies the positive Gaussian mixture as a tempered distribution. The source
hypotheses for the Riesz covariances below are thus satisfied.

\begin{lemma}[Static coefficient input]\label{L:static-holder-input}
  Let $U_n(t,x)$ be the ordered Stratonovich integrals for static Gaussian noise
  with covariance $\gamma_\chi$, where $d \in \{2,3\}$ and $1 < \chi < 2$. Put
  $U_0 = 1$ and
  \[
    B_k(t)=\sum_{m=0}^{2k}\E[U_m(t,x)U_{2k-m}(t,x)].
  \]
  Then $B_0 = 1$, and there is a finite constant $C_{\rm st} = C_{\rm
  st}(\chi,c,d) > 0$ such that
  \begin{equation}\label{E:static-holder-input}
    0\le B_k(t)\le
    \frac{C_{\rm st}^{\,k}t^{(4-\chi)k}}{(k!)^{3-\chi}},
    \qquad k\ge1,\quad t>0.
  \end{equation}
\end{lemma}

\begin{proof}
  We identify the full static coefficient sum using the imported fixed-order
  inputs. The imported upper coefficient rate controls all sufficiently large
  orders; the finitely many remaining orders are absorbed into one constant.

  \medskip
  \noindent\textit{Step 1.~}We verify the static input and pairing
  interpretation.
  The static solution has data $(1,0)$ and wave operator $\Delta$, as here. The
  degree satisfies $\chi < \min\{2,d\}$, so the Riesz covariance meets the
  source's Dalang condition. From
  Chen--Hu~\cite{chen.hu:25:hyperbolic}*{Theorem~3.8(i), p.~2165}, the
  fixed-order integrals exist in $L^2$. Their product rule and Wick expectation
  formula in that source's (2.17) and (2.19), pp.~2151--2152, identify $B_k$
  with the complete nonnegative sum over the $2k+1$ splits and every pairing of
  the $2k$ vertices. In particular, it includes within-chain pairs; no
  orthogonality of different orders is asserted. Empty chains contribute the
  factor one.

  \medskip
  \noindent\textit{Step 2.~}We convert the upper coefficient rate into an
  all-order bound.
  From~\cite{chen.hu:25:hyperbolic}*{(5.1), p.~2173, and (5.15), p.~2180},
  specialized to $p = 2$, we have
  \[
    B_k(t)=t^{(4-\chi)k}B_k(1),\qquad
    \limsup_{k\to\infty}\frac1k
    \log\left((k!)^{3-\chi}B_k(1)\right)<\infty.
  \]
  Only the finite upper coefficient rate from the latter formula is used. Choose
  $R > 0$ larger than this upper limit. For some $k_0 \ge 1$,
  $(k!)^{3-\chi}B_k(1) \le e^{Rk}$ for all $k \ge k_0$. Each of the remaining
  coefficients is finite by the fixed-order $L^2$ property. Increase $e^R$, if
  necessary, to a finite positive $C_{\rm st}$ that also dominates
  $\left((k!)^{3-\chi}B_k(1)\right)^{1/k}$ for $1 \le k < k_0$. Homogeneity
  gives~\eqref{E:static-holder-input}. The source results are imported inputs;
  this finite-order enlargement is the only deduction about their constant
  needed here. This proves Lemma~\ref{L:static-holder-input}.
\end{proof}

\begin{proposition}[Interpolation bound for Riesz covariance]
  \label{P:static-holder-upper}
  Suppose that $d \in \{2,3\}$, $c > 0$, $0 < \alpha_0 < 1$, and $1 < \alpha <
  2-\alpha_0$. Let the spatial covariance be $\gamma_\alpha(x) =
  c|x|^{-\alpha}$, and let $A_k(t)$ denote the full total-order-$2k$ coefficient
  of the second moment, as in Proposition~\ref{P:mixed}. Choose and fix
  \begin{equation}\label{E:static-holder-parameters}
    \alpha_0<\vartheta<2-\alpha,\qquad
    b_0=\frac{\alpha_0}{\vartheta},\qquad
    \chi=\frac{\alpha-\vartheta}{1-\vartheta}.
  \end{equation}
  Let $V_1$ be the spectral mass of the unit ball for $\gamma_1$, in the
  convention of Lemma~\ref{L:alpha-one-convolution}, and define
  \[
    Z_{b_0}=\frac2{(1-b_0)(2-b_0)},\qquad
    C_{\rm low}=2Z_{b_0}\pi V_1,\qquad
    C=C_{\rm low}^{\,\vartheta}C_{\rm st}^{\,1-\vartheta},
    \qquad q=3-\alpha,
  \]
  where $C_{\rm st}$ is supplied by Lemma~\ref{L:static-holder-input} for the
  displayed $\chi$. Then, for every $k \ge 1$, $t > 0$, and $x \in \R^d$,
  \begin{equation}\label{E:static-holder-upper}
    A_k(t)\le\frac{C^kt^{(4-\alpha-\alpha_0)k}}{(k!)^{3-\alpha}},
    \qquad
    \E u^2(t,x)\le\exp\left(qC^{1/q}t^\kappa\right).
  \end{equation}
  The finite constants may depend on $\alpha,\alpha_0,c,d$ and the fixed choice
  of $\vartheta$, but not on $k,t,x$.
\end{proposition}

\begin{proof}
  We interpolate the regularized covariance factors on one common positive
  measure containing all pairings and splits. We bound the low endpoint
  directly, pass fixed-order limits for the target and static endpoint, and
  only then use the static coefficient estimate.

  \medskip
  \noindent\textit{Step 1.~}We verify the endpoint parameters.
  The interval for $\vartheta$ is nonempty by $\alpha+\alpha_0 < 2$. Since
  $\alpha > 1$, we have $0 < \vartheta < 1$, $0 < b_0 < 1$, and
  \[
    \chi-1=\frac{\alpha-1}{1-\vartheta}>0,\qquad
    2-\chi=\frac{2-\alpha-\vartheta}{1-\vartheta}>0.
  \]
  Thus, both endpoint covariances are admissible in dimensions two and three.
  Moreover,
  \begin{equation}\label{E:static-holder-affine}
    \alpha=\vartheta+(1-\vartheta)\chi,
    \qquad \alpha_0=\vartheta b_0.
  \end{equation}

  \medskip
  \noindent\textit{Step 2.~}We place all three covariance profiles on the same
  measure.
  Fix $k \ge 1$ and $t > 0$. Use a common spatial noise width $\epsilon > 0$ and
  the target's temporal cutoff width $\delta > 0$. The covariance on each
  pairing edge has spatial width $2\epsilon$. In dimension three first use
  $G_{h,\eta} = p_\eta*G_h$ with $\eta > 0$; in dimension two take $\eta = 0$.
  For each split $m+n = 2k$ and pairing $D \in \Pi_k$, use the two separate
  roots, increments $h_l = s_l-s_{p(l)}$, and positions $x_l$ of
  Section~\ref{S:d1-upper}. Define a common positive measure $\nu_{k,t,\eta}$ by
  \[
    \int\Phi\ud\nu_{k,t,\eta}
    :=\sum_{m=0}^{2k}\sum_{D\in\Pi_k}
    \int_{[0,t]^m_<\times[0,t]^{2k-m}_<}\ud s
    \int\Phi(m,D,s,y)\prod_{l=1}^{2k}G_{h_l,\eta}(\ud y_l).
  \]
  Here $y_l$ are vector increments and each $x_l$ is their cumulative sum within
  its own chain. Empty blocks have empty products and their own root
  conventions. This measure is finite: the spatial mass at fixed times is
  $\prod_lh_l$, the time domains are bounded, and both discrete sums are finite.
  It is the same measure for all three covariance profiles.

  For $\lambda \in \{1,\alpha,\chi\}$ and $b > 0$, put
  \[
    H_\lambda=\prod_{(i,j)\in D}
    \gamma_{\lambda,2\epsilon}(x_i-x_j),\qquad
    T_b=\prod_{(i,j)\in D}|s_i-s_j|^{-b}.
  \]
  The pointwise identity $\gamma_\alpha =
  \gamma_1^{\vartheta}\gamma_\chi^{1-\vartheta}$ retains the amplitude $c$,
  since $c^\vartheta c^{1-\vartheta} = c$. H\"older's inequality on the Gaussian
  probability measure gives
  \[
    \gamma_{\alpha,2\epsilon}(z)
    \le\gamma_{1,2\epsilon}(z)^\vartheta
    \gamma_{\chi,2\epsilon}(z)^{1-\vartheta},
    \qquad z\in\R^d.
  \]
  These convolutions are finite because the degrees are below $d$. Let
  $Q_\delta$ denote the product of the target's truncated temporal covariances.
  From the regularization in Chen~\cite{chen:author2026:time-dependency}*{(2.9),
  (2.11)},
  \[
    0\le Q_\delta\le T_{\alpha_0}=T_{b_0}^{\,\vartheta}.
  \]
  This uses only domination by the original temporal power; no interpolation
  property of distinct temporal cutoffs is required.

  Let $A_k^{\epsilon,\delta,\eta}(t) = \int Q_\delta H_\alpha \ud\nu_{k,t,\eta}$
  be the target regularized coefficient. The Gaussian moment formula and the
  separate chain reversals in Proposition~\ref{P:mixed} identify it with the
  finite sum of regularized mixed moments. Also set
  \[
    L_k^{\epsilon,\eta}(t)=\int T_{b_0}H_1\ud\nu_{k,t,\eta},
    \qquad
    B_k^{\epsilon,\eta}(t)=\int H_\chi\ud\nu_{k,t,\eta}.
  \]
  The second expression is the regularized static coefficient, by the Gaussian
  moment formula with all pairings and splits. The first is a deterministic
  integral with original temporal covariance. We do not identify it with a
  stochastic integral having mixed regularizations. H\"older's inequality on the
  full measure, including its two discrete variables, yields
  \begin{equation}\label{E:static-holder-regularized}
    A_k^{\epsilon,\delta,\eta}(t)
    \le\left(L_k^{\epsilon,\eta}(t)\right)^\vartheta
    \left(B_k^{\epsilon,\eta}(t)\right)^{1-\vartheta}.
  \end{equation}

  \medskip
  \noindent\textit{Step 3.~}We bound the low endpoint before passing limits.
  The deterministic estimates in the proof of
  Proposition~\ref{P:low-alpha-upper}, at spatial degree one and temporal degree
  $b_0$, apply to $L_k^{\epsilon,\eta}$. In detail,
  Lemma~\ref{L:alpha-one-convolution} gives the closing bound $\pi V_1/2$,
  uniformly in both smoothing widths. Descending integration and the opening
  arithmetic--geometric mean give the common spatial bound $(\pi V_1)^kt^k/k^k$.
  The temporal pairing sum is symmetric on the full time cube. Its integral,
  after summing every split as in~\eqref{E:d1-time-symmetry}, is
  $2^k(Z_{b_0}t^{2-b_0})^k/k!$. Hence,
  \begin{equation}\label{E:static-holder-low}
    L_k^{\epsilon,\eta}(t)
    \le\frac{C_{\rm low}^{\,k}t^{(3-b_0)k}}{k!\,k^k}
    \le\frac{C_{\rm low}^{\,k}t^{(3-b_0)k}}{(k!)^2}.
  \end{equation}
  The original temporal integral is finite since $b_0 < 1$. The static
  regularized coefficient is finite because its covariance is bounded at fixed
  $\epsilon$. Thus, the H\"older step has finite right side.

  \medskip
  \noindent\textit{Step 4.~}We remove Green and noise regularization in order.
  We pass to fixed-order limits before using the static coefficient bound. In
  dimension three fix $\epsilon,\delta$ and let $\eta\downarrow0$. The target
  and static smoothed covariance factors are bounded continuous functions of the
  spatial positions. Products of the finite Green measures converge weakly with
  masses $\prod_lh_l$. A constant times this product dominates each time
  integrand at fixed noise smoothing, so dominated convergence applies on the
  bounded time simplexes. The same argument with two Green widths gives a common
  mixed covariance limit, and hence fixed-order $L^2$ convergence. This
  identifies the static limit with the measure convention of
  Chen--Hu~\cite{chen.hu:25:hyperbolic}*{(2.22), pp.~2152--2153}. The target
  uses the corresponding inner limit in
  Chen~\cite{chen:author2026:time-dependency}*{pp.~13--14}. We need no limit
  identification for $L_k^{\epsilon,\eta}$: its uniform upper bound
  in~\eqref{E:static-holder-low} suffices.

  Next let the noise-smoothing parameters tend to zero. From
  Chen--Hu~\cite{chen.hu:25:hyperbolic}*{Theorem~3.8(i), p.~2165}, each fixed
  static order converges in $L^2$ to $U_n$. For the target, the fixed-order
  $L^2$ convergence follows from Lemma~\ref{L:joint-fixed-order}. If $Y,Y',Z,Z'$
  are the corresponding approximations and limits, then
  \[
    \|YZ-Y'Z'\|_1
    \le\|Y-Y'\|_2\|Z\|_2+\|Y'\|_2\|Z-Z'\|_2.
  \]
  Therefore, every fixed split product converges in $L^1$. The finite sums over
  the $2k+1$ splits converge as well; a zero-order factor stays equal to one.
  This argument applies directly in dimension two, where no Green smoothing is
  needed. Passing to these limits in~\eqref{E:static-holder-regularized},
  using~\eqref{E:static-holder-low}, gives
  \[
    A_k(t)\le
    \left(\frac{C_{\rm low}^{\,k}t^{(3-b_0)k}}{(k!)^2}\right)^\vartheta
    B_k(t)^{1-\vartheta}.
  \]
  Only now apply Lemma~\ref{L:static-holder-input}. No monotonicity in
  $\epsilon$ and no uniform static large-time bound in that parameter have been
  assumed.

  \medskip
  \noindent\textit{Step 5.~}We combine the endpoint powers and sum.
  The powers combine without a counting loss. Indeed,
  by~\eqref{E:static-holder-affine},
  \[
    \begin{split}
      2\vartheta+(1-\vartheta)(3-\chi)&=3-\alpha,\\
      \vartheta(3-b_0)+(1-\vartheta)(4-\chi)
      &=4-\alpha-\alpha_0.
    \end{split}
  \]
  This proves the coefficient bound in~\eqref{E:static-holder-upper}. The
  constant term is $A_0 = 1$. Since $q = 3-\alpha > 1$, the scalar inequality
  $\sum_{k\ge0}z^k/(k!)^q \le e^{qz^{1/q}}$ proved in
  Proposition~\ref{P:d1-upper} applies with $z = Ct^{4-\alpha-\alpha_0}$.
  Proposition~\ref{P:mixed} identifies the nonnegative full coefficient series
  with $\E u^2(t,x)$. Here $(4-\alpha-\alpha_0)/q = \kappa$. This proves
  Proposition~\ref{P:static-holder-upper}.
\end{proof}

\begin{remark}[Calibrations and remaining scope]\label{R:static-holder-scope}
  The time exponent in~\eqref{E:static-holder-upper} tends to
  $(4-\alpha)/(3-\alpha)$ as $\alpha_0\downarrow0$, in agreement with the
  Chen--Hu static exponent. This does not interchange the static and large-time
  limits or assert uniformity of the constants. Cauchy--Schwarz gives $\E u(t,x)
  \le \exp(qC^{1/q}t^\kappa/2)$. If the first-moment lower rate is at least
  $L_*$, this upper bound implies $2L_* \le qC^{1/q}$. Thus, the first-moment
  and Jensen calibrations are retained, without identifying a sharp constant.
  The Skorohod comparison remains only the comparison recorded in
  Remark~\ref{R:higher-calibrations}; its exponent agrees, but its solution
  concept and covariance normalization differ. No Skorohod estimate is used in
  this proof. The argument covers isotropic Riesz covariance throughout the
  stated higher-degree range in dimensions two and three, including the
  previously uncovered two-dimensional interval $3/2 \le \alpha < 2-\alpha_0$.
  This argument does not cover arbitrary homogeneous covariances: anisotropic
  profiles need not factor into the admissible endpoints. The intrinsic
  second-moment formula is given separately in
  Theorem~\ref{T:intrinsic-second-moment}. No novelty is asserted here.
\end{remark}

\chapter{Covariance comparisons in higher dimensions}\label{S:appendix-covariance}

Suppose a bound is known for simple covariances, such as products of
one-dimensional ones. Can it be carried over to a general homogeneous
covariance in dimension two or three? Sometimes. A pointwise factorization is
not enough: each factor must itself be a covariance, and the comparison must
keep the amplitudes and the scaling. This matters most for measure-valued
and anisotropic covariances. This appendix has no figure.

We keep Hypothesis~\ref{H:homog} and the mixed pairing formula
of Chapter~\ref{S:G0}. Section~\ref{S:planar-holder}
builds planar endpoints and proves the transfer
(Proposition~\ref{P:planar-upper}). Section~\ref{S:three-endpoints}
gives product constructions in dimension three and a conditional transfer
(Proposition~\ref{P:three-transfer}). It does not claim a factorization
for every three-dimensional covariance. Section~\ref{S:energy-fidelity}
identifies the physical and spectral forms of the spatial energy for measure
covariances; this identity does not remove the assumptions of the transfers.

The upper bound of Chapter~\ref{S:book-independent-upper} and the
constant of Chapter~\ref{S:book-constant} do not rely on these
comparisons. On a first reading, the two transfer statements are
enough. Episode~\VideoRef{9}{eLz9Hd20ZT0} of the video series covers this
appendix.

\section{General homogeneous covariance in dimension two}
\label{S:planar-holder}

This section gives planar endpoint constructions and a transfer under its stated
hypotheses. It does not use the preceding Riesz interpolation argument as an
input. The full upper and intrinsic formula use separate methods.

Throughout this section $d = 2$. Write $\Gamma$ for the spatial covariance as a
nonnegative, even, positive-definite tempered Radon measure. For a locally
integrable covariance function, this means $\Gamma(\ud x) = \gamma(x)\ud x$. The
measure interpretation uses the Gaussian approximation framework of
Chen~\cite{chen:author2026:time-dependency}*{pp.~2--3, 13} and
Chen--Hu~\cite{chen.hu:25:hyperbolic}*{pp.~2145--2147, 2151--2153}. All products
of covariance factors below are taken after Gaussian smoothing; no product or
fractional power of distributions is used. Homogeneity means
\begin{equation}\label{E:planar-measure-scaling}
  \Gamma(rA)=r^{2-\alpha}\Gamma(A),\qquad r>0,
\end{equation}
for Borel sets $A$. This is equivalent to the covariance homogeneity in
Hypothesis~\ref{H:homog}. If $\Gamma = 0$, the noise vanishes and $u = 1$. We
first treat $\Gamma \ne 0$ and $1 < \alpha < 2$.

\begin{lemma}[Planar endpoint covariances]\label{L:planar-endpoints}
  Suppose that $\Gamma$ satisfies the preceding assumptions, and let $\alpha <
  \chi < 2$. Define the finite angular measure
  \[
    \sigma(E)=(2-\alpha)\Gamma
    \{x:0<|x|\le1,\ x/|x|\in E\},\qquad
    S=\sigma(\mathbb{S}^1).
  \]
  On $\R^2\setminus\{0\}$ define the weighted measures
  \begin{equation}\label{E:planar-endpoints}
    \Gamma_{\rm low}(\ud x)=|x|^{\alpha-1}\Gamma(\ud x),\qquad
    \Gamma_{\rm st}(\ud x)=|x|^{\alpha-\chi}\Gamma(\ud x),
  \end{equation}
  with zero mass at the origin. Both are nonnegative, even, positive-definite
  tempered Radon measures. Their measure scaling exponents are $1$ and $2-\chi$,
  respectively. For every $h > 0$, $\epsilon > 0$, $\eta \ge 0$, and $z \in
  \R^2$,
  \begin{equation}\label{E:planar-line-convolution}
    \int G_{h,\eta}(\ud y)(p_\epsilon*\Gamma_{\rm low})(y+z)
    \le\frac S4.
  \end{equation}
  If $\vartheta = (\chi-\alpha)/(\chi-1)$, then $0 < \vartheta < 1$ and
  \begin{equation}\label{E:planar-gaussian-holder}
    (p_\epsilon*\Gamma)(z)
    \le(p_\epsilon*\Gamma_{\rm low})(z)^\vartheta
    (p_\epsilon*\Gamma_{\rm st})(z)^{1-\vartheta}.
  \end{equation}
  No angular density or isotropy is assumed.
\end{lemma}

\begin{proof}
  We identify the physical polar measure and its two weighted endpoints. We
  prove their positive definiteness separately, estimate the low endpoint
  through line integrals, and apply Hölder on the original heat-weighted
  measure.

  \medskip
  \noindent\textit{Step 1.~}We obtain the polar measure and endpoint
  integrability.
  Since $2-\alpha > 0$, scaling and local finiteness force $\Gamma(\{0\}) = 0$.
  The unit-ball mass is finite, so $S < \infty$. It is positive: otherwise
  scaling would make every ball have zero mass, contradicting $\Gamma \ne 0$.
  Evenness of $\Gamma$ implies evenness of $\sigma$. For $0 < a < b$ and Borel
  $E \subset \mathbb{S}^1$,
  \[
    \Gamma\{x:a<|x|\le b,\ x/|x|\in E\}
    =\frac{b^{2-\alpha}-a^{2-\alpha}}{2-\alpha}\sigma(E)
    =\int_a^b r^{1-\alpha}\ud r\,\sigma(E).
  \]
  These polar rectangles form a generating family. Uniqueness of measures, first
  on finite annuli and then on their increasing union, identifies the polar
  pushforward with $r^{1-\alpha}\ud r\,\sigma(\ud\omega)$. Hence, for every
  nonnegative Borel function $f$, Tonelli gives
  \begin{equation}\label{E:planar-polar}
    \int f(x)\Gamma(\ud x)
    =\int_{\mathbb{S}^1}\int_0^\infty
    f(r\omega)r^{1-\alpha}\ud r\,\sigma(\ud\omega).
  \end{equation}
  Extended values are permitted. In particular, the formula applies to the
  nonnegative wave density used below.

  The endpoint polar measures are therefore $\ud r\,\sigma(\ud\omega)$ and
  $r^{1-\chi}\ud r\,\sigma(\ud\omega)$. Their radial integrals are finite near
  zero because $\chi < 2$. Rapid decay of a Schwartz function controls their
  integrals at infinity, uniformly in direction. They are thus locally finite
  and tempered, and their stated scaling and evenness follow directly. If
  $\Gamma$ has a density, the two weighted measures also have locally integrable
  densities.

  \medskip
  \noindent\textit{Step 2.~}We prove positive definiteness of the low
  endpoint.
  To prove positive definiteness of the low endpoint, evenness gives the
  full-sphere line representation
  \begin{equation}\label{E:planar-line-mixture}
    \int f(x)\Gamma_{\rm low}(\ud x)
    =\frac12\int_{\mathbb{S}^1}\int_\R f(r\omega)\ud r
    \,\sigma(\ud\omega).
  \end{equation}
  It holds for nonnegative Borel $f$, and for Schwartz $f$ by absolute
  integrability. For a complex Schwartz function $\varphi$, put
  $\widetilde{\varphi}(x) = \overline{\varphi(-x)}$ and $\psi =
  \varphi*\widetilde{\varphi}$. Fubini in line and perpendicular coordinates
  gives
  \[
    \int_\R\psi(r\omega)\ud r
    =\int_{\omega^\perp}
    \left|\int_\R\varphi(s\omega+v)\ud s\right|^2\ud v\ge0.
  \]
  Schwartz decay justifies all these integrals, uniformly for unit $\omega$.
  Integrating against the finite measure $\sigma$
  in~\eqref{E:planar-line-mixture} proves $\Gamma_{\rm low}(\psi) \ge 0$. The
  factor $1/2$ is essential because both orientations of each line occur on the
  full sphere.

  \medskip
  \noindent\textit{Step 3.~}We prove positive definiteness of the static
  endpoint.
  For the static endpoint put $\lambda = \chi-\alpha > 0$. For $x \ne 0$,
  \[
    |x|^{-\lambda}=\frac1{\Gamma(\lambda/2)}
    \int_0^\infty s^{\lambda/2-1}e^{-s|x|^2}\ud s.
  \]
  For each $s > 0$, the measure $e^{-s|x|^2}\Gamma(\ud x)$ is positive definite.
  Indeed, the Fourier transform of $e^{-s|x|^2}\psi(x)$ is a positive Gaussian
  convolution of $|\widehat{\varphi}|^2$. Its pairing with the original positive
  spectral measure is nonnegative; Schwartz decay makes the pairing finite.
  Moreover,
  \[
    \int |\psi(x)||x|^{-\lambda}\Gamma(\ud x)
    =\int |\psi(x)|\Gamma_{\rm st}(\ud x)<\infty
  \]
  by the endpoint polar integrability just proved. Tonelli with $|\psi|$,
  followed by Fubini, justifies the Gaussian mixture against $\psi$. Thus,
  $\Gamma_{\rm st}(\psi) \ge 0$. This uses smooth Gaussian multiplication of the
  original covariance and a convergent mixture, not a fractional power of that
  covariance.

  \medskip
  \noindent\textit{Step 4.~}We bound the planar wave convolution.
  We next prove the convolution estimate directly. Use the planar wave density
  \[
    G_h(y)=\frac{\mathbf1_{\{|y|<h\}}}{2\pi\sqrt{h^2-|y|^2}},
  \]
  with value zero on the circle. Write $z = a\omega+b\omega^\perp$ for a unit
  perpendicular vector. If $|b| < h$, put $A = \sqrt{h^2-b^2}$. Translation
  along the line gives
  \[
    \int_\R G_h(z-r\omega)\ud r
    =\frac1{2\pi}\int_{-A}^A\frac{\ud s}{\sqrt{A^2-s^2}}
    =\frac12.
  \]
  If $|b| \ge h$, the open support is empty or tangent and the line integral is
  zero. The tangency point has zero line Lebesgue measure.
  Formula~\eqref{E:planar-line-mixture} and nonnegative Tonelli yield
  \[
    (G_h*\Gamma_{\rm low})(z)
    =\frac14\int_{\mathbb{S}^1}
    \mathbf1_{\{|z\cdot\omega^\perp|<h\}}\sigma(\ud\omega)
    \le\frac S4.
  \]
  Evenness, Gaussian convolution, and Tonelli now give
  \[
    \int G_{h,\eta}(\ud y)(p_\epsilon*\Gamma_{\rm low})(y+z)
    =\left[p_{\epsilon+\eta}*(G_h*\Gamma_{\rm low})\right](z)
    \le\frac S4.
  \]
  For $\eta = 0$, use the unsmoothed Green measure. The Gaussian average has
  mass one, so the same constant works for all widths.

  \medskip
  \noindent\textit{Step 5.~}We compare the regularized endpoints by Hölder.
  Finally, the measure $p_\epsilon(z-x)\Gamma(\ud x)$ is finite. Gaussian decay
  and polynomial radial growth prove this directly; the endpoint convolutions
  are finite for the same reason. The identity
  \[
    \vartheta(\alpha-1)+(1-\vartheta)(\alpha-\chi)=0
  \]
  and absence of an origin atom permit H\"older's inequality on this
  unnormalized positive measure:
  \[
    \begin{split}
      (p_\epsilon*\Gamma)(z)
      &=\int |x|^{\vartheta(\alpha-1)}
      |x|^{(1-\vartheta)(\alpha-\chi)}
      p_\epsilon(z-x)\Gamma(\ud x)\\
      &\le(p_\epsilon*\Gamma_{\rm low})(z)^\vartheta
      (p_\epsilon*\Gamma_{\rm st})(z)^{1-\vartheta}.
    \end{split}
  \]
  This proves Lemma~\ref{L:planar-endpoints}.
\end{proof}

\begin{lemma}[General planar static coefficient input]
  \label{L:planar-static-input}
  Let $\Xi$ be a nonnegative, even, positive-definite tempered Radon measure on
  $\R^2$ with $\Xi(rA) = r^{2-\chi}\Xi(A)$, where $1 < \chi < 2$. For static
  Stratonovich noise with covariance $\Xi$, data $(1,0)$, and wave operator
  $\Delta$, let $U_n(t,x)$ be its ordered integrals, with $U_0 = 1$. There is a
  finite $C_{\rm st} = C_{\rm st}(\Xi) > 0$ such that, for $k \ge 1$ and $t >
  0$,
  \begin{equation}\label{E:planar-static-input}
    0\le B_k(t):=\sum_{m=0}^{2k}\E[U_m(t,x)U_{2k-m}(t,x)]
    \le\frac{C_{\rm st}^{\,k}t^{(4-\chi)k}}{(k!)^{3-\chi}}.
  \end{equation}
  Also $B_0 = 1$.
\end{lemma}

\begin{proof}
  We treat zero covariance directly. For nonzero covariance, homogeneity
  verifies the static spectral condition, and the cited finite coefficient
  rate is converted to a bound valid at every order.

  \medskip
  \noindent\textit{Step 1.~}We verify the static covariance hypotheses.
  If $\Xi = 0$, then $U_n = 0$ for $n \ge 1$, and hence $B_k = 0$ for $k \ge 1$.
  The choices $C_{\rm st} = 1$ and $B_0 = 1$ give the conclusion. We may
  therefore assume $\Xi \ne 0$.

  Fourier uniqueness transfers covariance homogeneity to its positive spectral
  measure $\mu_\Xi$: $\mu_\Xi(rA) = r^\chi\mu_\Xi(A)$. Its unit-ball mass is
  finite because the spectral measure is Radon. On dyadic annuli this gives
  \[
    \int_{\R^2}(1+|\xi|^2)^{-1}\mu_\Xi(\ud\xi)
    \le\mu_\Xi(B_1)
    +2^\chi\mu_\Xi(B_1)\sum_{j\ge0}2^{j(\chi-2)}<\infty.
  \]
  This verifies the static Dalang condition without an angular density. The
  present $\Xi$ also satisfies the source's covariance assumptions: it is a
  nonnegative positive-definite tempered distribution, and $1 < \chi < 2 = d$.
  The source expressly includes generalized covariances and uses their
  heat-regularized interpretation \cite{chen.hu:25:hyperbolic}*{p.~2152,
  (2.18)}. Thus, the coefficient estimate invoked below is a direct application
  within the source's stated class, not an extension from density kernels. From
  Chen--Hu~\cite{chen.hu:25:hyperbolic}*{Theorem~3.8(i), p.~2165}, every fixed
  order exists in $L^2$. The product and Wick formulas (2.17) and (2.19),
  pp.~2151--2152, identify $B_k$ with all pairings and all splits of the total
  order, including self-pairs and empty chains. Thus, every coefficient is
  finite and nonnegative.

  \medskip
  \noindent\textit{Step 2.~}We state the imported finiteness input used by the
  coefficient rate.
  We make explicit the finiteness input in the source's coefficient upper rate.
  The variational constant $\mathcal{M}$ in
  Chen--Hu~\cite{chen.hu:25:hyperbolic}*{p.~2147} is finite under the displayed
  Dalang condition. The source's invocation on p.~2175 identifies
  $2^{\chi/(4-\chi)}\mathcal{M}$ with the square-root variational quantity
  covered by Chen~\cite{chen:23:exponential}*{Theorem~1.1, p.~2}. That theorem
  explicitly states finiteness under Dalang's condition and permits generalized
  nonnegative positive-definite covariances. We import that finiteness, not a
  sharp value or its proof.

  \medskip
  \noindent\textit{Step 3.~}We turn the rate into an all-order constant.
  From Chen--Hu~\cite{chen.hu:25:hyperbolic}*{(5.1), p.~2173, and (5.15),
  p.~2180}, with $p = 2$, it follows that
  \[
    B_k(t)=t^{(4-\chi)k}B_k(1),\qquad
    \limsup_{k\to\infty}\frac1k
    \log\left((k!)^{3-\chi}B_k(1)\right)<\infty.
  \]
  Choose a finite $R > 0$ larger than this upper limit. For all sufficiently
  large $k$, $(k!)^{3-\chi}B_k(1) \le e^{Rk}$. Increase $e^R$ to a finite
  positive constant dominating the $k$th roots of the finitely many remaining
  quantities. The resulting $C_{\rm st}$ proves~\eqref{E:planar-static-input} by
  homogeneity. These are the source's general homogeneous-covariance statements,
  not a Riesz specialization. This proves Lemma~\ref{L:planar-static-input}.
\end{proof}

\begin{proposition}[Planar homogeneous-covariance upper bound]
  \label{P:planar-upper}
  Let $d = 2$, $0 < \alpha_0 < 1$, and $1 < \alpha < 2-\alpha_0$. Suppose that
  the nonzero spatial covariance $\Gamma$ is nonnegative, even, positive
  definite and tempered, with~\eqref{E:planar-measure-scaling}. Choose and fix
  \[
    \alpha_0<\vartheta<2-\alpha,\qquad
    b_0=\alpha_0/\vartheta,\qquad
    \chi=(\alpha-\vartheta)/(1-\vartheta).
  \]
  Use the endpoint measures and physical angular mass $S$ of
  Lemma~\ref{L:planar-endpoints}. Put
  \[
    Z_{b_0}=\frac2{(1-b_0)(2-b_0)},\qquad
    C_{\rm low}=Z_{b_0}S,\qquad
    C=C_{\rm low}^{\,\vartheta}C_{\rm st}^{\,1-\vartheta},
    \qquad q=3-\alpha,
  \]
  where $C_{\rm st}$ is from Lemma~\ref{L:planar-static-input} for $\Xi =
  \Gamma_{\rm st}$. For the Stratonovich solution with data $(1,0)$, every $k
  \ge 1$, $t > 0$, and $x \in \R^2$ satisfy
  \begin{equation}\label{E:planar-upper}
    A_k(t)\le\frac{C^kt^{(4-\alpha-\alpha_0)k}}{(k!)^{3-\alpha}},
    \qquad
    \E u^2(t,x)\le\exp\left(qC^{1/q}t^\kappa\right).
  \end{equation}
  Constants may depend on the covariance, $\alpha_0$, and the fixed choice of
  $\vartheta$, but not on $k,t,x$.
\end{proposition}

\begin{proof}
  We verify the planar endpoint parameters and use a common positive measure
  for all diagrams. We interpolate the regularized coefficients, bound the low
  endpoint by its line-convolution estimate, and then pass to the fixed-order
  limits before summing.

  \medskip
  \noindent\textit{Step 1.~}We verify the original and endpoint spectral
  conditions.
  The parameter interval is nonempty. We have $0 < \vartheta < 1$, $0 < b_0 <
  1$, and
  \[
    \chi-\alpha=\frac{\vartheta(\alpha-1)}{1-\vartheta}>0,
    \qquad
    2-\chi=\frac{2-\alpha-\vartheta}{1-\vartheta}>0.
  \]
  Thus, Lemma~\ref{L:planar-endpoints} applies, with its interpolation weight
  equal to $\vartheta$. The original spectral measure has ball growth of
  exponent $\alpha$, by the same Fourier homogeneity argument used in
  Lemma~\ref{L:planar-static-input}. A dyadic sum gives the time-dependent
  spectral condition (1.5) of
  Chen~\cite{chen:author2026:time-dependency}*{p.~3}, since $\alpha <
  2-\alpha_0$. The low endpoint has spectral ball exponent one and $b_0 < 1$;
  the static endpoint has exponent $\chi < 2$. Thus,
  Lemma~\ref{L:joint-fixed-order} applies to the target, and the imported static
  fixed-order theorem is used within its stated range.

  \medskip
  \noindent\textit{Step 2.~}We form the common measure and apply Hölder.
  Fix $k \ge 1$, $t > 0$, and common noise widths $\epsilon,\delta > 0$. The
  spatial covariance on each edge has width $2\epsilon$. No Green smoothing is
  needed in dimension two. On the disjoint union over $m = 0,\ldots,2k$ and $D
  \in \Pi_k$, define the finite positive measure $\nu_{k,t}$ by
  \[
    \int\Phi\ud\nu_{k,t}
    :=\sum_{m=0}^{2k}\sum_{D\in\Pi_k}
    \int_{[0,t]^m_<\times[0,t]^{2k-m}_<}\ud s
    \int\Phi(m,D,s,y)\prod_{l=1}^{2k}G_{h_l}(\ud y_l).
  \]
  Here each chain has its own root, $h_l = s_l-s_{p(l)}$, the $y_l$ are its
  increments, and $x_l$ is its cumulative position. A zero-length chain has
  empty products. At fixed times the spatial mass is $\prod_lh_l$; bounded time
  simplexes and finite discrete sums therefore make $\nu_{k,t}$ finite.

  For $J \in \{\Gamma,\Gamma_{\rm low},\Gamma_{\rm st}\}$ write
  \[
    H_J=\prod_{(i,j)\in D}(p_{2\epsilon}*J)(x_i-x_j),\qquad
    T_b=\prod_{(i,j)\in D}|s_i-s_j|^{-b}.
  \]
  The Gaussian covariance functions are bounded and continuous at each fixed
  width: their positive spectral measures are tempered, so Gaussian weighting
  makes them finite. Let $Q_\delta$ be the product of the target's truncated
  temporal covariances. The source regularization gives $0 \le Q_\delta \le
  T_{\alpha_0} = T_{b_0}^{\,\vartheta}$. Lemma~\ref{L:planar-endpoints}, at the
  common edge width, gives $H_\Gamma \le H_{\Gamma_{\rm low}}^\vartheta
  H_{\Gamma_{\rm st}}^{1-\vartheta}$.

  The Gaussian moment formula and the two separate chain reversals in
  Proposition~\ref{P:mixed} identify the target regularized coefficient with
  $\int Q_\delta H_\Gamma\ud\nu_{k,t}$. Define
  \[
    L_k^\epsilon(t)=\int T_{b_0}H_{\Gamma_{\rm low}}\ud\nu_{k,t},
    \qquad
    B_k^\epsilon(t)=\int H_{\Gamma_{\rm st}}\ud\nu_{k,t}.
  \]
  The latter is the regularized static coefficient, including every split and
  pairing. The former is only a deterministic integral with original temporal
  power. H\"older's inequality on the full common measure, including its
  discrete variables, yields
  \begin{equation}\label{E:planar-coefficient-holder}
    A_k^{\epsilon,\delta}(t)
    \le\left(L_k^\epsilon(t)\right)^\vartheta
    \left(B_k^\epsilon(t)\right)^{1-\vartheta}.
  \end{equation}
  No interpolation of different temporal cutoff kernels is asserted.

  \medskip
  \noindent\textit{Step 3.~}We estimate the low endpoint with the full
  temporal pairing sum.
  We bound $L_k^\epsilon$ directly. Orient each pair by increasing block label.
  Descending integration of increments, as in Lemma~\ref{L:d1-spatial}, removes
  at a closing label its only remaining covariance factor.
  Lemma~\ref{L:planar-endpoints} bounds that integral by $S/4$. At an opening
  label its paired covariance has already been removed, so the Green integral is
  its mass $h_l$. The $k$ opening increments sum to at most $2t$. Thus, the
  common spatial bound is
  \[
    (S/4)^k\prod_{i\in O_D}h_i
    \le(S/2)^k\frac{t^k}{k^k}.
  \]
  The temporal pairing sum is invariant under every permutation of its $2k$
  labels. Each pairing integrates on the full time cube to
  $(Z_{b_0}t^{2-b_0})^k$. Partitioning by the two within-block orders, and then
  summing all splits, gives
  \[
    \begin{split}
      \sum_{m=0}^{2k}\sum_{D\in\Pi_k}
      \int_{[0,t]^m_<\times[0,t]^{2k-m}_<}T_{b_0}\ud s
      &=\frac{(2k)!}{2^kk!}(Z_{b_0}t^{2-b_0})^k
      \sum_{m=0}^{2k}\frac1{m!(2k-m)!}\\
      &=\frac{2^k}{k!}(Z_{b_0}t^{2-b_0})^k.
    \end{split}
  \]
  Symmetry is used for the sum, not for an individual pairing. The integral is
  finite because $b_0 < 1$. Multiplication gives
  \begin{equation}\label{E:planar-low-coefficient}
    L_k^\epsilon(t)
    \le\frac{(Z_{b_0}S)^kt^{(3-b_0)k}}{k!\,k^k}
    \le\frac{C_{\rm low}^{\,k}t^{(3-b_0)k}}{(k!)^2}.
  \end{equation}
  At fixed $\epsilon$, bounded covariance and the finite common measure make
  $B_k^\epsilon$ finite as well. These observations also justify finiteness in
  the H\"older estimate.

  \medskip
  \noindent\textit{Step 4.~}We remove noise smoothing at fixed order.
  We next take fixed-order limits. The target orders converge in $L^2$ by
  Lemma~\ref{L:joint-fixed-order}. For the static orders the corresponding input
  is Chen--Hu~\cite{chen.hu:25:hyperbolic}*{Theorem~3.8(i), p.~2165}. For
  products of two fixed orders, the inequality
  \[
    \|YZ-Y'Z'\|_1
    \le\|Y-Y'\|_2\|Z\|_2+\|Y'\|_2\|Z-Z'\|_2
  \]
  gives $L^1$ convergence. There are only $2k+1$ split products, and a
  zero-order factor stays equal to one. Hence, $A_k^{\epsilon,\delta}(t) \to
  A_k(t)$ and $B_k^\epsilon(t) \to B_k(t)$ as the noise widths vanish. The
  uniform deterministic estimate~\eqref{E:planar-low-coefficient} needs no
  stochastic limit interpretation. Passing to limits
  in~\eqref{E:planar-coefficient-holder} gives
  \[
    A_k(t)\le
    \left(\frac{C_{\rm low}^{\,k}t^{(3-b_0)k}}{(k!)^2}\right)^\vartheta
    B_k(t)^{1-\vartheta}.
  \]
  Only now apply Lemma~\ref{L:planar-static-input}. No monotonicity in the
  spatial width or uniform static smoothed asymptotic was used.

  \medskip
  \noindent\textit{Step 5.~}We combine powers and sum the full moment series.
  The parameter identities $\alpha = \vartheta+(1-\vartheta)\chi$ and $\alpha_0
  = \vartheta b_0$ imply
  \[
    2\vartheta+(1-\vartheta)(3-\chi)=3-\alpha,\qquad
    \vartheta(3-b_0)+(1-\vartheta)(4-\chi)=4-\alpha-\alpha_0.
  \]
  This yields the coefficient bound in~\eqref{E:planar-upper} with $C = C_{\rm
  low}^\vartheta C_{\rm st}^{1-\vartheta}$. The constant coefficient is $A_0 =
  1$. For $q = 3-\alpha > 1$ and $z \ge 0$,
  \[
    \sum_{k\ge0}\frac{z^k}{(k!)^q}
    =\sum_{k\ge0}\left(\frac{z^{k/q}}{k!}\right)^q
    \le\left(\sum_{k\ge0}\frac{z^{k/q}}{k!}\right)^q
    =e^{qz^{1/q}}.
  \]
  Proposition~\ref{P:mixed} supplies the nonnegative full series for the target
  solution. Take $z = Ct^{4-\alpha-\alpha_0}$ and use $(4-\alpha-\alpha_0)/q =
  \kappa$. This proves Proposition~\ref{P:planar-upper}.
\end{proof}

\begin{corollary}[Full planar parameter range]\label{C:planar-range}
  In dimension two, under Hypothesis~\ref{H:homog} and the standing covariance
  positivity assumptions, the upper bound $\E u^2(t,x) \le \exp(C_{\rm dir}t^\kappa)$
  holds throughout $0 < \alpha < 2-\alpha_0$, with a finite constant $C_{\rm dir}$
  depending on $\alpha_0$ and the spatial covariance. A nonnegative generalized
  covariance is interpreted as the homogeneous measure and its Gaussian
  approximations described above.
\end{corollary}

\begin{proof}
  For $0 < \alpha \le 1$, apply Proposition~\ref{P:low-alpha-upper}. For $1 <
  \alpha < 2-\alpha_0$, apply Proposition~\ref{P:planar-upper}. If the
  covariance vanishes, $u = 1$ and $C_{\rm dir} = 0$ suffices. This proves
  Corollary~\ref{C:planar-range}.
\end{proof}

\begin{remark}[Planar calibrations and remaining goals]\label{R:planar-scope}
  The exponent tends to $(4-\alpha)/(3-\alpha)$ as $\alpha_0\downarrow0$,
  consistently with the static source, without an interchange of limits or
  uniformity claim for the constants. For the nonzero higher-exponent case,
  Cauchy--Schwarz gives $\E u(t,x) \le \exp(qC^{1/q}t^\kappa/2)$. If the
  first-moment lower rate is at least $L_*$, this implies $2L_* \le qC^{1/q}$.
  No sharp constant is identified. The existing comparisons with Skorohod
  exponents, in their stated Riesz subclasses, concern different solution
  concepts and are not inputs here. The planar endpoint construction is
  dimension-specific; no general higher-exponent result in dimension three
  follows from it. Section~\ref{S:interval-upper} treats the full
  homogeneous-covariance class by a separate argument. The intrinsic constant is
  identified separately in Theorem~\ref{T:intrinsic-second-moment}.
\end{remark}

\section{Three-dimensional bounds from covariance endpoints}
\label{S:three-endpoints}

The endpoint constructions and conditional moment transfer are retained as an
alternative route. No general three-dimensional factorization is asserted. The
full upper and intrinsic formula do not use these factorizations.

We first construct the product endpoints. Zero coordinate exponents below mean
constant factors, not one-dimensional white noise.

\LeanMarker{L:product-endpoints}
\begin{lemma}[Strict product endpoints]\label{L:product-endpoints}
  Let $0 \le a_j < 1$ for $j = 1,2,3$, let $\alpha = a_1+a_2+a_3 \in (1,2)$, and
  choose $0 < \vartheta < 2-\alpha$. Put $\chi =
  (\alpha-\vartheta)/(1-\vartheta)$ and define
  \begin{equation}\label{E:product-endpoint-vectors}
    \begin{split}
      L_j&=\max\{0,(a_j-1+\vartheta)/\vartheta\},\qquad
      U_j=\min\{1,a_j/\vartheta\},\\
      \Lambda&=\frac{1-\sum_jL_j}{\sum_jU_j-\sum_jL_j},\qquad
      b_j=L_j+\Lambda(U_j-L_j),\qquad
      v_j=\frac{a_j-\vartheta b_j}{1-\vartheta}.
    \end{split}
  \end{equation}
  Then $0 < \Lambda < 1$, $\sum_j b_j = 1$, $\sum_jv_j = \chi$, and
  \[
    a_j=\vartheta b_j+(1-\vartheta)v_j,\qquad
    a_j>0\Longrightarrow 0<b_j,v_j<1,\qquad
    a_j=0\Longrightarrow b_j=v_j=0.
  \]
  For any $c > 0$, the ordinary covariance densities
  \[
    \gamma_a(x)=c\prod_{j=1}^3|x_j|^{-a_j},\qquad
    \gamma_b(x)=c\prod_{j=1}^3|x_j|^{-b_j},\qquad
    \gamma_v(x)=c\prod_{j=1}^3|x_j|^{-v_j}
  \]
  are nonnegative, even, positive definite, locally integrable and tempered.
  Their covariance homogeneity exponents are $\alpha,1,\chi$, respectively, with
  $\alpha < \chi < 2$. For every $\epsilon > 0$ and $z \in \R^3$,
  \begin{equation}\label{E:product-gaussian-comparison}
    (p_\epsilon*\gamma_a)(z)
    \le(p_\epsilon*\gamma_b)(z)^\vartheta
    (p_\epsilon*\gamma_v)(z)^{1-\vartheta}.
  \end{equation}
  A factor with exponent zero is defined to be one everywhere.
\end{lemma}

\begin{proof}
  We first show that the coordinate endpoint intervals admit a strict interior
  choice with the required sums. We then verify the covariance properties of
  the resulting product densities and apply Gaussian Hölder.

  \medskip
  \noindent\textit{Step 1.~}We choose strict coordinate endpoints and verify
  their sums.
  We have $\vartheta < 1$. There are at least two positive coordinates, since
  every $a_j < 1$ and $\sum_j a_j > 1$. At a zero coordinate, $L_j = U_j = 0$.
  At a positive coordinate, $L_j < U_j$: if $a_j \ge \vartheta$, then $U_j = 1$
  and $L_j < 1$; if $a_j < \vartheta$, then $U_j = a_j/\vartheta > 0$, and
  either $L_j = 0$ or $U_j-L_j = (1-\vartheta)/\vartheta > 0$.

  Let $I = \{j:L_j>0\}$ and $r = |I|$. For $r = 0$ the sum of the lower
  endpoints is zero; for $r = 1$ it is less than one. For $r \ge 2$,
  \[
    \sum_{j\in I}a_j\le\alpha<2-\vartheta
    \le r-(r-1)\vartheta,
  \]
  because the last inequality is equivalent to $(r-2)(1-\vartheta) \ge 0$.
  Therefore,
  \[
    \sum_jL_j
    =\frac{\sum_{j\in I}a_j-r+r\vartheta}{\vartheta}<1.
  \]
  If some $a_j \ge \vartheta$, its upper endpoint equals one and a different
  positive coordinate has positive upper endpoint. Otherwise, $\sum_jU_j =
  \alpha/\vartheta > 1$. Hence, $\sum_jL_j < 1 < \sum_jU_j$, which proves $0 <
  \Lambda < 1$ and $\sum_jb_j = 1$. At every positive coordinate, $L_j < b_j <
  U_j$, so $0 < b_j < 1$. The inequalities $b_j < U_j \le a_j/\vartheta$ and
  $b_j > L_j \ge (a_j-1+\vartheta)/\vartheta$ give $0 < v_j < 1$. At a zero
  coordinate $b_j = v_j = 0$. The reconstruction identity is immediate, and
  summing it gives $\sum_jv_j = \chi$. Finally,
  \[
    \chi-\alpha=\frac{\vartheta(\alpha-1)}{1-\vartheta}>0,
    \qquad 2-\chi=\frac{2-\alpha-\vartheta}{1-\vartheta}>0.
  \]

  \medskip
  \noindent\textit{Step 2.~}We verify local integrability and positive
  definiteness.
  For any of the three exponent vectors $e$, local integrability is a product of
  one-dimensional integrals with exponents below one. For a Schwartz function
  $\psi$, choose $N > 1$ and use $|\psi(x)| \le C_N\prod_j(1+|x_j|)^{-N}$. This
  proves absolute integrability against $\gamma_e$ and also temperedness. For
  each active coordinate, $0 < e_j < 1$, use
  \[
    |u|^{-e_j}=\frac1{\Gamma(e_j/2)}
    \int_0^\infty s^{e_j/2-1}e^{-su^2}\ud s.
  \]
  The product over the active coordinates is a positive mixture of Gaussian
  functions in those coordinates, constant in all other coordinates. Each such
  Gaussian is positive definite. For a complex Schwartz $\varphi$, set $\psi =
  \varphi*\widetilde{\varphi}$, where $\widetilde{\varphi}(x) =
  \overline{\varphi(-x)}$. The absolute integrability just proved permits
  Tonelli with $|\psi|$ and then Fubini for the mixture against $\psi$. Its
  quadratic form is nonnegative. No Euler Gamma factor at exponent zero is
  introduced. Evenness and homogeneity follow from the displayed densities.

  \medskip
  \noindent\textit{Step 3.~}We compare the smoothed densities at common
  amplitude.
  The reconstruction of each exponent and the common amplitude imply $\gamma_a =
  \gamma_b^\vartheta\gamma_v^{1-\vartheta}$ almost everywhere; values on
  singular coordinate planes do not change these locally integrable covariances.
  H\"older's inequality on the Gaussian probability measure
  gives~\eqref{E:product-gaussian-comparison}. Its endpoint convolutions are
  finite by the product integrability and Gaussian decay. The amplitude remains
  $c^\vartheta c^{1-\vartheta} = c$. This proves
  Lemma~\ref{L:product-endpoints}.
\end{proof}

We now take $d = 3$, $0 < \alpha_0 < 1$, and $1 < \alpha < 2-\alpha_0$. Let
$\Gamma$ be a nonnegative, even, positive-definite tempered spatial covariance,
with homogeneity exponent $\alpha$. As a positive distribution it is a Radon
measure satisfying $\Gamma(rA) = r^{3-\alpha}\Gamma(A)$; a covariance function
is identified with its density measure. Gaussian convolution is used for all
covariance products. Choose and fix
\begin{equation}\label{E:three-parameters}
  \alpha_0<\vartheta<2-\alpha,\qquad
  b_0=\alpha_0/\vartheta,\qquad
  \chi=(\alpha-\vartheta)/(1-\vartheta).
\end{equation}
Then $0 < b_0 < 1$ and $\alpha < \chi < 2$, by the identities in the preceding
proof. The zero target covariance is trivial, since then $u = 1$.

\begin{hypothesis}[Gaussian covariance comparison]\label{H:three-endpoints}
  There are nonnegative, even, positive-definite tempered endpoint covariances
  $\Gamma_{\rm low}$ and $\Gamma_{\rm st}$ on $\R^3$ with measure scaling
  \[
    \Gamma_{\rm low}(rA)=r^2\Gamma_{\rm low}(A),\qquad
    \Gamma_{\rm st}(rA)=r^{3-\chi}\Gamma_{\rm st}(A),
  \]
  such that, for every $\epsilon > 0$ and $z \in \R^3$,
  \begin{equation}\label{E:three-gaussian-comparison}
    (p_\epsilon*\Gamma)(z)
    \le(p_\epsilon*\Gamma_{\rm low})(z)^\vartheta
    (p_\epsilon*\Gamma_{\rm st})(z)^{1-\vartheta}.
  \end{equation}
\end{hypothesis}

This is a condition on deterministic covariance measures and their Gaussian
convolutions. It assumes no moment estimate. We do not assert that every
homogeneous covariance in dimension three admits these endpoints.

\LeanMarker{L:plane-endpoints}
\begin{lemma}[Plane-supported endpoints]\label{L:plane-endpoints}
  Retain~\eqref{E:three-parameters}, let $c > 0$, and let $a_1,a_2 \ge 0$
  satisfy $s = a_1+a_2 = \alpha-1$. In coordinates $(y,z) \in \R^2\times\R$, put
  \[
    m(\ud y,\ud z)=\ud y\,\delta_0(\ud z),\qquad
    v_j=\frac{a_j}{1-\vartheta},\qquad
    \Gamma_{\rm pl}=c\prod_{j=1}^2|y_j|^{-a_j}m.
  \]
  Zero coordinate powers mean one. Then
  \[
    \Gamma_{\rm low}=cm,\qquad
    \Gamma_{\rm st}=c\prod_{j=1}^2|y_j|^{-v_j}m
  \]
  satisfy Hypothesis~\ref{H:three-endpoints} for $\Gamma_{\rm pl}$. In the
  spectral convention of Lemma~\ref{L:alpha-one-convolution}, this low endpoint
  has $V_{\rm low} = c/\pi$.
\end{lemma}

\begin{proof}
  We verify the exponents and covariance properties of the plane-supported
  measures. We then apply Hölder to their ordinary densities on the same plane
  measure and compute the low endpoint spectral mass.

  \medskip
  \noindent\textit{Step 1.~}We verify the endpoint parameters and measure
  scaling.
  Since $0 < s < 1-\alpha_0$ and $\vartheta < 2-\alpha = 1-s$,
  \[
    0\le a_j<1,\qquad 0\le v_j<1,\qquad
    v_1+v_2=\frac{s}{1-\vartheta}=\chi-1<1.
  \]
  A measure $f_e(y)\ud y\,\delta_0(\ud z)$, with $f_e(y) = \prod_j|y_j|^{-e_j}$
  and $0 \le e_j < 1$, is locally finite. For a Schwartz test function,
  product-weighted rapid decay proves absolute integrability and temperedness.
  Its scaling exponent as a measure is $2-e_1-e_2$. Thus, the three measures in
  the statement have exponents $3-\alpha$, $2$, and $3-\chi$. They are even and
  nonnegative.

  \medskip
  \noindent\textit{Step 2.~}We include the normal Dirac measure in the
  positivity argument.
  We verify positive definiteness, including the normal Dirac measure. Each
  positive coordinate power has the Gaussian mixture
  \[
    |u|^{-e}=\frac1{\Gamma(e/2)}
    \int_0^\infty r^{e/2-1}e^{-ru^2}\ud r,\qquad 0<e<1;
  \]
  a zero coordinate contributes the positive-definite constant one.
  Consequently, $f_e$ is a positive mixture of positive-definite tangential
  Gaussian functions. For $\varphi \in \mathcal{S}(\R^3)$, its quadratic form
  against the plane measure is
  \[
    \int_\R\ud z\int_{\R^2}\int_{\R^2}
    f_e(y-y')\varphi(y,z)\overline{\varphi(y',z)}
    \ud y\ud y'.
  \]
  Absolute Fubini follows by setting
  \[
    A(u)=\int_\R\int_{\R^2}
    |\varphi(y+u,z)\varphi(y,z)|\ud y\ud z.
  \]
  Rapid decay gives $A(u) \le C_N\prod_j(1+|u_j|)^{-N}$. Hence, $\int
  f_e(u)A(u)\ud u < \infty$ for $N > 1$, by one-dimensional integrability.
  Integrating the nonnegative Gaussian quadratic forms proves positive
  definiteness. This argument includes the low endpoint, where both coordinate
  powers are zero.

  \medskip
  \noindent\textit{Step 3.~}We compare the Gaussian convolutions on the plane.
  For $\epsilon > 0$ and $w \in \R^3$, apply H\"older on the finite positive
  measure $p_\epsilon(w-(y,z))m(\ud y,\ud z)$. Its mass is
  $p_\epsilon^{(1)}(w_3)$. If $f(y) = \prod_j|y_j|^{-a_j}$, the identity between
  ordinary densities
  \[
    cf=c^\vartheta
    \left(cf^{1/(1-\vartheta)}\right)^{1-\vartheta}
  \]
  holds $m$-almost everywhere. Both endpoint convolutions are finite, since $v_j
  < 1$ and the Gaussian decays at infinity. H\"older gives
  exactly~\eqref{E:three-gaussian-comparison}; no power of a distribution is
  taken.

  \medskip
  \noindent\textit{Step 4.~}We compute the spectral normalization of the low
  endpoint.
  Finally, the low spectral measure is
  \[
    \mu_{\rm low}(\ud\xi)
    =\frac{c}{2\pi}\delta_0(\ud\xi_1)\delta_0(\ud\xi_2)\ud\xi_3.
  \]
  This follows from the unit spectral atom of a constant and the one-dimensional
  Fourier identity for $\delta_0$ in the stated convention. Its unit-ball mass
  is $c/\pi$, so the uniform closing bound of
  Lemma~\ref{L:alpha-one-convolution} is $c/2$. This proves
  Lemma~\ref{L:plane-endpoints}.
\end{proof}

\LeanMarker{L:endpoint-closure}
\begin{lemma}[Linear images and finite positive sums]\label{L:endpoint-closure}
  Fix the parameters in~\eqref{E:three-parameters}. The covariance comparison of
  Hypothesis~\ref{H:three-endpoints} is preserved by simultaneously pushing the
  target and its two endpoints forward under any invertible linear map $T:\R^3
  \to \R^3$. It is also preserved by finite sums with nonnegative weights, when
  all component comparisons use the same $\alpha,\vartheta,\chi$.
\end{lemma}

\begin{proof}
  We first extend the comparison from isotropic to general nondegenerate
  Gaussian kernels. We then track the common determinant factor under linear
  maps and use discrete Hölder for finite positive sums.

  \medskip
  \noindent\textit{Step 1.~}We extend the Gaussian convolution comparison.
  First the comparison extends to every centered nondegenerate Gaussian
  probability density $q_\Sigma$ with covariance matrix $\Sigma$. Choose $a > 0$
  below the smallest eigenvalue of $\Sigma$. The Gaussian convolution identity
  gives $q_\Sigma = p_a*q_{\Sigma-aI}$. Put $F_J = p_a*J$ for each $J \in
  \{\Gamma,\Gamma_{\rm low},\Gamma_{\rm st}\}$. These functions are nonnegative,
  bounded, and continuous: Gaussian weighting makes their positive tempered
  spectral measures finite, and the covariance bound is $|F_J(x)| \le F_J(0) <
  \infty$. Tonelli and H\"older on the Gaussian probability measure yield
  \[
    q_\Sigma*\Gamma
    \le q_{\Sigma-aI}*
    \left(F_{\Gamma_{\rm low}}^\vartheta
    F_{\Gamma_{\rm st}}^{1-\vartheta}\right)
    \le(q_\Sigma*\Gamma_{\rm low})^\vartheta
    (q_\Sigma*\Gamma_{\rm st})^{1-\vartheta}.
  \]
  Every convolution here is finite.

  \medskip
  \noindent\textit{Step 2.~}We push all three measures through the same linear
  map.
  Write $J_T = T_\#J$ and $d_T = |\det T|$. Pushforward preserves nonnegativity,
  evenness, local finiteness, temperedness, and homogeneous scaling, since $T$
  is invertible and commutes with scalar dilation. For positive definiteness,
  take a complex Schwartz function $\varphi$ and put $\psi(y) = \varphi(Ty)$.
  Then
  \[
    (\varphi*\widetilde{\varphi})(Ty)
    =d_T(\psi*\widetilde{\psi})(y).
  \]
  The quadratic form of $T_\#J$ is therefore $d_T$ times a nonnegative quadratic
  form of $J$. With $\Sigma = \epsilon T^{-1}(T^{-1})^\top$, direct Gaussian
  change of variables gives, for every $\epsilon > 0$ and $x$,
  \[
    (p_\epsilon*J_T)(x)
    =d_T^{-1}(q_\Sigma*J)(T^{-1}x).
  \]
  Apply the preceding comparison. The two endpoint factors contribute
  $d_T^{-\vartheta}d_T^{-(1-\vartheta)} = d_T^{-1}$, which is the target factor.
  This proves the comparison for the three pushed-forward measures.

  \medskip
  \noindent\textit{Step 3.~}We take finite weighted sums of endpoint triples.
  For the finite-sum assertion, let $\Gamma_i,L_i,S_i$, $1 \le i \le N$, be
  component triples, and let $w_i \ge 0$. At a fixed width and translation put
  $l_i = (p_\epsilon*L_i)(x)$ and $s_i = (p_\epsilon*S_i)(x)$. Finite sums
  preserve all required covariance properties and scalings. Discrete H\"older
  gives
  \[
    \begin{split}
      p_\epsilon*\left(\sum_iw_i\Gamma_i\right)
      &\le\sum_i(w_il_i)^\vartheta(w_is_i)^{1-\vartheta}\\
      &\le\left(p_\epsilon*\sum_iw_iL_i\right)^\vartheta
      \left(p_\epsilon*\sum_iw_iS_i\right)^{1-\vartheta}.
    \end{split}
  \]
  Zero weights cause no ambiguity, since both powers are positive. The all-zero
  target also satisfies the comparison. This proves
  Lemma~\ref{L:endpoint-closure}.
\end{proof}

The linear maps in Lemma~\ref{L:endpoint-closure} act on the noise covariances.
All moment estimates below still use the original isotropic wave kernel.
Endpoint spectral masses and variational quantities are computed for the
transformed covariances; their values need not be unchanged by a nonorthogonal
map.

For the low endpoint let $V_{\rm low} = \mu_{\rm low}(B_1)$ in the spectral
convention of Lemma~\ref{L:alpha-one-convolution}, and put
\[
  Z_{b_0}=\frac2{(1-b_0)(2-b_0)},\qquad
  C_{\rm low}=2Z_{b_0}\pi V_{\rm low}.
\]
For the static endpoint, write
\[
  \mathcal{Q}_{\rm st}(g)
  =\int_{\R^3}\int_{\R^3}g(y+z)^2g(y)^2\ud y
  \,\Gamma_{\rm st}(\ud z),\qquad
  \mathcal{F}_3=\{g\in W^{1,2}(\R^3):\|g\|_2=1\},
\]
using real Sobolev functions, and define
\begin{equation}\label{E:three-static-variational}
  \begin{split}
    \mathcal{M}_{\rm st}
    &=\sup_{g\in\mathcal{F}_3}
    \left\{\mathcal{Q}_{\rm st}(g)^{1/2}-\|\nabla g\|_2^2\right\},\\
    Q&=3-\chi,\qquad \nu_{\rm st}=4-\chi,\qquad
    \Lambda_{\rm st}
    =2\left(\frac{2\sqrt{\mathcal{M}_{\rm st}}}{4-\chi}\right)^{4-\chi}.
  \end{split}
\end{equation}
The functional is precisely Chen--Hu's (1.11), p.~2147, for the endpoint
covariance, with unit gradient penalty. The positive-measure energy above agrees
with the spectral form on every real Sobolev profile by
Lemma~\ref{L:physical-spectral-energy}; thus no density is assigned to a
singular covariance. The lemma also supplies finiteness, since $\chi < 2$ gives
Dalang's condition as checked below. Nonnegativity follows by dilating a
unit-$L^2$ Schwartz function: its energy term is nonnegative and its gradient
penalty tends to zero.

\begin{proposition}[Conditional three-dimensional transfer]
  \label{P:three-transfer}
  Under the preceding assumptions and Hypothesis~\ref{H:three-endpoints}, with
  $\Gamma \ne 0$, there is a finite constant $C > 0$ such that
  \begin{equation}\label{E:three-transfer}
    A_k(t)\le\frac{C^kt^{(4-\alpha-\alpha_0)k}}{(k!)^{3-\alpha}}
    \quad(k\ge1),\qquad
    \E u^2(t,x)\le\exp\left(qC^{1/q}t^\kappa\right),
    \quad q=3-\alpha.
  \end{equation}
  The estimates hold for all $t > 0$ and $x \in \R^3$, with data $(1,0)$. The
  constant may depend on the fixed covariances, $\alpha_0$, and $\vartheta$, but
  not on $k,t,x$. In addition, the explicit variational upper constant satisfies
  \begin{equation}\label{E:three-variational-upper}
    \limsup_{t\to\infty}t^{-\kappa}\log\E u^2(t,x)
    \le q\left(C_{\rm low}^{\,\vartheta}
    \Lambda_{\rm st}^{\,1-\vartheta}\right)^{1/q}.
  \end{equation}
  This constant is not asserted sharp and is not substituted into the pure
  all-time exponential bound in~\eqref{E:three-transfer}.
\end{proposition}

\begin{proof}
  We verify the static input and establish a regularized coefficient
  comparison on one common measure. We identify the static and target Green
  limits separately before removing noise widths. We then sum the all-time
  coefficient bound and retain the finite low orders when deriving the sharper
  asymptotic constant.

  \medskip
  \noindent\textit{Step 1.~}We state and verify the imported static
  coefficient input.
  We first record the general static input needed here. For a static endpoint
  covariance $\Xi$ of exponent $\chi$, Fourier uniqueness gives $\mu_\Xi(rA) =
  r^\chi\mu_\Xi(A)$ for its positive spectral measure. Its unit-ball mass is
  finite, and
  \[
    \int_{\R^3}(1+|\xi|^2)^{-1}\mu_\Xi(\ud\xi)
    \le\mu_\Xi(B_1)
    +2^\chi\mu_\Xi(B_1)\sum_{j\ge0}2^{j(\chi-2)}<\infty.
  \]
  Thus, the static Dalang condition holds. Let $U_n$ be the source's ordered
  integrals for this covariance. Its nonnegative positive-definite
  tempered-distribution class is explicitly included by the generalized
  covariance convention in~\cite{chen.hu:25:hyperbolic}*{p.~2152, (2.18)}; its
  homogeneity degree satisfies $0 < \chi < 2 < 3$. No density or angular
  regularity is added here. Let $U_0 = 1$, and set $B_k(t) =
  \sum_{m=0}^{2k}\E[U_m(t,x)U_{2k-m}(t,x)]$. If $\Xi = 0$, then $B_k = 0$ for $k
  \ge 1$, and $C_{\rm st} = 1$ suffices in the bound below. For the remaining
  static calculation, assume $\Xi \ne 0$. From
  Chen--Hu~\cite{chen.hu:25:hyperbolic}*{Theorem~3.8(i), p.~2165, and (2.17),
  (2.19)--(2.22), pp.~2151--2153}, the expectation and product rules give finite
  nonnegative coefficients with every split and pairing included; the latter
  formulas also specify the two-chain and measure conventions. Finiteness of the
  variational input follows from Chen~\cite{chen:23:exponential}*{Theorem~1.1,
  p.~2}, which gives a finite square-root variational supremum under Dalang's
  condition. Chen--Hu identify that supremum with $2^{\chi/(4-\chi)}\mathcal{M}$
  on p.~2175. Hence, their coefficient-rate constant involving $\mathcal{M}$ is
  finite. For $\Xi = \Gamma_{\rm st}$ this functional is $\mathcal{M}_{\rm st}$
  defined above. Only finiteness, not a sharp value, is needed for the all-time
  coefficient bound.

  From Chen--Hu~\cite{chen.hu:25:hyperbolic}*{(5.1), p.~2173, and (5.15),
  p.~2180}, at $p = 2$,
  \[
    B_k(t)=t^{(4-\chi)k}B_k(1),\qquad
    \limsup_{k\to\infty}\frac1k
    \log\left((k!)^{3-\chi}B_k(1)\right)<\infty.
  \]
  For $\Xi \ne 0$, choose $R > 0$ above this upper limit, then enlarge $e^R$ to
  cover the $k$th roots of the finitely many remaining quantities
  $(k!)^{3-\chi}B_k(1)$. We obtain a finite $C_{\rm st} = C_{\rm st}(\Xi) > 0$
  such that
  \begin{equation}\label{E:three-static-coefficient}
    B_k(t)\le\frac{C_{\rm st}^{\,k}t^{(4-\chi)k}}{(k!)^{3-\chi}},
    \qquad k\ge1,\quad t>0.
  \end{equation}
  Also $B_0 = 1$. We apply this source bound only after taking the fixed-order
  noise limits, with $\Xi = \Gamma_{\rm st}$.

  \medskip
  \noindent\textit{Step 2.~}We verify the target hypotheses and form the
  common regularized measure.
  For the target, the same spectral scaling calculation has exponent $\alpha$.
  Its dyadic sum for the weight $(1+|\xi|^2)^{-(2-\alpha_0)/2}$ is finite
  because $\alpha < 2-\alpha_0$. The target fixed-order limit is supplied by
  Lemma~\ref{L:joint-fixed-order}. The endpoint $\Gamma_{\rm low}$ has
  covariance exponent one, with $V_{\rm low},C_{\rm low}$ defined above.
  Finiteness of $V_{\rm low}$ follows from the Radon spectral measure. It is
  positive: otherwise homogeneity gives $\mu_{\rm low} = 0$,
  contradicting~\eqref{E:three-gaussian-comparison} and $\Gamma \ne 0$.

  Fix $k \ge 1$, $t > 0$, and positive widths $\epsilon,\delta,\eta$. Noise
  width $\epsilon$ means edge covariance width $2\epsilon$; $\delta$ is the
  target temporal cutoff and $\eta$ is Green smoothing. Write
  $X_n^{\epsilon,\delta,\eta}(t,x)$ for the source's ordered integral against
  its smoothed field $W_{\epsilon,\delta}$, with every Green kernel replaced by
  $G_{h,\eta}$, and set $X_0^{\epsilon,\delta,\eta} = 1$. Explicitly, the target
  coefficient is
  \[
    A_k^{\epsilon,\delta,\eta}(t)
    =\sum_{m=0}^{2k}
    \E\left[X_m^{\epsilon,\delta,\eta}(t,x)
    X_{2k-m}^{\epsilon,\delta,\eta}(t,x)\right].
  \]
  On the full disjoint union over all splits $m+n = 2k$ and all $D \in \Pi_k$,
  define a common finite positive measure by
  \[
    \int\Phi\ud\nu_{k,t,\eta}
    :=\sum_{m=0}^{2k}\sum_{D\in\Pi_k}
    \int_{[0,t]^m_<\times[0,t]^{2k-m}_<}\ud s
    \int\Phi(m,D,s,y)\prod_{l=1}^{2k}G_{h_l,\eta}(\ud y_l).
  \]
  Each chain has its own root; $h_l = s_l-s_{p(l)}$, $y_l$ is an increment, and
  $x_l$ is the cumulative position in its chain. Empty chains use empty
  products. The spatial mass is $\prod_lh_l$, so bounded time simplexes and
  finite discrete sums prove finiteness. This measure does not depend on which
  covariance is inserted.

  For $J \in \{\Gamma,\Gamma_{\rm low},\Gamma_{\rm st}\}$ set
  \[
    H_J=\prod_{(i,j)\in D}(p_{2\epsilon}*J)(x_i-x_j),\qquad
    T_b=\prod_{(i,j)\in D}|s_i-s_j|^{-b}.
  \]
  If $Q_\delta$ is the product of the target temporal covariances, the source
  cutoff gives $0 \le Q_\delta \le T_{\alpha_0} = T_{b_0}^{\,\vartheta}$. The
  covariance hypothesis gives $H_\Gamma \le H_{\Gamma_{\rm low}}^\vartheta
  H_{\Gamma_{\rm st}}^{1-\vartheta}$ at the same edge width. Define the
  deterministic low endpoint and the regularized static coefficient by
  \[
    L_k^{\epsilon,\eta}(t)
    =\int T_{b_0}H_{\Gamma_{\rm low}}\ud\nu_{k,t,\eta},
    \qquad
    B_k^{\epsilon,\eta}(t)
    =\int H_{\Gamma_{\rm st}}\ud\nu_{k,t,\eta}.
  \]
  Gaussian moments and separate chain reversal identify the target regularized
  coefficient with $\int Q_\delta H_\Gamma \ud\nu_{k,t,\eta}$. The static
  identity has temporal covariance one and includes self-pairings. H\"older on
  the entire common measure yields
  \begin{equation}\label{E:three-regularized-holder}
    A_k^{\epsilon,\delta,\eta}(t)
    \le\left(L_k^{\epsilon,\eta}(t)\right)^\vartheta
    \left(B_k^{\epsilon,\eta}(t)\right)^{1-\vartheta}.
  \end{equation}
  No stochastic identity is assumed for $L_k^{\epsilon,\eta}$ with its original
  temporal power, and no relation between different temporal cutoff kernels is
  used.

  \medskip
  \noindent\textit{Step 3.~}We estimate the deterministic low endpoint.
  Lemma~\ref{L:alpha-one-convolution} applies to the low endpoint with no
  angular-density assumption. Its uniform closing bound is $\pi V_{\rm low}/2$.
  Descending increment integration removes the only remaining covariance at each
  closing label; at each opening the paired covariance is already removed and
  the Green mass is $h_l$. The $k$ opening increments sum to at most $2t$.
  Consequently, the common spatial bound is $(\pi V_{\rm low})^kt^k/k^k$,
  uniformly in $\epsilon,\eta$. After summing pairings, the temporal integrand
  is symmetric in all $2k$ labels. Its full cube integral is
  $P_k(Z_{b_0}t^{2-b_0})^k$, where $P_k = (2k)!/(2^kk!)$. The block-order
  factors and the binomial split sum give
  \[
    P_k\sum_{m=0}^{2k}\frac1{m!(2k-m)!}=\frac{2^k}{k!}.
  \]
  Therefore,
  \begin{equation}\label{E:three-low-coefficient}
    L_k^{\epsilon,\eta}(t)
    \le\frac{C_{\rm low}^{\,k}t^{(3-b_0)k}}{k!\,k^k}
    \le\frac{C_{\rm low}^{\,k}t^{(3-b_0)k}}{(k!)^2}.
  \end{equation}
  Here $b_0 < 1$ makes the original temporal cube integral finite. At fixed
  spatial width, the static covariance is bounded, so its regularized
  coefficient is finite too. This justifies the finite right side
  in~\eqref{E:three-regularized-holder}.

  \medskip
  \noindent\textit{Step 4.~}We identify the static source and its
  zero-Green-width integral.
  We now identify the static inner Green limit with the source's prescribed
  zero-Green-width approximant. Take the source's static isonormal noise $W^{\rm
  st}$ with covariance $\Gamma_{\rm st}$ and set $W_\epsilon(x) = W^{\rm
  st}(p_\epsilon(x-\cdot))$. All widths are therefore on the source's Gaussian
  space, and the covariance is $p_{2\epsilon}*\Gamma_{\rm st}$. This covariance
  is bounded and continuous, since Gaussian weighting makes its positive
  tempered spectral measure finite. It has a continuous modification. Indeed,
  with $\mu_{\rm st}$ the spectral measure,
  \[
    \E|W_\epsilon(x+h)-W_\epsilon(x)|^2
    \le |h|^2\int_{\R^3}|\xi|^2e^{-\epsilon|\xi|^2}
    \mu_{\rm st}(\ud\xi).
  \]
  The integral is finite by temperedness. Gaussian eighth moments are therefore
  bounded by a constant times $|h|^8$. Kolmogorov's continuity criterion in
  dimension three, applied on compact sets and then on a countable exhaustion,
  gives the claimed modification. We choose this version. For $n \ge 1$ define
  $\Phi_n(\mathbf{x}) = \prod_{j=1}^n W_\epsilon(x_j)$ as an $L^2$ random
  variable. Gaussian moments show that $\Phi_n$ is continuous as a map from
  $(\R^3)^n$ to $L^2$ and uniformly bounded there: if $v_\epsilon =
  (p_{2\epsilon}*\Gamma_{\rm st})(0)$, then
  \[
    \|\Phi_n(\mathbf{x})\|_2^2\le P_n v_\epsilon^n.
  \]
  The continuity follows by applying the same finite Gaussian moment formula to
  $\E|\Phi_n(\mathbf{x})-\Phi_n(\mathbf{x}')|^2$. Its separable continuous image
  is strongly measurable.

  Let $\mathfrak{g}_{n,t,x}^{\eta}$ be the finite positive ordered Green-chain
  measure, including its time variables, from the source's kernel (2.13), with
  each Green factor smoothed by $p_\eta$. For $\eta = 0$, its spatial marginal
  is exactly the measure $\mu_n^{t,x}$ specified immediately before (2.22),
  pp.~2152--2153. Its mass is finite because each spatial factor has mass its
  time gap. Thus,
  \[
    U_n^{\epsilon,0}(t,x)
    :=\int\Phi_n(\mathbf{x})
    \mathfrak{g}_{n,t,x}^{0}(\ud\mathbf{s},\ud\mathbf{x})
  \]
  exists as a Bochner integral in $L^2$. We verify its agreement with the
  source's pathwise integral, rather than identifying it only by notation. Put
  $\rho = \mathfrak{g}_{n,t,x}^{0}$ and $M_n = P_n^{1/2}v_\epsilon^{n/2}$, and
  denote the pathwise integral by $J_n = \int\Phi_n\ud\rho$. The continuous
  field version makes the integrand jointly measurable. Since
  \[
    \int\E|\Phi_n|\ud\rho\le M_n\rho(1)<\infty,
    \qquad
    \|J_n\|_2
    \le\int\|\Phi_n\|_2\ud\rho\le M_n\rho(1),
  \]
  the pathwise integral exists almost surely and belongs to $L^2$. For every $Y
  \in L^2$, Cauchy--Schwarz and Fubini give
  \[
    \E[YJ_n]
    =\int\E[Y\Phi_n]\ud\rho.
  \]
  Bochner duality gives the same identity for the $L^2$ integral. Taking $Y =
  J_n-U_n^{\epsilon,0}$ identifies the two random variables almost surely.
  Changing to another jointly measurable field modification leaves the integral
  unchanged by Fubini. Since the static integrand is independent of the time
  variables, this is the integral against $\mu_n^{t,x}$ in the source's
  convention and hence its smoothed approximant $S_{n,\epsilon}$ with common
  width, from Chen--Hu~\cite{chen.hu:25:hyperbolic}*{(3.1), p.~2153, and the
  convention preceding (2.22)}. Define $U_n^{\epsilon,\eta}$ by the same formula
  with $\mathfrak{g}_{n,t,x}^{\eta}$, and retain $U_0^{\epsilon,\eta} = 1$.

  \medskip
  \noindent\textit{Step 5.~}We identify both inner Green limits in mean
  square.
  In the second and mixed moments of these integrals, the Gaussian moment
  factors are bounded continuous functions of all spatial positions. For each
  pair of time vectors, weak convergence of the products of Green measures,
  including mixed products with one width already zero, applies to these
  functions. The masses are unchanged. A constant depending only on $n,\epsilon$
  times the product of all time gaps dominates the resulting time integrals.
  This is integrable on the two bounded time simplexes. Hence, as
  $\eta\downarrow0$,
  \[
    \E\left[(U_n^{\epsilon,\eta})^2\right]
    \longrightarrow\E\left[(U_n^{\epsilon,0})^2\right],
    \qquad
    \E\left[U_n^{\epsilon,\eta}U_n^{\epsilon,0}\right]
    \longrightarrow\E\left[(U_n^{\epsilon,0})^2\right].
  \]
  Subtracting twice the second limit from the two square terms proves
  $U_n^{\epsilon,\eta} \to U_n^{\epsilon,0}$ in $L^2$. This identifies the limit
  with the source approximant, rather than only showing that positive-width
  approximants are Cauchy. For the target we specify the zero-Green-width object
  as well. At fixed $\epsilon,\delta$, its field covariance is
  $\gamma_\delta^0(s-r)(p_{2\epsilon}*\Gamma)(y-z)$, where
  \[
    \gamma_\delta^0(u)
    =\Gamma(\alpha_0)^{-1}\int_0^{1/\delta}
    e^{-\lambda|u|}\lambda^{\alpha_0-1}\ud\lambda,\qquad
    \gamma_\delta^0(0)
    =\frac{\delta^{-\alpha_0}}{\alpha_0\Gamma(\alpha_0)}.
  \]
  Here $\Gamma(\alpha_0)$ is Euler's gamma function. The variance
  $v_{\epsilon,\delta} = \gamma_\delta^0(0)(p_{2\epsilon}*\Gamma)(0)$ is finite.
  This field has a jointly continuous modification. The temporal covariance
  satisfies $\gamma_\delta^0(0)-\gamma_\delta^0(h) \le C_\delta|h|$, by
  $1-e^{-\lambda|h|} \le \lambda|h|$. The spatial increment bound just proved
  applies to $p_{2\epsilon}*\Gamma$. The joint increment variance is therefore
  bounded by $C_{\epsilon,\delta}(|h_0|+|h'|^2)$. Gaussian twelfth moments are
  bounded by $C_{\epsilon,\delta}|h|^6$ for $|h| \le 1$. Kolmogorov's criterion
  in dimension four supplies the version.

  Let $\mathfrak{f}_{n,t,x}^{\eta}$ be the ordered Green measure of the target's
  kernel (2.16), with $s_{n+1} = t$, $x_{n+1} = x$: its factors are
  $G_{s_{j+1}-s_j,\eta}$. Define
  \[
    \begin{split}
      \Psi_n(\mathbf{s},\mathbf{x})
      &=\prod_{j=1}^nW_{\epsilon,\delta}(s_j,x_j),\\
      X_n^{\epsilon,\delta,0}(t,x)
      &=\int\Psi_n(\mathbf{s},\mathbf{x})
      \mathfrak{f}_{n,t,x}^{0}(\ud\mathbf{s},\ud\mathbf{x}).
    \end{split}
  \]
  The Gaussian moment bound is $\|\Psi_n\|_2^2 \le P_nv_{\epsilon,\delta}^n$.
  Joint measurability, finite measure mass, and the preceding Fubini argument
  show that this Bochner integral equals the pathwise integral in the target's
  prescribed convention \cite{chen:author2026:time-dependency}*{pp.~13--14}. The
  same definition at positive $\eta$ gives $X_n^{\epsilon,\delta,\eta}$ above.
  Its second and mixed $\eta/0$ moments converge by weak Green convergence
  against the bounded continuous Gaussian moment functions, with dominating
  bound a constant times the product of the time gaps. Thus,
  \[
    \|X_n^{\epsilon,\delta,\eta}
    -X_n^{\epsilon,\delta,0}\|_2\longrightarrow0.
  \]
  This proves the inner-limit identification with the target approximant. Both
  order-zero objects remain one.

  \medskip
  \noindent\textit{Step 6.~}We remove noise widths and sum the all-time
  coefficient bound.
  First take this Green limit in~\eqref{E:three-regularized-holder}, using the
  uniform bound~\eqref{E:three-low-coefficient}. Next send the noise widths to
  zero. The static fixed-order $L^2$ limit is from
  Chen--Hu~\cite{chen.hu:25:hyperbolic}*{Theorem~3.8(i), p.~2165}; the target
  fixed-order limit is from Lemma~\ref{L:joint-fixed-order}. Products of two
  fixed orders converge in $L^1$, since
  \[
    \|YZ-Y'Z'\|_1
    \le\|Y-Y'\|_2\|Z\|_2+\|Y'\|_2\|Z-Z'\|_2.
  \]
  The split sums are finite, and zero-order factors remain one. We obtain
  \[
    A_k(t)\le
    \left(\frac{C_{\rm low}^{\,k}t^{(3-b_0)k}}{(k!)^2}\right)^\vartheta
    B_k(t)^{1-\vartheta}.
  \]
  Only now use~\eqref{E:three-static-coefficient}. Neither monotonicity of
  spatial smoothing nor a uniform smoothed static asymptotic was needed. The
  identities
  \[
    \begin{split}
      2\vartheta+(1-\vartheta)(3-\chi)&=3-\alpha,\\
      \vartheta(3-b_0)+(1-\vartheta)(4-\chi)&=4-\alpha-\alpha_0
    \end{split}
  \]
  give the coefficient estimate with $C = C_{\rm low}^{\,\vartheta}C_{\rm
  st}^{\,1-\vartheta}$. Finally, $A_0 = 1$, $q = 3-\alpha > 1$, and
  \[
    \sum_{k\ge0}\frac{z^k}{(k!)^q}
    \le\left(\sum_{k\ge0}\frac{z^{k/q}}{k!}\right)^q
    =e^{qz^{1/q}},\qquad z\ge0.
  \]
  The nonnegative full series of Proposition~\ref{P:mixed}, with $z =
  Ct^{4-\alpha-\alpha_0}$ and $(4-\alpha-\alpha_0)/q = \kappa$, yields the
  moment estimate. It remains to prove~\eqref{E:three-variational-upper}. The
  precise upper coefficient rate in Chen--Hu's (5.15), p.~2180, at $p = 2$ is
  \[
    \begin{split}
      \limsup_{k\to\infty}
      \left((k!)^Q B_k(1)\right)^{1/k}
      &\le 2^{-Q}2^{\nu_{\rm st}}
      \left(\frac{2\sqrt{\mathcal{M}_{\rm st}}}{\nu_{\rm st}}\right)^{\nu_{\rm st}}\\
      &=\Lambda_{\rm st}.
    \end{split}
  \]

  \medskip
  \noindent\textit{Step 7.~}We retain the precise asymptotic coefficient rate
  and its low-order remainder.
  The factor two follows from $\nu_{\rm st}-Q = 1$; the logarithm in the source
  applies to the entire coefficient-rate product. Fix $r > 0$. There is $N_r \ge
  1$, depending on the fixed endpoint and $r$ but not on $t$, such that
  \[
    B_k(1)\le\frac{(\Lambda_{\rm st}+r)^k}{(k!)^Q},
    \qquad k\ge N_r.
  \]
  Static coefficient homogeneity and the already established coefficient
  H\"older inequality therefore give, simultaneously for every $t > 0$ and $k
  \ge N_r$,
  \[
    A_k(t)\le\frac{D_r^k t^{\beta k}}{(k!)^q},\qquad
    D_r=C_{\rm low}^{\,\vartheta}(\Lambda_{\rm st}+r)^{1-\vartheta},
    \qquad \beta=4-\alpha-\alpha_0.
  \]
  For each $1 \le k < N_r$ the same inequality with the exact static coefficient
  gives $A_k(t) \le d_kt^{\beta k}$, where
  \[
    d_k=\frac{C_{\rm low}^{\,\vartheta k}
    B_k(1)^{1-\vartheta}}{(k!)^{2\vartheta}}<\infty.
  \]
  Adding $A_0 = 1$ and the nonnegative series yields
  \[
    \E u^2(t,x)\le P_r(t)+e^{qD_r^{1/q}t^\kappa},\qquad
    P_r(t)=1+\sum_{k=1}^{N_r-1}d_kt^{\beta k}.
  \]
  This finite power sum has $\log P_r(t)/t^\kappa \to 0$; its exponents need not
  be integers. The logarithm of a sum of two positive terms is at most the
  larger logarithm plus $\log2$. Consequently, the required limsup is at most
  $qD_r^{1/q}$. Let $r\downarrow0$ to obtain~\eqref{E:three-variational-upper}.
  This also covers a zero static rate, without division by that rate. The exact
  static homogeneity makes $N_r$ independent of $t$; no target coefficient
  homogeneity or smoothed static asymptotic was invoked. The finitely many low
  orders were kept, so this argument does not replace the all-time constant by
  the asymptotic one. This proves Proposition~\ref{P:three-transfer}.
\end{proof}

\LeanMarker{C:three-product-riesz}
\begin{corollary}[Ordinary product and Riesz covariances in dimension three]
  \label{C:three-product-riesz}
  Let $d = 3$, $0 < \alpha_0 < 1$, and $1 < \alpha < 2-\alpha_0$. The
  coefficient and all-time moment estimates in~\eqref{E:three-transfer}, and the
  asymptotic estimate in~\eqref{E:three-variational-upper}, hold for either of
  the following spatial covariances, with $c > 0$:
  \begin{enumerate}
    \item $\gamma(x) = c\prod_{j=1}^3|x_j|^{-a_j}$, where
      $0 \le a_j < 1$ and $\sum_j a_j = \alpha$;
    \item $\gamma(x) = c|x|^{-\alpha}$.
  \end{enumerate}
  In the first case the constant depends only on $a_1,a_2,a_3,c,\alpha_0$; in
  the second it depends only on $\alpha,c,\alpha_0$. A zero coordinate exponent
  denotes a constant factor; the case $a_j = 1$ and any interpretation of it as
  coordinate white noise are excluded.
\end{corollary}

\begin{proof}
  In both cases choose the fixed midpoint $\vartheta = (\alpha_0+2-\alpha)/2$
  and define $b_0,\chi$ by~\eqref{E:three-parameters}. Then $\alpha_0 <
  \vartheta < 2-\alpha$, and all subsequent choices depend only on the data in
  the statement. For the product covariance, Lemma~\ref{L:product-endpoints}
  supplies $b,v$ by its explicit interval formula. Its same-amplitude endpoint
  densities $\gamma_b,\gamma_v$ are admissible covariance measures of exponents
  one and $\chi$. Their common-width Gaussian inequality is exactly
  Hypothesis~\ref{H:three-endpoints}.

  For the Riesz covariance take $\gamma_1(x) = c|x|^{-1}$ and $\gamma_\chi(x) =
  c|x|^{-\chi}$. Each degree $e \in \{1,\alpha,\chi\}$ is less than three, so
  the corresponding density is locally integrable and tempered. Positive
  definiteness follows from the Gaussian mixture
  \[
    c|x|^{-e}=\frac{c}{\Gamma(e/2)}
    \int_0^\infty s^{e/2-1}e^{-s|x|^2}\ud s.
  \]
  For a Schwartz autocorrelation, absolute integrability follows from $e < 3$
  near zero and rapid decay at infinity. Tonelli and Fubini therefore justify
  integrating the positive Gaussian quadratic forms. The affine identity $\alpha
  = \vartheta+(1-\vartheta)\chi$ gives $\gamma_\alpha =
  \gamma_1^\vartheta\gamma_\chi^{1-\vartheta}$ almost everywhere, with the same
  amplitude $c$. H\"older on the Gaussian probability measure proves the
  required common-width comparison for every $\epsilon,z$. Thus, in both cases
  all endpoint hypotheses are discharged. Proposition~\ref{P:three-transfer}
  applies. This proves Corollary~\ref{C:three-product-riesz}.
\end{proof}

\LeanMarker{C:three-finite-images}
\begin{corollary}[Finite sums of linear images]\label{C:three-finite-images}
  Let $d = 3$, $0 < \alpha_0 < 1$, and $1 < \alpha < 2-\alpha_0$. Let
  $\Gamma_i$, $1 \le i \le N$, belong to one of the following classes, always
  with homogeneity exponent $\alpha$:
  \begin{enumerate}
    \item the ordinary products in Corollary~\ref{C:three-product-riesz};
    \item the isotropic Riesz covariances in that corollary;
    \item the plane-supported measures of Lemma~\ref{L:plane-endpoints},
      with $a_1+a_2 = \alpha-1$.
  \end{enumerate}
  Let $T_i:\R^3 \to \R^3$ be invertible linear maps, and let $w_i \ge 0$ with
  $\sum_iw_i > 0$. Then the coefficient, all-time moment, and nonsharp
  variational estimates
  \eqref{E:three-transfer}--\eqref{E:three-variational-upper} hold for
  \[
    \Gamma=\sum_{i=1}^Nw_i(T_i)_\#\Gamma_i.
  \]
  The constants may depend on all fixed component data, maps and weights. The
  endpoints and their variational quantities are those constructed in the proof.
\end{corollary}

\begin{proof}
  Choose the common midpoint $\vartheta = (\alpha_0+2-\alpha)/2$, and let
  $b_0,\chi$ be given by~\eqref{E:three-parameters}.
  Lemma~\ref{L:product-endpoints} supplies the ordinary-product endpoints. For a
  Riesz component use $c_i|x|^{-1}$ and $c_i|x|^{-\chi}$: the Gaussian-mixture
  calculation in the preceding proof verifies their covariance properties, and
  $\alpha = \vartheta+(1-\vartheta)\chi$ gives their Gaussian comparison by
  H\"older. Lemma~\ref{L:plane-endpoints} supplies the plane-supported
  endpoints. Thus, every component has an endpoint triple with these same fixed
  parameters.

  Apply Lemma~\ref{L:endpoint-closure} first to $T_i$ on each triple, then to
  the finite weighted sums. The resulting endpoints are
  \[
    \Gamma_{\rm low}=\sum_iw_i(T_i)_\#\Gamma_{{\rm low},i},
    \qquad
    \Gamma_{\rm st}=\sum_iw_i(T_i)_\#\Gamma_{{\rm st},i}.
  \]
  They satisfy Hypothesis~\ref{H:three-endpoints}. The target is nonzero because
  all components are positive nonzero measures and some weight is positive.
  Proposition~\ref{P:three-transfer} now applies to these actual target and
  endpoint measures. No addition or coupling identity for random solutions is
  used. This proves Corollary~\ref{C:three-finite-images}.
\end{proof}

\begin{remark}[Scope and calibrations]\label{R:three-endpoint-scope}
  As $\alpha_0\downarrow0$, the time exponent tends to $(4-\alpha)/(3-\alpha)$,
  agreeing with the static source without interchanging limits or asserting
  uniform constants. Cauchy--Schwarz gives $\E u(t,x) \le
  \exp(qC^{1/q}t^\kappa/2)$. Conditionally on a first-moment lower bound with
  the intended constant $L_*$, both upper constants above are at least $2L_*$,
  by the same inequality and the argument of Proposition~\ref{P:Jensen}. The
  author version's printed variational definitions (1.11)--(1.13), with gradient
  penalties one and one-half, require normalization clarification. Their
  intended correction remains unspecified. This qualification concerns the
  source-dependent calibration; neither the first-moment upper constant nor that
  calibration is an input to the upper-bound proof here. Neither upper constant
  is asserted sharp. Comparisons with Skorohod exponents are restricted to the
  published Riesz and strictly positive product-exponent cases satisfying their
  stated hypotheses; see Balan--Conus~\cite{balan.conus:16:intermittency}*{
  Theorem~2.1(b), p.~1492}. Their different solution concept and covariance
  normalization are not inputs here. The zero-exponent product cases,
  plane-supported measures, linear images and finite sums above do not enlarge
  that published comparison. The product construction changes coordinate powers;
  it is not Euclidean radial lowering of the original covariance. For a general
  three-dimensional covariance, Hypothesis~\ref{H:three-endpoints} remains an
  additional condition, verified here for the stated families and their finite
  positive sums of invertible linear images. No representation or density
  theorem for arbitrary covariances is asserted, and no uniform estimate for a
  limiting infinite mixture is proved. Section~\ref{S:interval-upper} treats the
  full covariance class by a different argument. These endpoint comparisons are
  not inputs to the intrinsic constant theorem~\ref{T:intrinsic-second-moment}.
\end{remark}

Euclidean radial lowering can destroy positive definiteness, even for an
admissible coordinate-product covariance; see No-go~\ref{P:nogo-09}.
Two other substitutes fail for different reasons: a heat--Mellin square can
have infinite local mass, while angular smoothing of a singular covariance
cannot retain a uniform isotropic density envelope; see
No-gos~\ref{P:nogo-12} and~\ref{P:nogo-13}.

\section{Physical and spectral energy for measure covariances}
\label{S:energy-fidelity}

This section justifies the positive-measure energy used
in~\eqref{E:three-static-variational}. It does not assume
Hypothesis~\ref{H:three-endpoints} or the conditional moment bound. The imported
analytic input is the finiteness of the variation on the right of (1.5), stated
in the closing sentence of Chen~\cite{chen:23:exponential}*{Theorem~1.1, p.~2}.

This identification of variational functionals does not remove the additional
hypotheses of the separate moment argument. It does not resolve general
covariance factorization.

\begin{lemma}[Sobolev energy identification]\label{L:physical-spectral-energy}
  Let $\Xi$ be a nonnegative, even, positive-definite tempered Radon covariance
  on $\R^3$, with $\Xi(rA) = r^{3-\chi}\Xi(A)$ for $0 < \chi < 2$. Let $\mu_\Xi$
  be its positive even spectral measure in the convention $\Xi = \int
  e^{i\xi\cdot x}\mu_\Xi(\ud\xi)$, understood distributionally, and use the
  Fourier transform $\widehat{h}(\xi) = \int e^{-i\xi\cdot x}h(x)\ud x$. For
  every real $g \in W^{1,2}(\R^3)$,
  \begin{equation}\label{E:physical-spectral-energy}
    \mathcal{Q}_\Xi(g)
    :=\int_{\R^3}\int_{\R^3}g(y+z)^2g(y)^2\ud y\,\Xi(\ud z)
    =\int_{\R^3}|\widehat{g^2}(\xi)|^2\mu_\Xi(\ud\xi)<\infty.
  \end{equation}
  The physical expression is independent of the Borel representative of $g$.
  There is a finite $C_0 \ge 0$, depending on $\Xi$, such that
  \begin{equation}\label{E:energy-quartic-bound}
    \mathcal{Q}_\Xi(g)
    \le\frac12\|g\|_2^2\|\nabla g\|_2^2+C_0\|g\|_2^4.
  \end{equation}
\end{lemma}

\begin{proof}
  We use the imported finiteness only on smooth profiles. Two distinct
  Cauchy--Schwarz bounds control the physical and spectral approximation
  errors. We then identify both strong limits and pass the quartic bound to
  the full Sobolev class.

  \medskip
  \noindent\textit{Step 1.~}We obtain the smooth-profile energy bound.
  Fourier uniqueness gives $\mu_\Xi(rA) = r^\chi\mu_\Xi(A)$. Since its unit-ball
  mass is finite, the dyadic estimate with $\sum_{j\ge0}2^{j(\chi-2)}$ verifies
  Dalang's condition. For a real Schwartz function $h$ put
  \[
    \mathscr{E}(h)=\int|\widehat{h}(\xi)|^2\mu_\Xi(\ud\xi),
    \qquad Q(f)=\mathscr{E}(f^2).
  \]
  Schwartz Fourier pairing identifies $Q(f)$ with its ordinary positive-measure
  integral. From the finiteness assertion in
  Chen~\cite{chen:23:exponential}*{Theorem~1.1, closing sentence for the
  variation in (1.5), p.~2}, restricted to real unit-$L^2$ Schwartz profiles,
  there is $C_0 \ge 0$ such that $Q(g) \le \frac12\|\nabla g\|_2^2+C_0$ for
  $\|g\|_2 = 1$. For $f \ne 0$, quartic homogeneity gives
  \[
    Q(f)=\|f\|_2^4Q(f/\|f\|_2)
    \le\frac12\|f\|_2^2\|\nabla f\|_2^2+C_0\|f\|_2^4.
  \]
  The same inequality is immediate for $f = 0$. Only Schwartz pairings have been
  used at this point, so this import presupposes no energy identification for
  nonsmooth profiles. If $N(f)^2 = \|f\|_2^2+\|\nabla f\|_2^2$, it follows that
  $Q(f)^{1/2} \le D N(f)^2$ for some finite $D$.

  \medskip
  \noindent\textit{Step 2.~}We estimate the two types of bilinear pairing.
  We need two different Cauchy--Schwarz inequalities. For real Schwartz $f,h$,
  spectral Cauchy--Schwarz gives
  \[
    \begin{split}
      0\le P(f,h)
      &:=\int\Xi(\ud z)\int f(y+z)^2h(y)^2\ud y\\
      &=\int\widehat{f^2}(\xi)
      \overline{\widehat{h^2}(\xi)}\mu_\Xi(\ud\xi)
      \le Q(f)^{1/2}Q(h)^{1/2}.
    \end{split}
  \]
  Evenness makes the Fourier sign convention in this real pairing immaterial. On
  the other hand, Cauchy--Schwarz on the positive measure $\ud y\,\Xi(\ud z)$
  gives
  \[
    \mathscr{E}(fh)
    =\int\Xi(\ud z)\int
    f(y+z)f(y)h(y+z)h(y)\ud y
    \le Q(f)^{1/2}Q(h)^{1/2}.
  \]
  All these are absolutely convergent Schwartz pairings. In particular,
  $\mathscr{E}(fh)$ is the quadratic energy of $fh$, not the quartic expression
  $Q(fh)$.

  \medskip
  \noindent\textit{Step 3.~}We compare the physical and spectral approximation
  maps.
  Let $\mathfrak{m}(\ud y,\ud z) = \ud y\,\Xi(\ud z)$ and define
  \[
    A_f(y,z)=f(y+z)f(y),\qquad B_f(\xi)=\widehat{f^2}(\xi).
  \]
  The measure $\mathfrak{m}$ is positive and sigma-finite. The decomposition
  \[
    A_f-A_h=(f-h)(y+z)f(y)+h(y+z)(f-h)(y)
  \]
  and the bound for $P$ imply
  \[
    \|A_f-A_h\|_{L^2(\mathfrak{m})}
    \le D\bigl(N(f)+N(h)\bigr)N(f-h).
  \]
  Since $f^2-h^2 = (f-h)(f+h)$, the bound for $\mathscr{E}$ gives
  \[
    \|B_f-B_h\|_{L^2(\mu_\Xi)}
    \le D N(f-h)N(f+h)
    \le D\bigl(N(f)+N(h)\bigr)N(f-h).
  \]

  \medskip
  \noindent\textit{Step 4.~}We identify their strong limits for a Sobolev
  profile.
  Choose real Schwartz $f_n \to g$ in $W^{1,2}$. These estimates make both
  $A_{f_n}$ and $B_{f_n}$ Cauchy in their respective $L^2$ spaces. Pass to a
  subsequence with $f_n \to g$ Lebesgue-almost everywhere. Translation
  invariance of Lebesgue measure and sigma-finite Fubini give
  \[
    f_n(y+z)f_n(y)\longrightarrow g(y+z)g(y)
    \quad\text{for }\mathfrak{m}\text{-almost every }(y,z).
  \]
  Thus, the strong $L^2(\mathfrak{m})$ limit is $A_g$. The same Fubini argument
  shows that changing the Borel representative of $g$ does not change this
  limit.

  Also
  \[
    \|f_n^2-g^2\|_1
    \le\|f_n-g\|_2(\|f_n\|_2+\|g\|_2)\longrightarrow0.
  \]
  Hence, $\widehat{f_n^2} \to \widehat{g^2}$ uniformly on all frequencies. A
  subsequence of its strong $L^2(\mu_\Xi)$ limit converges $\mu_\Xi$-almost
  everywhere, and must agree with this uniform limit. Consequently, that strong
  limit is $B_g$. For every $n$, $\|A_{f_n}\|_{L^2(\mathfrak{m})}^2 = Q(f_n) =
  \|B_{f_n}\|_{L^2(\mu_\Xi)}^2$. Passing to the strong limits proves
  \eqref{E:physical-spectral-energy}, including finiteness. Passing the Schwartz
  quartic bound to the same limit proves \eqref{E:energy-quartic-bound}. This
  proves Lemma~\ref{L:physical-spectral-energy}.
\end{proof}

In particular, for $\|g\|_2 = 1$ the energy is at most $\|\nabla g\|_2^2/2+C_0$.
Therefore, the supremum of $\mathcal{Q}_\Xi(g)^{1/2}-\|\nabla g\|_2^2$ is
finite. Applied to $\Xi = \Gamma_{\rm st}$, the profile-by-profile identity
justifies the measure-form functional in~\eqref{E:three-static-variational}
without assigning a density to the covariance or borrowing a disputed
special-case constant.

\chapter{Static comparisons and the clock obstruction}\label{S:appendix-static-clock}

Two comparisons are natural here, and both stop short of the constant. As the
temporal correlation approaches the static case, how do the moments compare
with the static ones? And the first moment has a Brownian clock
representation: can it be copied for the second moment? The first question
gives useful inequalities. The second gets a negative answer at low order, for
the construction tried here. This appendix has no figure.

The appendix recalls Proposition~\ref{P:mixed}, the upper bound of
Chapter~\ref{S:book-independent-upper} and the source-version qualifications
in Section~\ref{S:known}. Section~\ref{S:near-static-lower} compares the
solution with its static counterpart
(Proposition~\ref{P:static-minorant-lower}) and shows, by temporal
interpolation, that the upper growth rates are log-convex under the stated
assumptions (Proposition~\ref{P:moment-growth-logconvex}).
Section~\ref{S:two-clock-screen} tests direct clock copying at small order
(Proposition~\ref{P:clock-copy-obstruction}). The obstruction concerns that
one representation, not every two-clock construction.

Neither comparison identifies the constant of
Theorem~\ref{T:intrinsic-second-moment}, and the main argument does not rely
on this appendix. Episode~\VideoRef{9}{eLz9Hd20ZT0} of the video series
covers it.

\section{A static lower comparison and a near-static mean gap}
\label{S:near-static-lower}

The comparison in this section does not identify the intrinsic constant of
Theorem~\ref{T:intrinsic-second-moment}; no novelty or numerical normalization
is claimed.

Fix a nonzero spatial covariance in Hypothesis~\ref{H:homog}, and put
\[
  \begin{gathered}
    \nu=3-\alpha>1,\qquad a_{\max}=\min(1,2-\alpha),\qquad
    \beta(a)=4-\alpha-a,\\
    \kappa(a)=\frac{\beta(a)}\nu,\qquad
    \kappa_{\rm st}=\frac{4-\alpha}\nu.
  \end{gathered}
\]
Here the temporal exponent varies: $u^{(a)}$ denotes our canonical solution with
covariance $|s-r|^{-a}\Gamma$, for $0 < a < a_{\max}$. At $a = 0$ define
$u^{(0)} = U$ to be the static canonical solution with the same spatial
covariance and initial data $(1,0)$. This endpoint means the temporal covariance
one, not a substitution into $\Gamma_{\mathrm{E}}(a)$. The notation
$\gamma_{\delta,a}^0$ records the temporal exponent $a$ in the same
Laplace-prefix regularization previously denoted by $\gamma_\delta^0$.

\paragraph{Imported static input.}
From Chen--Hu~\cite{chen.hu:25:hyperbolic}*{Theorem~1.2, (1.9), p.~2147}, the
static first-mean rate $L_0$ exists and
\begin{equation}\label{E:static-ratio-input}
  0<L_0:=\lim_{t\to\infty}t^{-\kappa_{\rm st}}\log\E U(t,x)<\infty,
  \qquad
  \lim_{t\to\infty}t^{-\kappa_{\rm st}}\log\E U(t,x)^2
  =2^{\kappa_{\rm st}}L_0.
\end{equation}
We use only the common $p^{\kappa_{\rm st}}$ factor for $p = 1,2$; no numerical
value of the static variational constant or its white-noise specialization is
adopted. Strict positivity in that formula follows by using a nonzero Gaussian
competitor: its interaction and gradient terms scale as $A c^{\alpha/2}$ and $B
c^2$, with $A,B > 0$, so their difference is positive for small $c$. Finiteness
and the moment law are the imported source result. The static kernel and
integral conventions are its (2.10)--(2.14), p.~2150; the fixed-order limits
agree with its Theorem~3.8(i), p.~2165. We supply the static series convergence
needed below from our regularized fractional-time upper bound, without using the
source's quantitative (4.3). The generalized-covariance convention is explicit
in its (2.18), p.~2152. Our spectral ball law gives
\[
  \int\frac{\mu(\ud\xi)}{1+|\xi|^2}
  \le\mu(B_1)+C_\alpha\mu(B_1)\sum_{j\ge0}2^{j(\alpha-2)}<\infty.
\]
Thus, the static Dalang condition holds in every admissible dimension. Measure
homogeneity of degree $d-\alpha$ is equivalent to distributional homogeneity of
degree $-\alpha$, so the source's (1.8) applies in its explicit
generalized-covariance convention. The static object is identified by its
canonical heat-limit expansion; no uniqueness assertion for other mild solutions
is used.

For $a > 0$ write $X_n^{(a)}$ for our fixed-order terms; at zero use the
corresponding static terms. Set
\[
  m_k(a)=\E X_{2k}^{(a)}(1,0),\qquad m_0(a)=1.
\]

\begin{lemma}[Mean coefficients for the power covariance]\label{L:mean-coefficient-holder}
  For every $0 \le a < a_{\max}$, the coefficients $m_k(a)$ are nonnegative, and
  \begin{equation}\label{E:mean-positive-series}
    \E u^{(a)}(t,x)=\sum_{k\ge0}m_k(a)t^{\beta(a)k}.
  \end{equation}
  The series $F_a(z) = \sum_{k\ge0}m_k(a)z^k$ is finite for every $z \ge 0$. If
  $0 < a < b < a_{\max}$, then
  \begin{equation}\label{E:mean-coefficient-holder}
    m_k(a)\le m_k(0)^{1-a/b}m_k(b)^{a/b},\qquad k\ge1.
  \end{equation}
\end{lemma}

\begin{proof}
  We construct the static endpoint series using a regularized positive
  comparison. We then interpolate single-chain coefficients on one common
  spatial measure, remove the spatial width, and recover the mean series by
  homogeneity.

  \medskip
  \noindent\textit{Step 1.~}We construct the static fixed-order limits and
  summable series.
  We first justify the static endpoint and its norm series. At positive spatial
  widths its ordered terms are the bounded-variance Gaussian Bochner integrals
  of Lemma~\ref{L:terminal-wick}, with a time-constant field.
  Lemma~\ref{L:spatial-wick-measures}, applied with temporal multiplier one to
  the diagonal and mixed-width second moments, gives a common finite diagram
  limit. The same Hilbert Cauchy argument as in Lemma~\ref{L:joint-fixed-order}
  therefore gives the static fixed-order limits $X_n^{(0)}$, with the source's
  heat-limit convention. This uses no infinite solution series.

  Fix any $b_0 \in (0,a_{\max})$ and $\delta_0 > 0$, and put $c =
  \gamma_{\delta_0,b_0}^0(1) > 0$. Let $U_N^\epsilon(r;1,x)$ be the static
  partial series at horizon one and noise amplitude $r > 1$. Compare its finite
  Wick sum with that of the fractional-time partial series at widths
  $\epsilon,\delta_0$, exponent $b_0$ and amplitude $r/\sqrt c$. Since
  $\gamma_{\delta_0,b_0}^0(s-r') \ge c$ for $s,r' \in [0,1]$, every pair factor
  on the second side dominates the static pair factor. Positivity, then the full
  regularized bound \eqref{E:regularized-scaled-moment}, gives
  \[
    \E|U_N^\epsilon(r;1,x)|^2
    \le \E|u^{\epsilon,\delta_0}_{b_0}(r/\sqrt c;1,x)|^2
    \le B_{b_0}(1;r/\sqrt c)<\infty,
  \]
  uniformly in $N,\epsilon$. The notation $B_{b_0}$ means the bound
  in~\eqref{E:regularized-scaled-moment} at temporal exponent $b_0$. Both
  positive-width full series exist by Lemma~\ref{L:terminal-wick}. Letting $N
  \to \infty$ at fixed $\epsilon$ preserves this bound.
  Lemma~\ref{L:uniform-series-tail} now applies to the static family (with only
  the spatial width varying). It proves absolute $L^2$ convergence of its
  limiting series at horizon one, with geometric coefficient bounds for every $r
  > 1$. Static fixed-order homogeneity gives
  \[
    \|X_n^{(0)}(t,x)\|_2
    =t^{(4-\alpha)n/2}\|X_n^{(0)}(1,0)\|_2.
  \]
  Choosing $r > \max\{1,t^{(4-\alpha)/2}\}$ in those bounds proves norm
  summability at every fixed $t > 0$. The resulting series has precisely the
  source's canonical static kernels and fixed-order limits.

  \medskip
  \noindent\textit{Step 2.~}We interpolate the temporal powers at fixed
  spatial smoothing.
  For a single chain of length $2k$, use the initial-chain variables of
  Proposition~\ref{P:mixed}. At common positive spatial width $\epsilon$, let
  $\omega_{\epsilon,D}$ be its finite positive spatial Green--Wick measure on
  the unit time simplex, with covariance $\gamma_{2\epsilon}$ on every pair. It
  is finite because that covariance is bounded and the Green masses are the time
  gaps. Define
  \[
    m_k^\epsilon(a)=\sum_{D\in\Pi_k}\int
    \prod_{\{j,h\}\in D}|s_j-s_h|^{-a}\ud\omega_{\epsilon,D},
    \qquad 0\le a<a_{\max}.
  \]
  These time integrals for the power covariance are finite. At fixed $\epsilon$
  the spatial factors and Green masses are bounded on the finite horizon; the
  disjoint pair endpoints give integrable powers on the full time cube since $a
  < 1$. Diagonals are null there. On the disjoint union of the pairing measures,
  the time multiplier at $a$ is exactly the $a/b$ power of its multiplier at
  $b$. H\"older's inequality therefore gives
  \[
    m_k^\epsilon(a)\le
    m_k^\epsilon(0)^{1-a/b}m_k^\epsilon(b)^{a/b}.
  \]
  This is an identity for the power covariance; no exponent-dependent Laplace
  cutoff or Euler Gamma normalization is being interpolated.

  \medskip
  \noindent\textit{Step 3.~}We remove the common spatial width.
  We justify removal of the common spatial width. For each positive exponent,
  the proof of Lemma~\ref{L:joint-fixed-order} gives this limit with the power
  covariance as well. Indeed, writing $F_{\epsilon,\eta,D}$ for the diagram with
  common temporal cutoff $\eta$, \eqref{E:joint-temporal-tail} gives
  \[
    F_{\epsilon,\eta,D}\le F_{\epsilon,0,D}
    \le F_{\epsilon,\eta,D}+r_a(\eta),\qquad r_a(\eta)\longrightarrow0,
  \]
  uniformly in $\epsilon$. First pass to the spatial limit at fixed $\eta$, then
  send $\eta$ to zero. The resulting limit is exactly the joint diagram value
  from that lemma. At zero exponent the multiplier is one, so
  Lemma~\ref{L:spatial-wick-measures} gives the spatial limit directly; the
  positive-width Wick calculation and the imported static fixed-order
  convergence identify it with the static expectation. Thus, $m_k^\epsilon(a)
  \to m_k(a)$ at all three fixed exponents $0,a,b$,
  proving~\eqref{E:mean-coefficient-holder} and nonnegativity. No convergence
  uniform in $a$ is claimed.

  \medskip
  \noindent\textit{Step 4.~}We identify the homogeneous positive mean series.
  Scaling all times and spatial displacements by $t$ in a regularized one-chain
  diagram gives the factor $t^{2(2k)-\alpha k-ak} = t^{\beta(a)k}$ and changes
  the spatial width to $\epsilon/t^2$. Passing to the just established limits
  proves $\E X_{2k}^{(a)}(t,x) = t^{\beta(a)k}m_k(a)$. Odd expectations vanish
  by Wick's formula. For positive $a$ the independently proved norm summability
  in Lemma~\ref{L:vertex-crude} justifies the mean-series interchange. At zero
  the static norm series just proved does so, with nonnegative even means. This
  proves \eqref{E:mean-positive-series}. Since $\beta(a) > 0$, every $z \ge 0$
  equals $t^{\beta(a)}$ for some $t \ge 0$, and the series is finite. This
  proves Lemma~\ref{L:mean-coefficient-holder}.
\end{proof}

\LeanMarker{L:positive-series-type}
\begin{lemma}[Type of a positive coefficient series]\label{L:positive-series-type}
  Let $\nu \ge 1$, $c_k \ge 0$, $c_0 = 1$, and suppose $F(z) =
  \sum_{k\ge0}c_kz^k < \infty$ for every $z \ge 0$. If
  \[
    \tau=\limsup_{z\to\infty}z^{-1/\nu}\log F(z)<\infty,
  \]
  then
  \begin{equation}\label{E:positive-series-type}
    K:=\limsup_{k\to\infty}\big((k!)^\nu c_k\big)^{1/k}<\infty,
    \qquad \tau=\nu K^{1/\nu}.
  \end{equation}
  Conversely, finite $K$ implies the same finite $\tau$ and equality.
\end{lemma}

\begin{proof}
  We first bound each coefficient by evaluating the positive series at an
  order-dependent argument. Conversely, a finite coefficient type gives an
  exponential series bound after keeping finitely many low orders separately.

  \medskip
  \noindent\textit{Step 1.~}We pass from growth type to coefficient type.
  If $\tau < \infty$, choose $L > \tau$. Since $c_0 = 1$, $\tau \ge 0$, and
  continuity of the power series on compact positive intervals gives $F(z) \le
  C_L e^{Lz^{1/\nu}}$ for all $z \ge 0$. For $k \ge 1$, the inequality $c_kz^k
  \le F(z)$ at $z = (\nu k/L)^\nu$ yields
  \[
    c_k\le C_L e^{\nu k}\left(\frac L{\nu k}\right)^{\nu k}.
  \]
  Monotonicity of $\log x$ and integration on $[1,k]$ give
  \[
    k\log k-k+1\le\log(k!)
    \le k\log k-k+1+\log k.
  \]
  Consequently, $K \le (L/\nu)^\nu$. Letting $L\downarrow\tau$ proves $K \le
  (\tau/\nu)^\nu$, including $K = 0$ when $\tau = 0$.

  \medskip
  \noindent\textit{Step 2.~}We recover growth type from the coefficient bound.
  Conversely, for any $A > K$ all sufficiently large coefficients obey $c_k \le
  A^k/(k!)^\nu$. Hence, for a polynomial $P$ with nonnegative coefficients,
  \[
    F(z)\le P(z)+\sum_{k\ge0}\frac{(Az)^k}{(k!)^\nu}
    \le P(z)+\exp\{\nu(Az)^{1/\nu}\}.
  \]
  The last inequality is the elementary $\ell^\nu$ bound $\sum a_k^\nu \le (\sum
  a_k)^\nu$, applied to $a_k = (Az)^{k/\nu}/k!$. The polynomial contributes zero
  to the normalized logarithmic growth. Therefore, $\tau \le \nu A^{1/\nu}$; let
  $A\downarrow K$. Combining both inequalities proves the equality when either
  quantity is finite. This proves Lemma~\ref{L:positive-series-type}.
\end{proof}

\begin{proposition}[Static temporal-minimum comparison]\label{P:static-minorant-lower}
  For every $0 < a < a_{\max}$, $T > 0$ and $p \in \{1,2\}$,
  \begin{equation}\label{E:static-minorant-moment}
    \E[u^{(a)}(T,x)^p]
    \ge\E\left[U\left(T^{1-a/(4-\alpha)},x\right)^p\right].
  \end{equation}
  In particular,
  \begin{equation}\label{E:static-minorant-rate}
    \liminf_{T\to\infty}T^{-\kappa(a)}\log\E[u^{(a)}(T,x)^p]
    \ge p^{\kappa_{\rm st}}L_0,\qquad p=1,2.
  \end{equation}
\end{proposition}

\begin{proof}
  We compare finite regularized moments with a static model using the minimum
  temporal covariance on the horizon. We remove regularizations at fixed
  order, sum, and then identify the static effective time.

  \medskip
  \noindent\textit{Step 1.~}We compare the finite expected powers at positive
  widths.
  Fix positive common spatial and temporal noise widths $\epsilon,\delta$. On
  $[0,T]$, the source's Laplace-prefix temporal covariance satisfies
  \[
    \gamma_{\delta,a}^0(s-r)\ge c_{\delta,T}:=
    \gamma_{\delta,a}^0(T)>0,\qquad
    c_{\delta,T}\longrightarrow T^{-a}\quad(\delta\downarrow0).
  \]
  Compare the order-$N$ partial series of $u^{(a)}$ with the same partial series
  of a static model having covariance $c_{\delta,T}\gamma_{2\epsilon}$. At
  positive widths, the mean or squared moment is a finite Wick sum of the same
  terminal Green kernels and pairings. Every spatial factor is the same
  nonnegative $\gamma_{2\epsilon}$; every time factor on the first side is at
  least $c_{\delta,T}$. Thus, the finite expected powers are ordered for $p =
  1,2$. For $d = 3$, remove Green smoothing first at these fixed positive noise
  widths using the bounded covariance calculation preceding
  Lemma~\ref{L:spatial-wick-measures}.

  \medskip
  \noindent\textit{Step 2.~}We pass fixed-order limits and the convergent
  series.
  Now send $(\epsilon,\delta)$ jointly to zero at fixed $N$.
  Lemma~\ref{L:joint-fixed-order} gives the fractional-time $L^2$ limits. Static
  fixed-order $L^2$ convergence gives all static split products in $L^1$, by
  Cauchy--Schwarz. Dependence on $c_{\delta,T}$ is a finite polynomial at fixed
  order, so its limit is the static model with noise amplitude $\sigma =
  T^{-a/2}$. This requires no uniform approximation in $T$ or in a varying
  infinite-order series. At every fixed $\sigma,T > 0$, static norm homogeneity
  gives
  \[
    \sigma^n\|X_n^{(0)}(T,x)\|_2
    =\|X_n^{(0)}(\sigma^{2/(4-\alpha)}T,x)\|_2.
  \]
  Thus, the amplitude-scaled static norm series is summable at the effective
  time already covered by Lemma~\ref{L:mean-coefficient-holder}. Finally, let $N
  \to \infty$ using the two $L^2$ expansions, with the static one supplied in
  Lemma~\ref{L:mean-coefficient-holder}. Means and squared moments converge,
  preserving the finite-order inequality. No covariance-order coupling of the
  two Gaussian fields is asserted.

  \medskip
  \noindent\textit{Step 3.~}We identify the static amplitude clock and its
  rate.
  Write $U_\sigma$ for the static solution with noise amplitude $\sigma$. In its
  $p$th moment, a total order $2k$ has amplitude factor $\sigma^{2k}$ and, by
  static diagram homogeneity, time factor $T^{(4-\alpha)k}$. For $p = 1,2$, the
  finite sums and then the convergent nonnegative moment series therefore give
  the exact identity
  \[
    \E[U_\sigma(T,x)^p]
    =\E[U(\sigma^{2/(4-\alpha)}T,x)^p].
  \]
  Only moments are identified; no spatially indexed random fields are coupled.
  At $\sigma = T^{-a/2}$ this proves \eqref{E:static-minorant-moment}. The
  effective static time tends to infinity and
  \[
    \left(T^{1-a/(4-\alpha)}\right)^{\kappa_{\rm st}}
    =T^{\kappa(a)}.
  \]
  The imported static rates~\eqref{E:static-ratio-input} now prove
  \eqref{E:static-minorant-rate}. This proves
  Proposition~\ref{P:static-minorant-lower}.
\end{proof}

\begin{proposition}[Near-static strict gap from the mean rate]
  \label{P:near-static-mean-gap}
  Define
  \[
    M(a)=\limsup_{T\to\infty}T^{-\kappa(a)}\log\E u^{(a)}(T,x),\qquad
    K(a)=\limsup_{k\to\infty}\big((k!)^\nu m_k(a)\big)^{1/k}.
  \]
  Then $K(0) = (L_0/\nu)^\nu > 0$, and $K(b) < \infty$ for every $0 < b <
  a_{\max}$. Fix such a $b$ and choose any finite $B_b \ge \max\{2K(0),K(b)\}$.
  For $0 < a < b$,
  \begin{equation}\label{E:near-static-mean-upper}
    M(a)\le L_0\left(\frac{B_b}{K(0)}\right)^{a/(b\nu)}.
  \end{equation}
  In particular, with
  \begin{equation}\label{E:near-static-window}
    a_*:=\frac{b\log2}{\log(B_b/K(0))}\in(0,b],
  \end{equation}
  every $0 < a < a_*$ satisfies
  \begin{equation}\label{E:near-static-mean-gap}
    \liminf_{T\to\infty}T^{-\kappa(a)}\log\E[u^{(a)}(T,x)^2]
    >2M(a).
  \end{equation}
\end{proposition}

\begin{proof}
  We convert mean growth rates into coefficient types. Coefficient
  interpolation then bounds the positive-exponent mean rate, and the static
  lower comparison gives the strict gap on the stated window.

  \medskip
  \noindent\textit{Step 1.~}We identify finite endpoint coefficient types.
  By Lemma~\ref{L:mean-coefficient-holder}, the substitution $z = T^{\beta(a)}$
  identifies the normalized logarithmic type of $F_a(z)$ with the mean limsup
  $M(a)$; at zero it equals $L_0$. The imported static law and
  Lemma~\ref{L:positive-series-type} give $K(0) = (L_0/\nu)^\nu$. For positive
  $b$, Cauchy--Schwarz and Theorem~\ref{T:full-upper}, applied at temporal
  exponent $b$, give $M(b) < \infty$. The same type lemma therefore gives $K(b)
  < \infty$. There is no use of a printed first-moment rate constant.

  \medskip
  \noindent\textit{Step 2.~}We interpolate the coefficient types and apply the
  static lower rate.
  Multiplying~\eqref{E:mean-coefficient-holder} by $(k!)^\nu$, taking $k$th
  roots and then limsups yields
  \[
    K(a)\le K(0)^{1-a/b}B_b^{a/b}.
  \]
  To justify this even if coefficients vanish, for any $\eta > 0$ both endpoint
  root sequences are eventually bounded by $K(0)+\eta$ and $B_b+\eta$; take
  their weighted product and then let $\eta\downarrow0$.
  Lemma~\ref{L:positive-series-type} now gives
  \[
    M(a)=\nu K(a)^{1/\nu}
    \le L_0\left(\frac{B_b}{K(0)}\right)^{a/(b\nu)},
  \]
  proving~\eqref{E:near-static-mean-upper}. Since $B_b/K(0) \ge 2$, the window
  in~\eqref{E:near-static-window} is positive and at most $b$. For $0 < a <
  a_*$,
  \[
    2M(a)<2L_0\cdot2^{1/\nu}
    =2^{\kappa_{\rm st}}L_0,
  \]
  because $\kappa_{\rm st} = 1+1/\nu$. Combine this with
  Proposition~\ref{P:static-minorant-lower} at $p = 2$. This proves
  Proposition~\ref{P:near-static-mean-gap}.
\end{proof}

\begin{remark}[Scope and calibrations]\label{R:near-static-scope}
  The static lower comparison holds at every admissible positive temporal
  exponent; the strict gap from twice the first-moment limsup is asserted only
  in the nonempty window~\eqref{E:near-static-window}. Its endpoint depends on
  the fixed covariance and the chosen reference exponent $b$; no uniformity over
  covariance shapes is claimed. As $a\downarrow0$, the comparison clock has
  exponent tending to one, the growth exponent tends to $\kappa_{\rm st}$, and
  the mean upper bound tends to $L_0$. The static second-moment constant is
  exactly $2^{\kappa_{\rm st}}L_0$.

  The gap concerns $\E u$, which is positive by its pairing series; no theorem
  for $\E|u|$ or all moment orders is asserted. If an intended source first-mean
  lower rate $L_*(a)$ is valid, then $L_*(a) \le M(a)$,
  so~\eqref{E:near-static-mean-gap} also improves the corresponding Jensen rate
  $2L_*(a)$ in this window. This does not establish the original lower target
  $2^{\kappa_{\rm st}}L_*(a)$, or either sharp bound for arbitrary $a$. The
  intrinsic formula in Theorem~\ref{T:intrinsic-second-moment} does not
  identify either printed source endpoint. The proof here uses the ordinary
  Stratonovich pairing expansion, including self-contractions; no Skorohod
  result or sibling-project development is made.
\end{remark}

\subsection{Temporal interpolation of the full moment}
\label{S:temporal-type-continuity}

The constant at a fixed positive temporal exponent is identified  in
Theorem~\ref{T:intrinsic-second-moment}.

Retain the fixed nonzero covariance and the notation of this section. For
\(p=1,2\), set \(c_{p,0}(a)=1\) and, for \(k\ge1\),
\begin{equation}\label{E:full-moment-coefficients}
  c_{1,k}(a)=m_k(a),\qquad
  c_{2,k}(a)=\sum_{m=0}^{2k}
  \E[X_m^{(a)}(1,0)X_{2k-m}^{(a)}(1,0)].
\end{equation}
The zero-exponent terms are the static limits constructed in
Lemma~\ref{L:mean-coefficient-holder}.

\begin{lemma}[Interpolation of the full moment coefficients]
  \label{L:full-moment-coefficient-holder}
  For \(0\le a<a_{\max}\) and \(p=1,2\), the coefficients are nonnegative and
  \begin{equation}\label{E:full-moment-positive-series}
    \E[u^{(a)}(t,x)^p]
    =\sum_{k\ge0}c_{p,k}(a)t^{\beta(a)k}<\infty,\qquad t>0.
  \end{equation}
  If \(0\le a_1<a<a_2<a_{\max}\) and \(\theta=(a-a_1)/(a_2-a_1)\), then
  \begin{equation}\label{E:full-moment-coefficient-holder}
    c_{p,k}(a)\le
    c_{p,k}(a_1)^{1-\theta}c_{p,k}(a_2)^\theta,\qquad k\ge0.
  \end{equation}
\end{lemma}

\begin{proof}
  We interpolate the full coefficient sum at common spatial smoothing,
  retaining all splits. We then identify each endpoint limit and use diagram
  homogeneity to obtain the full positive moment series.

  \medskip
  \noindent\textit{Step 1.~}We place the mean and second-moment diagrams on
  their common measures.
  Work first at a common positive spatial noise width \(\epsilon\) and at unit
  horizon. For each split \(m+n=2k\) and pairing \(D\), let
  \(\omega_{\epsilon,m,n,D}\) be the finite positive spatial Green--Wick measure
  in the initial-chain convention of Proposition~\ref{P:mixed}, with covariance
  \(\gamma_{2\epsilon}\). For \(p=1\) take the disjoint union over the single
  specified split \((m,n)=(0,2k)\) and all pairings. For \(p=2\), take the
  disjoint union over every \(m=0,\ldots,2k\) and all pairings. In particular,
  the two empty-chain orientations are distinct terms for \(p=2\), whereas
  \(p=1\) has just one. Wick's formula has coefficient one for each pairing on
  each ordered domain.

  \medskip
  \noindent\textit{Step 2.~}We interpolate the temporal power factors.
  With \(\omega_{\epsilon,p,k}\) denoting this disjoint-union measure, define
  its time integral for the power covariance by
  \[
    c_{p,k}^{\epsilon}(a)=
    \int \prod_{\{j,h\}\in D}|s_j-s_h|^{-a}
    \ud\omega_{\epsilon,p,k}.
  \]
  The pairing label in the multiplier is that of its component. These integrals
  are finite: the spatial covariances and Green masses are bounded on the unit
  horizons, and the disjoint pair endpoints give integrable factors on the full
  time cube since \(a<1\). The multiplier at \(a\) is exactly the product of the
  endpoint multipliers to powers \(1-\theta,\theta\), off null diagonals.
  Hölder's inequality on the disjoint union gives
  \[
    c_{p,k}^{\epsilon}(a)\le
    c_{p,k}^{\epsilon}(a_1)^{1-\theta}
    c_{p,k}^{\epsilon}(a_2)^\theta.
  \]

  \medskip
  \noindent\textit{Step 3.~}We identify all fixed-order limits.
  For each fixed positive endpoint, the all-chain bound
  \eqref{E:joint-temporal-tail} and the common temporal-cutoff sandwich used in
  Lemma~\ref{L:mean-coefficient-holder} identify the time-domain spatial limit
  with the joint diagram value. This applies to every unequal or empty split,
  not only to a mean diagram. At exponent zero, the multiplier is one and
  Lemma~\ref{L:spatial-wick-measures} gives the same limit. The static diagonal
  and mixed-width Gaussian limits identify those diagram values with the split
  expectations of the constructed static terms: products converge in \(L^1\) by
  Cauchy--Schwarz, and the terminal-to-initial chain reversal is the one in
  Proposition~\ref{P:mixed}, applied before removing widths. For positive
  exponents this identification is Proposition~\ref{P:mixed}. The union is
  finite at fixed \(k\); passing to the spatial limit proves
  \eqref{E:full-moment-coefficient-holder} and nonnegativity. No interpolation
  of exponent-dependent Laplace prefixes is used.

  \medskip
  \noindent\textit{Step 4.~}We use homogeneity and the convergent moment
  series.
  Every diagram has \(2k\) time Jacobians, \(2k\) Green scaling factors, \(k\)
  spatial covariance factors and \(k\) temporal factors. Thus, its homogeneous
  time degree is \(4k-\alpha k-ak=\beta(a)k\), including both empty-chain cases.
  For positive \(a\), the full second-moment series follows from
  Proposition~\ref{P:mixed}; the mean series is \eqref{E:mean-positive-series}.
  At zero, the static norm series gives convergence of the finite partial means
  and squared moments. Each finite moment expansion has nonnegative mixed terms,
  so its increasing limit may be regrouped by total order \(2k\). The
  homogeneity just proved then gives~\eqref{E:full-moment-positive-series}. This
  proves Lemma~\ref{L:full-moment-coefficient-holder}.
\end{proof}

For \(p=1,2\) define the normalized time growth rates
\begin{equation}\label{E:moment-growth-types}
  \begin{split}
    \overline{\Lambda}_p(a)&=\limsup_{T\to\infty}
    T^{-\kappa(a)}\log\E[u^{(a)}(T,x)^p],\\
    \underline{\Lambda}_p(a)&=\liminf_{T\to\infty}
    T^{-\kappa(a)}\log\E[u^{(a)}(T,x)^p].
  \end{split}
\end{equation}
Here \(p=1\) concerns \(\E u\), not \(\E|u|\).

\LeanMarker{P:moment-growth-logconvex}
\begin{proposition}[Log-convexity of the upper growth type]
  \label{P:moment-growth-logconvex}
  The functions \(\overline{\Lambda}_p(a)\), \(p=1,2\), are positive and finite
  on \(0\le a<a_{\max}\). For the endpoints and \(\theta\) in
  Lemma~\ref{L:full-moment-coefficient-holder},
  \begin{equation}\label{E:moment-growth-logconvex}
    \overline{\Lambda}_p(a)\le
    \overline{\Lambda}_p(a_1)^{1-\theta}
    \overline{\Lambda}_p(a_2)^\theta .
  \end{equation}
  In particular, \(\log\overline{\Lambda}_p\) is convex and
  \(\overline{\Lambda}_p\) is continuous on \((0,a_{\max})\).
\end{proposition}

\begin{proof}
  We establish that all rates are positive and finite, convert coefficient
  interpolation into growth-type interpolation, and then use convex secant
  slopes to obtain interior continuity.

  \medskip
  \noindent\textit{Step 1.~}We verify positivity and finiteness of the rates.
  Put \(C_p=p^{\kappa_{\rm st}}L_0>0\). The imported static law gives
  \(\overline{\Lambda}_p(0)=\underline{\Lambda}_p(0)=C_p\).
  Proposition~\ref{P:static-minorant-lower} gives \(\underline{\Lambda}_p(a)\ge
  C_p\) for positive \(a\). The full upper bound gives finiteness for \(p=2\).
  Cauchy--Schwarz gives finiteness for \(p=1\), whose mean is at least one by
  its positive series. Hence, all logarithms and rates used here are well
  defined, positive and finite.

  \medskip
  \noindent\textit{Step 2.~}We transfer coefficient interpolation to the
  logarithmic type.
  Let \(F_{p,a}(z)=\sum_{k\ge0}c_{p,k}(a)z^k\). It is finite for every
  \(z\ge0\), by Lemma~\ref{L:full-moment-coefficient-holder}. Its normalized
  logarithmic type equals \(\overline{\Lambda}_p(a)\): spatial stationarity
  makes the moment independent of \(x\), while \(z=T^{\beta(a)}\) and
  \(z^{1/\nu}=T^{\kappa(a)}\). Lemma~\ref{L:positive-series-type} therefore
  gives
  \[
    K_p(a):=\limsup_{k\to\infty}
    \big((k!)^\nu c_{p,k}(a)\big)^{1/k}
    =\left(\frac{\overline{\Lambda}_p(a)}{\nu}\right)^\nu.
  \]
  Multiplying the coefficient interpolation by \((k!)^\nu\), taking \(k\)th
  roots and then limsups gives
  \[
    K_p(a)\le K_p(a_1)^{1-\theta}K_p(a_2)^\theta.
  \]
  More explicitly, the two endpoint root sequences are eventually bounded by
  \(K_p(a_i)+\eta\); let \(\eta\downarrow0\) after taking their weighted
  product. This argument permits zero or oscillating coefficient subsequences.
  Taking the \(1/\nu\) power and multiplying by \(\nu\)
  proves~\eqref{E:moment-growth-logconvex}.

  \medskip
  \noindent\textit{Step 3.~}We obtain continuity from the convex slope bounds.
  For completeness, continuity follows directly from secant slopes. Write
  \(f=\log\overline{\Lambda}_p\). If \(L<l\le x<y\le r<R\) are in
  \((0,a_{\max})\), convexity gives
  \[
    \frac{f(l)-f(L)}{l-L}\le
    \frac{f(y)-f(x)}{y-x}\le
    \frac{f(R)-f(r)}{R-r}.
  \]
  Both outer slopes are finite. They bound the absolute difference quotient on
  each compact interior interval, proving continuity of \(f\), and hence of its
  exponential. This proves Proposition~\ref{P:moment-growth-logconvex}.
\end{proof}

\LeanMarker{C:static-moment-constant-continuity}
\begin{corollary}[Recovery of the static moment constant]
  \label{C:static-moment-constant-continuity}
  Put \(C_p=p^{\kappa_{\rm st}}L_0\), \(p=1,2\). For every fixed
  \(0<b<a_{\max}\) and \(0<a<b\),
  \begin{equation}\label{E:static-constant-envelope}
    C_p\le\underline{\Lambda}_p(a)\le\overline{\Lambda}_p(a)
    \le C_p\left(\frac{\overline{\Lambda}_p(b)}{C_p}\right)^{a/b}.
  \end{equation}
  Consequently,
  \begin{equation}\label{E:static-constant-continuity}
    \lim_{a\downarrow0}\underline{\Lambda}_p(a)
    =\lim_{a\downarrow0}\overline{\Lambda}_p(a)=C_p,\qquad p=1,2,
  \end{equation}
  and
  \begin{equation}\label{E:static-moment-ratio-continuity}
    \lim_{a\downarrow0}
    \frac{\underline{\Lambda}_2(a)}{\overline{\Lambda}_1(a)}
    =2^{\kappa_{\rm st}}.
  \end{equation}
\end{corollary}

\begin{proof}
  The lower estimate is Proposition~\ref{P:static-minorant-lower}. For the upper
  estimate, apply Proposition~\ref{P:moment-growth-logconvex} with endpoints
  \(0,b\), using \(\overline{\Lambda}_p(0)=C_p\). This
  proves~\eqref{E:static-constant-envelope}. The finite ratio
  \(\overline{\Lambda}_p(b)/C_p\) is at least one. With \(b\) and the covariance
  fixed, its \(a/b\) power tends to one. Squeezing both rates
  proves~\eqref{E:static-constant-continuity}. The numerator and denominator
  in~\eqref{E:static-moment-ratio-continuity} tend to \(2^{\kappa_{\rm st}}L_0\)
  and \(L_0>0\), respectively. Division gives the asserted ratio. This proves
  Corollary~\ref{C:static-moment-constant-continuity}.
\end{proof}

\begin{remark}[Order of limits and comparison with the intrinsic constant]
  \label{R:temporal-type-scope}
  The static calibration now concerns the constant, not only the time exponent.
  It is an iterated statement: first take a time liminf or limsup at each fixed
  \(a\), with its own normalization \(T^{-\kappa(a)}\), and then send
  \(a\downarrow0\). No uniformity over covariances or \(T\), no simultaneous
  limit, and no existence of the time limit at a fixed positive \(a\) follows.
  Only the upper rate is asserted to be log-convex on the interior.

  Neither \(L_*(a)\) nor \(L^*(a)\) enters this argument. The constant \(L_0\)
  is the static first-mean rate. Parameter continuity alone does not identify
  the conjectured curve. The separate intrinsic second-moment value is proved in
  Theorem~\ref{T:intrinsic-second-moment}. No all-\(p\), absolute-first-moment,
  Skorohod, novelty or circulation result is claimed.
\end{remark}

\section{A small-order test for direct Brownian clock copying}
\label{S:two-clock-screen}

The result concerns a specified gap-transform substitution. It does not rule out
ordering corrections, a diagonal-only representation, one-sided/asymptotic use,
or the original sharp conjecture.

Theorem~6.1 of~\cite{chen:author2026:time-dependency}*{pp.~37--42} represents
the time Laplace transform of a \emph{fixed-order} mean by a spatial Brownian
path and an independent scalar Brownian clock. The proof uses nonnegative cut
sums \(c_l\), obtained from one-chain time order on p.~38, before the
wave-to-heat substitution on p.~39. Remark~6.3, p.~42, singles out that order as
the obstruction for more copies. We test the direct replacement by two spatial
Brownian paths and two scalar Brownian clocks, without an additional ordering
weight.

All transforms below are at fixed total order. We begin with the original
admissible case \(d=\alpha=1\), \(\Gamma=\delta_0\), \(0<a<1\), and
\(G_s(y)=\frac12\mathbf1_{\{|y|<s\}}\). This is a falsification case for an
intended general representation, not a restriction of Theorem~\ref{T:intro-second-moment}.

\LeanMarker{L:white-clock-laplace}
\begin{lemma}[White-space low-order Laplace transforms]
  \label{L:white-clock-laplace}
  The \((1,1)\) cross-pair gap transform is
  \begin{equation}\label{E:white-cross-clock-transform}
    \begin{split}
      K_\times(\theta_1,\theta_2)
      &:=\int_0^\infty\int_0^\infty
      e^{-\theta_1s-\theta_2t}|s-t|^{-a}
      \frac{\min(s,t)}2\ud s\ud t\\
      &=\frac{\Gamma_{\mathrm{E}}(1-a)}
      {2(\theta_1+\theta_2)^2}
      \big(\theta_1^{a-1}+\theta_2^{a-1}\big),
      \qquad \theta_1,\theta_2>0.
    \end{split}
  \end{equation}
  In particular, \(K_\times(\theta,\theta)=
  \Gamma_{\mathrm{E}}(1-a)/(4\theta^{3-a})\). For the empty-plus-length-two
  self-pair, independent gap damping gives
  \begin{equation}\label{E:white-self-clock-transform}
    K_{\rm self}(\eta_1,\eta_2)
    :=\int_0^\infty\int_0^\infty
    e^{-\eta_1r_1-\eta_2r_2}r_2^{-a}\frac{r_1}2
    \ud r_1\ud r_2
    =\frac{\Gamma_{\mathrm{E}}(1-a)}
    {2\eta_1^2\eta_2^{1-a}}.
  \end{equation}
\end{lemma}

\begin{proof}
  We compute the cross-pair transform by splitting its time domain into two
  cones, then evaluate the independent gap integrals for the self-pair.

  \medskip
  \noindent\textit{Step 1.~}We compute the cross-pair Laplace transform.
  The spatial cross factor is
  \[
    \int G_s(y)G_t(z)\delta_0(y-z)\ud y\ud z
    =\frac14\int\mathbf1_{\{|y|<s\wedge t\}}\ud y
    =\frac{\min(s,t)}2.
  \]
  Split its time integral into \(s=t+v\) and \(t=s+v\), with \(v>0\). Each part
  is one half the product of
  \[
    \int_0^\infty t e^{-(\theta_1+\theta_2)t}\ud t
    =(\theta_1+\theta_2)^{-2}
  \]
  and \(\int_0^\infty e^{-\theta_i v}v^{-a}\ud v
  =\Gamma_{\mathrm{E}}(1-a)\theta_i^{a-1}\). All integrands are nonnegative and
  finite since \(a<1\), so Tonelli justifies the separation. This proves
  \eqref{E:white-cross-clock-transform}.

  \medskip
  \noindent\textit{Step 2.~}We compute the self-pair gap transform.
  For one length-two chain, its delta covariance forces the two positions equal.
  The second Green factor is \(G_{r_2}(0)=1/2\), and integration of the first
  gives mass \(r_1\). The time difference is \(r_2\). The two elementary Laplace
  integrals therefore give \eqref{E:white-self-clock-transform}. This proves
  Lemma~\ref{L:white-clock-laplace}.
\end{proof}

The source's analytic time covariance is
\[
  \gamma_a^{\rm an}(z)
  =\big(z\,\operatorname{sgn}\Re z\big)^{-a},\qquad \Re z\ne0,
\]
with the power taken on the right half-plane. For \(r>0,v\ge0\), put
\begin{equation}\label{E:averaged-clock-kernel}
  \Psi_a(r,v)=\frac1{\Gamma_{\mathrm{E}}(a)}
  \int_0^\infty \lambda^{a-1}e^{-r\lambda-v\lambda^2/2}\ud\lambda.
\end{equation}
From the source's complex Laplace identity
\cite{chen:author2026:time-dependency}*{(6.3), p.~41}, and the Gaussian
characteristic function, \(\Psi_a(r,v)=\E(r+iZ)^{-a}\) when \(Z\) is centered
normal with variance \(v\). Fubini at fixed \(r>0\) is justified by
\(\int_0^\infty\lambda^{a-1}e^{-r\lambda}\ud\lambda<\infty\). In particular,
\(\Psi_a\) is nonnegative, decreases with \(v\), and \(\Psi_a(r,v)\le r^{-a}\).

Return to an arbitrary nonzero covariance in Hypothesis~\ref{H:homog}, with
temporal exponent \(a\), so \(0<a<1\), \(0<\alpha<2\) and \(a+\alpha<2\). Define
\[
  h_\Gamma(w)=\int e^{-w|\xi|^2/2}\mu(\ud\xi)
  =(p_w*\Gamma)(0),\qquad
  \mathcal{R}_\Gamma(q)=\int\frac{\mu(\ud\xi)}{q^2+|\xi|^2}.
\]
Both are positive and finite for \(w,q>0\). Spectral homogeneity gives
\(h_\Gamma(w)=w^{-\alpha/2}h_\Gamma(1)\) and
\(\mathcal{R}_\Gamma(q)=q^{\alpha-2}\mathcal{R}_\Gamma(1)\). The cross-pair gap
transform at common damping is
\[
  K_\Gamma(\theta)=N_2(\theta)^2
\]
with temporal exponent \(a\) in~\eqref{E:damped-norm-definitions}. This is the
spatial-regularization limit of the \((1,1)\) diagram.

Let \(B_1,B_2\) be independent spatial standard Brownian motions. For auxiliary
standard Brownian clocks with \(\E[\beta_1(u)\beta_2(v)]=\rho\min(u,v)\),
\(-1\le\rho\le1\), independent of the spatial paths, the clock difference has
variance
\[
  v_\rho(u,v)=u+v-2\rho\min(u,v).
\]
At common damping \(\theta>0\), the direct unweighted, unrestricted two-clock
covariance substitution yields the averaged candidate
\begin{equation}\label{E:naive-clock-candidate}
  J_{\Gamma,\rho}(\theta)
  :=\frac14\int_0^\infty\int_0^\infty
  e^{-\theta^2(u+v)/2}h_\Gamma(u+v)
  \Psi_a(\theta|u-v|,v_\rho(u,v))\ud u\ud v .
\end{equation}
Values on \(u=v\) do not affect the integral. The factor \(1/4\) is one factor
\(1/2\) for each wave-to-heat gap transform; \(h_\Gamma(u+v)\) is the spatial
average of \(\Gamma(B_1(u)-B_2(v))\). This average is interpreted by spatial
heat regularization. For the white covariance it is \((2\pi(u+v))^{-1/2}\). The
analytic covariance, not its absolute-value bound, is being copied.

\LeanMarker{P:clock-copy-obstruction}
\begin{proposition}[Obstruction to unweighted constant-correlation clock copying]
  \label{P:clock-copy-obstruction}
  For every covariance and parameters just specified, put
  \begin{equation}\label{E:cross-clock-lambda-weight}
    w_\theta(\lambda)
    =\frac{\lambda^{a-2}}
    {2\theta\Gamma_{\mathrm{E}}(a)}
    \big(\mathcal{R}_\Gamma(\theta)
    -\mathcal{R}_\Gamma(\theta+\lambda)\big),\qquad\lambda>0.
  \end{equation}
  This is a positive integrable density, and
  \begin{equation}\label{E:clock-copy-obstruction}
    \begin{split}
      K_\Gamma(\theta)&=\int_0^\infty w_\theta(\lambda)\ud\lambda,\\
      0<J_{\Gamma,\rho}(\theta)
      &\le J_{\Gamma,1}(\theta)
      =\int_0^\infty
      \frac{2\theta}{2\theta+\lambda}w_\theta(\lambda)\ud\lambda
      <K_\Gamma(\theta).
    \end{split}
  \end{equation}
  When \(\alpha=1\),
  \begin{equation}\label{E:degree-one-clock-loss}
    J_{\Gamma,1}(\theta)=(2-2^a)K_\Gamma(\theta).
  \end{equation}
  Thus, this direct analytic clock-copying family, including independent clocks
  and a single shared clock, fails to reproduce the cross-pair transform
  throughout the original nonzero covariance class.
\end{proposition}

\begin{proof}
  We first compute the exact original cross transform. We bound the clock
  candidates and evaluate the largest one, corresponding to a shared clock.
  Comparing the two positive densities gives strict loss, and degree-one
  homogeneity then gives its exact factor.

  \medskip
  \noindent\textit{Step 1.~}We compute the density of the original cross
  transform.
  First evaluate the cross transform. Temporarily smooth the spatial covariance,
  so Fourier and time Fubini are absolute. For \(r=|\xi|\), direct integration
  over \(s=t+v\) and \(t=s+v\) gives
  \[
    \int_0^\infty\int_0^\infty
    e^{-\theta(s+t)-\lambda|s-t|}
    \widehat{G}_s(\xi)\overline{\widehat{G}_t(\xi)}\ud s\ud t
    =\frac{2\theta+\lambda}
    {2\theta(\theta^2+r^2)((\theta+\lambda)^2+r^2)} .
  \]
  For example, on the first cone the \(v\)-integral uses
  \[
    \int_0^\infty e^{-(\theta+\lambda)v}
    \sin(r(t+v))\ud v
    =\frac{(\theta+\lambda)\sin(rt)+r\cos(rt)}
    {(\theta+\lambda)^2+r^2}.
  \]
  Integrating the remaining sine-square and sine-cosine terms against
  \(e^{-2\theta t}\), and adding the second cone, gives the displayed identity.
  The value at \(r=0\) is its continuous limit.

  \medskip
  \noindent\textit{Step 2.~}We remove spatial smoothing in the nonnegative
  spectral expression.
  The resulting spectral integrand is nonnegative. Removing the spatial Gaussian
  multiplier and using the Laplace formula for the power covariance therefore
  gives
  \[
    \begin{split}
      K_\Gamma(\theta)
      &=\frac1{\Gamma_{\mathrm{E}}(a)}
      \int_0^\infty\lambda^{a-1}
      \frac{2\theta+\lambda}{2\theta}
      \int\frac{\mu(\ud\xi)}
      {(\theta^2+|\xi|^2)((\theta+\lambda)^2+|\xi|^2)}
      \ud\lambda\\
      &=\int_0^\infty w_\theta(\lambda)\ud\lambda.
    \end{split}
  \]
  The second equality is partial fractions, with
  \((\theta+\lambda)^2-\theta^2=\lambda(2\theta+\lambda)\). The smoothing
  passage is justified by monotone convergence in this nonnegative spectral
  expression and the finite quadratic form in
  Lemma~\ref{L:damped-homogeneous-norms}. Alternatively, homogeneity directly
  bounds the density near zero by a constant times \(\lambda^{a-1}\), and at
  infinity by a constant times \(\lambda^{a-2}\). Its strict positivity follows
  from \(\mathcal{R}_\Gamma(\theta) -\mathcal{R}_\Gamma(\theta+\lambda)>0\),
  since \(\mu\ne0\).

  \medskip
  \noindent\textit{Step 3.~}We prove integrability and compare the clock
  correlations.
  The candidate integrals are well defined and positive. The bound
  \(\Psi_a(\theta|u-v|,v_\rho)\le (\theta|u-v|)^{-a}\), together with
  \(h_\Gamma(w)=h_\Gamma(1)w^{-\alpha/2}\), supplies an absolute majorant. With
  \(w=u+v\), \(q=u-v\), the Jacobian is \(1/2\) and \(-w<q<w\). It is a constant
  times
  \[
    \int_0^\infty e^{-\theta^2w/2}w^{-\alpha/2}
    \int_{-w}^w|q|^{-a}\ud q\ud w
    =\frac2{1-a}\int_0^\infty
    e^{-\theta^2w/2}w^{1-\alpha/2-a}\ud w<\infty .
  \]
  Here \(a<1\) and \(a+\alpha/2<2\). The same bound dominates the spatial
  heat-regularized averages, whose spatial factor is \(h_\Gamma(w+2\epsilon)\le
  h_\Gamma(w)\). This is an averaged spectral bound, not pointwise covariance
  monotonicity. Since \(v_\rho(u,v)\ge|u-v|=v_1(u,v)\), monotonicity of
  \(\Psi_a\) gives \(J_{\Gamma,\rho}\le J_{\Gamma,1}\).

  \medskip
  \noindent\textit{Step 4.~}We evaluate the shared-clock candidate and compare
  its density.
  We compute the latter by nonnegative Tonelli. Write
  \(A=\theta\lambda+\lambda^2/2\). Integration in \(q\) gives
  \(2(1-e^{-Aw})/A\). Also
  \[
    \mathcal{R}_\Gamma(q)
    =\frac12\int_0^\infty e^{-q^2w/2}h_\Gamma(w)\ud w,
  \]
  by the elementary resolvent integral and spectral Tonelli. Consequently,
  \[
    \begin{split}
      J_{\Gamma,1}(\theta)
      &=\frac1{4\Gamma_{\mathrm{E}}(a)}
      \int_0^\infty\frac{\lambda^{a-1}}A
      \int_0^\infty h_\Gamma(w)
      \left(e^{-\theta^2w/2}
      -e^{-(\theta+\lambda)^2w/2}\right)\ud w\ud\lambda\\
      &=\frac1{\Gamma_{\mathrm{E}}(a)}
      \int_0^\infty
      \frac{\lambda^{a-2}
      (\mathcal{R}_\Gamma(\theta)-\mathcal{R}_\Gamma(\theta+\lambda))}
      {2\theta+\lambda}\ud\lambda\\
      &=\int_0^\infty
      \frac{2\theta}{2\theta+\lambda}w_\theta(\lambda)\ud\lambda.
    \end{split}
  \]
  The \(q\)-integral's factor two cancels its coordinate Jacobian \(1/2\); the
  resolvent identity supplies the other factor two. Since \(w_\theta>0\) and
  \(2\theta/(2\theta+\lambda)<1\) for every \(\lambda>0\), the strict inequality
  follows.

  \medskip
  \noindent\textit{Step 5.~}We compute the exact loss at spatial degree one.
  Finally, if \(\alpha=1\), homogeneity gives
  \(\mathcal{R}_\Gamma(q)=\mathcal{R}_\Gamma(1)/q\). Substitution in the two
  densities gives
  \[
    \begin{split}
      K_\Gamma(\theta)
      &=\frac{\mathcal{R}_\Gamma(1)}
      {2\theta^2\Gamma_{\mathrm{E}}(a)}
      \int_0^\infty\frac{\lambda^{a-1}}{\theta+\lambda}\ud\lambda,\\
      J_{\Gamma,1}(\theta)
      &=\frac{\mathcal{R}_\Gamma(1)}
      {\theta^2\Gamma_{\mathrm{E}}(a)}
      \int_0^\infty\lambda^{a-1}
      \left(\frac1{\theta+\lambda}
      -\frac1{2\theta+\lambda}\right)\ud\lambda .
    \end{split}
  \]
  The real Beta identity
  \[
    \int_0^\infty\frac{x^{a-1}}{x+c}\ud x
    =c^{a-1}\Gamma_{\mathrm{E}}(a)\Gamma_{\mathrm{E}}(1-a),
    \qquad c>0,
  \]
  follows by scaling and \((r,s)\mapsto(r+s,r/s)\) in the positive product of
  the two Gamma integrals. Their finiteness uses \(0<a<1\). Evaluation of these
  two integrals gives \(J_{\Gamma,1}/K_\Gamma=2(1-2^{a-1})=2-2^a\), proving
  \eqref{E:degree-one-clock-loss}. This proves
  Proposition~\ref{P:clock-copy-obstruction}.
\end{proof}

\begin{remark}[Normalization and scope of the obstruction]
  \label{R:clock-copy-scope}
  At the static endpoint define the time covariance directly as one. The
  candidate with \(\Psi_0=1\) gives
  \[
    J_{\Gamma,\rho}(\theta)
    =\int\frac{\mu(\ud\xi)}{(\theta^2+|\xi|^2)^2}
    =K_\Gamma(\theta).
  \]
  This follows by integrating \(w e^{-\theta^2w/2}h_\Gamma(w)/4\) spectrally. In
  the white case it equals \(1/(4\theta^3)\). Thus, static and
  homogeneous-degree checks alone miss the positive-\(a\) loss. Relative to the
  shared-clock value \(\rho=1\), the omitted cross-pair multiplier is exactly
  \(1+\lambda/(2\theta)\). For smaller correlations the deficit can be larger.
  This does not assert a correction rule at higher orders.

  The transforms \(K_\times\) use gap damping. For the auxiliary cutoff function
  \[
    F(T_1,T_2)=\int_0^{T_1}\int_0^{T_2}
    |s-t|^{-a}\frac{\min(s,t)}2\ud s\ud t,
  \]
  nonnegative Fubini adds \(1/(\theta_1\theta_2)\) to its double
  external-horizon transform. This function is not an unequal-time solution
  correlation. Reversing chains ending at different physical times instead
  shifts the cross covariance to \(|T_1-T_2-s+t|^{-a}\). Setting the two damping
  parameters equal is also not the single Laplace transform of a
  diagonal-horizon moment. For a self-pair at common gap damping, one external
  horizon instead adds \(1/\theta\), recovering the known first-mean
  coefficient.

  The proposition rules out exactly the analytic covariance substitution
  in~\eqref{E:naive-clock-candidate}, with product heat-chain measure,
  unrestricted operational-time rectangle and constant Brownian correlation.
  Additional relative-order weights,other derived clock couplings,matrix
  representations,modulus inequalities and their possible asymptotics are not
  covered. In particular, the inequality at this first cross pair leaves open a
  shared-clock lower comparison or asymptotic equivalence at higher orders. This
  finite clock calculation does not identify the intrinsic value proved  in
  Theorem~\ref{T:intrinsic-second-moment} with the source's all-moment
  conjecture. This is not a Skorohod result or a no-go for every two-clock
  method.
\end{remark}

No-go~\ref{P:nogo-17} summarizes the missing degree-two multiplier.
No-go~\ref{P:nogo-18} records a different degree-four obstruction: erasing
an additional shift makes one upper envelope too large to close its
comparison, without refuting the comparison itself.

\chapter{The AI work in numbers}\label{S:appendix-ai-statistics}

% Appendix: the AI work in numbers. The chapter heading is in the main file.
% Every number comes from notes/ai-statistics-data.tex, generated by
% scripts/ai-statistics/make_appendix.py.
\definecolor{aiClaude}{HTML}{4477AA}
\definecolor{aiDeepSeek}{HTML}{EE6677}
\definecolor{aiKimi}{HTML}{CCBB44}
\definecolor{aiGPT}{HTML}{228833}
\definecolor{aiCommits}{HTML}{555555}
\definecolor{aiMilestone}{HTML}{AA3377}

This appendix records how much of the work was done, by which kind of model,
and how often the checks failed. It counts what the project's own records
show. Every number is produced by a short script from files kept in the
repository, and the scripts are kept with it, so a reader can rebuild each
figure and each table. The record covers \aiDays{} days, from \aiFirstDay{}
to \aiLastDay{} (dates in UTC), with \aiCommits{} recorded changes (commits).
The retained reviews-folder census ends on \aiReviewLastDay{}. For \aiLaterFirstDay{}
to \aiLaterLastDay{} the figure also shows a separate count of recovered
review attempts, and the last section lists the reviews of the publication
closeout.

\section*{Who did the work}

The mathematics, the \LaTeX{} source and the bookkeeping were drafted by a
GPT-family model working through Codex. Its work was checked by reviewers from
other model families: Claude (Opus~5.5) at first, and from November~1 DeepSeek
V4 Pro as the main reviewer, with Kimi K3. A GPT model also reviewed the Lean
statements. Appendix~\ref{S:appendix-lean-checks} lists what the Lean companion
checks. A review counts only when the reviewer comes from a different family
from the author of the work it checks.

\begin{center}
\small
\begin{tabular}{@{}lrrr@{}}
  \toprule
  & Sessions or records & Input tokens & Output tokens \\
  \midrule
  Worker: GPT (Codex) & \aiCodexSessions{} sessions & \aiCodexInB{} billion
  & \aiCodexOutM{} million \\
  Claude Code logs & \aiClaudeCalls{} usage records
  & \aiClaudeCacheReadB{} billion (cache) & \aiClaudeOutM{} million \\
  DeepSeek V4 Pro & \aiDeepSeekSessions{} sessions & \aiDeepSeekInM{} +
  \aiDeepSeekCacheM{} million & \aiDeepSeekOutM{} million \\
  \bottomrule
\end{tabular}
\end{center}

Most of the worker's input, \aiCodexCachedPct\%, was text it had already read
and that was served again from a cache; \aiCodexReasonM{} million of its
output tokens were internal reasoning. Counting cached input on both sides,
the worker read about \aiWorkerToDeepSeekRatio{} times as much as the
DeepSeek reviews.
The Claude Code figures include Opus~5.5, Fable~5.2 and synthetic usage
records, covering editing and coordination as well as reviews. They count
logged assistant usage records, not verified API calls. Only the DeepSeek
figures carry a cost figure, a cost counter of \aiDeepSeekCost{} computed by
OpenCode, which is not an invoice and whose currency was not verified; the
other models were used on subscriptions, whose costs are not recorded in these
logs. The worker
machine's logs hold no token counts for Kimi, so these totals leave Kimi out;
the closeout records at the end of this appendix carry Kimi's own reported
counts.
These totals are a snapshot taken at \aiUsageSnapshot{} from those logs.
The DeepSeek line includes the recorded usage of the re-read calls of
\aiRereadDay{}, described at the end of this appendix, which ran on the same
machine; the token figure given there is part of these totals and must not be
added to them. The totals cover the sessions and usage records that one machine retained, not all
activity. Only recorded usage fields are summed; usage with no record,
including the abandoned attempt described in the closeout section, is not
included.

\section*{Reviews}

At that historical census, the repository held \aiReviewFolders{} review
folders. Of these,
\aiFoldersWithRun{} contain at least one review run; the others hold
preparation, decisions or diagnostics. A \emph{run} is one review task sent to
one reviewer model, with or without a usable report. There were \aiRuns{}
runs, and \aiFollowups{} further follow-up requests within the same reviewer
sessions.

\begin{center}
\small
\setlength{\tabcolsep}{5pt}
\begin{tabular}{@{}lrrrrrrr@{}}
  \toprule
  Reviewer & Runs & Completed & Not completed & Met bar & FAIL & MEDIUM
  & Hours \\
  \midrule
  \aiReviewTableRows
  \bottomrule
\end{tabular}

\smallskip
\begin{minipage}{0.92\textwidth}
  \footnotesize Hours are summed over the runs whose time was recorded:
  \aiHoursCoverage.
\end{minipage}
\end{center}

\emph{Met the bar} means a completed run graded PASS or PASS WITH NOTES at HIGH
confidence, which is what the project required before a result could be
accepted. Meeting the bar is not the same as acceptance: the final decision was
recorded separately, and only the Claude and GPT reports record it in a form
that can be counted. Of the \aiRunsFail{} FAIL grades, \aiFailGPT{} came from
the GPT reviews of the Lean companion; they reported defects in the Lean
statements or in their match with this book, which then had to be repaired. Runs
that did not complete stopped on usage limits, timeouts or reports in the wrong
format. Figure~\ref{Fig:ai-timeline} shows the work day by day and marks
\aiMilestones{} milestones, from setting up the project to closing the Lean
correspondence gates; their commits are listed in
\path{scripts/ai-statistics/milestones.json}.

The figure covers only this book, which began on \aiInfraBookStart{}. The
system that wrote it is older. The bibliography and the paper library have been
kept since \aiInfraBibStart{}. The collection of \LaTeX\ sources began in
\aiInfraSourcesStart{}, the wiki in its present form in \aiInfraWikiStart{}, and
the library of skills \aiInfraSkillsWhen{}. The book took a month because that
work was already done. The dates and the commits behind them are in
\path{scripts/ai-statistics/infrastructure.json}.

\section*{Approaches that did not work}

As of \aiNogoSnapshot{}, the record holds \aiNogoTotal{} failed approaches
and corrections: \aiNogoMath{} mathematical (\aiNogoRoutes{} approaches that
cannot work and \aiNogoCounter{} counterexamples or corrected claims),
\aiNogoSource{} defects in the source paper, \aiNogoMethod{} failures of a
method of checking, and \aiNogoProcess{} process record.
Chapter~\ref{S:book-no-go} has \aiNogoMath{} numbered mathematical entries
with written arguments at their printed scopes. NG-33 is a broader historical
group of refuted reviewer objections; its printed counterpart proves three
representative cases. The remaining items in that group are not claimed
proved there. In the first snapshot of the record, on November~10,
there were \aiNogoTotalFirst{} entries, and \aiNogoInBookFirst{} of the
mathematical ones were marked as proved in the book, \aiNogoReviewedFirst{} as
reviewed and adjudicated, and \aiNogoRecordedFirst{} as recorded only.
In the current snapshot, \emph{proved in book} denotes a written proof at
the record's stated scope, not scientific acceptance. NG-33 remains
\emph{reviewed and adjudicated} because only part of its historical group is
printed. The earlier review of NG-09 did not meet its acceptance threshold.
None of these records refutes the main theorems.

\begin{figure}[t]
  \centering
  \aiTimelineFigure

  \medskip
  {\scriptsize
  \begin{tabular}{@{}l@{}}
    \aiMilestoneKey
  \end{tabular}\par}
  \caption[The work day by day, with milestones]{Review runs per day by
    reviewer (bars, left scale) and commits per day (line, right scale), with
    the project's \aiMilestones{} milestones (numbered circles). Days are in
    UTC, and each circle sits on the UTC day of the milestone's commit; the
    key gives the date as recorded. Claude returned after November~1 for
    checks of rendered pages. The coloured bars stop on \aiReviewLastDay{},
    where the census of review runs ends; the grey bars count recovered
    review attempts, which are not part of the census. The figure starts with
    this book; the system behind it was built earlier, as the text says.}
  \label{Fig:ai-timeline}
\end{figure}

\section*{What the numbers do not show}

Counts of runs and tokens measure effort, not correctness: a thousand commits do
not prove a theorem, and a passed review is a check, not a proof. Some figures
are incomplete. The worker machine's logs hold no Kimi tokens, and only
\aiRunsKimiTimed{} of the \aiRunsKimi{} Kimi runs have a recorded time; some
DeepSeek runs recorded only their last request, so DeepSeek hours are too low;
and the dates of \aiRunsDatedByFolder{} runs are taken from folder names. The
token totals come from the logs of one machine. The \aiCommits{} commits are
counted in the working repository; the public history has fewer, because commits
that changed only private files are left out of it. A short list of review
folders could not be classified and is kept with the scripts.

The grey bars count \aiLaterAttempts{} review attempts recovered from
retained records for \aiLaterFirstDay{} to \aiLaterLastDay{}. They are
attempts, not completed reviews: some were blocked or returned an unusable
report. Their UTC days use recorded launch, client-session and panel-dispatch
evidence; some are bounded or decoded rather than exact invocation clocks.
Actual model identity is not verified for every attempt, so this series is not
split by reviewer. We show them separately from the \aiRuns{}-run census;
together they still do not establish a complete project-wide inventory. These
are counts of review tasks, not provider requests, and activity with no
retained record is missing. Of these attempts, \aiLaterInLedger{} are
represented in the closeout ledger below.

\section*{The publication closeout reviews of \aiCloseoutMonth}

This section covers only the closeout reviews before publication, in
\aiCloseoutMonth: \aiCloseoutRecords{} records in the experiment's run ledger,
whose snapshot is dated \aiCloseoutSnapshot. The figures come from the dated
correction of that ledger recorded on \aiCloseoutCorrection{}
(\aiCloseoutSource, commit \texttt{\aiCloseoutSourceCommit}, blob
\texttt{\aiCloseoutSourceBlob}). Earlier corrections added verified model
identities and accounting fields and separated one raw report's hash from
that of its parsed copy (C003). The current snapshot preserves all
\aiCloseoutHistoricalRecords{} older
rows unchanged and adds \aiCloseoutRecentRecords{} later cold-review records,
including the author's
bounded decision on the uniform-product review. Reported grades and token
objects remain unchanged. The experiment itself is the way of
doing mathematics described in the Preface, and its primary record is the whole
repository: its history, proofs, failed approaches, human and AI work, sources,
reviews and the Lean companion. The run ledger instruments only this closing
round, and only partly; its own coverage note says what it includes. Its
records are not added to the totals above, and its token counts, reported for
some records only, are not totals for the project. The table counts the outcome
recorded for each record, by the model that actually ran\aiCloseoutRunningNote.
No record accepts the whole book: each acceptance is limited to the scope its
review stated. The preserved older rows retain their historical status: a
row marked as running describes its snapshot date, not a live process now.

\aiCloseoutRecentWords{} later calls add scoped observations in this
November~15 update. A whole-source call returned inconclusive at medium
confidence;
a focused call returned pass with notes at high confidence but did not
complete the required core reconstruction. \aiCloseoutLocalWords{} local
calls returned \aiCloseoutLocalPassWords{} passes and
\aiCloseoutLocalNotesWords{} pass with notes, all at high confidence.
\aiCloseoutLocalAcceptedWords{} were accepted as scoped mathematical evidence
with notes and \aiCloseoutLocalPartialWords{} were partial. For the
uniform-product check, the author accepted a bounded
exception for one automatic native retry, with the incident disclosed. The
original JSON-format failure in one report remains recorded; a
copy in which only the JSON punctuation was corrected is not another review.
These calls do not replace a complete independent review of the book. The
combined critical review of Chapter~12 was not obtained: one attempt
exhausted its context and the next stopped at a safeguard with invalid partial
output. It was never marked passed. A separate high-confidence review of
whether Chapter~12 faithfully states its evidence and limitations does not
replace that missing assessment.
Author-side source comparisons and delivery checks are not additional
independent scientific votes.

The \aiCloseoutRecords{}-record ledger remains a dated snapshot. A separate
evidence addendum links later observations to its original identifiers,
including its missing launch clocks and C011's terminal report. The recovered
inventory contains \aiEvidenceAttempts{} documented review attempts:
\aiEvidenceBaseAttempts{} already represented in this ledger and
\aiEvidenceExtraAttempts{} additional records, with \aiEvidenceFollowups{}
same-session follow-ups listed separately. It is a timeline inventory, not a
new completion, pass-rate or acceptance table; it does not add tokens or costs
to the snapshots above. Its public metadata are in
\aiEvidenceInventoryPath.

\begin{center}
\small
\setlength{\tabcolsep}{5pt}
\aiCloseoutTable
\end{center}

Token counts are given for the \aiCloseoutDenominator{} original review
records, per model, reporting format and accounting scope, in the form in which
the review tool reported them. Only the \aiCloseoutUsageRecords{} records with
reported usage are summed; a record without figures is shown as having no
reported value, and nothing is estimated. Records separately marked as derived
recoveries, parser refusals, preflights or audit addenda do not represent new
provider calls in this ledger. Format continuations count toward usage only
when nonoverlap is verified.\aiCloseoutHeaderNote{} The formats, as recorded
in the ledger:
\begin{itemize}\small
  \setlength{\itemsep}{0pt}
  \aiCloseoutFormatKey
\end{itemize}

Accounting scopes are kept separate:
\begin{itemize}\small
  \setlength{\itemsep}{0pt}
  \aiCloseoutAccountingKey
\end{itemize}
These are sums of known reported components, not actual all-activity totals.
A review record counts a review task, not a provider-stream total. An automatic
retry within one later review record opened two streams of generated output;
it is separate from the recorded repairs to report formats. Usage and cost
for the earlier
abandoned stream are unavailable. The reported final-step total is not the
whole invocation. Actual whole-run usage and cost remain unknown, as does
whether the whole run stayed within its approved token limit.

\begin{center}
\small
\begin{tabular}{@{}>{\raggedright\arraybackslash}p{0.27\textwidth}>{\raggedright\arraybackslash}p{0.5\textwidth}r@{}}
  \toprule
  Model (form) & Reported totals & Reported on \\
  \midrule
  \aiCloseoutTokenRows
  \bottomrule
\end{tabular}
\end{center}

Known reported cost components of the closeout records (not actual
all-activity cost): \aiCloseoutCost.

\subsection*{Later evidence and the author's decision}

Missing cells in the original ledger mean that a value was not verified in
that snapshot, not that no record could ever be recovered. Original launch
records or CLI receipts now supply \aiEvidenceRecordedStarts{} of its
previously missing start clocks:

\begin{center}
\small
\begin{tabular}{@{}lll@{}}
  \toprule
  Record & Recorded clock (UTC) & Kind \\
  \midrule
  \aiEvidenceStartRows
  \bottomrule
\end{tabular}
\end{center}

These are client-runner or CLI receipt clocks, not independently measured
provider API boundaries. They were read from retained records, not calculated
by subtracting a duration from a later clock. The original ledger's nulls
remain unchanged.

For \aiEvidenceClientClocks{} other panel records, the retained clocks record
an OpenCode client session's creation and last update. These do not establish
the review invocation's exact start or finish. For the remaining panel
record, only its UTC day is bounded by the recorded panel creation and the
next role's receipt under the documented sequential workflow. Its exact
start, finish and actual model remain unknown. The evidence addendum keeps
recorded clocks, client clocks, decoded clocks and bounds distinct.

C011 was still running at the original snapshot. Its later receipt records
an invalid report; a punctuation-only recovery returned pass with notes at
medium confidence. The review is terminal and remains unaccepted under the
high-confidence standard. This is later evidence about the same invocation,
not another review or an upgraded grade.

On \aiEvidenceDecisionDate{}, the author chose the policy of an
\emph{honest-standing release with the grades stated plainly}. This decision
changes no grade, supplies no missing review and does not turn scoped
evidence into whole-book scientific acceptance. It selects how the record is
disclosed; it is not a public-release act or the final release accept file.
The dated observation and decision metadata are in \aiEvidenceUpdatesPath.

\section*{The re-reads after the readability pass}

After the reviews counted above, the proofs were rewritten so that their steps
are easier to follow. Statements and existing displayed formulas were kept, apart from the
separately approved corrections described below. Inside the
proofs, sentences were added, removed or reordered. Of the \aiRereadProofs{}
proofs in the \aiRereadFiles{} technical files, \aiRereadChanged{} changed and
\aiRereadUnchanged{} stayed identical byte for byte. In
\aiRereadStatementsChanged{} statements the wording was corrected under
separate approval. The
\aiRereadReaderEditionProofs{} proofs of the reader edition are not counted
here. The earlier reviews stay bound to the text they received, and they do not
cover the sentences added later.

Every one of the \aiRereadProofs{} proofs, changed or not, was then sent
again, the old text beside the new, to DeepSeek V4 Pro. The editor was a GPT
model, so the reader came from a different family. Each call was asked to read
both texts in full, to rebuild the deductions, and to say whether a step or a
qualification had been lost. All calls were made on \aiRereadDay{}.
Table~\ref{Tab:ai-rereads} gives the counts. This series is separate from the
\aiRuns{} review runs and the \aiLaterAttempts{} recovered attempts above, and
its numbers must not be added to theirs. It changes no earlier grade, and it is not an acceptance of the book. These
scoped re-reads do not upgrade the inconclusive whole-book cold read or
replace the missing skeptic synthesis of Chapter~\ref{S:book-no-go}.

\begin{table}[!ht]
  \centering
  \caption[Re-reads of the proofs after the readability pass]{Re-reads of the
    proofs after the readability pass. The last three rows are our own check
    of what the reports wrote, not grades given by the reader.}
  \label{Tab:ai-rereads}
  \begin{tabular}{@{}lr@{}}
    \toprule
    Calls made                                              & \aiRereadCalls \\
    \quad with a report                                     & \aiRereadReports \\
    \quad without a report (the output limit was reached)   & \aiRereadNoReport \\
    Reports graded PASS                                     & \aiRereadPass \\
    Reports graded PASS WITH NOTES                          & \aiRereadPassNotes \\
    Proofs sent, counting every repeat                      & \aiRereadPairsSent \\
    Scopes prepared but never sent                          & \aiRereadNeverSent \\
    \midrule
    Reports that cover the final text                       & \aiRereadCurrent \\
    \quad graded PASS                                       & \aiRereadCurrentPass \\
    \quad graded PASS WITH NOTES                            & \aiRereadCurrentPassNotes \\
    Earlier calls, on text that changed later or was sent again
    & \aiRereadEarlier \\
    Proofs of the final text covered, each by one report    & \aiRereadProofs \\
    \quad old and new text both reported as read in full    & \aiRereadComplete \\
    \midrule
    Report's reasoning checked for the proof                & \aiRereadChecked \\
    Checked, but resting on a parent result with a stated limit
    & \aiRereadConditional \\
    Report's reasoning partial for the proof                & \aiRereadPartial \\
    \bottomrule
  \end{tabular}
\end{table}

A grade belongs to a call, not to a proof. All \aiRereadReports{} reports
stated high confidence for their call, and none gave a confidence for a single
proof. The \aiRereadPairsSent{} proofs sent include repeats: some files were sent at
more than one frozen source revision, and some proofs appeared in more than
one assigned scope, so this is not a number of distinct proofs. The call without a report has no grade; it is not counted as a pass or
as an inconclusive review, and its proofs were read in later calls.

The last three rows of the table are the authors' side of the work. The
number \aiRereadPartial{} there counts proofs and the number
\aiRereadPassNotes{} above counts reports; that they are equal is a
coincidence. For each
proof we checked the reasoning that the report wrote down. \emph{Checked}
means that the deductions in the report hold as written for that proof. It is
a local check; it does not prove that everything the proof depends on is
sound. For \aiRereadConditional{} proofs the check rests on a parent result
whose own check has a stated limit. These fall into
\aiRereadConditionalGroups{} groups and repeat their parents' limits; they
are not further gaps. \emph{Partial} means that a deduction the report was
asked to give is missing or does not hold as written in the report. It does
not by itself show an error in the book. Our checks of the book's text found
no new mathematical defect in these proofs. Before the final calls, four added sentences that stated more than the proof
gives or misdescribed a calculation were corrected (the rank-one witness,
``canceling half'', the strictly positive profile, and the complex-vector
components), and sentence boundaries were repaired. For one of the partial proofs the central argument was
checked and only further derivations that the call had been asked for are
missing. In one report the requested word-for-word comparison of a corrected
statement was wrong, although its mathematics was checked; that comparison
was not obtained. This is a separate limit and is not one of the partial
proofs.

Some usage is unknown. In \aiRereadRestarted{} calls the client program
restarted an interrupted answer by itself, and the logs kept the usage of the
final answer only. Of these, \aiRereadRestartedHeld{} happened before the
author prospectively allowed automatic transport recovery under the later
policy; they keep the status \emph{held}, and
none of them covers the final text. The other one was allowed in advance, and
its report covers \aiRereadRestartedCurrentProofs{} proofs of the final text.
The available finish records report \aiRereadTokens{} tokens, of which
\aiRereadTokensCached{} were cached input. The client's own cost counter for
these records sums to \aiRereadCostCounter{} in the client's units; it is not an
invoice, and its currency was not verified. The usage of the
interrupted answers, the whole cost and the number of requests that reached
the provider are unknown, and no estimate is given. No whole-invocation token-limit
compliance is claimed.

Three kinds of help in this pass were not reviews and are counted nowhere
above. The editor, a GPT model, rewrote the proofs under written rules, and
each file had to pass mechanical checks that existing statements, displayed
formulas, labels and references were kept, with separately approved
original-text changes recorded as exceptions. The
\aiRereadStatementsChanged{} statement corrections, the other repairs of the
original text and the removal of unused setup in the proof of
Proposition~\ref{P:symplectic-block-mean} were separate edits of the original
text, made at the author's decision and not by the editor. In \aiRereadReaderRuns{} read-only runs,
Claude Fable~5.2 listed the places where a reader has to stop and the
passages that repeat; \aiRereadStopPoints{} such places were listed and given to the editor to
fill. These runs were editorial aids and were not counted as scientific
reviews. GPT helpers assisted in checking the reports'
reasoning, which is an author-side check and not an independent vote.

The data are in \path{reviews/2026-11-16-readable-final-source/}. They are
bound to the text of the proofs at development commit
\texttt{\aiRereadSource} of the private development repository. The
editorial record is
\path{scripts/ai-statistics/reviews/readability-pass-2026-11-16.json}.

\section*{The review of the Lean companion's references}

Before release v0.10.0 of the Lean companion, its references to lines of
this book were moved to the book's final text, and that work was reviewed in
\aiLeanReCalls{} calls to \aiLeanReModel{} at high effort, read-only, on
\aiLeanReDay{}. The companion's text was written by a Claude model, so the
reviewer came from a different family. The first call returned FAIL at high
confidence, with \aiLeanReFirstMajor{} major and \aiLeanReFirstMinor{} minor
findings about the companion's descriptions and quotation ranges. After the
corrections, a narrower confirmation returned FAIL at high confidence again,
with \aiLeanReSecondMajor{} major findings, because two corrections had not
reached every description. A second confirmation returned PASS WITH NOTES at
high confidence, with no fatal or major finding and \aiLeanReThirdMinor{}
minor note, and that report was accepted. The two earlier reports are kept.
The first call checked a sample of the references with broad mechanical checks;
the two confirmations checked only the corrections. Together the calls took
about \aiLeanReMinutes{} minutes. They used \aiLeanReInput{} input tokens, of
which \aiLeanReCached{} were cached, and \aiLeanReOutput{} output tokens, of
which \aiLeanReReasoning{} were reasoning; these counts include each other as
stated and must not be added again. No cost was recorded, which is not the same
as zero. This series is separate from every count above. It checked the
references and the repairs, not the mathematics of this book or the analytic
inputs outside Lean, and it is not a human check. The author accepted release
v0.10.0 on November~16, 2026. The record is
\path{scripts/ai-statistics/reviews/lean-reanchor-review-2026-11-16.json}.

\chapter{What Lean checks}\label{S:appendix-lean-checks}

% Appendix F: what Lean checks. The chapter heading is in the main file.
% The table and the numbers come from lean-coverage-data.tex, generated by
% scripts/lean-coverage/make_appendix_f.py from the Lean companion's
% coverage file and receipts (pinned by scripts/lean-coverage/pin_companion.py).

\hypertarget{lean-entry-top}{}%
Some steps of this book were also checked in Lean, a proof assistant. They
live in a separate Lean package, the Lean companion. This appendix describes its release
v0.10.0 (DOI~\href{https://doi.org/10.5281/zenodo.24295453}{10.5281/zenodo.24295453}),
which contains the formalization waves through L14 and refers to the lines
of this book's final text. The waves were checked at commit \texttt{8f51004}
of the private development repository; the release is cut from that commit,
its receipts are those of the release commit \texttt{\lcCommit}, and the
companion's \texttt{provenance/DEVELOPMENT.md} lists every change between
the two.

For each listed result, the entry identifies what Lean proves and what
remains outside that proof. Read the result reference and its title as one
entry. Where the conclusion is narrower, the entry states which part is
proved and which part is not yet formalized. Where a checked statement assumes an analytic fact, the entry names it as
a supplied input. Where the checked steps take no additional analytic
input, the entry states their assumptions and describes work left outside
Lean. A supplied input is a hypothesis of the Lean proof, not an additional
conclusion it establishes.

\section{How to read the entries}\label{S:lean-how-to-read}

Every result with a row here, except those listed as not formalized, carries
a small link in the left margin of the text, beside its heading or its
displayed equation: the word Lean and, under it, the class that this
appendix prints for the result. The link leads to the result's row.
\emph{Full} means that every conclusion is formalized under the book's
hypotheses.
\emph{Narrower} means that Lean works under the book's hypotheses but
proves only the conclusion specified in the entry; the remaining part is
not yet formalized.
\emph{Supplied} means that every conclusion is proved with analytic inputs
supplied as hypotheses.
\emph{Partial} means that Lean checks only the stated steps, parts or
narrower form, using the listed supplied inputs, if any; it does not prove
the whole result as stated in the book. The entry states which hypotheses
the checked steps take and which analytic work is left outside Lean.
Lean checks the formal statements. Their match with the statements of this
book was reviewed by AI models and has not been verified by a person. A
result without a link carries no Lean claim; this does not mean that the
result is in doubt.

\ifdefined\lcModelKey
Some entries name a registered modelling choice by a code. Lean proves the
statements for the model; the identification named in each line is not
proved in Lean.
\begin{description}\setlength{\itemsep}{0pt}
  \lcModelKey
\end{description}
\fi

Three points apply to every entry. Lean checks the statements as written
in Lean. Their correspondence with this book was reviewed by AI models
from a different family from the author of the formal statements; it has
not been verified by a person. A supplied input is a named mathematical
hypothesis: its proof and its application to the book's objects are not
certified merely by using it in Lean. The entries distinguish these inputs
from the conclusions Lean proves. They do not certify the book's main
stochastic moment theorems end to end.

\section{How the match with the book was reviewed}\label{S:lean-match-review}

The statements added in wave L9, and the changes that wave made to existing
entries of the coverage record, were checked against the book in eight
independent correspondence reviews, done by AI models from a different family
than the one that wrote the formal statements. The reviews found documentation
corrections only, and no error in a Lean proof. The corrections were made and
re-checked over ten confirmation rounds. During them the author fixed three
conventions for the record: a result counts as formalized in full only if every
conclusion is a Lean conclusion under the book's hypotheses, and otherwise the
missing part is named; a result about an object the companion models by a
registered standard choice counts in full, while a result that needs an object
the companion does not construct counts as narrower; and the record's index
lists, for each label, only the sections where it is treated. A final
confirmation passed with high confidence, and the author accepted wave L9 on
November~12, 2026, with a correction of the count of public definitions
recorded as a documentation-only successor.

The statements added in waves L10 to L14, and the changes those waves made to
existing entries of the coverage record, were checked against the manuscript in
correspondence reviews. The reviews were done by AI models from a different
family than the one that wrote the formal statements, under the conventions the
author fixed in wave L9. Each wave's reviews found documentation corrections
only, and no error in a Lean proof. The corrections were made, and a
confirmation review passed with high confidence.

The author accepted waves L10 and L11 on November~13, 2026:
\begin{itemize}
  \item In wave L10, one assembly step of
    Proposition~\ref{P:interval-coefficients} assumed a bound on the pairings of
    the cut-off coefficients that the manuscript does not state as a hypothesis.
    Wave L11 proves that bound from a supplied expansion of the coefficients.
  \item In wave L11, the Gaussian space is represented by Hermite coefficient
    vectors. Lemmas~\ref{L:fock-form} and~\ref{L:interval-contractions} hold
    with the manuscript's hypotheses but narrower conclusions, on a
    Fourier-side model not identified with the physical-space operators.
\end{itemize}
The author accepted wave L12 the same day, recording that the smoothed
covariance of Lemma~\ref{L:damped-homogeneous-norms} is modelled by its density;
with that choice the lemma is formalized in full.

Waves L13 and L14 were accepted together, with three further modelling choices
recorded:
\begin{itemize}
  \item the second moment's two-chain pairing integrals are represented by a
    deterministic diagram model, with the mixed-moment expansion of
    Proposition~\ref{P:mixed} supplied as an input;
  \item the occupation spaces are modelled by their $L^2$ completions;
  \item $H^1$ is modelled as the graph of the weak gradient.
\end{itemize}
On these models:
\begin{itemize}
  \item Propositions~\ref{P:spatial-four-linear}, \ref{P:occupation-physical-class},
    \ref{P:raw-occupation-energy} and~\ref{P:masked-auxiliary-norm}, and
    Lemma~\ref{L:variational-global-mean}, are formalized in full.
  \item Theorem~\ref{T:intrinsic-occupation-ceiling} is formalized in both forms
    with the exact constant $C_{\rm occ}$ over the occupation class itself,
    without its moment sentence $C_2\le U_2$.
  \item Propositions~\ref{P:d1-upper} and~\ref{P:low-alpha-upper} have their
    diagram values and coefficient bounds proved on the model, and their
    statements about the solution follow from the supplied
    Proposition~\ref{P:mixed}.
\end{itemize}
Open items in the record:
\begin{itemize}
  \item Theorems~\ref{T:full-upper} and~\ref{T:intrinsic-second-moment} are
    outside Lean.
  \item Theorem~\ref{T:white-moment-bounds} is formalized for the variational
    value $S_2^{\mathrm{pos}}$ rather than for $C_2$, with the white spatial
    energy supplied.
  \item Several results rest on named inputs, each listed with the result that
    supplies it. These include Propositions~\ref{P:mixed}
    and~\ref{P:three-transfer}, the first-moment series, and the signed
    multiplier form.
  \item No theorem of the manuscript is verified end to end.
\end{itemize}
These reviews are not a human verification, and they do not check the analytic
inputs that stay outside Lean.

For release v0.10.0, the companion's references to lines of this book were
moved to its final text. That work was reviewed by a model from a different
family; the first review and one confirmation failed on the companion's
descriptions and quotation ranges, the corrections were made, and a further
confirmation passed (Appendix~\ref{S:appendix-ai-statistics}). The author
accepted release v0.10.0 on November~16, 2026.

\Needspace{10\baselineskip}
\section{How the Lean package was checked}\label{S:lean-package-checks}

\ifdefined\lcHaveReceipts
\begin{wrapfigure}{r}{0.19\textwidth}
  \vspace{-\baselineskip}
  \centering
  \includegraphics[width=0.17\textwidth]{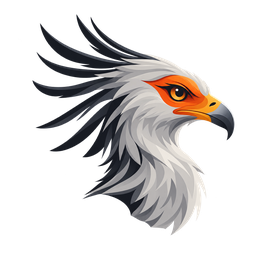}\\
  {\footnotesize\href{https://github.com/lean-pkg/lean-pkg}{lean-pkg}:
  registers Lean packages and records their checks.\par}
\end{wrapfigure}
The companion was built and checked with lean-pkg~\cite{chen:26:lean-pkg},
a tool for registering and checking Lean packages that the author developed
(version~\lcLeanPkgVersion). For commit \lcCommit, lean-pkg records that the
package builds from scratch with no unproved steps, and that every public
declaration, and each of the \lcClaimsChecked{} listed claims, uses only
Lean's three standard axioms (\texttt{propext}, \texttt{Classical.choice} and
\texttt{Quot.sound}). A second, separate check confirms that the proofs of
all \lcComparatorTheorems{} theorems prove their statements as restated in a
separate file.
\ifdefined\lcHaveLedger The list of claims has \lcVerified{} verified and
\lcOpen{} open.\fi{}
These checks concern the Lean statements themselves; whether those
statements say what the book says is the AI-reviewed match described above.
The private development verification record for commit \texttt{8f51004}
reports a separate kernel replay on a computing cluster: all 141 library
modules and the umbrella module passed \texttt{leanchecker --fresh}, and
all 142 targets also passed its default mode. The three libraries used only
for the separate statement check were not included in this replay.
That replay is of the gate commit. The release commit \texttt{\lcCommit}
drops two claims that name Lean-generated declarations and the one compared
theorem that named one of them; the receipts above are those of the release.
\fi

\section{The results}\label{S:lean-results}

In short: \lcSummary{} The \lcRows{} labelled results fall into \lcGroups{} groups.

\lcTables

\medskip
Besides these rows, the companion checks \lcUnlabelled{} further computations
that do not correspond to a single numbered result of the book.

\section{Counts and growth}\label{S:lean-counts}

\ifdefined\lcHaveCompletion
\Needspace{22\baselineskip}
\medskip
\noindent\textbf{Which numbered results the table counts.} The table covers
\lcNumberedTotal{} numbered theorems, propositions, lemmas and corollaries.
Besides the \lcOutsideNumbered{} main theorems recorded as outside Lean,
\lcNotFormalized{} further results (\lcNotFormalizedShare\% of the total) are
not formalized.
The separately numbered no-gos in Chapter~\ref{S:book-no-go} are not added to
this count; Appendix~\ref{S:appendix-check-status} counts them with the other
results. Some restate results listed here; the listed coverage retains
its own scope.
\ifnum\lcAddedSinceRecord>0 The companion's record was made against an earlier
revision of the book with \lcRecordTotal{} numbered results; the
\lcAddedSinceRecord{} added since then are counted here as not formalized.\fi
\begin{center}
\small\lcCompletionTable
\end{center}
The groups above also contain \lcOtherRows{} rows for displayed equations or
remarks, which are not numbered results. This is why the first group above has
\lcSameFullAll{} rows and the first line of this table \lcSameFullNumbered.

\Needspace{24\baselineskip}
\medskip
\noindent\textbf{How the companion grew.} Each stage is a tagged commit of the
companion, from its creation on \lcFirstStageDate{} to the version described in
this appendix. Claims counts the Lean declarations the companion lists for
checking, and Theorems counts those among them that are theorems (the others
are definitions); Results covered counts the rows of the lists above.
\begin{center}
\small\lcMilestoneTable
\end{center}
\fi

\chapter{How far each result has been checked}\label{S:appendix-check-status}

% Appendix G: how far each result has been checked. The chapter heading is in
% the main file. The counts come from verification-status-data.tex, generated
% by scripts/verification-status/make_appendix_g.py from the pinned results
% table of the verification ledger (scripts/verification-status/verification.md).

This appendix records how far the numbered results of this book have been
checked, in the three ways used in this project: by AI review, by Lean and by
people. It is generated from the verification ledger, a table in the book's
repository with one row for each of the \vsRows{} numbered theorems,
propositions, lemmas, corollaries and no-go results. The counts below are
those of the ledger on \vsDate{}. The ledger in the repository is the current
record, and it changes as results are checked. The ledger released with this
version of the book is at
\begin{center}\footnotesize
\href{https://github.com/chenle02/stratonovich-ham-second-moment/blob/v1.0.0/VERIFICATION.md}{\nolinkurl{github.com/chenle02/stratonovich-ham-second-moment/blob/v1.0.0/VERIFICATION.md}}
\end{center}

For credit, the ledger divides the results into five groups by chapter:
group~I is Chapters~\ref{S:book-setting}--\ref{S:book-independent-upper}
(the setting and the upper bound), group~II is
Chapters~\ref{S:book-gaussian-foundations}--\ref{S:book-lower} (Gaussian
approximation and the lower bound), group~III is
Chapters~\ref{S:book-means}--\ref{S:book-constant} (the matching upper bound
and the constant), group~IV is
Appendices~\ref{S:appendix-finite-order}--\ref{S:appendix-static-clock}, and
group~V is the no-go results of Chapter~\ref{S:book-no-go}.
Table~\ref{Tab:check-status} gives the counts.

\begin{table}[!ht]
  \centering
  \caption{How far the numbered results have been checked, by group
    (ledger of \vsDate).}\label{Tab:check-status}
  \begin{tabular}{@{}lrrrrrr@{}}
    \toprule
    & I & II & III & IV & V & All \\
    \midrule
    \vsTable
    \bottomrule
  \end{tabular}
\end{table}

\medskip
\noindent\textbf{AI review.} The repository keeps the AI reviews of this
project in dated folders, with the responses to their findings. A folder may
hold several review runs, or only preparation, decisions or diagnostics, and
not every run produced a usable report.
Appendix~\ref{S:appendix-ai-statistics} summarizes these records. No review
grade is assigned to an individual result, so the table has no row for them.
An AI review shows where to look; it is not a proof.

\medskip
\noindent\textbf{Lean.} The Lean rows come from the coverage file of the Lean
companion from which Appendix~\ref{S:appendix-lean-checks} is generated; the
script that generates this appendix checks the ledger's Lean entries and
numbers of claims against that file. Of the \vsRows{} results, \vsLeanAny{}
have some Lean claim, and \vsOutside{} more, all main theorems, are recorded as
outside Lean. Appendix~\ref{S:appendix-lean-checks} says, result by result,
what Lean checks and which inputs it takes as given. A Lean entry never
covers the analytic inputs, and whether the Lean statements say what the book
says was reviewed by AI, not by a person.

\medskip
\noindent\textbf{People.} To record a check, a mathematician opens a pull
request that changes exactly one row of the ledger and says what was checked
(the statement, the proof, the cited inputs and their hypotheses) and what was
not. A result restated from a cited source is checked against that source,
not reproved. To report a gap or an error, a reader opens an issue naming the
result's label and the step in question; the row is then marked as disputed,
and it returns to checked only through a new pull request that links the
closed issue.
\ifnum\vsHumanChecked=0 On the date above, no result had yet been checked in
this way. \fi
This is the column where readers can help most.

\section{What the main theorems rest on}\label{S:critical-path}

The four main theorems of this book are
Theorems~\ref{T:full-upper}, \ref{T:intrinsic-second-moment},
\ref{T:intrinsic-occupation-ceiling} and~\ref{T:white-moment-bounds}. This
section lists recorded citation links for these theorems in
Chapters~\cpFirstChapter{}--\cpLastChapter{}, as a reading aid.

A result \emph{rests directly on} another result when its statement or its
proof cites that result, cites an equation assigned to it, or uses a
preceding-result phrase recognized by the script, subject to the recorded
link corrections. A theorem
\emph{rests on} every result reached by following these links. The list is
made by a script from the source of the book
(\path{scripts/critical-path/make_critical_path.py}). It shows which results
a proof cites. It does not say that a proof is correct, and it is not a
review. It can err in both directions: a result cited only for comparison is
counted, and a result used without being cited by its label is missed. The
second happens often: a result is cited in the sentence that introduces a
statement, or it is used by name. For \cpReadByHand{} results, AI-assisted
author-side reading of the statements and proofs was recorded: the results
that most others depend on, and the results that the script did not reach. This added \cpLinksAdded{}
links and removed \cpLinksRemoved{}; each is recorded with the words that
justify it in \path{scripts/critical-path/link-corrections.json}. The other
results were not read in this way, and the third column of the list below may
miss links for them. In that column a number in parentheses is a result cited
only in the sentences that introduce the statement; it is shown for the reader
and is not counted as a link.

\begin{table}[!ht]
  \centering
  \caption[What each main theorem rests on]{What each main theorem rests on:
    the number of other results reached in the recorded graph, excluding the
    theorem named in the row, and their Lean classes as in
    Appendix~\ref{S:appendix-lean-checks} (the last column counts the results
    with no Lean claim or outside Lean).}\label{Tab:critical-path-theorems}
  \small\setlength{\tabcolsep}{4pt}
  \begin{tabular}{@{}l>{\raggedright\arraybackslash}p{0.24\textwidth}rrrrrr@{}}
    \toprule
    Theorem & Title & Results & Full & Narrower & Supplied & Partial & None \\
    \midrule
    \cpTheoremRows
    \bottomrule
  \end{tabular}
\end{table}

The recorded graph reaches \cpResults{} of the \cpProved{} results with
printed proofs in Chapters~\cpFirstChapter{}--\cpLastChapter{}, including
the four main theorems themselves. The second-moment theorem reaches
\cpSecondMoment{} other results in this graph. Table~\ref{Tab:critical-path-chapters}
gives the count by chapter. Of the \cpResults{} results, Lean formalizes
\cpLeanFull{} in full and \cpLeanNarrower{} with a narrower conclusion,
proves \cpLeanSupplied{} with supplied inputs and checks steps of
\cpLeanPartial{}; the other \cpLeanNone{} have no Lean claim or are outside
Lean. No result has been checked by a person.

\begin{table}[!ht]
  \centering
  \caption[The results the main theorems rest on, by chapter]{The results the
    four main theorems rest on, by chapter, with their Lean
    classes. Counts include the four main theorems.}\label{Tab:critical-path-chapters}
  \small
  \begin{tabular}{@{}rrrrrrrr@{}}
    \toprule
    Chapter & Proved results & Rested on & Full & Narrower & Supplied & Partial & None \\
    \midrule
    \cpChapterRows
    \bottomrule
  \end{tabular}
\end{table}

The list below gives all results reached by the recorded graph in
Chapters~\cpFirstChapter{}--\cpLastChapter{}, in book order. The third column
gives recorded direct links. The fourth counts other listed results that
reach the row through these links; a large count indicates many recorded
citation paths, not mathematical indispensability. The last column is its
Lean class; a dash means no Lean claim.

\medskip
\cpHead
\cpRows

\medskip
The recorded graph does not reach \cpNotReached{} results with printed
proofs in these chapters:
\cpNotReachedList.
This does not mean that they have no mathematical use.
Of Proposition~\ref{P:static-phase-variance}, the proof of
Proposition~\ref{P:interval-coefficients} uses one calculation, the
free-phase energy bound at the end of its proof, and none of its conclusions.

\backmatter
\bibliography{All,extra}
\chapter*{Working words and common misreadings}
\addcontentsline{toc}{chapter}{Working words and common misreadings}
\markboth{Working words and common misreadings}{Working words and common misreadings}
\label{S:glossary}

The proof uses a few everyday words in a fixed sense. This list collects them,
and then names the misreadings that are easiest to make. It adds no assumption
and no statement. For symbols see the List of notation, and for page numbers
the Index.

\section*{Working words}

\begin{description}
  \item[Chain.] The ordered $n$-vertex wave kernel $g_n$: wave Green factors
    at $n$ times in increasing order. The $n$-th term of the solution
    integrates $n$ noise factors against one chain (Section~\ref{S:known}).
  \item[Order.] The number of noise factors in a term. Terms of different
    orders need not be orthogonal here, and the total order $m+n$ of a product
    is not a Wiener-chaos order.
  \item[Pairing.] A matching of noise factors two by two, made when an
    expectation is taken. The second moment sums over all pairings of the
    vertices of two chains (Proposition~\ref{P:mixed}). A self-pairing matches
    two vertices of the same chain.
  \item[Canonical solution.] The sum of the fixed-order limits, obtained from
    the regularized classical equations
    (Proposition~\ref{P:mild-construction}).
  \item[Phase.] Displacement and velocity together. Localized perturbations
    lie in $E=H^{1/2}\times H^{-1/2}$; constant data are carried in
    $F=E\oplus\mathbb{R}^2$ (Section~\ref{S:constant-phase-young}).
    A restarted wave needs both components (No-go~\ref{P:nogo-14}).
  \item[Mask.] A smooth compactly supported spatial function that multiplies
    the physical potential. The lower bound takes a numerical supremum after
    fixed-mask bounds; the upper point comparison uses growing supports
    (Propositions~\ref{P:variational-remove-mask} and~\ref{P:original-point-transfer}).
  \item[Radius.] The size of the region where the mask equals one.
    The upper point comparison chooses a radius growing with its fast time
    horizon; lower-bound masks remain fixed before time asymptotics
    (Chapters~\ref{S:book-lower} and~\ref{S:book-growing-radius}).
  \item[Block.] One interval of a subdivision of the time axis. The noise in
    one block is correlated with the noise in the others, so bounds do not
    simply multiply across blocks.
  \item[Control.] A vector in the covariance Hilbert space that specifies a
    Gaussian mean shift. The lower bound uses nonzero nonnegative smooth
    compact tests $\Psi\in\mathcal{C}_{\mathrm{pos}}$; the upper bound also
    uses signed completed controls $\Phi\in H_I$.
    The shift is toward the induced mean.
  \item[Gaussian cost.] The price of that shift: $v_\Psi/2$, half the variance
    of $W(\Psi)$. It is paid once for the whole time interval, not once per
    block.
  \item[Variational value.] The nonnegative-control value
    $S_c^{\mathrm{pos}}=\sup_{\Psi\in\mathcal{C}_{\mathrm{pos}}}
    (cI_\Psi-v_\Psi/2)$ (Proposition~\ref{P:variational-positive-scores}).
    The constant of Theorem~\ref{T:intrinsic-second-moment} is $S_2^{\mathrm{pos}}$.
  \item[Ceiling.] An upper bound for the constant, computed from occupation
    and energy estimates (Theorem~\ref{T:intrinsic-occupation-ceiling}).
  \item[Clock.] An auxiliary scalar Brownian motion used with a spatial
    Brownian path to represent the time Laplace transform of the expectation
    of one expansion term (Section~\ref{S:two-clock-screen}). Copying it with
    two clocks, constant correlation and no extra weight fails for the first
    cross pairing (Proposition~\ref{P:clock-copy-obstruction}).
\end{description}

\section*{Common misreadings}

\begin{description}
  \item[``This is the Skorohod model.''] It is not. Stratonovich second
    moments retain self-pairings and mixed terms with $m\ne n$. In the
    Skorohod chaos expansion, unequal orders are orthogonal and there are no
    self-contractions. The exponent comparisons here require the cited
    covariance and data hypotheses; they identify neither solutions nor
    leading constants.
  \item[``Theorem~\ref{T:full-upper} gives the constant.''] It gives an upper
    bound at scale $T^\kappa$, with a constant not claimed sharp. The exact
    asymptotic constant is identified in Theorem~\ref{T:intrinsic-second-moment},
    by a separate argument.
  \item[``The white-noise formulas evaluate the constant.''] They give
    two-sided bounds (Theorem~\ref{T:white-moment-bounds}). The spatial energy
    is evaluated in Proposition~\ref{P:white-spatial-energy}; its relation to
    the moment constant gives bounds, not an evaluation of the full intrinsic
    supremum.
  \item[``The ceiling is the constant.''] The ceiling is not asserted to be
    attained or to equal the constant.
  \item[``$S_1^{\mathrm{pos}}$ is the first-moment constant.''] It is an
    auxiliary variational value. No identification with a first-moment
    constant is asserted.
  \item[``The arXiv version of the source is the one cited.''] Statements from
    Xia Chen's paper are cited from his public author version, posted on his
    homepage at the University of Tennessee and retrieved on October~26,
    2026~\cite{chen:author2026:time-dependency}. Its Theorem~1.2 differs from
    the statement in arXiv:2510.01412v1~\cite{chen:25:time-dependency}.
\end{description}

% Selected notation with explicit scope and definition links.
% Generated from notes/index/notation.json.
\chapter*{List of notation}
\addcontentsline{toc}{chapter}{List of notation}
\markboth{List of notation}{List of notation}
\label{S:notation}

The tables collect recurring notation and point to its definition or
principal use. Symbols introduced for one calculation are defined there.
Several letters are reused in different sections; the scope column records
the main distinctions. These tables do not add assumptions to the model.

\begingroup
\small
\setlength{\LTpre}{0pt}
\setlength{\LTpost}{12pt}

\Needspace{8\baselineskip}
\section*{Model, noise and scaling parameters}
\begin{longtable}{@{}>{\raggedright\arraybackslash}p{0.20\textwidth}
  >{\raggedright\arraybackslash}p{0.49\textwidth}
  >{\raggedright\arraybackslash}p{0.25\textwidth}@{}}
  \toprule
  Symbol & Meaning and scope & Definition or use \\
  \midrule
  \endfirsthead
  \toprule
  Symbol & Meaning and scope & Definition or use \\
  \midrule
  \endhead
  \midrule
  \multicolumn{3}{r}{\emph{Continued on the next page}} \\
  \endfoot
  \bottomrule
  \endlastfoot
  $d$ &
  Spatial dimension, $d \in \{1, 2, 3\}$. &
  Hypothesis~\ref{H:homog} \\[5pt]
  $\alpha$ &
  Exponent in the spatial covariance homogeneity relation. &
  Hypothesis~\ref{H:homog} \\[5pt]
  $\alpha_0$, $a$ &
  Temporal covariance exponent; $a = \alpha_0$ in the integrated arguments. &
  \eqref{E:cov} \\[5pt]
  $\Gamma$, $\gamma$ &
  Spatial covariance measure; $\gamma$ denotes its density when one exists and
  otherwise the corresponding distribution. &
  Hypothesis~\ref{H:homog} \\[5pt]
  $\mu$ &
  Spectral measure, in the Fourier convention fixed by the hypothesis. &
  Hypothesis~\ref{H:homog} \\[5pt]
  $\dot W$, $W(f)$ &
  Gaussian forcing and its isonormal realization on covariance test classes,
  respectively. &
  \eqref{E:cov}; Section~\ref{S:actual-operator-driver} \\[5pt]
  $\mathbb{E}$ &
  Expectation under the underlying probability measure. &
  \eqref{E:intro-second-moment-limit} \\[5pt]
  $\kappa$ &
  Normalization exponent $(4 - \alpha - \alpha_0)/(3 - \alpha)$ in the moment
  bounds and limit. &
  \eqref{E:kappa} \\[5pt]
  $\beta$, $\nu$ &
  In the interval-field estimates, $\beta = 4 - \alpha - \alpha_0$ and $\nu = 3
  - \alpha$. These letters have different local uses in path estimates. &
  Section~\ref{S:interval-upper} \\[5pt]
  $R_\alpha$ &
  Spectral-mass constant in the independent upper bound. &
  \eqref{E:intro-spectral-mass} \\[5pt]
  $\Gamma_E$, $\operatorname{Beta}$ &
  Euler Gamma and Beta functions; $\Gamma_E$ is distinct from the covariance
  measure $\Gamma$. &
  Lemma~\ref{L:equilibrium-gamma-beta} \\[5pt]
  $\chi$ &
  Degree of an auxiliary static covariance in the comparison arguments; it need
  not equal $\alpha$. &
  Section~\ref{S:three-endpoints} \\[5pt]
\end{longtable}

\Needspace{8\baselineskip}
\section*{Solutions, kernels and pairing expansions}
\begin{longtable}{@{}>{\raggedright\arraybackslash}p{0.20\textwidth}
  >{\raggedright\arraybackslash}p{0.49\textwidth}
  >{\raggedright\arraybackslash}p{0.25\textwidth}@{}}
  \toprule
  Symbol & Meaning and scope & Definition or use \\
  \midrule
  \endfirsthead
  \toprule
  Symbol & Meaning and scope & Definition or use \\
  \midrule
  \endhead
  \midrule
  \multicolumn{3}{r}{\emph{Continued on the next page}} \\
  \endfoot
  \bottomrule
  \endlastfoot
  $u(t,x)$ &
  Canonical Stratonovich solution with displacement $1$ and velocity $0$ at time
  zero. &
  \eqref{E:HAM}; Proposition~\ref{P:mild-construction} \\[5pt]
  $G_t$ &
  Wave Green measure, with mass $t$ and support in the closed radius-$t$ ball. &
  \eqref{E:expansion}; Lemma~\ref{L:constant-phase-free} \\[5pt]
  $g_n$, $S_n(g_n)$ &
  Ordered chain kernel and its multiple Stratonovich integral. &
  \eqref{E:expansion} \\[5pt]
  $X_n(t,x)$ &
  The $n$th Stratonovich term $S_n(g_n)$; its iteration order is not an
  orthogonal chaos order. &
  Proposition~\ref{P:mixed} \\[5pt]
  $I_1$, $I_2$, $\Pi_k$ &
  The two ordered chains and the pair partitions of their $2k$ vertices. &
  \eqref{E:mixed-diagram} \\[5pt]
  $P_k$ &
  Number of pair partitions, $(2k)!/(2^k k!)$. &
  \eqref{E:B-envelope} \\[5pt]
  $J_{m,n,D}$ &
  Pairing integral for chain lengths $m,n$ and pair partition $D$. &
  \eqref{E:mixed-diagram} \\[5pt]
  $N_0$, $N_1$, $N_2$ &
  Laplace-damped wave norms used in the construction estimates. &
  \eqref{E:damped-norm-definitions} \\[5pt]
  $B_K(t)$ &
  Numerical majorant obtained by estimating each pairing separately; it is not
  the actual second moment. &
  \eqref{E:B-envelope} \\[5pt]
\end{longtable}

\Needspace{8\baselineskip}
\section*{Covariance spaces and wave evolution}
\begin{longtable}{@{}>{\raggedright\arraybackslash}p{0.20\textwidth}
  >{\raggedright\arraybackslash}p{0.49\textwidth}
  >{\raggedright\arraybackslash}p{0.25\textwidth}@{}}
  \toprule
  Symbol & Meaning and scope & Definition or use \\
  \midrule
  \endfirsthead
  \toprule
  Symbol & Meaning and scope & Definition or use \\
  \midrule
  \endhead
  \midrule
  \multicolumn{3}{r}{\emph{Continued on the next page}} \\
  \endfoot
  \bottomrule
  \endlastfoot
  $H_W$ &
  Completion of smooth positive-time tests modulo covariance-null tests for the
  original noise. &
  Section~\ref{S:actual-operator-driver} \\[5pt]
  $I$ &
  The unit time interval; control tests have support in $(0,1)$, while compact
  parameter families also use $[0,1]$. &
  \eqref{E:variational-positive-class}; Section~\ref{S:compact-mean-family}
  \\[5pt]
  $H_t$, $H_{\mathrm{sp}}$ &
  Temporal and spatial covariance Hilbert spaces. &
  Section~\ref{S:actual-cameron-mean} \\[5pt]
  $H_I$ &
  Covariance space of tests supported in $I \times \mathbb{R}^d$, embedded in
  $H_W$. &
  Lemma~\ref{L:cameron-supported-tensor} \\[5pt]
  $\Lambda$ &
  Bessel multiplier $(1-\Delta)^{1/2}$. &
  Section~\ref{S:wave-variance} \\[5pt]
  $E_r$, $E$ &
  Wave phase $E_r = H^r \times H^{r-1}$; the balanced phase is $E = H^{1/2}
  \times H^{-1/2}$. &
  Section~\ref{S:actual-operator-driver} \\[5pt]
  $F$, $e_0$ &
  Extended phase $F = E \oplus \mathbb{R}^2$ and the vector encoding constant
  initial displacement. &
  \eqref{E:constant-phase-space} \\[5pt]
  $S_E(t)$, $S_F(t)$ &
  Free wave groups on the respective phase spaces. &
  Lemma~\ref{L:constant-phase-free} \\[5pt]
  $M_R(W)$ &
  Integrated Gaussian interaction operator path with spatial localization radius
  $R$. &
  Proposition~\ref{P:actual-operator-driver} \\[5pt]
  $U_D$ &
  Young evolution driven by an operator path $D$ in the interaction
  representation. &
  \eqref{E:constant-young-equation} \\[5pt]
  $P_D(t,s)$ &
  Physical wave evolution obtained by conjugating the interaction flow by the
  free wave group. &
  \eqref{E:constant-physical-flow}; Section~\ref{S:actual-symplectic-phase}
  \\[5pt]
  $B_E(V)$ &
  Potential insertion $(u,v) \mapsto (0,Vu)$ on $E$. &
  Lemma~\ref{L:uniform-E-lift} \\[5pt]
  $\omega$, $J_E$ &
  Symplectic form and its representing skew-adjoint operator on $E$. &
  Lemma~\ref{L:symplectic-free} \\[5pt]
\end{longtable}

\Needspace{8\baselineskip}
\section*{Regularization, spatial masks and time blocks}
\begin{longtable}{@{}>{\raggedright\arraybackslash}p{0.20\textwidth}
  >{\raggedright\arraybackslash}p{0.49\textwidth}
  >{\raggedright\arraybackslash}p{0.25\textwidth}@{}}
  \toprule
  Symbol & Meaning and scope & Definition or use \\
  \midrule
  \endfirsthead
  \toprule
  Symbol & Meaning and scope & Definition or use \\
  \midrule
  \endhead
  \midrule
  \multicolumn{3}{r}{\emph{Continued on the next page}} \\
  \endfoot
  \bottomrule
  \endlastfoot
  $p_\epsilon$ &
  Gaussian heat kernel used for spatial smoothing; its covariance normalization
  is specified locally. &
  Lemma~\ref{L:variational-heat-envelope} \\[5pt]
  $\rho$, $Y_\rho$ &
  A pair of noise-smoothing widths and the corresponding regularized source. In
  damped phase estimates, $\rho$ instead denotes a positive damping parameter. &
  Proposition~\ref{P:source-canonical-convergence};
  Section~\ref{S:masked-moments} \\[5pt]
  $R$, $\eta_R$, $m_R$ &
  Spatial radius, cutoff profile and physical mask. In the uniform-radius
  arguments, $\eta_R = \eta(\cdot/R)$ and $m_R = \eta_R^2$. &
  Section~\ref{S:uniform-mean-regulation} \\[5pt]
  $\delta$, $r_T$ &
  Long-time normalization $\delta = T^{-\kappa}$ and shrinking spatial scale
  $r_T = T\delta$. Smoothing parameters denoted by $\delta$ elsewhere are local.
  &
  Lemma~\ref{L:exact-canonical-scaling} \\[5pt]
  $g_\delta$ &
  Fast-noise amplitude $\delta^{(1 - a)/2}$. &
  \eqref{E:raw-block-clocks} \\[5pt]
  $B$, $N$, $s_m$, $\tau$ &
  Fast block length, number of full blocks, slow block starts and final block
  length. Here $\tau \in [0,B)$. &
  \eqref{E:raw-block-clocks} \\[5pt]
  $T_{\delta,m}$ &
  Isometry used to pull each fast-time test back to the original Gaussian noise.
  &
  \eqref{E:raw-block-pullback} \\[5pt]
  $X_m$, $H_m$ &
  Random block operator input and its deterministic mean shift; these local
  symbols differ from the Stratonovich terms $X_n(t,x)$. &
  \eqref{E:raw-block-drivers} \\[5pt]
  $R_\delta$ &
  Growing spatial radius $\delta^{-1} + 2$ in the point upper bound; it differs
  from the shrinking scale $r_T$. &
  Section~\ref{S:original-point-transfer} \\[5pt]
  $\gamma$, $\theta$, $\beta$ &
  Auxiliary spatial and temporal regularity exponents. Here $\gamma$ is not the
  covariance distribution and $\beta$ is not the coefficient-growth exponent. &
  \eqref{E:driver-parameters} \\[5pt]
\end{longtable}

\Needspace{8\baselineskip}
\section*{Controls, means and variational constants}
\begin{longtable}{@{}>{\raggedright\arraybackslash}p{0.20\textwidth}
  >{\raggedright\arraybackslash}p{0.49\textwidth}
  >{\raggedright\arraybackslash}p{0.25\textwidth}@{}}
  \toprule
  Symbol & Meaning and scope & Definition or use \\
  \midrule
  \endfirsthead
  \toprule
  Symbol & Meaning and scope & Definition or use \\
  \midrule
  \endhead
  \midrule
  \multicolumn{3}{r}{\emph{Continued on the next page}} \\
  \endfoot
  \bottomrule
  \endlastfoot
  $\mathcal{C}_{\mathrm{pos}}$ &
  Nonzero nonnegative smooth compact spacetime controls; separation into time
  and space factors is not required. &
  \eqref{E:intro-control-class} \\[5pt]
  $h_{\Psi,s}$ &
  Deterministic covariance mean induced by the control $\Psi$. &
  \eqref{E:intro-control-mean} \\[5pt]
  $v_\Psi$ &
  Original Gaussian variance, equal to the squared covariance-space norm of the
  control. &
  \eqref{E:intro-control-variance}; \eqref{E:compact-raw-cost} \\[5pt]
  $\lambda_{\Psi,s}$, $I_\Psi$ &
  Positive Rayleigh growth rate and its integral over the unit interval, for
  nonzero covariance. &
  \eqref{E:variational-real-rayleigh} \\[5pt]
  $\mathcal{R}_\Psi(s)$ &
  Squared Rayleigh rate used in the introductory constant formula. &
  \eqref{E:intro-control-rate} \\[5pt]
  $S_c^{\mathrm{pos}}$ &
  Supremum of $cI_\Psi - v_\Psi/2$. The coefficient-one value is auxiliary; the
  coefficient-two value identifies $C_2$. &
  Proposition~\ref{P:variational-positive-scores};
  Theorem~\ref{T:intrinsic-second-moment} \\[5pt]
  $C_2(\Gamma,a)$ &
  Intrinsic second-moment constant, including its zero-covariance value. &
  \eqref{E:intrinsic-second-moment-limit} \\[5pt]
  $k$, $q$, $m$ &
  In the spatial-energy and constant calculations, $k = 3 - \alpha$, $q = 4 -
  \alpha$ and $m = q/k$. These are local conventions. &
  Section~\ref{S:spatial-occupation-foundation};
  \eqref{E:variational-score-homogeneity} \\[5pt]
  $f_a$, $M_a$ &
  Temporal equilibrium density and its total mass. $M_a$ is distinct from the
  mean-growth rate $M(a)$. &
  \eqref{E:temporal-equilibrium-density}; \eqref{E:temporal-equilibrium-mass}
  \\[5pt]
  $K(s,t)$, $K_0(a)$ &
  Original temporal kernel $|s-t|^{-a}$ and its maximal row integral $2^a/(1 -
  a)$. &
  \eqref{E:positive-bochner-inclusion}; \eqref{E:white-moment-bounds} \\[5pt]
  $B_\Gamma(f,h)$ &
  Physical two-point covariance pairing. &
  \eqref{E:spatial-physical-pair} \\[5pt]
  $A_{\mathrm{sp}}(f)$, $b(f)$ &
  Spatial covariance energy of $f^2$ and gradient energy $\|\nabla f\|_2^2$. &
  Section~\ref{S:spatial-occupation-foundation} \\[5pt]
  $R_{\mathrm{sp}}$ &
  Best homogeneous spatial energy ratio, extended to real unit $H^1$ profiles. &
  Proposition~\ref{P:spatial-four-linear};
  Proposition~\ref{P:spatial-energy-map} \\[5pt]
  $\mathcal{E}_{\mathrm{sp}}$ &
  Supremum of $A_{\mathrm{sp}}(f) - b(f)/2$ over unit spatial profiles. &
  \eqref{E:spatial-energy-map} \\[5pt]
  $F_0(z)$ &
  Scalar map from spatial energy to a variational bound; it is distinct from the
  wave phase $F$. &
  \eqref{E:spatial-static-energy-score} \\[5pt]
  $\mathcal{O}$, $T(g)$ &
  Admissible time-dependent unit spatial profiles and their integrated gradient
  energy. &
  Proposition~\ref{P:occupation-physical-class} \\[5pt]
  $A(g)$, $C_{\rm occ}$ &
  Original temporal occupation energy and its homogeneous ratio supremum. &
  \eqref{E:raw-occupation-definitions} \\[5pt]
  $\mathcal{E}_\zeta^{\mathcal{O}}$ &
  Occupation energy with gradient penalty $\zeta T(g)$, using the original
  temporal kernel. &
  \eqref{E:raw-occupation-definitions} \\[5pt]
  $K_\star$, $\mathcal{E}_\star^{\mathcal{O}}$ &
  Auxiliary damped temporal kernel and its occupation energy; these differ from
  the original-kernel quantities. &
  \eqref{E:occupation-damped-kernel}; \eqref{E:occupation-own-energy} \\[5pt]
  $C_t$, $T_R$, $\mathcal{T}_R$ &
  Temporal, spatial and combined covariance-mean maps. Their domains and radius
  dependence are specified at each construction. &
  Propositions~\ref{P:cameron-temporal-map} and~\ref{P:cameron-full-mean};
  Section~\ref{S:uniform-mean-regulation} \\[5pt]
  $\mathcal{S}_R$, $V_{\rm up}$ &
  Signed-control variational value at radius $R$ and its limsup along the
  growing radii. &
  \eqref{E:signed-raw-score} \\[5pt]
\end{longtable}

\Needspace{8\baselineskip}
\section*{Gaussian geometry and finite optimization}
\begin{longtable}{@{}>{\raggedright\arraybackslash}p{0.20\textwidth}
  >{\raggedright\arraybackslash}p{0.49\textwidth}
  >{\raggedright\arraybackslash}p{0.25\textwidth}@{}}
  \toprule
  Symbol & Meaning and scope & Definition or use \\
  \midrule
  \endfirsthead
  \toprule
  Symbol & Meaning and scope & Definition or use \\
  \midrule
  \endhead
  \midrule
  \multicolumn{3}{r}{\emph{Continued on the next page}} \\
  \endfoot
  \bottomrule
  \endlastfoot
  $\mathcal{K}(E)$ &
  Compact operators on the Hilbert phase $E$. &
  Proposition~\ref{P:actual-operator-driver} \\[5pt]
  $\mathscr{B}_\delta$ &
  Banach space of block paths with its weighted temporal norm. &
  Proposition~\ref{P:tube-actual-law} \\[5pt]
  $C$, $\mathcal{M}_\delta$ &
  In Gaussian enlargement, $Ch = \mathbb{E}[XW(h)] =
  \sqrt\delta\,\mathcal{M}_\delta h$. The letter $C$ elsewhere often denotes an
  estimate constant. &
  \eqref{E:tube-covariance-scale} \\[5pt]
  $H_{\rm seen}$, $H_{\rm CM}$ &
  Visible isonormal subspace and Cameron--Martin image of the covariance map. &
  \eqref{E:tube-covariance-scale}; \eqref{E:tube-cameron-space} \\[5pt]
  $E_0$ &
  Closed support space of the Gaussian input; this is distinct from the wave
  phase $E$. &
  \eqref{E:tube-covariance-scale} \\[5pt]
  $M_\delta$ &
  Expected squared operator norm of the physical endpoint in the upper argument.
  &
  Section~\ref{S:original-signed-score-upper} \\[5pt]
  $C_{ij}$ &
  Finite Gram matrix entries $\langle k_i,k_j\rangle$, assumed nonnegative in
  the positive-cone optimization. &
  \eqref{E:finite-gram-score} \\[5pt]
  $\mathfrak{F}_C$, $\mathscr{S}(C)$ &
  Finite-control objective and its supremum. &
  \eqref{E:finite-gram-score} \\[5pt]
  $D_A^C$, $d_A^C$ &
  Dual objective and its infimum for a selected subset $A$. The empty subset has
  value zero. &
  \eqref{E:finite-gram-dual} \\[5pt]
  $w^A$, $\Phi_A$ &
  Positive dual minimizer and the corresponding positive combination of Gram
  vectors. &
  \eqref{E:gram-dual-stationarity}; Proposition~\ref{P:gram-positive-optimizer}
  \\[5pt]
\end{longtable}

\Needspace{8\baselineskip}
\section*{Static comparisons and local outlook notation}
\begin{longtable}{@{}>{\raggedright\arraybackslash}p{0.20\textwidth}
  >{\raggedright\arraybackslash}p{0.49\textwidth}
  >{\raggedright\arraybackslash}p{0.25\textwidth}@{}}
  \toprule
  Symbol & Meaning and scope & Definition or use \\
  \midrule
  \endfirsthead
  \toprule
  Symbol & Meaning and scope & Definition or use \\
  \midrule
  \endhead
  \midrule
  \multicolumn{3}{r}{\emph{Continued on the next page}} \\
  \endfoot
  \bottomrule
  \endlastfoot
  $L_0$, $\kappa_{\rm st}$ &
  Static first-mean rate and static time exponent; the comparison retains the
  source normalization stated there. &
  \eqref{E:static-ratio-input} \\[5pt]
  $M(a)$ &
  Normalized first-mean limsup in the near-static comparison; it is not the
  equilibrium mass $M_a$. &
  Proposition~\ref{P:near-static-mean-gap} \\[5pt]
  $\rho$ &
  In the clock comparison only, correlation of the two auxiliary Brownian
  clocks; neither a source width nor a damping parameter. &
  \eqref{E:naive-clock-candidate} \\[5pt]
  $p$ &
  In future-work moment statements, a fixed moment order. Elsewhere it may be an
  auxiliary integrability or Schatten index. &
  Conjecture~\ref{Conj:intrinsic-even-moments}; Section~\ref{S:compact-gaussian}
  \\[5pt]
  $p$, $\nu$ &
  Locally, $p = 2 - \alpha_0$ in the interval-field argument, while $\nu =
  1/\theta$ in Young growth estimates. &
  Section~\ref{S:interval-upper}; Lemma~\ref{L:young-upper-growth} \\[5pt]
  $A$, $z_*$ &
  In the conditional upper-tail discussion, $A = 2^{-m}C_2$ and $z_* = mA$. This
  $A$ is not the occupation energy $A(g)$. &
  Section~\ref{S:future-real-orders} \\[5pt]
  $\mathcal{I}(z)$ &
  Conditional upper-tail rate on the stated range of thresholds. &
  \eqref{E:future-tail-rate} \\[5pt]
  $H_a$ &
  Covariance Hilbert space for the spatial-white optimizer conjecture. &
  Section~\ref{S:future-optimizers} \\[5pt]
  $J_a$, $\tau_y$ &
  Control objective and spatial translation in the optimizer conjecture. &
  Conjecture~\ref{Conj:positive-control-compactness} \\[5pt]
  $h$, $b(t)$, $\varphi$ &
  In the smooth logarithmic-wave outlook, $h = \log|u|$ and $h = b(t)+\varphi$.
  This $b(t)$ is distinct from the spatial energy $b(f)$. &
  \eqref{E:future-log-wave}; \eqref{E:future-centered-log-wave} \\[5pt]
\end{longtable}
\endgroup

\printindex
\end{document}